\UseRawInputEncoding
\documentclass[11pt,a4paper,leqno]{article}

\usepackage{tikz-cd}
\usepackage{tikz}
\usepackage{dynkin-diagrams}

\usepackage{svg}

\usepackage{rotating}

\usepackage{graphicx,accents}

\makeatletter
\DeclareRobustCommand{\cev}[1]{%
  \mathpalette\do@cev{#1}%
}
\newcommand{\do@cev}[2]{%
  \fix@cev{#1}{+}%
  \reflectbox{$\m@th#1\vec{\reflectbox{$\fix@cev{#1}{-}\m@th#1#2\fix@cev{#1}{+}$}}$}%
  \fix@cev{#1}{-}%
}
\newcommand{\fix@cev}[2]{%
  \ifx#1\displaystyle
    \mkern#23mu
  \else
    \ifx#1\textstyle
      \mkern#23mu
    \else
      \ifx#1\scriptstyle
        \mkern#22mu
      \else
        \mkern#22mu
      \fi
    \fi
  \fi
}

\makeatother

\usepackage{psfrag}

\newenvironment{sproof}{%
  \proof}{\endproof}

\usepackage{titlesec}
\usepackage{bigints}

\titleformat{\subsubsection}[runin]
  {\normalfont\normalsize\bfseries}{\thesubsubsection}{1em}{}

\DeclareMathAlphabet{\mathpzc}{OT1}{pzc}{m}{it}

\usepackage{mathrsfs}

\usepackage{shuffle}

\usepackage{slashbox}

\usepackage{mathtools}

\usepackage{eepic}

\usepackage[hang]{footmisc}

\usepackage{graphicx,epsfig,amscd,graphics,xypic}

\usepackage{amsmath}
\usepackage{amssymb}
\usepackage{amsthm}
\usepackage{latexsym}
\usepackage{pstricks}
\usepackage{amsbsy}
\usepackage{amscd}
\usepackage{mathrsfs}
\usepackage{tipa}
\usepackage{txfonts}
\usepackage{amsfonts}
\usepackage{graphicx}
\usepackage{tabularx}
\usepackage{listings}
\usepackage{tikz}
\usepackage{hyperref}
\tikzstyle{dot} = [inner sep=0pt,thick,fill=black,circle,minimum size=2.5pt]
\tikzstyle{line} = [draw, -latex]
\newtheorem{thm}{Theorem}[section]
\newtheorem{cor}[thm]{Corollary}
\newtheorem{lem}[thm]{Lemma}

\newtheorem{exm}{Example}
\newtheorem{exms}{Examples}
\newtheorem{fact}[thm]{Fact}
\newtheorem{prop}[thm]{Proposition}
\newtheorem{defn}[thm]{Definition}
\newtheorem{rem}[thm]{Remark}
\newtheorem{defn-prop}[thm]{Definition-Proposition}
\newtheorem{conjecture}[thm]{Conjecture}
\newtheorem{warning}[thm]{Warning}
\newtheorem{question}[thm]{Question}
\newtheorem{questions}[thm]{Questions}

\def \l  #1 {{ {\rm {\bf L}i}_{#1}}}

\begin{document}

\begin{center} 
{\Large \bf  On  webs, polylogarithms and cluster algebras}

\bigskip

{\large Luc Pirio}

\bigskip


\bigskip


\end{center}

\begin{abstract}
In this text, we investigate webs which can be associated to cluster algebras from the point of view of the abelian functional equations these webs carry, focusing on the polylogarithmic ones.  We introduce a general notion of webs whose rank is `{\it As Maximal as Possible}' (AMP) and show that many webs associated  either to polylogarithmic functional equations or to cluster algebras are of this type. In particular, we prove a few results and state some conjectures about cluster webs associated to (pairs of) Dynkin diagrams.  Along the way,  we show that many of the classical functional equations  satisfied by low-order polylogarithms  
 (such as Spence-Kummer's equation of the trilogarithm or the tetralogarithmic one of Kummer) are of cluster type.
\end{abstract}


\newcommand{\sk}{\smallskip}
\newcommand{\mk}{\medskip}
\newcommand{\bk}{\bigskip}

\newcommand{\homeo}{\textup{Homeo}^+(S,M)}
\newcommand{\homeoo}{\textup{Homeo}_0(S,M)}
\newcommand{\mg}{\mathcal{MG}(S,M)}
\newcommand{\mmg}{\mathcal{MG}_{\bowtie}(S,M)}
\newcommand{\Z}{\mathbf{Z}}
\newcommand{\Q}{\mathbf{Q}}
\newcommand{\N}{\mathbf{N}}

\hyphenation{ap-pro-xi-ma-tion}


\tableofcontents

\newpage

\section{General introduction}
 \href{}{\bf Webs} are geometric objects formed (locally) by a finite number of foliations, whose  geometrical study was started at Hamburg in the 1930s, by a group of mathematicians led by \href{http://www-history.mcs.st-andrews.ac.uk/Biographies/Blaschke.html}{Blaschke}. 
\href{https://en.wikipedia.org/wiki/Polylogarithm}{\bf Polylogarithms} 
 are special functions of one variable,  the study of which can be traced back to Euler.  
   \href{https://en.wikipedia.org/wiki/Cluster_algebra}{\bf Cluster algebras} are commutative rings constructed inductively according to a combinatorial recipe, known as {\it mutation}. These have been introduced  at the beginning of the 2000's by \href{https://en.wikipedia.org/wiki/Sergey_Fomin}{Fomin} and \href{https://en.wikipedia.org/wiki/Andrei_Zelevinsky}{Zelevinsky}, and are nowadays studied from many points of view in mathematics. 
\mk 

Polylogarithms, webs and cluster algebras therefore are mathematical objects of three distinct natures but  despite this,    they are related. What links them  is a particular class of functional equations that we call {\bf Abelian functional equations} 
 and which play a crucial role in this text.\vspace{-0.2cm}
\begin{center}
$\star$
\end{center}

In the sequel of the present Introduction, we begin by introducing these three notions with more details before considering the kind of functional equations making a bridge between them. Then in subsection \S\ref{SS:Intro-Results}, we describe some of the main features of this text: main new notions introduced, our most appealing results and some interesting conjectures all this suggests.  In  \S\ref{SS:Intro3}, we describe the content of the text and comment briefly on the general framework of our study.

%

 \subsection{\bf The main objects: webs, polylogarithms and cluster algebras}
\label{SS:TheMainObjects}
We describe the main objects considered in this text and explain how they are related according to the study undertaken here. Everything below is presented and explained quite succinctly, the reason being that all this material will be introduced further on with much more details and explanations.

 \subsubsection{\bf Webs and web geometry.}
 \label{SSS:Webs}
 Let $n$  and  $d$ be fixed positive integers, with $d$ sufficiently big with respect to $n$ ({\it e.g.}\,$d\geq n+1$).
A `{\it $d$-web}' $\boldsymbol{\mathcal W}$ on a connected $n$-dimensional complex manifold $N$ is   a geometric object formed (locally) by a finite collection of foliations whose leaves intersect transversally (at least pairwise and at the generic point of $N$).  
We will mainly consider webs formed by global foliations on $N$, each of them  admitting a global meromorphic first integral $u_i : N\dashrightarrow \mathbf P^1$. In such a case, we will write 
$$\boldsymbol{\mathcal W}=\boldsymbol{\mathcal W}(u_1,\ldots,u_d)\, .$$  

Another web $\boldsymbol{\mathcal W}'$ on another connected complex manifold $M'$ is said to be `{\it equivalent}' to $\boldsymbol{\mathcal W}$ if there exists 
a biholomorphism 
$\varphi : U \rightarrow  U'$ between two open subsets $U,U'$ of $N$ and $N'$ respectively such that $\varphi^*(\boldsymbol{\mathcal W}')= \boldsymbol{\mathcal W}$ (possibly only as non-ordered webs). 
`{\it Web geometry}' is the subdomain of analytic geometry consisting in classifying or just studying (even global) webs up to this notion of (local) equivalence. \sk

An important notion to study webs in relation with classical algebraic geometry (see three paragraphs below) is the notion of `{\it abelian relation}' (abbreviated to AR).  For a web $\boldsymbol{\mathcal W}$ determined by 
$d$ first integrals $u_1,\ldots,u_d$, an abelian relation is a $d$-tuple of holomorphic functions $(F_1,\ldots,F_d)$ (each defined on a dense open subset of ${\rm Im}(u_i)\subset \mathbf P^1$, each considered up to the addition of a constant) such that the  following functional relation  holds true identically
$$
\sum_{i=1}^d F_i(u_i)=0\, .
$$ 
When the $u_i$'s are given, this is what we call an {\it Abelian functional equation} (in the $F_i$'s).\sk 

The space 
$\boldsymbol{\mathcal A}(\boldsymbol{\mathcal W})$ of the ARs of a web 
$\boldsymbol{\mathcal W}$ has a natural structure of complex vector space and its  dimension  
$\dim_{\mathbf C}\,\big( \boldsymbol{\mathcal A}(\boldsymbol{\mathcal W}) \big)$
is by definition the `{\it rank}' ${\rm rk}(\boldsymbol{\mathcal W})$ of $\boldsymbol{\mathcal W}$. It is an  invariant of this web, which can be proved to be always finite.
 \mk

The relevance of the notion of AR  becomes obvious if it is  considered in the case of a so-called `{\it algebraic web}' associated to a sufficiently generic degree $d$ projective curve $C$ in $\mathbf P^n$. This web, denoted by $\boldsymbol{\mathcal W}_C$ is the $d$-web  on the dual projective space 
$\check{\mathbf P}^n$ whose leaves are the hyperplanes tangent to the dual curve $C^*$.  
A nice consequence of Abel's addition theorem and of its converse in web geometry is that the ARs of
$\boldsymbol{\mathcal W}_C$ correspond to the abelian differentials of $C$ in a  very natural way. More precisely,  there is a  linear isomorphism $H^0(C,\omega_C^1)\rightarrow \boldsymbol{\mathcal A}(\boldsymbol{\mathcal W}_C)$ hence it follows that the arithmetic genus of $C$ 
coincides with the rank of $\boldsymbol{\mathcal W}_C$: one has $
p_a(C)=h^1(\mathcal O_C)={\rm rk}(\boldsymbol{\mathcal W}_C)$.  \sk 

It follows from classical results by Bol and Chern that
there exist universal bounds on the rank ${\rm rk}(\boldsymbol{\mathcal W})\leq \pi(d,n)$ satisfied for any $d$-web in dimension $n$ (with $d\geq n+1$) whose leaves satisfy a strong general position assumption denoted by (sGP)  (discussed in \S\ref{Par:GermOfWebs} below), where $ \pi(d,n)$ stands for \href{https://en.wikipedia.org/wiki/Castelnuovo_curve}{Castelnuovo's constant}.\footnote{By definition, for any $d,n$ such that $2\leq n<d$, $ \pi(d,n)$ is defined as the maximum of the genera of degree $d$ non-degenerate  irreducible projective curve in $\mathbf P^n$. Castelnuovo gave an explicit formula for this quantity.}  
The web $\boldsymbol{\mathcal W}$ is said to have `{\it maximal rank}' if ${\rm rk}(\boldsymbol{\mathcal W})=\pi(d,n)>0$. Determining whether a maximal rank web 
 is always `{\it algebraizable}', that is 
equivalent to an algebraic web 
such as described in the preceding paragraph, is an important problem in web geometry, as asserted by the following Chern and Griffiths' quote taken from the very end of \cite{CG2}: 
 \begin{quote}
{\it Before concluding, we cannot refrain from mentioning what we consider to be the fundamental problem on the subject, which is to determine the maximum rank non-linearizable\footnote{For webs of maximal rank, the properties of being non-algebraizable and non-linearlizable are equivalent (this follows from the converse to Abel's addition theorem, see \cite[Chap.\,4]{Coloquio}).} webs. The strong conditions must imply that there are not many. 
It may not be unreasonable to compare the situation with the exceptional simple Lie groups.}\sk 
\end{quote}

What makes this problem so interesting is the existence of the so-called `{\it exceptional webs}', 
that is non-algebraizable webs with  maximal rank.
Webs of this kind are very interesting since one can mimic for them several classical constructions of algebraic geometry,  which raises the following question: to which extent an `algebraic geometry of exceptional webs' might be developed. \sk

If exceptional webs do exist in dimension 2 (and there are many explicit examples of exceptional planar webs, see {\it e.g.}\,\cite{MPP,PP1,PirioTrepreau} and also \S\ref{SSS:Polylogs} just below), this is not the case for maximal rank webs satisfying (sGP) in dimension $n\geq 3$, as it follows from the remarkable algebraization theorem about them obtained by \href{https://en.wikipedia.org/wiki/Gerrit_Bol}{Bol} in dimension 3 and proved recently by Tr\'epreau in full generality (see \cite[Chap.\,5]{Coloquio} and the references therein). 
\begin{center}
\vspace{-0.2cm}$\star$
\end{center}

The purpose of the work undertaken here  is to show that  regarding the webs under scrutiny, if (1) one relaxes the general position hypothesis about how their leaves  are assumed to intersect (notion of `{\it weak general position}', {\it cf.}\,\S\ref{Par:GermOfWebs}) and (2) 
one considers a weaker but more flexible and still relevant condition regarding the maximality of their  rank (condition to have `{\it AMP rank}, see \S\ref{SS:Webs-AMP-Rank}), then one can define a less restrictive notion of exceptional webs which first of all coincides with the classical one in the case of planar webs but for which there are interesting examples in higher dimension as well. 
These examples are either webs naturally associated to functional equations in $k>3$ variables satisfied by higher polylogarithms or are webs which can be obtained  from cluster algebras and which also carry polylogarithmic ARs of weight 1 and 2 in most cases, but also of higher weight (namely 3 or 4) for some `cluster webs'. 

 \subsubsection{\bf Polylogarithms and their functional equations.}
 \label{SSS:Polylogs}

Polylogarithms are special functions which generalize the classical logarithm, this from several perspectives. For any integer $n\geq 1$, one defines the `{\it ${n}$-th (or weight $n$) polylogarithm}' $ \l {n} $ by means of the following series 
$$ 
\l {n} (z)= \sum_{k=1}^\infty \frac{z^k}{k^n}
$$
which converges for any $z$ in the unit disk $\mathbf D=\{\, \zeta \in \mathbf C\,  \hspace{0.1cm} \lvert \hspace{0.15cm} \lvert \zeta\lvert<1\, \}$.\sk

For $n=1$, one has $\l {1} (z)=- {\rm Log}(1-z)$ for $z\in \mathbf D$ hence the first polylogarithm essentially coincides with the classical logarithm. The latter extends as a multivalued holomorphic function on $\mathbf P^1$ with ramification at $0,1$ and $\infty$ and additive (hence unipotent) monodromy. And more interestingly for us, it satisfies the following classical identity 
$${\rm Log}(x)+{\rm Log}(y)-
{\rm Log}(xy)=0
$$ 
which gives an AR for the planar 3-web $\boldsymbol{\mathcal W}(x,y,xy)$ defined by the rational functions appearing as arguments of the logarithm in this functional relation.  Moreover, this web has maximal rank.\mk 

What makes polylogarithms interesting, in particular regarding web geometry,  is that the polylogarithms  of higher weight also satisfy nice functional identities which give rise to webs with interesting properties. The most emblematic case is that of the dilogarithm $ \l {2} $ which satisfies 
\begin{equation*}
 {}^{} \qquad 
\l{2} (x) -
\l{2} (y)
-\l{2} \left(\,\frac{x}{y}\, \right)
-\l{2} \left(\,\frac{1-y}{1-x}\,\right)+
\l{ 2 } \left(\,\frac{x(1-y)}{y(1-x)}\,\right)
={\rm Log}(y)\,{\rm Log}\left(\frac{1-y}{1-x}\right)
-\frac{{}^{}\, \, \pi^2}{6}
\end{equation*}
for any $x,y \in \mathbf R$ such that $0<x<y<1$, an identity which in the above form has been obtained by Abel, but whose equivalent forms were discovered by other authors in the XIXth and XXth centuries. This  identity  has a logarithmic second member but it admits an equivalent version without 
second member 
which gives rise to a dilogarithmic AR for the web defined by the five rational functions appearing as arguments
of  $ \l {2} $. This 5-web $\boldsymbol{\mathcal W}(x,y,x/y,(1-x)/(1-y),x(1-y)/(y(1-x)))$ was proved by Bol to be exceptional (of maximal rank 6 but not algebraizable) hence is now known as `{\it Bol's web}' and is accordingly denoted by $\boldsymbol{\mathcal B}$.\mk 

Things are similar for the trilogarithm $ \l {3} $ :  it satisfies the following `{\it Spence-Kummer identity}'
\begin{align} \qquad  2 \,\l {3} (\,x\,) &  \, +  \, 2  \,\l{3} (\,y\,)  \, -
  \,  \l{3} \left(\, \frac{x}{y}\, \right)\,     
+ \, 2  \, \l{3}  \left(\,  \frac{1-x}{1-y}\, \right) \,   + \,  2 \,\l{3} \left(\,
  \frac{x(1-y)}{y(1-x)}\, \right)    
- \, \l{3} (\, xy\,)  \qquad \qquad \qquad \qquad \nonumber \\
  \qquad \qquad \qquad   &  \, +   \,  2 \,
 \l{3} \left(\, -\frac{x(1-y)}{{\,}(1-x)}  \,         \right) 
+ \,  2 \,   \l{3} \left(\,  - \frac{{\,}(1-y)}{y(1-x)}\, \right)   -  \, \l{3}
  \left(\, \frac{x(1-y)^2}{y(1-x)^2}\, \right) 
 \qquad      \nonumber \\
  &=2\, 
  \zeta(3)
+
\frac{{}^{}\, \, \pi^2}{3}
\, {\rm Log} (y)
  -{\rm Log}(y)^2\,{\rm Log}\, \left( \frac{1-y}{1-x}\right) 
+\frac{1}{3}\, {\rm Log}(y)^3  \, ,  \nonumber 
\end{align} 
for all real numbers  $x$ and $y$ such that $0<x<y<1$. This identity gives rise to a trilogarithmic AR for the 
so-called `{\it Spence-Kummer web}  $\boldsymbol{\mathcal W}(x,y,x/y,(1-x)/(1-y), \ldots,
-(1-y)/(y(1-x)),  x(1-y)^2/(y(1-x)^2))$
admitting as first integrals the nine rational functions appearing as arguments of $ \l {3} $ in the above identity. In this case again, this web is exceptional.

\newpage
Considering the examples  in weight 1,2 and 3 above suggests to ask as Griffiths in \cite{Griffiths} 
 \begin{equation}
 \label{Griffiths-Question}
\begin{tabular}{l}
{\it whether or not for each $n$ there is an integer $d_n$ such that there is a ``new'' $d_n$-web of}\\
 {\it maximum rank one of whose abelian relations is a (the?) functional equation with $d_n$}\\ {\it   terms for the $n$-th polylogarithm $ \l {n}  $?}
\end{tabular} 
\end{equation}

The answer to this question is still not known, the reason explaining this being certainly that the theory of multivariable functional equations satisfied by higher polylogarithms  is still a mystery.  For instance, identities for $\l {n} $ similar to Abel's and Spence-Kummer's ones in weight 2 and 3 respectively, are known only for small values of $n$ ($n\leq 7$, see \S\ref{S:EFA-polylog} for more on this subject).   Moreover, the planar webs associated to the seemingly most basic and simple functional identities satisfied by $\l {4} $ and $ \l {5}  $, discovered by Kummer, do not have maximal rank hence are not exceptional. Even worse, no exceptional web carrying a  weight $n$ polylogarithmic  AR is known for $n\geq 4$, which makes us think that the answer to the question just above might well be negative.
\begin{center}
\vspace{-0.2cm}$\star$
\end{center}

In this memoir, we consider many functional identities in $k\geq 3$ variables satisfied by polylogarithms of weight $n =2,3, 4$ and we show that for some of them, the associated webs have AMP rank. 
We expect that most if not all of these webs are not algebraizable hence are generalized exceptional webs. We will prove the later property for an interesting series of cluster webs, one in $n$ variables for each $n\geq 2$.

 \subsubsection{\bf Cluster algebras and webs associated to them.}
 \label{SSS:ClusterAlgebras}
A '\href{https://en.wikipedia.org/wiki/Cluster_algebra}{\it Cluster algebra}'  of rank $m\geq 1$ is a subalgebra of ${\bf k}(\boldsymbol{u})={\bf k}(u_1,\ldots,u_m)$ (where $\bf k$ stands for a fixed field, {\it e.g.}\,$\mathbf Q$ or $\mathbf C$) generated by the union of certain subsets of cardinality $m$, called clusters, which all can be obtained inductively from a given initial cluster, via birational transformations called `mutations' given by a combinatorial recipe. Before being a bit more precise, let us say that there are actually two kind of clusters and mutations, say of type $\boldsymbol{\mathcal A}$ and $\boldsymbol{\mathcal X}$ respectively. We will not discuss the former here, first to simplify the exposition but mainly because it is the latter which is relevant regarding the construction of webs with polylogarithmic abelian relations.
\sk 

A ($\boldsymbol{\mathcal X}$-)seed is a pair $S=(\boldsymbol{x}, B)$ where $\boldsymbol{x}=(x_i)_{i=1}^m$ is a $m$-tuple of elements of ${\bf k}(\boldsymbol{u})$ (usually assumed to be algebraically independent) and $B=(b_{ij})_{i,j=1}^m$ is a $m$ skew-symmetrizable square matrix of size $m$ with integer  coefficients. The tuple $\boldsymbol{x}$ is the `{\it cluster}' associated to $S$, its components $x_i$'s are the corresponding `{\it cluster variables}' whereas $B$ is the `{\it exchange matrix}' of the seed. 
\sk

The notion of ($\boldsymbol{\mathcal X}$-)mutation is one of the most fundamental in the theory considered here, and allows to construct other seeds from $S$. More precisely, for any $k=1,\ldots,m$, the `{\it $k$-th mutation of $S$ in the $k$-th direction}' is the web seed $
\mu_k(S)=S'=
(\boldsymbol{x}', B')$ where the new  cluster $\boldsymbol{x}'=(x_i')_{i=1}^m$ and new exchange matrix $B'=(b'_{ij})_{i,j=1}^m$ are obtained from $x$ and $B$ by means of the following formulas
for $i=1,\dots,m$, where $[\cdot ]_+: \mathbf R \rightarrow \mathbf R^+$ stands for the function $x\mapsto \max ( 0,x)$: 
\begin{equation*}
\label{Eq:B-X-Mutation-formulae}
x_{i}'=\begin{cases}
{}^{} \hspace{0.1cm} {x_{k}}^{-1} \hspace{2.35cm} \mbox{ if }\,  \,  i=k   \\
{}^{} \hspace{0.1cm} 
x_i\,{\Big(1+{x_k}^{[- b_{ki}]_+}\Big)}^{-b_{ki}}  
 \hspace{0.1cm} \mbox{ if }\,  \,  i\neq k 
\end{cases}
 \mbox{and} \quad 
b_{ij}'=\begin{cases}
\, - b_{ij} \hspace{1.8cm} \mbox{ if }\,  \,  k \in \{i,j\} \, ,  \\
\, \hspace{0.2cm} b_{ij}\hspace{1.85cm} \mbox{ if } \,  k\not \in \{i,j\}  \, \mbox{  and   } \,   b_{ik}b_{kj}\leq 0\, , \\
\, \hspace{0.2cm} b_{ij}+\lvert b_{ik}\lvert \, b_{kj}\hspace{0.3cm} \mbox{ if } \,  k\not \in \{i,j\} \, \mbox{ and  }\,   b_{ik}b_{kj}> 0\, . 
\end{cases}
\end{equation*}
\newpage

 Mutations 
(in the same direction) are involutions since it can be verified that for any $k$, one has 
$\mu_k'(\mu_k(S))=S$ where $\mu_k'$ stands for the mutation of $S'=\mu_k(S)$ in the direction $k$, but now with respect to the muted exchange matrix $B'$.   \sk

Given an initial seed $S^0=( \boldsymbol{x}, B^0)$, by considering 
all the seeds obtained by successive mutations from it, one can construct the associated `{\it cluster exchange pattern}', which is the countable family of seeds  $S^t=( \boldsymbol{x}^t, B^t)$ indexed by the vertices $t$ of the $m$-regular tree $\mathbb T^m$, with the initial seed 
$S^0$ associated to the
root $t_0$ of this tree, the edges of  $\mathbb T^m$ being labeled by elements of $\{1,\ldots,m\}$ in the natural consistent way.  The components of the clusters $\boldsymbol{x}^t$ are the ($\boldsymbol{\mathcal X}$-)cluster variables of the cluster algebra defined by $S^0$. 
Each  cluster variable $x_i^t$ is a Laurent polynomial in the initial cluster variables
$x_i$'s, which moreover enjoys several other remarkable properties (separation formula, sign-coherence, positivity). \mk 

Given a finite set $\Sigma$ of cluster variables, one can consider the web 
${\boldsymbol{\mathcal W}}_\Sigma$ formed by the foliations admitting the elements of $\Sigma$ as first  integrals.  It is what we call a `{\it cluster web}' and webs of this kind   are among the main objects studied in this memoir.  A nice feature of the theory of cluster algebras is that it comes with several ways to get finite sets $\Sigma$ of cluster variables giving rise to webs carrying polylogarithmic ARs. 
The first way is given by the so called cluster algebras  of `{\it finite type}', which are by definition those admitting only a finite number of clusters. An early fundamental result of Fomin and Zelevinsky shows that the classification of such algebras is parallel to that of Dynkin diagrams. Thus given a Dynkin diagram $\Delta$ of rank $m\geq 2$, the set $\Sigma_\Delta$ of all the $\boldsymbol{\mathcal X}$-cluster variables of the corresponding cluster algebra is finite, which allows us to define the `{\it cluster web of Dynkin type $\Delta$}' as the cluster web
$$
{\boldsymbol{\mathcal X\hspace{-0.05cm} \mathcal W}}_{\hspace{-0.05cm} \Delta}=
{\boldsymbol{\mathcal W}}_{\Sigma_\Delta}\, . 
$$

The first  case to be considered is when $\Delta=A_2$:   
the ${\boldsymbol{\mathcal X}}$-cluster web of type $A_2$ is
$$
{\boldsymbol{\mathcal X\hspace{-0.05cm} \mathcal W}}_{\hspace{-0.05cm} A_2}
={\boldsymbol{ \mathcal W}}\bigg(\, 
x_1 \, , \, 
x_2 \, , \, 
\frac{1+x_2}{x_1}\, , \, 
\frac{1+x_1}{x_2}\, , \, 
\frac{1+x_1+x_2}{x_1x_2}\,\bigg)\, .
$$
 It is a planar 5-web carrying a dilogarithmic AR which is known to be equivalent  to Bol's web  associated to Abel's 5 terms identity.  Thus the $\boldsymbol{\mathcal X}$-cluster web of type $A_2$ is exceptional, which gives a strong motivation to study the ARs and the rank of the 
$\boldsymbol{\mathcal X}$-cluster web ${\boldsymbol{\mathcal X\hspace{-0.05cm} \mathcal W}}_{\hspace{-0.05cm} \Delta}$ 
for $\Delta$ arbitrary. 
\sk 

\vspace{-0.25cm}
To each seed $S^t$ is naturally attached a `{\it cluster torus}' $\mathbf T^t$. These tori can be glued birationally via the mutation transformations $\mu_k$ and give rise to the so-called `{\it cluster variety}', that we will denote by $\boldsymbol{\mathcal X}_\Delta$ in case of a finite type cluster algebra of Dynkin type $\Delta$.  The web ${\boldsymbol{\mathcal X\hspace{-0.05cm} \mathcal W}}_{\hspace{-0.05cm} \Delta}$ naturally lives on $\boldsymbol{ \mathcal X}_\Delta$ and restricting this web along a certain subvariety $\boldsymbol{ \mathcal U}_\Delta\subset 
\boldsymbol{ \mathcal X}_\Delta $ known as the `{\it secondary cluster variety}', one 
gets the so-called `{\it secondary cluster web}' of Dynkin type $\Delta$, denoted as follows
$$
{\boldsymbol{\mathcal U\hspace{-0.05cm} \mathcal W}}_{\hspace{-0.05cm} \Delta}=
\big( {\boldsymbol{\mathcal X\hspace{-0.05cm} \mathcal W}}_{\hspace{-0.05cm} \Delta}\big) \big\lvert_{
\boldsymbol{\mathcal U}_\Delta}
\,
$$ and which will be proven to be interesting, at least for certain Dynkin diagrams $\Delta$.
 \begin{center}
\vspace{-0.2cm} 
$\star$
\end{center}
Another way to construct interesting cluster webs, which is relevant even within cluster algebras of infinite type, relies on the notion of `{\it cluster period}' formalized by Nakanishi. This is a finite sequence $\boldsymbol{i}=(i_1,\ldots,i_k)\in \{ 1,\ldots,m\}^k$ for some $k\geq 2$, giving rise to a sequence of seeds defined inductively by $S_{\hspace{-0.04cm}\ell+1}=\mu_{i_\ell}(S_{\hspace{-0.04cm} \ell})$ for $\ell=0,\ldots,k-1$, with the crucial property that $S_k$ is isomorphic to the initial seed $S_0$.  A remarkable result proved in full generality by Nakanishi asserts that if $x_{\boldsymbol{i}}(\ell)$ stands for the $i_\ell$-th cluster variable of the $\ell$-th seed $S_{\hspace{-0.04cm}\ell}$, then an identity 
$$ \big({\sf R}_{\boldsymbol{i}} \big)\hspace{3cm} \sum_{\ell=1}^k d_\ell\,{\sf R}\big(x_{\boldsymbol{i}}(\ell)\big)=N_{\boldsymbol{i}}\, \pi^2/6
\hspace{4cm} {}^{}
$$
  holds true, for some positive  integers $d_1,\ldots,d_\ell $ and $N_{\boldsymbol{i}}$, where ${\sf R}$ stands for a cluster version of the dilogarithm. By definition, the `{\it cluster web associated to $\boldsymbol{i}$}', denoted by  ${\boldsymbol{\mathcal W}}_{\hspace{-0.05cm} \boldsymbol{i}}$, is the web admitting the $x_{\boldsymbol{i}}(\ell)$'s 
 as first integrals: one has ${\boldsymbol{\mathcal W}}_{\hspace{-0.05cm} \boldsymbol{i}}=
 {\boldsymbol{\mathcal W}}(x_{\boldsymbol{i}}(1),\ldots, x_{\boldsymbol{i}}(k))$. 
Clearly, the identity $\big({\sf R}_{\boldsymbol{i}} \big)$ gives rise to a dilogarithmic AR for this web, which therefore appears as worth studying from this point of view. Indeed, we show in this text that many cluster webs associated to cluster periods give rise to cluster webs with many polylogarithmic ARs, of weight 1 or 2. \mk 

An important class of examples of webs associated to cluster periods are those given by the 
 so-called $Y$-systems of mathematical physics. For a pair  of Dynkin diagrams $(\Delta,\Delta')$ of rank $m$ and $m'$ respectively, several authors  have proved the periodicity of the $Y$-system of bi-Dynkin type $(\Delta,\Delta')$, a result due to Keller in full generality. 
This is  equivalent to the fact that a certain sequence $\boldsymbol{i}_{\Delta,\Delta'}$ of mutations is a cluster period for the cluster algebra of bi-Dynkin type $(\Delta,\Delta')$.  
The associated cluster dilogarithmic identity is denoted by $({\sf R}_{\Delta,\Delta'})$ and 
 the $\boldsymbol{i}_{\Delta,\Delta'}$-cluster web, which we call the 
`{\it ${\boldsymbol{\mathcal Y}}$-cluster web of bi-Dynkin type $(\Delta,\Delta')$}' is denoted by ${\boldsymbol{\mathcal Y\hspace{-0.05cm} \mathcal W}}_{\hspace{-0.05cm} \Delta, \Delta'}$. It is a  web in $mm'$ variables which carries the complete dilogarithmic AR corresponding to the identity $({\sf R}_{\Delta,\Delta'})$. \sk 

When $\Delta'=A_1$, we drop it in the notation and  just write ${\boldsymbol{\mathcal Y\hspace{-0.05cm} \mathcal W}}_{\hspace{-0.05cm} \Delta}$
 and $({\sf R}_{\Delta})$ respectively. In this case, ${\boldsymbol{\mathcal Y\hspace{-0.05cm} \mathcal W}}_{\hspace{-0.05cm} \Delta}$ is a subweb of ${\boldsymbol{\mathcal X\hspace{-0.05cm} \mathcal W}}_{\hspace{-0.05cm} \Delta}$, with which it coincides when  $\Delta$ has rank 2. In particular, one has ${\boldsymbol{\mathcal Y\hspace{-0.05cm} \mathcal W}}_{\hspace{-0.05cm} A_2}={\boldsymbol{\mathcal X\hspace{-0.05cm} \mathcal W}}_{\hspace{-0.05cm} A_2}$ which is of maximal rank. One of the main results obtained in this memoir is that this extends to the whole family 
 of ${\boldsymbol{\mathcal Y}}$-cluster webs of type $A$: for any $n\geq 2$, although  it is non-linearizable, the 
 cluster web  ${\boldsymbol{\mathcal Y\hspace{-0.05cm} \mathcal W}}_{\hspace{-0.05cm} A_n}$ is AMP with only logarithmic ARs modulo  the dilogarithmic one associated to $({\sf R}_{A_n})$.


 \subsection{\bf 
 Main  results contained in and main features of this text}
\label{SS:Intro-Results}

We now discuss with a bit more detail 
 the main new results presented in this memoir.  
 
\subsubsection{}
In what concerns the general theory of webs, our main contribution is to define a new invariant attached  to quite many webs, the `{\it virtual rank}', which enjoys several nice properties, in particular that of giving a bound on the actual rank. More precisely, let $\boldsymbol{\mathcal W}$ be a $d$-web defined by first integrals $u_1,\ldots,u_d$, defined on a domain $\Omega$ of $\mathbf C^n$ say, with $n\geq 2$. We only require that   the associated foliations are pairwise transverse, that is analytically: one has $du_i\wedge du_j\neq 0$ on $\Omega$ if $i\neq j$. For any $\omega\in \Omega$, the first jet of $\boldsymbol{\mathcal W}$ at this point is the web  which can be identified with the linear web on the whole $\mathbf C^n$ defined by the differential forms $\ell_i=\ell_{i,\omega}=du_i(\omega)$ for $i=1,\ldots,d$.  For any $\sigma\geq 1$, taking $\omega$ generic in $\Omega$,  one defines the $\sigma$-th  virtual (resp.\,the total) rank $\rho^\sigma (\boldsymbol{\mathcal W})$  
 (resp.\,$\rho(\boldsymbol{\mathcal W})$) of $\boldsymbol{\mathcal W}$ by
$$
\rho^\sigma\big(\boldsymbol{\mathcal W}\big)=\dim_{\mathbf C}\, 
\Big\langle \, (c_i)_{i=1}^d\in \mathbf C^d 
\hspace{0.1cm}
 \big\lvert \hspace{0.1cm}
\sum_{i=1}^d  c_i \big( \ell_i\big)^\sigma = 0 
\hspace{0.1cm}
\Big\rangle  
\qquad \mbox{and}\qquad 
\rho(\boldsymbol{\mathcal W})=\sum_{\sigma\geq 1} \rho^\sigma\big(\boldsymbol{\mathcal W}\big)\, .
$$

These virtual ranks are well-defined, finite, and invariantly attached to $\boldsymbol{\mathcal W}$. They are easy to compute (when the $u_i$'s are given) and they satisfy many other nice properties, the main one being the following upper bound on the rank
\begin{equation}
\label{Eq:r=rho}
{\rm rk}\big(\boldsymbol{\mathcal W} \big) \leq  \rho(\boldsymbol{\mathcal W})\,,  
\end{equation}
a majoration which applies even to webs not satisfying the `strong general position assumption' classically required in web geometry.   When \eqref{Eq:r=rho} actually is an equality, the rank of  $\boldsymbol{\mathcal W}$ is said to be `{\it As Maximal as Possible}' and we will say that this web `{\it is AMP}'.  Albeit defined by elementary means, this notion of `{\it web with AMP rank}'  is new and  encompasses all the notions of `maximal rank webs' previously considered in web geometry. \mk 

We then state some general results about the `generalized webs' we are dealing with in this text, some in relation with the new notion of `being AMP' previously introduced, while some others are natural and easy generalizations of some well-known results holding true for classical (say planar) webs. \sk

First, for algebraic webs, it is not difficult to  relate the web-theoretic property of being AMP to a classical notion  of the theory of  algebraic curves, namely that  a curve is `Arithemically Cohen-Macaulay' (ACM): 
\begin{prop}
Given a reduced algebraic curve $C\subset \mathbf P^n$, the following equivalence holds true: \vspace{0.15cm}

\vspace{-0.4cm}
\begin{tabular}{r}
${}^{}$\hspace{2.6cm}
the algebraic web 
${\boldsymbol{\mathcal W}}_C$ is AMP
\hspace{0.2cm}
 $\Longleftrightarrow$
 \hspace{0.2cm}
the projective curve $C$ is ACM.
\end{tabular}
 \end{prop}
If we will not elaborate on this result in this text, we mention it here because it shows that the theory of AMP webs admits, as a special subcase, the theory of ACM projective curves, which has been and still is an important active field of research in algebraic geometry. Since there exist many examples of AMP webs which are not algebraizable, it follows that the theory of such webs is a strict generalization of that of ACM projective curves. 

%

Next we show that the main properties of the ARs  as well as the 
now classical methods to determine them and to compute the rank in the case of planar webs extend quite naturally to generalized webs in arbitrary dimension. 
To simplify the exposition below, we will deal with a web denoted by 
$\boldsymbol{\mathcal W}$,  defined by $d$ rational first integrals 
$u_1,\ldots,u_d\in \mathbf C(x_1,\ldots,x_n)$ in $n\geq 2$ variables. \sk

We denote by $\Sigma^c(\boldsymbol{\mathcal W})$ the union of the irreducible divisors in $\mathbf P^n$ which are invariant by at least two foliations of $\boldsymbol{\mathcal W}$ and we define finite subsets of $\mathbf P^1$ by setting  $\mathfrak B_i=
 u_i\big(\Sigma^c(\boldsymbol{\mathcal W})\big)$ for $i=1,\ldots,d$. 

 \begin{prop} 
 \label{P:AnalyticConinuation}
 For any $\omega\in \mathbf P^n\setminus \Sigma^c(\boldsymbol{\mathcal W})$ and any (germ of)  abelian functional relation $\sum_i F_i(u_i)=0$ at $\omega$, each germ $F_i$ extends as a  global multivalued holomorphic function on  $\mathbf P^1\setminus \mathfrak B_i$. 
\end{prop}

This a slight generalization of a well-known result holding true for planar webs  (with essentially the same proof) which allows to restrict the study of the ARs and of the rank of the web under scrutiny in the vicinity of any  point $\omega$ outside $\Sigma^c(\boldsymbol{\mathcal W})$. \sk 

After having discussed some ways to construct or even to determine the ARs of a given web (symbolic determination of ARs with iterated integrals as components and the so-called Abel's method for solving AFEs), we show that Pantazi-H\'enaut criterion 
 for characterizing planar webs of maximal rank 
  actually extends rather straightforwardly to webs in arbitrary dimension since the following proposition holds true: 
\begin{prop}
Let $\sigma({{\boldsymbol{\mathcal{ W}}}})$ be the biggest integer $\sigma\geq 1$ such that 
$\rho^\sigma({\boldsymbol{\mathcal W}})>0$.  Then 
$\boldsymbol{\mathcal{ W}}$ is AMP if and only if it  has rank $\rho(\boldsymbol{\mathcal W})$ at the order $\sigma({{\boldsymbol{\mathcal{ W}}}})+2$ at a generic point of $\Omega$. 
\end{prop}

Actually, this result admits a rather straightforward generalization which gives an effective tool to compute the rank of the web $\boldsymbol{\mathcal W}$ when the $u_i$'s are explicitly given. 

 \subsubsection{} Let us now discuss what we have obtained regarding polylogarithms.  
 \sk

In Part I, we discussed  the symbolic method to determine the ARs of webs whose components are iterated integrals. In the second part, we specialize to the polylogarithmic case,  that is the one with ramification points $0$, $1$ or $\infty$.  This material is essentially already known, but we give a more detailed treatment of it, in the suitable setting for our purpose, namely the holomorphic one.  

Let   $\zeta\in \mathbf P^1$ be a fixed base-point distinct from $0,1$ and $\infty$ and  set  $\omega_0=du/u$ and $\omega_1=du/(1-u)$. 
Then for any `symbol $w$', that is a word in the two letters $0$ and $1$, we set inductively 
${\mathfrak{L}}_{\emptyset}^{\zeta}=1$ and  ${\mathfrak{L}}_{w}^{\zeta}( \cdot)=
\int_{\zeta}^{\,\cdot} {\mathfrak{L}}_{w'}^{\zeta}(u)\,\omega_{\varepsilon}$  if $w=\varepsilon \, w'$ with $\varepsilon \in \{0,1\}$.  For $n\geq 2$ fixed,  we define  
${\mathfrak{L}}_{n}^{\zeta}$ as the modified polylogarithm  at $\zeta$ whose symbol is $0^{n-2}(01-10)$ and  for $k\leq n-2$, one sets 
${\mathfrak{L}}_{n,k}^{\zeta}={\mathfrak{L}}_{0^k}^{\zeta}\cdot {\mathfrak{L}}_{n-k}^{\zeta}$.
\mk 

As above,  let $u_1,\ldots,u_d$ be rational functions in $m\geq 2$ variables defining a (singular) $d$-web $\boldsymbol{\mathcal W}=\boldsymbol{\mathcal W}(u_1,\ldots,u_d)$ on $\mathbf C^m$.
For any word $w$ as above and any $i$, $u_i^*w$ stands for the  symbol corresponding to $w$ but where $\omega_s$  has been replaced by the pull-back $u_i^*(\omega_s)$  for $s=0,1$.  For $\zeta\in \mathbf C^n$ generic, we set $\zeta^i=u_i(\zeta)$ and we use  ${\rm P}_{\hspace{-0.06cm} <n}$ to denote any rational expression evaluated on tuples of polylogarithms $\l {n'} $ of weight $n'<n$,  themselves evaluated on rational functions. 
\begin{thm}
\label{Thm:PolylogAFE}
 For any non-zero tuple  $\boldsymbol{c}=(c_i)_{i=1}^d\in \mathbf C^d$,  the following assertions are equivalent: 
\begin{enumerate}
\item  a functional equation of the form  $\sum_{i=1}^d c_i \, { {\bf L}{\rm i}_{ n } }  (u_i)={\rm P}_{\hspace{-0.06cm} <n}$ holds true;
\item     the symbolic identity  $
\sum_{i =1}^d c_i\cdot 
 u_i^*\Big[ 0^{n-2}(01-10)  \Big]=0
$ is satisfied;
\item    the sum 
$ \sum_{i=1}^d c_i \,{\mathfrak{L}}_{n}^{\zeta_i}(u_i)$ vanishes identically on a neighbourhood of $\zeta$ ;
\item  for any $k=0,\ldots,n-2$,  one has $ \sum_{i=1}^d c_i \,{\mathfrak L}_{n,k}^{\zeta_i}(u_i)\equiv 0$ on a neighbourhood of $\zeta$. 
  \end{enumerate}
  When these assertions hold true for some $\boldsymbol{c} \in \mathbf  C^d\setminus \{0\}$, the $n-1$ identities mentioned in 4. furnish as many linearly independent polylogarithmic abelian relations for the web defined by the $u_i$'s.
\end{thm}
From this result, one deduces that webs associated to polylogarithmic identities 
$\sum_{i=1}^d c_i \, { {\bf L}{\rm i}_{ n } }  (u_i)={\rm P}_{\hspace{-0.06cm} <n}$ carry polylogarithmic ARs of weight $n'=1,\ldots,n$ hence a priori have high rank hence  are good candidates for being AMP.  From the fact that ${\mathfrak{L}}_{n}^{\zeta}$ has valuation $n$ at $\zeta$, one can also easily deduce the following general but interesting result: 
\begin{cor}
\label{Cor:cor}
 Let $
  \sum_{i=1}^d c_i\,  \l { n } (u_i)={\rm P}_{\hspace{-0.05cm}<n}$ be a non-trivial polylogarithmic FE  
 in several   variables. If the $u_i$'s appearing in the LHS define $d$ distinct foliations then 
%
  necessarily 	
 $d\geq n+3$.
\end{cor}

Actually, using not only  ${\mathfrak{L}}_{n}^{\zeta}$ but the $n-1$ functions 
 ${\mathfrak{L}}_{n,k}^{\zeta}$ for $k=0,\ldots,n$, 
   one should be able to improve the bound in this theorem (and obtain that  necessarily $d\geq 2n+1$ in the conclusion). 
\mk 

After these generalities about the polylogarithms and the AFEs they satisfy, we turn to considering the webs attached to several (one could say many) explicit  polylogarithmic functional identities, in weight $2,3$ and $4$. Using the notions and tools discussed in the first part on webs, we can study quite effectively these webs to get the 

\begin{thm}
\label{T:classical-cluster-webs}
Many of the webs associated to classical or more recent known AFEs satisfied by  polylogarithms, such as for instance the ones  given in Table 
 \ref{T:MRP-1923-2021} just below,  have AMP rank.
\end{thm}

\begin{table}[!h]
\begin{center}
\begin{tabular}{|l|c|c|}
\hline
${}^{}$  \hspace{0.2cm} {\bf Polylogarithmic identity}
& {\bf Weight}    &
{\bf Reference}
    \\
\hline \hline 
${}^{}$ \hspace{0.2cm} Abel's identity $\boldsymbol{({\mathcal A}b)}$ & 2  &   
\S\ref{Par:Abel'sAFE}
  \\ \hline  
${}^{}$ \hspace{0.2cm} Newman's identity $\boldsymbol{({\mathcal N}_6)}$& 2  &    
\S\ref{Par:Newman-equation}
     \\ \hline 
${}^{}$ \hspace{0.2cm} Mantel's identity $\boldsymbol{({\mathcal M})}$& 2  &    
\S\ref{Par:Mantel}
     \\ \hline 
${}^{}$ \hspace{0.2cm} Maier's identity $\boldsymbol{({\mathcal M}_8)}$& 2  &    
\S\ref{Par:Maier's-Web}
     \\ \hline  \hline
${}^{}$ \hspace{0.2cm} Spence-Kummer's identity $\boldsymbol{({\mathcal S \mathcal K})}$ & 3 & 
\S\ref{SS:Spence-Kummer} 
 \\ \hline 
 ${}^{}$ \hspace{0.2cm} Goncharov's identity $\boldsymbol{({\mathcal G}_{22})}$ & 3 & 
\S\ref{Par:Goncharov-22} 
 \\ \hline 
  ${}^{}$ \hspace{0.2cm} Gangl's identity $\boldsymbol{({\mathcal G an} _{21})}$ & 3 & 
\S\ref{Parag=Gan-21}
 \\ \hline  \hline
${}^{}$ \hspace{0.2cm} Tetralogarithmic identity $\boldsymbol{(
{\mathcal G \mathcal G \mathcal S \mathcal V\mathcal V})}$ & 4 & 
\S\ref{Par:GSVV}
\\ \hline 
\end{tabular}
\caption{Examples of polylogarithmic identities giving rise to AMP webs.} 
\label{T:MRP-1923-2021}
\end{center}
\end{table}

We also consider series of dilogarithmic identies, such as Zagier's identities 
${({\mathcal Z}(m,n))}$  of the recently discovered 
Bytsko-Volkov's identities ${({\mathcal B}{\mathcal V}(n))}$
(see \S\ref{SPar:Z(m,n)} and \S\ref{SPar:Bytsko-Volkov}  respectively).  
We have verified that the associated webs  are AMP for  $m$ and $n$ small and we conjecture that this holds true for the whole series.

 \subsubsection{} 
In a third direction, we study webs which can be constructed from cluster algebras.\sk

Our first results consist in linking some webs associated to classical polylogarithmic functional equations with cluster algebras.  As far as we know, this was previously known only  for Abel's five-terms dilogarithmic identity, which has been related  quite early to the finite type cluster algebra of type $A_2$.  We prove that a similar phenomenon occurs  for some other classical polylogarithmic identities as well. We shall say that a polylogarithmic identity is of `{\it cluster type}' if it is equivalent to a polylogarithmic identity 
 carried by a cluster web. 

\begin{thm}
\label{T:classical-cluster-webs}
For $n=2,3$ or $4$, the following classical identities satisfied by 
$ \l {n} $ 
 are of cluster type and can be constructed from finite type cluster webs associated to certain Dynkin diagrams: Abel's 5-terms and Newman's dilogarithmic identities,  Spence-Kummer trilogarithmic identity and Kummer's tetralogarithmic identity.  
The corresponding cluster webs associated to each of these identities are given in the table below: 
\begin{table}[!h]
\begin{center}
\begin{tabular}{|l|c|c|l|}
\hline
${}^{}$  \hspace{0.2cm} {\bf Polylogarithmic identity}
& {\bf Weight} $\boldsymbol{n}$    &
{\bf Cluster web} & {\bf Reference(s)}
    \\
\hline \hline 
${}^{}$ \hspace{0.2cm} Abel's identity $\boldsymbol{({\mathcal A}b)}$ & 2  &    ${\boldsymbol{\mathcal X \hspace{-0.03cm}\mathcal W}}_{ \hspace{-0.03cm} A_2}$  & 
 \hspace{0.14cm}
  \S\ref{SS:cluster-webs-type-A2}
 \\ \hline 
${}^{}$ \hspace{0.2cm} Newman's identity $\boldsymbol{({\mathcal N}_6)}$& 2  &    ${\boldsymbol{\mathcal X \hspace{-0.03cm}\mathcal W}}_{ \hspace{-0.03cm} B_2}$   &
 \hspace{0.14cm}  \S\ref{SS:cluster-webs:B2}
  \\ \hline 
${}^{}$ \hspace{0.2cm} Spence-Kummer's identity $\boldsymbol{({\mathcal S \mathcal K})}$ & 3 & ${\boldsymbol{\mathcal U \hspace{-0.03cm}\mathcal W}}_{ \hspace{-0.03cm} A_3}$ & 
 \hspace{0.14cm} \S\ref{SS:ClassicalPolylogarithmicIdentitiesAreOfClusterType}, 
\eqref{Eq:pXWA3}
 \\ \hline 
${}^{}$ \hspace{0.2cm} Kummer's identity $\boldsymbol{({\mathcal K}_4)}$ & 4 & ${\boldsymbol{\mathcal U \hspace{-0.03cm}\mathcal W}}_{ \hspace{-0.03cm} D_4}$ 
&  \hspace{0.14cm}  \S\ref{SS:ClassicalPolylogarithmicIdentitiesAreOfClusterType}, \eqref{Eq:pXWD4}
\\
\hline 
\end{tabular}
\caption{Some polylogarithmic identities of cluster type.} 
\label{T:Ex-of-cluster-webs}
\end{center}
\end{table}
\end{thm}

If the preceding result concerns specific cluster webs, we have tried to study 
 quite systematically some series of webs which can be constructed from the cluster algebras of finite type, namely the cluster webs 
 ${\boldsymbol{\mathcal X \hspace{-0.03cm}\mathcal W}}_{ \hspace{-0.03cm} \Delta}$, ${\boldsymbol{\mathcal Y \hspace{-0.03cm}\mathcal W}}_{ \hspace{-0.03cm} \Delta}$ or ${\boldsymbol{\mathcal Y \hspace{-0.03cm}\mathcal W}}_{ \hspace{-0.03cm} \Delta \square \Delta'}$ for any Dynkin diagrams $\Delta$ and $\Delta'$.   An extensive computational investigation, by means of the effective tools for studying webs mentioned above, led us to get quite precise, but in full generality still conjectural, statements about the members of these families of cluster webs, and more specifically regarding their ARs and their rank(s). \sk

The main result of \cite{Pereira} can be interpreted (modulo the  standard identification between the cluster variety $\mathcal X_{A_n}$ with $\mathcal M_{0,n+2}$) as the fact that for any $n\geq 2$, the cluster web  ${\boldsymbol{\mathcal X \hspace{-0.03cm}\mathcal W}}_{ \hspace{-0.03cm} A_n}$ is AMP with only polylogarithmic abelian relations of weight 1 or 2. We use this to prove that the same holds true for its $\mathcal Y$-subweb ${\boldsymbol{\mathcal Y \hspace{-0.03cm}\mathcal W}}_{ \hspace{-0.03cm} A_n}$. More precisely, we will prove the following result: 
 \begin{thm}
 \label{T:YW-An-intro}
 For any $n\geq 2$,  ${\boldsymbol{\mathcal Y\hspace{-0.05cm} \mathcal W}}_{\hspace{-0.05cm} A_n}$ is a ${n(n+3)}/{2}$-web in $n$ variables such that
$$
 \rho^\bullet\Big({\boldsymbol{\mathcal Y\hspace{-0.05cm}\mathcal W}}_{\hspace{-0.05cm}A_n}\Big)=
\left(\,  \frac{n(n+1)}{2} 
  \, ,\, n\, ,\, 1\, \right)  \qquad \mbox{and} \qquad 
{\rm polrk}^\bullet\Big(
{\boldsymbol{\mathcal Y\hspace{-0.05cm}\mathcal W}}_{\hspace{-0.05cm}A_n}
\Big)=\left( \, \frac{n(n+3)}{2}\, , \, 
1
 \right)\, . $$
 Consequently, this web is AMP with all its ARs logarithmic, except the one associated to the dilogarithmic identity  $({\sf R}_{A_n})$. Finally, 
 ${\boldsymbol{\mathcal Y\hspace{-0.05cm} \mathcal W}}_{\hspace{-0.05cm} A_n}$ is not linearizable hence not algebraizable.
\end{thm}

It is natural to expect that this result generalizes to the other Dynkin types but 
it turns out, a bit surprisingly, that this is not the case: for instance, although  the Dynkin diagrams of type $D$  are simply-laced, we have verified that ${\boldsymbol{\mathcal Y\hspace{-0.05cm} \mathcal W}}_{\hspace{-0.05cm} D_r}$ is not AMP for small values of $r$, a fact that we conjecture to hold true for all $r\geq 4$. However investigating many examples of $\boldsymbol{\mathcal Y}$-cluster webs of bi-Dynkin type leads us to make the following 
 \begin{conjecture}
 \label{Conjecture:YWDD'}
 Let $\Delta,\Delta'$ be two Dynkin diagrams of classical type, of  ranks 
$n$ and $n'$ and of 
 Coxeter numbers $h$ and $h'$ respectively.  
\begin{enumerate}
 \item   The web $\boldsymbol{\mathcal Y\hspace{-0.1cm} \mathcal W}_{\hspace{-0.1cm} \Delta,\Delta'}$ is a  
$d
_{\Delta , \Delta'}$-web in $nn'$ variables with 
$d_{\Delta ,\Delta'}
=nn'\big(h+h'\big)/2$. 
\item   If both $n$ and $n'$ are greater than 1 ({\it i.e.}\,both $\Delta$ and $\Delta'$ are distinct from $A_1$), then 
one has 
$$
\rho^\bullet\Big(\boldsymbol{\mathcal Y\hspace{-0.1cm} \mathcal W}_{\hspace{-0.06cm} \Delta,\Delta'}  \Big)= \Big(\, d_{\Delta , \Delta'}-nn'\, , \,  nn'\, , \, 
1\, \Big)\qquad \mbox{ and } \qquad 
{\rm polrk}^\bullet\Big(\boldsymbol{\mathcal Y\hspace{-0.1cm} \mathcal W}_{\hspace{-0.06cm}  \Delta,\Delta'}
\Big)=   \Big( \, d_{\Delta , \Delta'}\, , \, 1 \, \Big)\, .
$$ 
Consequently,  this web 
is AMP with all its ARs logarithmic, except the one associated to the dilogarithmic identity  $({\sf R}_{\hspace{-0.03cm} \Delta,\Delta'})$. 
 \item Moreover, 
 $\boldsymbol{\mathcal Y\hspace{-0.1cm} \mathcal W}_{\hspace{-0.06cm} \Delta,\Delta'}$ is not linearizable hence not algebraizable.
   \end{enumerate}
\end{conjecture}
Even computing the degree 
 of  $\boldsymbol{\mathcal Y\hspace{-0.1cm} \mathcal W}_{\hspace{-0.1cm} \Delta,\Delta'}$ does not seem to be  straightforward in full generality. We prove that it is equal to 
$nn'(h+h')/2$ in two cases, namely when one of the Dynkin diagram is $A_1$ or when 
both 
 are of type $A$ (thanks to some works by Sherman-Bennett and Volkov respectively). 
\begin{center}
$\star$
\end{center}

Our main goal when starting to study cluster webs was to better understand  their abelian relations and their rank. If the theorem and the conjecture just above show or suggest that many of these webs are interesting regarding the property of being AMP, we are far from having a complete understanding of them. Actually, our web-theoretic investigations led us to ask   basic but seemingly new questions about the algebraic and/or differential-geometric properties satisfied by $\boldsymbol{\mathcal X}$-cluster variables. The consideration of an important number of cases led us to make the following 
\begin{conjecture} Let $x,x'$ be two $\mathcal X$-cluster variables of a cluster algebra, considered as rational functions of some fixed initial cluster variables. 
\begin{itemize}
\item Both $x$ and $x'$ define the same foliation if and only if they coincide (possibly up to inversion).
\item  For any $\lambda\in \mathbf P^1$, if the fiber 
of $x$ over $\lambda$ 
is reducible then 
$\lambda\in \{0 , -1, \infty\}$. Moreover, its irreducible components are cut out by some $F$-polynomials of the considered cluster algebra. 
\item  For any $\lambda,\lambda'\in \mathbf P^1$, if the fibers 
$x^{-1}(\lambda)$ and ${x'}^{-1}(\lambda')$
have an irreducible component $H$ in common, then both $\lambda$ and $\lambda'$ belong to $ \{0 , -1, \infty\}$ and  $H$ is cut out by a $F$-polynomial.
\end{itemize}
(Note that it is necessary to consider the initial cluster variables as particular $F$-polynomials for this statement to have a chance to be satisfied).
\end{conjecture}

These properties that we conjecture for the $\boldsymbol{\mathcal X}$-cluster variables are reminiscent of some properties satisfied by the cross-ratios on the moduli spaces ${\mathcal M}_{0,n+3}$. This is not so surprising for the cluster variables of the finite cluster algebra of type $A$ since it is known in this case that these can be expressed in terms of cross-ratios.  If the above conjecture  is indeed satisfied, it should say roughly that in many ways, $\boldsymbol{\mathcal X}$-clusters variables behave  like cross-ratios  in full generality, which is more striking and, we believe, was not expected.

 \subsubsection{}
To conclude this section, we would like to say a few synthetic words about the main contributions of this text.   What we do in it can be roughly be resumed as follows:  
\begin{itemize}
\item First, regarding webs, we introduce the notion of `{\it virtual rank}' which is new despite its elementary character. It is also easy to compute in practice and the associated notion of '{\it web with AMP rank}' encompasses all the notions of 'maximal rank webs' previously considered in web geometry. It allows to generalize the classical problem of determining/studying webs of maximal rank to more webs than those classically considered in several variables. We also discuss some effective tools to study the abelian relations and the (classical) rank of webs. 
\item In the second stage, we discuss the polylogarithms and the abelian functional equations they satisfy from the point of view of web geometry. We consider many explicit cases that we analyse using the new web-theoretic perspective elaborated in the first part.  Our study shows that a great number of the webs naturally associated to a polylogarithmic AFE carries the maximal possible number of ARs, or more precisely, has `{\it AMP rank}'. 
In our opinion, this sheds an interesting new light on Griffiths' question \eqref{Griffiths-Question}, which probably needs to be rephrased using the new concepts that  have been introduced regarding webs. 
\item In a third and main part, we consider webs which can be obtained from cluster algebras and study their abelian relations, especially those with polylogarithmic components.  We show that some classical polylogarithmic identities are of cluster type. But we also get many new examples of cluster webs with AMP rank  whose ARs all are polylogarithmic.  A big number of explicit examples are studied, as well as several families of cluster webs (or some of the first elements of these families at least). We formulate also several conjectures about cluster variables and cluster webs, coming from natural considerations of web geometry.  We prove some of these conjectures in some cases. 
\item  An interesting feature of this work consists as well  in the important number of questions which are asked  regarding several distinct directions of research.  
\end{itemize}
\mk

What we have obtained in this text shows that although they are distinct 
mathematical objects, 
 webs, polylogarithms and cluster algebras are even more connected than it was thought before.  If we are happy with our results which contribute to make the connexions between these three fields more apparent, it seems to us that 
they do not offer any conceptual explanation of why these three fields are so deeply connected, which is not really satisfying in our opinion. 
We believe that what has been obtained in this text 
appeals further investigations,  and this from different perspectives.
%
%
%




 \subsection{\bf Organization of this text and a few general comments} 
\label{SS:Intro3}

First a word about the general setting: in the whole text, we work in the complex analytic category. Accordingly, except if it is explicitly mentioned, all the objects we will consider  (manifolds, maps, foliations, etc.) will be complex analytic and actually holomorphic, that is without any singularity. 

\subsubsection{}
Here we describe linearly the content of the paper, starting from the next section. \mk 

In {\bf Section 1}, we start by introducing basic material on webs in \S\ref{Par:GermOfWebs}. 
In particular, we define  the `weak general position  (wGP)' assumption   that the leaves of the webs we will deal with in the sequel are assumed to satisfy. Examples of webs are given in \S\ref{SS:Examples-of-Webs} whereas in \S\ref{SS:WebAttributes} we discuss several attributes for webs upon which our subsequent considerations rely. In particular, the new notions of `virtual rank(s)' and `AMP web' are introduced and discussed  in \S\ref{Para:Virtual-Rank} and \S\ref{SS:Webs-AMP-Rank} respectively. In \S\ref{SS:ACM-AMP}, we discuss the relations between the property of being AMP for webs and being ACM for a projective curve. In subsection \S\ref{Par:II-AR}, we review the symbolic algebraic approach for constructing abelian relations with iterated integrals as components for a given web. We use it on some explicit cases and give examples of AMP webs (see \S\ref{Par:Illustrate-Symbolic-Method} and 
\S\ref{SSub:Examples-AMP-Webs}).  In \S\ref{Para:DeterminingARs-and-Rank} we explain that the standard tools of planar web geometry that are Abel's method for determining the ARs and Pantazi-Henaut's criterion for characterizing maximal rank planar webs, actually generalize rather straightforwardly to the multivariable webs only satisfying  (wGP) that we are considering in this text. Finally in \S\ref{SS:Results-Planar-Web-Geometry}, we recall some fundamental results of web geometry that we will use in some proofs later on.
\bk

The {\bf second Section} concerns classical polylogarithms and more specifically the properties of the webs naturally attached to the functional identities these functions satisfy. In  \S\ref{SS:Definition-basic-properties} we review some basic material about polylogarithms and discuss several points such as the different modified versions of polylogarithms previously considered in the literature. Since it will be useful for our purpose (namely to construct polylogarithmic ARs of webs associated to polylogarithmic identities), 
in  \S\ref{SS:Modified-Higher-Polylogarithms}
we work out with some supplementary details, the material of \S\ref{Par:II-AR} in the polylogarithmic case (namely, for iterated integrals on $\mathbf P^1$ ramified at $0,1$ and $\infty$). It is there that we establish Theorem \ref{Thm:PolylogAFE} from which we deduce Corollary \ref{Cor:cor}. 
With the exception of this corollary, essentially all the content of \S\ref{SS:Definition-basic-properties} is well-known to experts.  This material has been incorporated here essentially for the sake of completeness, but also since the holomorphic setting, which is the one relevant for our purpose, is not explicitly considered in the existing literature about functional equations of polylogarithms.  To make this part more appealing, we have also added some material regarding the history of polylogarithms and their functional equations.  Following a common practice in historical papers, this has been incorporated into  the text (in \S\ref{SS:Definition-basic-properties})  by means of several footnotes. \sk 

In 
\S\ref{S:EFA-polylog}, we study the webs associated to many AFEs satisfied by polylogarithms (of weight up to 4, 4 included) that can be found in the classical or more recent literature. Classical dilogarithmic identities are considered in \S\ref{S:EFA-dilog}, 
trilogarithmic ones in \S\ref{SS:TrilogFE} and Kummer's  weight 4 and weight 5 identities in \S\ref{SS:Kummer-K(4)-K(5)}.  Recent (and more complex) polylogarithmic AFEs and the webs associated to them are discussed in \S\ref{SS:WebsRecentAFEs}. What emerges from these subsections is that, even if it is not systematically verified, the web associated to a polylogarithmic AFE is AMP in many cases, moreover quite often with only polylogarithmic ARs ({\it cf.}\,Theorem \ref{T:classical-cluster-webs} above). 
\bk

From {\bf Section \ref{S:Cluster}} to {\bf Section \ref{S:ClusterWebs-other-Dynkin-Type}}, we study webs which can be obtained from cluster algebras. \sk 

 {\bf Section \ref{S:Cluster}} consists in generalities about cluster algebras and the webs which can be obtained from them. After having recalled some basic notions of the theory of cluster algebras in \S\ref{S:Cluster-Algebra}, we introduce in \S\ref{SSub:ClusterWebs} the notion of `cluster web' which is one of the most important of the present memoir.  We also discuss here some basic objects associated to cluster algebras ($F$-polynomial, cluster varieties, etc.) which will prove to be useful for studying cluster webs. 
In the next subsection \S\ref{SS:Y-systems-Periods}, we discuss from the perspective of cluster algebras the webs which can be associated to the so-called $Y$-systems. We then explain that these webs actually are particular cases of a more general class of webs, namely those associated to a `cluster period'. Webs of this type are interesting  in what concerns their ARs hence their rank since each carries a dilogarithmic AR, according to a remarkable result by Nakanishi that we recall in \S\ref{Par:NakanishiDilogarithmicIdentity-Period}.  Finally in \S\ref{SS:Some-conjectures}, we discuss several conjectures about some differential and geometric properties of cluster variables suggested by our web-theoretic investigations about cluster webs. \sk

{\bf Section \ref{S:ClusterWebs-FiniteType}} is devoted to the general properties of cluster webs in finite type. Buildling on some previous works about the cluster variables of finite cluster algebras, we prove some of the conjectures of \S\ref{SS:Some-conjectures} in this case ({\it cf.}\,\S\ref{SS:DifferentialIndependenceClusterVariables}). From this, we deduce   closed formulas for the degrees of the $\boldsymbol{\mathcal X}$- and  $\boldsymbol{\mathcal Y}$-cluster webs associated to any cluster algebra of finite type  in \S\ref{SS:DegreesClusterWebsFiniteType}.  Finally in \S\ref{SS:Y-cluster-Web-non-Lin}, we establish that the $\boldsymbol{\mathcal Y}$-cluster web associated to any Dynkin diagram is  
non-linearizable. \sk
 
 Specific interesting examples of cluster webs are investigated in 
{\bf Section \ref{S:SpecialsCases}}.  First, in \S\ref{SS:cluster-webs:rank2}, we offer a careful study of the three cluster webs associated to the rank 2 cluster algebras, focusing on their abelian relations which are explicitly described. Then, in  \S\ref{SS:ClassicalPolylogarithmicIdentitiesAreOfClusterType}, we establish that several classical polylogarithmic identities actually are of cluster type, namely we prove Theorem \ref{T:classical-cluster-webs} stated above. 
\mk 

 {\bf Section \ref{S:ClusterWebs-type-A}} which follows is devoted to the cluster webs associated to type $A$ Dynkin diagrams.  After considering in \S\ref{SS:X-cluster-web-An} the $\boldsymbol{\mathcal X\hspace{-0.05cm}\mathcal W}_{\hspace{-0.05cm}A_n}$'s (for $n\geq 2$) which all are  AMP (which follows from an immediate reinterpretation of  the nice results of \cite{Pereira}), we focus on the associated $\boldsymbol{\mathcal Y}$-cluster webs in \S\ref{SS:Y-cluster-web-An}.  Our main result there is  Theorem \ref{T:YW-An} which is  a refined version of Theorem \ref{T:YW-An-intro} stated above. Our proof goes as follow: using the natural identification between $\mathcal X_{A_n}$ and $\mathcal M_{0,n+3}$ allows us to work on the latter moduli space and to use the same technics as in \cite{Pereira}. Arguing inductively on $n\geq 2$, one first proves the bound $\rho(\boldsymbol{\mathcal Y\hspace{-0.05cm}\mathcal W}_{\hspace{-0.05cm}A_n})\leq (n+1)(n+2)/2$ in \S\ref{SS:Virtual-rank-YW-An}, before obtaining in \S\ref{SS:log-ARs-YW-An} that the RHS $(n+1)(n+2)/2$ is actually equal to the polylogarithmic rank of the web under scrutiny, thus establishing that $\boldsymbol{\mathcal Y\hspace{-0.05cm}\mathcal W}_{\hspace{-0.05cm}A_n}$ is AMP.  The question of constructing explicit bases of the space of logarithmic ARs of $\boldsymbol{\mathcal Y\hspace{-0.05cm}\mathcal W}_{\hspace{-0.05cm}A_n}$ is discussed in \S\ref{SS:Bases-log-ARs-YW-An}, where two distinct approaches are considered, the first by considering the ARs corresponding to the algebraic identities defining the $Y$-system of type $A_n$ (in \S\ref{SS:Bases-log-ARs-YW-An}), the  second by derivating the dilogarithm identity $({\sf R}_{A_n})$ (see \S\ref{Par:Dilog-Generation-LogAR}). 
 We end the sixth section with subsection \S\ref{SubSec:ZigZag-Map} 
 which is about the birational identification
 $x^{T}: \mathcal M_{0,n+3}\dashrightarrow \boldsymbol{\mathcal X}_{A_n}$  
associated  to the zig-zag triangulation $T$ of the $(n+3)$-gon. We give explicit formulas for $x^T$ and its converse  from which we deduce that this map induces an isomorphism between $\mathcal M_{0,n+3}$ and the complement of the hypersurfaces cut out by the $F$-polynomials in the initial $\boldsymbol{\mathcal X}$-cluster torus ({\it cf.}\,Proposition \ref{P:Xn(Un)=XXn}). \mk 

The  cluster webs associated to cluster algebras of Dynkin type different from  $A$ are studied in {\bf Section \ref{S:ClusterWebs-other-Dynkin-Type}}.  We leave aside the exceptional Dynkin diagrams and mainly consider the cluster webs of type $B$, $C$ and $D$, which are studied in \S\ref{SS:ClusterWebs-type-B}, \S\ref{SS:ClusterWebs-type-C} 
and \S\ref{SS:ClusterWebs-type-D}
respectively. Contrarily to what occurs in type $A$, the cluster webs in the other classical cases do not seem to be AMP. However, using the effective methods to study the ARs and the rank of webs 
discussed in the first section, 
we studied the first webs of each of these three series in order to formulate in each case fairly precise conjectures about their ARs and their virtual and 
genuine ranks. For instance, see the precise statements pages \pageref{Conj:YWBn-page}, \pageref{Conj:YWCn-page}  and  \pageref{Conj:YWDn-page} regarding the $\boldsymbol{\mathcal Y}$-cluster webs of type $B_n$, $C_n$ and $D_n$ respectively.  In \S\ref{SS:YW-Delta-Delta'}, we discuss the 
 case of $\boldsymbol{\mathcal Y}$-cluster webs of bi-Dynkin type, which seem a bit surprisingly more  interesting than in the case of a single Dynkin type not of  type $A$. Indeed, when both Dynkin diagrams $\Delta$ and $\Delta'$ are of classical type and of rank at least 2, the cluster web  
$\boldsymbol{\mathcal Y\hspace{-0.05cm}\mathcal W}_{\hspace{-0.05cm}\Delta,\Delta'}$ seems to be AMP with only logarithmic ARs, plus the one associated to  the dilogarithmic  identity $({\sf R}_{\Delta,\Delta'})$. For a precise (but conjectural) statement, see 
\S\ref{SS:SomeConjecturalStatements} and the following subsections. 
\bk

The last section, {\bf Section \ref{S:Questions-Problems-Perspectives}}, is of a different nature from the previous ones since it offers no new result. 
 On the contrary, 
At the opposite, we discuss here a large number of questions and  problems 
suggested by the material presented in the preceding sections. The structure of this 
speculative section roughly mirrors the general structure of the rest of the text: we first discuss webs (starting from \S\ref{SS:InWebGeometry}), then polylogarithms in \S\ref{SS:AboutPolylog}  where among other things we wonder  about the cluster nature of some polylogarithmic identities (in \S\ref{SS:ClusterNature}) and also briefly discuss by means of interesting examples  the  webs (of codimension at least 2) which can be associated to functional identities satisfied by multivariable polylogarithms.  
Cluster algebras, more precisely cluster variables and cluster webs enter the dance in \S\ref{SS:AboutClusterVW} where basic questions about them are asked. The particularly interesting one of finding polylogarithmic identities of high weight (that is of weight 3 or higher) within the cluster webs is discussed in  
\S\ref{SS: Polylog-AFE-higher-weight-cluster-webs}. 
In many cases, cluster algebras can be quantized and it is natural to wonder about the possible consequences  with regard to  the cluster webs. This is discussed in \S\ref{SS:In-connection-with-quantization} where the consideration 
of some examples of quantum dilogarithmic identities leads us to dream about a possible quantized version of classical web geometry (notions of quantum webs, quantum abelian relations, etc). This subsection is highly speculative but it does suggest several interesting things, in particular regarding identities (with infinitely many terms) satisfied by classical polylogarithms and the infinite webs which may be associated to them. After wondering about a possible modular interpretation in terms of projective geometry for the $\boldsymbol{\mathcal X}$-cluster variety $\mathcal X_\Delta$ for any Dynkin diagram in \S\ref{SS:ClusterDelta-Configurations}, 
we end the memoir by discussing in \S\ref{SS:ClusterWebs-higher-codim} the curvilinear cluster webs associated to any cluster algebra of finite type.  
The consideration of the two webs $\boldsymbol{\mathcal X\hspace{-0.05cm} \mathcal W}_\Delta^{(1)}$ when $\Delta$ is $A_3$ or $B_3$ shows that 
some of the new basic notions about 1-codimensional webs introduced in this text (virtual ranks, property of being AMP) have natural analogues for curvilinear webs, are relevant and certainly deserve further study.
{\textcolor{gray}{
\subsubsection{Relations with other works (to be written).}
In this subsection (that we intend to develop in a subsequent version), we plan to discuss briefly, or at least to mention,  some recent works in which some interesting links between polylogarithms, cluster algebras, and other objects (mainly coming from mathematical physics) have appeared. }}

{\textcolor{gray}{
Here is a list (in chronological order)  of papers/preprints that we think   interesting   to mention considering the  approach taken in this memoir:
\begin{itemize}
\item  Several papers of mathematical physics about scattering amplitudes where  polylogarithms and cluster algebras come into the picture: the papers 
\cite{GGSVV}, \cite{GSVV} cited in the bibliography, but also the preprint J.\,Golden \& al,  \href{https://doi.org/10.1088/1751-8113/47/47/474005}{\it Cluster Polylogarithms for Scattering Ampli-} \href{https://doi.org/10.1088/1751-8113/47/47/474005}{\it tudes}, arXiv:1401.6446;
\sk 
\item A.B.\,Goncharov, D.\,Rudenko, \href{https://arxiv.org/abs/1803.08585}{\it Motivic correlators, cluster varieties and Zagier's conjecture}  
\href{https://arxiv.org/abs/1803.08585}{\it on zeta(F,4)}, arXiv:1803.08585;
\item  C.K.\,Zickert, \href{https://arxiv.org/abs/1902.03971}{\it Holomorphic polylogarithms and Bloch complexes}, arXiv:1902.03971;
\item R.\,de Jeu, \href{https://arxiv.org/abs/2007.11014}{\it Describing all multivariable functional equations of dilogarithms},  
arXiv:2007.11014;
\item  S.\,Charlton, H.\,Gangl, D.\,Radchenko.
\href{https://arxiv.org/abs/2012.09840}{\it Functional equations of polygonal type for multiple}
\href{https://arxiv.org/abs/2012.09840}{\it polylogarithms in weights 5, 6 and 7}, 
 arXiv:2012.09840;
\item S.\,Charlton, C.\,Duhr, H.\,Gangl. 
\href{https://arxiv.org/abs/2104.04344}{\it Clean single-valued polylogarithms},   arXiv:2104.04344.
\end{itemize}
}}

\subsection*{\bf Acknowledgements}
The author is indebted to \href{http://irma.math.unistra.fr/~chapoton/}{F.\,Chapoton} for sharing with him his Maple routines
to work with cluster algebras a couple of years ago.  
These routines were a great tool to explore the world of cluster algebras from the point of view of web geometry.
 As for polylogarithms and their functional equations,  the author has benefited 
from several very informative (e-mail) exchanges with \href{https://www.maths.dur.ac.uk/~dma0hg/}{H.\,Gangl}, circa the same period. We are very grateful to him for answering many of our questions and more generally for 
 sharing with us his deep knowledge on these matters.
\sk 

We are thankful to  \href{https://webusers.imj-prg.fr/~sophie.morier-genoud/}{S.\,Morier-Genoud} and \href{http://ovsienko.perso.math.cnrs.fr}{V.\,Ovsienko} for the early interest they both have shown in our work and for their invitations to give talks about it at Reims seminar  and Paris algebra seminar.  
When starting this project, we were almost entirely ignorant about cluster algebras. 
We are very grateful to \href{https://webusers.imj-prg.fr/~bernhard.keller/}{B.\,Keller} for numerous exchanges about this subject but also 
about many other topics related to it, such as the theory of $Y$-systems, the notion of scattering diagram, or the so-called quantum dilogarithm identities.
We thank  \href{https://math.berkeley.edu/~msb/}{M.\,Sherman-Bennett} as well, for having kindly answered to some of our questions about the cluster variables of finite type cluster algebras.\sk

I was very lucky to have the opportunity to talk with \href{https://smf.emath.fr/actualites-smf/deces-de-laurent-gruson}{L.\,Gruson} in the last few years, 
in particular about the notion of ACM projective curves that he was quite familiar with, and which plays an important role regarding  
the general problem addressed in this text.  More generally, I am very grateful to him for having always been available to enlighten me on questions of algebraic geometry. His mathematical knowledge was equalled only by his delicacy and kindness, I will miss all this.
\mk

Finally and as so often, we thank Brubru for her patient proofreading and the numerous English corrections she has  suggested. 

\newpage
\section{On webs and their geometry} 
\label{S:Webs}
We introduce below some basic notions of web geometry. Roughly, a web on a manifold $M$ is the geometric configuration formed locally by (the leaves of) a finite number of foliations on $M$. 
Actually, we will be only interested in a particularly simple kind of webs (namely the ones defined by rational functions on $\mathbf C^n$)  
hence we will essentially introduce the main notions of web geometry    
in this particularly simple case. 
We let the reader  get precise notions and definitions in full generality. 
\smallskip

For more material, details and perspectives, we refer to  \cite{Coloquio,PirioBook} and the references therein. 

 \subsection{\bf First definitions}
    In this text, we will only consider the 1-codimensional case hence all the foliations we will work with  will be holomorphic foliations of codimension 1, possibly with singularities.  \smallskip

\subsubsection{Germs of webs.} 
 \label{Par:GermOfWebs}
 A (germ of) 
 $\boldsymbol{d}${\bf-web}  at the origin of $ \mathbf C^n$
 is 
a $d$-tuple $\boldsymbol{\mathcal W}=(\mathcal F_i)_{i=1}^d$ of regular foliations by hypersurfaces  
on $(\mathbf C^n,0)$
which moreover  satisfy a general position assumption. 
\sk

\vspace{-0.3cm}
 There are several possibilities for such an hypothesis and  a convenient way for us to formulate them is by considering some first integrals $u_1,\ldots,u_d \in \mathcal O_{(\mathbf C^n,0)}$ for the foliations  which compose $\boldsymbol{\mathcal W}$ : one has $\mathcal F_i=\mathcal F(u_i)$ for $i=1,\ldots,d$. 
 The most classical general position assumption made in web geometry ({\bf strong general position} assumption, denoted by {(sGP)} in what follows) 
  requires that given $n_d=\max(d,n)$ foliations among the $\mathcal F_i$'s, these intersect as transversally as possible. Analytically, this means that  for every subset $\{i_1,\ldots,i_{n_d}\} \subset \{1,\ldots,d\}$ of cardinality $n_d$, the wedge product $du_{i_1}\wedge \cdots \wedge du_{i_{n_d}}$ does not vanish at the origin. \smallskip 
 
   The previous general position assumption is quite strong hence to deal with the webs we are interested in here, it is necessary to relax it and to allow webs satisfying the following `{\it weak general position assumption}', denoted by (wGP): 
$$
\begin{tabular}{lcr}
 \begin{tabular}{c}
${}^{}$\hspace{-0.6cm} {\bf (wGP) }
\end{tabular}
 &   &  \hspace{-0.6cm} \begin{tabular}{l}
 {\it  two foliations  of $\boldsymbol{\mathcal W}$ intersect transversally, that is in codimension 2, as soon  as} \\  {\it   they are distinct;  in terms of the  local first integrals $u_1,\ldots,u_d$, this means that} \\
  {\it for any $i,j=1,\ldots,d$ such that $i<j$, one has $du_i\wedge du_j\neq 0$ at the origin.}
\end{tabular}
\end{tabular}
   $$

This relaxed notion of general position already appears in some  papers on webs, such as \cite{Pereira}.  Note that it does not properly coincide with the notion bearing the same name in \cite{CavalierLehmann} where the authors make the (natural) supplementary assumption that  $n$ foliations among the ones of the considered web  are  in  strong general position.
 In what follows,  if the contrary is not explicitly mentioned, we will also assume that this latter hypothesis holds true. \mk 
 
The number of pairwise transverse foliations composing a web is called the {\bf degree} of this web.

\subsubsection{Webs.}  
A global $d$-web $\boldsymbol{\mathcal W}$ on a manifold $M$ is the data of germs of $d$-webs ${\boldsymbol{\mathcal W}}_m$  for any point $m$ of $M$ or for any $m$ in a dense open subset of $M$. \sk 
 
For instance, let $u_1,\ldots,u_d$ be $d$ non-constant meromorphic functions defined on a domain $U\subset \mathbf C^n$. Then for any $i$, the level subsets $u_i=\lambda$ organize themselves in a (possibly singular) foliation  on $U$, denoted by $\mathcal F_{u_i}$. Assuming that  $du_i\wedge du_j\neq 0$ at the generic point of $U$ as soon as $i\neq j$, one gets that the $\mathcal F_{u_i}$'s  satisfy (wGP) hence the 
  $d$-tuple of foliations $(\mathcal F_{u_1},\ldots,\mathcal F_{u_d})$ is a 
  web on $U$, which will be denoted by $\boldsymbol{\mathcal W}(u_1,\ldots,u_d)$.  Note that it can happen that some of the $\mathcal F_{u_i}$'s admit singularities on $U$ and/or that the general position assumption (wGP) is not satisfied everywhere on $U$ hence  there is a slight blurring about what precisely is the definition domain of the web  $\boldsymbol{\mathcal W}(u_1,\ldots,u_d)$.  However, 
we will consider properties of webs  which are generically satisfied (or not)  in the whole memoir, thus  this lack of precision regarding  the definition domains of the webs we will deal with will not cause any genuine problem hence we will not elaborate further on that in the sequel.
\mk 

%

A natural choice for  first integrals $u_i$  of a given web is to take some which are {\bf primitive}, that is such that for each $i$, the level-subset $u_i=\lambda$ is  connected  for any generic $\lambda\in {\rm Im}(u_i)$. \mk

In this text, we will mainly consider webs $\boldsymbol{\mathcal W}(u_1,\ldots,u_d)$ defined by rational first integrals $u_i\in \mathbf C(x_1,\ldots,x_n)$.  In this context, the fact that $u_i$ be a primitive first integral 
for  the foliation $\mathcal F_{u_i}$ it defines, 
 coincides with the fact that $u_i$ is {\bf noncomposite} as a rational function\footnote{This follows from a classical  theorem due to Bertini, see Theorem 37 in \cite{Schinzel} or Theorem 3.4.6 page 108 in \cite{Jouanolou}.}, {\it i.e.}\,there does not exist a pair of rational functions $(v_i,r)$ with 
$v_i  \in \mathbf C(x_1,\ldots,x_n)$ and $r\in \mathbf C(u)$ of  degree at least 2,  such that $u_i=r\circ v_i$.
By an obvious recurrence on the degree of $u_i$, it comes that  an algebraically integrable foliation such as $\mathcal F_{u_i}$ always admits a noncomposite (equivalently, a primitive) rational first integral, which moreover  is unique up to post-composition by a M\"{obius} transformation.\mk 

 In practice, when dealing with a web $\boldsymbol{\mathcal W}(u_1,\ldots,u_d)$ defined by rational first integrals $u_i$, it is 
 more convenient to assume that they are noncomposite.
\subsubsection{Equivalence of webs and web geometry.} 
\label{Parag:EquivalenceOfWebs}
Two germs of  $d$-webs $\boldsymbol{\mathcal W}$ and $\boldsymbol{\mathcal W}'$  on $(\mathbf C^n,0)$ are {\bf equivalent} if there exists 
$\varphi\in {\rm Aut}(\mathbf C^n,0)$ such that 
$$\boldsymbol{\mathcal W}'=\varphi^*(\boldsymbol{\mathcal W})\, ,$$
 possibly up to re-indexing the foliations.
 Analytically, this translates as follows once one introduces some first integrals $u_i$ and $u_i'$ for $\boldsymbol{\mathcal W}$ and $\boldsymbol{\mathcal W}'$ respectively: these webs are equivalent if and only if there exist $\varphi$ as above and $\sigma\in \mathfrak S_d$ as well as invertible  holomorphic functions $h_1,\ldots,h_d$ such that $u_{\sigma(i)}'=h_i\circ ( u_i\circ \varphi)$ for $i=1,\ldots,d$. \medskip 
 
Two global webs $\boldsymbol{\mathcal W}$ and $\boldsymbol{\mathcal W}'$ are said to be {\bf equivalent} if there are two points $x$ and $x'$ 
on their respective regular sets,  such that the two associated 
 regular germs ${\boldsymbol{\mathcal W}}_{\hspace{-0.05cm}x}$ and $ {\boldsymbol{\mathcal W}'}_{\hspace{-0.17cm}x'}$ are equivalent in the previous sense. 
 This is a notion of a local nature, but for global webs, it tends to become global, in the sense that the germ of biholomorphism $\varphi$ realising the local equivalence extends globally, possibly with ramification and/or as a multivalued map,  (at least) in any case we could consider until now.\mk

\begin{figure}
    \centering
    \begin{minipage}{0.45\textwidth}
        \centering
        \includegraphics[width=1\textwidth]{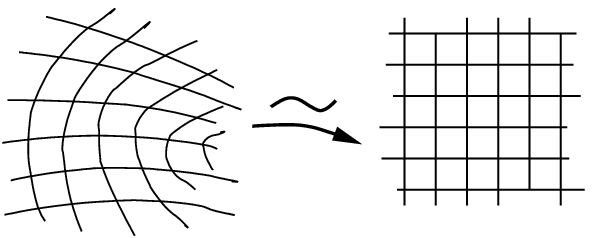} 
    \end{minipage}\hfill
    \begin{minipage}{0.45\textwidth}
        \centering
        \includegraphics[angle=0,origin=c,width=1\textwidth]{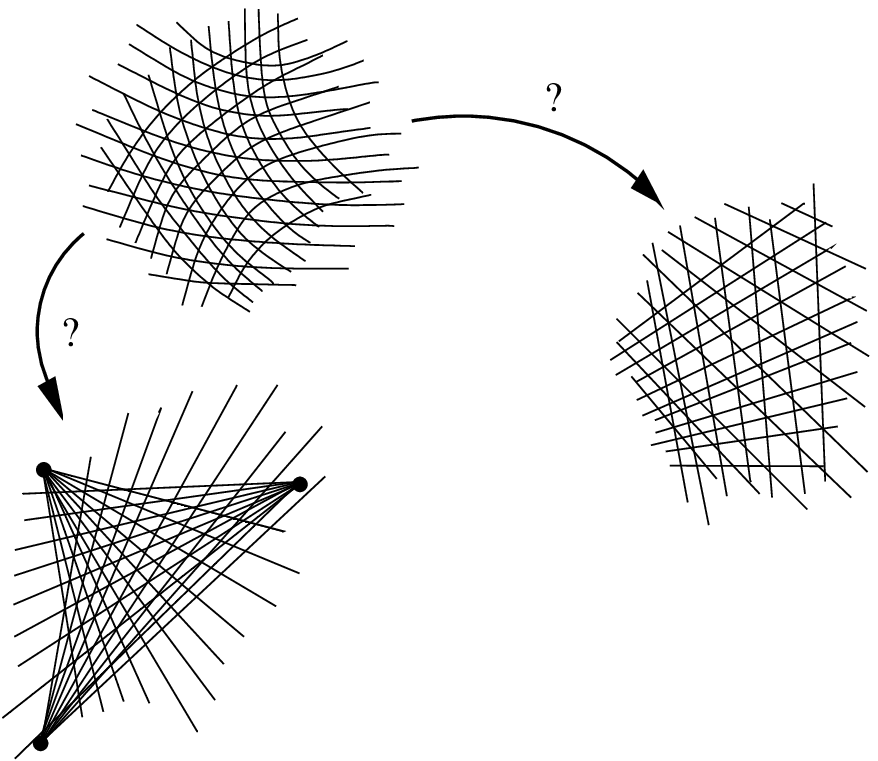} 
    \end{minipage}
        \caption{Two classical problems of web-geometry: if a planar 2-web is always locally trivial (left), even the local geometry of a planar 3-web is quite richer. For instance,  the questions of knowing whether such a web is  parallelizable or linearizable  (right) admit non-trivial answers.}
        \label{Fig:WebGeometry}
\end{figure}

 The {\bf geometry of webs} consists in the study of webs up to this notion of equivalence.  See Figure \ref{Fig:WebGeometry} for some very basic questions of planar web geometry.   The standard way to study webs from this point of view is to attach some invariants to them.   In the next paragraph, we describe a class of such invariants which is very specific to the study of webs.

\newpage

\subsection{\bf Examples of webs}  
\label{SS:Examples-of-Webs}
We review below some classical examples of webs. 

\subsubsection{\hspace{-0.2cm} Pencils of linear subspaces.}${}^{}$ \hspace{-0.4cm}
Given $d\geq 1$ pairwise distinct points  in the projective plane, the  associated pencils of lines form the foliations of a global (but singular) $d$-web on $\mathbf P^2$. 
\begin{figure}[h]
    \centering
    \begin{minipage}{0.45\textwidth}
        \centering
        \includegraphics[width=0.9\textwidth]{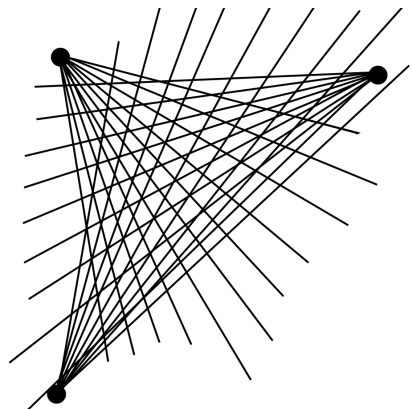}  
           \end{minipage}
        \caption{A planar 3-web formed by three pencils of lines.}
        \label{Fig:3-web-pencils}
\end{figure}

The  generalization in arbitrary dimension $n\geq 2$ is straightforward: $d$ pairwise distinct pencils of hyperplans in $\mathbf P^n$ define a $d$-web on the complement of the union of hyperplanes which are common to two  pencils among the $d$.



\subsubsection{Coordinate  and parallelizable webs.}
\label{SS:Coordinate-parallelizable-webs}
Let $n$ be an integer bigger than or equal to $2$. 
 By definition, a {\bf coordinate web} on $\mathbf C^n$ is  a  web which is `equivalent' ({\it cf.}\,\S\ref{Parag:EquivalenceOfWebs} below where this notion is properly defined) to the {\bf standard  coordinate web} $\boldsymbol{\mathcal W}(x_1,\ldots,x_n)$ where the $x_i$'s are the standard coordinates. 
 For instance, a planar 3-web formed by three pencils of lines is a coordinate web. 
\sk 

A {\bf parallelizable $d$-web} web on $\mathbf C^n$ (or on a domain in it) is a linear web formed by $d$-families of parallel hyperplanes.  Such a web can also be described as a web formed by $d$-pencils of hyperplanes, the vertices of which are all included in the hyperplane at infinity $\mathbf P^{n-1}_\infty$.  When $d=n$ and under the assumption that the vertices of the considered $n$ pencils are hyperplanes in general position in $\mathbf P^{n-1}_\infty$, one recovers the notion of coordinate web.

\subsubsection{\hspace{-0.2cm} Web associated to an AFE.}
   To a functional equation of the form $\sum_{i=1}^d F_i(u_i)=0$ where the $u_i$'s stand for  nontrivial functions of $n\geq 2$ variables, one can associate the web formed by the foliations $\mathscr F_{u_i}$ (for $i=1,\ldots,d$) on the intersection of the definition domains of the $u_i$'s.  
   
   The simplest concrete example of such a construction is certainly the one obtained by
   the functional equation of the logarithm: ${\rm Log}(x)-{\rm Log}(y)-{\rm Log}(x/y)=0$. 
  The associated web is $\boldsymbol{\mathcal W}(x,y,x/y)$ and is formed by 3 pencils of lines     whose vertices are not aligned ({\it cf.}\,Figure \ref{Fig:3-web-pencils}). 
\subsubsection{\hspace{-0.2cm} Bol's web.}
   To a functional equation of the form $\sum_{i=1}^d F_i(u_i)=0$ where the $u_i$'s stand for  nontrivial functions of $n\geq 2$ variables, one can associate the web formed by the foliations $\mathscr F_{u_i}$ (for $i=1,\ldots,d$) on the common intersection of the definition domains of the $u_i$'s. \sk

For instance, the web associated in this way to Abel's functional equation $(\boldsymbol{\mathcal A b})$ is the 5-web 
$$
\boldsymbol{\mathcal B}=\boldsymbol{\mathcal W}\bigg(x,y,\, \frac{x}{y}\, , \, \frac{1-y}{1-x}\, , \,\frac{x(1-y)}{y(1-x)}\, 
\bigg)
$$
which is known as {\bf Bol's web} in web geometry (since it  was first recognized by him as an  example of an `exceptional web' {\rm \cite{Bol}}; this explains the notation $\boldsymbol{\mathcal B}$ for this web). \sk

Bol's web can also be described geometrically given four points in general position in $\mathbf P^2$: it is the global 5-web formed by the four associated pencils of lines and the pencil formed by the conics passing through all the four points, see the picture juste below. 

\begin{figure}[h]
    \centering
    \begin{minipage}{0.45\textwidth}
        \centering
        \includegraphics[width=0.8\textwidth]{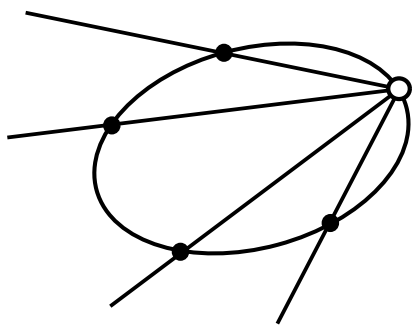}  
           \end{minipage}
        \caption{The leaves of Bol's web  passing through a generic point (white dot) 
               of $\mathbf P^2$.}
        \label{Fig:Bol}
\end{figure}
\subsubsection{Algebraic webs.}
\begin{figure}[h]
    \centering
    \begin{minipage}{0.45\textwidth}
        \centering
        \includegraphics[width=1\textwidth]{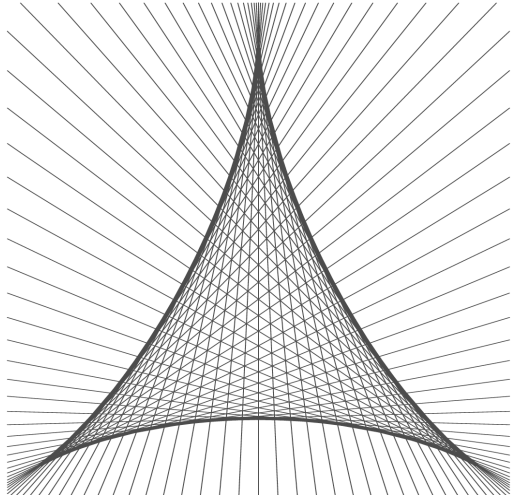}  
           \end{minipage}
        \caption{a planar algebraic 3-web.}
        \label{Fig:HHYYPP} 
\end{figure}

  To an irreducible and reduced algebraic curve $C\subset \mathbf P^n$ of degree $d$ is classically associated a dual web ${\boldsymbol{\mathcal W}}_C$ on $\check{\mathbf P}^n$. This web can be defined as follows: dually, the assumption that $C$ has degree $d$ translates into the fact that the dual hypersurface $C^*\subset \check{\mathbf P}^n$  of the considered curve\footnote{By definition, the {\bf dual hypersurface} $\boldsymbol{C^*}$  of the curve $C\subset \mathbf P^n$ is the closure  of the set of hyperplanes of $\mathbf P^n$ which are tangent to $C$ at one of its regular points.} is of `{\it class}' $d$: through a general point $H$ of $ \check{\mathbf P}$ pass $d$ hyperplanes tangent to $C^*$. These are the leaves of ${\boldsymbol{\mathcal W}}_C$  through $H$.  Note that when $C$ is degenerate in $\mathbf P^n$ (that is, when its span  $\langle C\rangle$ is a proper linear subspace of the ambiant projective space), the web ${\boldsymbol{\mathcal W}}_C$ as defined above is the pull-back, under a  linear projection $\pi: \mathbf P^n\setminus P\rightarrow \langle C\rangle$ from a subspace $P$ supplementary to $\langle C\rangle$, of the web associated to $C$ but now considered as a (non-degenerate) curve in $\langle C\rangle$.  From this, it follows that if $C$ is a line in $\mathbf P^n$, then ${\boldsymbol{\mathcal W}}_C$ is nothing else but the 1-web (that is, the foliation) formed by the pencil of hyperplanes containing $C^*$ which is a linear subspace of codimension 2 in $\mathbf P^n$. 
  Finally,  one generalizes 
straightforwardly  this construction to the case when $C$ is reduced but admits several irreducible components: if $C_1,\ldots,C_m$ denote them, then 
  $${\boldsymbol{\mathcal W}}_C={\boldsymbol{\mathcal W}}_{C_1}\boxtimes \cdots \boxtimes 
  {\boldsymbol{\mathcal W}}_{C_m}\, .$$ 
  The notation $\boxtimes$ means that we consider the web obtained by considering the union of the foliations of the webs involved.     Figure \ref{Fig:HHYYPP}  below is a picture of a real algebraic 3 web formed by the tangent lines to a hypocycloid curve with 3 cusps (which is the dual of a singular cubic in $\mathbf P^2$). 

\begin{figure}[h]
    \centering
    \begin{minipage}{0.45\textwidth}
        \centering
        \includegraphics[width=1.1\textwidth]{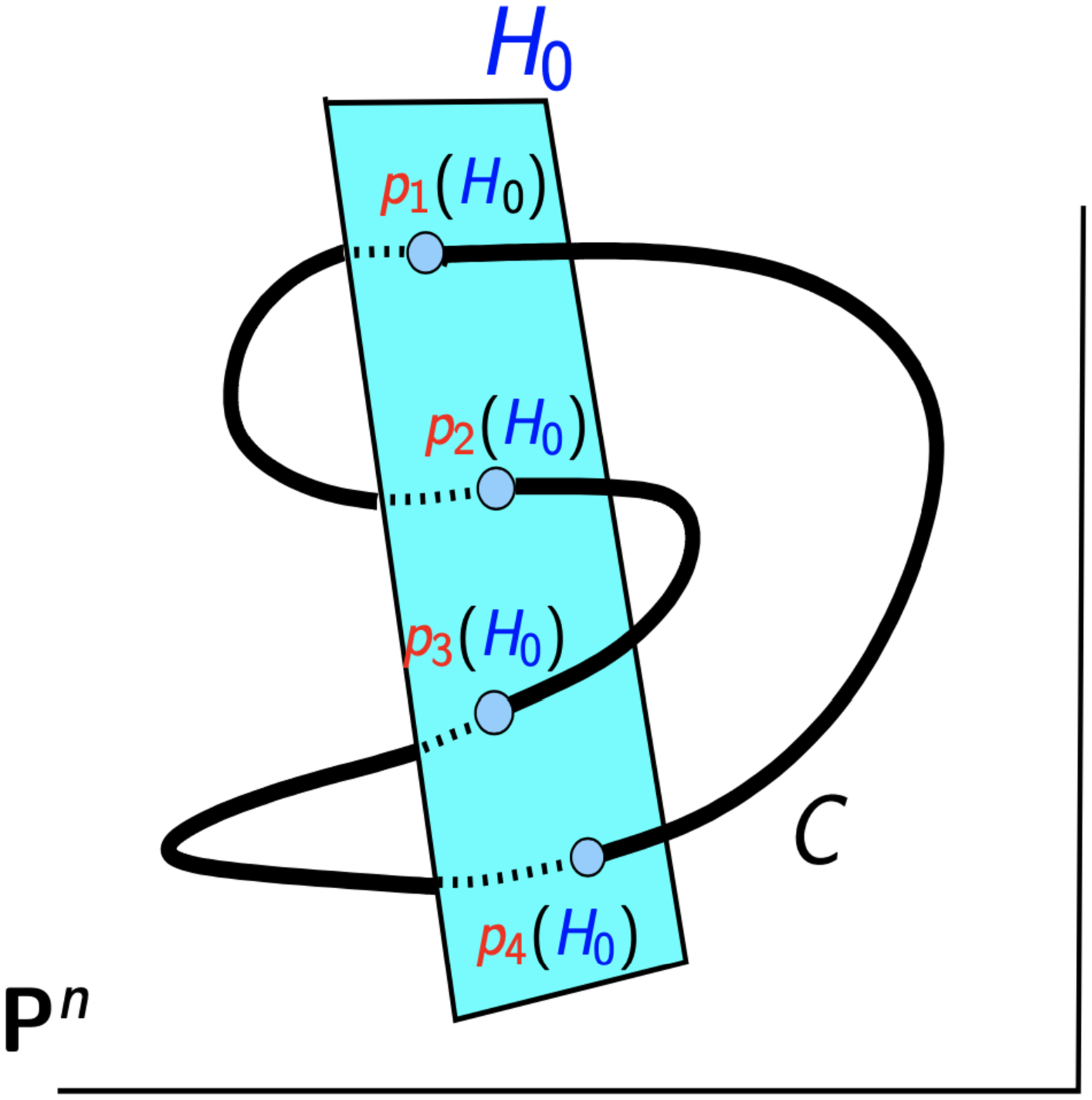}  
           \end{minipage}
       \caption{A generic hyperplane $H_0$ intersects $C$ in $d=\deg C$ distinct points.} 
        \label{Fig:WCWC} 
\end{figure}

  The web 
  ${\boldsymbol{\mathcal W}}_C$ can also be described locally by means of algebraic first integrals: let $H_0\in \check{\mathbf P}^n$ be a hyperplane which intersects $C$ in $d$ pairwise distinct points.  According to Bezout's Theorem, these points are smooth points $p_1(H_0),\ldots,p_d(H_0)$ of $C$ which intersects $H_0$ transversally 
   (see Figure \ref{Fig:WCWC}).
From this, it follows that there exist  $d$ germs of holomorphic 
  maps $p_i: (\,\check{\mathbf P}^n,H_0\,)\rightarrow C$ such that $H\cdot C=\sum_{i=1}^d p_i(H)$ as a 0-cycle on $C$, for any hyperplane $H$ sufficiently close to $H_0$.
 It is easily seen that the $p_i$'s are holomorphic submersions (this comes from the transversality assumption) and define a web which coincides with the germ of ${\boldsymbol{\mathcal W}}_C$ at $H_0$: one has 
  $$
  \big({\boldsymbol{\mathcal W}}_C\big)\lvert_{ (\, \check{\mathbf P}^n,H_0\,)}=
   {\boldsymbol{\mathcal W}}\big(\, p_1,\ldots,p_d 
    \, \big)\,.
  $$

By definition, a web is said to be {\bf algebraic} if it is of the form  ${\boldsymbol{\mathcal W}}_C$  for an algebraic curve $C$ as above.   Notice that such webs are global singular linear webs on $ \check{\mathbf P}$. Except when $C$ has several  irreducible components, ${\boldsymbol{\mathcal W}}_C$ is irreducible as a global  web on the dual projective space $\check{\mathbf P}^n$.\mk 

A web is said to be  {\bf algebraizable} if it is equivalent to an algebraic web  ${\boldsymbol{\mathcal W}}_C$ associated to a projective curve $C$. Since algebraic webs are  {\bf  linearizable} (that is equivalent to a {\bf linear web}, that is a web on a domain in an affine complex space whose leaves are pieces of hyperplanes), linearizability is a necessary condition for the algebraizability of a given web. 
%

\subsubsection{Webs induced by restriction along a subvariety.}
\label{SS:WebsRestriction}
One can obtain a new web from a given one by restricting the latter along a subvariety of the ambiant space. This is a very elementary and basic construction which, however, appears to be relevant in order to construct interesting webs (as Theorem \ref{T:classical-cluster-webs} will show for instance).
\sk

More precisely,  let  $\boldsymbol{\mathcal W}$ be a web  defined on a variety $M$ and let $Z\subset M$ be an irreducible and reduced subvariety (which may or may not have singularities).  The foliations of $\boldsymbol{\mathcal W}$ can be of two types regarding $Z$, depending on whether they
 leave it  invariant or not.  Considering only those 
 inducing a foliation of codimension 1 on $Z$ by restriction, one obtains a web on $Z$ that we will denote (again a bit abusively) $\boldsymbol{\mathcal W}\lvert_Z$.  Note that it may happen that two distinct foliations of $\boldsymbol{\mathcal W}$ generically transverse to $Z$ induce the same foliation when restricted along $Z$. So in full generality, 
it is not clear what might be the characteristics of $\boldsymbol{\mathcal W}\lvert_Z$, even some very basic ones such as its degree, even if those  of $\boldsymbol{\mathcal W}$ are known and well understood. 
  
%
%

\subsubsection{\hspace{-0.2cm} Webs on moduli spaces of projective configurations.}
\label{SubPar:Webs-spaces-configurations}
Let $n$ and $N$ be two positive integers. One defines the {\bf moduli space   of  projective configurations of $N$ points in $\mathbf P^n$}, denoted by $ \boldsymbol{{\bf Conf}_N(\mathbf P^n)}$,   as the quotient of the space of $N$-tuples of pairwise distinct points on $\mathbf P^n$ by the diagonal action of the group of projective automorphisms, {\it i.e.}  
$$
{\rm Conf}_N\big(\mathbf P^n\big)= \left(  \big(\mathbf P^n\big)^N\setminus  \Delta^{n,N}
\right)_{\big/ \, {\rm PGL}_{n+1}(\mathbf C)}
$$
where $\Delta^{n,N}$ stands for the set  of $N$-tuples $(x_1,\ldots,x_N)$ of elements  of the $n$-dimensional projective space $\mathbf P^n$ such that there exist two distinct indices $i,j\in \{1,\ldots,N\}$ such that $x_i=x_j$.  The moduli space ${\rm Conf}_N\big(\mathbf P^n\big)$ naturally carries the structure of a quasi-projective variety. It is finite for  $N\leq n+2$ and irreducible and of dimension $n(N-n-2)$ otherwise.

\paragraph{Webs on moduli spaces $\boldsymbol{\mathcal M_{0,n+3}}$.}
 For any $n\geq 2$, the moduli space $\mathcal M_{0,n+3}={\rm Conf}_{n+3}(\mathbf P^1)$ of projective configurations of $n+3$ points on the projective line carries a natural ${ n+3 \choose 4}$-web, denoted by ${\boldsymbol{\mathcal W}}_{\hspace{-0.05cm} \mathcal M_{0,n+3}}$ which is the one 
whose first integrals  are all the maps $\mathcal M_{0,n+3}\rightarrow \mathcal M_{0,4}\simeq \mathbf P^1\setminus \{0,1,\infty  \}$ obtained by forgetting $n-1$ points among the $n+3$ ones.  \sk 

One easily verifies  that ${\boldsymbol{\mathcal W}}_{\hspace{-0.05cm} \mathcal M_{0,5}}$ and Bol's web $\boldsymbol{\mathcal B}$ are equivalent planar 5-webs. It can be proved that these webs are not linearizable {\rm (}{\it cf.}\,{\rm \cite[\hspace{-0.08cm}\S6.1.4]{Coloquio}}{\rm )} hence they are not algebraizable, {\it i.e.}\,they are not equivalent to any  algebraic web associated to a plane quintic curve.\mk 

The webs ${\boldsymbol{\mathcal W}}_{\hspace{-0.05cm} \mathcal M_{0,n+3}}$ for any $n\geq 2$,   have been studied  in {\rm \cite{Pereira}}.

\paragraph{} ${}^{}$ \hspace{-0.4cm}
\label{Par:WebsOnConfugurationSpaces}
  The previous series of webs generalizes in a straightforward way as follows.  \sk 
%

For $m,n\geq 1$, let ${\rm Conf}_{m+n+2}(\mathbf P^n)$ be the space  of  projective  configurations of $m+n+2$ pairwise distinct points on $\mathbf P^n$: it is a rational variety of dimension $mn$.  For any subset $I\subset \{1,\ldots,m+n+2\}$ of cardinality $n-1$, let $\pi_I: {\rm Conf}_{m+n+2}(\mathbf P^n)\dashrightarrow {\rm Conf}_{m+3}(\mathbf P^1)=\mathcal M_{0,m+3}$  be the rational map 
obtained by projecting the points $x_s$'s for  $s\in \{1,\ldots,m+n+2\}\setminus I$   from the subspace $\langle x_i\, \lvert \, i\in I\, \rangle\subset \mathbf P^n$. Then for any $J\subset \{1,\ldots,m+n+2\}\setminus I$ of cardinality 4, one denotes by  $p_J: \mathcal M_{0,m+3}\rightarrow \mathcal M_{0,4}\simeq \mathbf P^1\setminus \{0,1,\infty \}$ the  forgetful map consisting in forgetting all the points except the $x_j$'s for $j\in J$. By composing the projections $\pi_I$ with the forgetful map $p_J$, one constructs non-constant rational maps 
\begin{equation}
\label{Eq:pi-I-J}
\xymatrix{
\pi_{I,J}=p_J\circ \pi_I\, :\,  {\rm Conf}_{m+n+2}(\mathbf P^n)\,  \ar@{-->}[r]& 
\, \mathcal M_{0,4}\simeq 
\mathbf P^1\setminus \{0,1,\infty \}
}
\end{equation}
for any pair $(I,J)$ of  disjoint subsets $I,J\subset  \{1,\ldots,m+n+2\}$ 
with $\lvert I\lvert=n-1$ and $\lvert J\lvert=4$.\footnote{Note that if $I$ is just a set,  $J$ has to be ordered in order to define 
$\pi_{I,J}$ as a rational function. But if one is only interested in the associated foliation, the ordering of the elements of $J$ does not matter hence can be forgotten.}   We denote the set of such pairs by $ \mathfrak P(m,n)$.
 Its cardinality is  $d_{m,n}= 
\lvert 
\mathfrak P(m,n) \lvert
= \scalebox{0.9}{${ m+n+2  \choose n-1} { n+3  \choose 4}$}\, . 
$ 

The proof of the following lemma is easy and left to the reader:
\begin{lem}
\label{Lem:pi-I-J--pi-I'-J'}
Given two pairs $(I,J)$ and $(I',J')$ elements of  $\mathfrak P(m,n)$, the 
  maps $\pi_{I,J}$ and $\pi_{I',J'}$ define the same foliation on 
${\rm Conf}_{m+n+2}(\mathbf P^n)$ if and only if  $I=I'$ and $J=J'$ (as non-ordered sets). 
\end{lem}
 By considering the  $\pi_{I,J}$'s for all the $(I,J)$'s in this set, one obtains a $d_{m,n}$-web in $mn$ variables: 
\begin{equation}
\label{Eq:W-Conf-n-m+2-P}
{\boldsymbol{\mathcal W}}_{\hspace{-0.05cm}  {\rm Conf}_{m+n+2}(\mathbf P^n) }= {\boldsymbol{\mathcal W}}\Big( \, \pi_{I,J}
\hspace{0.15cm } \lvert   \hspace{0.15cm }(I,J)\in \mathfrak P(m,n)\,  \Big)\, .
\end{equation}
These webs seem very interesting, especially  with regard to their abelian relations, but   have not been really studied until now, except in the (already interesting) case  $n=1$. 

\paragraph{Webs on strata of projective configurations.}
\label{Subpar:Strata}
The construction mentioned in \S\ref{SS:WebsRestriction} allows to construct many webs from those in \eqref{Eq:W-Conf-n-m+2-P} by considering the restrictions of the latter along strata of degenerate configurations.\sk 

 Recall that a configuration $(x_i)_{i=1}^M$ of $M>n+2$ points in $\mathbf P^n$ is {\bf degenerate} if there exists a non-empty subset $K \subset \{1,\ldots,M\}$ such that the dimension of the (projective) subspace spanned by $\{ x_k\, \lvert \, k\in K\}$  is strictly less than $\lvert K \lvert-1$.  The whole set of degeneracies between the elements of a degenerate configuration can be encoded by means of a combinatorial object associated to it, called a \href{https://en.wikipedia.org/wiki/Matroid}{\bf matroid} (on $M$ elements and of rank $n+1$ if the considered configuration belongs to ${\rm Conf}_{M}(\mathbf P^n) $). 
 Given such a matroid $\boldsymbol{\mathfrak M}$, 
the set  
 of all configurations in ${\rm Conf}_{M}(\mathbf P^n) $ whose associated matroid coincides with $\boldsymbol{\mathfrak M}$ form an algebraic subvariety of 
$   {\rm Conf}_{M}(\mathbf P^n) $ denoted by ${\rm Conf}_{\boldsymbol{\mathfrak M}}(\mathbf P^n)$, or sometimes just by $\boldsymbol{\mathfrak M}$ (abusively).\footnote{It should be noted that 
$ {\rm Conf}_{\boldsymbol{\mathfrak M}}(\mathbf P^n)$ 
may very well be empty, depending on whether  the matroid ${\mathfrak M}$ is `{\it realizable}' in the complex projective space $\mathbf P^n$ or not.}\sk 

To simplify our discussion, we assume that ${\rm Conf}_{\boldsymbol{\mathfrak M}}(\mathbf P^n)$ is irreducible. 
A rational map $\pi_{I,J}$ as in \eqref{Eq:pi-I-J} above is said to be {\bf $\boldsymbol{\mathfrak M}$-admissible} if its restriction to ${\rm Conf}_{\boldsymbol{\mathfrak M}}(\mathbf P^n)$ is 
well-defined\footnote{That is, if ${\rm Conf}_{\boldsymbol{\mathfrak M}}(\mathbf P^n)$ is not included in the indeterminacy locus of $\pi_{I,J}$.} and non-constant. By considering  only the maps of this type among all the $\pi_{I,J}$'s with $(I,J)\in \mathfrak P(M-n-2,n)$, one obtains after restriction a set of rational first integrals on 
${\rm Conf}_{\boldsymbol{\mathfrak M}}(\mathbf P^n)$  which define a web 
 called the {\bf restriction of ${\boldsymbol{\mathcal W}}_{\hspace{-0.05cm}  {\rm Conf}_{M}(\mathbf P^n) }$ to the $\boldsymbol{\mathfrak M}$-strata} and denoted 
 by $
 {\boldsymbol{\mathcal W}}_{\hspace{-0.05cm}  {\rm Conf}_{\boldsymbol{\mathfrak M}}(\mathbf P^n) }$: one has 
 \begin{equation} 
 \label{Eq=W-Conf_M(PN)} 
 {\boldsymbol{\mathcal W}}_{\hspace{-0.05cm}  {\rm Conf}_{\boldsymbol{\mathfrak M}}(\mathbf P^n) }=\Big(
  {\boldsymbol{\mathcal W}}_{\hspace{-0.05cm}  {\rm Conf}_{M}(\mathbf P^n) }
\Big)\big\lvert_{
\boldsymbol{\mathfrak M}
}
=
  {\boldsymbol{\mathcal W}}\Big( \, \big(\pi_{I,J}\big)\lvert_{\boldsymbol{\mathfrak M}}
\, \big\lvert \hspace{0.15cm } \pi_{I,J} \mbox{ is } 
{\boldsymbol{\mathfrak M}}\mbox{-admissible}
\,   \Big)\, .
 \end{equation}

 Given the matroid $\boldsymbol{\mathfrak M}$, it is not clear   
how to compute the dimension(s) of  the associated stratum/strata of degenerate configurations,  
 how to determine the  pairs 
$ (I,J)$ for which the map $\pi_{I,J}$ is $\boldsymbol{\mathfrak M}$-admissible, and even more complicated, given two distinct such admissible maps, to know when they define the same foliation when restricted to $ {\rm Conf}_{\boldsymbol{\mathfrak M}}(\mathbf P^n) $. 
 Consequently, determining the type of the web ${\boldsymbol{\mathcal W}}_{\hspace{-0.05cm}  {\rm Conf}_{\boldsymbol{\mathfrak M}}(\mathbf P^n) }$ from 
 $\boldsymbol{\mathfrak M}$ in full generality does not seem to be an easy task.\mk 
 
 For concrete examples of interesting webs obtained by means of this construction,  see \S\ref{Par:Goncharov-22} (Figure \ref{Fig:CSKG} and the paragraph just after it therein) and 
 also the web  ${\boldsymbol{{\mathcal U} \hspace{-0.04cm} {\mathcal X }
\hspace{-0.04cm}
{\mathcal W}}}_{A_3}$ 
 of Theorem 
 \ref{T:classical-cluster-webs} which 
  can be proved to be of the form \eqref{Eq=W-Conf_M(PN)}  thanks to Proposition  \ref{Prop:UA3-KV}.


\subsection{\bf Web attributes}
\label{SS:WebAttributes}
Here, we define and discuss briefly  some objects and invariants associated to webs which  are central in view of the approach undertaken in the sequel of this memoir. \mk

If everything below is elementary, there are some novelties, such as the notions of `virtual ranks' and of `web with AMP rank' which, as numerous results in this text show,  appear as being relevant to extend the study  of webs of maximal rank of classical web-geometry to a more general class of webs.

\subsubsection{Singular set and the union of common leaves.}  
\label{Par:SSet-UnionCommonLeaves}
Let $\boldsymbol{\mathcal W}$ be a web on an open domain  $U$ of 
$\mathbf C^n$, 
formed 
by $d$ distinct  foliations $\mathcal F_1,\ldots,\mathcal F_d$ (possibly with singularities).  By definition, the {\bf singular set} of $\boldsymbol{\mathcal W}$, denoted by $\Sigma(\boldsymbol{\mathcal W})$ is the set of points where at least two of the foliations of $\boldsymbol{\mathcal W}$ do not intersect transversally.   For  $u\in\Sigma(\boldsymbol{\mathcal W}) $, two possibilities can occur :  either (i) $u$ belongs to the singular set of one of the $\mathcal F_i$'s, or (ii)  it is not the case and there exist two distinct indices $i,j$ such that $T_{ \hspace{-0.04cm} u}   \mathcal F_i=T_{ \hspace{-0.04cm} u} \mathcal F_j$.  
We define the {\bf union of common leaves} of $\boldsymbol{\mathcal W}$ (on $U$), denoted by $\Sigma^c(\boldsymbol{\mathcal W})$ as the closure of the union of the irreducible subvarieties of  $U \setminus \big( \cup_{i} {\rm Sing}(\mathcal F_{ \hspace{-0.04cm} i }) \Big)$  which are invariant by (at least) two distinct foliations of $\boldsymbol{\mathcal W}$.  It is an analytic (in general proper) subset  of $\Sigma(\boldsymbol{\mathcal W})$.  
Finally, for each $i$, one denotes by  $\Sigma^c_i(\boldsymbol{\mathcal W})$ the union of the components of $\Sigma^c(\boldsymbol{\mathcal W})$ which are invariant by $\mathcal F_i$. Clearly, one has $\Sigma^c(\boldsymbol{\mathcal W})=\cup_i \Sigma^c_i(\boldsymbol{\mathcal W})$.
\medskip 

In the case of Bol's web, $\Sigma(\boldsymbol{\mathcal B})$ coincides with $\Sigma^c(\boldsymbol{\mathcal B})$ and is the arrangement in $\mathbf P^2$ formed by the six lines joigning two of the four base-points\footnote{Explicitly, these four base-points are $[0,0,1], \, [1,1,1], \, [1,0,0]$ and $[0,1,0]$ hence 
$\Sigma(\boldsymbol{\mathcal B})$ is cut out by $xyz(x-y)(x-z)(y-z)=0$  in homogeneous coordinates: it is the arrangement of hyperplanes of type $A_2$, {\it cf.} \cite{Pereira}.} of the pencil of conics $
x(1-y) = \lambda\cdot y(1-x)$, $\lambda\in \mathbf P^1$.
For a precise description of the singular set $\Sigma({\boldsymbol{\mathcal W}}_C)$ of the algebraic web ${\boldsymbol{\mathcal W}}_C$ associated to an 
algebraic curve $C\subset \mathbf P^n$, see \cite[\S2]{Nakai}.


\subsubsection{Intrinsic dimension.}  
\label{SS:Intrinsic-Dim}
Let $\boldsymbol{\mathcal W}$ be a (possibly singular) $d$-web defined on a connected variety $M$ of dimension $n\geq 2$.  For $m\in M\setminus \Sigma(\boldsymbol{\mathcal W})$, there exist $d$ germs of holomorphic submersions $u_i: (M,m)\rightarrow \mathbf C$ such that 
$\boldsymbol{\mathcal W}(u_1,\ldots,u_d)$ coincides with the germ of $\boldsymbol{\mathcal W}$ at $m$.  One defines the {\bf intrinsic dimension of $\boldsymbol{\mathcal W}$ at $\boldsymbol{m}$}, denoted by ${\bf IntDim}_m(\boldsymbol{\mathcal W})$, as the rank of the map 
$(u_i)_{i=1}^d : (M,m)\rightarrow \mathbf C^d$ at this point.  One verifies that this quantity does not depend on the choice of  the $u_i$'s but only on (the first order jet at $m$ of) $\boldsymbol{\mathcal W}$ hence it is well-defined. Moreover, one has $2\leq {\bf IntDim}_m(\boldsymbol{\mathcal W})\leq n$, the lower bound being an immediate consequence of  (wGP). 
One defines the {\bf intrinsic dimension ${\bf IntDim}(\boldsymbol{\mathcal W})$ of $\boldsymbol{\mathcal W}$},  as ${\bf IntDim}_m(\boldsymbol{\mathcal W})$ for $m$ generic in $M$. Clearly, one has 
${\bf IntDim}(\boldsymbol{\mathcal W})=\max_{m\in M\setminus \Sigma(\boldsymbol{\mathcal W}) } \{\, {\bf IntDim}_m(\boldsymbol{\mathcal W})
\,\}$.  \sk

An obvious remark is that  ${\bf IntDim}(\boldsymbol{\mathcal W})=\dim(M)=n$ is equivalent to the fact that  $\boldsymbol{\mathcal W}$ admits a $n$-subweb satisfying the `strong general position' hypothesis (sGP) ({\it cf.}\,\S\ref{Par:GermOfWebs}).  
The webs studied until now in classical web geometry were all assumed  to satisfy (sGP) hence were  all of maximal intrinsic dimension.   Therefore the notion of `intrinsic dimension' is  relevant only regarding  webs satisfying only (wGP) but not (sGP) which explains why it did not appear  in web geometry before the present work.\mk 

When $n'={\bf IntDim}(\boldsymbol{\mathcal W})<n$, the intrinsic dimension interprets itself more geometrically as follows: locally, the $d$-web $\boldsymbol{\mathcal W}$  is the pull-back of a $d$-web  
$\boldsymbol{\mathcal W}'$ in $n'$ variables and of maximal intrinsic dimension.  More precisely,  for any simply connected domain $U\subset M\setminus \Sigma(\boldsymbol{\mathcal W})$ such that the intrinsic dimension of $\boldsymbol{\mathcal W}$ at $u$ is $n'$ for any $u\in U$, there exists a submersion $\varphi_U: U\rightarrow \mathbf C^{n'}$ and a $d$-web $\boldsymbol{\mathcal W}'_U$ on $U'={\rm Im}(\varphi_U)$ of intrinsic dimension $n'$ (at any point of $U'$) such that 
$\boldsymbol{\mathcal W}\lvert _U=\varphi_U^*(\boldsymbol{\mathcal W}'_U)$: 
  $\boldsymbol{\mathcal W}'_U$ is called the {\bf intrinsic reduction of $\boldsymbol{\mathcal W}$ over $U$}. 
The fact that $(\varphi_U,\boldsymbol{\mathcal W}'_U)$ is essentially unique 
(as can easily be verified) implies first that defining $\boldsymbol{\mathcal W}'_U$ as the {\bf intrinsic reduction of $\boldsymbol{\mathcal W}$ over $U$} makes sense. Second, it allows to glue (in a certain sense which it is not difficult to make precise) the reductions  $\boldsymbol{\mathcal W}'_U$ obtained in this manner.  
Since all the  web-theoretic notions and properties of webs we are interested in in this text (linearizability/algebraizability, ARs, classical and virtual rank(s), being AMP (see below)) are of local nature, 
 the study of $ \boldsymbol{\mathcal W}$ from this perspective can be reduced to that of any of its intrinsic reductions $ \boldsymbol{\mathcal W}'_U$.


\subsubsection{Abelian relations and rank.}  ${}^{}$ \hspace{-0.4cm}
 Let $\boldsymbol{\mathcal W}=\boldsymbol{\mathcal W}(u_1,\ldots,u_d)$ be a germ of $d$-web on $(\mathbf C^n,0)$ defined by $d$ submersions $u_i\in \mathcal O(\mathbf C^n,0)$.  By definition, an {\bf abelian (functional) relation} for $\boldsymbol{\mathcal W}$ is a
 $d$-tuple of functions in one variable $(F_i)_{i=1}^d$ satisfying identically  a functional relation 
 of the following form:
\begin{equation}
\label{E:EFA}
\sum_{i=1}^d F_i(u_i)\equiv {\rm cst}\, .
\end{equation}

 It is natural to consider  such  a relation as {\bf trivial} when all the components $F_i$'s are constant.  The natural approach in web geometry consisting in dealing with abelian relations up to the trivial ones, we will identify  the 
 abelian relation \eqref{E:EFA} (which will be said in `functional form') with the following differential identity 
 \begin{equation}
\label{E:dEFA}
\sum_{i=1}^d F_i'(u_i) \, du_i\equiv 0 
\end{equation}
which will also be called an `abelian relation' (but in `classical (differential) form'). 
In the whole memoir, the expression `abelian relation' will be abbreviated to AR.  
\sk

%


   The space ${\mathcal A}(\boldsymbol{\mathcal W})$ of abelian relations of $\boldsymbol{\mathcal W}$ has  a natural structure of complex vector space. Thus one can define the {\bf rank} ${\bf rk}({\mathcal W})$ of $\boldsymbol{\mathcal W}$ as the dimension of this space: 
\begin{equation}
\label{E:rankW}
{\rm rk}\big(\boldsymbol{\mathcal W}\big)=\dim_{\mathbf C}\, {\mathcal A}\big(\boldsymbol{\mathcal W}\big) \in \mathbf N \cup \{ \infty \}
\, .
\end{equation}

Here is some terminology about 
abelian relations   $\boldsymbol{F}=(F_i)_{i=1}^d \in  {\mathcal A}(\boldsymbol{\mathcal W})  $:
\begin{itemize}
\item for $i=1,\ldots,d$, the {\bf $\boldsymbol{i}$-th component} of $\boldsymbol{F}$  is of course $F_i$, but 
considered up the addition of a complex constant; equivalently, it is the derivative $F_i'$;  \item the {\bf support}  
of $\boldsymbol{F}$, denoted by $\boldsymbol{{\bf Supp}(F)}$, is the subweb of $\boldsymbol{\mathcal W}$ formed by the foliations $\mathcal F_{u_i}$ associated to the  non-trivial component  $F_i$ of $\boldsymbol{F}$; equivalently, it is the smallest subweb of 
 $\boldsymbol{\mathcal W}$ (for the inclusion) admitting $\boldsymbol{F}$ as an abelian relation;
 \item the {\bf length} $\boldsymbol{\ell(F)}$ of $\boldsymbol{F}$ is the degree (that is the number of foliations) of its support; 
\item  $\boldsymbol{F}$ is  {\bf complete} if none of its components $F_i$ is trivial or, equivalently, if it has length $d$;
\item finally,  $\boldsymbol{F}$ will be said to be {\bf irreducible} if it 
cannot be written as a linear combination of elements of ${\mathcal A}(\boldsymbol{\mathcal W})  $ of length strictly less than that of $\boldsymbol{F}$.
\end{itemize}
\medskip

If we have defined all the notions above by means of some first integrals, one has to be aware (and this can be verified quite easily) that these notions only depend on 
$\boldsymbol{\mathcal W}$, and actually only on the analytic equivalence class of the latter. More generally, if $\boldsymbol{\mathcal W}$ is the pull-back of another web 
$\boldsymbol{\mathcal W}'$ under a holomorphic submersion $\varphi$, then the latter induces naturally a linear map  $\varphi^*: {\mathcal A}(\boldsymbol{\mathcal W}')\rightarrow {\mathcal A}(\boldsymbol{\mathcal W})$ which can easily be seen to be an isomorphism. In particular, it follows that both $\boldsymbol{\mathcal W}=\varphi^*(\boldsymbol{\mathcal W}')$ and $\boldsymbol{\mathcal W}'$ have the same rank.

\begin{exms}
{1.} The  case of (the germ of) an algebraic web is quite 
enlightening since its abelian relation can be interpreted in terms of well-known objetcs of classical algebraic geometry. 
\smallskip 

Let ${\boldsymbol{\mathcal W}}_{C,H_0}$ be the germ of ${\boldsymbol{\mathcal W}}_C$ at a hyperplane $H_0\in \check{\mathbf P}^n$ intersecting 
a given degree $d$ algebraic curve $C\subset \mathbf P^n$ transversely in $d$ points $p_1(H_0),\ldots,p_d(H_0)$.  There exist $d$ local submersions $p_i : ( \check{\mathbf P}^n,H_0)\rightarrow (C,p_i(H_0))$ such that as a 0-cycle on $C$, one has $C\cdot H=\sum_i p_i(H)$ for every hyperplane $H
 $ sufficiently close to $H_0$.  Then, as proved by Abel, given any abelian differential $\omega$ on $C$\footnote{These differential 1-forms are called `{\it differentials of the first kind}' in the classical literature.   In a more modern and rigorous way,  such a form can be defined as a global section of the sheaf $\omega_C^1$ of {\bf regular} (or {\bf Rosenlicht}) {\bf differentials} on $C$. 
When the latter is smooth,  one has $\omega_C^1=\Omega_C^1$ and such a 
 $\omega$ 
is nothing else but a global holomorphic 1-form on the curve $C$.}, 
   the {\bf abelian sum} $\sum_i \int_{p_i(H_0)}^{p_i(H)} \omega$ vanishes identically 
 as a function of $H$,  in the vicinity of $H_0$. \smallskip

 From this, it follows that the linear map $H^0\big(C,\omega_C^1\big)\rightarrow \mathcal A({\boldsymbol{\mathcal W}}_{C,H_0}):\, \omega\mapsto \big( \int^{p_i} \omega \,\big)_{i=1}^d$ is well-defined. It can be proved that it induces an isomorphism  which allows
 to interpret  the rank of the germ of ${\boldsymbol{\mathcal W}}_C$ at a generic hyperplane as the arithmetic genus 
  of $C$: 
one has
 $$
 H^0\Big(C,\omega_C^1\Big)\simeq  \mathcal A\Big({\boldsymbol{\mathcal W}}_{C,H_0}\Big) \qquad \quad \mbox{ and  }\qquad  \quad 
 p_a(C)=h^0(\omega_C^1)={\rm rk}\Big({\boldsymbol{\mathcal W}}_{C,H_0}\Big)\, .\mk 
 $$
 {2.} Another fundamental example is the identity $(\boldsymbol{\mathcal Ab})$ which can be seen as an AR for Bol's web $\boldsymbol{\mathcal B}$. 
 We will 
discuss in depth the 
 ARs
 of (another web equivalent to) this web further in \S\ref{SS:cluster-webs-type-A2}.
 \end{exms} 

Let us now consider the case of a global web $\boldsymbol{\mathcal W}$ defined by $d$ global foliations $\mathcal F_i$ on a domain $U\subset \mathbf C^n$.   Then it can be proved (see \cite[\S1.2.2]{PirioThese} for instance)  
that a germ of AR at any point $u$ in $ U\setminus \Sigma(\boldsymbol{\mathcal W})$ extends holomorphically along  any path in this set (possibly as a multivalued AR). Consequently the complex vector spaces of (germs of) abelian relations $\mathcal A({\boldsymbol{\mathcal W}}_u)$ of the germs of $\boldsymbol{\mathcal W}$ at points $u\in U\setminus \Sigma(\boldsymbol{\mathcal W})$ glue together to form a local system on  $U \setminus \Sigma(\boldsymbol{\mathcal W})$, denoted by $\mathcal A(\boldsymbol{\mathcal W})$.  The {\bf rank} of $\boldsymbol{\mathcal W}$ is the rank of this local system: it coincides with the rank of the germ ${\boldsymbol{\mathcal W}}_u$ at any point $u$  outside $\Sigma(\boldsymbol{\mathcal W})$. \medskip 

Finally, one  will notice that since the $i$th-component of a given abelian relation of $\boldsymbol{\mathcal W}$  is constant along the leaves of $\mathcal F_i$, one gets the stronger fact that the $i$th-component of a germ  of an  abelian relation of $\boldsymbol{\mathcal W}$ at a regular point  extends holomorphically along any path in $\Sigma_i^c(\boldsymbol{\mathcal W})$, this for any~$i$.

\subsubsection{Virtual ranks.}  ${}^{}$ \hspace{-0.4cm}
\label{Para:Virtual-Rank}
Set $V=\mathbf C^n$ and let $\boldsymbol{\mathcal W}=(\mathcal F_i)_{i=1}^d$ be a germ of $d$-web at the origin of $V$. 
 For any $i$, there exists  a non-trivial linear form $\ell_i \in V^* $  (uniquely determined up to multiplication by an element of $\mathbf C^*$), such that ${\rm ker}(\ell_i)=T_0\mathcal F_i
\subset T_0 V\simeq V$. \smallskip 

For any $i$, let $u_i$ be a local submersion such that $\mathcal F_i=\mathcal F(u_i)$ with $u_i(0)=0\in \mathbf C$. 
Then for any $\sigma\in \mathbf N_{>0}$, let $\mathcal A^{\sigma \leq}(\boldsymbol{\mathcal W})$ be the space of abelian relations of $\boldsymbol{\mathcal W}$ of valuation at least $\sigma$ (at $o\in V$): it is the subspace of $\mathcal A(\boldsymbol{\mathcal W})$ formed by the elements $(F_i)_{i=1}^d \in \mathcal O(\mathbf C,0)^d$ such that  $\sum_i F_i(u_i)\equiv 0$ and with  $F_k(z)={\rm O}_0(z^{\sigma})$ for every $k=1,\ldots,d$.   
 One gets a decreasing filtration 
  $$\mathcal A^{\bullet}(\boldsymbol{\mathcal W})\hspace{0.15cm} : \quad  \mathcal A^{}(\boldsymbol{\mathcal W})=\mathcal A^{1 \leq}(\boldsymbol{\mathcal W})\supset \mathcal A^{2 \leq}(\boldsymbol{\mathcal W})\supset \cdots  $$
 and for each $\sigma>0$, the 
 $\sigma$th-piece   $\mathcal A^{\sigma}(\boldsymbol{\mathcal W})=\mathcal A^{\sigma\leq }(\boldsymbol{\mathcal W})/ \mathcal A^{\sigma+1 \leq }(\boldsymbol{\mathcal W})$ of the associated graded vector space embeds linearly into the complex vector space 
 \begin{equation}
 \label{Eq:A-sigma(W)}
 A^\sigma(\boldsymbol{\mathcal W})=\left\{ \, \big(c_1,\ldots,c_d\big) \in \mathbf C^d 
 \hspace{0.1cm}
\Big\vert \hspace{0.1cm}
 \sum_{i=1}^d c_i\big(\ell_i\big)^\sigma=0  \,  
 {\mbox{ in }}  \,
  {\rm{Sym}}^\sigma(V^*)\,
  \right\}
  \end{equation}
 (this embedding is not canonical since the choice of the $\ell_i$'s is not, but it is  essentially unique). \sk 
 
A non-trivial element of $A^{\sigma}(\boldsymbol{\mathcal W})$ is called a {\bf virtual abelian relation of valuation $\sigma-1$} of $\boldsymbol{\mathcal W}$ (at the origin of $ V$) and by definition, the  associated {\bf $\boldsymbol{\sigma}$-virtual rank} of $\boldsymbol{\mathcal W}$  is 
$$
\rho^\sigma(\boldsymbol{\mathcal W})=\dim_{\mathbf C}\,  \Big(A^\sigma(\boldsymbol{\mathcal W}) \Big)=\max\Big( 0\, , \,  d-
\ell^\sigma(\boldsymbol{\mathcal W})
\Big)
$$
where $\ell^\sigma(\boldsymbol{\mathcal W})$ stands for the rank of the subfamily 
 $\{ \, \big(\ell_i\big)^\sigma \,\lvert \, i=1,\ldots,d \, \}$ of $ {\rm Sym}^\sigma(V^*)$. \mk 
 
 The {\bf virtual rank sequence} of $\boldsymbol{\mathcal W}$ is the sequence $\rho^\bullet(\boldsymbol{\mathcal W})$ of the $\rho^\sigma(\boldsymbol{\mathcal W})$'s and their possibly infinite sum is the {\bf virtual rank} $\rho(\boldsymbol{\mathcal W})$  of $\boldsymbol{\mathcal W}$: 
$$
\rho^\bullet(\boldsymbol{\mathcal W})=\big(\rho^\sigma(\boldsymbol{\mathcal W})\Big)_{\sigma>0} \qquad \mbox{ and }\qquad 
\rho(\boldsymbol{\mathcal W})=\sum_{\sigma>0} \rho^\sigma(\boldsymbol{\mathcal W})\in \mathbf N\cup \{\infty \}\, .
$$

We list below some properties of these objects\footnote{Everything can be verified quite easily and this is left to the reader.}: 
\begin{enumerate}
\item for any $\sigma>0$, $\ell^\sigma(\boldsymbol{\mathcal W})$ hence $\rho^\sigma(\boldsymbol{\mathcal W})$ does not depend on the $\ell_i$'s but actually  only 
on the equivalence class of the germ of $\boldsymbol{\mathcal W}$. It follows that 
$\rho^\bullet(\boldsymbol{\mathcal W})$ and $\rho(\boldsymbol{\mathcal W})$ are 
analytic invariants attached to $\boldsymbol{\mathcal W}$;
\item  if $\rho^{\sigma_0}(\boldsymbol{\mathcal W})=0$ for some ${\sigma_0}>0$, then $\rho^{\sigma}(\boldsymbol{\mathcal W})=0$ for any $\sigma\geq {\sigma_0}$; 
\item  one has  $\rho^{\sigma}(\boldsymbol{\mathcal W})=0$ for ${\sigma}>>0$ (actually for $\sigma\geq d-1$);  

\item the bound $\dim \big(\mathcal A^{\sigma}(\boldsymbol{\mathcal W}) \big) \leq \rho^\sigma(\boldsymbol{\mathcal W})$ holds true for any $\sigma\in \mathbf N_{>0}$ and thus ${\rm rk} (\boldsymbol{\mathcal W})\leq \rho(\boldsymbol{\mathcal W})$; 
\item  the previous definitions extend to any  general (germ of) web $\boldsymbol{\mathcal W}$ 
without  requiring that  (wGP) is satisfied. But as soon as two of the $\ell_i$'s 
are proportional, 
 one has $\rho^\sigma(\boldsymbol{\mathcal W})\geq 1$ for any $\sigma$ hence $\rho(\boldsymbol{\mathcal W})=+ \infty$. 
\item  Let $ \boldsymbol{\mathcal W}$ be a $d$-web defined on a  domain $U\subset \mathbf C^n$ 
  and satisfying (wGP) at any point of it. Then 
  denoting by $ {\boldsymbol{\mathcal W}}_u$ the germ of $ \boldsymbol{\mathcal W}$ at $u\in U$, one can set 
 $\ell^\sigma_u( \boldsymbol{\mathcal W})=\ell^\sigma( {\boldsymbol{\mathcal W}}_u)$, 
  $\rho_u^\sigma( \boldsymbol{\mathcal W})=\rho^\sigma( {\boldsymbol{\mathcal W}}_u)
$ (for any $\sigma>0$),  $\rho_u^\bullet( \boldsymbol{\mathcal W})=\rho^\bullet( {\boldsymbol{\mathcal W}}_u)
$ and $\rho_u( \boldsymbol{\mathcal W})=\rho( {\boldsymbol{\mathcal W}}_u)$. 
\item For any $\sigma>0$,  the function $u\mapsto \ell_u^\sigma( \boldsymbol{\mathcal W})$ is  integer-valued and lower semi-continuous on $U$.  Thus,  given $u_0\in U$,  one has $\rho^\sigma_u( \boldsymbol{\mathcal W})\leq \rho^\sigma_{u_0}( \boldsymbol{\mathcal W})$ as well as  $\rho_u( \boldsymbol{\mathcal W})\leq \rho_{u_0}( \boldsymbol{\mathcal W})$ for any element $u$ of $U$ sufficiently close to $u_0$.  Consequently, setting $\rho^\star( \boldsymbol{\mathcal W})=\min_{u\in U} \rho^\star_u( \boldsymbol{\mathcal W})$ where $\star$ stands either for nothing or 
for one of the two symbols $\sigma$ or $\bullet$, one gets well-defined invariants attached to $\boldsymbol{\mathcal W}$.  Clearly, 
one has $ \rho^\star( \boldsymbol{\mathcal W})= \rho^\star_u( \boldsymbol{\mathcal W})$ 
for $u\in U$ generic.
 \smallskip 

\vspace{-0.4cm} 
Clearly, the bound of  4. above extends for the global web $ \boldsymbol{\mathcal W}$:  ${\rm rk}( \boldsymbol{\mathcal W})\leq \rho( \boldsymbol{\mathcal W})$ holds true. 
\item Finally, given a global web $ \boldsymbol{\mathcal W}$ defined on a domain $U \subset \mathbf C^n$, we will write $\rho^\star( \boldsymbol{\mathcal W})$ for $\rho^\star( \boldsymbol{\mathcal W}\lvert_{U'})$ where $U'\subset U$ denotes the dense open subset of points where (wGP) is satisfied. 
\end{enumerate}
\smallskip


\begin{exm}
\label{Exm:Webs}
{\rm (1).} For a germ of $d$-web $\boldsymbol{\mathcal W}$ on $(\mathbf C^2,0)$, one has 
$\rho^\sigma(\boldsymbol{\mathcal W})=\max\big(0, (d-1-\sigma)\big)$ for any  $\sigma\geq 1$,  therefore 
$\rho(\boldsymbol{\mathcal W})=(d-1)(d-2)/2$ is independent of $\boldsymbol{\mathcal W}$.   Thus, in this case,  one cannot get  from the previously introduced invariants $\rho^\bullet(\boldsymbol{\mathcal W})$ and $\rho(\boldsymbol{\mathcal W})$ 
anything  else than the classical Bol's bound ${\rm rk}(\boldsymbol{\mathcal W})\leq \rho(\boldsymbol{\mathcal W})=(d-1)(d-2)/2$ on the rank.
\medskip 

{\rm (2).} In arbitrary dimension, if $\boldsymbol{\mathcal W}$ satisfies the strong general position assumption, then 
one has $\rho^\sigma(\boldsymbol{\mathcal W})\leq \max\big(0, (d-1-\sigma(n-1)\Big)$ for any positive integer $\sigma$. From this, it  follows that $\rho(\boldsymbol{\mathcal W})$ {\rm (}hence a fortiori $r(\boldsymbol{\mathcal W})${\rm )} is bounded from above by Castelnuovo constant $\pi(d,n)$ {\rm (}see {\rm \cite[\S2.3]{Coloquio}} for details and references{\rm )}.
\medskip 

    {\rm (3).} In {\rm \cite{CavalierLehmann}}, the authors consider the notion of {\it `ordinary webs'}. With the notions introduced above, 
 such a web $\boldsymbol{\mathcal W}$ is nothing but a web satisfying (wGP), with $\rho^\sigma(\boldsymbol{\mathcal W})=\max\Big(0 \, , \, d- {n-1+\sigma \choose n-1}\Big)$  for every $\sigma>0$.
\medskip 

{\rm (4).} For any integer $n\geq 2$,  the following equalities hold true {\rm (}{\it cf.} {\rm \cite[Proposition\,4.1]{Pereira}}\footnote{Beware that in several places in \cite{Pereira} (for instance in equation (2) 
 in \S2.1), $`\max'$ has been used instead of $`\min'$.}{\rm ):}
$$
\rho^\sigma\Big({\boldsymbol{\mathcal W}}_{ \hspace{-0.06cm}\mathcal M_{0,n+3}}\Big)=  \begin{cases}
 \,  \,  { n+3 \choose n-1} - {  n-1+\sigma  \choose n-1}  \hspace{0.4cm} \mbox{ for }  \,  \sigma=1,2,3\\ 
\hspace{1cm} 0 \hspace{1.7cm}   \mbox{ for }  \,  \sigma\geq 4\, , 
\end{cases}
$$
and consequently,  we have  
\begin{equation}
\label{E:Pereira-r(W)}
{\rm rk}\Big({\boldsymbol{\mathcal W}}_{ \hspace{-0.06cm}\mathcal M_{0,n+3}}\Big)\leq \rho\Big({\boldsymbol{\mathcal W}}_{ \hspace{-0.06cm}\mathcal M_{0,n+3}}\Big)=3 \, { n+3\choose n-1} 
-{  n+2  \choose n-1} 
- {  n+1  \choose n-1}-  {  n  \choose n-1}
\, .
\end{equation}
\end{exm}

\vspace{0.1cm}
To finish this paragraph, let us stress that the notion of virtual rank $\rho(\boldsymbol{\mathcal W})$ seems to be quite convenient. It is a non-trivial  invariant attached to any web $\boldsymbol{\mathcal W}$ which can be computed quite easily and which allows to state in the simple and uniform manner 
${\rm rk}(\boldsymbol{\mathcal W})\leq \rho(\boldsymbol{\mathcal W})$ 
all the already known bounds on the rank  of 1-codimensional webs  ({\it cf.}\,the statement of \cite[Corollary\,2.1]{Pereira} for instance).

\subsubsection{Webs of AMP rank.}  ${}^{}$ \hspace{-0.4cm}
\label{SS:Webs-AMP-Rank}
Let $\boldsymbol{\mathcal W}$ be a web defined on $U\subset \mathbf C^n$ (or a germ of such a web). 
 Then from the preceding paragraph, we know that 
 \begin{equation}
 \label{Eq:r(W)-leq-rho(W)}
 {\rm rk}(\boldsymbol{\mathcal W})\leq \rho(\boldsymbol{\mathcal W})<+\infty\, .
 \end{equation}
  The rank of $\boldsymbol{\mathcal W}$ is said to be {\bf As Maximal as Possible} (abbreviated AMP\footnote{We have chosen to use this abbreviation since it  also makes sense in French (`rang {\it Aussi Maximal que Possible}').})   if one has 
  ${\rm rk}(\boldsymbol{\mathcal W})=\rho(\boldsymbol{\mathcal W})>0$.  In this case, $\boldsymbol{\mathcal W}$ is said to be of {\bf AMP rank} or more concisely (but a bit abusively), that $\boldsymbol{\mathcal W}$ is an {\it `AMP web'}, or even shorter, that this web {\it `is AMP'}.  
  \medskip
  
  The interest of the notion of AMP rank is clear considering the two following points : \vspace{-0.1cm}
 \begin{enumerate}
 \item[$\bullet$] being of AMP rank 
 is clearly a property invariant up to local equivalence. Even better, this property is stable 
 under pullback: if  $\pi$ is a submersion whose image is contained in the definition domain of an AMP web $\boldsymbol{\mathcal W}$, then 
 $\pi^*(\boldsymbol{\mathcal W})$ is AMP as well.
 \item[$\bullet$] given a web defined by first integrals $u_1,\ldots,u_d: U\rightarrow \mathbf C$, it is quite easy to compute the values of the differential 1-forms $du_i$ at a given point $u\in U$.  From this,  the values $\rho^\sigma_u(\boldsymbol{\mathcal W})$  
 for $\sigma>0$ can be effectively determined by straightforward linear algebra (possibly with the help of a computer algebra system).  It follows that $\rho(\boldsymbol{\mathcal W})$ hence the property of having AMP rank  can be determined quite easily on concrete examples; \vspace{-0.1cm}
  \item[$\bullet$] as the list of examples below shows, the notion of web with AMP rank encompasses in a quite general, systematic and convenient way most of the already introduced notions of `web with maximal rank' previously introduced in web geometry. 
   \end{enumerate}
  
  \begin{exm}
  \label{Ex:toto}
  {\rm 1.}  Webs of maximal rank in the classical sense are of course webs with AMP rank. \mk
  
    {\rm 2.} The main objects of study  in 
    {\rm \cite{CavalierLehmann}} are ordinary webs, in particular those of maximal rank.  
  With the notions introduced above, such  a $d$-web  $\boldsymbol{\mathcal W}$ can be defined as  a web 
 of AMP rank with  
 $\rho^\sigma(\boldsymbol{\mathcal W})= \max\Big(0 \, , \, d- {n-1+\sigma \choose n-1}\Big)$  
  for every  posistive integer $\sigma$. 
\medskip  

 {\rm 3.} The main result of {\rm \cite{Pereira}} (namely Theorem 4.1 therein) can be stated as the fact that for any $n\geq 2$, 
   the web ${\boldsymbol{\mathcal W}}_{ \hspace{-0.06cm}\mathcal M_{0,n+3}}$  has AMP rank. Indeed, 
   in 
   {\rm \cite[\S4.1]{Pereira}} 
   (see also \S\ref{Par:W-M0n+3-all-AMP} below), it is proved that 
   this web carries at least $\rho\big({\boldsymbol{\mathcal W}}_{ \hspace{-0.06cm}\mathcal M_{0,n+3}}\big)$ linearly independent ARs which implies that   
 \eqref{E:Pereira-r(W)} is actually an equality.
%
\end{exm}

One of the main problems in web geometry is  to  classify and study  the so-called {\bf exceptional webs} namely, 
 the  $d$-webs  on a $n$-dimensional manifold,  
the foliations of which 
satisfy the strong general position assumption, 
which are not algebraizable and 
 whose  rank is maximal (that 
is equal to Castelnuovo constant $\pi(d,n)$, the latter being assumed to be positive\footnote{Note that requiring that $\pi(d,n)$ be positive implies that $d$ is sufficiently big.}).   But, according to some famous algebraization results in web geometry (due to Bol \cite{Bol1} in dimension $n=3$ and to Tr\'epreau  \cite{Trepreau} in arbitrary dimension $n\geq 3$), such a web is  necessarily planar and formed of $d\geq 5$ foliations.  This makes the class of such webs, although being rich of many interesting elements ({\it eg.}\,see \cite[Chapter\,6]{Coloquio}), very  particular and specific to the 2-dimensional case. \medskip 

For a web $\boldsymbol{\mathcal W}$ on a $n$-dimensional manifold with $\rho(\boldsymbol{\mathcal W})$ `big' (in a sense to be made precise), the fact  that ${\rm rk}(\boldsymbol{\mathcal W})=\rho(\boldsymbol{\mathcal W})$ means something very strong on the differential system the solutions of which are the AR's of $\boldsymbol{\mathcal W}$, hence on  $\boldsymbol{\mathcal W}$ itself. From this point of view, 
 the truly relevant generalization of the notion of `planar exceptional web' in dimension $n\geq 3$ could be, instead of the most straightforward one, which is empty according to  Bol-Tr\'epreau's algebraization theorem for webs of maximal rank (see  \cite{Bol1} and  \cite[Th\'eor\`eme 1.1]{Trepreau}),  rather that   of `{\it non-algebraizable web with AMP rank}'.\mk

The webs ${\boldsymbol{\mathcal W}}_{ \hspace{-0.06cm}\mathcal M_{0,n+3}}$ are of this kind 
(see Example \ref{Ex:toto}.3.\,just above), as well as numerous examples of webs we are going to consider in the sequel, such as Goncharov's trilogarithmic web 
discussed in \S\ref{Par:Goncharov-22} or 
the $\boldsymbol{\mathcal Y}$-webs of type $A_n$ we will study further on in Section \ref{SS:Y-cluster-web-An}, among many others.

This shows that this class of webs is non-empty and contains webs which, like many examples of planar exceptional webs, are related with functional equations of polylogarithms.
\begin{center}
$\star$
\end{center}

According to us,  the determination of webs with AMP rank is what corresponds, in dimension $n\geq 3$, to the problem of determining the exceptional planar webs, which is the main open question in  planar(=2-dimensional) web geometry,

\subsubsection{ACM curves and AMP webs.}  
\label{SS:ACM-AMP}
${}^{}$ \hspace{-0.4cm}
Considering the notions introduced above for algebraic webs is particularly interesting since it induces a connection with some projective curves of  particular interest.  \mk 

It is interesting to say a few words about the first basic results of Castelnuovo's theory of projective curves (see \cite{Harris,Ciliberto} for more details). Let $C\subset \mathbf P^n$ be a reduced  projective curve of degree $d$ and denote by $\Gamma$ 
 a general hyperplane section of it.  For  $k$ sufficiently big, the values of the Hilbert polynomial $P_C(k)$ and of the Hilbert function $h_C(k)$ coincide and for such a $k$, one has 
$$dk-p_a(C)+1=P_C(k)=h_C(k)=1+\sum_{\ell=1}^k\Big(h_C(\ell)-h_C(\ell-1) \Big)$$ and since  $h_C(\ell)-h_C(\ell-1)\geq h_\Gamma(\ell)$ for any integer $\ell$ 
(see \cite[(1.1)]{Ciliberto}, one deduces the majoration 
$$p_a(C)\leq \sum_{k=1}^\infty \Big(d-h_\Gamma(k)\Big)$$ 
(where the right-hand-side is actually a finite sum).\mk 

Recall that a ring is \href{https://en.wikipedia.org/wiki/CohenÐMacaulay_ring#The_unmixedness_theorem}{\it Cohen-Macaulay} if its \href{https://en.wikipedia.org/wiki/Depth_(ring_theory)}{\it depht} equals its 
\href{https://en.wikipedia.org/wiki/Krull_dimension}{\it Krull dimension}. A  curve $C$ as above is {\bf Arithmetically Cohen-Macaulay} ({\bf ACM} for short) if the coordinates ring $R/ I_C$ of its affine cone is CM\footnote{Here $R$ denotes  
the space of homogeneous polynomials in $n+1$ variables and $I_C$ stands for the saturated ideal of $C$.}.  It is well-known ({\it cf.}\,\cite[p.\,84]{Harris}) that $C$ is ACM if and only if  $h_C(\ell)-h_C(\ell-1)\geq h_\Gamma(\ell)$ for all $\ell\geq 1$ hence if and only if the equality $p_a(C)=\sum_{k=1}^\infty \Big(d-h_\Gamma(k)\Big)$ holds true.\mk 

Classically, $h_\Gamma(\ell)$ is defined as the number of linearly independent conditions imposed by $\Gamma$ on the hypersurfaces of degree $\ell$ in $\mathbf P^n$. Then $h_\Gamma(\ell)$ can be seen as equal to 1 plus the dimension of the linear space spanned by the image of $\Gamma$ by the $\ell$-th order Veronese embedding $v_\ell : \mathbf P^n\hookrightarrow  \mathbf P^{ { n+\ell \choose n}-1 }$: one has $h_\Gamma(\ell)=1+\dim \big(\langle v_\ell(\Gamma) \rangle\Big)$.  On the other hand, the ambiant  $\mathbf P^n$ can be identified with the projectivization $\mathbf P(V^*)$ of the dual of a complex vector space $V$ of dimension $n+1$.  Let $\gamma_1,\ldots,\gamma_d$ be $d$ non-trivial  linear forms on $V$ such that 
$\Gamma=[\gamma_1]+\cdots+[\gamma_d]$ as 0-cycles on $\mathbf P^n$. Then from the preceding remark, it follows that $h_\Gamma(\ell)$ coincides with the dimension of the vector subspace of ${\rm Sym}^\ell \big( V^*\big)$ spanned by the $\ell$-th power of the linear forms $\gamma_i$: one has 
$$h_\Gamma(\ell)=\dim \Big( \big\langle  
(\gamma_1)^\ell,\ldots,(\gamma_d)^\ell
\big\rangle \Big)\, .$$ 
Denoting by $H\in \check{\mathbf P}^n$ the linear span of $\Gamma$, one gets that the $\gamma_i$'s can be seen as linear forms defining the constant part of the algebraic web $\boldsymbol{\mathcal W}_C$ associated to $C$ at $H$, {\it i.e.}\,$[\boldsymbol{\mathcal W}_C]_H^0=
\boldsymbol{\mathcal W}(\gamma_1,\ldots,\gamma_d)$ 
which gives us  that $d-h_\Gamma(\ell)=
\rho^\ell_H\big(\boldsymbol{\mathcal W}_C \big)
$, and this for any $\ell\geq 1$.\mk 

One can gather  all the  facts considered above together in the following statement: 
\begin{prop} 
Let  $C$ and $\Gamma=C\cap H$ be 
as above: $C$ is a reduced projective curve of degree $d$ in $\mathbf P^n$ and $H$ stands for a general hyperplane. 
\begin{enumerate}
\item For any $\sigma \geq 1$, one has  $\rho^\sigma_H\big(\boldsymbol{\mathcal W}_C \big)=d-h_\Gamma(\sigma)$.
\sk 
\item  Thus the majoration $r\big(\boldsymbol{\mathcal W}_C\Big)\leq \sum_{\sigma\geq 1}\rho^\sigma_H\big(\boldsymbol{\mathcal W}_C\big)$ coincides exactly, but stated dually for 
 $\boldsymbol{\mathcal W}_C$, with the majoration $p_a(C)\leq \sum_{\sigma \geq 1}
\big(d - h_\gamma(\sigma) \big)
$
of the classical theory of projective curves.
\sk 
\item As an immediate corollary, one gets that the following equivalence holds true: 
\begin{equation*}
\mbox{ the curve }\, 
C\,  \mbox{ is ACM }\hspace{0.15cm}  \Longleftrightarrow 
\hspace{0.15cm} 
\mbox{ the web }\, 
\boldsymbol{\mathcal W}_C \,  \mbox{ is AMP}.
\end{equation*}
\end{enumerate}
%
\end{prop}

 Let us discuss an explicit example taken from \cite{GHL}:
\begin{exm} 
\label{Ex:Courbe-D-ACM}
The image of $\mathbf P^1\ni [\zeta:z]\mapsto [\zeta^8:\zeta^3z^5:\zeta z^7:z^8]\in \mathbf P^3$ is an  ACM singular rational curve of degree $d=8$ and genus $g=7$, noted by $D$. Its ideal is generated by 
$x_0x_3^3 -x_1^2x_2^2$, $ x_3^2x_1-x_2^3$ and $x_1^3-x_0x_2x_3$ from which one gets   
 (using \href{http://www2.macaulay2.com/Macaulay2/}{Macaulay2}) that its Hilbert series $h_D$ is given by 
 $$h_D(t)=(1-2t^3-t^4+2t^5)/(1-t^4)= 1+4 t+10 t^2+18 t^3+26 t^4+34 t^5+O(t^6)\, . $$
From this, we deduce that $\rho^\bullet(\boldsymbol{\mathcal W}_D)=(5, 2,0,\ldots)$ thus 
$\rho(\boldsymbol{\mathcal W}_D)=5+2=g$ as expected. 
 \end{exm}

A very classical question regarding the theory of projective curves is about the determination of the pairs $(d,g)$ such that there exists  a smooth curve 
$C\subset \mathbf P^n$ with 
 $d=\deg(C)$ and $g=g(C)$ (see \cite{GrusonPeskine} and the reference therein).   An (almost) complete answer to this problem for ACM curves is given  for curves in $\mathbf P^3$ in  \cite[\S2]{GrusonPeskine}\footnote{The {\bf postulation index} of a projective curve 
  $C$  is $s=s(C)=\min\{ \, k\in \mathbf Z\, \lvert \, h^0(\mathcal I_C(k))>0\, \}$ where $\mathcal I_C$ stands for the sheaf of ideals defining $C$.  The results of \cite[\S2]{GrusonPeskine} concern smooth ACM space curves. The only thing missing in  it regarding the determination of such curves is a description of a geometric construction of the general 
  smooth ACM curve $C\subset \mathbf P^3$ of degree $d$  when $d$ is small, namely when $d\leq s(s-1)$ (see the last sentence page 46 of \cite{GrusonPeskine}).} but as far as we know, the higher dimensional cases are still open, even if there have been some recent advances for curves in $\mathbf P^4$.\bk

As already mentioned before, Chern-Griffiths' problem of determining 
the $d$-webs of maximal rank in the classical sense is now solved in dimension $n\geq 3$. Indeed,  such a web is necessarily algebraic (for $d$ big enough) according to Bol-Tr\'epreau's theorem; and the algebraic webs of this kind are those attached to the projective curves of maximal genus which have been classified by Castelnuovo. In particular, this shows that except in dimension 2, the web-theoretic approach  of webs with maximal rank is not 
 really interesting and relevant since essentially everything can be reduced  to some considerations 
(moreover very specific) about  projective curves. \bk

We believe that the question of determining the (1-codimensional) AMP webs is  the appropriate generalization in arbitrary dimension $n\geq 2$ of the classical main problem in web geometry of determining the exceptional planar webs.  Indeed: 
\begin{itemize}
\item[$\bullet$] this new problem admits, as a particular subcase, the quite rich and difficult one of determining the ACM curves in $\mathbf P^n$; 
\item[$\bullet$]  as shown by many of the results of the present paper, there are many AMP webs which are not algebraizable, which are {\it `polylogarithmic'} (in a sense introduced below). 
Perhaps even more interestingly, via the theory of cluster algebras and from  Dynkin diagrams, one can construct several infinite families of {\it `cluster webs'} in arbitrary dimensions which are conjectured to be of AMP rank and also to be polylogarithmic.\footnote{For now, this has been proved only for 
two families: the one of type $A$ $\mathcal X$-cluster webs $\boldsymbol{\mathcal X \hspace{-0.1cm}\mathcal W}_{A_n}$ with $n\geq 2$ (in \cite{Pereira}) and the associated family of $\mathcal Y$-cluster subwebs $\boldsymbol{\mathcal Y\hspace{-0.1cm}\mathcal W}_{A_n}$ (see \S\ref{SS:Y-cluster-web-An} further).} 
\end{itemize}

 \begin{questions}
 \label{Question:1.3}
 Let the ambiant dimension $n\geq 3$ be fixed.\vspace{-0.15cm}
  \begin{enumerate}
 \item[\bf (1).] 
Are there pairs $(d, g)$ such that there is no ACM curve in $\mathbf P^n$ of degree $d$ and genus $g$ 
while  AMP $d$-webs in dimension $n$ of rank $g$ do exist ?
 \item[\bf (2).] If yes, determine such {\it `AMP pairs'} $(d,g)$, for instance when $n=3$. 
 \end{enumerate}
 \end{questions}


%
%

 Let $C$ be a reduced  ACM curve in $ \mathbf P^3$. We recall some attributes to $C$ and explain how these are interpreted in terms of the dual web $\boldsymbol{\mathcal W}_C$.  We follow \cite{GrusonPeskine} to which we refer  for more details.\sk

 Let $s=s(C){>0}$ be the postulation of $C$ and denote by 
$\boldsymbol{n}=
(n_i)_{i=0}^{s-1}$  its {\it `numerical character}' as defined in \cite[Def.\,2.4]{GrusonPeskine}: it is a decreasing sequence of integers satisfying several other properties\footnote{For instance, one has: $(i)$ $n_0=3+e$ with $e=e(C)={\max}\{ k\geq 0\, \lvert \,  \, h^1(\mathcal O_C(k))>0 \,\}$; $(ii)$ $n_{s-1}\geq s$; and $(iii)$ $\boldsymbol{n}$  is without gaps  ({\it i.e.}\,one has $n_i\leq n_{i+1}+1$ for $i=0,\ldots,s-2$) if $C$ is assumed to be irreducible.}  from which the arithmetic invariants of $C$ can all be explicitly recovered. Indeed, one has: 
\begin{itemize}
\item the degree $d=\deg(C)$ and the arithmetic genus $g=p_a(C)$ of $C$ satisfy 
   \begin{equation}
 \label{Eq:ACM-(d,g)}
 d=d_{\boldsymbol{n}}=
  \sum_{i=0}^{s-1}\big(n_i-i\big)\qquad  
\mbox{ and  } \qquad g=g_{\boldsymbol{{n}}}
=1+\frac{1}{2}\sum_{i=0}^{s-1}\big(n_i-i\big)\big(n_i+i-3\big)\, .
\end{equation} 
\item setting 
$h^1_{\boldsymbol{n}}(k) = \sum_{i=0}^{s-1}
\big(n_{i } - k - 1 \big)_+ -\big( k-i\big)_+$ for $k\geq 0$ (with $m_+=\max(0,m)$ for any $m\in \mathbf Z$),  the Hilbert function of $C$ coincides with the function $h_{\underline{n}}:\mathbf N\rightarrow \mathbf N$ defined by 
 \begin{equation}
 \label{Eq:ACM-H}
 h_{\boldsymbol{n}}(0)=1
  \qquad  
\mbox{ and  } \qquad 
 h_{\boldsymbol{n}}(k)
=k \cdot d_{\boldsymbol{n}}+1-
g_{\boldsymbol{n}} +\sum_{\ell\geq k+1} 
h^1_{\boldsymbol{n}}(\ell)
\quad \mbox{for }\, k\geq 1\, .
\end{equation} 
 \end{itemize}
 A direct consequence of the preceding formula is the following formula for the virtual ranks of the web $\boldsymbol{\mathcal W}_C$ associated to $C$: assuming that $C$ spans the whole ambiant space $\mathbf P^3$,   one has 
  \begin{equation}
 \label{Eq:ACM-Rho-sigma}
\rho^1\Big( \boldsymbol{\mathcal W}_C
\Big)=d-3 \qquad \mbox{ and }\qquad 
\rho^\sigma\Big( \boldsymbol{\mathcal W}_C
\Big)=h^1_{\boldsymbol{n}}(\sigma) \quad \mbox{for }\,\sigma\geq 2\, .
\end{equation} 


 Let us  consider the explicit example considered above from this point of view: 
 \begin{exm}[Continuation of Example \ref{Ex:Courbe-D-ACM}]
The postulation of $D \subset \mathbf P^3$ is 3 and $e=e(D)=1$. Therefore the numerical character of this rational curve is a 3-tuple $\boldsymbol{n}=(n_0,n_1,n_2)$ with $n_0=3+e=4$ thus one has $4=n_0\geq n_1\geq n_2\geq 3$.  From \eqref{Eq:ACM-(d,g)}, it follows that 
$\boldsymbol{n}=(4,4,3)$ and by direct computations, one gets that  $h^1_{\boldsymbol{n}}(2)=2$ and $h^1_{\boldsymbol{n}}(\ell)=0$ for $\ell\geq 3$.  From this, one retrieves that 
$\rho^\bullet(\boldsymbol{\mathcal W}_D)=(d-3,h_{\boldsymbol{n}}^1(2))=(5,2)$. 
  \end{exm}

 From the discussion above, one deduces the 

 \begin{prop} 
 \label{Prop:W-algeb-Ari-Conditions}
 A necessary condition for an 
 AMP $d$-web $\boldsymbol{\mathcal W}$ of intrinsic dimension 3  to be algebraizable is that there exist $s>1$ such that $s(s-1)\leq d$ and  
 a $s$-tuple $\boldsymbol{n}=(n_i)_{i=0}^{s-1}\in \mathbf N^s$ 
 such that $ n_0\geq n_1\geq \cdots n_{s-1}\geq s$ and verifying $d=d_{\boldsymbol{n}}$, ${\rm rk}(\boldsymbol{\mathcal W})=g_{\boldsymbol{n}}$ and $\rho^\sigma(\boldsymbol{\mathcal W})=h^1_{\boldsymbol{n}}(\sigma)$ for all $\sigma\geq 2$. 
 \end{prop}

The interest of this proposition lies in the fact that it can be used to prove 
 that  a given AMP web  of intrinsic dimension 3 is not algebraizable 
solely by means of arithmetic considerations. 
 For instance:
 \begin{itemize}
 \item  the web $\boldsymbol{\mathcal W}_{\hspace{-0.02cm}\mathcal M_{0,6}}$ is a 15-web in 3 variables  which is AMP and of rank 26 ({\it cf.}\,Example \ref{Ex:toto}.3). On the other hand, 
 one verifies that for any $s=2,3,4$, there does not exist any $s$-tuple 
 $\boldsymbol{n}=(n_i)_{i=0}^{s-1}$ with  
 $n_0\geq n_1\geq \cdots n_{s-1}\geq s$ and such that $d_{\boldsymbol{n}}=15$ and 
 $g_{\boldsymbol{n}}=26$. It follows that 
 $\boldsymbol{\mathcal W}_{\hspace{-0.02cm}\mathcal M_{0,6}}$ is not algebraizable; 
 \item  Goncharov's trilogarithmic web 
 $ {\boldsymbol{\mathcal W}}_{ \boldsymbol{\mathcal G}_{22}}$
 that we will consider further on in \S\ref{Par:Goncharov-22} is an AMP 22-web of intrinsic dimension 3 whose rank is 56. For $s=2,3,\ldots,5$, there is no $s$-tuple $\boldsymbol{n}=(n_i)_{i=0}^{s-1}$ as above such that $d_{\boldsymbol{n}}=22$ and 
 $g_{\boldsymbol{n}}=56$ hence  
 $ {\boldsymbol{\mathcal W}}_{ \boldsymbol{\mathcal G}_{22}}$  is not algebraizable.
 \end{itemize}
%
%
%
%
%
 
 Actually, since a linearizable web with a complete abelian relation is necessarily algebraizable (a well-known consequence of Abel-inverse theorem for webs), it follows that both $\boldsymbol{\mathcal W}_{\hspace{-0.01cm}\mathcal M_{0,6}}$ and $ {\boldsymbol{\mathcal W}}_{ \boldsymbol{\mathcal G}_{22}}$ are not linearizable. Thanks to Proposition \ref{Prop:W-algeb-Ari-Conditions} and to the fact that these two webs are AMP, we have obtained this result just by dealing with a finite number of arithmetic conditions whereas  the unique previous method we were aware of  to do so 
 ({\it cf.}\,\cite{PirioLin}) relied on the computations of some differential invariants.
  \mk
  
  However, the strategy sketched above to establish the non algebraizability of a web does not apply systematically, 
 for instance in the case of  the $\boldsymbol{\mathcal Y}$-cluster web of type $A_3$, which is the second member of a family of AMP webs associated to the Dynkin diagrams of type $A$ (see \S\ref{SS:Y-cluster-web-An}). 
This web, denoted by $\boldsymbol{\mathcal Y \hspace{-0.1cm}\mathcal W}_{A_3}$,  is an AMP  9-web in 3 variables, with rank 10 and such that 
 $\rho^\bullet(\boldsymbol{\mathcal Y \hspace{-0.1cm}\mathcal W}_{A_3})=(6,3,1)$.  It can be verified that $\boldsymbol{n}=(5,4,3)$ is the (unique) decreasing tuple $ \boldsymbol{n}=(n_i)_{i=0}^{s-1}$ with $s\in \{2,3\}$ such that $d_{\boldsymbol{n}}=9$, $g_{\boldsymbol{n}}=10$ and 
 $h_{\boldsymbol{n}}^1(\sigma)=\rho^\sigma(\boldsymbol{\mathcal Y \hspace{-0.1cm}\mathcal W}_{A_3})$ for $\sigma\geq 2$.  However  $\boldsymbol{\mathcal Y \hspace{-0.1cm}\mathcal W}_{A_3}$ is not linearizable hence not algebraizable (see Proposition \ref{Prop:YWD-non-linearizable} further on) but this cannot be deduced from Proposition \ref{Prop:W-algeb-Ari-Conditions}.
 \begin{center}
 $\star$
 \end{center}
 
 The examples above show that the answer to the first of the Questions \ref{Question:1.3} is affirmative which makes the second question the one that really matters.



\subsubsection{Branch loci.}  ${}^{}$ \hspace{-0.4cm}
We consider here the particular case of a global web $\boldsymbol{\mathcal W}$ determined by a $d$-tuple $(u_i)_{i=1}^d$ of rational functions in $n\geq 2$ variables satisfying (wGP) generically. \medskip 

Considering $\boldsymbol{\mathcal W}=\boldsymbol{\mathcal W}(u_1,\ldots,u_d)$ as a global singular web on $\mathbf P^n$, there is  a `space of common leaves' $\Sigma^c_i(\boldsymbol{\mathcal W})$ associated  to each first integral $u_i$ of $\boldsymbol{\mathcal W}$ (see above).  Then a fiber $u_i=\lambda$ with  $\lambda\in \mathbf P^1$ is said to be {\bf common} if it is entirely contained in $\Sigma_i^c(\boldsymbol{\mathcal W})$.  
From the assumption that (wGP) holds true, it follows immediately  
 that the set of $\lambda$'s giving rise to such a fiber is a finite subset of $\mathbf P^1$ denoted by $\mathfrak B_i=\mathfrak B_i(\boldsymbol{\mathcal W})$.   By definition, it is the {\bf $\boldsymbol{i}$-th  branch locus} of the web $\boldsymbol{\mathcal W}$.\footnote{More rigorously, 
$B_i$ depends on the pair $(\boldsymbol{\mathcal W},u_i)$ thus it would have been more accurate to call it the {\it `${u_i}$-branch locus'} of $\boldsymbol{\mathcal W}$.  Keeping this in mind, we will commit the abuse of not mentioning $u_i$ explicitly  in what follows.} 
The interest of considering  these finite subsets $\mathfrak B_i\subset \mathbf P^1$ comes  
from the   
 following proposition:
 \begin{prop} 
 \label{P:AnalyticConinuation}
 For any $x\in \mathbf P^n\setminus \Sigma(\boldsymbol{\mathcal W})$ and any (germ of)  abelian functional relation $\sum_i F_i(u_i)=0$ at $x$, each germ $F_i$ extends as a  global multivalued holomorphic function on the whole projective line $\mathbf P^1$ with (possible) ramification precisely at the points of $\mathfrak B_i$. 
\end{prop}
 \begin{proof} 
First we prove that,  for any (say smooth) path $\gamma:[0,1]\rightarrow \mathbf P^n\setminus \Sigma(\boldsymbol{\mathcal W})$ starting at $\gamma(0)=x$, each holomorphic germ  $\Psi_i= F_i\circ u_i\in \mathcal O_{x}$ extends analytically along $\gamma$, for any $i=1,\ldots,d$.  This is quite standard ({\it cf.}\,\cite[Th\'eor\`eme 1.2.2]{PirioThese}): let  $\zeta>0$ be the supremum of all $s\in[0,1]$ such that each $\Psi_i$ extends analytically along the restriction of $\gamma$ to the interval $[0,s]$. This means that one of the $\Psi_i$'s, say $\Psi_1$, admits a  holomorphically continuation along $\gamma\lvert _{[0,\zeta[}$ which does not extend at $y=\gamma(\zeta)$. Since $\Psi_1$ is constant along the leaves of $\boldsymbol{\mathcal F}_{\hspace{-0.1cm} u_1}$, this implies that $\gamma$ is tangent to this foliation at $y$, {\it i.e.}\,$\gamma'(\zeta)\in T_y  \boldsymbol{\mathcal F}_{\hspace{-0.1cm} u_1}$.\footnote{Note that $\gamma'(\zeta)$ is a real tangent vector whereas $T_y  \boldsymbol{\mathcal F}_{\hspace{-0.1cm} u_1}$ is a complex vector space; to make sense of the formula $\gamma'(\zeta)\in T_y  \boldsymbol{\mathcal F}_{\hspace{-0.1cm} u_1}$, both involved objects  have to be considered as included in the real tangent space to $\mathbf P^n$ at $y$.}.\sk 

Let 
$\pi:  \mathbf P^n \dashrightarrow P$ be a  generic linear projection 
onto a  generic 
2-plane $P$ passing through  $y$.  Then considering the restriction to $P$ of the $u_i$'s on the one hand, and the composition $\pi\circ \gamma$ on the other hand, it comes that one can assume that $n=2$. Then  for any $i=2,\ldots,d$,  one has 
$T_y  \boldsymbol{\mathcal F}_{\hspace{-0.1cm} u_i}\cap T_y  \boldsymbol{\mathcal F}_{\hspace{-0.1cm} u_1}=0$ hence 
$\gamma$ is transverse  to the 
foliation $\boldsymbol{\mathcal F}_{\hspace{-0.1cm} u_i}$ on a neighbourhood of $y$. This implies that $y$ belongs to the $\boldsymbol{\mathcal F}_{\hspace{-0.1cm} u_i}$-saturation of any open neighbourhood of $\gamma([0,\zeta[)$. According to the 
 definition of $\zeta$, the germ $\Psi_i\in \mathcal O_x$ admits precisely a holomorphic extension, denoted by $\widetilde \Psi_i$, on such an open neighborhood. Since $\widetilde \Psi_i$ is  constant along the leaves of $\boldsymbol{\mathcal F}_{\hspace{-0.1cm} u_i}$, it comes that  it extends analytically on an open neighbourhood of 
 $\gamma([0,\zeta])$. We have proved that $\Psi_2,\ldots,\Psi_{d}$ extend analytically along the restriction of $\gamma$ to $[0,\zeta]$. Using the functional relation $\sum_{i=1}^{d} \widetilde \Psi_i=0$ which is identically satisfied along $[0,\zeta[$, one obtains that the same holds true for $\widetilde \Psi_1$. This contradicts the initial assumption $\zeta<1$ hence gives us that all the $\Psi_i$'s extend analytically along $\gamma$. \sk 
 
 Set $\mathfrak B_i'=u_i(\Sigma(\boldsymbol{\mathcal W}))$ for $i=1,\ldots,d$: it is a finite subset of $\mathbf P^1$ containing (possibly properly) the branch locus $B_i$.  We have proved so far that for any $d$-tuple $(F_i)_{i=1}^d$ as in the statement of the proposition,   each $F_i$ extends holomorphically on $\mathbf P^1\setminus \mathfrak B_i'$ (as a multivalued function denoted by $\widetilde F_i$).  We now prove that $\widetilde F_1$ actually extends holomorphically at any point $b'_1\in \mathfrak B_1'\setminus \mathfrak B_1$, which would imply the proposition.  According to the  definitions of $\mathfrak B_1$ and $\mathfrak B_1'$, there exists an irreducible hypersurface $H_1\subset \mathbf P^n$  such that $u_1(H_1)=b'_1$ and $H_1\not \subset \Sigma^c_i(\boldsymbol{\mathcal W})$. This latter condition implies that $H_1$ is not invariant by $\boldsymbol{\mathcal F}_{\hspace{-0.1cm} u_k}$ for any $k>1$, hence that $u_k \lvert_{H_1}: H_1\dashrightarrow \mathbf P^1$ is dominant. This implies in particular that for $y\in H_1$ generic,  the rational map $u_k$ is defined at $y$ and the point $y_k=u_k(y)$
does not belong to  $ \mathfrak B_k'$. Given such a point $y$, let $\beta:[0,1]\rightarrow \mathbf P^n$ be a smooth injective  path joining $x$ to $y$ such that $\beta([0,1[)$ does not meet $\Sigma(\boldsymbol{\mathcal W})$. Using the very same arguments as above, one gets that all the $d$ germs $F_i\in \mathcal O_x$ extend analytically along $u_i(\beta([0,1[))\subset \mathbf P^1$ and because none of the $y_k$'s belong to $\mathfrak B_k'$ for $k>1$,  these prolongations can be extended to $u_k(\beta([0,1]))$ for any such $k$.  Then, using again the relation $\sum_{i=1}^dF_i(u_i)=0$, one deduces easily that the same holds true for $F_1$ as well, which hence extends holomorphically at $u_1(y)=b'_1$, which terminates the proof.  \end{proof}

\begin{rem}   The preceding proposition is a slight improvement 
of Th\'eor\`eme 1.2.2 of {\rm \cite{PirioThese}} since the branch locus $\mathfrak B_i$ can be a proper subset of 
the complement 
$\mathfrak B_i'$ in $\mathbf P^1$ of the image  of 
$\mathbf P^n\setminus \Sigma(\boldsymbol{\mathcal W})$ by $u_i$, 
 for some (or even for all) indices $i$.  \sk 
 
 For a concrete example, one can consider the Spence-Kummer web (see \S\ref{SS:Spence-Kummer} later on):
 $${\boldsymbol{\mathcal W}}_{\!{\cal S}{\cal K}}=\boldsymbol{\mathcal W}\left(\,  x \, , \, 
y  \, , \, 
 \frac{x}{y}\, , \, 
 \frac{1-x}{1-y} \, , \, 
 \frac{x(1-y)}{y(1-x)}\, , \, 
xy \, , \, 
 -\frac{x(1-y)}{{\,}(1-x)} \, , \, 
- \frac{(1-y)}{y(1-x)}  \, , \, 
  \frac{x(1-y)^2}{y(1-x)^2}
\,  \right)\, .
 $$
  For the rational first integrals above, 
 all the branch loci $\mathfrak B_i$ are equal to $\{ 0,1,\infty\}$ while one has
 $\mathfrak B_j'=\mathfrak B_j=\{ 0,1,\infty\}$ for $j\in J=\{3,6,9\}$ but 
 $\mathfrak B_i'=\{ 0,\pm 1,\infty\}  
 \supsetneq 
 \mathfrak B_i$ for $i\in \{1,\ldots,9\}\setminus J$.
 \end{rem}

The preceding result gives us immediately the
\begin{cor} {1.} Given a path $\gamma$ in $\mathbf P^n\setminus \Sigma(\boldsymbol{\mathcal W})$,  analytic continuation along $\gamma$ induces a linear isomorphism between the  spaces of germs of ARs of 
$\boldsymbol{\mathcal W}$ at the two extremities of $\gamma$.\sk 

{2.} 
Consequently, the spaces of (germs of) ARs of $\boldsymbol{\mathcal W}$ at its regular points organize themselves into a local system on $\mathbf P^n\setminus \Sigma(\boldsymbol{\mathcal W})$. 
\end{cor}

The  branch loci are interesting objects attached to a web 
$\boldsymbol{\mathcal W}=\boldsymbol{\mathcal W}(u_1,\ldots,u_d)$: 
they are rather easy to determine algebraically when the $u_i$'s are given and it follows from  Proposition \ref{P:AnalyticConinuation} that they can be quite useful regarding the question of how to find all the solutions of the associated AFE $\sum_{i=1}^d F_i(u_i)=0$, a fact that will be illustrated in the next paragraph. 
\mk

For any non-constant rational function $r\in \mathbf C(u)$,  the $\boldsymbol{\mathcal W}$-branch loci $\mathfrak B_{v_i}$ of the first integral $v_i=r(u_i)$ of 
$\boldsymbol{\mathcal F}_{\hspace{-0.1cm} u_i}$
  is clearly seen to be $r(\mathfrak B_{u_i})\subset \mathbf P^1$.  So taking suitable 
rational first integrals for   $\boldsymbol{\mathcal W}$, one could deal with a case when all the branch loci $B_i$ have only one point, say the point at infinity $\infty\in \mathbf P^1$.  Actually, this is not particularly relevant if interested by the ARs of $\boldsymbol{\mathcal W}$: indeed, any component $G_i$ of an AR $\sum_{i=1}^d G_i(v_i)$ would be written (locally) $G_i=F_i\circ  r^{-1}$ for a branch of the algebraic function $r^{-1}$ hence would be `more complicated' (we agree that the meaning of this should be made more precise) than the corresponding component $F_i$. \sk 

At the opposite, it is more convenient to deal with the case when the branch loci $\mathfrak B_i$ all are the biggest possible, when all the $u_i$'s have been assumed to be non-composite. First of all, it provides as much information as possible 
 on the ramification points that the components $F_i$ of an AR $\sum_{i=1}^dF_i(u_i)$ can have. And second, for such a first integral $u_i$, the associated $\boldsymbol{\mathcal W}$-branch locus $\mathfrak B_{u_i}$ is canonically attached to the corresponding foliation, up to projective automorphisms: indeed, if ${\tilde{u}}_i$ is another primitive rational first integral, then  
${\tilde{u}}_i =g(u_i)$ with $g\in {\rm PGL}_2(\mathbf C)$ hence $\mathfrak B_{{\tilde{u}}_i}=g(\mathfrak B_{u_i})$.  In this case, the corresponding branch loci will be said to be {\bf primitive}.
\sk 

\bk

For any $i$, we define the  {\bf $\boldsymbol{i}$-th ramification locus $\boldsymbol{\mathfrak B_i^r}$} of $\boldsymbol{\mathcal W}$: 
it is the set of points of $\mathbf P^1$ at which there exists an AR of $\boldsymbol{\mathcal W}$ whose $i$-th component ramifies non-trivially at $\lambda$ (thats is, does not extend holomorphically at $\lambda$).   Clearly,  $\mathfrak B_i^r=\emptyset $ means that all the $i$-th components 
 of the AR of $\boldsymbol{\mathcal W}$ are trivial. \mk 
 
 Proposition \ref{P:AnalyticConinuation} can be stated as the fact that 
the inclusion $\mathfrak B_i^r\subset \mathfrak B_i$ holds true for any $i$.
Note that this inclusion can be proper, {\it cf.}\,the example just below.  
An interesting property of the $\mathfrak B_i^r$'s is that they are invariantly attached to the web considered, contrarily to the $\mathfrak B_i$'s: 
let $\widetilde{\boldsymbol{\mathcal W}}$ be the web equivalent to 
$\boldsymbol{\mathcal W}$, 
defined by the first integrals $\tilde u_i=u_i\circ \varphi$ where $\varphi$ is a given 
 birational transformation of $\mathbf P^n$. We set $\tilde{\mathfrak B}_i=
{\mathfrak B_{\tilde v_i}}$ and $\tilde{\mathfrak B}_i^r=
{\mathfrak B_{\tilde v_i}^r}$ for any $i=1,\ldots,d$. Then  one has 
$\tilde{\mathfrak B}_i^r=\mathfrak B_i^r$  for any $i$,  whereas it can happen that 
$\tilde{\mathfrak B}_j$ and $\mathfrak B_j$ do not coincide for some (possibly for all) $j$.

\mk 

\begin{exm}  
By direct computations, one gets that for the planar 3-web ${\boldsymbol{\mathcal W}}_3$ defined by the first integrals 
$v_1=- {(x-y)}/{(xy)}$, $v_2=x(y-1)/{y}$ and $v_3= y(x-1)/x$, one has $\mathfrak B_1=\{ \pm 1, \infty\}$ and 
$\mathfrak B_i=\{0, \pm 1, \infty\}$ for $i=2,3$. However, this web admits a unique AR, the one associated to the following logarithmic identity 
$ F(v_1)-F(v_2)+F(v_3)=0$ where $F$ stands for the logarithmic function 
$F(u)={\rm Log}(1-u)-{\rm Log}(1+u)$. It follows that $\mathfrak B_i^r=\{\pm 1,\infty\}$ for $i=1,2,3$ hence the latter coincides with the corresponding branch locus $\mathfrak B_i$ only for $i=1$.\sk 

Regarding the invariance of the $\mathfrak B_i$'s and the $\mathfrak B_i^r$'s 
up to birational transformations,  we consider the relations $X=x(y-1)/y$ and $Y=y(x-1)/x$: they  define a birational change of variables and the web 
$\widetilde {\boldsymbol{\mathcal W}}_3$
equivalent to ${\boldsymbol{\mathcal W}}_3$ in the new  system  of coordinates $(X,Y)$ admits $\tilde v_1=(X-Y)/XY$, $\tilde v_2=X$ and $\tilde v_3=Y$ as first integrals.  Contrarily to what happens for ${\boldsymbol{\mathcal W}}_3$, one has 
$\tilde{\mathfrak B}_i^r=\tilde{\mathfrak B}_i=\{\pm 1,\infty\}$ for $i=1,2,3$ 
in what concerns $\widetilde {\boldsymbol{\mathcal W}}_3$. 
\end{exm}

If they are not necessarily invariant up to pre-composition by a birational map, it is not difficult (and left to the reader) to verify that $\mathfrak B_i$ and $\tilde {\mathfrak B}_i$ are related as follows for any $i=1,\ldots,d$: let $\mathcal E xc_{\varphi,i}$ be the set of irreducible divisors $Z\subset \mathbf P^n$ which are contracted by $\varphi$\footnote{That is, such that $\varphi(Z)$ be of codimension at least 2 in the target projective space.} and invariant by $u_i$ (that is such that $u_i(Z)$ be a point of $\mathbf P^1$).  Then 
$\mathfrak B_i$ is included (possibly properly of course) in the union 
of $ \tilde {\mathfrak B}_i $ with the image 
of $\mathcal E xc_{\varphi,i}$ by $u_i$: one has $\mathfrak B_i\subset 
 \tilde {\mathfrak B}_i \cup u_i( \mathcal E xc_{\varphi,i}) \subset \mathbf P^1$. 
 \begin{center}
 $\star$
 \end{center}
 
An important case we will deal with in this text is the one when all the branch loci of $\boldsymbol{\mathcal W}$ are of cardinality less than, or equal to 3. Since these are related to polylogarithms which are special multivalued functions on $\mathbf P^1$ with branch points $0,1$ and $\infty$, we set the 
\begin{defn} 
\label{Def:Polylogarithmic-Ramification}
The web $\boldsymbol{\mathcal W}=\boldsymbol{\mathcal W}(u_1,\ldots,u_d)$ 
has ({\bf primitive}) {\bf polylogarithmic ramification} if for any $i$, the $i$-th (primitive) branch locus of $\boldsymbol{\mathcal W}$  has cardinality equal to or less than 3.  
\end{defn}

Assuming that $\boldsymbol{\mathcal W}=\boldsymbol{\mathcal W}(u_1,\ldots,u_d)$ has polylogarithmic ramification, there are two natural choices 
 for a normal form for  the  $\boldsymbol{\mathcal W}$-branch loci of the $u_i$'s: 
 \begin{equation}
 \label{Eq:Normalisation-Branch-Locus}
 \begin{tabular}{l}
$\bullet$ \, $\mathfrak B_1=\{0, 1,\infty\}$, which is more suitable to deal with classical polylogarithms;  or \mk \\
  $\bullet$ \, $\mathfrak B_{-1}=\{0, -1,\infty\}$, which is more adapted to deal with the components of the abelian \\ ${}$ \hspace{0.2cm} relations of webs constructed from cluster algebras.
 \end{tabular}
 \end{equation}
Of course, one can pass from one of the normalizations to the other if needed. 
For instance, one has $\varphi(\mathfrak B_{-1})=\mathfrak B_{1}$ for $\varphi(u)=1/(1+u)$ or $\varphi(u)=u/(1+u)$. \mk 

To finish this paragraph, we mention that the use of $\mathfrak B_{-1}=\{0, -1,\infty\}$ to deal with polylogarithmic functions, if far less popular than that of  $\mathfrak B_{1}$  used quasi-systematically nowadays, can actually be considered as 
very classical, since according to \cite[\S1.5.3]{Lewin}, the function 
$ - \l {2} ( -x)=
\int_0^x 
{\log(1+u)}du/{u} 
$ was already considered by Spence in the first text \cite{Spence} (dating of 1809) entirely devoted to the study of polylogarithmic functions. 
In Remark \ref{Rem:II-dim>0} further,  we will briefly discuss  the ramification sets of some more general classes of hyper-logarithmic functions.

\subsection{\bf Iterated integrals and abelian relations.}  
\label{Par:II-AR}
Considering the fact that 
low-order polylogarithms  provide some particularly interesting abelian relations for webs defined by rational functions and because these functions are a specific kind of iterated integrals, it is natural to look more closely to the ARs of a web $\boldsymbol{\mathcal W}=\boldsymbol{\mathcal W}(u_1,\ldots,u_d)$ defined by rational functions whose all the components are iterated integrals. 
This subsection is devoted to this. We follow rather closely the presentation given in \cite[\S3]{Pereira} 

\subsubsection{\bf Abelian relations with logarithmic components.}  
\label{SSub:ARs-with-logarithmic-components}
We have defined an AR as the class of an abelian functional relation modulo the  constants. A more concrete and intrinsic way to consider a germ of  AR for $\boldsymbol{\mathcal W}$
at a point  $x\not \in \Sigma(\boldsymbol{\mathcal W})$ 
 is as a $d$-tuple $(\eta_i)_{i=1}^d$ of (germs of) holomorphic 1-forms 
  $\eta_i\in \Omega^1(\mathbf P^1,u_i(x))$ 
 such that $\sum_i u_i^*(\eta_i)=0$. 
  Then,  thanks to  the considerations above, we know that  each component $\eta_i$
 of such an AR extends as a global (but possibly multivalued) 1-form on $\mathbf P^1$, with ramification (if any) at the points of $\mathfrak B_i=\mathfrak B_{u_i}$.\smallskip

 Within this formalism, the reader would probably agree that 
some particularly elementary ARs  of $\boldsymbol{\mathcal W}$ are certainly the logarithmic ones, that is the 
$d$-tuple $(\eta_i)_{i=1}^d$ where  $\eta_i$ is a  global logarithmic 1-form on $\mathbf P^1$ with simple poles at the points of $\mathfrak B_i$, {\it i.e.}  
$${}^{} \qquad \eta_i\in H_i=H^0\Big(\Omega^1_{\mathbf P^1}(\mathfrak B_i)\Big)
\quad \mbox{ for }\quad 
i=1,\ldots,d.$$ 

Then one defines a sub-local system of $\mathcal A(\boldsymbol{\mathcal W})$  by considering the kernel of the map
$$
 \Psi_{\boldsymbol{\mathcal W}}^1  \, :\, \bigoplus_{i=1}^d 
H_i
 \longrightarrow H^0\Big(\mathbf P^n, \Omega^1_{\mathbf P^n}\big( \hspace{-0.05cm} {\rm Log} \hspace{-0.05cm} \Sigma(\boldsymbol{\mathcal W})\Big)\Big)\, , \, \, \,   (\eta_i)_{i=1}^d\longmapsto \sum_{i=1}^d u_i^*(\eta_i)
$$
where  $\Omega^1_{\mathbf P^n}\big( \hspace{-0.04cm}{\rm Log} \hspace{-0.04cm}\Sigma(\boldsymbol{\mathcal W})\Big)$ 
stands for the sheaf of \href{https://en.wikipedia.org/wiki/Logarithmic_form}{logarithmic 1-forms} on $\mathbf P^n$ with  poles along  $\Sigma(\boldsymbol{\mathcal W})$. Then 
$$
\mathcal I
\hspace{-0.08cm}
\mathcal I
\hspace{-0.08cm}
\mathcal A^1_{\boldsymbol{\mathcal W}}={\rm ker}\Big(
\Psi_{\boldsymbol{\mathcal W}}^1
\Big)$$ is a trivial local system on $U=\mathbf P^n\setminus \Sigma\big(\boldsymbol{\mathcal W}\big)$ which clearly embeds into 
$\mathcal A_{\boldsymbol{\mathcal W}}$. \mk

It is  (and it will be) convenient to consider  the space  $H_{\boldsymbol{\mathcal W}}$ 
of logarithmic 1-forms on $\mathbf P^n$ spanned by the pull-backs $u_i(\eta_i)$ for $\eta_i\in H_i$ for each $i$: 
\begin{equation}
\label{Eq:HW}
H_{\boldsymbol{\mathcal W}}=
\sum_{i=1}^d  u_i^*\big( H_i\big)
\subset H^0\Big(\mathbf P^n, \Omega^1_{\mathbf P^n}\big( \hspace{-0.05cm} {\rm Log} \hspace{-0.05cm} \Sigma\big(\boldsymbol{\mathcal W}\big)\Big)\Big)\, .
\end{equation}
Then by elementary linear algebra, 
$\mathcal I
\hspace{-0.08cm}
\mathcal I
\hspace{-0.08cm}
\mathcal A^1_{\boldsymbol{\mathcal W}}$ has rank equal to 
$\sum_{i=1}^d \dim (H_i)-\dim ( H_{\boldsymbol{\mathcal W}})$: by definition, this is the 
{\bf logarithmic rank}  of $\boldsymbol{\mathcal W}$. 

\subsubsection{\bf Reminders about iterated integrals.}  
\label{SSub:Reminder-Iterated-Integrals}
We now recall some classical facts about iterated integrals  in order to introduce an algebraic tensorial formalism to deal with  ARs whose components are iterated integrals of higher weight (in the subsection \S\ref{SSub:ARs-hyperlogarithmic-components} just following).   More precisely, we will only consider here iterated integrals associated to logarithmic rational 1-forms on $\mathbf P^1$, 
functions also known as `hyperlogarithms'. As classical references, we mention 
\cite{Poincare} and  \cite{Lappo}. For more modern ones, the interested reader can consult 
\cite{Wechsung1967}, 
\cite[\S3.6]{Brown} or \cite[\S2.1]{BanksPanzerPym}. 
\mk

We fix a finite subset $\mathfrak B\subset \mathbf P^1$ (hereafter,  it will be one of the branch loci $B_i$) and set $H=H_{\mathfrak B}=H^0\Big(\Omega^1_{\mathbf P^1}(\mathfrak B)\Big)$.
Given a path $\gamma : [0, 1]\rightarrow \mathbf P^1\setminus \mathfrak B$ and a collection of $w\geq 1$ non-zero elements $\eta_1,\ldots,\eta_w$ in $H$, the {\bf iterated integral} 
of   $\eta_1\otimes \cdots \otimes \eta_w$
along $\gamma$ is defined as the complex number 
\begin{equation}
\label{E:II}
\int_{\gamma} \eta_1\otimes \cdots \otimes \eta_w = \int_{\gamma_*(\Delta_w)} p_1^*(\eta_1)
\wedge \cdots \wedge p_w^*(\eta_w)
\end{equation}
where 
$p_s$ denotes the projection $\big(\mathbf P^1\setminus \mathfrak B\Big)^w\rightarrow \mathbf P^1\setminus \mathfrak B$ onto the $s$-th factor and $\gamma_*(\Delta_w)$ stands for the image by $\gamma\times \cdots \times \gamma: [0,1]^w\rightarrow \mathbf P^1\setminus \mathfrak B$ of the standard $w$-dimensional simplex  $\Delta_w\subset [0,1]^w$.
\sk

In the case under scrutiny, it can be proved\footnote{This follows easily from   the fact that we are considering iterated integrals on spaces of dimension 1.} that an iterated integral  \eqref{E:II}
 depends only on the homotopy class of $\gamma$ 
  hence, as soon as we assume that $x=\gamma(0)$ is a fixed base point,  can be viewed as a multivalued holomorphic function of $z=\gamma(1)\in \mathbf P^1\setminus \mathfrak B$.  
 Now we choose an affine coordinate $\zeta$ centered at $x$ on $\mathbf P^1$ (thus $\zeta(x)=0$) and for $z$ sufficiently close to $x$, we denote by $\gamma_z$ the path $[0,1]\ni t\mapsto \zeta^{-1}(t \zeta(z))\in \mathbf P^1$. 
 Then one obtains a well-defined $\mathbf C$-linear map 
 \begin{align*}
 \label{Eq:IIxw}
  {\rm II}^w_{x}=  {\rm II}^w_{\mathfrak B,x} \, : \hspace{0.2cm}
H^{\otimes w} \,  &  \longrightarrow \hspace{0.2cm} \mathcal O_{x} \\ 
\otimes_{i=1}^w \eta_i & \longmapsto \bigg( z \mapsto \int_{\gamma_z}  \eta_1\otimes \cdots \otimes \eta_w \bigg)
\end{align*}
which can be easily seen to be independent of the chosen affine coordinates $\zeta$. By definition, an element of  
$$\mathcal H_{x}^{ w}={\rm Im}\Big({\rm II}^w_{x} \Big)  \subset \mathcal O_{x}
$$ is (the germ at $x$ of) an {\bf iterated integral of weight $\boldsymbol{w}$} (on the projective line, with respect to $\mathfrak B$).  Given a 
 basis $\underline{\nu}=(\nu_1,\ldots,\nu_s)$ of $H$ (thus $s=\lvert \mathfrak B\lvert-1$) and for a word $\boldsymbol{{a}}=a_1\cdots a_w$ with $a_k\in \{1,\ldots,s\}$ for all $k$, we will denote 
 ${\rm II}^w_{x} (\nu_{a_1}\otimes \cdots \otimes \nu_{a_w})$  by  $L_{\boldsymbol{a},x}$. \sk 
 
 When $x$ is assumed to be fixed, no confusion can arise hence 
 we   drop it from all the notation in what follows (we will write   ${\rm II}^w$, $\mathcal H^w$  and $L_{\boldsymbol{a}}$ instead of ${\rm II}^w_x$, 
$\mathcal H_{x}^w$ and $L_{\boldsymbol{a},x}$, etc.)\mk

The iterated integrals satisfy several nice properties. To state these, it is useful to introduce the following objects constructed by taking direct sums: 
\begin{equation*}
H^{\otimes \bullet}=\oplus_{w\in \mathbf N} H ^{\otimes w}\, , \qquad \mathcal H^{\bullet}=\oplus_{w\in \mathbf N} \mathcal H^{ w}
\quad \mbox{ and } \quad 
{\rm II}^\bullet
=\oplus_{w} {\rm II}^w:\,  H^{\otimes \bullet} \rightarrow 
\mathcal H^{ \bullet}\subset \mathcal O_{x}\, , 
\end{equation*}
and since the \href{https://en.wikipedia.org/wiki/Shuffle_algebra}{`shuffle product'} plays an important role regarding the multiplicative structure of the algebra of iterated integrals, we need also to recall what  it is.  
Given a non-empty finite set $S$, let $\mathbf C\langle S \rangle = \oplus_{n\in \mathbf N} \big( \oplus_{ {\boldsymbol{\omega}} \in S^n} \mathbf C\cdot {\boldsymbol{\omega}}\big) $
 be the free non-commutative complex algebra with unity (denoted by $\boldsymbol{1}$) generated by it.  Given two words $\boldsymbol{\sigma}=s_1\cdots s_m$ and ${ \boldsymbol{\sigma}'}=s_1'\cdots s_n'$ on $S$ (that is $s_i, s_j'\in S$ for every $i$ and $j$), one denotes their concatenation by 
 $\boldsymbol{\sigma}\boldsymbol{\sigma}'=s_1\cdots s_m s_1'\cdots s_n'$ and  one sets $\boldsymbol{\sigma}^\nu=s_{\nu(1)}\cdots s_{\nu(m)}$ for any  permutation $\nu\in \mathfrak S_{m}$. By definition, the {\bf shuffle product $\boldsymbol{\shuffle}$} on 
 $\mathbf C\langle S \rangle$,  is the $\mathbf C$-linear commutative and associative product on this  algebra characterized by the fact that  for any two words $\boldsymbol{\sigma}$ and $\boldsymbol{\sigma}'$ as above (of length $m$ and $n$ respectively), one has 
 \begin{equation}
 \label{Eq:shuffle}
 \boldsymbol{\sigma}\shuffle \boldsymbol{\sigma}'=\sum_{\nu \in m\shuffle n} \big(\boldsymbol{\sigma}\boldsymbol{\sigma}'\big)^\nu
\end{equation}
 where $m\shuffle n$ stands for the set of \href{https://fr.wikipedia.org/wiki/(p,_q)-shuffle}{`$(m,n)$-shuffles'}, namely the set 
 of permutations $\nu\in \mathfrak S_{m+n}$ such that $\nu(1)<\cdots <\nu(m)$ and $\nu(m+1)<\cdots <\nu(m+n)$.
 \mk 
 
\label{Page:shuffle}
 

\newpage

\label{Page:Properties-Hyperlogarithms}
Here are the most interesting properties satisfied by iterated integrals: 
\begin{enumerate} 
\item  relatively to the algebra structure on $H^{\otimes \bullet}$ induced by the shuffle product on words in $\{1,\ldots,s\}$,
the map ${\rm II}^\bullet : H^{\otimes \bullet} \rightarrow \mathcal O_{x}$  is a  morphism of complex algebras, {\it i.e.} for any two words  $\boldsymbol{a}$ and $\boldsymbol{a}'$, 
as germs of holomorphic functions at $x$, 
one has 
$$
L_{\boldsymbol{a}}\, L_{\boldsymbol{a}'}= L_{ \boldsymbol{a} \shuffle \boldsymbol{a}'}\, ;$$
\item  the iterated integrals $L_{\boldsymbol{a}}$ for 
all words $\boldsymbol{a}$, are $\mathbf C[z]$-linearly independent. It follows that the map ${\rm II}^\bullet : H^{\otimes \bullet} \rightarrow \mathcal H_{x}^{\bullet}$ is an isomorphism of complex algebras; 
\item in particular,  the  word $\boldsymbol{{a}}$ (hence its length $w=w(\boldsymbol{a})$) is well-defined by 
 $L_{\boldsymbol{a}} \in \mathcal O_{x}$.  By definition,  $\boldsymbol{{a}}$ is the {\bf symbol}  of the latter and $w$ its {\bf weight} (with respect to the  basis $\underline{\nu}$).\sk 
 
More generally, for any iterated integral $L=\sum_{\boldsymbol{a}} \lambda_{\boldsymbol{a}} L_{\boldsymbol{a}}\in \mathcal H^{\otimes\bullet}$ (with 
$\lambda_{\boldsymbol{a}} \in \mathbf C$ non-zero for all but a finite number of words $\boldsymbol{a}$'s), one defines its {\bf weight} $\boldsymbol{w(L)}$ as $\max\{ \,   w(\boldsymbol{a})\, \lvert \, \lambda_{\boldsymbol{a}} \neq 0\,  \}$ and its {\bf symbol} (at $x$) by $\boldsymbol{\mathcal S}(L)=\boldsymbol{\mathcal S}_{{\hspace{-0.05cm}}x}(L)=\sum_{w(\boldsymbol{a})=w(L)} \lambda_{\boldsymbol{a}}\cdot {\boldsymbol{a}}$;
\item  a consequence of {1.} is that $\mathcal H^{\otimes \bullet} $ is a subalgebra of 
$\mathcal O_{x}$ which is graded by the weight: for any $w,w'\in \mathbf N$,  one has 
$$  \mathcal H^{ w} \cdot  \mathcal H^{ w'} \subset   \mathcal H^{w+w'} ;$$
\label{Page:Properties-Hyperlogarithms-II}
\item analytic continuation  along a path  $\gamma$ joining two points $x$ to $y$ in $\mathbf P^1\setminus \mathfrak B$   induces an injective morphism of  algebras  $\mathcal C^\gamma$ between $\mathcal H_{x}^{\bullet} $ and 
$\mathcal H_{y}^{\bullet} $ which is compatible with the filtrations associated to the gradings by the weight. In other terms,  for any weight $w\geq 0$, one has 
 $$ \mathcal C^\gamma\Big(\mathcal H_{x}^{w} \Big)\subset \mathcal H_{y}^{\leq w}=\bigoplus_{w'\leq w} \mathcal H_{y}^{ w'}\, .$$
 Moreover, the   map induced by $\mathcal C^\gamma$ between ${\rm Gr}^w \mathcal H_{x}^{\bullet}=\mathcal H_{x}^{w}$
and ${\rm Gr}^w \mathcal H_{y}^{\bullet}=\mathcal H_{y}^{w}$
is reduced to the identity, {\it i.e.} for any 
  $L_{\boldsymbol{a},x}\in \mathcal H_{x}^w$, one has 
 $$  \mathcal C^\gamma\Big(L_{\boldsymbol{a},x}\Big)-L_{\boldsymbol{a},y} \in 
  \mathcal H_{y}^{<w}=\oplus_{\tilde w<w}\mathcal H_{y}^{\tilde w}\, ; $$
\item from the preceding points, it follows  that \begin{enumerate} 
\item for any $x\not \in \mathfrak B$, iterated integrals elements of $\mathcal H_x$ extend as global (but multivalued) holomorphic function on $\mathbf P^1$, with $\mathfrak B$ as branch locus;
\item  the  symbol and the weight  of an iterated integral do not depend on the base point;
\item  iterated integrals have unipotent monodromy: for any $L\in \mathcal H_{x}^w$ and any loop $\gamma$ in 
$\mathbf P^1\setminus \mathfrak B$ 
centered at $x$, one has 
\begin{equation}
\label{Eq:C-gamma}
 \mathcal C^\gamma\big( L\big)-L\in \mathcal H_{x}^{< w}\, . 
\end{equation}
\end{enumerate}
\end{enumerate}

\label{Pg:Tangential-base-point}
The choice of the additional base-point $x$ outside $\mathfrak B$ is arbitrary in most cases hence may appear as unnatural. Another approach, encountered in many papers, is to take $x$ as one of the points of $\mathfrak B$ ({\it e.g.}\,one takes $x=0$ in the polylogarithmic case $\mathfrak B=\{\, 0,1\infty\,\}$) but this requires in addition to specify  a non-trivial  (real) tangent direction $\zeta$ at $x$: one then talks about of a {\it `tangential base point'}. By iterated integrations 
of 1-forms in $H_{\mathfrak B}$  along smooth paths $\gamma:[0,1]\rightarrow \mathbf P^1$ with $
\gamma(]0,1])\subset \mathbf P^1\setminus \mathfrak B$, $\gamma(0)=x\in \mathfrak B$ and $
\gamma'(0)=\zeta$,  one defines a class of multivalued holomorphic functions $L_{\boldsymbol{a},\zeta}$ on $\mathbf P^1\setminus \mathfrak B$ which enjoys similar versions of all the properties 1.\,to 6.\,listed above. 
However this approach produces genuine  multivalued functions  on the Riemann sphere and 
a determination of any such function must be specified when working locally. And above all, there is the technical problem that a symbol  $\boldsymbol{a}\in H_{\mathfrak {B}}^{\otimes \bullet}$ may be  {\it `divergent'} 
\footnote{Think of
the symbol $\nu_0=du/u$ in case $x=0$ and $\zeta=\partial/\partial z_1$ when $z_1$ stands for the real  part of the variable $z\in \mathbf C$: for any smooth path $\gamma$ emanating from the origin with tangent vector $\zeta$ at this point, the weight 1 integral $\int_\gamma \nu_0$ is divergent, precisely because the integrand has a logarithmic singularity at the origin.} and that in this case,  it is necessary to define $L_{\boldsymbol{a},\zeta}(\gamma)=\int_\gamma \nu_{\boldsymbol{a}}$ by means of a certain 'regularization process'.  All this is quite well-known and mentioned here for the sake of completeness. For more details, we refer 
to \cite{WojtkowiakGalois}
or  to \S2.1 of \cite{BanksPanzerPym} where everything is clearly explained. 
\mk 

In this text, when dealing with iterated integrals we will consider them when a base point outside the ramification locus $\mathfrak B$ has been fixed, essentially because it fits better with the geometric picture given by web geometry and also because their theory is somehow `cleaner'.  
However we will sometimes consider iterated integrals (actually polylogarithms) defined with respect to a tangential base point but mostly to connect with some results already existing in the literature. 
\mk

\begin{rem}
\label{Rem:II-dim>0}
{\bf 1.}{\,}Thanks to some of the preceding points,  when the base-point is fixed and for any weight $w\geq 1$, one can identify the iterated 
integral $L_{\boldsymbol{a}}$, the word $\boldsymbol{a}=\boldsymbol{\mathcal S}(L_{\boldsymbol{a}})$ and  the 
tensor product of logarithmic 1-forms $\nu_{\boldsymbol{a}}=\nu_{a_1}\otimes \cdots  \otimes \nu_{a_w}$. We will often do so in the sequel. 
\mk \\
{\bf 2.}{\,} About the notion of symbol of iterated integrals, one can consult  {\rm \cite{GSVV}}  
or {\rm \cite{DGR}}. Note that the notation  for the symbol  we use (namely $\boldsymbol{\mathcal S}(L_{a_1\cdot a_2\cdots a_n})=a_1a_2\cdots a_n$) is the opposite of the one used in the latter reference.
\mk \\
{\bf 3.}{\,} Actually, a notion of iterated integrals can be considered when working on a complex manifold $M$ of any dimension, relatively to 
any  $k$-tuple $(\omega_s)_{s=1}^k$ of 
 meromorphic 1-forms $\omega_s$ on $M$. Denote by 
 $Z$  the union of the polar sets of the $\omega_i$'s and by $H$ their linear span in $H^0(M,\Omega^1_M(*Z))$. \sk
 
Then for any path $\gamma$ in $M^*=M\setminus Z$ and any $\boldsymbol{\eta}=\eta_1\otimes \cdots \otimes \eta_w\in H^{\otimes w}$, one  defines  a complex value $\int_\gamma \boldsymbol{\eta}$ by means of formula  \eqref{E:II}.  However,  unless when $\boldsymbol{\eta}$ satisfies some {\it integrality conditions} due to K.T. Chen\footnote{These integrability conditions can be explicited but there is no point to do it  here.}, the iterated integral $ \int_\gamma \boldsymbol{\eta}$ does depend on the path $\gamma$ and is not constant on 
its homotopy class (with fixed extremities).  
A simple but interesting case when  
Chen's integrability conditions 
 are automatically satisfied is when all the $\omega_i$'s are closed (which happens for instance when $\dim(M)=1$). In this case, everything that has been said about dimension 1 iterated integrals generalizes  straightforwardly.  
 In particular, for any $x\in M^*$, there is a linear map
 ${\rm II}^w_x: \, H^w \rightarrow \mathcal O_x$, thus one can speak of iterated integrals on $M^*$, of the symbols and weights of such functions, etc.
  \mk \\
{\bf 4.}{\,}  Iterated integrals on the projective line have been considered a long time ago, in full generality by Poincar\'e (see p.\,215 in {\rm \cite{Poincare}}{\rm )} and more systematically by Lappo-Danilevski  who named 
 these functions  {\bf hyperlogarithms} (see p.\,104 of {\rm \cite{Lappo}}{\rm )}. But the specific case  when the ramification locus $\mathfrak B$ has 3 elements was considered even before, for instance by Kummer in 1840 (see his three papers cited in the bibliography). Another classical reference about the polylogarithmic iterated integrals is Nielsen's monograph {\rm \cite{Nielsen}}. 
 What we call `polylogarithmic iterated integrals' were previously called  `Nielsen polylogarithms' in several papers (see   {\rm \cite{CGR}} and some references therein).
 \mk \\
{\bf 5.}{\,}  Several authors have studied Lappo-Danilevski's hyperlogarithms from a more modern perspective, such as  Wechsung {\rm \cite{Wechsung1967}} (see also the book  {\rm \cite{LewinStruct}}{)}.  More recently, in a series of papers (see {\rm \cite{Wojtkowiak1989,Wojtkowiak,WojtkowiakGalois})},  Wojtkowiak  considered and studied  general iterated integrals on $\mathbf P^1$ (in particular, polylogarithms) and more generally on quasi-projective varieties, from the point of view of modern algebraic and arithmetic geometry. 
\mk

{\bf 6.}{\,}  Most of the papers dealing with iterated integrals on $\mathbf P^1$ concern the classical polylogarithmic case (that is when $\lvert \mathfrak B\lvert=3$)  but some 
other specific cases beyond that one have been considered recently. 
Stemming from questions arising in quantum electrodynamics,  
 the case of iterated integrals with ramification locus $\{0,\pm 1,\infty\}$ has been considered starting from the early 70's.  Such iterated integrals are now known as {\it `harmonic polylogarithms'}, a name coined in {\rm \cite{RemiddiVermaseren}}\footnote{This name comes from the fact that the ramification locus $\{\, 0\, ,- 1\, , \, 1 \, , \, \infty\, \}$ is formed by four points in `\href{https://en.wikipedia.org/wiki/Projective_harmonic_conjugate}{harmonic division}' on the projective line.}.   More recently, some authors have also considered the case of hyperlogarithms on $\mathbf P^1$ ramified in 4 points not necessarily in harmonic division (see {\rm \cite{HiroseSato1,HiroseSato2})}.  
\mk 

{\bf 7.}{\,}  
 Regarding concrete functional equations (more precisely, AFEs) satisfied by hyperlogarithms, most of the results until recently were about classical polylogarithms.  A slightly more general class of polylogarithms were considered from this perspective recently in \cite{CGR} where the authors  study (and first and foremost establish) an AFE satisfied by Nielsen polylogarithms whose existence was predicted in a conjecture by Goncharov in the realm of motives.  To our knowledge, at the time of writing there is no published paper in which 
 an interesting AFE satisfied by hyperlogarithms ramified in $p> 3$ points appears. However, we have discovered an interesting equation in the case $p=4$ that we plan to discuss in a future work.
%
\end{rem}
\subsubsection{\bf Abelian relations with hyper-logarithmic components.}  
\label{SSub:ARs-hyperlogarithmic-components}
Let us consider the web  $\boldsymbol{\mathcal W}=\boldsymbol{\mathcal W}(u_1,\ldots,u_d)$ again and let  $w\geq 1$ be an arbitrary weight.  
 Recall the notion introduced in the weight 1 case:  
 $$H_i=H^0\Big(\mathbf P^1,\Omega^1\big ({\rm Log}\mathfrak B_i\big)\Big)
 \qquad \mbox{and}\qquad 
 {H}_{\boldsymbol{\mathcal W}}=\Big\langle u_i^*\big( H_i\big)\, \lvert \,  i=1,\ldots,d\, \Big\rangle\, .$$   
 
\paragraph{}
 Since taking  weight $w$ tensor products of elements of $H_i$ and pull-backs under the $u_i$'s commute 
 (that is 
 $ u_i^*(\eta_1 \otimes \cdots \otimes \eta_w)=
 u_i^*(\eta_1)\otimes \cdots \otimes 
u_i^*(\eta_w)$ for any $\eta_1,\dots,  \eta_w\in H_i$), it comes that  $H_i$ is naturally embedded into  $({H}_{\boldsymbol{\mathcal W}})^{\otimes w}$
for each $i$. 
 Thus one has a well-defined 
  linear map 
\begin{equation}
\label{Eq:PsiWw}
 \Psi_{\boldsymbol{\mathcal W}}^w  \, :\, \bigoplus_{i=1}^d 
\Big(H_i\Big)^{\otimes w}
 \xrightarrow{\hspace*{1cm}} 
 \Big(H_{\boldsymbol{\mathcal W}}\Big)^{\otimes w}\, , \hspace{0.3cm}
 \Big( \boldsymbol{a}_i\Big)_{i=1}^d  
 \xmapsto{\phantom{Hello w}} 
 \sum_{i=1}^d u_i^*( \nu_{\boldsymbol{a}_i})
\end{equation}

Now we fix $x\in \mathbf P^n\setminus \Sigma(\boldsymbol{\mathcal W})$ and set  $x_i=u_i(x)\in \mathbf P^1\setminus \mathfrak B_i$ for every $i$. Since any element of $H_i$ is closed, the same holds true for any element of $u_i^*(H_i)$ hence for any element of 
$H_{\boldsymbol{\mathcal W}}$.  Hence ({\it cf.}\,Remark {\bf \ref{Rem:II-dim>0}.}{\bf 2}) there is an `iterated integration map' 
${\rm II}^w_x: \, H_{\boldsymbol{\mathcal W}}^{\otimes w} \rightarrow \mathcal O_x$ which 
 commutes with any ${\rm II}_{x_i}^w$: 
for any $\boldsymbol{\eta}\in (H_i)^{\otimes w}$, one has $L_{\boldsymbol{\eta},x_i}\circ u_i= L_{u_i^*(\boldsymbol{\eta}),x}$ 
as homolorphic germs at $x$ on $\mathbf P^n$.
\mk 

Consider now $K_{\boldsymbol{\mathcal W}}^w={\rm ker}\big(
\Psi_{\boldsymbol{\mathcal W}}^w\Big)$. 
Then from the previous paragraph, it follows that  $(u_i^*(\boldsymbol{\eta}_i))_{i=1}^d \in \oplus_i H_i^{\otimes w}$ belongs to $K_{\boldsymbol{\mathcal W}}^w$ if and only if ${\rm II}_x^w(\sum_i u_i^*(\boldsymbol{\eta}_i))=0$ in $H_{\boldsymbol{\mathcal W}}^{\otimes w}$, that is if and only if $\sum_{i} L_{\boldsymbol{\eta}_i}\circ u_i$ vanishes (as a germ at $x$).  Thus one gets a linear map 
from $
K_{\boldsymbol{\mathcal W}}^w $ into 
the space $\mathcal A_x({\boldsymbol{\mathcal W}})$ of ARs of $\boldsymbol{\mathcal W}$ at $x$. 
This map being obviously injective, its 
 image 
$$
\mathcal I
\hspace{-0.08cm}
\mathcal I
\hspace{-0.08cm}
\mathcal A^w_x({\boldsymbol{\mathcal W}})=
{\rm Im} \left(K_{\boldsymbol{\mathcal W}}^w 
\xhookrightarrow{\phantom{Hello}}
 \mathcal A_x(\boldsymbol{\mathcal W})\, 
 \right)
\subset   \mathcal A_x(\boldsymbol{\mathcal W})\, 
$$
is a linear subspace of $\mathcal A_x({\boldsymbol{\mathcal W}})$ isomorphic to $K_{\boldsymbol{\mathcal W}}^w$,  which is the stalk at $x$ 
of a (non-trivial) sub-local system of $\mathcal A_{\boldsymbol{\mathcal W}}$, that will be denoted by $\mathcal I
\hspace{-0.08cm}
\mathcal I
\hspace{-0.08cm}
\mathcal A^w_{\boldsymbol{\mathcal W}}$.\mk

\paragraph{}
  In the sequel of the text, we will use the following terminology: 
\begin{itemize}
\item    the adjective {\bf hyperlogarithmic} will be used to characterize anything constructed from iterated integrals on the projective line, with respect to a given/fixed  branch locus $\mathfrak B$;
\item the term {\bf polylogarithmic} will specially refer to the case when 
$\mathfrak B$ has cardinality 3, with the convention that when this set is not specified, then it has to be $\mathfrak B_1=
\{0,1,\infty\}$ or $\mathfrak B_{-1}=
\{0,-1,\infty\}$, with each of these two cases being clear depending on the context ({\it cf.}\,\eqref{Eq:Normalisation-Branch-Locus});\footnote{One verifies easily that  any `{\it $(\mathfrak B_{-1})-$polylogarithm}' (that is, any iterated integral 
constructed from 
the basis 
$\omega_0=du/u$ and $\omega_{-1}=du/(1+u)$ of $H^0(\mathbf P^1, \Omega^1({\rm Log}\, \mathfrak B_{-1}))$) admits a simple expression in terms of 
the classical ({\it i.e.}\,$(\mathfrak B_1)$-)polylogarithms.  
For instance, the standard `{\it $n$-th $(\mathfrak B_{-1})$-polylogarithm}' defined by 
${\bf L}{\rm i}_n^{\mathfrak B_{-1}}(x)=\int_0^x (\omega_0)^{\otimes (n-1)}\otimes \omega_{-1}$ coincides with $-\l {n} (-x)$. 
We note that $(\mathfrak B_{-1})$-polylogarithms have been considered at the very beginning of the study of higher polylogarithms. For instance,  
in his 1809 essay \cite{Spence}, Spence was already  considering  ${\bf L}{\rm i}_n^{\mathfrak B_{-1}}(x)$ (denoted $
\scalebox{0.9}{$\overset{{}^{}\hspace{0.1cm}n}{L}$}(1+x)$  by him, see \cite[p.\,4]{Spence}).} 
\item $K^w_{\boldsymbol{\mathcal W}}$ is the space of {\bf weight $\boldsymbol{w}$ symbolic hyperlogarithmic ARs} of $\boldsymbol{\mathcal W}$; 
\item $K_{\boldsymbol{\mathcal W}}=\oplus_{w\geq 1} 
K^w_{\boldsymbol{\mathcal W}}
$ is the space of {\bf symbolic hyperlogarithmic ARs} of $\boldsymbol{\mathcal W}$; 
\item $
\mathcal I
\hspace{-0.08cm}
\mathcal I
\hspace{-0.08cm}
\mathcal A^w({\boldsymbol{\mathcal W}})$  
(resp.\,$\mathcal I
\hspace{-0.08cm}
\mathcal I
\hspace{-0.08cm}
\mathcal A^w_x({\boldsymbol{\mathcal W}})$)
 is the local system (resp.\,the space of germs at $x$) of {\bf weight $\boldsymbol{w}$ hyperlogarithmic ARs} of $\boldsymbol{\mathcal W}$;
\item $
\mathcal I
\hspace{-0.08cm}
\mathcal I
\hspace{-0.08cm}
\mathcal A({\boldsymbol{\mathcal W}})
=\oplus_{w\geq 1} 
\mathcal I
\hspace{-0.08cm}
\mathcal I
\hspace{-0.08cm}
\mathcal A^w({\boldsymbol{\mathcal W}})
$  
(resp.\,$
\mathcal I
\hspace{-0.08cm}
\mathcal I
\hspace{-0.08cm}
{\mathcal A}_x({\boldsymbol{\mathcal W}})
=\oplus_{w\geq 1} 
\mathcal I
\hspace{-0.08cm}
\mathcal I
\hspace{-0.08cm}
{\mathcal A}_x({\boldsymbol{\mathcal W}})$)
 is the local system (resp.\,the space of germs at $x$) of {\bf hyperlogarithmic  ARs} of $\boldsymbol{\mathcal W}$;
\item ${\rm IIrk}^w(\boldsymbol{\mathcal W})={\rm rk}(\mathcal I
\hspace{-0.08cm}
\mathcal I
\hspace{-0.08cm}
\mathcal A^w_{\boldsymbol{\mathcal W}})=\dim(K^w_{\boldsymbol{\mathcal W}})$ is the {\bf weight $\boldsymbol{w}$   hyperlogarithmic rank} of $\boldsymbol{\mathcal W}$;
\item ${\rm IIrk}(\boldsymbol{\mathcal W}) =\sum_{w \geq 1} {\rm IIrk}^w(\boldsymbol{\mathcal W})=
\dim(K_{\boldsymbol{\mathcal W}})
$ is the {\bf (total)  hyperlogarithmic rank} of $\boldsymbol{\mathcal W}$;

%
%
%
\item in the polylogarithmic case,  
 we will write ${\rm polrk}(\boldsymbol{\mathcal W})$ and ${\rm polrk}^w(\boldsymbol{\mathcal W})$ 
instead of  
${\rm IIrk}(\boldsymbol{\mathcal W})$ and 
${\rm IIrk}^w(\boldsymbol{\mathcal W})$ respectively, and talk about the 
({\bf weight $\boldsymbol{w}$}) {\bf polylogarithmic rank} of $\boldsymbol{\mathcal W}$;
\item[${}^{}$\hspace{-0.9cm}$\bullet$] when $\mathcal I
\hspace{-0.08cm}
\mathcal I
\hspace{-0.08cm}
\mathcal A({\boldsymbol{\mathcal W}})
=\mathcal A({\boldsymbol{\mathcal W}})$,
 the considered web 
${\boldsymbol{\mathcal W}}$ is said to have 
 {\bf only hyperlogarithmic (or polylogarithmic) ARs}. Most of the time, the word `only' will be omitted.
\end{itemize}

\paragraph{}
\label{Par:Illustrate-Symbolic-Method}
To illustrate the symbolic method to find hyperlogarithmic ARs described above, let us consider the following planar 9-web
\begin{align}
\label{Eq:qqqXWA3}
\boldsymbol{{\mathcal W}}
=
\boldsymbol{\mathcal W}\Bigg( 
\, 
x_1\, , \, 
x_2\, , \, 
\frac{1+x_2}{x_1}
\, , \,  & 
\frac{1+x_1}{x_2}
\, , \, 
\frac{(1+x_1)^2}{x_2}
\, , \,  
\frac{1+x_1+x_2}{x_1x_2}
\, ,  \\ & 
\frac{(1+x_1)^2+x_2}{x_1x_2}
\, ,\,  
\frac{1+x_1+x_2}{x_1(1+x_1)}
\, , \, 
\frac{(1+x_1  + x_2)^2}{x_1^2x_2}
\hspace{0.15cm}
\Bigg)\, .
\nonumber
\end{align}
This web  is a particular example of a cluster web (namely, 
it is the `secondary cluster web of Dynkin type $A_3$', see \S\ref{Par:Cluster-varieties-cluster-ensembles} further on and in particular \eqref{Eq:pXWA3}) but this is not relevant for our current purpose  here hence we leave this aspect aside below in this paragraph. 
\sk

The web $\boldsymbol{{\mathcal W}}$ carries hyperlogarithmic ARs of weight 1, 2 and 3. We are going to focus on the case of weight 3 ARs, namely we are going to describe $K_{\boldsymbol{{\mathcal W}}}^3$  by giving an explicit basis of it (the case of $K_{\boldsymbol{{\mathcal W}}}^w$ for $w=1$ or $w=2$ could be treated in the same way). Denote by $u_1,\ldots,u_9$ the rational first integrals in \eqref{Eq:qqqXWA3}, labeled in their order of appearance (thus $u_1=x_1$, $u_2=x_2$, $u_3={1+x_2}/{x_1}$, $\ldots$ and $u_9={(1+x_1  + x_2)^2}/({x_1^2x_2})$).  Then one verifies that 
\begin{itemize}
\item  the divisor of common leaves is cut out (in 
  $\mathbf C^2$) 
by the following irreducible polynomials: $F_1=
x_1$, $ F_2=x_2$, $ 
F_3=1+x_1$, $ F_4=1+x_2$, $
F_5=1+x_1+x_2$ and $ F_6= (1+x_1)^2+x_2
 $. 
\item 
the corresponding ramifications loci of the $u_i$'s all coincide with $\mathfrak B_{-1}=\{ \, 0\, ,\,  -1\, ,  \infty\, \}$. 
\end{itemize}

It follows that for $i=1,\ldots,9$ and $\zeta\in \{0, 1\}$, the logarithmic 1-form $\Omega^i_{\zeta}=u_i^*(du/(u-\zeta))$ can be written as a linear combination (with integer coefficients) in the 1-forms $\eta_s=d{\rm Log}F_s=dF_s/F_s$ for $s=1,\ldots,6$, which form a basis of $H_{\boldsymbol{{\mathcal W}}}$. 
For any $i\in\{1,\ldots,9\}$, the two 1-forms $ \Omega^i_{0}$ and  $ \Omega^i_{1}$ form a basis of the subspace $H_i$ of $
H_{\boldsymbol{{\mathcal W}}}$ they span, hence the 8 tensors 
$$ {}^{} \qquad \Omega^i_{\zeta_1\zeta_2\zeta_3}= \Omega^i_{\zeta_1}
\otimes \Omega^i_{\zeta_2}
\otimes \Omega^i_{\zeta_3} \quad \mbox{ with }\quad  \zeta_1,\zeta_2,\zeta_3\in \{0, 1\}
$$
 form a basis of 
$(H_i)^{\otimes 3}$.  Elements 
of $K_{\boldsymbol{{\mathcal W}}}^3$ corresponds to 9-tuples $(\Omega^i)_{i=1}^9 \in \prod_{i=1}^9  (H_i)^{\otimes 3} $ summing up to zero, that is such that 
$ \sum_{i=1}^9 \Omega^i=0$ 
 in $(
H_{\boldsymbol{{\mathcal W}}})^{\otimes 3}$. Elementary (but lenghty)  linear algebra computations  in this vector space (which is of dimension $6^3=216$) lead us to the conclusion that $K_{\boldsymbol{{\mathcal W}}}^3$ is 2-dimensional and admits the following two 9-tuples as a basis
\begin{align*}
\Big( 
2 \,\Omega_{011-110} \, , \,  
\Omega_{001-100}\, , \,
2 \, \Omega_{011-110} \, , \,
-2\,  \Omega_{001-100}&\, , \, \\
 \Omega_{001-100}
\, , \,
-2\, \Omega_{001-100}
&\, , \,
-2\,\Omega_{001-100}\, , \,
2 \,\Omega_{011-110}
\, , \,
\Omega_{001-100}
\, \Big) \\
{}^{}\, \mbox{and} \quad 
 \Big( 
 2\, \Omega_{011 -101}
 \, , \, 
\Omega_{010-100} 
\, , \, 
2\, \Omega_{011 -101}
\, , \, 
-2\, \Omega_{010-100} 
 \, , \,  &\\
\Omega_{010-100} 
 \, , \, 
-2\, \Omega_{010-100} 
 & \, , \, 
-2\, \Omega_{010-100} 
   \, , \, 
 2\, \Omega_{011 -101}
    \, , \, 
\Omega_{010-100} 
\Big)\, , 
\end{align*}
where we have used the following notation: for any length 3 words $w,w'$ in  the two letters $0$ and $1$, 
$\Omega_{w -w'}$  is written for $\Omega^i_w-\Omega^i_{w'}$ when it appears at the  $i$-th place of the 9-tuple.



\paragraph{}
\label{Par:Generation-ARs-Monodromy}
From the material above, one can first deduce an interesting  general property of ARs with   iterated integrals as components.  Given $x\in 
\mathbf P^n\setminus  \Sigma(\boldsymbol{\mathcal W})$, 
 it follows from Proposition \ref{P:AnalyticConinuation} that analytic continuation of ARs along loops centered at $x$ induces a representation 
$$ \mu_{\boldsymbol{\mathcal W}}: \pi_1\Big(\mathbf P^n\setminus  \Sigma(\boldsymbol{\mathcal W}),x\Big)\rightarrow \mathcal A_x({\boldsymbol{\mathcal W}})\, .
$$

From what has been recalled above about the monodromy of iterated integrals, one deduces easily that  $\mu_{\boldsymbol{\mathcal W}}-{\rm Id}
 $ 
sends $
\mathcal I
\hspace{-0.08cm}
\mathcal I
\hspace{-0.08cm}
\mathcal A^w_x({\boldsymbol{\mathcal W}})$   into 
$
\mathcal I
\hspace{-0.08cm}
\mathcal I
\hspace{-0.08cm}
\mathcal A^{<w}_x({\boldsymbol{\mathcal W}})=\oplus_{\omega<w}
=
\mathcal I
\hspace{-0.08cm}
\mathcal I
\hspace{-0.08cm}
\mathcal A_x^\omega({\boldsymbol{\mathcal W}})$ for any $w>1$. This shows in particular that 
 the restriction of $\mu_{\boldsymbol{\mathcal W}}$ to 
 $
\mathcal I
\hspace{-0.08cm}
\mathcal I
\hspace{-0.08cm}
\mathcal A_x({\boldsymbol{\mathcal W}})$ is unipotent.  
An interesting consequence of this is that it is possible to construct  new hyperlogarithmic ARs starting from any given element of $\mathcal I
\hspace{-0.08cm}
\mathcal I
\hspace{-0.08cm}
\mathcal A_x^w({\boldsymbol{\mathcal W}})$, which is particularly useful when investigating webs from the point of view of their rank. 

\paragraph{}
The interest of looking specifically at iterated integrals ARs of webs defined by rational first integrals is clear: first many of the known ARs of such webs are of this kind; second, the determination of ARs of this type is equivalent to that of the vector space $K_{\boldsymbol{\mathcal W}}$, hence can be achieved using standard technics of linear algebra.   In order to illustrate the method, an explicit example is discussed in detail  in \S\ref{Par:YW-A3-is-AMP} just below. \sk 

As far as we know, the first occurrence of this symbolic method to solve an AFE is in \cite[Prop.\,(4.5)]{HainMacPherson}. The first systematic uses of this approach  for hunting ARs  in web geometry 
can be found  in the unpublished preprint \cite{Robert} then in the paper \cite{Pereira} (see also \S\ref{SSub:Examples-AMP-Webs} below) where the author determines the spaces of  ARs of the webs 
${\boldsymbol{\mathcal W}}_{\hspace{-0.1cm}\mathcal M_{0,n+3}}$, for $n\geq 2$. 
\mk 

How powerful the symbolic method is  can be illustrated by considering the following  elementary fact, which allows to get a priori many distinct 
hyperlogarithmic ARs for ${\boldsymbol{\mathcal W}}$ from any non-trivial 
element of $
\mathcal I
\hspace{-0.08cm}
\mathcal I
\hspace{-0.08cm}
\mathcal A_x({\boldsymbol{\mathcal W}})$, and which we will use  to prove Theorem \ref{T:FE-holomorphic} later on. \sk 

We continue to use the notation introduced above relatively to ${\boldsymbol{\mathcal W}}$. As is well-known, for any positive weight $w$, there is a natural linear action $\circ$ of the permutation  group $\mathfrak S_w$ on $(H_{{\boldsymbol{\mathcal W}}})^{\otimes w}$, characterized by the fact that, 
for any $\boldsymbol{v}=v_1\otimes \cdots \otimes v_w\in (H_{{\boldsymbol{\mathcal W}}})^{\otimes w}$ and any 
 $\sigma\in \mathfrak S_w$, one has: 
$${\boldsymbol{v}}^\sigma= \sigma\circ {\boldsymbol{v}}
=v_{\sigma^{-1}(1)}\otimes \cdots \otimes v_{\sigma^{-1}(w)}\, .$$

The  point which is relevant for our purpose is that this action commutes with pull-backs. 
More precisely, let $\sigma$ be a fixed permutation and assume that $\boldsymbol{\eta}_i
=\eta_{i,1}\otimes \cdots \otimes \eta_{i,w}\in (H_i)^{\otimes w}$ for some $i\in \{1,\ldots,d\}$. Then $u_i^*(\boldsymbol{\eta}_i)=u_i^*(\eta_{i,1})\otimes \cdots \otimes u_i^*(\eta_{i,w})$ belongs to $(H_{{\boldsymbol{\mathcal W}}})^{\otimes w}$ and one verifies that 
$$
u_i^*\Big(\boldsymbol{\eta}_i^\sigma\Big)=
u_i^*(\boldsymbol{\eta}_i)^\sigma\, .
$$
From this,  one immediately gets  the following 
\begin{lem} 
\label{L:KWw-Sw-equivariant}
The linear map \eqref{Eq:PsiWw} is $\mathfrak S_w$-equivariant and consequently, the space $K^w_{{\boldsymbol{\mathcal W}}}$ 
of weight $w$ symbolic hyperlogarithmic ARs of the web 
${{\boldsymbol{\mathcal W}}}$
is stable under the action of  
 $\mathfrak S_w$. \sk 

More concretely, for any  $\sigma\in \mathfrak S_w$ and  any $d$-tuple $(\boldsymbol{\eta}^i)_{i=1}^d\in \oplus_{i=1}^d (H_i)^{\otimes n}$, one has (in $(H_{{\boldsymbol{\mathcal W}}})^{\otimes w}$):
$$
\sum_{i=1}^d  u_i^*(\boldsymbol{\eta}_i)=0 \quad \Longleftrightarrow 
\quad 
\sum_{i=1}^d  u_i^*\Big(\boldsymbol{\eta}_i^\sigma\Big)=0\, .
$$
\end{lem}

This result, which is basically an elementary remark, allows to construct new hyperlogarithmic ARs starting from one, which is again  particularly useful when investigating webs from the point of view of their rank.  Below in this text, we will use it to get an easy proof of a (holomorphic version of a) classical result about functional equations of classical polylogarithms (see Theorem \ref{T:FE-holomorphic}). 
 However, note that in contrast with the generation of new hyperlogarithmic ARs by means of analytic continuation mentioned above in \S\ref{Par:Generation-ARs-Monodromy}, the ARs produced from a given one in $\mathcal I
\hspace{-0.08cm}
\mathcal I
\hspace{-0.08cm}
\mathcal A^w({\boldsymbol{\mathcal W}})$ by using the preceding lemma are of the same weight: for any $w\geq 2$, the symmetric group  $\mathfrak S_w$ acts on $\mathcal I
\hspace{-0.08cm}
\mathcal I
\hspace{-0.08cm}
\mathcal A^w({\boldsymbol{\mathcal W}})$. 
\sk 

It should be noted that Lemma 
\ref{L:KWw-Sw-equivariant} 
may well be quite far-reaching because it might well open the door to an approach to the hyperlogarithmic ARs of webs via methods derived from the representation theory of symmetric groups (see \S\ref{SS:II-Schur-Monodromy} further on for more in this direction).
\begin{center}
$\star$
\end{center}

\paragraph{}
\label{poupou}
Finally, we mention that in order to deal with the ARs with iterated integrals components of the rational webs considered later on in this text, we have implemented the symbolic method 
 in the computer algebra system Maple.\footnote{Any reader interested in these routines and who would like to use them to study webs is invited to contact the author.} We now have at disposal a series of Maple routines which, given 
 a rational web 
 ${\boldsymbol{\mathcal W}}(u_1,\ldots,u_d)$: 
 \begin{enumerate}
\item[$-$] first compute the ${\boldsymbol{\mathcal W}}$-branch loci $\mathfrak B_i$ of the $u_i$'s; \vspace{-0.25cm}
\item[$-$] then construct a basis $\boldsymbol{\beta}$ of the space $H_{{\boldsymbol{\mathcal W}}}$ defined in \eqref{Eq:HW}; \vspace{-0.25cm}
\item[$-$] next, given an element $u_i^*(\eta)\in H_{{\boldsymbol{\mathcal W}}}$ with $\eta\in H^0(\mathbf P^1,\Omega^1({\rm Log} \mathfrak B_i))$, give its decomposition in the basis $\boldsymbol{\beta}$;\vspace{-0.25cm}
\item[$-$]  next, construct the matrix corresponding to the map $\Psi^w_{{\boldsymbol{\mathcal W}}}$ defined in \eqref{Eq:PsiWw}; \vspace{-0.25cm}
\item[$-$] eventually, compute the dimension of $K_{\boldsymbol{\mathcal W}}^w$ (and possibly provides a basis for it if required).
\end{enumerate}

These routines constitute an efficient package to study ARs of rational webs. 
A great amount of the results concerning the polylogarithmic rank of the webs considered further on in the text have been established using these routines.  However we believe it is fair to say that this computational approach  actually has a rather limited scope since  it does not seem to allow to find  iterated integrals ARs of weight bigger than 4 within a reasonable time.  The reason is that such ARs have to be looked for as carried by (for instance cluster) $d$-webs with $d$ relatively big and so with spaces $H_{\boldsymbol{\mathcal W}}$ of  dimension of the order of several dozens or even more. Then starting from $w=4$ (and even for $w=3$ when $d$ is already big enough), the tensor space $(H_{\boldsymbol{\mathcal W}})^{\otimes w}$ is of quite high dimension which makes  performing linear algebra  in it time consuming hence out of reach of computational investigations.




\subsubsection{Some examples of AMP webs with only polylogarithmic ARs.}  
\label{SSub:Examples-AMP-Webs}
In this subsection, we give examples of AMP webs with only polylogarithmic abelian relations.  The first are the webs $\boldsymbol{\mathcal W}_{\mathcal M_{0,n+3}}$ for $n\geq 2$.  We recall the approach and the results of \cite{Pereira} about them. 
 Then we deal first with the so-called `{\it $\boldsymbol{\mathcal Y}$-cluster web of type $A_3$}' before considering examples of  webs defined on configuration spaces of points in the projective plane.

\paragraph{The webs $\boldsymbol{{\mathcal W}_{\mathcal M_{0,n+3}}}$ all are AMP.}  
\label{Par:W-M0n+3-all-AMP}
We recall some results of Pereira about the polylogarithmic ARs of the web  $\boldsymbol{\mathcal W}_{\mathcal M_{0,n+3}}$ for any fixed $n\geq 2$.  
  In \cite[\S4]{Pereira}, he deduces from general arguments 
  the two following lower bounds on the first two polylogarithmic ranks of  $\boldsymbol{\mathcal W}_{\mathcal M_{0,n+3}}$ (with all ramification loci equal to 
  $\mathfrak B_1=\{0,1,\infty\}$): 
\begin{align*}
{\rm polrk}^1 \Big( \boldsymbol{\mathcal W}_{\mathcal M_{0,n+3}}\Big)\geq  &\, 2 \, { n+3\choose 4}-h^1({\mathcal M_{0,n+3}}) \\
\mbox{ and }\quad 
{\rm polrk}^2 \Big( \boldsymbol{\mathcal W}_{\mathcal M_{0,n+3}}\Big)\geq  &\, 2^2\, { n+3\choose 4}-h^1\big({\mathcal M_{0,n+3}}\big)^2+h^2\big({\mathcal M_{0,n+3}}\big)\, . 
\end{align*}

On the other hand, thanks to a celebrated result of Arnold , one has 
$$
h^1\big({\mathcal M_{0,n+3}}\big)
= { n \choose 2}+2n\qquad \mbox{ and } \qquad 
h^2\big({\mathcal M_{0,n+3}}\big)=\frac{n(n-1)}{24}\big(3n^2+17n+26\big)
$$
from which one gets that the two following minorations hold true: 
$$
{\rm polrk}^1 \Big( \boldsymbol{\mathcal W}_{\mathcal M_{0,n+3}}\Big)
\geq  2{ n+3\choose 4 }-{ n\choose 2  }-2n\quad \mbox{ and } \quad  
{\rm polrk}^2 \Big( \boldsymbol{\mathcal W}_{\mathcal M_{0,n+3}}\Big)
\geq  { n+3\choose 4 }-{ n+2\choose 3  }={ n+2\choose 4  }\, .
$$

Combined with  \eqref{E:Pereira-r(W)}, we deduce from the previous bounds that 

$$3\,{ n+3\choose 4 }-{ n+2\choose 3  }-{ n\choose 2  }-2n\leq 
{\rm polrk} \Big( \boldsymbol{\mathcal W}_{\mathcal M_{0,n+3}}\Big)
\leq \rho\Big({\boldsymbol{\mathcal W}}_{ \hspace{-0.06cm}\mathcal M_{0,n+3}}\Big)=3\,{ n+3\choose n-1} 
-{  n+2  \choose n-1} 
- {  n+1  \choose n-1}-  {  n  \choose n-1}
$$
and since the arithmetic quantities at both extremities are equal, one obtains 
$$
{\rm polrk} \Big( \boldsymbol{\mathcal W}_{\mathcal M_{0,n+3}}\Big)=
{\rm rk} \Big( \boldsymbol{\mathcal W}_{\mathcal M_{0,n+3}}\Big)= \rho\Big({\boldsymbol{\mathcal W}}_{ \hspace{-0.06cm}\mathcal M_{0,n+3}}\Big)=3{ n+3\choose 4} 
-{  n+2  \choose 3} 
- {  n+1  \choose 2}- n \, .
$$
In particular, one gets  that $\boldsymbol{\mathcal W}_{\mathcal M_{0,n+3}}$ is AMP with only polylogarithmic ARs, of weigth 1 or 2. 

%
%


\paragraph{The $\boldsymbol{\mathcal Y}$-cluster web of type $\boldsymbol{A_3}$ is AMP.}  
\label{Par:YW-A3-is-AMP}

Further on we will give a  general definition of the notion of cluster web and in particular that of the $\boldsymbol{\mathcal Y}$-cluster web associated to a Dynkin diagram $\Delta$. For $\Delta=A_2$ one recovers (a model of) Bol's web whose all properties are well-known. 
We find it more interesting to consider the case when $\Delta=A_3$: the corresponding web has not been considered before and it is not too complicated to explain how to get all its ARs by using the algebraic approach described above in \S\ref{SSub:ARs-hyperlogarithmic-components}.
\mk 

We will explain later on how its first integrals are generated, but for the moment let us simply define the  `$\boldsymbol{\mathcal Y}$-cluster web of type $\boldsymbol{A_3}$', as the following 9-web in three variables $x_1,x_2,x_3$: 
\begin{align*}
\boldsymbol{\mathcal Y \hspace{-0.05cm}\mathcal W}_{\hspace{-0.05cm} A_3}= & \, 
\boldsymbol{\mathcal W}\Bigg(\, 
{x_1} \, , \, 
{x_2} \, , \, 
{x_3} \, , \, 
\frac{1+x_2}{x_1} \, , \, 
\frac{(1+x_1)(1+x_3)}{x_2} \, , \, 
\frac{1+x_2}{x_3} \, , \, 
 \\
& \hspace{1.5cm}
\frac{1+x_1+x_2+x_3+x_1x_3}{x_1x_2} \, ,
\, 
\frac{1+x_1+x_2+x_3+x_1x_3}{x_2x_3} \, , \, 
\frac{(1+x_1+x_2)(1+x_2+x_3)}{x_1x_2x_3} \,  \Bigg)\, . 
\end{align*}
We denote by $u_1,u_2,\ldots,u_9$ the rational first integrals appearing in this definition of 
$\boldsymbol{\mathcal Y \hspace{-0.05cm}\mathcal W}_{\hspace{-0.05cm} A_3}$, labeled in the corresponding order: one has 
$$u_1=x_1\, ,\quad u_2=x_2\, ,\quad u_3=x_3\, , 
\quad u_4=\frac{1+x_2}{x_1} \, , \, 
 \ldots  \quad \mbox{and }\quad u_9=\frac{(1+x_1+x_2)(1+x_2+x_3)}{x_1x_2x_3}\, .$$ 

This web has intrinsic dimension 3 and by direct computations, one verifies that:
\begin{itemize}
\item[{\bf 1.}] seen as a divisor in $\mathbf C^3$, its divisor of common leaves 
 has nine irreducible components; more precisely,  one has 
 $\Sigma^c( \boldsymbol{\mathcal Y \hspace{-0.05cm}\mathcal W}_{\hspace{-0.05cm} A_3}) =\cup_{s=1}^9 V_s$ with $V_s=(F_s=0)\subset \mathbf C^3$, 
 where $F_1,\ldots,F_9$ stand for the following polynomials: 
\begin{align*}
F_{1}=& \, x_1 && F_{4}= 1+x_1 &&  F_7=1+x_1+x_2
\\
F_{2}= &\, x_2 && F_5=1+x_2 && F_8=1+x_2+x_3
\\
F_{3}= & \, x_3  && F_6=1+x_3 && F_9=1+x_1+x_2+x_3+x_1x_3\, .
\end{align*}
\item[{\bf 2.}] for any $i=1,\ldots,9$, the $i$-th ramification locus  of 
$\boldsymbol{\mathcal Y \hspace{-0.05cm}\mathcal W}_{\hspace{-0.05cm} A_3}$
 is included in $\mathfrak B_{-1}$, {\it i.e.}\,one has  
 $$\mathfrak B_i=\mathbf P^1\setminus u_i\big(\Sigma^c( \boldsymbol{\mathcal Y \hspace{-0.05cm}\mathcal W}_{\hspace{-0.05cm} A_3})\big)=\mathfrak B_{-1}=\{\, -1\, ,\, 0\, ,\, \infty\, \}\, .$$
\end{itemize}

From these two points, one deduces that for any $i\in \{ 1,\ldots,9\}$ and any $\zeta\in \mathfrak B_{-1}=\{\, -1, 0 , \infty\, \}$, the fiber $u_i^{-1}(\zeta)$, seen as a divisor in $\mathbf C^3$ is an effective linear combination of the $V_s$ for $s=1,\ldots,9$. All the corresponding decompositions 
\begin{equation}
\label{Eq:Decomp}
{}^{}
\hspace{1cm}
u_i^{-1}(\zeta)=\sum_{s=1}^9 \lambda_{i,s}^\zeta\, V_s\qquad \mbox{\big(with }\, \lambda_{i,s}^\zeta\in \mathbf N\, \mbox{for every } i,s  \mbox{ and } \zeta\mbox{\big)}
\end{equation}
 can be entirely characterized by using logarithmic differential 1-forms as follows:  one sets  
\begin{align*}
\eta_s= & \, d{\rm Log}(F_s)=dF_s/F_s \quad \mbox{ as well as }\\
\Omega_i^0=&\, d{\rm Log}\,(u_i)=du_i/u_i\hspace{1cm}  \mbox{and} \qquad \Omega_i^{-1}=d{\rm Log}\,(1+u_i)=du_i/(1+u_i)\,
\end{align*}
for $i,s=1,\ldots,9$. Since  $H^0\big(\mathbf P^1,\Omega^1({\rm Log}\, \mathfrak B_{-1})\big)=\big\langle \, \omega^0\, , \, \omega^1\,\big\rangle$ with $\omega^\zeta=d{\rm Log}(\zeta+u)=du/(\zeta+u)$ for $\zeta\in \{-1,0\}$, it follows that for every  $i=1,\ldots,9$, the pair  $\big(\Omega_i^0,\Omega_i^1\big)$ form a basis of the subspace 
$$H_i=u_i^*\Big(H^0\big(\mathbf P^1,\Omega^1({\rm Log}\, \mathfrak B_{-1})\big)\Big)\qquad  \mbox{ of }\qquad  
H_{\boldsymbol{\mathcal Y \hspace{-0.05cm}\mathcal W}_{\hspace{-0.05cm} A_3}}
=H^0\Big(\mathbf C^3,\Omega^1({\rm Log}\,\Sigma^c\big( \boldsymbol{\mathcal Y \hspace{-0.05cm}\mathcal W}_{\hspace{-0.05cm} A_3}) \big)\Big)\, ,$$ the latter vector space admitting the 9-tuple $\boldsymbol{\eta}=(\eta_1,\ldots,\eta_9)$ as a basis.  In terms of the $\eta_s$ and the $\Omega_i^\zeta$'s, any decomposition \eqref{Eq:Decomp} is equivalent to $\Omega_i ^\zeta=\sum_{s=1}^9 \lambda_{i,s}^\zeta\,\eta_s$ and by direct (easy but a bit lengthy) computations, one obtains  
that  the expressions of the $\Omega_i^\zeta$'s in the basis $\boldsymbol{\eta}$ are as follows: 
\begin{align}
\label{aligno}
\Omega_1^0=&  \,  \eta_1     &&\Omega_1^{-1} = \eta_4\, ; \nonumber
\\ 
\Omega_2^0=& \, \eta_2     &&\Omega_2^{-1} = \eta_5\, ;  \nonumber
\\ 
\Omega_3^0=& \, \eta_3     &&\Omega_3^{-1} = \eta_6\, ; \nonumber
\\ 
\Omega_4^0=& \, \eta_5-\eta_1     &&\Omega_4^{-1} = \eta_7-\eta_1\, ;  \nonumber
\\ 
\Omega_5^0=&\,  \eta_4+\eta_6-\eta_2     &&\Omega_5^{-1} = \eta_9-\eta_2\, ;
\\ 
\Omega_6^0=& \, \eta_5-\eta_3     &&\Omega_6^{-1} = \eta_8-\eta_3\, ; \nonumber
\\ 
\Omega_7^0=&\,  \eta_9-\eta_2-\eta_1     &&\Omega_7^{-1} = \eta_8+\eta_4-\eta_1-\eta_2\, ; \nonumber
\\ 
\Omega_8^0=& \, \eta_9-\eta_2-\eta_3     &&\Omega_8^{-1} = \eta_6+\eta_7-\eta_2-\eta_3\, ; \nonumber
\\ 
\mbox{and }\qquad \Omega_9^0=& \, \eta_7+\eta_8-\eta_1-\eta_3     &&\Omega_9^{-1} = \eta_5+\eta_9-\eta_1-\eta_3\, .  \nonumber
\end{align}


It is then easy to determine the space of weight $w$ symbolic polylogarithmic ARs of ${\boldsymbol{\mathcal Y \hspace{-0.05cm}\mathcal W}_{\hspace{-0.05cm} A_3}}$
 for any $w\geq 1$. The determination of the logarithmic ARs amounts to finding a basis of the space of pairs of 9-tuples $\big((c^i_0)_{i=1}^9,(c^i_1)_{i=1}^9\big)\in \big(\mathbf C^9\big)^2$ such that $\sum_{i=1}^9 \big( c^i_0 \Omega_i^0+c^i_1 \,\Omega_i^{-1})=0$ (equality as a  logarithmic 1-form on $\mathbf C^3$). Using the formulas \eqref{aligno}, it is easy to decompose the RHS with respect to the basis $\boldsymbol{\eta}$: it gives us a linear system in the $c_i^\zeta$'s which is not difficult to solve. One obtains that 
$
\mathcal I
\hspace{-0.08cm}
\mathcal I
\hspace{-0.08cm}
\mathcal A^1\big({\boldsymbol{\mathcal Y \hspace{-0.05cm}\mathcal W}_{\hspace{-0.05cm} A_3}}\big)$ is 9-dimensional, in other terms : 
${\rm polrk}^1\big( 
{\boldsymbol{\mathcal Y \hspace{-0.05cm}\mathcal W}_{\hspace{-0.05cm} A_3}}\big)=9$. 
\sk 

The determination of weight 2 polylogarithmic ARs is  not much more complicated: it amounts to finding a basis of the space of scalars 
$c^i_{\zeta\zeta'}$ (with $i=1,\ldots,9$ and $\zeta,\zeta'\in \{0,-1\}$) such that 
\begin{equation}
\label{lololo}
\sum_{i=1}^9\sum_{\zeta,\zeta'\in \{-1,0\}}  c^i_{\zeta,\zeta'} \,\Omega_i^\zeta \otimes \Omega_i^{\zeta'}=0\, , 
\end{equation}
this equality holding in the double tensor product 
$H_{\boldsymbol{\mathcal Y \hspace{-0.05cm}\mathcal W}_{\hspace{-0.05cm} A_3}}\otimes 
H_{\boldsymbol{\mathcal Y \hspace{-0.05cm}\mathcal W}_{\hspace{-0.05cm} A_3}}$.
 Using the formulas \eqref{aligno} it is not difficult to express each term $\Omega_i^\zeta \otimes \Omega_i^{\zeta'}$ as a linear combination of the $\eta_i\otimes \eta_j$'s for $i,j=1,\ldots,9$. These tensors  forming a basis of $\big(H_{\boldsymbol{\mathcal Y \hspace{-0.05cm}\mathcal W}_{\hspace{-0.05cm} A_3}}\big)^{\otimes 2}$,   
determining the ARs with iterated integral components comes down to  solving an explicit linear system in the $c^i_{\zeta\zeta'}$'s.  Solving it gives us that any tensorial identity \eqref{lololo} is a multiple of the following one
  \begin{equation}
  \label{Eq:Lalolu}
\sum_{i=1}^9  \Big( \Omega_i^0 \otimes \Omega_i^{-1}- \Omega_i^{-1} \otimes \Omega_i^{0}\Big)=0\, , 
\end{equation}
from which it follows that $
\mathcal I
\hspace{-0.08cm}
\mathcal I
\hspace{-0.08cm}
\mathcal A^2\big({\boldsymbol{\mathcal Y \hspace{-0.05cm}\mathcal W}_{\hspace{-0.05cm} A_3}}\big)$ has dimensional one, {\it i.e.} 
${\rm polrk}^2\big( 
{\boldsymbol{\mathcal Y \hspace{-0.05cm}\mathcal W}_{\hspace{-0.05cm} A_3}}\big)=1$. 
 \mk 
 
 From the previous considerations, it comes that 
 $${\rm polrk}\big( 
{\boldsymbol{\mathcal Y \hspace{-0.05cm}\mathcal W}_{\hspace{-0.05cm} A_3}}\big)
\geq {\rm polrk}^1\big( 
{\boldsymbol{\mathcal Y \hspace{-0.05cm}\mathcal W}_{\hspace{-0.05cm} A_3}}\big)
+{\rm polrk}^2\big( 
{\boldsymbol{\mathcal Y \hspace{-0.05cm}\mathcal W}_{\hspace{-0.05cm} A_3}}\big)=9+1=10\, . $$
 On the other hand, direct computations give that 
$$ 
{\rho}^\bullet\big( 
{\boldsymbol{\mathcal Y \hspace{-0.05cm}\mathcal W}_{\hspace{-0.05cm} A_3}}\big)=\big( 
6,3,1
\big) \,. 
$$

We thus have $10= {\rm polrk}\big( 
{\boldsymbol{\mathcal Y \hspace{-0.05cm}\mathcal W}_{\hspace{-0.05cm} A_3}}\big)
\leq {\rm rk}\big( 
{\boldsymbol{\mathcal Y \hspace{-0.05cm}\mathcal W}_{\hspace{-0.05cm} A_3}}\big)
\leq 
{\rho}\big( 
{\boldsymbol{\mathcal Y \hspace{-0.05cm}\mathcal W}_{\hspace{-0.05cm} A_3}}\big)
=6+3+1=10$, hence all these inequalities actually are equalities.   In particular, we get that 
${\boldsymbol{\mathcal Y \hspace{-0.05cm}\mathcal W}_{\hspace{-0.05cm} A_3}}$ is AMP, with only logarithmic ARs, plus the extra polylogarithmic AR of weight 2 corresponding to the tensorial identity \eqref{Eq:Lalolu}.

\paragraph{The web $\boldsymbol{{\mathcal W}_{\hspace{-0.0cm}{\bf Conf}_{6}(\mathbf P^2)}}$ is AMP.}  
\label{Par:W-Conf6(IP2)-is-AMP}
We consider the web defined in \S\ref{Par:WebsOnConfugurationSpaces} above in the particular case when $m=n=2$.  \sk 

The web $\boldsymbol{\mathcal W}_{\hspace{-0.0cm}{\rm Conf}_{6}(\mathbf P^2)}$ is a 30-web in four variables and from explicit computations, one gets that 
$$
\rho^\bullet\Big( \boldsymbol{\mathcal W}_{\hspace{-0.0cm}{\rm Conf}_{6}(\mathbf P^2)}
\Big)=\big( 26, 20, 11\big) 
\qquad 
\mbox{ and } 
\qquad 
{\rm polrk}^\bullet\Big( \boldsymbol{\mathcal W}_{\hspace{-0.0cm}{\rm Conf}_{6}(\mathbf P^2)}
\Big)=\big( 46, 11\big)
$$
which  implies that $
 {\rm polrk}\big( \boldsymbol{\mathcal W}_{\hspace{-0.0cm}{\rm Conf}_{6}(\mathbf P^2)}
\big)
= {\rm rk}\big( \boldsymbol{\mathcal W}_{\hspace{-0.0cm}{\rm Conf}_{6}(\mathbf P^2)}
\big)=
{\rho}   \big( \boldsymbol{\mathcal W}_{\hspace{-0.0cm}{\rm Conf}_{6}(\mathbf P^2)}
\big)
=57
$ hence in particular that $\boldsymbol{\mathcal W}_{\hspace{-0.0cm}{\rm Conf}_{6}(\mathbf P^2)}$ is AMP as well, with only polylogarithmic ARs (with ramification loci $\{0,1,\infty\}$), of weight 1 or 2 at most. 
\begin{center}
\vspace{-0.3cm}$\star$
\end{center}

Considering that the webs $\boldsymbol{\mathcal W}_{\mathcal M_{0,n+3}}$ (for any $n\geq 2$)  and   $\boldsymbol{\mathcal W}_{\hspace{-0.0cm}{\rm Conf}_{6}(\mathbf P^2)}$ are AMP may suggest that being AMP is a property common to all webs 
$\boldsymbol{\mathcal W}_{\hspace{-0.0cm}{\rm Conf}_{m+n+2}(\mathbf P^m)}$ (for any $m,n\geq 1$ such that $m+n>4$). This is not the case. For instance, let us look at  
 the web $\boldsymbol{\mathcal W}_{\hspace{-0.0cm}{\rm Conf}_{7}(\mathbf P^2)}$ (case $m=2$ and $n=3$): it is a 105-web of intrinsic dimension 6, such that 
 $\rho^\bullet\big(\boldsymbol{\mathcal W}_{\hspace{-0.0cm}{\rm Conf}_{7}(\mathbf P^2)}\big) =(99,84,53)$. One can verify\footnote{This can be verified using the approach 
 described in \S\ref{SubPar:CharacterizationWebsMaximalRank} below.} that  ${\rm rk}\big(\boldsymbol{\mathcal W}_{\hspace{-0.0cm}{\rm Conf}_{7}(\mathbf P^2)}\big) =232$
  and ${\rm polrk}^\bullet\big(\boldsymbol{\mathcal W}_{\hspace{-0.0cm}{\rm Conf}_{7}(\mathbf P^2)}\big) =(182,50)$. We thus have  
 $$232={\rm polrk}\big(\boldsymbol{\mathcal W}_{\hspace{-0.0cm}{\rm Conf}_{7}(\mathbf P^2)}\big) ={\rm rk}\big(\boldsymbol{\mathcal W}_{\hspace{-0.0cm}{\rm Conf}_{7}(\mathbf P^2)}\big) <\rho\big(\boldsymbol{\mathcal W}_{\hspace{-0.0cm}{\rm Conf}_{7}(\mathbf P^2)}\big) =236 $$  hence  although it only carries  polylogarithmic ARs of weight 1 or 2, $\boldsymbol{\mathcal W}_{\hspace{-0.0cm}{\rm Conf}_{7}(\mathbf P^2)}$ is not AMP.



\subsection{\bf Determining the abelian relations and the rank.} 
\label{Para:DeterminingARs-and-Rank}
We want to discuss here some methods to determine the abelian relations and the rank of a given web. Everything presented below in this paragraph is already known but since we  have effectively used some of these methods  further in this paper, we think it  worth telling a few words about this here.\sk 

We first briefly describe how Abel's method  for solving functional equations specializes in the case of AFE's. We then discuss   a  characterization of AMP webs inspired by classical results related to this question for planar webs. We obtain a criterion characterizing AMP webs which can be used quite efficiently for webs defined by rational first integrals.
 Finally, we state a general conjecture about the  nature of the components of any AR of such a web.    \mk 
 
 In the sequel of the present paragraph, we use the following notation: 
  $\Omega$ is an open domain  in $\mathbf C^n$ with $n>2$ and 
  $\boldsymbol{\mathcal W}$ is a web on it, defined by $d$ holomorphic 
  submersions $v_1,\ldots,v_d\in \mathcal O(\Omega)$. 
  Moreover, to simplify,  we will assume that $\boldsymbol{\mathcal W}$ has no singularity on $\Omega$, {\it i.e.}\,for any $i,j$ distinct, the 2-form $dv_i\wedge dv_j$ does not vanish at any point of $\Omega$.
%
  %


\subsubsection{\bf Abel's method for solving AFEs.} In his first published paper \cite{Abel1923}, Abel describes a general method to 
determine  unknown fonctions $F_1,\ldots,F_d$ of one variable  when these are assumed to satisfy a multivariable functional equation $(F\hspace{-0.04cm}E)$. {\it `Abel's method'}\footnote{We attribute this method to Abel since it was first formalized
 in a general setting 
 by him  in \cite{Abel1923}, but  it was of course already known by many of his predecessors. For instance Pfaff already  used  this method in a systematic way to solve several functional equations in 
his  {\it Programma inaugurale}
  presented on the occasion
of his election as professor of mathematics at the University of Helmstedt 
in 1788 (see \cite{Dhombres}). \label{xxx}} is 
quite general and consists in performing differential elimination starting from 
the initial equation $(F \hspace{-0.04cm} E)$ in order to  construct some ODE's satisfied by each of the $F_i$ appearing in it. Then solving the latter allows to determine 
the solutions of the  original functional equation in several variables.\footnote{Rigorously speaking, in \cite{Abel1923} Abel only considers  FEs involving one unkonwn function. The remark that the very same method applies to solve FEs involving several unknown functions is due to Pexider \cite{Pexider}.}

\mk


 The point is that Abel's method, which is a bit vague in full generality, specializes in a precise and effective way when applied to AFE's. 
 \sk 
 
 First, note that by considering the restriction of the web along a sufficiently generic surface in $\Omega$, one can assume that $n=2$  without loss of generality.  
  Then Abel's method for AFEs in two variables can be formalized by an algorithm which  has been described in detail in \cite{PirioSelecta} hence we will not do so here.    We only mention that it can be decomposed into two steps: the first only requires standard algebraic operations on meromorphic functions (addition/substraction, multiplication/quotient) as well as, of course, differentiation with respect to both free variables, denoted here by $x$ and $y$.   For each $i$,  we obtain  a linear differential equation in two variables $K_i(F_i)=\sum_{k=1}^{N_i} C_k^i(x,y) F_i^{(k)}(v_i(x,y))=0$ (for some $N_i>0$ depending on $i$) constructed from the initial equation 
 $\sum_{i=1}^d F_i(v_i(x,y))=0$.\sk
  
  Only the second step (possibly) requires a transcendental operation: let $w_i=w_i(x,y)$ be any function such that $(x,y)\mapsto (v_i(x,y),w_i(x,y))$ be invertible (at a given base point).  Then $K_i(F_i)=0$
  is written $\sum_{k=1}^{N_i} \widetilde C_k^i(v_i,w_i) F_i^{(k)}(v_i)=0$
  in the coordinates $(v_i,w_i)$ 
   hence,  by successive differential eliminations, one reduces this equation in order to get a last one whose coefficients only depend on $v_i$, that is of the form 
   $L_i(F_i)=\sum_{k=1}^{M_i} D_k^i(v_i) F_i^{(k)}(v_i)=0$ (for some positive $M_i\leq N_i$). 
   \sk 
   
   Note that a priori the  differential operator $L_i
   =\sum_{k=1}^{M_i} D_k^i \cdot \partial^k
   $ obtained above is not unique: first the two-variables differential equation $K_i(F_i)=0$ depends on the ordering of the $v_k$'s regarding the first differential elimination process. Second, different choices for  the function $w_i(x,y)$ in order to get an ordinary differential equation  in the second step could give distinct differential operators  $L_i$.  In order to get a canonical ODE for $F_i$, a natural possibility would be to consider the greatest common divisors of the linear differential operators obtained by the two steps described above.   But in practice, knowing only one linear ODE $L_i(F_i)=0$ obtained by Abel's method  is sufficient to solve the initial AFE. \mk 
   
A case of interest for us is when the $v_i$'s are rational and 
such that for each $i$, there exists another rational   function $w_i\in \mathbf C(x,y)$ making of  $(x,y)\dashrightarrow (v_i(x,y),w_i(x,y))$ a Cremona transformation (the webs associated to cluster algebras that we will consider later in the text will all be of this kind).  In this case, one can  perform Abel's method in order to get ODEs 
$L_i(F_i)=\sum_{k=1}^{M_i} D_k^i(v_i) F_i^{(k)}(v_i)=0$ with rational coefficients, {\it i.e.}\,$D_k^i(z)\in \mathbf C(z)$ for all $i=1,\ldots, d$ and all $k=1,\ldots, M_i$.  
 \mk 
 
As an example, we consider a  4-web which is a subweb of a cluster web that we will consider later (namely, it is the subweb $\boldsymbol{\mathcal W}_{\hspace{-0.1cm} G_2}^{long}$ of the `cluster web of type $G_2$', see \S\ref{SS:cluster-webs:G2}): 
 \begin{exm} 
 \label{Ex:WG2-short-Abel'sMethod}
 Applied to the 4-web  
 $$  \boldsymbol{\mathcal W}\bigg(\, \frac{1}{x}\, ,
 \frac{(1+y)^3}{x}\, , \, \frac{(1+x  + y)^3}{xy^3}\, , \,  \frac{(x +(1+y)^2)^3}{x^2y^3}\, \bigg)\, , $$  
  Abel's method gives us that any component $F_1$ of a solution of the associated AFE must verify $L_1\big(F_1(z)\big)=0$ where $L_1$ stands for the following  composition of  four differential operators of the first order: 
$$
L_1= \bigg[\frac{\partial}{\partial z} + \frac{11z+8}{3z(1+z)}\bigg]\cdot \bigg[\frac{\partial}{\partial z} + 
 \frac{7z + 4}{3z(1+z)} \bigg]\cdot \bigg[
 \frac{\partial}{\partial z} + \frac{2z + 1}{z(1+z)}\bigg]\cdot \frac{\partial}{\partial z}\, 
 $$
 whose a  basis of the space of solutions is given by 
 $$
 {1}\, , \quad {\rm Log}(1+z)\, , \quad {\rm Log}\Big(1+z^{1/3}\Big)\quad \mbox{ and  } \quad 
 {\rm Log}\Big(1 - z^{1/3} + z^{2/3}\Big)\, .
 $$
 \end{exm}

\begin{rem}
Surprisingly, the fact that   Abel's method applies in a systematic way to find the solutions of any AFE does not seem to be well known, even by some 
people working in the field of functional equations. For instance, the paper  {\rm \cite{Ebanks}} concerns the solutions of 
$ F(x) + F(y) - F(x + y) = G(h(x, y))$ where  $h$  stands for a given function of two (real or complex) variables. The author mentions  that the problem of  finding a general method for solving
functional equations of this type was posed (by R. Ger and L. Reich) at  the 37th International Symposium on Functional Equations  in 2000.
 And this is not an isolated example: for some other recent papers, among many,  each dealing  with a specific AFE 
which can be easily solved by means of Abel's method, see 
{\it e.g.}\,{\rm \cite{KS}}  (which concerns 
 $F(x+y+\alpha xy)+G(xy)=H(x)+K(y)$), {\rm \cite{Wesolowski}} (about 
 $ G_1(\frac{1-x}{1-xy})+G_2(\frac{1-y}{1-xy})=\alpha_1(x)+\alpha_2(y)$) or 
%
{\rm \cite{GL}} (for the AFE $F(xy) + G(x + y) = H(xy + x) + K(y)$). 
%
%
\end{rem}



%
%
%
%


\subsubsection{\bf Characterization of webs with maximal rank.} 
  \label{SubPar:CharacterizationWebsMaximalRank}
A general  method to determine the rank of planar webs, and in particular to characterize such webs of maximal rank, 
  has been described in 1941 by Pantazi in a short note  \cite{Pantazi}\footnote{Pantazi's note }reproduced in \cite{PirioThese}) which has remained 
   ignored until quite recently. Unaware of  Pantazi's work and using a more modern approach,  H\'enaut has obtained 
    similar results in \cite{Henaut}.   
       A synthesis of both approaches has been given in \cite{Coloquio}, which can be referred to if more information is needed. 
    \mk

 The point is that Pantazi's approach generalizes quite straightforwardly to webs in arbitrary dimension, at least to characterize those whose rank is maximal.  
 We give below a concise exposition of such a generalization, following 
 quite closely \cite[\S6.3]{Coloquio}, in particular using notations similar to those used in it. Regarding the  web $ \boldsymbol{\mathcal W} =\boldsymbol{\mathcal W}(v_1,\ldots,v_d)$ on $\mathcal O(\Omega)$ considered here, 
 we will make the following hypotheses (just to simplify):  
  \sk
 
${}^{}$ \quad  $-$ for any $i=1,\ldots,d$,  the partial derivative $\partial_{x_1}(v_i)$ does not vanish identically on $\Omega$;\sk

${}^{}$ \quad  $-$ the map $\omega \mapsto \rho^\bullet_\omega(  \boldsymbol{\mathcal W})$ is constant on $\Omega$.
  \mk 
  
  Here are some notations used below:  for $\sigma\geq 1$, one writes $\rho^\sigma$ for $\rho^\sigma( \boldsymbol{\mathcal W})$, $\rho$ for the 
total virtual rank $\rho=\sum_{\sigma\geq 0} \rho^\sigma$ and one sets 
$s=\max\{\,\sigma\geq 1\, \lvert \, \rho^\sigma    >0\, \}$. 
   For each $i=1,\ldots,d$, we consider the vector fields 
 $X_{i,s}=\partial_{s}(v_i)\partial_1-\partial_{1}(v_i)\partial_s$ for $s=2,\ldots,n$. These are such that for any holomorphic germ $\Phi_i$ (at a point of $\Omega$), one has $\Phi_i=F_i(v_i)$ for a holomorphic germ $F_i$ in one variable  
  if and only if 
  $X_{i,s}(\Phi_i)=0$ for $s=2,\ldots,n$.\mk

  Clearly, the space of (functional) ARs of $\boldsymbol{\mathcal W}$ identifies with the space of solutions of the following 
 first order linear differential system 
 $$
\boldsymbol{S}_{\boldsymbol{\mathcal W}} \, : \qquad 
\begin{cases}
 \hspace{0.15cm} \sum_{i=1}^d \Phi_i=0 \, ,  \vspace{0.15cm} \\ 
  \hspace{0.15cm} X_{i,s}(\Phi_i)=0 \quad \mbox{ for }\,
 i=1,\ldots,d, \,  s=2,\ldots,n\, . 
 \end{cases}
 $$
 The latter can and will be
 identified to a sub-bundle of the bundle $J^1(E)$ of the first jets of sections of $E$, again denoted by the same notation.  \sk
 
 For any $k\geq 0$, denote by $\boldsymbol{S}_{\boldsymbol{\mathcal W}}^{(k)}\subset J^{k+1}(E)$ the $(k-1)$-th differential prolongation of $\boldsymbol{S}_{\boldsymbol{\mathcal W}}$ (defined inductively by  $\boldsymbol{S}^{(0)}=\boldsymbol{S}_{\boldsymbol{\mathcal W}}$ 
 and $\boldsymbol{S}_{\boldsymbol{\mathcal W}}^{(k+1)}= (\boldsymbol{S}_{\boldsymbol{\mathcal W}}^{(k)} )^{(1)}$ for any $k\geq 0$). 
 The differential system  $\boldsymbol{S}_{\boldsymbol{\mathcal W}}^{(k)}$ can be understood more concretely by realizing that its fiber at any $\omega\in \Omega$
  naturally identifies with the space of $d$-tuples of jets 
$(f_i)_{i=1}^d$ 
of order $k$ which satisfies $\sum_i f_i(v_i)=0$ at  the order $k$. We thus have an isomorphism of vector spaces: 
\begin{equation}
\label{Eq:SWx-k-Identifiction}
{\boldsymbol{S}}_{{\hspace{-0.1cm}}\boldsymbol{\mathcal W},\omega}^{(k)}\simeq 
\left\{\, \big(f_i\big)_{i=1}^d\in 
\Big( \mathbf C[t]_{\leq k} \Big)^d  \hspace{0.1cm} \Big\lvert 
\,  \sum_{i=1}^d f_i(v_i)\in   \mathbf C+\Big(\mathfrak M_{\omega}\Big)^{k+1}\, 
\right\}
\end{equation}
 (here $\mathfrak M_{\omega}$ stands for the maximal ideal of $\mathcal O_\omega$). 
 \bk

 The symbol $\mathfrak S_\omega(\boldsymbol{S}^{(k)}_{\boldsymbol{\mathcal W}})$ (at any point $\omega$ of $\Omega$) identifies with 
 the space $A_\omega(\boldsymbol{\mathcal W})$ hence has dimension $\rho^{k}$. Consequently, setting 
 $$
{\boldsymbol{\mathscr S}}_{{\hspace{-0.06cm}}\boldsymbol{\mathcal W}}=\boldsymbol{S}^{(s)}_{\boldsymbol{\mathcal W}}
\subset J^{s}(E)\, ,  $$
it comes that the symbol of ${\boldsymbol{\mathscr S}}_{{\hspace{-0.06cm}}\boldsymbol{\mathcal W}}^{(1)}$ 
is trivial. Hence, it follows from a result of the theory of formal integrability of differential systems  (see Theorem 6.3.1 in \cite{Coloquio}) that, assuming that ${\boldsymbol{\mathscr S}}_{{\hspace{-0.06cm}}\boldsymbol{\mathcal W}}^{(1)}$ is a sub-bundle of 
$J^{s+1}(E)$\footnote{This means that, the fibers of ${\boldsymbol{\mathscr S}}_{{\hspace{-0.06cm}}\boldsymbol{\mathcal W}}^{(1)}$, which are vector-subspaces of $J^{s+1}(E)$, all have the same dimension.},   then it  is formally integrable, hence completely integrable since it is holomorphic\footnote{For regular analytic differential systems, 
 formal integrability implies holomorphic integrability. It is 
a well-known  consequence of the classical {\it `method of majorants'}.}
if (and only if) the restriction of the projection map $J^{s+2}(E)\rightarrow J^s(E)$ induces an isomorphism ${\boldsymbol{\mathscr S}}_{{\hspace{-0.1cm}}\boldsymbol{\mathcal W}}^{(2)}\rightarrow {\boldsymbol{\mathscr S}}_{{\hspace{-0.1cm}}\boldsymbol{\mathcal W}}$. 
 From this, we deduce the 
 \begin{thm} 
 \label{T:two-ass-equiv}
 The two following  assertions are equivalent:
 \begin{enumerate}
\item the web $ \boldsymbol{\mathcal W}$ has  AMP rank;
 \item one has $\dim_{\mathbf C} \Big( {\boldsymbol{S}}_{{\hspace{-0.03cm}} \boldsymbol{\mathcal W},\omega}^{(s+2)} \Big)=
 \rho\big(\boldsymbol{\mathcal W}\big)+(d-1)$ for any generic point $\omega\in \Omega$.
 \end{enumerate}
  \end{thm}
 
 The latter criterion for complete integrability can be made  more concrete  by using the identification \eqref{Eq:SWx-k-Identifiction} for $k=s+2$.  The situation is even better when the $v_i$'s all are rational (and then are denoted by $u_i$). Indeed, in this case 
 the process of localizing at the generic point of $\mathbf C^n$ can be achieved quite simply in an algebraic effective way. 
\sk 
 
  Indeed, when $u_i\in \mathbf C(x)$ with $x=(x_1,\ldots,x_n)$ for $i=1,\ldots,d$, consider  a $n$-tuple $z=(z_1,\ldots,z_n)$ of new indeterminates and set $\mathbf k=\mathbf C(z)$. 
   Then write $\tilde u_i$ for $u_i(x+z)$ viewed as an element of 
   the ring of formal series $\mathbf k[[x]]$.  
   Then, denoting by $ \mathfrak M_{\mathbf k}$ the maximal ideal of  $\mathbf k[[x]]$ formed by the polynomials in the $x_i$'s with trivial constant term, we have that 
   \begin{equation}
\label{Eq:sur-C(z)}
\boldsymbol{S}_ {\boldsymbol{\mathcal W}}^{(s+2)}({\mathbf k})=
\left\{\, \big(f_i\big)_{i=1}^d\in 
\Big(
\mathbf k[t]_{\leq s+2}\Big)^d \hspace{0.1cm} \Big\lvert 
\,  \sum_{i=1}^d f_i(\tilde u_i)\in   {\mathbf k}+\Big(\mathfrak M_{\mathbf k}\Big)^{s+3}\, 
\right\}
\end{equation}
is isomorphic to the RHS of \eqref{Eq:SWx-k-Identifiction} for $\omega$ being  {\it `the generic point of $\mathbf C^n$'}.  We deduce the 

\begin{cor} When all  the $u_i$'s  are rational, the following equivalence holds true: 
$$\boldsymbol{\mathcal W} \, 
\mbox{ is AMP }\,  \hspace{0.2cm} \Longleftrightarrow \,  
  \hspace{0.3cm}
 \dim_{\mathbf k} \Big( \boldsymbol{S}_ {\boldsymbol{\mathcal W}}^{(s+2)}({\mathbf k})
 \Big)=\rho(\boldsymbol{\mathcal W})+(d-1) \, .
$$
\end{cor}
   
   The practical interest of this corollary is clear 
 for verifying if a web defined by  rational functions  is AMP or not:  
 it reduces this question to   computing the dimension of 
a finite dimensional vector space (over ${\mathbf k}$), which can be handled by standard techniques  of  linear algebra.
  \begin{center}
\vspace{-0.1cm}  
$\star$
\end{center} 

Building on Pantazi's characterization of planar webs with maximal rank, 
Mih\u{a}leanu's has sketched  an effective method to determine the rank of a given planar web in \cite{Mihaileanu}. It turns out that Mih\u{a}leanu's approach  generalizes as well, and rather straightforwardly, to any 1-codimensional web  $\boldsymbol{\mathcal W}$ as above ({\it i.e.}\,only assumed to satisfy  assumption (wGP)  of
 \S\ref{Par:GermOfWebs}).  We get the following result:


 \begin{prop}
 Let $\zeta$ be the smallest  integer bigger than or equal to $s$ such that 
 the map  ${\boldsymbol{S}}_{{\hspace{-0.03cm}} \boldsymbol{\mathcal W},\omega}^{(\zeta+2)} \rightarrow {\boldsymbol{S}}_{{\hspace{-0.03cm}} \boldsymbol{\mathcal W},\omega}^{(\zeta)}$ is surjective for  $\omega$ generic. Then 
 ${\rm rk}\big( \boldsymbol{\mathcal W} \big)= \dim_{\mathbf C} \big( {\boldsymbol{S}}_{{\hspace{-0.03cm}} \boldsymbol{\mathcal W},\omega}^{(\zeta)} \big)-(d-1)$.
  \end{prop}

Given a point $\omega$, determining the smallest  integer $\zeta$ bigger than or equal to $s$ such that 
 the map  ${\boldsymbol{S}}_{{\hspace{-0.03cm}} \boldsymbol{\mathcal W},\omega}^{(\zeta+2)} \rightarrow {\boldsymbol{S}}_{{\hspace{-0.03cm}} \boldsymbol{\mathcal W},\omega}^{(\zeta)}$ is surjective only depends on the 
 $(\zeta+2)$-th order jet  of  $ \boldsymbol{\mathcal W}$ at $\omega$. 
Consequently, the same approach as the one described just after Theorem \ref{T:two-ass-equiv} can be followed which allows us to get an effective general method for determining  the rank of webs from the preceding proposition.
\sk 

\vspace{-0.25cm}
We have implemented this effective method for determining the rank in the computer algebra system Maple and used it to compute the rank of many webs considered in the sequel.

 \subsubsection{\bf On the nature of components of rational AFEs.}  
 \label{SS:NatureComponentsRationalAFEs}
We now discuss briefly the 
 nature of the components $F_i$ of an AFE  $\sum_{i=1}^d F_i(u_i)=0$ when the 
 $u_i$'s all  are rational. As before, $\boldsymbol{\mathcal W}$ stands for the web defined by 
  the $u_i$'s and 
 $\mathfrak B_i$ denotes the corresponding $i$-th branch locus. 
 We call a {\bf generalized iterated  integral} any (multivalued) holomorphic function  on a given projective manifold $M$,  which can be written as a finite sum $\sum_k R_k\cdot {I}_k$ where the $R_k$'s are rational functions and the $I_k$'s iterated integrals (in the usual sense, {\it cf.}\,\S\ref{Par:II-AR}  and especially Remark {\bf \ref{Rem:II-dim>0}.2}).\sk 
 
 The author started thinking about such functional equation at the beginning of the 2000's, when working on his PhD about webs. Since then, he has considered and solved quite a lot of such rational AFEs. Considering all the solutions of these leads to state the following
\begin{conjecture}
\label{Conj:Nature-Fi}
Let  
 $\boldsymbol{F}=(F_i)_{i=1}^d$ be 
a germ of functional AR of $\boldsymbol{\mathcal W}$ at $ x\in \mathbf P^n\setminus \Sigma(\boldsymbol{\mathcal W})$, {\it i.e.}\,each  $F_i$ is a holomorphic germ at $u_i({x})$  and 
$\sum_{i=1}^d F_i(u_i)=0$ holds identically in the vicinity of $x$.\sk

  Then   $F_i=J_i\circ A_i$  for each $i$, 
  where $J_i$ is a generalized iterated integral (with poles and branching points in $\mathfrak B_i$) and 
 $ A_i$ an algebraic function.  
 Moreover, the $A_i$'s can be chosen  independently of $\boldsymbol{F}$  (hence they only depend  on the $u_i$'s).
\end{conjecture}

In addition to the fact that we are not aware of any counter-example to the previous statement,  our belief that it indeed holds true is supported by the following result by Aomoto, which has also inspired us the  statement of our conjecture:   {\it ``a global multivalued function on $\mathbf P^n$ is a generalized iterated integral if and only if it has unipotent monodromy and 
has moderate growth along its branching locus''}
 (see \cite{Aomoto}, Corollaire 1 p.\,154).
\mk 

Aomoto's result suggests an approach to prove Conjecture \ref{Conj:Nature-Fi}: 
prove that (1) possibly up to precomposition by an algebraic function, the monodromy of $F_i$ is unipotent; and (2)  each $F_i$ has moderate growth at any point of the branch locus $\mathfrak B_i$. 
 A way to get (2) could be to prove that  the  linear ODE's $L_i$ obtained through 
Abel's method have regular singularities (on the ramified coverings of $\mathbf P^1$ on which they live).  As for (1), we believe that it is equivalent to the fact that 
a suitable power of the monodromy operator $\mu_{\boldsymbol{\mathcal W}}: \pi_1(\mathbf P^n\setminus  \Sigma(\boldsymbol{\mathcal W}),x)\rightarrow \mathcal A_x({\boldsymbol{\mathcal W}})$ is unipotent, {\it i.e.}\,there exist two integers $N,M\in \mathbf N_{>0}$ such that $(\mu_{\boldsymbol{\mathcal W}}^N-{\rm Id})^M=0$. In order to prove that such a formula holds true, a first step would be  to have an idea of what the integers $N$ and $M$ might be.

\subsection{\bf Some results of planar web geometry}
\label{SS:Results-Planar-Web-Geometry}
We gather here some notions and results specific to planar webs that we will use later on in the text, essentially in \S\ref{SS:ClassicalPolylogarithmicIdentitiesAreOfClusterType}.  For more details and proofs concerning the material presented below, the general reference we refer to is \cite{Coloquio}.  
Another very valuable reference is Blaschke-Bol's classic book in german \cite{BB}. 

\subsubsection{\bf Hexagonality, curvature and flatness.}
A classical and important notion in planar web geometry is that of `hexagonality': a planar 3-web $\boldsymbol{\mathcal W}$ is {\bf hexagonal} if a walk along its leaves around any given base-point forms a  closed hexagon (see \cite[\S1.2.1]{Coloquio}).  \sk 

A basic result of the theory of webs is that the following conditions are equivalent\footnote{A partial version of this statement is due to Thomsen together with Blaschke (1927) and is considered as marking the birth of the study of webs for themselves (see the Appendix of  \cite{Coloquio}
 for a brief overview of web geometry from a historical perspective).}
\vspace{0.1cm} 

${}^{}$ \, $-$ $\boldsymbol{\mathcal W}$ is hexagonal;\vspace{0.1cm} 

${}^{}$ \, $-$ $\boldsymbol{\mathcal W}$ is parallelizable ({\it cf.}\,\S\ref{SS:Coordinate-parallelizable-webs});\vspace{0.1cm} 

${}^{}$ \, $-$ $\boldsymbol{\mathcal W}$ has maximal rank (namely ${\rm rk}(\boldsymbol{\mathcal W})=1$); \vspace{0.1cm} 

${}^{}$ \, $-$ $\boldsymbol{\mathcal W}$ is {\bf flat}, {\it i.e.} its curvature $\mathfrak K(\boldsymbol{\mathcal W})$  vanishes identically.\mk 

We recall that the {\bf curvature $\boldsymbol{\mathfrak K(\boldsymbol{\mathcal W})}$} of $\boldsymbol{\mathcal W}$ is a 2-form canonically attached to it, which is the most basic invariant attached to $\boldsymbol{\mathcal W}$ ({\it cf.}\,\cite[\S1.2.2]{Coloquio}).  When 
$x,y$ are local coordinates such that 
$\boldsymbol{\mathcal W}=\boldsymbol{\mathcal W}(x,y,U(x,y))$  for a function  of two variables $U$  (such that neither $\partial_x U$ nor $\partial_y U$ vanishes), one has $\mathfrak K(\boldsymbol{\mathcal W})=\partial^2_{xy}\big( {\rm Log}(\partial_x U/\partial_y U)\big) \,dx\wedge dy$  and there exists a more general but explicit formula for computing $\mathfrak K(\boldsymbol{\mathcal W})$  in terms of any 3-tuple of first integrals (see \cite[Theorem\,0.3]{PirioChern}). \begin{center}
$\star$
\end{center}

Interesting features of the notions of hexagonality and of curvature for 3-webs is that they furnish relevant characteristics of planar $d$-webs for any $d\geq 3$:  
a planar $d$-web $\boldsymbol{\mathcal W}_d$  is said to be {\bf hexagonal} if all its 
3-subwebs are hexagonal and one defines the {\bf curvature} of $\boldsymbol{\mathcal W}_d$, denoted by 
$\mathfrak K(\boldsymbol{\mathcal W}_d)$, as the sum of the curvatures of all its $3$-subwebs, that is $\mathfrak K(\boldsymbol{\mathcal W}_d)=\sum_{\boldsymbol{\mathcal W}_3< \boldsymbol{\mathcal W}_d} \mathfrak K(\boldsymbol{\mathcal W}_3)$. Since the curvature of a web can be computed effectively,   the flatness of $\boldsymbol{\mathcal W}_d$ can easily be checked in practice as soon as a $d$-tuple of first integrals for this web has been given. \sk



The notions  of hexagonality, of curvature and of flatness for $d$-webs with $d\geq 3$ are useful and relevant in many ways. First,  
as soon as a $d$-tuple of first integrals of a web $\boldsymbol{\mathcal W}_d$ has been given, its curvature can be computed quite easily (by means of a computer algebra system) hence the fact that $\boldsymbol{\mathcal W}_d$  is hexagonal or flat can be effectively checked. Moreover, these are invariant notions attached to $\boldsymbol{\mathcal W}_d$ and above all, they are meaningful regarding the characterization of planar webs of maximal rank thanks to the following theorems: 


\begin{thm}[\cite{Bol}]
Let $\boldsymbol{\mathcal W}$ be a planar hexagonal $d$-web. Then one of the two following possibilities occurs, depending on whether $\boldsymbol{\mathcal W}$ is linearizable or not: 
\begin{itemize}
\item either $\boldsymbol{\mathcal W}$ is equivalent to a web formed by d pencils of lines; or
\item the web $\boldsymbol{\mathcal W}$ is not linearizable, in which case $d=5$ and $\boldsymbol{\mathcal W}$ is equivalent to Bol's web $\boldsymbol{\mathcal B}$. 
\end{itemize}
\end{thm}
\sk

\begin{thm}
{1.} Planar webs of maximal rank are necessarily flat. \sk 

{2.} For linearizable planar webs,  flatness is sufficient for ensuring the maximality of the rank.
\end{thm}

In the preceding theorem, the first statement is a consequence of Pantazi's characterization 
of planar webs with maximal rank in  \cite{Pantazi} (discussed and generalized above in \S\ref{SubPar:CharacterizationWebsMaximalRank}).  For a linear web, the curvature admits a simple and nice expression, the vanishing of which allows to apply an Abel-inverse type theorem which in its turn implies  that the considered web actually is algebraic hence of maximal rank (see \S6.3.4 in \cite{Coloquio} for more details).

\subsubsection{\bf Arithmetic invariants.}
 Let us say that a property $\boldsymbol{\mathcal P}$ of planar webs is  `invariant' when 
the fact that it is satisfied or not only depends on the equivalence classes of webs which are considered. By their very definition, the linearizability or the  algebraizability of a web are examples of such properties.  Other examples are hexagonality, or flatness, denoted  by $\boldsymbol{\mathcal Hex}$ and  $\boldsymbol{\mathcal F \hspace{-0.1cm}lat}$ respectively. 
\sk

There is a very basic recipe to produce arithmetic invariants for webs from a 
given invariant property $\boldsymbol{\mathcal P}$ as above. Given a planar $d$-web 
$\boldsymbol{\mathcal W}_{\hspace{-0.07cm}d}$ (with $d\geq 3$) and for any $k\in \{3,\ldots,d\}$, 
 one  defines $\boldsymbol{\mathcal P}_{\hspace{-0.07cm}k}(\boldsymbol{\mathcal W})$ as the number of $k$-subwebs of $\boldsymbol{\mathcal W}_d$ satisfying $\boldsymbol{\mathcal P}$ and one denotes 
 by $\boldsymbol{\mathcal P}_{\hspace{-0.07cm}\bullet}(\boldsymbol{\mathcal W}_d)$ the tuple of the $\boldsymbol{\mathcal P}_{\hspace{-0.07cm}k}(\boldsymbol{\mathcal W}_d)$'s for $k$ ranging from 3 to $d$: one has 
 $ \boldsymbol{\mathcal P}_{\hspace{-0.07cm}\bullet}(\boldsymbol{\mathcal W}_{\hspace{-0.07cm}d})= \big( \boldsymbol{\mathcal P}_{\hspace{-0.07cm}k}(\boldsymbol{\mathcal W}_{\hspace{-0.07cm}d})\big)_{k=3}^d $ with 
 $$
\boldsymbol{\mathcal P}_{\hspace{-0.07cm}k}\big( \boldsymbol{\mathcal W}_{\hspace{-0.07cm}d} \big)
={\rm Card}\, \bigg(\, \Big\{\hspace{0.15cm}  \boldsymbol{\mathcal W}_{\hspace{-0.07cm}k}
\subset 
\boldsymbol{\mathcal W}_{\hspace{-0.07cm}d}\hspace{0.15cm} 
\big\lvert \hspace{0.15cm} 
\boldsymbol{\mathcal W}_{\hspace{-0.07cm}k} 
\mbox{ satisfies } \, 
\boldsymbol{\mathcal P}\hspace{0.1cm} 
\Big\}\, \bigg)\quad \mbox{for} \quad  k=3,\ldots,d\,.
$$

Given a web $\boldsymbol{W}$, we denote a bit abusively in the same way the property of being equivalent to it. 
 Then for any other web $\boldsymbol{\mathcal W}$, one denotes by $\boldsymbol{W}(\boldsymbol{\mathcal W})$ the cardinal (which can of course be zero) of the set of subwebs of $\boldsymbol{\mathcal W}$ which are equivalent to $\boldsymbol{W}$. 
For instance, in the case when  $\boldsymbol{W}=\boldsymbol{\mathcal B}$ (Bol's web),  
$\boldsymbol{\mathcal B}(\boldsymbol{\mathcal W})$ stands for the number of 5-subwebs of $\boldsymbol{\mathcal W}$ which are equivalent to $\boldsymbol{\mathcal B}$. 
\sk 

The interest of the invariants obtained as described just above is that they are arithmetic 
invariants which, for most of them, are easy to compute hence can be  useful when  dealing with the problem of determining whether two given webs are equivalent or not.

\begin{exm}
For instance, let us consider the case of the Spence-Kummer web ${\boldsymbol{\mathcal W}}_{\hspace{-0.1cm}{\cal S}{\cal K}}$ which is important in web geometry since it is exceptional and carries trilogarithmic ARs (it will be discussed more in depth further on in 
 \S\ref{SS:Spence-Kummer}). By elementary computations, one gets that 
 $$ \boldsymbol{\mathcal Hex}_{\hspace{-0.02cm}\leq 6}\Big({\boldsymbol{\mathcal W}}_{\hspace{-0.1cm}{\cal S}{\cal K}}\Big)=(48,30,9,1)  \, , 
\quad  \boldsymbol{\mathcal F \hspace{-0.1cm}lat}_{\hspace{-0.02cm}\bullet} \Big({\boldsymbol{\mathcal W}}_{\hspace{-0.1cm}{\cal S}{\cal K}}\Big)=(48, 48, 12, 11, 3, 0, 1)
 \quad \mbox{ and }\quad  
  {\boldsymbol{\mathcal B}}\big({\boldsymbol{\mathcal W}}_{\hspace{-0.1cm}{\cal S}{\cal K}}\big)=3\, . 
$$  
Any planar 9-web with the same invariants  can be suspected of being equivalent  to  ${\boldsymbol{\mathcal W}}_{\hspace{-0.1cm}{\cal S}{\cal K}}$. 
\end{exm}

\subsubsection{\bf Linearizability.}
We recall some results about the characterization of linearizable planar webs. This  is  classical and well understood for webs of degree bigger than or equal to 4 ({\it cf.}\,\cite[\S6.1]{Coloquio} and the references given therein).\sk 

Let ${\boldsymbol{\mathcal W}}$ be a $d$-web with $d\geq 3$.  Then the following equivalence is more or less obvious : \sk\\
\begin{tabular}{c}
{\it the web ${\boldsymbol{\mathcal W}}$  is linearizable} \quad $\Longleftrightarrow$ \quad 
\begin{tabular}{l}
{\it (1) \, ${\boldsymbol{\mathcal W}}$ is compatible with a projective connection; and} \\
{\it (2) \, this projective connection is flat.}
\end{tabular}
\end{tabular}
\sk 

\vspace{-0.2cm}
Condition {\it (1)} above is then necessary in order that 
${\boldsymbol{\mathcal W}}$   be linearizable and whether this condition holds true or not can be verified quite easily.  Indeed, assume that  the leaves of ${\boldsymbol{\mathcal W}}$
 are the integral curves of $d$ vector fields $X_1,\ldots,X_d$ on 
a domain  $U\subset \mathbf C^2$.  Without loss of generality (since being linearizable is a local condition), one can assume that each $X_i$ is written $X_i=\partial_x-b_i\, \partial_y$ for a certain holomorphic function $b_i$ (where $x,y$ stand for the standard coordinates on $ \mathbf C^2$).  The fact that ${\boldsymbol{\mathcal W}}$ is compatible with a projective connection is equivalent to the existence of a 4-tuple $(a_s)_{s=0}^3\in \mathcal O(U)^4$ 
 such  that for every $i=1,\ldots,d$, one has $
\sum_{s=0}^4 a_s (b_i)^s=X_i(b_i) $.  
 Then verifying if condition {\it (1)} above is satisfied is only a matter of linear algebra (with coefficients in the field of meromorphic functions on $U$) and one gets the following well-known and very useful consequence: 
\begin{prop}
\label{P:Unique-Lin}
If $d\geq 4$, the web ${\boldsymbol{\mathcal W}}$  admits at most one linearization (up to post-composition by a projective transformation). 
\end{prop}

\subsubsection{\bf Poincar\'e-Blaschke maps.}
\label{SS:Poincare-Blaschke-Maps}
 The Poincar\'e map of a planar web ${\boldsymbol{\mathcal W}}$ is a map constructed from its ARs (hence which exists when ${\rm rk}({\boldsymbol{\mathcal W}})\geq 3$) and which is very useful in order to study ${\boldsymbol{\mathcal W}}$ from a geometric and more canonical point of view.  Here we will only consider the classical  case of 4-webs, which goes back to Poincar\'e. 
   For more details and references, see \cite[\S4.3.4]{Coloquio}.
 \sk 
 
Let  $u_1,\ldots,u_4$ be four holomorphic submersions defining a 4-web 
${\boldsymbol{\mathcal W}}={\boldsymbol{\mathcal W}}(u_1,\ldots,u_4)$
 on a domain $U\subset \mathbf C^2$.  We assume that ${\boldsymbol{\mathcal W}}$ has maximal rank 3: there exist three 4-tuples of holomorphic functions $\boldsymbol{f}^\lambda=(f_i^\lambda)_{i=1}^4$ for $\lambda=1,2,3$ such that $\sum_{i=1}^4 f_i^\lambda(u_i)\, du_i\equiv 0 $ for any $\lambda$ and inducing a basis of ${\boldsymbol{\mathcal A}}({\boldsymbol{\mathcal W}})$.  
 For any $i=1,\ldots,4$, the map 
 $$\kappa_i=\Big[f_i^1(u_i) :f_i^2(u_i) : f_i^3(u_i)\Big] : U\rightarrow \mathbf P^2$$
  is well-defined and is a canonical first integral of the $i$-th foliation of ${\boldsymbol{\mathcal W}}$. Its image 
  $$C_i={\rm Im}(\kappa_i)\subset \mathbf P^2$$ 
  is a  smooth analytic  curve in $\mathbf P^2$, called the {\bf $\boldsymbol{i}$-th canonical curve of ${\boldsymbol{\mathcal W}}$}. Moreover, for any $x\in U$, the four points $\kappa_1(x),\ldots,\kappa_4(x)$ are on a same line  
 $\ell(x)\subset \mathbf P^2$, which can equivalently be seen as 
 a point, denoted by $[\ell(x)]$, in the dual projective space $\check{\mathbf P}^2$. By letting $x$ vary within $U$, one gets the {\bf Poincar\'e-Blaschke map} of ${\boldsymbol{\mathcal W}}$, namely  
$$
 {\rm P \! B}[{\boldsymbol{\mathcal W}}] \, :\, U\longrightarrow \check{\mathbf P}^2\, , \hspace{0.15cm} x\longmapsto \big[\ell(x)\big]\, .
 $$
 
 The interest of the objects introduced above is asserted by the following fundamental result, first obtained by Poincar\'e (but stated by him in another form\footnote{Note also that the arguments presented by Poincar\'e, if they were very interesting, were not very rigorous and complete.}): 
 \begin{thm} 
 \begin{enumerate}
\item $
 {\rm P \! B}[{\boldsymbol{\mathcal W}}]$ is a local biholomorphism which is invariantly attached to ${\boldsymbol{\mathcal W}}$;
 \item Therefore the push-forward 
 $
 {\boldsymbol{\mathcal W}}^{can}={\rm P \! B}[{\boldsymbol{\mathcal W}}]_*\big( {\boldsymbol{\mathcal W}}
 \big)
 $
 is a canonical model of ${\boldsymbol{\mathcal W}}$;
 \item Actually,   $
 {\boldsymbol{\mathcal W}}^{can}$ is linear hence algebraic: it is the web   associated to the 
 quartic curve $C\subset \mathbf P^2$ such that $\overline{C_i}^{Zar}=C$ for any $i=1,\ldots,4$. 
\end{enumerate}
\end{thm}


\newpage
\section{On polylogarithms and their functional equations}
\label{S:Polylogs-FEs}

In this section, we first introduce the classical polylogarithms in a single variable and discuss their main properties. Then we focus on the functional  equations in several variables  they satisfy. Even if it is not complete,  we give a  substantial review of the known such equations and for each of them, we discuss the associated web.\mk 

In order to make a (very) slight improvement to the existing literature on the subject and also because we find this interesting, we have mentioned many precise historical references. Following a common practice in the literature on history of sciences, these points as well as some associated historical references are included in many footnotes (contrarily to the mathematical references which  all appear in a list at the very end of this text).
\bk 

As general references about polylogarithms,  we refer readers to the books
 written or edited by Lewin \cite{Lewin} and \cite{LewinStruct}, to Zagier's papers \cite{Zagier1991,ZagierSpecial,Zagier}. Another text of general interest on the subject  is Zagier's joint paper with Gangl \cite{GanglZagier}. 
For historical details, an interesting source is \cite[\S10]{Maximon}. 
 For some modern references dealing more specifically with the functional equations satisfied by polylogarithms, 
  the reader can look at some  of the papers by Gangl (for instance 
  \cite{Gangl1} or \cite{Gangl}), at the two  dissertations 
\cite{Charlton} and \cite{Radchenko} or at the recent preprint \cite{Rudenko}. 
\sk

\vspace{-0.4cm}
As for hyperlogarithms, aka iterated integrals associated to logarithmic rational 1-forms on $\mathbf P^1$, we refer to  \cite{Poincare} and  \cite{Lappo} as classical references. For more modern ones, see \cite{Wechsung1967} and the recently published  \cite{BanksPanzerPym}.


 \subsection{\bf Definition and basic properties}
 \label{SS:Definition-basic-properties}
We follows a classical presentation:   after having recalled the definition and the main properties, we indicate how all this generalizes to higher polylogarithms.

\subsubsection{\bf \hspace{-0.2cm}The logarithm.}\hspace{-0.2cm}
The story of polylogarithms of course begins with that of the 
logarithm.\footnote{The birth of the notion of logarithm is attributed independently to 
\href{https://en.wikipedia.org/wiki/John_Napier}{J.\,Napier}
 and to \href{http://www-groups.dcs.st-and.ac.uk/~history/Biographies/Burgi.html}{J. B\"urgi}, who published tables of logarithms 
 \href{https://archive.org/details/mirificilogarit00napi/page/n8}{\it Mirifici logarithmorum canonis descriptio}
 (1614) and  
  \href{https://bildsuche.digitale-sammlungen.de/index.html?c=viewer&bandnummer=bsb00082065&pimage=00007&v=2p&nav=&l=fr}{\it Aritmetische und Geometrische Progress Tabulen}
 (1620) respectively.  Napier and B\"urgi both used (equivalents forms of) \eqref{Eq:EFA-Log}  to perform some computations and the corresponding property of logarithmic values was immediately incorporated to the `common knowledge' about  logarithms of that time, see {\it e.g.}\,Caput II  in 
\href{https://en.wikipedia.org/wiki/Henry_Briggs_(mathematician)}{H. Briggs'} 
 \href{https://archive.org/details/arithmeticalogar00brig/page/2}{\it Arithmetica Logarithmica} (1624) ({\it cf.}\,\href{http://www-history.mcs.st-andrews.ac.uk/Miscellaneous/Briggs/Chapters/Ch2.pdf}{here} for a modern  transcription of Briggs' essay). However, \eqref{Eq:EFA-Log} was not seen as a functional equation per se until about 25 years later since  the logarithm itself was not truly considered as a function at first.} This classical function has many nice very well-known features: indeed,  the logarithm ${\rm Log}$
\begin{itemize}
\item satisfies the  functional equation:\footnote{According to 
  \href{https://imsc.uni-graz.at/gronau/LogNeuhof89.pdf}{D. Gronau}, 
 the fact  that the logarithm satisfies \eqref{Eq:EFA-Log}
  can already  be found in \href{https://en.wikipedia.org/wiki/Johannes_Kepler}{J.\,Kepler's} 
  \href{https://www.e-rara.ch/zut/wihibe/content/titleinfo/1405839}{\it `Chilias}
   \href{https://www.e-rara.ch/zut/wihibe/content/titleinfo/1405839}{\it  logarithmorum'}
  (1624),  but the we are not convinced. 
  In 
   Propositio CIX of his  {\it Opus geometricum} (1647) (see also \href{https://link.springer.com/content/pdf/10.1007/BF01273373.pdf}{here} or \S2.4 \href{https://doi.org/10.1007/978-0-387-92154-9}{there} for modern expositions),  
  the belgian Jesuit
  \href{https://en.wikipedia.org/wiki/Gr\'egoire_de_Saint-Vincent}{G. de Saint-Vincent} establishes that,  considered as a function of the abscissa $x$,  the surface area $A(x)$ under the hyperbola $y=1/x$ satisfies the functional equation  $A(xy)=A(x)+A(y)$ which is formally similar to the one satisfied by the logarithm. 
%
 %
 This leaded G. de Saint-Vincent's student  \href{https://en.wikipedia.org/wiki/Alphonse_Antonio_de_Sarasa}{A.A.\,de Sarasa} to establish a little later the integral representation \eqref{Eq:IntRep-Log} for the logarithm in  \href{https://archive.org/details/bub_gb_TG-i3Dzr7QoC}{\it Solutio problematis a R.P. Marino Mersenno} (1649) (see \href{https://doi.org/10.1006/hmat.2000.2295}{here} for a modern exposition of de Sarasa's work). As far as we know, the first formal occurrence of the functional equation \eqref{Eq:EFA-Log} seems to be on the first page of W.\,Gardiner's \href{https://www.e-rara.ch/zut/content/titleinfo/2672212?lang=en}{\it Tables of logarithms} (1742).}
\begin{equation}
\label{Eq:EFA-Log}
{\rm Log}(xy)={\rm Log}(x)+{\rm Log}(y)\, ; \vspace{0.35cm}
\end{equation}

\item admits the following integral representation:\vspace{-0.25cm}
\begin{equation}
\label{Eq:IntRep-Log}
{\rm Log}(x)=\int_{1}^x \frac{du}{u} ; 
\qquad \quad {}^{}
\vspace{-0.05cm}
\end{equation}
\item is such that the following development in series 
holds true at the origin\footnote{Albeit it has been independently obtained by others scholars of that time such as G.\,de St.-Vincent or I.\,Newton, this development in series was first published by  \href{https://en.wikipedia.org/wiki/Nicholas_Mercator}{N.\,Mercator} in his  \href{https://doi.org/10.3931/e-rara-7457}{\it Logarithmotechnia} (1668).}:  
\begin{equation}
\label{Eq:Devel-Series-Log}
-{\rm Log}(1-x)=\sum_{k=1}^{+\infty} \frac{x^{k}}{k}
\, ; \quad  {}^{}
 \vspace{-0.15cm}
\end{equation}
\item extends as a multivalued holomorphic function on $\mathbf C^*$, with monodromy around the origin\footnote{The definition of the logarithms of negative integers and next of non-zero complex numbers  gave rise to a well-known controversy between famous mathematicians  of the 18th century such as G.W.\,Leibniz, J.\,Bernoulli, \href{https://en.wikipedia.org/wiki/Leonhard_Euler}{L.\,Euler}, solved by the latter in 
a \href{https://archive.org/details/euler-e168}{paper} written in 1749. 
The complete explanation of the multivalued nature of the logarithm of a complex argument by means  of complex integration and   use of  the integral representation \eqref{Eq:IntRep-Log}, was first encountered by \href{https://fr.wikipedia.org/wiki/Carl_Friedrich_Gauss}{C.F.\,Gauss} (see \href{https://archive.org/details/werkenachtraegez101cfga/page/n369}{its 1811 letter to Bessel}.)
}
\begin{equation}
\label{Eq:Monod-Log}
\mathcal M_0 \Big({\rm Log}(x)\Big)= {\rm Log}(x)+2i\pi\, . \quad {}^{}
\end{equation}
\end{itemize}

Among all these properties, the functional equation \eqref{Eq:EFA-Log}
 seems to be the most fundamental\footnote{This opinion has been explicitely formulated by J.\,Pfaff who writes  {\it `indoles logarithmorum hac aequatione fundamentali continetur'}   in his 1788 inaugural essay (refered  in footnote \ref{xxx} page \pageref{xxx})}: it can be stated  only by means 
of the two most basic operations of arithmetic (addition and multiplication)\footnote{
This contrasts with \eqref{Eq:IntRep-Log}   \eqref{Eq:Devel-Series-Log}   \eqref{Eq:Monod-Log}, which are only intelligible by using more sophisticated mathematical concepts.} and suffices to characterize 
 the logarithm  ${\rm Log}$ hence to recover all the properties this function satisfies.

\subsubsection{\bf The dilogarithm.}\hspace{-0.2cm} 
  \label{SS:TheDilogarithm}
 The {\bf dilogarithm} (or {\bf bilogarithm}) ${\boldsymbol{ \l {2} }}$  is a classical function\footnote{The dilogarithm, in the form of the integral $\int_{0}^z {\rm Log(1-u)}du/u$, is nowadays commonly attributed to Euler but was already considered at the 
end of the 17th century, apparently by Leibniz in some letters to the Bernoulli's in the first place.  The dilogarithm explicitly appears 
in (Sectio Prima, 
Caput VI, Exemplum II of) Euler's book \href{https://archive.org/details/institutionescal020326mbp/page/n129}{\it Institutionum calculi integralis} (1768) 
but was considered and studied a few years before by \href{https://en.wikipedia.org/wiki/John_Landen}{J.\,Landen} in his  paper 
\href{https://doi.org/10.1098/rstl.1759.0056}{\it A new method of computing the sums of certain series}  (1759).\label{footLanden}}
%
 %
 which 
 is nowadays usually defined as the sum 
 following series
\begin{equation}
\label{Eq:Devel-Series-Li2}
\l {2} (z)=\sum_{k=1}^{+\infty}\frac{z^k}{k^2}\, 
\end{equation}
which is abslolutely convergent on the unit  disk $\mathbf D
\subset \mathbf C$. It appears as  a generalization of the logarithm since, in addition of \eqref{Eq:Devel-Series-Li2} which has to be compared with the development in series \eqref{Eq:Devel-Series-Log}, it satisfies  similar properties to those of the logarithm listed above. \mk 

${}^{}$\quad Indeed, the dilogarithm: 
\begin{itemize}
\item admits the following integral representation \begin{equation}
\label{Eq:IntRep-Li2}
\l {{2}} (z)=- \int_0^z \frac{{\rm Log}(1-u) }{u}  du\, ; 
\qquad \quad {}^{}
\vspace{0.35cm}
\end{equation} 
\item extends as a multivalued holomorphic function on $\mathbf C\setminus \{0,1\}$, with monodromy around the origin and 1 along a small circle oriented in the direct order given by 
\begin{equation}
\label{Eq:Monod-Li1}
\mathcal M_0 \Big( 
\l { {2} } (z)
\Big)= \l { {2} } (z) \qquad \mbox{and}\qquad 
\mathcal M_1 \Big( 
\l {  {2} } (z)
\Big)= \l { {2} } (z) -2i\pi\,   {\rm Log}\, z
 \vspace{0.15cm}
\end{equation}
\item satisfies the following functional equation in two variables with five dilogarithmic terms:
 \vspace{-0.2cm}
\end{itemize}
\begin{equation}
\label{Eq:EFA-Li2} {}^{} \qquad 
\l{2} (x) -
\l{2} (y)
-\l{2} \left(\,\frac{x}{y}\, \right)
-\l{2} \left(\,\frac{1-y}{1-x}\,\right)+
\l{ 2 } \left(\,\frac{x(1-y)}{y(1-x)}\,\right)
={\rm Log}(y)\,{\rm Log}\left(\frac{1-y}{1-x}\right)
-\frac{{}^{}\, \, \pi^2}{6}\, .
\end{equation}

 The dilogarithm appears in  many domains of mathematics and physics (see \cite{Zagier} for more perspectives) and is  nowadays recognized as one of the most interesting special functions. As for the logarithm, the five-term functional equation above it satisfies seems to be the most important of all of its properties.
\mk 

Due to its importance and because it has been discovered independently  by many mathematicians throughout the XIX-th century (under different but equivalent forms), we think it  worthwhile to devote not just a simple footnote, but a whole paragraph to the history of this functional equation. For more mathematical and historical details on this, the reader can consult the first chapter of Lewin's book \cite{Lewin}.\sk

A functional equation for the dilogarithm equivalent to 
\eqref{Eq:EFA-Li2} seems to have been obtained first by Spence in 1809 (see p.\,9 of its essay \cite{Spence} ({\it cf.} also  \cite[\S1.5.3]{Lewin}).\footnote{The content of  Spence's essay on polylogarithms, in particular in what concerns the functional equations 
they satisfy, 
is discussed with many details in the recent historical paper \cite{Craik} .}
Spence's essay seems not to have been known to his successors who studied polylogarithms. Another form of \eqref{Eq:EFA-Li2}  (see formula (9) p.\,252  of \cite{Abel}) has been obtained in 1826 by Abel in a short  note which remained unpublished until Abel's works 
were published in 1881. 
Then,   Hill rediscovered this functional equation in 1830 
({\it cf.}\,equation X) p.\,9 of \cite{Hill2}). Hill's functional equation is mentioned in  1840 Kummer's paper in Crelle's journal \cite{Kummer1} (equation 2 p.\,83 {\it loc.\,cit.}). 
  A few years later in 1846, Schaeffer established \eqref{Eq:EFA-Li2} in a paper again in Crelle's    (see formula 38 p.\,288 in  \cite{Schaeffer}). Surprinsingly, Schaeffer does not mention Kummer's papers dating back to six years before and published in the very same journal. In his paper \cite{Rogers} dating of 1907, Rogers introduced a modified version  of the classical dilogarithm $\l {2} $ which is now named after him (and denoted by $\mathcal R$ in this text, see \S\ref{Par:NotationDilogarithmicFunctions}) and  he established that the latter posses the nice property of satisfying a five-term relation equivalent to 
  \eqref{Eq:EFA-Li2} but without  any logarithmic term (see formula (11) p.\,173 of Rogers' paper or just below). According to Hardy, Rogers' functional equation for the dilogarithm 
was also independently
discovered  by Ramanujan (see p.\,14 and 
  p.\hspace{0.05cm}21 in \cite{Hardy}).  
\vspace{-0.2cm}
\begin{center}
$\star$
\end{center}

The presence of a  logarithmic right hand-side in \eqref{Eq:EFA-Li2} 
is not (at least) very  aesthetic and satisfying and, as first remarked by Rogers in \cite{Rogers}, actually can be removed by considering some modified versions of the classical Euler's dilogarithm $\l {2} $. 
 The first example is  {\bf Rogers' dilogarithm} $\boldsymbol{R}$  defined 
 for $x\in ]0,1[$ 
 by
  \begin{equation}
\label{Eq:Def-R}
 \qquad R(x)=\l {2} (x) + \frac{1}{2}\, {\rm Log}(x)\, {\rm Log}(1-x)-\frac{{}^{}\hspace{0.2cm}  \pi^2}{6}
 \, .
  \end{equation}
(Actually, this function is a  translation (by 
   substraction of $\pi^2/6$) of the original dilogarithmic function considered in 1907 by Rogers in \cite[\S1]{Rogers}, this in order that the RHS of \eqref{Eq:EFA-R} be zero, see \S\ref{Par:NotationDilogarithmicFunctions}).   It is a dilogarithmic function  which satisfies 
 for any  $x,y\in \mathbf R$ such that $0<x<y<1$:
 \begin{equation}
\label{Eq:EFA-R}
R (x) -
R (y)
-R \bigg(\,\frac{x}{y}\, \bigg)
-R \bigg(\,\frac{1-y}{1-x}\,\bigg)+
R\bigg(\,\frac{x(1-y)}{y(1-x)}\,\bigg)
=0\, .
\end{equation}

There also exist global versions of Rogers' dilogarithm and of equation \eqref{Eq:EFA-R}:  following Dupont {\rm \cite{Dupont}},
one can extend $R$
 as a continuous global map 
$\overline{R}$ 
on $\mathbf R$ by setting 
$$
\overline{R} (0)=-{\pi^2}/{6}\, , \quad 
\overline{R} (1)=0
\quad \mbox{ and }\quad 
\overline{ R}(x)=\begin{cases}
\, \, \, \, \, R\Big(\frac{1}{1-x}\Big)- \frac{{}^{}\, \, \pi^2}{6} \hspace{0.3cm} \mbox{ if } \, x<0;\, 
 \\ 
\, \, \,  \, \, R(x) \hspace{1.55cm} \mbox{ if } \, 0<x<1;\, 
 \\ 
\, -R\Big(\frac{1}{x}\Big)\hspace{1.5cm} \mbox{ if } \, x> 1\,.  \footnotemark
\end{cases}
\footnotetext{See also {\rm \cite[\S1.4]{DupontZickert}}, p.\,580 in {\rm \cite{GanglZagier}} 
or {\rm \cite[\hspace{-0.1cm}II.1.A]{Zagier}}.}
$$
It gives rise to a   continuous map $\overline{R}: \mathbf R \mathbf P^1=\mathbf R\cup \{\infty \}\rightarrow \mathbf R/((\pi^2/2)\mathbf Z)$  which is real-analytic except at $0,1$ and $\infty$ and which satisfies a global real version (modulo $\pi^2/2$) of \eqref{Eq:EFA-R}.\bk 

Even more popular than $\overline{R}$, there exists another global version of the dilogarithm, also real-valued but depending on a complex argument, 
which is due to S.\,Bloch and D.\,Wigner.  
For $z\in \mathbf C\setminus \{0,1\}$, 
it can be proved that the real quantity 
\begin{equation}
\label{Eq:Def-D}
  D(z)={\bf I}{\rm m}\Big( \l {2} (z)+{\rm Log}\, \lvert z\lvert
\, {\rm Log}( 1-z)\Big)
  \end{equation}
  is well-defined, {\it i.e.}\,is defined unambiguously regarding the choices of the determinations of the complex quanities $\l {2} (z)$ and ${\rm Log}( 1-z)$ appearing in \eqref{Eq:Def-D}.
  Hence setting $D(0)=D(1)=D(\infty)=0$, one gets the 
so-called {\bf Bloch-Wigner's dilogarithm}\footnote{The notation $D$  is the one used by Bloch in \cite{Bloch}, where the role of this function in relation with K-theory  appears for the first time. It seems that this function has also been considered by D. Wigner about the same time, hence its name (but we haven't  been able to find any  paper by Wigner regarding this). Actually, an equivalent form of the function $D$  was considered long before Bloch and Wigner, by Lobachevsky \cite{Lob} and latter by Schl\"afli \cite{Schlafli} 
and Coxeter \cite{Coxeter35,Coxeter36}. } 
$\boldsymbol{D}: \mathbf P^1\rightarrow \mathbf R$ which is a single-valued global  continuous function on the Riemann sphere, real-analytic on  $\mathbf C\setminus \{0,1\}$, and  satisfies the same relation as  \eqref{Eq:EFA-R}, but  globally that is, for any $x$ and $y$ in $\mathbf P^1$ such that the five rational arguments in \eqref{Eq:EFA-R} are well-defined. Actually, thanks to the one variable functional equation  $D(1/z)=-D(z)$ which is identically satisfied, the five-terms AFE  satisfied by $D$ can be written in the following nice geometric form:  for any 5-tuple $(p_i)_{i=1}^5\in (\mathbf P^1)^5$ of pairwise distinct points, one has 
  \begin{equation}
  \label{Eq:Efa-(Ab)-symmetric}
  \sum_{i=1}^5 (-1)^i D\Big([p_0,\ldots,\widehat{p_i}, \ldots, p_4]\Big)=0
  \end{equation}
  (where $[\cdot,\cdot ; \cdot,\cdot ]$ stands for the classical cross-ratio of four points on the projective line).


\begin{rem} 
%
Rogers'\,dilogarithm  (or more rigorously,  its extension $\overline{R}$ considered just above) 
 as well as Bloch-Wigner's function  $D$ are two modified real-valued versions of the classical dilogarithm 
 which both enjoy the nice properties of being globally defined and of satisfying the suitable global homogeneous version  
(that is, without any logarithmic second member) of the 5-terms FE \eqref{Eq:EFA-Li2}.
\sk 

  
 Actually, these two classical modified dilogarithms  both can be obtained from a seemingly more fundamental 
 `extended (or enhanced)  dilogarithm' 
 \begin{equation}
 \label{Eq:Neumann's-dilog}
 \widehat{R}: \widehat{\mathbf C}\rightarrow \mathbf C/
 \mathbf Z(2)
 \end{equation}
  first considered by Neumann in {\rm \cite{Neumann}} and slightly improved in {\rm \cite[\hspace{-0.1cm}\S2]{GoetteZickert}} (see also 
  {\rm \cite[\hspace{-0.1cm}II.1.B]{Zagier}}). 
 Here $\widehat{\mathbf C}$ stands for the universal abelian covering of the 2-punctured complex plane $\mathbf C\setminus\{0,1\}$\footnote{The space $\widehat{\mathbf C}$ can equally be 
 defined as the Riemann surface of the multivalued map $z\mapsto({\rm Log(z)}, {\rm Log}(1-z))$. 
 Very explicitely, it can be seen  as the set of pairs $(u,v)\in \mathbf C^2$ such that $e^u+e^v=1$ and 
up this identification, the covering $\widehat{\mathbf C}\rightarrow \mathbf C\setminus\{0,1\}$ is simply given by $(u,v)\mapsto e^u$.\label{Chat}}
  and we use  the notation $\mathbf Z(2)=(2i\pi)^2\mathbf Z$ 
 hence 
 $\widehat{R}$ has values in 
 $ \mathbf C/(4\pi^2\mathbf Z)
 = \mathbf R/(4\pi^2\mathbf Z)\oplus i \mathbf R$.  
 The map  $\widehat{R}$ is holomorphic, satisfies a global version of 
  \eqref{Eq:EFA-Li2} ({\it cf.}\,{\rm \cite[Lemma 2.2]{GoetteZickert}}{\rm )}.  Moreover  it  is the complete holomorphic extension of the complexification of the real-analytic function 
  $R\hspace{-0.05cm}:\,]0,1[\rightarrow{\mathbf R}$;
and  its imaginary part coincides with  $D$ up to the addition of an elementary function (see formula (13) in {\rm \cite{GanglZagier}}{\rm )}.
 \sk 
 
 %

In addition to all these nice properties  which in our view already justify 
 of considering $\widehat{R}$ as being more fundamental than $R$ or $D$, 
it turns out, and this is the very reason why it has been considered
 by several authors, that $\widehat{R}$  can be used to represent some abstract conceptual objects, namely 
the universal second Cheeger-Chern-Simmons characteristic class or,  
equivalently,  of the second Beilinson's regulator of $K$-theory,  see {\it e.g.}\,{\rm \cite[Chap.\hspace{0.05cm}10]{DupontBook}} or the Introduction (and in particular diagram (1.13)) of {\rm \cite{Zickert}}.
\mk 

All the facts mentioned above support our opinion that the theory of functional equations of polylogarithms must concern complex analytic versions of  such equations, 
which is in perfect accordance with our approach of these equations, through the geometry of webs.
%
%
%
\end{rem}

\subsubsection{\bf The classical higher polylogarithms.}\hspace{-0.4cm} 
Considering the classical development in series \eqref{Eq:Devel-Series-Log} and 
 \eqref{Eq:Devel-Series-Li2}, 
it is natural  to define, for any positive integer $n$,  a function $ \boldsymbol{ \l {n} }$, named the {\bf $\boldsymbol{n}$-th} (or the {\bf weight $\boldsymbol{n}$}) {\bf  polylogarithm},   by setting
\begin{equation}
\label{Eq:Devel-Series-Lin}
\l {n} (z)=\sum_{k=1}^{+\infty}\frac{z^k}{k^n}\, 
\end{equation}
for any $z\in \mathbf C$ such that $\lvert z \lvert<1$. Then $\l {1} (z)= - {\rm Log}(1-z)$, and for $n=2$ one gets the dilogarithm.\footnote{As far we know,
the first appearance  of the $n$-th  polylogarithm for  $n$ arbitrary is in Landen's 1759 paper (cited in footnote \ref{footLanden}) 
where the development in series \eqref{Eq:Devel-Series-Lin} as well as  the integral representation \eqref{Eq:IntRep-Lin} are explicitly written down (in the second section, page 554). As for the first text essentially devoted to polylogarithms of any order, it seems to be \href{https://www-history.mcs.st-andrews.ac.uk/Biographies/Spence.html}{W. Spence}'s essay \cite{Spence}, published in 1809 but long ignored by 19th century mathematicians. 
 For a thorough discussion of Spence work on polylogarithms and especially on its  1809 essay, the reader can consult 
 \href{https://doi.org/10.1016/j.hm.2013.06.002}{\it Polylogarithms, functional equations and more: the elusive essays of William Spence} (2013) by D. Craik.
}
\sk


The higher polylogarithms 
 are generalizations of the logarithm, which satisfy similar properties to those listed above.  Indeed, for $n>1$, the $n$-th polylogarithm $\l {n} $: 
\begin{itemize}
\item admits the following integral representation \begin{equation}
\label{Eq:IntRep-Lin}
\l {{ n }} (z)=  \int_0^z \frac{  \l {{ n-1 }}  (u)  }{u}  du\, ; 
\qquad \quad {}^{}
\vspace{-0.15cm}
\end{equation} 
\item extends as a multivalued holomorphic function on $\mathbf C\setminus \{0,1\}$, with monodromy around the origin and 1 along a small circle oriented in the direct order given by \vspace{-0.25cm}
\begin{equation}
\label{Eq:Monod-Li1}
\mathcal M_0 \Big( 
\l { {n} } (z)
\Big)= \l { {n} } (z) \qquad \mbox{and}\qquad 
\mathcal M_1 \Big( 
\l {  {n} } (z)
\Big)= \l { {n} } (z) - 2i\pi  \frac{{}^{}\, \, \, {\Big(\rm Log} \, z \Big)^{n-1}    }{(n-1)!}
 \vspace{0.15cm}
\end{equation}
\item satisfies many functional equations and in particular the following ones  in one variable: 
\begin{align}
\label{Eq:EFA-Lin}
\l {n} \Big(z^r\Big)  = & \,  r^{n-1} \sum_{w^r=1} \l { n } (wz)  &&  \mbox{ for } z 
\in \mathbf C \mbox{ such that }
\lvert z\lvert <1  
\\ 
\l {n} (z) +  (-1)^n\, 
\l {n} \Big(z^{-1}\Big)
  = &  - \frac{(2i\pi)^n}{n!} \boldsymbol{B}_n\bigg( \frac{{\rm Log}(z)}{2i\pi}\bigg)  &&\mbox{ if }\,  z\in \mathbf C\setminus [0,+\infty[  
  \nonumber
\end{align}
(where $\boldsymbol{B}_n(x)$ stands for the $n$-th  \href{https://en.wikipedia.org/wiki/Bernoulli_polynomials}{Bernoulli polynomial}).
\end{itemize}

The two previous FEs are satisfied uniformly in $n$ by the polylogarithms. This contrasts with the case of FE in at least two variables: if many such equations are known (and will be given in Section \ref{S:EFA-polylog} hereafter), it is only for  small values of $n$, typically for $n\leq 5$.  \sk

A major problem about polylogarithms is to know whether or not they all satisfy FEs in several variables.  
 The actual record for a several variables FE satisfied by (a modified version of) the $n$-th polylogarithm is $n=7$ and has been discovered quite recently by  
Gangl \cite{Gangl} with the help of a computer. It is an AFE  in 2 variables 
with no less than 274  heptalogarithmic terms. 
  \sk 
\begin{center}
$\star$
\end{center}

 Any known functional equation involving $\l {n} $ can be reduced to a FE  of the following form 
\begin{equation}
\label{Eq-FE-Lin}
\sum_{j=1}^{N} c_j\, \l {n} \Big( u_j(x) \Big)= 
{\rm P}_{\hspace{-0.08cm}< n }(x)\, 
\end{equation}
where:  
\vspace{0.1cm}

-- $N$ is a certain positive integer; 
\vspace{0.1cm}

-- $x$ denotes a point of a certain non-void domain in $\mathbf C^s$, 
 for some  $s\in \mathbf N_{>0}$; 
 \vspace{0.1cm}

--  the $c_j$'s are non-zero  (rational) constants, for any $j=1,\ldots,N$;
\vspace{0.1cm}

--  the $u_j$'s are non-constant rational functions in $s$ variables 
(with rational coefficients)
for any $j$; 
\vspace{0.1cm}

--  ${\rm P}_{\hspace{-0.1cm}< n }(x)$ stands for a function of the form $P\Big( {\rm Li}_{ m_1 } (v_1(x)), \ldots, 
{\rm Li}_ { m_M } (v_M(x)) \Big)$  where $P$ is a complex 

\textcolor{white}{-} 
polynomial in $M$ variables, the $v_k$'s are of the same type as the $u_j$'s, and all the $m_k$'s are positive 

\textcolor{white}{-} integers strictly less than $n$.
\vspace{0.1cm}

\begin{rem} 
By means of composition and/or changes of variables, it is easy to construct many functional equations satisfied by $\l {n} $ starting from 
one (or even several) FEs of the form \eqref{Eq-FE-Lin}. A folkloric conjecture 
of the field (see conjecture $(i)$ in {\rm \cite [\S7]{Lewin}} for instance) 
is that any FE satisfied by $\l {n} $ can be obtained in this way.\vspace{0.1cm}

Almost nothing is known about this conjecture for $n\geq 2$.\footnote{For the logarithm case ({\it i.e.}\,when $n=1$), this is a folkloric result but we are not
aware of any place where a formal proof  is given in this case.}. The unique available result  is a partial one due to Wojtkowiak which concerns  dilogarithmic EFAs in  one variable satisfied (see 
\S\label{SubPar:Accessibility-Dilog} further).
\footnote{This was true when we were writing these lines. Since then, the 
case of dilogarithmic EFAs with rational function  in an arbitrary number of variables as arguments
 has been essentially settled by R. de Jeu in  \href{https://arxiv.org/pdf/2007.11014.pdf}{{\it Describing all multivariable functional equations of dilogarithms} (preprint arXiv: 2007.11014)}.}.
\end{rem}

As it will be clear from the classical (or not) examples that we will consider in 
Section \ref{S:EFA-dilog} below, the known FEs of the form \eqref{Eq-FE-Lin} satisfied by polylogarithms become more and more complicated as $n$ becomes bigger and because of the multivaluedness of these functions,  are not easy to handle:  it is necessary to specify on which domain in $\mathbf C^s$ and for which determination of the polylogarithms appearing in  it, the considered FE 
\eqref{Eq-FE-Lin} is truly satisfied.   This is a source of many complications and in order to bypass the technical difficulties this creates, several authors have introduced
(several distinct) modified versions of the polylogarithm ${ \l {n} } $ (for 
$n$ arbitrary), which enjoy the properties: \sk 
\begin{itemize}  
\item  of being defined as a real-valued univalued functions,  on the whole real line $\mathbf R$ (or on the whole Riemann sphere $\mathbf P^1$), continuous everywhere and even real-analytic
  outside $\{0,1,\infty\}$;
\item
 of satisfying clean versions of the FE of the form \eqref{Eq-FE-Lin} satisfied by the classical polylogarithms (where `{\it clean}' means: with a trivial, that is constant, second member, instead of the usually non-trivial polylogarithmic expression 
 ${\rm P}_{\hspace{-0.1cm}< n }(x)$ appearing in the RHS of \eqref{Eq-FE-Lin}).
\end{itemize}

In the paragraph below, we review the  modified polylogarithms  encountered in the existing literature.
 In the two following paragraphs after it, we introduce some holomorphic modified polyogarithms which are more suitable to deal with holomorphic AFEs.

\subsubsection{\bf Modified higher polylogarithms (short review).}\hspace{-0.4cm} 
\label{SS:Modified-Higher-Polylogarithms}
There are several versions of such modified polylogarithms. Our goal in this paragraph is to review them and to state some of their properties, essentially in regards with the AFE they satisfy.\mk

  It will be useful to use the following notation and terminology, for any $n\in \mathbf N_{>0}$: 
\begin{itemize}  
\item 
  we set $\boldsymbol{\mathfrak R}_n(z)={\rm Re}\Big((-1)^{\lfloor{n/2}+1\rfloor}i^{n+1}z\Big)$ for any $z\in \mathbf C$; more concretely, $\boldsymbol{\mathfrak R}_n$ stands for the real part $ {\bf R}{\rm e}$ if $n$ is odd and for the imaginary part $ {\bf I}{\rm m}$ if $n$ is even;
\item  by definition, a function $F$ of a real or complex argument $a$ 
 satisfies the {\bf `inversion formula at the order $\boldsymbol{n}$'} if  the relation $F(a^{-1})=(-1)^{n-1} F(a)$ holds true for $a$ generic (in $\mathbf R$ or $\mathbf C$ depending on the definition domain of the considered function).
\end{itemize}

There are two distinct modified dilogarithmic functions: Rogers dilogarithm $R$,  and Bloch-Wigner function $D$. The latter is a univalued global function on the Riemann sphere $\mathbf P^1$  while the former should rather be seen as a function of a real variable. 
Modified higher polylogarithms have been constructed for each of these two kinds. We first discuss the  higher order generalizations of  Bloch-Wigner function  and then we turn to the modified higher  polylogarithms of Rogers' type.

\paragraph{}  
The problem of constructing higher-order generalizations of Bloch-Wigner's function is first considered in \cite{Ramakrishnan} where Ramakrishnan establishes that single valued higher-order versions of $D$ exist.   These 
single-valued polylogarithms have been explicited  a few year later by  Zagier in \cite{Zagier1990}  who relates them in \cite{Zagier1991} 
 to some others  (only differing from Ramakrinshan ones by the addition of a power of the logarithm) obtained independently by Wojtkowiak:\mk 
\begin{itemize}
\item For $n\geq 2$, the {\it Bloch-Wigner-Ramakrishnan-Zagier's} polylogarithm  is the function $D_n$ 
 defined by 
%
$$
D_n(z)={\boldsymbol{\mathfrak R}}_n\Bigg( \sum_{s=0}^{n-1} \frac{(-1)^s }{s!} \Big({\rm Log}\,\lvert z\lvert\Big)^s \l { {n-s} } (z)
- \frac{(-1)^n}{(2n)!}\Big({\rm Log}\,\lvert z\lvert\Big)^n
\Bigg)
$$
for $\lvert z \lvert \leq 1$, and extended to $\mathbf C\setminus\{0,1\}$ by requiring that  the $n$-th order inversion formula  be satisfied at any point. 
 The function $D_n$ is real-analytic, extends continuously at $1$ and has a 
 singularity of type $\big({\rm Log}\,\lvert z\lvert\big)^n$ at $0$ (hence at $\infty$ as well thanks to the inversion formula); 
\mk 
\item   In \cite[\S7]{Zagier1991}, Zagier observes that  the logarithmic singularities of $D_n$ can be removed by adding a suitable logarithmic expression to it.  
 One recovers by this way (up to a minus sign) the modified  polylogarithms considered by Wojtkowiak  in \cite{Wojtkowiak1989} defined for any $n\geq 2$,   by 
 \begin{equation}
 \label{Eq:LWojtkowiak}
  \L{}_n(z)
= 
\boldsymbol{\mathfrak R}_n
\Bigg( \sum_{s=0}^{n-1} \frac{(-1)^s }{s!} ({\rm Log}\,\lvert z\lvert)^s \l { {n-s} } (z)
- \frac{(-1)^n}{n!}\Big({\rm Log}\,\lvert z\lvert \Big)^{n-1}{\rm Log}\,\lvert 1-z\, \lvert\, 
\Bigg)
\end{equation}
 for any $z\in \mathbf C\setminus \{0,1\}$.\footnote{Note that when $n$ is odd, 
 $\boldsymbol{\mathfrak R}_n={\bf I}{\rm m}$ thus the real logarithmic term 
 $- \frac{(-1)^n}{n!}({\rm Log}\,\lvert z\lvert )^{n-1}{\rm Log}\,\lvert 1-z\, \lvert$ 
 does not play any role 
  hence can be suppressed in the definition of $\L{}_n$, as it is the case in \cite{Wojtkowiak}.} 
Wojtkowiak's polylogarithm $\L{}_n$ is real-analytic and, as Bloch-Wigner dilogarithm,  extends as a continuous function on $\mathbf P^1$ by setting 
$\L{}_n(0)=\L{}_n(\infty)=0$ and $\L{}_n(1)=\l { n} (1)=\zeta(n)$ if $n$ is odd, and 
 $\L{}_n(1)=0$  otherwise.  
\mk 
\item Also in \cite{Zagier1991}, Zagier introduces another 
 real single-valued version of 
 $\l {n} $ 
  by setting 
$$
\mathcal L_n(z)=\boldsymbol{\mathfrak R}_n\Bigg( \sum_{s=0}^{n-1} \frac{2^s B_s}{s!} \big({\rm Log}\,\lvert z\lvert\Big)^s \l { {n-s} } (z)
\Bigg)
$$  for $z\in \mathbf C\setminus \{0,1\}$,
where the $B_s$'s stand for  the \href{https://en.wikipedia.org/wiki/Bernoulli_number}{Bernoulli numbers} ($B_0=1$, $B_1=-1/2$, $B_2=1/6$, $B_3=0$, $B_4=-1/30$, etc.). 
This function can be seen too as as a weight $n$ generalization of the classical Bloch-Wigner function $D$ (which coincides with $\mathcal L_2$).\footnote{As mentioned by Goncharov in \cite[\S5]{Goncharov2005}, there is also a modification due to  A.\,Levin of Zagier's $n$-th polylogarithm $\mathcal L_n$ which enjoys the property of admiting a global integral representation, over the whole complex projective space $\mathbf P^{n-1}$,  of an explicit real logarithmic differential form.}
\mk 
\end{itemize}

The previous functions $D_n,  \L{}_n$ and $\mathcal L_n$ are functions of a  parameter $z\in \mathbf P^1\setminus \{0,1,\infty\}$ hence as such, have to be seen as generalizations of Bloch-Wigner's function $D$  to any arbitrary order $n$. 
 \mk

    
    In \cite[(16)]{Lewin1986}, Lewin 
gave a recursive definition     
of  some modified polylogarithms  $L_n$ which depend on a real argument and which can  be seen as higher-order generalizations  of 
  Rogers' dilogarithm $R$. These functions have been made explicit by Zagier in \cite{Zagier1991} who proved that for $n\geq 2$, the value of $L_n$ at  $x\in ]-1,1[$ is given by 
  \begin{equation}
  \label{Eq:Lewin-Dilog}
L_n(x) =\sum_{s=0}^{n-1} \frac{(-1)^s}{s!} \Big({\rm Log}\, \lvert x \lvert \Big)^s
 \l { {n-s} } (x) -\frac{(-1)^n}{n!}{\rm Log}\, \lvert 1-x \lvert \cdot \big({\rm Log}\, \lvert x \lvert \Big)^{n-1}\, .
\end{equation}
 By requiring that the relation $L_n(1/x)=(-1)^{n-1}L_n(x)$ holds true for any $x\in \mathbf R$ and setting  $L_n(\infty)=L_n(0)=0$, one gets a function defined on the whole projective line $L_n:  \mathbf P^1_{\hspace{-0.07cm}\mathbf R}\rightarrow \mathbf R$.    For  $n=2$,  one has $L_2(x)=\l { 2 } (x)+\frac{1}{2}{\rm Log}(x){\rm Log}(1-x)$ when $x$ belongs to $]0,1[$  hence on this interval and up to addition of $\pi^2/6$, 
Lewin's bilogarithm   $L_2$ coincides with Rogers' one $R$.\mk 
 
   For $n$ odd, one gets a continuous function $L_n: \mathbf P^1_{\mathbf R} \rightarrow  [-\zeta(n)\, , \, \zeta(n)\,]$ which is real-analytic, except at $-1, 0,1$ and $\infty$. In this case,  
$L_n$ agrees with the restriction to the real axis of Wojtkowiak's modified polylogarithm $\L{}_n$.  For $n$ even,  $L_n$ is real-analytic on $\mathbf R\setminus \{\pm 1\}$ but has a  gap at $\pm 1$ equal to  $2\l {n} (\pm 1)$.  It follows that, similarly to Rogers' dilogarithm,  $L_n$ induces 
a continuous map 
 $\overline{L}_n : \mathbf P^1_{\hspace{-0.1cm}\mathbf R}\rightarrow \mathbf R/\Big(  \frac{2\zeta(n)}{2^{n-1}-1} \mathbf Z\Big)$.\mk

Given $n\geq 2$, one interesting property of $L_n$  ({\it cf.}\,\cite[\hspace{-0.1cm}(19)]{Lewin1986} or \cite[\hspace{-0.1cm}p.\,412]{Zagier1991}) is that its derivative admits a quite simple explicit  expression involving only rational terms and powers of the logarithm. Indeed,  one has 
$
L_n'(x)=\frac{(-1)^{n-1}}{n(n-2)!}\big({\rm Log}\,x\big)^{n-2}\Big( 
{{\rm Log}(x)}/{(1-x)} +
{{\rm Log}(1-x)}/{x}
\Big)
$ for any $x\in ]0,1[$. 
\bk

Since we are most interested in holomorphic objects, we have to mention here 
 that Lewin's polylogarithms 
 naturally  give rise  to 
some 
modified holomorphic  polylogarithms, the so-called 
{\it `enhanced polylogarithms'}, denoted by $\mathscr{L}_n$ here, 
first considered by 
Gangl and Zagier in \cite[\S4]{GanglZagier} and  studied again recently  by 
Zickert in \cite[\S5.1]{Zickert}. \sk 

 Let $\Omega$ be the domain obtained by removing from the complex plane $\mathbf C$ the two cuts $]-\infty,0]$ and $[1,+\infty[$. Then the value at  
$z\in \Omega$ of the $n$-th order enhanced polylogarithm is given by 
$$
\mathscr{L}_n(z)=\sum_{s=0}^{n-1}\frac{(-1)^s}{s!} \Big({\rm Log}\,z\Big)^s \l { {n-s} } (z) 
+\frac{(-1)^n}{n!} \Big({\rm Log}\, z\Big)^{n-1}\l { 1 } (z)\, .
$$
%
%

Of course, one has $\mathscr{L}_n(x)=L_n(x)$  for every $x\in ]0,1[$ hence $\mathscr{L}_n$ can also be defined as the holomorphic extension to $\Omega$ of the complexification of the (real-analytic) restriction of 
 $ L_n$ to $]0,1[$.\mk 
 
 By analytic continuation,  $\mathscr{L}_n$ gives rise to a global holomorphic map  $\widehat{\mathscr{L}}_n$ defined on the Riemann surface  $\widehat{\mathbf C}$  considered above regarding the case $n=2$ (see footnote \footref{Chat}).  Considered as a multivalued function on $\mathbf C$ with ramification at $0$ and $1$, $\widehat{\mathscr{L}}_n$ has monodromy in $\frac{(2i\pi)^n}{(n-1)!}\mathbf Z$
  (see \cite[Thm.\,2.2]{Zickert})  hence induces a well defined holomorphic map
   denoted in the same way
   $\widehat{\mathscr{L}}_n: \mathbf C\setminus \{0,1\}\rightarrow \mathbf C/(\mathbf Z(n)/(n-1)!)$ generalizing Neumann's dilogarithm \eqref{Eq:Neumann's-dilog} to any $n\geq 2$.

\begin{rem}
%
%
%
 Zagier's polylogarithms appear as  the most interesting univalued versions  of the classical polylogarithms, not only thanks to the nice properties they satisfy 
 (Theorem \ref{T:FE-Zagier-Polylogs} being by far the most important)
 but  also because they have a motivic origin ({\it cf.}\,{\rm \cite[p.\,415]{Zagier1991}} or  {\rm \cite[\S1.5]{BeilinsonDeligne}} for more details). 
Actually, Zagier's  polylog $\mathcal L_n$ can be explicited  in terms of the 
real or imaginary part (depending whether $n$ is odd or not) of the holomorphic enhanced polylogarithm  $\widehat{\mathscr{L}}_n$ and it has been speculated  that the latter could give an explicit description of  a complex motivic regulator  (see 	\S1.2.1 and \S1.3.1 in {\rm \cite{Zickert})} which would be the most fundamental object in all this story. 
 More generally, many of the recent works on polylogarithms are within the framework of motives.  
 One may think that a motivic  approach possibly could  lead to a better conceptual understanding 
  of the polylogarithmic webs and their ARs.
\end{rem}

  \begin{center}
  \vspace{-0.1cm}
  $\star$
  \end{center}

 The modified polylogarithms mentioned above enjoy some nice properties similar to those satisfied by Bloch-Wigner's or Rogers' dilogarithms: they are globally defined (either on the Riemann sphere, or on the real projective line $\mathbf P_{\hspace{-0.03cm}\mathbf R}^1$ or on the Riemann surface $\widehat{\mathbf C}$), they are real and then single-valued, or holomorphic with values well-defined up to a rational multiple of a power of $2i\pi$, they satisfy simple and explicit first order differential equations (which are expressed in terms of the logarithm and of the modified polylogarithms of the same kind but of a lower order).  \mk 
 
 But what is in our view the most interesting property these functions verify is that they satisfy the clean version of any FE of the type \eqref{Eq-FE-Lin} satisfied by the classical $n$-th order polylogarithm $\l {n} $. 
  More precisely, for any $n\geq 2$,  assuming that  $\sum_{j=1}^{N} c_j\, \l {n} \Big( u_j(x) \Big)= 
{\rm P}_{\hspace{-0.08cm}< n }(x)$ holds identically for $x$ in a domain $U$ of $\mathbf C^s$ 
 and if $\boldsymbol{L}_n$ stands for one of the modified polylogarithms considered above, it is expected that $\sum_{j=1}^{N} c_j\, \boldsymbol{L}_n\big( u_j \big)= {\rm cst.}$ on a Zariski open set  of $\mathbf C^s$ intersecting  $U$.  These  are kind of  folkloric results  in what concerns polylogarithms and   can be found, sometimes in some particular  case ({\it e.g.}\,one-variable case) in several references, see the following table: \sk


\begin{table}[!h]
\begin{center}
\begin{tabular}{|c|c|c|}
\hline
{\bf Modified polylogarithm(s)}
& {\bf Reference(s)}    &
{\bf Remark(s)}
    \\
\hline
$\L{}_n$, $L_n$  &  \cite[Proposition 1]{Zagier1991}  &    ${}^{}$ \hspace{-0.8cm}One-variable case       \\
\hline
$D_n$  &  \cite[Proposition 2]{Zagier1991}  &   ${}^{}$ \hspace{-0.8cm}  One-variable case         \\
\hline
$\mathcal L_n$ &\begin{tabular}{l}
 \cite[Proposition 3]{Zagier1991}\\
  \cite[Th\'eor\`eme 2]{Oesterle}
\end{tabular}
   &      \begin{tabular}{l}
One-variable case \\
Several-variables case
\end{tabular}     \\
\hline
$\mathscr{L}_n$  &  \cite[Proposition 6.1]{Zickert}  &    Several-variables case         \\
\hline
\end{tabular}
\end{center}
\end{table}


\paragraph{}  
We think it is interesting to state a precise result regarding the clean AFE satisfied by modified polylogarithms. We deal with the most popular such polylog, namely Zagier's version $\mathcal L_n$.\footnote{In view of the main purpose of this text, namely the construction of holomorphic AFEs, it would have been more natural and  relevant to consider the holomorphic case dealt with in \cite{Zickert} (compare the proposition and the corollary of the sixth section of Zickert's paper with Theorem \ref{T:FE-Zagier-Polylogs} and  Theorem \ref{T:FE-holomorphic} below). But we became aware of Zickert's text quite late when the writing of the present memoir was already well advanced so we have decided to advertise only the case of $\mathcal L_n$ which, moreover, is by far the most popular.\label{Why}} \mk 

First we recall  the formalism of \S\ref{Par:II-AR}\hspace{0.2cm} relatively to the logarithmic pair $(\mathbf P^1,\mathfrak B_1)$: a word $\boldsymbol{w}=w_1\cdots w_n$ in the  letters $0$ and $1$ 
 corresponds to the logarithmic $n$-differential form $d\log(z-w_1)\otimes \cdots \otimes d\log(z-w_n)
$ viewed as an  element of the $n$-tensorial product of $\boldsymbol{H}_{\mathfrak B_1}=H^0(\mathbf P^1,\Omega^1_{\mathbf P^1}({\rm log} \mathfrak B_1))
$.  Now given  an 
 irreducible algebraic variety $V$ with function field $\mathbf K=\mathbf C(V)$, we define $\Omega^1_{\mathbf K}(\log)
 $ as  the space of rational 1-forms on $V$ spanned by the logarithmic differentials $d\log f=df/f$ for $f\in \mathbf K$. Given  $u\in \mathbf K$ and a word $\boldsymbol{w}$
 as above, we 
 denote by $u^*( \boldsymbol{w})$ the element 
$
\otimes_{i=1}^n u^*(d\log(z-w_i))=
d\log(u-w_1)\otimes \cdots \otimes d\log(u-w_n)$  of  $\Omega^1_{\mathbf K}(\log)^{\otimes n}$. 
\mk

Now assume that $n\geq 2$ and let $\{(c_i,u_i)\}_{i\in I}$ be a finite collection of pairs of complex numbers $c_i$ and of non-constant rational functions $u_i\in \mathbf K$.  Denote by $\boldsymbol{H}_I$ the vector subspace (of finite dimension) of $\Omega^1_{\mathbf K}(\log)$ spanned by the logarithmic 1-forms $ d\log(u_i)=du_i/u_i$ and $-d\log(1-u_i)=du_i/(1-u_i)$ for $i$ in $I$. Then for any  word $\boldsymbol{w}$ of length $n$ as above, one has $u_i^*(\boldsymbol{w})\in (\boldsymbol{H}_I)^{\otimes n}\subset \Omega_{\mathbf K}^1(\log)^{\otimes n}$.

\begin{thm}
\label{T:FE-Zagier-Polylogs}
{\bf 1.} The following two assertions are equivalent: 
\begin{enumerate}
\item[$(a).$]  a functional equation of the form  $\sum_{i\in I} c_i \, { {\bf L}{\rm i}_{ n } }  (u_i)={\rm P}_{\hspace{-0.06cm} <n}$ holds true;
\item[$(b).$]   in the tensorial  
vector space $(\boldsymbol{H}_I)^{\otimes n}$,  one has  $
\sum_{i \in I} c_i\cdot 
 u_i^*\Big[ 0^{n-2}(01-10)  \Big]=0
$.\sk 
\end{enumerate}
{\bf 2.} If $(a)$ and $(b)$ are satisfied, then 
$ \sum_{i \in I} c_i \,({\rm Log}\lvert u_i\lvert)^k \mathcal {L}_{n-k}(u_i)$ 
is constant  for any $k=0,\ldots,n-2$.
\end{thm}

In the case when $V=\mathbf P^1$ (and therefore $\mathbf K=\mathbf C(t)$), the second part of this result (moreover only for $s=0$) has been obtained first by 
 Zagier \cite[Prop.\,3]{Zagier1991} who proved it by induction via derivation with respect to the variable $t$ and its conjugate. 
  In full generality, the statement above appears in Oesterl\'e's survey \cite[\S4.2]{Oesterle} 
 where only the main lines of a proof are sketched, the arguments used to establish {\bf 2.} being variants of the one used in \cite{Zagier1991} in  the case of $\mathbf P^1$.
\mk 

This theorem is very helpful to construct ARs for polylogarithmic webs
 since it indicates that 
any FE 
$\sum_i c_i \l {n} (u_i)={\rm P}_{{}\hspace{-0.07cm}<n}$ 
of the form \eqref{Eq-FE-Lin}  satisfied by the $n$th-order classical polylogarithm $\l {n} $ gives rise to $n-1$ linearly independent `global and real abelian relations' for the polylogarithmic web $\boldsymbol{\mathcal W}(u_i)$ defined by the $u_i$'s.    When the latter have real coefficients (as it is the case for all the known polylogarithmic functional equations), since all the functions ${\mathcal L}_{n,k}: z\mapsto ({\rm Log}\vert z\lvert)^k {\mathcal L}_{n-k}(z)$ (for $k=0,\ldots,n-2$) are real-analytic, one can {\it (i)} 
  first restrict the (`real-valued') AFE 
$\sum_i c_i  \mathcal L_{n,k}(u_i)={\rm cst.}$ 
 to the reals;  and then {\it (ii)} consider the 
  complexification of the latter real-analytic AFE, in order to get 
   a genuine (that is a `holomorphic') AR for the aforementioned web $\boldsymbol{\mathcal W}(u_i)$. However, this procedure is obviously not injective since 
   ${\mathcal L}_{n,k}$   
  vanishes identically on 
  ${\mathbf R}$
  as soon as $n-k$ is even, hence one  gets a trivial AFE already after  step {\it (i)} mentioned above.  \mk 
   
   The preceding considerations show that, 
if  Theorem   \ref{T:FE-Zagier-Polylogs} is interesting to construct ARs for polylogarithmic webs, its real nature makes that it is not suitable for this purpose as such. Moreover,  no proof of it is available in the existing literature. For these reasons, we give (and prove)  below a   holomorphic (but local) version of this theorem which is better suited for dealing with polylogarithmic webs.


\subsubsection{\bf A holomorphic version of Theorem \ref{T:FE-Zagier-Polylogs}.}
\label{Par:HolomSymbolicRA}
We use the same notation as above: 
$\mathbf K=\mathbf C(V)$ is the function field  of $V$  which is irreducible 
and $\{(c_i,u_i)\}_{i \in I}$ stands for  a finite collection as in Theorem \ref{T:FE-Zagier-Polylogs}. We set 
$\boldsymbol{{u}}_I=(u_i)_{i\in I}$ and 
$\boldsymbol{H}_I=\Big\langle \,
du_i/u_i\, , \, du_i/(1-u_i)
\, 
 \lvert \, i \in I \,  \Big\rangle \subset \Omega_{\mathbf K}^1(\log)$.
 \mk 

We fix $x\in V$,  a regular point for the web $\boldsymbol{\mathcal W}(\boldsymbol{{u}}_I)$ defined by the $u_i$'s
\footnote{By this, we mean that  $v$ is a regular point for all the $u_i$'s  and that the 2-form $du_i\wedge du_{i'}$ does not vanish at $v$ as soon as 
the rational functions $u_i$ and $u_{i'}$ define two distinct foliations on $V$.} and for every   $i\in I$, we  set 
$$x_i=u_i(x) \in \mathbf P^1\, . 
$$ For any word $\boldsymbol{w}=a_1\cdots a_n$ in the two letters 0 and 1, one has 
\begin{equation}
\label{Eq:togodoc}
 {\rm II}_x^n \Big(u_i^*(\boldsymbol{w}) \Big)= u_i^*\Big({\rm II}_{x_i}^n (\boldsymbol{w})\Big)=
u_i^*(L_{\boldsymbol{w}}^i)=
L_{\boldsymbol{w}}^i(u_i)
\end{equation}
where $L_{\boldsymbol{w}}^i=L_{\boldsymbol{w},x_i}$ stands for the (germ at $x_i$ of the) iterated integral 
$$
L_{\boldsymbol{w}}^i(\bullet)= \int_{x_i}^\bullet 
\nu_{a_n}\otimes \cdots \otimes
  \nu_{a_1}
\in  \mathcal I \hspace{-0.07cm} \mathcal I_{\hspace{-0.07cm}x_i}^n
$$
with 
$\nu_0=du/u$ and $\nu_1=
du/(1-u)$.  To state the result we have in mind, we need to consider special polylogarithmic iterated integrals: 
for any $i\in I$, any $n\geq 2$ and $k\in \{ 0,\ldots,n-2\}$, one sets 
$$
\mathfrak L_n^i=\frac{1}{n}\,L_{0^{n-2}(01-10) }^i\qquad \mbox{ and } \qquad 
\mathfrak L_{n,k}^i=
\frac{k!}{n-k} \,
L_{0^k\shuffle 0^{n-k-2}(01-10) }^i= (k!) \,
L_{0^k}^i  \cdot \mathfrak L_{n-k}^ i
$$
(more explicitly, one has $\mathfrak L_{n,k}^i(z)={\rm Log}(z/x_i)^k \hspace{0.02cm}
\mathfrak L_{n-k}^i(z)$ for $z$ sufficiently close to $x_i$ in $\mathbf P^1$).

\begin{thm}
\label{T:FE-holomorphic}
{\bf 1.} Given a non-zero $d$-tuple $(c_i)_{i=1}^d\in \mathbf C^d$, the 
 following assertions are equivalent: 
\begin{enumerate}
\item[$(a).$]  a functional equation of the form  $\sum_{i\in I} c_i \, { {\bf L}{\rm i}_{ n } }  (u_i)={\rm P}_{\hspace{-0.06cm} <n}$ holds true;
\sk 
\item[$(b).$]     one has  $
\sum_{i \in I} c_i\cdot 
 u_i^*\Big[ 0^{n-2}(01-10)  \Big]=0
$ in $(\boldsymbol{H}_I)^{\otimes n}$;
\item[$(c).$]     the sum 
$ \sum_{i \in I} c_i \,{\mathfrak{L}}_{n}^i(u_i)$ vanishes identically on a neighbourhood of $x$ in $V$.
\end{enumerate}
{\bf 2.} If the conditions of\hspace{0.1cm}{\bf 1.}\,are satisfied, then 
$ \sum_{i \in I} c_i \,{\mathfrak L}_{n,k}^i(u_i)\equiv 0$ 
 for any $k=0,\ldots,n-2$ and consequently these $n-1$ 
 AFE's 
 provide as many  linearly independant abelian relations  for 
  $\boldsymbol{\mathcal W}(\boldsymbol{{u}}_I)$. 
\end{thm}
This result, which is a holomorphic counterpart of Theorem  \ref{T:FE-Zagier-Polylogs}, is not difficult to prove and is certainly well-known by the specialists of the subject\footnote{See \cite[\S6]{Zickert} and also 
footnote \footref{Why} above for some comments.}. But since we aren't aware of any written proof of it  in the literature, we describe below the main lines of an algebraic proof.\mk 
 
{\it Proof (sketch).}  Before entering into the proof, let us discuss a bit the point $(a)$: it means that a FE of the form \eqref{Eq-FE-Lin} is satisfied locally at some point of $V$. By analytic continuation, this implies that this point can be assumed to be our chosen base-point $x$. \mk

To prove {\bf 1.}, we are going to use standard material about the shuffle product 
 that we recall here (see p.\,\pageref{Page:shuffle} above and also \cite[\S5.4]{DGR} for some details and references). Given 
a complex vector space $\boldsymbol{H}$ and  $n\geq 2$,  denote by ${\rm Im}(\shuffle)^n$ the  subspace of 
$\boldsymbol{H}^{\otimes n}$ 
spanned by all the non-trivial shuffle products (that is by shuffles $v_1\shuffle v_2$ where $v_i\in \boldsymbol{H}^{\otimes n_i}$ with $n_i>0$ for $i=1,2$ such that $n_1+n_2=n$).  One defines inductively an endomorphism  $\Pi_n\in {\rm End}(\boldsymbol{H}^{\otimes n})$ for $n\geq 1$,  by setting $\Pi_1={\rm Id}_{\boldsymbol{H}}$  and 
$$
\Pi_n\Big(a_1\otimes \cdots \otimes a_n\Big)=\frac{n-1}{n} \bigg[\,a_1\otimes \Pi_{n-1}\Big(a_2\otimes\cdots \otimes a_n\Big) -a_n \otimes \Pi_{n-1}\Big( a_1\otimes \cdots \otimes a_{n-1}\Big) \,\bigg]
$$ 
for $n\geq 2$ and any $a_1\otimes \cdots \otimes a_n\in  \boldsymbol{H}^{\otimes n}$.\footnote{Our formula for $\Pi_n$ does not exactly coincide with formula (5.20) in 
\cite{DGR} on the grounds of our convention for the symbol which is opposite to  the one used in {\it loc.\,cit.},\,see Remark \ref{Rem:II-dim>0}.2 above.} 
For any $n\geq 1$, the endomorphism $\Pi_n$ enjoys two properties: first it is a projector, that is $\Pi_n\circ \Pi_n=\Pi_n$;  second, its kernel precisely coincides with the subspace of non-trivial shuffle products, {\it i.e.}\,one has ${\rm Ker}(\Pi_n)={\rm Im}(\shuffle)^n$ as subspaces of $\boldsymbol{H}^{\otimes n}$. \sk

\vspace{-0.25cm}
A key fact to prove the first part of the preceding theorem 
 is that the following identity is satisfied for every $u,v\in \boldsymbol{H}$  (the verification of which is left to the reader)
\begin{equation}
\label{Eq:Pi-n}
\Pi_n\Big(u^{n-1}v\Big)=\frac{1}{n}\, u^{n-2}\big(uv-vu \big)
\end{equation}
(where $u^{n-1}v$ stands for the symbol $u^{\otimes (n-1)}\otimes v$, etc.) 
\sk


We are going to use the preceding algebraic considerations in the case when 
$\boldsymbol{H}=\boldsymbol{H}_I$. Since this space is formed by closed 1-forms, many of the properties of
one-variable hyperlogarithms  discussed page \pageref{Page:Properties-Hyperlogarithms} 
 admit direct analogues for the (germs of) iterated integrals constructed from $\boldsymbol{H}_I$, see Remark \ref{Rem:II-dim>0}.3. In particular, for any sufficiently generic $x\in V$,  there exists a linear map ${\rm II}_x : \boldsymbol{H}_I^{\otimes \bullet}\rightarrow \mathcal O_x$ which induces an isomorphism of 
$\mathbf C$-algebras from  $\boldsymbol{H}_I^{\otimes \bullet}$ onto its image $\mathcal H_x={\rm Im}\big( {\rm II}_x \big)$. 
\sk 

With the preceding considerations at hand, it is then easy to prove {\bf 1}:
first,  since $\boldsymbol{\mathcal  S}(\l {n} )=0^{n-1}1$ and because 
 ${\rm II}_x$ is an isomorphism of algebras, we get that  
an identity of the form $\sum_{i\in I} c_i \, { {\bf L}{\rm i}_{ n } }  (u_i)={\rm P}_{\hspace{-0.06cm} <n}$ holds true identically at $x$ if and only if 
$\sum_{i\in I} c_i \,  u_i^*(0^{n-1}1)$ belongs to $ {\rm Im}(\shuffle)^n$. 
When this indeed holds true,  it follows from the properties of $\Pi_n$ recalled just above (namely $ {\rm Im}(\shuffle)^n={\rm Ker}(\Pi_n)$ and \eqref{Eq:Pi-n} that one has  $\sum_{i\in I} c_i \,  u_i^* \big(\Pi_n(0^{n-1}1)\big)=\sum_{i\in I} c_i \,  u_i^* \big(0^{n-2}(01-10))\big)=0$ in $\boldsymbol{H}_I^{\otimes n}$, which gives us that  {\it (a)} implies {\it (b).} 
Conversely,  for any $u,v\in \boldsymbol{H}$,  one has 
 $u^{n-1}v-u^{n-2}(u v-vu)=u^{n-1}v-\Pi_n(u^{n-1}v)\in {\rm Im}(\shuffle )^n$ 
since $\Pi_n$ is a projector onto ${\rm Im}(\shuffle )^n$, 
 which allows to deduce easily {\it (a)} from {\it (b)} (details left to the reader). 
That {\it (b)} and {\it (c)} are equivalent follows immediately   from the fact that ${\rm II}_x : \boldsymbol{H}_I^{\otimes n}\rightarrow \mathcal H_x^n$ is an isomorphism.  This proves that the three points of the first part of Theorem \ref{T:FE-holomorphic} are equivalent. 
\sk

We now assume that the three 
 equivalent conditions of {\bf 1.}\,hold true. 
 Then from \eqref{Eq:togodoc}, it follows that $\sum_{i\in I} c_i  \mathfrak L_n^i(u_i)\equiv 0$ at $x$.   Set $\boldsymbol{w}_n=0^{n-2}(01-10)$ (whose pull-back under $u_i$ coincides with $ \boldsymbol{\mathcal S}\big(\mathfrak L_n^i(u_i)\big)$ for any $i\in I$). 
 Then, according to Lemma \ref{L:KWw-Sw-equivariant},   for any permutation $\sigma\in \mathfrak S_n$, since $u_i^*(\boldsymbol{w}_n)^\sigma=u_i^*(\boldsymbol{w}_n^\sigma)$, 
one has 
$\sum_i c_i  u_i^*(\boldsymbol{w}_n^\sigma)=0$ in $\boldsymbol{H}_I^{\otimes n}$ hence 
the sum 
$\sum_{i\in I} c_i  L_{\boldsymbol{w}_n^\sigma}^i(u_i)$ vanishes identically as well on a neighborhood of $x$.
 Now for $k\in \{0,\ldots,n-2\}$, since $0^{k}\shuffle 0^{n-2-k}(01-10)$ is the sum of the 
permuted symbols $ (0^{n-2}(01-10))^\sigma$ for $\sigma \in  k\shuffle (n-k)$ ({\it cf.}\,\eqref{Eq:shuffle} and the whole paragraph it is in), it follows that  $\mathfrak L_{n,k}^i=\sum_{\sigma \in k\shuffle (n-k)} L_{\boldsymbol{w}_n^\sigma}^i$ for each $i\in I$, from which it comes that 
$\sum_{i\in I} c_i  \mathfrak L_{n,k}^i(u_i)\equiv 0$ for any $k\leq n-2$.
 From the fact  that    the $n-1$ 
  functions $\mathfrak L^i_{n,k}$'s are linearly independent for any $i\in I$ (according to Proposition  \ref{P:Lnk-Lin-Indep} below),  it follows that the ARs corresponding to these 
    functional identities 
  are  linearly independent as well, which 
concludes the proof of the theorem.  
%
 \hfill $\square$


\begin{rem} 
Assuming that the equivalent assertions of the first point of the previous theorem hold true, we deduce from Theorem \ref{} that Zagier's single valued polylogarithm $\mathcal L_n$ satisfies the functional relation   $\sum_{i\in I} c_i\,\mathcal L_n(u_i)={\rm cst.}$ globally on (a certain Zariski open subset of) $V$.  It is natural to wonder whether this last condition does not imply the symbolic relation  $\sum_{i\in I} c_i\,u_i^*[0^{n-1}(01-10)]=0$ in return.  In the unpublished text {\rm \cite{Rudenko}}, Rudenko proves that it is indeed the case  when $V=\mathbf P^1$.
We expect that this holds true in full generality.
\end{rem}


\subsubsection{\bf The modified holomorphic polylogarithms $\boldsymbol{\mathfrak L_n}$.}
We now discuss the modified holomorphic polylogarithms 
 appearing in the statement of Theorem \ref{T:FE-holomorphic}. \sk

Let $\zeta$ be a generalized base point on $\mathbf P^1\setminus \{0,1,\infty\}$: either $\zeta$ belongs to this set or $\zeta$ is the tangent vector  $\partial/\partial x$ at the origin, denoted by $\vec{0}$ (see the last paragraph page \pageref{Pg:Tangential-base-point}). \mk 
 
By definition, $ {\mathfrak L}_n^{\zeta}$  is the (germ of) iterated integral at $\zeta$ which is the 
image of the symbol $\frac{1}{n}0^{n-2}(0\wedge 1)$ by the map ${\rm II}_{\zeta}^n $.  It is not difficult to get an explicit closed formula for $\mathfrak L_n^\zeta$ in terms of the classical polylogarithms (with coefficients being classical polylogarithms evaluated at $\zeta$).  To simplify our formulas, we will work on the open domain $\Omega=\mathbf C\setminus \Big(
]-\infty,0]\cup [1,+\infty[\Big)$,  with the principal determinations of the standard (poly)logarithmic functions on it.  \mk

Let us start with the simplest case, namely when $\zeta=\vec{0}$.  
For any $z \in \Omega$, we denote by $\gamma_{\vec{0}}^z: [0,1]\rightarrow \mathbf C$ a smooth path  from 0 to $z$ with derivative $\vec{0}$ at 0 and such that 
$\gamma_{\vec{0}}^z(t)\in \Omega$ for any $t\in ]0,1]$. 
Then, for any $n\geq 2$,  since 
$0^{n-2}\Big(01-10 \Big)=n\cdot 0^{n-1}1-0^{n-2}1\shuffle 0 $, 
  one has for any $n\geq 2$ and any $z \in \Omega$:
\begin{align*}
{\mathfrak L}_n^{\vec{0}}(z)=\int_{\gamma_{\vec{0}}^z} \frac{1}{n}\cdot 0^{n-2}\Big(01-10\Big)=\, &  {\rm {\bf L}i}_{ n } (z)-\frac{1}{n} \, 
\l{ n-1} (z) \cdot {\rm Log}(z) \, 
\end{align*}
(note that the second equality follows after the regularization process used to deal with iterated integrals with base point a tangent vector, {\it cf.}\,page \pageref{Pg:Tangential-base-point}).\mk

More interesting, especially in view of Theorem \ref{T:FE-holomorphic}, is the case when $\zeta$ belongs to $ \mathbf P^1\setminus \{0,1,\infty\}$.  To simplify, we assume that $\zeta$ belongs to $\Omega$ too.
 Then for  $z$ sufficiently close to the latter, one set 
$$
\widetilde \gamma_\zeta^{\,z}=\gamma_{\vec{0}}^z\cdot \gamma^{\vec{0}}_\zeta
$$
where $\gamma^{\vec{0}}_\zeta$ stands for 
$(\gamma_{\vec{0}}^\zeta)^{-1}$. 
 Thanks to the concatenation formula (see  (21) in  \cite{BanksPanzerPym}), it comes: 
\begin{align*}
n\cdot \mathfrak L_n^\zeta(z)=& \, 
\int_{\widetilde \gamma_\zeta^{\,z}}
0^{n-2}\Big(01-10\Big)\\
=& \,n\cdot \mathfrak L_n^{\vec{0}}(z)
+L_{0^{n-2}0}^{\vec{0}}(z) \int_{\gamma^{\vec{0}}_\zeta}1- 
L_{0^{n-2}1}^{\vec{0}}(z) \int_{\gamma^{\vec{0}}_\zeta}0+
\sum_{k=2}^n L_{0^{n-k}}^{\vec{0}}(z)
\int_{\gamma^{\vec{0}}_\zeta} 0^{k-2}(01-10)
\\
=& \,n\cdot \mathfrak L_n^{\vec{0}}(z)
+ {\rm {\bf L}i}_{n-1}(z) \,{\rm Log}(\zeta)
+\frac{{}^{}\,{\rm Log}(z)^{n-1}}{(n-1)!}\, {\rm Log}(1-\zeta)
+\sum_{s=0}^{n-2} 
\frac{{}^{}\,{\rm Log}(z)^{s}}{s!}\, 
\int_{\gamma^{\vec{0}}_\zeta} 0^{n-s-2}(01-10)\, .
\end{align*}

On the other hand, from  \cite[(13)]{BanksPanzerPym}, it follows that 
for any $k\geq 2$:
$$\int_{\gamma^{\vec{0}}_\zeta} 0^{k-2}(01-10)= 
\int_{(\gamma_{\vec{0}}^\zeta)^{-1}} {0^{k-2}(01-10)}
 = (-1)^{k} \int_{\gamma_{\vec{0}}^\zeta} \overline{0^{k-2}(01-10)}$$ 
 (where we note $\overline{\boldsymbol{w}}=w_n\cdots w_1$ for any word $\boldsymbol{w}=w_1\cdots w_n$), hence 
\begin{align*}
\int_{\gamma^{\vec{0}}_\zeta} 0^{k-2}(01-10)
 = (-1)^{k}\cdot  \Bigg(\ell_{k}(\zeta)- 
\int_{\vec{0}}^\zeta 
 \frac{\ell_{k-1}(u)}{u} du
 \Bigg)
\end{align*}
with $\ell_k(z)=L_{10^{k-1}}^{\vec{0}}(z)=\int_{\vec{0}}^z
\frac{{\rm Log}(u)^{k-1}}{(1-u)(k-1)!} du$ for every $k\geq 1$ and every $z\in \Omega$. 
  By direct computations,  for any such $k$ and $z$, one gets: 
\begin{align*}
\ell_k(z)
=& \, \sum_{s=0}^{k-1} (-1)^{k-1-s} \frac{1}{s!}\, {\rm Log(\zeta)^{s}} \cdot {\rm {\bf L}i}_{k-s}(\zeta) \\
\mbox{and }\, 
\int_{\vec{0}}^\zeta 
\ell_{k-1}(u) \frac{du}{u}
=& \, 
\sum_{s=0}^{k-2}(-1)^{k-s}\frac{ (k-1-s)}{s!} \, {\rm Log(\zeta)^{s}} \cdot {\rm {\bf L}i}_{k-s}(\zeta)\, , 
\end{align*}
thus  
\begin{equation*}
\int_{\gamma^{\vec{0}}_\zeta} 0^{k-2}(01-10)= \sum_{s=0}^k(-1)^{s+1}\frac{(k-s)}{s!} \,  {\rm Log(\zeta)^{s}} \cdot {\rm {\bf L}i}_{k-s}(\zeta)\, .
\end{equation*}
\bk 

Collecting everything, we get: 
  \begin{align*}
 \mathfrak L_n^\zeta(z)=& \, 
{\rm {\bf L}i}_{ n } (z)-\frac{1}{n} \, 
\l{ n-1} (z) \, {\rm Log}\left(\frac{z}{\zeta}\right)
+\frac{1}{n}\sum_{k=0}^{n-1}\frac{{}^{}\,{\rm Log}(z)^{k}}{k!}\, \kappa_{n}^k(\zeta)
  \end{align*}
with
$$
\kappa_{n}^{k}(\zeta)=\begin{cases} {}^{}\, \,\,
 \sum_{s=0}^{n-k}(-1)^{s+1}\frac{(n-k-s)}{s!} \,  {\rm Log(\zeta)^{s}} \, {\rm {\bf L}i}_{n-k-s}(\zeta) \qquad \mbox{if }\, k<n-1\\
\,\,  {}^{}\hspace{0.15cm}{\rm Log}(1-\zeta)  \hspace{4.8cm} \mbox{ if }\, k=n-1.
 \end{cases}
$$

It is then straightforward to obtain that, for $n\geq 2$ and any $z\in \Omega$, one has: 
 \begin{equation}
 \label{Eq:mathfrakLn}
 \mathfrak L_n^\zeta(z)= \, 
{\rm {\bf L}i}_{ n } (z)
- \frac{1}{n}\, \l { {n-1} } (z) \cdot 
{\rm Log}\bigg(\frac{z}{\zeta}\bigg) 
-\frac{1}{n}\sum_{k=0}^{n-1} \frac{n-k}{k!} 
{\rm Log}\bigg(\frac{z}{\zeta}\bigg)^{k} \cdot 
\l { {n-k} } (\zeta)\, . 
  \end{equation}


For instance, for $n=2$ and $n=3$ and $z\in \Omega$, one has : 
\begin{align}
\label{Eq:mathfrak-L2-zeta}
\mathfrak L_2^\zeta(z)=& \, \l {2} (z)- \l {2} (\zeta)
+\frac{1}{2}{\rm Log}\bigg(\frac{z}{\zeta}\bigg) \cdot {\rm Log}\Big((1-\zeta)(1-z)\Big)\\
\mathfrak L_3^\zeta(z)=& \, \l {3} (z)- \l {3 } (\zeta)-
\frac{1}{3}\, {\rm Log}\bigg(\frac{z}{\zeta}\bigg) \cdot 
\Big(  \l {2} (z) +2 \l {2} (\zeta)
\Big)
+\frac{1}{6}\, {\rm Log}\bigg(\frac{z}{\zeta}\bigg)^2\cdot {\rm Log}(1-\zeta)\, .
\nonumber
\end{align}


From formula \eqref{Eq:mathfrakLn}, it follows that $\mathfrak L_n^\zeta$ can be  characterized  as the unique weight $n$ iterated integral at $\zeta$ written 
$$
\mathfrak L_n^\zeta(z)=\l {n} (z)-\frac{1}{n}\, \l {n-1} (z)\cdot {\rm Log}\bigg(\frac{z}{\zeta}\bigg) +\sum_{k=0}^{n-1} A_k^\zeta \, {\rm Log}\bigg(\frac{z}{\zeta}\bigg)^k
$$
for some complex constant $A_k^\zeta$ uniquely determined by the fact that 
$\mathfrak L_n^\zeta(z)=O((z-\zeta)^{n+1})$ at $\zeta$. \mk 

Finally, it is not difficult to establish that $\mathfrak L_n^\zeta$ admits a development in series at $\zeta$ of the form
\begin{equation}
\label{Eq:DevelSeries-mathfrak Ln}
\mathfrak L_n^\zeta(z)=\sum_{k \geq n+1} \frac{P_{n,k}(\zeta)}{N_{n,k}} \frac{(z-\zeta)^k}{\zeta^{k-1}(\zeta-1)^{k-n+1}}\, , 
\end{equation}
{for some universal constants $N_{n,k}\in \mathbf Z^*$ and some universal polynomials $P_{n,k}(\zeta)\in \mathbf Z[\zeta]$, both of  which can  be made explicit easily (this is left to the reader).}
\begin{center}
\vspace{-0.1cm}
$\star$
\end{center}

Theorem \ref{T:FE-holomorphic} shows that, with regard to the AFEs it satisfies, the function   $\mathfrak L_{n}$ does not come alone but with the other weight $n$ iterated integrals $\mathfrak L_{n,k}^\zeta$ (for $k=0,\ldots,n-2$) defined by 
$$\mathfrak L_{n,k}^\zeta(z)={\rm Log}\Big(z/\zeta\Big)^k \, \mathfrak L_{n-k}^\zeta(z)
$$ 
for   $z$ sufficiently close to $\zeta$.

\begin{prop}
\label{P:Lnk-Lin-Indep}
The functions $\mathfrak L_{n,k}^\zeta$'s (for $k=0,\ldots,n-2$) 
 are linearly independant (over $\mathbf C$). 
\end{prop}
\begin{proof}  For $k=0,\ldots,n-2$ and $z$ sufficiently close to $\zeta$ in $\mathbf P^1$, one sets 
$$\mathfrak L_{n,k}^{\vec{0}}(z)={\rm Log}(z)^k  \cdot \mathfrak L_{n-k}^{\vec{0}}(z) = 
\bigg( 
\l { n-k } (z) -\frac{1}{n-k}  \l { n-k-1 } (z) \cdot {\rm Log}(z)
\bigg)\cdot {\rm Log}(z)^k$$ 
for some fixed  determinations of the classical polylogarithms at $\zeta$.  
 Since the $\mathfrak L_{n,k}^\zeta$'s are weight $n$ iterated integrals, the proposition is equivalent to the fact that their symbols $\boldsymbol{\mathcal S}(\mathfrak L_{n,k}^\zeta)=(k!) \cdot  0^k\shuffle   0^{n-2-k}(01-10)$ are linearly independent.  But as 
$\mathfrak L_{n,k}^{\vec{0}}$ and $\mathfrak L_{n,k}^\zeta$ have the same symbol for any $k$, it thus suffices to prove the corresponding statement for the former functions, which we are going to  establish   by induction below. \mk

Let $n$ be bigger than 2 and 
assume that, as a holomorphic germ at $\zeta$,  one has $\sum_{k=0}^{n-2} c_{k} \mathfrak L_{n,k}^{\vec{0}}\equiv 0$ for some constant $c_0,\ldots,c_{n-2}\in \mathbf C$.  Since all the functions involved extend as multivalued holomorphic functions on the complex plane with ramifications at $0$ and $1$, the same equation holds true in the space of such functions.  Expanding $\sum_{k=0}^{n-2} c_{k} \mathfrak L_{n,k}^{\vec{0}}$ in powers of ${\rm Log}(z)$ with holomorphic functions at the origin as coefficients, it comes immediately that $c_0=0$. Hence one has
$$0\equiv \sum_{k=1}^{n-2} c_{k} \mathfrak L_{n,k}^{\vec{0}}= {\rm Log}(z)\sum_{k=0}^{n-3}c_{k+1}\mathfrak L_{n-1,k}^{\vec{0}}
\qquad \mbox{ which implies  }\qquad 
0\equiv \sum_{k=0}^{n-3}c_{k+1}\mathfrak L_{n-1,k}^{\vec{0}}\,.$$ Using the induction hypothesis, one deduces that all the $c_k$'s are trivial 
 which ends the proof. 
\end{proof}

Regarding FEs,  Theorem  \ref{T:FE-holomorphic} applies  in particular when the $u_i$ are rational functions in a single variable. For instance, for any $n\geq 2$, we  get that  the classical  functional equations in one variable \eqref{Eq:EFA-Lin} satisfied by $\l { n} $ take the following form for the modified polylogarithms $\mathfrak L_n$: for  any $\zeta\in \mathbf P^1\setminus\{0,1,\infty\}$, any $z$ sufficiently close to $\zeta$ and for any positive integer $r$,  one has 
\begin{equation}
\label{Eq:Ln-inversion-distribution-formula}
\mathfrak L_{n}^\zeta(z)+(-1)^n\mathfrak L_{n}^{1/\zeta}\Big(z^{-1}\Big)=0 \qquad 
\mbox{ and } \qquad 
\mathfrak L_{n}^{\zeta^r}\Big(z^r\Big)=r^{n-1}\sum_{w^r=1}
\mathfrak L_{n}^{w  \zeta}(w z)\, 
\end{equation}
and these relations are satisfied as well when replacing $\mathfrak L_{n}$ by  $\mathfrak L_{n,k}$,  for any $k=0,\ldots,n-2$.
\begin{center}
\vspace{-0.1cm}
$\star$
\end{center}

The fact that $\mathfrak L_n^\zeta$ has valuation $n+1$ at $\zeta$ has an easy but interesting consequence regarding the functional equations of  type  \eqref{Eq-FE-Lin} satisfied by the classical polylogarithm $\l { n } $.

\begin{prop}
\label{P:d-geq-n+2}
 Let $(\mathcal E):\,
  \sum_{i=1}^d c_i\,  \l { n } (u_i)={\rm P}_{\hspace{-0.05cm}<n}$ be a non-trivial polylogarithmic FE  
 in several   variables. If the $u_i$'s appearing in it define $d$ distinct foliations then 
%
  necessarily 	
 $d\geq n+3$.
\end{prop}
\begin{proof} 
Let $x$ be a regular point for the $d$-web $\boldsymbol{\mathcal W}_I$ defined by the $u_i$'s: this means that all the $u_i$'s are defined at $x$ and that $du_i\wedge du_j(x)\neq 0$ whenever  the indices $i$ and $j$ are distinct (in other terms, $\boldsymbol{\mathcal W}_I$ satisfies (wGP)  at $x$). From  a property satisfied by the virtual rank (see 3.\,in the list of properties listed in {\S\ref{Para:Virtual-Rank}}\, ), it follows that $\rho^\sigma_x(\boldsymbol{\mathcal W}_I)=0$ for any integer $\sigma \geq d-1$. \sk

Now assuming that the functional equation $\sum_{i=1}^d c_i\,  \l { n } (u_i)={\rm P}_{\hspace{-0.05cm}<n}$ holds true, it follows from Theorem \ref{T:FE-holomorphic} that $\sum_{i=1}^d c_i\,  \mathfrak L_{ n }^i (u_i)=0$  (as a holomorphic germ at $x$). This functional relation can be seen as an AR for $\boldsymbol{\mathcal W}_I$ which is of valuation $n+1$ at $x$ since each  $\mathfrak L_{ n }^i $ has valuation $n+1$ at $u_i(x)$ for any $i=1,\ldots,d$ and because not all  the $c_i$'s are equal to 0 (otherwise the considered functional equation 
 $(\mathcal E)$  would be trivial). It follows that $\rho^{n+1}_x(\boldsymbol{\mathcal W}_I)>0 $.  Combining this with the result stated at the end of the preceding paragraph, we get that $n+1<d-1$, which proves the 
  proposition. 
\end{proof}

We think that,  under the assumptions of the preceding proposition, one can actually get the better bound $d\geq 2n+1$, as suggested by the following argument, which we believe is very likely.\mk 

Consider the following statement
for $n\geq 2$ fixed and any $\zeta\in \mathbf P^1\setminus \{0,1,\infty\}$: 
\begin{equation}
\label{Eq:jets-Ln-k,k}
\begin{tabular}{l}
{\it  
there exists 
 a non-trivial linear combination 
$
{\mathpzc{ L}}\hspace{-0.24cm}{\mathpzc{ L}}_n^\zeta=\sum_{k=0}^{n-2} \delta_{n,k}
 \, \mathfrak L_{n-k,k}
$,  } \\ {\it   for some coefficients $\delta_{n,k}\in \mathbf C$, 
which is of valuation $2n-1$ at $\zeta$.
 }
\end{tabular}
\end{equation} 


Assume that this statement holds true for $\zeta$ generic. 
Then given a generic point $x$:  (1) this point is  a regular point for the web $\boldsymbol{\mathcal W}_I$;  and  (2) the 
 points $\zeta_i=u_i(x)$ (for $i=1,\ldots,d$) are generic  points of $\mathbf P^1$. 
It follows that  the functions 
${\mathpzc{ L}}\hspace{-0.24cm}{\mathpzc{ L}}_n^i={\mathpzc{ L}}\hspace{-0.24cm}{\mathpzc{ L}}_n^{\zeta_i}$ have valuation $2n-1$ at $\zeta_i$ for each $i$,  hence the functional equation $\sum_{i=1}^d c_i\, {\mathpzc{ L}}\hspace{-0.24cm}{\mathpzc{ L}}_n^i(u_i)=0$ provides an AR of valuation $2n-1$ at $x$ for $\boldsymbol{\mathcal W}_I$. Arguing as in the proof of Proposition \ref{P:d-geq-n+2}, one deduces from this that one necessarily has $d\geq 2n+1$. \mk

We have verified that \eqref{Eq:jets-Ln-k,k} is satisfied (for any $\zeta$ in 
$\mathbf P^1\setminus\{0,1,\infty\}$) 
 for $n$ from 2 to 15, hence for those values,  one has $d\geq 2n+1$ in the conclusion of Proposition \ref{P:d-geq-n+2}.  Actually, for 
 each explicit value of $n$ we have dealt with, 
 normalizing the $\delta_{n,k}$'s by requiring that $\delta_{n,n-2}=(-1)^{n}$, we have found that the other  coefficients $\delta_{n,k}$'s are uniquely determined integers which 
moreover do not depend on $\zeta$. This is a strong evidence for believing that \eqref{Eq:jets-Ln-k,k} indeed holds true for any $n\geq 2$ and any $\zeta$.  
  By way of examples,
  we give below the sequences $\boldsymbol{\delta}_{n}=(\delta_{n,k})_{k=0}^{n-2}$ 
for $n$ from 3 to 9: 
\begin{align*}
\boldsymbol{\delta}_{3}= &\, (6, -1)\\
\boldsymbol{\delta}_{4}= &\, (60, -15, 1) \\
\boldsymbol{\delta}_{5}= &\, (840, -252, 27, -1) \\
\boldsymbol{\delta}_{6}= &\, (15120, -5040, 672, -42, 1) \\
\boldsymbol{\delta}_{7}= &\, (332640, -118800, 18000, -1440, 60, -1) \\
\boldsymbol{\delta}_{8}= &\, (8648640, -3243240, 534600, -49500, 2700, -81, 1) \\
\boldsymbol{\delta}_{9}= &\, (259459200, -100900800, 17657640, -1801800, 115500, -4620, 105, -1)\, .
\end{align*}

\begin{rem}
Using the \href{http://oeis.org}{OEIS}, we have verified that the following formulas hold true for $n\leq  15$: 
\begin{itemize}
\item[$\bullet$] $\delta_{n,0}=(2n-3)!/(n-2)! $
\hspace{1.7cm}
 ({\it cf.}\, sequence \href{https://oeis.org/A000407}{A000407 of the OEIS}) ;
\item[$\bullet$] $\delta_{n,n-3}=(-1)^{n-1}3(n-2)(n+1)/2$ \quad  ({\it cf.}\,sequence 
\href{https://oeis.org/A140091}{
A140091 of the OEIS}).
\end{itemize}
We conjecture that they hold true for  $n\geq 3$ arbitrary. More generally,  
it would be interesting to understand better the other integer values $\delta_{n,k}$'s and ideally, to have  closed formulas for them. \mk
\end{rem}

Since all the $\mathfrak L_{n-k,k}$'s have valuation $n+1$ at $\zeta$, 
\eqref{Eq:jets-Ln-k,k} is equivalent  to the fact that 
\begin{quote}
{\it  
${}^{}$ \hspace{1cm}`the $(2n-2)$-th  order 
  jets  at $\zeta$ of the  
   holomorphic 
    polylogarithmic \\  
    ${}^{}$ \hspace{1cm} 
    functions 
   $\mathfrak L_{n-k,k}^\zeta$ for $k=0,\ldots,n-2$ are linearly dependent'}.  
   \end{quote}
   That this holds true seems very likely but 
   the direct approach we tried to follow using 
  the power  
   series representations  \eqref{Eq:DevelSeries-mathfrak Ln}  led us to a combinatorial problem which requires some work to be solved. 
   We plan to come back to this in a future paper. 
 \subsection{\bf Some polylogarithmic functional equations and the webs associated to them}
\label{S:EFA-polylog}

We now consider some explicit polylogarithmic functional equations and study the webs associated to them from the point of view of their ARs and of their rank.   Some of the polylogarithmic FE we consider are quite classical and well-known, some others are either quite less popular or more recent.  The number of such equations appearing in the literature is huge 
and it is not the point here to make a complete survey on the topic. We have considered either the classical 
 polylogarithmic FEs or  some others, either classical too but not very well known, or quite recent, but which in any case give rise to interesting webs in terms of their rank. 
 \mk

Our presentation below mainly follows the chronological order of appearance  of these polygarithmic identities in the literature.
Our main  references here are the books by Lewin on polylogarithms, Zagier's and Gangl's papers already cited above, and  the two recent  dissertations \cite{Charlton,Radchenko}. 
\mk

\vspace{-0.3cm}
 Polylogarithmic AFEs  in several variables are known only when the weight $w$ is not too big. There are plenty of examples for the dilogarithm ($w=2$), some of them very classical, others quite more recent. Then the number of known  polylogarithmic identities decreases drastically as the weight increases.   As we write these lines, the current record is in weight $w=8$ and is due to Radchenko.
 Classical examples of polylogarithmic identities are known up to weight $5$: 
 the classical record is due to Kummer in 1840 who gave an identity involving 
 33 
  pentalogarithmic terms.  Starting from weight 6, the few known polylogarithmic AFEs involve far more terms (several hundreds typically) and  
 have been discovered much more recently (by Gangl and Radchenko),  
  via 
 a formal algebraic approach  implemented on a computer. 
 \vspace{-0.2cm}
 \begin{center}
$\star$
\end{center}
 
 
 Due to the complexity of the recently discovered  polylogarithmic AFEs (in particular all those in weight $w\geq 6$), we will not consider most of them below. Rather, we will mainly focus on  polylogarithmic AFEs  involving only some dozens of terms, typically 
in weight 2 and 3: the webs associated to those AFEs  are relatively easy to study, 
we have a sufficient number of them and  these are already interesting from the point of view of web geometry. 
\mk

 Almost all the results stated about the polylogarithmic webs we consider below have been obtained by brute force computations, using the formal algebraic methods discussed above (see the last paragraph of \S\ref{Par:II-AR}\hspace{0.14cm}page \pageref{poupou}  and \S\ref{Para:DeterminingARs-and-Rank}\hspace{0.14cm}as a whole).

\subsubsection{Notation for functional equations}
\label{SS:NotationForFEs}
Abelian functional equations satisfied by polylogarithms can be quite complicated, especially when the weight of the polylogarithm involved is high (already for $n\geq 3$, say).  Everytime it is possible and reasonable, we have chosen to write down the considered AFE as explicitly as possible.  
When this will not be  possible anymore, we will use some standard  notation to write down functional identities in a concise way, that we review here.\mk

In what follows, $n$ stands for a fixed integer bigger than or equal to 2. We denote by $\mathfrak S_n$ the permutation group on $\{1,\ldots,n\}$ and by 
$\varepsilon: \mathfrak S_n\rightarrow \{\pm 1\}$  the associated signature morphism.
\mk 

To simplify, we assume here that $\boldsymbol{F}$ is the field of rational functions in $n$ variables.
Then $\mathbf Z[\boldsymbol{F}]$ stands for the free $\mathbf Z$-module freely spanned by the set of symbols $[f]$ with $f\in \mathbf F$. Set theoretically, it is the set of  linear combinations $\sum_{i\in I} n_i [u_i]$ with  $n_i\in \mathbf Z$ and $u_i\in \boldsymbol{F}$ for any $i$ in $I$, the latter set being finite by assumption.  Then given a function $F$ (defined on a certain domain of $\mathbf P^1$), we will say that $\sum_{i\in I} n_i [u_i] $ is an AFE satisfied by $F$ if the sum 
$\sum_{i\in I} n_i F(u_i)$, considered as a function on a non-empty domain of $\mathbf P^n$ on which it is well-defined, is constant.  
If a bit unprecise, since it does not require to specify either 
the domain on which this identity holds true or the corresponding constant, 
this notation is very useful and will be sufficient for us.  With this notation in mind, setting $$\boldsymbol{[Ab]}=[x]-[y]-\left[\frac{x}{y}\right]-\left[\frac{1-y}{1-x}\right] + 
\left[\frac{x(1-y)}{y(1-x)}\right]
\in \mathbf Z\Big[\mathbf C (x,y)\Big]\, ,$$ the 5-term dilogarithmic identity \eqref{Eq:EFA-R} will be denoted by $R(\boldsymbol{[Ab]})=0$ or even just by $\boldsymbol{[Ab]}=0$ itself if we keep in mind that this identity concerns Roger's dilogarithm $R$. 
\begin{center}
$\star$
\end{center}


We recall the definition of the operators of (anti)symmetrization of a function of several variables (considered {\it e.g.}\,in \cite{Radchenko}) which prove to be useful for writing down in concise form certain functional identities  satisfied by polylogarithms involving an important number of terms. 
\sk

Let $\Phi$ be  a function of $n$ variables. Then for any 
 $\sigma\in \mathfrak S_n$, we denote by $\Phi^\sigma$ the function 
 $x \mapsto \Phi(x^\sigma)$ where we have set  $x^\sigma=(x_{\sigma(i)})_{i=1}^n $ for any 
  $n$-tuple  of (real or complex) numbers $x=(x_i)_{i=1}^n$. Then, 
 one denotes by ${\rm Sym}_n(\Phi)$ (resp.\,${\rm Alt}_n(\Phi)$) the symmetrization (resp.\,the anti-symmetrization) of $\Phi$:
$$ {\rm Sym}_n(\Phi)=\sum_{\sigma \in \mathfrak S_n} \Phi ^\sigma\qquad 
\mbox{ and } \qquad 
{\rm Alt}_n(\Phi) =\sum_{\sigma \in \mathfrak S_n} \varepsilon(\sigma)\Phi^\sigma\, . 
$$


With this notation, the symmetric form \eqref{Eq:Efa-(Ab)-symmetric} 
on ${\rm Conf}_5(\mathbf P^1)\simeq \mathcal M_{0,5}$ (which is birational to the projective plane $\mathbf P^2$) 
of the 5-terms relation for the dilogarithm 
has the very condensed form 
$$
0= {\rm Alt}_{5}\Big( [r_1] \Big)\, , 
$$ where  $r_1$ stands for  
the rational function on $\mathcal M_{0,5}$ given by $[(p_1,\ldots,p_5)]\mapsto [ p_1,p_2 ; p_3,p_4]$ where $[\cdot  ; \cdot]: \mathcal M_{0,4}\rightarrow \mathbf P^1\setminus \{0,1,\infty\}$ is the usual cross-ratio of four points on the projective line.

\subsubsection{Functional equations of the dilogarithm.}
\label{S:EFA-dilog}

If some FEs for weight $n$ polylogarithms are known for $n$ up to 8, it is for the case $n=2$, that is for the dilogarithm, that the theory as well as a important number of examples are the most developed. 
Here we give a short review of many interesting dilogarithmic identities.    References are in numbers for such AFEs:
in addition to the general references given above, one can mention 
 \cite{Kirillov}.  \mk 
 
  It is by now well known that the theory of cluster algebras is a source of many interesting dilogarithmic identities. We will not touch upon this 
  here
  since it will be discussed in more detail, from a web-theoretic perspective, in Section \ref{S:Cluster}. In it, the reader will find many series of examples and references about dilogarithmic identities within the theory of cluster algebras. 
%
%
%

\paragraph{\hspace{-0.2cm}Notation for dilogarithmic functions.}\hspace{-0.2cm}
\label{Par:NotationDilogarithmicFunctions}
Several distinct (but related) versions  of Roger's dilogarithm are used in this text and in particular in the present section. We first give below their definition and then recall some properties of these functions and also how they are related.
 \mk

Here are the notations for the different variables we use in this paragraph: \sk

${}^{}$\quad $-$ $x$ stands for a real number in $]0,1[$; \sk 

${}^{}$\quad $-$ $u$ denotes an element of $\mathbf R_{>0}=]0,+\infty[$;\sk 

${}^{}$\quad $-$ $\zeta$ is a fixed arbitrary element of $\mathbf P^1$, distinct from $0,1$ or $\infty$; and \sk 

${}^{}$\quad $-$ $z$ refers to a complex number sufficiently close to $\zeta$.
\bk 

By definition, the original {\bf Roger's dilogarithm} $\boldsymbol{{\mathcal R}}$ is the function of $x$ defined by
$${\mathcal R}(x)= \l {2} (x) +\frac{1}{2} {\rm Log}(x)\, {\rm Log}(1-x)
\, ; 
$$

It is sometimes  useful to use the following translated version $R$ of this function, defined by 
$$ {}^{}\hspace{0.8cm}
R(x)=\l {2} (x) 
+\frac{1}{2} {\rm Log}(x)\, {\rm Log}(1-x)- \frac{{}^{}\hspace{0.15cm}\pi^2}{6}\, ;
$$

 The {\bf $\boldsymbol{(\mathfrak B_{-1})}$-} or {\bf cluster Rogers' dilogarithm $\boldsymbol{\mathsf R}$} is the function of $u$ defined 
 by: 
$$
 {}^{}\hspace{0.4cm}
{\mathsf R}(u)=
\frac{1}{2} \bigintsss_{0}^u \left( \frac{{\rm Log}(1+v)}{v}-
 \frac{{\rm Log}(v)}{1+v}\right)\,  dv\, ;
$$

Finally, the `{\bf $\boldsymbol{\zeta}$-localized Roger's dilogarithm} 
$\boldsymbol{{\mathfrak R}^\zeta}$ is the function of $z$ defined by 
$$ {}^{}\hspace{0.7cm}
{\mathfrak R}^\zeta(z)= 
-\frac{1}{2}
\bigintsss_{\zeta}^z 
\left(
\frac{{\rm Log}\left(\frac{1-\sigma}{1-\zeta}\right)}{\sigma} +
\frac{{\rm Log}\left(\frac{\sigma}{\zeta}\right)}{1-\sigma}
\right)\, d\sigma\, .
$$

   The three functions $\mathcal R$, ${\mathsf R}$ and ${\mathfrak R}^\zeta$  
   (we leave $R$ out since it is just a translation of $\mathcal R$) 
   are weight 2 iterated integrals with formally the same symbol. Indeed, 
setting 
$\omega_{\xi_0}=d\xi/(\xi-\xi_0)$ for any $\xi_0\in \mathbf C$, 
and writing for short $\xi_0\xi_1$ for $\omega_{\xi_0}\omega_{\xi_1}$ for any 
$\xi_0,\xi_1\in \mathbf C$, one verifies that:\sk 

${}^{}$\quad $-$ $\mathcal R$ is the iterated integral with symbol $(01-10)/2$ and tangential base point $\vec{0}$ ({\it cf.}\,page \pageref{Pg:Tangential-base-point}); 
\sk 

${}^{}$\quad $-$ $\mathsf R$ is the iterated integral with symbol $(0(-1)-(-1)0)/2$ and tangential base point $\vec{0}$; and \sk

${}^{}$\quad $-$ ${\mathfrak R}^\zeta$ is the iterated integral with symbol $(01-10)/2$ and base point $\zeta$.\mk 

We now indicate in which way the previous functions are related. First, we note that generalizing the definition of ${\mathfrak R}^\zeta$ by allowing tangential base-points, one can write $\mathcal R={\mathfrak R}^{\vec{0}}$. The functions $\mathcal R$, $R$ and ${\mathsf R}$ are related by simple functional identities: for any $x$ and $u$ as above, one has 
$$
R(x)={\mathcal R}(x)-\pi^2/6
\qquad \mbox{ and }\qquad 
{\mathsf R}(u)=\mathcal R\left(\frac{u}{1+u} \right) .
$$
For any $x\in ]0,1[$, the following formula (already known by Abel, see \cite{Abel}) also holds true: 
$${\mathsf R}(x)=- \l {2} (-x) - \frac{1}{2} \,{\rm Log}(x)\,{\rm Log}(1+x).$$ 
Finally, we recall the following expression of $\mathfrak R^\zeta$ in terms of classical polylogarithms (see \eqref{Eq:mathfrak-L2-zeta}): 
$$
\mathfrak R^\zeta(z)=  \l {2} (z)- \l {2} (\zeta)
+\frac{1}{2}\,{\rm Log}\big({z}/{\zeta}\big) \cdot {\rm Log}\big((1-\zeta)(1-z)\big)\, .
$$

%
%
%
%

\paragraph{\hspace{-0.2cm}The five-term functional equation of the dilogarithm and Bol's web.}
\label{Par:Abel'sAFE}
This paragraph is about the by far most famous and well-known polylogarithmic functional equation, the so-called  `five-terms' or `Abel's' functional equation.  It has already been written down above (when introducing the dilogarithm, see \S\ref{SS:TheDilogarithm}) but due to its importance, we discuss it again   below. \mk 

Even if the AFE below satisfied by $R$ has been obtained by Rogers, we will call it `{\bf Abel's relation}' for the dilogarithm since this name is  one of the two which are the most commonly used nowadays (the other being the `{\bf five terms relation}'): 
$$
(\mathcal A{\rm b})\qquad  \qquad\qquad
R (x) -
R (y)
-R\bigg(\,\frac{x}{y}\, \bigg)
-R \bigg(\,\frac{1-y}{1-x}\,\bigg)+
R\bigg(\,\frac{x(1-y)}{y(1-x)}\,\bigg)
=0\, .\qquad\qquad\qquad\qquad
$$
This identity holds true for any real numbers $x,y$  such that $0<x<y<1$.
\mk

%
%
%

To $(\mathcal A{\rm b})$ one associates the web defined by the  rational functions appearing as arguments of $R$ in it. One gets a planar 5-web on $\mathbf P^2$, which could be denoted by $\boldsymbol{\mathcal W}_{\hspace{-0.05cm}\mathcal A{\rm b}}$ and named 
{\bf Abel's 5-web} $\boldsymbol{\mathcal W}_{\hspace{-0.05cm}\mathcal A{\rm b}}$, but which is better known in web geometry as 
{\bf Bol's web} and noted accordingly by $\boldsymbol{\mathcal B}$: 
$$
\boldsymbol{\mathcal W}_{\hspace{-0.05cm}\mathcal A{\rm b}}=
\boldsymbol{\mathcal B}=\boldsymbol{\mathcal W}\left(\,
x    \, , \, 
y         \, , \, 
    \frac{x}{y}      \, , \, 
  \frac{1-y}{1-x}
       \, , \, 
  \frac{x(1-y)}{y(1-x)}             
\,\right)\, . 
$$

This web has been named after Bol since he was the first 
to recognize it (in 1936) as an example of an exceptional (that is, of maximal rank but non algebraizable) web.  For the anecdote, 
 $\boldsymbol{\mathcal W}_{\hspace{-0.05cm}\mathcal A{\rm b}}$ was actually first considered by Blaschke in \cite[\S4]{Blaschke1933}, but as an example of a planar 5-web of almost maximal rank 5. This was wrong, as observed three years later by Bol in \cite{Bol}:  Blaschke had precisely forgotten to consider the   AR of  $\boldsymbol{\mathcal W}_{\hspace{-0.05cm}\mathcal A{\rm b}}$ 
associated to 
 $(\mathcal A{\rm b})$!
\mk 

Bol's web is hexagonal  and all his ARs are logarithmic, except the one associated to Abel's relation $(\mathcal A{\rm b})$.  In other terms, one has
$$ \boldsymbol{\mathcal Hex}_3\big(\boldsymbol{\mathcal B}\big)=10 
\qquad \mbox{ and }
 \qquad {\rm polrk}^\bullet\big(\boldsymbol{\mathcal B}\big)=(5 ,1)\,. 
$$

Bol's web is of great importance in web geometry.  Historically, first of all,  it has been the single  known example of an exceptional web for almost 70 years. Second, mathematically, in addition to being exceptional, it is characterized by a nice and simple  criterion\footnote{Namely: ``{\it Bol's web is the unique hexagonal planar web which is not linearizable}''.} and it can be obtained  in several natural geometric ways: it is the 5-web on the moduli space $\mathcal M_{0,5}$ admitting as first integrals the five forgetful maps $f_i : \mathcal M_{0,5}\rightarrow 
\mathcal M_{0,4}\simeq \mathbf P^1\setminus \{0,1,\infty\}$; moreover it is Segre's 5-web of the anticanonical model in $\mathbf P^5$ of the compactification 
$\overline{\mathcal M}_{0,5}$ which is a Del Pezzo surface of degree 5. Finally, Bol's web has a kind of `universality  property' for planar 5-webs.\footnote{To any  planar web $\boldsymbol{\mathcal W}=(\mathcal F_1,\ldots,\mathcal F_5)$ defined on a complex surface $S$ is invariantly attached the map $T_{\boldsymbol{\mathcal W}}: 
S\rightarrow \mathcal M_{0,5}
$ associating to $s\in S$  the configuration 
$([T_s \mathcal F_1], \ldots,[T_s \mathcal F_5])\in \mathcal M_{0,5}$
obtained from the tangent directions $[T_s \mathcal F_i]\in \mathbf P(T_s S)\simeq \mathbf P^1$ of the foliations of $\boldsymbol{\mathcal W}$ at this point. Then except when $ T_{\boldsymbol{\mathcal W}}$ is tangentially degenerate everywhere, $T_{\boldsymbol{\mathcal W}}^*(\boldsymbol{\mathcal B})$ is a model of Bol's web on $S$ which is canonically attached to the considered web ${\boldsymbol{\mathcal W}}$.} 
\mk 

All these properties make  Bol's web  
 $\boldsymbol{\mathcal B}$ 
 a fundamental object in web geometry, which is very likely related 
 to the fact that its most important abelian relation (the dilogarithmic one)
  seems to be one of the most fundamental polylogarithmic functional equations.  What is still a mystery is what a conceptual explanation about why and
  how these two facts are related could be. 

\subparagraph{}
\label{Subpar:(AB)Bol}
\hspace{-0.7cm}
Due to its importance in what concerns dilogarithmic functional equations, it is natural to look for a conceptual understanding of Abel's identity. This has been obtained by many authors, from several points of views. Among others, we list  some of these below: 
\begin{itemize}
\item[$-$] In \cite{GelfandMcPherson}, Gelfand and MacPherson have constructed geometrically Bol's web and the real version of Abel's relation by means of integration along the fibers (corresponding to quotienting by the diagonal action of Cartan subgroups) of  invariant differential forms on real grassmannians representing some characteristic classes (namely the Pontryagin classes);
\item[$-$] In their {\it `quest for higher logarithms'}, Hain and MacPherson 
have recovered in \cite{HainMacPherson}  (an holomorphic and multivalued version of) $(\mathcal A{\rm b})$. They used  multivalued holomorphic forms on (zariski open-subsets of) complex grassmannians and were motivated by the possibility of constructing this way kinds of {\it `universal Chern classes'};  
\sk 
\item[$-$] In \cite{Goncharov1995}, motivated by the construction of `regulators' in $K$-theory 
(for  function fields of algebraic varieties) and by means of integral transforms, 
Goncharov constructs currents on spaces 
of effective cycles in any projective space $\mathbf P^n$ in general position with respect to a simplex  $L\subset \mathbf P^n$.  
In particular,  he gets a real-analytic function $\mathcal P_2$ defined on the space of generic effective 1-cycles  in $\mathbf P^3$, that he calls the {\it `Chow dilogarithm'}. The latter  satisfies a certain geometric functional equation with 5 terms on the space of generic effective 2-cycles in $\mathbf P^4$ which gives (a global real-analytic version of) Abel's five-terms relation when restricted to  the subvariety formed by linear 2-cycles (that is, 2-planes)  in $\mathbf P^4$.  In the 7-th section
of \cite{Goncharov1995} is  sketched the construction of a corresponding (complex) analytic Chow dilogarithm which  satisfies the corresponding 
(holomorphic version) of Abel's relation $(\mathcal A{\rm b})$;
%
%
%
%
%
\item[$-$]  In some of their papers, Kerr and Lewis return on Goncharov's work, correcting and clarifying some points.  Recently with 
Lopatto 
in \cite{KerrLewisLopatto}, they  have translated into the simplical setting some complex analytic constructions 
 (namely, Abel-Jacobi maps for higher Chow groups) 
previously given in a `cubical setting' by the first two authors (also with M\"uller-Stach).
 The main results obtained in \cite{KerrLewisLopatto} are two  {\it `Reciprocity laws'} for algebraic cycles on smooth (quasi-)projective varieties that are both used 
 in the fifth section to derive two holomorphic versions of the five-terms relation for the dilogarithm (see formulas (5.2) and (5.5) in {\it loc.\,cit.}).
 \sk 
\end{itemize}

It would be very interesting to better understand the links between the approaches listed above and the webs associated to polylogarithmic AFEs. It is not impossible that any such functional equation (or at least, many of these) could be explained by means of a family of  cycles lying in the kernel of a certain Abel-Jacobi map. 
As we are writing this, this has been only obtained for Abel's dilogarithmic identity
$(\mathcal A{\rm b})$ and also, but in a less satisfactory way, for 
 Spence-Kummer equation of the trilogarithm  (see the end of \S\ref{SS:Spence-Kummer} further). 
If true in general, this would place (at least some) polylogarithmic webs and the classical algebraic webs associated to algebraic curves on a very similar footing. This is a research direction that we find particularly worthwhile.

\paragraph{\hspace{-0.2cm}Newman's six-terms functional equation $\boldsymbol{({\cal N}_6)}$.}
\label{Par:Newman-equation}
In \cite{Newman} (see also \cite[\S1.6]{Lewin}),  
Newman establishes that  the bilogarithm $\l{ 2 } $ satisfies the following clean functional equation\footnote{This nice functional  equation is nowadays attributed to Newman but an equivalent identity 
was actually  obtained by Kummer more than fifty years before (see formula (13) p.\,86 of \cite{Kummer1}  and compare with 
$({\cal N}_6)$).}
$$
({\cal N}_6)
\qquad \quad 
2\, \l{2} ( x ) 
+2\, \l{2} ( y ) 
+
2\, \l{2} ( z ) 
-
 \l{2} \big(x+y-xy \big) 
- \l{2} \big(x+z-xz \big) 
- \l{2} \big( y+z-yz \big) =0\, \qquad \qquad  {}^{}
$$ which holds true 
as soon as the real variables  $x,y,z$ have absolute value in $]0,1[$ and 
satisfy the relation 
\begin{equation}
\label{Eq:xyz1} \qquad \qquad 
\frac{1}{x}+\frac{1}{y}+\frac{1}{z}=1 
\quad \Longleftrightarrow \quad xy+xz+yz-xyz=0\, .
\end{equation}

Then one defines Newman's web $
\boldsymbol{\mathcal W}_{\hspace{-0.1cm}{\cal N}_6}$ as the planar 6-web defined by the six rational first integrals appearing as arguments of $\l{ 2} $ in $({\cal N}_6)$ (up to substraction of 1 for the last three, in order to get nice formulas for these): 
$$ 
\qquad \qquad 
\boldsymbol{\mathcal W}_{\hspace{-0.1cm}{\cal N}_6}=\boldsymbol{\mathcal W}\bigg( \, x \, , \, y 
\, ,  \,  
z  \, ,  \, 
(1-x)(1-y)
\, ,  \, 
(1-x)(1-z)
\, ,  \, 
(1-y)(1-z)
\bigg)\, . \qquad \quad 
$$
Of course here and in what follows,  the variables $x,y$ and $z$ are assumed to verify the  algebraic relation \eqref{Eq:xyz1}. 
Solving the latter with respect to $x$ and $y$, and considering new variables $u,v$ defined by the relations 
$ x=u$ and $y=\frac{u(1-v)}{v(1-u)}$, one has  $z=-\frac{u(1-v)}{1-u}$ and 
one gets a presentation of (a web equivalent to) Newman's 6-web defined by rational  first integrals expressed as products of powers of simple affine functions of $u$ and $v$: 
\begin{equation}
\label{Eq:N6-Kummer}
\boldsymbol{\mathcal W}_{\hspace{-0.1cm}{\cal N}_6}\simeq 
 \boldsymbol{\mathcal W}\Bigg(\, u\, ,\, uv, \frac{u}{v} \, ,\, \frac{u(1-v)}{v(1-u)}\, ,\,  -\frac{u(1-v)}{1-u}\, ,\,   \frac{u(1-v)^2}{v(1-u)^2} \, \Bigg) \, . 
\end{equation}

Newman's web has maximal rank: it admits a basis of its space of abelian relations  formed by six linearly independent logarithmic ARs; two non-colinear dilogaritmic ARs (one of which being the one associated to  
$({\cal N}_6 )$); the rational AR corresponding to the relation
\eqref{Eq:xyz1}
  which is assumed to be identically satisfied; and a last AR which is the one  associated to the following functional relation 
  $$A\left(\frac{u}{v}\right)-A(uv)- A\left(\frac{u(1-v)^2}{v(1-u)^2}\right)=0$$ which is identically satisfied by the function 
  $A:  u\mapsto {\rm Arctanh}(\sqrt{u})$ 
 on the set of pairs $(u,v)$  
of real numbers $u$ and $v$ such that 
  $0<u<v<1$.  We will see later on in 
  \S\ref{SS:cluster-webs:B2}  
  that $\boldsymbol{\mathcal W}_{\hspace{-0.1cm}{\cal N}_6}$  is equivalent to  the cluster web of type $B_2$.\footnote{An explicit basis of the space of abelian relations of $\boldsymbol{\mathcal W}_{\hspace{-0.1cm}{\cal N}_6}$ can also be obtained from the one of 
$\boldsymbol{\mathcal X\hspace{-0.05cm}\mathcal W}_{\hspace{-0.05cm}B_2}$ given in \S\ref{SS:cluster-webs:B2} (see especially the decomposition in direct sum \eqref{Eq:AR-XWB2} therein).}  \sk

Here is a list of some  invariants of Newman's web: 
\begin{equation}
\label{Eq:Invariants-WN6} 
\mathcal Hex\Big(\boldsymbol{\mathcal W}_{\hspace{-0.1cm}{\cal N}_6}\Big)=\big(8,0,0,0\big) 
\qquad  \mathcal F \hspace{-0.1cm}lat \Big(\boldsymbol{\mathcal W}_{\hspace{-0.1cm} {\cal N}_6}\Big)=\big(8,6,0,1\big) 
 \qquad {\rm polrk}^\bullet\Big(\boldsymbol{\mathcal W}_{\hspace{-0.1cm}{\cal N}_6}\Big)=\big(6 ,2\big)\,. 
\end{equation}


\subparagraph{}
\hspace{-0.6cm}
\label{Subpar:N6-lambda}
Newman's web is a specific element of an interesting 1-dimensional family of webs which all carry a dilogarithmic AFE which is a deformation of Newman's equation $({\cal N}_6)$. 

Indeed, for any complex parameter $\lambda\in \mathbf C$, we consider the web 
$$ 
\qquad \qquad 
\boldsymbol{\mathcal W}
_{\hspace{-0.05cm}{N}_{6,\lambda}}=\boldsymbol{\mathcal W}\bigg( \, X \, , \, Y
\, ,  \,  
Z  \, ,  \, 
XY
\, ,  \, 
XZ
\, ,  \, 
YZ\, 
\bigg) \qquad \quad 
$$
where the variables $X,Y$ and $Z$ are assumed to satisfy $X+Y+Z-XYZ=2 \lambda$ (note that the latter relation is preserved by the symmetry $(X,Y,Z)\mapsto (-X,-Y,-Z)$ from which it comes  that    
$\boldsymbol{\mathcal W}_{\hspace{-0.05cm}{N_{6,-\lambda}}}$ is isomorphic to $\boldsymbol{\mathcal W}_{\hspace{-0.05cm}{N}_{6,\lambda}}$ for any parameter $\lambda$). 
 For $\lambda=0$ and $\lambda=\pm\, 1$, one gets   webs which are isomorphic to the initial 
Newman's web:  one recovers $\boldsymbol{\mathcal W}_{\hspace{-0.1cm}{\cal N}_6}$ from 
$\boldsymbol{\mathcal W}_{\hspace{-0.05cm}{N}_{6,1}}$ by means of the change of variables $X=x+1$, $Y=y+1$ and $Z=z+1$, whereas it requires a transcendental change of variables to get $\boldsymbol{\mathcal W}_{\hspace{-0.1cm}{\cal N}_6}$ from 
$\boldsymbol{\mathcal W}_{\hspace{-0.05cm}{N}_{6,0}}$, see \cite[p.\,134]{Lewin}.  More generally, it can be verified that $\boldsymbol{\mathcal W}_{\hspace{-0.05cm}{N}_{6,\lambda}}$ has maximal rank  (that is, rank equal to 10) only for $\lambda =0,-1, 1$ whereas it has (almost maximal) rank 9 for any complex 
parameter $\lambda$ distinct from any of these three particular values. \mk   

For any $\lambda\in \mathbf C$, one denotes by $\lambda_\pm= \lambda\pm \sqrt{\lambda^2-1}$  the 
 two roots  of the trinomial $u^2-2\lambda u+1=0$ (which coincide if and only if $\lambda=\pm 1$).  Then 
\begin{equation}
\label{Eq:L2lambda}
{\bf L}_{2,\lambda}(\bullet)=\int_0^{\,\bullet} \frac{\log\Big(u^2-2\lambda u+1\Big)}{u} du
=\int^{\,\bullet}_0  \frac{\log(u-\lambda_+)}{u} du
+
\int^{\,\bullet}_0  \frac{\log(u-\lambda_-)}{u} du
\end{equation}
is a weight 2 iterated integral on 
 the punctured complex plane  $\mathbf C\setminus \{0, \lambda_+, \lambda_-\}$. 
  Then as proved by Newman (see again \cite[p.\,134]{Lewin}),  there is a 1-dimensional family of dilogarithmic AFEs
 $$
 ({N}_{6,\lambda})\qquad \quad 
2 \, {\bf L}_{2,\lambda}
 (X) + 
2 \, {\bf L}_{2,\lambda} (Y) 
+
2 \, {\bf L}_{2,\lambda} (Z) 
- \l {2}  ( XY ) 
- \l {2 } (XZ) 
-\l { 2 } (YZ)=2\,\bigg( \frac{\pi}{2}-\arccos(\lambda)\bigg)^2\, , 
$$
\label{Page:N6lambda}
(at least for any $\lambda$ in the real interval $ [0,1]$, but this certainly holds true for any parameter $\lambda \in \mathbf C$). The identity  $ ({N}_{6,\lambda})$ can be seen as a deformation of Newman's equation $({\cal N}_6)$ since the latter is equivalent to $ ({N}_{6,1})$ (and to $ ({N}_{6,1})$ as well).
Finally, here is a list of some  invariants of deformed Newman's  web $\boldsymbol{\mathcal W}_{\hspace{-0.1cm}{N}_6,\lambda}$ when $\lambda\in \mathbf C\setminus\{ \pm 1\, , \, 0\, \}$: 
\begin{equation}
\label{Eq:Invariants-WN6-lambda}
 \mathcal Hex\Big(\boldsymbol{\mathcal W}_{\hspace{-0.05cm}{N}_6,\lambda}\Big)=\big(7,0,0,0\big) 
\qquad  \mathcal F \hspace{-0.1cm}lat \Big(\boldsymbol{\mathcal W}_{\hspace{-0.05cm} {N}_6,\lambda }\Big)=\big(7,3,0,0,0\big) 
 \qquad {\rm polrk}^\bullet\Big(\boldsymbol{\mathcal W}_{\hspace{-0.05cm}{N}_6,\lambda}\Big)=\big(6 ,2\big)\,. 
\end{equation}
\begin{center}
$\star$
\end{center}

The function ${\bf L}_{2,\lambda}$ has been considered for the first time by Hill (compare \eqref{Eq:L2lambda} with the function denoted by $D_{\hspace{-0.05cm}x}^{\hspace{0.02cm}\alpha}$ in {\rm \cite{Hill})}.  It is a particular case of a class of functions recently introduced by Nakanishi in {\rm \cite{Nakanishi2015}} (see also {\rm \cite[\S3]{Nakanishi2018}}), the so called `{generalized (Euler) dilogarithms}'.  This author has proved ({\it cf.}\,{\rm \cite[Theorem\,4.5]{Nakanishi2018})} that to any period in a `{\it generalized cluster algebra}' (a notion introduced by Chekov and Shapiro in {\rm \cite{ChekhovShapiro})} is associated a functional equation satisfied by the corresponding  generalized Rogers dilogarithm. Since $({\cal N}_6)$ corresponds to the case of the (classical) cluster algebra of type $B_2/C_2$, 
 one can wonder whether or not, for any parameter $\lambda\in \mathbf C$, the deformation $(N_{6,\lambda})$ of $({\cal N}_6)$ is  a generalized dilogarithmic identity associated to a period in a generalized cluster algebra {(see \S\ref{Par:Periods} below for a definition of webs associated to a period of a classical cluster algebra)}.  This will be discussed again at the end of this memoir,  in  \S\ref{Par:Generalized-cluster-Algebras}.

\subparagraph{}
\hspace{-0.6cm}
In addition to being a specific element of an interesting 1-dimensional family of AFEs, 
Newman's equation $( {\mathcal N}_6)$ has another interesting property:  
it admits an equivalent formulation in terms of a two-variables hypergeometric function which is 
quite  intriguing and has been unnoticed before (as far as we know).  \mk

Denote by $F$ the \href{https://en.wikipedia.org/wiki/Appell_series}{third Appell's hypergeometric function  $F_3(x,y)=F_3(a_1,a_2,a_3,a_4;c ,x,y)$} when the parameters are specialized as  $a_1=\cdots=a_4=1$ and $c=3$. The main result of the note \cite{Sanchis-Lozano} is the following formula relating  the two-variables function  $F$ to a dilogarithmic expression
$$
 \frac{1}{2}xy\, F( x,y)={\rm Li}_2(x)+{\rm Li}_2(y)-{\rm Li}_2(x+y-xy)\, , 
$$
which holds true identically assuming that the absolute value of the arguments of $1-x$, $1-y$ and $(1-x)(1-y)$ be less than $\pi$.  Then, up to the previous relation 
and still assuming that \eqref{Eq:xyz1} holds true,   Newman's equation  
  $( {\mathcal N}_6)$   is equivalent  to the following functional identity in two variables
  \begin{equation}
\label{Eq:FE-F3-hypergeom}
 xy\,F(x,y)+xz\,F(x,z)+yz\,F(y,z)=0
  \end{equation}
 for the hypergeometric function $F$.  We didn't see any functional equations  of this kind in the literature about multivariable hypergeometric functions we have looked at. It would be interesting to know if   
\eqref{Eq:FE-F3-hypergeom} is truly new and to have a conceptual understanding of it within the theory of hypergeometric functions.

\paragraph{\hspace{-0.2cm}Mantel's dilogarithmic functional equation.}
\label{Par:Mantel}
In his 1898 paper \cite{Mantel} (see also \cite[\S1.6]{Kirillov} or  p.\,17 and p.\,394 in \cite{LewinStruct}), Mantel 
establishes that  the following identity 
\begin{align}
\label{E:Mantel}
\l{2} (x)
 + \l{2} (y)
&\,  - \l{2} (v)
- \l{2} (w)
 + \l{2} \left(\frac{v}{x}\right)
+ \l{2} \left(\frac{v}{y}\right) \\ 
&\, + \l{2} \left(\frac{w}{x}\right)
+ \l{2} \left(\frac{w}{y}\right) - \l{2} \left(\frac{vw}{xy}\right)=
- \frac{1}{2}\, {\rm Log}^2\bigg(- \frac{x}{y} \bigg) \nonumber
\end{align}
holds true as soon as the four real variables 
$x,y,v,w$ are assumed to satisfy identically the relation $(1-v)(1-w)=(1-x)(1-y)$.
 Under this condition, one gets {\bf Mantel's web} 
$$\boldsymbol{\mathcal W}_{\hspace{-0.05cm}\mathcal M}=\boldsymbol{\mathcal W}\Bigg( \, x \, , \, y 
\, , \, v
\, , \, w 
\, , \,  \frac{v}{x}
\, , \,  \frac{v}{y}
\, , \,  \frac{w}{x}
\, , \, \frac{w}{y}
\, , \,\frac{vw}{xy}\, 
\Bigg)$$
that we see as a 9-web in the three variables $x,y$ and $v$ (thus 
with $w=1-(1-x)(1-y)/(1-v)$).  \mk 

By direct computations, one verifies that 
$$
\rho^\bullet\Big(\boldsymbol{\mathcal W}_{\hspace{-0.05cm}\mathcal M}\Big)=(6,3,1)
\qquad \mbox{ and }
\qquad 
{\rm polrk}^\bullet\Big(\boldsymbol{\mathcal W}_{\hspace{-0.05cm}\mathcal M}\Big)=(9,1)
$$
 from which it comes that this web is AMP and has all his ARs polylogarithmic.
 
 Some geometric properties of Mantel's web (its symmetries, its definition in terms of configurations of points on the projective line) are discussed by Zagier in \cite[\S2.3.3]{ZagierSpecial}. From this, one deduces easily that Mantel's web coincides with the $Y$-cluster web of type $A_3$ which will be introduced and studied in the third part of this memoir. 
\mk

Since the term appearing in the RHS of \eqref{E:Mantel} is also a weight 2 iterated integral (but a trivial one in  a certain sense, since it is a square of a weight 1 such integral),  it is natural to wonder about the `complete web' 
$
\boldsymbol{\mathcal W}_{\hspace{-0.05cm}\mathcal M}^{\,c}=
\boldsymbol{\mathcal W}_{\hspace{-0.05cm}\mathcal M}
\cup 
\boldsymbol{\mathcal W}(x/y)
$ 
associated to \eqref{E:Mantel},  obtained from $\boldsymbol{\mathcal W}_{\hspace{-0.05cm}\mathcal M}$ by adjoining to it the foliation defined by the 
the first integral appearing as argument of the logarithmic term.\sk 

Again by direct computations, one gets 
$$
\rho^\bullet\Big(\boldsymbol{\mathcal W}_{\hspace{-0.05cm}\mathcal M}^{\,c}\Big)=(7,4,1)
\qquad  \mbox{ and } \qquad 
{\rm polrk}^\bullet\Big(\boldsymbol{\mathcal W}_{\hspace{-0.05cm}\mathcal M}^{\,c}\Big)=(10,2)
$$
thus  this web is AMP with polylogarithmic ARs as well. 

%
%
%
%
%
%

\paragraph{\hspace{-0.2cm}Rogers' multivariable dilogarithmic functional equations.}
\label{Par:RogersMultivariableAFEs}
In his paper \cite{Rogers}, Rogers did not only obtain the clean version \eqref{Eq:EFA-R} of the dilogarithmic five-terms relation but actually has established a whole series of such functional equation, one in $m$ variables with $m^2+1$ dilogarithmic terms, for any $m\geq 1$.   We follow the concise presentation of \cite{Zagier}  below (see \cite[\S1.7]{Lewin} for a more explicit treatment). \mk 

Let $f(t)=
\prod_{i=1}^m \big( 1- x_i \, t\big)-1$ 
be a degree $n$ polynomial without constant term (thus none of  the $x_i$'s is equal to 0). 
Given an indeterminate $y$, one denotes by $\zeta_1,\ldots,\zeta_n$ the $n$ roots of $f(t)=y$ viewed as algebraic functions of the $x_i$'s and of $y$.  Then there exists a constant $C$ such that the following identity holds true: 
$$
\Big(R(m)\Big) \hspace{4cm}
\sum_{i=1}^m \sum_{j=1}^m \l {2} \Big(x_i \, \zeta_j\Big)= \l { 2}  (y)+C\, .
\hspace{5cm} {}^{}
$$

To this identity is attached the  web defined by
the $m^2+1$ algebraic first  integrals appearing in it:
$$
\boldsymbol{\mathcal W}_{R(m)}=\boldsymbol{\mathcal W} \left( 
\, y\, ,\,  
 x_i \hspace{0.1cm} \zeta_j \hspace{0.1cm} \Big\lvert
 \, i,j=1,\ldots,m\, 
\right)\, .
$$

Actually, there is a subtlety in the definition of this web. Indeed, 
 the constant $C$ in the preceding identity 
is equal to $-\pi^2/6+\sum_{i,j} \l {2} (x_i/x_j)$ 
({\it cf.}\,\cite[\S1.7]{Lewin}) 
hence the quantities $x_i/x_j$ appearing as arguments of the bilogarithm
    in it should a priori  be considered as  first integrals as well for the web naturally associated to $(R(m))$. However, since the reciprocal $x_j/x_i$ of each such term also appears in the summation formula for $C$, these  terms compensate each others thanks to the inversion formula $\eqref{Eq:EFA-Lin}$. But each such  compensation actually holds true up to 
     linear combinations 
 of terms of the form $\log(x_i/x_j)^k$ with $k=0,1,2$. Hence, thanks to 
 Theorem \ref{T:FE-holomorphic}, the genuine dilogarithmic AFE truly associated to $ \boldsymbol{\mathcal W}_{R(m)}$ is not $(R(m))$ itself, but the associated identity for Rogers' dilogarithm
\begin{equation}
\label{Eq:(RogersR(m))}
\sum_{i=1}^m \sum_{j=1}^m  R \Big(x_i \, \zeta_j\Big)-R(y)={\rm cst.}
\end{equation}
(For another remark on this, see the case $m=3$ discussed below).\bk 

For  $m=2$, $R(2)$ is (equivalent  to) the five terms  relation thus $\boldsymbol{\mathcal W}_{R(2)}$ is (equivalent to) Bol's web. \mk 

The case $m=3$ has been explicited in \cite{GordonMcIntosh} (see p.\,434 in it, where it is proved that the associated dilogarithmic functional equation can be obtained from several copies of Abel's five terms relation). The
web  associated to $R(3)$ is 
$$\boldsymbol{\mathcal W}_{R(3)}=\boldsymbol{\mathcal W}\bigg( \, a, b, c, u, v, abc, \frac{ac}{u}, \frac{bc}{v}, \frac{av}{u}, \frac{bu}{v}\, \bigg)
$$
where the variables $a,b,c$ and $u,v$ are assumed to satisfy the two following relations: 
$$
av(1-bc)+bu(1-ac) = uv(1-ab)  \qquad \mbox{and}  \qquad v(1-a)+u(1-b) = 1-abc\, .
$$ 
After some computations, 
one gets 
$$\rho^\bullet(\boldsymbol{\mathcal W}_{R(3)})=(7,4,1)\, , \qquad {\rm polrk}^\bullet(\boldsymbol{\mathcal W}_{R(3)})=(10,1)
\qquad \mbox{ and } \qquad 
{\rm rk}(\boldsymbol{\mathcal W}_{R(3)})=11.$$ 
Thus, if  
this web 
has 
 polylogarithmic ARs, it is not AMP since 
its rank  is 
$\rho(\boldsymbol{\mathcal  W}_{R(3)})$  minus 1
 (one could say that it is of `{\it almost AMP rank}'). 
 Note that the equality $ {\rm polrk}^2(\boldsymbol{\mathcal W}_{R(3)})=1$ shows that 
 the AFE \eqref{Eq:(RogersR(m))} is the unique weight 2 polylogarithmic AR of 
$\boldsymbol{\mathcal W}_{R(3)}$ and that this web does not carry any dilogarithmic AR a component of which is the classical bilogarithm $\l { 2} $.  Indeed, otherwise we would have $ {\rm polrk}^2(\boldsymbol{\mathcal W}_{R(3)})>1$ according to  Lemma \ref{L:KWw-Sw-equivariant}.



\paragraph{\hspace{-0.2cm}Maier's dilogarithmic 8-web.}
\label{Par:Maier's-Web}
In \cite{M}, Maier has established  a functional equation which turns out to correspond 
to a dilogarihtmic identity, see \cite{MW}.  Indeed, in the latter paper, the authors explain how to get the following functional identity\footnote{The identity $(\mathcal M_8)$ holds true for $(x,y,z)$ in a certain open subset of $\mathbf C^3$. Replacing $R$ by Lewin's dilogarithm $L_2$ (see \eqref{Eq:Lewin-Dilog}) in it, one obtains a real functional relation  which is identically satisfied for any $x,y,z\in \mathbf R$.} in three independent variables $x,y$ and $z$ for 
Rogers dilogarithm 
$R(x)=\l {2} (x)+
{\rm Log}(x){\rm Log}(1-x)/2-\pi^2/6$
from two copies of the five-term relations:  
$$
(\mathcal M_8)\qquad 
R(x)+R( y)+ R\left(\frac{y(x-1)}{1-y}\right) +R \left(\frac{x(y-1)}{1-x}\right)
 - R(z)  - R\left(\frac{xy}{z}\right)  - R\left(\frac{xy(z-1)}{z-xy}\right)
- R\left(\frac{ z-xy}{z-1}\right) =0\,. 
$$
To this identity, 
one associates  the following  8-web in three variables
$$\boldsymbol{\mathcal W}_{\hspace{-0.05cm}\mathcal M_8}=
 \boldsymbol{\mathcal W}\left( x \, , \,  y \, , \, 
\frac{y(x-1)}{1-y}\, , \,
\frac{x(y-1)}{1-x}\, , \,
z \, , \,  
\frac{xy}{z}\, , \,
\frac{xy(z-1)}{z-xy}\, , \,
\frac{ z-xy}{z-1}\, 
\right) 
$$
whose invariants are:  $\rho^\bullet(\boldsymbol{\mathcal W}_{\hspace{-0.05cm}\mathcal M_8})=(5,3,1)$  (thus $\rho(\boldsymbol{\mathcal W}_{\hspace{-0.05cm}\mathcal M_8})=9$)  and 
${\rm polrk}^\bullet\Big(\boldsymbol{\mathcal W}_{\hspace{-0.05cm}\mathcal M_8}\Big)=(8,1)$. It follows that 
$\boldsymbol{\mathcal W}_{\hspace{-0.05cm}\mathcal M_8}$ 
is AMP with only weight 1 or 2 polylogarithmic ARs.

\subparagraph{}
\hspace{-0.6cm}
From \cite{MW}, one can deduce a nice geometric (that is coordinate free) construction of 
$\boldsymbol{\mathcal W}_{\hspace{-0.05cm}\mathcal M_8}$ which 
 will prove to be useful later on page 
 \pageref{P:WM8-WBV(4)}
 to compare this web with another one. \sk 

Let $\boldsymbol{\mathcal B}$ be Bol's web viewed geometrically as the web $\boldsymbol{\mathcal W}(f_1,\ldots,f_5)$ on the moduli space 
$\mathcal M_{0,5}$ defined by the five forgetful maps $f_i:\mathcal M_{0,5}\rightarrow \mathcal M_{0,4}\simeq \mathbf P^{1}\setminus \{0,1,\infty\}$.  Taking other but independent copies of these objects, denoted by $\boldsymbol{\mathcal B}'$, $\mathcal M_{0,5}'$ and $f_i'$ for $i=1,\ldots,5$, one can define the `{\it fiber product of $\boldsymbol{\mathcal B}$ and $\boldsymbol{\mathcal B}'$ over $\mathcal M_{0,4}$ or more precisely, with respect to  a pair $\psi=(f_i,f_{i'}')$ of two rational first integrals, one for each web}'. This fiber product will be denoted by 
$\boldsymbol{\mathcal B}\times_\psi \boldsymbol{\mathcal B}'$.
\sk 

 To explain the construction, we assume that $\psi=(f_5,f_5')$ and  one denotes by 
$\pi$ (resp.\,by $\pi'$) the projection of the product 
$\mathcal M_{0,5}\times \mathcal M_{0,5}'$ onto its first (resp.\,its second) factor.   Then by definition, the fiber product $ \boldsymbol{\mathcal B}\times_{\psi} \boldsymbol{\mathcal B}'$ is the restriction to the hypersurface $\Sigma_\psi$ cut out by $f_5\circ \pi=f_5' \circ \pi'$,  of the product web 
$\boldsymbol{\mathcal B} \times \boldsymbol{\mathcal B}'=\pi^*(\boldsymbol{\mathcal B})\cup (\pi')^*(\boldsymbol{\mathcal B}')$. 
In other terms, one has 
$$ \boldsymbol{\mathcal B}\times_{\psi} \boldsymbol{\mathcal B}'=
\boldsymbol{\mathcal W}\Big( \hspace{0.1cm} f_i\circ \pi \, , \,  f_{i'}'\circ \pi' 
\, , \,  f_5\circ \pi \hspace{0.1cm}
\, \Big\lvert 
\,\hspace{0.1cm} i,i'=1,\ldots,4\,\hspace{0.1cm} \Big)\, \bigg\lvert_{\Sigma_\psi}\, .
$$

Since two foliations of the 10-web $\boldsymbol{\mathcal B} \times \boldsymbol{\mathcal B}'$ have been identified in order to define $\boldsymbol{\mathcal B}\times_{\psi} \boldsymbol{\mathcal B}'$, it comes that the latter is a 9-web defined on a space of dimension 3. It is easy to show that (wGP) is satisfied, and also that  $ \boldsymbol{\mathcal B}\times_{\psi} \boldsymbol{\mathcal B}'$ does not depend  on $(i)$ the chosen pair $\psi= (f_i,f_{i'}')$, 
and on $(ii)$ the identification between the target spaces of the components of $\psi$ neither.  Then denoting by $\mathcal F(\psi)$ the foliation defined by the restriction of 
$f_i\circ \pi$ on $\Sigma_\psi$ (which,  by definition, coincides with that of  $f_{i'}'\circ \pi'$), we let the reader verify that (up to equivalence) Maier's  web is nothing else but the 8-web obtained by removing 
$\mathcal F(\psi)$ from $ \boldsymbol{\mathcal B}\times_{\psi} \boldsymbol{\mathcal B}'$, {\it i.e.}\,one has: 
$$
\boldsymbol{\mathcal W}_{\hspace{-0.05cm}\mathcal M_8} \simeq 
 \Big(\boldsymbol{\mathcal B}\times_{\psi} \boldsymbol{\mathcal B}'\Big)\setminus 
\mathcal F(\psi)=
\boldsymbol{\mathcal W}\Bigg(\, \hspace{0.1cm} \Big(f_j\circ \pi\Big)
\Big\lvert_{\Sigma_\psi}
 \, , \, \hspace{0.1cm}  \Big(f_{j'}'\circ \pi' \Big)\Big\lvert_{\Sigma_\psi}
\hspace{0.1cm} \Big\lvert 
\, j,j'=1,\ldots,5\,, j\neq i, \, j'\neq i' \,  \Bigg)
 \, .
 $$

\paragraph{\hspace{-0.2cm}Some recent dilogarithmic identities.}
Starting from the late 1980s, it has been gradually realized that dilogarithmic identities naturally arise 
in several domains of mathematics, such as number theory (the $K$-theory of number fields more precisely), some integrable systems coming from mathematical physics, 
  hyperbolic geometry, and also and rather systematically within the theory of cluster algebras.
 \mk

We present below some dilogarithmic identities obtained recently within these different contexts, with the exception of that of the cluster algebras, which will be studied more systematically further on. We only discuss some identities which we have found interesting from the point of view of web geometry. By no means the lines below constitute a   complete overview of the subject. 
 



\subparagraph{}
\hspace{-0.6cm}
\label{Par:Dilog_varphi}
The  algebraic symbolic approach to deal with functional equations of polylogarithms 
({\it cf.}\,\S\ref{Par:II-AR} and more specifically 
\S\ref{Par:HolomSymbolicRA} above) 
can be traced back to Bloch (see \cite{Bloch} and also \cite[Cor.\,6.2.2]{BlochBook}) for the dilogarithm and has been formalized and generalized to any weight  by Zagier,  in relation with a conjecture of its own about the K-theory of number fields, see {\it e.g.}\,\cite{Zagier1991,Zagier}.  Among many other things, this led  to the following generalization of Rogers' multivariables dilogarithmic functional equation (see \cite{Zagier1989})\footnote{This generalization of Rogers' multivariable dilogarithmic identity has been obtained independently almost at the same time by 
 Wechsung and 
  Wojtkowiak,  see the  eighth and  tenth  chapters of  \cite{LewinStruct} respectively. 
 See also the  seventh chapter of this book for an intermediary version, between Rogers' and Zagier's one,  due to Ray.}, which is more convenient to state using Bloch-Wigner dilogarithm $D$ (which,  it should be recalled, is defined at every point of the Riemann sphere $\mathbf P^1$). 
   Let $\varphi \in \mathbf C(x)$ be a sufficiently generic rational function
 of positive degree  $m>0$.  
  Then for $z$ generic in $\mathbf P^1$,  setting 
 $\xi=\varphi(z)$, one  has 
 \begin{equation}
  \label{Eq:Zagier-varphi}
  \sum_{
  \scalebox{0.7}{
 \begin{tabular}{l}
 $ z_1 \in \varphi^{-1}(1) $\\
  $ z_0 \in \varphi^{-1}(0) $\\
   $ z_\infty \in \varphi^{-1}(\infty) $
 \end{tabular}  }}
D\Big([z,z_1,z_0,z_\infty ] \Big)= m\,
D(\xi)\, . 
\end{equation}
  This formula being identically satisfied with respect to $z$ and 
the coefficients of  
 $\varphi$, it  can be seen as a (global and real) AR on the $(m^3+1)$-web defined on the space of pairs 
 $(z,\varphi)$  (which is of complex dimension $2m+2$) by 
 the projection 
 $(z,\varphi)\mapsto 
 z$ and  
 the $m^3$ algebraic first integrals $(z,\varphi)\mapsto [z,z_1,z_0,z_\infty ] \in \mathbf P^1 $ with 
  $z_1\in  \varphi^{-1}(1)$, $z_0\in  \varphi^{-1}(0)$ and   $z_\infty\in  \varphi^{-1}(\infty)$. 
 \mk 
 
 In the  case when $\varphi$ is assumed to be a polynomial,  $\varphi^{-1}(\infty)$  is reduced to $\infty$ (but with multiplicity $m$), the summation over $\varphi^{-1}(\infty)$ 
 in \eqref{Eq:Zagier-varphi}  is trivial and  one recovers Rogers' equation $R(m)$.  
 Since the web associated to the latter is not AMP for $m=3$, we suspect that the same holds true for the web associated to \eqref{Eq:Zagier-varphi} for any $m\geq 3$. 
 %
%
\subparagraph{}
\hspace{-0.6cm}
Another class of dilogarithmic identities are related to 
certain integrable systems coming from mathematical physics
 ($Y$-systems associated to the `{\it Thermodynamic Bethe Ansatz}' in CFT/QFT, see \cite{FrenkelSzenes} or \cite{Kirillov} for more background and references).   \mk 
 
The first important result  within this field from our perspective has been obtained  in \cite{FrenkelSzenes} (and also independently and almost at the same time in \cite{GliozziTateo}) where Frenkel and Szenes have proved 
 in type $A_n$ for any positive integer $n$, 
some previously conjectured results.  In particular, they have established that   (an equivalent form of) the following identity  with $n(n+3)/2$ dilogarithmic terms
\begin{equation}
\label{Eq:Type-An}
\sum_{\substack{i,j=1 \\ \lvert i-j\lvert \geq 2}}^n R\Big([x_i,x_{i+1},x_j,x_{j+1}]\Big)=(n-3)\frac{{}^{}\hspace{0.15cm} \pi^2}{6}
\end{equation}
holds true for any  $x_1, \ldots,x_{n+3}\in \mathbf R$ listed modulo  $n+3$ ({\it i.e.}\,$x_{n+4}=x_1$, $x_{n+5}=x_2$, etc.)
\mk

For $n=2$,  \eqref{Eq:Type-An} corresponds to the five terms relation \eqref{Eq:Efa-(Ab)-symmetric}  
whereas when $n=3$, it is equivalent to the dilogarithmic identity   for Rogers dilogarithm $R$ associated to Mantel's functional equation \eqref{E:Mantel}.
For $n$ arbitrary, the identity \eqref{Eq:Type-An} has been obtained 
again, 
 with no reference to any previous work, first by Bridgeman in  \cite{Bridgeman} (see also \cite[\S5.2]{BridgemanTan}) in relation with the theory of hyperbolic surfaces and more recently in \cite{Souderes} where   Soud\`eres  proves that for any $n\geq 2$, the identity \eqref{Eq:Type-An} can be obtained as a linear combination of the classical five-terms relation.
\sk 

The relation \eqref{Eq:Type-An} can be seen as an AR for a web in $n$ variables defined by the $n(n+3)/2$ cross-ratios appearing as arguments of $R$ in it.  This web will be interpreted in the next part as the `{\it $\mathcal Y$-cluster web of type $A_n$}', denoted by 
$\boldsymbol{\mathcal Y\hspace{-0.06cm}\mathcal W}_{\hspace{-0.1cm} A_n}$, 
 and  will be proved to be AMP with only polylogarithmic abelian relations of weight less than or equal to 2 {({\it cf.}\,Theorem \ref{T:YW-An} below)}.

\subparagraph{}
\hspace{-0.6cm}
\label{SPar:Z(m,n)}
The fact that   volumes in hyperbolic geometry  can 
in some cases be written in terms of values of polylogarithmic functions is quite classical (at least in the case of surfaces) and has been recognized  as a source of dilogarithmic identities (see \cite{BridgemanTan} for a nice survey). \mk 

The known such identities in the case of hyperbolic surfaces, one due to Bridgeman \cite{Bridgeman}, the other to Luo and Tan \cite{LuoTan}, 
 are not of `finite type' in general, in the sense that they involve an infinite summation of dilogarithmic terms. Hence these identities can be seen as  kind of ARs for some  webs defined on moduli spaces of hyperbolic surfaces of certain types,  but webs constitued of 
  an infinite number of foliations.    As far as we know, the unique dilogarithmic identity with a finite number of terms which can be deduced in the case of hyperbolic surfaces is \eqref{Eq:Type-An} which appears as the specialization of Bridgeman's 
    identity (which is of infinite type in general) to the case of an ideal hyperbolic $n$-gon (compare equations (1) and (6) in \cite{Bridgeman}).\mk

\begin{rem}
\label{Rem:Bridgeman}
Actually in {\rm \cite{Bridgeman}}, Bridgeman obtains \eqref{Eq:Type-An} as a very particular case of a more general dilogarithmic identity with infinitely many terms.  Namely, 
given any hyperbolic surface  $S$  with  finite area and geodesic boundary
endowed with 
a geodesic lamination $\lambda$, one has 
\begin{equation}
\label{Eq:Bridgeman-infinite-sum}
\sum_{i} \mathcal L\left( \frac{1}{{\rm cosh}^2(\ell_i/2)}\right)=\frac{\pi^2}{12}\Big( 6 \chi(S)-N_{\lambda}\Big)\, , 
\end{equation}
where the summation is carried over the `$\lambda$-orthospectrum',  $N_\lambda$ stands for the number of so-called `$\lambda$-cusps', with $\mathcal L$ denoting the function defined by $\mathcal L(x)={\bf L}{\rm i}_2(x)+\frac{1}{2}{\rm Log}\lvert x\lvert\,{\rm Log}(1-x)$ for $x\leq 1$.
\sk

 In general, the orthospectrum is a countable set hence  the sum in 
  \eqref{Eq:Bridgeman-infinite-sum} is infinite. 
This identity, called `Bridgeman orthospectrum identity', can be seen as an  abelian relation for the infinite web 
 defined by the $\ell_i$'s, these quantities being considered as analytic functions on a suitable moduli space of hyperbolic surfaces endowed with a geodesic lamination. 
 \sk

Several other identities of this kind exist within the theory of hyperbolic surfaces, e.g.\,`Mac-Shane' or `Luo-Tan' identities (see {\rm \cite{BridgemanTan}} for a general overview) 
and all might be interpreted as ARs for some infinite webs.  It would be interesting to begin a  more systematic study of these identities from the perspective of web geometry but we leave this aside for future works (possibly by others). We only mention that identities as \eqref{Eq:Bridgeman-infinite-sum} are similar by several aspects, to conjectural  identities which might be obtained  as semiclassical limits of certain quantum dilogarithmic identities we will say a few words about in the very last part of this memoir (see \S\ref{Par:QdilogIdentInfiniteNumberOfTerms} more precisely).

\end{rem}

\begin{center}
$\star$
\end{center}

  The case of hyperbolic 3-folds is more interesting than the case of surfaces for the purpose 
  of getting 
   (finite) dilogarithmic identities.   We refer to \cite{GliozziTateo,NeumannZagier} and to \cite{Zagier} for more details on what follows.\mk

  We can trace back to \href{https://en.wikipedia.org/wiki/Nikolai_Lobachevsky}{Lobatchevsky} the fact that the volume of an ideal tetrahedron $\tau$ in the hyperbolic 3-space $\mathbf H^3$ is expressed as the value of a dilogarithmic function (Bloch-Wigner function $D$ more precisely) evaluated at the cross-ratio of the four ideal vertices  $v_1(\tau),\ldots,v_4(\tau)$ of $\tau$ viewed as four elements of $\boldsymbol{\partial} \mathbf H^3\simeq \mathbf P^1$.  Now given a hyperbolic threefold $M$ of finite volume (possibly with cusps  and/or a geodesic boundary) endowed with an ideal triangulation $\mathcal T$, its volume is expressed as 
  $$
  {\rm Vol}(M)=\sum_{i=1}^n D\Big( [v_1(\tau_i),\ldots,v_4(\tau_i)]  \Big)
  $$
  where the summation is  over all the tetrahedra 
 $\tau_1,\ldots,\tau_n$  
  of $\mathcal T$.  \mk 
  
  Then,  considering this identity on the moduli space of triangulated hyperbolic 3-folds $(M',\mathcal T')$ with ${\rm Vol}(M')={\rm Vol}(M)$ and $\mathcal T'$ combinatorially equivalent to $\mathcal T$ gives a dilogarithmic AFE  (with only a finite number of terms).  An interesting explicit example which has been worked out by Zagier in \cite[II.2.C]{Zagier},  is the join of a $m$-gon with a $n$-gon for any two integers $m,n\geq 3$.  This case gives rise to the following identity 
$$ 
\Big(\mathcal Z(m,n)\Big) \qquad 
\qquad \quad 
\sum_{\substack{j =1,\ldots,m\\ k=1,\ldots,n}}
R\left( \frac{(y_k-y_{k+1})(x_j-x_{j+1})}{({y_k}-x_j)(y_{k+1}-x_{j+1})} 
\right)=  N_{m,n}\,  \frac{{}^{}   \, \pi^2}{6}
\qquad \qquad \qquad \qquad \qquad   {}^{}
$$
satisfied for any real numbers $x_1,\ldots, x_{m}$ and $y_1,\ldots,y_n$ (with $x_{m+1}=x_1$ and $y_{n+1}=y_1$), where $N_{n,m}$ stands for an integer which depends on $m,n$ and on the relative position of the $x_i$'s and the $y_j$'s on the real axis. 
   To ${\mathcal Z}(m,n)$ is naturally associated the following web 
$$ 
\boldsymbol{\mathcal W}_{\hspace{-0.05cm}\mathcal Z(m,n)}=
\boldsymbol{\mathcal W}\left(\, 
\,   \frac{(y_k-y_{k+1})}{({y_k}-x_j)} \, \frac{(x_j-x_{j+1})}{(y_{k+1}-x_{j+1})} 
\, 
\hspace{0.15cm}\bigg\lvert \, 
\scalebox{0.94}{\begin{tabular}{l}
$j=1,\ldots,m$
\vspace{-0.0cm}
\\
$ k=1,\ldots,n$
\end{tabular}}
\right)
$$
which is then a 
$mn$-web in $m+n$ variables.\footnote{Let us mention that the relation $\mathcal Z(m,n)$ is also satisfied when $m=2$ for any $n\geq 2$ but in this case there is no threefold involved  and the  foliations defined by the 
rational first integrals  appearing in the definition of  $\boldsymbol{\mathcal W}_{\hspace{-0.05cm}\mathcal Z(m,n)}$ do not satisfy the `weak general position' assumption discussed in \S\ref{S:Webs} hence there is no interesting web to deal with. Note also that  thanks to the projective invariance of the cross-ratio, the web $\boldsymbol{\mathcal W}_{\hspace{-0.05cm}\mathcal Z(m,n)}$ (for any $m,n\geq 3$) can easily be seen as the pull-back of a $mn$-web in $m+n-3$ variables.}\sk 

By means of direct computations, we have verified that the following equalities  
$$\rho^\bullet\Big(
\boldsymbol{\mathcal W}_{\hspace{-0.05cm}\mathcal Z(m,n)}\Big)=
\Big((m-1)(n-1)+2, m+n-3,1\Big)\, , 
\hspace{0.2cm}
\rho\Big(
\boldsymbol{\mathcal W}_{\hspace{-0.05cm}\mathcal Z(m,n)}\Big)=
mn+1\, , 
\hspace{0.2cm}
{\rm polrk}^\bullet\Big(
\boldsymbol{\mathcal W}_{\hspace{-0.05cm}\mathcal Z(m,n)}\Big)=
\Big(mn,1\Big)
$$
hold true for any $m,n\leq 12$ and we conjecture that it is actually the case in full generality.  If true, the 
$\boldsymbol{\mathcal W}_{\hspace{-0.05cm}\mathcal Z(m,n)}$'s form another family of 
AMP webs with only weight 1or  2 polylogarithmic ARs. \mk

\subparagraph{}
\hspace{-0.6cm}
\label{SPar:Bytsko-Volkov}
Finally, we mention the new family of dilogarithmic identities obtained recently by Bytsko and Volkov in  \cite{BV} and discuss some properties of the web associated to each of them. \mk 

In the lines below:  $n$ stands for an integer bigger than or equal to 4,  $T_n$ denotes the set  of 3-tuples $(a,b,c)\in \mathbf N^3$ such that $1\leq a < b<c \leq n$ and one sets $d_n=2\lvert T_n\lvert=2 \scalebox{0.9}{${ n \choose 3}$}=
n(n-1)(n-2)/3$.
\sk

Let $N_n^+$ be the space of   totally positive upper-triangular $n\times n$ unipotent matrices.  There exists an explicit rational map   (see \cite[Lemma 1]{BV}) 
$$
\mathbf R^{\scalebox{0.8}{${ n \choose 2}$}} \longrightarrow M_n(\mathbf R)\, , \quad \boldsymbol{x}=(x_{ij})_{1\leq i<j\leq n}\longmapsto M(x)
$$
 which induces a bijection between the set of {\it `positive'} $ \boldsymbol{x}\in \mathbf R^
 \frac{n(n-1)}{2}
 $ (that is $ \boldsymbol{x}=(x_{ij})$ with $x_{ij}>0$ for any $i<j$) and $N_n^+$ ({\it cf.}\,\cite[Prop.\hspace{0.05cm}1.7]{BFZ}). 
 Thus $N_n^+$ is a smooth submanifold of $N_n$, and the $x_{ij}$'s, 
named the {\it Jacobi coordinates}',  form a global system  of coordinates on it. 
 Then Bytsko and Volkov consider two families  $Y_{abc}(M)$ and $Y^{abc}(M)$ indexed by the elements of $T_n$, of  ratios of products of 
  generalized minors of $M\in N_n^+$, which behave nicely with respect to two involutions $M\mapsto M'$ and $M\mapsto M''$ defined on $N_n^+$. \mk 
 

One of the main results of \cite{BV} is that the following identity holds true uniformly on $N_n^+$: 
$$
\Big({\mathcal B}{\mathcal V}(n)\Big)\quad :  \hspace{2cm}
\sum_{(a,b,c)\in T_n} 
{\mathsf R}\Big( Y_{abc}\Big) \, - \hspace{-0.3cm}\sum_{(a,b,c)\in T_n}
{\mathsf R}
\Big( Y^{abc}\Big)=0\, , \hspace{5cm} {}^{}
$$
where ${\sf R}$ stands for the cluster dilogarithm defined above  in \S\ref{Par:NotationDilogarithmicFunctions}.
\mk


This can be seen as a dilogarithmic AR for what we call the {\bf $\boldsymbol{n}$-th Bytsko-Volkov web}
$$
\boldsymbol{\mathcal W}_{\hspace{-0.05cm}{\mathcal B}{\mathcal V}(n)}=
\boldsymbol{\mathcal W}\left( \, Y_{abc}\, , \, Y^{abc}\, \Big\lvert \hspace{0.1cm} 
1\leq a < b < c \leq n\, 
\right)
$$
which  can be seen as a 
 web  in the Jacobi coordinates $x_{ij}$, hence a web in $\scalebox{0.85}{${ n \choose 2}$}=n(n-1)/2$  
 variables.\mk 

The $\boldsymbol{\mathcal W}_{\hspace{-0.05cm}{\mathcal B}{\mathcal V}(n)}$'s form a family of webs which seem to be interesting in what concerns their rank and their abelian relations. Indeed, some explicit computations for  the small values of  $n$ 
 lead us to conjecture that the following formulae hold true for any $n \geq 4$: 
$$
\rho^\bullet\Big(\boldsymbol{\mathcal W}_{{\mathcal B}{\mathcal V}(n)} \Big)
= 
\left(  
\, 
d_n
- \scalebox{1}{${ n -1\choose 2}$}  \, , \, 
\scalebox{1}{${ n -1\choose 2}$}
 \, , \, 1\right)
\qquad \mbox{ and } \qquad 
{\rm polrk}^\bullet\Big( 
\boldsymbol{\mathcal W}_{{\mathcal B}{\mathcal V}(n)}
\Big)=\Big( d_n,1\Big)\,.
$$
If true (and we have verified that it is indeed the case for $n\leq 9$),
 this would imply that the $\boldsymbol{\mathcal W}_{{\mathcal B}{\mathcal V}(n)}$'s form a series of AMP webs with only logarithmic and dilogarithmic ARs.  \bk 

The case $n=4$ is interesting and completely explicited in \cite[\S2]{BV}.  
The authors prove that 
   ${\mathcal B}{\mathcal V}(4)$ can actually be written as the following identity involving only three variables $x,y$ and $z$ (which can be written explicitly in terms of the Jacobi coordinates $x_{ij}$): 
\begin{align*}
\Big({\mathcal B}{\mathcal V}'(4)\Big)\hspace{0.1cm}:  \hspace{0.2cm}
{\mathsf R}\Bigg(
 \frac{x}{1+y}  \Bigg)& \, +  {\mathsf R}\Bigg( \frac{y}{1+x}   \Bigg)
+ 
{\mathsf R}\Bigg(
\frac{z(1+x+y)}{(1+x)(1+y)}   \Bigg)
+{\mathsf R}\Bigg( \frac{xy}{(1+x+y)(1+z)}  \Bigg)
\\
\nonumber 
& \, -{\mathsf R}\Bigg(  \frac{x}{1+z}  \Bigg)
-{\mathsf R}\Bigg( \frac{z}{1+x}   \Bigg)
-{\mathsf R}\Bigg( \frac{y(1+x+z)}{(1+x)(1+z)}   \Bigg)
-{\mathsf R}\Bigg(  \frac{xz}{(1+x+z)(1+y)} 
\Bigg)=0\, .
\end{align*}
This shows that $\boldsymbol{\mathcal W}_{\hspace{-0.05cm}{\mathcal B}{\mathcal V}(4)}$ actually is the pull-back of the web in three variables, denoted by 
$\boldsymbol{\mathcal W}_{\hspace{-0.05cm}{\mathcal B}{\mathcal V}(4)}'$, 
 defined by the  rational functions appearing as arguments of ${\mathsf R}$ in the previous identity. 
  It is natural to wonder if something similar holds true for any $n\geq 4$.
 The conjectural equality   $\rho^1(\boldsymbol{\mathcal W}_{{\mathcal B}{\mathcal V}(n)})
= d_n
-\scalebox{0.8}{${ n -1\choose 2}$}$ and 
the fact that this is indeed the case when $n=4$  leads us to believe that,  for any $n\geq 4$, $\boldsymbol{\mathcal W}_{{\mathcal B}{\mathcal V}(n)}$ might be a pull-back of a $d_n$-web 
$\boldsymbol{\mathcal W}_{\hspace{-0.05cm}{\mathcal B}{\mathcal V}(n)}'$
in only  
\scalebox{0.8}{${ n -1\choose 2}$} variables. 
For 
 $n\leq 10$, we have computed the intrinsic dimension of 
$\boldsymbol{\mathcal W}_{{\mathcal B}{\mathcal V}(n)}$ and verified  that 
 it indeed coincides with $\scalebox{0.8}{${ n -1\choose 2}$}$.
\begin{center}
\vspace{-0.2cm}$\star$
\end{center}

Another interesting property of the webs $\boldsymbol{\mathcal W}_{\hspace{-0.05cm}{\mathcal B}{\mathcal V}(n)}$ (and possibly of the $\boldsymbol{\mathcal W}_{\hspace{-0.05cm}{\mathcal B}{\mathcal V}(n)}'$'s if they exist) is that they all are defined by first integrals which are positive (or substraction-free) real rational functions.  Since this is a property that the cluster webs we will study in the next part all satisfy as well, 
it would be interesting to know whether or not  Bytsko-Volkov's webs are webs associated to periods of some cluster algebras (one period and one clusters algebra for each $n$). 
\subparagraph{}
\hspace{-0.6cm}
 \label{P:WM8-WBV(4)}
We remark  that the functional equations  $({\mathcal B}{\mathcal V}'(4))$ just above and Maier's 
$(\mathcal M_8)$ are formally similar since both are written 
as the sum of four dilogarithmic terms minus the sum of four others. 
In addition to this, the fact that both associated webs 
$\boldsymbol{\mathcal W}_{\hspace{-0.05cm}{\mathcal B}{\mathcal V}(4)}'$
and 
$\boldsymbol{\mathcal W}_{\hspace{-0.05cm}{\mathcal M_{8}}}$
are AMP with only polylogarithmic ARs and have the same invariants $\rho^\bullet$ and ${\rm polrk}^\bullet$ naturally makes us wonder whether these two webs are equivalent or not. \mk

 As it happens, these two webs are equivalent.  To verify this, one starts by noticing that 
  Bytsko-Volkov web when $n=4$ can be decomposed into the union of two disjoint  4-subwebs, each of intrinsic dimension 2: one has 
  $\boldsymbol{\mathcal W}_{\hspace{-0.05cm}{\mathcal B}{\mathcal V}(4)}'
  =
  \boldsymbol{\mathcal W}_4
  ' \sqcup \boldsymbol{\mathcal W}_4
   ''$ 
  with 
\begin{align*}
\boldsymbol{\mathcal W}_4
 '= &\, \boldsymbol{\mathcal W}
 \Bigg(
\, \frac{x}{y+1}, \frac{z(1+x+y)}{(x+1)(y+1)}\, , \, \frac{z}{1+x}\, , \, \frac{xz}{(1+x+z)(y+1)}\,
  \Bigg) \\
  \mbox{ and }\quad 
  \boldsymbol{\mathcal W}_4
  ''= &\, \boldsymbol{\mathcal W}
  \Bigg(
  \frac{xy}{(1+x+y)(z+1)} \, , \, 
  \frac{y}{1+x} \, , \, 
    \frac{x}{z+1} \, , \, 
      \frac{y(1+x+z)}{(1+x)(z+1)} \,\Bigg)\, .
\end{align*}

Now introducing to pairs of new coordinates,  $(a,b)=(-x/(y + 1), -(1 + x)/z)$ for $\boldsymbol{\mathcal W}_4
'$
and 
$(\alpha ,\beta)= (-x/(z + 1), -(1 + x)/y)$
for $\boldsymbol{\mathcal W}_4''$, one has 
$$
\boldsymbol{\mathcal W}_4
 '=
\boldsymbol{\mathcal W}\left(\, -a\, ,\,  \frac{a-1}{b}\, , \,  \frac{a}{b-1}\,  ,\,  -\frac{1}{b} \, 
\right)
\qquad \mbox{ and }
\qquad 
\boldsymbol{\mathcal W}_4
 ''=
\boldsymbol{\mathcal W}\left(\frac{\alpha}{\beta-1}, -\frac{1}{\beta}\, , \, -\alpha \, , \frac{\alpha-1}{\beta}\, \right)$$
After some easy  polynomial eliminations, one arrives to the following relation between the four variables $a,b,\alpha,\beta$: 
$$\frac{ab}{a+b-1}=\frac{\alpha\beta}{\alpha+\beta-1}\, .$$

Since $\boldsymbol{\mathcal W}_4' \sqcup \boldsymbol{\mathcal F}({ab}/({a+b-1}))= \boldsymbol{\mathcal W}(-a, (a-1)/b, a/(b-1) , 1/b, {ab}/({a+b-1}))$ is equivalent to Bol's web (and similarly for $\boldsymbol{\mathcal W}_4'' \sqcup \boldsymbol{\mathcal F}({\alpha\beta}/({\alpha+\beta-1}))$, one can deduce from this the same geometric construction of 
$\boldsymbol{\mathcal W}_{\hspace{-0.05cm}{\mathcal B}{\mathcal V}(4)}'$ 
by means of a fiber product of two copies of Bol's web as that of Maier's web $\boldsymbol{\mathcal W}_{\hspace{-0.05cm}{{\mathcal M}_8}}$ given at the end of \S\ref{Par:Maier's-Web}. This proves our claim that 
$\boldsymbol{\mathcal W}_{\hspace{-0.05cm}{\mathcal B}{\mathcal V}(4)}'$  and Maier's 8-web are equivalent. In particular, this implies that the functional equation $({\mathcal B}{\mathcal V}'(4))$ is just Maier's one but written in some other coordinates. 


{{
\paragraph{Some remarks about dilogarithmic AFEs.}
\label{Par:Remarks-Dilog-AFE}
To conclude our discussion of dilogarithmic functional equations, we would like two discuss them from the perspective of being AMP which is the new (albeit elementary) concept of web geometry we introduce in this memoir.}
\sk

%


When elaborating this text, we have considered a fairly large number of webs 
carrying dilogarithmic ARs. This led us to observe a certain phenomenon that could well be verified in a rather general way and that we are going to discuss quickly now. \mk 

It is convenient here to use the formalism of webs. Recall that an 
  AR of a web $\boldsymbol{\mathcal W}$ is said to be
  \begin{itemize}
  \item  {\bf proper} 
 if it cannot be obtained as a linear combination of ARs of proper subwebs of $\boldsymbol{\mathcal W}$;  In particular, such an AR is necessarily {\bf complete}, {\it i.e.}\,none of its components is trivial; 
\item {\bf dilogarithmic} if when you choose certain first integrals  $\boldsymbol{\mathcal W}$,  its components are dilogarithmic functions, 
{\it i.e.}\,multiples of  weight 2 iterated integrals with symbol $
 \boldsymbol{\mathcal S}(R)= \frac{1}{2}(01-10)$.
 \end{itemize}

 
 The study of many webs carrying dilogarithmic ARs led us to state the following 
 \mk 


{\bf 
 ``$\boldsymbol{\rho^3}$-dilogarithmic 
 conjecture'':} {\it Let $ \boldsymbol{\mathcal W}$ be a web 
 carrying a proper dilogarithmic AR. Then}
  \begin{itemize}
\item[{\bf (1).}] \quad {\it 
${\rm IntrDim}(\boldsymbol{\mathcal W})\geq 3$ 
\, $\Longrightarrow $  \hspace{0.2cm}
$\rho^3( \boldsymbol{\mathcal W})>0$ \quad  and  \quad $\rho^\sigma( \boldsymbol{\mathcal W})=0$ \quad  for  any $\sigma\geq 4$.}
 \item[{\bf (2).}]  \quad {\it ${\rm polrk}^2(\boldsymbol{\mathcal W})=1$ \hspace{0.3cm} $\Longrightarrow $  \hspace{0.2cm}
 $\rho^3( \boldsymbol{\mathcal W})=1$ \quad and  \quad $\rho^\sigma( \boldsymbol{\mathcal W})=0$ \quad for  any $\sigma\geq 4$.}
\mk 
\end{itemize}


We first remark that the hypothesis that $\boldsymbol{\mathcal W}$ carries a non-trivial dilogarithmic AR (not necessarily assumed to be proper) implies quite directly that $\rho^3(\boldsymbol{\mathcal W})>0$. Indeed, assuming that this web is defined by first integrals $u_i:\Omega\rightarrow \mathbf C$ on a domain $\Omega$ of $V=\mathbf C^m$ with $m\geq 2$, our assumption is that for $\omega\in \Omega$ generic (in particular such that  $\omega_i=u_i(\omega)$ be distinct from 0 or 1 for every $i=1,\ldots,d$), there 
exist some complex constants $c_1,\ldots,c_d\in \mathbf C$ 
(not all trivial) 
such that 
the 
AFE  $(E)$ : 
$\sum_{i=1}^d c_i\, {\mathsf R}^i(u_i)\equiv 0$ holds true 
 in the vicinity of $\omega^*$ ({\it cf.}\,Theorem \ref{T:FE-holomorphic}.2). Here for each $i$,  ${\mathsf R}^i$ stands for the iterated integral centered at $\omega_i$ with symbol $
 \boldsymbol{\mathcal S}({\mathsf R}^i)=
 (01-10)/2$, that is explicitly
$${\mathsf R}^i(u)=\frac{1}{2} \bigintss_{\omega_i}^u \left( \, 
\frac{\Big(\log(1+\tau)-\log(1+\omega_i)\Big)}{\tau}-\frac{\Big(\log(\tau)-\log(\omega_i)\Big)}{1+\tau}
\, \right)d\tau 
\,. 
$$ 
For any $i$, let $\ell_i$ be the linear form on $V$ corresponding to the differential 
$du_i(\omega)$ of $u_i$ at $\omega$ via the standard identification $V^*\simeq T_\omega^* V$.  Now since for any $i=1,\ldots,d$, one has
$$\hspace{0.5cm} {\mathsf R}^i(u)=\frac{-(u-\omega_i)^3}{12\,\omega_i^2\,(1+\omega_i)^2}+O\Big( \,(u-\omega_i)^4\Big)\, 
\hspace{0.15cm}
\quad 
\mbox{at}\hspace{0.25cm} \omega_i\, ,  
$$
 it comes   that  ${\mathsf R}^i(u_i)=
 {-(\ell_i)^3}/({12\,\omega_i^2\,(1+\omega_i)^2})+O( (\ell_i)^4)$ at $\omega$ for every $i$ which, together with the AFE $(E)$ gives us that the $c_i$'s are such that 
 the identity 
  $\sum_{i=1}^d {c_i}/({12\,\omega_i^2\,(1+\omega_i)^2})\,(\ell_i)^3
 =0$
  holds true in ${\rm Sym}^3(V^*)$.  Since some of the $c_i$'s are non zero, 
  this gives us that $\rho^{3}(\boldsymbol{\mathcal W})\geq 1$ as claimed. 
  \bk 
  
  It follows that what the 
  ``$\boldsymbol{\rho^3}$-dilogarithmic conjecture'' actually is about  is the fact that 
  $\rho^\sigma(\boldsymbol{\mathcal W})=0$ for $\sigma\geq 4$ in both cases, and that $\rho^3(\boldsymbol{\mathcal W})$ is precisely equal to 1 in  case {\bf (2)}. 
  \sk 
  
    Let us discuss briefly both cases of this conjecture: 
  \begin{itemize}
  \item The condition 
${\rm IntrDim}(\boldsymbol{\mathcal W})\geq 3$ of {\bf (1)} is assumed to disregard the case of planar webs which seems to be quite particular. 
 For instance, Newman's 6-web $\boldsymbol{\mathcal W}_{\hspace{-0.05cm}\mathcal N_6}$ (which is equivalent to the cluster web 
$\boldsymbol{\mathcal W}_{B_2}$
of type $B_2$, see \S\ref{SS:cluster-webs:B2} below) 
carries a proper dilogarithmic AR but is planar hence  
 $\rho^\bullet(\boldsymbol{\mathcal W}_{B_2})=(4,3,2,1)$ so in particular  
 $\rho^4(\boldsymbol{\mathcal W}_{B_2})=1>0$. If the latter web 
 has maximal rank ({\it i.e.}\,is AMP), one can consider the cluster web $\boldsymbol{\mathcal W}_{G_2}$ of type $G_2$ to get one which is not: it is a planar 8-web (thus $\rho^\bullet(\boldsymbol{\mathcal W}_{G_2})=(6,5,4,3,2,1)$) carrying a proper dilogarithimic AR. However, it has rank 
${\rm rk}(\boldsymbol{\mathcal W}_{G_2})=14$, strictly less than 
$\rho(\boldsymbol{\mathcal W}_{G_2})=21$ ({\it cf.}\,\S\ref{SS:cluster-webs:G2}) hence 
this web is not AMP; 
  \sk 
  \item As for {\bf (2)},  we must admit that we don't really have any idea why this could be true in full generality.  We  have just remarked that it is satisfied on all  the examples of webs considered in this text which satisfy the aforementioned conditions.  Thus, it seems to be satisfied by all the dilogarithmic webs considered above (moreover, all these webs seem to be AMP, at the likely exception of the webs $\boldsymbol{\mathcal W}_{R(m)}$ for $m\geq 3$), 
but it is also satisfied by  many of the $\boldsymbol{\mathcal Y}$-cluster webs that we will consider in the next section.
  \end{itemize}

To conclude our discussion of the above conjecture, 
it seems to us that, if true, it 
 accounts for some existing links between two kinds of `(multi)linear relations' associated to the first integrals defining the considered dilogarithmic web.
 To make the discussion simpler, we specialize as follows 
the notation introduced in the paragraph immediately following the statement of the conjecture:  now $\Omega$ is a Zariski open subset of $V=\mathbf C^m$,  the $u_i's$ are assumed to be rational 
and $\omega$ stands for the generic point of $V$. 
 Then, to the $d$-tuple $(u_i)_{i=1}^d$ of elements of 
 $\mathbf C(V)=\mathbf C(x_1,\ldots,x_m)$ defining $\boldsymbol{\mathcal W}$
  are associated two kinds of relations, which can be said of (multi)linear nature since they all  hold (or not) in some complex vector spaces: 
\begin{itemize}
\item the first are relations in $\mathbf C(V)^\times\wedge
 \mathbf C(V)^\times $
of the following form, for some scalar constants $c_i$:
  \begin{equation}
  \label{Eq:WedgeEquations}
  {}^{}\qquad \qquad 
  \sum_{i=1}^d c_i\, u_i\wedge (1-u_i)=0\,
   ;
   \end{equation}
\item the second are identities in ${\rm Sym}^\sigma(V^*)$, for 
 each $\sigma\geq 1$,  of the following type
  \begin{equation}
  \label{Eq:PolynomialEquations}
\sum_{i=1}^d \kappa_{i}^\sigma \,\Big(\ell_i\Big)^\sigma=0\,  
\end{equation}
 for some scalars $\kappa_i^\sigma \in \mathbf C$ (we recall 
that for any $i=1,\ldots,d$, $\ell_i$ stands for the linear form on $V$ corresponding 
 to the 
the differential of $u_i$ at $\omega$ (via the natural identification $V^*\simeq T^*_{\omega} V$)). 
\end{itemize}

We have seen above in the paragraph just following the statement of the conjecture, that there are links between these two kinds of relations. 
It would be interesting, and the `${\rho}^3$-dilogarithmic conjecture' goes into this direction, to understand better how are identities \eqref{Eq:WedgeEquations} and \eqref{Eq:PolynomialEquations} related for any $d$-tuple  of functions 
$(u_1,\ldots,u_d)$ 
defining a web carrying a proper dilogarithmic AR.
}

%

\subsubsection{Trilogarithmic functional equations.}
\label{SS:TrilogFE}

We discuss now some of the AFEs satisfied by trilogarithmic functions that one can find in the literature.  Few such identities  have been discovered by classical authors: as far as we now, in this category one finds only the Spence-Kummer equation (in many different but equivalent forms) and a family due to Sandham. After discussing these functional equations and the webs associated to them, we turn to an interesting and  more recent identity discovered by Goncharov and to a few others which can be obtained from it.

\paragraph{\hspace{-0.2cm}Spence-Kummer functional equation of the trilogarithm.}
\label{SS:Spence-Kummer}

Spence ({\it cf.}\,p.\,33 in \cite{Spence} and also \S6.7 in \cite{Lewin}) and independently Kummer (see formula (93) p.\,336 in \cite{Kummer3}) have independently established that the 
 trilogarithm $\l{3 } $ satisfies the  following  functional equation
\begin{align} \qquad  2 \,\l {3} (\,x\,) &  \, +  \, 2  \,\l{3} (\,y\,)  \, -
  \,  \l{3} \left(\, \frac{x}{y}\, \right)\,     
+ \, 2  \, \l{3}  \left(\,  \frac{1-x}{1-y}\, \right) \,   + \,  2 \,\l{3} \left(\,
  \frac{x(1-y)}{y(1-x)}\, \right)    
- \, \l{3} (\, xy\,)  \qquad \qquad \qquad \qquad \nonumber \\
 ({\cal S}{\cal K})
  \qquad \qquad \qquad   &  \, +   \,  2 \,
 \l{3} \left(\, -\frac{x(1-y)}{{\,}(1-x)}  \,         \right) 
+ \,  2 \,   \l{3} \left(\,  - \frac{{\,}(1-y)}{y(1-x)}\, \right)   -  \, \l{3}
  \left(\, \frac{x(1-y)^2}{y(1-x)^2}\, \right) 
 \qquad      \nonumber \\
  &=2\, 
  \zeta(3)
+
2\,\zeta(2)
\, {\rm Log} (y)
  -{\rm Log}(y)^2\,{\rm Log}\, \left( \frac{1-y}{1-x}\right)
+\frac{1}{3}\, {\rm Log}(y)^3  \, ,  \nonumber 
\end{align} \vspace{0.15cm}
for every real numbers  $x,y$ such that $0<x<y<1$.
\mk

In accordance with Theorem \ref{T:FE-Zagier-Polylogs}, it follows that 
 Zagier's trilogarithm $\mathcal L_3$ identically satisfies the following identity  (this time for any pair $(x,y)\in \mathbf C^2$  distinct 
 from $(0,0)$ or from $(1,1)$): 
\begin{align*} 
\qquad  2 \,\mathcal L_3 (\,x\,) &  \, +  \, 2  \,\mathcal L_3 (\,y\,)  \, -
  \, \mathcal L_3 \left(\, \frac{x}{y}\, \right)\,     
+ \, 2  \, \mathcal L_3 \left(\,  \frac{1-x}{1-y}\, \right) \,   + \,  2 \,\mathcal L_3 \left(\,
  \frac{x(1-y)}{y(1-x)}\, \right)    
 \qquad \qquad \qquad \qquad \nonumber \\
  \qquad \qquad \qquad  - \, \mathcal L_3 (\, xy\,)   &  \, +   \,  2 \,
\mathcal L_3 \left(\, -\frac{x(1-y)}{{\,}(1-x)}  \,         \right) 
+ \,  2 \,   \mathcal L_3 \left(\,  - \frac{{\,}(1-y)}{y(1-x)}\, \right)   -  \, \mathcal L_3 
  \left(\, \frac{x(1-y)^2}{y(1-x)^2}\, \right) 
=2\, \zeta(3)\, . 
\end{align*}

\subparagraph{\hspace{-0.2cm}Spence-Kummer web.}
\hspace{-0.2cm}
By definition, {\bf Spence-Kummer web $
\boldsymbol{{\boldsymbol{\mathcal W}}_{\hspace{-0.05cm}{\cal S}{\cal K}}}$}  is the 
planar 9-web defined by the rational functions appearing as arguments of 
the trilogarithm 
$\l{3} $ in $ ({\cal S}{\cal K})$, that is: 
$$
\boldsymbol{{\boldsymbol{\mathcal W}}_{\hspace{-0.05cm}{\cal S}{\cal K}}}
=\boldsymbol{\mathcal W}\left(\, 
  x  \, , \, 
  y   \, , \, 
    \frac{x}{y}    \, , \, 
    \frac{1-x}{1-y}        \, , \, 
     \frac{x(1-y)}{y(1-x)}            \, , \, 
     xy \, , \, 
         -\frac{x(1-y)}{{\,}(1-x)}           \, , \, 
                   - \frac{{\,}(1-y)}{y(1-x)}      \, , \, 
                  \frac{x(1-y)^2}{y(1-x)^2}   \hspace{0.1cm} \right)\, .
$$
We denote by $U_i$ the rational first integrals of 
$\boldsymbol{{\boldsymbol{\mathcal W}}_{\hspace{-0.05cm}{\cal S}{\cal K}}}$ given above: 
\begin{equation}
\label{Eq:U_i-SK}
U_1=x\, , \quad  U_2=y   \, , \quad 
  U_3 = \frac{x}{y}    \, , \quad 
 U_4=  \frac{1-x}{1-y}        \, , \hspace{0.1cm} 
     \ldots \hspace{0.1cm} , \,  
    U_8=  - \frac{(1-y)}{y(1-x)}
     \, , \quad 
 U_9=  \frac{x(1-y)^2}{y(1-x)^2} \, .
\end{equation}

By direct computations, one obtains that $\mathfrak B_i^s=\{ 0, \pm 1, \infty \}$ except for $i=3,6,9$, in which case one has $\mathfrak B_i^s=\{ 0, 1, \infty \}$.  In what concerns the ramification loci (common leaves), one has  $\mathfrak B_i=\{0 , 1, \infty\}$  for any $i=1,\ldots,9$.\mk

Spence-Kummer web  is of a certain importance in web geometry: indeed,  it is  the second discovered web (independently in \cite{PirioSelecta} and \cite{Robert},  65 years after Bol's example!)  which is exceptional, that is of maximal rank but non-algebraizable.  A complete and explicit basis with 28 elements of  $\mathcal A({\boldsymbol{\mathcal W}}_{\hspace{-0.1cm}{\cal S}{\cal K}})$  
has been given in \cite[\S3.2]{PirioSelecta}:  
5 elements of this basis are  non-polylogarithmic, the 23 others being  polylogarithmic ({\it i.e.}\,iterated integrals with ramification at $0,1$ or $\infty$).  
 In particular, one gets that 
 $${\rm polrk}^\bullet\Big({\boldsymbol{\mathcal W}}_{\hspace{-0.1cm}{\cal S}{\cal K}}\Big)=(12 ,9,2)\,.$$

Since ${\boldsymbol{\mathcal W}}_{\hspace{-0.1cm}{\cal S}{\cal K}}$ is defined by quite simple rational first integrals, it is not difficult to compute all its web-theoretic invariants, which can be  useful to identify this web up to a (possibly unknown) change of coordinates (as it has been the case to identify it with a cluster web, see \S\ref{SS:ClassicalPolylogarithmicIdentitiesAreOfClusterType} below).  

For instance, one has 
$$ \boldsymbol{\mathcal Hex}\Big({\boldsymbol{\mathcal W}}_{\hspace{-0.1cm}{\cal S}{\cal K}}\Big)=(48,30,9,1)  \, , 
\quad  \boldsymbol{\mathcal F \hspace{-0.1cm}lat} \Big({\boldsymbol{\mathcal W}}_{\hspace{-0.1cm}{\cal S}{\cal K}}\Big)=(48, 48, 12, 11, 3, 0, 1)
 \quad \mbox{ and }\quad  
  {\boldsymbol{\mathcal B}}\Big({\boldsymbol{\mathcal W}}_{\hspace{-0.1cm}{\cal S}{\cal K}}\Big)=3\, . 
$$

Regarding the  (birational) symmetries of ${\boldsymbol{\mathcal W}}_{\hspace{-0.1cm}{\cal S}{\cal K}}$, we refer to  \cite{PirioThese} and 
to  Example 3 in \cite[\S7]{ZagierSpecial} where they are discussed. 

\subparagraph{\hspace{-0.2cm}Some geometric descriptions of Spence-Kummer web.} 
\label{SubPar:GeomSK-Web}
\hspace{-0.2cm}
The fact that $ \boldsymbol{\mathcal Hex}_6({\boldsymbol{\mathcal W}}_{\hspace{-0.1cm}{\cal S}{\cal K}})$ is equal to 1 is interesting: it means that Spence-Kummer web admits a unique hexagonal 6-subweb. It is the subweb, denoted by 
${\boldsymbol{\mathcal W}}_6$, 
 formed by the first integrals in \eqref{Eq:U_i-SK}, except $U_3=x/y$, $U_6=xy$ and $U_9=x(1-y)^2/(y(1-x)^2)$. According to Bol's classification of hexagonal planar webs, ${\boldsymbol{\mathcal W}}_6$ is equivalent to a linear web formed by six pencils of lines. 
 Considering the birational map 
  $$\Phi:(X,Y)\rightarrow (x,y)=\left(\, \frac{1+X}{X}\, , \, \frac{1+Y}{Y}\, \right)\, ,$$
  one verifies that $\Phi^*({\boldsymbol{\mathcal W}}_6)$ 
  is the web formed by the pencils whose vertices $p_i$ form a
  degenerate 
   configuration in $\mathbf P^2$,  
  unique up to  ${\rm PGL}_3(\mathbf C)$, pictured just below. 
  We call it {\bf Spence-Kummer configuration} and denote it by $C_{\hspace{-0.05cm}{\cal S}{\cal K}}$.\vspace{0.2cm}
\begin{figure}[h]
\begin{center}
\resizebox{1.4in}{1.4in}{
 \includegraphics{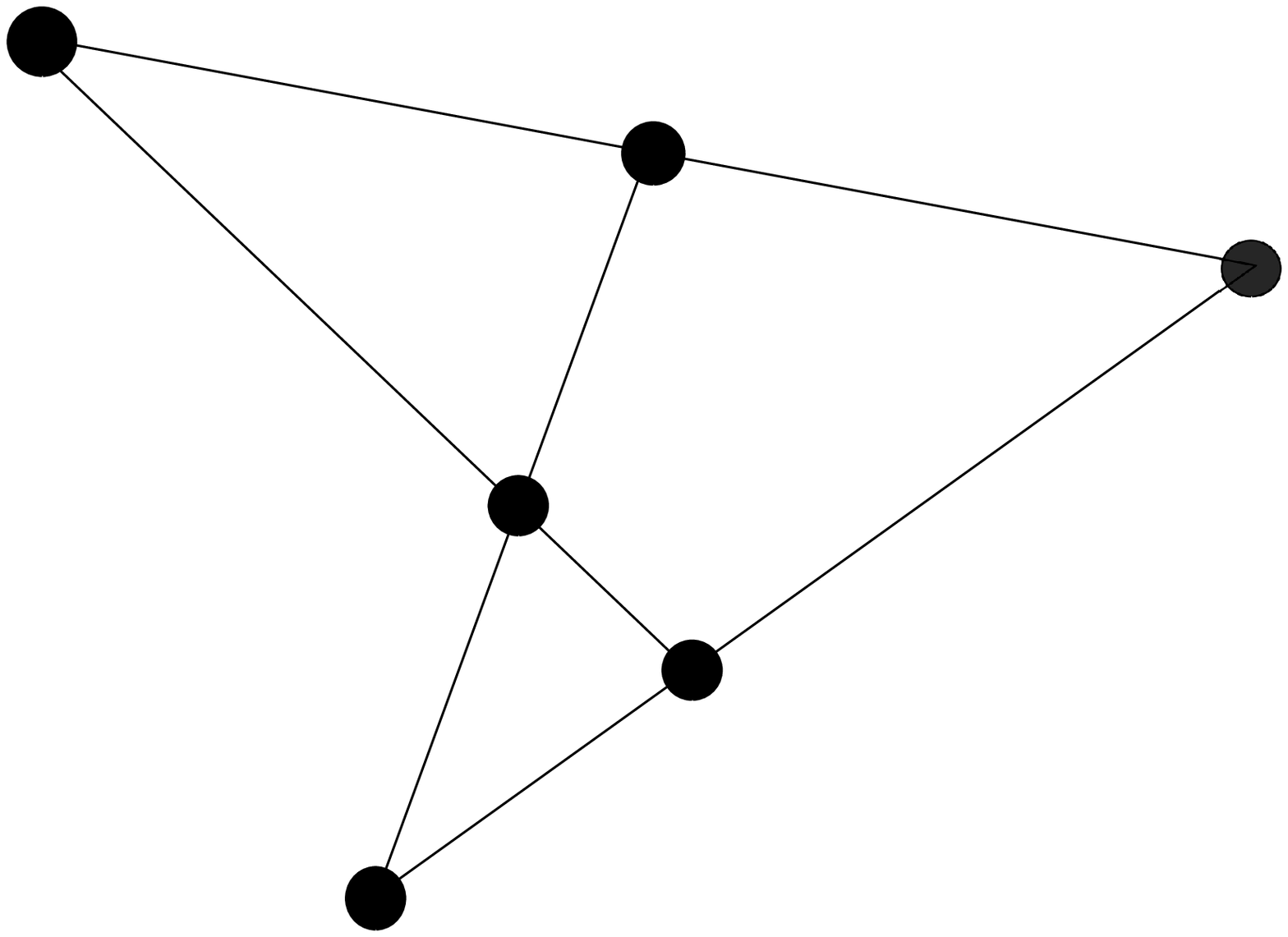}}
 \vspace{-0.1cm}
\caption{Spence-Kummer configuration $C_{\hspace{-0.05cm}{\cal S}{\cal K}}\in {\rm Conf}_6(\mathbf P^2)$.}  \label{Fig:Spence-Kummer configuration}
\end{center}
\end{figure}
  
Since a linearizable planar $d$-web admits at most one linearization (up to composition with a projective transform, see Proposition \ref{P:Unique-Lin}),  $ C_{\hspace{-0.1cm}{\cal S}{\cal K}}=[p_1,\ldots,p_6]\in {\rm Conf}_6(\mathbf P^2)$ is canonically attached to ${\boldsymbol{\mathcal W}}_{\hspace{-0.0cm}{\cal S}{\cal K}}$ and 
 $$\Phi^*({\boldsymbol{\mathcal W}}_{\hspace{-0.1cm}{\cal S}{\cal K}})=
 {\boldsymbol{\mathcal W}}\left(  \, X 
  \, , \,  Y 
   \, , \,  \frac{X(1+Y)}{Y(1+X)}
  \, , \,  \frac{X}{Y}
  \, , \,   \frac{1+X}{1+Y}
  \, , \,  \frac{XY}{(1+X)(1+Y)}
  \, , \,  \frac{X}{1+Y}
  \, , \,  \frac{1+X}{Y}
    \, , \,  \frac{X(1+X)}{Y(1+Y)}\, 
  \right)
$$ is a canonical model of Spence-Kummer web (unique, up to projective transformations). This web can be described as the 9-web on the projective plane, formed by the pencils whose vertices are the points of 
$ C_{\hspace{-0.1cm}{\cal S}{\cal K}}$ together with the three pencils of conics associated to the three subconfigurations of $ C_{\hspace{-0.0cm}{\cal S}{\cal K}}$ constituted by four points in general position. From this geometric description of the `canonical model' of Spence-Kummer web and the knowledge of many of its invariants, one deduces a quite efficient way to decide whether  a given 9-web is equivalent to it or not.\mk

The previous description of ${\boldsymbol{\mathcal W}}_{\hspace{-0.1cm}{\cal S}{\cal K}}$ by means of a projective configuration makes us wonder whether  this web is equivalent to 
  a web associated to a stratum of projective configuration, through the construction presented in \S\ref{Subpar:Strata} above. An affirmative answer to this can be found in \cite[\S1.5]{Goncharov1995}:  let $\boldsymbol{S}({\hspace{-0.0cm}{\cal S}{\cal K}})$ be the set  formed by degenerate configurations of seven points on $\mathbf P^2$  of the following type: 
\begin{figure}[h]
\begin{center}
\resizebox{1.5in}{1.1in}{
 \includegraphics{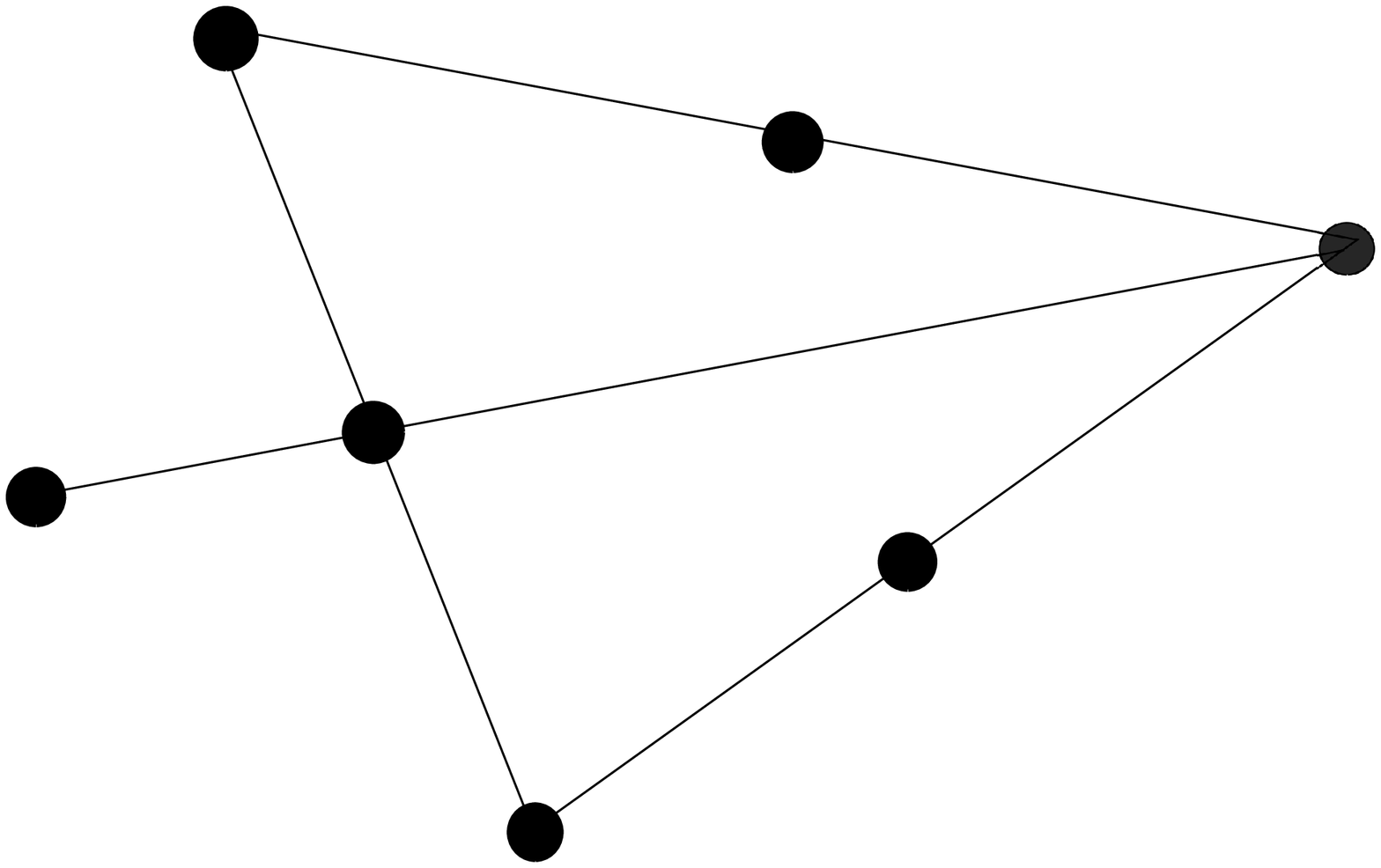}}
 \vspace{-0.0cm}
\caption{An element of $\boldsymbol{S}({\hspace{-0.0cm}{\cal S}{\cal K}})$.} 
\label{Fig:CSKG}
\end{center}
\end{figure}

  It is a stratum  of ${\rm Conf}_7(\mathbf P^2)$ of dimension 2. 
One verifies that, up to equivalence of webs, one has 
 $$ {\boldsymbol{\mathcal W}}_{\boldsymbol{S}({\hspace{-0.0cm}{\cal S}{\cal K})}}= \Big( {\boldsymbol{\mathcal W}}_{{\rm Conf}_7(\mathbf P^2)} \Big)
 \Big\lvert_{\boldsymbol{S}({\hspace{-0.0cm}{\cal S}{\cal K})}}
 = {\boldsymbol{\mathcal W}}_{\hspace{-0.1cm}{\cal S}{\cal K}}\, .$$

The two types of projective configurations pictured 
in the two figures 
above are related to Spence-Kummer web but we are not aware of  any direct and natural geometric link between them, which is a bit surprizing. However, there is a way to obtain ${\boldsymbol{\mathcal W}}_{\hspace{-0.1cm}{\cal S}{\cal K}}$ by means of another sub-stratum in 
${\rm Conf}_7(\mathbf P^2)$ constructed from the Spence-Kummer configuration $ C_{\hspace{-0.1cm}{\cal S}{\cal K}}$ as follows: let $\varphi_i: {\rm Conf}_7(\mathbf P^2)\rightarrow {\rm Conf}_6(\mathbf P^2)$ be the forgetful map of the $i$-th point for some $i$, say $i=7$. Then 
$\boldsymbol{\Sigma}({\hspace{-0.0cm}{\cal S}{\cal K}})=\varphi_7^{-1}( C_{\hspace{-0.1cm}{\cal S}{\cal K}})$ is the 2-dimensional stratum  whose generic element is obtained from $C_{\hspace{-0.1cm}{\cal S}{\cal K}}$ by adjoining to it a seventh point 
not lying on any line between two of its points (see Figure \ref{Fig:AgenCSK} below).
\begin{figure}[h]
\begin{center}
\resizebox{1.7in}{1.7in}{
 \includegraphics{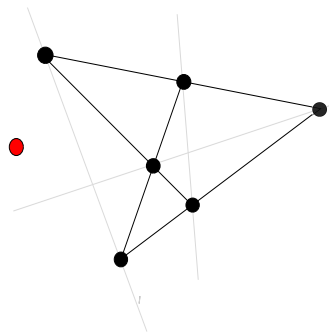}}
 \vspace{-0.3cm}
\caption{A generic element of $\boldsymbol{\Sigma}({\hspace{-0.0cm}{\cal S}{\cal K}}) 
$
 (the seventh point added to 
$C_{\hspace{-0.05cm}{\cal S}{\cal K}}$ is pictured in red).} 
 \label{Fig:AgenCSK}
\end{center}
\vspace{-0.4cm}
\end{figure}

Then, despite the fact that $\boldsymbol{\Sigma}({\hspace{-0.0cm}{\cal S}{\cal K}})$ is a  stratum of  ${\rm Conf}_7(\mathbf P^2)$ distinct from  $\boldsymbol{S}({\hspace{-0.0cm}{\cal S}{\cal K}})$ (their 
 combinatorial types are clearly not the same), 
we also have (up to equivalence): 
  $$ {\boldsymbol{\mathcal W}}_{\boldsymbol{\Sigma}({\hspace{-0.0cm}{\cal S}{\cal K})}}=
   \Big( {\boldsymbol{\mathcal W}}_{{\rm Conf}_7(\mathbf P^2)} \Big)
 \, \big\lvert_{\,\boldsymbol{\Sigma}({\hspace{-0.0cm}{\cal S}{\cal K})}}
 = {\boldsymbol{\mathcal W}}_{\hspace{-0.1cm}{\cal S}{\cal K}}\, .$$
 \sk 
 
  We have indicated above a couple of geometric constructions of Spence-Kummer web, by means of projective configurations of points.    In the third section further, 
 we will give  
   another way to obtain ${\boldsymbol{\mathcal W}}_{\hspace{-0.1cm}{\cal S}{\cal K}}$, but of a more algebraic nature. Indeed,  we will prove that it is equivalent to a cluster web constructed from the $\boldsymbol{\mathcal X}$-cluster web of type $A_3$ (see Theorem \ref{T:classical-cluster-webs} below). 
\sk

\subparagraph{}\hspace{-0.6cm}
To finish this paragraph about the Spence-Kummer equation of the trilogarithm, we would like to  mention the very 
interesting paper \cite{KerrLewisLopatto} (especially the sixth section therein) where,  using the Abel-Jacobi maps for higher Chow groups they have previously constructed, the authors are able to deduce (an holomorphic version of) $({\cal S}{\cal K})$ via a `{\it reprocity law}' arising from subvarieties of projective space. \sk 

What is especially interesting within this approach is that the  framework used by these authors to get the functional equation $(\mathcal S\mathcal K)$ is conceptually the same as the much more classical one, that allows to obtain the abelian relations of algebraic webs. Since something similar also holds true for Bol's web (see  
\S\ref{Subpar:(AB)Bol} above), this puts on the same level, the standard algebraic webs on the one hand, and the exceptional polylogarithmic webs on the other.  It is a somewhat `revolutionary' change of perspective in web geometry  since it has long been thought that these webs were of fundamentally different natures. \sk

However, in contrast with what has been obtained for Abel's equation of the dilogarithm (see \S\ref{Subpar:(AB)Bol} above and also  the fifth section of {\it loc.\,cit.}), 
the approach followed by Kerr, Lewis and Lopatto only allows them to construct (the holomorphic version of) the Spence-Kummer functional equation, but not the nine foliations 
forming Spence-Kummer web ({\it cf.}\,Remark 6.3.(b) in \cite{KerrLewisLopatto}).  
Moreover, their approach shows several features which are geometrically distinct from the ones mentioned above: 
\begin{itemize}
\item[\bf (i).] first they do not work on a space of projective configurations
of seven points in $\mathbf P^2$ 
 but rather over it, on  a certain open domain $\mathcal U$ in 
 the complex grassmannian of 2-planes in $\mathbf P^6$; 
 \item[\bf (ii).] second, and more importantly, the main result they use (namely what they call {\it `Reprocity law B'}) involves not only (natural lifts to the corresponding grassmannians of the) projection maps from $\mathbf P^2$ to $\mathbf P^1$ with center the points of the considered configurations, but also 
(the corresponding natural grassmannians lifts) of  the forgetful maps consisting in forgetting a point of the configurations. 
\end{itemize}

As explained above, Spence-Kummer web (that is, the nine rational fibrations which compose it) is  obtained as the restriction along a stratum in ${\rm Conf}_7(\mathbf P^2)$, of the ${ 7 \choose 1}{ 6 \choose 4 }=105$ projection maps ${\rm Conf}_7(\mathbf P^2)\rightarrow \mathbf P^1$ hence no forgetful maps are involved.  Note that {\it `Reprocity law A'} in \cite{KerrLewisLopatto} only involves (lifts of) projection maps (and no forgetful map) hence one can ask whether  it is possible to get the Spence-Kummer equation, and even better, the associated  web ${\boldsymbol{\mathcal W}}_{\hspace{-0.1cm}{\cal S}{\cal K}}$ from this simpler reciprocity law or not. 
\sk

In any case, and more generally, understanding better the content of  \cite{KerrLewisLopatto} could be very interesting: it could provide conceptual explanations for the (or at least for some of the ) functional equations satisfied by polylogarithms which are, in some way, similar to the ones responsible for the abelian relations of the algebraic webs (Abel's type theorems, that is vanishing results  for traces of abelian differentials).

\paragraph{\hspace{-0.2cm} Sandham's trilogarithmic functional equations.}
 In \cite{Sandham} (see also \S6.9 in \cite{Lewin}), Sandham establishes a series of trilogarithmic 
 functional equations, one for each positive integer, which is formally very similar to Rogers' series which has been discussed  \S\ref{Par:RogersMultivariableAFEs} above.\mk 
 
Sandham's functional equations  are expressed in terms of a non-standard trilogarithmic function\footnote{Namely the function $M: x\mapsto 
\mu(x)-\frac{1}{3}\log x\, \log^2(1-x)$ with $\mu(x)=\int_0^x\log^2(1-u) du/u$,  for $x\in \mathbf R$.} and are presented in symmetric form. This has the consequence that some terms in each functional equation are functionally dependent, which is not the most convenient approach if one is interested in the web associated to it.  For our purpose, it is more interesting: first, not to use Sandham's original identity  but rather the one satisfied by Zagier's trilogarithm $\mathcal L_3$ corresponding to it; second, to clean the equation 
by grouping together the terms which are functionally dependent.  
\mk 

Let $N>0$ be a fixed integer.  One considers $N+1$ indeterminates $\lambda_1,\ldots,\lambda_N$ and $y$ 
and one sets $\lambda=\lambda=(\lambda_1,\ldots,\lambda_N)$. 
 The $N$ roots $x_s=x_s(y,\lambda)$,  with $s=1,\ldots,N$,  of the polynomial equation $\prod_{i=1}^N(1-\lambda_i x)=y$ in the variable $x$ are elements of 
 the field of algebraic functions of $(y,\lambda)$, denoted by 
$\boldsymbol{F}$.  One considers the following element of 
$\mathbf Z[\boldsymbol{F}^{**}]$: 
\begin{align*}
\boldsymbol{\mathcal S}_N= -\,  \Big[y\Big]  
    + & \, 
   N\sum_{i,s=1}^N \Big[ 1- \lambda_i x_s  \Big]
+ N \hspace{-0.25cm}\sum_{1\leq i < j \leq N}   \left[ \frac{\lambda_i }{\lambda_j} \right]    
     + & \, 
     \hspace{-0.25cm} \sum_{1\leq i < j \leq N}
     \sum_{s=1}^N  
    \Bigg(\, 
     \left[ \frac{\lambda_j(\lambda_i x_s-1) }{\lambda_i(\lambda_j x_s-1)} \right] -   \left[ \frac{\lambda_i x_s-1 }{\lambda_j x_s-1} \right]
  \,   \Bigg) 
  \, . 
\end{align*}

Then the $N$-th Sandham functional identity is equivalent to 
 the relation 
\begin{equation}
\label{Eq:Sandham}
\mathcal L_3(\boldsymbol{\mathcal S}_N)\equiv \frac{1}{2}N^2(N-1)\, \zeta(3)\, .
\end{equation}
 By definition, the {\bf $\boldsymbol{N}$-th Sandham web $\boldsymbol{\mathcal W}_{ {}^{}\hspace{-0.1cm}
\boldsymbol{\mathcal S}_N}$} is the web associated to this functional identity: 
\begin{equation}
\label{Eq:WS_N}
\boldsymbol{\mathcal W}_{ {}^{}\hspace{-0.1cm}
\boldsymbol{\mathcal S}_N}=
\boldsymbol{\mathcal W}\left( 
 \, y 
 \, , \,  \lambda_i x_s
 \, , \,  \frac{\lambda_i}{\lambda_j}
 \, , \,  \frac{\lambda_j(\lambda_ix_s-1)}{\lambda_i(\lambda_j x_s-1)}
 \, , \,   \frac{\lambda_ix_s-1}{\lambda_j x_s-1} 
 \hspace{0.2cm}
  \bigg\lvert 
   \begin{tabular}{l}
  $1\leq i< j  \leq N  $\\
  $s=1,\ldots,N$ 
 \end{tabular}
 \, 
\right)\, . 
\end{equation}
It is a web  defined by $N(N+1)(2N-1)/2+1$ algebraic first integrals in $(y,\lambda)$, that is in $N+1$ variables. When $N>2$, it admits 
$\boldsymbol{\mathcal W}\left(  \, y 
 \, , \,  {\lambda_i}/{\lambda_j} \, \lvert \, 1\leq i< j  \leq N   \, \right)$ as a subweb, which ensures that its intrinsic dimension 
is $N+1$ (this constrasts with the case when $N=2$, see just below).  In this case, we believe that  the algebraic first integrals in \eqref{Eq:WS_N} define pairwise 
 distinct foliations hence  that 
$\boldsymbol{\mathcal W}_{{}^{}\hspace{-0.1cm}
\boldsymbol{\mathcal S}_N}$ is a genuine $(N(N+1)(2N-1)/2+1)$-web in $N+1$ variables.\footnote{We have verified this only for $N=3$.} \mk 

The $N=2$ case is particularly interesting. Indeed, in this case,  exactly two of the 10 foliations defined by the functions appearing in \eqref{Eq:WS_N} coincide thus  $\boldsymbol{\mathcal W}_{ {}^{}\hspace{-0.1cm}
\boldsymbol{\mathcal S}_2}$ is a 9-web; second, the intrinsic dimension of the latter is 2 hence it is the pull-back of a planar $9$-web that we denote by 
$\boldsymbol{\mathcal W}'_{ {}^{}\hspace{-0.1cm}
\boldsymbol{\mathcal S}_2}$.  In this case, Sandam  has established 
(see \cite[\S8]{Sandham} or \cite[\S6.8]{Lewin}) that $
\mathcal L_3(\boldsymbol{\mathcal S}_2)= 2\, \zeta(3)$ is equivalent to the following identity for  the classical trilogarithm $\l {3} $: 
\begin{align*}
\l {3} \left( -X^2  \right)
 + & \, 
\l {3} \left(   -Y^2\right)
+
\l {3} \left(   -Z^2 \right)
-  2\, 
\l {3} \left( XY  \right)
-2\, 
\l {3} \left(  XZ \right) \\
-& 2\, 
\l {3} \left(  YZ \right)
-2\, 
\l {3} \left(  -X/Y  \right)
-2 
\l {3} \left(  -Y/Z  \right)
-2 
\l {3} \left(  -X/Z  \right) \\
 =& \,  \frac{1}{3}\pi^2\log(-Z/Y)+
\frac{1}{3}\log^3(-Z/Y)+\log^2(-Z/Y)\log(-X/Y)-2\,\zeta(3)\, , 
\end{align*}
which is identically verified under the assumption 
that the variables $X,Y$ and $Z$ satisfy the algebraic relation
$X+Y+Z=XYZ$.   Up to a birational change of coordinates,  it can be seen that the previous identity  is equivalent to the Spence-Kummer equation $({\cal S}{\cal K})$ ({\it cf.}\,\cite[\S6.8.2]{Lewin}) from which it comes that 
 the second Sandham's web $\boldsymbol{\mathcal W}_{ {}^{}\hspace{-0.1cm}
\boldsymbol{\mathcal S}_2}$ actually is a pull-back of Spence-Kummer web ${\boldsymbol{\mathcal W}}_{\!{\cal S}{\cal K}}$. The latter being AMP, one can wonder whether  or not the same holds true for all Sandham's webs 
$\boldsymbol{\mathcal W}_{ {}^{}\hspace{-0.1cm}
\boldsymbol{\mathcal S}_N}$, $N\geq 2$. \vspace{-0.2cm}
\begin{center}
$\star$
\end{center}

In the same way as Rogers' equations $(R(m))$ are in fact special cases of a functional identity that applies to each rational fraction  $\varphi \in \mathbf C(z)$ (see \S\ref{Par:RogersMultivariableAFEs} ans \S\ref{Par:Dilog_varphi}), Sandham's trilogarithmic functional equations are special cases of a family of such identities, one for  each sufficiently generic rational function $\varphi$. For more details, see  the chapters by Ray, Weschung and Wojtkowiak in \cite{LewinStruct} or the beginning of the third section in \cite{Gangl}.

%

\paragraph{\bf \hspace{-0.3cm} Goncharov's 22-terms equation of the trilogarithm.}
\label{Par:Goncharov-22}
\hspace{-0.2cm} In \cite{Goncharov} (see also \S3.1 in \cite{Gangl}), Goncharov shows that the trilogarithm satisfies an interesting functional equation in three variables that he interprets geometrically in terms of (degenerate) configurations of 7 points on $\mathbf P^2$. \sk 

\vspace{-0.3cm}
More precisely, for generic complex numbers $a,b,c$, Goncharov introduces the following element 
\begin{align*}
g(a,b,c)=-\big[  1\big]+ \big[ c \big]+ 
\big[  ca-a+1\big]  &\, + 
\left[    \frac{  ca-a+1}{ca} \right]  \, +
\left[    \frac{bc-c+1}{  (ca-a+1)b} \right]\\ & \, 
- 
\left[    \frac{  ca-a+1}{c} \right] 
-
\left[    \frac{bc-c+1}{  (ca-a+1)bc} \right]
+
\left[    \frac{(bc-c+1)a}{  ca-a+1} \right]
\end{align*}
of $ \mathbf Z[ \mathbf Q(a,b,c)  ]$ 
 and then he proves that 
$$
   \boldsymbol{\mathcal G}_{22}
   =g(a,b,c)+g(c,a,b)+g(b,c,a)+\big[ -abc \big]$$
is a functional equation for the trilogarithm, {\it i.e.} one has  $\mathcal L_3( \boldsymbol{\mathcal G}_{22})=0$ identically in $a,b,c$.
 \mk

Then one can consider the associated Goncharov's web: it is the 22 web in 3 variables (satisfying (PGw)), denoted by 
$ {\boldsymbol{\mathcal W}}_{\hspace{-0.07cm} \boldsymbol{\mathcal G}_{22}}$, defined by  the  22 rational primitive first integrals appearing as arguments of $\mathcal L_3$ in Goncharov's equation. Explicitly,  
 one has
\begin{align}
\label{Def:G22}
 {\boldsymbol{\mathcal W}}_{ \boldsymbol{\mathcal G}_{22}}= \boldsymbol{\mathcal W}\bigg( \, a\, , \, & b\, ,\,  c\, ,\,  -abc\, ,\,  ab-b+1\, , \, bc-c+1\, ,\,  ac-a+1\, , \, \frac{ab-b+1}{a}\, ,\,  
\frac{bc-c+1}{b}\, ,  \nonumber \\ 
\nonumber
& \frac{ac-a+1}{c}\, ,\,  -\frac{(bc-c+1)a}{ac-a+1}\, , \, 
-\frac{(ac-a+1)b}{ab-b+1}\, ,\,  -\frac{(ab-b+1)c}{bc-c+1}\, , \, \frac{ab-b+1}{ab}\, , \\ 
&  \frac{ab-b+1}{a(bc-c+1)}\, , \,  \frac{bc-c+1}{bc}\, ,\,  \frac{bc-c+1}{b(ac-a+1)}\, , \, \frac{ac-a+1}{ac}\, ,\, 
 \frac{ac-a+1}{(ab-b+1)c}\, ,  \\ &  \hspace{4.25cm}
\frac{ab-b+1}{(bc-c+1)ab}\, ,\,  \frac{bc-c+1}{(ac-a+1)bc)}\, , \, 
\frac{ac-a+1}{(ab-b+1)ac}\,  \bigg)\, . \nonumber
\end{align}





By direct explicit computations, it can verified that this web has polylogarithmic ramification ({\it i.e.}\,one has $\mathfrak B_i=\{ 0 ,  1 ,  \infty\}$ for every $i=1,\ldots,22$), and that 
\begin{equation}
\label{Eq:rhoG22}
\rho^\bullet\Big( {\boldsymbol{\mathcal W}}_{  \boldsymbol{\mathcal G}_{22} }\Big)=\Big(19,16,12,7,2\Big)\qquad \mbox{hence } \qquad \rho\Big( {\boldsymbol{\mathcal W}}_{  \boldsymbol{\mathcal G}_{22} }\Big)=56\, .
\end{equation}

On the other hand, using the method described in \S\ref{Par:II-AR} 
 and considering iterated integrals abelian relations of  $ {\boldsymbol{\mathcal W}}_{ \boldsymbol{\mathcal G}_{22}}$ with ramification at $0$, $1$ and $\infty$ for each of the first integrals appearing in  \eqref{Def:G22}, one obtains after some computations: 
$$
{\rm polrk}^\bullet\Big( {\boldsymbol{\mathcal W}}_{ \boldsymbol{\mathcal G}_{22} }\Big)=\big(\, 34\, , \, 20\,  , \, 2\, \big) \qquad  
\mbox{thus}\qquad 
{\rm polrk}\Big( {\boldsymbol{\mathcal W}}_{  \boldsymbol{\mathcal G}_{22} }\Big)=56\, . 
$$
Comparing with \eqref{Eq:rhoG22}, it follows that  Goncharov's web $ {\boldsymbol{\mathcal W}}_{ \boldsymbol{\mathcal G}_{22} }$ has AMP rank. Moreover,   all its ARs are polylogarithmic which makes of it a web quite similar to Bol's web, but in weight 3. 
\begin{center}
\vspace{-0.15cm}
$\star$
\end{center}




In \cite{Goncharov1995Adv}, Goncharov gave a geometric description of   ${\boldsymbol{\mathcal W}}_{\hspace{-0.08cm} \boldsymbol{\mathcal G}_{22}}$  in terms of projective configurations of seven points  on the projective plane. 
  Let $
\boldsymbol{\Sigma}(\boldsymbol{{\mathcal G}}_{22})$ be the 3-dimensional subvariety of the space of configurations ${\rm Conf}_7(\mathbf P^2)$ formed by the degenerate configurations as in the picture below: 
\begin{figure}[h]
\begin{center}
\resizebox{1.5in}{1.1in}{
 \includegraphics{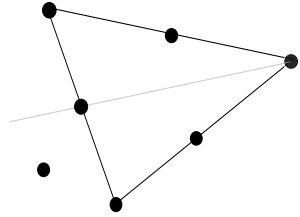}}
 \vspace{-0.0cm}
\caption{A generic  element of Goncharov's stratum $ \boldsymbol{\Sigma}(\boldsymbol{{\mathcal G}}_{22}) \subset {\rm Conf}_7(\mathbf P^2)
$.} \vspace{-0.3cm}
\label{Fig:CSKG}
\end{center}
\end{figure}

Then one verifies that  Goncharov's 22-web 
coincides with 
 the web obtained by taking the restriction of 
 ${\boldsymbol{\mathcal W}}_{\hspace{-0.08cm}
 {\rm Conf}_7(\mathbf P^2)}$ on 
 $\boldsymbol{\Sigma}(\boldsymbol{{\mathcal G}}_{22})$: one has 
$$
{\boldsymbol{\mathcal W}}_{ \hspace{-0.08cm}\boldsymbol{{\mathcal G}}_{22}}=\bigg(\, {\boldsymbol{\mathcal W}}_{ \hspace{-0.08cm}{\rm Conf}_7(\mathbf P^2)}
\bigg)\, 
\Big\lvert_{\, \boldsymbol{\Sigma(\boldsymbol{{\mathcal G}}_{22})}}\, .
$$

For other descriptions of $\boldsymbol{\Sigma}(\boldsymbol{{\mathcal G}}_{22})$ as well as a discussion of its symmetries,  see \S3.1 in \cite{Gangl}. 
\mk

It has been  conjectured\footnote{We believe that this conjecture has to be attributed to Goncharov but we have not be able to locate a precise statement of it in the existing literature. Any information on this subject would be of interest to the author.} that Goncharov's identity $(\boldsymbol{{\mathcal G}}_{22})$ could be the `{\it fundamental functional equation}' for the trilogarithm, in the sense that any AFE $\sum_{i=1}^N c_i \mathcal L_3(u_i)={\rm cst.}$ with rational coefficients $c_i$ and 
rational functions $u_i$  as arguments (of a finite but arbitrary number of variables), can be obtained formally as a linear combination of a finite number of copies of Goncharov's equation. 
As far as we know, this conjecture is still wide open.


\paragraph{Gangl's 21-terms equation of the trilogarithm.}
\label{Parag=Gan-21}
In \cite[\S3.4]{Gangl}, Gangl explains that a certain sum of four copies of 
$\boldsymbol{{\mathcal G}}_{22}$ gives rise to a functional identity with 21 terms for $\mathcal L_3$ which corresponds to the following element\footnote{Actually, there is a misprint in \cite[\S3.4]{Gangl}, the formula  
$\Gamma(x_1,x_2,z_1)+\Gamma(x_2,x_1,z_1)$  given there  is not 
annihilated by $\mathcal L_3$ hence is not correct. The expression given here for $\boldsymbol{\mathcal G an} _{21}$ has been communicated to us by Gangl (private exchange).}  of $\mathbf Z[\mathbf Q(x_1,x_2,z_1)]$ (with $z_2=x_1x_2z_1$): 
\begin{align*}
\boldsymbol{\mathcal G an} _{21}=   & 
[x_1] +  2\,[x_2] - 2\,[x_1x_2]
- 2\,[x_2z_1] + 2\,[x_1x_2z_1] + 2\,\left[z_1^{-1}\right] - 2\,[x_1z_1] \\ 
& + 
2 \,\left[-\frac{x_1( x_2-1)}{x_1-1}\right]  + 
 2\left[\frac{x_2(x_1-1)}{x_2-1}\right]   - 
 2\,\left[ -\frac{x_1(x_2-1)z_1}{x_1-1}\right]  
 - 2\, \left[\frac{x_2(z_1-1)}{z_2-1}\right] 
  \\ 
 & + 2\,\left[\frac{z_1-1}{z_2-1}\right] - 
 2\,\left[\frac{x_1( z_1-1)}{z_2-1}\right]   -  
 2\,\left[-\frac{x_1(x_2-1)(z_1-1)}{(x_1-1)(z_2-1)}\right] 
 -  2\left[-\frac{(x_1-1)x_2z_1}{x_2-1}\right] 
  \\
 & + 2\, \left[\frac{x_1x_2(z_1-1)}{z_2-1}\right] - 
 2\, \left[- \frac{(x_1-1)x_2(z_1-1)}{(x_2-1)(z_2-1)}\right] + 
 \left[-\frac{(x_1-1)(z_1-1)}{(x_2-1)z_1(z_2-1)}\right] \\ &  + 
 \hspace{0.3cm}
 \left[-\frac{(x_2-1)(z_1-1)}{(x_1-1)z_1(z_2-1)}\right] + 
\left[-\frac{x_1^2(x_2-1)(z_1-1)z_1}{(x_1-1)(z_2-1)}\right] + 
\left[-\frac{(x_1-1)x_2^2(z_1-1)z_1}{(x_2-1)(z_2-1)}\right] \,  .
\end{align*}

Then one denotes by $ {\boldsymbol{\mathcal W}}_{ \hspace{-0.1cm}\boldsymbol{\mathcal G an} _{21} }$ the 21-web in three variables defined by the rational functions appearing in the definition of $\boldsymbol{\mathcal G an} _{21}$. 
By direct computations, one verifies that 
$$
\rho^\bullet\Big( {\boldsymbol{\mathcal W}}_{ \hspace{-0.1cm}\boldsymbol{\mathcal G an} _{21} } \Big)=(18,15,11,6,1) \qquad \mbox{hence}
\qquad  \rho\Big( {\boldsymbol{\mathcal W}}_{ \hspace{-0.1cm}\boldsymbol{\mathcal G an} _{21} } \Big)=51\, . 
$$
On the other hand, looking at harmonic polylogarithmic ARs for this web (that is, whose components are harmonic polylogarithms, that is ramified at $\pm1$, $0$ and $\infty$), one gets
$$
{\rm polrk}^\bullet\Big( {\boldsymbol{\mathcal W}}_{ \hspace{-0.1cm}\boldsymbol{\mathcal G an} _{21} } \Big)=(30,18,2) \qquad \mbox{hence}
\qquad {\rm polrk}\Big( {\boldsymbol{\mathcal W}}_{ \hspace{-0.1cm}\boldsymbol{\mathcal G an} _{21} } \Big)=50\, . 
$$
Since  
${\rm polrk}\big( {\boldsymbol{\mathcal W}}_{ \hspace{-0.1cm}\boldsymbol{\mathcal G an} _{21} } \big)=50<51=\rho\big( {\boldsymbol{\mathcal W}}_{ \hspace{-0.1cm}\boldsymbol{\mathcal G an} _{21} } \big)$, one might think that Gangl's 21-web 
is not AMP (which could be an indication that  $\boldsymbol{\mathcal G an} _{21} $ might be less fundamental 
than Goncharov's functional equation $\boldsymbol{\mathcal G on} _{22}$), but actually it is.  
Indeed, setting 
$A(u)={\rm Arctan}(\sqrt{u})$, then the following identity holds true and constitutes an AR for Gangl's 21-web: 
\begin{align*}
A\left( \frac{(x_1-1)(z_1-1)}{(x_2 - 1)z_1(z_2 - 1)}\right) 
- & A\left( \frac{(x_2 - 1)(z_1-1)}{(x_1-1)z_1(z_2 - 1)}\right)  \\
+ &A\left( \frac{x_1^2(x_2 - 1)(z_1-1)z_1}{(x_1-1)(z_2 - 1)}\right) 
-A\left( \frac{(x_1-1)x_2^2(z_1-1)z_1}{(x_2 - 1)(z_2 - 1)}\right) =0\, .
\end{align*}
(Note that this AR is not polylogarithmic since its 
components are not iterated integrals\footnote{However, see  Conjecture \ref{Conj:Nature-Fi}  above regarding this.}).  
%
It follows that ${\rm rk}\big( {\boldsymbol{\mathcal W}}_{ \hspace{-0.07cm}\boldsymbol{\mathcal G an} _{21} } \big)=\rho\big( {\boldsymbol{\mathcal W}}_{ \hspace{-0.06cm}\boldsymbol{\mathcal G an} _{21} } \big)=51$ hence 
${\boldsymbol{\mathcal W}}_{ \hspace{-0.06cm}\boldsymbol{\mathcal G an} _{21} }$ is  AMP as well. 
\mk 

To finish our discussion of $\boldsymbol{\mathcal G an} _{21}$, 
we remark that restricting Gangl's 21-web on the hyperplane  of $\mathbf C^3$
  cut out by $x_1-x_2=0$, one obtains the following 9-web in two variables
$${\boldsymbol{\mathcal W}}\bigg(\, x_1\, ,\,  z_1\,, \, x_1z_1\, , \, x_1^2z_1\, ,\,  \frac{z_1-1}{x_1^2z_1 - 1}\, , \,x_1 \frac{ z_1-1}{x_1^2z_1 - 1}\, , \, x_1^2 \frac{z_1-1}{x_1^2z_1 - 1}\, ,\,  \frac{z_1-1}{z_1(x_1^2z_1 - 1)}\, ,\, x_1^2z_1\frac{z_1-1}{x_1^2z_1 - 1}\, \bigg)
$$
which, as is easily verified,  turns out to be a model of Spence-Kummer's web
$\boldsymbol{\mathcal W}_{ \hspace{-0.07cm} \boldsymbol{{\mathcal S}{\mathcal K}}}$.
\subsubsection{Kummer's functional equations for the tetra- and penta-logarithm.}
\label{SS:Kummer-K(4)-K(5)}
In a series of papers published in 1840,  Kummer  has  
established that the polylogarithms of weight 4 and 5 satisfy AFEs which look similar, in some respect, to the five-terms  and the Spence-Kummer equations of the dilogarithm and the trilogarithm respectively. 
 \mk 
 
Kummer has worked with the polylogarithmic functions defined for any 
$x\in \mathbf R$ by 
 $$
\boldsymbol{\Lambda}_n(x)=\int_{\vec{0}}^x 
 \frac{\log^{n-1}\lvert u\lvert}{1+u}du
$$
(with $n\in \mathbf N^*$). 
This function is related to the  classical 
 $n$-th polylogarithm  $\l {n} $ through the following formula (see (3.12) p.\,27 in \cite{LewinStruct}): 
\begin{equation}
\label{Eq:Lin-Lambdan}
\l {n} (x) =\frac{(-1)^{n}}{(n-1)!}  
\boldsymbol{\Lambda}_n(-x) - \sum_{k=1}^{n-1}\frac{(-1)^{k}}{k!} \log^k\lvert x\lvert \, \l {n-k} (x) \, .
\end{equation}
Consequently,  to any (possibly with a logarithmic second member) AFE 
satisfied by Kummer's polylogarithm $\boldsymbol{\Lambda}_n$ corresponds an identity identically satisfied by Zagier's polylogarithm $\mathcal L_n$.
 \begin{center}
 $\star$
 \end{center}

In what follows,  $x$ and $y$  stand for real numbers such that $0<x<y<1$ and one sets $$
 \zeta=1-x \qquad  \quad \mbox{and} \qquad  \quad  \eta=1-y\, .$$

 
 \paragraph{Kummer's functional equation for the tetralogarithm.} 
 \label{Parag:Kummer_tetralog} 
 Kummer has established ({\it cf.} formula (144) p.\,366 in \cite{Kummer3}) that the tetralogarithmic function $\boldsymbol{\Lambda}_4$ satisfies a certain functional equation with 18 terms. 
The corresponding AFE satisfied by  the classical tetralogarithm $\l { 4} $ has been made explicit by Lewin (see formula (7.90) p.\,211 in  \cite{Lewin}) and is the following: 
\begin{align}
&\l{4}  \left(\,-\frac{x^2y\, \eta}{\zeta} \, \right)+  \l {4} \left(\,-\frac{y^2x\, \zeta}{\eta}\,\right) +
\l{4} \bigg(\,\frac{x^2y}{ \zeta \eta^2} \, \bigg) +\l {4} \bigg(\,\frac{y^2x}{\zeta^2\eta}\,
\bigg) \nonumber \\ 
 & -6\, \l {4} \Big(\,xy  \,\Big)  -6\,\l {4} \bigg( \,\frac{xy}{\eta \zeta} \,\bigg) 
-6\,\l {4} \bigg(\,-\frac{xy}{\eta} \, \bigg) 
-6\, \l {4} \bigg(\,-\frac{xy}{\zeta} \, \bigg)  
\nonumber \\ 
\big(\boldsymbol{\cal K}_4\Big) \qquad \qquad \quad 
& -3\,\l {4} \Big(\,x \,\eta \,\Big)  -3\, \l {4} \Big(\,y\,\zeta\,
\Big) -3\,\l {4} \Big(\,\frac{x}{\eta} \, \Big) 
-3 \,\l {4} \Big(\,\frac{y}{\zeta}\,  \Big)  
 \nonumber \\ 
 & -3\, \l{4} \Big(-\frac{x \,\eta}{\zeta}\,  \Big)  
-3\,\l{4} \bigg(-\frac{y \,\zeta}{\eta}  \bigg)  
-3\,\l {4} \bigg(\,  -\frac{x}{\eta \, \zeta} \,\bigg) 
-3 \, \l {4} \bigg(\,  -\frac{y}{ \eta \, \zeta  }\,  \bigg)  \nonumber
\\ 
 & +6\,  \l {4} \big(\, x\, \big)  +6\, \l{4}
\Big(-{x}/{\zeta} \Big)   +6\, \l{4} \big(y\big) 
+6\, \l {4} \Big(-{y}/{\eta} \Big) 
=  \frac{3}{2} \, {\log}(\zeta)^2 \, {\log}(\eta)^2\, .   \qquad \qquad 
\nonumber 
\end{align}

      Then the following element of $\mathbf Z\big[\mathbf C(x,y)\big]$ is annihilated by $\mathcal L_4$: 
      \begin{align}
&  \left[\,-\frac{x^2y\, \eta}{\zeta} \, \right]+ \left[\,-\frac{y^2x\, \zeta}{\eta}\,\right]+
 \bigg[\,\frac{x^2y}{ \zeta \eta^2} \, \bigg] +  \bigg[\,\frac{y^2x}{\zeta^2\eta}\,
\bigg] \nonumber \\ 
 &+6\, \left(\,  \big[\, x\, \big]  + \Big[-{x}/{\zeta} \Big]   + \big[y\big] 
+ \Big[-{y}/{\eta} \Big] \,  \right)   -6\, \left(\, \Big[\,xy  \,\Big]  +\, \bigg[ \,\frac{xy}{\eta \zeta} \,\bigg] 
+\, \bigg[\,-\frac{xy}{\eta} \, \bigg]
+\,  \bigg[\,-\frac{xy}{\zeta} \, \bigg]\, \right)
\nonumber \\ 
& -3\, \left(\, \Big[\,x \,\eta \,\Big]  +   \Big[\,y\,\zeta\,
\Big] + \Big[\,\frac{x}{\eta} \, \Big] 
+ \Big[\,\frac{y}{\zeta}\,  \Big]
%
  +  \Big[-\frac{x \,\eta}{\zeta}\,  \Big]  
+ \bigg[-\frac{y \,\zeta}{\eta}  \bigg]
+ \bigg[\,  -\frac{x}{\eta \, \zeta} \,\bigg]
+   \bigg[\,  -\frac{y}{ \eta \, \zeta  }\,  \bigg]  \, \right) \, .
\nonumber
%
\end{align}
%
%

      We consider the following 18 primitive rational functions, denoted by  $V_1,\dots,V_{18}$,  which appear as arguments of $\l{4} $ in Kummer's equation $(\boldsymbol{\cal K}_4)$: 
     $${x,\,y,\,-{\frac {y}{\eta\, \zeta }},\,-{\frac {x}{\eta\, \zeta  }},\,-{\frac {y \,\zeta }{\eta}},\,-{\frac {x \,\eta }{\zeta}},\,{\frac {y}{\zeta}},\,{\frac {x}{\eta}},\,y \,\zeta  ,\,x\, \eta  ,\,-{\frac {xy}{\zeta}},\,-{\frac {xy}{\eta}},\,{\frac {xy}{ \eta\,\zeta  }},\,xy,\,{\frac {x{y}^{2}}{\zeta^{2}\eta }},\,{\frac {{x}^{2}y}{\eta ^{2}\, \zeta }},\,-{\frac {{y}^{2}x \,\zeta }{\eta}},\,-{\frac {{x}^{2}y\, \eta }{\zeta}}
}\, .
$$

Then  {\bf Kummer's tetralogarithmic web} $\boldsymbol{{\boldsymbol{\mathcal W}}_{{\cal K}_4}}$ is the planar 18-web
defined by the  functions $V_i$'s.  
\mk 

\vspace{-0.3cm}
Since we will show later that a certain cluster web is equivalent to $\boldsymbol{{\boldsymbol{\mathcal W}}_{{\cal K}_4}}$, it is worth giving several of invariants of the latter: 
\begin{itemize}
\item Kummer's tetralogarithmic web 
 is not of maximal rank (since its total curvature does not vanish identically);  Some computations show that ${\rm rk}(\boldsymbol{{\boldsymbol{\mathcal W}}_{{\cal K}_4}})\leq 127$ (strictly less than ${\rho}(\boldsymbol{{\boldsymbol{\mathcal W}}_{{\cal K}_4}})= 136$); 
\sk 
\item For any $i=1,\ldots,18$, one has $\mathfrak B_i=\{0,1,\infty\}$ so this web has polylogarithmic ramification;
\sk 
\item Considering only iterated integrals ramified at $0,1$ or at infinity on $\mathbf P^1$, one gets: 
$$
{\rm polrk}^\bullet\Big({\boldsymbol{\mathcal W}}_{\hspace{-0.1cm}{\cal K}_4}\Big)=(28,28,16,3)\, ;
$$
\item Regarding the hexagonal and flat subwebs of $\boldsymbol{{\boldsymbol{\mathcal W}}_{{\cal K}_4}}$, one has: 
$$ \boldsymbol{\mathcal Hex}\Big({\boldsymbol{\mathcal W}}_{\hspace{-0.0cm}{\cal K}_4}\Big)=(186,112,36,4)   
\qquad \mbox{ and }\qquad  
\quad  \boldsymbol{\mathcal F \hspace{-0.1cm}lat} \Big({\boldsymbol{\mathcal W}}_{\hspace{-0.0cm}{\cal K}_4}\Big)=(186, 206, 60, 51, \ldots)
 \, . 
$$
\end{itemize}

 \paragraph{Kummer's pentalogarithmic functional equation.}
Kummer   proved  
({\it cf.}\,formula (148) p.\,368 in \cite{Kummer3},  or formula (7.104) p.\,216 in  \cite{Lewin}) that the pentalogarithmic 
 function $\Lambda_5$ satisfies the following functional identity
 for any reals $x$ and $y$ such that $0<x<y<1$: 
\begin{align*}
&\hspace{0.2cm} \textcolor{white}{+} 
\boldsymbol{\Lambda}_5 \left(\,  \frac{x^2y\eta}{\zeta} \right)+
\boldsymbol{\Lambda}_5 \left(\,-\frac{x^2y}{\zeta\eta^2} \,  \right)+
\boldsymbol{\Lambda}_5 \left(\,-\frac{x^2\eta}{ y^2 \zeta}  \, \right)+
\boldsymbol{\Lambda}_5 \left(\,-\frac{x\zeta}{y\eta}  \, \right)+
\boldsymbol{\Lambda}_5 \left(\,\frac{x\zeta\eta^2}{y}  \, \right)  \\ 
 & \,  +
\boldsymbol{\Lambda}_5 \left(\,\frac{x y^2\zeta }{\eta}  \, \right)+
\boldsymbol{\Lambda}_5 \left(\,\frac{x}{y \zeta^2 \eta} \right)+
\boldsymbol{\Lambda}_5 \left(\,-\frac{x\eta^2}{y \zeta^2} \right)+ 
\boldsymbol{\Lambda}_5 \left(\,-\frac{xy^2}{\zeta^2\eta} \right)
 -9  \,  \boldsymbol{\Lambda}_5 \Big(  -xy   \Big)  \\ &    -9  \,  
  \boldsymbol{\Lambda}_5 \left(\,  -\frac{x}{y}  \, \right)   -9  \,   
   \boldsymbol{\Lambda}_5 \Big(  -x\eta   \Big) 
    -9  \,  
    \boldsymbol{\Lambda}_5 \left(\,   -\frac{x}{\eta}   \, \right)
    -9  \,  
    \boldsymbol{\Lambda}_5 \left(\,   \frac{xy}{\eta}   \, \right)  -9  \,   
     \boldsymbol{\Lambda}_5 \left(\,   \frac{x\eta}{y} \, \right)  \\ &   -9  \,   
      \boldsymbol{\Lambda}_5 \left(\,   \frac{xy}{\zeta}   \, \right)   -9  \,   
       \boldsymbol{\Lambda}_5 \left(\,   \frac{x}{y\zeta} \, \right)    -9  \,  
        \boldsymbol{\Lambda}_5 \left(\, \frac{x\eta}{\zeta}   \, \right)  -9  \,  
         \boldsymbol{\Lambda}_5 \left(\,   \frac{x}{\zeta \eta}  \, \right)  -9  \,  
          \boldsymbol{\Lambda}_5 \left(\,  -\frac{xy}{\zeta \eta}  \, \right)    -9  \,  
           \boldsymbol{\Lambda}_5 \left(\,  -\frac{x\eta}{y\zeta}   \, \right)    \\  \,  &   -9  \,  
          \boldsymbol{\Lambda}_5 \Big(  - y\zeta   \Big)    -9  \,  
             \boldsymbol{\Lambda}_5 \Big( -\zeta\eta  \Big)  -9  \,  
              \boldsymbol{\Lambda}_5 \left(\,  \frac{y\zeta }{\eta} \, \right)   -9  \,   
               \boldsymbol{\Lambda}_5 \left(\,   -\frac{y}{\zeta}   \, \right)    -9  \,  
                \boldsymbol{\Lambda}_5 \left(\,     -\frac{\eta}{\zeta} \right)  
                 -9  \,  
                 \boldsymbol{\Lambda}_5 \left(\,  \frac{y}{\zeta\eta}  \, \right) 
                 \\ & \,
                  +18\, 
  \boldsymbol{\Lambda}_5\Big(
 -x\Big) +18\, \boldsymbol{\Lambda}_5\Big(-\zeta\Big) 
 +18\, \boldsymbol{\Lambda}_5
  \left(\frac{x}{\zeta}\right)  +
  18  \,\boldsymbol{\Lambda}_5\Big(-y\Big)+ 18 \,\boldsymbol{\Lambda}_5\Big( -\eta\Big) +18\,  \boldsymbol{\Lambda}_5\left(\frac{ y}{\eta}\right)
   \\ &\, 
     = 18 \log\bigg(\frac{x}{y}\bigg) \log(\zeta)^4
     - 108 \, \log(y) \log(\zeta)^2\log(\eta)^2 - \frac{36}{5}\log( \zeta)^5 +36 \log(\eta)^3 \log(\zeta)^2 - 18 \, \zeta(5)\, . 
\end{align*}

\begin{rem}
The  functional equation for the classical pentalogarithm  $\l {5} $ corresponding to the one just above  has been worked out in {\rm \cite{ALewin}} (see also (3.12) p.\,30 
 in {\rm \cite{LewinStruct}}). If the sum of the 33 pentalogarithmic terms 
 (each multiplied with the suitable integer coefficient)
 given there is the correct one (since identity $(\boldsymbol{\mathcal K}_5)$ below holds true identically), it does not seem to be the case for the logarithmic second member $18 \,\zeta(5)+18\,\zeta(4)\,\log(\zeta)+3\,\zeta(2) \,\log(\zeta)^2\log(\zeta/\eta^3)+\frac{3}{2}\log(\zeta)^2\log(\eta)^2\log(y^3/\eta)+\frac{3}{20} \log(\zeta)^4\log(y^5\zeta^2/x^5)
$. There is probably a typographical error~here.
\end{rem}
\sk

In the light of  \eqref{Eq:Lin-Lambdan}, if one considers the opposite of the arguments of $\boldsymbol{\Lambda}_5$ in the  equation above, one gets a functional identity satisfied by the  classical pentalogarithm. Indeed, setting 
\begin{align*}
(U_i)_{i=1}^9= & \, \left( \, -\frac{x^2y\eta}{\zeta} 
  \, , \, - \frac{x y^2\zeta }{\eta} 
\, , \,   \frac{x^2y}{\zeta\eta^2}
  \, , \,\frac{xy^2}{\zeta^2\eta}
\, , \,\frac{x^2\eta}{ y^2 \zeta} 
   \, , \,\frac{x\eta^2}{y \zeta^2}
  \, , \,-\frac{x}{y \zeta^2 \eta}
 \, , \,-\frac{x\zeta\eta^2}{y}
  \, , \,\frac{x\zeta}{y\eta}\, 
   \right)\\ 
(V_j)_{j=1}^{18}=   & 
\left( \, 
xy \, , \frac{x}{y}\, ,\,   x\,\eta\,  , \frac{x}{\eta}\,  ,  -\frac{xy}{\eta}\, ,\,   - \frac{x\eta}{y}\, , \, 
- \frac{xy}{\zeta}\, , \, - \frac{x}{y\zeta}\, , \,  - \frac{x\eta}{\zeta}\, , - \frac{x}{\zeta\eta}\, ,\,   \frac{xy}{\zeta\eta}\, ,\,   \frac{x\eta}{y\zeta}\, , \, y\zeta, \zeta\eta, - \frac{y\zeta}{\eta}\, ,
  \frac{y}{\zeta}\, ,\,  \frac{\eta}{\zeta}, \, -\frac{ y}{\zeta\eta}\, 
\right)
\\
 (W_k)_{k=1}^6= & \,  \left( \, x 
  \, , \,  \zeta
  \, , \,-\frac{x }{\zeta } 
    \, , \,y
      \, , \,  \eta 
        \, , \,-\frac{y }{\eta } 
   \right)\, , 
\end{align*}
one verifies that Zagier's pentalogarithm $\mathcal L_5$ satisifies identically  the following functional identity: 
\begin{align*}
(\boldsymbol{\mathcal K}_5) \hspace{2cm}
\sum_{i=1}^9 \mathcal L_5 (U_i ) 
-9 \sum_{j=1}^{18} \mathcal L_5(V_j)
+18 \sum_{k=1}^6 
\mathcal L_5(
W_k)= 18 \, \zeta(5)\, .  \hspace{3cm} \textcolor{white}{.}
\end{align*}

This equation, that we call {\bf Kummer's 
      functional equation of the pentalogarithm}, involves  33 pentalogarithmic terms, but the arguments $W_k$ for $k=1,\ldots,6$ altogether define only two foliations, since obviously $\mathcal F_x=\mathcal F_\zeta=\mathcal F_{x/\zeta}$,  and similarly
$\mathcal F_y=\mathcal F_\eta=\mathcal F_{y/\eta}$. Consequently, the web naturally associated to $(\boldsymbol{\mathcal K}_5)$ is the following 
29-web 
$$
\boldsymbol{\mathcal W}_{\boldsymbol{\mathcal K}_5}=\boldsymbol{\mathcal W}\left( \, x\, , \,  y \, , 
 \, U_i \, , \, V _j\,  \Big\lvert  \, \begin{tabular}{l}
$ i=1,\ldots,9$ \vspace{-0.1cm}
\\  ${}^{}$\hspace{-0.2cm} $ j=1,\ldots,18$
\end{tabular}
\right)\, . 
$$

Here are some features of 
$\boldsymbol{\mathcal W}_{\boldsymbol{\mathcal K}_5}$
 which have been established  by direct computations: 
\begin{itemize}
\item 
 The total curvature of 
 Kummer's pentalogarithmic web   does not vanish identically hence 
 this web is not of maximal rank:  ${\rm rk}(\boldsymbol{{\boldsymbol{\mathcal W}}_{{\cal K}_5}})$ 
is strictly less than ${\rho}(\boldsymbol{{\boldsymbol{\mathcal W}}_{{\cal K}_5}})= 378$; 
\sk 
\item For any $i=1,\ldots,29$, one has $\mathfrak B_i=\{0,1,\infty\}$ so as  Bol's web, Spence-Kummer's web and Kummer's tetralogarithmic web,  
$\boldsymbol{\mathcal W}_{\boldsymbol{\mathcal K}_5}$  has polylogarithmic ramification as well;
\sk 
\item Considering iterated integrals ramified at $0,1$ or at infinity on $\mathbf P^1$, one gets: 
$$
{\rm polrk}^\bullet\Big({\boldsymbol{\mathcal W}}_{\hspace{-0.1cm}{\cal K}_5}\Big)=\big( 48, 53, 42, 22, 4\big)\, ;
$$
\item Regarding the count of hexagonal subwebs of $\boldsymbol{{\boldsymbol{\mathcal W}}_{{\cal K}_5}}$, one has: 
$$ \boldsymbol{\mathcal Hex}^\bullet\Big({\boldsymbol{\mathcal W}}_{\hspace{-0.0cm}{\cal K}_5}\Big)=\big( 474, 288, 105 , \ldots \big)   
 \, . 
$$
\end{itemize}

Because Bol's web $\boldsymbol{\mathcal B}$, Spence-Kummer's trilogarithmic one $\boldsymbol{{\boldsymbol{\mathcal W}}_{\hspace{-0.05cm}{\cal S\cal K}}}$ and  Kummer's web $\boldsymbol{{\boldsymbol{\mathcal W}}_{{\cal K}_4}}$ are of cluster type (see Theorem \ref{T:classical-cluster-webs} below), it is natural to ask whether the same holds true for $\boldsymbol{{\boldsymbol{\mathcal W}}_{{\cal K}_5}}$. That this web has polylogarithmic ramification  let us think that it might be the case.  This is the reason why we have given the few numerical invariants of $\boldsymbol{{\boldsymbol{\mathcal W}}_{{\cal K}_5}}$ above: any (cluster?) planar 29-web 
sharing the same invariants would be a good candidate for being equivalent to $\boldsymbol{{\boldsymbol{\mathcal W}}_{{\cal K}_5}}$.

\subsubsection{Webs associated to some recent polylogarithmic AFEs.}
\label{SS:WebsRecentAFEs}
The problem of finding and better understanding  functional equations of higher polylogarithms is an active area of current research. In this subsection, we briefly discuss some equations of this type recently  discovered by several authors. 
The corresponding webs are  typically formed by an important number  of foliations (several hundreds) hence are out of reach  by the formal/algebraic computations (of their ARs, their rank, etc.) which would allow to understand them better. So we do not have so much to say about them. They are essentially mentioned here to show that there are 
many interesting polylogarithmic webs in  weight higher than 2 or 3 which would be worth studying.

%
%
%
%

\paragraph{Remarks about a series of tetralogarithmic AFEs discovered by Gangl.}  
In \cite{Gangl}, Gangl exhibits a whole series 
 of AFEs for the tetralogarithm $\mathcal L_4$ that we are going to discuss succinctly from the perspective of webs.\mk 
 
Let $n>0$ be a fixed integer, set  $\phi(t)=t^{n-1}(t-1)$ and given two indeterminates $x$ and $y$, one denotes by $x_1,\ldots,x_n$ 
(resp.\,by $y_1,\ldots,x_n$) the 
roots of the polynomial equation $\phi(t)=x$ (resp.\,$\phi(t)=y$).
 Then, according to Theorem 4.1 in \cite{Gangl}, the following element is annihilated by $\mathcal L_4$: 
\begin{align*}
\boldsymbol{\mathcal Gan}(n)=n(n-2)\left[ \frac{\prod_{i=1}^n x_i}{\prod_{j=1}^n y_j}   \right]
&\, -(n-1)^2\sum_{i,j=1}^n\left[ \frac{y_j(1-x_i)}{x_i(1-y_j)}   \right]+n^2 
\sum_{i,j=1}^n
\left[\frac{1-x_i}{1-y_j}    \right] \\
&\, -  n^2(n-1)^2 \sum_{i,j=1}^n
\left[\frac{x_i}{y_j}    \right]
+n(n-1)^2\sum_{i=1}^n
\bigg( \, \Big[1-x_i^{-1}\Big]  - \Big[1-y_i^{-1}\Big]\,  \bigg)
\, . \sk
\end{align*}

It is necessary to assume $n>1$ otherwise the statement is not correct: indeed, it would be equivalent to the relation $\mathcal L_4((1+x)/(1+y))+\mathcal L_4(x/y)=0$, which obviously is not true. Nevertheless, the case $n=1$ is still interesting from the point of view of webs, as will be discussed below.  Note that  assuming $n>1$, 
the first term of $\boldsymbol{\mathcal Gan}(n)$ simplifies since 
 ${\prod_{i=1}^n x_i}/{\prod_{j=1}^n y_j}$ is equal to $x/y$. \mk 
 
The case $n=2$  is also particular and should be disregarded. Indeed, in this case one has $x_1=(1+\sqrt{1+4x})/2 $ and $x_2=(1-\sqrt{1+4x})/2$ with the same formulas,  but where $x$ has been replaced by  $y$, for $y_1$ and $y_2$.  Thus thanks to the inversion formula 
$\mathcal L_4([x]+[x^{-1}])=0$, 
it comes that each of the terms $\sum_{i=1}^2 [1-x_i^{-1}]$ and $\sum_{i=1}^2 [1-y_i^{-1}]$ is already annihilited by $\mathcal L_4$.  The same inversion formula also implies that  
$ \sum_{i,j=1}^n
[(1-x_i)/(1-y_j)]$ vanishes when evaluated by $\mathcal L_4$. 
 Finally,  each term of the sum 
$\sum_{i,j=1}^n
[(1-x_i)/(1-y_j)]   $ is the inverse of one of the terms appearing in 
$\sum_{i,j=1}^n [{x_i}/{y_j} ]  $.  
Thus ${\boldsymbol{\mathcal Gan}(2)}$ is  a linear combination of elements of the form  $ [R]+[1/R]$ for some rational function $R\in \mathbf C(x,y)$. The associated tetralogarithmic functional equation is accessible from the inversion relation, which makes it 
not so interesting: it is not a genuine two-variables functional equation for the tetralogarithm.  Nevertheless,  again, this case is still interesting  from the point of view of webs (see below). 
\mk 

Since it is indeed the case when $n=2$, one may wonder if the identity $\mathcal L_4({\boldsymbol{\mathcal Gan}(n)})=0$ is not accessible from the inversion identity for any $n\geq 2$.  We believe that it is not the case and that one gets genuine tetralogarithmic identities  starting from $n=3$.  We are brought to this conclusion by looking at the web whose first integrals are the   algebraic functions appearing in the definition of $\boldsymbol{\mathcal Gan}(n)$: 
\begin{equation}
\label{Eq:Gan_n}
\boldsymbol{\mathcal W}_
{\boldsymbol{\mathcal Gan}(n)}=\boldsymbol{\mathcal W}
\left( \, G_n \, , \, \frac{y_j(1-x_i)}{x_i(1-y_j)}
\, , \, \frac{1-x_i}{1-y_j}  
\, , \, \frac{x_i}{y_j}   
 \, , \, x_i \,  , \, y_j
\hspace{0.15cm}
\Big\lvert 
\hspace{0.15cm}
\begin{tabular}{l}
$i=1,\ldots,n$\vspace{-0.1cm} \\
${}^{}$\hspace{-0.1cm}$j=1,\ldots,n$
\end{tabular}
\right)
\end{equation}
with $ G_1= (x+1)/(y+1)$ and $G_n= {x}/{y}$  for $n\geq 2$. 
 Since, as functions of the variables $x,y$, the $x_i$'s only depend on $x$ whereas the $y_j$'s only on $y$, it follows that 
 $\mathcal F_{x_i}=\mathcal F_x$ and $\mathcal F_{y_i}=\mathcal F_y$ for any $i$. Consequently,  Gangl's web $\boldsymbol{\mathcal W}_
{\boldsymbol{\mathcal Gan}(n)}$ actually is formed by 
at most  $3(n^2+1)$ foliations:  one has 
$$
\boldsymbol{\mathcal W}_
{\boldsymbol{\mathcal Gan}(n)}=\boldsymbol{\mathcal W}
\left( \, x\, , \, y\, ,  G_n\, , \, \frac{y_j(1-x_i)}{x_i(1-y_j)}
\, , \, \frac{1-x_i}{1-y_j}  
\, , \, \frac{x_i}{y_j}   
\hspace{0.15cm}
\Big\lvert 
\hspace{0.15cm}
\begin{tabular}{l}
$i=1,\ldots,n$\vspace{-0.1cm} \\
${}^{}$\hspace{-0.1cm}$j=1,\ldots,n$
\end{tabular}
\right)\, . \sk 
$$
By direct computations, we have verified that when $n=3$, the corresponding 30 foliations are pairwise distinct hence $
\boldsymbol{\mathcal W}_
{\boldsymbol{\mathcal Gan}(3)}$ is  a genuine 30-web.  Thus, in the  identity $\mathcal L_4({\boldsymbol{\mathcal Gan}(3)})=0$, there are  30 purely tetralogarithmic terms which are pairwise functionally independent. This contrasts with the case $n=2$ for which precisely the opposite holds true.\footnote{Namely, given any of the $17$ elements appearing in \eqref{Eq:Gan_n} when $n=2$, there exists another one which is its inverse.} 
We suspect that,  
 for any $n\geq 3$,  $
\boldsymbol{\mathcal W}_
{\boldsymbol{\mathcal Gan}(n)}$ is a planar $3(n^2+1)$-web which carries a genuine tetralogarithmic AR not accessible from the inversion relation, but we have not any  proof of that.\footnote{The first question regarding that is to know whether  or not the functional equation 
$\mathcal L_4({\boldsymbol{\mathcal Gan}(3)})=0$ is accessible from the inversion relation. It seems to us that this is not an easy question to answer. 
} \mk

The $n=1$ and $n=2$ cases, although not really compelling regarding  tetralogarithmic functional equations as seen just above, are quite interesting from the  perspective of web geometry.  First, we remark that the definition of $
\boldsymbol{\mathcal W}_
{\boldsymbol{\mathcal Gan}(n)}$ makes sense when $n=1$ or $n=2$ as well. For $n=1$, one has 
$$
\boldsymbol{\mathcal W}_
{\boldsymbol{\mathcal Gan}(1)}=\boldsymbol{\mathcal W}
\left( \, x\, , \, y\, , \frac{x}{y}\, , \, \frac{1+x}{1+y}\, , \,  \frac{x(1+y)}{y(1+x)}\,\right)\, , 
$$

which is nothing else but Bol's web, associated to Abel's 5-terms relation of the dilogarithm. \mk 

As for the case $n=2,$ it is no less interesting:  one verifies that the web  $
\boldsymbol{\mathcal W}_
{\boldsymbol{\mathcal Gan}(2)}$ is actually a 9-web defined by 
algebraic first integrals, whose explicit formulas are not difficult to give:  
\begin{align*}
\frac{x}{y}
\, , \,  
\frac{ \Big( \sqrt{1 + 4x}-1\Big)\Big(\sqrt{1 + 4y}+1\Big)}{(\sqrt{1 + 4x}+1) (\sqrt{1 + 4y}-1)}  & 
 \, , \, 
\frac{  \Big(\sqrt{1 + 4x}+1\Big)\Big(\sqrt{1 + 4y}+1\Big)}{\Big(\sqrt{1 + 4x}-1\Big)\Big(\sqrt{1 + 4y}-1\Big)}
  \, , \, \ldots \\
  \ldots & \,  
    , \,   \frac{\sqrt{1 + 4x}-1}{\sqrt{1 + 4y}+1} 
     \, , \, 
       \frac{\sqrt{1 + 4x}+1}{\sqrt{1 + 4y}+1}
      \, , \, 
        \frac{\sqrt{1 + 4x}-1}{\sqrt{1 + 4x}+1}
       \, , \, 
       \frac{\sqrt{1 + 4y}-1}{\sqrt{1 + 4y}+1}\, . \qquad \textcolor{white}{ooo}
 \end{align*}


Then, setting $x=(X^2 - 1)/4$ and $y=(Y^2 - 1)/4$, one defines a rational  map $\varphi : (X,Y)\mapsto (x,y)$ which is such that $\varphi^*(
\boldsymbol{\mathcal W}_
{\boldsymbol{\mathcal Gan}(2)} ) $ 
is defined by the following  rational first integrals: 
  \begin{align*}
  \frac{ X+1}{Y+1}
       \hspace{0.2cm} , \hspace{0.2cm}
        & \frac{ X-1}{X+1}
       \hspace{0.2cm} ,  \hspace{0.2cm}
   \frac{     Y-1}{Y+1}  \hspace{0.2cm} , \hspace{0.1cm}
      -  \,  \frac{ X+1}{Y-1}
   \hspace{0.2cm} , \hspace{0.1cm}
   - \, \frac{  X-1}{Y+1}
      \hspace{0.1cm} , \hspace{0.2cm}
  \vspace{0.2cm} \\
   &  \frac{X-1}{Y-1}
        \hspace{0.2cm} , \hspace{0.2cm} 
 \frac{X^2 - 1}{Y^2 - 1}
 \hspace{0.2cm} , \hspace{0.2cm}
 \frac{(X-1)(Y+1)}{(X+1)(Y-1)}
 \hspace{0.4cm} \mbox{ and }    \hspace{0.4cm} 
  \frac{ (X+1)(Y+1)}{(X-1)(Y-1)}
 \hspace{0.2cm} .
   \end{align*}
It is not difficult to recognize the nine rational first integrals defining (the projective canonical model of) Spence-Kummer web 
${\boldsymbol{\mathcal W}}_{\!{\cal S}{\cal K}}$
associated to the famous functional equation $(\mathcal S\mathcal K)$ satisfied by the trilogarithm (see \S\ref{SS:Spence-Kummer} above). \mk 

It is quite surprising that the first two of the webs associated to a 
 series of tetralogarithmic functional equations specialize to the two webs (namely, Bol's web and Spence-Kummer web) 
  which actually do not carry any tetralogarithmic functional equation, but  the simplest truly multivariable AFEs of the dilogarithm and the trilogarithm respectively.\footnote{Since both 
  $\boldsymbol{\mathcal W}_{\boldsymbol{\mathcal Gan}(1)}\simeq \boldsymbol{\mathcal B}$ and 
$ \boldsymbol{\mathcal W}_
{\boldsymbol{\mathcal Gan}(2)} \simeq {\boldsymbol{\mathcal W}}_{\!{\cal S}{\cal K}}$ are webs with maximal rank, a natural question about the whole series of webs $
\boldsymbol{\mathcal W}_{\boldsymbol{\mathcal Gan}(n)}$'s is whether  or not this still holds true for $n\geq 3$. It doesn't seem easy to answer: even with the help of a computer, we haven't been able to verify if  the curvature of the 30-web 
$\boldsymbol{\mathcal W}_{\boldsymbol{\mathcal Gan}(3)}$  vanishes identically...} 
We find this  very intriguing  
  and have absolutely no idea of any explanation of this fact (if there is one, of course).


\paragraph{A   functional equation with 40 tetralogarithmic terms.}  
\label{Par:GSVV}
In  \cite{GGSVV}, the authors study what they call {\it `motivic amplitudes'}, that is motivic objects containing all the mathematical content of some scattering amplitudes (as they describe it). They find that these objetcs can be better understood by means of the `cluster structures' of spaces of projective configurations ${\rm Conf}_n(\mathbf P^3)$. In particular, as a byproduct of their approach, they `{\it find and prove the first known functional equation for the
trilogarithm in which all 40 arguments are cluster $\boldsymbol{\mathcal X}$-coordinates of a single algebra}'. 
\mk 

The functional equation discussed here lives on the space 
${\rm Conf}_6(\mathbf P^2)$ of configurations of six points in 
the projective plane.  We denote here by $K$ the field of rational functions on this space (which is rational, hence $K$ is isomorphic to $\mathbf C(t_1,\ldots,t_4)$ but not in a canonical way). 
For any $i\in \{1,\ldots,6\}$, one denotes by 
 $\pi_i: {\rm Conf}_6(\mathbf P^2)\dashrightarrow {\rm Conf}_5(\mathbf P^1)$ 
 (resp.\,$f_i : {\rm Conf}_2(\mathbf P^6)\dashrightarrow {\rm Conf}_5(\mathbf P^2)$) the rational map induced by 
 projecting from (resp.\,forgetting) 
 the $i$-th point of a configuration. Let  also $(2,4,3, 25\cap 36)$  be the map associating to $(x_i)_{i=1}^6\in {\rm Conf}_6(\mathbf P^2)$ the configuration of four points $(x_2,x_4,x_3, x_{25\cap 36}) 
  \in {\rm Conf}_4(\mathbf P^2)$
  where $x_{25\cap 36}$ stands for the intersection of the two lines $\langle x_2,x_5\rangle$ and $\langle x_3,x_6\rangle$ in $\mathbf P^2$; 
  Then one defines two rational functions 
 on ${\rm Conf}_6(\mathbf P^2)$ by setting:
 \begin{equation}
 \label{Eq:r1-and-r3}
  r_1={\rm Cr}\circ \pi_1\circ f_6
 \qquad \mbox{ and }\qquad  
 r_3={\rm Cr}\circ \pi_1(2,4,3, 25\cap 36)\, .
 \end{equation}
 (The rational function $r_3$ is the {\it triple ratio} introduced 
  in \cite{Goncharov1995Adv}, see also \S\ref{Par:ChowPol-Webs} below.)\sk
 
 For any $\Upsilon\in \mathbf Q[K^{**}]$, let  $\boldsymbol{{\rm Cyc}_6(\Upsilon)}$ (resp.\,$\boldsymbol{{\rm Acyc}_6(\Upsilon)}$) be the sum $\sum_{\tau } \Upsilon^{\tau}$  where $\tau$ ranges in the set of cyclic (resp.\,anticyclic)  permutations of $\{1,\ldots,6\}$.\footnote{{\it I.e.} $\tau\in \{ \alpha^k\}_{k=0}^5$ with $\alpha=(612345)$ in the cyclic case  and 
 $\tau\in \{ \alpha^k\beta\}_{k=0}^5$ with $\beta=(654321)$ in the anticyclic one.} 
 Denoting respectively by $a,b$ and $c$  the permutations $(3456)$ and $(4326)$ and $(235)(46)$,  one defines elements of $\mathbf Q[K^{**}] $  
 by setting   
$$
\Gamma=\Gamma_1+\Gamma_3\qquad \mbox{ with } 
\qquad \Gamma_1=[r_1]+\big[ r_1^{a}\big]\quad  \mbox{ and } \quad 
\Gamma_3=\big[r_3^{b}\big]+
 (1/3)\big[r_3^{c}\big]\, .
$$
Then according to Theorem B.1 in \cite{GGSVV}, 
$$
\boldsymbol{\mathcal G \mathcal G \mathcal S \mathcal V\mathcal V}
={\rm Cyc}_6\Big(   \Gamma\Big) - 
{\rm Acyc}_6\Big(   \Gamma\Big)\in \mathbf Q[K^{**}]
$$
 is annihilated by $\mathcal L_3$. 
 Moreover, if each sum ${\rm Cyc}_6(\gamma) - 
{\rm Acyc}_6(\gamma)$  involves 12 terms for $\gamma=[r_1],[r_1^a]$ and $[r_3^b]$, 
this is not the case for    ${\rm Cyc}_6([r_3^c]) - 
{\rm Acyc}_6([r_3^c])$ which is a linear combination of only four primitive elements of 
$ \mathbf Q[K^{**}]$. Consequently, the functional equation $\mathcal L_3(\boldsymbol{\mathcal G \mathcal G \mathcal S \mathcal V\mathcal V})=0$ involves only $3\times 12+4=40$ trilogarithmic terms and the associated web
$$
{\boldsymbol{\mathcal W}}_{\boldsymbol{\mathcal G \mathcal G \mathcal S \mathcal V\mathcal V}}=
{\boldsymbol{\mathcal W}}\left(
\, r_1^\tau \, , \, r_1^{a\tau} \, , \, r_3^{b\tau}\, , \, 
 r_3^{c\tau} \, \Big\lvert \, \tau \, \mbox{ cyclic or anticyclic}\, 
\right)
$$
is a 40-web on ${\rm Conf}_6(\mathbf P^2)$. By some direct computations, one establishes that 
$$ 
\rho^\bullet\Big({\boldsymbol{\mathcal W}}_{\boldsymbol{\mathcal G \mathcal G \mathcal S \mathcal V\mathcal V}}\Big)= \Big(36, 30, 20, 8, 1,0\Big)  \qquad \mbox{ and }
\qquad 
{\rm polrk}^\bullet\Big({\boldsymbol{\mathcal W}}_{ \hspace{-0.05cm} \boldsymbol{\mathcal G \mathcal G \mathcal S \mathcal V\mathcal V}}\Big)=(64,29,2)
$$
from which it comes that 
$ 
\rho\big({\boldsymbol{\mathcal W}}_{\hspace{-0.05cm}\boldsymbol{\mathcal G \mathcal G \mathcal S \mathcal V\mathcal V}}\big)
= {\rm polrk}\big({\boldsymbol{\mathcal W}}_{ \hspace{-0.05cm} \boldsymbol{\mathcal G \mathcal G \mathcal S \mathcal V\mathcal V}}\big)=95$. It follows that ${\boldsymbol{\mathcal W}}_{ \hspace{-0.05cm} \boldsymbol{\mathcal G \mathcal G \mathcal S \mathcal V\mathcal V}}$ is AMP, with only polylogarithmic ARs, of weight 1,2 and 3. 
\mk 

The web ${\boldsymbol{\mathcal W}}_{ \hspace{-0.05cm} \boldsymbol{\mathcal G \mathcal G \mathcal S \mathcal V\mathcal V}}$ contains interesting subwebs. Indeed,  
setting
$$
{\boldsymbol{\mathcal W}}\boldsymbol{\Gamma}_1=
{\boldsymbol{\mathcal W}}\left(
\, r_1^\tau \, , \, r_1^{a\tau} \, 
\right)
\qquad \mbox{ and }
\qquad 
{\boldsymbol{\mathcal W}}\boldsymbol{\Gamma}_3=
{\boldsymbol{\mathcal W}}\left(
\, r_3^{b\tau} \, , \, r_1^{c\tau} \, 
\right)
$$
where, as above,  $\tau$ ranges in the set of cyclic or anticyclic permutations of $\mathfrak S_6$, one verifies that 
\begin{align*}
\rho^\bullet\big({\boldsymbol{\mathcal W\Gamma}_1}\big)
= &\, (20,14,5) &&  {\rm polrk}^\bullet\big({\boldsymbol{\mathcal W\Gamma}_1}\big)=(34,5)\\
\mbox{ and }
\quad 
\rho^\bullet \big({\boldsymbol{\mathcal W\Gamma}_3}\big)=& \, (12,6) &&  {\rm polrk}^\bullet\big({\boldsymbol{\mathcal W\Gamma}_3}\big)=(16,2)\, , 
\end{align*}
from which it comes that these two subwebs of ${\boldsymbol{\mathcal W}}_{ \hspace{-0.05cm} \boldsymbol{\mathcal G \mathcal G \mathcal S \mathcal V\mathcal V}}$ are AMP as well, with only polylogarithmic ARs, of weight at most 2.
\mk 

In \S\ref{Par:Remarks-trilog-AR-XWD4} below, we will recognize and describe the web ${\boldsymbol{\mathcal W}}_{ \hspace{-0.05cm} \boldsymbol{\mathcal G \mathcal G \mathcal S \mathcal V\mathcal V}}$ as a subweb of the $\boldsymbol{\mathcal X}$-cluster web of type $D_4$, which is a 52-web in four variables.

\paragraph{Chow polylogarithms and associated webs.}
\label{Par:ChowPol-Webs}
In \cite{Goncharov1995}, for each integer $n\geq 2$, by restricting to linear subspaces the  {\it $n$-th Chow polylogarithm} $\mathcal P_n$ he has constructed by means of an integral transform of a current,  on the space of algebraic $n$-cycles in $\mathbf P^{2n-1}$ intersecting properly the faces of a fixed simplex in it, Goncharov obtains the $n$-th {\bf grassmannian polylogarithm} 
$\mathcal L_n^G$.  \sk

It is a real analytic function  defined on a Zariski open-subset   
 $G_n^0(\mathbf C^{2n})$ of the grassmannian 
 of $n$-planes in $\mathbf C^{2n}$ which enjoys interesting properties. 
First,  $\mathcal L_n^G$ is constant along the orbit of  a  $n$-plane 
belonging to $G_n^0(\mathbf C^{2n})$ 
under the standard (diagonal) action of the maximal torus $H_{2n-1}$ of ${\sf SL}_{2n}(\mathbf C)$.  Hence, thanks to the natural identification between the quotient $
G_n^0(\mathbf C^{2n})/H_{2n-1}$ and the space  of {\it generic} configurations ${\rm Conf}_{2n}^0(\mathbf P^{n-1})$ (see \cite[\S2.2.2]{GelfandMcPherson}), $\mathcal L_n^G$ can be seen as a function on the latter. Viewed this way, it satisfies the following  functional equations  
\begin{equation}
\label{Eq:Chow-Polylog-AFE}
\sum_{i=1}^{2n+1}
 (-1)^i
 \mathcal L_n^G\big(f_i\big)=0 \qquad 
\qquad 
\mbox{ and } \qquad \qquad 
\sum_{j=1}^{2n+1}
 (-1)^j
 \mathcal L_n^G\big(\pi_j\big)=0
\end{equation}
which hold true on the configuration spaces ${\rm Conf}_{2n+1}^0(\mathbf P^{n-1})$ and ${\rm Conf}_{2n+1}^0(\mathbf P^{n})$ respectively. \big(Here, as in the preceding paragraph, $f_i: {\rm Conf}_{2n+1}(\mathbf P^{n-1})\rightarrow {\rm Conf}_{2n}(\mathbf P^{n-1})$ stands for the map given by forgetting  the $i$-th point of a configuration  and $\pi_i:  {\rm Conf}_{2n+1}(\mathbf P^{n-1})\rightarrow {\rm Conf}_{2n}(\mathbf P^{n-1})$ for the one induced by the linear projection $\mathbf P^{n}\dashrightarrow \mathbf P^{n-1}$ from the $i$-th point.\big)\sk 

Actually, as proved by Goncharov, $\mathcal L_n^G$ can be defined on the whole 
 space ${\rm Conf}_{2n}(\mathbf P^{n-1})$ but as such, is not smooth (and even continuous?) on it. Anyway,  Goncharov shows that there exists a stratum of non-generic configurations $c_z$ depending on  $z\in \mathbf C\setminus \{0,1\}$ such that for any such complex number, one has $\mathcal L_n^G(c_z)=\mathcal L_n(z)$ which indicates that $\mathcal L_n^G$ is related to the classical $n$-th polylogarithm.

\begin{rem}
Associated to the two functional equations \eqref{Eq:Chow-Polylog-AFE}, we can consider the two webs
with first integrals, the forgetful maps $\varphi_i$ on the one hand, and the projections $\pi_j$ on the  other hand. These are webs with $n$-dimensional leaves  which turn out to be the same up to the Gelfand-MacPherson correspondence (or Gale Transform)  ${\rm Conf}_{2n+1}(\mathbf P^{n-1})\simeq {\rm Conf}_{2n+1}(\mathbf P^{n})$.   The identity $\sum_{i=1}^{2n+1}  (-1)^i \mathcal L_n^G(f_i)=0$ can be seen as a real-analytic abelian relation of order 0 for this web. It would be interesting to know more about these webs, in particular about their $k$-abelian relations and their $k$-rank, for $k=0,\ldots,n$.\footnote{For the notions of {\it $k$-abelian relation} and {\it $k$-rank} of a web of codimension $c\geq 1$ (with $k\leq c$), see \cite{GriffithsHopkins} and \cite{Henaut1}.}  
\end{rem}

For $n=2$, one recovers Bol's web:  indeed for $\boldsymbol{x}=[x_1,\ldots, x_4]\in {\rm Conf}_4(\mathbf P^1)$, 
 one has $\mathcal L^G_2(\boldsymbol{x})=\mathcal L_2({\rm Cr}(x_1,\ldots, x_4))$ hence \eqref{Eq:Chow-Polylog-AFE} coincides with identity \eqref{Eq:Efa-(Ab)-symmetric} satisfied by Bloch-Wigner function $\mathcal L_2=\mathcal D$.   It has been a question regarding grassmannian (and even Chow) polylogarithms to know whether the general value $\mathcal L_n^G(\boldsymbol{x})$ for a generic configuration $\boldsymbol{x}\in {\rm Conf}_{2n}(\mathbf P^{n-1})$ can be expressed as a finite sum $\sum_i \mathcal L_n(r_i(\boldsymbol{x}))$ for some rational functions $r_i$.   It does not seem to be the case already for $n=3$  but in \cite{GoncharovZhao}, 
 Goncharov and Zhao express $\mathcal L_3^G$ in terms of a {\it motivic grassmannian trilogarithm} $L_3^G$ which, in turn, can be expressed by means of $\mathcal L_3$ since $L_3^G={\rm Alt}_6(\mathcal L_3(r_3))$, that is 
 $$L_3^G(\boldsymbol{x})=\sum_{\sigma\in \mathfrak S_6}
 {{\rm sgn}(\sigma)}\, 
 \mathcal L_3\Big(r_3\big(x_{\sigma (1)}, \ldots,x_{\sigma (6)}\big)\Big)
 $$ 
 for a generic  configuration $\boldsymbol{x}=[x_1,\ldots,x_6]\in {\rm Conf}_6(\mathbf P^{2})$, where $r_3: {\rm Conf}_6(\mathbf P^2)\dashrightarrow \mathbf P^1$ stands for the triple-ratio (defined in \eqref{Eq:r1-and-r3} above).  Moreover, the trilogarithm $L_3^G$ satisfies the same  identity 
as  $\mathcal L_3^G$ 
on ${\rm Conf}_7(\mathbf P^2)$, namely 
$\sum_{i=1}^7 (-1)^i      L_3^G(f_i)=0$ ({\it cf.}\,Theorem 3.11 in \cite{Goncharov1994}),   which,  expressed in terms of 
Zagier's trilogarithm 
$\mathcal L_3$, takes the following nice (anti)symmetric  form:
\begin{equation}
 \label{Eq:Eq-L3-r3*ui}
\sum_{i=1}^7 \sum_{\sigma \in \mathfrak S_6} 
(-1)^i \,{{\rm sgn}(\sigma)}\, \mathcal L_3\Big( \,r_3^\sigma\circ f_i  \, \Big)
=0\, .
\end{equation}
\big(We recall that for any $\sigma\in \mathfrak S_6$, 
one has  $r_3^\sigma(\boldsymbol{x})=r_3(x_{\sigma(1)},\ldots,x_{\sigma(6)})$ for $\boldsymbol{x}=(x_i)_{i=1}^6\in {\rm Conf}_6(\mathbf P^2)$.\big)
\mk 

  The previous identity leads us to consider the 
web on ${\rm Conf}_7(\mathbf P^2)$, denoted by 
$\boldsymbol{{\mathcal W} \hspace{-0.05cm} L}_3^G$, 
 defined by the rational functions 
 $r_3^\sigma\circ f_i:  {\rm Conf}_7(\mathbf P^2)\rightarrow {\rm Conf}_6(\mathbf P^2)\dashrightarrow \mathbf P^1$:  
 \begin{equation}
 \label{Eq:Web-r3*ui}
 \boldsymbol{{\mathcal W}\hspace{-0.05cm} L}_3^G = 
  \boldsymbol{{\mathcal W}}\Big(  \, r_3^\sigma\circ f_i\hspace{0.2cm} \big\lvert \hspace{0.2cm}
 i=1,\ldots,7, \, 
 \, \sigma \in \mathfrak S_6\, 
 \Big)\, .
 \end{equation}
 
One verifies that, thanks to some invariance properties of the triple-ratio,  
some of the $r_3^\sigma$'s with $\sigma\in \mathfrak S_6$ 
coincide (as rational functions on ${\rm Conf}_6(\mathbf P^2)$) and that  the set of such functions is of cardinality 120. Consequently, there are only 840 trilogarithmic terms in \eqref{Eq:Eq-L3-r3*ui} ({\it cf.}\,Remark p.\,66 of  \cite{Goncharov1994})  and one verifies that $\boldsymbol{{\mathcal W} \hspace{-0.05cm} L}_3^G$ is indeed a 840-web in 6 variables. 
 After some lengthy computations on a computer, one gets
$$ 
\rho^\bullet\Big(\boldsymbol{{\mathcal W} \hspace{-0.05cm} L}_3^G
\Big)=
\big( \, 834\, , \,  819\, , \, 784\, , \, 714\, , \,588\, , \,399\, , \, 168\, \big) \qquad 
\mbox{ and } \qquad {\rm polrk}^1\Big(\boldsymbol{{\mathcal W} \hspace{-0.05cm} L}_3^G
\Big)= 1547\,.
$$
We have been unable to obtain the polylogarithmic rank of higher weight
but considering the previous numbers leads us to doubt that 
$\boldsymbol{{\mathcal W} \hspace{-0.05cm} L}_3^G$
is AMP.
\begin{center}
\vspace{-0.2cm}
$\star$
\end{center}

%

The general case has been considered in \cite{CGR} where the authors study the $n$-th grassmannian  polylogarithm 
$\mathcal L_n^G$
for any $n\geq 2$: 
they prove that it can be expressed as a linear combination of weight 4 iterated integrals (in one or two variables) and give a formula for the corresponding symbol ({\it cf.}\,their Theorem 6).  Then, specializing 
their formula 
to the case when $n=4$ (see \cite[Theorem 6]{CGR}),  they 
are able to find an explicit element $q$ of $\mathbf Z[K_8(3)^{**}]$ (where $
K_8(3)$ stands for the field of rational functions on ${\rm Conf}_8(\mathbf P^3)$) such that the functional equation 
\begin{equation}
 \label{Eq:Eq-L4-q*ui}
\sum_{i=1}^9 (-1)^{i}\mathcal L_4\Big( {\rm Alt}_8(q)\circ f_i\Big)=0
\end{equation}
is satisfied on ${\rm Conf}_9(\mathbf P^3)$. They see this  abelian functional equation in 12 variables  as the tetralogarithmic analogue of the 5-terms relation for the dilogarithm and of Goncharov's 
840-terms relation \eqref{Eq:Eq-L3-r3*ui} for the trilogarithm. It would be interesting to know more about the web associated to \eqref{Eq:Eq-L4-q*ui}.

%

\subsubsection{Radchenko's polylogarithmic functional equations and associated webs.}
In his PhD dissertation \cite{Radchenko}, Radchenko elaborates a new approach to construct functional equations for higher polylogarithms.  The main new conceptual ingredient he uses is the notion of {\it `$S$-cross-ratio`'}. 
According to Radchenko (we refer to his dissertation for more details), this notion is a  
 {\it `generalization of the classical cross-ratio which, roughly speaking,   is a function of $n$ points in a projective
space that is invariant under the change of coordinates and satisfies an arithmetic
condition somewhat similar to the Pl\"ucker relation'}. In particular, the usual cross-ration ${\rm Cr}$ on ${\rm Conf}_4(\mathbf P^1)$ as well as Goncharov's triple ratio $r_3$ on ${\rm Conf}_6(\mathbf P^2)$ are just particular cases of the generalized cross-ratios which can be obtained following Radchenko's approach. \mk 
 
 After describing the generalized cross-ratios on ${\rm Conf}_m(\mathbf P^k)$ when $m$ and $k$ are small enough (see \cite[Appendix A]{Radchenko}), he studies the AFE of the form $\sum_{i=1}^N c_i \mathcal L_n(r_i)$ when the $r_i$'s are generalized cross-ratios and the $c_i$'s are rational numbers. 
 In the  final  chapter of his dissertation, Radchenko gives  an important number of such new AFEs  obtained using this approach, in weight 3, 4 and 5 
 (see the tables 5.1 and 5.2, and Appendix B in \cite{Radchenko}).
 \mk

We briefly discuss below a few of the new polylogarithmic AFEs obtained by Radchenko from the point of view of web geometry. We have very little to say about these webs: there are mentioned here to advertise future researches. 
The notation used for the $r_i$'s in the lines below are those of \cite[Appendix A]{Radchenko}. 
\begin{itemize}
\item[$\bullet$] The  following functional equation  
$$
 \mathcal L_3\Big({\rm Sym}_5\big([r_5]-2[r_3]\big)
 \Big)=
\sum_{\sigma\in \mathfrak S_5} \mathcal L_3\big(r_5^\sigma\big)-2\, 
\mathcal L_3\big(r_3^\sigma\big)=20 \, \zeta(3)
$$  
holds true on  $\mathcal M_{0,5}$
 and provides a trilogarithmic AR for 
$
\boldsymbol{\mathcal W}_{{\rm Sym}_5([r_5]-2[r_3])}
=
\boldsymbol{\mathcal W}\big(  
r_3^\sigma\, , \, r_5^\sigma
\, \lvert 
\, \sigma\in \mathfrak S_5\, 
\big)
$.  The latter  is a planar 45-web with 
${\rm polrk}^\bullet(\boldsymbol{\mathcal W}_{{\rm Sym}_5([r_5]-2[r_3])}
)=(85, 86, 36, \ldots)$. 
\sk
\item[$\bullet$]  Looking at Tables  5.1 and 5.2 in \cite{Radchenko}, 
we get that the following identities hold true
\begin{align}
\label{Eq:Radchenko-r3r4r7}
 \mathcal L_3\Big(\,  {\rm Sym}_5\big([r_7]-3[r_4]-5[r_3] \, \Big)
 = & -40 \, \zeta(3) \\
 \mbox{ and } \qquad 
 \mathcal L_3\Big( \, {\rm Alt}_5\big([r_7]+3[r_4]-3[r_3]\, \Big)=& \,\,  0 \nonumber 
 \end{align}
and provide trilogarithmic ARs for  the web  
$
{\boldsymbol{\mathcal W}}\big( 
\, r_3^\sigma\, , \, r_4^\sigma 
\, , \, r_7^\sigma \, \lvert \, \sigma\in \mathfrak S_5\, 
\big) 
$ on  $\mathcal M_{0,5}$.  After cleaning, we obtain a 135-web in two variables with (harmonic) polylogarithmic rank  
%
${\rm polrk}^\bullet=(245, 276, \ldots)$. That the second polylogarithmic rank be higher than the first is a bit surprising. It would be interesting to understand this phenomenon better.
\sk
\item[$\bullet$] 
From our perspective, it might be interesting to consider the {\it `complete Radchenko's web'} on $\mathcal M_{0,5}$, namely the web 
$
\boldsymbol{\mathcal W}\big(  \, 
r_i^\sigma
\,  \, \lvert \,  \, 
i=1,\ldots,7\, , \, 
\sigma\in \mathfrak S_5\, 
\big)
$
defined by all the generalized cross-ratios on $\mathcal M_{0,5}$, up to permutations: this web clearly carries polylogarithmic ARs up to weight 3 at least. 
It can be verified that it is a 170-web in two variables. It would be interesting to describe all its ARs and to know its rank.
\sk
\item[$\bullet$]  
The symmetric  identity  on $\mathcal M_{0,6}$ $$
\mathcal L_3\Big(\, {\rm Sym}_6\big(3[r_{11}]-4[r_8]\big)\, \Big)= 120\, \zeta(3)
$$ 
 gives rise to a trilogarithmic AR for   
$\boldsymbol{\mathcal W}\big(  \, r_8^\sigma
\, ,  \, r_{11}^\sigma \,  \lvert \,  \, 
\sigma\in \mathfrak S_6\, 
\big) $. After cleaning, one sees that the latter is a 105-web in three variables 
such that $\rho^\bullet=(102, 99, 95, 90, 84, 77, 69, 60, 50, 39, 27,$ $ 14)$ and 
$ {\rm polrk}^\bullet=(186,170,\ldots)$. We remark that   $\rho^{\sigma-1}=\rho^{\sigma}+\sigma+1$ for $\sigma=2,\ldots,12$: this nice behavior for its  
virtual ranks  might indicate that this web  could have other interesting features. 
\item[$\bullet$]  Looking at 
 \cite[Table\,5.2]{Radchenko},  
we get that the following identity holds true on $\mathcal M_{0,6}$ 
\begin{align*}
 \mathcal L_5\Big({\rm Alt}_6\big([r_7]+9[r_4]-15[r_3]\big)\Big)=0
\end{align*}
and provides a pentalogarithmic AR for the web  
${\boldsymbol{\mathcal W}}_{r_3,r_4,r_7}={\boldsymbol{\mathcal W}}\big( 
\, r_3^\sigma\, , \, r_4^\sigma 
\, , \, r_7^\sigma \, \lvert  \, \sigma\in \mathfrak S_6\, 
\big)$ which can be verified to be a 810-web in three variables. 
Since for any $i=1,\ldots,6$, it admits,  as a subweb,  the pull-back under the forgetful map $\varphi_i$ of the 135-web on $\mathcal M_{0,5}$ associated to the functional equations \eqref{Eq:Radchenko-r3r4r7}, we see that this web carries many trilogarithmic ARs as well. 
Some computations give us 
 $
{\rho}^\bullet\Big({\boldsymbol{\mathcal W}}_{r_3,r_4,r_7}\Big)=(807, 804, 800, 795, 789, 783, 777, 771,$ $ 765, 759, 753, 747,\ldots)$ 
hence it seems that $\rho^{\sigma+1}=\max( 0, \rho^\sigma-6)$ starting from $\sigma=5$. These  nice relations between its virtual ranks and the fact 
 that it carries many polylogarithmic ARs of weight $\leq 5$ make 
 this web particularly interesting. 
\end{itemize}

\newpage
\section{Cluster algebras and associated webs}
\label{S:Cluster}

In this section, we discuss some interesting webs constructed from cluster algebras.  We will first of all consider webs associated to cluster algebras of finite type but in a second step, we will have to work with more general cluster algebras. \mk

We first give a review of the theory of cluster algebras before 
introducing the webs associated to them that we will study afterwards. 
We will essentially consider two kinds of cluster webs: the {\it $\boldsymbol{\mathcal X}$-cluster webs'}  that  are the webs defined by all the ${\mathcal X}$-cluster variables of a given cluster algebra of finite type, and the 
{\it $\boldsymbol{\mathcal Y}$-cluster webs} which are the webs  associated to a {\it `period'} of a non-necessarily finite cluster algebra. \sk

Our web-theoretic motivations will  make us ask some questions about cluster variables which already make senss and are interesting within the theory of cluster algebras (without any references to webs) and we will mention some conjectural answers to them.  \mk 

As a (already interesting) warm-up, we will consider the webs associated 
to finite type cluster algebras of rank 2. Then, we will turn  to the cases of $\boldsymbol{\mathcal X}$-cluster webs or $\boldsymbol{\mathcal Y}$-cluster webs associated to a finite type cluster algebra. We will discuss many of these webs and state several conjectures about them (more precisely about their (virtual or standard) rank(s) and their abelian relations). 
\begin{center}
$\star$
\end{center}

We denote by $[\cdot ]_+ : \mathbf R\rightarrow \mathbf R_{\geq 0}$ the function defined by 
$[ x ]_+=\max (0, x)$ 
for any  real  $x$. 

 \subsection{\bf Cluster algebras}
\label{S:Cluster-Algebra}
The notion of cluster algebra has been introduced by Fomin and Zelevinsky at the beginning of the 2000's and has been recognized since as a very fertile one.  Indeed, 
cluster algebras have been intensively studied from several points of views: it appears now that these objects are related to many other subjects of mathematics and mathematical physics, among which one can mention the following ones: 
 Poisson geometry, positivity and dual canonical bases for (quantum) algebraic groups, higher Teichm\"uller theory, Lie Theory and representation theory of quivers, discrete (integrable)  dynamical systems, (quantum and classical) dilogarithmic identities, 
algebraic geometry (more precisely geometry of log-Calabi-Yau varieties and  
Donaldson-Thomas invariants), mathematical physics (BPS states, scattering amplitudes).\footnote{To get a glimpse of the vivacity of the subject these days, the reader can also consult the \href{http://www.math.lsa.umich.edu/~fomin/cluster.html}{\it `Cluster algebras portal'}.} 
\mk

There are nowadays several introductions to or survey papers about cluster algebras  to which we refer the reader for more details and background: to name a few, we mention  the seminal papers by Fomin and Zelevinsky ({\it e.g.}\,\cite{FZ,FZII,FZIV} cited in the bibliography), the books \cite{GSV}, \cite{FWZ} or \cite{Marsh}, the short(er) introductions 
\cite{Williams,GR}, and the survey papers \cite{FominC,KellerBourbaki}.  Among the research 
papers   on cluster algebras   which are important regarding the  material to come, one can mention  \cite{FockGoncharovENS} by Fock and Goncharov, the papers by Keller cited in the bibliography (especially \cite{KellerCTQDI} and 
\cite{KellerAnnals}) and  those by Nakanishi such as \cite{NakanishiNagoya,Nakanishi,NakanishiTrop}. 
\mk 


The theory of cluster algebras has developed in a very dynamic way over the last twenty years and 
 several deep and important results in the field have been obtained, some quite recently.  Although we will state some of the most important results in concise form, we will not discuss them that much since we will essentially use them as black boxes for our purpose, which can be summarized as that of `{\it getting new webs defined by 
 cluster variables (as first integrals) which are interesting from the point of view of web geometry since carrying many (mostly polylogarithmic) abelian relations'.} 
\begin{center}
$\star$
\end{center}



Before going more seriously into the details, let us say a few general words about the kind of cluster algebras we will work with in the sequel. 
There are several versions,  more or less general,  of cluster algebras. 
Those we will consider in this text all are cluster algebras without frozen variables.  Another ingredient used is a certain `semi-field of coefficients' $\mathbb P$. 
The simplest case is when  this semi-field is trivial, that is when $\mathbb P={\bf 1}$.  It is the one we will consider almost everywhere. 
 Another important case encountered in the literature is the one of `principal coefficient', when $\mathbb P$ is the tropical semi-field on $n$ generators, where $n$ stands for the rank of the considered cluster algebra. Considering our purpose, this case is important from a theoretical point of view, regarding  the key notion of `$F$-polynomial' 
of the theory of cluster algebras. It will be briefly discussed  in \S\ref{SubPar:Fpoly-finite-type}.  
Everywhere else, the underlying semi-fields of coefficients of the cluster algebras we will consider will be implicitly assumed to be trivial.


\subsubsection{Basics on cluster algebras}
\label{SS:Basics-Cluster-Algebras}
We will mainly consider cluster algebras over the field of rational or complex numbers. Accordingly, we will work within a field $F=\mathbf k(z_1,\ldots,z_n)$  of rational functions in $n$ indeterminates  with coefficients in $\mathbf k=\mathbf Q$ or $\mathbf C$, for some integer $n\geq 2$. \mk 

\paragraph{Skew-symmetrizable matrices, quivers and mutations.}
Let $B=(b_{ij})_{i,j=1}^n$ be a square $n\times n$ matrix with integers coefficients $b_{ij}$. It is {\bf skew-symmetrizable} if there exists a diagonal matrix $D={\rm diag}(d_1,\ldots,d_n)$ with diagonal entries  $d_i\in \mathbf Z_{>0}$ such that $DB=(d_ib_{ij})_{i,j=1}^n$  be skew-symmetric. In this case, such a  matrix $D$ is a {\bf skew-symmetrizer} of $B$.
\mk 


By definition, for $k\in \{1,\ldots,n\}$, the matrix {\bf mutation}  
$\mu_k(B)$
of $B$ in the direction $k$ is the matrix 
$B'=(b_{ij}')_{i,j=1}^n$ whose coefficients $b_{ij}'$ are given by the following formulas: 
$$
b_{ij}'=\begin{cases}
\, - b_{ij} \hspace{1.95cm} \mbox{ if }\,  \,  k \in \{i,j\} ;  \\
\, \hspace{0.2cm} b_{ij}\hspace{2cm} \mbox{ if } \,  k\not \in \{i,j\}  \, \mbox{  and   } \,   b_{ik}b_{kj}\leq 0;\\
\, \hspace{0.2cm} b_{ij}+\lvert b_{ik}\lvert \, b_{kj}\hspace{0.5cm} \mbox{ if } \,  k\not \in \{i,j\} \, \mbox{ and  }\,   b_{ik}b_{kj}> 0.
\end{cases}
$$
It can be verified that: (1) $DB'$ is antisymmetric hence $B'=\mu_k(B)$ is again skew-symmetrizable; (2) this operation is involutive, that is $\mu_k(\mu_k(B))=B$.  The point is that the different mutations $\mu_i$'s do not commute all together and this is the source of all the complexity and all the  richness of the cluster algebras. 
\mk

Now let $Q$ be a finite quiver with $n$ vertices, that is an oriented graph together with two maps, the source $s$ and the target $t$, defined on the set of oriented edges (or {\it `arrows'}) $Q_1$ of $Q$, with values in its set of vertices $Q_0$.  We will always assume that there is no loop nor 2-cycle in $Q$.\sk

A {\bf valued quiver} (see \cite[\S2.1.2]{Dupont2008} for more details) is
a quiver $Q$ endowed with a function $v : Q_1\rightarrow  \mathbf N^2$  such that $(i)$ there is at most one arrow between any two vertices of $Q$; and $(ii)$ there is a function $d : Q_0 \rightarrow  \mathbf N_{>0}$ such that  for each arrow
$a:   i \longrightarrow j$  (with $i,j\in Q_0$) one has 
 $ d(i)\, v(a)_1 = d(j)\,v(a)_2$ where $v(a) = (v(a)_1, v(a)_2)$.
\sk

On the graph of a valued quiver we
write the weights on the top or bottom of the arrows. Note that 
$ i\stackrel{ 
 {(v_1,v_2)}    }{\longrightarrow} j$ 
is the same weighted
arrow as $ j\stackrel{(v_2,v_1)}{\longleftarrow} i$. 
In the case when $v_1=v_2=w$, we replace the weighted arrow  by $ i\stackrel{w}{\longrightarrow} j$ or by $w$ arrows from $i$ to $j$. We also adopt the following convention: 
\begin{equation}
\label{Eq:Convention}
\begin{tabular}{l}
\begin{tabular}{l}
$\bullet $
we  allow edge $e\in Q_1$ with $v(e)=(0,0)$ with the convention that   
$ i\stackrel{(0,0)}{\longrightarrow} j$ 
 (or equiva\- \\  
 ${}^{}$\hspace{0.25cm} -lently  
 $ i\stackrel{0}{\longrightarrow} j$) means that there is no  arrow from  $i$ to $j$;  and\vspace{0.4cm}
\end{tabular}
\\
\begin{tabular}{l}
$\bullet $
 $ i\stackrel{(-v_1,-v_2)}{\longrightarrow} j$  with both $v_1$ and $v_2$ non-negative 
stands for $ i\stackrel{(v_1,v_2)}{\longleftarrow} j$.
\end{tabular}
\end{tabular}\sk
\end{equation}

Now to  a valued quiver $Q$ with vertex set $Q_0=\{1,\ldots,n\}$, one associates an integer matrix $B = B_Q=(b_{ij})_{i,j\in I}$ as follows: $
b_{ij}=0$ if there is no arrow between $i$ and $j$; $b_{ij}=v(a)_1$ if there is an arrow $a:i\longrightarrow j$; and $b_{ij}=-v(\alpha)_2$ if  there is an arrow $\alpha$ from $j$ to $i$. Then denoting by $D$ the $I$-diagonal matrix with diagonal entries $d(i)$ with $i\in I$, one verifies that $DB$ is 
antisymmetric thus the matrix $B$ is skew-symmetrizable. 
Conversely, to a skew-symmetrizable matrix $B$,  one associates the valued quiver $Q=Q(B)$ with set of vertices $Q_0=\{1,\ldots,n\}$ and valued arrow $ a_{ij}:i{\longrightarrow} j$ for any distinct $i,j\in Q_0$, 
with $v(a_{ij})= (\lvert b_{ij}\lvert, \lvert b_{ji}\lvert)$.
One obtains this way 
a bijection 
\begin{equation}
\label{Eq:Bijection-B<->Q}
B=B_Q 
\qquad \xleftrightarrow{\hspace*{1cm}}
\qquad 
Q=Q_B
\end{equation}
 between the set of skew-symmetrizable $I$-square matrices with integer entries and that of valued quivers with vertex set $I$. Using this bijection, 
one defines the mutation $\mu_k(Q)$ in the direction $k\in I$ of a valued quiver $Q$ via the formula 
$$\mu_k(Q)=Q_{\mu_k(B_Q)}.$$
(There are well-known rules 
 allowing to  
 construct $\mu_k(Q)$ directly from $Q$, without using the bijection \eqref{Eq:Bijection-B<->Q}, see \cite[Prop.\,8.1]{FZII} for instance). 
 \mk 
 
 Here is a skew-symmetrizable example (of Dynkin type $B_3$, see \S\ref{Par:Cluster-Algebras-Finite-Type} further): 
 \begin{equation*}
B_{B_3}=\begin{bmatrix}
0 & 1 & 0 \\
-1 & 0 & -2 \\
0 & 1 & 0
\end{bmatrix}
\qquad \xleftrightarrow{\hspace*{1cm}}
\qquad 
Q_{B_3}\, : 
\hspace{0.3cm}
  {\bf 1}  \hspace{-0.1cm}
\begin{tabular}{c} 
$ \stackrel{1}{
\xrightarrow{\hspace*{1cm}}
}$ \vspace{-0.15cm}\\  ${}^{}$
 \end{tabular} 
 \hspace{-0.2cm} {\bf 2} 
 \hspace{-0.1cm}
\begin{tabular}{c} 
$ \stackrel{(1,2)}{
\xleftarrow{\hspace*{1.6cm}}
}$ \vspace{-0.15cm}\\  ${}^{}$
 \end{tabular} 
\hspace{-0.2cm} {\bf 3}\, .
\end{equation*}



\paragraph{Seeds, mutations and cluster algebras}
\label{Par:Seeds}

We now describe a central notion of the theory of cluster algebras, that of seeds. There are two parallel notions of seeds, the $\mathcal A$-seeds, and the $\mathcal X$-seeds. The  latter is the most relevant for applications to web theory, but since the former is a fundamental one, we introduce it as well.
\mk

A {\bf (labeled) seed} is a triple $S=(\boldsymbol{a},\boldsymbol{x},B)$ where $\boldsymbol{a}=(a_1,\ldots,a_n)$ is a $n$-tuple of indeterminates, 
$\boldsymbol{x}=(x_1,\ldots,x_n)$ is another one (independant of the former) and  $B = (b_{ij})_{i,j=1}^n$ is a skew-symmetrizable $n \times n$  integer matrix, called the {\bf exchange matrix} of the seed.  
The associated 
{\bf (labeled) $\boldsymbol{\mathcal A}$-seed}  (resp.\,{\bf $\boldsymbol{\mathcal X}$-seed})  is the pair $(\boldsymbol{a},B)$ (resp.\,$(\boldsymbol{x},B)$).  The $n$-tuples $\boldsymbol{a}$ and $\boldsymbol{x}$ are called {\bf `$\boldsymbol{\mathcal A}-$} (resp.\,{\bf $\boldsymbol{\mathcal X}-$}) {\bf clusters'} and their elements  $a_i$ and $x_i$ 
 {\bf `$\boldsymbol{\mathcal A}-$} (resp.\,{\bf $\boldsymbol{\mathcal X}-$}) {\bf cluster variables'}.  Note that all these objects are a priori associated to the considered seed. \mk 
 
 By definition, for $k\in \{1,\ldots,n\}$, the {\bf (cluster) mutation} in the $k$-th direction  of the seed $S$ is the new seed $$S'=\big(\boldsymbol{a}',\boldsymbol{x}',B'\big)=\Big(\mu_k\big(\boldsymbol{a},B\big), \mu_k\big(\boldsymbol{x},B\big),\mu_k(B)\Big)=\mu_k\big(S\big)$$ 
  with $B'=\mu_k(B)$ being the matrix mutation of $B$ in direction $k$ 
  and where the new clusters $\boldsymbol{a}'=(a_1',\ldots,a_n')=\mu_k(\boldsymbol{a},B)$ and $\boldsymbol{x}'=(x_1',\ldots,x_n')=\mu_k(\boldsymbol{x},B)$ are defined by the following formulas for the corresponding cluster variables: 
\begin{align}
\label{Eq:A-X-Mutation-formulae}
a_{j}'= & \, a_j  \hspace{-0.5cm} &\mbox{if }\, j\neq k   \qquad \quad    &\mbox{ and }\quad \qquad
a_k'=  \frac{1}{a_k}\left( \,   \prod_{
b_{\ell k}>0} a_\ell^{b_{\ell k}}+
 \prod_{
   b_{lk }<0} a_l^{-b_{l k  }}
\, \right) \hspace{0.4cm}  \mbox{if }\, j= k \sk \\
x_{j}'= & \, {x_j}^{-1}  \hspace{-0.5cm} &\mbox{ if }\, j= k   \qquad \quad    &\mbox{ and }\quad \qquad
x_j'=  x_j\,\left(1+x_k^{[- b_{kj}]_+}\right)^{-b_{kj}}
\hspace{1.1cm}  \quad  \mbox{if }\, j\neq  k\, .\nonumber
\end{align}

Again, one can verify that this mutation is an involution: one has $\mu_k(\mu_k(S))=S$ as a seed from which it comes that the 
the {\bf `$\boldsymbol{\mathcal A}-$  mutation'} 
in the $k$-th direction  $\boldsymbol{a}\mapsto \boldsymbol{a}'=\mu_k(\boldsymbol{a},B)$ as well as the associated {\bf `$\boldsymbol{\mathcal X}-$  mutation'} $\boldsymbol{x}\mapsto \boldsymbol{x}'=\mu_k(\boldsymbol{x},B)$ both are birational maps (whose inverse respectively are  
$\boldsymbol{a}'\mapsto \boldsymbol{a}=\mu_k(\boldsymbol{a}',B')$  
and $\boldsymbol{x}'\mapsto \boldsymbol{x}=\mu_k(\boldsymbol{x}',B')$,  
with $B'=\mu_k(B)$). 
\mk

Given a seed $S=(\boldsymbol{a},\boldsymbol{x},B)$ as above and  a permutation
$\sigma\in \mathfrak S_n$, one sets $S^\sigma=(\boldsymbol{a}^\sigma,\boldsymbol{x}^\sigma,B^\sigma)$ with 
\begin{equation}
\label{Eq:Equiv-of-Seeds}
\boldsymbol{a}^\sigma=\big(\, a_{\sigma(i)}\, \big)_{i=1}^n\, , \qquad 
\boldsymbol{x}^\sigma=\big(\, x_{\sigma(i)}\, \big)_{i=1}^n
\qquad \mbox{ and }\qquad B^\sigma=\Big(\, b_{\sigma(i)\sigma(j)}\, \Big)_{i,j=1}^n\, .
\end{equation}
  Then another seed $S'=(\boldsymbol{a}',\boldsymbol{x}',B')$ is said to be {\bf equivalent to $S$} if is of the form $S^\sigma$ for a certain $\sigma$.  The notion of equivalence for two 
$\boldsymbol{\mathcal A}-$ or $\boldsymbol{\mathcal X}-$seeds  can be stated similarly in an obvious way. We will write $S\simeq S'$ whenever  these two seeds are equivalent and $[S]$ will stand for the set of $\boldsymbol{\mathcal S}-$seeds equivalent to $S$. 
\mk  


Now let $\mathbb T_n$ stand for the $n$-regular tree, that is the connected acyclic simple $n$-regular graph with a labeling by integers from 1 to $n$ of the edges of $\mathbb T_n$, requiring moreover that the $n$ edges emanating from a vertex $t$ of $\mathbb T_n$ receive pairwise different labels.  By definition, a {\bf cluster pattern} $\boldsymbol{\mathcal C\mathcal P}$ 
 (possibly of {\bf $\boldsymbol{\mathcal S}-$seeds} with 
$\boldsymbol{\mathcal S}=\boldsymbol{\mathcal A}$ or $\boldsymbol{\mathcal S}=\boldsymbol{\mathcal X}$) is an assignment of a $\boldsymbol{\mathcal S}-$seed $S_t$ to any vertex $t$  of $\mathbb T_n$ such that for any edge, say between $t$ and $t'$ and labeled by $k$,  $S_{t'}$ is obtained from $S_t$ by means of the mutation of the latter seed in the direction $k$ (and equivalently for $S_t$ from $S_{t'}$), {\it i.e.}: 
$$
 t
 \hspace{-0.1cm}
\begin{tabular}{c} 
$ \stackrel{k}{\line(1,0){25}}$ \vspace{-0.35cm}\\  ${}^{}$
 \end{tabular} 
\hspace{-0.2cm}
 t'
 \qquad 
 \Longrightarrow \qquad S_{t'}=\mu_k\Big(S_t\Big) \hspace{0.2cm} \Longleftrightarrow 
 \hspace{0.2cm} S_{t}=\mu_k\Big(S_{t'}\Big) \, .
$$
We define the {\bf $\boldsymbol{\mathcal S}-$exchange graph} of a $\boldsymbol{\mathcal S}-$cluster pattern as the quotient of ${\mathbb T}_n$  by the relation $t\simeq t'$ for two vertices if the two corresponding seeds $S_t$ and $S_{t'}$ are equivalent $\boldsymbol{\mathcal S}-$seeds.  This exchange graph can be finite or infinite. 
\mk

Given a cluster pattern, for any vertex $t$ of $\mathbb T_n$,  we denote by $a_i^t$ (resp.\,$x_i^t$) the corresponding $\boldsymbol{\mathcal A}-$ (resp.\,$\boldsymbol{\mathcal X}-$) cluster variables and by 
$b_{ij}^t$ the coefficients of the exchange matrix $B^t$ of the seed $S_t=
(\boldsymbol{a}^t,\boldsymbol{x}^t,B^t)$. When a vertex $t_0$ is chosen, one refers to the associated seed $S_{t_0}$ as the {\bf initial seed}, and to the corresponding cluster variables $a_i^{t_0}$ (or $x_i^{t_0}$) as the {\it `initial cluster variables'}, and to $B^{t_0}$ as the {\it `initial exchange matrix'}. Note that an initial seed $S_{t_0}$ entirely determines the full cluster pattern.  And since this 
seed is completely determined by $B^{t_0}$, actually everything rely on the initial matrix.
To simplify the notation, we will often drop the suscript $t_0$ for the initial objects, that is the `{\it cluster objects}' associated to $S_{t_0}$. 
\mk 

Assume that $\boldsymbol{\mathcal S}$ is empty or is $\boldsymbol{\mathcal A}$.  
The cluster algebra $\boldsymbol{A}=\boldsymbol{A}(\boldsymbol{\mathcal C\mathcal P})$ (also denoted by $\boldsymbol{A}(S_{t_0})$ or $\boldsymbol{A}(B^{t_0})$ assuming that an initial seed has been chosen) associated to a $\boldsymbol{\mathcal S}$-cluster pattern $\boldsymbol{\mathcal C\mathcal P}$ is the $\boldsymbol{\bf k}$-subalgebra of $F=\boldsymbol{\bf k}(a_1,\ldots,a_n)$ generated by all the 
$\boldsymbol{\mathcal A}$-cluster variables $a_i^t$ of all the cluster $(\boldsymbol{a}^t,B^t)$ of the considered cluster pattern.  In other terms, 
 given an {initial seed} $(\boldsymbol{a}, B)$, $\boldsymbol{A}(B)$ is   the 
subalgebra of  $F$ generated by all the $\boldsymbol{\mathcal A}$-cluster variables 
which can be constructed from $(\boldsymbol{a}, B)$
 by means of a finite but arbitrary number of seed mutations. 
 \mk

At some points, it will be useful to use the following notation:  given an initial seed $(\boldsymbol{a},\boldsymbol{x},B)=((a_i)_{i=1}^n,(x_i)_{i=1}^n,B)$, one sets
\begin{align*}
\boldsymbol{\mathcal A var}\big( B\big)= 
\boldsymbol{\mathcal A var}\big( \boldsymbol{A}(B)\big)=
& \, \Big\{  \hspace{0.1cm} 
  \boldsymbol{\mathcal A}-\mbox{cluster variables of  }\boldsymbol{A}(B)\, 
\Big\} 
\subset {\bf k}(a_1,\ldots,a_n)
\\
\mbox{ and }\quad 
\boldsymbol{\mathcal X var}\big( B\big)= 
\boldsymbol{\mathcal X var}\big( \boldsymbol{A}(B)\big)= 
& \, \Big\{  \hspace{0.1cm}
 \boldsymbol{\mathcal X}- \mbox{cluster variables of }
\boldsymbol{A}(B)\, 
\Big\} 
\subset {\bf k}(x_1,\ldots,x_n)\, . 
\end{align*}
A priori, both sets are infinite sets of rational functions in the initial variables $a_i$'s and $x_i$'s respectively.
\bk

\label{Pageref:Warning-Apology}
{\bf A warning and an apology:} {\it there are  two distinct conventions for the terminology  in the theory of cluster algebra, one by Fomin and Zelevinsky, the other by Fock and Goncharov:  the $\boldsymbol{\mathcal A}$- and 
$\boldsymbol{\mathcal X}$-cluster variables of the two latter authors are respectively called 
$X$- and 
$Y$-cluster variables by the former two. Moreover, the mutations  of a seed $(\boldsymbol{a},\boldsymbol{x},B)$ considered by Fock and Goncharov correspond to those of Fomin and Zelevinsky, but with the transpose of $B$  instead of $B$. }\sk

{\it Here,  due to a clumsiness that occurred at an early stage in the development of this work, we have accomplished the feat of using the `$\boldsymbol{\mathcal A}/\boldsymbol{\mathcal X}$-cluster terminology' of Fock and Goncharov while working with the formula of Fomin and Zelevinsky for the mutations.    If this doesn't cause any error if working carefully, this can make it uneasy to connect with the papers of the authors mentioned above, 
which are basic references on cluster algebras. We apologize to the reader for any inconvenience this may cause.}
\bk

In this text, we will consider several examples of cluster algebras constructed from Dynkin diagrams, to which we refer in the sequel. 
As for another  classical and interesting type  of cluster algebras, we mention  the ones associated to triangulated surfaces 
introduced by Fomin, Shapiro and D. Thurston in \cite{FominShapiroThurston}. 

\begin{exm}
Let $(S,P)$ be a {\it `bordered surface with marked points'}: $S$ is an oriented surface with (possibly empty)  boundary $\partial S$ and 
$P$ is a finite set of points in $S$ which has non-empty intersection with any connected component of $\partial S$.  
Let $\mathcal T=(\tau_1,\ldots,\tau_n)$ be a set of arcs in $S\setminus \partial S$, with endpoints in $P$ which together with the boundary segments (namely, the closure of the connected components of $\partial S\setminus P$) form a triangulation 
$\overline{\mathcal T}$
 of $S$.\footnote{Since we will not actually use this construction in this text, we only describe it in a simple case and refer to \cite{FominShapiroThurston} for the general construction.} One associates to such a $\mathcal T$ a $n\times n$ `{\it exchange matrix}'  $B({\mathcal T})$   whose 
$(i,j)$-th coefficient is given by the difference 
between the number of triangles in $\overline{\mathcal T}$ with sides 
$\tau_i$   and 
 $\tau_j$  
with $\tau_j$ following $\tau_i$ in clockwise order, minus the number of triangles with sides $\tau_i$ and $\tau_j$, but this time when $\tau_j$ follows 
$\tau_i$ in counterclockwise order.\sk

For any $k\in \{1,\ldots,n\}$, one sets $\tau_i=\tau_i$ for $i\neq k$ 
and one defines  $\tau_k'$ as the ideal arc obtained by flipping $\tau_k$ in the topological quadrilateral formed by the two adjacent triangles  of $\mathcal T$ having $\tau_k$ as common edge (see Figure above).
  Then as $\mathcal T$, $\mu_k(\mathcal T)=(\tau_j')_{j=1}^n$ is a maximal collection of ideal arcs on $S$ and one can verify that the exchange matrix 
associated to it coincides with the one obtained by mutation of $B^{\mathcal T}$ in the $k$-th direction: i.e.\,$\mu_k(B({\mathcal T}))=B({\mu_k(\mathcal T)})$.  One defines the cluster algebra $\boldsymbol{A}(S,P)$ associated to the pair $(S,P)$ as the one with exchange matrix $B(\mathcal T)$ for a maximal collection of arcs $\mathcal T$ as above. Since any other such collection can be obtained from $\mathcal T$ after a finite number of flips,  
 $\boldsymbol{A}(S,P)$ does not depend on the choice of $\mathcal T$ (up to mutations) hence this definition makes sense.  
\end{exm}



 \subsection{\bf Cluster webs.}
\label{SSub:ClusterWebs}
We have reached the point where we can finally define the main objects of study in this text. As many good and relevant definitions in mathematics, the general definition we first give is quite simple,  general and not very precise.\sk 

Let $\boldsymbol{A}=\boldsymbol{A}(B)$ be a cluster algebra defined by a skew-symmetrizable $n\times n$ matrix $B$, for some $n\geq 2$.  Let $\Sigma$ be a non-empty finite set of ($\boldsymbol{\mathcal A}-$ or $\boldsymbol{\mathcal X}-$)cluster variables. By definition, the {\bf cluster web associated} to it is the web
$\boldsymbol{\mathcal W}_{\Sigma}$ admitting as first integrals the elements of $\Sigma$: 
$$
\boldsymbol{\mathcal W}_{\Sigma}=\boldsymbol{\mathcal W}\big( \,  x \hspace{0.15cm}  \lvert  \hspace{0.15cm}  x\in  \Sigma\, 
\big)\, .
$$
Such a `{\it $\boldsymbol{\mathcal A}-$ or $\boldsymbol{\mathcal X}-$cluster web'} (the terminology depending on $\Sigma$ is included in $\boldsymbol{\mathcal A var}( B)$ or in $\boldsymbol{\mathcal X var}( B)$) is a web on $\mathbf C^n$ defined by rational first integrals.  We will give further some general properties of such webs but for the moment we agree that this definition is a bit elliptic. \mk 

{ The notion of `cluster web' is one of the most important of this text and a few preliminary basic remarks  about it are in order: 
\begin{itemize}
\item  First, in general the set of cluster variables is infinite and 
 the potential relevance of the notion of cluster web resides precisely in considering finite subsets  of cluster variables giving rise to interesting webs. 
\sk
\item Second, there is a slight inaccuracy in the above definition of what is a `cluster web' since the space on which it has to be considered is not clearly indicated. We will  settle this later on in \S\ref{Par:Cluster-tori}.
\sk
\item  The way they are defined above, cluster webs are unordered webs a priori since there is no natural ordering of the cluster variables. However,  we will always work with ordered webs in practice, after having ordered the cluster first integrals  elements of $\Sigma$ 
in one way or another.\footnote{In most of the cases we will deal with, we will order the cluster variables from the less to the more {\it complex} one, being aware that the notion of `complexity' we are considering here is subjective and not rigorously defined mathematically.}
\sk  
\item 
A basic but important thing 
that 
might happen is that two distinct cluster variables belonging to $\Sigma$ define the same foliation on $\mathbf C^n$ and that even knowing $\Sigma$ quite explicitly, it may be rather difficult to describe precisely the foliations of $\boldsymbol{\mathcal W}_{\Sigma}$, or even just to enumerate them.  This point, regarding which we have precise conjectures,  
will be discussed more in length further on in \S\ref{SS:Some-conjectures}.
\sk 
\item
Finally, we  have to say that from a web-theoretic perspective, the interesting cluster webs are the ones of $\boldsymbol{\mathcal X}$-type, that is those defined by rational first integrals which are $\boldsymbol{\mathcal X}$-cluster variables.  This is  not surprising since it is known that such cluster variables are of a much more geometric nature than the $\boldsymbol{\mathcal A}$-cluster variables.  
In this text, we will  focus on $\boldsymbol{\mathcal X}$-cluster webs.  
\end{itemize}
}

 At this point, it is time to  look at a  first example, 
 which is already interesting although very basic.

\subsubsection{An example: the $\boldsymbol{A_2}$ case.} 
\label{Par:A2-case}
To the Cartan matrix 
\scalebox{0.8}
{$\bigg[\hspace{-0.15cm}
\begin{tabular}{c}
\, 2 \, -1\\ 
-1 \, \,  2
\end{tabular}\hspace{-0.15cm}\bigg]$
}
of type $A_2$  is associated the following exchange matrix and  quiver:
$$ 
 {}^{} \hspace{0.2cm}
 B_{A_2}=\scalebox{0.9}{$\begin{bmatrix}
 0 & 1 \\
 -1 & 0
 \end{bmatrix}$} \qquad \qquad \mbox{ and } \qquad   \qquad 
  Q_{A_2}\, :  \,   \,
 {\bf 1} \longrightarrow  {\bf 2}\, . 
 $$

In this case, and this contrasts with what happens in general,  it is not difficult to write down explicit formulas for the seeds obtained by successive mutations.   In order to get nicer formulas for the $\boldsymbol{\mathcal X}$-variables, one takes $S_0=\big( \, (1/x_1,x_2)\, ,  B_{A_2}\big)$ for the   initial $\boldsymbol{\mathcal X}$-seed. Then the seeds which can be constructed from $S_0$ by means of  mutations are given in Figure \ref{Fig:ExchangeGraph-A2} below. \mk 


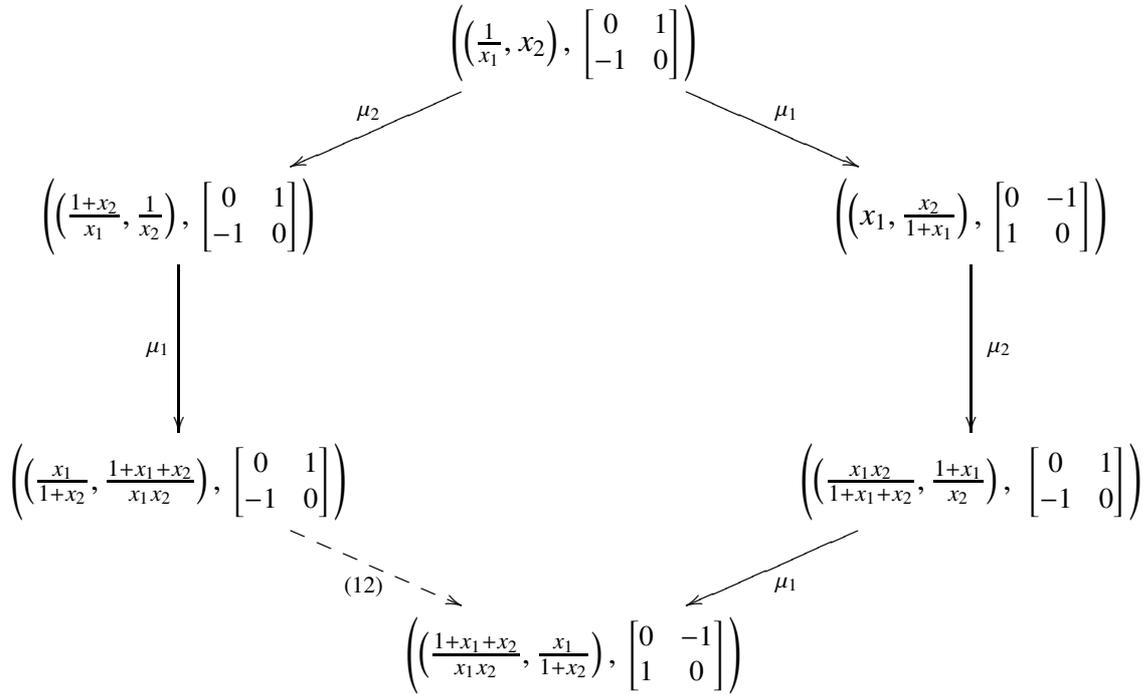
\begin{figure}[!h]
   \begin{center}
  \scalebox{1.1}{
    \xymatrix@R=0.9cm@C=0.5cm{    
    &  {{ \Bigg( \Big(\frac{1}{x_1},x_2\Big)\, , 
\scalebox{0.9}{ 
$\begin{bmatrix}
 0 & 1 \\
 -1 & 0 
 \end{bmatrix}$}\,  \Bigg)}}
 \ar@{->}[rd]^{ \hspace{0.1cm}\mu_1} 
  \ar@{->}[dl]_{ \hspace{0.1cm} \mu_2} &  \\
   \Bigg( \Big(\frac{1+x_2}{x_1},\frac{1}{x_2}\Big)\, , 
\scalebox{0.9}{ 
$\begin{bmatrix}
 0 & 1 \\
 -1 & 0 
 \end{bmatrix}$} \, \Bigg)  \ar@{->}[dd]_{ \hspace{0.1cm}\mu_1}     
      &    & 
     \ar@{->}[dd]^{ \hspace{0.1cm}\mu_2}    \Bigg( \Big({x_1},\frac{x_2}{1+x_1}\Big)\, , 
\scalebox{0.9}{ 
$\begin{bmatrix}
 0 & -1 \\
 1 & 0 
 \end{bmatrix}$} \, \Bigg)\\ 
 & & 
 \\
 \Bigg( \Big(\frac{x_1}{1+x_2},\frac{1+x_1+x_2}{x_1x_2}\Big)\, , 
\scalebox{0.9}{ 
$\begin{bmatrix}
 0 & 1 \\
 -1 & 0 
 \end{bmatrix}$}\,  \Bigg) 
  \ar@{-->}[dr]_{ \hspace{0.0cm} (12)} 
  &      &  
   \ar@{->}[dl]^{ \hspace{0.1cm}\mu_1} 
    \Bigg( \Big(\frac{x_1x_2}{1+x_1+x_2},\frac{1+x_1}{x_2}\Big)\, ,\, 
\scalebox{0.9}{ 
$\begin{bmatrix}
 0 & 1 \\
 -1 & 0 
 \end{bmatrix}$}\,  \Bigg)\\
 &       \Bigg( \Big(\frac{1+x_1+x_2}{x_1x_2},\frac{x_1}{1+x_2}\Big)\, , 
\scalebox{0.9}{ 
$\begin{bmatrix}
 0 & -1 \\
 1 & 0 
 \end{bmatrix}$} \, \Bigg)  &
 }}
 \end{center}
 \caption{The $\boldsymbol{\mathcal X}$-exchange graph of the cluster algebra of type $A_2$}
 \label{Fig:ExchangeGraph-A2}
 \end{figure}

We see that the two $\boldsymbol{\mathcal X}$-seeds $S=\mu_2\circ \mu_1(S_0)$ and $S'=\mu_1\circ \mu_2\circ \mu_1(S_0)$  are equivalent by means of the transposition $(1,2)$ ({\it i.e.}\,one has $S^{(12)}=S'$, see \eqref{Eq:Equiv-of-Seeds} for this notation). This implies that the $\boldsymbol{\mathcal X}$-exchange graph of the cluster algebra defined by $B_{A_2}$ is a pentagon and that there are precisely five $\boldsymbol{\mathcal X}$-clusters and just as many $\boldsymbol{\mathcal X}$-custer variables (up to inversion $x\leftrightarrow x^{-1}$):\footnote{It is well-known (and easy to verify) that the complete exchange graph $\mathcal E\Gamma_{A_2}$ associated to $B_{A_2}$ is a pentagon as well.}
 \begin{align*}
\boldsymbol{\mathcal X var}\Big( B_{A_2}\Big)= &\, 
\left\{  \, {x_1}  \,  , \,  {x_2} \,  ,  \, {\frac{1+x_2}{x_1}} \,  ,  \, {\frac{1+x_1}{x_2}} \,  ,  \, {\frac{1+x_1+x_2}{x_1x_2}}\right\} 
 \end{align*} 
  
 By definition (see \eqref{Eq:AWD-XWD} further), 
   the `{\it $\boldsymbol{\mathcal X}$-cluster web of type $A_2$}', denoted by $\boldsymbol{\mathcal X}\hspace{-0.01cm}\boldsymbol{\mathcal W}_{\hspace{-0.01cm} A_2}$,  
 is the 5-web on $\mathbf C^2$ defined by all the $\boldsymbol{\mathcal X}$-custer variables, {\it i.e.}
$$
\boldsymbol{\mathcal X}\hspace{-0.01cm}\boldsymbol{\mathcal W}_{\hspace{-0.01cm} A_2}
=
\boldsymbol{\mathcal W}\left( \, {x_1}  \,  , \,  {x_2} \,  ,  \,  {\frac{1+x_2}{x_1}} \,  ,  \, 
 {\frac{1+x_1}{x_2}} \,  ,  \,  {\frac{1+x_1+x_2}{x_1x_2}} \, 
 \right) \, . 
$$
 
 The main feature of this web, and this was our main motivation to consider and study the {\it `cluster webs'}, is that the following functional identity  holds true for any $x_1,x_2\in \mathbf R_{>0}$
 $$
 {\sf R}({x_1})  +{\sf R}(  {x_2} )+ {\sf R}\left( {\frac{1+x_2}{x_1}}\right) 
 +
{\sf R}\left(  {\frac{1+x_1}{x_2}}\right) 
+{\sf R}\left(  {\frac{1+x_1+x_2}{x_1x_2}}\right)= {\pi^2}/{2}\, , 
$$ where ${\sf R}$ stands for the 
  {\it `cluster dilogarithm'} defined by 
${\sf R}(u)
=-\l {2} (-u)-\frac{1}{2} {\rm Log}(u){\rm Log}(1+u)$ for  $u\in \mathbf R_{>0}$  (introduced and discussed in \S\ref{Par:NotationDilogarithmicFunctions} above).\mk 

It is not difficult to verify that $\boldsymbol{\mathcal X}\hspace{-0.01cm}\boldsymbol{\mathcal W}_{\hspace{-0.01cm} A_2}$ is equivalent to Bol's web $\boldsymbol{\mathcal B}$. The latter being so important regarding web-geometry, this motivates for studying cluster webs more systematically. This is the main theme of the present section of this text. 
\begin{center}
$\star$
\end{center}
In the next paragraph, using  a basic result of the theory of cluster algebras, we are going to introduce several series of webs generalizing the $\boldsymbol{\mathcal X}$-cluster web of type $A_2$ considered above, and 
 which look quite interesting in what concerns their abelian relations and their rank(s).

\paragraph{Cluster algebras of finite type.}
\label{Par:Cluster-Algebras-Finite-Type}
%
We now describe the  classification of cluster algebras of {\it `finite type'}. This is  a fundamental result  of the theory due to Fomin and Zelevinsky \cite{FZII}  which allows us to define several families of cluster webs, one for each Dynkin diagram $\Delta$.\mk

Let $\mathcal E\Gamma_{\hspace{-0.05cm}Q}$ be the exchange graph associated to the cluster pattern defined by (the skew-symmetri\-zable matrix $B=B_Q$ associated to)   a valued cluster quiver 
  $Q$.    Then the associated cluster algebra $\boldsymbol{A}_Q$ is said to be of {\bf finite type} if $\mathcal E\Gamma_{\hspace{-0.05cm}Q}$ is finite. It is clearly equivalent to the facts that $(i)$ $B$ is mutation finite;  and $(ii)$ there exists only a finite number of $\boldsymbol{\mathcal A}$-clusters or equivalently, of $\boldsymbol{\mathcal A}$-cluster variables.   
  
 \subparagraph*{Dynkin diagram and valued quivers.} 
  Let $(Q,v)$ be a valued quiver. 
  By \eqref{Eq:Convention}, one can assume that $v_1(e)>0$
  for any edge 
   $e $ in it.  Then,  by definition, the diagram  
 associated to it is the quiver   $D(Q,v)$ with the same vertices as $Q$ but where each oriented edge $ e:  i{\longrightarrow} j$ is replaced by $v_1(e)$ distinct oriented edges from $i$ to $j$ if $v_1(e)>1$, and just by 
 $ i \stackrel{\line(1,0){15}}{} j$ otherwise.   Then $(Q,v)$ is said to be {\bf of Dynkin type $\boldsymbol{\Delta}$} if $D(Q,v)=\Delta$ for $\Delta$ being one of the Dynkin diagrams in Figure \ref{Fig:DynkinDiagram} below.
 \mk 
 \begin{figure}[h!]
\begin{center}
\resizebox{6in}{3in}{
 \includegraphics{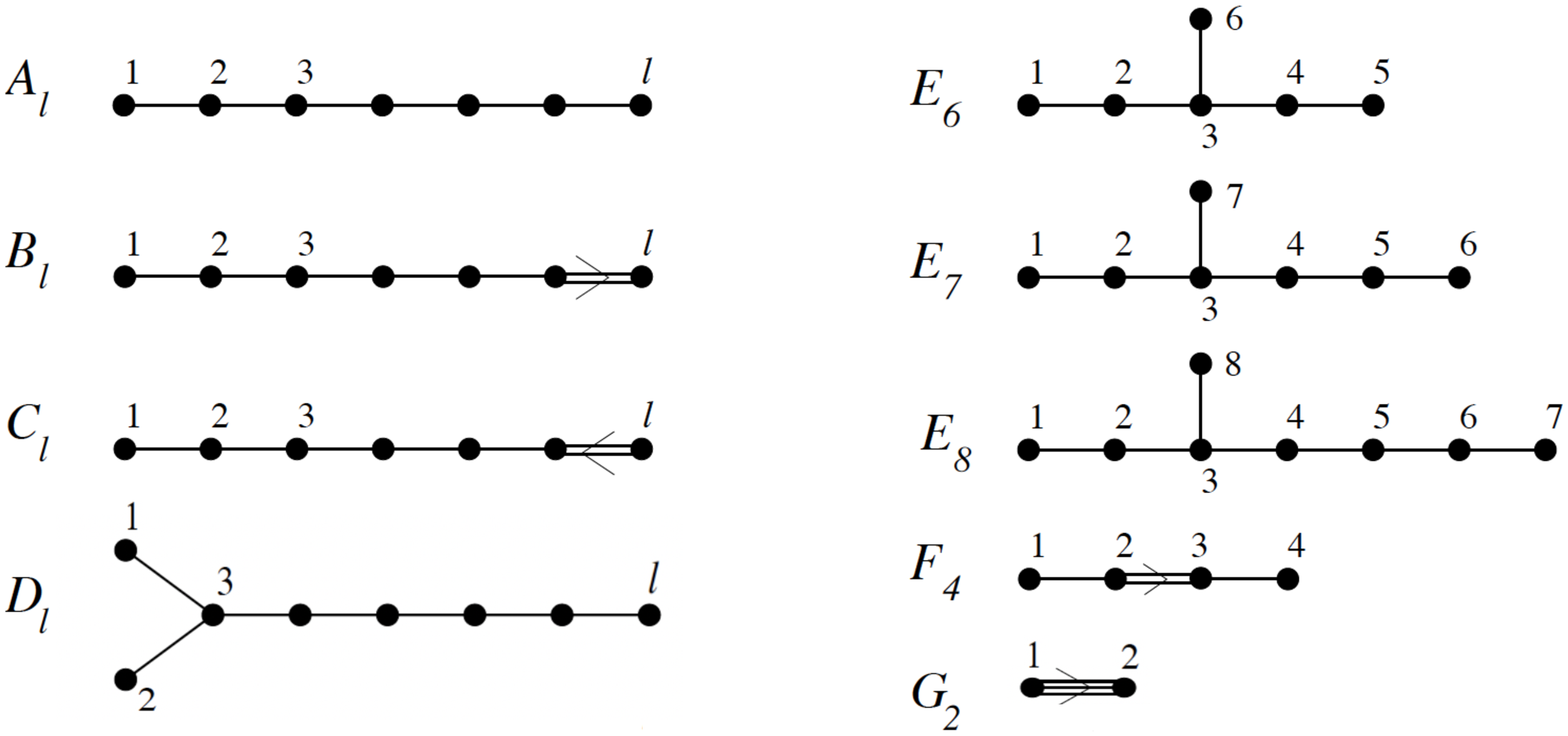}}
 \vspace{-0.5cm}
\caption{Dynkin diagrams} 
\label{Fig:DynkinDiagram}
\end{center}
\end{figure}

%
%
%
%
%
%
%
%
%
%
%
%
%
%
%
%
%
%

Conversely, let us associate a valued quiver $\vec{\Delta}$ attached to any 
Dynkin diagram $\Delta$. First, we remind that a {\bf source} (resp.\,a {\bf sink}) in a quiver $Q$ is a vertice 
which is not the endpoint  (resp. the initial point) of any oriented edge $e\in Q_1$. A quiver is bipartite if any of its vertices is either a source or a sink. 
In such a quiver, sinks will be indicated by black dots 
\scalebox{1.4}{$\bullet$}, whereas the sources by white ones 
\scalebox{0.7}{$\boldsymbol{\bigcirc}$}. 

\begin{figure}[h!]
\begin{center}
\resizebox{3in}{3in}{
 \includegraphics{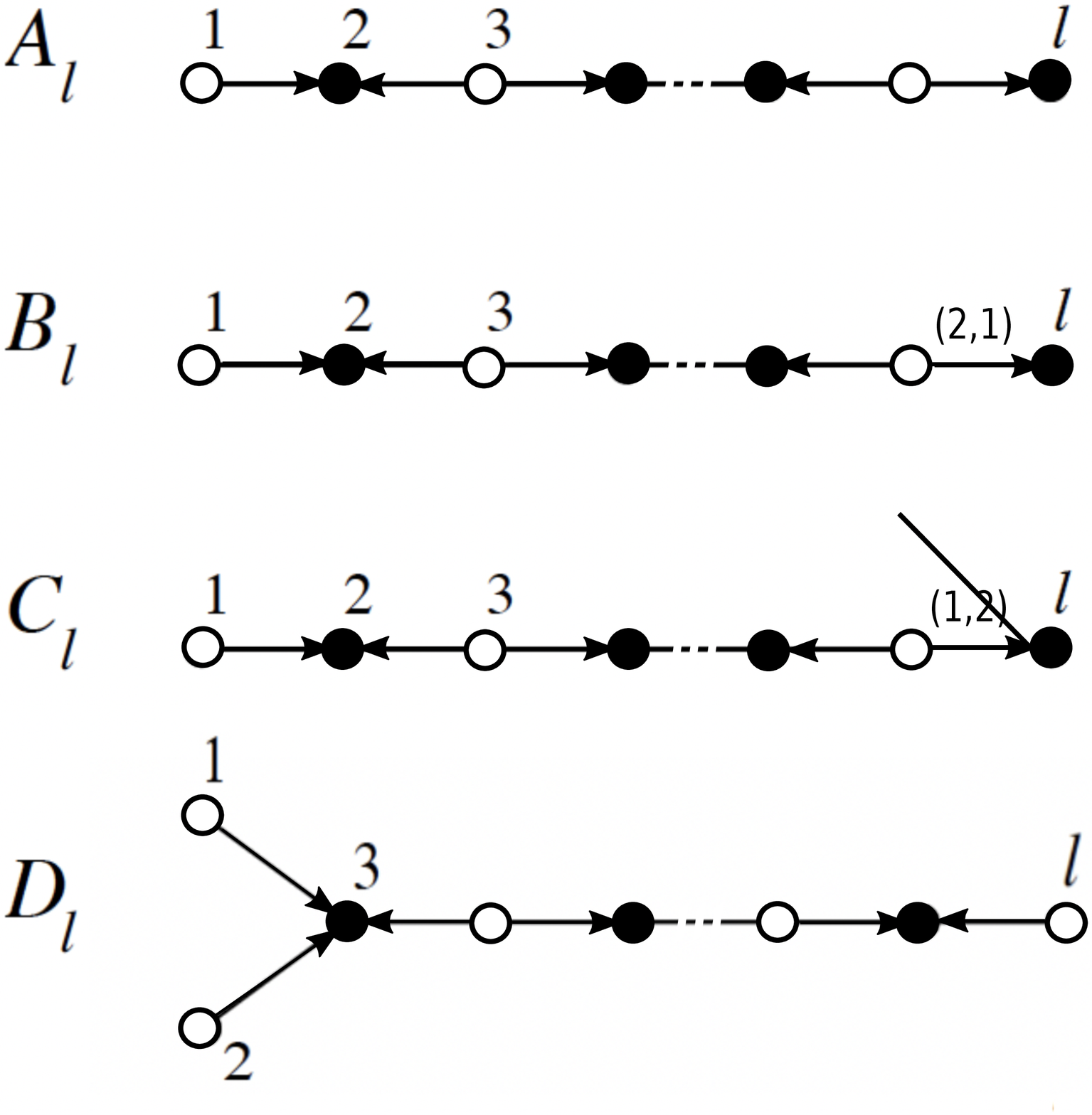}}
 \vspace{-0.5cm}
\caption{Dynkin quivers $\vec{\Delta}$ for $\Delta$ of classical type (of rank  $l$ even).} 
\label{Fig:DynkinQuiverDiagram}
\end{center}
\end{figure}

 Let $\Delta$ be a connected Dynkin diagram.  Then there exists a unique coloring of the vertices of $\Delta$ as sources of sinks in order to get a bipartite graph $\Delta^b$ with the first vertex being a source. 
When $\Delta$ is simply-laced ($A,D,E$ case), one  
gets the valued quiver $\vec{\Delta}$ by orienting the edges of $\Delta^b$ from sources to sinks as follows 
$
\scalebox{0.7}{$\boldsymbol{\bigcirc}$} 
 \hspace{-0.1cm}
\begin{tabular}{c} 
$\xrightarrow{\hspace*{0.6cm}}$ 
 \end{tabular} 
\hspace{-0.1cm} \scalebox{1.4}{${\bullet}$}
$
 and setting $v(e)=(1,1)$ for any such oriented edge $e$.  
 In case when  $\Delta$ is multi-laced,  for two vertices $i$ and $j$ in $\Delta^b$ with $k\in \{1,2,3\}$ edges from $i$ to $j$, one associates an edge in $\vec{\Delta}$ which goes from $i$ to $j$ and has valuation $(k,1)$ if $i$ is a source, 
and which goes from $j$ to $i$ with valuation $(1,k)$ if $i$ is a sink.  
 Some Dynkin quivers are pictured in  Figure \ref{Fig:DynkinQuiverDiagram}. 
\mk

Given $\Delta$, we denote by $B_{{\Delta}}$ the bipartite skew-symmetrizable matrix associated to $\vec{\Delta}$.  For instance, these matrices  for $\Delta$ of rank 2 and 3 are the following: 
\begin{align*}
B_{A_2}=& \begin{bmatrix}
0 & 1\\
-1 & 0
\end{bmatrix}
 &&B_{A_3}= \begin{bmatrix}
0 & 1 & 0\\
-1 & 0 & -1\\
0 & 1 & 0
\end{bmatrix}
 \\ 
    B_{B_2}=& \begin{bmatrix}
0 & 2\\
-1 & 0
\end{bmatrix}
&& B_{B_3}= \begin{bmatrix}
0 & 1 & 0\\
-1 & 0 & -2 \\
0 & 1 & 0 \\
\end{bmatrix}
 \\ 
B_{C_2}=& \begin{bmatrix}
0 & 1\\
-2 & 0
\end{bmatrix}
&& B_{C_3}= \begin{bmatrix}
0 & 1 & 0\\
-1 & 0 & -1 \\
0 & 2 & 0 \\
\end{bmatrix}
 \\ 
\mbox{and } \qquad  B_{G_2}=& \begin{bmatrix}
0 & 3\\
-1 & 0 
\end{bmatrix}
&& 
\end{align*}
By definition, the {\bf cluster algebra of Dynkin type $\boldsymbol{\Delta}$} is the one obtained by taking $B_{{\Delta}}$ as initial exchange matrix. 
  One can then state one of the main basic results of the theory of cluster algebras, due to Fomin and Zelevinsky: 
\begin{thm}[\cite{FZII}] 
\begin{enumerate}
\item The cluster algebra  $\boldsymbol{A}_Q$ associated to a valued quiver
$Q$ 
 is of finite type if and only if  $Q$ is mutation equivalent to a valued quiver of Dynkin type $\vec{\Delta}$.
\item  Two valued cluster quivers $Q_1$ and $Q_2$ of Dynkin type are mutation equivalent if and only if they define the same Dynkin (or `Cartan-Killing') type. 
\end{enumerate}
\end{thm}

To any Dynkin Diagram $\Delta$, since the total number of associated cluster variables of 
$\boldsymbol{A}_\Delta$
is finite,  one can consider the following two webs:
\begin{align}
\label{Eq:AWD-XWD}
\boldsymbol{\mathcal A}\hspace{-0.05cm}\boldsymbol{\mathcal W}_\Delta=& \, \boldsymbol{\mathcal W}\Big(   \, \mbox{ all } \boldsymbol{\mathcal A}-\mbox{cluster variables of }\, 
\boldsymbol{A}_\Delta\, 
\Big)\\ 
\mbox{ and }\quad 
\boldsymbol{\mathcal X}\hspace{-0.05cm}\boldsymbol{\mathcal W}_\Delta=& \, \boldsymbol{\mathcal W}\Big(   \, \mbox{ all } \boldsymbol{\mathcal X}-\mbox{cluster variables of }\, 
\boldsymbol{A}_\Delta\, 
\Big)\, .  \nonumber 
\end{align}

By definition, these two webs are respectively the {\bf $\boldsymbol{\mathcal A}$-} and {\bf $\boldsymbol{\mathcal X}$-cluster web} of type $\boldsymbol{\Delta}$.   \mk

From the second point of the preceding theorem, it follows that the Dynkin type of a finite type cluster  algebra $\boldsymbol{A}$ is well-defined hence  for 
$\boldsymbol{\mathcal R}$ standing for $\boldsymbol{\mathcal A}$ or 
$\boldsymbol{\mathcal X}$, 
we will  call (a bit abusively) `$\mathcal R$-cluster web'  and will denote by the same notation $\boldsymbol{\mathcal R}\hspace{-0.05cm}\boldsymbol{\mathcal W}_\Delta$,  any cluster web whose the first integrals are all the $\boldsymbol{\mathcal R}$-cluster variables of a finite type cluster algebra of Dynkin type $\Delta$.    
In other words, 
$\boldsymbol{\mathcal R}\hspace{-0.05cm}\boldsymbol{\mathcal W}_\Delta$ will stand for any webs in the equivalence class (up to isomorphisms of webs induced by cluster mutations) of the web with initial matrix $B_\Delta$.  Since the objetcs/properties of  webs  we are dealing with (RAs, rank(s), being AMP, etc.) make sense and behave well  up to equivalence, this will not be the source of any real problem.

\subparagraph*{Some examples of cluster webs associated to finite type cluster algebras.}
 Since $B_2\simeq C_2$, the cluster webs of these two types are equivalent. Thus there exist only three  Dynkin diagrams of rank 2  to consider: $A_2$, already studied in \S\ref{Par:A2-case}, $B_2$ and $G_2$.   
 The 
 $\boldsymbol{\mathcal A}$-cluster webs 
 for the initial $\boldsymbol{\mathcal A}$-seeds  $((a_1,a_2),B_\Delta)$ with $B_{A_2}=
 \Big[ \hspace{-0.1cm}
\scalebox{0.8}{\begin{tabular}{cc}
0 &  \hspace{-0.3cm}1  \vspace{-0.1cm}\\
-1 & \hspace{-0.45cm} 0
\end{tabular}}
   \hspace{-0.1cm}  \Big]$, $B_{B_2}=\Big[ \hspace{-0.1cm}
\scalebox{0.8}{\begin{tabular}{cc}
0 &  \hspace{-0.3cm}2  \vspace{-0.1cm}\\
-1 & \hspace{-0.44cm} 0
\end{tabular}}
   \hspace{-0.1cm}  \Big]$ and 
   $B_{G_2}=\Big[ \hspace{-0.1cm}
\scalebox{0.8}{\begin{tabular}{cc}
0 &  \hspace{-0.3cm}-1  \vspace{-0.1cm}\\
3 & \hspace{-0.3cm} 0
\end{tabular}}
   \hspace{-0.1cm}  \Big]$ are 
   \begin{align*}
\boldsymbol{\mathcal A}\hspace{-0.01cm}\boldsymbol{\mathcal W}_{\hspace{-0.01cm} A_2}
=& \, 
\boldsymbol{\mathcal W}\left( \, \frac{1}{a_{{1}}} \,,\, \frac{1}{a_{{2}}}\, ,\, 
{\frac {1+a_{{2}}}{a_{{1}}}},{\frac {1+a_{{1}}}{a_{{2}}}},{\frac {1+a_{{1}}+a_{{2}}}{a_{{1}}a_{{2}}}}
\right) \, , 
\\
 \boldsymbol{\mathcal A}\hspace{-0.02cm}\boldsymbol{\mathcal W}_{\hspace{-0.01cm} B_2}
=& \, 
\boldsymbol{\mathcal W}\left( \, \frac{1}{a_{{1}}} \,,\, \frac{1}{a_{{2}}}\, ,\, 
{\frac {1+a_{{2}}}{a_{{1}}}},{\frac {1+{a_{{1}}}^{2}}{a_{{2}}}},{\frac {1+a_{{2}}+{a_{{1}}}^{2}}{a_{{1}}a_{{2}}}},{\frac {{a_{{1}}}^{2}+{(1+
a_{{2}})}^{2}}{{a_{{1}}}^{2}a_{{2}}}} 
 \right)  \, , 
 \\
 \mbox{and}\quad 
 \boldsymbol{\mathcal A}\hspace{-0.01cm}\boldsymbol{\mathcal W}_{\hspace{-0.01cm} G_2}
=& \, 
\boldsymbol{\mathcal W}\Bigg( \, \frac{1}{a_{{1}}} \,,\, \frac{1}{a_{{2}}}\, ,\, 
{\frac {1+a_{{1}}}{a_{{2}}}},{\frac {1+{a_{{2}}}^{3}}{a_{{1}}}} \, ,\, 
{\frac {1+a_{{1}}+{a_{{2}}}^{3}}{a_{{1}}a_{{2}}}}\, ,\, 
{\frac {{(1+a_1)^2}+{a_{{2}}}^{3}}{a_{{1}} {a_{{2}}}^{2}}} \, ,\,
\\ 
&\hspace{3.5cm}  {\frac {{(1+a_{{1}})}^{3}+{a_{{2}}}^{3}}{a_{{1}} {a_{{2}}}^{3}}} 
\, , \, 
{\frac {{(1+a_{{1}})}^{3}+  3\,a_1 {a_{{2}}}^{3} +{a_{{2}}}^{6}}{{a_{{1}}}^{2}{a_{{2}}}^{3}}} \, 
 \Bigg) \, , 
\end{align*}
 whereas,   when  taking $(x_1^{-1},x_2)$ as initial $\boldsymbol{\mathcal X}$-cluster    in the cases ${A_2}$ and $B_2$,  and  $(x_1,x_2^{-1})$  in the ${G_2}$-case, the associated $\boldsymbol{\mathcal X}$-cluster webs
are the following ones:
\begin{align*}
\boldsymbol{\mathcal X}\hspace{-0.01cm}\boldsymbol{\mathcal W}_{\hspace{-0.01cm} A_2}
=& \, 
\boldsymbol{\mathcal W}\left( \, \frac{1}{x_1}  \,  , \,  \frac{1}{x_2} \,  ,  \,  {\frac{1+x_2}{x_1}} \,  ,  \, 
 {\frac{1+x_1}{x_2}} \,  ,  \,  {\frac{1+x_1+x_2}{x_1x_2}} \, 
 \right) \\
 \boldsymbol{\mathcal X}\hspace{-0.02cm}\boldsymbol{\mathcal W}_{\hspace{-0.01cm} B_2}
=& \, 
\boldsymbol{\mathcal W}\left( \, \frac{1}{x_1}  \,  , \,  \frac{1}{x_2} \,  ,  \,  {\frac{1+x_2}{x_1}} \,  ,  \, 
 {\frac{(1+x_1)^2}{x_2}} \,  ,  \,  {\frac{(1+x_1)^2+x_2}{x_1x_2}} \, , 
 \,  {\frac{(1+x_1+x_2)^2}{x_1^2x_2}} \,
 \, 
 \right) \\
 \mbox{and}\quad 
 \boldsymbol{\mathcal X}\hspace{-0.01cm}\boldsymbol{\mathcal W}_{\hspace{-0.01cm} G_2}
=& \, 
\boldsymbol{\mathcal W}\Bigg( \, \frac{1}{x_1}  \,  , \,  \frac{1}{x_2} \,  ,  \,  {\frac{1+x_1}{x_2}} \,  ,  \, 
{\frac{(1+x_2)^3}{x_1}} \,  ,  \, 
\frac{x_1+(1+x_2)^3}{x_1x_2}\,  ,  \\ 
& {}^{} \hspace{1.4cm}
\frac{(1+x_1)^2+ 3(1+x_1)x_2
+3x_2^2+x_2^3}{x_1x_2^2}  \, , \, 
\frac{ (1+x_1+x_2)^3}{x_1x_2^3} \,  ,  \, \frac{(x_1 + (1+x_2)^2)^3}{x_1^2x_2^3}\,   \Bigg) \, .
\end{align*}



For an example in  rank 3, 
 the two cluster webs in type $A_3$ 
 with 
initial seeds $( (a_1,a_2,a_3), B_{A_3})$
and $\big(  (\, x_1^{-1}\, ,\, x_2\, ,\, x_3^{-1})\, , B_{A_3}\big)$ are respectively: 
\begin{align}
\label{Eq:XWA3-Definition}
\boldsymbol{\mathcal A \hspace{-0.05cm}\mathcal W}_{\hspace{-0.05cm} A_3}= &\, 
\boldsymbol{\mathcal W}
 \Bigg( \, \frac{1}{a_{{1}}} \,,\, \frac{1}{a_{{2}}}\, ,\, \frac{1}{a_{{3}}}  \,,\,   {\frac {1+a_{{2}}}{a_{{1}}}}  \,,\,  {\frac {1+a_{{2}}}{a_{{3}}}}  \,,\,  {\frac {1+a_{{1}}a_{{3}}}{a_{{2}}}}  \,,\,   {\frac {
 1+a_{{2}}+
 a_{{1}}a_{{3}}}{a_{{1}}a_{{2}}}}  \,,\,  \nonumber 
 \\ &    \hspace{4.8cm}
  {\frac {
  1+a_{{2}}+
 a_{{1}}a_{{3}}
 }{a_{{2}}a_{{3}}}}  \,,\,   
 {\frac {
 1+2\,a_{{2}}+{a_{{2}}}^{2}+
 a_{{1}}a_{{3}}}{a_{{1}}a_{{2}}a_{{3}}}} \Bigg) \nonumber 
\\
\mbox{and } \quad 
\boldsymbol{\mathcal X \hspace{-0.05cm}\mathcal W}_{\hspace{-0.05cm} A_3}= & \, 
\boldsymbol{\mathcal W}\Bigg(\, 
\frac{1}{x_1} \, , \, 
\frac{1}{x_2} \, , \, 
\frac{1}{x_3} \, , \, 
\frac{1+x_2}{x_1} \, , \, 
\frac{1+x_1}{x_2} \, , \, 
\frac{(1+x_1)(1+x_3)}{x_2} \, , \, 
\frac{1+x_3}{x_2} \, , \, 
\frac{1+x_2}{x_3} \, , \, 
\nonumber  \\
& \hspace{1.5cm}
\frac{1+x_1+x_2}{(1+x_1)x_3} \, ,  \, 
  \frac{1+x_2+x_3}{(1+x_3)x_1} \, , \, 
\frac{1+x_2+x_3}{x_2x_3} \, , \, 
\frac{1+x_1+x_2+x_3+x_1x_3}{x_1x_2} \, ,
\, 
 \\
& \hspace{3.5cm}
\frac{1+x_1+x_2+x_3+x_1x_3}{x_2x_3} \, , \, 
\frac{(1+x_1+x_2)(1+x_2+x_3)}{x_1x_2x_3} \,  \Bigg)\, . \nonumber 
\end{align}

As it appears clearly when considering the examples above (but this is a general fact which follows easily from the `separation formulas'  
\eqref{Eq:X-var:separation-formula} and \eqref{Eq:A-var:separation-formula} given below), the $\boldsymbol{\mathcal A}$-cluster webs are defined by cluster variables which are formally simpler than thoses defining the 
$\boldsymbol{\mathcal X}$-cluster webs. However, in all the examples we have studied, the latter webs seem far more interesting than the former from a web-theoretic perspective, especially  in what concerns their ARs and their rank (see for instance  \S\ref{subpar:AWA3} for an illustration of this fact). 
 Accordingly, we will only consider $\boldsymbol{\mathcal X}$-cluster webs in this text and will  
study  these webs more in depth further.

\begin{rem}
\label{Rem:wGP-cluster-Webs}
The cluster variables $1/x_1$, $1/x_2$ and $({1+x_2})/{x_1}$ are first integrals of a 3-subweb of $\boldsymbol{\mathcal X \hspace{-0.05cm}\mathcal W}_{\hspace{-0.05cm} A_3}$ and obviously only depend on  $x_1$ and $x_2$. Consequently, this 3-subweb has intrinsic dimension 2 hence the $\mathcal X$-cluster web of type $A_3$ does not satisfy the `{\it strong general position assumption}' for webs discussed in \S\ref{Par:GermOfWebs}.  This seems to be a general fact regarding cluster webs: they only satisfy the assumption (wGP).  
\end{rem}

\paragraph{Some fundamental objects and results.}
\label{Par:Fundamental-Results-Cluster-Algebras}
We discuss here some important properties of cluster variables and what these properties imply for cluster webs. Since only the $\boldsymbol{\mathcal X}$-side is relevant for webs, we will focus on the  cluster variables  of this type even though the properties discussed below admit versions for the $\boldsymbol{\mathcal A}$-variables.


\subparagraph{$\boldsymbol{F}$-polynomials and $\boldsymbol{c}$-vectors.}
\label{SubPar:Fpoly-and-c-vectors}
We describe some objects associated to each seed of a given cluster algebra, introduced very early in the development of the theory, by Fomin and Zelevinsky in \cite{FZIV}. \mk 

Let $\big\{S^t=(\boldsymbol{x}^t,\boldsymbol{a}^t,B^t) \big\}_{t\in \mathbb T_n}$ 
be the cluster pattern of a cluster algebra ${\bf A}$ of rank $n$, with initial seed corresponding to a vertex $t_0$. To it is associated another pattern indexed by the vertices of $ \mathbb T_n$, formed by {\bf $\boldsymbol{FC}$-seeds} $\Gamma^t=(\boldsymbol{F}^t, C^t)$ where for any $t$: 
\begin{itemize}
\item $\boldsymbol{F}^t$ is a $n$-tuple $(F^t_{i}(u))_{i=1}^n$ of rational function 
 in $n$ indeterminates $u_1,\ldots,u_n$;
\item $C^t=(c^t_{ij})_{i,j=1}^n$ is a $n\times n$ square matrix with integer entries.
\end{itemize}
The $\Gamma^t$'s are obtained recursively from 
$\Gamma^{t_0}$ by requiring that for the $FC$-seed, one has 
$$
F^{t_0}_i(u)=1\quad \mbox{ for }\, i=1,\ldots,n\, \qquad \mbox{ and }
\qquad C^{t_0}={\rm Id}_n\, ; 
$$
and assuming that, for any edge  $t
 \hspace{-0.1cm}
\begin{tabular}{c} 
$ \stackrel{k}{\line(1,0){25}}$ \vspace{-0.35cm}\\  ${}^{}$
 \end{tabular} 
\hspace{-0.2cm}
 t'$ in $\mathbb T_n$, one has  for any $i,j=1,\ldots,n$ with $j\neq k$: 
\begin{align}
 c^{t'}_{ik}= &
-c^{t}_{ik}
&
 \mbox{ and } \qquad
 & \hspace{0.5cm}
 c^{t'}_{ij}=   c^{t}_{ij}+c^t_{ik}\cdot \big[b_{ij}^t\big]_+  + \big[ -c_{ik}^t\big]_+\cdot b_{kj}^t\, 
\,  \\
\nonumber
F^{t'}_{k}(u)= &\, {M_k^t(u)}/{F_k^t(u)}
&
 \mbox{ and } \qquad
 & F^{t'}_{\hspace{-0.05cm}j}(u)= {F_{\hspace{-0.05cm}j}^t(u)}\, , 
\end{align}
  where $M_k^t(u)=\prod_{l=1}^n u_l^{\big[c^t_{lk}\big]_+} F_{\hspace{-0.05cm}l}^t(u)^{\big[b^t_{lk}\big]_+}
+\prod_{l=1}^n u_l^{\big[-c^t_{lk}\big]_+} F_{\hspace{-0.05cm}l}^t(u)^{\big[-b^t_{lk}\big]_+}$.\mk 
 
 From the previous formulas, it is straightforward that $F_{\hspace{-0.05cm}i}^t$ belongs to $\mathbf Q_{\sf sf}(u)=\mathbf Q_{\sf sf}(u_1,\ldots,u_n)$  for any $i=1,\ldots,n$ and seed vertex $t$.\footnote{Here  $\mathbf Q_{\sf sf}(u)$ stands for the 
 universal semifield  of {\it `substraction free'}  rational functions in $n$ indeterminates $u_1,\ldots,u_n$ with coefficients in the field of rational numbers $\mathbf Q$. We recall that a rational function $F\in {\mathbf Q}(u)$ is `substraction free' if it can be written $F=P/Q$ with $P,Q\in \mathbf Q_{>0}[u]$. For instance, $G=
{(x^2-x+y+1)}/{(x^2-x+1)}$ belongs to $ \mathbf Q_{\sf sf}(x,y)$ since it can also be written  $G= {(y(1+x)+x^3+1)}/{(x^3+1)}$ as a rational function.} Actually, much more is true
 since it can be proved that each $F_{\hspace{-0.05cm}i}^t$ actually is a polynomial with integer coefficients, that is 
 $F_{\hspace{-0.05cm}i}^t\in \mathbf Z[u_1,\ldots,u_n]$ (see \cite[Proposition\,3.6]{FZIV}).\sk 
 
  For any $t$, the $F_{\hspace{-0.08cm}i}^t$'s  with $i=1,\ldots,n$ are called the {\bf ${\boldsymbol{ F}}$-polynomials} of the seed $S_t$. We will see below that the $F$-polynomials of a cluster algebra actually enjoy quite stronger properties.  
The interest of considering the $F$-polynomials $F_{\hspace{-0.08cm}i}^t$  as well as the integer matrices $C^t$ is made clear once considering the following fundamental result: 
 \begin{thm}[{\bf Separation formula for $\boldsymbol{\mathcal X}$-cluster variables, \cite[Prop.\,3.13]{FZIV}}] 
 \label{Thm:Separation-Formula}
 For any $i=1,\ldots,n$ and  any $t\in \mathbb T_n$, the corresponding 
 $\boldsymbol{\mathcal X}$-cluster variable $x_{i}^t$ 
  expresses in terms of the initial cluster variables $x_1,\ldots,x_n$ by means of the following formula:
 \begin{equation}
 \label{Eq:X-var:separation-formula}
 x_{i}^t=\left(\prod_{j=1}^n x_j^{c^t_{\hspace{-0.05cm}ji}}\right) 
\left( \prod_{k=1}^n F_{\hspace{-0.08cm}k}^t({x})^{b_{\hspace{-0.05cm}ki}^t}
\right)
\, . 
\end{equation}
 \end{thm}
 The matrix $C^t$ is the 
 {\bf $\boldsymbol{C}$-matrix} of the seed $S^t$ and by definition, its $i$-th column vector $c_i^t=\big[c_{1i},\ldots,c_{ni}\big]^t$ is the {\bf $\boldsymbol{c}$-vector} of the cluster variable $x_i^t$. 
 \mk

  There is also a  'separation formula' for the $\mathcal A$-cluster variable (see  \cite[Cor\,6.3]{FZIV}):  there exists $G^t=(g_{ij}^t)_{i,j=1}^n\in M_n(\mathbf Z)$ such that, 
  if $\hat{a}=(\hat{a}_i)_{i=1}^n$ with $\hat{a}_i=\prod_{j=1}^n  a_{j}^{b_{ji}}$ for any $i=1,\ldots, n$, one has 
 \begin{equation}
 \label{Eq:A-var:separation-formula}
  a_i^t= \left(\prod_{j=1}^n a_j^{g^t_{\hspace{-0.05cm}ji}}\right) \cdot 
   F_{\hspace{-0.08cm}i}^t(\hat{a})\, .
\end{equation}
  \mk

 The notions of $F$-polynomial and of $c$-vector are fundamental in the theory of cluster algebras and have been intensively studied in recent times. Some deep conjectures about them have been proved recently by means of deep and conceptual approaches (categorification of cluster algebras, scattering diagrams). Among the most important results we are referring to here, let us mention the two following theorems:

 \begin{thm}[{\bf Positivity for $\boldsymbol{F}$-polynomials}]
 \label{Thm:Positivity}
 For  $t\in \mathbb T_n$ and  $i=1,\ldots, n$, the 
 polynomial $F_{\hspace{-0.08cm}i}^t(u)$:
 \begin{enumerate}
  \item has constant term equal to 1, that is $F_{\hspace{-0.08cm}i}^t(0)=1$;
 \sk 
  \item has only one monomial  of maximal degree and the latter has coefficient 1 and  is   
a (proper monomial) multiple of any other   monomial appearing in $F_{\hspace{-0.08cm}i}^t(u)$;
 \sk 
 \item has  positive integer coefficients, 
 {\it i.e.} $F_{\hspace{-0.08cm}i}^t(u)\in \mathbf Z_{>0}[u_1,\ldots,u_n]$.
 \sk  
 \end{enumerate}
 \end{thm}
 (The first two  statements in the preceding theorem have been conjectured in \cite[\S5]{FZIV} where it is proven that they actually are equivalent. 
In finite type, one can be much more precise about the monomial of maximal degree of the second point, see the paragraph just before Remark \ref{Rem:coco}.)
 
 \begin{thm}[{\bf Sign-coherence for $\boldsymbol{c}$-vector}]
 \label{Thm:Sign-Coherence}
 For any $t\in \mathbb T_n$ and  any $i=1,\ldots, n$, the 
 $c$-vector $c_i^t$  is non-zero
and has its non-zero entries either all positive or all negative. 
\end{thm}

We mention that the presentation above does not reflect accurately how the different statements appearing in these two fundamental theorems are related. 
Indeed, each of the two last points of Theorem \ref{Thm:Positivity} actually is equivalent to Theorem \ref{Thm:Sign-Coherence} ({\it cf.}\,Proposition 5.6 in 
\cite{FZIV}).  Moreover, despite their importance regarding the theory of cluster algebras, we will  not really use these results in the current text, but only a few easy, but already important for our purpose, consequences.\mk 

$F$-polynomials  are important objects attached to cluster algebras and they have been the subject of various studies from different points of view ({\it e.g.}\,see the references given in \cite[\S1]{CK}). When studying cluster webs, and more specifically their polylogarithmic ARs (see \ref{SS:Some-conjectures} further for more about this), it becomes relevant to know whether  the hypersurfaces $H_F=\{ F_{{}^{}\hspace{-0.08cm}i}^t(u)=0\}$ cut out by a non constant $F$-polynomial is irreducible or not.  The consideration of many concrete examples led us to conjecture that it is always the case. Motivated by a question of ours, Cao and Keller prove that it indeed holds true in the recent preprint \cite{CK}: 
\begin{thm}
 \label{Thm:CK}
Each non constant $F$-polynomial $F_{{}^{}\hspace{-0.08cm}i}^t$ 
 is irreducible in $\mathbf Z[u_1,\ldots,u_n]$.
 \end{thm}
\mk

A remarkable feature of the results presented above is that they hold true for any choice of initial seed.  This highly non trivial fact has very pleasant/useful  consequence for the cluster webs we are going to study. 
\begin{center}
$\star$
\end{center}

 By definition, the {\bf $\boldsymbol{F}$-polynomials associated} to 
 a given a cluster variable $x_i^t$
are the $F_{\hspace{-0.05cm}j}^t$ in  \eqref{Eq:X-var:separation-formula} appearing with a non-vanishing exponent $b_{ij}^t$.  Given a cluster web $\boldsymbol{\mathcal W}_\Sigma$ defined by 
a finite subset of $\boldsymbol{\mathcal X}$-cluster variables $\Sigma$, ones denotes by 
\begin{equation}
\label{Eq:F(W)}
\boldsymbol{F}\Big(\boldsymbol{\mathcal W}_\Sigma\Big)
 \end{equation}
 (or just by $\boldsymbol{F}\big(\Sigma\big)$) 
  the union of the non-constant $F$-polynomials associated to the elements of $\Sigma$. Remark that they are polynomials (with positive integers as coefficients) in the cluster variables $x_1,\ldots,x_n$ of the initial seed.   \mk

Following \cite{NakanishiTrop} (where this notion has been formally introduced for the first time and to which we refer for more details and perspective), we define  the {\bf $\boldsymbol{i}$-th tropical sign} $\boldsymbol{\varepsilon_i^t}\in \{-1,+1\}$
 of the cluster associated to a vertex $t$ in the exchange graph of a cluster algebra, as the sign of the non-zero coefficients of the $c$-vector $c_i^t$. Note that, in this case again, this notion a priori depends on the choice of the initial seed. In the present text, 
 we will use it to formulate some dilogarithmic identities associated to cluster algebras (for instance, see identity  $\big(\boldsymbol{{\sf R}_{\boldsymbol{i}}}^{\hspace{-0.05cm}\epsilon}\big)$ p.\,\pageref{kooko} further). 
\begin{center}
$\star$
\end{center}

As an example, we consider the case of finite type $A_2$ (see \cite{FZII} and the different tables in \cite{FZIV}).  In Figure 
\ref{Fig:FCC-ExchangeGraph-A2} below, we give the associated {\it `${FC}$-exchange graph'} whose vertices $\boldsymbol{FC}^t$ corresponding to a vertex $t$ of $\mathbb T_2$ is $({F}^t,C^t)$ where ${F}^t=(F_1^t,F_2^t)$ is the pair of associated $F$-polynomials and  where $C^t$ stands for the matrix of $c$-vectors. We note $\boldsymbol{FC}^0$ the initial `$FC$-seed' and 
for any $i,j$, two distinct elements of $\{1,2\}$, we set $\boldsymbol{FC}^{i}$ and $\boldsymbol{FC}^{ij}$ for the `$FC$-seeds' after mutation $\mu_i$ and $\mu_j\circ \mu_i$ respectively).

 
  \begin{figure}[!h]
   \begin{center}
  \scalebox{0.9}{
    \xymatrix@R=0.9cm@C=-1.5cm{    
    &  {{\boldsymbol{FC}^0 = \Bigg( \big(1,1\big)\, , 
\scalebox{0.9}{ 
$\begin{bmatrix}
 1 & 0 \\
 0 & 1 
 \end{bmatrix}$}\, \Bigg)}}
 \ar@{->}[rd]^{ \hspace{0.1cm}\mu_1} 
  \ar@{->}[dl]_{ \hspace{0.1cm} \mu_2} &  \\
 \boldsymbol{FC}^2 =  \Bigg( \big(\, 1 , 
  \, {1+x_2}\big)\, , 
\scalebox{0.9}{ 
$\begin{bmatrix}
 1 & 0 \\
 0 & -1 
 \end{bmatrix}$} \, \Bigg)  \ar@{->}[dd]_{ \hspace{0.1cm}\mu_1}     
      &    & 
     \ar@{->}[dd]^{ \hspace{0.1cm}\mu_2} 
{}^{} \hspace{-1cm}   \boldsymbol{FC}^1 =
        \Bigg( \, \big( {1+x_1}\, ,\,  1\big)\, , 
\scalebox{0.9}{ 
$\begin{bmatrix}
 -1 & 1 \\
 0 & 1 
 \end{bmatrix}$} \, \Bigg)\\ 
 & & 
 \\
 \boldsymbol{FC}^{12}= \Bigg( \Big({1+x_1+x_1x_2}\, , \,  {1+x_2}\Big)\, , 
\scalebox{0.9}{ 
$\begin{bmatrix}
 -1 & 0 \\
 0 & -1 
 \end{bmatrix}$}\,  \Bigg) 
  \ar@{-->}[dr]_{ \hspace{0.0cm} (12)} 
  &      &  
 \boldsymbol{FC}^{21}=   \Bigg( \Big(\,  1+x_1 \, , \, 1+x_1+x_1x_2\Big)\, , 
\scalebox{0.9}{ 
$\begin{bmatrix}
 0 & -1 \\
 1 & -1 
 \end{bmatrix}$} \, \Bigg)  \ar[dl]^{\mu_1} \\
 &     \boldsymbol{FC}^{121} =   \Bigg( \Big(1+x_2 \, , \, 1+x_1+x_1x_2\Big)\, , 
\scalebox{0.9}{ 
$\begin{bmatrix}
 0 & -1 \\
 -1 & 0 
 \end{bmatrix}$} \, \Bigg)  &
 }}
 \end{center}
 \caption{The $FC$-exchange graph of the cluster algebra of type $A_2$ (see also \cite[Table\,2]{FZIV}).}
 \label{Fig:FCC-ExchangeGraph-A2}
 \end{figure}
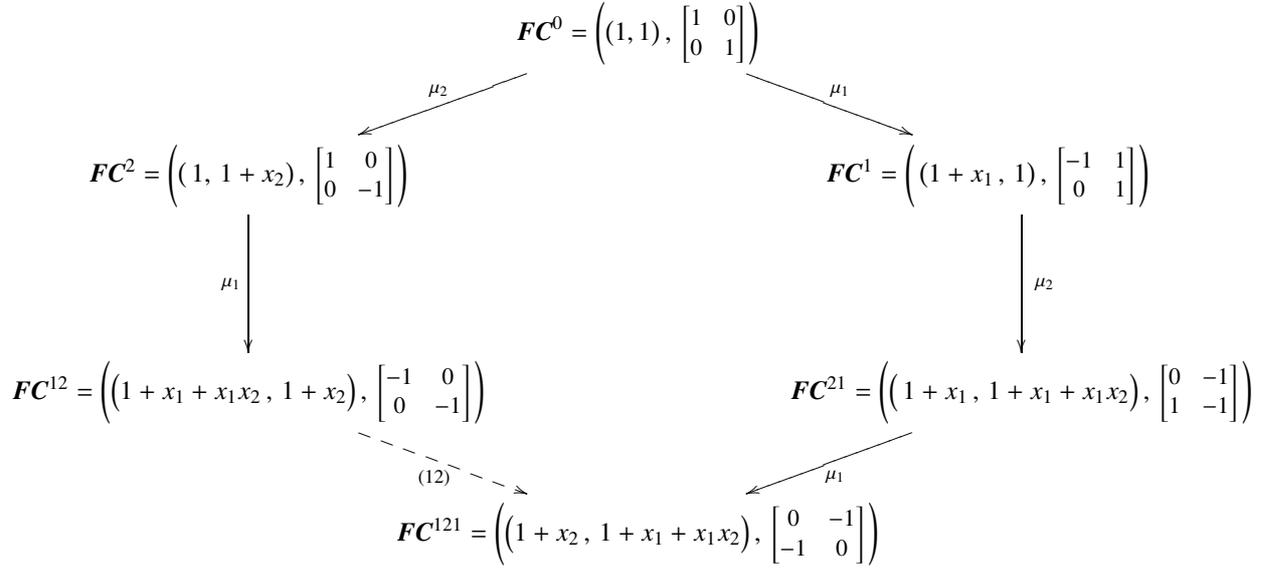

 \mk

From Figure \ref{Fig:FCC-ExchangeGraph-A2}, it follows 
that 
 $$ \boldsymbol{F}(A_2)=\big\{ \, F_\alpha \, \lvert \, \alpha\in (A_2)_{>0} \,  \big\}=
 \Big\{ 1+x_1,1+x_2,1+x_1+x_1x_2\Big\}\, .$$  
   If the expressions for the elements of $\boldsymbol{F}(A_2)$ are quite simple, we remark that there is a lack of symmetry since exchanging 
    $x_1$ with $x_2$ in the   $F$-polynomial $1+x_1+x_1x_2$ gives a polynomial which does not belong to  $\boldsymbol{F}(A_2)$. \mk 
   
   We will see just below that, up to inverting some of the initial cluster variables, one obtains some variants of the $F$-polynomials which, in case of cluster algebras of  finite type, have a nice Lie-theoretic and combinatorial description.

\subparagraph{Some features of $\boldsymbol{F}$-polynomials in finite type.}
\label{SubPar:Fpoly-finite-type}
 Since we are going to look at the cluster webs  associated to finite type cluster algebras  
in a significant way, we describe below some particular features of $F$-polynomials in this case that  will prove to be useful later on. \mk

Because  finite type cluster algebras are intimately related to Lie theory, it is not surprising that many things about the 
$F$-polynomials of a cluster algebra associated to a given  Dynkin diagram 
$\Delta$ can be described in terms of the root system 
$\Phi=\Phi(\Delta)$  
associated to the latter. 
\sk

    In what follows, we will use the following notations:
    \begin{itemize}
 \item  $n$ stands  for the rank of $\Delta$; 
 \item a bit abusively, we write again $\Delta$ for the set of simple roots 
 $\delta_1,\ldots,\delta_n$ of $\Phi$;
 \item  then $\Phi_{>0}$ denotes the set of positive roots, that is the intersection of $\mathbf Z_{>0}\Delta$ with $\Phi$,  and  $\Phi_{\geq -1}$ stands  the set of {\bf almost-positive roots}, {\it i.e.} 
  $$ 
  \Phi_{\geq -1}=\Phi_{>0} \sqcup \big( - \Pi \big)\, .
  $$
We will frequently denote $\Phi_{>0}$ (resp.\,$\Phi_{\geq -1}$) by $\Delta_{>0}$ or 
$\Delta^+$ (resp.\,by $\Delta_{\geq -1}$);
\end{itemize}
\sk 
 
Most of the facts recalled below are satisfied `up to mutations', that is  hold true for any choice of an exchange matrix 
$B$ mutation equivalent to the bipartite matrix $B_\Delta=B_{\vec{\Delta}}$ defined above.\footnote{A classical example of a matrix mutation equivalent to $B_\Delta$  is the matrix $B(\Delta)=(b_{ij})_{i,j=1}^n\in {\rm M}_n(\mathbf Z)$ defined as follows:  if $C=C_{\Delta}=(c_{ij})_{i,j=1}^n$ is the Cartan matrix associated to $\Delta$, then  for all $i,j=1,\ldots,n$ one has $b_{ii}=0$ and $b_{ij}=\epsilon_{ij} c_{ij}$ when $i\neq j$, with $\epsilon_{ij}=1$ if $i>j$, and $\epsilon_{ij}=-1$ otherwise}.  
   For simplicity, we assume that $B$ stands precisely for  $B_\Delta$ in what follows. 
\sk

 In order to  remind the definition of the $F$-polynomials of the cluster algebra associated to $B$ given in \cite[\S3]{FZIV}, 
one needs to recall that, as explained in  \cite{FZIV}  (and in many other places as well),  one can consider cluster algebras with coefficients in any  semi-field $\mathbb P$, the case when the latter is trivial (that is $\mathbb P=\{1\}$) corresponding to the case considered in the present text.   An important case is that of {\it `principal coefficients'}, that is when $\mathbb P$ coincides (up to isomorphism) with the `tropical semi-field'
${\rm Trop}(x)={\rm Trop}(x_\delta, \delta\in \Delta)$ ({\it cf.}\,Definition 2.2 in 
\cite{FZIV}).  Indeed, the cluster variables of the associated cluster algebra, denoted by $A_{B}^{\rm prin}$, can be used to define the $F$-polynomials and also to give a separation formula (analogous to \eqref{Eq:A-var:separation-formula}) for the 
$\boldsymbol{\mathcal A}$-cluster variables 
of the cluster algebra $A_B^{\mathbb P}$ associated to any coefficient semi-field ({\it cf.}\,\cite[Thm.\,3.7]{FZIV}).\mk

Let $(a_\delta)_{\delta\in \Delta}$ be the cluster of the initial $\boldsymbol{\mathcal A}$-seed of 
$A_{B}^{\rm prin}$. 
%
%
%
For any $\alpha=\sum_{\delta\in \Delta } \alpha_\delta \delta$ in $\Phi$ (with $\alpha_\delta \in \mathbf Z$ for every $\delta\in \Delta$), one sets $a^\alpha=\prod_\delta a_\delta^{\alpha_\delta}$.  
The $\boldsymbol{\mathcal A}$-cluster variable of $A_B^{\rm prin}$
associated to any pair $(t,\lambda)\in \mathbb T^n\times \Delta$ is denoted by $A_{t,\lambda}$. By construction, it is a substraction-free element of 
the ring $ \mathbf Q\mathbb P(a)$ of polynomials in the $a_\delta$'s with coefficients in the fraction field of $ \mathbf Z\mathbb P$. Since 
this fraction field $ \mathbf Q\mathbb P$ naturally identifies with $\mathbf Q_{\rm sf}(x_\delta,\delta \in \Delta)$, any cluster variable can  be seen as a substraction-free rational function in the $a_\delta$'s and the $x_\delta$'s:  one has
$$
A_{t,\lambda}\in  \mathbf Q\mathbb P(a)
\subset \mathbf Q_{\rm sf}(a,x)\, .
$$

Then by definition,  the {\bf $\boldsymbol{F}$-polynomial} $\boldsymbol{F_{t,\lambda}}$ associated to the pair $(t,\lambda)$ is the element of $\mathbf Q_{\rm sf}(y)$ obtained by specializing $A_{t,\lambda}$ by setting  $a_\delta=1$ for every $\delta$: one has 
$$
F_{t,\lambda}=
\mathbf S(A_{t,\lambda })
$$
where $\mathbf S$ stands for the specialization map 
$\mathbf Q_{\rm sf}(a,x)\longrightarrow \mathbf Q_{\rm sf}(x)$,
$R(a,x)\longmapsto 
R(a=1,x)
 $.
In other terms : the previous map induces a surjection from the set of $\boldsymbol{\mathcal A}$-cluster variables of $A_B^{\rm prin}$ onto the set of associated $F$-polynomials.   On the  other hand, thanks to a fundamental  result (namely, Theorem 1.9 of \cite{FZII}), it is known that  
 there exists a (unique) bijection 
 \begin{equation}
 \label{Eq:Phi-Avar-bij}
\Phi_{\geq -1}\longrightarrow 
 \boldsymbol{{\mathcal A}var}\big(A_B^{\rm prin}\big), \, \alpha\longmapsto A[{\alpha}]
 \end{equation}
such that, in terms of the initial cluster variables $a_\delta$, 
 one has 
$
A[\alpha]= 
{P_\alpha(a)}/{a^\alpha}\, , 
$  where  $P_\alpha$ is an element of $\mathbf Z\mathbb P[a]$
with non-trivial constant term (moreover, under this
bijection, one has $A[-\delta] = a_\delta$   for any $-\delta \in (-\Pi)=\Phi_{\geq -1}\setminus \Phi_{>0}$).  
 Composing $\boldsymbol{S}$ with \eqref{Eq:Phi-Avar-bij}, one gets  a surjective map between the set of almost-positive roots and that of $F$-polynomials of the cluster algebra considered: 
  \begin{equation}
 \label{Eq:Phi-Fpolynomials}
\Phi_{>0} \longrightarrow 
 \boldsymbol{F}\big(A_B\big)\, , \, \alpha\longmapsto f[{\alpha}]\, .
 \end{equation}

 Finally, it is known that for any $\alpha=\sum_{\delta \in \Delta}a_\delta \delta \in \Phi_{>0}$, the monomial of highest degree of $f[\alpha]$ is precisely 
 $x^\alpha=\prod_\delta (x_\delta)^{\alpha_\delta}$.\footnote{This is a non trivial result. It is known that any cluster algebra $\mathcal A_{\Delta}$ of finite type admits a categorification $\mathscr C_{\Delta}$. This allows to express any $F$-polynomial $f[\alpha]$ 
as the value taken by the so-called  {\it `Caldero-Chapoton character'} when evaluated at a specific object of $\mathscr C_{\Delta}$. That the highest degree monomial of $f[\alpha]$ is $x^\alpha$ follows at once from this.  See 
\cite{KellerCat} for a general overview of the interplay between cluster algebras and their categorifications.} From this, it comes immediately that \eqref{Eq:Phi-Fpolynomials} is bijective.

 %
%

 \begin{rem}
 \label{Rem:coco}
The exchange matrix $B$ being given, 
the bijection  \eqref{Eq:Phi-Avar-bij} actually exists for any cluster algebra of `{\it geometric type}', that is for any cluster algebras $A_B^{\mathbb P}$ with coefficients in a tropical semi-field ${\mathbb P}$. The case when the latter is trivial (that is $\mathbb P=\{1\}$) corresponds to the cluster algebras considered here and to which we associate some webs. Consequently, one deduces from this the following more explicit description of the $\boldsymbol{\mathcal A}$-cluster web associated to 
$\Delta$:  one has 
\begin{align*}
\boldsymbol{\mathcal A} \boldsymbol{\mathcal W}_\Delta=& \, 
\, \boldsymbol{\mathcal W}
\left(  \, a_\alpha \, \big\lvert \hspace{0.05cm}  \alpha\in \Phi_{\geq -1}\, 
\right)
=
\boldsymbol{\mathcal W}
\left(\,  a_\delta \, ,   
{P_\alpha(a)}/{a^\alpha} 
\hspace{0.1cm}
 \Big\lvert \hspace{0.1cm}\delta \in \Delta\, , \, \alpha\in \Phi_{>0}\, 
\right)\, .  
\end{align*}
 \end{rem}

 For  the case of $A_2$    considered above ({\it i.e.}\,$B=B_{A_2}$), if one denotes by $\alpha_i$ the principal root attached to the $i$th vertex of the Dynkin diagram 
  \scalebox{0.7}{\begin{dynkinDiagram}[edge length=.5cm,labels*={1,...,2}]{A}{2}\end{dynkinDiagram}}, then the polynomial $f[\alpha]$s associated to 
  the  positive roots are given in the second column of Table
  \ref{T:FpolyA2} 
   below (extracted from Table 7 of \cite{FZIV}).  \sk 
   
  \begin{table}[h]
\begin{center}
\begin{tabular}{|l|l|l|}
\hline 
  \begin{tabular}{c}   \vspace{-0.15cm} \\
$\boldsymbol{\alpha \in A_2^+}$ 
 \vspace{0.2cm} \\ 
\end{tabular}
& 
 \hspace{0.3cm} $\boldsymbol{f[\alpha]}$ 
 &  \hspace{0.3cm} $\boldsymbol{F[\alpha]}$ 
 \\
\hline
\hspace{0.3cm}
 \begin{tabular}{c}   \vspace{-0.2cm} \\
$\alpha_1$ 
 \vspace{0.2cm} \\ 
\end{tabular}
& \hspace{0.1cm} $1+x_1$ & \hspace{0.1cm} $1+x_1$
\\
\hline 
\hspace{0.3cm} 
 \begin{tabular}{c}   \vspace{-0.2cm} \\
$\alpha_2$ 
 \vspace{0.2cm} \\ 
\end{tabular}
 & \hspace{0.2cm}$1+x_2$ & \hspace{0.1cm} $1+x_2$
\\
\hline
\hspace{0.3cm} 
\begin{tabular}{c}   \vspace{-0.2cm} \\
$\alpha_1+\alpha_2$  \vspace{0.2cm} \\ 
\end{tabular}
&  \hspace{0.1cm}
$1+x_1+x_1x_2$ 
  &  \hspace{0.1cm}
$1+x_1+x_2$ 
\\ \hline
\end{tabular}  
\caption{$F$-polynomials in type $A_2$}
\label{T:FpolyA2}
\end{center}
\end{table}
%
%
%
%
%
  We see that the polynomials $f[\alpha]$s are not symmetric in the variables  $x_1$ and $x_2$, 
   whilst everything regarding the roots is.     This lack of symmetry 
for the $f[\alpha]$'s  holds true for any Dynkin diagram and 
tends to make artificially a bit too complicated  the formulas for the cluster variables, which  can be the source of difficulties (for instance when performing explicit computations in order to determine the ARs of a given cluster web).
 \sk 
 
 There is a way to get more symmetric formulas for some variants of the polynomials $f[\alpha]$'s, via a rational change of variables described in \cite[\S11]{FZIV} ({\it cf.}\,Definition 11.5 therein).  
Let $E$ be the automorphism $E$ of $\mathbf Q(x_1,\ldots,x_n)$ characterized by the following relations: 
 \begin{equation}
 \label{Eq:Bir-change-E}
 E(x_i)=\begin{cases}
 \hspace{0.18cm}
 {x_i}^{-1}
 \hspace{0.3cm} \mbox{if }\, i\hspace{0.1cm}
 \mbox{is a source in }\, \vec{\Delta}; \vspace{0.15cm}\\
  \hspace{0.18cm} x_i\hspace{0.5cm} \mbox{if }\, i\hspace{0.1cm}
 \mbox{is a sink  in }\, \vec{\Delta}\,.
  \end{cases} 
  \end{equation}
   
 Then it follows from \cite[\S11]{FZIV}, that setting for any almost-positive root $\alpha\in \Delta_{\geq -1}$: 
 $$
 F[\alpha]=E\big(f_\alpha(x)\big)\prod_{
 \substack{
 i\,{source}
 }
  }
  {x_i}^{[\alpha_i]_+}\, , 
 $$
 one recovers the $F$-polynomials introduced in \cite{FZ}.
 \sk

 What makes the $F[\alpha]$s interesting compared to the $f[\alpha]$s is 
first that the former behave  in a better way than the latter from a Lie theoretic point of view 
 (see Theorem 1.6 or Proposition 2.9 of \cite{FZ} for instance).  
Second, the $F[\alpha]$s are formally simpler than the $f[\alpha]$s\footnote{For a simple but  concrete example, compare the last two columns of Table \ref{T:FpolyA2} above.}, which makes them more useful when performing some computations for instance. In classical simply-laced type ($A$ and $D$), there even exists a closed  (but complicated) formula for $F[\alpha]$ in terms of some combinatorial objects associated to any positive root $\alpha$ ({\it cf.}\,Proposition 2.10 in \cite{FZ}).\sk 

\begin{exm}
 In type $A$, everything is particularly elementary and there are nice descriptions of the $F$-polynomials:  a closed combinatorial one, 
and another  inductive one, which is simpler and which, moreover, has the merit of explaining the terminology.  
For $n\geq 1$ fixed, one denotes by $\alpha_1,\ldots,\alpha_n$ the principal roots associated in the following way to the nodes of 
 the Dynkin diagram: 
\begin{center}
$\Delta(A_n)$\hspace{0.15cm}: \hspace{0.3cm}
\begin{tabular}{c}
\vspace{-0.4cm}
 \begin{dynkinDiagram}[text style={scale=1.3,blue},
edge length=.75cm,
labels={1,2,n-1,n},
label macro/.code={\alpha_{\drlap{#1}}}
]{A}{}
\end{dynkinDiagram}
\end{tabular}
\, .
 \end{center}
 
 It is well known that, 
for  any non-trivial  interval $[a,b]\subset [1,n]$, the sum 
 $\alpha_{[a,b]}= \sum_{a \leq k \leq b} \alpha_k$ is a positive root and conversely, any positive root is written in this way. For any sub-interval as above, one denotes by $F[a,b]$ the $F$-polynomial associated to 
 $\alpha_{[a,b]}\in (A_n)_{>0}$.  Then, according to Example 2.15 of \cite{FZ}, writing $x_S=\prod_{s\in S}x_s$ for any $S\subset [1,n]$, one has  
 $$
 F[a,b]=\sum_{\Omega \subset [a,b]} x_\Omega\, 
 $$
 where the sum is carried over all the totally disconnected subsets $\Omega$ of $[a,b]$. \sk 
 
 For an inductive formula for $F[a,b]$, we recall the definition of what 
it is natural to call  the `Fibonacci's polynomials' $(F_m)_{m\geq -1}$:
one sets $F_{-1}=1$,   $F_0=1$  and
  $$
 F_m(s_1,\ldots,s_m)=F_{m-1}(s_1,\ldots,s_{m-1})
+ F_{m-2}(s_1,\ldots,s_{m-2}) \, s_m
  \quad \mbox{for } \, m\geq 1\, 
 $$ 
  (thus 
 $ F_1(s_1)=1+s_1$, $F_2(s_1,s_2)=1+s_1+s_2$, 
 $F_3(s_1,s_2,s_3)=1+s_1+s_2+s_3+s_1s_3$, etc.) \sk 
 
 For each $m\geq 1$, $F_m$ is a polynomial in $m$ variables with positive integers as coefficients and constant term equal to 1. And for any $n$ and any subinterval $[a,b]\subset  
 [1,n]$, it can be verified that,  as a polynomial in the variables $x_1,\ldots,x_n$, one has 
 $$
 F[a,b](x_1,\ldots,x_n)=F_{b-a+1}(x_a,\ldots,x_b)\, .
 $$
\end{exm}

\begin{rem}
\label{Rem:birat-E}
Since any $\boldsymbol{\mathcal X}$-cluster variable $x$ is expressed in  terms 
of the initial cluster variables and of the $F$-polynomials in a rather simple way (see \eqref{Eq:X-var:separation-formula}), and because the polynomials $F[\alpha]$s are given by simpler formulas than the $f[\alpha]$s, it is more convenient in practice to consider $x\circ E$ instead of $x$, 
where $E$ stands for the rational change of coordinates \eqref{Eq:Bir-change-E}.   Since $E$ is  a birational involution, this doesn't change anything  when considering cluster webs from a geometric point of view.
\end{rem}

\paragraph{Cluster tori and cluster webs.}
\label{Par:Cluster-tori}
For more on the material discussed in this paragraph and in the next one, the reader can consult 
  \cite{FockGoncharovIHES}, 
 \cite[\S2]{FockGoncharovENS} and the fifth and sixth chapters of \cite{GSV}.\mk

Let $\boldsymbol{A}_B$ be a cluster algebra of rank $n$ with fixed initial seed
$(\boldsymbol{a},\boldsymbol{x},B)$ associated to a vertice $t_0$ of $\mathbb T^n$. Below, we use $\boldsymbol{\mathcal Z}$ to denote either $\boldsymbol{\mathcal A}$ or $\boldsymbol{\mathcal X}$ since our discussion covers both cases. \mk 

To a $\boldsymbol{\mathcal Z}$-cluster $\boldsymbol{z}^t=(z_1^t,\ldots,z_n^t)$ corresponding to a
vertex $t$ of $\mathbb T_n$,   one associates the {\bf $\boldsymbol{\mathcal Z}$-cluster torus $\boldsymbol{\mathbf T^{\boldsymbol{z}^t}}$} which is the complex torus defined by
$$  
 \mathbf T^{\boldsymbol{z}^t}={\rm Spec} \left(\, \mathbf C\Big[(z_1^t)^{\pm 1},\ldots,(z_n^t)^{\pm 1}\Big]\, \right) \simeq \big( \mathbf C^*\big)^n\,
$$
and that we will also denote sometimes by $ \boldsymbol{\mathcal Z}\mathbf T^t$ (or even by $\mathbf T^t$ if no ambiguity is possible) for simplicity.  We will denote by $ \mathbf T^{\boldsymbol{z}^t}_{>0}$ or 
$ \boldsymbol{\mathcal Z}\mathbf T^t_{>0}$ its positive part, that is its points with values 
in the positive semi-field $\mathbf R_{>0}$, that is: $
 \mathbf T^{\boldsymbol{z}^t}_{>0}= \mathbf T^{\boldsymbol{z}^t}(\mathbf R_{>0})\simeq (\mathbf R_{>0})^n\subset (\mathbf C^*)^n$. 
\mk

Since $\boldsymbol{\mathcal Z}$ is fixed, we drop it most of the time in what follows. 
 Let $\boldsymbol{\mathcal W}$  be a cluster web 
 defined by cluster variables $x_\sigma$ indexed by a finite set $\Sigma$.   
We will see it  as the 
 web on the initial cluster tori $ \mathbf T^{t_0}$ 
 whose first integrals are the $x_\sigma$'s seen as rational functions on it. It is not difficult to see that  this definition does not really depend on the choice of the initial seed. \sk 
 
 Indeed, let $\mu_P:  \mathbf T^{t_0} \dashrightarrow \mathbf T^{t}$ be the rational map of tori given by the composition of the mutations associated to the edges of a path $P: 
  t_0\rightarrow t_1 \rightarrow \cdots \rightarrow t_{n-1}\rightarrow t_n=t$ 
 in $\mathbb T^n$ joining the initial vertex $t_0$ to another one $t$. 
 Then it is
 a birational map whose inverse is the rational map  $\mu_{P^{-1}}:\mathbf T^{t} \dashrightarrow \mathbf T^{t_0}$ associated to 
 the inverse path $P^{-1}
 : t=t_n \rightarrow \cdots  \rightarrow t_{0}$. Moreover, since  $\mu_P$ and $\mu_{P^{-1}}$ are composed by mutations, both are positive maps hence are defined at any point of $\mathbf T^{t_0}_{>0}$ and  $\mathbf T^{t}_{>0}$ respectively. 
 Since $\mu_{P^{-1}}\circ \mu_P={\rm Id}$, it follows that $\mu_P$ induces a real-analytic diffeomorphism from  $\mathbf T^{t_0}_{>0}$ onto $\mathbf T^{t}_{>0}$. Now, let  $(\zeta,\xi)$ be a pair   of $\boldsymbol{\mathcal Z}$-cluster variables viewed as rational functions on $  \mathbf T^{t_0}$ and denote by $\zeta^t$ and $\xi^t $ the corresponding cluster variables but seen as rational functions on $\mathbf T^{t}$: then $\zeta=\mu_P^*(\zeta^t)$ and similarly 
  $\xi=\mu_P^*(\xi^t)$. Since $\mu_P$ is birational, it follows that 
  the rational 2-form $d\zeta\wedge d\xi$ on $\mathbf T^{t_0}$  vanishes identically if and only if the same holds true for $d\zeta^t\wedge d\xi^t$ on $\mathbf T^{t}$. Equivalently:  the two foliations $\mathcal F_{\zeta}$ and 
  $\mathcal F_{\xi}$ coincide on $\mathbf T^{t_0}$ if and only if the same holds true for $\mathcal F_{\zeta^t}$ and 
  $\mathcal F_{\xi^t}$ on $\mathbf T^{t}$. \sk 
  
  From the discussion above, one immediately gets the 
  \begin{prop} 
  \label{P:The notion of cluster web does not depend on the choice of an initial seed}
  {1.}  
  The notion of cluster web does not depend on the choice of an initial seed. \sk
  
  {2.} The foliations of a cluster web are globally defined on the positive part  
   of any cluster tori.\footnote{However, it may not be excluded that (some of) these foliations have singularities on some positive part $\mathbf T^{t}_{>0}$.}
  \end{prop}
 
Since we are interested only in studying cluster webs up to equivalence, we can  consider them only on the initial cluster torus $\mathbf T^{t_0}$, and even to consider only their real versions given by considering their restrictions on the positive part $\mathbf T^{t_0}_{>0}=\mathbf T^{t_0}(\mathbf R_{>0})\simeq (\mathbf R_{>0})^n$ since it appears that  both share exactly the same invariants (their virtual, polylogarithmic, or total rank in particular). \sk

A notable feature of cluster  webs is that each one of the foliations composing such a web is defined by a cluster variable, which in particular is a substraction-free hence a positive rational function on $\mathbf T^{t_0}$. A natural question is whether  such a first integral obtained  through the inductive algebraic machinery of cluster algebras are also primitive in the sense of foliation theory, that is have connected generic fiber. We will address  this  point further on in \S\ref{Par:ClusterVarNotComposite}.

\paragraph{Cluster varieties, cluster ensembles and $\boldsymbol{\mathcal U}$-cluster webs.}
\label{Par:Cluster-varieties-cluster-ensembles}
We now introduce the notion of cluster varieties,  which is of great importance since 
it  connects the theory of cluster algebras with many fields of geometry. \sk 

We continue to use the notation introduced in the preceding pargraph, with $\boldsymbol{\mathcal Z}$ denoting either $\boldsymbol{\mathcal A}$ or $\boldsymbol{\mathcal X}$ and  $B$ standing for the exchange matrix of the initial seed. One denotes by $D$ a skew-symmetrizer for  $B$: one has $D={\rm diag}(d_1,\ldots,d_n)$ with $d_1,\ldots,d_n\in \mathbf N^*$  such that $DB$ is skew-symmetric. We assume  that the $d_i$'s are relatively prime and 
  denote by $r$ 
   the rank of $B$.
\mk

Let 
$(t,t')$ be a pair of vertices in $\mathbb T^n$. Since the latter is a tree, there exists a unique injective continuous path $p_{t}^{t'}: [0,1]\rightarrow \mathbb T^n$ joining $t$ to $t'$ to which corresponds the composition, denoted by   
$\mu_{t}^{t'}$, of the mutations 
associated to the edges  of $\mathbb T^n$ encountered  when one goes along ${\rm Im}(p_t^{t'})$.   Then $\mu_{t}^{t'}$ is a positive birational map from $\mathbf T^t$ to $\mathbf T^{t'}$ with inverse $\mu_{t'}^t$. On the disjoint union of cluster tori $\mathbf T^t$, one defines an equivalence relation `$\boldsymbol{\mathcal Z}{\sf -mut}$'
by setting that $z\in \mathbf T^t$ and $z'\in \mathbf T^{t'}$ are equivalent if 
$\mu_{t}^{t'}$ is defined at $z$ and similarly for $\mu_{t'}^t$ at $z'$ and $z'=\mu_{t}^{t'}(z)$. Then one defines the 
{\bf $\boldsymbol{\mathcal Z}$-cluster variety $\boldsymbol{\mathcal Z}_B$} associated to $B$ as the quotient of $\sqcup_t \mathbf T^t$ by this equivalence relation: 
\begin{equation}
 \label{Eq:Cluster-variety}
\boldsymbol{\mathcal Z}_B=\Bigg(\,  \bigsqcup _{t\in \mathbb T^n}
\boldsymbol{\mathcal Z} \mathbf T^t\, \Bigg)\, \Big/_{\boldsymbol{\mathcal Z}{\sf -mut}}\, .
 \end{equation}
It can be proved that  $\boldsymbol{\mathcal Z}_B$ is a smooth rational positive manifold.  In particular, the positive parts $\mathbf T^t(\mathbf R_{>0})$ of the $\boldsymbol{\mathcal Z}$-cluster tori glue together and form the {\bf `positive part $\boldsymbol{\mathcal Z}_B(\mathbf R_{>0})$'} of the cluster variety 
whose  importance  in geometry is well known.  For us, the interest of considering this positive part of the cluster variety comes from the fact that it is topologically trivial (since it is isomorphic to any $\mathbf T^t_{>0}$, each 
of them being in turn isomorphic to the positive orthant $(\mathbf R_{>0})^n$) and that the restriction on it of any cluster variable is a real-analytic positive function (with no indeterminacy points in particular). 


\begin{warning} 
Actually, the definitions of the $\boldsymbol{\mathcal A}$- and 
$\boldsymbol{\mathcal X}$-cluster varieties associated to a cluster algebra are more subtile and much more involved than as presented just above.  Note also 
there are several variations of the constructions described above ({\it e.g.}\,the varieties $\boldsymbol{\mathcal A}^\star$ and $\boldsymbol{\mathcal X}^\star$
considered in {\rm \cite{MandelCox}} for instance). 
\sk

For a rigorous definition/construction in the realm of algebraic geometry, of the cluster varieties under scrutiny here, we refer to the second section of {\rm \cite{GrossHackingKeel}} and in particular to Proposition 2.4 therein. This proposition allows to give a rigorous meaning to \eqref{Eq:Cluster-variety}
from which it follows that  the varieties $\boldsymbol{\mathcal A}_B$ and $\boldsymbol{\mathcal X}_B$ are well-defined schemes. However,
it appears that these varieties, if well defined, can be a bit wilder than those usually considered: for instance, 
the $\boldsymbol{\mathcal X}$-cluster variety   is usually not separated, contrarily to the $\boldsymbol{\mathcal A}$-version, and it is not clear whether both are noetherian in general, see  Remark 2.6 of 
{\rm \cite{GrossHackingKeel}}. However, consequently to what is explained in the latter remark, this will not cause any particular problem when dealing with webs on these varieties, since in this case there is no loss of information by  working  locally, in particular,  on any cluster torus in a cluster variety.  This is what we are going to do systematically (but implicitly) in the sequel hence we will not elaborate further about any possible pathology (like non-separatedness)
of the cluster varieties we will deal with.
\end{warning}

To a cluster algebra ${\bf A}_B$, we thus have associated two cluster varieties of dimension $n$, namely $\boldsymbol{\mathcal A}_B$ and $\boldsymbol{\mathcal X}_B$. It turns out that these varieties carry supplementary structures.  Indeed, for 
any seed coresponding to a vertex $t$ of $\mathbb T_n$, one considers:
\begin{itemize}
\item[{\bf (i).}]
  the following monomial map from the $\boldsymbol{\mathcal A}$-cluster torus to the $\boldsymbol{\mathcal X}$-cluster torus: 
\begin{align*}
p^t : \boldsymbol{\mathcal A}\mathbf T^t     \longrightarrow 
\boldsymbol{\mathcal X}\mathbf T^t \,  , \quad 
(a_i)_{i=1}^n  \longmapsto \Big( \prod_{j=1}^n a_j^{b_{ji}} \Big)_{i=1}^n\, .
\end{align*}
\item[{\bf (ii).}] for any $n$-tuple of integers $w=(w_1,\ldots,w_n)$ whose transpose is in the kernel of $B^t$, one has a 1-dimensional torus action $\mathbf C^*\times 
\boldsymbol{\mathcal A}\mathbf T^t\rightarrow 
\boldsymbol{\mathcal A}\mathbf T^t :\, (t,(a_i)_{i=1}^n)\mapsto (t^{w_i} a_i)_{i=1}^n$. 
Summing up on a basis of ${\rm ker}(B^t)$, these actions give rise to a torus  action $ q^t : (\mathbf C^*)^{n-r}\times 
\boldsymbol{\mathcal A}\mathbf T^t\rightarrow\boldsymbol{\mathcal A}\mathbf T^t$. 
\item[{\bf (iii).}]  the closed holomorphic logarithmic 2-form $\Omega^t=\sum_{i,j=1}^n (d_i b_{ij}^t) \cdot \frac{da_i^t}{a_i^t} \wedge  \frac{da_j^t}{a_j^t} $ on $\boldsymbol{\mathcal A}\mathbf T^t$; 
\item[{\bf (iv).}] the Poisson bracket $\{\,\cdot \, , \, \cdot\,\}^t$ on $\boldsymbol{\mathcal X}\mathbf T^t$, associated to the bivector $\omega^t=\sum_{i,j=1}^n (b_{ij}^t /d_j)\cdot 
\frac{\partial}{\partial x_i^t}\wedge 
\frac{\partial}{\partial x_j^t}$.
\end{itemize}
Then it can be proved that all these objects behave well modulo mutations\footnote{More precisely, given an edge with label $k$ between $t$ and an adjacent vertex $t'$ in $\mathbb T^n$ with associated mutation $\mu_k^t$, 
it can be verified that  $\mu_k^t\circ p^t=\mu_k^t \circ p^{t'}$, 
$(\mu_k^t)^*(\Omega^{t'})=\Omega^{t'}$, 
$(\mu_k^t)_*(\omega^t)=\omega^{t'}$, etc (see \cite{GSV} and \cite{FockGoncharovENS}).} and globalize to the whole cluster varieties and give rise to 
a global map 
\begin{equation}
\label{Eq:Map-p}
p=p_B: \boldsymbol{\mathcal A}_B\longrightarrow \boldsymbol{\mathcal X}_B
\end{equation}
 whose rank is equal to that of $B$\footnote{We recall that the rank of a skew-symmetrizable matrix is stable by mutations.} and whose fibers are precisely the orbits of the global toric action $q=q_B : (\mathbf C^*)^{n-r}\times 
\boldsymbol{\mathcal A}_B\rightarrow\boldsymbol{\mathcal A}_B$ obtained as the gluing of the local torus actions $q^t$ of {\bf (ii).} Moreover, the $\Omega^t$'s (resp.\,the $\omega^t$'s) give rise to a global pre-symplectic 2-form $\Omega=\Omega_B$ (resp.\,a Poisson bivector $\omega=\omega_B$)  on the corresponding cluster variety.\footnote{This explains the 
terminology, introduced by Fock and Goncharov, of  {\it `cluster Poisson manifold'} to designate 
$\boldsymbol{\mathcal X}_B$.} The leaves of the null-foliation defined by $\Omega$  on the $\boldsymbol{\mathcal A}$-cluster variety 
$\boldsymbol{\mathcal A}_B$ 
are precisely the orbits of $q$ hence the fibers of $p$. As a result, the map $p$ can be seen as the pre-symplectic reduction of $(\boldsymbol{\mathcal A}_B,\Omega_B)$. \mk

We denote by $\boldsymbol{\mathcal U}_B$ (or just by 
$\boldsymbol{\mathcal U}$ if no confusion can arise)  
the image of
$\boldsymbol{\mathcal A}_B$
 by 
 the map \eqref{Eq:Map-p}:
\begin{equation}
\boldsymbol{\mathcal U} = 
\boldsymbol{\mathcal U}_B = {\bf Im}\big(p_B\big)=p_B\big(\boldsymbol{\mathcal A}_B\big) \subset 
\boldsymbol{\mathcal X}_B\, .
\end{equation} 
It can be proved that  it is a smooth rational subvariety of $\boldsymbol{\mathcal X}_B$ which moreover is positive (in the sense of `positive geometry', {\it cf.}\,\cite[\S1.1]{FockGoncharovENS}). It carries the non-degenerate closed 2-form $p_*(\Omega)$ hence is naturally a holomorphic symplectic variety. Moreover, one verifies that 
the restriction of the Poisson structure of $\boldsymbol{\mathcal X}_B$ along 
 $\boldsymbol{\mathcal U}_B$ coincides 
with the one induced by $p_*(\Omega)$. \sk 

The subvariety 
$\boldsymbol{\mathcal U}_B$ of $\boldsymbol{\mathcal X}_B$ 
 is called the {\bf `secondary cluster manifold'} in \cite{GSV}.\footnote{
 Note that our cluster variety $\boldsymbol{\mathcal A}_B$ 
 is called the (first) {\it `cluster manifold'} (associated to $A_B$) in \cite{GSV}, where it is denoted by $\boldsymbol{\mathcal X}_B$. And $\boldsymbol{\mathcal Y}_B$ is used there in place of our $\boldsymbol{\mathcal U}_B$. 
 }
 Following a terminology introduced by Fock and Goncharov, one now refers to the 3-tuple 
 $(\boldsymbol{\mathcal A}_B , \boldsymbol{\mathcal X}_B, p_B)$ as the {\bf `cluster ensemble'} associated to the cluster algebra determined by $B$. 
\mk 

Here are some classical concrete examples of cluster varieties or of cluster ensembles: 
\begin{exm} 
\begin{enumerate}
\item[{\bf 1.}] {\bf Surface type.} Let $\overline{S}$ be a compact oriented smooth surface and  $p_1,\ldots,p_n$ be $n$ distinct points  on it. We denote by $S=\overline{S}\setminus\{ p_i\}$ the associated punctured surface.  Let $G$ be a semi-simple algebraic group with trivial center and denote by $G'$ the associated adjoint group 
({\it e.g.}\,in type $A_m$, one has $G={\rm SL}_m$ and $G'={\rm PGL}_m$).

In {\rm \cite{FockGoncharovIHES}}, Fock and Goncharov define two moduli spaces $ \boldsymbol{\mathcal A}_{G,S}$ and  
$ \boldsymbol{\mathcal X}_{G',S}$: the first para-  -metrizes  
 `twisted' $G$-local systems on $S$, the second  `framed' $G'$-local systems on $S$, which are in both cases local systems equipped with some extra data near the punctures.  These two moduli spaces come with a natural map 
    $p_{G,S} : \boldsymbol{\mathcal A}_{G,S}\rightarrow \boldsymbol{\mathcal X}_{G,S}$ and an important result regarding these objects is that $(\boldsymbol{\mathcal A}_{G,S}, \boldsymbol{\mathcal X}_{G',S},p_{G,S})$ is a cluster ensemble.\footnote{For $G={\rm SL}_m$, this is due to Fock and Goncharov (see \cite[\S9]{FockGoncharovIHES}). For classical groups, this has been obtained by Le in \cite{Le} by an impressive case by case treatment. The general result follows from a more conceptual and general approach developed recently by Goncharov and Shen \cite{GoncharovShen}.}
\item[{\bf 2.}] {\bf Type $\boldsymbol{A_n}$.}    In type $A_n$, when the considered decorated surface is the unit disk in $\mathbf C$ with $n+3$ pairwise distinct points on its 
 boundary, one obtains a $\boldsymbol{\mathcal X}$-cluster variety  which can be proved to identify in a natural way to the moduli space $\mathcal M_{0,n+3}$. From this, one gets that the latter identifies with the cluster variety $\boldsymbol{\mathcal X}_{A_n}$ (see  
 {\rm \cite[Appendix B]{FockGoncharovAmalgation}}, 
  {\rm \cite[ \S1.2]{FockGoncharovXInfinity}} and 
  {\rm \cite{King}}
   for more details).  
\item[{\bf 3.}] {\bf Type $\boldsymbol{A_m\boxtimes A_n}$.} By a classical construction, 
to a pair $(\Delta,\Delta')$ of Dynkin diagrams,  one can 
construct a certain product quiver  $\Delta\boxtimes \Delta'$ (see \S\ref{Par:Delta-box-Delta'} below). 
When $\Delta=A_m$ and $\Delta'=A_n$, it is known  ({\it cf.}\,{\rm \cite{Weng}} or 
{\rm \cite[\S6]{GGSVV}}) that the associated $\boldsymbol{\mathcal X}$-cluster variety $\boldsymbol{\mathcal X}_{A_m \boxtimes  A_n}$ identifies naturally with the moduli space 
$ {\rm Conf}_{m+n+2}(\mathbf P^m)$ of projective  configurations 
of $m+n+2$ points in $\mathbf P^m$.
\item[{\bf 4.}] {\bf Finite type $\boldsymbol{\Delta}$.}
In a work in progress {\rm \cite{Pirio-XDelta}}, we endeavour to show that quite similarly to what happens in type $A$, 
the  $\boldsymbol{\mathcal X}$-cluster variety $\boldsymbol{\mathcal X}_\Delta$ associated to any Dynkin diagram $\Delta$ might possibly be interpreted as a  moduli space of polygons. 
%
%
%
\end{enumerate}
\end{exm}
\begin{center}
$\star$
\end{center}

Now given a $\boldsymbol{\mathcal X}$-cluster web $\boldsymbol{\mathcal W}_\Sigma$ defined by a finite set $\Sigma$ of $\boldsymbol{\mathcal X}$-cluster variables, we define a new web by considering its restriction along $\boldsymbol{\mathcal U}$, that we denote by $\boldsymbol{\mathcal U\hspace{-0.05cm}{\mathcal X} \hspace{-0.05cm} \mathcal W}_\Sigma$:
\begin{equation*}
\boldsymbol{\mathcal U\hspace{-0.05cm}{\mathcal X} \hspace{-0.05cm}
 \mathcal W}_\Sigma:= 
\big(\boldsymbol{\mathcal W}_\Sigma\big)\big\lvert_{\boldsymbol{\mathcal U}}\, . 
\end{equation*} 
We will call a web obtained in this way a {\bf `$\boldsymbol{\mathcal U}$-cluster web'} when $\boldsymbol{\mathcal U}$ is a proper subvariety of $\boldsymbol{\mathcal X}$ (otherwise this web coincides with $\boldsymbol{\mathcal W}_\Sigma$ hence one gets nothing new).  As for their $\boldsymbol{\mathcal A}$ or -$\boldsymbol{\mathcal X}$-versions ({\it cf.}\,Proposition \ref{P:The notion of cluster web does not depend on the choice of an initial seed}), it is not difficult to show that the notion of $\boldsymbol{\mathcal U}$-cluster web does not depend on the choice of the initial seed.
\mk 

An interesting case we will study further on, is the finite type case.  
Given a Dynkin diagram, we use the notation 
$\boldsymbol{\mathcal A}_\Delta$, 
$\boldsymbol{\mathcal X}_\Delta$ and $\boldsymbol{\mathcal U}_\Delta$ to denote the corresponding cluster varieties associated to the initial bipartite Dynkin quiver $\vec{\Delta}$.  When the inclusion $\boldsymbol{\mathcal U}_\Delta\subset \boldsymbol{\mathcal X}_\Delta$ is proper, one obtains a new cluster web by considering the restriction of 
$\boldsymbol{\mathcal X\hspace{-0.05cm}\mathcal W}_\Delta$ along $\boldsymbol{\mathcal U}_\Delta$. One gets what we call the {\bf $\boldsymbol{\mathcal U}$-cluster web of type $\boldsymbol{\Delta}$}, denoted as follows
\begin{equation*}
\boldsymbol{\mathcal U\hspace{-0.05cm}{\mathcal X} \hspace{-0.05cm}
 \mathcal W}_\Delta:= 
\big(\boldsymbol{\mathcal X\hspace{-0.05cm}\mathcal W}_\Delta\big)\big\lvert_{\boldsymbol{\mathcal U}_\Delta}\, . 
\end{equation*} 

We will use a similar notation for the $\boldsymbol{\mathcal Y}$-cluster webs:  we will denote by $\boldsymbol{\mathcal U\hspace{-0.05cm}{\mathcal Y} \hspace{-0.05cm}
 \mathcal W}_\Delta$ the restriction of the $\boldsymbol{\mathcal Y}$-cluster web along the secondary cluster variety, {\it i.e.} 
\begin{equation*}
\boldsymbol{\mathcal U\hspace{-0.05cm}{\mathcal Y} \hspace{-0.05cm}
 \mathcal W}_\Delta:= 
\big(\boldsymbol{\mathcal Y\hspace{-0.05cm}\mathcal W}_\Delta\big)\big\lvert_{\boldsymbol{\mathcal U}_\Delta}\, . 
\end{equation*} 

The dimension of $\boldsymbol{\mathcal U}_\Delta$ coincides with the rank  of the initial exchange matrix $B_\Delta$, which is  not difficult to determine.  This rank 
is maximal (that is equal to the rank of $\Delta$)  
in the following cases: $A_{2n}$, $B_{2n}$, $C_{2n}$,  $E_6$, $E_8$, $F_4$  and $G_2$.   For   $\Delta$ of another type,  the rank of $B_{{\Delta}}$ is not maximal and is given in the following table:

  \begin{table}[h]
\begin{center}
\begin{tabular}{|c|c|c|c|c|c|}
\hline 
 \begin{tabular}{c}   \vspace{-0.15cm} \\
 $\boldsymbol{{\Delta}}$
 \vspace{0.2cm} \\ 
\end{tabular}
  & 
 $\boldsymbol{A_{2n+1}}$ & $\boldsymbol{B_{2n+1}}$ 
  & $\boldsymbol{C_{2n+1}}$ 
 &  
  $\boldsymbol{D_n}$ &  $\boldsymbol{E_7}$
    \\ \hline
     \begin{tabular}{c}   \vspace{-0.15cm} \\
 $\boldsymbol{{\bf rk}\big(B_{{\Delta}}\big)}$
 \vspace{0.2cm} \\ 
\end{tabular}
  & 
    \begin{tabular}{l}   
 $2n$  
  \end{tabular} &   
      \begin{tabular}{l}   
 $2n$  
  \end{tabular} & 
  \begin{tabular}{l}   
  $2n$ 
  \end{tabular}
    & 
  \Bigg\{ \hspace{-0.35cm}
   \begin{tabular}{l}   
 \hspace{0.06cm}  $n-2$ \hspace{0.06cm}  for $n$ even \vspace{-0.01cm} \\
 \hspace{0.06cm}   $n-1$  \, for $n$ odd
  \end{tabular} & 6
\\
 \hline 
\end{tabular}
\end{center}
\caption{Rank of the exchange matrix $B_{
 {\Delta}}$ when it is not maximal.}
\end{table}


  
   
The first example of $\boldsymbol{\mathcal U}$-cluster web in finite type to be considered is certainly the one associated to $A_3$: in this case $p$ has rank 2 hence $\boldsymbol{\mathcal U\hspace{-0.05cm}{\mathcal X} \hspace{-0.05cm}
 \mathcal W }_{A_3}$ is a web in two-variables. 
  In the $\boldsymbol{\mathcal X}$-cluster torus associated to the initial seed 
$((x_1^{-1},x_2,x_3^{-1}), B_{A_3})$\footnote{For the reason why we use $(x_1^{-1},x_2,x_3^{-1})$ instead of 
$(x_1,x_2,x_3)$ for the initial cluster, see Remark \ref{Rem:birat-E} above.}, 
the image of $p$ is cut out by $x_1=x_3$. Injecting this in \eqref{Eq:XWA3-Definition} and after a few elementary computations, one deduces  the following explicit expression for $\boldsymbol{\mathcal U\hspace{-0.05cm}{\mathcal X} \hspace{-0.05cm}
 \mathcal W }_{\hspace{-0.05cm} A_3}$ in the coordinates $x_1,x_2$: 
\begin{align}
\label{Eq:pXWA3}
\boldsymbol{\mathcal U\hspace{-0.05cm}{\mathcal X} \hspace{-0.05cm}
 \mathcal W}_{\hspace{-0.05cm} A_3}
=
\boldsymbol{\mathcal W}\Bigg( 
\, 
x_1\, , \, 
x_2\, , \, 
\frac{1+x_2}{x_1}
\, , \,  & 
\frac{1+x_1}{x_2}
\, , \, 
\frac{(1+x_1)^2}{x_2}
\, , \,  
\frac{1+x_1+x_2}{x_1x_2}
\, ,  \\ & 
\frac{(1+x_1)^2+x_2}{x_1x_2}
\, ,\,  
\frac{1+x_1+x_2}{x_1(1+x_1)}
\, , \, 
\frac{(1+x_1  + x_2)^2}{x_1^2x_2}
\hspace{0.15cm}
\Bigg)\, .
\nonumber
\end{align}
We will study this planar 9-web defined by positive rational functions further  in \S\ref{SS:ClassicalPolylogarithmicIdentitiesAreOfClusterType}, where we will prove that it is equivalent to Spence-Kummer trilogarithmic web ${\boldsymbol{\mathcal W}}_{\!{\cal S}{\cal K}}$ of  \S\ref{SS:Spence-Kummer}. 
\mk 

Another interesting example is the case of $D_4$:  
 since $B_{D_4}$ has rank 2, 
$\boldsymbol{\mathcal U\hspace{-0.05cm}{\mathcal X} \hspace{-0.05cm}
 \mathcal W}_{\hspace{-0.1cm}D_4}$ is a web in two-variables. 
  In the $\boldsymbol{\mathcal X}$-cluster torus associated to the initial seed 
$((x_1^{-1},x_2,x_3^{-1},x_4^{-1}), B_{D_4})$, the map $p$ is given by 
$(a_i)_{i=1}^4\mapsto(a_2^{-1}, a_1a_3a_4,a_2^{-1},a_2^{-1})$ hence $
\boldsymbol{\mathcal U}_{D_4}$  is cut out by $x_1=x_3=x_4$. After some computations, one gets the following explicit expression for 
$\boldsymbol{\mathcal U\hspace{-0.05cm}{\mathcal X} \hspace{-0.05cm}
 \mathcal W}_{\hspace{-0.1cm}D_4}$ in the coordinates $x_1$ and $x_2$ on the intersection 
of $
\boldsymbol{\mathcal U}_{D_4}$ with the $\boldsymbol{\mathcal X}$-cluster tori associated to the initial seed: 
\begin{align}
\label{Eq:pXWD4}
\nonumber
\boldsymbol{\mathcal U\hspace{-0.05cm}{\mathcal X} \hspace{-0.05cm}
 \mathcal W}_{\hspace{-0.1cm}D_4}=
\boldsymbol{\mathcal W}
 \Bigg( \, {x_{{1}}}
 \, , \, 
 {x_{{2}}}
  \, , \,  &
 {\frac{1+x_{{2}}}{x_{{1}}}}
  \, , \, 
  {\frac {1+x_{{1}}}{x_{{2}}}}
   \, , \, 
 {\frac { \left( 1+x_{{1}} \right) ^{2}}{x_{{2}}}}
   \, , \, 
 {\frac{ \left( 1+x_{{1}} \right) ^{3}} {x_{{2}}}}
  \, , \, 
{\frac {1+x_{{1}}+x_{{2}}}{x_{{1}}x_{{2}}}}
 \, , 
 \\ 
 \nonumber
 & 
{\frac {1+x_{{1}}+x_{{2}}}{x_{{1}} \left( 1+x_{{1}} \right) }}
 \, , \, 
{\frac { \left(1+ x_{{1}}+x_{{2}} \right) ^{2}}{{x_{{1}}}^{2}x_{{2}}}}
 \, , \, 
{\frac { \left( x_{{2}}+1+x_{{1}} \right) ^{2}}{  (1+{x_{{1}}})^{2}+x_{{2}} }}
 \, , \, 
{\frac { \left( 1+x_{{1}}+x_{{2}} \right) ^{3}}{{x_{{1}}}^{3}x_{{2}}}}
 \, , \, \\  
 &
{\frac { (1+x_{{1}})^{2}+x_{{2}}}{x_{{1}}}}
 \, , \, 
{\frac {(1+{x_{{1}}})^{2}+x_{{2}}}{x_{{1}} \left( 1+x_{{1}} \right) ^{2}}}
\, , \,
{\frac { \left( {(1+x_{{1}})}^{2}+x_{{2}}\right) ^{3}}{{x_{{2}}}^{2}{x_{{1}}}^{3}}}
 \, , \,
{\frac {(1+x_{{1}})^{3}+x_{{2}}}{x_{{1}}x_{{2}}}}
 \, , \, 
 \\ 
 \nonumber
 &  
 {\frac { (1+x_{{1}})^{3}+2\,x_{{2}}+{x_{{2}}}^{2}
 +3\,x_{{1}}x_{{2}}
  }{{x_{{1}}}^{2}x_{{2}}}}
\, ,\, 
{\frac  { \left( (1+x_{{1}})^{2}+x_{{2}} \right) ^{2}}
{x_{{1}}x_{{2}}
(1+x_1+x_2)}}
\, ,\, 
{\frac { \left( (1+x_{{1}})^{2}+x_{{2}} \right) ^{2}}{{x_{{1}}}^{2} x_2 \left( 1+x_{{1}} \right) }} \, \,  \Bigg)\, .
\end{align}
One obtains this way a planar 18-web, again defined by positive rational functions. We will also study it further and will prove further ({\it cf.}\,Theorem \ref{T:classical-cluster-webs} below) that it is equivalent to Kummer's tetralogarithmic web 
$\boldsymbol{{\boldsymbol{\mathcal W}}_{{\cal K}_4}}$ 
considered in  \S\ref{Parag:Kummer_tetralog}.
\mk

The fact that the webs in the two preceding examples are equivalent to webs associated to classical  polylogarithmic identities (of weight 3 and weight 4 respectively) shows that the notion of `$\boldsymbol{{\mathcal U}}$-cluster web' deserves to be studied more in depth.

  \subsection{\bf \bf  Cluster webs associated to $\boldsymbol{Y}$-systems and periods}
\label{SS:Y-systems-Periods}
The several notions of cluster webs introduced in the previous subsections are very specific to the finite type case since 
to define theses webs, we 
consider all the cluster variables we have at disposal.  
In order to consider interesting webs when the set of all cluster variables is not finite, one has to be more subtle and find particular finite subsets of such variables.  \sk

From this perspective, a relevant notion of the theory of cluster algebras is that  of `{\it period}', which has been introduced then studied  by Nakanishi  in several papers devoted to understand better, in terms of cluster algebras, some periodicity phenomena satisfied by certain families of rational identities coming from mathematical physics, the so-called {\it `$Y$-systems'}.  The webs associated to some $Y$-systems of Dynkin type being particularly interesting from the point of view of their ARs and their rank, we first give a very concise account of the theory of $Y$-systems before introducing the more general notion of (cluster) period.


\subsubsection{$\boldsymbol{Y}$-systems and webs associated to them.}
\label{SS:Generalities-Y-systems}
A {\it `${Y}$-system'} is a set of countable rational identities concerning the solutions of the {\it `Thermodynamic Bethe Ansatz'}, a method  designed to  better understand  some integrable systems coming from conformal fields theory. 
 There are many such $Y$-systems and there is now a rich literature (in physical  as well as in  mathematical journals) on the subject. In the sequel,  we will only consider some particular examples and focus on the mathematical side of the theory, saying almost nothing of the physical one.  For more details about this notion, from a mathematical point of view,  the interested reader can consult \cite{FZ},   
\cite{Kuniba-Al}, \cite{KellerAnnals}, \cite{NakanishiNagoya},  \cite{Nakanishi} and the references therein. 


\paragraph{$\boldsymbol{Y}$-systems.}
\hspace{-0.2cm} 
Let $I$ be a fixed finite set of indices. 
We consider 
$A=(a_{ij})_{i,j\in I}$, a 
 matrix 
with integer coefficients
and $\epsilon=(\epsilon_{ij})$ a matrix of signs ({\it i.e.}\,$\epsilon_{ij}\in \{1,-1\}$ for any $i,j\in I$). By definition, the {\bf $\boldsymbol{Y}$-system}  associated to these data is  the family of rational identities 
$${}^{} 
\hspace{1.5cm}
\Big( \boldsymbol{Y}_{\hspace{-0.05cm} A,\epsilon}\Big) \hspace{3cm}
Y_i(t+1)Y_i(t-1)=\prod_{j\in I} \Big( 1+Y_{\hspace{-0.1cm}j}(t)^{\epsilon_{ij}}\Big)^{a_{ij}}\qquad i\in I, \, t\in \mathbf Z\, .
  \hspace{4cm} {}^{} 
$$

In full generality, it is certain that not much can be told about such a set of identities, but for some particular cases of physical importance, and most of the time related to Lie theory, two important 
properties are conjectured to hold true for any set of variables $\{ Y_i(t) \}_{i\in I,t\in \mathbf Z}$ satisfying $( \boldsymbol{Y}_{\hspace{-0.05cm}A,\epsilon})$: 
\begin{itemize}
\item[$\bullet$] {\bf  [Periodicity]:} 
for a certain positive integer $m$, one has 
$$
\big({\bf Per}\boldsymbol{( \boldsymbol{Y}_{\hspace{-0.05cm}A,\epsilon})}
\big)\,\hspace{2.3cm}
Y_i(t+m)=Y_i(t) \qquad \mbox{for all }\, (i,t)\in I\times \mathbf Z\, ;
  \hspace{3cm} {}^{} 
$$
\item[$\bullet$] {\bf  [Dilogarithmic identity]:}  
 one sets $u_i=Y_i(0)$ for any $i\in I$ and one considers the $u_i$'s as initial variables. Then there exist:
 \begin{itemize}
 \item  a certain finite set $\mathcal S( \boldsymbol{Y}_{A,\epsilon})$ of pairs 
$(i,t)\in 
I\times \mathbf Z$;
 \item some rational functions $y_i(t)$ in the $u_i$'s, for $(i,t)\in \mathcal S( \boldsymbol{Y}_{A,\epsilon})$;
 \item some non-zero integers $d_{i,t}$ for $(i,t)\in \mathcal S( \boldsymbol{Y}_{A,\epsilon})$;
 \item a rational constant $c( \boldsymbol{Y}_{A,\epsilon})$, 
  \end{itemize}
 such that 
 the following dilogarithmic identity is satisfied
$$ 
\big(\boldsymbol{\mathcal R( \boldsymbol{Y}_{\hspace{-0.05cm} A,\epsilon})}
\big)\,\hspace{2.3cm}
\sum_{(i,t)\in \mathcal S(
 \boldsymbol{Y}_{A,\epsilon} 
)} d_{i,t}\cdot \mathcal R\left(  
\frac{y_i(t)}{1+y_i(t)}
\right)= c(\boldsymbol{Y}_{A,\epsilon})\,
\frac{{}^{}\hspace{0.1cm}\pi^2}{6}\, 
  \hspace{3cm} {}^{} 
$$
(where $\mathcal R$ stands for the original Rogers' dilogarithm, {\it cf.}\,\S\ref{Par:NotationDilogarithmicFunctions}). 
\end{itemize}

As far as  we know, these two properties (the latter can only be formulated if the former  holds true) are still conjectural in full generality. But they have been established in many interesting cases, some of which are discussed more in detail  in the next subsection.  
\sk 

But  two important remarks, on which we will elaborate further, are already in order  for a given $Y$-system $(\boldsymbol{Y})$, under the assumption that the periodicity holds true. 
\begin{rem}
\label{Rem:Y-systems}
\begin{enumerate}
\item There are formulae for all the quantities (the set of indices $\boldsymbol{\mathcal S}(\boldsymbol{Y})$, the rational functions $y_i(t)$'s, the integers $d_{i,t}$s, etc.) involved in the dilogarithmic identity 
$({\mathcal R( \boldsymbol{Y})})$, in terms  of the data used to construct $\boldsymbol{Y}$;\footnote{See 
\S\ref{Par:NakanishiDilogarithmicIdentity-Period} for some formulas in the more general case of the dilogarithmic identity associated to a cluster period.}
\sk
\item 
 The  $y_i(t)$'s involved in the dilogarithmic identity 
 $({\mathcal R( \boldsymbol{Y})})$ 
 are of course related to the $Y_i(t)$'s appearing in the very definition of 
 the considered $Y$-system  $\boldsymbol{(Y)}$ and in many classical cases the  two corresponding sets of rational functions in the $u_i$'s, namely  
 $\{ \, y_i(t) \, \lvert \,  (i,t)\in \boldsymbol{\mathcal S}(\boldsymbol{Y}) \}$ and $\{ \, Y_i(t)\, \lvert \,   i\in I,\, 0\leq t<m  \}$, 
  are the same (possibly up to inversion $x\leftrightarrow x^{-1}$).  However, this does not hold true in full generality. 
  \end{enumerate}
\end{rem}

\paragraph{Webs associated to $Y$-systems.} 
\hspace{-0.4cm}
\label{parag:Webs-associated-to-Y-systems}
What makes   the  notion of $Y$-system interesting for us is that 
as soon as the two properties above are satisfied for such a system $(\boldsymbol{Y})$, then it is natural to associate to it the web formed by the $y_i(t)$'s appearing in the dilogarithmic 
identity  $(\mathcal R(\boldsymbol{Y}))$, viewed as rational functions of the initial $Y$-variables $u_i=Y_i(0)$'s  (with $i\in I$): 
$$
\boldsymbol{\mathcal{W}}(\boldsymbol{Y})=\boldsymbol{\mathcal{W}}\left(
\, y_i(t)\hspace{0.15cm} \Big\lvert \hspace{0.05cm} 
\begin{tabular}{l}
$y_i(t)$ {\it appears in the dilog- }
\vspace{-0.0cm}
\\
{\it -arithmic identity}  $ \mathcal R(\boldsymbol{Y})$  
\end{tabular} \right)\, .
$$ 
This is a web defined by rational functions on $\mathbf C^n$ where $n=\lvert I\lvert$.   Assuming that $\mathcal R(\boldsymbol{Y})$ holds true gives that  $\boldsymbol{\mathcal{W}}(\boldsymbol{Y})$ carries a non-trivial dilogarithmic AR that will be denoted (a bit abusively)  again by $\mathcal R(\boldsymbol{Y})$.\sk 

We will say that $(\boldsymbol{Y})$ is {\bf $\boldsymbol{\mathcal W}$-regular} (as a $Y$-system) if it can also be defined as the web admitting the $Y_i(t)$'s as first integrals ({\it cf.}\,Remark \ref{Rem:Y-systems}.2 just above).  Assuming that 
$(\boldsymbol{Y})$ is of this type, then 
 taking the logarithm
of any multiplicative relation 
$Y_i(t+1)Y_i(t-1)=\prod_{j\in I} \big( 1+Y_{\hspace{-0.1cm}j}(t)^{\epsilon_{ij}}\big)^{a_{ij}}$ appearing in the definition of $Y$ 
gives  the following additive relation
$$ 
\Big({Log}AR_i(t)\Big)
\hspace{1.5cm}
{\rm Log}\big( Y_i(t+1)\big)+ {\rm Log} 
\big(
Y_i(t-1)\big)=\sum_{j \in I} a_{ij}{\rm Log}\big( 1+Y_{\hspace{-0.1cm}j}(t)^{\epsilon_{ij}}\big)\, , 
\hspace{2.5cm}{}^{}
$$ 
which can be seen as a logarithmic AR for $ \boldsymbol{\mathcal{W}}(\boldsymbol{Y})$. Thus one gets a linear inclusion
\begin{equation}
\label{Eq:AR-Y-web}
\Big\langle \,  {Log}AR_i(t) \, \big\lvert \, i\in I,\, t\in \mathbf Z\, 
\Big\rangle  \oplus 
\big\langle \,  \mathcal R(\boldsymbol{Y})
\big\rangle \subset  \boldsymbol{\mathcal A}(\boldsymbol{\mathcal{W}}(\boldsymbol{Y}))
\end{equation}
from which it follows that $\boldsymbol{\mathcal{W}}(\boldsymbol{Y})$ carries many polylogarithmic ARs  (of weight 1 and 2 at least). \sk 

The preceding considerations make us ask the following questions:
\begin{questions}
\label{Quest:YWeb-AMP?}
   Let $(\boldsymbol{Y})$ be a $Y$-system. \mk \\
{\bf 1.} Is the web $\boldsymbol{\mathcal{W}}(\boldsymbol{Y})$ AMP with only polylogarithmic ARs (of weight 1 or 2)?\mk \\
{\bf 2.}  Does the linear inclusion \eqref{Eq:AR-Y-web} induce an isomorphism ?
\mk \\
{\bf 3.} Do the ${Log}AR_i(t) $'s form a basis of the space of logarithmic ARs of $\boldsymbol{\mathcal{W}}(\boldsymbol{Y})$? 
\mk \\
{\bf 4.} When the answer to  {\bf 3.} is negative, what are the relations between the  ${Log}AR_i(t)$'s and how can  a supplementary space of 
$\big\langle \,  {Log}AR_i(t) \, \big\lvert \, i\in I,\, t\in \mathbf Z\, 
\big\rangle$ in $\boldsymbol{{Log}AR}\big( 
\boldsymbol{\mathcal{W}}(\boldsymbol{Y})\big)$ be described? 
\end{questions}

What makes  the notion of $Y$-system interesting and the previous questions relevant  from the point of view of web geometry is that the answer  to the first  seems 
 to be affirmative, at least for many of the most classical types of $Y$-systems, that we will now introduce.


\paragraph{Some classical examples of $\boldsymbol{Y}$-systems associated to Dynkin diagrams.} 
\label{Par:Classical-Y-systems}
In this paragraph, we use the following notations: 
for a given Dynkin diagram $\Delta$, we denote by $n$ its rank, 
$h$ its Coxeter number, $\check{h}$ the dual one (see
 Table \ref{Table:CoxeterNumbers} below); $C=C(\Delta)=(c_{ij})_{ij=1}^n$ stands for the associated Cartan matrix and we set $A=(a_{ij})=2\,{\rm Id}_n-C$. If $\Delta'$ denotes another Dynkin diagram, we use prime notations $n'$, $h'$, $A'=(a'_{ij})$ etc.  to denote the corresponding objects.

\begin{table}[!h]
\begin{center}
\begin{tabular}{|c|c|c|c|c|c|c|c|c|c|c|} 
\hline
\begin{tabular}{c}
\vspace{-0.3cm} \\
${\boldsymbol{\Delta}}$ 
\vspace{0.15cm}
\end{tabular}
&   
$\boldsymbol{A_n}$ 
  &   $\boldsymbol{B_n}$  & $\boldsymbol{C_n}$ & $\boldsymbol{D_n}$  & $\boldsymbol{E_6}$ & $\boldsymbol{E_7}$ & $\boldsymbol{E_8}$  &$\boldsymbol{F_4}$ &  $\boldsymbol{G_2}$     
   \\ 
    \hline 
    \begin{tabular}{c}
\vspace{-0.3cm} \\
 ${\boldsymbol{h}}$
\vspace{0.15cm}
\end{tabular}
   & $n+1$ &  $2n$ & $2n$ & $2n-2$ & 12 & 18 & 30 & 12 & 6
       \\ 
    \hline 
    \begin{tabular}{c}
\vspace{-0.3cm} \\
\, 
${\boldsymbol{{h}^\vee}}$
\vspace{0.15cm}
\end{tabular}
& $n+1$ &  $2n-1$ & $n+1$ & $2n-2$ & 12 & 18 & 30 & 9 & 4 
   \\ 
    \hline 
\end{tabular}
\mk
\caption{Coxeter numbers and dual Coxeter numbers of Dynkin diagrams}
\label{Table:CoxeterNumbers}
\end{center}
\end{table}

 We now introduce some of the most classical $Y$-systems considered in the literature. For each of them, we give precise forms for the periodicity and for the associated dilogarithmic identity when these latter  do already appear in the existing literature. 
\begin{enumerate}
\item[$\bullet$]  {\bf [Type $\boldsymbol{Y (\Delta)}$].}  
The $Y$-system  associated to a Dynkin diagram $\Delta$  of rank $n$ is
$$
\big(\boldsymbol{Y (\Delta)}\big)\hspace{3cm}
Y_i(t-1)Y_i(t+1)=\prod_{j=1}^n\Big(1+Y_j(t) \Big)^{a_{ij}}\, \hspace{0.5cm} i=1,\ldots,n, \, t\in \mathbf Z
\,;
\hspace{4cm}
$$
In this case, the periodicity $({\bf Per}\boldsymbol{( \boldsymbol{Y}_{\hspace{-0.05cm}\Delta})}$) is 
expressed as $Y_i(t+2(h+2))=Y_i(t)$ for all $i$ and $t\in \mathbf Z$, and the associated dilogarithmic  identity
$(\mathcal R(\boldsymbol{Y}_{\hspace{-0.05cm}\Delta}))$
 takes the following form 
  $$
\big(\boldsymbol{\mathcal R (\Delta)}\big)
\hspace{3cm}
\sum_{i=1}^n \sum_{t=0}^{2(h+2)-1} \mathcal R\bigg( \frac{Y_{i}(t)}{
  1+Y_{i}(t) }     \bigg)=2 hn\frac
  {{}^{}\hspace{0.1cm}\pi^2}{6} \,;
\hspace{4.3cm}
  $$
\item[$\bullet$]  {\bf [Type $\boldsymbol{Y (\Delta,\Delta')}$].}   
Let $\Delta,\Delta'$ be two Dynkin diagrams of rank $n$ and $n'$.   One denotes by $I$ the set of pairs $(i,i')$ with $i\in \{1,\ldots,n\}$ and $i'\in \{1,\ldots,n'\}$.  
 The associated $Y$-system is 
$$
\big(\boldsymbol{Y (\Delta,\Delta')}\big)\hspace{1cm}
  Y_{i,i'}(t+1) Y_{i,i'}(t-1)=\frac
  { \prod_{j=1}^n \Big( 1+Y_{j,i'}(t)\Big)^{a_{ij} }}
  {\prod_{j'=1}^{n'} \Big( 1+Y_{i,j'}^{-1}(t)\Big)^{a'_{i'j'}}}
    \qquad  
   { \begin{tabular}{l}
    $i=1,\ldots,n$\vspace{-0.0cm}\\
    $i'\!=1,\ldots,n'$, \, 
      $t\in \mathbf Z$
      \end{tabular}}
\hspace{4cm}
$$
The periodicity is $Y_i(t+2(h+h'))=Y_i(t)$ for all $(i,i')\in I$ and $t\in \mathbf Z$.  When both $\Delta$ and $\Delta'$ are simply-laced, the associated dilogarithmic identity  is (see \cite[Conjecture 1.6]{NakanishiNagoya}):
     $$
\big(\boldsymbol{\mathcal R (\Delta,\Delta')}\big)
\hspace{2.7cm}
 \sum_{(i,i')\in I\times I'} 
 \hspace{-0.3cm}
 \sum_{t=0}^{2(h+h')-1}\mathcal R\left( \frac{Y_{ii'}(t)}{  1+Y_{ii'}(t)}     \right)=2hnn' \frac
  {{}^{}\hspace{0.3cm}\pi^2}{6} \,.
  \hspace{4cm}
  $$ 
  When at least one of the Dynkin diagrams is not simply-laced,  the dilogarithmic identity $({\mathcal R (\Delta,\Delta')})$ 
does not appear in explicit form   
   in the existing literature (except if one of those is of type $A$, see just below).   However, such an identity is a particular case of the theory of dilogarithmic identities associated to cluster periods developed by Nakanishi, that will be succinctly discussed in the next subsection. \sk 
   
 Finally, we mention the fact that exchanging $\Delta$ and $\Delta'$ corresponds to inverting the variables $Y_{i,i'}(t)\leftrightarrow 1/Y_{i'i}(t)$ in the corresponding $Y$-systems, hence is essentially not relevant from the web geometric perspective we are interested in. 
\item[$\bullet$]  {\bf [Type $\boldsymbol{Y_\ell (\Delta)}$].}   
Let $\Delta$ be a Dynkin diagram of rank $n$ and let $\ell$ be a {\it `level'}, that is an integer bigger than or equal to 2\footnote{When $\Delta$ is not simply-laced, one can even take $\ell=0,1$ but we will not consider these particular  cases here.}. 
Then one can define the {\bf $Y$-system of Dynkin type $\Delta$ and level $\ell$}, denoted by   $({Y_\ell (\Delta)})$.  
If  $\Delta$ is of type $ADE$, this $Y$-system is equivalent to 
$({Y (A_{\ell-1},\Delta)})$ hence one gets something new only when $\Delta$ is not simply-laced.  In the latter case, there is no uniform way (with respect to the Cartan matrix of $\Delta$) to write down the algebraic relations characterizing ${Y_\ell (\Delta)}$. To save space, we will not write these here but refer to \cite[\S2.2]{IIKNS} where they are given in explicit form.\sk 

The $Y$-system ${Y_\ell (\Delta)}$ involves a set of 
variables $Y_{i,m}(t)$ with $(i,m)\in I\times \mathbf N_{>0}$ and $t\in \mathbf Z$ and the associated periodicity 
condition is : `$Y_{i,m}(t+2(h^\vee+\ell))=Y_i(t)$ for all $(i,m)\in I$ and $t\in \mathbf Z$'.  As for the associated dilogaritmic identity, 
it  is written   ({\it cf.}\,\cite{Nakanishi} and \cite{Inoue2})
%
%
%
   $$
\big(\boldsymbol{\mathcal R_\ell (\Delta)}\big)
\hspace{2cm}
 \sum_{ (i,m,t)\in \mathcal Y_\ell(\Delta)} 
  \mathcal R\left( \frac{Y_{i,m}(t)}{  1+Y_{i,m}(t)}     \right)=
2 \tau  n(\ell h-h^\vee)
   \frac
  {{}^{}\,\pi^2}{6} 
  \hspace{4cm}
  $$
where $\mathcal Y_\ell$ is a finite set (depending on $\Delta$ and $\ell$) and $\tau$ is an integer  equal to 1,2 or 3.\footnote{More precisely, one has $t=1$ for $\Delta$ simply-laced, $\tau=3$ for $\Delta=G_2$ and $\tau=2$ otherwise.}.\sk
%
\end{enumerate}
 
 \begin{rem} 
 {1.} Among the $Y$-systems described above, it seems that 
 those which are  $\boldsymbol{\mathcal W}$-regular according to the terminology introduced above are: the $Y(\Delta)$'s for any $\Delta$ and 
the $Y(\Delta,\Delta')$  when both $\Delta$ and $\Delta'$ are simply laced (note that this covers the level $\ell$ $Y$-systems for any Dynkin diagram of type $ADE$).
 \mk \\
  {\it 2.} 
   When $\Delta$ or $\Delta'$, say $\Delta$,  is not simply-laced, then 
 it seems that the dilogarithmic identities $(\mathcal R(\Delta,\Delta'))$ or $\mathcal R_\ell(\Delta)$ (with $\ell\geq 2$) 
  can be expressed in terms of rational functions $y_i(t)$'s, some of which define foliations which do not admit any function $Y_i(t)$ as a first integral. Consequently,  the corresponding $Y$-system is not $\boldsymbol{\mathcal W}$-regular.
 \end{rem}
 \sk 
 
 As already mentioned before, the notion of $Y$-system
 comes from mathematical physics, as well as the expectation that 
  such a system always satisfies the two  fundamental properties 
  mentioned above (periodicity and dilogarithmic identity).  These properties were only conjectured at the end of the 80s/beginning of the 90s\footnote{As for early references on the subject, in which the periodicity and/or the corresponding dilogarithmic identity are conjectured (and even established in some cases), we mention 
  \cite{Zamolodchikov,GliozziTateo0} for the periodicity and more specifically \cite{KNS1,Kirillov1,KM,Kirillov} in what concerns the corresponding dilogarithmic identities.} and  have given rise  since  to a great deal of mathematical papers 
  in which the periodicity and the  associated dilogarithmic identities are established:
%
 %
%
%
\begin{thm}  For any three of the $Y$-systems 
of Dynkin type presented  above, the periodicity property  as well as the corresponding dilogarithmic identity hold true. 
\end{thm}
\sk

\begin{table}[!h]
\begin{center}
\begin{tabular}{|c|c|l|} 
\hline
\begin{tabular}{c}
\vspace{-0.3cm} \\
${\boldsymbol{Y}}${\bf -system} 
\vspace{0.15cm}
\end{tabular}
&   
{\bf Periodicity}
  &    ${}^{}$\,  {\bf Dilogarithmic identity}
    \\ 
    \hline 
\begin{tabular}{c}
\vspace{-0.3cm} \\
${\boldsymbol{Y}}(A_m
 )$  
\vspace{0.15cm}
\end{tabular}
    &    \cite{FrenkelSzenes}, \cite{GliozziTateo}   &       \cite{FrenkelSzenes}, \cite{GliozziTateo}  \\
\hline
\begin{tabular}{c}
\vspace{-0.3cm} \\
${\boldsymbol{Y}}(A_m,A_n)$ 
\vspace{0.15cm}
\end{tabular}
     &   \cite{Szenes}, \cite{Volkov}   &  
     \cite{NakanishiNagoya}
     \\
\hline
\begin{tabular}{c}
\vspace{-0.3cm} \\
 ${\boldsymbol{Y}}(\Delta
  )$  
\vspace{0.15cm}
\end{tabular}     &    \cite{FZ}   &       \cite{Chapoton}
\hspace{0.38cm}
 ($\Delta$ ADE) \\
\hline
\begin{tabular}{c}
\vspace{-0.3cm} \\
 ${\boldsymbol{Y}}(\Delta,\Delta')$  
\vspace{0.15cm}
\end{tabular}     &    \cite{KellerAnnals}   &    
\hspace{-0.3cm}   
\begin{tabular}{l}
\vspace{-0.3cm} \\
 \cite{NakanishiNagoya} \hspace{0.05cm} ($\Delta\,  \& \,\Delta'$  $ADE$)
 \vspace{0.05cm}
 \\
 \cite{Nakanishi}
\vspace{0.15cm}
\end{tabular}  \\
\hline
\begin{tabular}{c}
\vspace{-0.3cm} \\
 ${\boldsymbol{Y}}_\ell(\Delta)$  
\vspace{0.15cm}
\end{tabular}  &  
\begin{tabular}{c}     \vspace{-0.3cm} \\  
\cite{KellerAnnals}
\vspace{0.05cm}
\\
\cite{Inoue1} \hspace{0.05cm} \vspace{0.05cm}
\\ 
\cite{Inoue2} \hspace{0.05cm}
\vspace{0.15cm}
\end{tabular}
 &  
\hspace{-0.3cm}
\begin{tabular}{l}     \vspace{-0.3cm} \\  
\cite{NakanishiNagoya} \hspace{0.15cm}($\Delta$  $ADE$) \vspace{0.05cm}
\\
\cite{Inoue1} \hspace{0.05cm}($B_n$)  \vspace{0.05cm}
\\ 
\cite{Inoue2} \hspace{0.05cm}($C_n,F_4,G_2$) 
\vspace{0.15cm}
\end{tabular}
  \\
\hline
\end{tabular}\mk
\caption{Periodicity and dilogarithmic identities for $Y$-systems  of Dynkin types}
\label{Table:Y-systemsDD'}
\end{center}
\end{table}


The corresponding references for the different cases are given in  Table \ref{Table:Y-systemsDD'} above.  It is interesting to make a few comments on how these results have been established

\begin{itemize}
\item[$\bullet$]
The crucial point to establish seems to be  the periodicity property of a  given type of $Y$-system. Once it has been established,  the proof that the corresponding dilogarithmic identities hold true always consists in verifying that 
the algebraic criterion given by Theorem \ref{T:FE-Zagier-Polylogs}
is satisfied. Note however that establishing the latter fact is far from being trivial in general. But this has been established in 
the quite more general context of {\it `periods'}  of the theory of cluster algebras by Nakanishi in  \cite{Nakanishi} that will be discussed further below.

\item[$\bullet$]
To establish the periodicity property of a given type of $Y$-system, several distinct methods have been used so far. 
For instance, rather direct methods are used in  \cite{GliozziTateo} (in relation with 3-dimensional differential geometry), in \cite{Volkov} (via projective geometry) or in \cite{Szenes} (using the flatness of certain connexions on some graphs). 
Each of these approaches seems however quite specific and ad hoc to the  considered cases.   
 There are also  other, more general methods of establishing periodicity ({\it cf.}\,\cite[\S3.4]{IIKNS} for instance).  Those of greater scope and which seem the most powerful are the {\it `cluster algebra/category method'} and  the one by means of {\it tropicalization} (see \cite{KellerAnnals} and \cite{NakanishiTrop} respectively). Both  rely in a crucial way on the fact that the $Y$-systems of Dynkin type can be formalized by means of cluster algebras. 
\end{itemize}

The two previous remarks indicate that their interpretation in terms of cluster algebras allows to  better understand the $Y$-systems of Dynkin type. This 
gives us grounds for looking at the webs associated to $Y$-systems introduced above  from the same point of view.

\paragraph{Interpretation of the periodicity of $\boldsymbol{Y}$-systems of Dynkin type
in terms of cluster algebras.} 
\label{Par:Interpretation-cluster-Y-systems-Dynkin-type}

We only deal here with the $Y$-systems associated 
to pairs of Dynkin diagrams $(\Delta,\Delta')$. This covers all the cases, except those of the form $\boldsymbol{Y}_\ell(\Delta'')$ with $\ell\geq 2$ and  $\Delta''$ not simply-laced. These  can also be interpreted by means of cluster algebras, but this will be discussed briefly below, in relation to the notion of period.
\mk 

\subparagraph{}\hspace{-0.6cm}
\label{SubPar:BipartiteBelt}
We start with the case when $\Delta'$ coincides with $A_1$. Our main references are \cite{FZ} and more specifically the eighth section  of \cite{FZIV}. 
 Let $n>1$ be the rank of $\Delta$ and denote by $B=B_{\vec{\Delta}}=(b_{ij})$ the exchange matrix associated with the 
bipartite quiver of type $\Delta$ (see \S\ref{Par:Cluster-Algebras-Finite-Type}).

This matrix is bipartite: $b_{ij}>0$ if and only if  $i$ is a source and $j$ a sink. 
We use $\circ$ (resp.\,$\bullet$) to denote an arbitrary source (resp.\,sink) of $\vec{\Delta}$ and $\mu_k$ stands for the mutation at $k$. 
Then  the two 'composite' mutations
\begin{equation}
\label{Eq:mu-circ--mu-bullet}
\mu_\circ=\prod_{k \mbox{ source}} \mu_k  
\qquad \mbox{ and } \qquad 
\mu_\bullet=\prod_{\ell \mbox{ sink}} \mu_\ell
\end{equation}
are well defined, involutive and $\mu_\circ(B)=\mu_\bullet(B)=-B$.  
In particular $\mu_\circ$ transforms the initial bipartite $\boldsymbol{\mathcal X}$-seed, denoted by $(\boldsymbol{y}_0,B)$ with $\boldsymbol{y}(0)=(y_1(0),\ldots,y_n(0))$, into a seed 
$(\boldsymbol{y}(1),-B)$ of the same type (that is bipartite) to which the same construction applies. This way, one defines the {\bf bipartite belt} 
as the family of $\boldsymbol{\mathcal X}$-seeds $S_t=(\boldsymbol{y}({t}),(-1)^{t}B)$ with 
$\boldsymbol{y}({t})=({y}_1({t}),\ldots,{y}_n({t}))$  for $t\in \mathbf Z$, defined inductively by 
$$
S_{t+1}=\Big(\boldsymbol{y}({t+1}),(-1)^{t+1}B\Big)= \mu_t\Big(\big(\boldsymbol{y}({t}),(-1)^{t}B\big)\Big)\qquad 
\mbox{with }\qquad  \mu_t=\begin{cases} {}^{}\hspace{0.2cm} \mu_\bullet \hspace{0.2cm}\mbox{if } t \mbox{ is even}\\
 {}^{}\hspace{0.2cm}
\mu_\circ \hspace{0.2cm}\mbox{otherwise}. 
\end{cases}
$$
Within this formalism, the half-periodicity property is that for a certain involution $\sigma$ of $\{1,\ldots,n\}$ depending on $\Delta$\footnote{\label{note1}The permutation $\sigma$ is given by   
$k\mapsto n+1-k$ in type $A_n$, is the transposition exchanging $1$ and $2$ in case $D_n$ with $n$ odd, is the product of transposistions $(16)(25)$ for $E_6$  and is trivial otherwise, see \cite[\S3.3]{IIKNS}.} the seed $S_{t+h+2}$ is $\sigma$-isomorphic to $S_t$ for any $t\in \mathbf Z$, and consequently in terms of cluster variables, one has $y_i(t+h+2)=y_{\sigma i}(t)$  for any $i=1,\ldots,n$.\sk 

As a consequence, there is only a finite number of cluster variables appearing in the bipartite belt which leads to consider the cluster web defined by them.  Let $\epsilon: \{1,\ldots,n\}\rightarrow \{\pm 1\}$ be the map  
given by $\epsilon(i)=1$ if $i$ is a source and -1 if it is a sink. Because 
$y_i(t+1)=y_i(t)^{-1}$ if $\epsilon(i)(-1)^t=1$, one can only consider the $y_i(t)$'s with $\epsilon(i)(-1)^t=-1$ in order to construct the web associated to the bipartite belt.  One constructs that way the {\bf bipartite belt cluster web} of type $\Delta$
\begin{equation}
\label{Eq:Y-web}
 \boldsymbol{ B \hspace{-0.05cm}B
 \hspace{-0.05cm}{\mathcal W}}_{\Delta}=
\boldsymbol{\mathcal W}\bigg( \hspace{0.1cm}
y_i(t)\hspace{0.15cm} \Big\lvert \, 
\begin{tabular}{l}
$i=1,\ldots,n$\\
$t=0,\ldots,h+1$, 
\end{tabular}
 \mbox{ with }\,\epsilon(i)(-1)^t=-1\, \bigg)
\end{equation}
that will be shown further equivalent to the  $\boldsymbol{\mathcal Y}$-cluster web $\boldsymbol{{\mathcal Y} \hspace{-0.05cm}{\mathcal W}}_{\Delta}$. 
Remark that $h+2$ is always even, except  when $\Delta=A_n$ with $n$ even, in which case $h=n+1$ is odd, hence so is $h+2$.  But in any case, $n(h+2)$ is  even as well and considering the condition $\epsilon(i)(-1)^t=-1$ appearing in the RHS of  \eqref{Eq:Y-web}, we deduce that  
$ \boldsymbol{{\mathcal Y}\!{\mathcal W}}_{\Delta}$ is defined by 
$n(h+2)/2$ cluster variables hence is  a $d'_{\Delta}$-web in $n$ variables, with $d'_\Delta\leq d_\Delta=n(h+2)/2$. Note that the majoration may be strict since one cannot exclude,  at this point,   that two of the $n(h+2)/2$ cluster first integrals of 
$ \boldsymbol{{\mathcal Y}\!{\mathcal W}}_{\Delta}$
do define the same foliation.
\begin{center}
$\star$\sk 
\end{center}

There is an equivalent way to construct the bipartite belt web, by means of a birational map of finite order, which may be useful  in practice.
 The periodicity (resp.\,half-periodicity) property of the $Y$-system of type $\Delta$ is equivalent to the fact that 
the birational map $u \mapsto 
F_{\Delta}(u)$ of $u=(u_i)_{i=1}^n\in \mathbf C^n$ induced by the action on the $\boldsymbol{\mathcal X}$-cluster variables by the successive  composition of  mutations $\mu_{\bullet}$ and $\mu_{\circ }$ 
 (thus one can write $
F_{\Delta}= \mu_{\circ } \circ \mu_{\bullet}
$ a bit abusively), is of order $m_\Delta=2(h+2)/\gcd (h,2)$ (resp. $m'_{\Delta}=(h+2)/\gcd (h,2)$). Then the web $
 \boldsymbol{ B \hspace{-0.05cm}B
 \hspace{-0.05cm}{\mathcal W}}_{\Delta}$ can be seen as the web on 
 $\mathbf C^n$ 
 defined by the components of the iterations $(F_{\Delta})^{\circ \ell}$ for $\ell=0,\ldots,m'_\Delta$, seen as rational functions in the 
 $u_i$'s.

\subparagraph{}\hspace{-0.6cm}
\label{Par:Delta-box-Delta'}
 We now discuss  the case when both $\Delta$ and $\Delta'$ have rank 2 or more. We recall some constructions of \cite{KellerAnnals} first when both considered Dynkin diagrams are simply laced, then we say a few words about the (more involved) general case. \sk 
 
 We assume that both $\Delta$ and $\Delta'$ are simply laced. 
 To simplify, we denote in the same way the associated bipartite quivers $\vec{\Delta}$ and $\vec{\Delta}'$.  
 In \cite[\S3.3]{KellerAnnals}, Keller defines the tensor product quiver $\Delta\otimes \Delta'$, from which he constructs two others, the square product $\Delta\square \Delta'$ and the triangle product $\Delta\boxtimes \Delta'$.  Each of these last two quivers could be used and the resulting webs are equivalent. However, typically $\Delta\square \Delta'$ is simpler than 
 $\Delta\boxtimes \Delta'$ (the former carries less arrows than the latter, see {\it e.g.}\,\cite[Fig.\,2]{KellerAnnals}) hence we will leave aside the cross product in the lines below, except for a few examples. \sk 
 
 The {\bf tensor product} $\boldsymbol{\Delta\otimes \Delta'}$ is the quiver whithout loops whose set of vertices is the set of pairs $(i,i')\in \Delta\times \Delta'$ where the number of arrows from $(i,i')$ to $(j,j')$ is equal to the number of arrows from $j$ to $j'$ if $i=i'$; equal to the number of arrows from $i$ to $i'$ if $j=j'$; and  equal to 0 otherwise.  The 
 {\bf square product} $\boldsymbol{\Delta\square \Delta'}$ is obtained from ${\Delta\otimes \Delta'}$ by reversing all arrows in the full subquivers of the form $\{i\}\times \Delta'$ and $\Delta\times \{i'\}$ where $i$ is a sink  and $i'$ a source (of  $\Delta$ and $\Delta'$ respectively).  The vertices of 
 $\Delta\square \Delta'$ are of two kinds: even vertices (pair of two sources or of two sinks), denoted with $\circ$, or odd vertices (mixed pairs), denoted by $\bullet$. An example of crossed product is given in Figure \ref{Fig:crossed-A4*D4} below. 
 \sk 
 
 \begin{figure}[h]
 \begin{center}
 \resizebox{3in}{1.5in}{
 \includegraphics{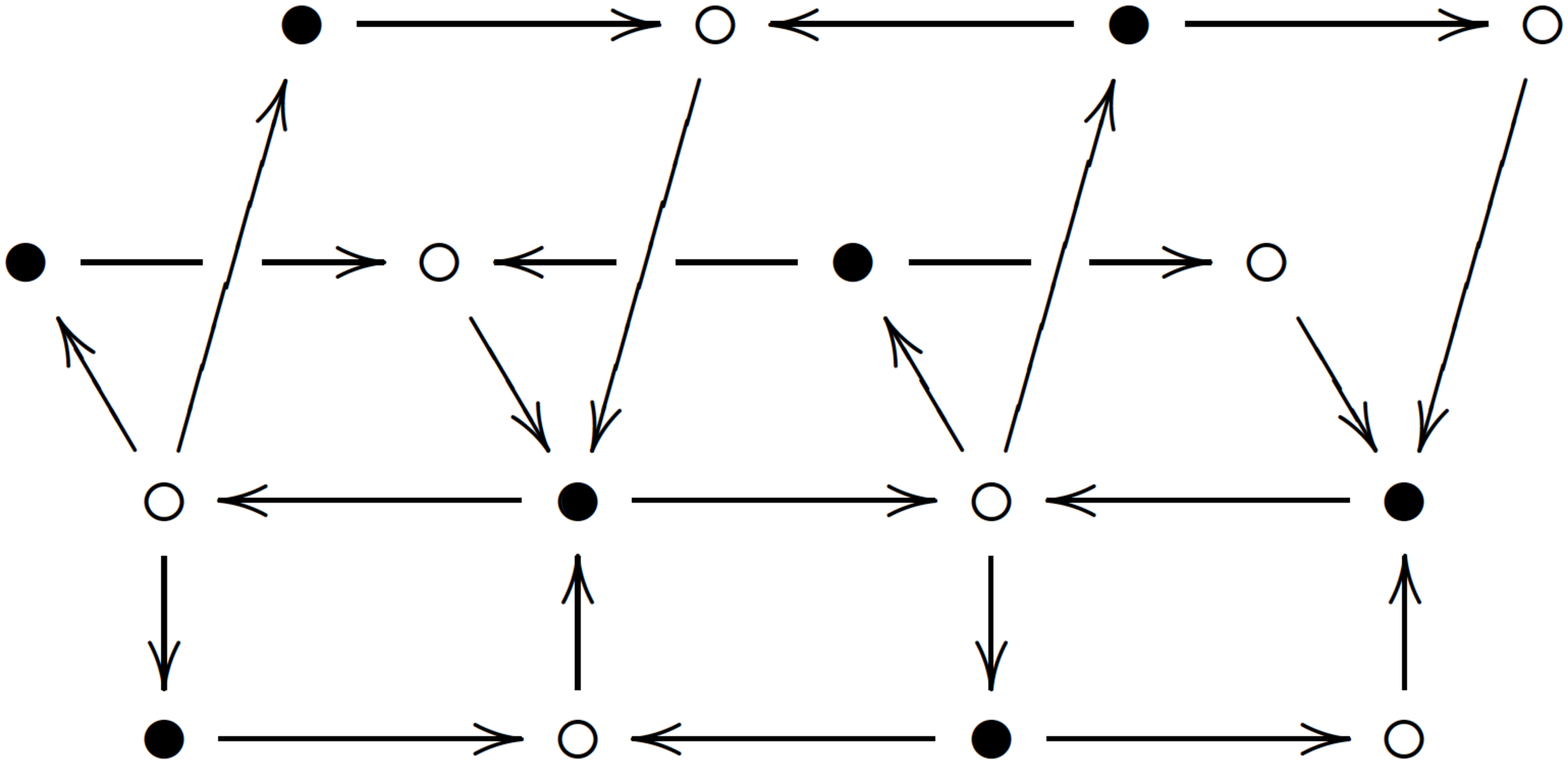}}  
 \vspace{-0.3cm}
 \end{center}
 \caption{The square product ${A}_4\square D_4$.}
  \label{Fig:crossed-A4*D4}
 \end{figure}


Similarly as in the previous case when $\Delta'=A_1$,  one denotes by 
 $\mu_\circ$ (resp.\,$\mu_\bullet$) the compositions of the mutations with respect to all even (resp.\,odd) vertices of $\Delta\square \Delta'$ (the orders are not important) and one considers 
 $
 \mu_{\Delta\square \Delta'}= \mu_\circ \circ \mu_\bullet
 $. The action on the $\boldsymbol{\mathcal X}$-cluster variables of  this composition of mutations corresponds to a birational map $F_{\Delta \square \Delta'}$ of $\mathbf C^{\Delta\times \Delta'}\simeq \mathbf C^{nn'}$.  Then as proved by Keller, the periodicity property for 
 the $Y$-system $Y(\Delta,\Delta')$ is equivalent to the fact that 
 $F_{\Delta \square \Delta'}$ is $m_{\Delta,\Delta'}=2(h+h')/\gcd(h,h')$-periodic as a rational map. As for the half-periodicity, 
  it corresponds to the fact that, setting 
  $m'_{\Delta,\Delta'}=(h+h')/\gcd(h,h')$, the birational map  $(F_{\Delta \square \Delta'})^{m'_{\Delta,\Delta'}}$ acts on $\mathbf C^{\Delta\times \Delta'}$ as the permutation 
$u=(u_{ii'}) \mapsto(u_{\sigma(i)\sigma'(i')})$, where $\sigma$ and $\sigma'$  stand respectively for  the permutations associated to $\Delta$ and $\Delta'$  described in footnote \ref{note1}.\sk 

 Then, although the term `bipartite' (and therefore the notation with two $B$'s) is no longer warranted we define the web 
\begin{equation}
\label{Eq:YDD'-web}
 \boldsymbol{ B \hspace{-0.05cm}B
 \hspace{-0.05cm}{\mathcal W}}_{\Delta\square\Delta'}=
\boldsymbol{\mathcal W}\Bigg( \hspace{0.1cm}
f(u)\hspace{0.15cm} \Big\lvert \, 
\begin{tabular}{l}
 $f(u)$ is a cluster variables appearing as a com- \\ ponent
 of $(F_{\Delta \square \Delta'})^\ell$ for some $\ell=0,\ldots, 
 m'_{\Delta,\Delta'}-1$
\end{tabular}
\Bigg)\, .
\end{equation}
 Because both  $\Delta$ and $\Delta'$ have been assumed to be simply laced,  it can be demonstrated that this web precisely coincides with the 
 $ \boldsymbol{ \mathcal Y}$-cluster web 
$ \boldsymbol{ \mathcal Y
 \hspace{-0.05cm}{\mathcal W}}_{\Delta\square\Delta'}$ associated to 
 the dilogarithmic identity $\mathcal R(\Delta,\Delta')$. 
 \vspace{-0.2cm}
\begin{center}
$\star$
\end{center}

We now quickly  discuss the general case when the considered Dynkin diagrams can have multiplicities higher than 1. We define two integers 
$m_{\Delta,\Delta'}$ and $m'_{\Delta,\Delta'}$ thanks to the same formulas than above. Following \cite[\S9.6]{KellerAnnals}, one  extends  the products of quivers discussed above to the case of valued quivers and one can then play the same game as  in the simply-laced case. Considering a certain composition of some mutations, one obtains a birational map $F_{\Delta\square  \Delta'}$ of order $m_{\Delta,\Delta'}$, which is equivalent to the periodicity of the associated $Y$-system. One can thus consider a web $\boldsymbol{ B \hspace{-0.05cm}B
 \hspace{-0.05cm}{\mathcal W}}_{\Delta\square\Delta'}$ defined as in 
\ref{Eq:YDD'-web} but it does not necessarily coincides with the suitable 
 $ \boldsymbol{ \mathcal Y}$-cluster web 
$ \boldsymbol{ \mathcal Y
 \hspace{-0.05cm}{\mathcal W}}_{\Delta\square\Delta'}$ in general (see \S\ref{SubPar:A2*B2} below for an explicit example).

\paragraph{Some explicit examples of  cluster webs associated to $\boldsymbol{Y}$-systems.} 
\label{Par:Exampels-Y-cluster-webs}

For any Dynkin diagram $\Delta$, it is clear that 
$\boldsymbol{\mathcal Y}\boldsymbol{\mathcal W}_{\Delta}$ is a subweb of the associated $\boldsymbol{\mathcal X}$-cluster web $\boldsymbol{\mathcal X}\boldsymbol{\mathcal W}_{\Delta}$. It is always a proper subweb, except in rank 2 where both webs coincide. For that reason, we do not consider $\boldsymbol{\mathcal Y}$-cluster webs associated to Dynkin diagrams of rank 2 here: these will be studied in depth further on ({\it cf.}\,\S\ref{SS:cluster-webs:rank2}). \sk 

Below, we discuss a few examples of $\boldsymbol{\mathcal Y}$-cluster webs from the point of view of their ARs and of their rank. As we will see, many of them, but not all, are AMP.

\subparagraph{The cluster web $\boldsymbol{\mathcal Y}\boldsymbol{\mathcal W}_{A_3}$.}
The  bipartite  Dynkin diagram of type $A_3$ is the following: 
%
%
$$\scalebox{0.4}{ \includegraphics{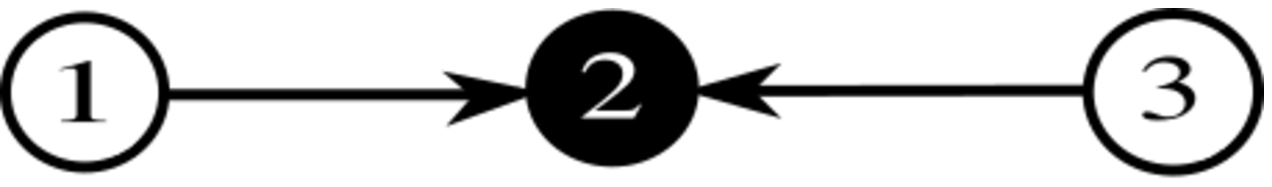}}$$

Using  $u=(u_1,u_2,u_3)$ as initial coordinates,  we obtain that, 
up to composition at the source and at the target by the 
birational involution, $u\mapsto (u_1^{-1},u_2,u_3^{-1})$ (see \eqref{Eq:Bir-change-E}),   the composition of the mutations of the  $\boldsymbol{\mathcal X}$-seeds corresponding  to the sink then followed by the ones associated to  the two sources is given by: 
\begin{align*}
F_{A_3}(u)=F_\circ \big( F_\bullet(u)\big)=& 
\left( \frac{1+u_2}{u_1}, \frac{(1+u_1+u_2)(1+u_2+u_3)}{u_1u_2u_3}  ,
\frac{1+u_2}{u_3}
\right)\, .
\end{align*}
One verifies that $(F_{A_3})^{\circ 3}(u)=(u_3,u_2,u_1)$ (half-periodicity) from which one deduces that 
 the $\boldsymbol{\mathcal Y}$-cluster web of type $A_3$ is  the following 9-subweb of the 
cluster 15-web 
  $\boldsymbol{\mathcal X \mathcal W}_{A_3}$  
\begin{align*}
 \boldsymbol{\mathcal Y \mathcal W}_{\hspace{-0.01cm}A_3}=
\boldsymbol{ \mathcal W}\Bigg( \, 
u_1\, , \, u_2\,  , \, u_3\, , \, 
{\frac {1+u_2}{u_1}}\, , \,   &
{\frac { \left( 1+u_1 \right)  \left( 1+u_3 \right) }{u_2}} \, , \, 
{\frac {1+u_2}{u_3}}\, , \, 
{\frac {1+u_1+u_2+u_3+u_1u_3}{u_1u_2}}\, , \, \\
 &
{\frac {1+u_1+u_2+u_3+u_1u_3}{u_2u_3}}\, , \, 
{\frac { \left( 1+u_1+u_2 \right)  \left( 1+u_2+u_3 \right) }{u_1u_2u_3}}
\hspace{0.15cm}
\Bigg)\, .
\end{align*}
  After some computations, one gets 
  $$\rho^\bullet( \boldsymbol{\mathcal Y \mathcal W}_{\hspace{-0.01cm}A_3})=(6,3,1)
  \qquad \mbox{ and }\qquad 
  {\rm polrk}^\bullet( \boldsymbol{\mathcal Y \mathcal W}_{\hspace{-0.01cm}A_3})=(9,1)$$ 
  from which it comes that this web is AMP with only polylogarithmic ARs of weight 1 and 2. 

\subparagraph{An aside: the cluster web $\boldsymbol{\mathcal A \hspace{-0.05cm}\mathcal W}_{\hspace{-0.05cm} A_3}$.}
\label{subpar:AWA3}
It is interesting to compare $\boldsymbol{\mathcal Y \hspace{-0.05cm}\mathcal W}_{\hspace{-0.05cm} A_3}$  with the $\boldsymbol{\mathcal A }$-cluster web of type $A_3$, which is the web 
\begin{align*}
\boldsymbol{\mathcal A \hspace{-0.05cm}\mathcal W}_{\hspace{-0.05cm} A_3}= 
\boldsymbol{\mathcal W}
 \Bigg( \, {a_{{1}}} \,,\, {a_{{2}}}\, ,\, {a_{{3}}}  \,,\,   {\frac {1+a_{{2}}}{a_{{1}}}}  \,,\,  {\frac {1+a_{{2}}}{a_{{3}}}}  \,,\,  {\frac {1+a_{{1}}a_{{3}}}{a_{{2}}}}  \,,\,  &  {\frac {
 1+a_{{2}}+
 a_{{1}}a_{{3}}}{a_{{1}}a_{{2}}}} \,,\,    \\ 
 &  {\frac {
  1+a_{{2}}+
 a_{{1}}a_{{3}}
 }{a_{{2}}a_{{3}}}}  \,,\,  
 {\frac {
 (1+{a_{{2}}})^{2}+
 a_{{1}}a_{{3}}}{a_{{1}}a_{{2}}a_{{3}}}} \Bigg) \nonumber 
\, . \nonumber 
\end{align*}
The two webs $\boldsymbol{\mathcal Y \hspace{-0.05cm}\mathcal W}_{\hspace{-0.05cm} A_3}$ and $\boldsymbol{\mathcal A \hspace{-0.05cm}\mathcal W}_{\hspace{-0.05cm} A_3}$ are very similar in some ways: both have degree 9 and are defined by Laurent polynomials the monomial denominators of which can be put in bijection with the set of almost-positive roots 
of the root system of type $A_3$ ({\it cf.}\;\S\ref{SS:log-ARs-YW-An} or 
\S\ref{SS:Cluster notation} further). Moreover, one has 
$$\rho^\bullet\big( \boldsymbol{\mathcal A \hspace{-0.05cm}\mathcal W}_{\hspace{-0.05cm} A_3}\big)=\rho^\bullet\big( \boldsymbol{\mathcal Y \hspace{-0.05cm}\mathcal W}_{\hspace{-0.05cm} A_3}\big)=(6,3,1)\, .$$

However, despite these similarities, both webs are quite different in what concerns their ARs and their ranks. Indeed, by direct computations, we have established that contrarily to $\boldsymbol{\mathcal Y \hspace{-0.05cm}\mathcal W}_{\hspace{-0.05cm} A_3}$, the  web $\boldsymbol{\mathcal A \hspace{-0.05cm}\mathcal W}_{\hspace{-0.05cm} A_3}$: 
\begin{itemize}
\item    does not carry any dilogarithmic AR but only logarithmic ones. More precisely, one has $$ {\rm polrk}^\bullet\big(\boldsymbol{\mathcal A \hspace{-0.05cm}\mathcal W}_{\hspace{-0.05cm} A_3}\big)=(5,0)\, ;$$ 
\item does carry a rational (hence non-polylogarithmic) abelian relation, which is the one corresponding to the following functional identity satisfied by the function ${\sf J}: u\mapsto 1/(1+u)$: 
\begin{equation}
\label{Eq:JOKO}
{\sf J}(a_2)+{\sf J}\left( {\frac {1+a_{{1}}a_{{3}}}{a_{{2}}}} \right)
+{\sf J}\left( {\frac {
 (1+{a_{{2}}})^{2}+
 a_{{1}}a_{{3}}}{a_{{1}}a_{{2}}a_{{3}}}}
\right)\equiv 1\, ; 
\vspace{-0.3cm}
\footnotemark
\end{equation}
\footnotetext{Remark that the 3-web defined by the first integrals appearing in \eqref{Eq:JOKO}  has intrinsic dimension 2 hence this identity is the pull-back (under the map $(a_1,a_2,a_3)\mapsto (a_2,a_1a_3)$ to be explicit) of a functional relation in two variables.}
\item 
has rank ${\rm rk}\big(\boldsymbol{\mathcal A \hspace{-0.05cm}\mathcal W}_{\hspace{-0.05cm} A_3}\big)=6<10=\rho\big(\boldsymbol{\mathcal A \hspace{-0.05cm}\mathcal W}_{\hspace{-0.05cm} A_3}\big)$ hence is not AMP.
\end{itemize}
\sk

To be honest, we have no idea why $\boldsymbol{\mathcal Y \hspace{-0.05cm}\mathcal W}_{\hspace{-0.05cm} A_3}$ is much more interesting than 
$\boldsymbol{\mathcal A \hspace{-0.05cm}\mathcal W}_{\hspace{-0.05cm} A_3}$ regarding their ARs, their ranks and the property of being AMP. However other investigations have led us to be quite convinced that webs defined by 
$\boldsymbol{\mathcal A}$-cluster variables do not carry many abelian relations 
and are not those to be considered in view of the study undertaken in this text.  This explains why, up to very few exceptions, we are not dealing with webs defined by 
$\boldsymbol{\mathcal A}$-cluster variables
in the whole text.


\subparagraph{The cluster web $\boldsymbol{\mathcal Y}\boldsymbol{\mathcal W}_{D_4}$.}
The  quiver $\vec{D}_4$ and the associated exchange matrix 
are: \sk 

\begin{tabular}{ccc} ${}^{}$ \qquad 
\scalebox{0.38}{ \includegraphics{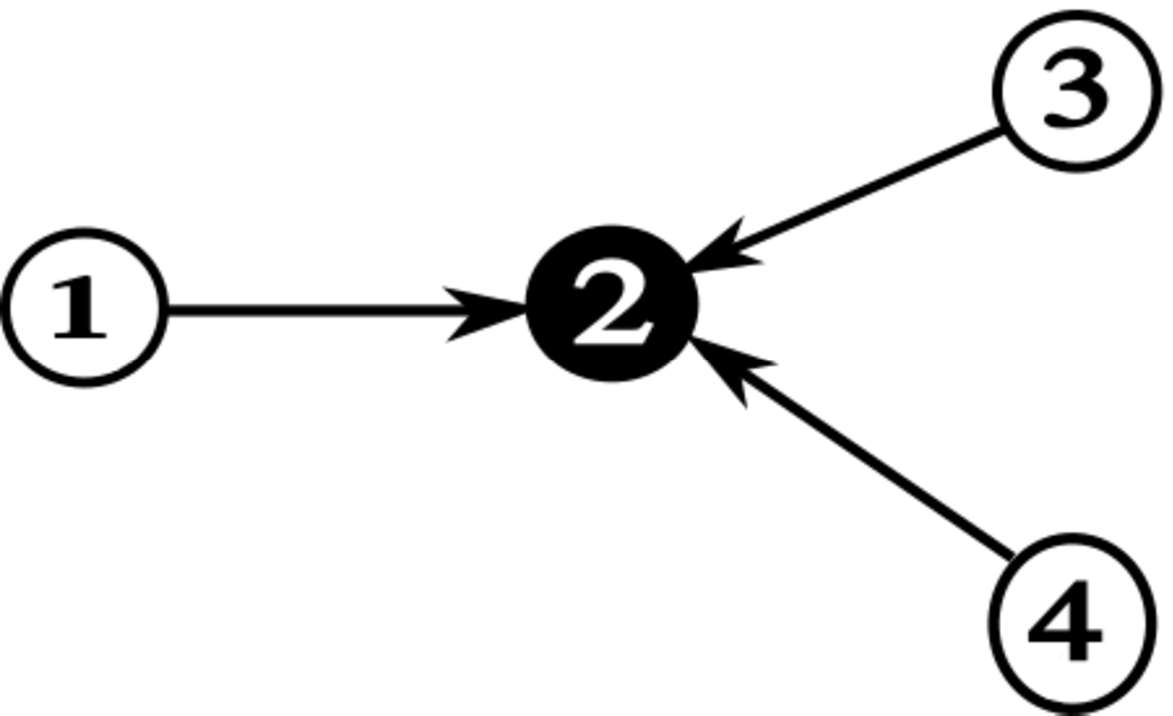}} &  $\qquad$  &  
$\begin{tabular}{c} 
${}^{}$ \vspace{-3cm}\\ 
 and  $\quad \qquad$ $B_{D_4}=\begin{bmatrix}
0 & 1 & 0 & 0 \\
0 & 1 & 0 & 0\\ 
-1 & 0 & -1 & -1 \\
0 & 1 & 0 & 0
\end{bmatrix}$\, .
\end{tabular}$\sk 
\end{tabular}

Using  $u=(u_i)_{i=1}^4$ as initial coordinates,  we obtain that, 
up to conjugation by the 
birational involution $u\mapsto (u_1^{-1},u_2,u_3^{-1},u_4^{-1})$, 
 the composition of the mutations of the  $\boldsymbol{\mathcal X}$-seeds corresponding  to the sources followed by the one associated to  the sink is given by: 
$$
F_{D_4}(u)=F_\circ \big( F_\bullet(u)\big)=\left(
{\frac {1+u_{{2}}}{u_{{1}}}},{\frac { \left( 1+u_{{2}}+u_{{3}} \right)  \left(1+ u_{{1}}+u_{{2}} \right) \\
\mbox{} \left( 1+u_{{2}}+u_{{4}} \right) }{u_{{1}}u_{{2}}u_{{3}}u_{{4}}}},{\frac {1+u_{{2}}}{u_{{3}}}},{\frac {1+u_{{2}}}{u_{{4\\
\mbox{}}}}}\right)
$$

One verifies that $(F_{D_4})^{\circ 4}={\rm Id}_{\mathbf C^4}$ (half-periodicity) and  that the associated web $ \boldsymbol{B\hspace{-0.04cm}B\hspace{-0.04cm}\mathcal W}_{\hspace{-0.05cm}D_4}$, which coincides with 
 the $\boldsymbol{\mathcal Y}$-cluster web of type $D_4$,  is  a 16-subweb of the 
$\boldsymbol{\mathcal X}$-cluster web 
  $\boldsymbol{\mathcal X \mathcal W}_{D_4}$ which can be explicited. One has:
  \begin{align*}
 \boldsymbol{\mathcal Y \mathcal W}_{\hspace{-0.01cm}D_4}=
\boldsymbol{ \mathcal W}\Bigg( \, 
& u_1\, , \, u_2\,  , \,  u_3\, , \, u_4\, , \, 
{\frac {1+u_{{2}}}{u_{{1}}}}
\, , \,  
{\frac {(1+u_1)(1+u_3)(1+u_4)}{u_{{2}}}}
\, , \,  
{\frac {1+u_{{2}}}{u_{{3}}}}
\, , \,  
{\frac {1+u_{{2}}}{u_{{4}}}}
\, , \,  \ldots
 \\ 
& \ldots\, , \, 
\frac { \left(  \left( 1+u_{{1}} \right) u_{{4}}+u_{{1}}+u_{{2}}+1 \right)  \left(  \left( 1+u_{{1}} \right) u_{{3}}+u_{{1}}+u_{{2}}+1 \right)  \left(  \left( 1+u_{{3}} \right) u_{{4}}+u_{{3}}+u_{{2}}+1 \right) }{u_1u_2^{2} u_3u_4}
\Bigg)
\end{align*}
Some computations give us  
$\rho^\bullet( \boldsymbol{\mathcal Y \mathcal W}_{\hspace{-0.01cm}D_4})=(12,6,1)$, ${\rm polrk}^\bullet( \boldsymbol{\mathcal Y \mathcal W}_{\hspace{-0.01cm}D_4})=(16,1)$ and 
\begin{equation}
\label{Eq:YWD4-invariants}
 {\rm rk}\Big( \boldsymbol{\mathcal Y \mathcal W}_{\hspace{-0.01cm}D_4}\Big)=
{\rm polrk}\Big( \boldsymbol{\mathcal Y \mathcal W}_{\hspace{-0.01cm}D_4}\Big)=17 < 19 =  \rho( \boldsymbol{\mathcal Y \mathcal W}_{\hspace{-0.01cm}D_4}) 
\end{equation}
hence this web only carries polylogarithmic ARs (of weight 1 and 2) but is not AMP (contrarily to what one may naively expect).

\subparagraph{The web $\boldsymbol{\mathcal Y}\boldsymbol{\mathcal W}_{A_2\square A_2}$.} The 
quiver $A_2\square A_2$ is
$$\scalebox{0.4}{ \includegraphics{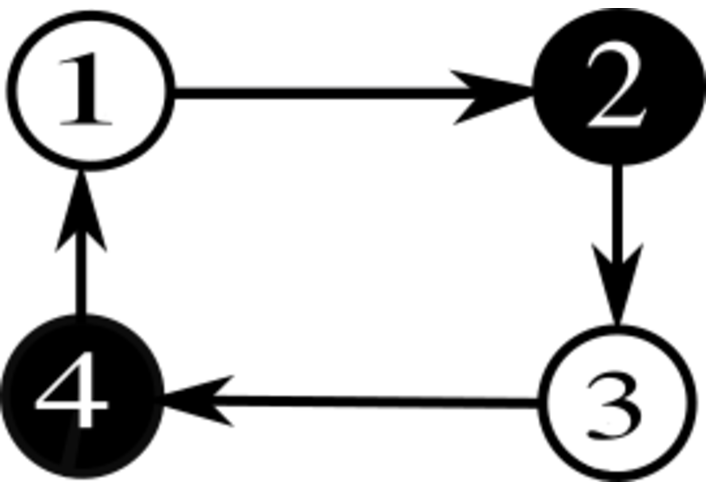}}$$
from which it is easy to construct the associated exchange matrix $B_{A_2\square A_2}$ hence the birational map $F_{A_2\square A_2}$ which is such that $(F_{A_2\square A_2})^{\circ 3}$ acts as a permutation on $\mathbf C^4$.  One verifies that  the web $ \boldsymbol{B\hspace{-0.04cm}B\hspace{-0.04cm}\mathcal W}_{\hspace{-0.05cm}{A_2\square A_2}}$  
coincides with $\boldsymbol{\mathcal Y}\boldsymbol{\mathcal W}_{A_2\square  A_2}$ and is a 12-web in four variables. 
By direct computations, one obtains that  
$$\rho^\bullet\Big( \boldsymbol{\mathcal Y \mathcal W}_{\hspace{-0.01cm}A_2 \square A_2}\Big)=(8,4,1)\qquad 
\mbox{ and }\qquad 
{\rm polrk}^\bullet\Big( \boldsymbol{\mathcal Y \mathcal W}_{\hspace{-0.01cm}A_2\square  A_2}\Big)=(12,1) \, . $$ 
It follows that $\boldsymbol{\mathcal Y}\boldsymbol{\mathcal W}_{A_2 \square A_2}$ is AMP with all its ARs polylogarithmic (all logarithmic, except the dilogarithmic one associated to $\mathcal R(A_2,A_2)$). \sk

\begin{rem} It is known that the $D_4$ quiver is mutation equivalent to the quiver 
$A_2\square A_2$ hence the two corresponding cluster algebras are mutation equivalent as well. Despite this, the two $\mathcal Y$-cluster webs $\boldsymbol{\mathcal Y}\boldsymbol{\mathcal W}_{\hspace{-0.05cm} D_4}$ and 
$\boldsymbol{\mathcal Y}\boldsymbol{\mathcal W}_{\hspace{-0.05cm}A_2\square  A_2}$
 are not equivalent (they are not of the same degree) which might be a  bit surprising at first. 
Actually,  there is no inconsistency here and this can be understood in terms of the notion of `cluster period' which will be discussed below {\rm (}{\it cf.}\,\S\ref{Par:Periods}{\rm )}. Indeed,  both webs can be understood within any one of the two involved (isomorphic) cluster algebras: within the one associated to $D_4$ say,  the cluster web $\boldsymbol{\mathcal Y}\boldsymbol{\mathcal W}_{\hspace{-0.05cm} D_4}$ corresponds  to the cluster period of length 16 associated to the $Y$-system of type $D_4$ whereas $\boldsymbol{\mathcal Y}\boldsymbol{\mathcal W}_{\hspace{-0.05cm} A_2\square  A_2}$ corresponds to a different cluster period for this cluster algebra, of length 12. 
\end{rem}

\subparagraph{The web $\boldsymbol{\mathcal Y}\boldsymbol{\mathcal W}_{A_2\square D_4}$.} The 
square product quiver $A_2\square D_4$ is
$$\scalebox{0.4}{ \includegraphics{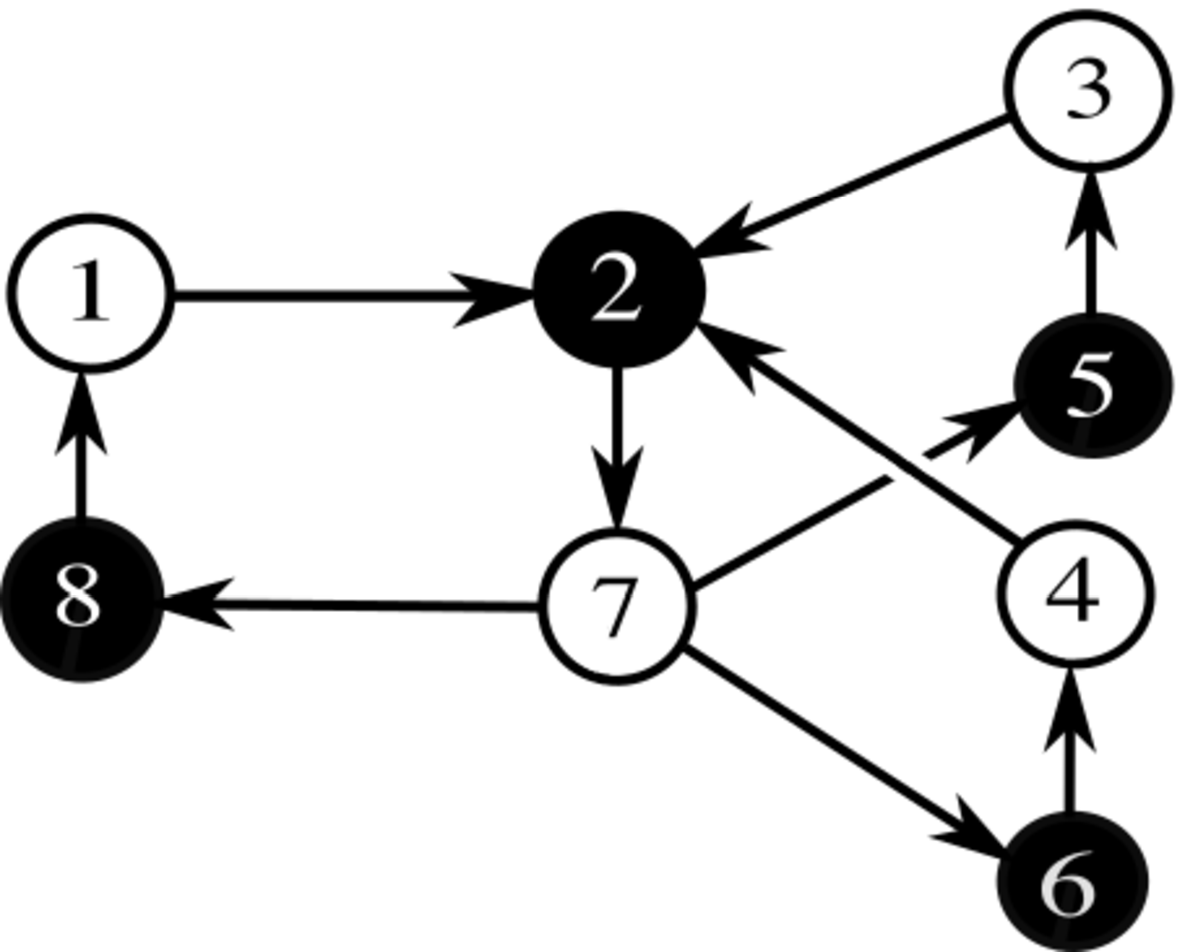}}$$
from which it is easy to construct the associated exchange matrix $B_{A_2\square D_4}$ hence the birational map $F_{A_2\square D_4}: \mathbf C^8\dashrightarrow \mathbf C^8$ which can be verified of 
order $h(A_2)+h(D_4)=3+(2\times 4-2)=9$. \sk 

One verifies that the web $ \boldsymbol{B\hspace{-0.04cm}B\hspace{-0.04cm}\mathcal W}_{\hspace{-0.05cm}{A_2\square D_4}}$  
coincides with $\boldsymbol{\mathcal Y}\boldsymbol{\mathcal W}_{A_2\square  D_4}$. It 
is a 36-web in 8 variables and by direct computations, one obtains that  
$$\rho^\bullet( \boldsymbol{\mathcal Y \mathcal W}_{\hspace{-0.01cm}A_2 \square D_4})=(28,8,1)\qquad 
\mbox{ and }\qquad 
{\rm polrk}^\bullet( \boldsymbol{\mathcal Y \mathcal W}_{\hspace{-0.01cm}A_2\square  D_4})=(36,1) \, . $$ 
It follows that $\boldsymbol{\mathcal Y}\boldsymbol{\mathcal W}_{A_2 \square D_4}$ is AMP with all its ARs polylogarithmic (all logarithmic, except a dilogarithmic one). It is surprising that $ \boldsymbol{\mathcal Y \mathcal W}_{\hspace{-0.01cm}A_2
\square
 D_4}$ is AMP whereas the a priori simplest web 
$ \boldsymbol{\mathcal Y \mathcal W}_{\hspace{-0.01cm}D_4}$ is not. If we believe that it is a general 
phenomenon (see  \S\ref{SS:BothDynkinDiagramsClassical} further and in particular our conjecture there), we do not have any explanation of that  to offer for now.

\subparagraph{The web $\boldsymbol{B\hspace{-0.04cm}B\hspace{-0.04cm}\mathcal W}_{A_2\boxtimes B_2}$.}
\label{SubPar:A2*B2}
In this case and in the next one, we deal with quivers obtained by means of the triangle product $\boxtimes$. Using $\square$ instead would have led to similar results.\mk

The weighted quiver of type $A_2\boxtimes B_2$ 
and the associated exchange matrix  are the following: 

\vspace{-1.7cm}
\begin{tabular}{lcr}
${}^{}$ \quad 
  \scalebox{1.1}{
    \xymatrix@R=1.5cm@C=2cm{
 \scalebox{1.4}{{\textcircled{\raisebox{-0.9pt}{\scalebox{0.8}{2}}}}}_{-+}
  \ar@{->}[r]^{ (1,1)}
& 
\ar@{->}[ld]^{ (1,2)}
 \scalebox{1.4}{{\textcircled{\raisebox{-0.9pt}{\scalebox{0.8}{4}}}}}_{--}
 \\  \ar@{->}[u]^{ (2,1)}
  \ar@{->}[r]_{ (1,1)}
\scalebox{1.4}{{\textcircled{\raisebox{-0.9pt}{\scalebox{0.8}{1}}}}}_{++} &
\scalebox{1.4}{{\textcircled{\raisebox{-0.9pt}{\scalebox{0.8}{3}}}}}_{+-}
 \ar@{->}[u]_{ (2,1)}
}
} &  $\quad$  &  
$\begin{tabular}{c} 
${}^{}$ \vspace{2cm}\\ 
 and  $\quad \qquad$ $B_{A_2\boxtimes B_2}=\begin{bmatrix}
0 & 2 & 1 & -2 \\
-1 & 0 & 0 & 1\\ 
-1 & 0 & 0 & 2 \\
1 & -1 & -1 & 0
\end{bmatrix}$\, .
\end{tabular}$\sk 
\end{tabular}

The associated map $F_{A_2\boxtimes B2}$ is given 
by $u=(u_i)_{i=1}^4\mapsto (f_i(u))_{i=1}^4$ with  
\begin{align*}
f_1(u)=& 
 \, {\frac {u_{{1}}u_{{2}} \left( u_{{4}}{u_{{3}}}^{2}+{u_{{3}}}^{2}+2\,u_{{3}}+1 \right) }{ \left( 1+u_{{3}} \right) \Phi_{A_2\boxtimes B2}(u) }}\,
&& f_2(u)=
 {\frac { \left( u_{{1}}u_{{3}}+u_{{1}}+1 \right) ^{2}}{{u_{{1}}}^{2}u_{{2}} \left( u_{{4}}{u_{{3}}}^{2}+{u_{{3}}}^{2}+2\,u_{{3}}+1 \right) }} 
 \\
f_3(u)=&   {\frac {u_{{4}}u_{{3}} \left( u_{{1}}u_{{3}}+u_{{1}}+1 \right) }{u_{{4}}{u_{{3}}}^{2}+{u_{{3}}}^{2}+2\,u_{{3}}+1}}
&& f_4(u)=
    {\frac { \left( 1+u_{{3}} \right) ^{2} \Phi_{A_2\boxtimes B2}(u) }{u_{{4}}{u_{{3}}}^{2} \left( u_{{1}}u_{{3}}+u_{{1}}+1 \right) ^{2}}}
\end{align*}
where $\Phi_{A_2\boxtimes B2}(u)$ stands for the following polynomial expression in the $u_i$'s: 
%
 $$\Phi_{A_2\boxtimes B2}(u)= {u_{{1}}}^{2}u_{{2}}{u_{{3}}}^{2}u_{{4}}+{u_{{1}}}^{2}u_{{2}}{u_{{3}}}^{2}+2\,{u_{{1}}}^{2}u_{{2}}u_{{3}}+{u_{{1}}}^{2}{u_{{3}}}^{2}+{u_{{1}}}^{2}u_{{2}}+2\,{u_{{1}}}^{2}u_{{3}}+{u_{{1}}}^{2}+2\,u_{{1}}u_{{3}}+2\,u_{{1}}+1 \, .$$

One verifies that the birational map $F_{A_2\boxtimes B_2}$ has 
order $h(A_2)+h(B_2)=3+4=7$  and that the associated  cluster web 
$\boldsymbol{B\hspace{-0.04cm}B\hspace{-0.04cm}\mathcal W}_{A_2\boxtimes B_2}$
is a 28-web in four variables such that 
$$\rho^\bullet( \boldsymbol{B\hspace{-0.04cm}B\hspace{-0.04cm}\mathcal W}_{\hspace{-0.01cm}A_2\boxtimes B_2})=   (24,18,8),\quad  
{\rm polrk}^\bullet( \boldsymbol{B\hspace{-0.04cm}B\hspace{-0.04cm}\mathcal W}_{\hspace{-0.01cm}A_2\boxtimes B_2})=  (28,2) \qquad \mbox{ and } \qquad  
{\rm rk}(\boldsymbol{B\hspace{-0.04cm}B\hspace{-0.04cm}\mathcal W}_{\hspace{-0.01cm}A_2\boxtimes B_2})
 =   30.$$ 
 Thus $\boldsymbol{B\hspace{-0.04cm}B\hspace{-0.04cm}\mathcal W}_{A_2\boxtimes B_2}$ is not AMP but all its ARs are polylogarithmic, of weight 1 or 2.  We will see further that this web is distinct from the corresponding $\boldsymbol{\mathcal Y}$-cluster web 
 $\boldsymbol{\mathcal Y\hspace{-0.04cm}\mathcal W}_{A_2\boxtimes B_2}$ which is AMP.

\subparagraph{The web $\boldsymbol{B\hspace{-0.04cm}B\hspace{-0.04cm}\mathcal W}_{B_2\boxtimes B_2}$.}
\label{SubPar:B2*B2}
The weighted quiver of type $B_2\boxtimes B_2$ 
and the associated exchange matrix  are the following (see p.\,162 in \cite{KellerAnnals}): \vspace{-2cm}
\begin{equation}
\label{Eq:kotok}
\begin{tabular}{lcr}
  \scalebox{1.1}{
    \xymatrix@R=1.5cm@C=2cm{
 \scalebox{1.4}{{\textcircled{\raisebox{-0.9pt}{\scalebox{0.8}{2}}}}}_{-+}
  \ar@{->}[r]^{ (2,1)}
& 
\ar@{->}[ld]^{ (1,4)}
 \scalebox{1.4}{{\textcircled{\raisebox{-0.9pt}{\scalebox{0.8}{4}}}}}_{--}
 \\  \ar@{->}[u]^{ (2,1)}
  \ar@{->}[r]_{ (2,1)}
\scalebox{1.4}{{\textcircled{\raisebox{-0.9pt}{\scalebox{0.8}{1}}}}}_{++} &
\scalebox{1.4}{{\textcircled{\raisebox{-0.9pt}{\scalebox{0.8}{3}}}}}_{+-}
 \ar@{->}[u]_{ (2,1)}
}
} &   &  
$\begin{tabular}{c} 
${}^{}$ \vspace{2cm}\\ 
 and  $ \qquad$ $B_{B_2\boxtimes B_2}=\begin{bmatrix}
0 & 2 & 2 & -4 \\
-1 & 0 & 0 & 2\\ 
-1 & 0 & 0 & 2 \\
1 & -1 & -1 & 0
\end{bmatrix}$\, .
\end{tabular}$\sk 
\end{tabular}
\end{equation}

As in the preceding case of $A_2\boxtimes B_2$, the map $F_{B_2\boxtimes B_2}$ (corresponding to the composition of mutations $\mu_\boxtimes=\mu_{-+}\circ \mu_{++} \circ \mu_{--} \circ \mu_{+-} $, {\it cf.}\,formula (3.6.1) in \cite{KellerAnnals}), can be made explicit and one verifies that, as a birational map of $\mathbf C^4$, it is of order 
$(h(B_2)+h(B_2))/2=h(B_2)=4$
 (half-periodicity) and that the associated 
$\boldsymbol{\mathcal Y}$-cluster web 
$\boldsymbol{B\hspace{-0.04cm}B\hspace{-0.04cm}\mathcal W}_{B_2\boxtimes B_2}$ is a 16-web in four variables such that 
$$\rho^\bullet\Big( \boldsymbol{B\hspace{-0.04cm}B\hspace{-0.04cm}\mathcal W}_{\hspace{-0.01cm}B_2\boxtimes B_2}\Big)=   (12,6),\quad  
{\rm polrk}^\bullet\Big(\boldsymbol{B\hspace{-0.04cm}B\hspace{-0.04cm}\mathcal W}_{\hspace{-0.01cm}A_2\boxtimes B_2}\Big)=  (10)\qquad  \mbox{and}\qquad  
{\rm rk}\Big(\boldsymbol{B\hspace{-0.04cm}B\hspace{-0.04cm}\mathcal W}_{\hspace{-0.01cm}B_2\boxtimes B_2}\Big)
=
    10.$$ 
Thus this web is not AMP,  all its ARs are logarithmic and it does not carry any dilogarithmic AR.  This contrats with the suitable $\boldsymbol{\mathcal Y}$-cluster of type $(B_2, B_2)$: we will see in \S\ref{Subpar:YWB2*B2} below that 
$\boldsymbol{\mathcal Y\hspace{-0.04cm}\mathcal W}_{\hspace{-0.01cm}B_2\boxtimes B_2}$ is AMP with only logarithmic ARs, except one which corresponds to   $\mathcal R(B_2,B_2)$.
\subsubsection{Periods of cluster algebras and webs associated to them.}
\label{Par:Periods}
We now turn to the notion of period within the theory of cluster algebras. 
It has been formally introduced by Nakanishi who has proved very interesting results associating a dilogarithmic identity to any cluster period. We briefly review this material here before explaining 
how it leads to consider new cluster webs of interest, in particular to  the $\boldsymbol{\mathcal Y}$-cluster webs associated to $Y$-systems previously discussed. \mk 

Our main reference for this subsection is Nakanishi's paper 
\cite{Nakanishi}. 

 \paragraph{Periods.}
 Let $S=(\boldsymbol{a},\boldsymbol{x},B)$ be a seed of a cluster algebra $\mathcal A$ of rank $n\geq 2$.  Given an integer $k\geq 2$, we consider 
 a $k$-tuple $\boldsymbol{i}=(i_k,\ldots,i_1)$ of consecutive distinct elements of $I= \{1,\ldots,n\}$ and, unsurprisingly, we define the associated (composition of) mutation(s) by $\mu_{\boldsymbol{i}}=\mu_{i_k}\circ \cdots \circ \mu_{i_1}$. For a permutation $\nu\in \mathfrak S_n$, we will say that $\boldsymbol{i}$ is a {\bf $\boldsymbol{\nu}$-period} for the seed $S$ if 
 $$
 \mu_{\boldsymbol{i}}(S)=S^\nu\, .
 $$

The interpretation in terms of cluster algebras of the periodicity property of a  $Y$-system considered in the preceding subsection gives us some examples of periods. But this notion is more general, as we will see below by considering  some examples of different kinds.\mk

 \paragraph{Nakanishi's dilogarithmic identity associated to a period.}
 \label{Par:NakanishiDilogarithmicIdentity-Period}
The main interest of the notion  of period is that, thanks to a very nice result of Nakanishi (fully proved by him in \cite{Nakanishi} only in the particular case of seeds with a skew-symmetric exchange matrix), there is a dilogarithm identity, hence a possibly interesting web,  associated to any cluster period. 
\sk 

We assume that $\boldsymbol{{i}}=(i_k,\ldots,i_1)$ is a $\nu$-period of length $k$ 
for a seed $S_0=(\boldsymbol{a}_0,\boldsymbol{x}_0,B_0)$.

In order to state Nakanishi's result, we need to introduce some notations.\footnote{Beware that the notation we use differ slightly from the one used in \cite[\S6.5]{Nakanishi}.} 
\begin{itemize}
\item For $\ell=1,\ldots,k$, one defines $S_\ell=\mu_{i_\ell}\circ \cdots \circ \mu_{i_1}(S_0)=(\boldsymbol{a}_\ell,\boldsymbol{x}_\ell,B_\ell)$
 and one denotes by $x_\ell(i)$ (with $i=1,\ldots,n$) the $\boldsymbol{\mathcal X}$-cluster variables of this seed, 
  seen as rational functions in the cluster variables  $x_0(1),\ldots,x_0(n)$ of the initial seed $S_0$.

Then one defines the {\bf $\boldsymbol{\ell}$-th cluster coordinates 
 $\boldsymbol{x_\ell({\boldsymbol{i}})}$} associated to the $\nu$-period $\boldsymbol{i}$ as 
the $i_\ell$ component of the $\boldsymbol{\mathcal X}$-cluster $\boldsymbol{x}_\ell$, {\it i.e.}
\begin{equation}
\label{Eq=l-th-cluster-variable-period}
x_\ell({\boldsymbol{i}})=x_\ell(i_\ell)\, .
\end{equation}
 \item Let $D$ be the right skew symmetrizer of $B_0$: it is the diagonal matrix 
$D={\rm diag}(d_1,\ldots,d_n)$ with $d_i\in \mathbf N^*$ for $i=1,\ldots,n$ and ${\rm lcm}(d_i)$ minimal, such that $B_0D$ is skew-symmetric. 
\sk
\item  A $\boldsymbol{\mathcal X}$-cluster variable $x$ is said to be {\bf negative} (notation $x<0$)  if 
all the non-zero entries of its  $c$-vector are negative (beware that this notion makes sense only with respect to a fixed initial seed).\footnote{One defines the notion of {\it positive} $\boldsymbol{\mathcal X}$-cluster variable similarly, by replacing `negative' by positive. Then the `sign-coherence' mentioned in \S\ref{SubPar:Fpoly-and-c-vectors} can be stated as the fact that any $\boldsymbol{\mathcal X}$-cluster variable is either negative or positive.}  Then one defines 
$$
N_{{\boldsymbol{i}}}=\sum_{
\substack{\ell=1,\ldots,k,  \\  
x_\ell({\boldsymbol{i}})
<0
 }
} d_{i_\ell} \in \mathbf N_{>0}\, . 
$$
\sk
\end{itemize}
Then with the preceding notation, Nakanishi's result states as follows: 
\begin{thm}
\label{Thm:Nakanishis-identity}
The following dilogarithmic identity holds true
identically: 
 $$
\big(\boldsymbol{\mathcal R_{\boldsymbol{i}}}\big)
\hspace{3.5cm}
\sum_{\ell=1}^k d_{i_\ell} \, \mathcal R\left( 
\frac{x_\ell({\boldsymbol{i}})}{1+x_\ell({\boldsymbol{i}})}
\right)= \frac{{}^{}\hspace{0.09cm}\pi^2}{6}N_{{\boldsymbol{i}}}\hspace{0.1cm}.
  \hspace{7cm}
  $$
\end{thm} 
 
 A few remarks are in order about this statement: 
 \begin{enumerate}
 \item  In \cite{Nakanishi}, only the case when the initial exchange matrix  is skew-symmetric is fully proved. The more general statement above 
is only established (at the very end, in \S6.5) under the assumption that a certain conjecture (Conjecture 2.1 in  \cite{Nakanishi}) is satisfied. 
At the time when Nakanishi wrote his article, this conjecture (which is essentially equivalent to the `sign-coherence' of the $c$-vectors of 
$\boldsymbol{\mathcal X}$-cluster variables) was known to hold true only in the skew-symmetric case. However, thanks to 
recent advances in the theory of cluster algebras, mainly obtained  in 
\cite{GHKK}, we now know that the aforementioned conjecture holds true in full generality, which allows us to state the preceding result in all cases. 
 \item  Actually, the results given in Nakanishi's paper corresponding to the  theorem just above (namely Theorem 6.4 and Theorem 6.8 in \cite{Nakanishi}) are stated in a more general form, 
 with respect to what Nakanishi's calls a `{\it slice}' of the considered period $\boldsymbol{i}$. Since this basically does not change the associated dilogarithmic identity (hence the web associated to it, see \eqref{Eq:PWeb-i} below),   we have chosen to simplify the presentation as much as possible  
  by only considering the maximal slice of the period $\boldsymbol{i}$. 
\item In the identity $({\mathcal R_{\boldsymbol{i}}})$, the $x_\ell(\boldsymbol{i})$'s are seen as (positive) rational functions in the initial cluster variables $x_1(0),\ldots,x_n(0)$ and the theorem asserts that this  functional identity is satisfied on the positive part  (isomorphic to the positive orthant 
$(\mathbf R_{>0})^n$) 
 of the $\boldsymbol{\mathcal X}$-cluster torus 
 $ \mathbf T_{S_0} $
 associated to the initial seed $S_0$.
  \item   In terms of the cluster dilogarithm ${\sf R}$ considered in \S\ref{Par:NotationDilogarithmicFunctions}, the dilogarithmic identity 
$(\mathcal R_{\boldsymbol{i}})$ can be  written   in the simpler but equivalent following form:
\begin{center}
$\big(\boldsymbol{{\sf R}_{\boldsymbol{i}}}\big)
\hspace{3.5cm}
\sum_{\ell=1}^k d_{i_\ell} \, {\sf R} \Big( 
x_\ell({\boldsymbol{i}})\Big)= (\pi^2/6)\, N_{{\boldsymbol{i}}}
\hspace{0.1cm}.
  \hspace{5cm}{}^{}$
\end{center}
 \item 
 When $B_0$ is skew-symmetric, one has $d_{i_\ell}=1$ for every $\ell$ and there is a equivalent version of $\big(\boldsymbol{{\sf R}_{\boldsymbol{i}}}\big)$ with trivial constant second member: if $\epsilon_\ell$ stands for the tropical sign of $x_\ell(\boldsymbol{i})$ for any $\ell=1,\ldots,k$, then one has ({\it cf.}\,\cite[Theorem 2.9]{KN}):
 \begin{center}
 \label{kooko}
$\big(\boldsymbol{{\sf R}_{\boldsymbol{i}}}^{\hspace{-0.05cm}\epsilon}\big)
\hspace{3.3cm}
\sum_{\ell=1}^k \epsilon_{\ell} \, {\sf R} \Big( 
{x_\ell({\boldsymbol{i}})}^{\epsilon_{\ell}}\Big)= 0
\hspace{0.1cm}.
  \hspace{6cm}{}^{}$
\end{center}
 \end{enumerate}

The identity $(\mathcal R_{\boldsymbol{i}})$ (or equivalently, 
$(\boldsymbol{{\sf R}_{\boldsymbol{i}}})$ or $(\boldsymbol{{\sf R}_{\boldsymbol{i}}}^{\hspace{-0.05cm}\epsilon})$)  induces a dilogarithmic AR, denoted again (a bit abusively) by $\mathcal R_{\boldsymbol{i}}$ (or $({{\sf R}_{\boldsymbol{i}}})$, it does not really matter), for the cluster web 
 defined by the $x_\ell({\boldsymbol{i}})$'s for $\ell$ ranging from 1 to $k$.  By definition, the latter is the 
{\bf $\boldsymbol{\mathcal P}$-cluster (dilogarithmic) web} associated to the period $\boldsymbol{i}$ (the ${\mathcal P}$ 
  referring of course to the word `period')  and it is denoted by 
\begin{equation}
\label{Eq:PWeb-i}
\boldsymbol{\mathcal P\hspace{-0.05cm}\mathcal W}_{\boldsymbol{i}}=
\boldsymbol{\mathcal W}\left(  \, x_\ell({\boldsymbol{i}})\hspace{0.1cm} \big\lvert \hspace{0.1cm}  \ell=1,\ldots,k\hspace{0.1cm}
\right)\, . 
\end{equation}

Such a web carries a complete dilogarithmic AR, hence many logarithmic ARs as well.  It is then natural to ask   the following question about it: 
\begin{question}
\label{Quest:PWeb-AMP?} Is a $\boldsymbol{\mathcal P}$-cluster  web necessarily AMP with only polylogarithmic ARs? 
\end{question}

In view of the examples of different kinds considered in the paragraphs below, if the answer to this question is not always positive, it is positive in sufficiently many cases to make it appear particularly relevant  hence deserving further studies.\sk 

As by now, almost nothing is known about $\boldsymbol{\mathcal P}$-cluster  webs in general (except that each carries a dilogarithmic AFE of course). For instance, even what is the degree of such a web is not clear. Indeed, 
 given a cluster period $\boldsymbol{i}$, it is not clear yet whether  all the associated cluster variables
 $x_\ell({\boldsymbol{i}})$ define pairwise distinct foliations or not.    Of course, it could not be  the case if ${\boldsymbol{i}}$ is not irreducible\footnote{For instance, the cluster web associated to $\boldsymbol{i}\lvert {\boldsymbol{i}}$ is easily seen to be the same as the one defined by ${\boldsymbol{i}}$.}, but when 
 ${\boldsymbol{i}}$ is assumed to be irreducible, we believe  that  this is indeed the case and state it as the following 
\begin{conjecture}
\label{Conj:Diff-Indep-in-periods}
Let $\boldsymbol{i}$ be a period of a cluster algebra, of length $k$. If it is irreducible, then the associated cluster variables $x_{\boldsymbol{i}}(\ell)$ for $\ell=1,\ldots,k$ 
  are pairwise differentially distinct, {\it i.e.} two of them define distinct foliations. Equivalently:  the cluster web $\boldsymbol{\mathcal P\mathcal W}_{\boldsymbol{i}}$ has degree~$k$.
\end{conjecture}

\paragraph{Examples of periods I: $\boldsymbol{Y}$-systems} 
 \label{Par:ExamplesPeriods-I}

%
%

Any of the $Y$-systems considered in \S\ref{Par:Classical-Y-systems} above can be formalized in terms of cluster algebras and from this point of view,  its (half-)periodicity property corresponds to the fact that a certain finite sequence of mutations is a cluster period. 
\sk

For the  $Y(\Delta,\Delta')$'s, this is evident from \S\ref{Par:Interpretation-cluster-Y-systems-Dynkin-type}.
 In the case when both 
Dynkin diagrams $\Delta,\Delta'$ are simply-laced, let $\boldsymbol{j}_\bullet $ (resp.\,$\boldsymbol{j}_\circ$) be a tuple of the numbers of the odd (resp.\,even) vertices (in any order, this does not matter for what follows). Then one sets 
$\boldsymbol{j}_{\Delta,\Delta'}= 
\boldsymbol{j}_\bullet
 \lvert 
\boldsymbol{j}_\circ $ where the short vertical line stands for the concatenation between the tuples on both sides.  
Setting $h_{\Delta,\Delta'}=h+h'$, it comes that 
 the full-periodicity for  the $Y$-system $Y(\Delta,\Delta')$ is equivalent to the fact that 
the concatenation 
\begin{equation}
\label{Eq:period-slice}
\boldsymbol{i}_{\Delta,\Delta'}=
\underbrace{\boldsymbol{j}_{\Delta,\Delta'} \lvert \, \cdots \, \lvert \, \boldsymbol{j}_{\Delta,\Delta'}}_{(h_{\Delta,\Delta'})-times} 
\end{equation}
of $h_{\Delta,\Delta'}$
 copies of  $\boldsymbol{j}_{\Delta,\Delta'}$ is a period 
  for the initial $\boldsymbol{\mathcal X}$-seed with exchange matrix $B_{\Delta\square \Delta'}$. \mk 
  
\begin{rem}
\label{Rem:Half-Periods}
{1.}  Actually, in many cases (and possibly in all, see below), 
  $\boldsymbol{i}_{\Delta,\Delta'}$ is not irreducible as a period.  And 
 there is a {\it `half-periodicity'} property satisfied by the $Y$-system $Y(\Delta,\Delta')$ which corresponds to the fact that $\boldsymbol{i}_{\Delta,\Delta'}$ admits a non-trivial factorization as a period. \sk 
 
{2.}  When both $\Delta$ and $\Delta'$ are simply laced, there are permutations  
  $\omega=\omega_\Delta$ and $\omega'=\omega'_{\Delta'}$ of the sets of vertices of $\Delta$ and $\Delta'$ respectively, such that, according to {\rm \cite[Theorem 4.27]{IIKNS}}, one has: 
  \begin{itemize}
  \item[$-$] 
   if both $h$ and $h'$ are even or if both are odd, then 
\begin{equation}
\label{Eq:i-Delta-Delta'-half}
  \boldsymbol{i}^{half}_{\Delta,\Delta'}=\big(\boldsymbol{j}_{\Delta,\Delta'}\big)^{(h+h')/2}
  \end{equation}
   is a $(\omega \otimes \omega')$-period for the initial $\boldsymbol{\mathcal X}$-seed with quiver ${\Delta\square \Delta'}$, such that $\boldsymbol{i}_{\Delta,\Delta'}=\big(\boldsymbol{i}_{\Delta,\Delta'}^{half}\big)^2$;
    \item[$-$] if  $h$ is even but  $h'$ is odd, then $\Delta'=A_r$ for an even integer $r$.  Then setting 
     $\boldsymbol{j}_{\bullet,\circ}=\boldsymbol{j}_{\Delta,A_r}=
 \boldsymbol{j}_{\bullet}\lvert \boldsymbol{j}_{\circ}$
and  $\boldsymbol{j}_{\circ,\bullet}=\big( 
     \boldsymbol{j}_{\Delta,A_r}
     \big)^{-1}=
 \boldsymbol{j}_{\circ}\lvert \boldsymbol{j}_{\bullet}$, it comes that 
\begin{equation}
\label{Eq:PeriodAn-n-even}
\boldsymbol{i}_{\bullet,\circ}^{half}=\boldsymbol{j}_\circ\big \lvert \,  \big(
\boldsymbol{j}_{\bullet,\circ}
\big)^{\frac{h+h'-1}{2}}
\qquad \mbox{or equivalently}\qquad 
\boldsymbol{i}_{\circ,\bullet}^{half}=\boldsymbol{j}_\bullet \big\lvert \,  \big(
\boldsymbol{j}_{\circ,\bullet}
\big)^{\frac{h+h'-1}{2}}
\end{equation}
 are periods 
  for the initial $\boldsymbol{\mathcal X}$-seed with quiver ${\Delta\square A_r}$, which are such that 
$ \boldsymbol{i}_{\Delta,A_r}= \big(\boldsymbol{i}_{\circ,\bullet}^{half} 
\big)^{-1}
\, 
\lvert  
\, \boldsymbol{i}_{\bullet,\circ}^{half}$. \sk
  \end{itemize}
 
 \vspace{-0.4cm}
{3.} Much less has been established in the case when at least one of the two involved Dynkin diagrams is not simply-laced. When one of them is not simply laced, say $\Delta$,  the half-periodicity is known in full-generality only when $\Delta'=A_1$ (and this is due to Fomin and Zelevinsky, see Theorem 4.4 in 
{\rm \cite{IIKNS}} and the references therein). \mk 

{4.} The period $  \boldsymbol{i}_{\Delta,\Delta'}$ can be seen to be of length $(h+h')nn'$  where 
$n$ and $n'$ stand for the rank of $\Delta$ and $\Delta'$ respectively.
When the half-periodicity property holds true, $  \boldsymbol{i}_{\Delta,\Delta'}$ can be factorized as the concatenation of two periods, each  of length $(h+h')nn'/2$ (which  makes sense since this number is always an integer).  In all the explicit examples we have considered (see below regarding a few of them),  the half-periodicity was verified. So we conjecture that it holds true in full generality. 
 \mk 

{5.} The half-periodicity of the $Y$-system $Y_\ell(\Delta)$ has also been established in many cases {\rm (}{\it cf.}\,{\rm \cite{IIKNS}, \cite{Inoue1},  \cite{Inoue2})}. 
 \end{rem}

 At this point, it is interesting to discuss the links between the cluster  webs 
 $\boldsymbol{B\hspace{-0.04cm}B\hspace{-0.04cm}\mathcal W}_{\Delta\square \Delta'}$ considered above and the 
 cluster web  associated to $\boldsymbol{i}_{\Delta,\Delta'}$: 
 \begin{itemize}
 \item  writing $\boldsymbol{j}_{\Delta,\Delta'}=(j_{nn'},\ldots,j_{1})$, the 
  $\boldsymbol{\mathcal P}$-cluster web $\boldsymbol{\mathcal P\mathcal W}_{\boldsymbol{i}_{\Delta,\Delta'}} $   is by definition the one associated to $\boldsymbol{i}_{\Delta,\Delta'}$, the latter being  considered as a $nn'm_{\delta,\Delta'}$-tuple of elements of  $\{1,\ldots,nn'\}$.
 \footnote{Actually, since a  certain `power' of $\boldsymbol{i}^{half}_{\Delta,\Delta'}$ is equal to $\boldsymbol{i}_{\Delta,\Delta'}$, it comes that the two $\boldsymbol{\mathcal P}$-cluster webs 
  associated to these periods coincide. In practice, it is more convenient to work with $\boldsymbol{i}^{half}_{\Delta,\Delta'}$ which is shorter than $\boldsymbol{i}_{\Delta,\Delta'}$.}
   It follows from Theorem  \ref{Thm:Nakanishis-identity} that this web coincides with the $\boldsymbol{\mathcal Y}$-cluster web $\boldsymbol{\mathcal Y\mathcal W}_{{\Delta,\Delta'}} $ defined before in  \S\ref{parag:Webs-associated-to-Y-systems}. 
 \item the web $\boldsymbol{B\hspace{-0.04cm}B\hspace{-0.04cm}\mathcal W}_{\Delta\square \Delta'}$ is rather 
 associated to its writing  as the concatenation \eqref{Eq:period-slice}
 of some copies of $\boldsymbol{j}_{\Delta,\Delta'}$
 than to   $\boldsymbol{j}_{\Delta,\Delta'}$ per se.   This is related to the notion of {\it `slice'} of a period considered in \cite{Nakanishi}  about which there is no real need to  elaborate with some details here. 
 We just say that when the period considered is {\it `regular'} ({\it cf.}\,Definition 5.1 in \cite{Nakanishi}), then the $Y$-systems (and as it happens, the associated $\boldsymbol{\mathcal Y}$-cluster webs) associated to its slices, all are equivalent.  
 \item  When both $\Delta$ and $\Delta'$ are simply-laced, the period $\boldsymbol{i}_{\Delta,\Delta'}$ is 
 regular hence the  three webs $\boldsymbol{B\hspace{-0.04cm}B\hspace{-0.04cm}\mathcal W}_{\Delta\square \Delta'}$, 
 $\boldsymbol{\mathcal P\mathcal W}_{\boldsymbol{i}_{\Delta,\Delta'}} $
 and $\boldsymbol{\mathcal Y\mathcal W}_{{\Delta,\Delta'}} $  
  coincide.  This is not the case when (at least) one of the Dynkin diagram $\Delta$ of $\Delta'$ is not simply-laced (see \S\ref{Subpar:YWB2*B2} below for the case $\Delta=\Delta'=B_2$ for instance). 
 \end{itemize}

 For further details, see \cite[\S3] {KellerAnnals} and \cite{Nakanishi}.
  The case when one (or two)  of the Dynkin diagrams involved  is (or are) not simply-laced is treated in \cite[\S9] {KellerAnnals}. 
The cluster formulations of the $Y$-systems   $Y_\ell(\Delta)$'s when $\Delta$ is not simply-laced\footnote{We remind that $Y_\ell(\Delta)$ is equivalent to $Y(A_{\ell-1},\Delta)$ when $\Delta$ is simply-laced.} is 
 given in the second section of \cite{Inoue1} in type $B$ and in \cite{Inoue2} for the other cases (type $C$, $F_4$, $G_2$).  For some explicit examples, one can also look at the third section of \cite{Nakanishi}.

 \subparagraph{The $\boldsymbol{\mathcal Y}$-cluster web of type $\boldsymbol{A_2\boxtimes B_2}$.} 
 \label{Subpar:YWA2*B2}
We work with the same initial exchange matrix as in \S\ref{SubPar:A2*B2}. 
From \cite{KellerAnnals} (in particular formula (3.6.1) therein),  the half-periodicity property of the $Y$-system is equivalent to the fact that 
\begin{equation}
\label{Eq:period-i-B2*B2}
\boldsymbol{i}_{A_2\boxtimes B_2}=\big(3, 4, 1, 2,  3, 4, 1, 2,3, 4, 1, 2, 3, 4\big)
\end{equation}
(obtained by taking the 14 first entries of the concatenation of $h(A_2)+h(B_2)=7$ copies of 
$(3, 4, 1, 2)$) is a period for the initial seed with exchange matric $B_{A_2\boxtimes A_2}$.  We deduce that $\boldsymbol{\mathcal Y 
 \hspace{-0.05cm} 
 \mathcal W}_{\hspace{-0.05cm}  A_2\square B_2}$ is a 14-web in four variables and 
 by direct computations, one gets  that 
 $$
 \rho^\bullet\Big( \boldsymbol{\mathcal Y 
 \hspace{-0.05cm} 
 \mathcal W}_{\hspace{-0.05cm}  A_2\boxtimes B_2}
 \Big)=(10,4,1)
 \qquad 
 \mbox{ and }
  \qquad
 {\rm polrk}^\bullet
 \Big( \boldsymbol{\mathcal Y 
 \hspace{-0.05cm} 
 \mathcal W}_{\hspace{-0.05cm}  A_2\boxtimes B_2}
 \Big)=(14,1)\,.
 $$
It follows that this web is AMP with only logarithmic ARs plus the dilogarithmic one ${\sf R}_{A_2\boxtimes B_2}$. 
 \subparagraph{The $\boldsymbol{\mathcal Y}$-cluster web of type $\boldsymbol{B_2\boxtimes B_2}$.} 
 \label{Subpar:YWB2*B2}
 We have seen above that the web 
 $ \boldsymbol{B\hspace{-0.04cm}B\hspace{-0.04cm}\mathcal W}_{\hspace{-0.01cm}B_2\boxtimes B_2}$ is not AMP. 
 We now explain that this web actually differs from 
 $ \boldsymbol{\mathcal Y\hspace{-0.04cm}\mathcal W}_{\hspace{-0.01cm}B_2\boxtimes B_2}$. \sk 
 
 The map $F_{B_2\boxtimes B_2}$ considered in \S\ref{SubPar:B2*B2}
corresponds to the sequence of mutations $\boldsymbol{j}_{B_2,B_2}=(2,1,4,3)$ hence the half-periodicity is equivalent to the fact that 
\begin{equation}
\label{Eq:period-i-B2*B2}
\boldsymbol{i}_{B_2\boxtimes B_2}=\big(2,1,4,3,2,1,4,3,2,1,4,3,2,1,4,3\big)
\end{equation}
 is a period for the initial seed with exchange matrix $B_{B_2\boxtimes B_2}$ (see \eqref{Eq:kotok}). The 
$\boldsymbol{\mathcal P}$-cluster 16-web 
 associated to the period of length 16 above, and denoted here by 
 $\boldsymbol{\mathcal Y 
 \hspace{-0.05cm} 
 \mathcal W}_{\hspace{-0.05cm}  B_2\boxtimes B_2}$, 
   is equivalent to the $\boldsymbol{\mathcal Y}$-cluster web $\boldsymbol{\mathcal Y 
 \hspace{-0.05cm} 
 \mathcal W}_{\hspace{-0.05cm}  B_2\square B_2}$ 
 ({\it cf.}\,\S\ref{parag:Webs-associated-to-Y-systems}). By direct computations, one gets 
 $$
 \rho^\bullet\Big( \boldsymbol{\mathcal Y 
 \hspace{-0.05cm} 
 \mathcal W}_{\hspace{-0.05cm}  B_2\boxtimes B_2}
 \Big)=(12,6,1)
 \qquad 
 \mbox{ and }
  \qquad
 {\rm polrk}^\bullet
 \Big( \boldsymbol{\mathcal Y 
 \hspace{-0.05cm} 
 \mathcal W}_{\hspace{-0.05cm}  B_2\boxtimes B_2}
 \Big)=(16,1)\,
 $$  
 thus 
 $
 {\rm polrk}
 \big( \boldsymbol{\mathcal Y 
 \hspace{-0.05cm} 
 \mathcal W}_{\hspace{-0.05cm}  B_2\boxtimes B_2}
 \big)=17<
 19= \rho^\bullet\big( \boldsymbol{\mathcal Y 
 \hspace{-0.05cm} 
 \mathcal W}_{\hspace{-0.05cm}  B_2\boxtimes B_2}
 \big)
 $. 
Denoting by $x_1,\ldots,x_{16}$ the cluster first integrals associated to the period \eqref{Eq:period-i-B2*B2} and setting ${\sf A}(x)=\arctan(\sqrt{x}\,)$ for any $x>0$, one verifies that the  two following functional  identities hold true identically: 
 \begin{align*}
 2 =&\hspace{0.1cm}   \frac{1}{1+x_3}
 +\frac{1}{1+x_7}
  +\frac{1}{1+x_{11}}
  + \frac{1}{1+x_{15}}\,, \vspace{0.5cm} \\
\pi    =&\hspace{0.1cm} 
{\sf A}(x_2) + {\sf A}(x_6) + {\sf A}(x_{10}) + {\sf A}(x_{14})\, .
 \end{align*}
 The two ARs associated to these two identities provide a supplementary space to the space of logarithmic and dilogarithmic ARs of $ \boldsymbol{\mathcal Y 
 \hspace{-0.05cm} 
 \mathcal W}_{\hspace{-0.05cm}  B_2\boxtimes B_2}$ which ensures that $\rho( \boldsymbol{\mathcal Y 
 \hspace{-0.05cm} 
 \mathcal W}_{\hspace{-0.05cm}  B_2\boxtimes B_2})=
 {\rm rk}^\bullet
 ( \boldsymbol{\mathcal Y 
 \hspace{-0.05cm} 
 \mathcal W}_{\hspace{-0.05cm}  B_2\boxtimes B_2}
)$. Thus this web is AMP, unlike  $\boldsymbol{B\hspace{-0.04cm}B\hspace{-0.04cm}\mathcal W}_{\hspace{-0.01cm}B_2\boxtimes B_2}$. Note however that  not all the  ARs
$\boldsymbol{\mathcal Y 
 \hspace{-0.05cm} 
 \mathcal W}_{\hspace{-0.05cm}  B_2\boxtimes B_2}$ are polylogarithmic (but they are not far from being almost of this kind, see the notion 
 of `{generalized iterated integral}' considered at the end of 
  \S\ref{SubPar:CharacterizationWebsMaximalRank}).


 \subparagraph{The $\boldsymbol{\mathcal Y}$-cluster web of type $\boldsymbol{A_2\boxtimes G_2}$.} 
 \label{Subpar:YWA2*G2}
 The exchange skew-diagonalizable matrix associated to the valued quiver ${A}_2\boxtimes G_2$ is  the one obtained from that of ${A}_2\boxtimes B_2$ by replacing the $2$'s  by $3$' in it (and keeping the corresponding signs). 
 From \cite{KellerAnnals} (see formula (3,6,1) therein in particular), we know that  that the concatenation $\mu_{A_2,G_2}$ of $h(A_2)+h(G_2)=3+6=9$ copies of $(2,4,1,3)$ is a cluster period. Actually, $\mu_{A_2,G_2}$ is not irreducible as its 18 first terms already form a cluster period. We deduce that 
 $\boldsymbol{\mathcal Y 
 \hspace{-0.05cm} 
 \mathcal W}_{\hspace{-0.05cm}  A_2\boxtimes G_2}$
  is a 18-web in four variables and by direct computations, we have obtained that 
 $$
 \rho^\bullet\Big(\boldsymbol{\mathcal Y 
 \hspace{-0.05cm} 
 \mathcal W}_{\hspace{-0.05cm}  A_2\boxtimes G_2}\Big)=(14, 8, 1)
 \, , \qquad 
{\rm polrk}^\bullet\Big(\boldsymbol{\mathcal Y 
 \hspace{-0.05cm} 
 \mathcal W}_{\hspace{-0.05cm}  A_2\boxtimes G_2}\Big)=(18, 1)
 \qquad \mbox{ and }\qquad 
 {\rm rk}^\bullet\Big(\boldsymbol{\mathcal Y 
 \hspace{-0.05cm} 
 \mathcal W}_{\hspace{-0.05cm}  A_2\boxtimes G_2}\Big)=19\, .
 $$

 In particular, one has 
${\rm polrk}
\big(\boldsymbol{\mathcal Y 
 \hspace{-0.05cm} 
 \mathcal W}_{\hspace{-0.05cm}  A_2\boxtimes G_2}\big)=
{\rm rk}
\big(\boldsymbol{\mathcal Y 
 \hspace{-0.05cm} 
 \mathcal W}_{\hspace{-0.05cm}  A_2\boxtimes G_2}\big)=19<27=\rho^\bullet\big(\boldsymbol{\mathcal Y 
 \hspace{-0.05cm} 
 \mathcal W}_{\hspace{-0.05cm}  A_2\boxtimes G_2}\big)$ hence all the ARs of this web 
 are polylogarithmic (of weight $\leq 2$)  but it is not AMP.

 \subparagraph{The $\boldsymbol{\mathcal Y}$-cluster web of type $\boldsymbol{B_2\boxtimes G_2}$.} 
 \label{Subpar:YWB2*G2}

The exchange matrix in the case $B_2\boxtimes G_2$
is 
$$
B_{B_2\boxtimes G_2}=\begin{bmatrix}
0& 3& 2& -6 \\
-1& 0& 0& 2 \\
-1& 0& 0& 3\\ 
1& -1& -1& 0
\end{bmatrix}\, . 
$$

The associated period has length 20, and is obtained by 
concatenating  $5$ copies of $(2, 4, 1, 3)$.  The associated cluster web
$\boldsymbol{\mathcal Y 
 \hspace{-0.05cm} 
 \mathcal W}_{\hspace{-0.05cm}  B_2\boxtimes G_2}$ is a 20-web in four variables, such that 
$$
\rho^\bullet\big(\boldsymbol{\mathcal Y 
 \hspace{-0.05cm} 
 \mathcal W}_{\hspace{-0.05cm}  B_2\boxtimes G_2}\big)=(16, 10, 1)\, , 
\qquad 
{\rm polrk}^\bullet
\big(\boldsymbol{\mathcal Y 
 \hspace{-0.05cm} 
 \mathcal W}_{\hspace{-0.05cm}  B_2\boxtimes G_2}\big)=(20, 1) 
 \qquad \mbox{and}\qquad 
{\rm rk}
\big(\boldsymbol{\mathcal Y 
 \hspace{-0.05cm} 
 \mathcal W}_{\hspace{-0.05cm}  B_2\boxtimes G_2}\big)=22\, .$$

In particular, one has 
${\rm rk}
\big(\boldsymbol{\mathcal Y 
 \hspace{-0.05cm} 
 \mathcal W}_{\hspace{-0.05cm}  B_2\boxtimes G_2}\big)=22<27=\rho^\bullet\big(\boldsymbol{\mathcal Y 
 \hspace{-0.05cm} 
 \mathcal W}_{\hspace{-0.05cm}  B_2\boxtimes G_2}\big)$ hence this web is not AMP.

 \subparagraph{The $\boldsymbol{\mathcal Y}$-cluster web associated to $\boldsymbol{Y_2(B_2)}$.} 
We refer to \cite[\S2.3]{Inoue1} for a description of the quiver 
 $Q_\ell(B_n)$ associated to the $Y$-system of level $\ell\geq 2$ and of type $B_n$.  In the case $\ell=n=2$ that we consider here, this quiver and the associated exchange matrix are
%
%
\begin{equation*}
\label{Eq:kotok}
\begin{tabular}{lcr}
 \begin{tabular}{c}  \\
  \scalebox{0.36}{ \includegraphics{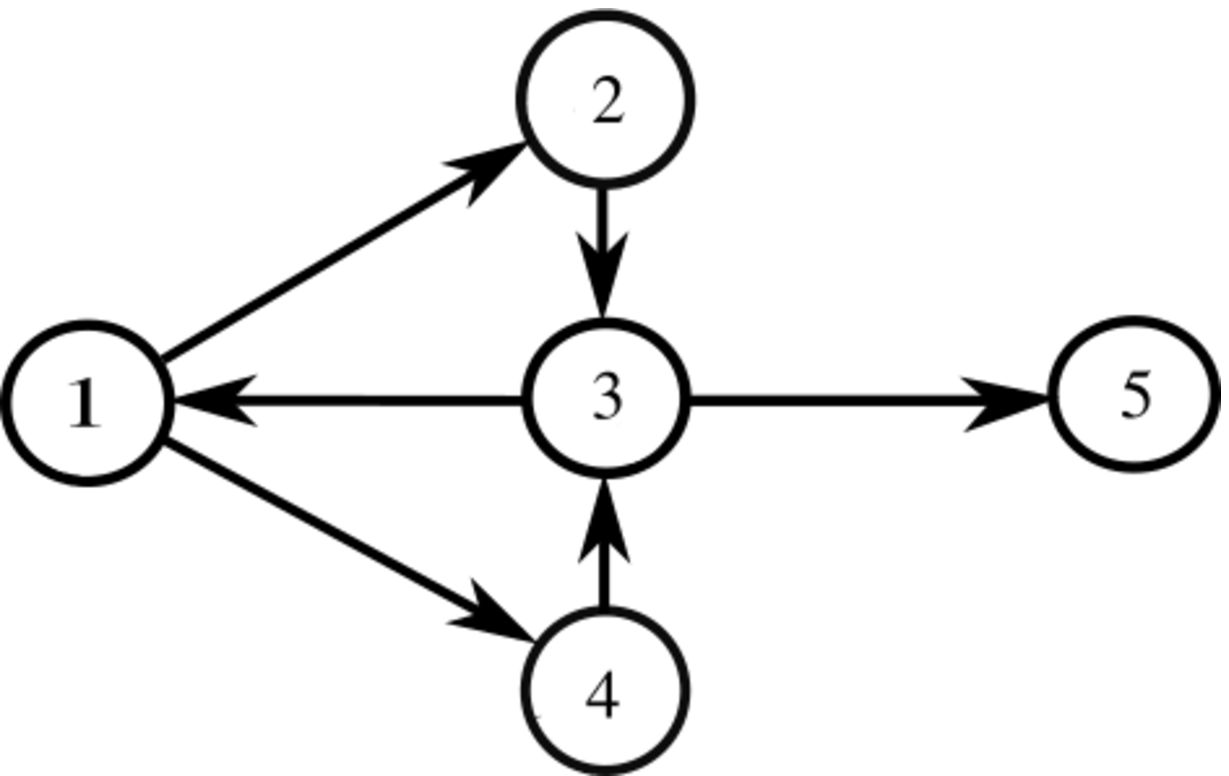}}
  \end{tabular}
     &   &  
$\begin{tabular}{c} 
 and  $ \qquad$ $B_{Q_2(B_2)}=\begin{bmatrix}
0 & 1 & -1 & 1 & 0 \\
-1 & 0& 1& 0& 0\\ 
1& -1& 0& -1& 1 \\
-1& 0& 1& 0& 0 \\
0& 0& -1& 0& 0
\end{bmatrix}$\, .
\end{tabular}$\sk 
\end{tabular}
\end{equation*}
 
 Then the concatenation of 5 copies of $(3, 1, 2, 4, 3, 5, 2, 4)$ gives a cluster period  $\boldsymbol{i}_{Q_2(B_2)}=(i_{40},\ldots, i_1)$
  for the initial $\boldsymbol{\mathcal X}$-seed $S_0=((x_i)_{i=1}^5 \, B_{Q_2(B_2)})$. 
   Actually, half of $\boldsymbol{i}_{Q_2(B_2)}$, namely  $(i_{20},\ldots,i_1)$, is already a period for $S_0$, but relatively to the permutation $i\mapsto 6-i$ of $\{1,\ldots,5\}$.  The associated $\boldsymbol{\mathcal P}$-cluster web is
   $\boldsymbol{\mathcal Y\hspace{-0.05cm}\mathcal W}_{Q_2(B_2)}$: it is  a 20-web in five variables which, after some direct computations, can be proved to be such that 
 $$
 \rho^\bullet\Big( \boldsymbol{\mathcal Y\hspace{-0.05cm}\mathcal W}_{\hspace{-0.05cm} Q_2(B_2)}
 \Big)=(15,5,1)
 \qquad 
 \mbox{ and }
  \qquad
 {\rm polrk}^\bullet
 \Big( \boldsymbol{\mathcal Y\hspace{-0.05cm}\mathcal W}_{\hspace{-0.05cm} Q_2(B_2)}
 \Big)=(20,1)\,. 
 $$
 Thus the $\boldsymbol{\mathcal Y}$-cluster web of type $Q_2(B_2)$ is AMP with only polylogarithmic  ARs (all of weight 1, except the dilogarithmic one $\mathcal R_2(B_2)$).

 \paragraph{Examples of periods II: Zamolodchikov periodic quivers.} 
The general  machinery of $Y$-systems can be generalized to a more general 
framework considered in \cite{GP}. \sk 
 
 Let $Q$ be a  quiver with bipartite underlying graph: its vertices are   either black  $\bullet$ or white $\circ$.  Let $\boldsymbol{j}_\bullet$ (resp.\, $\boldsymbol{j}_\circ$) be a tuple consisting of all the black (resp.\,white) vertices of $Q$ and denote by  $\mu_\bullet$ (resp.\,$\mu_\circ$) the associated composition of mutations. 
 Since mutations associated to vertices of the same color pairwise commute, 
 both  
 $\mu_\bullet$  and $\mu_\circ$ are well-defined.
The quiver $Q$ is said to be 
\begin{itemize}
\item[$-$] {\it `recurrent'} if $\mu_Q=\mu_\circ \circ \mu_\bullet$ transforms $Q$ into   its {\it `dual'} $-Q$, which is the quiver with the same underlying bipartite graph but with all its arrows reversed in it;  
\item[$-$] {\it `Zamolodchikov periodic'} if it is 
recurrent and if the concatenation, denoted by  $\boldsymbol{i}_Q$, of a certain number $N_Q\in \mathbf N^*$ of copies of 
$\boldsymbol{j}_Q=\boldsymbol{j}_\circ\lvert \boldsymbol{j}_\bullet$ is a 
cluster period for the initial 
$\boldsymbol{\mathcal X}$-seed with $Q$ as its associated quiver. 
 \end{itemize}

In  \cite{GP}, the authors give a complete classification of the 
{\it `Zamolodchikov periodic quivers'}, a result which can be interpreted as the classification of certain cluster periods. They prove 
 that a bipartite recurrent quiver $Q$ is Zamolodchikov periodic if and only if a certain bipartite bigraph $G(Q)$,  associated to $\pm Q$  in a one-to-one way,  is of a certain kind, namely `{\it admissible and of ADE type}' (see \S1.3.2 in \cite{GP} for more details and more specifically Theorem 1.10 therein).  Since bigraphs of this kind have been classified before ({\it cf.}\,\cite[Thm\,3.1]{Stembridge}), one obtains a complete list of Zamolodchikov periodic quivers up to duality: 
there are five infinite families as well as eight exceptional cases. The tensor products $\Delta\boxtimes\Delta'$ (where both Dynkin diagrams are  
 simply-laced) constitute of course one of these families, and the  four other families 
as well as the exceptional cases
can be described quite explicitly  (see respectively \S1.3.2 and Appendix A in \cite{GP}).  

Given a Zamolodchikov periodic quiver $Q$, one denotes by $\boldsymbol{\mathcal Z\hspace{-0.05cm}\mathcal W}_Q$ the $\boldsymbol{\mathcal X}$-cluster web associated to the period $\boldsymbol{i}_Q=\boldsymbol{j}_Q^{N_Q}$. Theorem \ref{Thm:Nakanishis-identity}  ensures that this web carries a complete dilogarithmic AR.\sk 

 Thanks to the classification of  Zamolodchikov periodic quivers, we have now at disposal five families and eight exceptional cases of cluster webs carrying polylogarithmic ARs (of weight 1 and 2) 
 which all are good candidates for being AMP webs with only logarithmic and dilogarithmic ARs.
  Two explicit examples are discussed below, which support our belief that all the cluster webs
associated to Zamolodchikov periodic quivers are indeed AMP.

 \subparagraph{The cluster web associated to the twist $\boldsymbol{A_3\times A_3}$.} 
Another infinite family of Zamo\-lodchikov periodic quivers is the one formed by the {\it `twists'} $\Delta\times \Delta$ for any $ADE$ Dynkin diagram $\Delta$ ({\it cf.}\,Example 1.4 in \cite{Stembridge} or p.\,453 in \cite{GP}).  For instance, the corresponding Zamolodchikov periodic quiver and the associated exchange matrix in the case when $\Delta=A_3$ are the following: 
\begin{equation*}
\label{Eq:kotok}
\begin{tabular}{lcr}
 \begin{tabular}{c}  \\
  \scalebox{0.36}{ \includegraphics{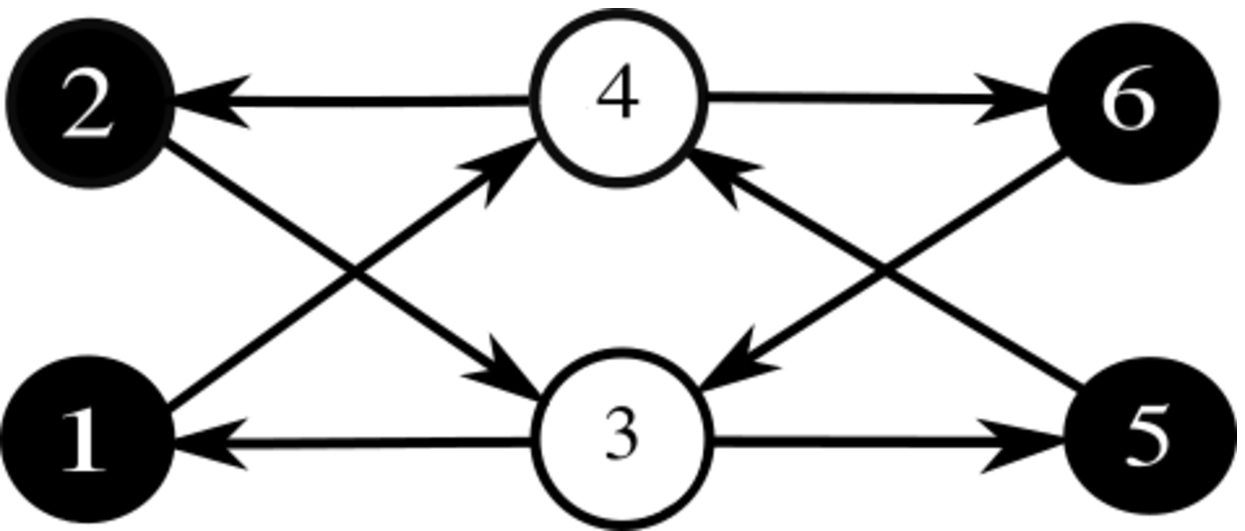}}
  \end{tabular}
     &   &  
$\begin{tabular}{c} 
 and  $ \qquad$ $B_{A_3\times A_3}=\begin{bmatrix}
0 & 0 & -1 & 1 & 0 & 0  \\
0 & 0& 1& -1& 0&0\\ 
1& -1& 0& 0& 1 &-1\\
-1& 1& 0& 0& -1& 1\\
0& 0& -1& 1& 0& 0 \\
0& 0& 1& -1& 0&  0
\end{bmatrix}$\, .
\end{tabular}$\sk 
\end{tabular}
\end{equation*}

In this case one has $\boldsymbol{j}_{A_3\times A_3}=(4,3,6,5,2,1)$ and one verifies that $\boldsymbol{i}_{A_3\times A_3}=(\boldsymbol{j}_{A_3\times A_3})^4$ is a cluster period for the initial $\boldsymbol{\mathcal X}$-seed 
whose quiver is the one just above. 
The associated cluster web $
\boldsymbol{\mathcal Z\hspace{-0.05cm}\mathcal W}_{A_3\times A_3}$ is a 24-web in 6 variables  such that $$
\rho^\bullet\Big( \boldsymbol{\mathcal Z\hspace{-0.05cm}\mathcal W}_{A_3\times A_3}\Big)=
\big(18, 6, 1\big) 
\qquad {\rm polrk}^\bullet
\Big( \boldsymbol{\mathcal Z\hspace{-0.05cm}\mathcal W}_{A_3\times A_3}\Big)=(24,1)\,. 
$$
So $\boldsymbol{\mathcal Z\hspace{-0.05cm}\mathcal W}_{A_3\times A_3}$ is AMP with only polylogarithmic ARs (all of weight 1, except the dilogarithmic one).

 \subparagraph{The cluster web associated to the sporadic  Zamolodchikov periodic quiver $\boldsymbol{D_6 *  D_6}$.} 
The simplest of the eight sporadic/exceptional  Zamo\-lodchikov periodic quivers, denoted by $D_6*D_6$,  is the following one: 
$$
  \scalebox{0.36}{ \includegraphics{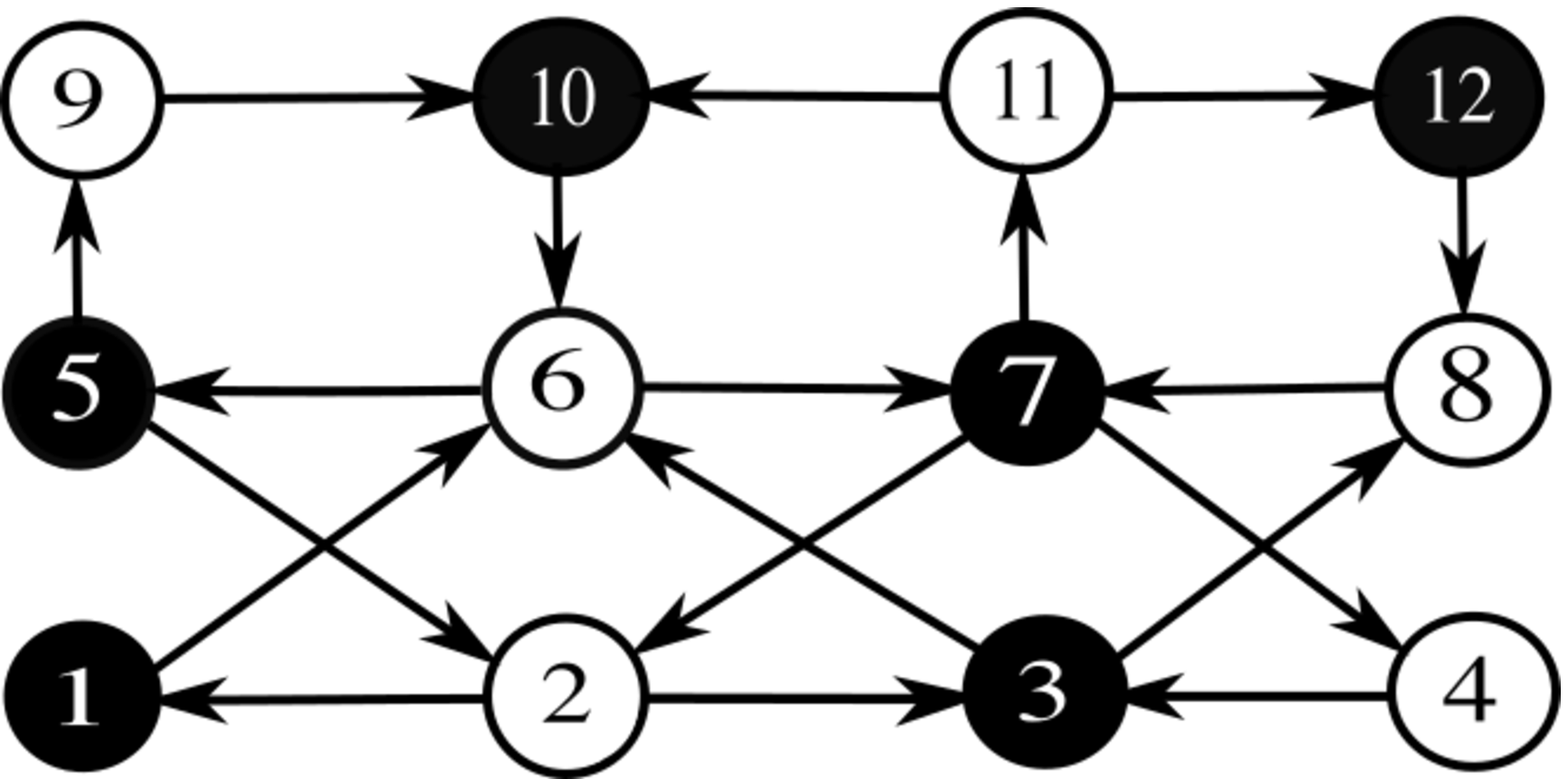}}
$$

In this case one has $\boldsymbol{j}_{D_6* D_6}=(11,9,8,6,4,2,12,10,7,5,3,1)$ and one verifies that $\boldsymbol{i}_{D_6* D_6}=(\boldsymbol{j}_{D_6*  D_6})^{15}$ is a cluster period for the initial $\boldsymbol{\mathcal X}$-seed whose associated quiver is the one  just above.\footnote{Note that 
only the fact that $(\boldsymbol{i}_{D_6* D_6})^2= (\boldsymbol{j}_{D_6*  D_6})^{30}$ is a cluster period is established in \cite{GP}. That it is already the case for the concatenation of only 15 copies of $\boldsymbol{j}_{D_6*  D_6}$ could be seen as a half-periodicity property. It would be interesting to know whether  something similar (that is, a certain half-periodicity property) holds true for any 
 Zamo\-lodchikov periodic quiver or not. As far as we are aware of,  this is not known.}
The associated cluster web $
\boldsymbol{\mathcal Z\hspace{-0.05cm}\mathcal W}_{\hspace{-0.05cm} D_6*D_6}$ is a 180-web in twelve variables.  It would be interesting to know whether  or not this web is AMP with only polylogarithmic ARs. 
We have tried to verify this by direct computations but these were too heavy to give us a definitive answer.  Preliminary computations suggest  that this web is not AMP.

\begin{question}
Determine the Zamolodchikov periodic quivers $Q$ such that $
\boldsymbol{\mathcal Z\hspace{-0.05cm}\mathcal W}_{\hspace{-0.05cm} Q}$ is AMP.
\end{question}


 \paragraph{Some topological cluster periods.} 
 \label{Par:KimYamazaki}
Given a cluster algebra $\mathcal A$ (of rank $n\geq 1$, say), there is a one-to-one correspondance between between cluster periods and non-trivial (finite and closed) loops in the  associated exchange graph $\Gamma_{\hspace{-0.05cm}\mathcal A}$.  \mk

In \cite{KY}, H. Kim and M. Yamazaki discuss a certain property of   $\Gamma_{\hspace{-0.05cm}\mathcal A}$ that  they call `Property $\star$'. In order to state it, 
we introduce the following subgroup, denoted by $\pi'_{\hspace{-0.05cm}\mathcal A}$ (or just by $\pi'$ when the cluster algebra we are dealing with is fixed), 
of the fundamental group of the exchange graph of $\mathcal A$: 
it is 
the subgroup 
 generated by elements of the form $PLP^{-1}$, where
$L$ is a closed loop (with a base point) obtained by mutating
2 of the $n$ cluster variables while keeping the remaining
other $n-2$ variables fixed, and $P$ is an arbitrary path originating at
the base point, and $P^{-1}$ is its inverse.\sk 

 With this definition at hand, the aforementioned property can be stated concisely as follows: \mk \\
{\bf Property $\star$\,:} {\it 
$\pi'_{\hspace{-0.05cm}\mathcal A}$  coincides with 
the full fundamental group 
$\pi_1(\Gamma_{\hspace{-0.05cm}\mathcal A})$ 
of the exchange graph 
of ${\mathcal A}$.}
\mk 


 The main motivation of Kim Yamazaki to consider this property is the accessability of the dilogarithmic identity $(\mathcal R_{\mathcal P})$ associated to any period/loop $\mathcal P$ in $\Gamma_{\hspace{-0.05cm}\mathcal A}$.  
  Indeed,  considering the equivalence `cluster period$\leftrightarrow$loop in $\Gamma_{\hspace{-0.05cm}\mathcal A}$',  then it is intuitively natural to think that, given a period $\mathcal P$, if 
the associated dilogarithmic identity 
$(\mathcal R_{\mathcal P})$ 
 is accessible from 
 Abel's 5-terms relation, then  seen as a loop in exchange graph, $\mathcal P$ belongs to the subgroup 
 $\pi'_{\hspace{-0.05cm}\mathcal A}$ 
 of 
$\pi_1(\Gamma_{\hspace{-0.05cm}\mathcal A})$ 
 spanned by elements of the form $PLP^{-1}$ as above  (and it is also natural to require additionally that the loops $L$ in such elements be,  in some way, isomorphic to the unique length 5 cluster period of type $A_2$ corresponding to Abel's identity).  Thus, one might think that Property $\star$ might be related (or even better, might be equivalent) to the following statement: 
 {\it `for any cluster period $\mathcal P$ in ${\mathcal A}$, the dilogarithmic identity $(\mathcal R_{\mathcal P})$ is accessible from dilogarithmic identities in two variables'.}
  \mk 
 
 It is known that Property $\star$ holds true in several interesting cases\footnote{See \cite{KY} for more details and for some references.}, such as
 \begin{itemize}
 \item[$(1)$] for any seed of a finite type cluster algebra. Moreover, in this case, the loops $L$ in the statement of Property $\star$ correspond to sequences of mutations of rank two
skew-symmetrizable matrices of type either $A_1 \times 
A_1$, $A_2$, $B_2$ or $G_2$;
 \item[$(2)$]  for any seed generated from a signed adjacency
matrix of an ideal triangulation of a bordered surface with
marked points, except when the surface is a closed surface
with exactly two punctures.
 \end{itemize}
 
Since the dilogarithmic identities associated to cluster periods of rank 2 finite cluster algebras all are accessible from Abel's identity (which is associated to $A_2$), it might follow
 from (1) above that the identity $(\mathcal R_\mathcal P)$ for any period $\mathcal P$ of a cluster algebra of finite type is accessible from $(\boldsymbol{\mathcal A b})$.

\begin{questions}
 {\bf (1).} 
If the preceding statement holds true, give a rigorous proof of it.\sk 

{\bf (2).} Under the reasonable assumption that it indeed holds true, give 
 a way (or better, an algorithm) allowing to express any dilogarithmic identity  associated to a period of a finite type cluster algebra as a (finite) linear combination of 5-terms Abel's identities.
\end{questions}

An interesting result ({\it cf.}\,\cite[Proposition 4.1]{KY}) is that 
 Property $\star$
  does not hold for the seed associated to a certain ideal triangulation of a closed surface (of genus
1 or higher) with two punctures. And explicit examples, one for each genus $g\geq 1$, are given in \cite[\S4]{KY}. \sk
  
  We consider here the simplest example which corresponds to the case when $g=1$.  
  In Example 4.1 of \cite{KY}, the authors describe a 
  triangulation of the 2-punctured torus to which is associated the following exchange matrix
%
%
$$
\begin{bmatrix} 
0&0&-1&1&-1&1\\ 
0&0&1&-1&1&-1\\
1&-1&0&1&0&-1\\
-1&1&-1&0&1&0\\
1&-1&0&-1&0&1\\
-1&1&1&0&-1&0
\end{bmatrix}
$$

Then one verifies that the following 32-tuple  
\begin{equation}
\label{Eq:PeriodKY}
(5, 6, 4, 3, 6, 5, 1, 2, 4, 3, 6, 5, 3, 4, 2, 1, 6, 5, 3, 4, 5, 6, 1, 2, 3, 4, 5, 6, 4, 3, 2, 1)
\end{equation}
is a period for the $\boldsymbol{\mathcal X}$-seed 
having the matrix above for exchange matrix. The web associated to it, denoted by 
$\boldsymbol{\mathcal W}_{KM,1}$ (where the suscript 1 referers to the genus $g=1$), is a  32-web in six variables  which can be verified to be such that 
$$\rho^\bullet\Big(\boldsymbol{\mathcal W}_{KM,1}\Big)
=\big(26,11,1 \big)
\, , 
\qquad {\rm polrk}^\bullet
\Big(\boldsymbol{\mathcal W}_{KM,1}\Big)=(32,1)
\qquad \mbox{and}
\qquad {\rm rk}\Big(\boldsymbol{\mathcal W}_{KM,1}\Big)
 =33\, . $$  
Thus if all the  ARs of $\boldsymbol{\mathcal W}_{KM,1}$ are polylogarithmic (of weight 1 or 2), this web is not AMP. \sk

On the other hand, it should be noted that  \eqref{Eq:PeriodKY} is precisely an example of a period which, when seen as a loop in the exchange graph, does not belong to the subgroup $\pi'$  of the fundamental group of the exchange graph described above.  Whether it is concomitant with the fact that the web $\boldsymbol{\mathcal W}_{KM,1}$ associated to   \eqref{Eq:PeriodKY} is not AMP makes us wonder how these two facts are related. \sk

These considerations suggest the following questions which, because of the interest they have for us, are worth formulating  for any cluster period: 
\begin{question} Let $\mathcal P$ be a period of a cluster algebra.
 How are related the following properties?
\begin{enumerate}
\item  The cluster web $\boldsymbol{\mathcal W}_{\mathcal P}$ is AMP with only polylogarithmic ARs (of weight 1 and 2).
\item As a loop, $\mathcal P$ belongs to the subgroup $\pi'$ of  the fundamental group of the exchange graph.
 \item  The dilogarithmic identity $(\mathcal R_{\mathcal P})$ is accessible from  cluster dilogarithmic identities in one or two variables of type 
 $A_1$, $A_1\times A_1$, $A_2$, $B_2$ and $G_2$, 
  in which all the arguments of the dilogarithm functions involved are $\boldsymbol{\mathcal X}$-cluster variables (of the cluster algebra considered). 
 \end{enumerate}
\end{question}

If we do think that these three properties are related, some examples we are aware of show that they are not equivalent. Such an example is provided by the $\boldsymbol{\mathcal Y}$-cluster web of type $D_4$. 
The associated period is $\boldsymbol{i}^{half}_{D_4,A_1}$ ({\it cf.}\,\eqref{Eq:period-slice})  for which 
property {\it 2.}\,above is known to hold true  (because of finite type) as well as the third (according to \cite{KY}\footnote{{\it Cf.}\,the proof of Lemma 6.2 in \cite{KY}, where it is said that the accessibility of (the quantum version) of 
the $\boldsymbol{\mathcal Y}$-cluster dilogarithmic identity $(\mathcal R_{D_4})$  can be proven easily. It would be more satisfying to have an explicit proof of this.}). However, it follows from 
\eqref{Eq:YWD4-invariants} that the associated cluster web 
$\boldsymbol{\mathcal Y\hspace{-0.05cm}\mathcal W}_{D_4}$ is not AMP. 

  \subsection{\bf \bf Some (conjectural) properties of cluster variables}
 \label{SS:Some-conjectures}
 Considering  finite collections of cluster variables from the perspective of web geometry naturally leads to look at a given cluster variable on its own from a new perspective,  more differential geometric in nature,
 which hasn't been done until now.\sk 
 
 Indeed, 
 given a cluster variable $x$, 
we are actually more interested in the foliation $\mathcal F_x$  that it defines, than in $x$ itself.  And given another such variable $x'$, it is crucial, for instance in order  to get a better understanding of the iterated integrals ARs of a web containing both $ \mathcal F_x$ and $\mathcal F_{x'}$, to know more about 
the algebraic subsets of the ambiant space  which are invariant by these two foliations. \sk

More specifically, 
given a web $\boldsymbol{\mathcal W}=\boldsymbol{\mathcal W}(x_1,\ldots,x_d)$ defined by $d$ cluster variables $x_i$ (of a same fixed cluster algebra),  the following questions naturally arise  and appear to us as very relevant from a web-theoretic perspective: 
\begin{enumerate}
\item Is $x_i$ a primitive first integral for the foliation $\mathcal F_{x_i}$ it defines? If yes, what can be said about the singular fibers $x_i^{-1}(\lambda)$ and the corresponding singular values $\lambda\in \mathbf P^1$?
\sk 
\item What can be said when two foliations $\mathcal F_{x_i}$ and $\mathcal F_{x_j}$ coincide?  Does this imply that the two cluster variables $x_i$ and $x_j$ defining them   coincide or are related in a simple way?
\sk 
\item Assuming that $\mathcal F_{x_i}$ and $\mathcal F_{x_j}$  are distinct foliations, what can be said about the invariant subsets $\Sigma$ for both? And for the corresponding values $x_i(\Sigma)$ and $x_j(\Sigma)$?
\sk 
\end{enumerate}

During the preparation of this text, we have considered an important number of cluster webs (those of the preceding subsection  forming only a part of it) the study of which gave us rather precise ideas about what could be the answers to the above questions.  In order to discuss our  insights more in depth, we first recall some classical notions about (complex) rational functions. 
\mk

\subsubsection{Reminder about rational functions.}
Let $R$ be a rational function in $n$ variables: for a fixed $n$-tuple $t=(t_i)_{i=1}^n$ of indeterminates, one has $R(t)=P(t)/Q(t)$ where $P,Q\in \mathbf C[t]$ are two polynomials without common factors.  
 Then the following facts hold true\footnote{Regarding this matter, 
 for a glimpse of the general context, statements, proofs and references, we refer to \cite{BCN}.}:  
\begin{itemize}
\item[$-$] the general fiber 
under 
  $R$ is not connected if and only is $R$ is {\bf composite}, {\it i.e.}\,there exists  $\tilde R\in \mathbf C(t)$ and a non-invertible rational morphism $g: \mathbf P^1\rightarrow \mathbf P^1$  such that $R=g\circ \tilde R$; 
\item[$-$] thus  $R$ is a primitive first integral of the foliation $\mathcal F_R$ 
if and only if 
  it is non-composite; 
\item[$-$]  an element $\lambda $ of $\mathbf P^1$ is a {\bf remarkable value} for $R$ if the associated fiber $R^{-1}(\lambda)$  is not irreducible.  By definition, the {\bf spectrum $\boldsymbol{{\bf Spect}(R)}$} of $R$  is the set of its remarkable values;
\item[$-$] it is well-known that $R$ is non-composite if and only if 
${\rm Spect}(R)\subset \mathbf P^1$ is finite;
\item[$-$] there exists a homogeneous polynomial  in two variables 
${\rm spect}_{P,Q}(\lambda_0,\lambda_1)$ such that

${}^{}$ \hspace{2cm} ${\rm Spect}(R)=\Big\{\,  \big[\lambda_0:\lambda_1\big]\in \mathbf P^1\, \big\lvert \, 
{\rm spect}_{P,Q}(\lambda_0,\lambda_1)=0\, \Big\}$.

Moreover, ${\rm spect}_{P,Q}$ 
 can be constructed from $P$ and $Q$ in an effective way.
\end{itemize}

We stress that everything above is stated and makes sense relatively to a fixed set of indeterminates $t_i$. Given a birational map $\varphi: \mathbf P^n\dashrightarrow \mathbf P^n$, one can wonder how is related the spectrum of $R$ (relatively to the $t_i$'s) with the one of $R\circ \varphi$ (relatively to the indeterminates $\tilde t_i=\varphi^*(t_i)$ in this case).
It seems to us that it is difficult to give a general answer to this question. 
However, for $R$  fixed and $\varphi$ generic, it is easy to see that 
the  following inclusion holds true
$$ {\rm Spect} (R\circ \varphi)\subset {\rm Spect} (R\circ \varphi)\cup R\big( S_\varphi \setminus I_{R}\big)
\subset \mathbf P^1
\, ,  $$
where $S_\varphi$ and $I_R$ 
stand for the set of singular values  of $\varphi$ 
and the indeterminacy locus of $R$  respectively.\footnote{We recall that  $S_\varphi$ is defined as the image by $\varphi$ of the jacobian divisor $J_\varphi$ of this rational map. It is an algebraic subset of the target projective space, of codimension at least 2.}

\subsubsection{The case of cluster variables.}
\label{SS:case-cluster-var}
In what follows, we fix a certain cluster algebra (of rank $n\geq 2$) defined by means of an initial seed with initial quiver $Q^0$.  We will mainly focus on the case of $\boldsymbol{\mathcal X}$-cluster variables, while that  of $\boldsymbol{\mathcal A}$-cluster variables, which is similar in many ways,  will only be mentioned briefly.   
\mk 

We use below the following notation: we denote by 
$u_1,\ldots,u_n$ the initial $\boldsymbol{\mathcal X}$-cluster variables 
and accordingly, $\boldsymbol{u}=(u_i)_{i=1}^n$ stands for the initial cluster. 
We set $\mathbf C[\boldsymbol{u}]$ (resp.\,$\mathbf C(\boldsymbol{u})$) 
 for the polynomial ring $\mathbf C[{u}_1,\ldots,{u}_n]$ (resp.\,for
its fraction field $\mathbf C({u}_1,\ldots,{u}_n)$) and $\mathbf C[\boldsymbol{u}^{\pm 1}]$ denotes the ring of Laurent polynomials in the $u_i$'s.  To these complex algebras are associated the following geometric objects:\vspace{0.2cm} \\
${}^{}$\hspace{0.3cm} $-$ $\mathbf T_0={\rm Spec}\big(\mathbf C[\boldsymbol{u}^{\pm 1}]\big)$ : the initial 'cluster torus';\sk \\
${}^{}$\hspace{0.3cm} $-$ 
$\mathbf C^n_0={\rm Spec}\big(\mathbf C[\boldsymbol{u}]\big)$ : the initial 'cluster affine space';\sk 
\\
${}^{}$\hspace{0.3cm} $-$ $\mathbf P^n_0={\rm Proj}\big(\mathbf C^h[\boldsymbol{u},u_0]\big)$ : the initial `cluster projective space' ($u_0$ is an extra variable).\sk 

We have of course the following inclusions of these spaces 
$\mathbf T_0\subset  \mathbf C^n_0 \subset 
\mathbf P^n_0$.\mk 

 In the  initial projective space $\mathbf P^n_0$, we consider the following  divisors: 
the initial `cluster hyperplanes divisors' $H_j=\{u_j=0\}$ for $j=1,\ldots,n$, and the divisor at infinity $H_0=\{u_0=0\}$.  Next,  for any $F$-polynomial $F$ of the considered cluster algebra, 
one denotes by $H_F$ the effective divisor defined as the closure in $\mathbf P^n_0$ of  the irreducible affine hypersurface $\{F(\boldsymbol{u})=0\} \subset \mathbf C^n_0$. We recall the following important property of such a $F$: one has $F(0)=1$ and $F$ admits a unique monomial of the highest degree, which is moreover a multiple of all the other monomials appearing in it.  Finally, for any vertex $t\in \mathbb T^n$, we denote by $\boldsymbol{H}^t$ the union of the divisors $H_s$'s for $s=0,1,\ldots,n$ with the set 
of those 
cut out by the $F$-polynomials $F_{\hspace{-0.05cm}1}^t,\ldots,F_{\hspace{-0.05cm}n}^t$ at $t$: 
\begin{equation}
\label{Eq:Ht}
\boldsymbol{H}^t=\Big\{  \, H_0\, , \, H_i\, , \,  H_{F_{\hspace{-0.05cm}i}^t}\hspace{0.1cm} \big\lvert \hspace{0.2cm} i=1,\ldots,n\,
\Big\} \subset {\rm Div}\big(\mathbf P^n_0\big)\, .
\end{equation}

\paragraph{} 
\label{Par:ClusterVarNotComposite}

Let $x$ be a ($\boldsymbol{\mathcal X}$-)cluster variable that we consider as a rational function in the initial cluster variables $u_i$'s.  By definition, $x=x_{i}^t$ for some vertex $t$ of $\mathbb T^n$ and some $i\in \{1,\ldots,n\}$.  The rational map $\boldsymbol{u}\dashrightarrow 
\boldsymbol{x}^t=(x_{k}^t)_{k=1}^n$ being invertible (birational), it comes that $x\in \mathbf C(\boldsymbol{u})$ is not composite.  Hence the following lemma holds true: 
 \begin{lem}  
\label{L:primitive-first-integral} 
The cluster variable $x$ is a primitive first integral 
 of the foliation $\mathcal F_{\hspace{-0.1cm}x}$ it defines. 
 Moreover, it is homaloidal thus  its generic fiber is rational. 
 \end{lem}

An immediate remark coming after this is  that, regarding the foliations they define, cluster variables always come in pairs.  
 Indeed, from the basic formulae   \eqref{Eq:A-X-Mutation-formulae} for mutations, it comes that 
 the new cluster variable $x'$ obtained by mutating $x=x_{i}^t$ in the $i$-th direction is nothing else but $1/x$ as a rational map. 
  If $x$ and $x'$ are (formally) distinct cluster variables, they define the same foliation.   Consequently, from a web-theoretic perpective, we will deal with cluster variables up to inversion $x\leftrightarrow {1}/{x}$. 
\sk 

Before turning on the general case, let us consider the particular case when the cluster algebra we are dealing with is the one of finite type $A_n$.  \sk 

\paragraph{} We will see further in \S\ref{SSub:Type-An}  that 
 the $\boldsymbol{\mathcal X}$-cluster variety $\boldsymbol{\mathcal X}_{A_n}$ identifies birationaly to the moduli space $\mathcal M_{0,n+3}$ (in some ways which can be made precise, see \ref{SSub:Type-An} further) and that, up to such an identification, any 
$\boldsymbol{\mathcal X}$-cluster variable 
can be written  as a cross-ratio
\begin{equation}
\label{Eq:xij-cross-ratio}
x_{ijkl}(c)=\frac{(\xi_i-\xi_{j})( \xi_k-\xi_{l})}{
(\xi_i-\xi_{l})( \xi_{j}-\xi_{k})}
\end{equation}
 for some indices $i,j,k,l\in \{1,\ldots,n+3\}  $ cyclically ordered modulo $n+3$\footnote{\label{Footnote:coco}{\it I.e.}\,such that their images by the map $z\mapsto \exp({2\sqrt{-1}\pi\, z}/{(n+3)})$ appear clockwise on the unit circle in $\mathbf C$.}, where $c=[0,-1,\xi_1,\ldots,$ $\xi_n,\infty]_{i=1}^n
 \in \mathcal M_{0,n+3}$ with $(\xi_s)_{s=1}^n\in \mathbf C^n$. 
 Using the explicit (quadratic) formula \eqref{Eq:xij-cross-ratio} in terms of the 
 $\xi_s$'s which form a system of rational coordinates on $ \mathcal M_{0,n+3}$, it is then easy to study the foliation $\mathcal F_{x_{ijkj}}$ from a differential geometric point of view, when considered as a rational fibration on the projective space $\mathbf P^n_\xi$ associated to the  $\xi_i$'s. For instance,  it is easy to verify that  
\begin{equation}
 \label{Eq:prop-cross-ratios}
 \begin{tabular}{l}
 {\it  ${\bf 1.}$\,  The remarkable values of $x_{ijkl}$ are  $0$, $-1$ and $\infty$, {\it i.e.}\,${\rm Spect}(x_{ijkl})=\{0,-1,\infty\}$;}\vspace{0.1cm} \\
{\it ${\bf 2.}$\, Given another cluster variable $x_{i'j'k'l'}$, the two following facts hold true:}
\vspace{0.1cm}
\\
${}^{}$\hspace{0.1cm} ${\bf i.}$\hspace{0.15cm}{\it 
$\mathcal F_{x_{ijkl}}$ and $\mathcal F_{x_{i'j'k'l'}}$ 
coincide  
$\Leftrightarrow$   \,$(i,j,k,l)=(i',j',k',l')$  \hspace{0.1cm} or   \hspace{0.1cm}$(i,j,k,l)=(k',l',i',j')$}\\
{\it ${}^{}$\hspace{4.6cm}   $\Leftrightarrow$ \hspace{0.6cm} $x_{ijkl}=x_{i'j'k'l'}$\hspace{0.95cm} or \hspace{0.7cm} $x_{ijkl}=1/x_{{i'j'k'l'}}$;
} 
\vspace{0.1cm}\\
 ${}^{}$\hspace{0.1cm} ${\bf ii.}$\hspace{0.15cm}{\it if $\mathcal F_{x_{ijkl}}$ and $\mathcal F_{x_{i',j',k',l'}}$  are distinct then 
any biinvariant irreducible hypersurface}
\vspace{0.0cm}\\ 
${}^{}$\hspace{0.45cm} {\it $H$ of \,$\mathbf P^n_\xi$ 
is a component of reducible fibers of $x_{ijkl}$ and $x_{i'j'k'l'}$. Thus 
 $x_{ijkl}(H)$ and}
  \vspace{0.0cm}\\
${}^{}$\hspace{0.45cm} {\it $x_{i'j'k'l'}(H)$ are remarkable values hence both belong to
  $\{0,-1,\infty\}$.}
 \end{tabular}\hspace{0.0cm} {}^{}
 \end{equation}



 Regarding to what we have in mind, 
 a warning is in order about  the properties 
 listed above: 
 these are properties of cross-ratios on $\mathcal M_{0,n+3}$, stated when these 
 cross-ratios are considered as rational functions in the standard coordinates $\xi_i$'s on this moduli space. But obviously  it also follows  from \eqref{Eq:xij-cross-ratio} that the $\xi_i$'s, are not cluster coordinates, but only depend birationally from any $\boldsymbol{\mathcal X}$-cluster of type $A_n$. 
As a result,  it is a priori unclear whether or not  \eqref{Eq:prop-cross-ratios} admits an analog when $x_{ijkl}$ and $x_{i'j'k'l'}$ are now truly considered as cluster variables,
 that is as rational functions of the initial cluster variables $u_i$'s.    
 The point {\bf 2.i.} above can be obtained rather easily  ({\it cf.}\,Corollary \ref{Cor:X-var-web-An}.{\it 2}) and we believe that,  with some supplementary work, it  should be possible to prove statements similar to   {\bf 1.} and {\bf 2.ii.} for the $\boldsymbol{\mathcal X}$-cluster variables in type $A$   above (see 
 Remark \ref{Rem:X-cluster-var-properties-yoyo} further and more generally Conjecture \ref{Conj:Spectrum-Cluster-Variables}
 and Conjecture  \ref{Conj:Inv-divisor-Fx-Ftildex} below).

 
\paragraph{}
\label{Par:x0-(-1)-infty}
 We now come back to the general case  considered at the very beginning  of  this subsection: $x$ is a cluster variable  of the cluster algebra with $(\boldsymbol{u},Q^0)$ as initial $\boldsymbol{\mathcal X}$-seed.  One has $x=x_i^t$ for  some $i\in \{1,\ldots,n\}$ and a certain vertex $t\in \mathbb T^n$. 
 \mk

  Since the general fiber of $x$ is  irreducible (by Lemma \ref{L:primitive-first-integral}), it is natural to look at its  remarkable values.  
The  separation formula \eqref{Eq:X-var:separation-formula} for $x=x_{i}^t$ 
 is  inspiring in this respect.  Indeed, since it is meaningless to deal with its inverse $1/x$, one can assume that the tropical sign 
 $\epsilon_x\in \{\pm 1\}$ 
 of $x$ is positive: all the coefficients 
$c_{\hspace{-0.05cm} ji}^t$ of the  associated    $c$-vector $c_x=(c_{\hspace{-0.05cm} ji}^t)_{j=1}^n$ are non-negative and 
at least one of them is non-zero hence positive (by `sign-coherence'). 
Consequently, it follows that the fiber $x^{-1}(\lambda_0)$ of the pencil $\mathcal F_{\hspace{-0.05cm}x}$
 passing through the origin of $\mathbf C^n_0$ corresponds to the value $\lambda_0=0$ and, 
 up to the addition of a non-negative multiple of the divisor at infinity $H_0$,  is the effective divisor $\sum_{ j=1 }^n  c_{\hspace{-0.05cm}ji}^t  \,    H_j=
\sum_{ c_{\hspace{-0.05cm}ji}^t>0 }  c_{\hspace{-0.05cm}ji}^t  \,    H_j$.  
Hence  any irreducible component of $x^{-1}(0)$ belongs to $\boldsymbol{H}^t$ and except  in the particular case when $\lvert c_x\lvert=\sum_j c_{\hspace{-0.05cm}ji}^t=1=\deg(x)$,  we have that $0$ belongs to the spectrum of $x$. 
\mk 

We now consider the fiber $x^{-1}(\infty)$.
Still assuming that $\epsilon_x=1$, 
it follows again from    \eqref{Eq:X-var:separation-formula} that  in the affine chart $\mathbf C^n_0$, this fiber is cut out  by $\prod_{b_{ki}<0} (F_{\hspace{-0.05cm}k}^t)^{-b_{ki}}=0$.  Hence 
all the irreducible components of $x^{-1}(\infty)$ belong to $\boldsymbol{H}^t$ as well  and except in another particular case, it follows that $\infty$ is also a remarkable value of $x$. \mk 
 
 Eventually, we deal with the fiber of $x$ over $-1$ which is cut out by the equation $x+1=0$.  Since the initial quiver $Q^0$ is connected, the same holds true for for $Q^t$ (this property being invariant by mutations) hence there exists $j\in \{1,\ldots,n\}$ distinct from $i$ such that $
  b_{\hspace{-0.03cm}ij}=
 b_{\hspace{-0.03cm}ij}^t\neq 0$. Assume that $ b_{\hspace{-0.03cm}ij}$ is negative (the case when it is assumed positive can be treated along the same arguments) and 
 let $x_j=x_j^t$ be the $k$-th cluster variable of the cluster $\boldsymbol{x}^t$ to which $x=x_i^t$ belongs.  Then according to \eqref{Eq:A-X-Mutation-formulae}, 
 the mutation  $x_j'=\mu_i(x_j)$ of $x_j$ in the $i$-th direction is given by $x'_j=x_j\big(1+x\big)^{-b_{\hspace{-0.02cm}ij}}$ 
hence one has $1+x=(x'_j/x_j)^{-1/b_{\hspace{-0.02cm}ij}}$ in the field of algebraic functions in the $u_i$'s.\footnote{This elementary argument is taken from the proof of the `Theorem on $1+\mathcal X$ variables' of 
 \cite[Appendix\,B]{PSSV}. 
For the cluster variables of the bipartite belt of a cluster algebra of finite type, a more precise result can be deduced from formula (2.18) of \cite{FZ}.} 
 Thus $1+x=0$ implies that $x_j'=0$ or $x_j=\infty$ and considering the two preceding paragraphs (regarding $x_j'$ and $x_j$ respectively), we deduce that any irreducible component of 
 $x^{-1}(-1)$ 
  belongs to $\boldsymbol{H}^t\cup \boldsymbol{H}^{t'}$, where $t'$ stands for the vertex in $\mathbb T^n$ obtained from $t$ by the $i$-th mutation.  Thus $-1$ as well is a good candidate for being a remarkable value for the considered cluster variable $x$.

\paragraph{A conjectural description of reducible fibers of cluster variables.} 
 The above considerations, first about the finite case of type $A_n$, then in the case of an arbitrary cluster algebra, have led us to state several conjectures about some differential and geometric properties of the foliations defined by cluster variables. We begin by stating these conjectures before giving a few comments regarding them. \sk 
 
 In order to be precise,  we recall the nature of the foliations we are dealing with: $x$ is a cluster variable hence is considered as a rational function  
 in the initial cluster variables $u_1,\ldots,u_n$ and consequently, we see 
   $\mathcal F_{\hspace{-0.03cm}x}$  as a rational fibration on the initial cluster projective space $\mathbf P^n_0$.   To make some statements shorter and
   although
    this is not standard, 
  any one of the $u_i$'s will be also called a  `$F$-polynomial' below. 
   \mk 
   
   Our first conjecture is about the decomposable fibers of $\mathcal F_{\hspace{-0.03cm}x}$: \sk
 \begin{conjecture}  
 \label{Conj:Spectrum-Cluster-Variables}
\vspace{-0.2cm}
\begin{itemize}
\item[(a).] 
The  set of remarkable values of $x$  is included in $\{ \, 0\, , \,  -1\, , \, \infty\, \}$. 
\item[(b).]  The irreducible components of the reducible fibers of $x$ are cut out by some $F$-polynomials.
\end{itemize}
\end{conjecture}
\sk

The statements in this conjecture deserve to be commented. 
\begin{itemize}
\item A first remark is that only an inclusion can hold true in {\it (a)}: taking for $x$ one of the initial cluster variables $u_i$, we see that there is no singular or remarkable value at all. But one can get a proper inclusion ${\rm Spect}(x)\subset 
\{ \, 0\, , \,  -1\, , \, \infty\, \}$ even for a  non-initial cluster variable $x$.

 For instance,  $u_1u_2/(1+u_2)$ is a  $\boldsymbol{\mathcal X}$-cluster variable of the cluster algebra of finite type $A_2$ and it has only two remarkable values. Indeed, one has  ${\rm Spect}\big(u_1u_2/(1+u_2)\big)=\{0,\infty\}$.
 \item  The previous remark naturally leads to wonder about the case when $x$ has no remarkable value. This will be discussed just below ({\it cf.}\,Conjecture \ref{Conj:Initial-Cluster-Variables}).
  \item  
 As it follows from  the considerations in paragraph \S\ref{Par:x0-(-1)-infty} above, any of the fibers $x^{-1}(0)$, $x^{-1}(-1)$ or $x^{-1}(\infty)$
 can be written as a sum with positive integer coefficients of divisors $H_F$ for some $F$-polynomials $F$ of the considered cluster algebra.  Since any divisor $H_F$ is irreducible according to Theorem \ref{Thm:CK},  we obtain that  {\it (a)} implies {\it (b)}.
%
%
 \item 
It is natural to wonder if there exists a version of the above conjecture  for 
$\boldsymbol{\mathcal A}$-cluster variables.  It might be so, but the statements should certainly be modified, as the sfollowing examples (concerning assertion {\it (b)}) show:  
${(1+y+x^2 )}/{xy}$ is a $\boldsymbol{\mathcal A}$-cluster variable in finite type $B_2$ and it can be verified that its remarkable values are $i$, $-i$ and $\infty$: one has ${\rm Spect}\big( {(1+y+x^2)}/{xy}\big)
=\big\{ \, \pm i \, , \, \infty \big\}$.  Another example is given by 
${(1+y+x^3)}/{xy}$ which is a $\boldsymbol{\mathcal A}$-cluster variable
 of type $G_2$: one has ${\rm Spect}\big( {(1+y+x^3)}/{xy}\big)
=\big\{ \,\infty \, , -1 \, , \, ({1\pm i \sqrt{3}})/{2}\, \big\}$. 
%

The two  preceding examples suggest a generalization of the previous conjecture: 
 let $\mathscr Q$ be a weighted quiver. Then for any $\mathcal A$-cluster variable $a$ of $A_{\mathscr Q}$, one has: ${\rm Spect}(a)\subset \{0,\infty\} \cup  \boldsymbol{\mu}_q $ where $q$ stands for a certain positive integer obtained from the weights in $\mathscr Q$ (for instance, it could be the least common multiple of these weights). 
\end{itemize}

\paragraph{A conjectural characterization of initial cluster variables by  their spectrum.} 
A question that emerges naturally from the above conjecture and comments 
  is that of characterizing the cluster variables in terms of the cardinality of their spectra.  Thinking about that and looking at numerous explicit examples led us to state the following conjecture about  the characterization of initial cluster variables (where, as above, $x$ stands for a $\boldsymbol{\mathcal X}$-cluster variable): 
 \begin{conjecture}  
 \label{Conj:Initial-Cluster-Variables}
 The following statements are 
equivalent: 
\begin{itemize}
\item[$-$] ${\rm Spect}(x)=\emptyset$, 
 {\it i.e.}\,$x$ has no remarkable value;
 \item[$-$] $\big \lvert\, {\rm Spect}(x)\, \big\lvert \leq 1$, 
 {\it i.e.}\,$x$ has at most one remarkable value;
\item[$-$] $\deg(x)=1$, {\it i.e.}\,$x$ is fractional linear in the 
$u_i$'s;
\item[$-$] $x$ or its inverse $1/x$ coincides with one of the initial  cluster variables.
\end{itemize}
\end{conjecture}
 
 Several authors have introduced the notion of `{\it compatibility degree'} for cluster variables which could be relevant to investigate the validity of the above conjecture (see \cite{CKQ} for an overview and some references).   Assuming this conjecture  holds true as well as point {\it (a)} in Conjecture \ref{Conj:Spectrum-Cluster-Variables} naturally leads to wonder about the cluster variables whose spectrum has cardinality only 2: what can be said 
 about them in addition of this  single hypothesis? 

\paragraph{A conjecture about hypersurfaces invariant two cluster fibrations.}
The material above is about some (conjectural) geometric or algebraic properties of the fibration 
 $\mathcal F_{\hspace{-0.05cm}x}$ 
 defined by a  cluster variable $x$, considered on its own.  We now discuss and propose some interesting but conjectural statements about how 
 two such fibrations
 are related.  \sk 
 
 
 Considering the  description of the fibers of cluster variables over $0, -1$ and $\infty$ in \S\ref{Par:x0-(-1)-infty},  it comes that elements of 
 $\boldsymbol{H}^t$ ({\it cf.}\,\eqref{Eq:Ht}) for some $t\in \mathbb T^n$ are good candidates for being invariant by two distinct cluster fibrations 
 $\mathcal F_{\hspace{-0.05cm}x}$  and $\mathcal F_{\hspace{-0.05cm}\tilde x}$  and that the corresponding values for $x$ and $\tilde x$ are elements of $\{0, -1 , \infty\}$.   
 In view of this, the question arises as to whether there exist fibers $x^{-1}(\lambda)$ and $\tilde x^{-1}(\tilde \lambda)$ 
 sharing an irreducible component, 
 for $\lambda $ or $\tilde \lambda$ not in $\{0, -1 , \infty\}$. About this, an extensive investigation of many explicit cases led us to make the following
\begin{conjecture}
\label{Conj:Inv-divisor-Fx-Ftildex}
Let  $\tilde{x}$  be another $\boldsymbol{\mathcal X}$-cluster variable (of  the same
 cluster algebra as $x$).  
\begin{enumerate}
\item The 
 fibrations 
$\boldsymbol{\mathcal F}_{{\hspace{-0.1cm}}x}$ and $\boldsymbol{\mathcal F}_{{\hspace{-0.1cm}}\tilde{x}}$ coincide if and only if $\tilde{x}=x$ or $\tilde{x}=x^{-1}$.
\item   Assuming  that $\boldsymbol{\mathcal F}_{{\hspace{-0.1cm}}x}$ and $ \boldsymbol{\mathcal F}_{{\hspace{-0.1cm}}\tilde{x}}$ are distinct, the following holds true:   \sk \\ 
given $\lambda,\tilde{\lambda}\in \mathbf P^1$, if the two fibers 
$\{x=\lambda\}$ and  $\{\tilde{x}=\tilde{\lambda}\}$ have an  irreducible component $H$ in common then both $\lambda$ and $\tilde{\lambda}$ belong to $\{ \, 0 ,  -1 ,  \infty\, \}$ and $H$ is cut out by a $F$-polynomial.
\end{enumerate}
\end{conjecture}

It is worth commenting these two statements:
\begin{itemize}
\item  The first statement being birationally invariant, one can assume that $x$ is precisely one of the initial cluster variables $u_i$. On the other hand, since $\tilde x$ is non composite (by Lemma \ref{L:primitive-first-integral}), one has $\tilde x=g(x)$ for a certain projective transformation 
$g\in {\rm PSL}_2(\mathbf C)$ if $\boldsymbol{\mathcal F}_{{\hspace{-0.1cm}}x}$ and $ \boldsymbol{\mathcal F}_{{\hspace{-0.1cm}}\tilde{x}}$ coincide.   Then {\it 1.}\,follows from 
our above Conjecture \ref{Conj:Inv-divisor-Fx-Ftildex}; 
\item  We recall  that the initial cluster variables $u_i$ are considered as $F$-polynomials in {\it 2}.
\item 
 A consequence of the second statement is that actually $\boldsymbol{\mathcal F}_{{\hspace{-0.1cm}}x}$ and $ \boldsymbol{\mathcal F}_{{\hspace{-0.1cm}}\tilde{x}}$ necessarily coincide as soon as they share a single leaf of the type $x^{-1}(\lambda)=\tilde{x}^{-1}(\tilde \lambda)$ for some 
 $\lambda,\tilde \lambda \not \in \{ 0,-1,\infty\}$.  
  This is quite stronger than statement {\it 1.}\,which can also be read as the fact that if the generic leaf of $\boldsymbol{\mathcal F}_{{\hspace{-0.1cm}}\tilde x}$ is also invariant by $\boldsymbol{\mathcal F}_{{\hspace{-0.1cm}}x}$ 
  then $\tilde x$  coincides with $x$ or $x^{-1}$.
\item  The first statement in the above conjecture can be proven to hold true for any cluster algebra of finite type, see Proposition \ref{Prop:propro} further (see also 
Corollary   \ref{Cor:Y-web-type-(Am,An)} for some other results in this direction).
\item   
 We don't see how to prove 
the second statement 
in full generality, even when restricting to the better understood case of cluster algebras of finite type. 
\end{itemize}
\begin{center}
$\star$
\end{center}

The conjectural properties of cluster variables stated in the three preceding conjectures  have to be compared with those  of the cross-ratio stated in \eqref{Eq:prop-cross-ratios}.  From this, one can state the following general (but rather vague hence non truly mathematical) motto:  
 \begin{center}
{\it  in many ways,  $\boldsymbol{\mathcal X}$-cluster variables behave like cross-ratios.}
 \end{center}
 
 It would be interesting to make this more rigorous and to know 
   to which  extent it holds true.

\paragraph{Consequences of the preceding conjectures for cluster webs.}
\label{Par:ClusterWeb-Polylogarithmic?}
 The previous conjectures have interesting consequences regarding the ARs with iterated integrals components of  cluster webs. 
We shall now discuss these consequences.   \sk 
 
    Let $\boldsymbol{\mathcal W}$ be a cluster web defined by a finite set $\{ x_s\}_{ s\in S}$
     of cardinality $d$ of $\boldsymbol{\mathcal X}$-cluster variables.  We assume that 
    the corresponding cluster variables $x_s$  and $x_{s'}$ associated to  two distinct elements $s,s'\in S$ are not the same (as rational functions of the fixed initial cluster coordinates), even up to inversion. Then $\boldsymbol{\mathcal W}$ is a $d$-web in $n$ variables (by Conjecture \ref{Conj:Inv-divisor-Fx-Ftildex}.{1}). \sk
    
     Let $\boldsymbol{H}_{\boldsymbol{\mathcal W}}$ be the set of $F$-polynomials such that the principal divisor $H_F$ appears with positive multiplicity in a fiber $x_s^{-1}(\lambda)$ for some points $s$ in $S$ and some values $\lambda \in \{0, -1,\infty\}$.  Then it follows from the second point of  Conjecture \ref{Conj:Inv-divisor-Fx-Ftildex} that the 
    irreducible components of the set  
      of common leaves $\Sigma^c({\boldsymbol{\mathcal W}})$  of 
${\boldsymbol{\mathcal W}}$ (discussed in \S\ref{Par:SSet-UnionCommonLeaves} above), 
 is included in  $\boldsymbol{H}_{\boldsymbol{\mathcal W}}$. 
By Proposition \ref{P:AnalyticConinuation}, this implies that for any $s\in S$, the $x_s$-ramification locus of $\boldsymbol{\mathcal W}$ 
is included in $\{0,-1,\infty\}$ hence is polylogarithmic ({\it cf.}\,Definition \ref{Def:Polylogarithmic-Ramification}). 
Thus, assuming than the conjecture stated in the preceding paragraphs are satisfied,  one can deduce that the same holds true for the following one: 
\begin{conjecture}
A $\boldsymbol{\mathcal X}$-cluster web has polylogarithmic ramification: all its ramification locus are included in $\{0,-1,\infty\}$ and consequently,  all its iterated integral ARs are polylogarithmic.
\end{conjecture}

Assuming that this conjecture holds true, the next step would be to have a more precise guess about the (virtual/polylogarithmic/standard) rank and the (polylogarithmic) abelian relations of a $\boldsymbol{\mathcal X}$-cluster web.  
In the next sections, we are going to investigate webs associated to Dynkin diagrams from this perspective.

\newpage

\section{Some general properties of cluster webs 
associated in finite type}
\label{S:ClusterWebs-FiniteType}
This section is devoted to the $\boldsymbol{\boldsymbol{\mathcal X}}$- and $\boldsymbol{\boldsymbol{\mathcal Y}}$-cluster webs of finite type. We first recall or introduce some notations about them that will be used in the sequel, and state several questions regarding their rank and ARs which are worth studying according to us.  After that in \S\ref{SS:DifferentialIndependenceClusterVariables}, we establish some very basic properties of these webs.

 \subsection{\bf Some notation and several basic questions}
\label{SS:SomeNotations}
First we recall/collect some notation that will be used in the sequel.  Next we make explicit the basic questions about cluster webs in finite type we are interested in.

\subsubsection{}
\label{SS:???}
{\bf Root systems:}
\sk

${}^{}$\quad $-$  $\Delta$ stands for a given irreducible root system of rank $n\geq 2$ (thus $A_1$ is excluded); \sk

${}^{}$  \quad $-$ we also use the same notation $\Delta$ to denote either the associated Dynkin diagram or the

${}^{}$\hspace{0.65cm}  corresponding (weighted bipartite) quiver (denoted by ${\vec{\Delta}}$ in \S\ref{Par:Cluster-Algebras-Finite-Type} above). This abuse of

${}^{}$\hspace{0.65cm} notation will not cause any lack of understanding;

${}^{}$\quad $-$ $\alpha_1,\ldots,\alpha_n$ denote fixed 
simple roots of $\Delta$;\sk

${}^{}$\quad $-$ 
$\Delta_{>0}$ (or $\Delta^+$) stands for the set of positive roots;  
\sk

${}^{}$\quad $-$ by definition, 
$\Delta$ is simply laced if all its roots have same length. Otherwise (that is when $\Delta$ 

${}^{}$\hspace{0.65cm} is of type $B_n,C_n, F_4$ or $G_2$) there are precisely two kinds of roots, depending on their length: 

${}^{}$\hspace{0.65cm}
the short ones and the long ones, whose sets are denoted by $\Delta^{short}$ and $\Delta^{long}$ respectively;
\sk

${}^{}$\quad $-$
the map $\alpha\mapsto \alpha^\vee=2\alpha/\langle \alpha,\alpha\rangle $ induces a duality between $\Delta$ and the dual root system $\Delta^\vee$. We 

${}^{}$\hspace{0.65cm} recall that $\Delta$ and $\Delta^\vee$ are isomorphic, except for $B_n$ and $C_n$ when $n\geq 3$;\footnote{On has $B_2\simeq C_2$ but $B_n^\vee=C_n$  is not isomorphic  with $B_n$ as a root system as soon as $n\geq 3$.}
\sk 

${}^{}$\quad $-$ for $\Delta$ simply laced ({\it i.e.}\,of type $A,D$ or $E$), one has $\alpha= \alpha^\vee$ but otherwise the duality is not 

${}^{}$\hspace{0.65cm} trivial (for instance it exchanges the short roots with the long ones {\it i.e.}\,$(\Delta^{short})^\vee=(\Delta^\vee)^{long}$);
\sk

${}^{}$\quad $-$ given a root $\alpha=\sum_{i=1}^n n_i\alpha_i\in \Delta$ (with $n_i\in \mathbf Z$ for any $i=1,\dots,n$), one sets $u^\alpha=\prod_{i}u_i^{n_i}$.\mk

\subsubsection{\bf Cluster notation:} ${}^{}$ \sk
\label{SS:Cluster notation}

${}^{}$\quad $-$ $A_\Delta$ denotes the cluster algebra with initial seed $(\boldsymbol{a},\boldsymbol{x},B_\Delta)$ where $
\boldsymbol{a}=(a_i)_{i}^n$ and $B_\Delta$ is the ex- 

${}^{}$\hspace{0.65cm}
change matrix associated to  $\vec{\Delta}$ ({\it cf.}\,page \label{page:BDelta}). As for the initial $\boldsymbol{\mathcal X}$-cluster, we will take 
either  

${}^{}$\hspace{0.65cm} $\boldsymbol{x}=(u_i)_{i=1}^n$ or $\boldsymbol{x}=(u_i^{\epsilon(i)})_{i=1}^n$, where in the second case $\epsilon:\{1,\ldots,n\}\rightarrow \{\pm 1\}$ is the function 

${}^{}$\hspace{0.65cm} corresponding to the bipartition of $\Delta$: one has $\epsilon(i)=-1$ if and only if the $i$-th vertex of the 

${}^{}$\hspace{0.65cm} bipartite quiver $\Delta$ is a source; 
\sk

${}^{}$\quad $-$  for $\boldsymbol{Z}\in \{ 
\, \boldsymbol{\mathcal A} \, 
, \, \boldsymbol{\mathcal X}\,  \}$, we denote by   
${\boldsymbol{Z \hspace{-0.05cm}\mathcal W}}_{\Delta}$  the cluster web 
(on the initial 
$\boldsymbol{Z}$-cluster torus $\boldsymbol{Z}\mathbf T^n$)   

${}^{}$\hspace{0.65cm}  defined  by  all   
 the $\boldsymbol{Z}$-cluster variables of $A_{\Delta}$;\sk 

${}^{}$\quad $-$  ${\boldsymbol{\mathcal Y\mathcal W}}_{\Delta}$ stands for  the subweb  of ${\boldsymbol{\mathcal X\mathcal W}}_{\Delta}$ defined by all the $\boldsymbol{\mathcal Y}$-cluster variables, which are the 

${}^{}$\hspace{0.65cm}  $\boldsymbol{\mathcal X}$-cluster variables
appearing in  
the bipartite belt of the cluster algebra $A_\Delta$ (see \S\ref{SubPar:BipartiteBelt});\sk


${}^{}$\quad $-$  the set 
$\boldsymbol{\mathcal Y}_{\Delta}$ of  $\boldsymbol{\mathcal Y}$-cluster variables 
up to inversion ($y\leftrightarrow y^{-1}$) is in bijection with $\Delta_{\geq -1}$. More 

${}^{}$\hspace{0.65cm}precisely, for any $\alpha\in \Delta_{\geq -1}$, there exists a unique $\boldsymbol{\mathcal Y}$-cluster variable, denoted by $y[\alpha]$, which 

${}^{}$\hspace{0.65cm}  can be  written $y[\alpha]=P[\alpha]/u^{\alpha^{\vee}}$ relatively to the initial $\boldsymbol{\mathcal X}$-cluster 
$\boldsymbol{x}=(u_i^{\epsilon(i)})_{i=1}^n$, where $P[\alpha]$

${}^{}$\hspace{0.65cm} is
a polynomial with positive integer coefficients (in the $u_i$'s) 
with constant term equal to 1

${}^{}$\hspace{0.65cm}  ({\it cf.}\cite[Theorem 1.5]{FZ}). Consequently, one has 
\begin{equation}
\label{Eq:YWDelta}
{\boldsymbol{\mathcal Y\mathcal W}}_{\hspace{-0.08cm}\Delta}=
{\boldsymbol{\mathcal W}}\left( \, 
 y[\alpha] \hspace{0.15cm} \big\lvert  
 \hspace{0.15cm} 
  \alpha\in \Delta_{\geq 1} \, 
\right)
=  {\boldsymbol{\mathcal W}}\left( \, u_i\, , \, 
 y[\alpha] \hspace{0.15cm} \big\lvert  
 \hspace{0.15cm}  i=1,\ldots,n,\, 
  \alpha\in \Delta_{>0} \, 
\right)
\, ; 
\end{equation}

${}^{}$\quad $-$  when $\Delta$ is not simply laced, 
${\boldsymbol{\mathcal Y\mathcal W}}_{\hspace{-0.08cm}\Delta^{short}}$ (resp.\,${\boldsymbol{\mathcal Y\mathcal W}}_{\hspace{-0.08cm}\Delta^{long}}$) is the subweb of 
${\boldsymbol{\mathcal Y\mathcal W}}_{\hspace{-0.08cm}\Delta}$ 
formed by the 

${}^{}$\hspace{0.65cm}  ${\boldsymbol{\mathcal Y}}$-cluster variables  indexed by the short roots  (resp.\,by the 
long roots) of $ \Delta$:  one has 
$$
{\boldsymbol{\mathcal Y\mathcal W}}_{\hspace{-0.08cm}\Delta^{short}}=
{\boldsymbol{\mathcal W}}\left( \, 
 y[\alpha] \hspace{0.15cm} \big\lvert\, \alpha\in \Delta^{short}_{\geq -1} \, 
\right)\qquad 
\mbox{ and }
\qquad 
{\boldsymbol{\mathcal Y\mathcal W}}_{\hspace{-0.08cm}\Delta^{long}}=
{\boldsymbol{\mathcal W}}\left( \, 
 y[\alpha] \hspace{0.15cm} \big\lvert\, \alpha\in \Delta^{long}_{\geq -1} \, 
\right)
$$
${}^{}$\hspace{0.65cm} with $\Delta^{-\hspace{-0.05cm}-}_{\geq -1}=\Delta_{\geq -1}\cap \Delta^{-\hspace{-0.07cm}-}$, where  $-\hspace{-0.05cm}$ stands for anyone of the two words $short$ or $long$;\mk

${}^{}$\quad $-$  for $\boldsymbol{Z}\in \{ 
\, \boldsymbol{\mathcal A} \, 
, \, \boldsymbol{\mathcal X}\,  , \, \boldsymbol{\mathcal Y} \, \}$, we denote by   
$d_\Delta^{\boldsymbol{Z}}$
  the degree  of the web ${\boldsymbol{Z\mathcal W}}_{\Delta}$ 
  (that is, the number of 
  
${}^{}$\hspace{0.65cm}   distinct foliations appearing in it).  
\mk

${}^{}$\quad $-$
  since $
\lvert \Delta_{\geq -1}\lvert=\lvert \Delta_{>0}\lvert+n=hn/2+n=h(n+2)/2$,\footnote{This follows from the  classical formula for the Coxeter number $h=h(\Delta)=2\lvert \Delta^+\lvert/n$, see \cite[\S3.18]{Humphreys}.}, it follows from \eqref{Eq:YWDelta} that 
\begin{equation}
\label{Eq:Majo-dYDelta}
d_\Delta^{\boldsymbol{Y}}\leq
 n(h+2)/2\, ;
\end{equation}

${}^{}$\quad $-$   the ${\boldsymbol{\mathcal Y}}$-cluster variables 
can also be indexed by integers, as it follows from the paragraphs 

${}^{}$\hspace{0.65cm} 
\S\ref{Par:NakanishiDilogarithmicIdentity-Period}
and \S\ref{Par:ExamplesPeriods-I}. Denote by 
$\boldsymbol{i}_\Delta$ the period $\boldsymbol{i}^{half}_{\Delta,A_1}$ of \eqref{Eq:PeriodAn-n-even}, except when 
$\Delta$ is of  type 
 $A_n$

${}^{}$\hspace{0.65cm} with $n$ even, in which case 
$\boldsymbol{i}_{A_n}$ stands for the period $\boldsymbol{i}_{\bullet,\circ}$ defined in \eqref{Eq:period-slice}. 
In any case, $\boldsymbol{i}_\Delta$ is of

${}^{}$\hspace{0.65cm} length $n(h+2)/2$ and its 
$\ell$-cluster coordinates $y_\ell=x_\ell(\boldsymbol{i}_{\Delta})$ for $\ell=1,\ldots,n(h+2)/2$ 
 (defined

${}^{}$\hspace{0.65cm} in  \eqref{Eq=l-th-cluster-variable-period})  form a set of first integrals for ${\boldsymbol{\mathcal Y\mathcal W}}_{\Delta}$.

\subsubsection{}
\label{SS:???}
{\bf $\boldsymbol{F}$-polynomials and the cluster arrangement.}
For the finite type cluster algebra $A_\Delta$, there are nice descriptions of the geometric objects  attached to ${\boldsymbol{\mathcal X\mathcal W}}_{\Delta}$ in \S\ref{Par:ClusterWeb-Polylogarithmic?}.
\mk 

${}^{}$\quad $-$ We denote by $\boldsymbol{F}_\Delta$ the set of $F$-polynomials of $A_\Delta$.  There is a bijection 
$\alpha\mapsto f[\alpha]$ between 

${}^{}$\hspace{0.65cm} the  set of almost-positive roots $\Delta_{\geq -1}$ and the set $\boldsymbol{F}_\Delta$
of $F$-polynomials of $A_\Delta$  which moreover 

${}^{}$\hspace{0.65cm} is such that $f[-\alpha_i]=u_i$ for any simple root $\alpha_i\in \Delta_{>0}$  ({\it cf.}\,\eqref{Eq:Phi-Fpolynomials}).
\mk

${}^{} $ \quad $-$ This gives us a way to enumerate the elements of $\boldsymbol{H}_\Delta=\boldsymbol{H}_{\boldsymbol{\mathcal X\mathcal W}_\Delta}$: its elements are the  

${}^{}$\hspace{0.65cm} 
hyperplanes  $H_i=(u_i=0)$'s for $i=0,\ldots,n$ ($H_0$ being the divisor at infinity) and the 

${}^{}$\hspace{0.65cm} 
irreducible hypersurfaces $H_\alpha=H_{f[\alpha]}$'s for $\alpha\in \Delta_{>0}$. 
\mk 

${}^{} $ \quad $-$ We define $\boldsymbol{Arr}_\Delta$,  the {\bf `cluster arrangement  of type $\boldsymbol{\Delta}$'} as the 
one  in $\mathbf P^n_0$ whose irreducible 

${}^{}$\hspace{0.65cm} components are the elements of  $\boldsymbol{H}_\Delta$: 
 $$\boldsymbol{Arr}_\Delta=\Big(\cup_{i=0}^n H_i \Big) \cup \Big(\cup_{\alpha\in \Delta_{>0}} H_\alpha \Big)\, .
$$
${}^{}$\hspace{0.65cm} It is an arrangement of  $\big(n+1+\lvert \Delta_{>0}\lvert\big)$ irreducible hypersurfaces in $\mathbf P^n_0$.\mk 

${}^{} $ \quad $-$ We denote by $\boldsymbol{U_\Delta}$ the complement of the cluster arrangement of type $\Delta$: 
 \begin{equation}
 \label{Eq:U--Delta}
U_\Delta=  
\mathbf P^n_0 \setminus 
\boldsymbol{Arr}_\Delta\, .
 \end{equation}

${}^{} $ \quad $-$ Assuming that (point {\it 2.}\,of) Conjecture \ref{Conj:Inv-divisor-Fx-Ftildex} holds true for all pairs of $\boldsymbol{\mathcal X}$-variables, 
   we get 
  
   ${}^{}$\hspace{0.8cm}that  the
divisor of common leaves 
$\Sigma^c_\Delta=\Sigma^c(\boldsymbol{\mathcal X\mathcal W}_\Delta)$
 is included in the cluster arrangement: 
 \begin{equation}
 \label{Eq:SigmacDelta-AADelta}
 \Sigma^c_\Delta\subset \boldsymbol{Arr}_\Delta\, .
 \end{equation}
 \mk

\subsubsection{\bf Some basic questions about cluster webs.}
\label{SS:???}
With the notations introduced above, we can now state and comment
on 
 the main questions we are interested about the cluster webs of the form $\boldsymbol{Z\mathcal W}_{\Delta}$ for $\boldsymbol{Z}=
\boldsymbol{\mathcal X}$ or $\boldsymbol{\mathcal Y}$\footnote{The $\boldsymbol{\mathcal A}$-cluster webs 
${\boldsymbol{\mathcal A\mathcal W}}_{\Delta}$ 
are far less interesting 
in what concerns 
 their rank and ARs hence are left aside.}:   
\mk

${}^{}$\quad $-$ $ {\bf   [a_{\boldsymbol{\mathcal X}}].}$  Is there a closed formula for 
the number
$d_\Delta^{\boldsymbol{\mathcal X}}$  foliations appearing effectively in ${\boldsymbol{\mathcal X\mathcal W}}_{\Delta}$? \sk

\mk 


${}^{}$\quad $-$ 
$ {\bf   [b_{\boldsymbol{\mathcal Y}}].}$
 Is the majoration \eqref{Eq:Majo-dYDelta} an equality? 
\mk


${}^{}$\quad $-$ 
$ {\bf   [c_{\boldsymbol{\mathcal X}}].}$
 Are there closed formulae for 
 $\rho^\bullet({\boldsymbol{\mathcal X\mathcal W}}_{\Delta})$, 
 ${\rm IIrk}^\bullet(\mathcal X{\boldsymbol{\mathcal W}}_{\Delta})$ and ${\rm rk}(\mathcal X{\boldsymbol{\mathcal W}}_{\Delta})$
?\mk 
 
${}^{}$\quad $-$ 
$ {\bf   [d_{\boldsymbol{\mathcal X}}].}$
 Is the web ${\boldsymbol{\mathcal X\mathcal W}}_{\Delta}$ AMP
?\mk

${}^{}$\quad $-$ 
$ {\bf   [e_{\boldsymbol{\mathcal X}}].}$
Is there a (combinatorial?) description of the space of iterated integral ARs 
$\boldsymbol{
\mathcal I
\hspace{-0.08cm}
\mathcal I
\hspace{-0.08cm}
{\mathcal A}}({\boldsymbol{\mathcal X\mathcal W}}_{\Delta})$?

${}^{}$\hspace{1.4cm} Same question for the full space $\boldsymbol{\mathcal A}({\boldsymbol{\mathcal X\mathcal W}}_{\Delta})$ of all ARs.\mk

${}^{}$\quad $-$ 
$ {\bf   [f_{\boldsymbol{\mathcal X}}].}$ Is the inclusion \eqref{Eq:SigmacDelta-AADelta} proper?
\mk

${}^{}$\quad $-$ 
$ {\bf   [c_{\boldsymbol{\mathcal Y}}]}$, $ {\bf   [d_{\boldsymbol{\mathcal Y}}]}$, $ {\bf   [e_{\boldsymbol{\mathcal Y}}]}$ 
and $ {\bf   [f_{\boldsymbol{\mathcal Y}}].}$
 Same questions as  the  ones above, but for the web ${\boldsymbol{\mathcal Y\mathcal W}}_{\Delta}$.
\mk

The sequel of this text can be seen mainly as an attempt to answer these questions. Those about the rank and the ARs will be investigated in the next sections. In the subsections below, we 
answer to $ {\bf   [a_{\boldsymbol{\mathcal X}}]}$ and $ {\bf   [b_{\boldsymbol{\mathcal Y}}]}$. Note that if these two questions sound quite basic, it is necessary to refer to some non-trivial previous results in order to provide them with some answers.

 \subsection{\bf Differential independence of cluster variables and consequences}
\label{SS:DifferentialIndependenceClusterVariables}
In this section, we prove the first point of Conjecture \ref{Conj:Inv-divisor-Fx-Ftildex} for the $\boldsymbol{\mathcal X}$-cluter variables in finite type and deduce from this a general formula for the 
degree of any  cluster web $\boldsymbol{\mathcal X\mathcal W}_{\Delta}$.  All this follows rather easily from some results by Shermann-Bennet in   \cite{Sherman-Bennett2} (see also \cite{Sherman-Bennett1}) where she enumerates the number of cluster variables in terms of some combinatorial objects associated to the finite type cluster algebra considered. \mk 

We get that 
the first point of Conjecture 
\ref{Conj:Inv-divisor-Fx-Ftildex} always holds true for   a finite type cluster algebra: 
%
%
\begin{prop} 
\label{Prop:propro}
For any $\boldsymbol{\mathcal X}$-cluster 
variables $x$ and $\tilde x$ of a cluster algebra of finite type: 
$x $ and $\tilde x$ define the same foliation if and only if $\tilde x=x$ or $\tilde x=x^{-1}$ (as rational functions).
\end{prop}

From this, we will answer to questions ${\bf   [a_{\boldsymbol{\mathcal X}}]}$ and ${\bf   [b_{\boldsymbol{\mathcal Y}}]}$ by providing 
explicit closed formulas for the degrees of the ${\boldsymbol{\mathcal X}}$- and ${\boldsymbol{\mathcal Y}}$-cluster webs of Dynkin type (see Corollary \ref{Coro:d-X-Delta--d-Y-Delta}).
\begin{center}
$\star$
\end{center}


The five exceptional cases (namely $E_6,E_7,E_8,F_4$ and $G_2$) are settled via brute force computations, by generating  all the clusters using a computer algebra system (Maple\textsuperscript{\textregistered}  as it happens) and verifying that the statement of Proposition \ref{Prop:propro} holds true for each pair of cluster variables. \mk

%


Our approach for proving Proposition 
\ref{Prop:propro} when $\Delta$ is a  Dynkin diagram of type $A,B,C$ or $D$ consists just in pushing a bit forward the one used in \cite{Sherman-Bennett1} to enumerate the cluster variables of the corresponding cluster algebra.  In each case, there is a (classical) combinatorial description of the ($\boldsymbol{\mathcal A}$- or) $\boldsymbol{\mathcal X}$-clusters and associated cluster variables in terms of some triangulations of a (possibly punctured) polygon $\boldsymbol{P}_{\Delta}$.  From this, one can explicitly compute the differential $dx$ of a given cluster variable $x$ (in some suitable coordinates) and verify that, given another cluster variable $\tilde x$, one necessarily has $\tilde x\in \{\, x\, ,\, x^{-1}\, \}$ if the rational 2-form $dx\wedge d\tilde x$  vanishes identically.   
However, a problem with this approach is that we are not aware of a uniform description of $\boldsymbol{P}_{\Delta}$ and of the associated  cluster  combinatorics  which cover all the types.  
Indeed, the combinatorics associated to  $\boldsymbol{P}_{\Delta}$ 
is specific to each  type $A,B,C$ or $D$  considered hence one has to deal with each type independently of the others, even if all these are handled along the same lines.  
\sk

Below, we consider the cases of type $A$ and $D$ but leave aside types $B$ and $C$. These can be treated in similar ways (and this is left to the reader). 

As for some general references (some more specific are given below), we mention 
\cite{FominShapiroThurston} (for type $A$ and $D$) and  
the fourth and fifth chapters of the book \cite{FWZ} (for type $B$ and $C$).  
All the material  we need can be found in it which henceforth is a handy reference regarding our purpose.

%
%
%
%
%

%
%
%
%
%
%
%
%
\subsubsection{Triangulations, clusters and cluster variables in type $\boldsymbol{A}$.}
\label{SSub:Type-An}
Since this case  will be considered again later, we give many details about  it, even if it is very well understood.
\sk 

For some additional references specific to this case, we mention among many others:  the PhD dissertation \cite{King},  section \S1.3 of \cite{FockGoncharovENS} and  the appendix of \cite{FockGoncharovAmalgation}.
\begin{center}
$\star$
\end{center}

For $n\geq 2$, let $\boldsymbol{P}_{\hspace{-0.05cm} n+3}=\boldsymbol{P}_{\hspace{-0.05cm}A_n}$ by the standard $(n+3)$-gon, with vertices labeled cyclically by elements of $\mathbf Z/(n+3)\mathbf Z$ ({\it i.e.}\,$n+3+i=i$ for any $i$).  We consider triangulations of $\boldsymbol{P}_{\hspace{-0.05cm} n+3}$ 
whose set of vertices is the one of the latter polygon.  Let $\mathcal T_{n+3}$ stand for the set of all such triangulations, each of them being identified with 
 the set of its interior edges (which is of cardinality $n$). 
   \sk 
 
 For four pairwise distinct points $z_1,\ldots,z_4\in \mathbf C$, one sets
 \begin{equation}
 \label{Eq:cluster-r}
  r\big(z_1,z_2,z_3,z_4\big)=  \frac{(z_1-z_2)(z_3-z_4)}{(z_1-z_4)(z_2-z_3)}\in  \mathbf C\setminus\{0,-1\}
 \end{equation}
 and one verifies that this expression  enjoys the following invariance properties: 
 \begin{equation}
 \label{Eq:r-invariance-property}
 r(z_1,z_2,z_3,z_4)= r(z_3,z_4,z_1,z_2) 
 \qquad \mbox{ and }\qquad 
 r(z_2,z_3,z_4,z_1)= r(z_1,z_2,z_3,z_4)^{-1}\, .
  \end{equation}
 This map  induces a  rational map  
\begin{align*}
  r\, : \quad
   \mathcal M_{0,4}& \longrightarrow \mathbf C\setminus\{0,-1\} \\ 
 [p_i]_{i=1}^4& \longmapsto  r(p_1,p_2,p_3,p_4)\, , 
\end{align*} 
 which is  a cross-ratio normalized such that $  r(\infty,-1,0,z)=z$ for any $z\in \mathbf C\setminus \{0,-1\}$.  \mk
 
 Now let $(i,j,k,l)$ be a quadruple of elements of $\mathbf Z/(n+3)\mathbf Z$. 
 Then post-composing by $r$ the forgetful map $ \mathcal M_{0,n+3}\longrightarrow  \mathcal M_{0,4} :  [p_s]_{s=1}^{n+3)\mathbf Z}\longmapsto 
 [p_i,p_j,p_k,p_l]$  (with the convention that the $p_s$'s are labeled modulo $n+3$), one gets a surjective  morphism 
 \begin{align}
 \label{Eq:r-ijkl}
  r_{i,j,k,l}\, : \quad
   \mathcal M_{0,n+3}& \longrightarrow \mathbf C\setminus\{0,-1\} \\ 
 [p_i]_{i=1}^{n+3}& \longmapsto  r(p_i,p_j,p_k,p_l)\, .  \nonumber
\end{align}

\begin{figure}[!h]
\begin{center}
\includegraphics[scale=0.6]{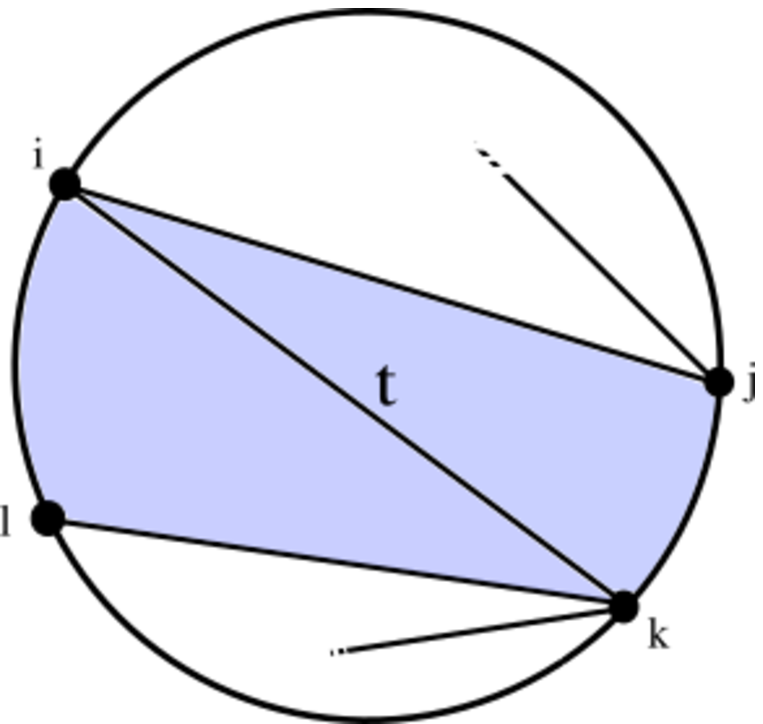}
 \hspace{2.5cm}
 \includegraphics[scale=0.6]{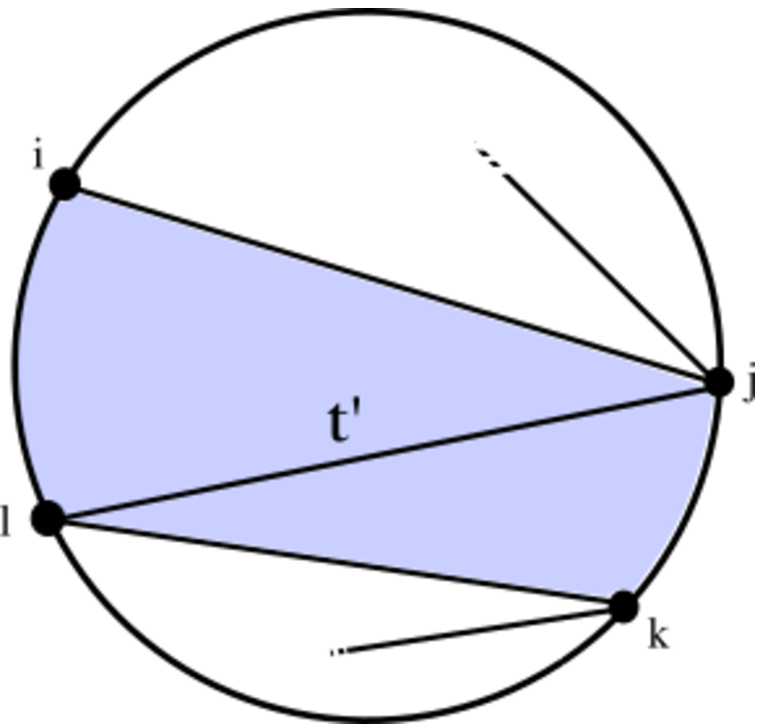}
\caption{on the left, the quadrilateral $q^T_t$ associated to the edge $t$ of $T$ (in blue) and its four vertices, noted here i,j,k and l. On the right, the flip $T'$ of $T$ at $t$ (with $q^{T'}_{t'}=q^T_t$).}
\label{Fig:ZigZag}
\end{center}
\end{figure}

 Let $T\in \mathcal T_{n+3}$ be a triangulation of $\boldsymbol{P}_{\hspace{-0.05cm} n+3}$. Any edge $t\in T$ is a diagonal of a unique and well-defined quadrilateral $q^T_{t}$ with sides in $\overline{T}$. 
Denote by $(i_{t}^T,j_{t}^T,k_{t}^T,l_{t}^T)\in \{1,\ldots,n+3\}^{4}$ the 4-tuple whose components are the vertices of $q^T_{t}$, labeled as they appear cyclically on the boundary $\partial  \boldsymbol{P}_{\hspace{-0.05cm} n+3}$, with the supplementary assumption that $i_t^T$ (hence $k^T_t$ as well) be an extremity of $t$, see Figure \ref{Fig:ZigZag}.  

 Then one defines a rational morphism by setting
\begin{align}
x_{t}^T =r_{i_{t}^T,j_{t}^T,k_{t}^T,l_{t}^T}\, : \quad  \mathcal M_{0,n+3} & \longrightarrow \mathbf C\setminus\{0,-1\}\\
\big[c(s)\big]_{s=1}^{n+3}& \longmapsto  r\Big(
c\big({i_{t}^T}\big),c\big({j_{t}^T}\big),c\big({k_{t}^T}\big),c\big({l_{t}^T}\big)
\Big) \, .
\nonumber
\end{align}
Note that although the 4-tuple $(i_{t}^T,j_{t}^T,k_{t}^T,l_{t}^T)$ is actually only defined  by $q^T_t$ up to a cyclic shift of length 2, the map $x_{t}^T$ only depends on the pair $(T,t)$ (thanks to the first invariance property  in \eqref{Eq:r-invariance-property}). 
Considering all the maps  $x^T_t$ for all the edges $t\in T$, we 
get a well-defined morphism
\begin{align}
\label{Eq:xT}
x^T : \mathcal M_{0,n+3}   & \longrightarrow
 \mathbf C^{T} \\
c  & \longmapsto  \big( x^T_t(c)
\big)_{t\in T} \, , 
\nonumber
\end{align}
which can be proved to be birational. \footnote{This is well-known. See \S\ref{SubSec:ZigZag-Map} further where we explicit $x^t$ in a particular case.} \mk

Given $t$ a edge of $T$, one defines the flip of $T$ at $t$ as the new triangulation $T'=T'_t\in \mathcal T_{n+3}$ defined as follows: 
$T'$ shares with $T$ all its edges except $t$  which is replaced by the other diagonal  of $q^T_t$, denoted by $t'$. Of course, $T$ can be recovered from $T'$ by flipping $t'$. Then one defines a graph 
 $\boldsymbol{\Gamma}_{n+3}$ as follows: it has 
 ${\mathcal T}_{n+3}$ as 
 set of vertices and there is an edge between two triangulations $T_1,T_2$  if and only if one is obtained from the other by a flip.  It is well known that $\boldsymbol{\Gamma}_{n+3}$ is connected.  

\begin{figure}[!h]
\begin{center}
\psfrag{n+3}[][][1]{$n+3 $}
\includegraphics[scale=0.5]{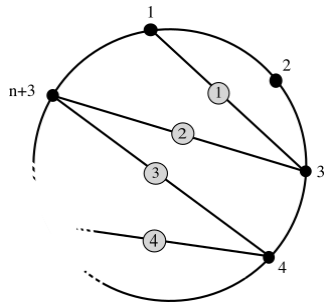}
\caption{the Zig-zag triangulation $T_0$ of the $(n+3)$-gon $ \boldsymbol{P}_{\hspace{-0.05cm} n+3}$.}
\label{Fig:ZigZag-n+3}
\end{center}
\end{figure}

Let $T_0$ be the {\bf `zig-zag triangulation} of $\boldsymbol{P}_{\hspace{-0.05cm} n+3}$ defined as the one whose interior edges are the segments $[1,3]$, $[3,n+3]$, $[n+3,4]$, $\ldots$. We label these edges according to their order of appearance when going along the interior zig-zag starting from the first vertex of $\boldsymbol{P}_{\hspace{-0.05cm} n+3}$: $[1,3]$ is labeled by 1, $[3,n+3]$ by 2, etc (see Figure \ref{Fig:ZigZag-n+3})
\mk

Since $\boldsymbol{\Gamma}_{n+3}$ is connected, this gives a labeling by $\{1,\ldots,n\}$ of the interior edges of any $T\in 
{\mathcal T}_{n+3}$, from which we deduce a well-defined 
 identification $\mathbf C^{T}\simeq \mathbf C^n$ for any triangulation. 
 Accordingly, we denote by $x^T_1,\ldots,x^T_n$ the components of the map $x^T$ in \eqref{Eq:xT}.

 For two (interior) edges labeled by $i,j\in \{1,\ldots,n\}$, 
 one defines an integer $b_{ij}$ as follows: it is zero if these two edges are not adjacent; otherwise, both are sides of a same 2-face 
 $\tau_{ij}$ of $T$ (a triangle) and one sets $b_{ij}=1$ (resp.\,$b_{ij}=-1$) if $j$ follows (resp.\,precedes) $i$ when circulating 
 clockwise 
 along the boundary of  $\tau_{ij}$. Then $B^T=(b_{ij})_{i,j=1}^n$ is a skew-symmetric matrix, the {\bf exchange matrix} associated to $T$. 
\begin{figure}[!h]
\begin{center}
\includegraphics[scale=0.6]{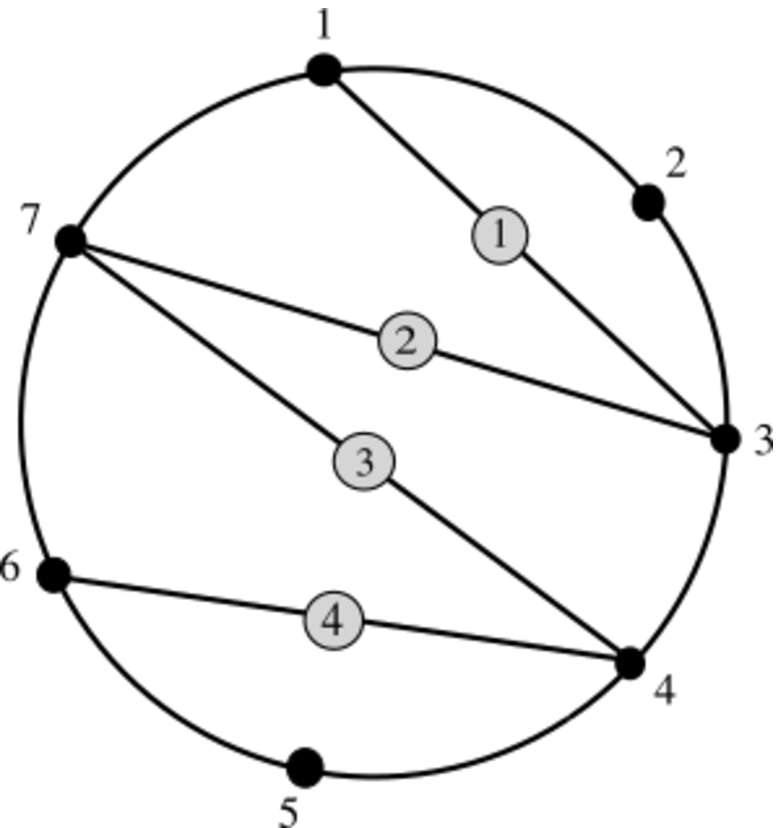} \qquad \qquad 
\begin{tabular}{c} \vspace{-5cm}\\
$B^{T_0}= \begin{bmatrix}
 0 & 1 & 0 & 0\\
 -1 & 0 & -1 & 0 \\
 0 & 1 & 0 & 1\\
 0 &0&-1 & 0
 \end{bmatrix}\, .$
 \end{tabular}
\caption{ the zig-zag triangulation and the associated exchange matrix in case $n=4$.}
\end{center}
\end{figure}

The key result  in order to relate the material considered above (of combinatorial and/or geometric nature) to the world of cluster algebras is the following one (which is well-known): 
\begin{lem}
\label{Lem:Flip-Mut-An}
Assume that $T,T'\in \mathcal T_{n+3}$ are related by the flip with respect to their $k$-th edge, for some $k$.  Then,     as rational functions of the $x^T_i$'s for $i=1,\ldots,n$, the $x^{T'}_j$'s are given by means of the $\boldsymbol{\mathcal X}$-cluster mutation formulas  \eqref{Eq:A-X-Mutation-formulae}, relatively to the exchange matrix $B^T$.
\end{lem}
%

Then, since $B^{T_0}$ coincides with the exchange matrix associated to the bipartite quiver $Q_{\vec{A_n}}$ for any $n\geq 2$ (as it is easy to verify), one obtains the following fundamental result: 
\begin{thm}
\label{Thm:Type-An}
The map associating $(\boldsymbol{x}^T,B^T)$
to $T\in \mathcal T_{n+3}$  
 is a $\boldsymbol{\mathcal X}$-cluster pattern of type $A_n$. \sk 
\end{thm}

Actually, many interesting things come along with this result since it 
 provides some geometric and combinatorial interpretations of the objects attached to the cluster algebra $A_{A_n}$:
\begin{itemize}
 \item  the set ${\mathcal T}_{n+3}$ of triangulations of $\boldsymbol{P}_{n+3}$ identifies with the set of clusters of $A_{A_n}$;
 \item  the `triangulation graph  $\boldsymbol{\mathcal T}_{n+3}$' corresponds to the exchange graph $\boldsymbol{\Gamma}_{A_n}$;
 \item   the cross-ratios $x^T_t$ for $T\in {\mathcal T}_{n+3}$ and $t\in T$ identify with the $\boldsymbol{\mathcal X}$-cluster variables;
 \item  any  map \eqref{Eq:xT} induces a birational identification  
 $\mathcal M_{0,n+3}\simeq  \boldsymbol{\mathcal X}_{A_n}$;
 \end{itemize}

Regarding the question we are interested in here, we get the 
%
\begin{cor}
\label{Cor:X-var-web-An}
\begin{enumerate}
\item Up to pull-back on $\mathcal M_{0,n+3}$ under $x^{T_0}$, 
the $\boldsymbol{\mathcal X}$-cluster variables of type $A_n$ are the 
cross-ratios $r_{i,j,k,l}$ for all cyclically oriented $4$-tuples $(i,j,k,j)\in (\mathbf Z/(n+3)\mathbf Z)^4$. Thus, up to inversion ({\it i.e.}\,modulo $x\leftrightarrow x^{-1}$),  there are precisely ${n+3 \choose 4}$ such cluster variables.
\sk
\item Two $\boldsymbol{\mathcal X}$-cluster variables define the same foliation if and only if  both coincide (possibly up to inversion). 
\sk
\item The  pull-back of 
$\boldsymbol{\mathcal X\mathcal W}_{A_n}$   under $x^{T_0}$ coincides with the web 
$\boldsymbol{\mathcal W}_{\hspace{-0.05cm}\mathcal M_{0,n+3}}$
 considered above in  \S\ref{SubPar:Webs-spaces-configurations}. Hence $\boldsymbol{\mathcal X\mathcal W}_{A_n}$ is AMP with only polylogarithmic ARs (of weight 1 or 2). 
\end{enumerate}
\end{cor}
\begin{proof}
The first point follows from all what has been said just before. 
When combined with Lemma \ref{Lem:pi-I-J--pi-I'-J'}, it gives {\it 2}. Finally, 
 the third point comes from the first combined together with  the results of \cite{Pereira}
mentioned in \S\ref{Par:W-M0n+3-all-AMP} above.
\end{proof}


\begin{rem} 
\label{Rem:X-cluster-var-properties-yoyo}
From the first point of the previous corollary, one can deduce easily that  relatively to the rational coordinate system 
$ z=(z_1,\ldots,z_n)\dashrightarrow  [\infty,0,1,z_1,\ldots,z_n]$ 
on $\mathcal M_{0,n+3}$ ({\it cf.}\,\S\ref{SS:map-Un} further), 
the spectrum of the pull-back of any  $\boldsymbol{\mathcal X}$-cluster variable on $\mathcal M_{0,n+3}$ under $x^{T_0}$     is exactly 
$\{0,-1,\infty \}$.   Moreover, given two such cluster variables, one verifies that the statement corresponding to the second point of Conjecture \ref{Conj:Inv-divisor-Fx-Ftildex} holds true for their pull-back under $x^{T_0}$. 
%
One should be able to deduce from these  facts that the 
conjectures \ref{Conj:Spectrum-Cluster-Variables} and \ref{Conj:Inv-divisor-Fx-Ftildex} hold true for the $\boldsymbol{\mathcal X}$-cluster variables of Dynkin type $A_n$. 
\end{rem}

%
%


\subsubsection{Triangulations, clusters and cluster variables in type $\boldsymbol{D}$.}
Everything in this case is very similar to the preceding $A$ case, modulo slight differences. The material below is  taken from subsections \S3.1 and \S5.4 of \cite{Sherman-Bennett1}. We refer to this paper  for more details, pictures  and additional references. 
\mk

%
%

We assume $n\geq 4$.  Let $\boldsymbol{P}^\circ_n$
be the  (topological) $n$-gon with a puncture $\circ$ in its interior.  
One considers arcs of  two different types: arcs between two exterior vertices of 
$\boldsymbol{P}^\circ_n$ or radii joining such a vertex with the puncture $p$\footnote{Actually, by 'arc' one has to understand 'isotopy class with fixed extremities of arc' in $\boldsymbol{P}^\circ _n$.}. These arcs are labeled in two ways: `notched' or `plain'.   We also need to consider a {\it 'cut'} which is the radius  between the middle of the boundary arc between $1$ and $n$ and the puncture $\circ$.  We will use it  to specify some rational functions associated to arcs in $\boldsymbol{P}^\circ_n$ below.
\sk 

One considers the set ${\mathcal T}^\circ_n$ of  `tagged triangulations'  of $\boldsymbol{P}^\circ_n$, these 
being defined as 
 maximal collections of pairwise compatible tagged arcs in $\boldsymbol{P}^\circ_n$ (for a certain compatibility about which we will not elaborate here).   A triangulation $T\in {\mathcal T}^\circ_n$ is identified with the set of its edges, which is of cardinality $n$. 
 Each  $t\in T$ 
   can be flipped, which allows to construct a graph with ${\mathcal T}^\circ_n$ as set of vertices, denoted by $\boldsymbol{\Gamma}^\circ_{\hspace{-0.05cm}n}$. 
It can be proved that it is connected, hence any $T\in {\mathcal T}^\circ_n$ can be obtained by a finite sequence of flip from any other one, for instance from 
the  triangulation $T^0$ consisting only of radii,   pictured on the left in  Figure \ref{Fig:T0-Radii}  below. \mk

 \begin{figure}[!h]
\begin{center}
\includegraphics[scale=0.75]{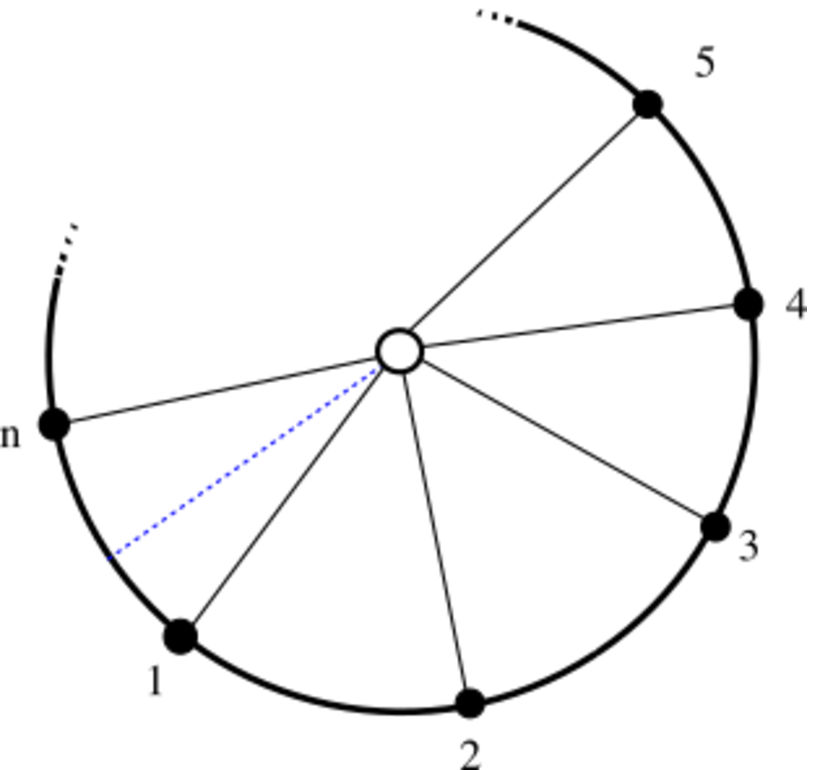}
 \hspace{1.5cm}
 \includegraphics[scale=0.7]{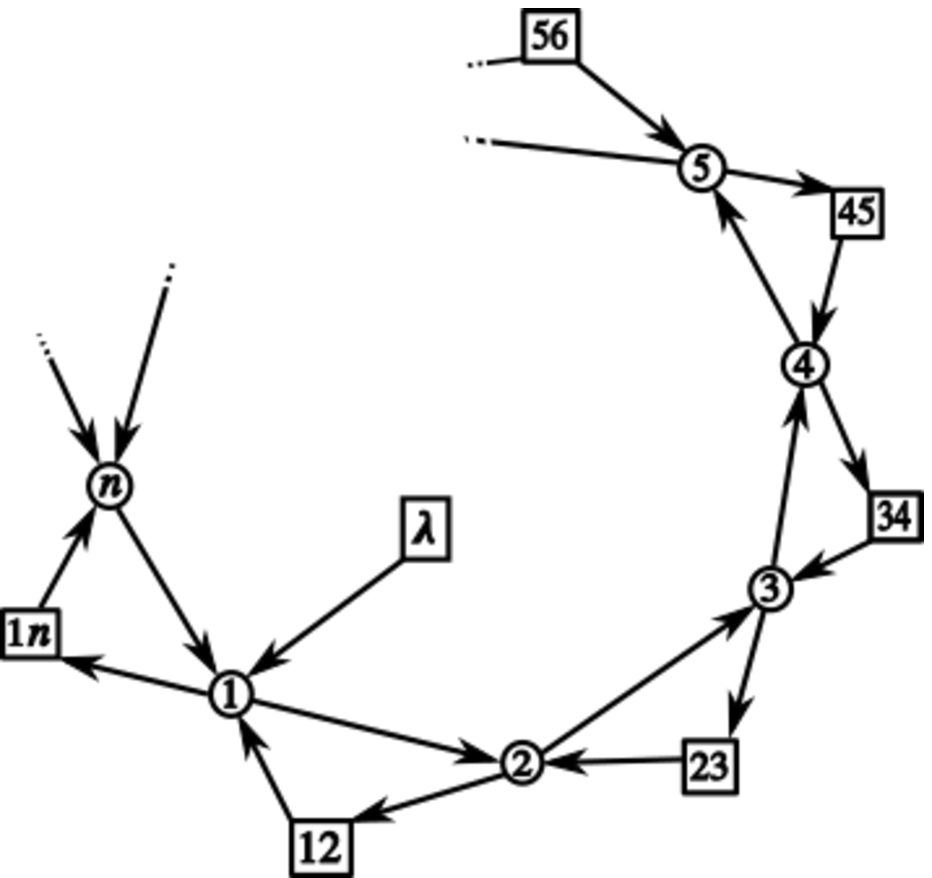}
\caption{On the left, the tagged triangulation $T_0$ of $\boldsymbol{P}^\circ_n$ whose all edges are radii. The cut is the dotted line in blue. On the right, the associated ice quiver $\overline{Q}^{T_0}$ 
 (its frozen vertices are boxed, the non frozen ones are circled). 
The $i$-th vertex is associated to the radius from $i$ to the puncture $\circ$ and the boxed(=frozen) one labeled by $i(i+1)$ corresponds to the arc between $i$ and $i+1$ on the boundary $\partial \boldsymbol{P}^\circ_{n}$ (with $1n$ standing for $n(n+1)$). The unfrozen part of $\overline{Q}^{T_0}$ is the  anticlockwise oriented $n$-cycle.}
\label{Fig:T0-Radii}
\end{center}
\end{figure}

 Our goal is to  construct rational functions $x_t^T$ for all $(T,t)$ with $T\in  
 {\mathcal T}^\bullet_n$ and $t\in T$, which will turn out to be the 
 $\boldsymbol{\mathcal X}$-cluster variables of a cluster algebra of finite type $D_n$. In order to do so,  we first define an extension $\overline{T}$ of  $T$ by adjoining to it the $n$ circular arcs of $\partial \boldsymbol{P}^\circ_n$.   Using the same rule as the one used to construct 
 an iced quiver from a triangulation of $\boldsymbol{P}_{\hspace{-0.05cm}n+3}$ in type $A$, one associates an iced  quiver to
 $\overline{T}$: its non frozen vertices correspond to the 
the edges of $T$ (interior edges)
and there is one frozen vertex for each boundary circular 
 arc   joining $i$ to $i+1$ for $i=1,\ldots,n$ (with $n+1$ identified with $1$). 
 However in type $D$, the construction of the cluster variables requires to  introduce an additional frozen vertex, denoted by $\lambda$, to get a quiver $\overline{Q}^T$ with $2n+1$ vertices, $n+1$ of which are frozen. 
In the case of $T_0$, there is a unique arrow 
adjacent to $\lambda$ in 
$\overline{Q}^{T_0}$, which goes from $\lambda$ to $1$ 
(see the picture  on the right in Figure  \ref{Fig:T0-Radii}). Since the graph $\boldsymbol{\Gamma}^\circ_{\hspace{-0.05cm}n}$ is connected, this defines unambiguously the arrow(s) adjacent to $\lambda$ in $\overline{Q}^{T}$
for any tagged  triangulation $T\in \mathcal T^\circ_{\hspace{-0.05cm}n}$. 
 Finally, we denote by $Q^T$ the unfrozen part of $\overline{Q}^T$. 
\mk

Let $\boldsymbol{\mathcal K}$ be the field of rational functions in $n+1$ indeterminates $v_1,\ldots,v_n$ and $\lambda$.  Each $v_i$ has to be thought of as associated to the $i$-th vertex of $\boldsymbol{P}^\circ_{n}$ and $\lambda$ is an extra indeterminate which is needed (about the reason why this supplementary variable is needed, see the discussion in the paragraph just after \eqref{Eq:xT-Dn} below). 
For $1\leq i<j\leq n$, one defines some elements of $\boldsymbol{\mathcal K}$ by setting 
\begin{equation}
\label{Eq:Pij-and-co}
P_{i,j}= (v_i-v_j) \, , \qquad P_{\hspace{-0.05cm}i,\overline{j}}=v_j(2-v_i)-v_i
\, , \qquad
P_{i,\circ}=v_i-1
\quad \mbox{ and } \quad 
P_{\hspace{-0.05cm} i,\circ^{\bowtie}}=v_i+1
\end{equation}
 and for any arc $t$ in $\overline{\boldsymbol{P}}^\circ_{n}$, one sets
\begin{equation}
\label{Eq:Pt}
P_t=\begin{cases}\hspace{0.15cm} 
 P_{i,j}\hspace{0.5cm} \mbox{if } t \mbox{ has endpoints } i,j\, \mbox{ with } i<j \mbox{ and does not cross the cut;} \\
 \hspace{0.15cm} P_{\hspace{-0.05cm}i,\overline{j}} \hspace{0.55cm} \mbox{if } t \mbox{ has endpoints } i,j\, \mbox{ with } i<j \mbox{ and does cross the cut;} \\
\hspace{0.15cm}  P_{i,\circ} \hspace{0.5cm} \mbox{if } t \mbox{ is the plain radius from } i \mbox{ to the puncture } \circ; \\
\hspace{0.15cm}  P_{\hspace{-0.05cm} i,\circ^{\bowtie}} \hspace{0.4cm} \mbox{if } t \mbox{ is the notched radius from } i   \mbox{ to the puncture } \circ. \end{cases}
\end{equation}

Setting in addition $P_{\hspace{-0.05cm} \lambda}=\lambda$,  the previous definitions allow to define a 
$\boldsymbol{\mathcal K}$-valued map on the set of vertices of  $\overline{Q}^T$, this for any tagged triangulation $T$ of $\boldsymbol{P}^\circ_{n}$. Hence for  
$T\in \mathcal T^\circ_{\hspace{-0.05cm}n}$ and any edge $t\in T$, 
 one defines an element of $\boldsymbol{\mathcal K}$ by setting 
\begin{equation}
x^T_t =\frac{ \prod_{ \tau   \rightarrow t}  P_\tau}{ \prod_{t \rightarrow \tau }  P_\tau}
\end{equation}
where the numerator (resp.\,denominator) is the product of the $P_\tau$'s for all vertices $\tau$ such that there is an arrow from $\tau$ to $t$ (resp.\,from $t$ to $\tau$) in $\overline{Q}^T$.


\begin{exm}
\label{Ex:Cluster-Dn-T0}
For the triangulation $T_0$ of Figure \ref{Fig:T0-Radii},  if  $x_i$ stands for the rational expression associated to the radius between $i$ and $\circ$, one gets easily that 
$$ x_1=  \frac{\lambda \,w_{1,2} \,w_n}{
w_2\, w_{1\overline{n}} }\, , \qquad 
x_i= \frac{ w_{i-1} \, w_{i,{i+1}}}{ w_{i+1} \, w_{i-1,i}} \quad \mbox{ for } i=2,\ldots,n-1
\qquad \mbox{and}\qquad 
x_n= \frac{w_{n-1} w_{1\overline{n}}  }{w_1 w_{1\overline{n}} }
$$
where we have set $w_i=P_{i,\circ}=v_i-1$ and  $w_{i,j}=P_{i,j}=w_i-w_j$  
for $i,j=1,\ldots,n$ with $i<j$ 
and $w_{1\overline{n}}=w_n-w_1(2+w_n)$. 
It is then easy to verify that the rational map $(w_1,\ldots,w_n,s)\mapsto (x_i)_{i=1}^n$ is dominant. This implies that the functions 
$x_1,\ldots,x_n$ are algebraically independent  in 
$\boldsymbol{\mathcal K}$
hence generate a subfield isomorphic to the field of rational functions in $n$ indeterminates. 
\end{exm}


It follows from general results  about cluster algebra associated to (tagged) triangulations of bordered surfaces ({\it cf.}\,\cite{FominShapiroThurston})  that the  
 statement of Lemma \ref{Lem:Flip-Mut-An}  applies to any 
 $T,T'\in \mathcal T^\circ_{n}$ related by a flip with respect to an edge $\tau \in T$: as rational functions of the $x^T_t$'s (for $t\in T$), the $x^{T'}_{t'}$'s (with $t'\in T'$) are given by means of the $\boldsymbol{\mathcal X}$-cluster mutation formulas  \eqref{Eq:A-X-Mutation-formulae} (in the direction corresponding to $\tau$), relatively to the quiver $Q^T$.     
 It then follows from Example \ref{Ex:Cluster-Dn-T0} just above that 
for any $T\in \mathcal T^\circ_{n}$, $\{x^T_t\, \lvert \, t\in T\, \}$ is a set of $n$  algebraically independent elements of $\boldsymbol{\mathcal K}$. Finally, 
 one verifies that  the  quiver $\mu_{n-2}\circ \ldots \mu_2\circ \mu_1(Q^{T_0})$  is of Dynkin type $D_n$.\footnote{We recall that  $Q^{T_0}$ is the unfrozen component of the quiver 
 $\overline{Q}^{T_0}$
  in Figure \ref{Fig:T0-Radii}: it is the $n$-cycle quiver.} 
 \sk 
 
  From all the preceding observations, we deduce the  following analog for $D_n$, of Theorem \ref{Thm:Type-An}:
\begin{thm}
\label{Thm:Type-Dn}
The map associating $(\boldsymbol{x}^T,B^T)$
to $T\in \mathcal T_{\hspace{-0.05cm}n}^\circ$  
 is a $\boldsymbol{\mathcal X}$-cluster pattern of type $D_n$. \sk 
\end{thm}

As in the $A_n$-case, many interesting things come along with this theorem  since it 
 provides some geometric and combinatorial interpretations of the objects attached to the cluster algebra $A_{D_n}$: for instance, 
  the triangulation graph  $\boldsymbol{\Gamma}\mathcal T_{n}^\circ $ identifies with the exchange graph $\boldsymbol{\Gamma}_{D_n}$, the 
set $\boldsymbol{\mathcal Xvar}_{D_n}$ of $\boldsymbol{\mathcal X}$-cluster variables in  type $D_n$ is in bijection with   $\{\, x_t^T\,  \lvert \, T\in 
 \mathcal T_{n}^\circ ,\,  t\in T\}$, etc. 
 \mk 
 
 The whole picture in type $D_n$ is analogous to the one in type $A_n$ but  there is a difference as well, which requires to be discussed.  Let $\mathcal M_{D_n}$ stand for the 
 affine space ${\rm Spec}(\mathbf C[v_1,\ldots,v_n,\lambda]\big)$ whose function field is $\boldsymbol{\mathcal K}$. 
 As they have been defined, the $x_t^T$'s are rational functions 
 on $\mathcal M_{D_n}$ which are the components of dominant rational maps 
\begin{equation}
\label{Eq:xT-Dn}
x^T= \big( x^T_t
\big)_{t\in T} : \mathcal M_{D_n}    \dashrightarrow
 \mathbf C^{T} \, , 
\end{equation}
 one for each $T\in \mathcal T_{n}^\circ$, which is similar to 
 \eqref{Eq:xT-Dn}. 
  But a notable difference with the $A_n$-case
 is that the source $\mathcal M_{D_n} $ has dimension $n+1$, which is
 one more than the dimension of the `cluster affine spaces $\mathbf C^T$'. The reason for that lies of course in the use of the additional 
  indeterminate $\lambda$  to construct the cluster variables $x_t^T$'s, and not only the $v_i$'s, as it would have be more natural (but naive) to do considering the $A_n$-case.  But the use of such an extra variable with no natural geometric interpretation actually seems  unavoidable:  one verifies that any specialization of $\lambda$ in a map \eqref{Eq:xT-Dn} gives rise to a map which is not dominant any more.\footnote{This is well illustrated in Example \ref{Ex:Cluster-Dn-T0} : with the notation used there, one verifies easily that $\prod_{i=1}^n x_i=\lambda$.} 
  \mk
  
  Despite the preceding fact, in order to compare the foliations attached to cluster variables which is the problem we are interested in, it will be sufficient to specialize $\lambda$ to 1, {\it i.e.}\,to  work with the restrictions, denoted by $\tilde x_t^T$, of the 
   $x^T_t$'s to the affine hyperplane $\widetilde{ \mathcal M}_{D_n}=\{\lambda = 1\}   \subset \mathcal M_{D_n}$.\sk 
   
   In \cite{Sherman-Bennett1}, the author gives explicit expressions for the 
  $\boldsymbol{\mathcal X}$-cluster variables. To  
  any  $x_t^T$ is associated a quadrilateral $q^T_t$  whose edges are in $\overline{T}$ and admitting $t$ as one of its diagonals. By definition,   
  the {\bf type    $\boldsymbol{\tau(x^T_t)}$} of  $x^T_t$ is the set of vertices of $q^T_t$ lying on the boundary $\partial \boldsymbol{P}^\circ_n$.   \sk 
  
The different combinatorial types of quadrilateral $q_t^T$ are described and pictured p.\,286 in \cite{Sherman-Bennett1}: there are four possibilities, one for each type of cardinality 2 or 4, but two when the cardinality is 3. 
The vertices being naturally ordered, the types can be written as a $k$-tuples of strictly increasing numbers in $\{1,\ldots,n\}$ for $k=2,3$ or $4$. Accordingly, the four aforementioned combinatorial types for $q_t^T$ are denoted by $(i,j)$, $(i,j,k)$ and $(i,j,k)'$ (for $k=3$) and  $(i,j,k,l)$ where $i,j,k,l$ are elements of $\{1,\ldots,n\}$. For each type, it is not difficult to 
 give a finite  list of possible explicit formulas  for the cluster variable $x_{t}^T$ : up to inversion and multiplication by $\lambda$ or its inverse, $x_t^T$ has the following expression according to its type ({\it cf.}\,formulas (5.7) to (5.10) in  \cite{Sherman-Bennett1}):  
\begin{align}
\label{Align:tildex-t-T}
\nonumber
{\rm Type }\, (i,j)\, :&   \qquad 
\frac{P_{i,j}}{P_{i,\overline{j} }}
\vspace{0.3cm} \\
{\rm Type }\, (i,j,k) :&   \qquad 
\frac{P_{i,j}\,P_{k,\circ }}{P_{i,\circ}P_{j,k}}\, , \quad 
\frac{P_{ {\hspace{-0.05cm}i,\overline{k}  }}   \, P_{j,\circ}}
{P_{j,k }\,P_{\hspace{-0.05cm} i,\circ }}\, , \quad 
\frac{P_{ {\hspace{-0.05cm}i,\overline{k}  }}   \, P_{j,\circ}}
{P_{i,j }\,P_{\hspace{-0.05cm} k,\circ }}\, ; 
\mk\\
\nonumber
{\rm Type }\, (i,j,k)' :&   \qquad 
\frac{P_{i,j}\,P_{k,\circ }\,P_{\hspace{-0.05cm} k,\circ^{\bowtie}}}{
P_{ {\hspace{-0.05cm}i,\overline{k}  }}
\, P_{\hspace{-0.05cm} j,k}\,}\, , \quad 
\frac{P_{\hspace{-0.05cm}i,\overline{k} }\,P_{j,\circ }\,P_{\hspace{-0.05cm} j,\circ^{\bowtie}}}
{P_{ i,j  } \, P_{\hspace{-0.05cm} j,k}}\, , \quad
\frac{P_{j,k}\,P_{i,\circ }\,P_{\hspace{-0.05cm} i,\circ^{\bowtie}}}{
P_{ {\hspace{-0.05cm}i,\overline{k}  }}
\, P_{\hspace{-0.05cm} i,j}\,}\, ;
\mk\\
\nonumber
{\rm Type }\, (i,j,k,l)\, :&   
\qquad 
\frac{P_{i,j}\,P_{k,l}}{P_{i,l }P_{\hspace{-0.05cm} j,k }}\, , \quad 
\frac{P_{\hspace{-0.05cm} i,j}\, P_{\hspace{-0.05cm}kl}}{P_{\hspace{-0.05cm}i,\overline{l} }\, P_{\hspace{-0.05cm} j,\overline{k} }}
\, , \quad 
\frac{P_{\hspace{-0.05cm} i,\overline{j}}\,P_{\hspace{-0.05cm}kl}}{P_{\hspace{-0.05cm} i,\overline{l} }\, P_{\hspace{-0.05cm} j, k }}
\, , \quad 
\frac{P_{i,j}\,P_{\hspace{-0.05cm}k,\overline{l} }}{P_{\hspace{-0.05cm} i,\overline{l} }\, P_{\hspace{-0.05cm} j, k }}
\qquad 
\end{align}

These expressions can also be seen as formulas for $\tilde x_t^T$ (up to inversion) since the  latter is obtained  by specializing  $x_t^T$ at $\lambda=1$.  Since we have explicit expressions  
\eqref{Eq:Pij-and-co} for the components 
$P_{i,j}$, $P_{\hspace{-0.05cm}i,\overline{j} }$, 
$P_{i,\circ }$ and $P_{\hspace{-0.05cm} i,\circ^{\bowtie}}$
from  which the $\tilde x_t^T$'s are built, we deduce the
\begin{lem}
The type $\tau(x^T_t)$ coincides with the set of indices 
$s\in \{1,\ldots,n\}$ such that $\partial \tilde{x}^T_t/\partial v_s\neq 0$.
\end{lem}
\begin{proof} 
It suffices to verify that for any of the eleven expressions in 
\eqref{Align:tildex-t-T}, its partial derivative with respect to  $v_s$ is not  trivial if and only if $s$ appears in the corresponding type. 
\end{proof}

If two cluster variables $x_{t}^T$ and $x_{t'}^{T'}$ define the same foliation, then the same holds true for $\tilde x_{t}^T$ and $\tilde x_{t'}^{T'}$ hence 
their differentials $d\tilde x_{t}^T$ and $d\tilde x_{t'}^{T'}$ are colinear  (over $\mathbf  C(v_1,\ldots,v_n)$), thus both  necessarily share the same type. On the other hand, one verifies by direct explicit computations
  that the differentials of the eleven explicit rational expressions of four variables in \eqref{Align:tildex-t-T} are pairwise non colinear.  This gives us  the following version in type $D_n$ of Corollary \ref{Cor:X-var-web-An}:
\begin{cor}
\label{Cor:X-var-web-Dn}
\begin{enumerate}
\item Up to pull-back on $\mathcal M_{D_n}$,  
the $\boldsymbol{\mathcal X}$-cluster variables of type $D_n$ are the 
$x_t^T$'s for $T\in \mathcal T^\circ _n$ and $t\in T$.  Thus, up to inversion ({\it i.e.}\,modulo $x\leftrightarrow x^{-1}$),  there are precisely $n(n-1)(n^2+4n-6)/6$ such cluster variables.
\sk
\item Two such cluster variables define the same foliation if and only if  both coincide possibly up to inversion. Consequently, one has    $d^{\boldsymbol{\mathcal X}}_{D_n}=n(n-1)(n^2+4n-6)/6$. 
\end{enumerate}
\end{cor}
 \begin{proof}
 Everything has already be proved above, except the count for the number of quadrilaterals $q_t^T$ which has been computed in \cite[Proposition\,4.8]{Sherman-Bennett1}. 
 \end{proof}

\subsubsection{Types $\boldsymbol{B}$ and $\boldsymbol{C}$.}
The $\boldsymbol{\mathcal X}$-cluster variables in type $B$ and $C$ can be constructed following an approach 
 similar  to the one used in cases 
 $A$ and $D$, namely in terms of some triangulations of some polygons, 
the combinatorial type of these latter being specific to each type.  
 As in the $D_n$-case,  it is possible to define a notion of type for a $\boldsymbol{\mathcal X}$-cluster variable $x$ and for each type, there is a finite list of possible explicit rational expression in some coordinates for $x$: 
 see the formulas analogous to \eqref{Align:tildex-t-T}
  in \S5.2 and in \S5.3 of \cite{Sherman-Bennett1} (for the cases  $B_n$ and $C_n$ respectively). 
  Then arguing as in the paragraph just above Corollary \ref{Cor:X-var-web-Dn}, one deduces that the first point of Conjecture \ref{Conj:Inv-divisor-Fx-Ftildex} holds true as well for any finite  type cluster algebra of Dynkin type  $B$ or $C$. 
 \begin{center}
 $\star$
 \end{center}

 This finishes the (partially sketched) proof of Proposition \ref{Prop:propro}.

\subsubsection{Degrees of cluster webs in finite type.}
\label{SS:DegreesClusterWebsFiniteType}
From above, we deduce the following 
\begin{cor}
\label{Coro:d-X-Delta--d-Y-Delta}
 Let $\Delta$ be any Dynkin diagram (of rank $n\geq 2$). 
\begin{enumerate}
\item  The degree $d_\Delta^{\boldsymbol{\mathcal X}}$  of $\boldsymbol{\mathcal X\mathcal W}_\Delta$ is equal to the cardinal of the set $
\boldsymbol{\mathcal X var}_\Delta^{\pm}$ and is given in Table 
\ref{Table:Degree-d-X}.\sk 
\item 
The degree $d_\Delta^{\boldsymbol{\mathcal Y}}$  of $\boldsymbol{\mathcal Y\mathcal W}_\Delta$ is equal to   $\lvert 
\Delta_{\geq -1}\lvert=
n+\lvert \Delta^+\lvert$ (thus  \eqref{Eq:Majo-dYDelta} is indeed an equality).

Closed 
formulas for $\lvert \Delta^+\lvert$ and $d_\Delta^{\boldsymbol{\mathcal Y}}$ are given in Table 
\ref{Table:Degree-d-Y}. 
\end{enumerate}
\end{cor}
\begin{proof}
The first point follows from \cite{Sherman-Bennett1} combined with the verifications to be made described before  (our Table \ref{Table:Degree-d-X} below 
is extracted from the one at the end of \cite[\S1]{Sherman-Bennett1}).\mk

 Regarding the second point, we recall that for each
 $\alpha\in \Delta_{\geq -1}$, one has $y[\alpha]=P_\alpha/u^{\alpha^\vee}$ for a polynomial $P_\alpha$ with constant term 1 hence $\alpha$ can be recovered from the fiber through the origin of the fibration $\mathcal F_{y[\alpha]}$.  Consequently all the $y[\alpha]$'s are pairwise distinct (even up to inversion) hence  the first point  
  gives us that $d_\Delta^{\boldsymbol{\mathcal Y}}$ is equal to the cardinality of $\Delta_{\geq 1}$.   The computation of $\lvert \Delta^+ \lvert $ in terms  of the  rank of $\Delta$  according to its type  is classical, see  \cite[p.\,216]{Schiffler} for instance.
\end{proof}

\begin{table}[!h]
\begin{center}
\begin{tabular}{|c|c|c|c|c|c|c|c|c|} 
\hline
$\boldsymbol{\Delta}$ & $A_n$& $B_n\, , \, C_n$  & $D_n$ & $E_6$ & $E_7$  & 
$E_8$ & $F_4$ & $G_2$  \\ \hline
\begin{tabular}{c} \vspace{-0.2cm}\\
$
\boldsymbol{d_\Delta^{\boldsymbol{\mathcal X}}}$ \vspace{0.15cm}
\end{tabular}\vspace{0.03cm}
 & 
${n+3 \choose 4 }$
& $\frac{1}{6}n(n+1)(n^2+2)$ &$\frac{1}{6}n(n-1)(n^2+4n-6)$  & 335 & 1050 & 3120 & 98 &  8 \\
\hline
\end{tabular}\mk
\caption{Degrees of ${\boldsymbol{\mathcal X}}$-cluster webs of finite type.}
\label{Table:Degree-d-X}
\end{center}
\end{table}


\begin{table}[!h]
\begin{center}
\begin{tabular}{|c|c|c|
 c|c|c|c|c|c|} 
\hline
$\boldsymbol{\Delta}$ & $A_n$& $B_n\, , \, C_n$ & $D_n$ & $E_6$ & $E_7$  & 
$E_8$ & $F_4$ & $G_2$
 \\ \hline
\begin{tabular}{c} \vspace{-0.2cm}\\
$\boldsymbol{\big\lvert \, \Delta^+\, \big\lvert}$\vspace{0.15cm}
\end{tabular}\vspace{0.03cm}
 & 
$\frac{1}{2}n(n+1)$
& $n^2$ 
 &  $n(n-1)$ & 36 & 63 & 120 & 24 &  6 
  \\ \hline
\begin{tabular}{c} \vspace{-0.2cm}\\
$
\boldsymbol{d_\Delta^{\boldsymbol{\mathcal Y}}}$ \vspace{0.15cm}
\end{tabular}\vspace{0.03cm}
 & 
$\frac{1}{2}n(n+3)$
& $n(n+1)$ 
 &  $n^2$ & 42 & 70 & 128 & 28 &  8  \\
\hline
\end{tabular}\mk
\caption{Number of positive roots in root systems  and degrees of ${\boldsymbol{\mathcal Y}}$-cluster webs of Dynkin type.}
\label{Table:Degree-d-Y} 
\end{center}
\end{table}

We thus have determined the very first basic invariant of cluster webs :  their degree. Remark that it was not trivial to get the answer in  full generality.
 In the next subsection, we consider the case of ${\boldsymbol{\mathcal Y}}$-cluster webs associated to pairs of Dynkin diagrams: we conjecture a nice simple formula for their degrees that however are able to be proveb only in type $(A,A)$. 

\subsubsection{Degrees of ${\boldsymbol{\mathcal Y}}$-cluster webs of  type ${\boldsymbol{(A,A)}}$ and some conjectures.} ${}^{}$ 
\label{SS:Y-cluster-Web-type-AA}


We consider the case of the ${\boldsymbol{\mathcal Y}}$-cluster web of type 
$(A_m,A_n)$ where $m,n$ are two positive integers. Since the case when one of them is 1 has already be treated before, we assume here that both are bigger than or equal to 2. \sk 

 The Dynkin diagrams are simply-laced, the half-periodicity holds true  thus  
${\boldsymbol{\mathcal Y\mathcal W}}_{A_m,A_{n}}$ is the web associated to the period   $\boldsymbol{i}^{half}_{A_m,A_n}$ if $m$ and $n$ have the same parity, or $\boldsymbol{i}^{half}_{\bullet,\circ}$ otherwise 
(see respectively \eqref{Eq:i-Delta-Delta'-half} and \eqref{Eq:PeriodAn-n-even} in Remark \ref{Rem:Half-Periods} above). In both cases,  it is a period of length $mn(h+h')/2$ with $h=h(A_n)=n+1$ and $h'=h(A_m)=m+1$ which gives us $d^{\boldsymbol{\mathcal Y}}_{A_m,A_n}\leq (m+n+2)mn/2$. \sk 

Our goal is to prove that this majoration actually is an equality. To this end, we are going to use the explicit formulas, in terms of cross-ratios, given by Volkov in \cite{Volkov} for the $Y$-variables of the $Y$-system $\boldsymbol{Y}(A_m,A_n)$. As in the $A_n$-case, there are natural birational identifications between ${\rm Conf}_{m+n+2}(\mathbf P^{m})$ and the $\boldsymbol{\mathcal X}$-cluster variety of type ${A_m\boxtimes A_n}$ (see 
\cite[\S6.3]{GGSVV} and 
\cite[\S2.3]{Weng}): given a graph $\Gamma$ of a certain type (namely, a `minimal bipartite graph' ), one constructs a birational map $x^\Gamma: {\rm Conf}_{m+n+2}(\mathbf P^{m})\dashrightarrow \boldsymbol{\mathcal X}_{A_m\boxtimes A_n}$. For instance, one can (and we will) consider the graph $\Gamma_0$ of \cite[Example 2.11]{Weng}.   Composed with the natural projection 
\begin{align*}
\big(\mathbf C^{m+1}\big) & \dashrightarrow {\rm Conf}_{m+n+2}(\mathbf P^{m})\\
 (x_s)_{s=1}^{n+m+2}& \longmapsto  
 \boldsymbol{ \big[\hspace{-0.07cm}\big[} \, [x_1],\ldots,[x_{m+n+2}]\,\boldsymbol{\big]\hspace{-0.07cm}\big]}
\end{align*}
where $[\cdot]$ stands for the projectivization $\mathbf C^{m+1}\setminus \{0\}\rightarrow \mathbf P^r$ and $ \boldsymbol{ [\hspace{-0.05cm}[} - 
\boldsymbol{]\hspace{-0.05cm}]}$ means `modulo ${\rm PSL}_{m+1}(\mathbf C)$', one obtains a dominant rational map $\widehat x^{\Gamma_0}: \big(\mathbf C^{m+1}\big)^{m+n+2}\dashrightarrow \boldsymbol{\mathcal X}_{A_m\boxtimes A_n}$.
\sk

We now recall  formula(1.10) of \cite{Volkov} which is an explicit expression 
for the (pull-back under $\widehat x^{\Gamma_{0}}$ of) the $Y$-variable $Y_{i,j}(t)$ of the $Y$-system $Y_{A_m,A_n}$, in terms of the normalized (one could say `cluster') cross-ratio  $r$  defined in \eqref{Eq:cluster-r}: 
for any $(i,j)\in \{1,\ldots,m\}\times \{1,\ldots,n\}$ and any $t\in \mathbf Z$ such that $i+j+t$ is even, one has 
\begin{equation}
\label{Eq:Volkov-Formula}
\Big(\widehat x^{\Gamma_{0}}\Big)^*Y_{i,j}(t)= r\Big( \, \big[x_a\big]_{i,j}^t\, , 
\, \big[x_b\big]_{i,j}^t\, , 
\, \big[x_c\big]_{i,j}^t\, , 
\, \big[x_d\big]_{i,j}^t\,  \Big)
\end{equation}
as a rational function on $\big(\mathbf C^{m+1}\big)^{m+n+2}$, where 
\begin{itemize}
\item the $x_s$'s  are now cyclically indexed by  $s\in \mathbf Z/(m+n+2)\mathbf Z$
 ({\it i.e.}\,$x_{s+(m+n+2)}=x_s$ for any $s$);
\item $a,b,c,d $ are the four integers 
given by 
\begin{equation}
\label{Eq:ijt->ABCD}
A=
(-i-j-t)/2\, , \qquad 
B= A+i
 \, , \qquad C=B+j
 \quad \mbox{ and } \quad D=C+m+1-i\, ; 
\end{equation}
\item $[ \cdot ]_{i,j}^t$ stands from the linear projection $\mathbf P^m\dashrightarrow \mathbf P^1$ from the linear subspace (generically of dimension $m-1$) spanned   by  the $[x_\sigma]$'s for 
$\sigma \in ]\,A,B\,[\cup ]\,C,D\,[$. \sk 
\end{itemize}
Actually, one has $Y_{i,j}\big(t+2(m+n+2)\big)=Y_{i,j}(t)$ thanks to the full-periodicity of the $Y$-system of type $(A_m,A_n)$ hence only the class of $t$ modulo $2(m+n+2)\mathbf Z$ matters.\mk

Let $(a,b,c,d)$ be a 4-tuple of integers such that 
\begin{equation}
\label{Eq:a,b,c,d}
0<  b-a < m+1\, , \qquad  
0< c-b  < n+1 
\qquad \mbox{ and }\qquad  
(b-a)+( d-c) =\, m+1\, .
\end{equation}
These conditions 
being invariant up to translations (applied componentwise),  they make  sense for the class $[a,b,c,d]_{m,n}$ of $(a,b,c,d)$ modulo $(m+n+2)\mathbf Z$ hence the following set is well-defined: 
\begin{equation}
\boldsymbol{\mathscr Y}_{\hspace{-0.05cm}m,n}=\left\{ \hspace{0.05cm}
\big[ a,b,c,d \big]_{m,n} \hspace{0.05cm} \Big\lvert \hspace{0.05cm}
(a,b,c,d) \in {\mathbf Z^4}  \mbox{ satisfies } \eqref{Eq:a,b,c,d} \hspace{0.05cm}
\right\}\, . 
\end{equation}

Now, given $\big[ a,b,c,d \big]_{m,n}\in \boldsymbol{\mathscr Y}_{\hspace{-0.05cm}m,n}$, one defines a rational map by setting 
\begin{equation}
\label{Eq:Piabcd}
\Pi_{[a,b,c,d]_{m,n}}=\pi_{]\,a,b\,[\,\cup \,]\,c,d\,[, (a,b,c,d)} : {\rm Conf}_{m+n+2}(\mathbf P^{n})\dashrightarrow \mathbf P^1\, ,  
\end{equation}
where we use the notation  \eqref{Eq:pi-I-J} (with  $I=]a,b[\cup ]c,d[$,  $J=(a,b,c,d)$ in the case under scrutiny and where the $x_s$'s are labeled cyclically modulo $m+n+2$): it obviously only depends on $\big[ a,b,c,d \big]_{m,n}$
hence  formula \eqref{Eq:Piabcd} makes sense.
\mk

With this notation Volkov's formula \eqref{Eq:Volkov-Formula}, 
can be paraphrazed by saying that 
the pull-back on ${\rm Conf}_{m+n+2}(\mathbf P^{n})$ of the foliation induced  by $Y_{i,j}(t)$  on $\boldsymbol{\mathcal X}_{A_m\boxtimes A_n}$ coincides with the one admitting 
$\Pi_{[A,B,C,D]_{m,n}}$ as first integral, where $A,B,C,D$ are given by the formulas in \eqref{Eq:ijt->ABCD}.  Conversely, given $\tau=[a,b,c,d]_{m,n}\in \boldsymbol{\mathscr Y}_{\hspace{-0.05cm}m,n}$, setting $i_\tau=b-a$, $j_\tau=c-b$ and $t_\tau=-a-c$, one has $i_\tau+j_\tau+t_\tau=-2a$ ${\rm mod}$ $2(m+n+2)$ which thus is even and one obtains a triplet $(i_\tau,j_\tau,t_\tau)\in \{1,\ldots,\}\times \{1,\ldots,m\}\times \mathbf Z/(2(m+n+2)\mathbf Z)$ which, first, only depends on $\tau$ and  not on the choice of the representative $(a,b,c,d)\in \mathbf Z^4$; and, second,  is such that 
$(x^{\Gamma_0})^*(Y_{i_{\tau},j_\tau}(t_\tau))=\Pi_{\tau}$. 

From the above, we first deduce the
\begin{prop}
The pull-back of $\boldsymbol{\mathcal Y\mathcal W}_{A_m,A_n}$ 
under $x^{\Gamma_{0}}$ is the web on ${\rm Conf}_{m+n+2}(\mathbf P^{n})$ whose first integrals are the maps \eqref{Eq:Piabcd}, for all classes 
$[a,b,c,d]_{m,n}$ in $\boldsymbol{\mathscr Y}_{\hspace{-0.05cm}m,n}$: 
\begin{equation}
\label{Eq:Pull-Back-YW-Am-An}
\big( x^{\Gamma_0}\big)^*\Big( \boldsymbol{\mathcal Y\mathcal W}_{A_m,A_n}\Big)= 
\boldsymbol{\mathcal W}\Big(
\, \mathcal F_{\Pi_\tau}
\hspace{0.1cm}\big\lvert \hspace{0.1cm} \tau \in \boldsymbol{\mathscr Y}_{\hspace{-0.05cm}m,n}\,  \Big)\, . 
\end{equation}
\end{prop}

From this proposition, we deduce the result we were looking for: 
\begin{cor}
\label{Cor:Y-web-type-(Am,An)}
\begin{enumerate}
\item
 The $\boldsymbol{\mathcal Y}$-cluster variables of type $(A_m,A_n)$ satisfy the first point of Conjecture \ref{Conj:Inv-divisor-Fx-Ftildex}: two of them define the same foliation iff they coincide (possibly up to inversion). 
  \item As a consequence of 1., one obtains that $d_{A_m,A_n}^{\boldsymbol{\mathcal Y}}\hspace{-0.1cm} = 
  mn\big(m+n+2\big)/2$.
\end{enumerate}
\end{cor}
\begin{proof}
From Lemma \ref{Lem:pi-I-J--pi-I'-J'}, it follows that for $[a,b,c,d]_{m,n},[a,b,c,d]_{m,n}\in \boldsymbol{\mathscr Y}_{\hspace{-0.05cm}m,n}$, the maps $\Pi_{[a,b,c,d]_{m,n}}$ and $\Pi_{[a',b',c',d']_{m,n}}$  define the same 
foliations on ${\rm Conf}_{m+n+2}(\mathbf P^m)$ if and only if the two sets of intervals $\{\, ]a,b[_{m,n}\, , \,  ]c,d[_{m,n}\, \}$ and $\{ ]a',b'[_{m,n}\, , \,  ]c',d'[_{m,n}\, \}$ coincide, where  the notation ${I}_{m,n}$ means that we are considering the image of $I$ in  $\mathbf Z/(m+n+2)\mathbf Z$ (for any $I\subset \mathbf Z$).  So there are only two possibilities: either (1) 
$ ]a,b[_{m,n}=]a',b'[_{m,n}$ and $ ]c,d[_{m,n}=]c',d'[_{m,n}$ or (2) 
$ ]a,b[_{m,n}=]c',d'[_{m,n}$ and $ ]c,d[_{m,n}=]a',b'[_{m,n}$. 
Since both 4-tuples $(a,b,c,d)$ and $(a',b',c',d')$ satisfy \eqref{Eq:a,b,c,d}, we get  that  $[a,b,c,d]_{m,n}=[a',b',c',d']_{m,n}$  in the first case, whereas 
$[a,b,c,d]_{m,n}=[c',d',a'',b'']_{m,n}$ in the second one, with $a''=a'+(m+n+2)$ and similarly $b''=b'+(m+n+2)$.  In the second case, one deduces 
immediately from \eqref{Eq:Piabcd} 
that  $\Pi_{[a',b',c',d']_{m,n}}=1/\Pi_{[a,b,c,d]_{m,n}}$, which  proves {\it 1}. \mk 

It is then easy to get the second point from the arguments above. Indeed, the (easily seen well-defined and involutive)  map  
$\boldsymbol{\mathscr Y}_{\hspace{-0.05cm}m,n}\rightarrow 
\boldsymbol{\mathscr Y}_{\hspace{-0.05cm}m,n},\, 
[a,b,c,d]_{m,n}\mapsto [c,d,a+(m+n+2),b+(m+n+2)]_{m,n}$ has no fixed point 
(otherwise there would exist $\Pi_{[a,b,c,d]_{m,n}}$ such that 
$\Pi_{[a,b,c,d]_{m,n}}=1/\Pi_{[a,b,c,d]_{m,n}}$). Hence 
from \eqref{Eq:Pull-Back-YW-Am-An}, it comes that the number of distinct foliations of the $\boldsymbol{\mathcal Y}$-cluster web of type $(A_m,A_n)$ is half the cardinality of  $\boldsymbol{\mathscr Y}_{\hspace{-0.05cm}m,n}$. 
But computing the latter is easy: there are $m+n+2$ choices for $a$, $m$ choices for  $b$, $n$ for the one of $c$ and that is all since then $d$ is determined by the rightest relation in \eqref{Eq:a,b,c,d}.  Thus one has  $\lvert \boldsymbol{\mathscr Y}_{\hspace{-0.05cm}m,n} \lvert =(m+n+2)mn$  from which 
{\it 2.} follows. 
\end{proof}

The determination of the degree of ${\boldsymbol{\mathcal Y}\mathcal W}_{A_m,A_n}$ is not really difficult, but not trivial either.  
Note however that we used an ad hoc method for each case. 
\sk 

Since ${\boldsymbol{\mathcal Y}\mathcal W}_{A_m,A_n}$ is the cluster web 
 associated to a period of length $mn(m+n+2)/2$ (half-periodicity), we knew for free that its degree is less than or equal to this number.  If Conjecture \ref{P:Y-web-type-Delta,Delta')} was known to hold true, it would have given directly the second point of  Corollary \ref{Cor:Y-web-type-(Am,An)}.  \sk

It still makes sense to replace $(A_m,A_n)$ by any pair $(\Delta,\Delta')$ of Dynkin diagrams in the preceding paragraph. Assuming that any $Y$-system $Y(\Delta,\Delta')$ satisfies half-periodicity ({\it cf.}\,Remark \ref{Rem:Half-Periods}))  and that Conjecture \ref{P:Y-web-type-Delta,Delta')} holds true, we  would get a proof of the following 
\begin{conjecture}
\label{P:Y-web-type-Delta,Delta')}
Let $\Delta,\Delta'$ be two Dynkin diagrams, of rank $n,n'$ and Coxeter numbers $h,h'$ respectively. 
\begin{enumerate}
\item
 The $\boldsymbol{\mathcal Y}$-cluster variables of type $(\Delta,\Delta')$ satisfy the first point of Conjecture \ref{Conj:Inv-divisor-Fx-Ftildex}. 
 
 \item One has $d_{\Delta,\Delta'}^{\boldsymbol{\mathcal Y}}=
 n n'\big(h+h'\big)/2$.
\end{enumerate}
\end{conjecture}

This conjecture is the first of a series of conjectural statements about $\boldsymbol{\mathcal Y}$-cluster webs of type $(\Delta,\Delta')$ 
that will be discussed further on in \S\ref{SS:YW-Delta-Delta'}.


\subsubsection{Non-linearizability of ${\boldsymbol{\mathcal Y}}$-cluster webs of Dynkin type.} ${}^{}$ 
\label{SS:Y-cluster-Web-non-Lin}
For  any  $n\geq 2$, the $\mathcal X$-cluster web of type 
$A_n$ 
 contains (pull-backs under forgetful maps of) several Bol's webs (aka cluster webs of type $A_2$) henceforth is known to be non-linearizable.  This does not apply to its subweb $\boldsymbol{\mathcal Y\hspace{-0.03cm}\mathcal W}_{\hspace{-0.03cm}A_n}$ which  does not contain any such pull-back hence another argument is needed to show that this web cannot be linearized.
\sk 

In this subsection, arguing by induction on the ranks of Dynkin diagrams (of arbitrary type), we are going to prove the 
\begin{prop}
\label{Prop:YWD-non-linearizable}
 For any Dynkin diagram $\Delta$, the web 
$\boldsymbol{\mathcal Y\hspace{-0.03cm}\mathcal W}_{\hspace{-0.03cm}\Delta}$
is not linearizable.
\end{prop}
When $\Delta$ has rank 2, one can use the classical characterization of linearizable planar webs (see \cite[\S6.1.3]{Coloquio} for instance) to verify that none of the cluster webs of type $A_2,B_2$ and $G_2$ is linearizable.  \mk

We now assume that $n={\rm rk}(\Delta)\geq 3$.  Let 
$\boldsymbol{\mathcal Y\hspace{-0.03cm}\mathcal W}^1_{\hspace{-0.03cm}\Delta}$
be
 the $2n$-subweb  
of $\boldsymbol{\mathcal Y\hspace{-0.03cm}\mathcal W}_{\hspace{-0.03cm}\Delta}=\boldsymbol{\mathcal W}\big(\, y[\alpha]\, \lvert \, \alpha\in \Delta_{\geq -1}\, \big)$ ({\it cf.}\,\eqref{Eq:YWDelta}) whose first integrals are the $y[\pm\alpha_i]$'s for $i=1,\ldots,n$. It can easily be verified (using the formalism explained in \cite[\S2]{KellerAnnals} if needed) that one has   
\begin{equation}
\label{Eq:z[alpha]}
\boldsymbol{\mathcal Y\hspace{-0.03cm}\mathcal W}_{\hspace{-0.03cm}\Delta}^1=
\boldsymbol{\mathcal W}\Big( 
\, u_i\, , \, z[\alpha_i]\hspace{0.1cm} \big\lvert \hspace{0.1cm} i\leq n\, 
\Big) 
\qquad \mbox{with}\qquad  z[\alpha_i]=\frac{\prod_{j=1 }^n \big(1+u_j\big)^{a_{ij}}}{u_i}
\quad \big(\, i=1,\ldots,n\,\big)\, ,  
\end{equation}
 where the integers $a_{ij}$'s are the coefficients of the matrix
$A_\Delta=(a_{ij})_{i,j=1}^n=2{\rm Id}_n-C_{\Delta}$ where $C_\Delta$ stands for the Cartan matrix of type  $\Delta$. We are going to prove that $\boldsymbol{\mathcal Y\hspace{-0.03cm}\mathcal W}_{\hspace{-0.03cm}\Delta}^1$ is not linearizable which 
 will imply Proposition \ref{Prop:YWD-non-linearizable} immediately.\sk 
 
We first consider the case when $n=3$. Even if the web 
 $
\boldsymbol{\mathcal Y\hspace{-0.03cm}\mathcal W}_{\hspace{-0.03cm}\Delta}^1$ 
for $\Delta=A_3,B_3$ or $C_3$  does not satisfy the `{\it Strong general position assumption}' discussed in \S\ref{Par:GermOfWebs} (see Remark 
\ref{Rem:wGP-cluster-Webs} in case $A_3$), one can adapt for it the approach of \cite{PirioLin}  (especially the arguments used in the proof of Proposition 4.3 therein) and verify by means of explicit computations that $
\boldsymbol{\mathcal Y\hspace{-0.03cm}\mathcal W}_{\hspace{-0.03cm}\Delta}^1$ is not compatible with a projective connexion and consequently, cannot be linearized.
\sk

To handle the higher rank case, we notice that 
 for any $k\leq n$ labeling an exterior vertex of $\Delta$ (that is a vertex adjacent with only one edge of the Dynkin diagram), then the graph 
$\Delta_k$  obtained by removing the $k$-th vertex (and the edge adjacent to it) in $\Delta$ is still a Dynkin diagram. The crucial fact which will allow us to argue by induction is given by the following
\begin{lem} 
For any complex number $\lambda$ distinct from $0$ and $-1$, the trace of 
$\boldsymbol{\mathcal Y\hspace{-0.03cm}\mathcal W}_{\hspace{-0.03cm}\Delta}^1$ along the affine hyperplane cut out by $u_k=\lambda$ coincides with the 
$2(n-1)$-web $\boldsymbol{\mathcal Y\hspace{-0.03cm}\mathcal W}_{\hspace{-0.03cm}\Delta_{k}}^1$ : one has 
$$
\Big(\boldsymbol{\mathcal Y\hspace{-0.03cm}\mathcal W}_{\hspace{-0.03cm}\Delta}^1\Big)\big\lvert_{u_k=\lambda }=\boldsymbol{\mathcal Y\hspace{-0.03cm}\mathcal W}_{\hspace{-0.03cm}\Delta_{k}}^1\, .
$$
\end{lem}
\begin{proof}
This follows easily from the explicit formulas given in \eqref{Eq:z[alpha]} for the first integrals $z[\alpha_i]$ combined with the fact that the Cartan matrix 
 of $\Delta_k$ is the one obtained by removing the $k$-th line and the $k$-th column of $C_\Delta$. 
\end{proof}

With the preceding lemma, we can argue by induction as follows: assume that 
$\boldsymbol{\mathcal Y\hspace{-0.03cm}\mathcal W}_{\hspace{-0.03cm}\Delta}^1$ is linearizable and that  $k$ labels an exterior vertex of $\Delta$. 
Hence for $\lambda $ generic, the fiber $u_k=\lambda$ is a leaf of $\boldsymbol{\mathcal Y\hspace{-0.03cm}\mathcal W}_{\hspace{-0.03cm}\Delta}^1$ hence the trace of this web on it is linearizable as well. Hence so is $\boldsymbol{\mathcal Y\hspace{-0.03cm}\mathcal W}_{\hspace{-0.03cm}\Delta_{k}}^1$ according to the previous lemma.
Since  $\boldsymbol{\mathcal Y\hspace{-0.03cm}\mathcal W}_{\hspace{-0.03cm}\Delta}^1$ is not linearizable when  $n=3$, we obtain that the same holds true for $n$ arbitrary which in turn proves Proposition \ref{Prop:YWD-non-linearizable}.

\newpage 
\section{Some specific examples of cluster webs}
\label{S:SpecialsCases}

In this somehow intermediary section, we look at some webs associated to some of the most classical polylogarithmic identities and show that these are of cluster type. 
\sk

We first discuss carefully the cluster webs associated to rank 2 Dynkin diagrams. In particular, we give an explicit basis of the space of their ARs. 
 Then we turn to more sophisticated cluster webs and/or cluster algebras and are able to recover several webs associated to polylogarithmic FEs, for weights from 2 to 4. Our results are gathered in Theorem \ref{T:classical-cluster-webs}.

 \subsection{\bf Cluster webs in rank 2}
\label{SS:cluster-webs:rank2}

Here we consider the case when $\Delta$ is a (bipartite weighted) Dynkin quiver of rank 2, namely one of the following type $A_2$, $B_2=C_2 $ and $G_2$ respectively:
$$
 \scalebox{0.35}{ \includegraphics{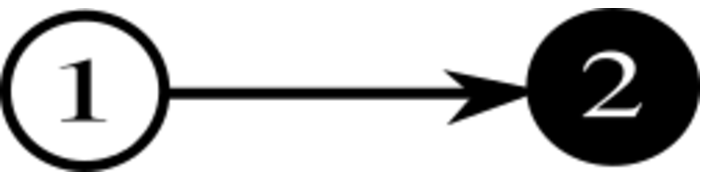}}
\qquad \quad  \scalebox{0.35}{ \includegraphics{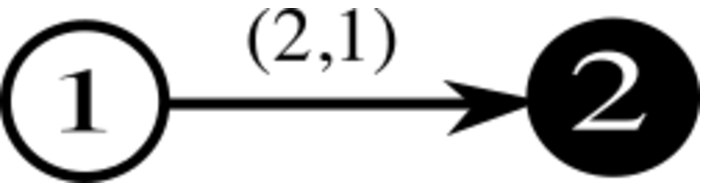}}
\qquad \quad \scalebox{0.35}{ \includegraphics{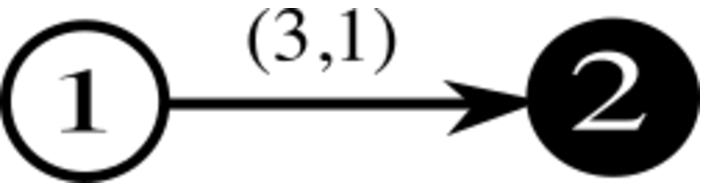}}\, .
$$
The  associated 
exchange matrices are 
$$
B_{A_2}=\begin{bmatrix} 0& 1 \\
-1 & 0 \end{bmatrix}\, , \qquad
B_{B_2}=\begin{bmatrix} 0& 2 \\
-1 & 0 \end{bmatrix} \qquad
 \mbox{and }  \qquad 
 B_{G_2}=\begin{bmatrix} 0& 3 \\
-1 & 0 \end{bmatrix}\, .
$$

Taking in each case  $((a_1,a_2),(x_1^{-1},x_2),B_{\Delta})$ as initial cluster seed,  
we get that the associated cluster webs are :   
\begin{itemize}
\item[$\bullet$]  {\bf {$\big[$ Type $\boldsymbol{A_2}$ $\big]$} }  : 
\begin{align*}
\vspace{-0.5cm}
\boldsymbol{\mathcal A\mathcal W}_{\! A_2}=&  \,  
\boldsymbol{\mathcal W}\left( a_1\, , \, a_2
\, , \, \frac{1+a_1}{a_2} 
\, , \, \frac{1+a_2}{a_1} 
\, , \, \frac{1+a_1+a_2}{a_1a_2} \, 
\right)\, , 
\\
\boldsymbol{\mathcal X\mathcal W}_{\! A_2}=&  \,  
\boldsymbol{\mathcal W}\left( x_1\, , \, x_2
\, , \, \frac{1+x_1}{x_2} 
\, , \, \frac{1+x_2}{x_1} 
\, , \, \frac{1+x_1+x_2}{x_1x_2} \, 
\right)\, ;
\end{align*}
\item[$\bullet$]   {\bf {$\big[$ Type $\boldsymbol{B_2}$ $\big]$} }  :   
\begin{align*}
\boldsymbol{\mathcal A\mathcal W}_{B_2}=&  \,  
\boldsymbol{\mathcal W}\left(\, 
a_{{1}},
\,a_{{2}},
\,{\frac {1+a_{{1}}}{a_{{2}}}},
\,{\frac {1+{a_2}^2}{a_1}},
\,{\frac {1+a_{{1}}+{a_{{2}}}^{2}}{a_{{1}}a_{{2}}}},
\,{\frac {(1+a_1)^{2}+a_{{2}}^{2}}{{a_{{1}}a_{{2}}}^{2}}} \, 
\right)\, , 
\\
\boldsymbol{\mathcal X\mathcal W}_{B_2}=&  \,  
\boldsymbol{\mathcal W}\left(\, x_{{1}}\, ,
  \,{\frac { \left( 1+x_{{1}} \right) ^{2}}{x_{{2}}}}\, ,
   \,{\frac {(1+x_{{1}})^2+x_{{2}}}{x_{{1}}x_{{2}}}}\, , 
      \,{\frac { \left(1+x_{{1}}+ x_{{2}} \right)^{2}}{{x_{{1}}}^{2}x_{{2}}}}\,  , 
  \,{\frac {1+x_{{2}}}{x_{{1}}}}\,  , \,x_{{2}}\, 
\right)\, ; 
\end{align*}
\item[$\bullet$]   {\bf {$\big[$ Type $\boldsymbol{G_2}$ $\big]$} }  : 
    (there is no need to write down  $\boldsymbol{\mathcal A\mathcal W}_{G_2}$, see  Remark \ref{Rem:riri} below) : 
\begin{align*}
\boldsymbol{\mathcal X\mathcal W}_{G_2}=&  \, 
\boldsymbol{\mathcal W}\Bigg(\, x_{{1}}\, ,\,
{\frac { \left( 1+x_{{1}} \right) ^{3}}{x_{{2}}}}\, ,\,
{\frac {(1+x_{{1}})^3+x_{{2}}}{x_{{1}}x_{{2}}}}\, , \, 
{\frac { \left( \big(1+{x_{{1}}}\big)^{2}+x_{{2}} \right) ^{3}}{
{x_{{1}}}^{3}{x_{{2}}}^{2}}}  \, , \,
\\
& {}^{} \hspace{3cm}
{\frac {\big( 1+x_{{1}}\big)^{3}+
2x_2+3x_1x_2+{x_{{2}}}^{2}
}{{x_{{1}}}^{2}x_{{2}}}}\, , \,  
{\frac { \big( 1+x_{{1}} +x_{{2}} \big) ^{3}}{{x_{{1}}}^{3}x_{{2}}}}\, ,\,
{\frac { 1+x_{{2}}}{x_{{1}}}}\, ,\,
x_{{2}}
\,\Bigg)\, .
\end{align*}
\end{itemize}

All the cluster webs above are ordered.  The first integrals of the $\boldsymbol{\mathcal A}$-cluster webs are ordered with respect to the natural order on the weights of their denominators. As for the $\boldsymbol{\mathcal X}$-versions, their first integrals are  labeled according to  \eqref{Eq=l-th-cluster-variable-period}, with respect to the following periods: 
\begin{equation}
\boldsymbol{i}_{A_2}=(1,2)^2\lvert 1 \, , \qquad 
\boldsymbol{i}_{B_2}=(1,2)^3
\qquad \mbox{ and } \qquad 
\boldsymbol{i}_{G_2}=(1,2)^4\, , 
\end{equation}
where for any $k\in \mathbf N^*$, $(1,2)^k$ stands for the concatenation of $k$ copies of $(1,2)$.
\begin{center}
$\star$
\end{center}

Below, we are going to  study these webs by determining some  of their invariants, more specifically those related to their ARs. 

\begin{rem}
\label{Rem:riri}
An interesting (or possibly just curious) fact about the cluster webs above is that 
$$p^*(\boldsymbol{\mathcal X\mathcal W}_{\Delta})=
\boldsymbol{\mathcal A\mathcal W}_{\Delta} $$
for any Dynkin diagram $\Delta$ of rank 2. This implies that in this case, the $\boldsymbol{\mathcal A}$-cluster webs carry interesting polylogarithmic ARs (of weight 1 or 2) as  well. This seems to be very specific to the rank 2 case since all the other cluster webs 
defined by $\boldsymbol{\mathcal A}$-cluster variables that we have been able to consider do not carry any dilogarithmic ARs.   It would be interesting to know a conceptual explanation of this (if there is one). 
\end{rem}

\subsubsection{The cluster web of type $\boldsymbol{A_2}$.}
\label{SS:cluster-webs-type-A2}


It is well-known that the  dilogarithmic identity 
\begin{equation}
\label{Eq:R-A2}
{\mathsf R}(x_1) + {\mathsf R}(x_2) + {\mathsf R}\left(\frac{1 + x_2}{x_1}\right) + 
{\mathsf R}\left(\frac{1 + x_1+x_2}{x_1x_2}\right) +{\mathsf R}\left(\frac{1 + x_1}{x_2}\right)=
\frac{{}^{}\hspace{0.1cm}\pi^2}{2}
\end{equation}
(which holds  true indentically for $x_1,x_2>0$) is equivalent to (Rogers' version \eqref{Eq:EFA-R} of) Abel's functional equation $(\boldsymbol{\mathcal A b})$. It gives us that $\boldsymbol{\mathcal X\mathcal W}_{\hspace{-0.05cm}A_2}$  is equivalent to Bol's web.  Another (more geometric) way to prove this comes by noticing that the first four  rational arguments of ${\mathsf R}$ in 
 \eqref{Eq:R-A2} define four pencils of lines 
  the vertices of which form the base locus of the pencil of conics defined by the fifth argument of  ${\mathsf R}$  in this identity.  In more mathematical terms, we have an equivalence of webs: 
  \begin{equation}
  \label{Eq:XWA2--Bol}
  \boldsymbol{\mathcal B}\, \simeq \, \boldsymbol{\mathcal X\mathcal W}_{\!A_2}\, .
  \end{equation}
%

We denote by $X_1,\ldots,X_5$  the rational functions appearing in 
\eqref{Eq:R-A2} in the corresponding order: $X_1=x_1$, $X_2=x_2$, $X_3=(1+x_2)/x_1$, etc.\footnote{{\it I.e.}\,one has $X_\ell=\boldsymbol{x}_{\boldsymbol{i}_{A_2}}(\ell)$
 for $\ell=1,\ldots,5$.}    Setting $X_0=X_5$ and $X_6=X_{1}$, it is known that the following relation holds true for any $\ell=1,\ldots,5$ (as it is well-known, these relations can be used to define the $X_\ell$'s inductively): 
 \begin{equation}
 \label{Eq:Xell}
  {}^{} \quad X_{\ell+1}X_{\ell-1}=1+X_\ell \quad 
 \Longleftrightarrow \quad {\rm Log}\big(X_{\ell-1}\big)-{\rm Log}\big(1+X_\ell\big)+{\rm Log}\big(X_{\ell+1}\big) =0\, . 
 \end{equation}
 For each $\ell$, the logarithmic identity on the right in \eqref{Eq:Xell} can be seen as an  AR for $\boldsymbol{\mathcal X\mathcal W}_{\!A_2}$, denoted by $LogAR_\ell$.   We denote by ${\sf R }_{\! A_2}$ the dilogarithmic AR associated to \eqref{Eq:R-A2}.  Then one has
  \begin{equation}
  \label{Eq:AR-XWA2}
 \boldsymbol{\mathcal A}\big( \boldsymbol{\mathcal X\mathcal W}_{\hspace{-0.07cm}A_2} \big)= \left\langle \, LogAR_\ell\, \lvert \, \ell=1,\ldots,5\, \right\rangle \oplus 
 \big\langle \, {\sf R }_{\! A_2} \, \big\rangle\,   
 \end{equation}
 from which it can be deduced that $ \boldsymbol{\mathcal X\mathcal W}_{\!A_2} $ is of maximal rank, hence is AMP.\sk 
 
 Since the period ${\boldsymbol{i}_{A_2}}=(1,2,1,2,1)$ used to get $\boldsymbol{\mathcal X\mathcal W}_{\!A_2}$ is unchanged by a cyclic shift of length 2, it follows that $ F_{A_2}: (x_1,x_2)\rightarrow (( 1+x_1)/x_2, x_1)$ is a birational automorphism of order 5 of this web such that  $\big(F_{A_2}\big)^*(X_\ell)=X_{\ell+1}$ for $\ell=1,\ldots,5$. Thus 
 $\big(F_{A_2}\big)^*(LogAR_\ell)=LogAR_{\ell+1}$ for any $\ell\in \mathbf Z/5\mathbf Z$
and $\big(F_{A_2}\big)^*({\sf R }_{\! A_2})={\sf R }_{\! A_2}$ from which we deduce that 
the decomposition as a direct sum \eqref{Eq:AR-XWA2} is invariant with respect to the natural action of the (birational) automorphism $F_{A_2}$. \mk 


Another nice feature of choosing the cluster variables $X_\ell$ as first integrals for the $\boldsymbol{\mathcal X}$-cluster web of type $A_2$ 
 is that one can easily construct 
the dilogarithmic identity ${\sf R }_{\! A_2}$ from the five logarithmic abelian relations $LogAR_{\ell}$'s. Indeed, considering the differential of 
$LogAR_{\ell}$, it follows that 
$$
\frac{dX_\ell}{1+X_\ell}=\frac{dX_{\ell-1}}{X_{\ell-1}}+\frac{dX_{\ell+1}}{X_{\ell+1}}
$$
for any $\ell$. Then, summing for $\ell $ in $ \mathbf Z/5\mathbf Z$, it comes 
\begin{align}
\label{Eq:RA2-from-LogARl}
-\sum_{\ell}  LogAR_{\ell}\,  \frac{dX_\ell}{X_\ell}= & \, 
\sum_{\ell }
\bigg( -{\rm Log}\big(X_{\ell-1}\big)+{\rm Log}\big(1+X_\ell\big)-{\rm Log}\big(X_{\ell+1}\big) \bigg)
\,
  \frac{dX_\ell}{X_\ell} \nonumber \\ 
  = & \,  
 \sum_{\ell}
 {\rm Log}\big(1+X_\ell\big)\,
  \frac{dX_\ell}{X_\ell} 
  -
  \sum_{\ell}
  {\rm Log}\big(X_{\ell}\big)\, \bigg( \frac{dX_{\ell-1}}{X_{\ell-1}}+\frac{dX_{\ell+1}}{X_{\ell+1}} 
  \bigg) 
 \nonumber   \\
   = & \,  \sum_{\ell}
\bigg( \frac{{\rm Log}\big(1+X_{\ell}\big)}{X_\ell}
- \frac{{\rm Log}\big(X_{\ell}\big)}{1+X_\ell}
\bigg) \, 
{dX_\ell}\\
  = & \, \sum_{\ell} 2 \,{\sf R}'\big(X_\ell\big) \,dX_\ell=
  2 \,d\bigg(   \sum_{\ell} {\sf R}\Big(X_\ell \Big) 
  \bigg)\,. \nonumber
\end{align}
Since all the $LogAR_{\ell}$'s vanish, the same holds true for 
 the total derivative of $ \sum_{\ell=1}^5 {\sf R}(X_\ell)$ hence this sum is  identically equal to a constant. We thus have recovered very symmetrically the $A_2$-dilogarithmic identity $({\sf R}_{A_2})$ from the logarithmic identities $LogAR_{\ell}$'s. 
\begin{center}
$\star$
\end{center}
All the material discussed above is well-known. 
The reason which motivated us to recall all that is that everything can be formulated  in a very nice and effective way 
within the cluster algebra formalism.

\subsubsection{The cluster web of type $\boldsymbol{B_2}$.}
\label{SS:cluster-webs:B2}
Here, we denote by $X_\ell$ the $\boldsymbol{\mathcal X}$-cluster variables appearing in the definition  of $\boldsymbol{\mathcal X\mathcal W}_{\hspace{-0.05cm}B_2}$ above: we have $X_1=x_1$ and 
$$
  X_2={\frac { \left( 1+x_{{1}} \right) ^{2}}{x_{{2}}}}\,, \quad 
  X_3={\frac {1+x_{{2}}}{x_{{1}}}}\,  , \quad 
  X_4= {\frac { \left(1+x_{{1}}+ x_{{2}} \right)^{2}}{{x_{{1}}}^{2}x_{{2}}}}\, , \quad     X_5={\frac {(1+x_{{1}})^2+x_{{2}}}{x_{{1}}x_{{2}}}}
  \quad 
  \mbox{and}
   \quad
X_6=x_2  \, .
$$

The associated $B_2$-cluster dilogarithmic identity is 
$$
\big({\sf R}_{B_2}\big)\quad : \hspace{1.6cm}
2\,{\sf R}(X_1)+
{\sf R}(X_2)+
2\,{\sf R}(X_3)+
{\sf R}(X_4)+
2\,{\sf R}(X_5)+
{\sf R}(X_6)=\pi^2\, .
\hspace{1.6cm} {}^{}
$$

On the other hand, one verifies that the following rational relation between the cluster variables $X_1,X_3$ and $X_5$ holds true identically: 
\begin{equation}
\label{E:EqXYZ=1} 
\frac{1}{1+X_1}+ \frac{1}{1+X_3}+\frac{1}{1+X_5}=1\, .
\end{equation}

Then one defines the following three rational quantities 
$$
X=X_1+1=1+x_1 \,  , \quad 
Y=X_5+1= \frac{1+x_1+x_2}{x_1}
 \quad \quad \mbox{ and} \quad \quad
Z=X_3+1=\frac{(1+x_1)(1+x_1+x_2)}{x_1x_2}
$$
 in such a way that  $\varphi : (x_1,x_2)\mapsto (X,Y)$ induces a 
 (birational) change of coordinates such that 
 \begin{equation}
 \label{Eq:XWB2--WN6}
 \varphi_*\Big(  \boldsymbol{\mathcal X\mathcal W}_{\hspace{-0.05cm}B_2}
 \Big)=\boldsymbol{\mathcal W}\bigg( \,    X\, , \, X+Y-XY
 \, , \, Y \, , \,  YZ-Y-Z \, , \, Z \, , \, XZ-X-Z
  \bigg)
  \end{equation}
  (equality between ordered webs), 
where $X,Y$ and $Z$ satisfy  the relation $X^{-1}+Y^{-1}+Z^{-1}=1$ ({\it cf.} 
 \eqref{E:EqXYZ=1}).   In this way, one obtains that 
  $ \boldsymbol{\mathcal X\mathcal W}_{\hspace{-0.05cm}B_2}$
   is equivalent to Newman's web $\boldsymbol{\mathcal W}_{\hspace{-0.03cm} \boldsymbol{\mathcal N}}$  associated to Newman's bilogarithmic identity  discussed in paragraph \S\ref{Par:Newman-equation} above.  Then setting 
$$ X=x\, , \qquad Y=\frac{x(1-y)}{y(1-x)} \qquad \mbox{ and } \qquad Z=-\frac{x(1-y)}{1-x} \,  , $$ 
one deduces that Newman's 6-web is equivalent to 
$\boldsymbol{\mathcal W}\Big(x,xy, \frac{x}{y} ,\frac{x(1-y)}{y(1-x)}, \frac{x(1-y)}{1-x}, \frac{x(1-y)^2}{y(1-x)^2} \Big)$. This is a  subweb of Spence-Kummer web
 ${\boldsymbol{\mathcal W}}_{\boldsymbol{\hspace{-0.05cm}{\mathcal S}{\mathcal K}}}$,   denoted by ${\boldsymbol{\mathcal W}}_{\widehat{248}}$ in  \cite{PirioThese}. In particular, we know that this web, hence consequently 
$\boldsymbol{\mathcal X\mathcal W}_{\hspace{-0.05cm}B_2}$, has maximal rank. 
\mk

We now turn to describing a nice basis of the space of ARs of 
$\boldsymbol{\mathcal X\mathcal W}_{\hspace{-0.05cm}B_2}$. 
To deal with dilogarithmic ARs, one considers the following two iterated integrals (of weight 2), defined for any $u>0$: 
\begin{equation}
\label{Eq:L01-L10}
L_{10}(u)=\int_0^u \frac{\log(s)}{1+s} ds
\qquad \mbox{and}
\qquad 
L_{01}(u)=\int_0^u\frac{\log(1+s)}{s} ds\, .
\end{equation}
Then one verifies that the symmetric iterated integral of weight two  
\begin{equation}
\label{Eq:Dilog-Symmetric-S}
{\sf S}=\frac{1}{2}\Big( L_{10}+L_{01}\Big)
\end{equation}
  satisfies identically the following identity: 
\begin{align*}
\big({\sf S}_{B_2}\big)\quad  \hspace{2.6cm}
2\, {\sf S}(X_1) -{\sf S}(X_2)+2\, {\sf S}(X_3) -{\sf S}(X_4)+ 2\, {\sf S}(X_5) -{\sf S}(X_6)=&0 \hspace{2.6cm}{}^{}
\end{align*}
\label{(S-B2)}
and that together with $({\sf R}_{B_2}\big)$, it forms a basis of the space of dilogarithmic abelian relations of $\boldsymbol{\mathcal X\mathcal W}_{\hspace{-0.05cm}B_2}$. Note that $({\sf S}_{B_2})$ is symmetric with respect to the non-trivial involution $\sigma=(12)$ on weight two polylogarithmic iterated integrals, whereas 
$({\sf R}_{B_2})$ is antisymmetric for the action of $\sigma$, {\it i.e.}  
$$
\big({\sf R}_{B_2}\big)^{\sigma}=-\big({\sf R}_{B_2}\big)\qquad \mbox{ and }
\qquad 
\big({\sf S}_{B_2}\big)^{\sigma}=\big({\sf S}_{B_2}\big)\, . 
$$

As for the logarithmic ARs, recall that the cluster variables $X_\ell$'s can also be defined recursively by means of the following polynomial relations (where $\ell$ is taken modulo 6): 
\begin{align*}
 \quad X_{\ell+1}X_{\ell-1}= 
(1+X_\ell) \quad ( \mbox{for}\, \ell \, \mbox{ even})\, \qquad  
\quad \mbox{ and }\qquad 
 \quad X_{\ell+1}X_{\ell-1}= 
(1+X_\ell)^2 \quad (\mbox{for}\, \ell \, \mbox{ odd})\,  .
\end{align*}

To these relations are associated the following logarithtmic identities which  can be considered as as many ARs for $\boldsymbol{\mathcal X\mathcal W}_{\hspace{-0.05cm}B_2}$: 
\begin{align}
{\rm Log}\big(X_{\ell-1}\big)-{\rm Log}\big(1+X_\ell\big)-{\rm Log}\big(X_{\ell+1}\big)=0&   \quad ( \mbox{for}\, \ell \, \mbox{ even})\\ 
 \mbox{ and }
 \qquad {\rm Log}\big(X_{\ell-1}\big)-2\,{\rm Log}\big(1+X_\ell\big)-{\rm Log}\big(X_{\ell+1}\big)=0& \quad (\mbox{for}\, \ell \, \mbox{ odd})\,  . \nonumber
\end{align}

We get six logarithmic abelian relations, denoted by $LogAR_\ell$ for $\ell=1,\ldots,6$, which are 
 linearly independant and which, together with the dilogarithmic ARs $({\sf R}_{B_2}\big)$ and $({\sf S}_{B_2})$ and the rational one \eqref{E:EqXYZ=1}, span a subspace of dimension 9 of $\boldsymbol{\mathcal A}(\boldsymbol{\mathcal X\mathcal W}_{\hspace{-0.05cm}B_2})$. 
 One last AR is missing, in order to have a basis of this space.  
 One verifies that it is the one associated to the following identity : 
\begin{align*}
\big({\sf A}_{B_2}\big)\quad  \hspace{3cm}
{\sf A}\big(X_2\big)+{\sf A}\big(X_4\big)+{\sf A}\big(X_6\big)=0 \hspace{7cm}{}^{}
\end{align*}
which holds true for any $x_1,x_2>0$, where ${\sf A}$ stands for the function defined by 
\begin{equation}
\label{Eq:Arctan-sqrt}
{\sf A}(u)={\rm Arctan}\big(\sqrt{u}\, \big)\, , \quad u>0\, .
\end{equation}

Thus, we end with the following decomposition of the space of ARs: 
\begin{equation}
\label{Eq:AR-XWB2}
\boldsymbol{\mathcal A}\Big(\boldsymbol{\mathcal X\mathcal W}_{\hspace{-0.05cm}B_2}\Big) = \underbrace{\Big\langle 
\scalebox{0.97}{$LogAR_\ell$} \, \lvert \, 
\scalebox{0.83}{$\ell=1,\ldots,6$}
\, \Big\rangle}_{logarithmic}  \, \oplus \, 
\underbrace{
\Big\langle
\, {\sf R}_{B_2} \, , \, {\sf S}_{B_2}\, 
 \Big\rangle  }_{dilogarithmic} \, \oplus  \,
 \underbrace{
  \big\langle
\, \eqref{E:EqXYZ=1}\,  
 \big\rangle}_{rational}\, 
 \oplus \, \big\langle
\, {\sf A}_{B_2}\, 
 \big\rangle\, . 
\end{equation}

\begin{rem}
\label{Rk:TildeXWB2}
At first sight, the abelian relation $({\sf A}_{B_2})$ looks a bit particular  compared to the others since  it doesn't look (and rigorously speaking, is not) of the same nature as the others, which all are polylogarithmic (here we see the rational AR 
\eqref{E:EqXYZ=1} as polylogarithmic, but of weight 0).   But there is a way to make  this go away. Indeed, considering the expression 
$ \arctan(u)=\frac{i}{2}\log\big( (u+i)/(u-i)\big)+\pi/2$  (valid for any $u>0$) 
of the arctangent function by means of logarithms,   it comes that $({\sf A}_{B_2})$ is not far from being  a logarithmic AR, the obstruction for this being the square root appearing in the definition \eqref{Eq:Arctan-sqrt} of ${\sf A}$. But denoting by $\tilde X_\ell$ the pull-back of $X_\ell$ under $\Psi : (u_1,u_2)\mapsto (u_1,u_2^2)=(x_1,x_2)$, we obtain that $\tilde X_\ell$ is a square for $\ell$ even: for any such $\ell$,  there exists $Z_\ell\in \mathbf Q_{sf}(u_1,u_2)$ such that $\tilde X_\ell=(Z_\ell)^2$, from which it  follows that $\Psi^*({\sf A}_{B_2})$ is written  $${\rm Arctan}(Z_2)+{\rm Arctan}(Z_4)+{\rm Arctan}(Z_6)=0\,,  $$ which is an identity 
actually of logarithmic type since it corresponds to the vanishing of the logarithmic derivative of  the following expression
which can be verified to be  identically equal to -1:
 $$\frac{(Z_2+i)}{(Z_2-i)}\frac{(Z_4+i)}{(Z_4-i)}\frac{(Z_6+i)}{(Z_6-i)}\, .$$ 

From the above considerations, it follows that  
the pull-back of $\boldsymbol{\mathcal X\mathcal W}_{\hspace{-0.05cm}B_2}$ under $\Psi$ is a model of this web which is  defined by positive (substraction-free) rational first integrals (with coefficients in $\mathbf Z$) and whose ARs are all of iterated integral type, of weight 0 (rational), 1 or 2. Note however that  $\Psi^*\big(\boldsymbol{\mathcal X\mathcal W}_{\hspace{-0.05cm}B_2}\big)$ is no longer polylogarithmic since for instance the ramification locus associated to its first integral $Z_\ell$  is $\{\, 0,-1, \, \pm i\, ,\, \infty\, \}$ for any $\ell \in \{2,4,6\}$.\mk

A completely similar phenomenon holds true for the cluster web of type $G_2$ 
({\it cf.}\,Remark \ref{Rk:TildeXWG2}).  Both could be two particular cases of a  general phenomenon (see Conjecture \ref{Conj:Nature-Fi} above).\mk 
\end{rem}

The following map $F_{\hspace{-0.06cm}B_2}: (x_1,x_2)\mapsto \big(X_3,X_4\big) =
\big((1+x_2)/{x_1}\, , \, 
{(1+x_1+x_2)^2}/{(x_1^2x_2)}\big) $ 
is a birational automorphism of order 3 such that $F_{ \hspace{-0.06cm} B_2}^*(X_\ell)=X_{\ell-2}$ for any $\ell$.  The decomposition as a direct sum \eqref{Eq:AR-XWB2} is invariant by this map.   More precisely, one has $\big(F_{\hspace{-0.06cm} B_2}\big)^*( LogAR_\ell)=LogAR_{\ell-2}$ for any $\ell$ whereas all the other ARs appearing in 
\eqref{Eq:AR-XWB2} are invariant by $F_{\hspace{-0.06cm} B_2}$. \sk

It is known that Newman's identity $({\cal N}_6)$ is accessible from three copies of Abel's five terms relation ({\it cf.}\,\S16.3 in \cite{Lewin}). Another interesting feature of $F_{ \hspace{-0.06cm} B_2}$ is that it can be used to obtain this quite neatly.  
 Indeed, modulo the inversion relation ${\mathsf R}(1/x)={\mathsf R}(x)+\pi^2/6$ (satisfied for any $x>0$), Abel's identity is equivalent to the fact that  the function of two variables $x,y>0$ 
$$AB(x,y)= 
{\mathsf R}\bigg(\frac{1}{x}\bigg) -{\mathsf R}(y) + {\mathsf R}
\left(\frac{x}{1+y}\right)
 - {\mathsf R}\left(\frac{1+x+y}{xy}\right)
+ {\mathsf R}\left(\frac{y}{1+x}\right)
$$
vanishes identically.  So the same holds true for the combination 
$AB+AB\circ F_{\hspace{-0.05cm}B_2}+AB\circ \big(F_{\hspace{-0.05cm}B_2}\big)^2$ 
and one verifies that this is  formally equivalent to 
$$ 
2\,{\mathsf R}\bigg(\frac{1}{x}\bigg) 
-{\mathsf R}(y)
+2\,{\mathsf R}\bigg(\frac{x}{1+y}\bigg) 
-{\mathsf R}\bigg(\frac{(1+x+y)^2}{x^2y}\bigg) 
+ 2\,{\mathsf R}\bigg(\frac{xy}{(1+x)^2+y}\bigg) 
-{\mathsf R}\bigg(\frac{(1+x)^2}{y}\bigg) 
\equiv 0\, ,  $$
 an identity which is nothing but $({\sf R}_{B_2})$ (up to several uses of the inversion relation for ${\sf R}$).
 \begin{center}
 $\star$
 \end{center}

 The presentation above of the cluster web of type $B_2$ and of its main features  is neat and well formalized, in a much  better way than what was known before about Newman's web.  It is another illustration of the value of having a description of a web as a web of cluster type.


%
%


%
%
\subsubsection{The cluster web of type $\boldsymbol{G_2}$.}
\label{SS:cluster-webs:G2}
The cluster web associated to the initial seed 
$(1/x_1,x_2) , B_{G_2})$  with respect to the period $(1,2)^4$ is 
\begin{align}
\label{Eq:XWG2}
\boldsymbol{\mathcal X \hspace{-0.05cm}\mathcal W}_{G_2}=&  \, 
\boldsymbol{\mathcal W}\Bigg(\, x_{{1}}\, ,\,
{\frac { \left( 1+x_{{1}} \right) ^{3}}{x_{{2}}}}\, ,\,
{\frac {(1+x_{{1}})^3+x_{{2}}}{x_{{1}}x_{{2}}}}\, , \, 
{\frac { \left( \big(1+{x_{{1}}}\big)^{2}+x_{{2}} \right) ^{3}}{
{x_{{1}}}^{3}{x_{{2}}}^{2}}}  \, , \, \nonumber
\\
& {}^{} \hspace{3cm}
{\frac {\big( 1+x_{{1}}\big)^{3}+
2x_2+3x_1x_2+{x_{{2}}}^{2}
}{{x_{{1}}}^{2}x_{{2}}}}\, , \,  
{\frac { \big( 1+x_{{1}} +x_{{2}} \big) ^{3}}{{x_{{1}}}^{3}x_{{2}}}}\, ,\,
{\frac {  1+x_{{2}}}{x_{{1}}}}\, ,\,
x_{{2}}
\,\Bigg)\, .
\end{align}
We denote by $X_1,\ldots,X_8$ the cluster variables defining 
$\boldsymbol{\mathcal X \hspace{-0.05cm}\mathcal W}_{G_2}$: $ X_1=x_1$, $X_2={\frac { \left( 1+x_{{1}} \right) ^{3}}{x_{{2}}}}, \ldots$, $X_8=x_2$. 
\mk 

\begin{figure}[h]
\begin{center}
\resizebox{3.2in}{3.2in}{
 \includegraphics{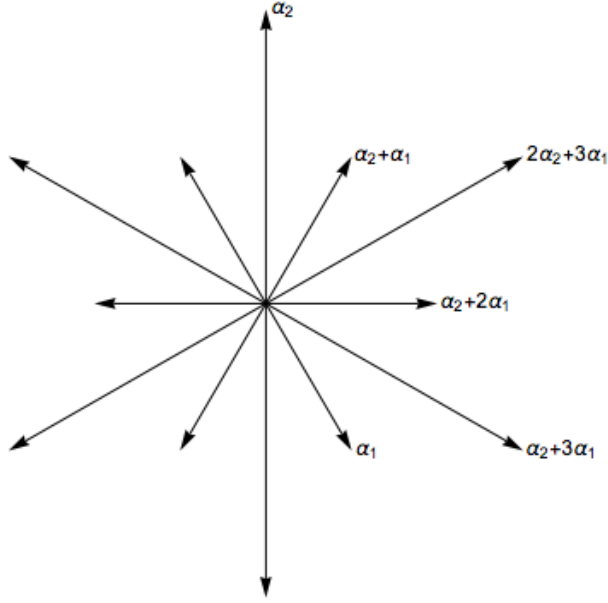}}
 \vspace{-0.3cm}
\caption{The root system of type $G_2$} 
\label{Fig:RG2}
\end{center}
\end{figure}

\vspace{-0.3cm}
Recall the root system of type $G_2$, denoted by $R$: see Figure \ref{Fig:RG2} where the principal roots are denoted by $\alpha_1$ and $\alpha_2$.   
There are eight elements in $R_{\geq -1}$:  four short roots $\pm \alpha_1$, $\alpha_1+\alpha_2$ and $2\alpha_1+\alpha_2$, and as many 
long roots, which are $\pm \alpha_2$, $3\alpha_1+\alpha_2$ and $3\alpha_1+2\alpha_2$.


\mk

In rank 2, the $\boldsymbol{\mathcal X}$- and $\boldsymbol{\mathcal Y}$-cluster webs coincide, hence one can label the $\boldsymbol{\mathcal Y}$-variables which are first integrals for $\boldsymbol{\mathcal X\mathcal W}_{\hspace{-0.05cm}}$ by the elements $\alpha\in R_{\geq 1}$ or equivalently by the pair of coordinates $(a_1,a_2)\in \mathbf Z^2$ of $\alpha=a_1\alpha_1+a_2\alpha_2$ 
with respect to $\alpha_1$ and $\alpha_2$.  
 In \eqref{Eq:XWG2}, the coordinate $x_1$ is  related to the short root $\alpha_1$, whereas $x_2$ is related to  $\alpha_2$ hence the first integrals of $\boldsymbol{\mathcal X \hspace{-0.05cm}\mathcal W}_{G_2}$ can be labeled as follows (where $U_{a_1,a_2}$ stands for $U_{a_1\alpha_1+a_2\alpha_2}$ for any $a_1,a_2\in \mathbf Z$): 
 \begin{align*}
U_{-1,0}=& \,x_1 \hspace{0.1cm} && U_{1,0}= \frac{1+x_2}{x_1} &&  \hspace{0.1cm} U_{1,1}= {\frac {(1+x_1)^3+x_2}{x_1x_2}}      &&U_{3,1}= {\frac { \big( 1+x_1 +x_2 \big) ^{3}}{{x_1}^{3}x_2}}
\\ 
U_{0,-1}=&\, x_2 \hspace{0.1cm} && U_{0,1}=  \frac{(1+x_1)^3}{x_2} && \hspace{0.1cm} U_{2,1}={\frac {\big( 1+x_1\big)^{3}+
x_2(2+3x_1+{x_2})
}{{x_1}^{2}x_2}}  &&U_{3,2}= {\frac { \left( \big(1+{x_1}\big)^{2}+x_2 \right) ^{3}}{
{x_1}^{3}{x_2}^{2}}}\, .
&&
 \end{align*}


In terms of the $X_\ell$'s or of the $U_\alpha$'s, the dilogarithmic identity of type $G_2$ is written: 
$$
\big( {\sf R}_{G_2}\big)\quad :  \hspace{1.5cm} \; 
3\sum_{\ell \, odd} {\sf R}\big(X_{\ell}\big)
+\sum_{\ell\, even} {\sf R}\big(X_{\ell}\big)
=
3\hspace{-0.2cm}  \sum_{\alpha\in R_{\geq -1}^{short}}{\sf R}\big(U_{\alpha}\big)+ \sum_{\beta\in R_{\geq -1}^{long}}{\sf R}\big(U_{\beta}\big)=2\,\pi^2\, .  
\hspace{1cm}{}^{}
$$

As in the $B_2$-case, there is another polylogarithmic functional identity satisfied by the symmetric dilogarithm ${\sf S}$ defined in \eqref{Eq:Dilog-Symmetric-S}. Indeed, the following identity is  satisfied for every  $x_1,x_2>0$: 
$$
\big( {\sf S}_{G_2}\big)\quad :  \hspace{1.5cm} \; 
3\sum_{\ell \, odd} {\sf S}\big(X_{\ell}\big)
-\sum_{\ell\, even} {\sf S}\big(X_{\ell}\big)
=
3\hspace{-0.2cm} \sum_{\alpha\in R_{\geq -1}^{short}}{\sf S}\big(U_{\alpha}\big)- \sum_{\beta\in R_{\geq -1}^{long}}{\sf S}\big(U_{\beta}\big)=0\, .  
\hspace{1.5cm}{}^{}
$$

The two ARs associated to $( {\sf R}_{G_2})$ and $( {\sf S}_{G_2})$ (denoted in the same way) form a basis of the space of weight two iterated integrals ARs of the cluster web of type $G_2$.
\begin{rem}
  Contrarily to $( {\sf R}_{\Delta})$ when $\Delta=A_2$ or $B_2$, the identity $( {\sf R}_{G_2})$ does not seem to have been known by classical authors.  To the best of our knowledge, it only appeared rather recently in the literature once a link between Lie algebras and dilogarithmic indentities was uncovered by math-physicists. 
  For some references in which (some variants of) this identities appear, see    {\rm \cite[\S B.5.3]{AlexandrovPioline}} and {\rm \cite{KY}}. 
\end{rem}

Regarding  the logarithmic ARs, what happens is similar to the previous cases: the  $X_\ell$'s can also be defined recursively by means of the following polynomial relations (where $\ell$ is taken modulo 8): 
\begin{align*}
 \quad X_{\ell+1}X_{\ell-1}= 
(1+X_\ell) \quad ( \mbox{for}\, \ell \, \mbox{ even})\, \quad  
\quad \mbox{ and }\quad 
 \quad X_{\ell+1}X_{\ell-1}= 
(1+X_\ell)^3 \quad (\mbox{for}\, \ell \, \mbox{ odd})\,  .
\end{align*}

To these relations are associated the following logarithtmic identities which  can be considered as as many ARs for $\boldsymbol{\mathcal X\mathcal W}_{\hspace{-0.05cm}G_2}$: 
\begin{align}
{\rm Log}\big(X_{\ell-1}\big)-{\rm Log}\big(1+X_\ell\big)-{\rm Log}\big(X_{\ell+1}\big)=0&   \quad \mbox{for }\, \ell \, \mbox{ even}\,, \\ 
 \mbox{ and}
 \quad {\rm Log}\big(X_{\ell-1}\big)-3\,{\rm Log}\big(1+X_\ell\big)-{\rm Log}\big(X_{\ell+1}\big)=0& \quad \mbox{for }\, \ell \, \mbox{ odd}. \nonumber
\end{align}

We get eight  logarithmic abelian relations, denoted by $LogAR_\ell$ for $\ell=1,\ldots,8$, which are 
 linearly independant and which, together with the dilogarithmic ARs $({\sf R}_{B_2}\big)$ and $({\sf S}_{B_2})$, span a subspace of dimension 10 of $\boldsymbol{\mathcal A}(\boldsymbol{\mathcal X\mathcal W}_{\hspace{-0.05cm}G_2})$.  \sk

 It turns out that the web ${\boldsymbol{\mathcal X \mathcal W}}_{G_2}$ is not of maximal rank 21 but of  rank 14 (this can be established using the method sketched in \S\ref{SubPar:CharacterizationWebsMaximalRank}).    In addition of its polylogarithmic ARs, ${\boldsymbol{\mathcal X \mathcal W}}_{G_2}$ carries four other ARs which are not of this kind.   These are ARs of the two following 4-subwebs  of the web under scrutiny: 
\begin{align*}
{\boldsymbol{\mathcal W}}_{\hspace{-0.05cm}G_2}^{short}=& \, \boldsymbol{\mathcal W}\Big(\, X_{\ell}\, \lvert \,   \ell=1,3,5,7\, 
 \Big)=\boldsymbol{\mathcal W}\Big(\, U_{\alpha}\, \lvert \,   \alpha \in  R
 _{\geq -1}^{ short }\, 
 \Big)=  \boldsymbol{ \mathcal W}\Big( U_{-1,0}\, , \, U_{1,0}\, , \, U_{1,1}\, , \, U_{2,1}\, \Big)\\
 \mbox{ and }\quad 
 {\boldsymbol{\mathcal W}}_{\hspace{-0.05cm} G_2}^{long}=& \,  \boldsymbol{\mathcal W}\Big(\, X_{\ell}\, \lvert \,   \ell=2,4,6,8\, 
 \Big)=\boldsymbol{ \mathcal W}\Big(\, U_{\alpha}\, \lvert \,   \alpha \in   R
 _{\geq -1}^{ long}\, 
 \Big)=\boldsymbol{\mathcal W}\Big( U_{0,-1}\, , \, U_{0,1}\, , \, U_{3,1}\, , \, U_{3,2}\, \Big)\, .
\end{align*} 
 
 The subweb associated to short roots carries two linearly independent ARs with rational components, namely the ones associated to the following two rational identities: 
\begin{align}
\label{Eq:RA-rat-WG2}
1=&\,\sum_{\ell\, odd} \frac{1}{1+X_\ell} =  \sum_{\alpha\in R_{\geq -1}^{short}} \frac{1}{1+U_\alpha}
\vspace{0.1cm}\\
\mbox{and }\quad 1 =  & \,
 \frac{1}{\big(1+X_{1}\big)^2}
 + \frac{1+2\, X_3}{\big(1+X_{3}\big)^2}
  +  \frac{1}{\big(1+X_{5}\big)^2}
  +   \frac{1+2\, X_7}{\big(1+X_{7}\big)^2} \,  . \nonumber 
 \end{align}

The 4-web ${\boldsymbol{\mathcal W}}_{\hspace{-0.05cm} G_2}^{long} $ has rank 3. It carries one logarithmic AR (which belongs to the linear span of the $LogAR_\ell$'s) and two other new ARs, whose type is different from all the types of the previously considered abelian relations.  Indeed, it can be verified ({\it cf.}\,Example \ref{Ex:WG2-short-Abel'sMethod} above) that 
the following identity holds true identically
\begin{align}
\label{Eq:RA-exotic-WG2}
0= & F\big(X_2   \big) -F\big(X_4   \big) +F\big(X_6   \big) -F\big(X_8   \big)\, 
 \end{align}
where for any  $u\in \mathbf R_{>0}$,    $F(u)$ stands for 
${\rm Log}\big(1+u^{\frac{1}{3}}\big)$ or 
${\rm Log}\big(1-u^{\frac{1}{3}}+u^{\frac{2}{3}}\big)$. 
%
%
\mk

Thus, we end with the following decomposition of the space of ARs: 
\begin{equation}
\label{Eq:AR-XWG2}
\boldsymbol{\mathcal A}\Big(\boldsymbol{\mathcal X\mathcal W}_{\hspace{-0.05cm}G_2}\Big) = \underbrace{\Big\langle \, 
\scalebox{0.97}{$LogAR_\ell$} \, \lvert \, 
\scalebox{0.83}{$\ell=1,\ldots,8$}
\, \Big\rangle}_{logarithmic}  \, \oplus \, 
\underbrace{
\Big\langle
\, {\sf R}_{G_2} \, , \, {\sf S}_{G_2}\, 
 \Big\rangle  }_{dilogarithmic} \, \oplus  \,
 \underbrace{
  \big\langle
\, \eqref{Eq:RA-rat-WG2}\,  
 \big\rangle}_{rational}\, 
 \oplus \, \big\langle
\, \eqref{Eq:RA-exotic-WG2}\, 
 \big\rangle\, . 
\end{equation}

\begin{rem}
\label{Rk:TildeXWG2}
All the ARs in \eqref{Eq:AR-XWG2} are polylogarithmic, except the two 
'exotic' ones associated to the two identities \eqref{Eq:RA-exotic-WG2}.
But as in the $B_2$ case ({\it cf.}\,Remark \ref{Rk:TildeXWG2} above), there is a formally natural way to remedy this and to get a model of $\boldsymbol{\mathcal X\mathcal W}_{\hspace{-0.05cm}G_2}$ whose all ARs are of iterated integral type.\sk

Indeed,  denoting by $\tilde X_\ell$ the pull-back of $X_\ell$ under $\Psi : (u_1,u_2)\mapsto (u_1,u_2^3)=(x_1,x_2)$, we obtain that $\tilde X_\ell$ is a cube for $\ell$ even: for any such $\ell$,  there exists $Z_\ell\in \mathbf Q_{sf}(u_1,u_2)$ such that $\tilde X_\ell=(Z_\ell)^3$. In terms of these $Z_\ell$'s, the ARs associated to \eqref{Eq:RA-exotic-WG2} are  logarithmic since one verifies that they 
are obtained by taking the logarithmic derivatives of the following multiplicative identities:
\begin{align}
\frac{(1+Z_2)(1+Z_6)}{(1+Z_4)(1+Z_8)}=1
 \qquad 
\mbox{and}\qquad 
\frac{\big(Z_2+\varrho\big)\big(Z_2+\overline{\varrho}\big)\cdot 
\big(Z_6+\varrho\big)\big(Z_6+\overline{\varrho}\big)
}{\big(Z_4+\varrho\big)\big(Z_4+\overline{\varrho}\big)\cdot \big(Z_8+\varrho\big)\big(Z_8+\overline{\varrho}\big)}=1 
\end{align}
with $\varrho=(1+i\sqrt{3})/2$  (thus $\varrho^3=-1$ and 
$1-X+X^2=\big(X+{\varrho}\big) \big(X+\overline{\varrho}\big)$ as polynomials in  $X$).\sk

It follows that  $\Psi^*\big(\boldsymbol{\mathcal X\mathcal W}_{\hspace{-0.05cm}G_2}\big)$ is a model of the $G_2$-cluster web which is again defined by positive (substraction-free) rational first integrals (with coefficients in $\mathbf Z$) and whose ARs are all of iterated integral type, of weight 0 (rational), 1 or 2. 
Note however that  it is no longer polylogarithmic since for any even $\ell$,  the ramification locus associated to its first integrals $Z_\ell$  is $\{\, 0,-1, \, \varrho\, , 
\, \overline{\varrho} \, , \, 
\, \infty\, \}$.\mk

All this is completely similar to what holds true in the $B_2$-case  
({\it cf.}\,Remark \ref{Rk:TildeXWB2})  and we believe that it is just a particular case of a more general phenomenon.
\end{rem}

Finally we mention that, as far as we know, the accessibility of $({\sf R}_{G_2})$ from Abel's five terms relation has not been worked out in detail yet in the published literature. It is mentioned in \cite{KY} where the authors claim that it can 
be obtained,   by folding,  from the accessibility of the $\boldsymbol{\mathcal Y}$-cluster dilogarithm identity  of type $D_4$, which can be established easily.\footnote{Strictly speaking, the authors in \cite{KY} consider the quantum versions of the classical functional identities we are talking about here. But it is our interpretation that their claims can also be applied to the classical dilogarithmic identities.}  
For all we know, the description of  the $D_4$-identity as a combination of several 5-terms identities $({\sf R}_{A_2})$ has not been worked out in explicit form anywhere.

 



\subsection{\bf Some classical polylogarithmic identities are of cluster type}
\label{SS:ClassicalPolylogarithmicIdentitiesAreOfClusterType}
If it is now common knowledge that cluster algebras provide many interesting dilogarithmic identities, the same was not clear concerning the case of higher polylogarithms.  In the result below, we show that some of the most well-known polylogarithmic AFEs of weight 2,3 and 4 are of cluster type.

\begin{thm}
\label{T:classical-cluster-webs}
Up to some birational automorphisms (which can be make explicit in each case), one has the following interpretations in terms of cluster webs of the planar webs associated to the classical functional equations  of low-order polylogarithms presented in Section \S\ref{S:EFA-polylog}: 
\begin{itemize}
\item[] $\bullet$  {\rm \big[\,Rank 2\,\big] :}   \hspace{0.7cm} $
 \boldsymbol{\cal B}\, \simeq \, {\boldsymbol{\mathcal X\mathcal W}}_{\! A_2}$\,; \hspace{0.65cm} 
 ${\boldsymbol{\mathcal W}}_{\!{\cal N}_6} \, 
\simeq \, {\boldsymbol{\mathcal X\mathcal W}}_{B_2}$; 
\smallskip 
\item[] $\bullet$ {\rm \big[\,Rank 3\,\big] :}
\hspace{0.05cm}
 ${\boldsymbol{\mathcal W}}_{\hspace{-0.03cm}{\cal S}{\cal K}}\, \simeq\,  {\boldsymbol{{\mathcal U} \hspace{-0.04cm} {\mathcal X }
\hspace{-0.04cm}
{\mathcal W}}}_{A_3}$\,; \hspace{0.35cm} ${\boldsymbol{\mathcal W}}_{\!{\cal K}_4} \, \simeq \, {\boldsymbol{{\mathcal U} \hspace{-0.03cm} {\mathcal X }
\hspace{-0.03cm}
{\mathcal W}}}_{ \! C_3}$\,; \hspace{0.1cm}
$ {\boldsymbol{\mathcal W}}_{\!{\cal N}_6}
 \,  \simeq  \, 
{\boldsymbol{{\mathcal U} \hspace{-0.03cm} {\mathcal Y}
\hspace{-0.05cm}
{\mathcal W}}}_{A_3}$;
\smallskip 
\item[] $\bullet$  {\rm \big[\,Rank 4\,\big] :} \,  ${\boldsymbol{\mathcal W}}_{\!{\cal K}_4} \, \simeq \, {\boldsymbol{{\mathcal U} \hspace{-0.03cm} {\mathcal X }
\hspace{-0.03cm}
{\mathcal W}}}_{ \! D_4}$\,; \hspace{0.1cm}
$ {\boldsymbol{\mathcal X\mathcal W}}_{ \! G_2} \,  \simeq\,  {\boldsymbol{{\mathcal U} \hspace{-0.04cm} {\mathcal Y}
\hspace{-0.04cm}
{\mathcal W}}}_{\! D_4}$.
\smallskip 
\end{itemize}

\end{thm}

Before entering into a proof, let us comment briefly this result which 
was the one which motivated the author to undertake a more systematic study of cluster webs.  \sk 

First, the most interesting equivalences in this theorem are the ones 
 relating the polylogarithmic webs associated to the classical AFEs
 $(\boldsymbol{\mathcal Ab})$, $(\boldsymbol{\mathcal S\mathcal K})$ and $(\boldsymbol{\mathcal K}_4)$ to some cluster webs of Dynkin types. The equivalences between two cluster webs in the statements above  are just manifestations for cluster webs of the well-known folding property of Dynkin diagrams, with the unique exception of 
 ${\boldsymbol{{\mathcal U} \hspace{-0.03cm} {\mathcal X }\hspace{-0.03cm}
{\mathcal W}}}_{\hspace{-0.05cm} C_3}$ and ${\boldsymbol{{\mathcal U} \hspace{-0.03cm} {\mathcal X }\hspace{-0.03cm}
{\mathcal W}}}_{\hspace{-0.05cm} D_4}$
 which both are equivalent to 
 ${\boldsymbol{\mathcal W}}_{\!{\cal K}_4}$.  This is a bit surprising since  the Dynkin types $D_4$ and $C_3$ are not related in a natural way. This contrasts with  $D_4$ and $B_3$, the latter being a $2:1$ folded version of the former. However, ${\boldsymbol{{\mathcal U} \hspace{-0.03cm} {\mathcal X }\hspace{-0.03cm}
{\mathcal W}}}_{\hspace{-0.05cm} B_3}$ is not equivalent to ${\boldsymbol{\mathcal W}}_{{\boldsymbol{\cal K}}_4}$ (for instance because one has ${\bf Hex}_3({\boldsymbol{{\mathcal U} \hspace{-0.03cm} {\mathcal X }\hspace{-0.03cm}
{\mathcal W}}}_{\hspace{-0.05cm} B_3})=96 < 186= 
{\bf Hex}_3({\boldsymbol{\mathcal W}}_{{\boldsymbol{\cal K}}_4})$.\footnote{We wonder 
whether this strange fact does not come from a possible confusion we made at some point: here it could be the case that we actually have worked with $B_3$ instead of $C_3$ and vice versa ({\it cf.}\,the warning page \pageref{Pageref:Warning-Apology}).}
\sk
  
%
Finally,  it is expected (at least by the author of these lines) that other  previously known polylogarithmic identities  can also be formulated using cluster variables.  For instance, since it holds true in weight 3 and 4, 
one would not be surprised whether  the answer to the following question were  affirmative: 
\begin{question}
\label{Question:K5-cluster-type?}
Is Kummer's pentalogarithmic web  ${\boldsymbol{\mathcal W}}_{{\boldsymbol{\cal K}}_5}$ equivalent to a cluster web?  
\end{question}
\begin{proof}
Since all the webs involved are completely explicit, in order to prove the theorem it suffices to give for each case, an appropriate explicit birational map. 
\sk

For Bol's and Newman's webs, this has been done above, see  \eqref{Eq:XWA2--Bol} and \eqref{Eq:XWB2--WN6} respectively. 
\sk


Regarding Spence-Kummer's 9-web,  setting 
 $\varphi: (x,y)\dashrightarrow \big(\,  1/y-1\, , \, -{1}/{(xy)}\big)$,  
 it is easy to verify by direct easy computations that 
$\varphi^*\big( {\boldsymbol{\mathcal U\hspace{-0.03cm}\mathcal X \hspace{-0.03cm} \mathcal W}}_{A_3})= {\boldsymbol{\mathcal W}}_{\hspace{-0.03cm} \boldsymbol{{\cal S}{\cal K}}}$ hence proving the theorem in this case. 
However, the simple analytical form of $\varphi$ could let the reader think that finding an appropriate birational equivalence  between the two webs involved 
  is straightforward 
and that proving the statement we are interested in only amounts to some elementary computations.  We claim that it is not the case. 
 Below, we explain  how we have been able to  find explicitly   the appropriate birational equivalence  
 in two cases, namely 
 $\varphi
$ between $
{\boldsymbol{{\mathcal U} \hspace{-0.03cm}{\mathcal X} \hspace{-0.03cm} {\mathcal W}}}_{\! A_3}
$ and ${\boldsymbol{\mathcal W}}_{\hspace{-0.03cm} \boldsymbol{{\cal S}{\cal K}}}$.  This case will show the reader that finding a suitable $\varphi$ relies on some non-trivial (but classical indeed) results and constructions of web geometry. 
 In our opinion, treating the  case under scrutiny 
is sufficiently enlightening  so that we do not give details about the other cases, since they can be treated using similar arguments.    
    \sk

The key point for finding an explicit equivalence between 
${\boldsymbol{{\mathcal U} \hspace{-0.03cm}{\mathcal X} \hspace{-0.03cm} {\mathcal W}}}_{\! A_3}
$ and ${\boldsymbol{\mathcal W}}_{\hspace{-0.03cm} \boldsymbol{{\cal S}{\cal K}}}$ is to realize that both contain a unique hexagonal 6-subweb (as explained in \S\ref{}, testing if a given web is hexagonal amounts to compute the  curvatures of all of its 3-subwebs, and this can be done easily with the help of a computer algebra system), 
which will be denoted by $
{\boldsymbol{{\mathcal U} \hspace{-0.03cm}{\mathcal X} \hspace{-0.03cm} {\mathcal W}}'}_{\hspace{-0.03cm}A_3}$ and 
$
{\boldsymbol{\mathcal W}'}_{\hspace{-0.15cm} \boldsymbol{{\cal S}{\cal K}}}$ 
 respectively.\footnote{The 6-web $
{\boldsymbol{\mathcal W}'}_{\hspace{-0.15cm} \boldsymbol{{\cal S}{\cal K}}}$ is the subweb of 
${\boldsymbol{\mathcal W}}_{\hspace{-0.03cm} \boldsymbol{{\cal S}{\cal K}}}$
denoted by ${\boldsymbol{\mathcal W}}_{\hspace{-0.03cm}\widehat{369}}$ in \cite{PirioSelecta}.}
 Any (even local)  analytic equivalence $\phi$ between ${\boldsymbol{{\mathcal U} \hspace{-0.03cm}{\mathcal X} \hspace{-0.03cm} {\mathcal W}}}_{\! A_3}$ and ${\boldsymbol{\mathcal W}}_{\hspace{-0.03cm} \boldsymbol{{\cal S}{\cal K}}}$ induces an equivalence between 
${\boldsymbol{{\mathcal U} \hspace{-0.03cm}{\mathcal X} \hspace{-0.05cm} {\mathcal W}}}_{ \hspace{-0.05cm} A_3}'$ and ${\boldsymbol{\mathcal W}}_{\hspace{-0.03cm} \boldsymbol{{\cal S}{\cal K}}}'$  and this is going to simplify everything.  Indeed, according to Bol's theorem ({\it cf.} \S\ref{S:Webs}), these 6-webs both are equivalent to webs 
denoted by 
${\boldsymbol{{\mathcal U} \hspace{-0.03cm}{\mathcal X} \hspace{-0.03cm} {\mathcal W}}}_{ \hspace{-0.03cm}A_3}^{lin}$ and ${\boldsymbol{\mathcal W}}_{\hspace{-0.03cm} \boldsymbol{ \hspace{-0.03cm} {\cal S}{\cal K}}}^{lin}$ respectively, each formed by six pencils of lines. Then,  since the linearization of a linear $k$-web is projectively unique as soon as $k\geq 4$ ({\it cf.}\,Proposition \ref{P:Unique-Lin}), one obtains that 
${\boldsymbol{{\mathcal U} \hspace{-0.03cm}{\mathcal X} \hspace{-0.03cm} {\mathcal W}}}_{ \hspace{-0.03cm}A_3}^{lin}$ and ${\boldsymbol{\mathcal W}}_{\hspace{-0.03cm} \boldsymbol{ \hspace{-0.03cm} {\cal S}{\cal K}}}^{lin}$ 
are  necessarily  equivalent by means of a projective transformation,  which can be determined using elementary techniques of plane projective geometry\footnote{By projective duality, this is equivalent to  finding a projective equivalence between two configurations of 6 pairwise distinct points on $\mathbf P^2$.}. 
\mk

We now briefly discuss the case of ${\boldsymbol{\mathcal W}}_{\hspace{-0.03cm} {\cal K}_4} $ and ${\boldsymbol{{\mathcal U} \hspace{-0.03cm} {\mathcal X }
\hspace{-0.03cm} {\mathcal W}}}_{ \! D_4}$.  
By straightforward computations (done with the help of a computer algebra system), one 
can determine explicitly the hexagonal 6-subwebs of each of these two 18-webs.  Each contains precisely four such subwebs, the union of which forms a 14-subweb.  By considering the complement of these latter, one gets  two 4-subwebs  
${\boldsymbol{\mathcal W}}_{\hspace{-0.03cm} {\cal K}_4}' $ and 
${\boldsymbol{\mathcal W}}_{\hspace{-0.04cm}D_4}' $
of  
${\boldsymbol{\mathcal W}}_{\hspace{-0.03cm} {\cal K}_4} $ and ${\boldsymbol{{\mathcal U} \hspace{-0.03cm} {\mathcal X }
\hspace{-0.03cm} {\mathcal W}}}_{ \! D_4}$
 respectively.  Explicitly, one has 
\begin{align*}
{\boldsymbol{\mathcal W}}_{\hspace{-0.03cm} {\cal K}_4}'=& \, \boldsymbol{\mathcal W}\bigg(\,  \frac{xy^2}{(1-y)(1-x)^2} \, , \,  \frac{ x^2y}{(1-y)^2(1-x)}\, ,\,  -
\frac{xy^2(1-x)}{1-y}\, , \,  -\frac{x^2y(1-y)}{1-x} \, \bigg) \\
\mbox{ and }
\quad {\boldsymbol{\mathcal W}}_{\hspace{-0.04cm}D_4}' = & \,  \boldsymbol{\mathcal W}\bigg( \,  \frac{1}{y}\, , \,  \frac{(x+1)^3}{y}\, , \, \frac{(x+1+y)^3}{x^3y}\, , \,\frac{(x^2+2x+y+1)^3}{x^3y^2}\, \bigg)\, .
\end{align*}

One verifies easily  that both these 4-subwebs have maximal rank hence are linearizable in a unique way (up to projective equivalence).  Determining explicitly their abelian relations is not difficult using Abel's method. Thus one can compute the associated Poincar\'e-Blaschke maps ${\rm P \! B}[{{\boldsymbol{\mathcal W}}_{\!{\cal K}_4}'}]$ and ${\rm P \! B}[{\boldsymbol{\mathcal W}}_{\hspace{-0.04cm}D_4}' ]$ (see \S\ref{SS:Poincare-Blaschke-Maps}). From  Proposition \ref{P:Unique-Lin}, it comes that any (even local) analytic equivalence between 
 ${\boldsymbol{\mathcal W}}_{\hspace{-0.03cm} {\cal K}_4} $ and ${\boldsymbol{{\mathcal U} \hspace{-0.03cm} {\mathcal X }
\hspace{-0.03cm} {\mathcal W}}}_{ \! D_4}$ necessarily induces an equivalence 
 between   ${\boldsymbol{\mathcal W}}_{\hspace{-0.03cm} {\cal K}_4}' $ and 
${\boldsymbol{\mathcal W}}_{\hspace{-0.04cm}D_4}' $  
from which one deduces that   ${\rm P \! B}[{{\boldsymbol{\mathcal W}}_{\!{\cal K}_4}'}]^{-1}\circ g\circ {\rm P \! B}[{p{\cal X} \! {\boldsymbol{\mathcal W}}_{D_4}'}]$ for some $g\in {\rm PGL}_3(\mathbf C)$ which can be determined explicitly. 
One eventually gets the following  birational map  
$$
\Phi \, : \, (x,y)\longmapsto \bigg(   \frac{x+y-1}{(x-1)(y-1)} \, ,\, \frac{x^2y}{(x-1)(y-1)^2}       \bigg)
$$
which can easily be verified to be  such that  $\Phi^*\big( {\boldsymbol{{\mathcal U} \hspace{-0.03cm} {\mathcal X }
\hspace{-0.03cm} {\mathcal W}}}_{ \! D_4} \big)  = {\boldsymbol{\mathcal W}}_{\!{\cal K}_4}$.   Using the same arguments, one finds an explicit birational maps $\psi$  
inducing an equivalence between 
$ {\boldsymbol{{\mathcal U} \hspace{-0.03cm} {\mathcal X }
\hspace{-0.03cm} {\mathcal W}}}_{ \! C_3} $ and  $ {\boldsymbol{\mathcal W}}_{\!{\cal K}_4}$. \sk 

All the other equivalences in Theorem \ref{T:classical-cluster-webs} are obtained using a similar approach. \end{proof}

Considering the previous result naturally leads to wonder about 
other polylogarithmic identities. The two first explicit  examples coming to author's mind 
are Kummer's pentalogarithmic identity ({\it cf.}\,Question \ref{Question:K5-cluster-type?} just above) or  Goncharov's functional equation
$(\boldsymbol{\mathcal G on} _{22}$) 
 of the trilogarithm (see \S\ref{Par:Goncharov-22} above and also \S\ref{SS:ClusterNature} in the last section).  At the other extreme, the most general question Theorem \ref{T:classical-cluster-webs}   
 suggests  is that  of knowing whether  there are some polylogarithmic AFEs which are not accessible from polylogarithmic identities  of cluster type.  Answering this seems completely out of reach at the time of writing.

\newpage



\newpage 
\section{Cluster webs of type $\boldsymbol{A}$.}
\label{S:ClusterWebs-type-A}
In this section, we study the ${\boldsymbol{\mathcal X}}$- and ${\boldsymbol{\mathcal Y}}$-cluster webs in type $A$, again regarding their ARs and their rank(s).  For the ${\boldsymbol{\mathcal X}}$-cluster webs of type $A$, everything is known thanks to the results in \cite{Pereira}. We first briefly review  them, before focusing on the main objects we are going to study in this section, namely the webs ${\boldsymbol{\mathcal Y\mathcal W}}_{\hspace{-0.05cm}A_n}$. 
\sk 

In what follows, $n$ will stand for a  fixed integer bigger than or equal to 2.

 \subsection{\bf The ${\boldsymbol{\mathcal X}}$-cluster web of type $\boldsymbol{A_n}$.} 
\label{SS:X-cluster-web-An}
The cluster web 
${\boldsymbol{\mathcal X\hspace{-0.05cm}\mathcal W}}_{\hspace{-0.05cm}A_n}$ is a web in $n$ variables, of degree $d_{A_n}^{\mathcal X}$.  We have seen above that, up to pull-back under a map  \eqref{Eq:xT}, 
this web can (and will) be identified to the web 
${\boldsymbol{\mathcal W}}_{\hspace{-0.05cm} \mathcal M_{0,n+3}}$ on the moduli space $\mathcal M_{0,n+3}$ defined by all the forgetful maps $\mathcal M_{0,n+3}\rightarrow \mathcal M_{0,4}\simeq \mathbf P^1\setminus \{0,1,\infty\}$. Thus $d_{A_n}^{\mathcal X}={n+3  \choose 4}$ and 
according to the results of 
\cite[\S4.2]{Pereira} recalled above in \S\ref{Par:W-M0n+3-all-AMP}, the following assertions hold  true: 
\begin{itemize}
\item[$\boldsymbol{(i).}$] 
$ \rho^\bullet\Big({\boldsymbol{\mathcal X\mathcal W}}_{A_n}\Big)=\left(d^{\mathcal X}_{A_n}-n,  d^{\mathcal X}_{A_n}-{n+1 \choose 2},   d^{\mathcal X}_{A_n}- {n+2 \choose 3}  \right)
$; in particular  $\rho^\nu({\boldsymbol{\mathcal X\hspace{-0.05cm}\mathcal W}}_{\hspace{-0.05cm}A_n})=0$ for $\nu>3$; 
\sk 
\item[$\boldsymbol{(ii).}$]    ${\rm polrk}^\bullet\Big(
{\boldsymbol{\mathcal X\hspace{-0.05cm}\mathcal W}}_{\hspace{-0.05cm}A_n}
\Big)=\left( 
 2\cdot d^{\mathcal X}_{A_n}-{ n \choose 2 }-2n\, , \, 
 { n+2\choose 4}\, 
 \right) $; thus ${\rm polrk}^w\big(
 {\boldsymbol{\mathcal X\hspace{-0.05cm}\mathcal W}}_{\hspace{-0.05cm}A_n}
 \big)=0$ for $w\geq 3$;
\sk 
\item[$\boldsymbol{(iii).}$] 
  ${\boldsymbol{\mathcal X\hspace{-0.03cm}\mathcal W}}_{\hspace{-0.05cm}A_n}$ has AMP rank; moreover, all its ARs are polyogarithmic,  of weight 1 or 2;
\end{itemize}

We are going to use these results, in particular the point $\boldsymbol{(ii)}$,  to study the ${\boldsymbol{\mathcal Y}}$-cluster web
${\boldsymbol{\mathcal Y\hspace{-0.05cm}\mathcal W}}_{\hspace{-0.05cm}A_n}$, thanks to the fact that it is  a subweb of ${\boldsymbol{\mathcal X\hspace{-0.05cm}\mathcal W}}_{\hspace{-0.05cm}A_n}$

 \subsection{\bf The ${\boldsymbol{\mathcal Y}}$-cluster web of type $\boldsymbol{A_n}$.} 
\label{SS:Y-cluster-web-An}
 
 Our goal in this subsection is to prove the following result:
 \begin{thm}
 \label{T:YW-An}
 For any $n\geq 2$,  the six following assertions hold true:
\begin{enumerate}
\item  ${\boldsymbol{\mathcal Y\mathcal W}}_{A_n}$ is a $d_{A_n}^{\boldsymbol{\mathcal Y}}$-web in $n$ variables with $d_{A_n}^{\boldsymbol{\mathcal Y}}=n(n+3)/2$; 
\sk 
\item $ \rho^\bullet\Big({\boldsymbol{\mathcal Y\hspace{-0.05cm}\mathcal W}}_{\hspace{-0.05cm}A_n}\Big)=
\left( 
d_{A_n}^{\boldsymbol{\mathcal Y}}-n,n,1
\right)=\big(\, n(n+1)/2\, ,\, n\, ,\, 1\, \big) $; in particular  $\rho^\nu({\boldsymbol{\mathcal Y\hspace{-0.05cm}\mathcal W}}_{\hspace{-0.05cm}A_n})=0$ for $\nu>3$; 
\sk 
\item    ${\rm polrk}^\bullet\Big(
{\boldsymbol{\mathcal Y\hspace{-0.05cm}\mathcal W}}_{\hspace{-0.05cm}A_n}
\Big)=\left( d_{A_n}^{\boldsymbol{\mathcal Y}}\, , \, 
1
 \right) $; thus  ${\rm polrk}^w\big({\boldsymbol{\mathcal Y\hspace{-0.05cm}\mathcal W}}_{\hspace{-0.05cm}A_n} \big)=0$ for any weight $w\geq 3$;
\sk 
\item
  ${\boldsymbol{\mathcal Y\hspace{-0.05cm}\mathcal W}}_{\hspace{-0.05cm}A_n}$ is AMP; moreover, all its ARs are logarithmic, except the dilogarithmic one ${\sf R}_{A_n}$;
   \sk 
\item The  abelian relation  ${\sf R}_{A_n}$ is complete and accessible from Abel's five terms identity ${\sf R}_{A_2}$;
   \sk 
\item The web ${\boldsymbol{\mathcal Y\hspace{-0.03cm}\mathcal W}}_{\hspace{-0.02cm}A_n}$ is non linearizable.
\end{enumerate}
\end{thm}
 \begin{proof}
 That $d_{A_n}^{\boldsymbol{\mathcal Y}}=n(n+3)/2$ has been established above 
 (see Corollary \ref{Coro:d-X-Delta--d-Y-Delta}). 
  The second and third assertions are proved below, 
  {\it cf.}\,\eqref{Eq:rk=polrk=rho}. 
 That ${\boldsymbol{\mathcal Y\hspace{-0.05cm}\mathcal W}}_{\hspace{-0.05cm}A_n}$ is AMP follows immediately.  
The fifth point is a consequence of one of the main results of \cite{Souderes} (namely Th\'eor\`eme 4.1 therein) combined with the explicit description of the pull-back of ${\boldsymbol{\mathcal Y\hspace{-0.05cm}\mathcal W}}_{\hspace{-0.05cm}A_n}$ on $\mathcal M_{0,n+3}$ under the map 
\eqref{Eq:xT} for  the zig-zag triangulation $T_0$ (see also the third point of Corollary \ref{Cor:X-var-web-An}).  Finally, the sixth point has been proved before in \S\ref{SS:Y-cluster-Web-non-Lin}.
 \end{proof}
  \sk

We recall some notation we are going to use below:  $\Delta$ stands for the root system of type $A_n$ for a fixed $n\geq 2$ and we use 
 $\alpha>0$ (resp.\,$\alpha\geq -1$) as a  shorthand   for $\alpha\in \Delta_{>0}$ (resp.\,for $\alpha\in \Delta_{\geq -1}$). 
 By $({\sf{ R}_{A_n}})$ we refer to any one of these two (equivalent) dilogarithmic identities: 
 \begin{align}
 \label{Eq:RA-An-¤6}
 \sum_{\alpha\geq -1 } {\sf R}\Big(Y[\alpha] \Big)= \frac{n(n+1)}{2}\cdot \frac{{}^{} \hspace{0.1cm}\pi^2}{6}  
 \qquad 
 \Longleftrightarrow \qquad 
 \sum_{i=1}^n {\sf R}\big(u_i\big)- 
 \sum_{\alpha>0 } {\sf R}\Big(1/{Y[\alpha]}\Big)= 0\, . 
  \end{align}

  \subsubsection{The virtual ranks of $\boldsymbol{{\mathcal Y\hspace{-0.05cm}\mathcal W}}_{\hspace{-0.06cm}\boldsymbol{A_n}}$.} 
  \label{SS:Virtual-rank-YW-An}
 We set $d=n(n+3)/2$ and we denote by $V$ the complex vector space $\mathbf C^n$. 
Keeping in mind that $n$ is fixed and to simplify the notation,  we set here $\rho^\sigma= \rho^\sigma({\boldsymbol{\mathcal Y\hspace{-0.05cm}\mathcal W}}_{\hspace{-0.05cm}A_n})$ for any $\sigma \geq 1$.
\mk 

Since $\boldsymbol{{\mathcal Y\hspace{-0.05cm}\mathcal W}}_{\hspace{-0.06cm}{A_n}}$ is a subweb of $\boldsymbol{{\mathcal X\hspace{-0.05cm}\mathcal W}}_{\hspace{-0.06cm}{A_n}}$, it follows 
immediately from point $\boldsymbol{(i)}$ of \S\ref{SS:X-cluster-web-An} that $\rho^\sigma=0$ for any $\sigma\geq 3$, hence  we just have to prove that  $\rho^1=n(n+1)/2$, $\rho^2=n$ and $\rho^3=1$.\mk 

Let $\zeta$ be a generic point of the initial cluster torus $\mathbf T^n_0$ and
for $i=1,\ldots,N$,  denote by $\ell_i$ the differential of $y_i$ at $\zeta$, viewed as a linear form on $V$.  \sk 

Clearly, the $\ell_i$'s span $V^\vee$ (actually, so do $\ell_1,\ldots,\ell_n$ already) and 
there are $d$ of them. So obviously, there are $d-n=n(n+1)/2$ linear relations between them, thus $\rho^1=n(n+1)/2$. \mk 

The cases of $\rho^2$ and $\rho^3$ cannot be settled so easily. Actually, we are going first to prove the following bounds 
(in \S\ref{Par:rho2} and \S\ref{Par:rho3} respectively): 
\begin{equation}
\label{Eq:rho1-rho2}
 \rho^2\leq n \quad\quad \mbox{ and }\quad \quad
\rho^3\leq 1\, .
\end{equation}
Combined with $\rho^1=n(n+1)/2$, it gives 
\begin{equation}
\label{Eq:tokou}
\rho\Big(\boldsymbol{{\mathcal Y\hspace{-0.05cm}\mathcal W}}_{\hspace{-0.06cm}{A_n}}\big)=\rho^1+\rho^2+\rho^3 \leq d+1\, .
\end{equation}

On the other hand, in Subsection \S\ref{SS:log-ARs-YW-An} just below, one establishes that 
\begin{equation}
\label{Eq:polrk1-polrk2}
{\rm rk}\big(\boldsymbol{{\mathcal Y\hspace{-0.05cm}\mathcal W}}_{\hspace{-0.06cm}{A_n}}\big)\geq 
{\rm polrk}\big(\boldsymbol{{\mathcal Y\hspace{-0.05cm}\mathcal W}}_{\hspace{-0.06cm}{A_n}}\big)\geq  
{\rm polrk}^1+{\rm polrk}^2=d+1\, . 
\end{equation}

Since $ {\rm rk}\big(\boldsymbol{{\mathcal Y\hspace{-0.05cm}\mathcal W}}_{\hspace{-0.06cm}{A_n}}\big)\leq
 {\rho}\big(\boldsymbol{{\mathcal Y\hspace{-0.05cm}\mathcal W}}_{\hspace{-0.06cm}{A_n}}\big)$, it follows from 
  \eqref{Eq:rho1-rho2} and 
 \eqref{Eq:polrk1-polrk2} together 
 that 
\begin{equation}
\label{Eq:rk=polrk=rho}
{\rm rk}\big(\boldsymbol{{\mathcal Y\hspace{-0.05cm}\mathcal W}}_{\hspace{-0.06cm}{A_n}}\big)
=
{\rm polrk}\big(\boldsymbol{{\mathcal Y\hspace{-0.05cm}\mathcal W}}_{\hspace{-0.06cm}{A_n}}\big)
=
 {\rho}\big(\boldsymbol{{\mathcal Y\hspace{-0.05cm}\mathcal W}}_{\hspace{-0.06cm}{A_n}}\big)= d+1\, 
 \end{equation}
which in  turn implies that all the majorations in \eqref{Eq:rho1-rho2}, \eqref{Eq:tokou} and 
\eqref{Eq:polrk1-polrk2} actually are equalities.

 \paragraph{The pull-back of $\boldsymbol{{\mathcal Y\hspace{-0.05cm}\mathcal W}}_{\hspace{-0.06cm}{A_n}}$ on $\mathcal M_{0,n+3}$.}
\label{}

Instead of working with $\boldsymbol{{\mathcal Y\hspace{-0.05cm}\mathcal W}}_{\hspace{-0.06cm}{A_n}}$, we are going to work with its pull-back on $\mathcal M_{0,n+3}$, 
 denoted by $\boldsymbol{{\mathcal Y\hspace{-0.02cm}\mathcal W}}_{\hspace{-0.06cm}\hspace{-0.06cm}{\mathcal M_{0,n+3}}}$. We will then take benefit of  the geometric description of the latter to prove the majorations   \eqref{Eq:rho1-rho2}.  \mk

First, combining Corollary \ref{Cor:X-var-web-An} above with the description in terms of projective configurations given in \cite{Volkov} (see formula (1.10) and \S2 in this paper) of the $Y$-variables $Y_{i}(t)$ (with $i$ ranging from 1 to $n$ and $t\in \mathbf Z$) of the $Y$-system of type $A_n$, we obtain the

\begin{prop} 
The pull-back 
$\boldsymbol{{\mathcal Y\hspace{-0.02cm}\mathcal W}}_{\hspace{-0.06cm}{\mathcal M_{0,n+3}}}$ 
of $\boldsymbol{{\mathcal Y\hspace{-0.05cm}\mathcal W}}_{\hspace{-0.06cm}{A_n}}$ under the map $x^{T_0}$ is the subweb of 
$\boldsymbol{\mathcal X\hspace{-0.05cm} {\mathcal W}}_{\hspace{-0.03cm}{\mathcal M_{0,n+3}}}$ admitting as first integrals the cross-ratios $r_{i-1,i,j,j+1}: {\mathcal M_{0,n+3}}\rightarrow \mathbb C\setminus\{0,1\}$ ({\it cf.}\,\eqref{Eq:r-ijkl}) for the cyclically ordered quadruples $(i-1,i,j, j+1)\in (\mathbf Z/(n+3)\mathbf Z)^4$ associated to  all pairs $(i,j)$ with $i=2,\ldots,n+2$ and  $j=i+1,\ldots,n+3$, 
with the exception of $(2,n+3)$.
\end{prop}


%
%

\begin{rem} \begin{enumerate}
\item  The cross-ratios  $r_{i-1,i,j,j+1}$  essentially coincide\footnote{More precisely,  one has $[i,j\, \lvert k , l]=r_{i,j,k,l}+1$
where $[i,j\, \lvert k , l]$ stands for the cross-ratio defined in \cite[\S2.1]{Brown} and $r_{i,j,k,l}$ the one   in \eqref{Eq:r-ijkl}.}
with some others considered in  {\rm \cite{Brown}}, where they are denoted by $[i-1,i\,\lvert \,j,j+1]$ and called `dihedral coordinates' (see \S2.2 therein).  
As Brown's results show, these appear as 
very convenient to compute  periods on $\overline{\mathcal M}_{0,n}$. 
\sk 
\item   In relation with dilogarithmic identities, the dihedral coordinates have 
alreadt been used in {\rm \cite{Souderes}} to give a constructive proof that the identity 
$(\boldsymbol{\sf R}_{A_n})$ is accessible from Abel's five terms equation for any $n\geq 2$.

\item 
Except in the bibliography, no mention of cluster algebras is made in {\rm \cite{Souderes}}. However, the combinatoric used in the proof of Th\'eor\`eme 4.1 therein bears strong similarities with the combinatoric used to describe 
$\boldsymbol{\mathcal X}$-cluster objects (clusters, cluster variables, etc.) for the  cluster algebra of finite Dynkin type $A_n$. 
For this reason, we may think it is possible to rephrase the whole fourth section 
of Soud\`eres' paper 
 in terms of cluster algebras, which could be very interesting in order  to prove in an  effective way 
 the accessibility of the cluster dilogarithmic identity $(\boldsymbol{\sf R}_{\Delta})$ from Abel's identity  
 $(\boldsymbol{\sf R}_{A_2})$  for any Dynkin diagram $\Delta$. 
 \end{enumerate}
\end{rem} 

Let  $\mathscr Y_n$ be the set of pairs  $(i,j)$ defined in the previous proposition: 
$$\mathscr Y_n=\Big\{ \hspace{0.1cm}(i,j) \in \mathbf N^2\hspace{0.1cm} \big\lvert \hspace{0.1cm} 2\leq i <j\leq n+3\, , 
\hspace{0.1cm} (i,j) \neq (2,n+3)
\hspace{0.1cm}
 \Big\}
  \, .$$  
For $(i,j)$ in this set, we  use $r_{i,j}$  as a shorthand for 
the cross-ratio $r_{i-1,i,j,j+1}$ (to simplify). 
Since $[i-1,i\,\lvert \,j,j+1]=r_{i,j}+1$ is also a primitive first integral of the foliation of  $\boldsymbol{{\mathcal Y\hspace{-0.02cm}\mathcal W}}_{\hspace{-0.06cm}{\mathcal M_{0,n+3}}}$
defined by $r_{i,j}$, we will sometimes use it instead of $r_{i,j }$ at some places
below.

\paragraph{}${}^{}$\hspace{-0.6cm}
\label{Par:Wa-Wa}
Transporting the problem on the moduli space $\mathcal M_{0,n+3}$ allows us to work geometrically and to use the same powerful (albeit elementary) techniques used in the proof of \cite[Proposition 4.1]{Pereira}. We first recall below some notations taken from this paper which are needed for our purpose. \mk

Let $\boldsymbol{\mathcal W}$ be a subweb of $\boldsymbol{{\mathcal X\hspace{-0.02cm}\mathcal W}}_{\hspace{-0.03cm}{\mathcal M_{0,n+3}}}$. Then for any 
$a\in \{1,\ldots,n+3\}$, we denote by 
\begin{itemize}
\item $\boldsymbol{\mathcal W}_a$ the  subweb  of 
$\boldsymbol{\mathcal W}$  the first integrals of which are the first integrals $r_{i,j,k,l}$ of $\boldsymbol{\mathcal W}$ with $a\in \{i,j,k,l\}$;
\item  $\boldsymbol{\mathcal W}^a$ the subweb of $\boldsymbol{\mathcal W}$ 
the first integrals of which are the first integrals $r_{i,j,k,l}$ of $\boldsymbol{\mathcal W}$ with $a\not \in \{i,j,k,l\}$.
\end{itemize}

Obviously, one has $\boldsymbol{\mathcal W}=\boldsymbol{\mathcal W}_a\sqcup \boldsymbol{\mathcal W}^a$ and one verifies easily that 
$\boldsymbol{\mathcal W}^a$ is  the pull-back of a (uniquely determined) subweb 
of $\boldsymbol{{\mathcal X\hspace{-0.02cm}\mathcal W}}_{\hspace{-0.05cm}{\mathcal M_{0,n+2}}}$, denoted by $\boldsymbol{W}^a$
 under the  forgetful map $f_a : {\mathcal M_{0,n+3}}\rightarrow {\mathcal M_{0,n+2}}$ 
corresponding to forgetting the $a$-th point: one has 
 \begin{equation}
 \label{Eq:Wa=fa*Wa}
 \boldsymbol{\mathcal W}^a=f_a^*\big( \boldsymbol{W}^a\big).
 \end{equation}

Let $I\subset \{1,\ldots,n+3\}$ be a subset of cardinality $n$ and set $V_I=\oplus_{i\in I} \mathbf C\simeq \mathbf C^n$. Denoting by $p,q,r$ the elements of $\{1,\ldots,n+3\}\setminus I$ labeled in increasing order, one sets $x_{I}(p)=0$, $x_{I}(q)=1$ and $x_{I}(r)=\infty$ (as elements of $\mathbf P^1$) and one defines a birational map  $V_I\dashrightarrow
{\mathcal M_{0,n+3}}$ by associating to $(x_i)_{i\in I}\in V_I$ the 
projective class of $(x_{I,s})_{s=1}^{n+3}\in (\mathbf P^1)^{n+3}$ where $x_{I,s}=x_s$ if $s\in I$ and $x_{I,s}=x_I(s)$  as defined above otherwise.  We obtain that way a system of affine coordinates $(x_i)_{i\in I}$ on ${\mathcal M_{0,n+3}}$. 
\mk 

We assume that $I$ as above is fixed, we write $V$ instead of $V_I$ to simplify  and consider the cross-ratios $r_{a,b,c,d}$ of $\boldsymbol{{\mathcal X\hspace{-0.02cm}\mathcal W}}_{\hspace{-0.03cm}{\mathcal M_{0,n+3}}}$ as  rational functions in the $x_i$'s for $i\in I$.  Let $N$ be the degree of the considered subweb $\boldsymbol{\mathcal W}$ and denote by $\ell_1,\ldots,\ell_N$ be the differentials (viewed as linear forms on $V$) of its  first integrals  at a generic point of ${\mathcal M_{0,n+3}}$ (birationally identified with $V$). 
\sk

We recall that $A^\sigma(\boldsymbol{\mathcal W})$ stands for 
the space of homogeneous  ARs of degree $\sigma$ of the linear web 
whose first integrals are the $\ell_s$  for $s=1,\ldots,N$ 
 ({\it cf.}\,\eqref{Eq:A-sigma(W)} above). 
 For $a\in I$, the partial derivative $\partial/\partial x_a$ acts on ${\rm Sym}(V^*)={\mathbf C}[x_i\, \lvert \, i\in I]$  and, for any $\sigma\geq 2$, induces a linear map 
 $A^\sigma(\boldsymbol{\mathcal W})\rightarrow A^{\sigma-1}(\boldsymbol{\mathcal W}_a)$ whose kernel can easily be seen to be $A^\sigma(\boldsymbol{\mathcal W}^a)$. In other terms, there is an exact sequence of complex vector spaces
 \begin{equation}
 \label{Eq:Asigma-exact-sequence}
0\rightarrow A^\sigma\big(\boldsymbol{\mathcal W}^a\big)
\longrightarrow 
A^\sigma\big(\boldsymbol{\mathcal W}\big)
\longrightarrow 
A^{\sigma-1}\big(\boldsymbol{\mathcal W}_a\big)
\end{equation}
 from which one deduces that the following inequality holds true
  \begin{equation}
  \label{Eq:ririri}
 \rho^\sigma\big(\boldsymbol{\mathcal W}\big)\leq 
  \rho^\sigma\big(\boldsymbol{\mathcal W}^a\big)+ \rho^{\sigma-1}\big(\boldsymbol{\mathcal W}_a\big)
 \end{equation}
 (note that  since  \eqref{Eq:Asigma-exact-sequence} is a priori not exact at right, 
 this majoration has no reason to be sharp). 
 \mk 
 
 \paragraph{}${}^{}$\hspace{-0.6cm} 
 \label{Par:rho2}
We are going to apply the material above to a certain subweb of $\boldsymbol{{\mathcal Y\hspace{-0.02cm}\mathcal W}}_{\hspace{-0.06cm}{\mathcal M_{0,n+3}}}$ then argue by  induction on $n$ to finally get  $\rho^2\leq n$. 
Let $\mathscr Y_n'$ be the subset of $\mathscr Y_n$ formed by the pairs $(i,j)$ not of the form $(i,n+3)$, {\it i.e.} 
$\mathscr Y_n'=\mathscr Y_n\setminus \big\{ (3,n+3),\ldots,(n+2,n+3)\big \}
$ and consider the associated web 
$$
\boldsymbol{{\mathcal Y\hspace{-0.03cm}\mathcal W}}'_{\hspace{-0.06cm}{\mathcal M_{0,n+3}}}=
\boldsymbol{{\mathcal W}}\Big( \, r_{i,j}\hspace{0.06cm} \lvert
\hspace{0.06cm}
 (i,j)\in \mathscr Y_n'\hspace{0.06cm}
\Big)\, .
$$
It is a $n(n+1)/2$-subweb of $\boldsymbol{{\mathcal Y\hspace{-0.02cm}\mathcal W}}_{\hspace{-0.06cm}{\mathcal M_{0,n+3}}}$, whose complement web $
 \boldsymbol{{\mathcal W}}\big( \, 
 [2,3\, \lvert \, n+3,1]\, , \, \ldots\, , \, [n+1,n+2\, \lvert \, n+3,1]
 \,\big)$  can easily be seen as being a coordinate $n$-web, that is equivalent to the web defined by the $n$ standard coordinates on $\mathbf C^n$. 
 
\begin{prop}
\label{P:rho2(W')=0}
 For any $n\geq 2$, one has 
 $\rho^2(\boldsymbol{{\mathcal Y\hspace{-0.03cm}\mathcal W}}'_{\hspace{-0.06cm}{\mathcal M_{0,n+3}}})=0$ and consequently $\rho^2(\boldsymbol{{\mathcal Y\hspace{-0.02cm}\mathcal W}}_{\hspace{-0.06cm}{\mathcal M_{0,n+3}}})\leq n$.
 \end{prop}
 \begin{proof}
 Since the complement of $\boldsymbol{{\mathcal Y\hspace{-0.03cm}\mathcal W}}'_{\hspace{-0.06cm}{\mathcal M_{0,n+3}}}$ in $\boldsymbol{{\mathcal Y\hspace{-0.02cm}\mathcal W}}_{\hspace{-0.06cm}{\mathcal M_{0,n+3}}}$ is a $n$-web, the second majoration for 
 $\rho^2(\boldsymbol{{\mathcal Y\hspace{-0.02cm}\mathcal W}}_{\hspace{-0.06cm}{\mathcal M_{0,n+3}}})$
 is an immediate consequence of  
 $\rho^2(\boldsymbol{{\mathcal Y\hspace{-0.03cm}\mathcal W}}'_{\hspace{-0.06cm}{\mathcal M_{0,n+3}}})=0$.  \mk

 The proof that the latter equality holds true goes by induction on $n$.   The initial case is obvious ($\boldsymbol{{\mathcal Y\hspace{-0.03cm}\mathcal W}}'_{\hspace{-0.06cm}{\mathcal M_{0,5}}}$ is a planar 3-web) hence 
we assume  that $n\geq 3$   and that  $\rho^2(\boldsymbol{{\mathcal Y\hspace{-0.03cm}\mathcal W}}'_{\hspace{-0.06cm}{\mathcal M_{0,n+2}}})=0$.   Let $I\subset \{1,\ldots,n+3\}$ be a subset of cardinality $n$ as above,
and let us assume moreover that 
it contains $n+3$.  \sk 

 We then apply the results of \S\ref{Par:Wa-Wa} to 
$\boldsymbol{\mathcal W}
=\boldsymbol{{\mathcal Y\hspace{-0.03cm}\mathcal W}}'_{\hspace{-0.06cm}{\mathcal M_{0,5}}}$ and $a=n+3$.  With these notations:
\begin{itemize}
\item $\boldsymbol{\mathcal W}^{n+3}$ is the web  defined by the cross-ratios 
$[i-1,i \, \lvert \, j,j+1]$ for $i=2,\ldots,n+1$ and $j=i+1,\ldots, n+2$. Hence it is the pull-back 
 of the web $\boldsymbol{{\mathcal Y\hspace{-0.03cm}\mathcal W}}'_{\hspace{-0.06cm}{\mathcal M_{0,n+2}}}$ 
under the 
$(n+3)$-forgetful map $f_{n+3}: \mathcal M_{0,n+3}\rightarrow \mathcal M_{0,n+2}$: one has $\boldsymbol{\mathcal W}^{n+3}=(f_{n+3})^*\big(\boldsymbol{{\mathcal Y\hspace{-0.03cm}\mathcal W}}'_{\hspace{-0.06cm}{\mathcal M_{0,n+2}}}\big)$;
\item  $\boldsymbol{\mathcal W}_{n+3}$ is the web  defined by the cross-ratios $[i-1,i\, \lvert \, n+2,n+3]$ for $i=2,\ldots,n+1$. It is easily seen to be equivalent to a coordinate $n$-web.
\end{itemize}

Since the virtual rank is invariant by pull-back, one has 
$\rho^\sigma(\boldsymbol{\mathcal W}^a)=\rho^\sigma(\boldsymbol{{\mathcal Y\hspace{-0.03cm}\mathcal W}}'_{\hspace{-0.06cm}{\mathcal M_{0,n+2}}})$ for any $\sigma>0$ and because  all the virtual ranks of a coordinate web are trivial, 
 \eqref{Eq:ririri} is written as follows in the case when $\sigma=2$: 
$$
\rho^2\big(\boldsymbol{\mathcal W}\big)\leq 
  \rho^2\big(
  \boldsymbol{{\mathcal Y\hspace{-0.03cm}\mathcal W}}'_{\hspace{-0.06cm}{\mathcal M_{0,n+2}}}
  \big)+ 0.
$$
From the induction hypothesis, it follows that $
\rho^2\big(
  \boldsymbol{{\mathcal Y\hspace{-0.03cm}\mathcal W}}'_{\hspace{-0.09cm}{\mathcal M_{0,n+3}}}
\big)=0$ which concludes the proof.
 \end{proof}

 \paragraph{}${}^{}$\hspace{-0.6cm} 
 \label{Par:rho3}
 We now turn to $\rho^3$, which we are going to prove to be 1 at most, 
using similar methods as in the preceding paragraph. 
 \mk 
 
 Setting $ \boldsymbol{{\mathcal W}}= \boldsymbol{{\mathcal Y\hspace{-0.03cm}\mathcal W}}_{\hspace{-0.09cm}{\mathcal M_{0,n+3}}}$, $a=n+3$ and $\sigma=3$ in \eqref{Eq:ririri}, we have 
 \begin{equation}
 \label{Eq:rho3-blabla}
 \rho^3\leq \rho^3\big( \boldsymbol{{\mathcal W}}^{n+3} \big) +
 \rho^2\big( \boldsymbol{{\mathcal W}}_{n+3} \big)\, .
 \end{equation}
The web $ \boldsymbol{{\mathcal W}}^{n+3} $ admits as first integrals the cross-ratios 
$r_{i,j}$ for any pairs $(i,j)\in \mathscr Y_n$ such that $2\leq i<j\leq n+1$ hence can be seen as the pull-back under the $(n+3)$-th forgetful map  of $\boldsymbol{{\mathcal Y\hspace{-0.03cm}\mathcal W}}'_{\hspace{-0.09cm}{\mathcal M_{0,n+2}}}$, whose third virtual rank is zero (since  its second virtual rank already vanishes according to  Proposition \ref{P:rho2(W')=0}): one has 
 \begin{equation}
 \label{Eq:rho3-bloblo}
\rho^3\big( \boldsymbol{{\mathcal W}}^{n+3} \big) =
\rho^3\Big( \boldsymbol{{\mathcal Y\hspace{-0.03cm}\mathcal W}}'_{\hspace{-0.09cm}{\mathcal M_{0,n+2}}} \Big)=0\, . 
 \end{equation}

 We now apply the techniques of paragraph \S\ref{Par:Wa-Wa} 
to $\boldsymbol{{W}}=\boldsymbol{{\mathcal W}}_{n+3}$, but now with $a=1$. 
Since it is the $2n$-web defined by the $r_{i,j}$ for  
pairs $(i,n+2)$ with  $i=2,\ldots,n+1$ and pairs 
$(i,n+3)$ with $i=3,\ldots,n+2$, it is not difficult to describe the two corresponding   webs: 
\begin{itemize}
\item $\boldsymbol{{W}}^1$ is the $(n-1)$-web defined by the cross-ratios 
$[i-1,i\, \lvert \, n+2,n+3]$ for $i=3,\ldots, n+1$;
\item $\boldsymbol{{W}}_1$  is the web with 
the following  cross-ratios as first integrals: $[1,2\, \lvert \,n+2,n+3]$ and $[i-1,i\, \lvert \,n+3,1]$ for $i=3,\ldots,n+3$. 
\end{itemize}
 Clearly, $\boldsymbol{{W}}^1$ is equivalent to a $(n-1)$-subweb of a coordinate $n$-web in dimension $n$, thus $\rho^2(\boldsymbol{{W}}^1)=0$. As for the $(n+1)$-web $\boldsymbol{{W}}_1$, its subweb defined by the $n$ cross-ratios $[i-1,i\, \lvert \,n+3,1]$ for $i=3,\ldots,n+3$ is a coordinates web as well. Thus $\rho^1(\boldsymbol{{W}}_1)=1$ and when $\sigma=2$, \eqref{Eq:ririri} is written as follows in the case under scrutiny: 
 $$
 \rho^2\big( \boldsymbol{{W}}
 \big)
 \leq \rho^2\big(\boldsymbol{{W}}^1 \big) +
 \rho^1\big( \boldsymbol{{W}}_1 \big)=0+1=1\, . 
 $$
Combined with \eqref{Eq:rho3-blabla} and 
 \eqref{Eq:rho3-bloblo}, this gives us the majoration 
 $\rho^3\leq 1$ 
 we were looking for.

 \subsubsection{The polylogarithmic ARs of $\boldsymbol{{\mathcal Y\hspace{-0.05cm}\mathcal W}}_{\hspace{-0.06cm}\boldsymbol{A_n}}$.} 
  \label{SS:log-ARs-YW-An}
Our main goal here is to prove that ${\rm polrk}^1(\boldsymbol{{\mathcal Y\hspace{-0.05cm}\mathcal W}}_{\hspace{-0.06cm}\boldsymbol{A_n}})$ is bigger than or equal to $ n(n+3)/2$.    Combined with the upper bounds obtained above for the virtual ranks $\rho^\sigma$ for $\sigma=1,2,3$,  it will show that actually the equality holds true. 
Our main tools to achieve this are some nice formulas of \cite{FZ}. We will also discuss the problem of describing an explicit basis of the space of $\{0,-1,\infty\}$-logarithmic ARs  of $\boldsymbol{{\mathcal Y\hspace{-0.05cm}\mathcal W}}_{\hspace{-0.06cm}\boldsymbol{A_n}}$, a space that we denote just by 
$\boldsymbol{Log\hspace{-0.05cm}A\hspace{-0.05cm}R}$  below. 
\sk

In this subsection, we work with holomorphic  functions $f$ on the 
positive orthant $(\mathbf R_{>0})^n$ (that is, with germs of 
functions a representative of which is a holomorphic function defined on an open neighborhood of $(\mathbf R_{>0})^n$), up to addition of a constant\footnote{Disregarding the constant is useful and standard when considering abelian relations.}. It is the same as working indiscriminately  
 with the differential $df$ of such a function $f$.
 \begin{center}
 $\star$
\end{center}
We assume that $n\geq 2$ is  fixed and we denote  by $\Delta$ the root system of type $A_n$. We recall that $\Delta_{>0}$ (resp.\,$\Delta_{\geq -1}$) stands for the subset of positive (resp.\,almost-positive) roots.  
To deal with the logarithmic ARs of $\boldsymbol{{\mathcal Y\hspace{-0.05cm}\mathcal W}}_{\hspace{-0.06cm}\boldsymbol{A_n}}$, we use the following model of this web, in terms of (some of) the $Y$-variables appearing in the $Y$-system $\boldsymbol{Y}_{A_n}$, but indexed by the almost-positive roots: 
$$
\boldsymbol{{\mathcal Y\hspace{-0.05cm}\mathcal W}}_{\hspace{-0.06cm}\boldsymbol{A_n}}=
\boldsymbol{{\mathcal W}}\Big(\, Y[\alpha]\hspace{0.1cm} \big\lvert \hspace{0.1cm} \alpha \in \Delta_{\geq -1}\hspace{0.1cm}\Big)\, .
$$
We see the $Y[\alpha]$'s as rational functions of the initial variables $Y[-\alpha_i]=u_i$, $i=1,\ldots,n$.  For $\alpha\in \Delta$, we  recall that $[\alpha,\alpha_i]\in \mathbf Z$ stands for its $i$-th coordinates  with respect to the set  $(\alpha_1,\ldots,\alpha_n)$ of simple roots (thus $\alpha=\sum_{i=1}^n [\alpha,\alpha_i]\cdot \alpha_i$), one denotes  the product of the $u_i^{[\alpha,\alpha_i]}$'s by $u^\alpha$ and one sets $\tilde u^\alpha=u^\alpha$ if $\alpha\in \Delta_{>0}$ and $\tilde u^\alpha=1$ otherwise.
\mk

In \cite{FZ}, Fomin and Zelevinsky have introduced 
the so-called `$F$-polynomials' $F[\alpha]$, indexed by $\alpha\in \Delta_{\geq -1}$ (with $F[-\alpha_i]=1$ for any $i$) which 
appear to be very useful for investigating the logarithmic relations  between the first 
integrals  $Y[\alpha]$'s of $\boldsymbol{{\mathcal Y\hspace{-0.05cm}\mathcal W}}_{\hspace{-0.06cm}\boldsymbol{A_n}}$, thanks to two nice formulas :

\begin{prop}[formulas (2.14) and (2.18) in \cite{FZ}] 
\label{P:???}
For any $\alpha\in \Delta_{\geq -1}$, one has
\begin{equation}
\label{Eq:Yalpha-Yalpha+1}
Y[\alpha]=\frac{\prod_{(\beta,d)\in \Psi(\alpha)} F[\beta]^d}{u^{\alpha}}
\qquad\quad \mbox{ and } \qquad \quad
1+Y[\alpha]=\frac{F[\tau_+\alpha]F[\tau_-\alpha]}{\tilde u^{\alpha}}\, 
\end{equation}
for a certain set $\Psi(\alpha)\subset \Delta_{\geq -1}\times \mathbf N_{>0}$ and two permutations $\tau_\pm$ of  $\Delta_{\geq -1}$ 
(defined in {\rm \cite[\S1]{FZ}}{\rm )}.\footnote{Actually, the monomial denominators of the formulas (2.14) and (2.18) in \cite{FZ} are more involved than those appearing in \eqref{Eq:Yalpha-Yalpha+1}  
 but this is because they concern the general ({\it i.e.}\,possibly multi-laced) case. For a simply laced Dynkin diagram such as $A_n$, these  formulas simplify and can be written as in \eqref{Eq:Yalpha-Yalpha+1}.}
\end{prop}

	We consider the following complex vector space of functions on $(\mathbf R_{>0})^n$:
$$
\boldsymbol{Log}\big(u ,  \boldsymbol{F}_{\hspace{-0.05cm}\Delta}\big)=
\Big\langle \hspace{0.1cm} {\rm Log}(u_i) \, \big\lvert \, i=1,\ldots,n\hspace{0.1cm} \Big\rangle \oplus  
\Big\langle 
\hspace{0.1cm}
{\rm Log}\big(F[\alpha]\big) \hspace{0.05cm} \big\lvert 
\hspace{0.05cm} \alpha\in \Delta_{>0}\hspace{0.05cm} 
\Big\rangle\, . 
$$

From \eqref{Eq:Yalpha-Yalpha+1}, it follows that 
${\rm Log}\big(Y[-\alpha_i]\big)={\rm Log}(u_i)$ and 
${\rm Log}\big(1+Y[-\alpha_i]\big)={\rm Log}\big(F[\alpha_i]\big)$ for $i=1,\ldots,n$ whereas for any 
 $\alpha\in \Delta_{>0}$, one has 
\begin{align*}
{\rm Log}
\big(Y[\alpha] \big)=& \, \sum_{(\beta,d)\in \Psi(\alpha)}  d\,{\rm Log}
\big(F[\beta] \big)-\sum_{i=1}^n [\alpha;\alpha_i]\cdot  {\rm Log}(u_i) \\
\mbox{ and }\quad  
{\rm Log} \big(1+Y[\alpha] \big)=& \,{\rm Log}
\big(F[\tau_+ \alpha] \big)+
{\rm Log}
\big(F[\tau_- \alpha] \big)
-\sum_{i=1}^n [\alpha;\alpha_i]\cdot  {\rm Log}(u_i)\, . 
\end{align*}

We thus have a well-defined $\mathbf C$-linear map $\big\langle \hspace{0.05cm}  {\rm Log}\big(Y[\alpha] \big) \hspace{0.05cm}, \hspace{0.05cm} {\rm Log}\big(1+Y[\alpha] \big) \hspace{0.05cm} \big\rangle \rightarrow \boldsymbol{Log}\big(u ,  \boldsymbol{F}_{\hspace{-0.05cm}\Delta}\big)$ for any $\alpha$,    
from which we deduce the following exact sequence of vector spaces: 
\begin{equation}
\label{Eq:ExactSequence}
0 \longrightarrow
\boldsymbol{Log\hspace{-0.05cm}A\hspace{-0.05cm}R}
\longrightarrow \oplus_{\alpha \in \Delta_{\geq -1}} 
\Big\langle \hspace{0.05cm}  {\rm Log}\big(Y[\alpha] \big) \hspace{0.05cm}, \hspace{0.05cm} {\rm Log}\big(1+Y[\alpha] \big) \hspace{0.05cm} \Big\rangle 
  \longrightarrow 
\boldsymbol{Log}\big(u ,  \boldsymbol{F}_{\hspace{-0.05cm}\Delta}\big)\, . 
\end{equation}
		Since 
		$\oplus_{\alpha \in \Delta_{\geq -1}} 
\big\langle \hspace{0.05cm}  {\rm Log}\big(Y[\alpha] \big) \hspace{0.05cm}, \hspace{0.05cm} {\rm Log}\big(1+Y[\alpha] \big) \hspace{0.05cm} \big\rangle $ 
is obviously of dimension $2d$ (since each space in the direct sum is 2-dimensional)  and because $\dim  \big(\boldsymbol{Log}\big(u ,  \boldsymbol{F}_{\hspace{-0.05cm}\Delta}\big)\big)\leq \lvert \Delta_{\geq -1}\lvert =d$, it comes that one necessarily has
\begin{equation}
\label{Eq:dim-LogAR}
{\rm polrk}^1\Big(\boldsymbol{{\mathcal Y\hspace{-0.05cm}\mathcal W}}_{\hspace{-0.06cm}\boldsymbol{A_n}}  \Big)=\dim(\boldsymbol{Log\hspace{-0.05cm}A\hspace{-0.05cm}R})\geq d=n(n+3)/2\, . 
\end{equation}

Since  $(\boldsymbol{\sf R}_{A_n})$ provides a non trivial dilogarithmic AR, one has 
also ${\rm polrk}^2\big(\boldsymbol{{\mathcal Y\hspace{-0.05cm}\mathcal W}}_{\hspace{-0.06cm}\boldsymbol{A_n}}  \big)\geq 1$ thus 
$$
{\rm rk}\Big(\boldsymbol{{\mathcal Y\hspace{-0.05cm}\mathcal W}}_{\hspace{-0.06cm}\boldsymbol{A_n}}  \Big)
\geq  {\rm polrk}^1\Big(\boldsymbol{{\mathcal Y\hspace{-0.05cm}\mathcal W}}_{\hspace{-0.06cm}\boldsymbol{A_n}}  \Big)+
{\rm polrk}^2\Big(\boldsymbol{{\mathcal Y\hspace{-0.05cm}\mathcal W}}_{\hspace{-0.06cm}\boldsymbol{A_n}}  \Big)\geq d+1\, . 
$$

As explained above, this gives us that the following equalities hold true for 
$\boldsymbol{{\mathcal Y\hspace{-0.05cm}\mathcal W}}_{\hspace{-0.06cm}\boldsymbol{A_n}}$:
\begin{align*}
\rho^1=  n(n+1)/2 \, , \qquad  \rho^2=n 
\, , \qquad 
\rho^3=1\, , \qquad 
{\rm polrk}^1=   n(n+3)/2 
 \quad \mbox{ and } \qquad 
{\rm polrk}^2=1\, .
\end{align*}

In particular, it follows from these equalities that $\boldsymbol{Log}\big(u ,  \boldsymbol{F}_{\hspace{-0.05cm}\Delta}\big)$ has dimension $d$, which is easy to establish directly, but also that \eqref{Eq:ExactSequence} is exact at right, 
what can be formulated as the fact that the linear map  $\oplus_{\alpha \in \Delta_{\geq -1}} 
\big\langle \hspace{0.05cm}  {\rm Log}\big(Y[\alpha] \big) \hspace{0.05cm}, \hspace{0.05cm} {\rm Log}\big(1+Y[\alpha] \big) \hspace{0.05cm} \big\rangle 
  \longrightarrow 
\boldsymbol{Log}\big(u ,  \boldsymbol{F}_{\hspace{-0.05cm}\Delta}\big)$ is surjective. \mk 

One interesting feature of the above  approach is that it makes sense when $\Delta$ stands for any Dynkin diagram (provided that one considers 
formulas (2.14) and (2.18) of \cite{FZ} instead of those of \eqref{Eq:Yalpha-Yalpha+1} when $\Delta$ is not simply-laced). In particular,  there exists an exact sequence similar to \eqref{Eq:ExactSequence} for any $\Delta$. 

We thus can ask the following questions: 
\begin{question} Given a Dynkin diagram $\Delta$ (not necessarily of type $A$ any more): 
\begin{enumerate}
\item is the space $\boldsymbol{Log}\big(u ,  \boldsymbol{F}_{\hspace{-0.05cm}\Delta}\big)$ of dimension $\lvert \Delta_{\geq -1}\lvert$?
\item is the exact sequence \eqref{Eq:ExactSequence} exact at right?
\end{enumerate}
\end{question}

We have proved above than the answers to both questions is affirmative in type $A$. We believe that this is the case in full generality. Let us sketch a proof that {\it 1.} indeed holds true in full generality. First, from the fact that $F[\alpha]$ is a positive polynomial with integer coefficients and constant term equal to 1, it comes that 
$\boldsymbol{Log}\big(u ,  \boldsymbol{F}_{\hspace{-0.05cm}\Delta}\big)=
\boldsymbol{Log}(u) \oplus\boldsymbol{Log}(\boldsymbol{F}_{\hspace{-0.05cm}\Delta})$ 
where 
$\boldsymbol{Log}(u)=\big\langle \hspace{0.1cm} {\rm Log}(u_i)  \hspace{0.1cm} \lvert  \hspace{0.1cm} i=1,\ldots,n\hspace{0.1cm} \big\rangle$
and 
$\boldsymbol{Log}(\boldsymbol{F}_{\hspace{-0.05cm}\Delta})=\big\langle 
\hspace{0.1cm}
{\rm Log}\big (F[\alpha]\big) \hspace{0.05cm} \lvert 
\hspace{0.05cm} \alpha\in \Delta_{>0}\hspace{0.05cm} 
\big\rangle$, hence the question is reduced to the one of proving that the latter space has dimension $\lvert \Delta_{>0}\lvert$. Taking for granted the result announced by Cao and Keller that the $F[\alpha]$'s are irreducible ({\it cf.}\,Theorem \ref{Thm:CK} above), it then suffices to prove that two polynomials $F[\alpha]$ and $F[\alpha']$ are distinct if $\alpha,\alpha'\in \Delta_{>0}$ are.   But this follows from the fact that  \eqref{Eq:Phi-Fpolynomials} is bijective (which in turn follows from  the Caldero-Chapoton formula and relies on the fact that the cluster algebra associated to $\Delta$ admits a categorification). \mk 

Regarding the second question, we think that this is true as well, but 
possibly not in a completely uniform manner. 
 Indeed, denoting by $h$ the Coxeter number of $\Delta$, 
some investigations lead us  to believe that the two following assertions are true: 
\begin{itemize}
\item {\it if $h$ is even, then 
$\big\langle \hspace{0.1cm}  {\rm Log}\big(1+Y[\alpha] \big) \hspace{0.1cm} \big\lvert \hspace{0.1cm}  \alpha\in \Delta_{>0}\hspace{0.1cm}  \big\rangle 
= 
\boldsymbol{Log}(\boldsymbol{F}_{\hspace{-0.05cm}\Delta})$ 
modulo $\boldsymbol{Log}(u)$
;}
 \item {\it in any case, one has 
$\big\langle \hspace{0.1cm} {\rm Log}\big(1+Y[\alpha] \big) 
\hspace{0.1cm} , \hspace{0.1cm} {\rm Log}(u_i)\hspace{0.1cm}
\hspace{0.1cm}  \big\lvert \hspace{0.1cm} \alpha\in \Delta_{>0}\, i=1,\ldots,n\hspace{0.1cm}  \big\rangle 
=
\boldsymbol{Log}\big(u ,  \boldsymbol{F}_{\hspace{-0.05cm}\Delta}\big)
 $.}
 \end{itemize}

What has led us to consider statements such as the two ones above is the 
problem of describing an explicit basis of the space of logarithmic ARs $\boldsymbol{Log\hspace{-0.05cm}A\hspace{-0.05cm}R}$ which is uniform in $n$ (possibly according to its parity as well).  We discuss this in the next subsection.

 \subsubsection{About bases of the space of logarithmic ARs of $\boldsymbol{{\mathcal Y\hspace{-0.05cm}\mathcal W}}_{\hspace{-0.06cm}\boldsymbol{A_n}}$.} 
  \label{SS:Bases-log-ARs-YW-An}
In the previous subsection, we have proved that 
${\rm polrk
}^1(\boldsymbol{{\mathcal Y\hspace{-0.05cm}\mathcal W}}_{\hspace{-0.06cm}\boldsymbol{A_n}})=n(n+3)/2$ but in a non-constructive way, which is not fully satisfying. Our primary approach was more ambitious:   it consisted in building an explicit basis of the corresponding space of ARs, a task that unfortunately could not be completed.  However, this leads us to identify two possible approaches to tackle this problem. We discuss both below, starting with the case of $A_2$ in both cases as a motivation. 

 \paragraph{Building logarithmic abelian relations via the `$\boldsymbol{Y}$-method'.}  
 \label{Par:LogAR-Y-method}
 We have seen above in \S\ref{SS:cluster-webs-type-A2}
 that the set  of  ARs associated to the multiplicative relations  
 \eqref{Eq:Xell} constitutes a basis of the space $\boldsymbol{Log\hspace{-0.05cm}A\hspace{-0.05cm}R}(A_2)$
  of logarithmic ARs of $\boldsymbol{{\mathcal Y\hspace{-0.05cm}\mathcal W}}_{\hspace{-0.06cm}\boldsymbol{A_2}}$ (see \eqref{Eq:AR-XWA2}). \mk

An immediate but interesting remark is that the relations  \eqref{Eq:Xell} are just the specialization, in the $A_2$-case under consideration, of the general identities   of the $Y$-system of type $A_n$
\begin{equation*}
\big(\boldsymbol{Y}_i(t)\big) \hspace{2cm} \qquad 
Y_i(t+1)Y_i(t-1)=\prod_{j\neq i}
  \Big(1+Y_j(t)\Big)^{-a_{ij}}
  \qquad \qquad \hspace{2cm} {}^{}
   \end{equation*}
for $i=1,\ldots,n$ 
 and $ t\in \mathbf Z$ (where $(a_{ij})_{i,j=1}^n$ stands for the corresponding Cartan matrix). To any identity $(\boldsymbol{Y}_i(t))$ 
is associated the followin logarithmic abelian relation: 
\begin{equation*}
\big(\boldsymbol{LogAR}_i(t)\big) \hspace{0.2cm} \qquad 
d{\rm Log}\Big(Y_i(t-1)\Big) +
 \sum_{j\neq i}a_{ij} d{\rm Log}\Big(1+Y_i(t)\Big) 
 +
d{\rm Log}\Big(Y_i(t+1)\Big) 
=0\, .  
  \qquad \qquad {}^{}
   \end{equation*}
This suggests to wonder whether or  not the  $\boldsymbol{LogAR}_i(t)$'s span the whole space of logarithmic ARs  of $\boldsymbol{{\mathcal Y\hspace{-0.05cm}\mathcal W}}_{\hspace{-0.06cm}\boldsymbol{A_n}}$,  as in the case when $n=2$.
 \sk

 We set the following notations: $\Delta=A_n$, $h=h(\Delta)$ (Coxeter number; here $h=n+1$), $I=\{1,\ldots,n\}$ and $\sigma_\Delta : i\mapsto \overline{\imath}$  denotes 
 the involution of $I$ involved in  the  half-periodicity property of the $Y$-system of type $\Delta$ (for $\Delta=A_n$, one has 
$\overline{\imath} =n+1-i$ for any $i\in I$). Then, denoting by 
$\epsilon: I\rightarrow \{\pm 1\}$  the map associating 1 or -1 to $i$, depending on whether the $i$-th node of 
the bipartite quiver $\vec{\Delta}$ 
is a source or  a sink, one sets
 $$
 \mathfrak I_\Delta=\left\{  \hspace{0.1cm} (i,t) \hspace{0.1cm}  \big\lvert \hspace{0.1cm} 
 \begin{tabular}{l}
$i\in I$, \, 
$t=0,\ldots, h+1$
\end{tabular} 
 \right\} 
 \qquad
  \mbox{ and }\qquad 
  \mathfrak I_\Delta^\pm=\left\{  \hspace{0.1cm}(i,t)\in  \mathfrak I_\Delta
  \hspace{0.1cm}  \lvert \hspace{0.1cm}  \epsilon(i)(-1)^{t}=\pm 1
  \hspace{0.1cm}\right\}\, , 
 $$
%
 and as a model for  $\boldsymbol{{\mathcal Y\hspace{-0.05cm}\mathcal W}}_{\hspace{-0.06cm}\boldsymbol{\Delta}}$ in this paragraph, one takes the one admitting as first integrals the $Y$-variables $Y_i(t)$ for $i$ and $t$ such that 
 $\epsilon(i)(-1)^t= -1$ (see \eqref{Eq:Y-web}
  above and (8.10) in \cite{FZIV}): 
$$
 \boldsymbol{{\mathcal Y\hspace{-0.05cm}\mathcal W}}_{\hspace{-0.06cm}\boldsymbol{\Delta}}=
  \boldsymbol{{\mathcal W}}
  \left( \hspace{0.1cm} Y_i(t) \hspace{0.1cm} \big\lvert \hspace{0.1cm}
(i,t) \in    \mathfrak I_\Delta^-
\hspace{0.1cm}   \right)
 $$
(with  the $u_i$'s defined by $u_i=Y_i(0)^{-{\epsilon(i)}}$ for any $i\in I$ as initial variables). \sk

 In the lines below, it is useful to set \sk
 
 ${}^{}$ \quad $-$   $Y_i(-1)=Y_{\overline{\imath}}(n+2)$ and $Y_i(n+3)=Y_{\overline{\imath}}(0)$,  
  for any $i\in I$;
 \sk 
 
 ${}^{}$ \quad $-$ $Y_0(t)=Y_{n+1}(t)=0$ for any $t\in \mathbf Z$.
  \sk\\
  (Note: the first convention is in accordance with the half-periodicity property in type $A_n$ (see Remark \ref{Rem:Half-Periods} above) and the second is introduced  to write the ARs $LogAR_i(t)$ in an uniform manner regarding $i$).\mk

%

 Now for $i$ and $t$ in the same ranges of indices but such that $\epsilon(i)(-1)^t= 1$, is associated to  the multiplicative identity $\big(\boldsymbol{Y}_i(t)\big)$  the following logarithmic abelian relation 
 \begin{equation}
\label{Eq:LogARi(t)}
 {\rm Log}\Big( Y_i(t-1)\Big) - {\rm Log}\Big( 1+Y_{i-1}(t)\Big)-{\rm Log}\Big( 1+Y_{i+1}(t)\Big)
 +
  {\rm Log}\Big( Y_i(t+1)\Big)
 =0
 \end{equation}
 for  $\boldsymbol{{\mathcal Y\hspace{-0.05cm}\mathcal W}}_{\hspace{-0.06cm}\boldsymbol{\Delta}}$, denoted by $LogAR_i(t)$. This abelian relation has 4 or 3 logarithmic terms, the latter occurring only in the `extremal cases'  when $i=1$ or  $i=n$.  We denote by  $\boldsymbol{\mathcal Y\hspace{-0.05cm}Log\hspace{-0.05cm}A\hspace{-0.05cm}R}(\Delta)$ the space of ARs spanned by the  $LogAR_i(t)$'s:
 \begin{equation}
\label{Eq:Y-LogAR(An)-LogAR(An)}
 \boldsymbol{\mathcal Y\hspace{-0.05cm}Log\hspace{-0.05cm}A\hspace{-0.05cm}R}(\Delta)=
 \left\langle \hspace{0.1cm}
LogAR_i(t) \hspace{0.1cm} \big\lvert \hspace{0.1cm}
\hspace{0.1cm}
(i,t)\in  \mathfrak I_\Delta^+ 
\hspace{0.1cm}
\right\rangle \subset 
\boldsymbol{Log\hspace{-0.05cm}A\hspace{-0.05cm}R}(\Delta)
\, .
\end{equation}
The main questions considered here are first  to know whether  this inclusion is strict or not, and second,  to determine the linear relations between the $LogAR_i(t)$'s and in particular the dimension of 
 $\boldsymbol{\mathcal Y\hspace{-0.05cm}Log\hspace{-0.05cm}A\hspace{-0.05cm}R}(A_n)$ (we recall that 
$
\boldsymbol{Log\hspace{-0.05cm}A\hspace{-0.05cm}R}(A_n)$ has dimension $d$, as proved above).

\begin{prop}
\label{Prop:LogARit-An-n-even}
When $n$ is even, the $LogAR_i(t)$'s for $(i,t)\in \mathfrak I_n^+$ 
 constitute a basis of the space of logarithmic ARs of 
 $\boldsymbol{{\mathcal Y\hspace{-0.05cm}\mathcal W}}_{\hspace{-0.06cm}\boldsymbol{A_n}}$ 
and consequently, 
the   inclusion  \eqref{Eq:Y-LogAR(An)-LogAR(An)} actually is an equality.
 \end{prop}
 The proof use the following elementary fact, the proof of which is left to the reader: 
 
 \begin{fact}
 \label{Fact:foco}
{\it Let $(a_s)_{s\in \mathbf Z}$ be a $k$-periodic sequence of complex numbers such that $a_s+{a}_{s+1}=0$ for any $s\in \mathbf Z $. If it is $k$-periodic with $k$ odd then all the $a_s$ vanish.} 
 \end{fact} 
\begin{proof}
Assume that $\sum_{(i,t)\in \mathfrak I_\Delta^+ } \lambda_{i,t} \, LogAR_i(t)=0$ as an abelian relation, for some complex constants $\lambda_{i,t} $.
 Considering the coefficients of ${\rm Log}(Y_i(t))$ for $(i,t)\in \mathfrak I_n^{-}$ in this abelian relation, we get that the following relations hold true: 
\begin{align}
\label{Eq:lambda-i-t}
\lambda_{i,t-1}+\lambda_{i,t+1}=&\, 0 \quad \mbox{ for } t\in \{1,\ldots,n+1\}\, \\
\mbox{ and }
\quad 
\lambda_{\overline{\imath},n+2}+\lambda_{i,1}=\lambda_{i,n+1}+\lambda_{\overline{\imath},0}=&\, 0 \quad \mbox{ (case } t\in \{ 0, n+2\} 
\mbox{)}\, .\nonumber
\end{align}
Then, we set $m=\lfloor n/2\rfloor\in \mathbf N^*$ and for $i\in \{1,\ldots,m\}$ such that $\epsilon(i)=1$ (which is equivalent to $i$ being odd), we define a $(n+3)$-periodic sequence of complex numbers  $(\ell_{i,s})_{s\in \mathbf Z}$ by setting 
$$
\ell_{i,s}=\begin{cases}
\hspace{0.1cm}
\lambda_{\, i\, ,\, 2s} \hspace{1.1cm} \mbox{ for }\, s=0,\ldots,m+1\,, \\
\hspace{0.1cm} \lambda_{\, \overline{\imath}\, ,\, 2s-n-3} \hspace{0.4cm} \mbox{ for }\, s=m+1,\ldots,n+2\,.
\end{cases}
$$
Because $n$ is even, the sequence $(\ell_{i,s})_{s\in \mathbf Z}$ is odd-periodic and from \eqref{Eq:lambda-i-t} it follows that the sums of two consecutive terms all vanish.  From the Fact above, it follows that all the terms $\lambda_{i,t}$ 
and $\lambda_{\overline{\imath},t}$
with 
$(i,t) \in \mathfrak I_\Delta^{+}$ and $(\overline{\imath},t)\in 
\mathfrak I_\Delta^{+}$ vanish. \mk 

For $i\in \{1,\ldots,m\}$ such that $\epsilon(i)=-1$ ({\it i.e.}\,$i$ is even), one proceeds similarly to obtain the vanishing of all the coefficients $\lambda_{i,t}$ and $\lambda_{\overline{\imath},t}$. This shows that the $LogAR_i(t)$'s for $(i,t)\in \mathfrak I_{n}^+$ form  a basis of   $\boldsymbol{\mathcal Y\hspace{-0.05cm}Log\hspace{-0.05cm}A\hspace{-0.05cm}R}(A_n)$. Therefore this space has dimension $\lvert \mathfrak I_{\Delta}^+\lvert=d={\rm polrk}^1(\boldsymbol{{\mathcal Y\hspace{-0.05cm}\mathcal W}}_{\hspace{-0.06cm}\boldsymbol{A_n}})$ which finishes the proof of the proposition. 
\end{proof}

Since the statement of Fact \ref{Fact:foco} obviously does not hold when the period $k$ is even, the above proof does not apply when $n$ is odd. Actually, in this case and contrarily to what 
one could (too naively) expect at first, it appears that  the inclusion 
 \eqref{Eq:Y-LogAR(An)-LogAR(An)}  is strict as the $LogAR_i(t)$'s do no longer  form a free family of abelian relations. 
Indeed, setting $n=2m+1$ with $m\geq 1$ and using the same 
 approach  as in the proof above, it is not difficult to determine the possible form  that a linear relation between the $LogAR_i(t)$'s can have.  
One sets
\begin{align*}
LogAR_i=&\, \sum_{\tau=0}^{m+1}(-1)^\tau \Big(  LogAR_i(2\tau)+(-1)^{m} 
 LogAR_{ 2n+2-i
 }(2\tau) \Big)\quad \mbox{ for }\hspace{0.1cm} i=1,\ldots,m\, ; \\
\mbox{ and}\quad 
LogAR_{m+1}=&\, \sum_{\tau=0}^{m+1}(-1)^\tau LogAR_{m+1}(2\tau)\, .
\end{align*}
 
 Then by a direct computation, one verifies that the following proposition holds true:
\begin{prop}
\label{Prop:}
For  $\Delta=A_{2m+1}$,  the $LogAR_i(t)$'s for $(i,t)\in \mathfrak I_{\Delta}^+$ 
satisfy the  linear relation: 
\begin{equation}
\sum_{j=0}^{\lfloor m/2\rfloor} (-1)^j LogAR_{2j+1}=0\,.
\end{equation}
Moreover, this is the unique linear relation between the $LogAR_i(t)$'s with $(i,t)\in \mathfrak I_{2m+1}^+$ and 
consequently, the space $ \boldsymbol{\mathcal Y\hspace{-0.05cm}Log\hspace{-0.05cm}A\hspace{-0.05cm}R}(A_{2m+1})$ they span  has corank 1 in 
$\boldsymbol{Log\hspace{-0.05cm}A\hspace{-0.05cm}R}(A_{2m+1})$.
 \end{prop}
 
 
 When $n$ is odd, this result tells us that one has to find a new logarithmic AR not belonging to $ \boldsymbol{\mathcal Y\hspace{-0.05cm}Log\hspace{-0.05cm}A\hspace{-0.05cm}R}(A_{n})$ in order to get a spanning family for the whole space of logarithmic ARs. We do think that this is not a difficult task but we let this for future investigation, possibly by others. 
 \begin{center}
 $\star$
 \end{center}
 
 Let us comment about the approach presented in this paragraph to build a basis of $\boldsymbol{Log\hspace{-0.05cm}A\hspace{-0.05cm}R}(A_{n})$.
\sk 

 First, it shows again, but from a new perspective, the interest of presenting a web algebraically, by means of a $Y$-system as here, or more generally, of a period of a cluster algebra. In the case under scrutiny here, the relations $\big(\boldsymbol{Y}_i(t)\big) $ for some pairs $(i,t)$, first serve to define the web $\boldsymbol{{\mathcal Y\hspace{-0.05cm}\mathcal W}}_{\hspace{-0.06cm}\boldsymbol{A_n}}$ itself, whereas other pairs correspond in a quite immediate way to simple logarithmic ARs, and sufficiently many of them   span a space which is almost equal to $\boldsymbol{Log\hspace{-0.05cm}A\hspace{-0.05cm}R}(A_{n})$. 
Second, this approach is not specific to the  type $A$ but also makes sense for an arbitrary Dynkin diagram, as very briefly discussed  at the beginning of the current paragraph. Our third point goes in the opposite direction, since it concerns the limitations of this approach: the dichotomy according the parity of the rank  already in type $A_n$, shows that this approach  does not allow 
to build a basis of the space of logarithmic ARs of a web 
$\boldsymbol{{\mathcal Y\hspace{-0.05cm}\mathcal W}}_{\hspace{-0.06cm}\boldsymbol{\Delta}}$ (see also just after Proposition \ref{P:cococo}, where the cases of $D_4,D_5$ and $E_6$ are briefly discussed).  \mk 

As interesting and powerful it may be to build logarithmic ARs, the `$Y$-method'  presented in this paragraph carries its own limitations 
 hence one can wonder whether  another more systematic approach  exists or not. We believe that it is indeed the case and this is the subject of the next paragraph. 
 
 \paragraph{Dilogarithmic generation of logarithmic abelian relations.}  
 \label{Par:Dilog-Generation-LogAR}
The idea behind the approach discussed below also comes from the consideration of the $A_2$-case.  It is well-known (at least by the author of these lines) that the dilogarithmic abelian relation ${\sf R}_{A_2}$ (which is the 
one naturally associated to 
the 5-terms identity denoted by $({\rm R})$ in \cite{PirioSelecta}) can be seen as the most important of the ARs of $\boldsymbol{{\mathcal Y\hspace{-0.05cm}\mathcal W}}_{\hspace{-0.06cm}\boldsymbol{A_2}}$ (which is equivalent to Bol's web) since  it generates its whole space of ARs 
by analytic continuation (see \cite{PirioSelecta} and in particular \S3.1.1 therein).   \mk 

One can wonder whether  something similar holds true in general or not. More precisely, for $\Delta=A_n$ with $n\geq 2$, the ARs of 
$\boldsymbol{{\mathcal Y\hspace{-0.05cm}\mathcal W}}_{\hspace{-0.06cm}\boldsymbol{\Delta}}$ extend holomorphically (as multivalued ARs) on the complement 
 $U_\Delta$ of the set of common leaves $\Sigma_{\Delta}^c
 \subset \mathbf P^n$.   Given a base point $\star$, in $(\mathbf R_{>0})^n$ say,  analytic continuation of the dilogarithmic abelian relation  ${\sf R}_{\Delta}$ along any 
 $\gamma\in \pi_1(U_\Delta,\star)$ gives rise to a new AR, denoted by  ${\sf R}_{A_n}^\gamma$, such that $L_{\Delta}^\gamma={\sf R}_{\Delta}^\gamma-{\sf R}_{\Delta}$ be (the germ at $\star$ of) a logarithmic AR of $\boldsymbol{{\mathcal Y\hspace{-0.05cm}\mathcal W}}_{\hspace{-0.06cm}\boldsymbol{\Delta}}$. 
Since the answer is affirmative when $n=2$, it is natural to ask the 

\begin{question}
Is the inclusion of vector spaces $\big\langle 
\hspace{0.1cm} L_{\Delta}^\gamma 
\hspace{0.1cm} \lvert 
\hspace{0.1cm} 
\gamma \in \pi_1(U_\Delta,\star)
\hspace{0.1cm}
\big\rangle
\subset 
\boldsymbol{Log\hspace{-0.05cm}A\hspace{-0.05cm}R}(\Delta)$ an equality?\sk

If not, what are the linear relations between the $L_{\Delta}^\gamma $'s and what  is $\dim \big(\langle 
\hspace{0.1cm} \mathcal L_{\Delta}^\gamma 
\hspace{0.1cm} \lvert 
\hspace{0.1cm} 
\gamma \in \pi_1(U_\Delta,\star)
\hspace{0.1cm}
\rangle\big)$?
\end{question}

Note that the question makes perfectly  sense when replacing $\Delta$ by any other Dynkin diagram.  We conjecture that the answer to this question is `yes' in full generality.  But an inherent difficulty to handle it is that it relies on monodromy arguments hence of a global nature. In particular, it is necessary to know more about $U_\Delta$ and in particular about  its fundamental group (a set of some generators for instance), which is not known in general (however, see 
 Proposition \ref{P:Xn(Un)=XXn} just below in type $A$). \mk

There is a more formal approach to construct logarithmic ARs from ${\sf R}_{\Delta}$, which  relies on differentiation. It consists in taking the total derivative of the identity $({\sf R}_{\Delta})$ and isolating its components. More precisely, since  ${\sf R}(u)=(1/2)\big( \int_{0}^u \log(1+s)ds/s-
 \int_{0}^u \log(s)ds/(1+s)
 \big)$, one has $2dR=\log(1+u)/u-\log(u)/(1+u)$ for any $u>0$ 
 from which  it follows that 
 $$
 2\,d{\sf R}\big(Y[\alpha]\big)= 2\cdot Y[\alpha]^*\big(d{\sf R}\big)={\rm Log}\Big(1+Y[\alpha]\Big) d 
 {\rm Log}\Big(Y[\alpha]\Big) -{\rm Log}\Big(Y[\alpha]\Big) d 
 {\rm Log}\Big(1+Y[\alpha]\Big)
 $$
 holds true for any first integral $Y[\alpha]$ of 
 $\boldsymbol{{\mathcal Y\hspace{-0.05cm}\mathcal W}}_{\hspace{-0.06cm}\boldsymbol{\Delta}}$. \mk 
 
 We set $\omega_{i}=\omega_{-\alpha_i}=du_i/u_i$ for $i=1,\ldots,n$ and $\omega_\alpha= d{\rm Log}(F[\alpha])=dF[\alpha]/F[\alpha]$ for any $\alpha\in \Delta_{>0}$.  With this notation, it follows from 
the relations \eqref{Eq:Yalpha-Yalpha+1} that  for $i=1,\ldots,n$ and $\alpha\in \Delta_{>0}$, one has  
  \begin{align*}
  d{\rm Log}\Big(Y[-\alpha_i]\Big)=& \hspace{0.2cm}
\omega_i\, , && 
d{\rm Log}\Big(1+Y[-\alpha_i]\Big)= \omega_{\alpha_i}
 \\
    d{\rm Log}\Big(Y[\alpha]\Big)=& \hspace{-0.2cm}
 \sum_{(\beta,q)\in \Psi(\alpha)} \hspace{-0.2cm} q\,  \omega_\beta 
 -\sum_{i=1}^n [\alpha,\alpha_i]\,  \omega_i
\, \quad  && 
   d{\rm Log}\Big( 1+Y[\alpha]\Big)=  \omega_{\alpha_-}+\omega_{\alpha_+}-\sum_{i=1}^n [\alpha,\alpha_i] \, \omega_i\, , 
  \end{align*}
 where for simplicity $\alpha_\pm$ stands for the root denoted by $\tau_\pm(\alpha)$  in \cite[\S1]{FZ}.
 From the formulas above, it follows that for any $\alpha\in \Delta_{\geq -1}$, there exists an expression 
$$
2\cdot Y[\alpha]^*\big(d{\sf R}\big)=\sum_{\beta \in \Delta_{\geq -1}} \mathcal L_{\alpha,\beta}\hspace{0.1cm} \omega_{\beta}
$$
 where  for any $\beta$, $ \mathcal L_{\alpha,\beta}$ stands for 
 a function of the form 
 \begin{equation}
 {}^{} \qquad 
  \mathcal L_{\alpha,\beta}=\mathcal L _{\alpha,\beta}^0\cdot {\rm Log}\big( Y[\alpha]\big)
  + \mathcal L _{\alpha,\beta}^1 \cdot  {\rm Log}\big( 1+Y[\alpha]\big)
  \qquad \mbox{ with }\hspace{0.2cm} \mathcal L_{\alpha,\beta}^0,\mathcal L_{\alpha,\beta}^1\in \mathbf Z\, . 
 \end{equation}
   Taking the total derivative of the identity \eqref{Eq:RA-An-¤6}, one gets 
 \begin{align}
 \label{Eq:dR-An}
 0= &\, \sum_{\alpha \in \Delta_{\geq -1}} 
2\,  d\big({\sf R}(Y[\alpha] )\big)
 = \sum_{\alpha ,\beta \in \Delta_{\geq -1}} 
 \mathcal L_{\alpha,\beta} \hspace{0.1cm} \omega_{\beta}  = 
 \sum_{\beta \in \Delta_{\geq -1}} \bigg[ \sum_{\alpha  \in \Delta_{\geq -1}} \mathcal L_{\alpha,\beta}\bigg] \hspace{0.1cm}  \omega_{\beta} \,.
 \end{align}

The proof of the following lemma is easy and left to the reader
\begin{lem}
The $\mathcal O((\mathbf R_{>0})^n)$-module generated by the $\omega_\beta$'s with $\beta\in \Delta_{\geq -1}$ is free: there is no non trivial identity $\sum_{\beta\geq -1} c_\beta\, \omega_\beta=0$ with $c_\beta\in \mathcal O((\mathbf R_{>0})^n)$ for every $\beta\in \Delta_{\geq -1}$. 
\end{lem}

Since the coefficients in 
\eqref{Eq:dR-An} 
all belong to the subspace 
$\boldsymbol{Log(u,F_{\Delta})}=\langle \, {\rm Log}(Y[\alpha])\, , \, 
{\rm Log}(1+Y[\alpha])\, \lvert \, \alpha\in\Delta_{\geq -1}\, \rangle$ of 
$\mathcal O((\mathbf R_{>0})^n)$, this lemma implies that 
for any $\beta\in \Delta_{\geq -1}$, 
the following relation holds true  identically  on $(\mathbf R_{>0})^n$: 
$$
\sum_{\alpha  \in \Delta_{\geq -1}} \mathcal L_{\alpha,\beta}
=
\sum_{\alpha  \in \Delta_{\geq -1}}  
\Big( \mathcal L _{\alpha,\beta}^0\cdot {\rm Log}\big( Y[\alpha]\big)
  + \mathcal L _{\alpha,\beta}^1 \cdot  {\rm Log}\big( 1+Y[\alpha]\big)
  \Big)
=0 \, . 
$$
This identity can be seen as a logarithmic abelian relation for 
 $\boldsymbol{{\mathcal Y\hspace{-0.05cm}\mathcal W}}_{\hspace{-0.06cm}\boldsymbol{\Delta}}$, denoted by $\mathcal L _{\Delta,\beta}$.  
 \begin{question}
 Is the inclusion of vector spaces $\big\langle 
\hspace{0.1cm} \mathcal L_{\Delta,\beta}
\hspace{0.1cm} \lvert 
\hspace{0.1cm} 
\beta \in \Delta_{\geq -1}
\hspace{0.1cm}
\big\rangle
\subset 
\boldsymbol{Log\hspace{-0.05cm}A\hspace{-0.05cm}R}(\Delta)$ an equality?\sk 

 If not, what are the linear relations between the $\mathcal L_{\Delta,\beta}$'s and what  is $\dim \big(\langle 
\hspace{0.1cm} \mathcal L_{\Delta,\beta}
\hspace{0.1cm} \lvert 
\hspace{0.1cm} 
\beta \in \Delta_{\geq -1}
\hspace{0.1cm}
\rangle\big)$?
\end{question} 
Note that since  the cardinality of $\Delta_{\geq -1}$ coincides with $d=
\dim \big(\boldsymbol{Log\hspace{-0.05cm}A\hspace{-0.05cm}R}(\Delta)\big)$, an affirmative answer to the first question would imply that 
 $\big\{
\hspace{0.1cm} \mathcal L_{\Delta,\beta}
\hspace{0.1cm} \lvert 
\hspace{0.1cm} 
\beta \in \Delta_{\geq -1}
\hspace{0.1cm}
\big\}$ is a basis of $\boldsymbol{Log\hspace{-0.05cm}A\hspace{-0.05cm}R}(\Delta)$. \mk 

\begin{exm} 
The cluster variables $u_1=x$, $u_2=({1+x+y})/({xy})$, $u_3=y$, $u_4=({1+x})/{y}$ and $u_5=({1+y})/{x}$ can be taken as first integrals 
for (a model of) 
$\boldsymbol{{\mathcal Y\hspace{-0.05cm}\mathcal W}}_{\hspace{-0.06cm}\boldsymbol{A_2}}$.  The associated $F$-polynomials are $x$, $y$, $1+x$, $1+y$ and $1+x+y$.  Writing $\omega_F$ for the logarithmic derivative $d{\rm Log}(F)=dF/F$ for any $F$-polynomials, It can be verified that the total derivative of the LHS of the (identically satisfied) identity $\sum_{i=1}^5 {\sf R}(u_i)=\pi^2$ can be written $\sum \mathcal L_F\cdot \omega_F$, the summation being over the set of $F$-polynomials, with 
\begin{align*}
\mathcal L_x= & \,{\rm Log}  \big( 1+u_{{1}} \big) 
 +{\rm Log}  \, \bigg( \frac{ u_{{2}}  } {  1+u_{{2}} } \bigg)
  +{\rm Log} \bigg( \frac{ u_{{5}}  } {  1+u_{{5}} } \bigg)
 \\
\mathcal L_y= &\, {\rm Log}  \bigg( \frac{ u_{{2}}  } {  1+u_{{2}} } \bigg)+{\rm Log}  \left( 1+u_{{3}} \right)+{\rm Log}  \bigg( \frac{ u_{{4}}  } {  1+u_{{4}} } \bigg) \\
\mathcal L_{1+x}= &
-{\rm Log}  \left( u_{{1}} \right) -{\rm Log}  \left( u_{{2}} \right) +{\rm Log}  \left( 1+u_{{4}} \right)  \\
\mathcal L_{1+y}= & -{\rm Log}  \left( u_{{2}} \right) -{\rm Log}  \left( u_{{3}} \right) +{\rm Log}  \left( 1+u_{{5}} \right) \\
\mbox{and } \quad \mathcal L_{1+x+y}= & \, {\rm Log}  \left( 1+u_{{2}} \right) -{\rm Log}  \left( u_{{4}} \right) -{\rm Log}  \left( u_{{5}} \right) \, .
\end{align*}
One verifies easily that these five ARs  form a basis of the space $\boldsymbol{Log\hspace{-0.05cm}A\hspace{-0.05cm}R}(A_2)$. 
\end{exm}

By direct computations, we have verified that it is indeed the case for $n\leq 8$ and we conjecture that it actually holds true for every $n\geq 2$. We discuss this in a more general setting  in  the next paragraph.


 \paragraph{Generalizations to the case of an arbitrary Dynkin diagram.}  
The approaches discussed in  the two preceding paragraphs both generalize, each in a straightforward way, when replacing $A_n$ by an arbitrary Dynkin diagram $\Delta$.  We discuss very briefly these two  generalizations below. 
\mk 

The $Y$-method of paragraph \S\ref{Par:LogAR-Y-method} applies as such in general.   Everything in it makes sense when $\Delta$ stands for an arbitrary Dynkin diagram.  The question whether  there are non trivial linear relations between the $LogAR_i(t)$'s with $(i,t)\in \mathfrak I_\Delta^+$ can be investigated using the same approach as in the proof of Proposition  \ref{Prop:LogARit-An-n-even}, by first taking into consideration the parts ${\rm Log}(Y_i(t-1))+
{\rm Log}(Y_i(t+1))$ of the $LogAR_i(t)$'s. From this, one can deduce rather precisely the general form that any relations between the $LogAR_i(t)$'s can take. \mk 

For any $i\in I$, let $\mathfrak I_\Delta^+(i)$ be the set of $t$'s such that  $(i,t)\in  \mathfrak I_\Delta^+$: if $\epsilon(i)=1$ (resp.\,$\epsilon(i)=-1$), it is the set of even (resp.\,odd) elements of $\{0,\ldots,h+1\}$.  One sets 
\begin{align*}
m_i=\big\lvert \mathfrak I_\Delta^+(i)\big\lvert 
\qquad \mbox{ and }\qquad 
logAR_{i}=
 \sum_{ \tau  \in \mathfrak I_\Delta^+(i) } (-1)^{\lfloor \tau/2 \rfloor}  LogAR_i(\tau )\, . 
\end{align*}
(Note that $m_i=h/2+1$ for all $i$ if $h$ is even, whereas 
$m_i\in \{ \lfloor h/2\rfloor+1,  \lfloor h/2\rfloor+2\}$ if $h$ is odd).\sk

Now let $J=J_\Delta\subset I$ be the set of representatives of $\sigma_{\Delta}$-orbits in $I$ such that $j\leq \overline{\jmath}$ for any $j\in J$. 
\begin{equation*}
LogAR_{j}= \begin{cases}\hspace{0.1cm} logAR_j \hspace{3cm} \mbox{if }\, j= \overline{\jmath}\\
\hspace{0.1cm}
logAR_j+(-1)^{m_j} logAR_{\overline{\jmath}}
\hspace{0.35cm} \mbox{if }\, j< \overline{\jmath}\, . 
\end{cases}
\end{equation*}
Finally, one defines $K=K_\Delta \subset J$ as the subset of the element $k$ of $J$ 
such that 
 an even number of terms $LogAR_l(t)$ are involved in the 
 sum  defining $ LogAR_{k}$ just above, {\it i.e.} 
$$
K=\Bigg\{ \hspace{0.2cm} k\in J\hspace{0.2cm} \Big\lvert \hspace{0.0cm}  
\begin{tabular}{c}
$k=\overline{k} \, \mbox{ and }\, m_k \, \mbox{ is even; or  }$\hspace{0.17cm} ${}^{}$\\
$k<\overline{k} \, \mbox{ and }\, m_k+m_{\overline{k}}\,  \mbox{ is even}$
\end{tabular}
\hspace{0.0cm}
\Bigg\}\, .
$$

Then using similar arguments as those used in the proof of Proposition  \ref{Prop:LogARit-An-n-even}, one can prove the 
\begin{prop}
\label{P:cococo}
{1.} Any linear relation between the  logarithmic abelian relations $LogAR_i(t)$'s with $(i,t)\in \mathfrak I_\Delta^+$ can be written as a linear combinations of the elements $
LogAR_{k}$ for $k\in K_\Delta$. \sk 

{2.} Consequently,  the dimension of 
$ \boldsymbol{\mathcal Y\hspace{-0.05cm}Log\hspace{-0.05cm}A\hspace{-0.05cm}R}(\Delta)=
\big\langle \, 
LogAR_i(t)\, \lvert \, (i,t)\in \mathfrak I_\Delta^+
\, 
\big\rangle 
$
is at least $d_{\Delta}^{\boldsymbol{\mathcal Y}}-\lvert K_\Delta\lvert$.
\end{prop}

The bound  $\dim\big(\boldsymbol{\mathcal Y\hspace{-0.05cm}Log\hspace{-0.05cm}A\hspace{-0.05cm}R}(\Delta)\leq d_{\Delta}-\lvert K_\Delta\lvert$ is not sharp, as the case of $\Delta=A_{n}$ with $n$ odd considered above already shows. This occurs for other types as well, as the three following cases,  which have been studied by means of explicit computations, show: 
\begin{itemize}
\item $\boldsymbol{D_4:}$ one has  $\sigma_{D_4}={\rm Id}$, $h(D_4)=6$, 
$m_i=4$ for any $i\in I=\{1,\ldots,4\}$, from which it comes that $K=J=I$ thus $\lvert K\lvert =4$. 
However one verifies that $\dim\big(\boldsymbol{\mathcal Y\hspace{-0.05cm}Log\hspace{-0.05cm}A\hspace{-0.05cm}R}(D_4)\big)=14=16-2$. In particular, in this case the  relations $
LogAR_{k}$ for $k=1,\ldots,4$ between the $LogAR_i(t)$'s are not linearly independant. 
\sk 
\item $\boldsymbol{D_5:}$  $\sigma_{D_5}$ is the transposition $(45)$, one has  $h(D_5)=8$, and 
$m_i=4$ for any $i\in \{1,\ldots,5\}$, from which it comes that $K$ is the singleton $\{4\}$.  It turns out that $LogAR_4$ corresponds to a genuine linear relation between the $LogAR_i(t)$ for $i\in \{4,5\}$ with $(i,t)\in \mathfrak I_{D_5}^+$ hence 
 one has  $\dim\big(\boldsymbol{\mathcal Y\hspace{-0.05cm}Log\hspace{-0.05cm}A\hspace{-0.05cm}R}(D_5)\big)=24=25-1$. \sk 
\item $\boldsymbol{E_6:}$ 
 $\sigma_{E_6}$ is the product of the two transpositions $(16)$ and $(25)$, 
 $h(E_6)=12$ and 
$m_i=7$ for every $i$. It follows that $K=\{1, 2\}$.  
But it turns out that there is actually no linear relation at all between the $LogAR_i(t)$'s which form a basis of the space of logarithmic ARs of 
$\boldsymbol{\mathcal Y\hspace{-0.05cm} \mathcal W}_{\hspace{-0.05cm} E_6}$: one has 
 $\dim\big(\boldsymbol{\mathcal Y\hspace{-0.05cm}Log\hspace{-0.05cm}A\hspace{-0.05cm}R}(E_6)\big)=42$. 
 \sk 
\end{itemize}

Thus the $Y$-method furnishes an important number of linearly independent logarithmic ARs, but not always a basis of $\boldsymbol{Log\hspace{-0.05cm}A\hspace{-0.05cm}R}(\Delta)$. It would be interesting to describe another family of elements of this space, whose union with $\{ LogAR_i(t)\, \lvert \, (i,t)\in \mathfrak I_{\Delta}^+\,\}$ forms a spanning set.   
%
%
%
%
%
%
\begin{center}
$\star$
\end{center}  

The generation of a basis of 
$\boldsymbol{Log\hspace{-0.05cm}A\hspace{-0.05cm}R}(\Delta)$ 
by means of the  dilogarithmic identity 
$({\sf R}_\Delta)$ discussed in \S\ref{Par:Dilog-Generation-LogAR}
allows to get closed formulas for the $\mathcal L_{\Delta,\beta}$, but in terms of the involutions $\tau_{\pm}: \Delta_{\geq -1}$ and of the sets $\Psi(\alpha)$'s for $\alpha\in \Delta_{\geq -1}$. But these formulas are not as explicit and simple as those of the $LogAR_i(t)$'s given by the `$Y$-method'. However, 
the construction of the $\mathcal L_{\Delta,\beta}$'s 
 is effective, in the sense that given $\Delta$, one can find the $\mathcal L_{\Delta,\beta}$ explicitly, a procedure that we have implemented on a computer algebra system. Playing with it leads us to state the following 

 
 
  
  

  
  
  
  

  \begin{conjecture} 
  For any Dynkin diagram $\Delta$,  
  the logarithmic abelian relations $\mathcal L_{\Delta,\beta}$ (with $\beta\in \Delta_{\geq -1}$) obtained from 
  $({\sf R}_{\Delta})$  by taking its total derivative, 
   form a basis of the space $\boldsymbol{Log\hspace{-0.05cm}A\hspace{-0.05cm}R}(\Delta)$.
\end{conjecture}
  
  By direct computations, we have verified that this conjecture is satisfied for all Dynkin diagrams of rank less than or equal to 8 (hence in particular for all 
  the exceptional cases).  \mk

  Finally, we mention a last natural question about the two families $\{ \, LogAR_i(t)  \hspace{0.1cm} \lvert   \hspace{0.1cm} (i,t)\in \mathfrak I_{\Delta}^+\}$ and $\{ \, \mathcal L_{\Delta,\beta}   \hspace{0.1cm} \lvert   \hspace{0.1cm}  \beta \in {\Delta}_{\geq -1}\}$ of logarithmic abelian relations of $\boldsymbol{{\mathcal Y\hspace{-0.05cm}\mathcal W}}_{\hspace{-0.06cm}\boldsymbol{\Delta}}$: 
  \begin{question}
{\it  How are these two families  related? In particular, assuming that the conjecture just above holds true, 
  how is each $LogAR_i(t)$ expressed in terms of the 
  basis
$\{ \hspace{0.1cm} \mathcal L_{\Delta,\beta}\hspace{0.1cm}\}_{\beta \in {\Delta}_{\geq -1}}$ of 
$\boldsymbol{Log\hspace{-0.05cm}A\hspace{-0.05cm}R}(\Delta)$?}
\end{question}


  \subsection{\bf The zig-zag map $\boldsymbol{x^{T_0}}$ and its inverse.}
  \label{SubSec:ZigZag-Map}
 Our goal here is to give an explicit formula for the birational map 
 $x^{T}: \mathcal M_{0,n+3}\dashrightarrow \boldsymbol{\mathcal X}_{A_n}$ associated to a certain zig-zag triangulation $T$ of the $(n+3)$-gon ${\bf P}_{n+3}$ and to get from this interesting results about the geometry of $\boldsymbol{\mathcal X}_{A_n}$. 
 \bk 
 
 
  In what follows, $n\geq 2$ stands for  an integer bigger than 1 and one sets:  
  \begin{equation} 
  m=\lfloor n/2 \rfloor\geq 1\, , \qquad m^*=m+1
  \qquad \mbox{ and } \qquad 
  \widetilde m=\lfloor (n-1)/2 \rfloor
   \, .
  \end{equation}
 One verifies that, independently of the parity of $n$,  $m^*$ and $\widetilde m$ are always such that 
 $m^*+\widetilde m=n$.
 

  \subsubsection{The zig-zag.}
    We consider the zig-zag triangulation $\boldsymbol{Z}_{n}=\boldsymbol{Z}_{A_n}$ of $\mathbf P_{n+3}$ starting by 
  $3\rightarrow 1 \rightarrow 4 \rightarrow n+3 \rightarrow \cdots  $ and ending as follows, depending on the parity of $n$:
\begin{align*}
  \cdots \longrightarrow m+2\longrightarrow m+5 \longrightarrow m+3&  \hspace{0.4cm} \mbox{for } \, n\, \mbox{ even};
\\
   \cdots \longrightarrow m+6\longrightarrow m+3\longrightarrow m+5 &\hspace{0.4cm} \mbox{for } \, n\, \mbox{ odd}.
\end{align*}  
This triangulation is pictured below in Figure  \ref{Fig:zzag}.
\begin{figure}[h]
\begin{center}
\scalebox{0.5}{
 \includegraphics{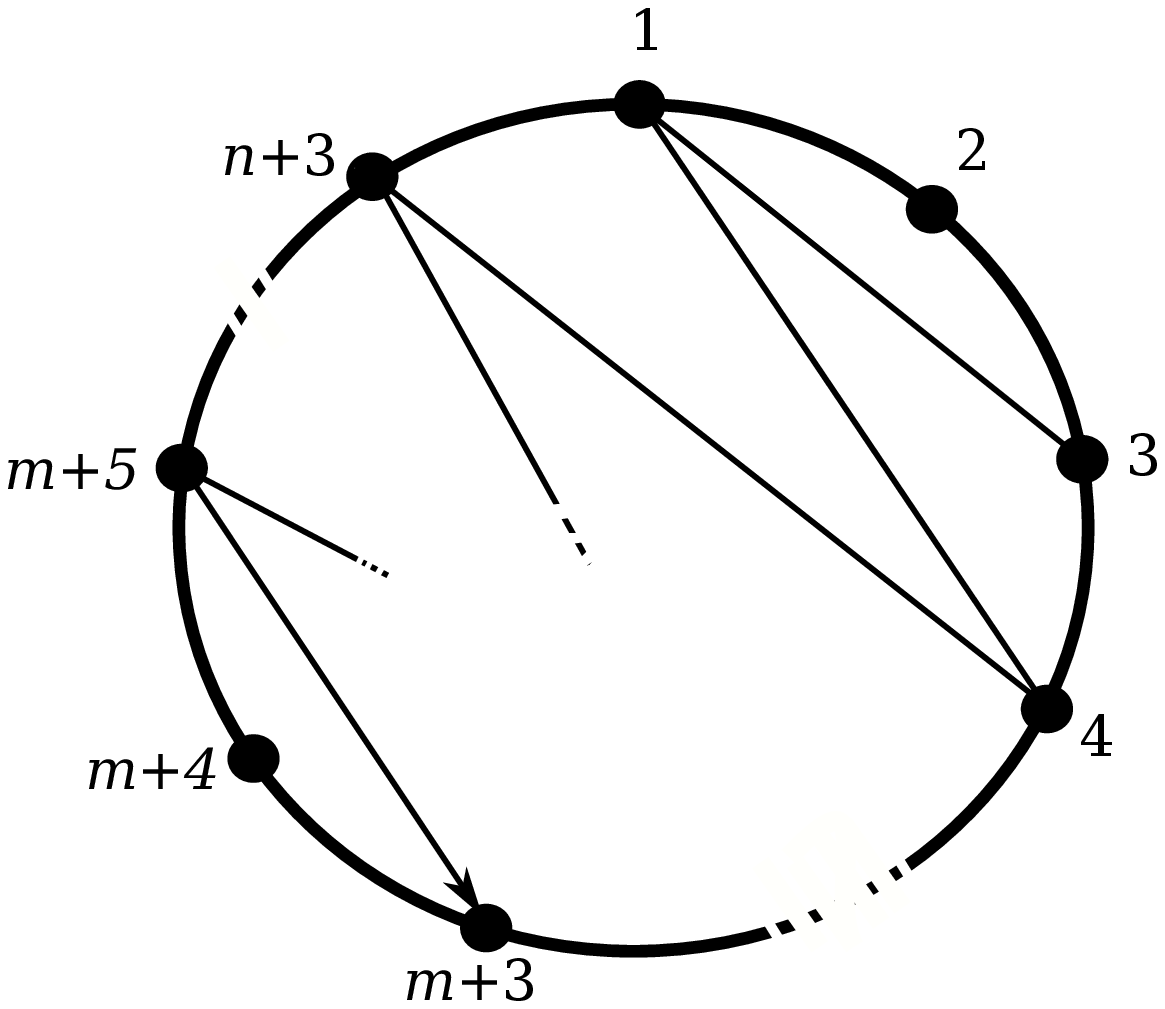}}
 \qquad \qquad \scalebox{0.5}{
 \includegraphics{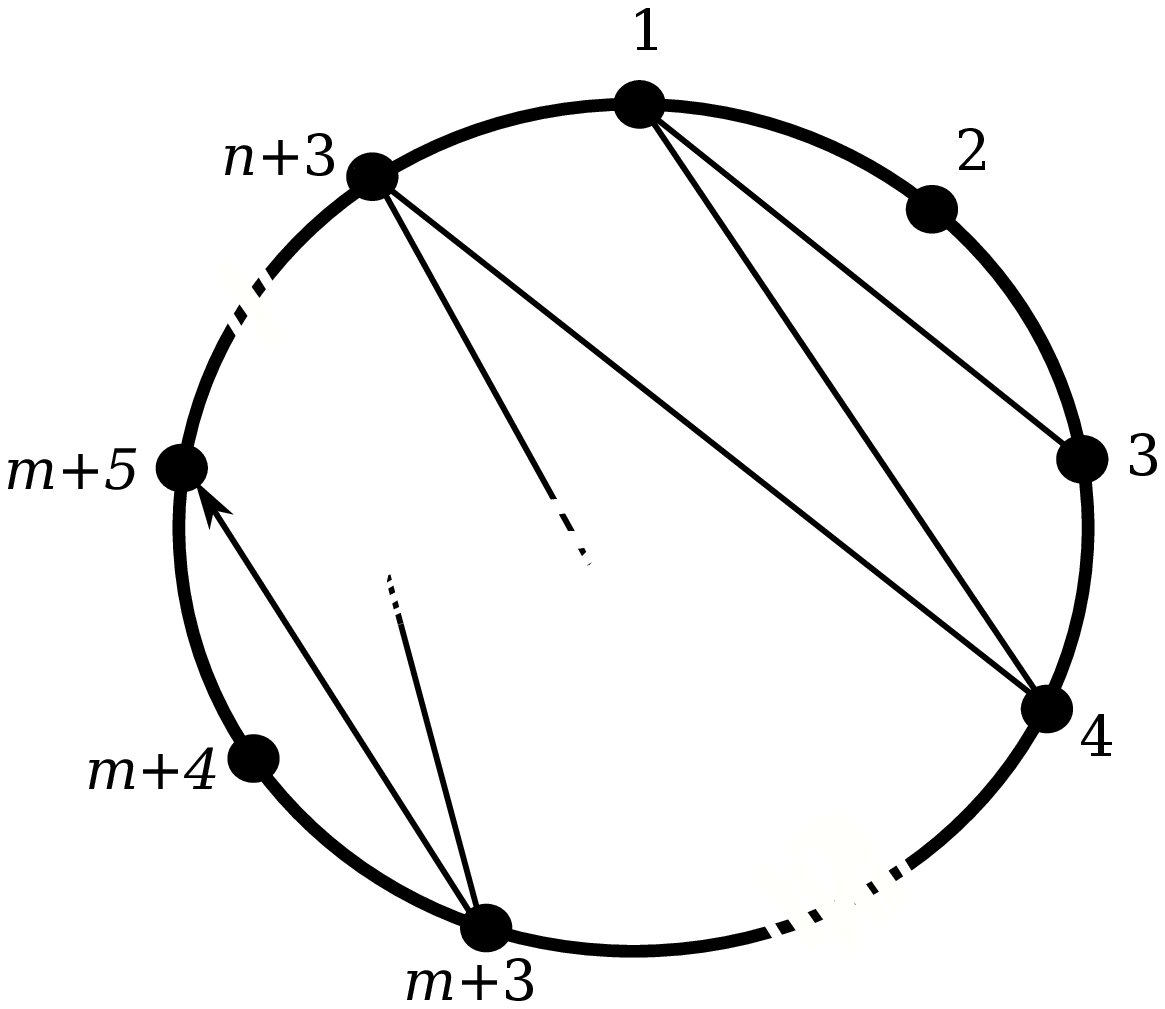}}
 \vspace{0.3cm}
\caption{The zig-zag $\boldsymbol{Z}_{n}$ for $n$ even (on the left) and $n$ odd (on the right)} 
\label{Fig:zzag}
\end{center}
\end{figure}

 In order to get simpler formulas for the $F$-polynomials we will deal with below, we will not work exactly with the cluster variables associated to the triangulation $\boldsymbol{Z}_{n}$ of $\mathbf P_{n+3}$ but will invert the ones associated to the sources of the quiver associated to $\boldsymbol{Z}_{n}$ ({\it cf.}\,Remark \ref{Rem:birat-E} above).  To this end, we consider 
 the birational involution $E=E_n$ characterized by $E^*(v_i)=v_i$ when $i$ is even and $E^*(v_i)=1/v_i$ when $i$ is odd ({\it cf.}\,\eqref{Eq:Bir-change-E}). 

 \subsubsection{The map $\boldsymbol{U_n}$.}
 \label{SS:map-Un}
To deal with a concrete explicit model for the birational map $x^{\boldsymbol{Z}_{n} }: 
 \mathcal M_{0,n+3} \dashrightarrow \mathcal X _{ A_n}$ associated to the zig-zag $\boldsymbol{Z}_{n}$ (defined in \eqref{Eq:xT}), we consider the birational parametrization $ \mu: 
 \mathbf C^n\dashrightarrow  \mathcal M_{0,n+3}, (x_i)_{i=1}^n
\longmapsto   [\infty,0,1,x_1,\ldots,x_n]$ and consider the affine birational map 
$$\boldsymbol{{\sf U}}_n= E\circ x^{\boldsymbol{Z}_{n} } \circ \mu : \mathbf C^n\dashrightarrow \mathbf C^n\, .$$

To make things as clear as possible, we denote by $\mathbf C_x^n$ the source affine space (with coordinates $x_1,\ldots,x_n$) and by $\mathbf C_u^n$ the target one (with coordinates $u_1,\ldots,u_n$). 
  Then it is not difficult to get explicit expressions for the components ${\sf U}_n^k$ of  $\boldsymbol{{\sf U}}_n= \big({\sf u}_n^k\big)_{k=1}^n : \mathbf C_x^n\dashrightarrow \mathbf C_u^n$: one has  
  \begin{align}  \nonumber
{\sf u}_n^1= & \, 
   \frac{1}{x_1-1} \, , \qquad 
   {\sf u}_n^2= 
   \frac{x_1-x_n}{1-x_1} \, ,  \qquad 
   {\sf u}_n^3=
   \frac{x_2-x_{n}}{x_1-x_2} \, ,  \qquad  
     {\sf u}_n^4=
   \frac{(x_2-x_{n-1})(x_1-x_{n})}{(x_{n-1}-x_n)(x_1-x_2)} \, ,  
\vspace{0.15cm}    \\ 
   {\sf u}_n^5= & \,  
   \frac{(x_2-x_{n})(x_3-x_{n-1})}{(x_{n-1}-x_n)(x_2-x_3)} \, ,  \, 
 \qquad   \ldots
     \ldots  \ldots
      \hspace{1.4cm}   
  {\sf u}_n^{n}=  \frac{(x_{\widetilde m+1}-x_{\widetilde m+2})(x_{\widetilde m}-x_{\widetilde m+3})}{
  (x_{\widetilde m}-x_{\widetilde m+1})(x_{\widetilde m+2}-x_{\widetilde m+3})}
  \label{F:unn}
       \end{align}

Our goal here is to study the birational map $\boldsymbol{\sf U}_n$ and its converse. 
Aiming at arguing inductively, we consider the two following linear projections (defined for any $n\geq 3$): 
\begin{itemize}
\item $p_n: \mathbf C^n_u\longrightarrow \mathbf C^{n-1}_u$ stands for the  projection given by forgetting the $n$-th coordinate: one has 
$$
p_n(u)=(u_1,\ldots,u_{n-1})  \hspace{3.15cm}  \mbox{ for any }\, u=(u_i)_{i=1}^n\in \mathbf C^n_u\, ; 
$$
\item 
 $\pi_n : \mathbf C^n_{x}  \longrightarrow  \mathbf C^{n-1}_{x}$ denotes the  projection given by forgetting the $m^*$-th coordinate, {\it i.e.}
 \begin{align*}
\pi_n(x)= \big( \, x_1,\ldots,x_{m}\, ,\,  \widehat{x_{m^*}}\, , \, 
  x_{m+2},\ldots,x_n\, 
  \big)
 \quad \mbox{ for any }\, x=(x_i)_{i=1}^n\in \mathbf C^n_x\, .
  \end{align*}
\end{itemize}

The basic result which will allow to study  $\boldsymbol{\sf U}_n$ inductively is 
  the following 
\begin{lem}
\label{L:Un-Un-1}
For any $n\geq 3$, as rational maps from $\mathbf C_x^{n}$ to $\mathbf C_u^{n-1}$, one has 
\begin{equation}
\label{Eq:Un-Un-1}
p_n\circ \boldsymbol{{\sf U}}_n= \boldsymbol{{\sf U}}_{n-1}\circ \ \pi_n
\, .
\end{equation}
\end{lem}
\begin{sproof}
One sets $\chi=(\chi_s)_{s=1}^{n+3}=(\infty,0,1,x_1,\ldots,x_n)$ (thus $\chi_{ i+3}=x_i$ for $i=1,\ldots,n$). 
Denote by $\boldsymbol{Z}_{n}^x$, the zig-zag $\boldsymbol{Z}_{n}$ decorated by attaching  $\chi_s$ to its $s$-th vertex for  $s=1,\ldots,n+3$.  Given $x=(x_i)_{ i=1}^n$, one sets for simplicity $x^\circ=\pi_n(x)=(x_1,\ldots,x_m, \widehat{x_{m^*}}, x_{m+2}, \ldots, x_n)$. Then the lemma follows easily by noticing that,  $\boldsymbol{Z}_{n-1}^{x^\circ}$ is obtained from $\boldsymbol{Z}_{n}^{x}$ by first removing the last arrow of the zig-zag (from $m+5$ to $m+3$ when $n$ is even, from   $m+3$ to $m+5$ when it is odd) and its $(m+4)$-th vertex (decorated by $x_{m^*}$) then re-labeling the remaining $n+2$ vertices cyclically (looking at Figure \ref{Fig:zzag} should make everything cristal clear). 
\end{sproof}

In view of given a closed formula for the converse map to $\boldsymbol{\sf U}_n
$, we need first to discuss some nice properties of $F$-polynomials in the cluster type we are dealing with.

\subsubsection{Some properties of $\boldsymbol{F}$-polynomials in type $\boldsymbol{A}$.}
\label{SS:SomeProp-F-polynomials}
 Below, we deal with rational functions in the $u_i$'s.  
 Being aware of this, we will sometimes denote an element $R(u)$ of $\mathbf C(u)$  just by $R$ in order to simplify the writing. 
 We will also use the following notations:  given any $n$-tuple
   $u=(u_1,\ldots,u_n)$, one sets 
$u'=(   u_1,\ldots,u_{n-1})$, $u''=(u')'=(   u_1,\ldots,u_{n-2})$;  and $'\hspace{-0.06cm}u=(   u_2,u_3,\ldots,u_{n})$, $''\hspace{-0.06cm}u=\,'\hspace{-0.04cm}('\hspace{-0.06cm}u)=( u_3, \ldots,u_{n})$, etc. 
(with the convention that a $k$-tuple is the empty set for any $k\leq 0$). \mk 

%
%

   One defines recursively  polynomials  $F_n=F_n(u)\in \mathbf Z_{>0}[u_1,\ldots,u_n]$  for every $n\in \mathbf N$ by setting  
   \begin{equation}
   \label{Eq:Fn-induction}
F_0=1,\quad  F_1=1+u_1 \qquad \mbox{ and } \qquad 
   F_n(u)=F_{n-1}(u')+u_n F_{n-2}(u'') \quad \mbox{ for any }\, n\geq 2\, .
   \end{equation}
   For instance, one has  $F_2=1+u_1+u_2$, $F_3=1+u_1+u_2+u_3+u_1u_3$, etc.  
   It also  follows immediately from the recurrence relations above that the non zero coefficients of any $F_n$ are positive integers.

   \begin{prop}
   \label{P:Prop-F-poly-An}
   \begin{enumerate}
   \item For any $n$ and any $u=(u_i)_{i=1}^n\in \mathbf C^n$, one has $F_n(u)=F_n(u_n,u_{n-1},\ldots,u_1)$.
   \item Consequently,  for $n\geq 2$, the identity  $F_n(u)=F_{n-1}('\hspace{-0.06cm}u)+u_1 
   F_{n-1}(''\hspace{-0.06cm}u)$ is satisfied as well. 
   \sk 
   \item  Let $a,b\in \mathbf N$ be such that $1\leq a \leq b \leq n$ and set $m=b+1-a\geq 1$. The $F$-polynomial  $F_{[a,b]}$ associated to the root $\alpha_{[a,b]}=\sum_{s} \alpha_s$ of the root system of type $A_n$ coincides with  $F_m$ evaluated on the $m$-tuple 
   $u_{[a,b]}=(u_a,u_{a+1},\ldots,u_b)$: as elements of $\mathbf Z[u_1,\ldots,u_n]$,  one has
     \begin{equation}
     F_{[a,b]}=F_m\big(u_{[a,b]}\big)\, .
     \end{equation}
\end{enumerate}
   \end{prop}
   \begin{proof}
For $(a,b)\in \mathbf N^2$ with $1\leq a\leq b$, let $\Omega_{[a,b]}$ be the set of totally disconnected subsets of $\big[\hspace{-0.07cm}\lvert a,b \lvert\hspace{-0.07cm}\big]$: elements of $\Omega_{[a,b]}$ can be written as increasing tuples $I=(i_1,\ldots,i_\ell)$ of integers $i_s\in  \big[\hspace{-0.07cm}\lvert a,b \lvert\hspace{-0.07cm}\big]$ such that no $i_s$ be adjacent to another element of $I$, {\it i.e.} $i_{s+1}> i_s+1$ for every $s=1,\ldots,\ell-1$. 
   We also set $\Omega_{\emptyset}=\{ \emptyset \}$.  For 
  such a tuple $I$, one  defines a monomial in the $u_i$'s (for $i=1,\ldots,n$), by setting  $u^I=\prod_{s=1}^\ell u_{i_s}$ (with $u^\emptyset=1$) 
   and summing  over all the elements of $\Omega_{[a,b]}$, one gets a polynomial with positive integer coefficients: 
   \begin{equation}
   \label{Eq:F[a,b]}
   {\sf  F}_{[a,b]}(u)=\sum_{I \in \Omega_{[a,b]}} u^I\, .
   \end{equation}
   
Obviously, one has  ${\sf  F}_{\emptyset}(u)=1$, ${\sf  F}_{[1]}(u)=1+u_1$ ($\Omega_{[1]}$ is the union of the empty set with the singleton 
   $\{1\}$) and for any $n\geq 2$, $\Omega_{[1,n]}$ is the disjoint union of 
   $\Omega_{[1,n-1]}$ with the set of concatenations $(J,n)$ for all elements $J\in 
   \Omega_{[1,n-2]}$. Consequently, for $u\in \mathbf C^n$, one has  ${\sf  F}_{[1,n]}(u)={\sf  F}_{[1,n-1]}(u')+u_n{\sf  F}_{[1,n-2]}(u'')$, from which it follows that 
   $F_n={\sf  F}_{[1,n]}$ for every $n\geq 0$. \mk 

Since $(i_1,\ldots,i_\ell)\mapsto (n+1-i_\ell, \ldots,n+1-i_1)$ induces an involution 
of $\Omega_{[1,n]}$, we get immediately that 
   the polynomial $F_n={\sf  F}_{[1,n]}$ is invariant by $(u_s)_{s=1}^n\mapsto (u_{n+1-s})_{s=1}^n$, which proves the first point of the proposition. 
   When combined with \eqref{Eq:Fn-induction}, this gives us the second one.\mk 
   
   As for point {\it 3.},  it follows from 
\eqref{Eq:F[a,b]} combined with Proposition 2.14 of \cite{FZ}. 
   \end{proof}

We are now ready to define  an explicit rational map 
  which will turn out to be the converse of $\boldsymbol{\sf U}_n$.

\subsubsection{The map $\boldsymbol{{\sf X}_n}$.}

Let   $n^*$ be the biggest even integer less than or equal to $n$, 
{\it i.e.}\,one has  $n^*=n$ if it is even and $n^*=n-1$ otherwise (mathematically: $n^*=2 \lfloor n/2\rfloor=2m$).\sk 
  
    Then  one defines an affine rational map of $u_1,\ldots,u_n$ by setting
  \begin{equation}
  \label{Eq:Formula-Xn}
   \boldsymbol{\sf X}_n= 
\left(
 \left(\,   \frac{F_{[1, 1+2k]}}{u_1 F_{[3, 1+2k]}}\right)_{k=0}^{\widetilde m}
\,   , \, 
  \left(  \frac{F_{[1, n^*-2\ell  ]}}{u_1 F_{[3,n^*-2\ell )]}}\right)_{\ell=0}^{m-1}
  \, \right) 
  \, . 
  \end{equation}
    Since $(\widetilde m+1)+m=n$,  $\boldsymbol{\sf X}_n$ can (and will) be considered as a rational map from $\mathbf C^n_u$ to $ \mathbf C^n_x$.  
  \sk   
    
    The rational map $\boldsymbol{\sf X}_n : \, \mathbf C^n_u\dashrightarrow  \mathbf C^n_x$
 also admits a nice inductive definition.  For any $n\geq 2$, one defines a 
   pair  $({\sf Y}_n,{\sf Z}_n)$ of rational maps in the variables $u_1,\ldots,u_n$, by first setting 
   $$ \big(\, {\sf Y}_2 \, , \, 
{\sf Z}_2
 \big) = 
\left(\,  \frac{F_{[1]}}{u_1}\, ,\, \frac{F_{[1,2]}}{u_1}
\, \right) 
=
\left(\,  \frac{1+u_1}{u_1}\, ,\, \frac{1+u_1+u_2}{u_1}\, 
\right) 
$$
 and for any $n\geq 2$, we define the pair 
 $( {\sf Y}_{n+1} ,{\sf Z}_{n+1})$ by setting:
\begin{align}
\label{Eq:Yn-Zn-Rec}
 {\sf Y}_{n+1}=\, & \big({\sf Y}_n\, , \, \varphi_{n+1}\big)  \quad  \mbox{and } \quad  {\sf Z}_{n+1}={\sf Z}_n \hspace{1.5cm } \mbox{for } \, n \,  \mbox{ even;  }  \quad\\
 {\sf Y}_{n+1}=\, & {\sf Y}_n  \hspace{1.6cm}
\mbox{and } \quad  {\sf Z}_{n+1}=\big(\varphi_{n+1}\, , \, {\sf Z}_n\big) 
\hspace{0.3cm }\mbox{for } \, n\,   \mbox{ odd, }  \quad
\nonumber
\end{align}
where for any $n\geq 3$, $ \varphi_n$ stands for the following rational function in $u_1,\ldots,u_n$: 
 \begin{equation}
 \label{Eq:varphi-n}
   \varphi_n= \varphi_n(u)= \frac{ F_{[1,n]}}{u_1F_{[3,n]}}\in \mathbf C(u)\, .
  \end{equation}
Then it is not difficult (and left to the reader) to verify that the tuple of rational functions obtained by concatenating the two corresponding to $ {\sf Y}_{n} $ and ${\sf Z}_n$ (in this order) coincides with the one formed by the components of $ \boldsymbol{\sf X}_{n} $, a relation which can be written  a bit abusively (for any $n\geq 3$) as 
  \begin{equation}
  \label{Eq:Xn=(Yn,Zn)}
   \boldsymbol{\sf X}_{n} =\big( \,{\sf Y}_{n} \, , \, {\sf Z}_n\,\big)
   =\big( \,{\sf Y}_{n-1} \, , \,
   \varphi_n 
    \, , \, {\sf Z}_n-1\,\big)
   \, .
    \end{equation}

One can notice that the rational function $ \varphi_{n}$ inserted between ${\sf Y}_{n-1}$ and ${\sf Z}_{n-1}$ to get $ \boldsymbol{\sf X}_{n}$ (for $n\geq 3$) can be characterized as the $m^*$-component of the latter. On the other hand, one verifies that   
${\sf Y}_{n-1}$ and ${\sf Z}_{n-1}$ only depend on 
$u'=(u_1,\ldots, u_{n-1} )$ 
 and that 
  $\varphi_{n}$ is the single component of $ \boldsymbol{\sf X}_{n}$ depending on the 
   variable 
   $u_{n}$ as well.  From these observations, one deduces the following counterpart of Lemma \ref{L:Un-Un-1}: 
\begin{lem}
\label{L:Xn-Xn-1}
For any $n\geq 3$, as rational maps from $\mathbf C_u^{n}$ to $\mathbf C_x^{n-1}$, one has 
\begin{equation}
\label{Eq:Xn-Xn-1}
\boldsymbol{{\sf X}}_{n-1}\circ p_n=  \pi_n\circ \boldsymbol{{\sf X}}_n
\, .
\end{equation}
\end{lem}


Using the relations \eqref{Eq:Un-Un-1} and \eqref{Eq:Xn-Xn-1} and the properties of $F$-polynomials discussed in the preceding paragraph, we are going to prove our first main result about $\boldsymbol{\sf U}_n$ and $\boldsymbol{\sf X}_n$: 

  \begin{prop}
  \label{P:Un-Xn=Id}
  As affine rational maps, $\boldsymbol{\sf X}_n$ and $\boldsymbol{\sf U}_n$ are inverse to each other, {\it i.e.}\,one has 
  \begin{equation} 
  \label{Eq:Un-Xn=Id}
 \boldsymbol{\sf X}_n\circ \boldsymbol{\sf U}_n={\rm Id}_{\mathbf C^n_x}
 \qquad \quad  \mbox{ and }  \quad 
   \qquad 
     \boldsymbol{\sf U}_n\circ\boldsymbol{\sf X}_n={\rm Id}_{\mathbf C^n_u}\, .
     \end{equation}
  \end{prop}
\begin{proof}
 Since $\boldsymbol{\sf U}_n$ is known to be birational, it suffices to prove that $\boldsymbol{\sf X}_n\circ \boldsymbol{\sf U}_n$ coincides with the identity of $\mathbf C^n_x$ as an affine rational map, which we  are going to prove by induction on $n\geq 2$. \sk
   
  For $n$ small enough (say $n$ less than 10), one verifies that \eqref{Eq:Un-Xn=Id} indeed holds true by explicit direct computations.  Let $n$ be big enough for what is to come ($n\geq 6$ is ok)  and  let us assume that the proposition holds true for all $n'\leq n$, in particular for $n'=n-1$. 
  Then, thanks to \eqref{Eq:Xn=(Yn,Zn)} and \eqref{Eq:Xn-Xn-1}, establishing the proposition for $n$ reduces to only proving that the following rational identity
  in the variables $x_1,\ldots,x_n$ holds true : 
  $$
 \varphi_n\big(\boldsymbol{\sf U}_n(x)\big)
  =  x_{m^*}
  \, . 
  $$
   
  From \eqref{Eq:Fn-induction}  combined 
with formulas    \eqref{F:unn}  and  \eqref{Eq:varphi-n} for ${\sf u}_n^1$ and $\varphi_n$, we deduce that this equality between rational functions is equivalent to the following one 
\begin{equation*}
\big(\star_n\big) \hspace{5cm}
\frac{F_{[2,n]}\big( \boldsymbol{\sf U}_n(x)\big)}{F_{ [ 3,n]}\big( \boldsymbol{\sf U}_n(x)\big)} =
\frac{ x_{ m^*}-1}{x_1-1} \, ,
\hspace{5cm} {}^{}
\end{equation*}
which we are going to prove by induction.  To simplify the writing, 
below we write  $ \boldsymbol{\sf U}$,  ${\sf u}_n^k$ instead of  $\boldsymbol{\sf U}_n(x)$, ${\sf u}_n^k(x)$, etc.\,(but we remain aware that all the functions we are dealing with  are actually  rational functions in the $x_i$'s). 
Using \eqref{Eq:Fn-induction}, one obtains that 
\begin{align}
\label{Eq:popopo}
\frac{F_{[2,n]}\big( \boldsymbol{\sf U}\big)}{F_{ [ 3,n]}\big( \boldsymbol{\sf U}\big)}=\frac
{ F_{ [ 2,n-1]}\big(\boldsymbol{\sf U}'\big)+
{\sf u}_n^n\cdot F_{ [ 2,n-2]}\big(\boldsymbol{\sf U}''\big)
}
{ F_{ [ 3,n-1]}\big(\boldsymbol{\sf U}'\big)+
{\sf u}_n^n\cdot F_{ [ 3,n-2]}\big(\boldsymbol{\sf U}''\big)
}
\end{align}
   where $\boldsymbol{\sf U}'$ stands for  $({\sf u}_n^k(x))_{k=1}^{n-1}$ viewed as  a $(n-1)$-tuple of rational functions in $x_1,\ldots,x_n$ (and similarly for $\boldsymbol{\sf U}''$). Using  $(\star_{n-1})$ and $(\star_{n-2})$  to express 
   $F_{ [ 2,n-1]}\big(\boldsymbol{\sf U}'\big) $ and 
   $F_{ [ 2,n-2]}\big(\boldsymbol{\sf U}''\big)$ in terms of 
   $F_{ [ 3,n-1]}\big(\boldsymbol{\sf U}'\big)$, $F_{ [ 3,n-2]}\big(\boldsymbol{\sf U}''\big)$ (respectively) and the $x_i$'s in  \eqref{Eq:popopo}, one gets that 
 \begin{align*}
\label{Eq:popopo2}
\frac{F_{[2,n]}\big( \boldsymbol{\sf U}\big)}{F_{ [ 3,n]}\big( \boldsymbol{\sf U}\big)}=\frac
{ x_{ \lfloor (n-1)/2 \rfloor+1}\cdot F_{ [ 3,n-1]}\big(\boldsymbol{\sf U}'\big)+
x_{ \lfloor (n-2)/2 \rfloor+1} \cdot {\sf u}_n^n\cdot F_{ [ 3,n-2]}\big(\boldsymbol{\sf U}''\big)
}
{ F_{ [ 3,n-1]}\big(\boldsymbol{\sf U}'\big)+
{\sf u}_n^n\cdot F_{ [ 3,n-2]}\big(\boldsymbol{\sf U}''\big)
}\, .
\end{align*}  
   
  Thanks to the technical lemma below, the  $F$-polynomial 
   $F_{ [ 3,n-1]}\big(\boldsymbol{\sf U}'\big)$ is equal to $F_{ [ 3,n-2]}\big(\boldsymbol{\sf U}''\big)$ times a certain explicit rational function  (in the $x_i$'s) of degree 1. Using this and the  explicit formula for ${\sf u}_n^n$ given in  \eqref{F:unn}, it is then easy to obtain that $(\star_n)$ indeed holds true, which concludes the proof. \end{proof}

%
%

   \begin{lem}
 For $m$ big enough, the following equalities hold true (as rational functions of $x$\footnote{Actually, all the expressions in \eqref{Eq:kiki} only depend on 
 $x^\circ=\pi_n(x)=(x_1,\ldots,x_m, \widehat{x_{m^*}}, x_{m+2}, \ldots, x_n)$.}): 
 \begin{equation}
 \label{Eq:kiki}
\frac{F_{[3,2m-1]}\big( \boldsymbol{\sf U}'\big)}{F_{ [ 3,2m-2]}\big( \boldsymbol{\sf U}''\big)}=
\frac{x_{m-1}-x_{m+2}}{x_{m-1}-x_{m}\hspace{0.3cm}{}^{}}
\qquad \quad \mbox { and } \qquad  \quad 
\frac{F_{[3,2m]}\big( \boldsymbol{\sf U}'\big)}{F_{ [ 3,2m-1]}\big( \boldsymbol{\sf U}''\big)}=
\frac{{}^{} \hspace{0.3cm}x_{m}-x_{m+3}}{x_{m+2}-x_{m+3}}\, . 
\end{equation}
    \end{lem}
   \begin{proof}
  The proof goes  by recurrence as well, the first cases being easily verified by means of direct explicit computations. \sk 
  
  Assume that  $2m-1$ is big enough and that the  identities \eqref{Eq:kiki} have been established for all the previous cases.  Considering the inductive property of the $F$-polynomials, it comes: 
\begin{align}
\nonumber 
\frac{F_{[3,2m-1]}\big( \boldsymbol{\sf U}'\big)}{F_{ [ 3,2m-2]}\big( \boldsymbol{\sf U}''\big)}=& \, 
\frac{    F_{[3,2m-2]}\big( \boldsymbol{\sf U}''\big) 
+({\sf u}')_{2m-1}^{2m-1} \cdot 
F_{[3,2m-3]}\big( \boldsymbol{\sf U}'''\big)     }
{ F_{ [ 3,2m-2]}\big( \boldsymbol{\sf U}''\big)  }  \\ 
= & \, 1+({\sf u}')_{2m-1}^{2m-1} \cdot \frac{F_{[3,2m-3]}\big( \boldsymbol{\sf U}'''\big)}{F_{ [ 3,2m-2]}\big( \boldsymbol{\sf U}''\big)}\, .
\label{lala}
\end{align}   
All the terms in the RHS of the second equality  can be made explicit. Indeed,  one has
\begin{equation*}
({\sf u}')_{2m-1}^{2m-1}= \frac{(x_m-x_{m+2 })(x_{ m-1}-x_{ m+3}) }{
(x_{m+2 }-x_{m+3 })(x_{ m-1}-x_{ m})
} \qquad 
\mbox{ and } \qquad \frac{F_{[3,2m-3]}\big( \boldsymbol{\sf U}'''\big)}{F_{ [ 3,2m-2]}\big( \boldsymbol{\sf U}''\big)}= 
\frac{x_{m+2 }-x_{ m+3} }{x_{m-1 }-x_{m+3 }}\, , 
\end{equation*}
the former formula following easily from  \eqref{F:unn} and the latter being given by the induction hypothesis.  Injecting both formulas in the RHS of equality \eqref{lala} and after some trivial simplifications, one obtains the left formula of \eqref{Eq:kiki}. 
 The proof that the right formula holds true is similar. \end{proof}

%

\subsubsection{Some properties of the maps $\boldsymbol{{\sf U}_n}$ and $\boldsymbol{{\sf X}_n}$.}
 In view of establishing that $\boldsymbol{{\sf U}_n}$ and $\boldsymbol{{\sf X}_n}$ induce isomorphisms between nice Zariski open domains of $\mathbf C_x^n$ and $\mathbf C_u^n$, we first establish a few basic properties of these maps, considered as affine rational maps. 
\mk 

 \paragraph{} ${}^{}$\hspace{-0.6cm} 
 \label{Par:BraidArr}
 We begin with the map $\boldsymbol{{\sf X}_n}$.  Considering its inductive definition discussed above, it comes that its indeterminacy locus, denoted by ${\rm Ind}(\boldsymbol{{\sf X}_n})$, is exactly the union of the hyperplane $H_{u_1}=\{ u_1=0\}$ with the  irreducible hypersurfaces $H_{F_{ [3,k]}}=\{ \,F_{ [3,k]}=0\, \}$ for $k=3,\ldots,n$: 
one has
  $$
  {\rm Ind}(\boldsymbol{{\sf X}_n})= 
  H_{u_1}\bigcup \Big( \cup_{k=3}^n H_{F_{ [3,k]}} \Big) \subset \mathbf C^n_u\, .
  $$
  
 Recall that except the $m^*$-component $\varphi_m$ of   $\boldsymbol{{\sf X}_n}$, all the others actually only depend on $u'=(u_1,\ldots,\widehat{ u_{ m^*}},\ldots,u_n)$. Thus if $ \boldsymbol{\mathcal J{\sf X}}_n$ stands for the jacobian determinant of ${\sf X}_n$ (in the variables $u_1,\ldots,u_n$), it follows from 
\eqref{Eq:Xn=(Yn,Zn)} and \eqref{Eq:varphi-n}
 that for any $n\geq 3$, one has 
    $$
   \boldsymbol{\mathcal J{\sf X}}_n
   =(-1)^{n+m^*}
\frac{\partial }{\partial u_n}  \bigg(  \frac{F_{[1, n]}}{u_1 F_{[3, n]}}
   \bigg)\cdot  \big(  \boldsymbol{\mathcal J{\sf X}_{n-1}}\circ p_n\big)
  $$
Since ${\partial }  \big( 
{F_{[1, n]}}/{(u_1 F_{[3, n]})}
\big)/{\partial u_n}=(-1)^n  \big( \prod_{j=2}^{n-1} u_j \big)/\big(u_1 F[3,n]^2 \big)
$ for any $n\geq 3$ (as it can easily be proven, still by induction),  we deduce that 
  $$
   \boldsymbol{\mathcal J{\sf X}}_2=-\frac{1}{u_1^3} \qquad \qquad  \mbox{ and } \qquad  \qquad 
   \boldsymbol{\mathcal J{\sf X}}_n=- \frac{\prod_{j=2}^{n-1} u_j^{n-j}}{u_1^{n+1}
  \big(\prod_{k=3}^{n} F[3,k]\big)^2
  }\quad \mbox{ for any } n\geq 3\, .
  $$
  
  Denoting by $ \boldsymbol{\sf J{\sf X}}_n$ the union (in 
  $\mathbf C^n_u$) 
  of ${\rm Ind}(\boldsymbol{{\sf X}_n})$  with $\cup_{k=3}^n H_{u_j}$, it follows from the formula above for $  \boldsymbol{\mathcal J{\sf X}}_n$ that 
   $\boldsymbol{\sf X}_n$ induces an isomorphism from the complement
   of $ \boldsymbol{\sf J{\sf X}}_n$ 
    in
  $\mathbf C^n_u$  onto its image in $\mathbf C_x^n$. 
  Since the complement $U_\Delta$ of the arrangement of type $\Delta=A_n$ (see 
  \eqref{Eq:U--Delta} above)  is included in  $\mathbf C^n_u\setminus   \boldsymbol{\sf J{\sf X}}_n$, we get that the map under scrutiny induces a rational affine isomorphism 
  \begin{equation}
 \label{Eq:Xn-isom}
 \boldsymbol{{\sf X}_n}\lvert_{U_\Delta} \, :\, U_\Delta \stackrel{\sim} {\longrightarrow}
  \boldsymbol{{\sf X}_n}(U_\Delta)\subset \mathbf C_x^n\,. 
\end{equation}
  In what follows, to simplify,  we will denote again  this isomorphism by $\boldsymbol{{\sf X}_n}$.

 \paragraph{} ${}^{}$\hspace{-0.6cm}
  \label{Par:BraidArr-2}
  We now consider the map $\boldsymbol{{\sf U}_n}$. To state some of its properties, we recall the definition of the `Braid arrangement' of type $A_n$, denoted here by ${\sf A}_n$ (see also 
 \cite[\S4]{Pereira}, where it is 
 denoted by   ${\sf A}_{0,n+3}$): by definition, it is the arrangement of $n(n+3)/2$ affine hyperplanes in $\mathbf C^n_x$ cut out by 
 \begin{equation}
\label{Eq:cut-An}
\Big( \prod_{i=1}^n x_i\Big)
\Big( \prod_{i=1}^n \big(x_i-1\big)\Big)\Big( \prod_{1\leq i <j \leq n}^n \big(x_i-x_j\big)\Big)\, . 
 \end{equation}

We denote by $X_\Delta=\mathbf C^n_x\setminus {\sf A}_n$ the complement of  ${\sf A}_n$ in $\mathbf C^n_x$: it is known that this   Zariski open subset of 
$\mathbf C^n_x$ 
 can naturally (although not in a unique way) be identified with the moduli space $\mathcal M_{0,n+3}$. \sk

  From the formulas \eqref{F:unn} and using Lemma \ref{Eq:Un-Un-1}, one gets that  
  the indeterminacy set ${\rm Ind}(\boldsymbol{\sf U}_n)$  of $\boldsymbol{\sf U}_n$
is included in $X_\Delta$ and that the Jacobian determinant 
$\boldsymbol{\mathcal J{\sf U}}_n$
 of the induced morphism  
satisfies
\begin{equation}
\label{Eq:JUn}
   \boldsymbol{\mathcal J{\sf U}}_2=-\frac{1}{(x_1-1)^3} \qquad  \mbox{ and } \qquad  
  \boldsymbol{\mathcal J{\sf U}}_n=(-1)^{n+m^*}\bigg(\frac{\partial {\sf U}^n_n}{\partial x_{m^*}}\bigg) \,
\Big(\boldsymbol{\mathcal J{\sf U}}_{n-1}\circ \pi_n\Big)
\quad \mbox{for any } n\geq 3\, .
  \end{equation}

Let   $\mathcal S_{n,x}$  be the multiplicative subgroup of 
the field of rational functions in the $x_i$'s 
 spanned by the affine functions $ x_i$, $x_i-1$ and $ x_i-x_j$ for all 
 $ i,j=1,\ldots,n$ (with $ i\neq j $) appearing in \eqref{Eq:cut-An}:
 \begin{equation}
\label{Eq:Set-Sn}
\mathcal S_{n,x}=\Big\langle 
\hspace{0.1cm}
x_i\, , \, x_i-1\, , \, x_i-x_j
\hspace{0.1cm} \big\lvert 
\hspace{0.1cm}
\begin{tabular}{l}
$i,j=1,\ldots,n$, 
$i\neq j$
\end{tabular}
\Big\rangle \subset \mathbf C(\boldsymbol{x})\, .
\end{equation}

From  \eqref{F:unn}  and \eqref{Eq:JUn}, one deduces easily that $\boldsymbol{\mathcal J{\sf U}}_n$ belongs to $\mathcal S_{n,x}$ for any $n\geq 2$.\footnote{Actually, it is not difficult to get an explicit closed formula for $\boldsymbol{\mathcal J{\sf U}}_n$ but we will not elaborate on this since it is not needed for our purpose.} Consequently, 
$ \boldsymbol{{\sf U}}_n$ induces 
 a rational affine isomorphism 
 \begin{equation}
 \label{Eq:Un-isom}
 \boldsymbol{{\sf U}_n}\lvert_{X_\Delta} \, :\, X_\Delta \stackrel{\sim} {\longrightarrow}
  \boldsymbol{{\sf U}_n}(X_\Delta)\subset \mathbf C_u^n 
\end{equation}
 that we will denote again by $\boldsymbol{{\sf U}_n}$  in what follows to simplify the writing. 
\sk 

The result we aim to prove regarding  the two morphisms considered above  
is the following 
\begin{prop} 
\label{P:Xn(Un)=XXn}
1. One has $\boldsymbol{{\sf X}_n}(U_\Delta)=X_\Delta$ and 
$\boldsymbol{{\sf U}_n}(X_\Delta)=U_\Delta$.\mk 

2. Consequently,  the affine morphism  $\boldsymbol{{\sf X}_n}$ in \eqref{Eq:Xn-isom}
induces an isomorphism 
from $X_\Delta$ onto $U_\Delta$ whose converse  is induced by the isomorphism $\boldsymbol{{\sf U}_n}$ in \eqref{Eq:Un-isom}.
\mk 

3.  It follows that the fundamental group of 
$U_\Delta$ is isomorphic to the one of $X_\Delta\simeq \mathcal M_{0,n+3}$, 
namely the pure braid group on $n+2$  strands : $\pi_1(U_\Delta)\simeq {\rm PB}_{n+2}$. 
\end{prop}
The second and third points of this proposition follow rather straightforwardly from the first, so we will only focus on it. Since both \eqref{Eq:Xn-isom}
and \eqref{Eq:Un-isom} are isomorphisms onto their respective images, it suffices to prove that the two following inclusions hold true: 
\begin{equation}
\label{Eq:Un-Xn-Inclusions}
\boldsymbol{{\sf U}_n}\big(X_\Delta\big)\subset U_\Delta
\qquad \qquad \mbox{ and } \qquad \qquad
\boldsymbol{{\sf X}_n}\big(U_\Delta\big)\subset X_\Delta
\, .
\end{equation}
We prove successively both inclusions 
in the next two paragraphs.

 \paragraph{Proof that $\boldsymbol{{\sf U}_n(X_\Delta)}$ is included in  $\boldsymbol{U_\Delta}$.}
    One has to prove that given any $x\in X_\Delta$, 
    that is $x=(x_i)_{i=1}^n\in \mathbf C^n$ with none of the linear expressions $x_i$, $x_i-1$ and $x_i-x_j$ for $i,j=1,\ldots,n$ with $i<j$, being equal to zero, then $\boldsymbol{{\sf U}_n}(x)$ 
  does not belong to any of the coordinate hyperplanes $H_{u_k}$ (for $k=1,\ldots,n$) 
  nor to any hypersurface 
  $H_\alpha$ with $\alpha\in \Delta_{>0}$.     
  This is implied by the fact that 
  \begin{equation}
  \label{Eq:Ze-Fact}
  \begin{tabular}{l}
  {\it for any $k$ or any pair $(a,b)\in \mathbf N^2$ such that $1\leq a\leq b\leq n$, as rational functions } \\ 
  {\it  of the $x_i$'s,
  the $k$-th component 
   ${\sf U}_n^k$ of $\boldsymbol{\sf U}_n$ and the composition
 $P_{a,b}=F_{[a,b]}\big( 
 \boldsymbol{{\sf u}_n}\big)$ }\\ 
{\it   belong to the set  $\mathcal S_{n,x}$ defined in \eqref{Eq:Set-Sn}.}
  \end{tabular}
  \end{equation}
  
For the ${\sf u}_n^k$'s, this follows immediately from  their formulas \eqref{F:unn} hence we only have to prove that 
 $P_{a,b}\in \mathcal S_{n,x}$ for any pair $(a,b)$ as above, which we are going to do by induction on $n$. \mk

 For $n$ small,  everything can be established by means of explicit computations. So we assume that the statement under scrutiny has been established for  $n-1\geq 2$: 
setting $x^\circ=\pi_n(x)=(x_1,\ldots , \widehat{x_{m^*}}, \ldots,x_n)$, one has 
 $P_{a,b}(x)=P_{a,b}(x^\circ)\in \mathcal S_{n-1,x^\circ}$ for any $(a,b)$ with $1\leq a\leq b\leq n-1$.   Since $\mathcal S_{n-1,x^\circ}$ is naturally a subset of $ \mathcal S_{n,x}$, this answers the question for all these cases and we only have to deal with the $P_{k,n}$'s for $k=1,\ldots,n$.  
 The case of $P_{n,n}$ is easy to deal with: since $F_{[n,n]}=1+u_n$, it follows from  the formula for ${\sf u}_n^{n}$ in  \eqref{F:unn}   that 
 $$P_{n,n}= 1+ {\sf u}_n^{n}= 
 \frac{  (x_{\widetilde m + 1} - x_{\widetilde m + 3})( x_{\widetilde m}-x_{\widetilde m + 2})}{(x_{\widetilde m + 2} - x_{\widetilde m + 3})(x_{\widetilde m} - x_{\widetilde m + 1})} 
 \in \mathcal S_{n,x}\, .$$

To deal with the other cases, we set: 
\begin{itemize}
\item $\rho_{k,n}=
  P_{k,n}/P_{k,n-1} \in \mathbf C(x)$ for $ k=1,\ldots,n-1$;

  \item $\epsilon=(-1)^n\in \{ \, \pm 1\, \}$ ({\it i.e.}\,$\epsilon$ is equal to 1 if $n$ is even and to $-1$ when it is odd);
  \item $c=(c_s)_{s=1}^{n+3}=(0,1,\infty,x)$, that is $c_1=0$, $c_2=1$, $c_3=\infty$ and $c_k=x_{k-3}$ for $k=4,\ldots,n+3$;
  \item $K=K'=m-2$ when $n$ is odd and $K=K'+1=m-2$  ({\it i.e.}\,$K'=m-3$) when $n$ is even; 
  \item  $\sigma=\sigma_n: \{1,n-1\}\rightarrow \{1,n\}$ for the map 
 defined by  
  \begin{align*}
    \sigma\big(\ell\big)=&\, \ell \hspace{1.75cm}\mbox{ for }\, \ell=1,\ldots,4\, ;\\
  \sigma\big(3+2k\big)=&\, n+1-k \hspace{0.5cm}\mbox{ for }\, k=1,\ldots,K\, ; \\
 \mbox{and} \quad  \sigma\big(4+2k'\big)=&\,4+k' \hspace{1.0cm}\mbox{ for }\, k'=1,\ldots,K'\, .
   \end{align*}
\end{itemize}
                                          
For each $k$, $\rho_{k,n}$ is considered as an element of $\mathbf C(x)$ with $x=(x_1,\ldots,x_n)$.  But the reader will have in mind that in the formula
$\rho_{k,n}= P_{k,n}/P_{k,n-1}$ defining it, 
the denominator $P_{k,n-1}$ coincides with $F_{[k,n-1]}(\boldsymbol{\sf U}_{n-1})$  evaluated on the $(n-1)$-tuple of indeterminates $x^\circ$. 
\sk 

With these notations at hand, we can state the technical result which we will use  to prove \eqref{Eq:Ze-Fact} for the remaining cases: 
 \begin{lem}
  For $n\geq 3$ and $k=1,\ldots,n-1$, one has 
\begin{equation}
\label{Eq:rho-k-n}
  \rho_{k,n}=\frac{\big(x_{m^*-\epsilon}-x_{m^*+\epsilon}\big)}{
  \big(\hspace{0.2cm}x_{m^*}\hspace{0.15cm} -x_{m^*+\epsilon}\big)} \,
\frac{\big(x_{m^*}\hspace{0.15cm}-c_{\sigma(k)}\big)}{
  \big(x_{m^*-\epsilon}-c_{\sigma(k)}\big)} \, .
\end{equation}
(Note that for $k=3$, one has $c_{\sigma(3)}=c_3=\infty$ hence $\frac{(x_{m^*} \, -c_{\sigma(k)})}{(x_{m^*-\epsilon}-c_{\sigma(k)})}=\frac{(x_{m^*}\, -\infty)}{(x_{m^*-\epsilon}-\infty)}=1$ by convention). 
\end{lem}
\begin{proof}
The proof goes by induction as well. 
From  \eqref{Eq:Un-Un-1} and \eqref{Eq:Fn-induction}, 
it comes that for $n$ big enough, one has $F_{[k,n]}(\boldsymbol{\sf U}_n(x))=F_{[k,n-1]}(\boldsymbol{\sf U}_{n-1}(x^\circ))+{\sf u}_n^n\cdot  F_{[k,n-2]}(\boldsymbol{\sf U}_{n-2}(x^{\circ\circ}))$ where $x^{\circ\circ}$ stands for the 
$(n-2)$-tuple of indeterminates given by $\pi_{n-1}(x^\circ)$.  Thus, one gets 
$$
  \rho_{k,n}=\frac{
  F_{[k,n-1]}(\boldsymbol{\sf U}_{n-1}(x^\circ))+{\sf u}_n^n\cdot  F_{[k,n-2]}(\boldsymbol{\sf U}_{n-2}(x^{\circ\circ}))
  }{F_{[k,n-1]}(\boldsymbol{\sf U}_{n-1}(x^\circ))}
  =1+
 \frac{{\sf u}_n^n}{ \rho_{k,n-1}} .
$$ 
Injecting  formula   \eqref{F:unn}  for ${\sf u}_n^n$ and the one for 
$\rho_{k,n-1}$ given by the 
 induction hypothesis, it is then not difficult to verify that formula 
 \eqref{Eq:rho-k-n} for $\rho_{k,n}$ is satisfied as well.
\end{proof}

Using the preceding lemma, the fact  that \eqref{Eq:Ze-Fact} holds true for  $P_{k,n}$ becomes obvious (for any $k\in \{1,\ldots,n-1\}$):  \eqref{Eq:rho-k-n} gives us that $\rho_{k,n}\in  \mathcal S_{n,x}$ and because $P_{k,n-1}\in  \mathcal S_{n,x^\circ}\subset \mathcal S_{n,x}$ according to the induction hypothesis, we get immediately that $
P_{k,n}=\rho_{k,n}\, P_{k,n-1} \in  \mathcal S_{n,x}$ as well, which finishes the proof  that the first inclusion 
in \eqref{Eq:Un-Xn-Inclusions} indeed holds true.

 \paragraph{Proof of the inclusion of $\boldsymbol{{\sf X}_n(U_\Delta)}$ into $\boldsymbol{X_\Delta}$.} 
Our proof is very similar in spirit to the one given in the preceding paragraph hence many details are left to the reader. 
\sk 

For $n\geq 2$, we recall that $u$ stands for the $n$-tuple of indeterminates $(u_i)_{i=1}^n$ and we denote by $\mathcal F_{n,u}$ the multiplicative stable subset of $\mathbf C(u)$ generated by all the $u_i$'s and the $F$-polynomials $F_{[a,b]}$ and their inverses: 
 \begin{equation}
\label{Eq:Set-Fn}
\mathcal F_{n,u}=\Big\langle 
\hspace{0.1cm}
u_i^{\pm 1} \, , \, \big(F_{[a,b]}\big)^{\pm 1}
\hspace{0.1cm} \big\lvert 
\hspace{0.1cm}
\begin{tabular}{l}
$i=1,\ldots,n$\, ,\,   $1\leq a\leq b \leq n$
 \, 
\end{tabular}
\Big\rangle \subset \mathbf C(\boldsymbol{u})\, .
\end{equation}

For $i=1,\ldots,n$,  we denote by ${\sf x}_n^i$ the $i$-th component of 
$\boldsymbol{\sf X}_n$.
  The inclusion under scrutiny follows from the following fact: 
 \begin{equation}
  \label{Eq:Ze-Fact-Xn}
  \begin{tabular}{l}
  {\it for any $i,j,k=1,\ldots,n$ such that $1\leq j<k\leq n$,  the rational  } \\ 
  {\it 
 functions (of the $u_i$'s) ${\sf x}_n^i$, ${\sf x}_n^i-1$ and ${\sf x}_n^j-{\sf x}_n^k$  all 
  belong to $\mathcal F_{n,u}$.}
  \end{tabular}
  \end{equation}
  
    That each ${\sf x}_n^i$ belongs to $\mathcal F_{n,u}$ follows immediately from \eqref{Eq:Formula-Xn}.  As for the ${\sf x}_n^i-1$'s, this follows from the fact that for any $k\leq n$, thanks to  the second point of Proposition  \ref{P:Prop-F-poly-An}, one has  
    $$
    \frac{F_{[1,k]}}{ u_1 F_{[3,k]}}-1=\frac{F_{[2,k]}+u_1 F_{[3,k]}}{ u_1 F_{[3,k]}}-1=
    \frac{F_{[2,k]}}{ u_1 F_{[3,k]}}\in \mathcal F_{n,u}\, , 
    $$
    the second equality following from the second point of Proposition  \ref{P:Prop-F-poly-An}.\mk

 Lemma \ref{L:Xn-Xn-1} allows us 
to argue by induction  in what concern the differences ${\sf x}_n^j-{\sf x}_n^k$ and for $n$ given, to only consider the case of those for which $j$ is equal to $m^*$.  To this end we set 
\begin{align}
\label{Al:pototo}
\alpha_n^l= \frac{{\sf x}_n^{m^*}-{\sf x}_n^l}{F_{[2+2l\,,\,n]}} \hspace{0.3cm} \mbox{for }\, l=1,\ldots,m
\qquad \mbox{ and } \qquad 
  \beta_n^\ell=  \frac{ {\sf x}_n^{m^*}-{\sf x}_n^{m^*+\ell} }{ F_{[n^\square+3-2\ell\, ,\, n]} } 
  \hspace{0.3cm} \mbox{for }\, \ell=1,\ldots,\widetilde{m}\,
  \end{align}
  where 
$n^\square$ stands for the smallest even integer bigger than or equal to $n$
(that is  $n^\square$ is equal to $n$ if the latter is even but to  $n+1$ otherwise; 
mathematically: $n^\square=2 \lfloor{(n+1)/2}\rfloor$).\footnote{In \eqref{Al:pototo} as everywhere else, we use the following convention: 
 we set $F_{[p,q]}=1$ for any $p,q\in \mathbf N$ such that $p>q$.} Then, using 
\eqref{Eq:Yn-Zn-Rec}
and some properties of the $F$-polynomials 
discussed in \S\ref{SS:SomeProp-F-polynomials}, 
one can verify that for any $n\geq 3$,  both tuples $\alpha_n=(\alpha_n^l)_{l=1}^m$ and 
  $\beta_n=(\beta_n^\ell)_{\ell=1}^{\widetilde m}$ satisfy the following recursive relations: 
\begin{align}
\label{Eq:Alpha-n-Beta-n-Rec}
 {\alpha}_{n}=\, & \big({\alpha}_{n-1}\, , \, \psi_{n}\big)  \quad  \mbox{and } \quad  {\beta}_{n}={\beta}_{n-1} \hspace{1.5cm } \mbox{for } \, n \,  \mbox{ even;  }  \quad
 \\
  {\alpha}_{n}=\, &  {\alpha}_{n-1}  \hspace{1.35cm}
\mbox{and } \quad  {\beta}_{n}=\big(\psi_{n}\, , \, {\beta}_{n-1}\big) 
\hspace{0.55cm }\mbox{for } \, n\,   \mbox{ odd, }  \quad
\nonumber
\end{align}
where we have set  $\psi_n=(-1)^n u_3\cdots u_n/ 
F_{[3,n-1]}$. From this and still arguing inductively, one deduces closed formulas 
  for the differences 
  ${\sf x}_n^{m^*}-{\sf x}_n^j$: 
  one has
\begin{equation*}
{\sf x}_n^{m^*}-{\sf x}_n^l =\bigg(\frac{\prod_{s=1}^{2l} u_s}{u_1u_2}\bigg)\, 
\frac{F_{[2+2l,n]}}{F_{[3\, ,\, 2k-1]}}
\qquad \mbox{ and } \qquad 
{\sf x}_n^{m^*}-{\sf x}_n^{m^*+\ell}=- \bigg(u_3 \hspace{-0.15cm}\prod_{s=2}^{\widetilde{m}+1-\ell} 
 \hspace{-0.15cm} u_{2s}u_{2s+1}\bigg)\,
\frac{F_{[n^\square+3-2\ell,n]}}{F_{[3\, ,\, n^\square -2\ell]}} \, .
\end{equation*}
   for any $l=1,\ldots,m$ and  any $\ell=1,\ldots,\widetilde{m}$. In particular, one has  ${\sf x}_n^{m^*}-{\sf x}_n^j\in \mathcal F_{n,u}$ for any $j$.  As said above, this implies that $\boldsymbol{\sf X}_n(U_\Delta)$ is  included in ${X_\Delta}$.

 \paragraph{For other Dynkin types.} 
   The arguments of the two preceding paragraphs taken together give a proof of the first point of Proposition 
  \ref{P:Xn(Un)=XXn}, from which the last two can be deduced easily.\mk 

  It is natural to wonder whether  this proposition admits versions for the other Dynkin types. One can say that the affine birational map $\boldsymbol{\sf X}_n: \mathbf C^n_u\dashrightarrow \mathbf C^n_x$ allows to `linearize' the arrangement of hypersurfaces $\boldsymbol{Arr}_{\Delta}$ when $\Delta$ has type $A_n$ in the sense that the complement of the image  of the complement of $\boldsymbol{Arr}_{\Delta}$ by $\boldsymbol{\sf X}_n$ is an arrangement of hyperplanes (namely, it is ${\sf A}_n$). Does a similar statement holds true for  Dynkin diagram of other types? More precisely, one asks the 
\begin{question} 
Given  a rank $n$ Dynkin diagram of arbitrary type $\Delta$, is there an (explicit?) affine birational map  $\boldsymbol{\sf X}_\Delta: \mathbf C^n_u\dashrightarrow \mathbf C^n_x$ inducing an isomorphism when restricted to $U_\Delta= 
\mathbf C^n_u \setminus \boldsymbol{Arr}_{\Delta}
$ and such that the complement of $\boldsymbol{\sf X}_\Delta(U_\Delta)$ in $\mathbf C^n_x$ by a linear arrangement?
\end{question}
  If the answer is yes (and is explicit) in type $A$, some considerations when $\Delta=D_4$ suggest that it may not be the case for other types.  The fact that the fundamental group of $U_\Delta$ in type $A$ is a well-known group in geometry  suggests to look at the $\pi_1(U_\Delta)$'s for $\Delta$ of arbitrary type.  
\begin{questions} 
Given a Dynkin diagram $\Delta$ not of type $A$, is there a presentation of 
$\pi_1(U_\Delta)$ by generators and relations similar to the classical one of the pure braid group? More generally, is $\pi_1(U_\Delta)$ isomorphic to a group already known and studied in topology and/or geometry? 
\end{questions}
  
%

For other considerations and questions in the same vein, 
 see \S\ref{Parag:ArrangementsResonanceWebs} further on.

  \newpage


\section{Cluster webs of other Dynkin types}
\label{S:ClusterWebs-other-Dynkin-Type}

The rather satisfying results obtained concerning cluster webs in type $A$ so far naturally lead to consider the cases of other Dynkin types. 
\mk

In this section, we discuss the cluster webs in type $B,C$ and $D$ 
(mainly their ${\mathcal Y}$-versions, but we say a few words about their 
${\mathcal X}$-versions as well).  These webs does not seem to be as nice as the corresponding ones in type $A$ (not AMP in general) but still seem interesting regarding  their Abelian relations and their rank.  We do not have well proved results concerning them but several precise conjectural statements which are interesting according to us.\sk

After that we turn to the case of bi-Dynkin type cluster webs.  We do not have any conclusive results concerning them, but only a fairly precise  conjecture that has been verified in a large number of the particular cases we have considered.  If true, this conjecture would give us several families of AMP cluster webs.  

 \subsection{\bf Cluster webs of type $\boldsymbol{B}.$}
\label{SS:ClusterWebs-type-B}
Below, we first introduce some notation specific to the $B_n$ case for $n\geq 2$ arbitrary, then discuss the dilogarithmic identity $({\sf R}_{B_n})$ in \S\ref{SSS:ClusterWebs-type-B:dilog-identity}, before dealing with the web $\boldsymbol{\mathcal Y\hspace{-0.1cm}\mathcal W}_{\hspace{-0.05cm} B_n}$. 
  By way of illustration, we first treat the case of $\boldsymbol{\mathcal Y\hspace{-0.1cm}\mathcal W}_{\hspace{-0.05cm} B_3}$ (the $B_2$ case has been considered before in \S\ref{SS:cluster-webs:B2}) before discussing the rank and the ARs when the rank is arbitrary. 
\mk 

In what follows, $n$ stands for an integer bigger than 1.

 \subsubsection{\bf Some notations in type $\boldsymbol{B}.$}
\label{SSS:ClusterWebs-type-B:notations}
There is a well-defined bipartite valued quiver $\vec{B}_n$  obtained from  from  Dynkin diagram 
$
\dynkin B{}
$ of type $B_n$ by requiring that it ends by an arrow
 from $n-1$ to $n$ valued by $(2,1)$, {\it i.e.}\,$\vec{B}_n$ ends like this: 
\begin{figure}[h]
\begin{center}
 \vspace{-0.2cm}
\scalebox{0.4}{
 \includegraphics{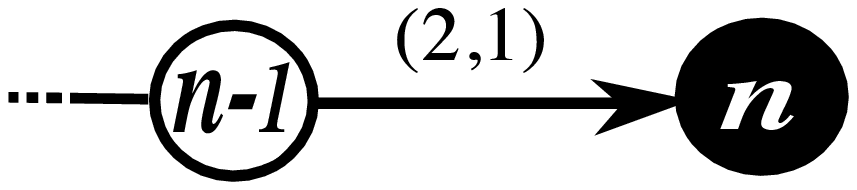}}
 \vspace{-0.6cm}
\end{center}
\end{figure}\sk\\
The corresponding exchange matrix is :
$$B_{\vec{B}_n}=
\begin{bmatrix}
0 & (-1)^n & 0   &\cdots &  \cdots& 0\\
(-1)^{n+1} & 0 & \ddots   & &  & \vdots\\
0&   (-1)^n&   \ddots &  \ddots & &  \vdots\\
\vdots & \ddots &  \ddots  &\ddots &-1 & 0 \\
\vdots&  &   \ddots &1 & {}^{}\hspace{0.2cm}0 & 2\\
0 & \cdots &  0  &0 & -1&0
\end{bmatrix} \, .
$$

The set of sources  (resp.\,of sinks)  in $\vec{B}_n$ is formed by the integers in $\{1,\ldots,n\}$ with opposite (resp.\,the same) parity to the one of $n$.  
As initial ($\mathcal X$-)seed, we take $\boldsymbol{\mathcal S}_0=\big(\boldsymbol{v}_0 , {\vec{B}_n}\big)$ where $\boldsymbol{v}_0$ stands for the $n$-tuple of variables 
\begin{equation}
\label{Eq:nu-0}
\boldsymbol{v}_0=\Big(u_i^{(-1)^{n-i}}\Big)_{i=1}^n=\Big( (u_1)^{(-1)^{n-1}} , \,   (u_2)^{(-1)^{n-2}}    , \,   \ldots , \,   u_{n-2} \,  , \, 
1/u_{n-1}
\, , \, u_n\Big)\,.
\end{equation}

We denote by $\mu_\circ$ (resp.\,$\mu_\bullet$) the composition of mutations on the set of sources (resp.\,of sinks) in 
$\vec{B}_n$ (for each, the order does not matter) and we define tuples of integers by setting 
%
%
%
\begin{align*}
\nu_\bullet=& \, \big( \, n\, , \, 
n-2\, , \,  \ldots\, , \,  4,2
\big)\, ,  &&  \nu_\circ=\big(\,  n-1\, , \, n-3\, , \,  \ldots \, , \,  3\, , \, 1\big)  && \mbox{ for } \, n\, \mbox{ even}\,;\\
\mbox{and }\quad 
\nu_\bullet=& \,\big(\,  n\, , \, n-2\, , \,  \ldots \, , \,  3\, , \, 1\big)\, ,  &&  \nu_\circ=\big( \, n-1\, , \, 
n-3\, , \,  \ldots\, , \,  4,2
\big)   && \mbox{ for } \, n\, \mbox{ odd }\,.
\end{align*}

Then the $n$-tuple $\nu_{\bullet,\circ}=\nu_\bullet \lvert \nu_\circ$ obtained by concatenating  $\nu_\bullet $ and $\nu_\circ$ is such that 
$\mu_\square=(\nu_{\bullet,\circ})^{n+1}$ (which is the $n(n+1)$-tuple obtained by concatenating $n+1$ copies of $\nu_{\bullet,\circ}$) is a period  for $\boldsymbol{\mathcal S}_0$. 

Indeed,  setting 
$$\mu_\square=\Big(i_{n(n+1)},i_{n(n+1)-1},\ldots,i_2,i_1\Big)
$$ 
with $i_s\in \{1,\ldots,n\}$ for all $s=1,\ldots,n(n+1)$, we define  new seeds inductively  as follows: 
\begin{align*}
\boldsymbol{\mathcal S}_1= \mu_{i_1}\big(\boldsymbol{\mathcal S}_0\big)=\begin{cases}  \hspace{0.2cm}\mu_1\big(\boldsymbol{\mathcal S}_0\big)\hspace{0.2cm} \mbox{ for } \, n\, \mbox{ even}\\
\hspace{0.2cm}
\mu_2\big(\boldsymbol{\mathcal S}_0\big)\hspace{0.2cm} \mbox{ for } \, n\, \mbox{ odd}
\end{cases}
\qquad \mbox{and}\qquad 
\boldsymbol{\mathcal S}_{s+1}=\mu_{i_s}(\boldsymbol{\mathcal S}_{s})  \mbox{ for all }\, s\geq 1.
\end{align*}
Then one has $
\boldsymbol{\mathcal S}_{n(n+1)} = 
\mu_\square\big(\boldsymbol{\mathcal S}_{0}\big)=\boldsymbol{\mathcal S}_{0}$ (as $\mathcal X$-seeds).
\mk 

Given $s\in \{1,\ldots,n(n+1)\}$, one denotes by 
$
Y_s
$ the $i_s$-th cluster variable of the $s$-th seed $\boldsymbol{\mathcal S}_s$, that we see as an element of $\mathbf C(u) $.  The model of the $\mathcal Y$-cluster web of type $B_n$ we will deal with is the $n(n+1)$-web defined by the 
all the $Y_s $: 
\begin{equation}
\label{Eq:Web-YWBn}
\boldsymbol{\mathcal Y\hspace{-0.1cm}\mathcal W}_{\hspace{-0.05cm} B_n}
=\boldsymbol{\mathcal Y\hspace{-0.1cm}\mathcal W}\Big( \, 
Y_s\, \big\lvert \, s=1,\ldots,n(n+1)\, \Big)\, .
\end{equation}

Finally, we notice that a right skew-symmetrizer for $B_{\vec{B}_n}$ is given by the diagonal matrix whose diagonal coefficients form the $n$-tuple 
$d_{B_n}=(2,\ldots,2,1)$ since one has
%
%
$$
B_{\vec{B}_n}\cdot 
\begin{pmatrix}
2 & 0 & \cdots & 0 \\
0 & \ddots &  \ddots  &  \vdots\\
 \vdots& \ddots & 2 &0 \\
0& \cdots  & 0& 1\\
\end{pmatrix}=
\begin{bmatrix}
0 & (-2)^n & 0   &\cdots &  \cdots& 0\\
(-2)^{n+1} & 0 & \ddots   & &  & \vdots\\
0&   (-2)^n&   \ddots &  \ddots & &  \vdots\\
\vdots & \ddots &  \ddots  &\ddots &-2 & 0 \\
\vdots&  &   \ddots &2 & {}^{}\hspace{0.2cm}0 & 2\\
0 & \cdots &  0  &0 & -2&0
\end{bmatrix} \in {\rm Asym}_n(\mathbf C) \, .
$$

We denote by $d_1,\ldots,d_n$  (with no reference to $B_n$ to simplify the notations) the coefficients of $d_{B_n}$, {\it i.e.}
 one has $d_i=2$ for $i=1,\ldots,n-1$ and $d_n=1$.

 \subsubsection{\bf The dilogarithmic identity $\boldsymbol{\big({\sf R}_{B_n}\big)}$.}
\label{SSS:ClusterWebs-type-B:dilog-identity}
According to Nakanishi's theorem ({\it cf.}\,Theorem \ref{Thm:Nakanishis-identity} above),  the $\mathcal Y$-cluster dilogarithm identity 
of type $B_n$ is written as follows: 
$$
{\Big({\sf R}_{B_n}\Big)}
\hspace{5cm}
\sum_{s=1}^{n(n+1)} d_{i_s} {\sf R}\big( Y_s\big) =N_{B_n}\, \frac{{}^{}\hspace{0.02cm}\pi^2}{6}
\hspace{7cm}{}^{} 
$$
for a certain non-negative integer $N_{B_n}$, which can verified to be equal to  $n(2n-1)$. 
In order to make this (a bit) more explicit, we introduce a new labelling of the $Y_s$'s by setting $Y_i^{[\kappa]}=Y_{i+\kappa n} $ 
for $i=1,\ldots,n$ and $\kappa=0,\ldots,n$. Then the dilogarithmic identity above can also be written  as
$$
\sum_{\kappa=0}^{n} \left ( \Bigg(\sum_{i=1}^{n-1}\,
2\,{\sf R}\Big( Y_i^{[\kappa]}\Big)\Bigg)
+  {\sf R}\Big( Y_n^{[\kappa]}\Big)\,
\right)= n(2n-1)
\, \frac{{}^{}\hspace{0.06cm}\pi^2}{6}\, . 
$$

The complete dilogarithmic AR  of $\boldsymbol{{\mathcal Y\hspace{-0.1cm}\mathcal W}_{\hspace{-0.05cm} B_3}}$
associated to ${\big({\sf R}_{B_n}\big)}$ will be denoted in the same way. \mk 

Before discussing the other ARS of this web for $n$ arbitrary,
it will be enlightening to first consider some explicit examples.

 \subsubsection{\bf The webs $\boldsymbol{{\mathcal Y\hspace{-0.1cm}\mathcal W}_{\hspace{-0.05cm} B_n}}$ for $\boldsymbol{n}$ small ($\boldsymbol{n\leq 4}$).}
\label{SSS:ClusterWebs-type-B:Y-cluster web-B234}

Since the $\mathcal Y$-cluster web of type $B_2$ coincides with 
the corresponding $\mathcal X$-cluster web
$\boldsymbol{{\mathcal X\hspace{-0.1cm}\mathcal W}}_{\hspace{-0.05cm} B_2}$
 which has been studied above in \S\ref{SS:cluster-webs:B2}, we refer the reader 
 to this subsection for a detailed description of its main features.

 \paragraph{\bf $\boldsymbol{{\mathcal Y\hspace{-0.1cm}\mathcal W}_{\hspace{-0.05cm} B_3}}$.}
\label{SSS:ClusterWebs-type-B:Y-cluster web-B3}

We can give $\boldsymbol{{\mathcal Y\hspace{-0.1cm}}\mathcal W}_{\hspace{-0.05cm} B_3}$ in explicit form: as an ordered web, one has
\begin{align*}
\boldsymbol{{\mathcal Y\hspace{-0.1cm}}\mathcal W}_{\hspace{-0.05cm} B_3}
=\boldsymbol{{\mathcal Y\hspace{-0.1cm}}\mathcal W}
\bigg( \,&   u_{{2}}
\, , \,  
{\frac {1+u_{{2}}}{u_{{1}}}}
\, , \, 
{\frac { \left( 1+u_{{2}} \right) ^{2}}{u_{{3}}}}
\, , \, 
{\frac { \left( u_{{1}}+1+u_{{2}} \right)  \left( {u_{{2}}}^{2}+2\,u_{{2}}+u_{{3}}+1 \right) }{u_{{1}}u_{{2}}u_{{3}}}}\, , 
\\ & 
{\frac {u_{{2}}u_{{1}}+u_{{1}}u_{{3}}+{u_{{2}}}^{2}+u_{{1}}+2\,u_{{2}}+u_{{3}}+1}{u_{{2}}u_{{3}}}} 
\, , \, 
{\frac { \left( u_{{2}}u_{{1}}+u_{{1}}u_{{3}}+{u_{{2}}}^{2}+u_{{1}}+2\,u_{{2}}+u_{{3}}+1 \right) ^{2}}{{u_{{1}}}^{2}{u_{{2}}}^{2}u_{{3}}}} 
\, , \\
&\, 
{\frac { \left( {u_{{1}}}^{2}u_{{3}}+{u_{{1}}}^{2}+2\,u_{{2}}u_{{1}}+2\,u_{{1}}u_{{3}}+{u_{{2}}}^{2}+2\,u_{{1}}+2\,u_{{2}}+u_{{3}}+1 \right)  \left( u_{{2}}+u_{{3}}+1 \right)
}{u_{{1}}{u_{{2}}}^{2}u_{{3}}}}
\, , \, \\ &\, 
{\frac {u_{{1}}u_{{3}}+u_{{1}}+u_{{2}}+u_{{3}}+1}{u_{{1}} u_{{2}}}}
\, , \,  
{\frac { \left( u_{{1}}u_{{3}}+u_{{1}}+u_{{2}}+u_{{3}}+1 \right) ^{2}}{{u_{{2}}}^{2}u_{{3}}}}
\, , \, 
{\frac { \left( u_{{3}}+1 \right)  \left( u_{{1}}+1 \right) }{u_{{2}}}}
\, , \,  
u_{{1}}
\, , \,  
u_{{3}}
\, \bigg)\, .
\end{align*}

By direct explicit computations, one obtains that the following assertions hold true: 
\begin{itemize}
\item The union of common leaves of this web coincides exactly with the cluster arrangement $\boldsymbol{Arr}_{B_3}$ and its associated ramification loci all are included in $\mathfrak B$ (hence $\boldsymbol{{\mathcal Y\hspace{-0.1cm}}\mathcal W}_{\hspace{-0.05cm} B_3}$ has   `polylogarithmic ramification' according to the terminology introduced in Definition \ref{Def:Polylogarithmic-Ramification});
\sk 
\item  One has $\rho^\bullet\big(\boldsymbol{{\mathcal Y\hspace{-0.1cm}}\mathcal W}_{\hspace{-0.05cm} B_3}\big)= \big(9,6,2\big) $ hence its virtual rank is 
$ \rho\big(\boldsymbol{{\mathcal Y\hspace{-0.1cm}}\mathcal W}_{\hspace{-0.05cm} B_3}\big)
  =17$;
\sk 
\item  ${\rm polrk}^\bullet\big(\boldsymbol{{\mathcal Y\hspace{-0.1cm}}\mathcal W}_{\hspace{-0.05cm} B_3}\big)=\big(12,1 \big)$ hence 
${\rm polrk}\big(\boldsymbol{{\mathcal Y\hspace{-0.1cm}}\mathcal W}_{\hspace{-0.05cm} B_3}\big)=
13$. In particular,   ${\big({\sf R}_{B_3}\big)}$ is the unique weight 2 polylogarithmic AR of $\boldsymbol{{\mathcal Y\hspace{-0.1cm}}\mathcal W}_{\hspace{-0.05cm} B_3}$ (up to  multiplication by a non zero scalar);
\sk 
\item  This web has rank  ${\rm rk}\big(\boldsymbol{{\mathcal Y\hspace{-0.1cm}}\mathcal W}_{\hspace{-0.05cm} B_3}\big)=14$;
\sk 
\item Thus  $ {\rm polrk}\big(\boldsymbol{{\mathcal Y\hspace{-0.1cm}}\mathcal W}_{\hspace{-0.05cm} B_3}\big)<{\rm rk}\big(\boldsymbol{{\mathcal Y\hspace{-0.1cm}}\mathcal W}_{\hspace{-0.05cm} B_3}\big)< \rho\big(\boldsymbol{{\mathcal Y\hspace{-0.1cm}}\mathcal W}_{\hspace{-0.05cm} B_3}\big)$ hence $\boldsymbol{{\mathcal Y\hspace{-0.1cm}}\mathcal W}_{\hspace{-0.05cm} B_3}$ is not AMP;
\sk 
\item Because  $ {\rm polrk}\big(\boldsymbol{{\mathcal Y\hspace{-0.1cm}}\mathcal W}_{\hspace{-0.05cm} B_3}\big)=13<14={\rm rk}\big(\boldsymbol{{\mathcal Y\hspace{-0.1cm}}\mathcal W}_{\hspace{-0.05cm} B_3}\big)$, the space of polylogarithmic ARs of $\boldsymbol{{\mathcal Y\hspace{-0.1cm}}\mathcal W}_{\hspace{-0.05cm} B_3}$ has codimension 1 in 
$\boldsymbol{\mathcal A}\big( 
\boldsymbol{{\mathcal Y\hspace{-0.1cm}}\mathcal W}_{\hspace{-0.05cm} B_3}
\big)$; 
\sk 
\item  If ${\sf A}$ stands for the function defined by ${\sf A}(u)={\rm Arctan}(\sqrt{u})$ for any $u>0$, then the following functional relation holds true indentically on the positive orthant: 
$$
{\Big({\sf A}_{B_3}\Big)}
\hspace{2cm}
{\sf A}\Big( Y_3\Big)+
 {\sf A}\Big( Y_6\Big)+ {\sf A}\Big( Y_9\Big)+ {\sf A}\Big( Y_{12}\Big)
=\sum_{\nu=1}^4 {\sf A}\Big( Y_{3\nu }\Big)=\frac{3}{2}\,\pi\, ;
\hspace{7cm}{}^{} 
$$
\item  Still denoting by ${\big({\sf A}_{B_3}\big)}$ the non-polylogarithmic AR associated to the preceding functional 
 identity, we thus have the following decomposition in direct sum of the space of ARs of 
$\boldsymbol{{\mathcal Y\hspace{-0.1cm}}\mathcal W}_{\hspace{-0.05cm} B_3}$ (to be compared with \eqref{Eq:AR-XWB2}): 
\begin{equation}
\label{Eq:AR-YWB3}
\boldsymbol{\mathcal A}\Big(\boldsymbol{\mathcal Y\hspace{-0.1cm}\mathcal W}_{\hspace{-0.05cm}B_3}\Big) = 
\boldsymbol{LogAR}\Big(\boldsymbol{{\mathcal Y\hspace{-0.1cm}}\mathcal W}_{\hspace{-0.05cm} B_3}\Big) 
 \, \oplus \, 
 {
\Big\langle
\, {\sf R}_{B_3} \,
 \Big\rangle  }
 \oplus  \,
 \big\langle
\, {\sf A}_{B_3}\, 
 \big\rangle\, . 
\end{equation}
\end{itemize}

Comparing the results above to the corresponding ones for $\boldsymbol{\mathcal Y\hspace{-0.1cm}\mathcal W}_{\hspace{-0.05cm}B_2}$, we notice that: 
\begin{itemize}
\item[$-$] contrarily to $\boldsymbol{\mathcal Y\hspace{-0.1cm}\mathcal W}_{\hspace{-0.05cm}B_2}$ which is of maximal rank, $\boldsymbol{\mathcal Y\hspace{-0.1cm}\mathcal W}_{\hspace{-0.05cm}B_3}$ is not AMP; 
\sk
\item[$-$] the space of weight two dilogarithmic ARs of $\boldsymbol{\mathcal Y\hspace{-0.1cm}\mathcal W}_{\hspace{-0.05cm}B_3}$ is spanned 
by $({\sf R}_{B_3})$ whereas the one in case $B_2$  is of dimension 2, spanned by  $({\sf R}_{B_2})$ and $({\sf S}_{B_2})$ ({\it cf.}\,\S\ref{SS:cluster-webs:B2} above).
\end{itemize}



 \paragraph{\bf $\boldsymbol{{\mathcal Y\hspace{-0.1cm}\mathcal W}_{\hspace{-0.05cm} B_4}}$.}
\label{SSS:ClusterWebs-type-B:Y-cluster web-B4}

The  $\mathcal Y$-cluster web of type $B_4$  is a 20-web in four variables which 
shares some properties of both 
$\boldsymbol{\mathcal Y\hspace{-0.1cm}\mathcal W}_{\hspace{-0.05cm} B_2}$ and 
$\boldsymbol{\mathcal Y\hspace{-0.1cm}\mathcal W}_{\hspace{-0.05cm} B_3}$. 
\mk 

By direct explicit computations, one obtains that the following assertions hold true: 
\begin{itemize}
\item The union of common leaves of 
$\boldsymbol{{\mathcal Y\hspace{-0.1cm}}\mathcal W}_{\hspace{-0.05cm} B_4}$ 
 coincides exactly with the cluster arrangement $\boldsymbol{Arr}_{B_4}$ ({\it i.e.}\,one has $\Sigma^c(\boldsymbol{\mathcal Y\hspace{-0.1cm}\mathcal W}_{\hspace{-0.05cm} B_4})=\boldsymbol{Arr}_{B_4}$) thus this web  has   `polylogarithmic ramification';
\sk 
\item The usual, virtual or polylogarithmic rank(s) of this web are given by 
\begin{align*}
  \rho^\bullet\big(\boldsymbol{{\mathcal Y\hspace{-0.1cm}}\mathcal W}_{\hspace{-0.05cm} B_4}\big)=& \, \big(16,10,2\big)  && 
  \hspace{0.65cm}\rho\big(\boldsymbol{{\mathcal Y\hspace{-0.1cm}}\mathcal W}_{\hspace{-0.05cm} B_4}\big)
  =28   \\ 
   {\rm polrk}^\bullet\Big(\boldsymbol{{\mathcal Y\hspace{-0.1cm}}\mathcal W}_{\hspace{-0.05cm} B_4}\Big)=&\, \Big(20,2 \Big)
&&{\rm polrk}\Big(\boldsymbol{{\mathcal Y\hspace{-0.1cm}}\mathcal W}_{\hspace{-0.05cm} B_4}\Big)=
22
\qquad \mbox{ and } \qquad   {\rm rk}\Big(\boldsymbol{{\mathcal Y\hspace{-0.1cm}}\mathcal W}_{\hspace{-0.05cm} B_4}\Big)=23\, ;
\end{align*}

\item  Thus  $\boldsymbol{{\mathcal Y\hspace{-0.1cm}}\mathcal W}_{\hspace{-0.05cm} B_4}$ is not AMP ;
\sk 
\item  
One has ${\rm polrk}^2\big(\boldsymbol{{\mathcal Y\hspace{-0.1cm}}\mathcal W}_{\hspace{-0.05cm} B_4}\big)=2$  hence the space  of weight 2 polylogarithmic ARs is 
  not solely spanned by $\big({\sf R}_{B_4}\big)$. 

Let $\delta=\delta_{B_4}=(\delta_s)_{s=1}^{20}$ be the 20-tuple obtained by 
 concatenating in a row $5$ copies of 
$(2, 2, -2, -1)$.  Then the following relation holds true identically on $(\mathbf R_{>0})^4$: 
$$
{\Big({\sf S}_{B_4}\Big)}
\hspace{2cm}
\sum_{s=1}^{20} \delta_{s}\, {\sf S}\big( Y_s\big) =
\sum_{\kappa=0}^{4} \Bigg( 
2\,{\sf S}\Big( Y_1^{[\kappa]}\Big)+
2\,{\sf S}\Big( Y_2^{[\kappa]}\Big)
-2 \,{\sf S}\Big( Y_3^{[\kappa]}\Big)
-\,{\sf S}\Big( Y_4^{[\kappa]}\Big)
\,
\Bigg)=
0
\hspace{7cm}{}^{} 
$$ 
where ${\sf S}$ stands for the 
symmetric dilogarithm defined in \eqref{Eq:Dilog-Symmetric-S} (see also \eqref{Eq:sf-S} below).  Then the two weight 2 polylogarithmic ARs  $\big({\sf R}_{B_4} \big)$ and 
${\big({\sf S}_{B_4}\big)}$ form a basis of $\boldsymbol{DilogAR}\big(\boldsymbol{{\mathcal Y\hspace{-0.1cm}}\mathcal W}_{\hspace{-0.05cm} B_3}\big) $;

\item  The space of the polylogarithmic ARs of 
$\boldsymbol{{\mathcal Y\hspace{-0.1cm}}\mathcal W}_{\hspace{-0.05cm} B_4}$ 
is 1-codimensional in  $\boldsymbol{\mathcal A}\big(\boldsymbol{\mathcal Y\hspace{-0.1cm}\mathcal W}_{\hspace{-0.05cm}B_4}\big)$. As a supplementary space, one can take the line 
spanned by the non-polylogarithmic AR associated to the following functional relation which holds true identically on the positive orthant:
$$
{\Big({\sf A}_{B_4}\Big)}
\hspace{2cm}
{\sf A}\Big( Y_{4}\Big)+{\sf A}\Big( Y_{8 }\Big)+{\sf A}\Big( Y_{12 }\Big)+{\sf A}\Big( Y_{16 }\Big)+{\sf A}\Big( Y_{20 }\Big)=
\sum_{k=1}^5 {\sf A}\Big( Y_{4k }\Big)=2\,\pi\,;
\hspace{7cm}{}^{} 
$$
\item  Consequently,  we have the following decomposition in direct sum of the space of ARs:
\begin{equation}
\label{Eq:AR-YWB4}
\boldsymbol{\mathcal A}\Big(\boldsymbol{\mathcal Y\hspace{-0.1cm}\mathcal W}_{\hspace{-0.05cm}B_4}\Big) = 
\boldsymbol{Lo\hspace{-0.03cm}g\hspace{-0.03cm}A\hspace{-0.03cm}R}\Big(\boldsymbol{{\mathcal Y\hspace{-0.1cm}}\mathcal W}_{\hspace{-0.05cm} B_4}\Big) 
 \, \oplus \, 
 {
\Big\langle
\, {\sf R}_{B_4} \,\,  \, {\sf S}_{B_4}
 \Big\rangle  }
 \oplus  \,
 \big\langle
\, {\sf A}_{B_4}\, 
 \big\rangle\, . 
\end{equation}
\end{itemize}
\begin{center}
$\star$
\end{center}

Comparing the properties of the webs $\boldsymbol{{\mathcal Y\hspace{-0.1cm}}\mathcal W}_{\hspace{-0.05cm} B_n}$ for $n\leq 4$ allows us to make the following observations:  
\begin{itemize}
\item[$-$] compared to the $\mathcal Y$-cluster webs associated to $B_3$ and $B_4$, the one associated to $B_2$ seems a bit particular. Indeed, contrarily to the two former, (1) it is of maximal rank (hence AMP) and (2) it carries a rational AR;\sk 

As  will be clear further in \S\ref{SS:ClusterWebs-type-C} where we will  discuss the $\mathcal Y$-cluster webs of type $C$, (2) reflects the fact that $\boldsymbol{{\mathcal Y\hspace{-0.1cm}}\mathcal W}_{\hspace{-0.05cm} B_2}$ is also the $\mathcal Y$-cluster web of type $C_2$ (which follows from the fact that  $B_2$ and $C_2$ are isomorphic as Dynkin diagrams). And the fact that $\boldsymbol{{\mathcal Y\hspace{-0.1cm}}\mathcal W}_{\hspace{-0.05cm} B_2}$ is of maximal rank comes from this extra rational AR. 
\sk

\item[$-$] Considering only the weight 2 polylogarithmic ARs, we observe that, contrarily to $\boldsymbol{{\mathcal Y\hspace{-0.1cm}}\mathcal W}_{\hspace{-0.05cm} B_3}$ for which  
$\boldsymbol{Dilo\hspace{-0.03cm}g\hspace{-0.03cm}A\hspace{-0.03cm}R}\Big(\boldsymbol{{\mathcal Y\hspace{-0.1cm}}\mathcal W}_{\hspace{-0.05cm} B_3}\Big) $ is spanned by $({\sf R}_{B_3})$, the webs
$\boldsymbol{{\mathcal Y\hspace{-0.1cm}}\mathcal W}_{\hspace{-0.05cm} B_n}$ 
for $n=2,4$ both    carry two linearly independent dilogarithmic ARs,  
$({\sf R}_{B_n})$ and the symmetric one $({\sf S}_{B_n})$. 
\end{itemize}

These two observations suggest a general pattern for the ranks and the ARs of the webs of the family $\boldsymbol{{\mathcal Y\hspace{-0.1cm}}\mathcal W}_{\hspace{-0.05cm} B_n}$ that we have explored (by means of direct computations) in many cases (say $n\leq 10$) and that we are going to present in the next subsection. 

The case when $n=2$ being particular since an extra AR induces the maximality of its rank (this presumably thanks to the coincidence between $B_2$ and $C_2$), this case will left aside from now on. We thus assume that $n\geq 3$ in what follows.

 \subsubsection{\bf A  conjectural pattern 
 about the ranks and  the abelian relations of $\boldsymbol{{\mathcal Y\hspace{-0.1cm}\mathcal W}_{\hspace{-0.05cm} B_n}}$.}
\label{SSS:ClusterWebs-type-B:Y-cluster web}
One has  $d^{\mathcal Y}_{B_\bullet}=n(n+1)$ so  $\boldsymbol{{\mathcal Y\hspace{-0.1cm}\mathcal W}}_{\hspace{-0.05cm} B_n}$ is a $n(n+1)$-web in $n$ variables. Below we use the definition \eqref{Eq:Web-YWBn} and use the $Y_s$ (for $s=1,\ldots,n(n+1)$) as first integrals for it.   We will also use 
the notation $Y_i^{[\kappa]}=Y_{i+\kappa n}$ with $i=1,\ldots,n$ and $\kappa=0,\ldots,n$.
\mk

It is by studying many cases through explicit computations that we have been able to draw the general outline we are  going to present about the generalization for $n\geq 3$ arbitrary.  We give below a list of facts verified for the webs 
$
\boldsymbol{{\mathcal Y\hspace{-0.1cm}\mathcal W}}_{\hspace{-0.05cm} B_n}
$ for $n$ ranging from 3 to 8. 
\mk

{\bf Space of common leaves.} One has $\Sigma^c(\boldsymbol{{\mathcal Y\hspace{-0.1cm}\mathcal W}}_{\hspace{-0.05cm} B_n})=\boldsymbol{Arr}_{B_n}$ for all $n\leq 8$. 
\mk 

{\bf Ranks.} Using the technics discussed in \S\ref{Para:DeterminingARs-and-Rank}, we obtain the following values for the (virtual, polylogarithmic and usual) ranks  of the webs  $
\boldsymbol{{\mathcal Y\hspace{-0.1cm}\mathcal W}}_{\hspace{-0.05cm} B_n}$'s 
considered:  
 \begin{align*}
  \rho^\bullet(\boldsymbol{{\mathcal Y\hspace{-0.1cm}\mathcal W}}_{\hspace{-0.05cm} B_3})=& \big(9,6,2\big)  && \rho(\boldsymbol{{\mathcal Y\hspace{-0.1cm}\mathcal W}}_{\hspace{-0.05cm} B_3})=17 && 
   {\rm polrk}^\bullet\big(\boldsymbol{{\mathcal Y\hspace{-0.1cm}\mathcal W}}_{\hspace{-0.05cm} B_3}\big)=\big(12,1 \big)
   &&  {\rm rk}\big(\boldsymbol{{\mathcal Y\hspace{-0.1cm}\mathcal W}}_{\hspace{-0.05cm} B_3}\big)=14
  \\
    \rho^\bullet(\boldsymbol{{\mathcal Y\hspace{-0.1cm}\mathcal W}}_{\hspace{-0.05cm} B_4})=& \big(\, 16\, , \, 10\, , \, 2 \, \big)  && \rho(\boldsymbol{{\mathcal Y\hspace{-0.1cm}\mathcal W}}_{\hspace{-0.05cm} B_4})= 28&& 
   {\rm polrk}^\bullet\big(\boldsymbol{{\mathcal Y\hspace{-0.1cm}\mathcal W}}_{\hspace{-0.05cm} B_4}\big)=\big(\,   20\, , \, 2      \,  \big)
   &&  {\rm rk}\big(\boldsymbol{{\mathcal Y\hspace{-0.1cm}\mathcal W}}_{\hspace{-0.05cm} B_4}\big)= 23
  \\
      \rho^\bullet(\boldsymbol{{\mathcal Y\hspace{-0.1cm}\mathcal W}}_{\hspace{-0.05cm} B_5})=& \big(\, 25\, , \, 15\,,\,  2\, \big)  && \rho(\boldsymbol{{\mathcal Y\hspace{-0.1cm}\mathcal W}}_{\hspace{-0.05cm} B_5})= 42 && 
   {\rm polrk}^\bullet\big(\boldsymbol{{\mathcal Y\hspace{-0.1cm}\mathcal W}}_{\hspace{-0.05cm} B_5}\big)=\big(\,  30\, , \, 1       \,  \big)
   &&  {\rm rk}\big(\boldsymbol{{\mathcal Y\hspace{-0.1cm}\mathcal W}}_{\hspace{-0.05cm} B_5}\big)= 32
  \\
      \rho^\bullet(\boldsymbol{{\mathcal Y\hspace{-0.1cm}\mathcal W}}_{\hspace{-0.05cm} B_6})=& \big(\, 36\, , \, 21\, , \, 2 \, \big)  && \rho(\boldsymbol{{\mathcal Y\hspace{-0.1cm}\mathcal W}}_{\hspace{-0.05cm} B_6})= 59&& 
   {\rm polrk}^\bullet\big(\boldsymbol{{\mathcal Y\hspace{-0.1cm}\mathcal W}}_{\hspace{-0.05cm} B_6}\big)=\big(\,     42\, , \, 2     \,  \big)
   &&  {\rm rk}\big(\boldsymbol{{\mathcal Y\hspace{-0.1cm}\mathcal W}}_{\hspace{-0.05cm} B_6}\big)= 45
  \\
      \rho^\bullet(\boldsymbol{{\mathcal Y\hspace{-0.1cm}\mathcal W}}_{\hspace{-0.05cm} B_7})=& \big(\, 49\, ,\, 28\, 2\, \big)  && \rho(\boldsymbol{{\mathcal Y\hspace{-0.1cm}\mathcal W}}_{\hspace{-0.05cm} B_7})= 79&& 
   {\rm polrk}^\bullet\big(\boldsymbol{{\mathcal Y\hspace{-0.1cm}\mathcal W}}_{\hspace{-0.05cm} B_7}\big)=\big(\,  56\, , \, 1        \,  \big)
   &&  {\rm rk}\big(\boldsymbol{{\mathcal Y\hspace{-0.1cm}\mathcal W}}_{\hspace{-0.05cm} B_7}\big)= 58
  \\
      \rho^\bullet(\boldsymbol{{\mathcal Y\hspace{-0.1cm}\mathcal W}}_{\hspace{-0.05cm} B_8})=& \big(\,64\, ,\, 36\, , \, 2 \, \big)  && \rho(\boldsymbol{{\mathcal Y\hspace{-0.1cm}\mathcal W}}_{\hspace{-0.05cm} B_8})= 102&& 
   {\rm polrk}^\bullet\big(\boldsymbol{{\mathcal Y\hspace{-0.1cm}\mathcal W}}_{\hspace{-0.05cm} B_8}\big)=\big(\,  72\, , \, 2       \,  \big)
   &&  {\rm rk}\big(\boldsymbol{{\mathcal Y\hspace{-0.1cm}\mathcal W}}_{\hspace{-0.05cm} B_8}\big)= 75 \, . 
\end{align*}

We remark that the following relations are satisfied:
\begin{itemize}
\item[$-$] the non trivial virtual ranks are given by 
\begin{align*}
\rho^1\big(\boldsymbol{{\mathcal Y\hspace{-0.1cm}\mathcal W}}_{\hspace{-0.05cm} B_n}\big)= &\, d^{\mathcal Y}_{B_n}-n=n^2\, ,\\ 
\rho^2\big(\boldsymbol{{\mathcal Y\hspace{-0.1cm}\mathcal W}}_{\hspace{-0.05cm} B_n}\big)=& \, {n(n+1)}/{2}\, ,  \\
\mbox{and} \quad 
 \rho^3\big(\boldsymbol{{\mathcal Y\hspace{-0.1cm}\mathcal W}}_{\hspace{-0.05cm} B_n}\big)=& \,2\,  ; 
\end{align*}
\item[$-$] in what concerns the non zero polylogarithmic ranks, these are 
\begin{align*}
{\rm polrk}^1\big(\boldsymbol{{\mathcal Y\hspace{-0.1cm}\mathcal W}}_{\hspace{-0.05cm} B_n}\big)=  &\, d^{\mathcal Y}_{B_n} =n(n+1) \\
 \mbox{and} 
\quad 
{\rm polrk}^2\big(\boldsymbol{{\mathcal Y\hspace{-0.1cm}\mathcal W}}_{\hspace{-0.05cm} B_n}\big)= &\, 
1\,  (n \,\mbox{odd});  \mbox{ or }\, 2 
\,   (n \,\mbox{even})\, ;
  \end{align*} 
\item[$-$] finally, the total rank and the polylogarithmic one are related by the following formula:  
\begin{equation}
\label{Eq:Bn-rkYWBn-polrkYWBn}
{\rm rk}\big(\boldsymbol{{\mathcal Y\hspace{-0.1cm}\mathcal W}}_{\hspace{-0.05cm} B_n}\big)={\rm polrk}\big(\boldsymbol{{\mathcal Y\hspace{-0.1cm}\mathcal W}}_{\hspace{-0.05cm} B_n}\big)+1\, .
\end{equation}
\end{itemize}

{\bf Polyogarithmic ARs of weight 2.} 
As explained above, the  dilogarithmic identity $({\sf R}_{B_n})$ always gives rise to a dilogarithmic AR (denoted in the same way) for $\boldsymbol{{\mathcal Y\hspace{-0.1cm}\mathcal W}}_{\hspace{-0.05cm} B_n}$. \sk 

When $n$ is odd, this is the unique dilogarithmic AR (up to multiplication by a non zero scalar).  This is not true when $n$ is even (that is for $n=2,4,6,8$) since in this case 
there exists another such AR which is not a multiple of $({\sf R}_{B_n})$.  For each such $n$, this AR corresponds to a functional identity, denoted by 
$({\sf S}_{B_n})$, 
 satisfied by the 
{\bf symmetric (cluster) dilogarithm} ${\sf S}$ the definition of which is recalled below: it is the weight 2 cluster iterated integral, with symmetric symbol $\frac{1}{2}(01+10)$ with $0$ as chosen  base point: one has 
\begin{equation}
\label{Eq:sf-S}
{\sf S}(x)=\frac{1}{2}\left( \, \int_{0}^x \frac{{\rm Log}(1+u)}{u}du+
 \int_{0}^x \frac{{\rm Log}(u)}{1+u}du\, 
\right)=\frac{1}{2}{\rm Log}(x){\rm Log}(1+x)
\end{equation}
for any $x\in \mathbf R_{>0}$. \sk 

For $n$ even, set $m=n/2$ and denote by $\delta=\delta_{B_n}=\big(\delta_s\big)_{s=1}^{n(n+1)}$ 
 the $n(n+1)$-tuple obtained by concatenating in a row $n+1$ copies of 
the $n$-tuple  $d'=d'_{B_n}=(d_i')_{i=1}^n$ defined by 
$$
d'=d'_{B_n}=(d_i')_{i=1}^n=
 \Big(\,\overbrace{\, 2\, ,\, \ldots\, ,\, 2\, }^{m\,\mbox{ copies}}\, , \, 
\underbrace{\, -2\, ,\, \ldots\, , -2\, }_{m-1\,\mbox{ copies}}\,,\,-1\,\Big)\, . 
$$

The generalization to $B_n=B_{2m}$ of the symmetric dilogarithmic identities
$({\sf S}_{B_2})$ and $({\sf S}_{B_4})$ 
is 
$$
{\Big({\sf S}_{B_n}\Big)}
\hspace{5cm}
\sum_{s=1}^{n(n+1)} \delta_{s} \,{\sf S}\big( Y_s\big) =0\, , 
\hspace{7cm}{}^{} 
$$
a functional identity on $(\mathbf R_{>0})^n$ 
which can as well be written  (a bit) more explicitly as 
$$
\sum_{\kappa=0}^{n} \left[ 2 \,\Bigg( \sum_{i=1}^{m}{\sf S}\big( Y_i^{[\kappa]}\Big)
\Bigg) 
-2\,\Bigg(
\sum_{j=1}^{m-1}{\sf S}\Big( Y_{m+j}^{[\kappa]}\Big)\Bigg)
- {\sf S}\Big( Y_n^{[\kappa]}\Big)\,
\right]= 0\, . 
$$

{\bf The non polylogarithmic AR.} 
 The fact that the rank and the polylogarithmic rank 
of $\boldsymbol{{\mathcal Y\hspace{-0.1cm}\mathcal W}}_{\hspace{-0.05cm} B_n}$ differ by one ({\it cf.}\,\eqref{Eq:Bn-rkYWBn-polrkYWBn}) is the manifestation of the existence of an extra non-polylogarithmic AR, whose components are multiples of the function $A$ defined in \eqref{Eq:Arctan-sqrt} and 
which is the one associated to the following identity which 
holds true identically on $(\mathbf R_{>0})^4$: 
$$
{\Big({\sf A}_{B_n}\Big)}
\hspace{5cm}
\sum_{k=1}^{n+1} A\Big( Y_{kn }\Big)=n\,\frac{\pi}{2}\, . 
\hspace{7cm}{}^{} 
$$
\begin{center}
$\star$
\end{center}
All the assertions above are satisfied for any $n$ less than or equal to 8 (and bigger than 2). But we have taken care to state things  in order that everything makes sense, at least formally, for $n\geq 3$ arbitrary.  We believe that everything indeed holds true in full generality: \mk
%

\label{Conj:YWBn-page}
{\bf Conjecture $\boldsymbol{{\mathcal Y\hspace{-0.1cm} \mathcal W}_{B_n}}$.} 
{\it 
Let $n$ be an integer bigger than or equal to 3. 

\begin{enumerate}
\item {\bf ${\big[}$Ramification\big].}
 The ramification of $\boldsymbol{{\mathcal Y\hspace{-0.1cm}\mathcal W}}_{\hspace{-0.05cm} B_n}$ is polylogarithmic and 
the  variety of common leaves $\Sigma^c(\boldsymbol{{\mathcal Y\hspace{-0.1cm}\mathcal W}}_{\hspace{-0.05cm} B_n})$ 
 coincides with the corresponding cluster arrangement:  one has 
 $$
 \Sigma^c\Big(\boldsymbol{{\mathcal Y\hspace{-0.1cm}\mathcal W}}_{\hspace{-0.05cm} B_n}\Big)=\boldsymbol{Arr}_{B_n}\, .
 $$
\item  {\bf ${\big[}$Virtual rank(s)\big].} The virtual ranks  
of $\boldsymbol{{\mathcal Y\hspace{-0.1cm}\mathcal W}}_{\hspace{-0.05cm} B_n}$ are  given by  
\begin{align*}
\rho^\bullet\Big(\boldsymbol{{\mathcal Y\hspace{-0.1cm}\mathcal W}}_{\hspace{-0.05cm} B_n}\Big)= & \,  \Big(\,  n^2\, , \, {n(n+1)}/{2}\, , \, 2\,  \Big) \\
\mbox{and} \quad 
\rho\Big(
\boldsymbol{{\mathcal Y\hspace{-0.1cm}\mathcal W}}_{\hspace{-0.05cm} B_n}
\Big)= & \, 
{\big(n-1\big)\big(3n+4\big)}/{2 }+4 \, . 
\end{align*}
 \item {\bf ${\big[}$Polylogarithmic rank(s)\big].}  Regarding the polylogarithmic ranks of $\boldsymbol{\mathcal Y\hspace{-0.1cm} \mathcal W}_{B_n}$, one has 
\begin{align*}
{\rm polrk}^\bullet\Big(
\boldsymbol{{\mathcal Y\hspace{-0.1cm}\mathcal W}}_{\hspace{-0.05cm} B_n}
\Big)= &\,  \Big( \, n(n+1)\, , \, \beta_n \, \Big) 
\qquad \mbox{with}\qquad \beta_n=\begin{cases} \hspace{0.1cm} 1 \quad \mbox{for }\, n\, \mbox{ odd}\\
\hspace{0.1cm}  2 \quad \mbox{for }\, n\, \mbox{ even} 
\end{cases}, \\
\mbox{therefore}\quad 
{\rm polrk}\Big(
\boldsymbol{{\mathcal Y\hspace{-0.1cm}\mathcal W}}_{\hspace{-0.05cm} B_n}
\Big) =  & \, n(n+1)+\beta_n\, .
\end{align*}
\item  {\bf ${\big[}$Dilogarithmic AR(s){\big]}.} 
Whether $n$ is even or odd, the web $\boldsymbol{\mathcal Y\hspace{-0.1cm} \mathcal W}_{B_n}$ admits an antisymmetric dilogarithmic AR 
which corresponds to the  functional identity $({\sf R}_{B_n})$.\sk 

 If $n$ is odd, this AR is a basis of $\boldsymbol{DilogAR}\big(\boldsymbol{\mathcal Y\hspace{-0.1cm} \mathcal W}_{B_n}\big)$ which is 1-dimensional. When $n$ is even, $\boldsymbol{\mathcal Y\hspace{-0.1cm} \mathcal W}_{B_n}$ carries another  (symmetric) dilogarithmic AR corresponding to the  functional equation $({\sf S}_{B_n})$.  
 In this case, the two ARs $({\sf R}_{B_n})$ and $({\sf S}_{B_n})$
  form a basis of the space of dilogarithmic ARs of $\boldsymbol{\mathcal Y\hspace{-0.1cm} \mathcal W}_{B_n}$. Thus,  one has 
  \begin{equation*}
\boldsymbol{Dilo\hspace{-0.03cm}g\hspace{-0.03cm}A\hspace{-0.03cm}R}\Big(\boldsymbol{{\mathcal Y\hspace{-0.1cm}}\mathcal W}_{\hspace{-0.05cm} B_n}\Big) =\begin{cases} 
\hspace{0.2cm}  \Big\langle \, 
 {\sf R}_{B_n}\,
 \Big\rangle  \hspace{1.1cm} \mbox{when } n  \mbox{ is odd,}
 \vspace{0.1cm}\\
 \hspace{0.2cm}  \Big\langle \, 
 {\sf R}_{B_n} \,  , \, 
 \, {\sf S}_{B_n}\,
 \Big\rangle \hspace{0.2cm} \mbox{when } n  \mbox{ is even}\, .
 \end{cases} 
 \end{equation*}
%
%
\item  {\bf ${\big[}$The non-polylogarithmic AR\big].}
 The following functional equation
\begin{equation}
\label{Eq:A(YWBn)}
\sum_{k=1}^{n+1} {\sf A}\Big(Y_{kn}\Big)=
\sum_{\nu=0}^{n} {\sf A} \Big(Y_{n}^{[\nu]}\Big)
=n\,\frac{{}\, \pi}{2}
\end{equation}
is identically satisfied on $(\mathbf R_{>0})^n$. The associated AR of $\boldsymbol{\mathcal Y\hspace{-0.1cm}\mathcal W}_{\hspace{-0.05cm}B_n}$ is not of polylogarithmic type hence is linearly independent from the ARs described in 4. just above.
\mk
\item {\bf ${\big[}$Rank and basis of ARs{\big]}.} 
  One has 
    $$
  {\rm rk}\big(\boldsymbol{\mathcal Y\hspace{-0.1cm} \mathcal W}_{\hspace{-0.05cm}B_n}\big)={\rm polrk}\big(\boldsymbol{\mathcal Y\hspace{-0.1cm} \mathcal W}_{\hspace{-0.05cm}B_n}\Big)+1=n(n+1)+\beta_n+1
  $$
 hence the following decomposition in direct sum holds true:
\begin{equation}
\label{Eq:AR-YWBn}
\boldsymbol{\mathcal A}\Big(\boldsymbol{\mathcal Y\hspace{-0.1cm}\mathcal W}_{\hspace{-0.05cm}B_n}\Big) = 
\boldsymbol{Lo\hspace{-0.03cm}g\hspace{-0.03cm}A\hspace{-0.03cm}R}\Big(\boldsymbol{{\mathcal Y\hspace{-0.1cm}}\mathcal W}_{\hspace{-0.05cm} B_n}\Big) 
 \, \oplus \, 
\boldsymbol{Dilo\hspace{-0.03cm}g\hspace{-0.03cm}A\hspace{-0.03cm}R}\Big(\boldsymbol{{\mathcal Y\hspace{-0.1cm}}\mathcal W}_{\hspace{-0.05cm} B_n}\Big) 
\,  \oplus  \,
 \big\langle
\, {\sf A}_{B_n}\, 
 \big\rangle\, . 
\end{equation}
\end{enumerate}}

Assuming that this conjecture is true for any $n\geq 3$, we obtain that contrarily to the corresponding webs in type $A$, the webs $\boldsymbol{\mathcal Y\hspace{-0.1cm} \mathcal W}_{\hspace{-0.05cm}B_n}$ are not AMP. If this 
 may seem a little disappointing at first sight, we really believe that the list of statements of this conjecture are interesting since they give a precise description of the ARs (most of these (if not all, see below) being of polylogarithmic type) of a any member of a whole countable series of webs.   \mk 
 
 Given $n\geq 3$, the non logarithmic ARs of $\boldsymbol{\mathcal Y\hspace{-0.1cm} \mathcal W}_{\hspace{-0.05cm}B_n}$, if they are as indicated in the conjecture (which is known to hold true for $n\leq 8$) are interesting for themselves. This is the subject of the next subsection.

 \subsubsection{\bf A few remark about  Conjecture $\boldsymbol{{\mathcal Y\hspace{-0.1cm}\mathcal W}_{\hspace{-0.05cm} B_n}}$.}
\label{SSS:Remarks-Conj-YWBn}
Here we discuss some interesting  AFEs which would be satisfied assuming that the supplementary conjectured AR(s) of  $\boldsymbol{\mathcal Y\hspace{-0.1cm}\mathcal W}_{\hspace{-0.05cm} B_n}$ ($({\sf A}_{B_n})$ for any $n$, and $({\sf S}_{B_n})$ when $n$ is even) hold true. 

 \paragraph{The case when $\boldsymbol{n}$ is even.} An interesting feature of the 
cluster web of type $B_2$ is that, modulo equivalence with Newman 6-web ({\it cf.}\,\S\ref{SS:cluster-webs:B2}), it carries
 an AR whose components all are multiple of the classical bilogarithm 
 ${\bf L}{\rm i}_2$ (see \S\ref{Par:Newman-equation}). Denoting (a bit abusively)  this AR by $\mathcal N_6$, it can be verified that it coincides with the difference between the antisymmetric and the symmetric ARs of 
 $\boldsymbol{\mathcal Y\hspace{-0.1cm}\mathcal W}_{\hspace{-0.05cm} B_2}$: one has 
 $$ 
 \mathcal N_6={\sf R}_{B_2}-{\sf S}_{B_2}\in \boldsymbol{\mathcal A}\big( 
 \boldsymbol{\mathcal Y\hspace{-0.1cm}\mathcal W}_{\hspace{-0.05cm} B_2}
 \big)\, .
 $$

It is natural to wonder whether  a similar phenomenon does occur when the space of dilogarithmic ARs has dimension 2, that is (conjecturally) when $n$ is even: in this case, is the difference 
${\sf R}_{B_2n}-{\sf S}_{B_n}$ equivalent to an AFE whose components all are multiples of the classical dilogarithm?  We discuss explicitly the case $B_4$ below, guessing that things go the same way for any $B_{2m}$ with $m\geq 1$.\mk 

The difference ${\sf R}_{B_4}-{\sf S}_{B_4}$
 is written as follows in explicit form: 
$$
\sum_{k=0}^{4}  \Big(\, 2\, L_{01}\big( Y_{1+4k}\big)  + 2\, L_{01} \big( Y_{2+4k}\big) 
-2\, L_{10}\big( Y_{3+4k}\big) 
\, -\, L_{10}\big( Y_{4+4k}\big) 
\Big)= 28\,\frac{{}^{}\,\pi^2}{6}\, .
$$
where $L_{10}$ and $L_{01}$ stand for the two weight 2 iterated integrals defined in 
\eqref{Eq:L01-L10}. 
Since 
$$
L_{01}(x)=- {\bf L}{\rm i}_2(-x) \qquad \mbox{and}\qquad 
 L_{10}(x)=- {\bf L}{\rm i}_2 \big(1+x\big)-i\pi\,{\rm Log}(1+x)+\pi^2/{6}$$
for any positive $x\in \mathbf R$, we deduce that setting 
$Z_{i+4k}=-Y_{i+4k}$ for $i=1,2$ 
and $Z_{j+4k}=1+Y_{j+4k}$ for $j=3,4$ (and $k=0,\ldots,4$), then
the previous functional identity can be written 
\begin{equation}
\label{Eq:oko?}
\sum_{k=0}^{4}  \Big(\, 2\, {\bf L}{\rm i}_2\big( Z_{1+4k}\big)  +2\, {\bf L}{\rm i}_2 \big( -Y_{2+4k}\big) 
-2\, {\bf L}{\rm i}_2\big( Z_{3+4k}\big)  - {\bf L}{\rm i}_2\big( Z_{4+4k}\big) 
\Big)= -43\,\frac{{}^{}\,\pi^2}{6}+i\pi \cdot \boldsymbol{\bf L}_{B_4}
\end{equation}
with $\boldsymbol{\bf L}_{B_4}=
\sum_{k=0}^{4}  \big(\, 
2\, {\rm Log}\big( Z_{3+4k}\big)  + {\rm Log}( Z_{4+4k}\big) 
\big)$. Noticing that all the arguments of the logarithms in the latter sum are real and strictly greater than 1, one deduces that $\boldsymbol{\bf L}_{B_4}\in \mathbf R_{>0}$ hence that the real part of the LHS of equality \eqref{Eq:oko?} is identically equal to $-43\pi^2/6$ on $(\mathbf R_{>0})^n$.  It would be interesting to know whether  it is possible to modify the argument of the bilogarithm  in \eqref{Eq:oko?} in order to make the non-constant imaginary part disappear and thus obtain an holomorphic identity  analogous to Newman's identity $(\mathcal N_6)$ only involving  the classical bilogarithm ${\bf L}{\rm i}_2$. 

 \paragraph{About the identity $\boldsymbol{({\sf A}_{B_n})}$.}
 We have verified that identity $({\sf A}_{B_n})$ is satisfied for $n$ ranging from 2 to 12. We are quite confident that it holds true for any $n$. We make two remarks about this identity below which could be useful to establish the latter fact. 
The first  is a nice Lie-theoretic way to write down this identity. After that, we explain that $({\sf A}_{B_n})$ is 
 equivalent to an explicit algebraic   
identity which can be expressed in terms of certain $F$-polynomials. \mk 

To make the notations simpler below, we use in some places   $\Delta$ to denote the root system $B_n$.
\begin{center}
$\star$
\end{center}

To begin with, we write down in explicit form  the identity 
 under scrutiny 
 for $n=2,3$ and $4$: 
\begin{align*}
\big({\sf A}_{B_2}\big)\,& \hspace{1cm} \frac{2\pi}{2}= {\sf A} \left( {\frac { \left( u_{{1}}+1 \right) ^{2}}{u_{{2}}}} \right) +{\sf A} \left( {\frac { \left( u_{{1}}+1+u_{{2}} \right) ^{2}}{{u_{{1}}}^{2}u_{{2}}}  } \right) +{\sf A} \left( u_{{2}} \right)\, ;
\mk \\
\big({\sf A}_{B_3}\big)\,& \hspace{1cm} \frac{3\pi}{2}={\sf A} \left( {\frac { \left( 1+u_{{2}} \right) ^{2}}{u_{{3}}}} \right) +
{\sf A} \left( {\frac { \left( u_{{1}}u_{{3}}+u_{{1}}+u_{{2}}+u_{{3}}+1 \right) ^{2}}{{u_{{2}}}^{2}u_{{3}}}} \right)
 \\ 
& \hspace{4cm}+ {\sf A} \left( {\frac { \left( u_{{1}}u_{{2}}+u_{{1}}u_{{3}}+{u_{{2}}}^{2}+u_{{1}}+2\,u_{{2}}+u_{{3}}+1 \right) ^{2}}{{u_{{1}}}^{2}{u_{{2}}}^{2}u_{{3}}}} \right) +{\sf A}\left( u_{{3}} \right)\, ;\sk \\
\big({\sf A}_{B_4}\big)\,& \hspace{1cm} \frac{4\pi}{2}=
{\sf A} \left( {\frac { \left( u_{{3}}+1 \right) ^{2}}{u_{{4}}}} \right) +
{\sf A} \left( {\frac { \left( u_{{2}}u_{{4}}+u_{{2}}+u_{{3}}+u_{{4}}+1 \right) ^{2}}{{u_{{3}}}^{2}u_{{4}}}} \right)
 \\ 
& \hspace{4cm}
+ {\sf A} \left( {\frac { \left( {u_{{3}}}^{2}u_{{1}}+2\,u_{{1}}u_{{3}}+\cdots+u_{{4}}+1 \right) ^{2}}{{u_{{2}}}^{2}{u_{{3}}}^{2}u_{{4}}}} \right) \\ 
& \hspace{4cm}
+ {\sf A} \left( {\frac { \left( u_{{2}}u_{{1}}u_{{3}}+u_{{1}}u_{{2}}u_{{4}}+\cdots+u_{{4}}+1 \right) ^{2}}{{u_{{1}}}^{2}{u_{{2}}}^{2}{u_{{3}}}^{2}u_{{4}}}} \right)
 +{\sf A} \left( u_{{4}} \right)\, . \hspace{1cm} {}^{}
\end{align*}

Considering the non trivial $g$-vectors of the arguments of ${\sf A}$ in each case, we obtain the following list: $(0,1)$, $(2,1)$ for $n=2$; $(0,0,1)$, $(0,2,1)$ and 
$(2,2,1)$ when $n=3$; and $(0,0,0,1)$, $(0,0,2,1)$ and 
$(0,2,2,1)$ and $(2,2,2,1)$ if $n=4$. We recognize, for each of these three cases,  the coordinates (relatively to the usual principal roots) of the $n$ long roots of the root system of type $B_n^\vee=C_n$.  Since the long roots of $C_n$ correspond to the short ones of $B_n$ via duality $\alpha\leftrightarrow \alpha^\vee$ and considering the labelling of the first integrals of 
$\boldsymbol{\mathcal Y\hspace{-0.1cm}\mathcal W}_{\hspace{-0.05cm} B_n}$
by the elements $\alpha$ of $\Delta_{\geq -1}$ discussed in \S\ref{SS:Cluster notation},  it comes that the $Y[\alpha]$'s appearing as arguments of the function $A$ in each of the three functional identities above are the ones associated with the short roots elements of $\Delta_{\geq -1}$.\footnote{Remark that viewed the orientation of the last arrow in the Dynkin diagram (from the $(n-1)$-th vertex to the last one), the short root among the $n$ principal ones  $\alpha_1,\ldots,\alpha_n$ is $\alpha_n$ and  therefore its opposite $-\alpha_n$ is short as well.}
\sk

From the discussion above, it comes that one has 
\begin{equation}
\label{Eq:Akn-Ashort}
{}^{}\hspace{-3.5cm}
\sum_{k=1}^{n+1} {\sf A}\big( Y_{kn}\big)= \sum_{
  \substack{\alpha \in \Delta_{\geq -1} \\  \alpha \mbox{\scriptsize{ short} }  }}
  {\sf A} \Big(     Y[\alpha]   \Big)
\end{equation}
from which it follows immediately that  
$({\sf A}_n)$ is equivalent to the following identity 
$$
{\Big({\sf A}^{short}_{B_n}\Big)}
\hspace{4.35cm}
 \sum_{\alpha \in \Delta^{short}_{\geq -1}}
  {\sf A} \Big(     Y[\alpha]   \Big)
=n\,\frac{\pi}{2}\,
\hspace{7cm}{}^{} 
$$
where $\Delta^{short}_{\geq -1}$ stands for the intersection $\Delta_{\geq -1}\cap \Delta^{short}$ (see \S\ref{SS:Cluster notation} again). \mk

Actually, all that is not at all specific to the cases $n=2,3$ and $4$ since one has the 
\begin{prop}
\label{P:tokou}
For any $n\geq 2$, the set 
$\{\, Y_{kn}\,\}_{k=1}^{n+1}$
 of cluster first integrals
of $\boldsymbol{\mathcal Y\hspace{-0.1cm}\mathcal W}_{\hspace{-0.05cm} B_n}$ appearing in $({\sf A}_n)$ 
 coincides with the one formed by the $Y[\alpha]$'s labeled by the elements of $\Delta^{short}_{\geq -1}$.\sk 
 
Consequently,  the equality \eqref{Eq:Akn-Ashort} holds true identically and 
 ${\big({\sf A}_{B_n}\big)}$  can be written as $\big({\sf A}^{short}_{B_n}\big)$.
\end{prop}
\begin{proof}
One has $Y_n=u_n$ hence it coincides with $Y[-\alpha_n]$ and $-\alpha_n$ is 
a short root ($\alpha_n$ is the unique short principal root in type $B_n$). 
Let $\mu_{\bullet,\circ}$ be the (composition of)  mutation(s) corresponding to the $n$-tuple $\nu_{\bullet,\circ}$ considered in \S\ref{SSS:ClusterWebs-type-B:notations} above. Then it follows from Theorem 1.4 of \cite{FZ} that 
for $k=1,\ldots,n+1$, one has $Y_{kn}=Y[\tau^{k-1}(-\alpha_n)]$ where $\tau$ stands for the composition $\tau_-\circ \tau_+$ where $\tau_\pm$ are the automorphisms of $\Delta_{\geq -1}$ considered in \cite{FZ}.  \sk

According to formula (2.7) therein, $\tau_+$ (resp.\,$\tau_-$) is the product of the  piecewise-linear modifications $\sigma_i$ of the simple reflexion $s_i$'s for all integers $i$ corresponding to a source (resp.\,to a skink) of the Dynkin quiver $\vec{B}_n$. Since each $\sigma_i$ transforms a short root into a root of the same type (thanks to \cite[(2.5)]{FZ}), we deduce that $\Delta_{\geq -1}^{short}$ 
is stable by $\tau$  which implies that $\tau^{k-1}(-\alpha_n)$ is short for any $k=1,\ldots,n$.  Since $\Delta_{>0}^{short}$ has cardinality $n$, this proves the proposition.
\end{proof}

For $k\in \{1,\ldots,n+1\}$, we denote by $\alpha_k $ the element of $\Delta_{\geq -1}^{short}$ such that $Y_{kn}=Y[\alpha_k]$. Note that in terms of the terminology introduced in \cite{FZ}, one has $\alpha_k=\alpha(k-1;n)$ for every $k$. 
\begin{center}
$\star$
\end{center}

The starting point of  our second remark concerning the identity ${\big({\sf A}_{B_n}\big)}$ is the following classical formula for the arctangent function (satisfied for any $x\in \mathbf R_{>0}$): 
\begin{equation}
\label{Eq:A-Arctan-sqrt}
{\rm Arctan}(x)=\frac{1}{2i}\,{\rm Log}\left(
\frac{x-i}{x+i}
\right)+\frac{\pi}{2}
 \,
\end{equation}
which shows that it can be seen as a logarithmic ({\it i.e.}\,of weight 1) iterated integral on $\mathbf P^1$ ramified at $\pm i$ and $\infty$.  Since ${\sf A}(x)={\rm Arctan}(\sqrt{x})$ for $x>0$,  it appears that the fact that this function is not an iterated integral comes from the square root.  \sk 

At this point, considering the identities $({\sf A}_{B_n})$ (for $n=2,3,4$) written above in explicit form, 
an interesting remark is that replacing $u_n$ by its square 
have the result of transforming  any of the $n+1$ arguments of $A$ into squares of positive rational functions.\sk 

 More formally, considering the elementary monomial map defined by 
\begin{equation}
\label{Eq:eta}
\eta=\eta_{B_n}\, :\, \big(\mathbf R_{>0}\big)^n\longrightarrow \big(\mathbf R_{>0}\big)^n\, , \hspace{0.15cm}  \big(u_i\big)_{i=1}^n\longmapsto \overline{u}= \big(u_1,\ldots,u_{n-1}, u_n^2\big)\, , 
\end{equation}
then we have the following result (where we use the notations of \cite{FZ}):
%
%
%
%


\begin{prop} For any  $\kappa\in  \{1,\ldots,n\}$: 
\begin{enumerate}
\item  there exists $\beta_\kappa\in \Delta_{>0}$ such that $Y_{\kappa n}={F[\beta_\kappa]^2}/{u^{\alpha_\kappa^\vee}}
$;
\item there exists $\overline{\alpha}_\kappa\in \{0,1\}^n $ such that $\eta^*\big({u^{\alpha_\kappa^\vee}}\big)=\big(u^{\overline{\alpha}_\kappa}\big)^2$ and 
 consequently, one has  
$$\eta^*\big( Y_{\kappa n}
\big)=\left( \frac{ \overline{F}[\beta_\kappa]}{ u^{\overline{\alpha}_\kappa}}\right)^2\,
$$
where $\overline{F}[\beta_\kappa]$ stands for ${F}[\beta_\kappa]$ evaluated on $\overline{u}$ ({\it i.e.}\,$\overline{F}[\beta_\kappa]=\eta^*\big({F}[\beta_\kappa]\big)${\rm )}.
\end{enumerate}
\end{prop}
\begin{proof}
From \cite[\S2]{FZ} and because $Y_{\kappa n}=Y[\alpha_\kappa]$ 
with 
$\alpha_\kappa=\alpha(\kappa-1;n)\in \Delta_{>0}^{short}$ 
(according to the proof of Proposition \ref{P:tokou}), it comes that 
$Y_{\kappa n}=N[\alpha_\kappa]/{u^{\alpha_\kappa^\vee}}$ where the numerator is a product of powers of certain $F$-polynomials. More precisely,  
one has $N[\alpha_\kappa]=\prod_{(\beta,d)\in \Psi(\alpha((\kappa-1);n)} F[\beta]^d$ for a certain subset 
$\Psi(\alpha_k)$ 
of $\Delta_{\geq -1}\times \mathbf N_{>0}$. 
Since $\alpha_\kappa=\alpha(\kappa-1;n)$,  we deduce  from the formula for 
$\Psi(\alpha(k;i))$ given p.\,911 of \cite{FZ} that $\Psi(\alpha_\kappa)$ has only one element, namely $\big( \alpha(1-\kappa;n-1) ,  2 \big)$.  Thus setting $\beta_\kappa= 
 \alpha(1-\kappa;n-1) $, one has 
 $Y_{\kappa n}={F[\beta_\kappa]^2}/{u^{\alpha_\kappa^\vee}}
$, proving the first point.\sk 

Since the duality $\alpha\leftrightarrow \alpha^\vee$ 
sends the short roots of $\Delta=B_n$ onto the long ones of 
the dual root system $\Delta^\vee=C_n$, we obtain that the $\alpha_\kappa^\vee$'s (for $\kappa=1,\ldots,n$) are exactly the $n$ positive long roots of $C_n$. Relatively to the standard simple roots of $\Delta^\vee$ (which form the dual basis of the one formed by the standard simple roots of $\Delta$), the coordinate vectors of the $n$ long roots of $\Delta^\vee$ are exactly the 
$n$-tuples  of the form $a[\ell]=(0^\ell,2^{n-1-\ell},1)$ for $\ell=0,\ldots,n-1$, where  $u^\ell$ stands here for the string formed by $\ell$ $u$'s in a row for any $u$ (with the convention that $a^0$ is the empty string).  If $\alpha_\kappa=a[\ell]$, then setting $\overline{\alpha}_\kappa=(0^\ell, 1^{n-\ell})$ one indeed has $\eta^*(u^{\alpha_\kappa})=(u^{\overline{\alpha}_\kappa})^2$  hence the proposition.
\end{proof}

We define $n+1$ elements of $\mathbf C(u)$ by setting 
\begin{equation*}
{\sf y}_{\kappa n}= { \overline{F}[\beta_\kappa]}/{ u^{\overline{\alpha}_\kappa}}\quad \mbox{for }\, \kappa=1,\ldots,n\, , 
\qquad \mbox{ and } \qquad 
{\sf y}_{(n+1)n}=u_n\, .
\end{equation*}
 Then $\eta^*(Y_{kn})=({\sf y}_{k n})^2$  and consequently, $\eta^*\big({\sf A}(Y_{kn})\big)={\rm Arctan}({\sf y}_{kn})$  for every $k=1,\ldots,n+1$. Since \eqref{Eq:eta}
 induces an isomorphism of the positive orthant, it follows that $({\sf A}_{B_n})$ is equivalent (by pull-back under $\eta$) to the identity $\sum_{k=1}^{n+1}{\rm Arctan}({\sf y}_{kn})=n\pi/2$. Then using \eqref{Eq:A-Arctan-sqrt}, one obtains 
that the following chain of equivalences holds true: 
 \begin{align*}
 \sum_{k=1}^{n+1} A\big( Y_{kn}\big)=n\frac{\pi}{2}  \qquad   & \Longleftrightarrow  \qquad 
 \sum_{k=1}^{n+1} {\rm Log}\left( 
 \frac{{\sf y}_{kn}-i}{{\sf y}_{kn}+i}
 \right)=-i\pi\\
 & \Longleftrightarrow \qquad 
 \prod_{k=1}^{n+1}\left( 
 \frac{{\sf y}_{kn}-i}{{\sf y}_{kn}+i}
 \right)=-1\, . 
 \\
 & \Longleftrightarrow  \qquad 
 \prod_{k=1}^{n+1}\left( 
{{\sf y}_{kn}-i}
 \right)
 + \prod_{k=1}^{n+1}\left( 
{{\sf y}_{kn}+i}
 \right)=0
 \, . 
 \end{align*}
 For $\nu=0,\ldots,n+1$, let $\sigma_\nu$ be the $\nu$-th elementary symmetric polynomial on 
$n+1$ indeterminates.  Since $\boldsymbol{\sf y}=({\sf y}_{n},{\sf y}_{2n}\ldots, {\sf y}_{(n+1)n})$ is a $(n+1)$-tuple of rational functions with positive integer coefficients, 
 we deduce from the previous chain of equivalences the 
 \begin{lem}
 The AFE $\big({\sf A}_{B_n}\big)$ holds true if and only if the ${\sf y}_{kn}$'s defined above satisfy the following algebraic relation:
 $$
{\Big({\sf A}^{alg}_{B_n}\Big)}
\hspace{5cm}
  \sum_{ \substack{ \nu=0,\ldots,n+1 \\ \nu\,\equiv \,n+1\,[2]}}
 (-1)^{\lfloor \nu/2\rfloor}\cdot 
\sigma_\nu\big(\boldsymbol{\sf y}\big)
=0\, .
\hspace{7cm}{}^{} 
$$
 \end{lem}

Actually, this is the algebraic identity   ${\big({\sf A}^{alg}_{B_n}\big)}$ 
  that we have verified to hold true for  many values of $n$ ($n\leq 12$) and not strictly speaking $\big({\sf A}_{B_n}\big)$ which is a bit more delicate to manipulate.\mk 
  
  To finish, let us remark that although it is not formally 'polylogarithmic', 
 it follows from \eqref{Eq:A-Arctan-sqrt}
 that $\big({\sf A}_{B_n}\big)$ is indeed of the type predicted by Conjecture
\ref{Conj:Nature-Fi}. Actually, the situation is even better since by pull-back under the algebraic (actually monomial) map $\eta$, we obtain the AR 
  $$
\eta^*\big({\sf A}_{B_n}\big)
\hspace{4cm}
  \sum_{k=1}^{n+1}
\Big({\rm Log}\big({\sf y}_{kn}-i\big)-
{\rm Log}\big({\sf y}_{kn}+i
\big)\Big)
=-i\pi \, . 
\hspace{7cm}{}^{} 
$$
which is the most logarithmic.  Since this is very similar to what holds true in the $G_2$-case (see Remark \ref{Rk:TildeXWG2}), one can wonder whether  this is a general phenomenon or not and ask the 
\begin{question}
Let $\boldsymbol{\mathcal W}$ be a $\mathcal X$-cluster web (that is, a web defined by $\mathcal X$-cluster variables as first integrals).  Does it exist a rational (or even a monomial) dominant map $\chi$ (a priori with non trivial ramifications) such that all the ARs of $\chi^*\big( \boldsymbol{\mathcal W}\big)$ be polylogarithmic?
  \end{question}

 \subsubsection{\bf About the numerology of  $\boldsymbol{{\mathcal X\hspace{-0.1cm}\mathcal W}_{\hspace{-0.05cm} B_n}}$.}
\label{SSS:ClusterWebs-type-B:X-cluster web}
We now say a few words about the $\mathcal X$-cluster webs of type $B$. 
After looking at the invariants of the webs $\boldsymbol{{\mathcal X\hspace{-0.1cm}\mathcal W}}_{\hspace{-0.05cm} B_n}$ for $n\leq 6$ (computed by direct methods), we guess what the ranks for $n$ arbitrary might be. Our approach here is by entering the series of numbers we have obtained  for $n\leq 6$ into the \href{https://oeis.org}{OEIS} in order to see what comes up. Needless to say (but we write it anyway), this subsection is mainly exploratory and our formulas below for $n\geq 2$ arbitrary are conjectural.
\mk 

In what follows, $n$ stands for an integer bigger than 1. 
According to Corollary \ref{Coro:d-X-Delta--d-Y-Delta}, 
 the $\mathcal X$-cluster web of type $B_n$ is a  $d^\mathcal X_{B_n}$-web in $n$ variables  with $d^\mathcal X_{B_n}={n(n+1)(n^2+2)}/{6}$.  Here are the ranks (virtual, polylogarithmic, etc) that we have been able to compute for $n=2,\ldots,8$: 
\begin{align*}
  \rho^\bullet\big(\boldsymbol{{\mathcal X\hspace{-0.07cm}\mathcal W}}_{\hspace{-0.05cm} B_2}\big)=& \, (4,3,2,1) && 
   {\rm polrk}^\bullet(\boldsymbol{{\mathcal X\hspace{-0.07cm}\mathcal W}}_{\hspace{-0.05cm} B_2}\big)=\big(6,2 \big) 
   && 
   {\rm rk}(\boldsymbol{{\mathcal X\hspace{-0.07cm}\mathcal W}}_{\hspace{-0.05cm} B_2}\big)=10
   \\
 \rho^\bullet\big(\boldsymbol{{\mathcal X\hspace{-0.07cm}\mathcal W}}_{\hspace{-0.05cm} B_3}\big)=&\, (19,16,12,7,1) && 
  {\rm polrk}^\bullet(\boldsymbol{{\mathcal X\hspace{-0.07cm}\mathcal W}}_{\hspace{-0.05cm} B_3}\big)=\big(32,11 \big)
  && 
  {\rm rk}(\boldsymbol{{\mathcal X\hspace{-0.07cm}\mathcal W}}_{\hspace{-0.05cm} B_3}\big)=50
  \\
  \rho^\bullet\big(\boldsymbol{{\mathcal X\hspace{-0.07cm}\mathcal W}}_{\hspace{-0.05cm} B_4}\big)=& \,(56,50,40,25,5)
  &&   {\rm polrk}^\bullet(\boldsymbol{{\mathcal X\hspace{-0.07cm}\mathcal W}}_{\hspace{-0.05cm} B_4}\big)=\big(100,37 \big) && {\rm rk}(\boldsymbol{{\mathcal X\hspace{-0.07cm}\mathcal W}}_{\hspace{-0.05cm} B_4}\big)=161
   \\ 
  \rho^\bullet\big(\boldsymbol{{\mathcal X\hspace{-0.07cm}\mathcal W}}_{\hspace{-0.05cm} B_5}\big)=& \, (130,120,100,65,15)
  &&     {\rm polrk}^\bullet(\boldsymbol{{\mathcal X\hspace{-0.07cm}\mathcal W}}_{\hspace{-0.05cm} B_5}\big)=\big(240,101 \big)
     \\ 
\rho^\bullet\big(\boldsymbol{{\mathcal X\hspace{-0.07cm}\mathcal W}}_{\hspace{-0.05cm} B_6}\big)=&\,  (260,245,210,140,35)
  &&     {\rm polrk}^\bullet(\boldsymbol{{\mathcal X\hspace{-0.07cm}\mathcal W}}_{\hspace{-0.05cm} B_6}\big)=\big(490 , 226  \big)
  \\
  \rho^\bullet\big(\boldsymbol{{\mathcal X\hspace{-0.07cm}\mathcal W}}_{\hspace{-0.05cm} B_7}\big)=&\,  (469, 448 , 392 , 266, 70)
  &&     {\rm polrk}^\bullet(\boldsymbol{{\mathcal X\hspace{-0.07cm}\mathcal W}}_{\hspace{-0.05cm} B_7}\big)=\big(896 , 442  \big)
   \\
  \rho^\bullet\big(\boldsymbol{{\mathcal X\hspace{-0.07cm}\mathcal W}}_{\hspace{-0.05cm} B_8}\big)=&\,  (784, 756, 672, 462, 126)
  &&     {\rm polrk}^\bullet(\boldsymbol{{\mathcal X\hspace{-0.07cm}\mathcal W}}_{\hspace{-0.05cm} B_8}\big)=\big(1512, 785  \big)\, .
\end{align*}


First, one observes that the following relations are satisfied for $n=2,\ldots, 8$: 
 $$ \rho^1\big(\boldsymbol{{\mathcal X\hspace{-0.07cm}\mathcal W}}_{\hspace{-0.05cm} B_n}\big)= d^{\mathcal X}_{B_n}-n
 \qquad \mbox{ and } \qquad 
   \rho^5\big(\boldsymbol{{\mathcal X\hspace{-0.07cm}\mathcal W}}_{\hspace{-0.05cm} B_n}\big)=\rho^4\big(\boldsymbol{{\mathcal X\hspace{-0.07cm}\mathcal W}}_{\hspace{-0.05cm} B_n}\big)-\rho^3\big(\boldsymbol{{\mathcal X\hspace{-0.07cm}\mathcal W}}_{\hspace{-0.05cm} B_n}\big)/2\, .$$
 Considering these relations and injecting the explicit values above into the \href{https://oeis.org}{OEIS}, one can extrapolate most of 
 the ranks of $\boldsymbol{{\mathcal X\hspace{-0.07cm}\mathcal W}}_{\hspace{-0.05cm} B_n}$ for $n$ arbitrary and state the following conjecture (where any sequence of the \href{https://oeis.org/}{OEIS} given on the right corresponds (possibly up to a shift) to the formula on the left of the same line): 
 \mk
 
{\bf Conjecture about the ranks of  $\boldsymbol{\boldsymbol{{\mathcal X\hspace{-0.07cm}\mathcal W}}_{\hspace{-0.05cm} B_n}}$.} 
{\it Assume that  $n$ is bigger than or equal to $2$.}
  \begin{itemize}
 \item {\bf [Virtual ranks].} {\it One has}
 \begin{align*}
 \rho^1\big(\boldsymbol{{\mathcal X\hspace{-0.07cm}\mathcal W}}_{\hspace{-0.05cm} B_n}\big)=& \, d^{\mathcal X}_{B_n}-n
 = n(n-1)(n^2+2n+4)/6 &&  
  \mbox{\Big({\it cf.}\,\href{https://oeis.org/A332697}{A332697} \Big)}
  \\
 \rho^2\big(\boldsymbol{{\mathcal X\hspace{-0.07cm}\mathcal W}}_{\hspace{-0.05cm} B_n}\big)= & \, 
 (n+1)\, { n+1 \choose 3}
 && 
  \mbox{\Big({\it cf.}\,\href{https://oeis.org/A004320}{A004320} \Big)}
\\
\rho^3\big(\boldsymbol{{\mathcal X\hspace{-0.07cm}\mathcal W}}_{\hspace{-0.05cm} B_n}\big)=&\,  
n\, { n+1 \choose 3}
&&
  \mbox{\Big({\it cf.}\,\href{https://oeis.org/A00891}{A00891} \Big)}
\\
\rho^4\big(\boldsymbol{{\mathcal X\hspace{-0.07cm}\mathcal W}}_{\hspace{-0.05cm} B_n}\big)=&\,  {(3n-2)}\, { n+1 \choose  3}/{4}
&&
  \mbox{\Big({\it cf.}\,\href{https://oeis.org/A001296}{A001296} \Big)}
  \\
\rho^5\big(\boldsymbol{{\mathcal X\hspace{-0.07cm}\mathcal W}}_{\hspace{-0.05cm} B_n}\big)=&  \, {n+1 \choose 4} &&
 \mbox{\Big({\it cf.}\,\href{https://oeis.org/A000332}{A000332} \Big)}
\vspace{0.5cm} \\
 \mbox{ and }\quad 
\rho^\sigma\big(\boldsymbol{{\mathcal X\hspace{-0.07cm}\mathcal W}}_{\hspace{-0.05cm} B_n}\big)=&  \, 0 \quad \mbox{ for any }\, \sigma\geq 6\,; &&
 \end{align*}
 \item   
 {\bf [Polylogarithmic ranks].} {\it One has}
  \begin{align*}
{\rm polrk}^1\big(\boldsymbol{{\mathcal X\hspace{-0.07cm}\mathcal W}}_{\hspace{-0.05cm} B_n}\big)=& \, 2(n+1) {n+1 \choose 3 }  \hspace{3.1cm} \quad \Big(\mbox{\it cf.}\, 
 \href{https://oeis.org/A288961}{A288961}\,\Big)\\ 
\mbox{ \it and } \quad {\rm polrk}^w\big(\boldsymbol{{\mathcal X\hspace{-0.07cm}\mathcal W}}_{\hspace{-0.05cm} B_n}\big)= & \, 0  \quad \mbox{for any } \, w\geq 3\, . 
 \end{align*}
 \end{itemize}
  
It would be quite interesting to know whether  the above formulas indeed hold true  but for the moment, we are not aware of any approach allowing to do so.  Here are a few comments on the formulas above and on what should come just next: 

\begin{itemize}
\item The sequence of weight 2 polylogarithmic ranks is 
\begin{equation}
\label{Eq:polrk2XWBn}
{\rm polrk}^2\Big(\boldsymbol{{\mathcal X\hspace{-0.07cm}\mathcal W}}_{\hspace{-0.05cm} B_\bullet}\Big)=\big(2,11,37,101,226,442, 785 , \ldots\big)\,.
\end{equation}
We have no idea of how  this sequence continues and the OEIS is not helping in this regard.  Having a (even conjectural) general  formula for this quantity would be interesting.  The same remark applies to the sequence of total ranks ${\rm rk}^2\big(\boldsymbol{{\mathcal X\hspace{-0.07cm}\mathcal W}}_{\hspace{-0.05cm} B_\bullet}\big)=(10,50, 161,\ldots)$ as well.
\item However, if one disregards the case $n=3$, then there is a (unique) entry in the OEIS corresponding to \eqref{Eq:polrk2XWBn}, namely the sequence \href{https://oeis.org/A154323}{A154323},  the first height values of which are
	1, 2, 10, 37, 101, 226, 442 and  785.  Actually, the sequence A000537 whose general term is 
	$A000537(n)=A154323(n)-1$ 
	 seems even more interesting: it appears in various contexts, and in particular each term can be combinatorially characterized as the  number of hexagons with vertices in an hexagonal grid. 
	 There is certainly more to find out about the ranks of the $\mathcal X$-cluster webs 
	from a combinatorial perspective.
\item It is known that numerology of cluster algebras is related, at least in many ways, to combinatoric ({\it e.g.}\,the enumeration of $\mathcal X$-clusters and $\mathcal X$-cluster variables in finite type in terms of certain polygons, see \cite{Sherman-Bennett1}).  This might again be the case for the diverse ranks of 
$\boldsymbol{{\mathcal X\hspace{-0.07cm}\mathcal W}}_{\hspace{-0.05cm} B_n}$.

 For instance, the OEIS provides the following combinatorial interpretation of the $k$-th term ${k \choose 4}$ of the sequence A000332  which might coincide with 
 $\rho^5\big(\boldsymbol{{\mathcal X\hspace{-0.07cm}\mathcal W}}_{\hspace{-0.05cm} B_\bullet}\big)=(0,1,5,15,35,\ldots)$ (up to a shift by 1) :  the binomial coefficient ${n+1 \choose 4}$ which is conjectured to coincide with $\rho^5\big(\boldsymbol{{\mathcal X\hspace{-0.07cm}\mathcal W}}_{\hspace{-0.05cm} B_n}\big)$ is equal to the 
number of intersection points of diagonals of a $(n+1)$-gon where no more than two diagonals intersect at any point inside. We believe that the fact that combinatorial objects of essentially the same type serve also to enumerate the cluster and the cluster variables of the cluster 
algebra of type $B_n$ ({\it i.e.}\,see \ref[\S12.3]{FZII} or \cite[\S5.5]{FWZ}) is not a coincidence.  
\sk 
\item  Our true interest actually concerns abelian relations and we do not have any guess about a possible description of (a basis of)  the space $\boldsymbol{\mathcal A}
\big(\boldsymbol{{\mathcal X\hspace{-0.07cm}\mathcal W}}_{\hspace{-0.05cm} B_n}
\big)$ for $n$ arbitrary.\sk 

Since the $\mathcal X-$ and $\mathcal Y$-cluster webs coincide in type $B_2$,  the space of ARs in this case has been described above in \ref{SS:cluster-webs:B2}: it admits a basis formed by logarithmic ARs, two dilogarithmic ARs (${ ({\sf R}_{B_2})}$ and 
${ ({\sf S}_{B_2})}$), the AR ${ ({\sf A}_{B_2})}$ plus an extra rational AR.
\sk

We have verified that 
there is a basis of 
$\boldsymbol{\mathcal A}\big(\boldsymbol{\mathcal X\hspace{-0.1cm}\mathcal W}_{\hspace{-0.05cm}B_3}\big)$ formed by logarithmic ARS,   eleven dilogarithmic ARs, four rational ARs (all the components of which are multiples of ${\sf J}:x\mapsto 1/(1+x)$) and three others whose components are multiples of  the function ${\sf A}$.

Extrapolating from the two preceding cases, one may 
wonder whether   a decomposition in direct sum of the following from holds true 
for any $n\geq 2$:
\begin{equation*}
\label{Eq:AR-XWBn}
\boldsymbol{\mathcal A}\big(\boldsymbol{\mathcal X\hspace{-0.1cm}\mathcal W}_{\hspace{-0.05cm}B_n}\big) = 
\boldsymbol{Lo\hspace{-0.03cm}g\hspace{-0.03cm}A\hspace{-0.03cm}R}\big(\boldsymbol{{\mathcal X\hspace{-0.1cm}}\mathcal W}_{\hspace{-0.05cm} B_n}\big) 
 \, \oplus \, 
\boldsymbol{Dilo\hspace{-0.03cm}g\hspace{-0.03cm}A\hspace{-0.03cm}R}\big(\boldsymbol{{\mathcal X\hspace{-0.1cm}}\mathcal W}_{\hspace{-0.05cm} B_n}\big) 
\,  \oplus  \,
\boldsymbol{ {{\sf J}}\hspace{-0.01cm}A\hspace{-0.01cm}R}\big(\boldsymbol{{\mathcal X\hspace{-0.1cm}}\mathcal W}_{\hspace{-0.05cm} B_n}\big)  \, \oplus \, 
\boldsymbol{ {\sf A}\hspace{-0.00cm}A\hspace{-0.03cm}R}\big(\boldsymbol{{\mathcal X\hspace{-0.1cm}}\mathcal W}_{\hspace{-0.05cm} B_n}\big) \, 
\end{equation*}
where $ \boldsymbol{{\sf J}\hspace{-0.01cm}A\hspace{-0.03cm}R}\big(\boldsymbol{{\mathcal X\hspace{-0.1cm}}\mathcal W}_{\hspace{-0.05cm} B_n}\big) $
 (resp.\,$\boldsymbol{ {\sf A}\hspace{-0.00cm}A\hspace{-0.03cm}R}\big(\boldsymbol{{\mathcal X\hspace{-0.1cm}}\mathcal W}_{\hspace{-0.05cm} B_n}\big)$) 
stands for the subspace of  $\boldsymbol{\mathcal A}\big(\boldsymbol{\mathcal X\hspace{-0.1cm}\mathcal W}_{\hspace{-0.05cm}B_3}\big)$ 
formed by ARs whose components are multiples of  ${\sf J}$ (resp.\,of ${\sf A}$).

\end{itemize}


 \subsection{\bf Cluster webs of type $\boldsymbol{C}$}
\label{SS:ClusterWebs-type-C}

The $\mathcal Y$-cluster webs of type $C$ share many properties with the same webs in type $B$. In particular, for any $n\geq 2$, 
$\boldsymbol{\mathcal Y\hspace{-0.1cm}\mathcal W}_{\hspace{-0.05cm}B_n}$ and 
$\boldsymbol{\mathcal Y\hspace{-0.1cm}\mathcal W}_{\hspace{-0.05cm}C_n}$
have the same (virtual, polylogarithmic, total) rank(s).  However, these webs probably are not equivalent for $n\geq 3$. \sk 

For this reason, on the one hand the structure of the present subsection is similar to that of  the previous one, and on the second hand, we give less details and 
go straight to the point.\mk 

As before, below $n$ stands for an integer bigger than 1.

 \subsubsection{\bf Some notations in type $\boldsymbol{C}.$}
\label{SSS:ClusterWebs-type-C:notations}
The bipartite initial exchange matrix associated to the Dynkin diagram 
$
\dynkin C{}
$ of type $C_n$ we will work with is
$$B_{\vec{C}_n}=
\begin{bmatrix}
0 & (-1)^n & 0   &\cdots &  \cdots& 0\\
(-1)^{n+1} & 0 & \ddots   & &  & \vdots\\
0&   (-1)^n&   \ddots &  \ddots & &  \vdots\\
\vdots & \ddots &  \ddots  &\ddots &-1 & 0 \\
\vdots&  &   \ddots &1 & {}^{}\hspace{0.2cm}0 & 1\\
0 & \cdots &  0  &0 & -2&0
\end{bmatrix} \, .
$$

The other notations in type $C_n$ coincide with those of $B_n$: 
the corresponding sets of sources  and of sinks are the same as in the case of $B_n$;  
as initial $\mathcal X$-seed, we take $\boldsymbol{\mathcal S}_0=\big(\boldsymbol{v}_0 , B_{\vec{C}_n}\big)$ where $\boldsymbol{v}_0$ denotes 
a $n$-tuple of variables 
\eqref{Eq:nu-0};  $\mu_\circ$ (resp.\,$\mu_\bullet$) stands for the composition of mutations on the set of sources (resp.\,of sinks), etc.  Then the $n$-tuple $\nu_{\bullet,\circ}$ (defined as in the $B_n$ case) is a period  for $\boldsymbol{\mathcal S}_0$ and denoting here by 
$Y_s $ the $i_s$-th cluster variable of the $s$-th seed, we get 
\begin{equation}
\label{Eq:Web-YWCn}
\boldsymbol{\mathcal Y\hspace{-0.1cm}\mathcal W}_{\hspace{-0.05cm} C_n}
=\boldsymbol{\mathcal Y\hspace{-0.1cm}\mathcal W}\Big( \, 
Y_s\, \big\lvert \, s=1,\ldots,n(n+1)\, \Big)\, .
\end{equation}

Finally, one verifies that a right skew-symmetrizer for $B_{\vec{C}_n}$ is given by the diagonal matrix whose diagonal coefficients form the $n$-tuple 
$d_{C_n}=(1,\ldots,1,2)$ whose coefficients are denoted simply by  $d_1,\ldots,d_n$  (with no reference to $C_n$) in the lines below ($d_i=1$ for $i=1,\ldots,n-1$ and $d_n=2$).
\begin{center}
$\star$
\end{center}

We now discuss the non-logarithmic ARs of $\boldsymbol{\mathcal Y\hspace{-0.1cm}\mathcal W}_{\hspace{-0.05cm} C_n}$: conjecturally, such ARs are either dilogarithmic, or rational. We describe the ARs of these two types below before stating a general conjecture about the ranks and the ARs of the $\mathcal Y$-cluster web of type $C_n$ for $n\geq 2$ arbitrary. \mk 

Using the labelling of the $Y_s$'s given by  $Y_i^{[\kappa]}=Y_{i+\kappa n} $ 
for $i=1,\ldots,n$ and $\kappa=0,\ldots,n$, it follows from Nakanishi's theorem 
 that the following identity holds true: 
$$
{\Big({\sf R}_{C_n}\Big)}
\hspace{3.5cm}
\sum_{\kappa=0}^{n} \left( \Bigg( \sum_{i=1}^{n-1}\,
{\sf R}\Big( Y_i^{[\kappa]}\Big) \Bigg)
+ 2\,  {\sf R}\Big( Y_n^{[\kappa]}\Big)\,
\right)= n(n+1)
\, \frac{{}^{}\hspace{0.06cm}\pi^2}{6}\, . 
\hspace{7cm}{}^{} 
$$

Similarly to what does occur in type $B$, when $n=2m$ is even, it is expected that the web
$\boldsymbol{\mathcal Y\hspace{-0.1cm}\mathcal W}_{\hspace{-0.05cm} C_n}$ 
carries a supplementary symmetric dilogarithmic AR associated to the following identity  
$$
{\Big({\sf S}_{C_n}\Big)}
\hspace{2.5cm}
\sum_{\kappa=0}^{n} \left[ \Bigg( \sum_{i=1}^{m}{\sf S}\big( Y_i^{[\kappa]}\Big)
\Bigg) 
-\Bigg(
\sum_{j=1}^{m-1}{\sf S}\Big( Y_{m+j}^{[\kappa]}\Big)\Bigg)
- 2\,{\sf S}\Big( Y_n^{[\kappa]}\Big)\,
\right]= 0\, 
\hspace{7cm}{}^{} 
$$
which is conjectured to hold true identically on $(\mathbf R_{>0})^n$. \mk

Finally, one conjectures that the following rational identity holds true: 
$$
{\Big({\sf J}_{C_n}\Big)}
\hspace{2.5cm}
{\sf J}\Big( Y_{n}\Big)\,+\,\cdots\,+\,{\sf J}\Big( Y_{(n+1)n}\Big)=
\sum_{\kappa=1}^{n+1} {\sf J}\Big( Y_{\kappa n}\Big) = 1\, 
\hspace{9cm}{}^{} 
$$
 where ${\sf J}$ stands for the rational function ${\sf J}: x\mapsto 1/(1+x)$.  As it can be verified that the cluster variables $Y_{\kappa n}$'s for $\kappa=1,\ldots,n+1$ coincide with the $Y[\alpha]$'s with $\alpha\in \Delta_{\geq -1}^{long}$ (where $\Delta$ now stands for the  root system of type $C_n$), the previous algebraic identity can also be written as 
 $$\sum_{\alpha \in 
 \Delta_{\geq -1}^{long}}  {\sf J}\big( Y[\alpha]\big) = 1\, .$$
 
 Then one states the following version in type $C$ of the 
 `Conjecture $\boldsymbol{{\mathcal Y\hspace{-0.1cm} \mathcal W}_{B_n}}$' stated above: 
 \mk

\label{Conj:YWCn-page}
{\bf Conjecture $\boldsymbol{{\mathcal Y\hspace{-0.1cm} \mathcal W}_{C_n}}$.} 
{\it 
Let $n$ be an integer bigger than or equal to 3. 

\begin{enumerate}
\item The first  five points of {\bf Conjecture $\boldsymbol{{\mathcal Y\hspace{-0.1cm} \mathcal W}_{B_n}}$} hold true if one replaces  $B_n$ by $C_n$ everywhere. 

\item  {\bf ${\big[}$The rational AR\big].}
 The following rational identity 
\begin{equation*}
\label{Eq:J(YWCn)}
\sum_{k=1}^{n+1} {\sf J}\Big(Y_{kn}\Big)=
\sum_{\alpha \in 
 \Delta_{\geq -1}^{long}}  {\sf J}\big( Y[\alpha]\big) = 1
\end{equation*}
is identically satisfied. The associated AR of $\boldsymbol{\mathcal Y\hspace{-0.1cm}\mathcal W}_{\hspace{-0.05cm}C_n}$ is not of polylogarithmic type.
\mk
\item {\bf ${\big[}$Rank and basis of ARs{\big]}.} 
  One has 
    $$
  {\rm rk}\Big(\boldsymbol{\mathcal Y\hspace{-0.1cm} \mathcal W}_{\hspace{-0.05cm}C_n}\Big)={\rm polrk}\Big(\boldsymbol{\mathcal Y\hspace{-0.1cm} \mathcal W}_{\hspace{-0.05cm}C_n}\Big)+1=n(n+1)+\beta_n+1
  $$
 hence the following decomposition in direct sum holds true:
\begin{equation}
\label{Eq:AR-YWCn}
\boldsymbol{\mathcal A}\Big(\boldsymbol{\mathcal Y\hspace{-0.1cm}\mathcal W}_{\hspace{-0.05cm}C_n}\Big) = 
\boldsymbol{Lo\hspace{-0.03cm}g\hspace{-0.03cm}A\hspace{-0.03cm}R}\Big(\boldsymbol{{\mathcal Y\hspace{-0.1cm}}\mathcal W}_{\hspace{-0.05cm} C_n}\Big) 
 \, \oplus \, 
\boldsymbol{Dilo\hspace{-0.03cm}g\hspace{-0.03cm}A\hspace{-0.03cm}R}\Big(\boldsymbol{{\mathcal Y\hspace{-0.1cm}}\mathcal W}_{\hspace{-0.05cm} C_n}\Big) 
\,  \oplus  \,
 \big\langle
\, {\sf J}_{C_n}\, 
 \big\rangle\, . 
\end{equation}
\end{enumerate}}

By direct computations, we have verified the above conjecture for $n\leq 10$.

 \subsubsection{\bf About the numerology of  $\boldsymbol{{\mathcal X\hspace{-0.1cm}\mathcal W}_{\hspace{-0.05cm} C_n}}$.}
\label{SSS:ClusterWebs-type-C:X-cluster web}
We won't say much about the ranks of the cluster webs 
$\boldsymbol{\mathcal X\hspace{-0.1cm}\mathcal W}_{\hspace{-0.05cm} C_n}$
 because all these seem to precisely coincide with the corresponding ones of 
 $\boldsymbol{\mathcal X\hspace{-0.1cm}\mathcal W}_{\hspace{-0.05cm} B_n}$. \sk 
 
 Replacing $B_n$ by $C_n$ everywhere in \S\ref{SSS:ClusterWebs-type-B:X-cluster web} gives us a series of statements which all make sense and which all have been checked 
to be valid  for  every $n\leq 10$. 
\sk 
 
The $\boldsymbol{\mathcal Y}$-cluster webs of type $C$ share many properties with the same webs in type $B$. In particular, for any $n\geq 2$, 
$\boldsymbol{\mathcal Y\hspace{-0.1cm}\mathcal W}_{\hspace{-0.05cm}B_n}$ and 
$\boldsymbol{\mathcal Y\hspace{-0.1cm}\mathcal W}_{\hspace{-0.05cm}C_n}$
have the same (virtual, polylogarithmic, total) rank(s).  However, these webs probably are not equivalent for $n\geq 3$.

\subsection{\bf Cluster webs of type $\boldsymbol{D}$.}
\label{SS:ClusterWebs-type-D}
The structure of the present subsection is similar to that of  the previous one, therefore we give less details and go straight to our points. 
 Since the Dynkin diagrams $D_n$ are simply laced, one might expect that the 
 $\boldsymbol{\mathcal Y}$-cluster webs of type $D$ satisfy the same nice properties as those of the corresponding webs in type $A$. Unfortunately and somewhat surprisingly, this is not the case.
\mk 

In the sequel $n$ stands for an integer bigger than 3.
 We set $m=\lfloor n/2\rfloor$ and 
$ \widetilde m=\lfloor (n-1)/2\rfloor$.

 \subsubsection{\bf Numerology of cluster webs of type $\boldsymbol{D}$.}
\label{SSS:NumerologyClusterWebs-type-D}
The bipartite initial exchange matrix associated to the Dynkin diagram 
$
\dynkin[backwards] D{}
$ of type $D_n$ we will work with is ({\it cf.}\,Figure \ref{Fig:DynkinQuiverDiagram}): 
$$B_{{D}_n}=
\begin{bmatrix}
0 &1 & 0  &0 &0 &     \cdots & 0\\
-1 & 0 & -1   & -1&   0& \cdots  & \vdots \\
0&   1&  0 & 0&  0 &    \\
0 & 1 & 0  &0 &1 & 0 & \vdots \\
0& 0 & 0   &-1 &    0 & \ddots  & 0\\
\vdots &  &    & &\ddots & \ddots & (-1)^{n-1} \\
0 & \cdots &    &  \cdots & 0&(-1)^n &  0 \\
\end{bmatrix} \, .
$$

The  set of sources  (resp.\,of sinks) is $\{1,3\} \cup \{ 2s\, \lvert \, s=2,\ldots , m\, \} $ (resp.\,$\{2\}\cup \{ 1 + 2s \, \lvert \,s=2,\ldots , \widetilde m \, \} $).  The $n^2$-tuple $\mu_{D_n}$ obtained by concatenating $n$ copies of the $n$-tuple 
$(1,3,4,\ldots, 2m, 2,5,7,\ldots, 1+2\widetilde m)$ is a period for the initial seed $\big(\boldsymbol{v},  B_{\vec{D}_n}\big)$ where $\boldsymbol{v}$ stands for the $n$-tuple 
of variables 
$$\Big( (u_1)^{-1}\, ,\, u_2\, , \, (u_3)^{-1}\, ,\, (u_4)^{-1}\, , \, u_5\, , \, \ldots, (u_{n-1})^{(-1)^{n-2}}\, ,\, 
(u_{n})^{(-1)^{n-1}} \Big)\, .$$

Denoting by 
$Y_s $ the $i_s$-th cluster variable of the $s$-th seed (which is the one obtained from the initial one by applying the mutations corresponding to the $s$ first elements of $\mu_{D_n}$), we have
\begin{equation}
\label{Eq:Web-YWDn}
\boldsymbol{\mathcal Y\hspace{-0.1cm}\mathcal W}_{\hspace{-0.05cm} D_n}
=\boldsymbol{\mathcal Y\hspace{-0.1cm}\mathcal W}\Big( \, 
Y_s\, \big\lvert \, s=1,\ldots,n^2\, \Big)\, .
\end{equation}

It is thus a $n^2$-web in $n$ variables which carries the abelian relation associated to the 
 following dilogarithmic identity, which holds true identically 
 according to the main result of \cite{Chapoton}: 
 $$
{\Big({\sf R}_{D_n}\Big)}
\hspace{3.5cm}
\sum_{i=1}^{n^2} {\sf R}(Y_i)= n(n-1)
\, \frac{{}^{}\hspace{0.06cm}\pi^2}{6}\, . 
\hspace{7cm}{}^{} 
$$

\paragraph{\bf Numerology of the $\boldsymbol{{\mathcal Y\hspace{-0.1cm}\mathcal W}_{\hspace{-0.05cm} D_n}}$'s.}
By direct computations, one obtains that the ranks of $\boldsymbol{\mathcal Y\hspace{-0.1cm}\mathcal W}_{\hspace{-0.05cm} D_n}$ for $n=4,\ldots,8$ are as follows: 
 \begin{align*}
    \rho^\bullet(\boldsymbol{{\mathcal Y\hspace{-0.1cm}\mathcal W}}_{\hspace{-0.05cm} D_4})=& \big(\, 12\, , \, 6\, , \, 1 \, \big)  && \rho(\boldsymbol{{\mathcal Y\hspace{-0.1cm}\mathcal W}}_{\hspace{-0.05cm} D_4})= 19&& 
   {\rm polrk}^\bullet\big(\boldsymbol{{\mathcal Y\hspace{-0.1cm}\mathcal W}}_{\hspace{-0.05cm} D_4}\big)=\big(\,   16\, , \, 1      \,  \big)
   &&  {\rm rk}\big(\boldsymbol{{\mathcal Y\hspace{-0.1cm}\mathcal W}}_{\hspace{-0.05cm} D_4}\big)= 17
  \\
      \rho^\bullet(\boldsymbol{{\mathcal Y\hspace{-0.1cm}\mathcal W}}_{\hspace{-0.05cm} D_5})=& \big(\, 20\, , \, 10\,,\,  1\, \big)  && \rho(\boldsymbol{{\mathcal Y\hspace{-0.1cm}\mathcal W}}_{\hspace{-0.05cm} D_5})= 31 && 
   {\rm polrk}^\bullet\big(\boldsymbol{{\mathcal Y\hspace{-0.1cm}\mathcal W}}_{\hspace{-0.05cm} D_5}\big)=\big(\,  25\, , \, 1       \,  \big)
   &&  {\rm rk}\big(\boldsymbol{{\mathcal Y\hspace{-0.1cm}\mathcal W}}_{\hspace{-0.05cm} D_5}\big)=   26
  \\
      \rho^\bullet(\boldsymbol{{\mathcal Y\hspace{-0.1cm}\mathcal W}}_{\hspace{-0.05cm} D_6})=& \big(\, 30\, , \, 15\, , \, 1 \, \big)  && \rho(\boldsymbol{{\mathcal Y\hspace{-0.1cm}\mathcal W}}_{\hspace{-0.05cm} D_6})= 46&& 
   {\rm polrk}^\bullet\big(\boldsymbol{{\mathcal Y\hspace{-0.1cm}\mathcal W}}_{\hspace{-0.05cm} D_6}\big)=\big(\,     36\, , \, 1     \,  \big)
   &&  {\rm rk}\big(\boldsymbol{{\mathcal Y\hspace{-0.1cm}\mathcal W}}_{\hspace{-0.05cm} D_6}\big)= 37
  \\
      \rho^\bullet(\boldsymbol{{\mathcal Y\hspace{-0.1cm}\mathcal W}}_{\hspace{-0.05cm} D_7})=& \big(\, 42\, ,\, 21 \, ,  \, 1\,  \big)  && \rho(\boldsymbol{{\mathcal Y\hspace{-0.1cm}\mathcal W}}_{\hspace{-0.05cm} D_7})= 64&& 
   {\rm polrk}^\bullet\big(\boldsymbol{{\mathcal Y\hspace{-0.1cm}\mathcal W}}_{\hspace{-0.05cm} D_7}\big)=\big(\,  49\, , \, 1        \,  \big)
   &&  {\rm rk}\big(\boldsymbol{{\mathcal Y\hspace{-0.1cm}\mathcal W}}_{\hspace{-0.05cm} B_7}\big)= 50\\
      \rho^\bullet(\boldsymbol{{\mathcal Y\hspace{-0.1cm}\mathcal W}}_{\hspace{-0.05cm} D_8})=& \big(\,56\, ,\, 28\, , \, 1 \, \big)  && \rho(\boldsymbol{{\mathcal Y\hspace{-0.1cm}\mathcal W}}_{\hspace{-0.05cm} D_8})= 102&& 
   {\rm polrk}^\bullet\big(\boldsymbol{{\mathcal Y\hspace{-0.1cm}\mathcal W}}_{\hspace{-0.05cm} D_8}\big)=\big(\,  64\, , \, 1       \,  \big)
   &&  {\rm rk}\big(\boldsymbol{{\mathcal Y\hspace{-0.1cm}\mathcal W}}_{\hspace{-0.05cm} D_8}\big)= 65 \, . 
\end{align*}

From these values, one can extrapolate the values of the ranks for $n$ arbitrary and state the \mk 

\label{Conj:YWDn-page}
{\bf Conjecture $\boldsymbol{{\mathcal Y\hspace{-0.1cm} \mathcal W}_{D_n}}$.}
{\it 
Let $n$ be an integer bigger than or equal to 4. 
\begin{enumerate}
\item 
 The ramification of $\boldsymbol{{\mathcal Y\hspace{-0.1cm}\mathcal W}}_{\hspace{-0.05cm} D_n}$ is polylogarithmic and 
the  variety of common leaves $\Sigma^c\big(\boldsymbol{{\mathcal Y\hspace{-0.1cm}\mathcal W}}_{\hspace{-0.05cm} D_n}\big)$
 coincides with the corresponding 
 cluster arrangement:  one has 
 $
 \Sigma^c\Big(\boldsymbol{{\mathcal Y\hspace{-0.1cm}\mathcal W}}_{\hspace{-0.05cm} D_n}\Big)=\boldsymbol{Arr}_{D_n}\, .
 $
\item The virtual and polylogarithmic ranks  
of $\boldsymbol{{\mathcal Y\hspace{-0.1cm}\mathcal W}}_{\hspace{-0.05cm} D_n}$ are  given by  
\begin{align*}
\rho^\bullet\Big(\boldsymbol{\mathcal Y\hspace{-0.1cm}\mathcal W}_{\hspace{-0.05cm} D_n}\Big)= &\,  \Big(\,  n(n-1)\, , \, {n(n-1)}/{2}\, , \, 1\,  \Big) 
\qquad 
\mbox{and } \qquad {\rm polrk}^\bullet\Big(
\boldsymbol{{\mathcal Y\hspace{-0.1cm}\mathcal W}}_{\hspace{-0.05cm} D_n}
\Big)=  \Big( \, n^2\, , \,1 \, \Big) 
\, . 
\end{align*}
Thus 
$\rho\big(
\boldsymbol{{\mathcal Y\hspace{-0.1cm}\mathcal W}}_{\hspace{-0.05cm} D_n}
\big)=  3
 n (n-1)/{2 }+1 $
and 
${\rm polrk}\big(
\boldsymbol{{\mathcal Y\hspace{-0.1cm}\mathcal W}}_{\hspace{-0.05cm} D_n}
\big) =   n^2+1$.
\item 
  One has 
    $
  {\rm rk}\big(\boldsymbol{\mathcal Y\hspace{-0.1cm} \mathcal W}_{\hspace{-0.05cm}D_n}\big)={\rm polrk}\big(\boldsymbol{\mathcal Y\hspace{-0.1cm} \mathcal W}_{\hspace{-0.05cm}D_n}\big)= n^2+1
  $ 
   hence the following decomposition: 
\begin{equation}
\label{Eq:AR-YWDn}
\boldsymbol{\mathcal A}\Big(\boldsymbol{\mathcal Y\hspace{-0.1cm}\mathcal W}_{\hspace{-0.05cm}D_n}\Big) = 
\boldsymbol{Lo\hspace{-0.03cm}g\hspace{-0.03cm}A\hspace{-0.03cm}R}\Big(\boldsymbol{{\mathcal Y\hspace{-0.1cm}}\mathcal W}_{\hspace{-0.05cm} D_n}\Big) 
 \, \oplus \, 
  \Big\langle \, 
 {\sf R}_{D_n}\, 
 \Big\rangle\, . 
\end{equation}
\end{enumerate}}

We have verified that this conjecture is satisfied for $n$ ranging from 4 to 8.  
If our conjecture is true, this would imply that the $\boldsymbol{{\mathcal Y\hspace{-0.1cm} \mathcal W}_{D_n}}$'s do not form a series of AMP webs, 
which would be a little disappointing, especially considering that the associated Dynkin diagrams are simply laced.
\mk

We now turn to the study of the $\mathcal X$-cluster webs of type $D$.\mk

\paragraph{\bf Numerology of the $\boldsymbol{{\mathcal X\hspace{-0.1cm}\mathcal W}_{\hspace{-0.05cm} D_n}}$'s.}
\label{Par:Numerology-XWDn}
We have seen that $\boldsymbol{{\mathcal X\hspace{-0.1cm}}\mathcal W}_{\hspace{-0.05cm} D_n}$ is a $d^\mathcal X_{D_n}$-web in $n$ variables, with  $d^\mathcal X_{D_n}={n(n-1)(n^2+4n-6)}/{6} $. By direct computations, one obtains that the ranks of $\boldsymbol{\mathcal X\hspace{-0.1cm}\mathcal W}_{\hspace{-0.05cm} D_n}$ for $n=4,\ldots,8$ are as follows (where we write $\boldsymbol{{\mathcal W}}_{D_n}$ instead of 
$\boldsymbol{{\mathcal X\hspace{-0.1cm}}\mathcal W}_{\hspace{-0.05cm} D_n}$ to save space): 
  \begin{align}
  \nonumber
   \rho^\bullet\big(\boldsymbol{{\mathcal W}}_{ D_4}\big)=& \, (48, 42, 32, 17, 1)
  &&     {\rm polrk}^\bullet\big(\boldsymbol{{\mathcal W}}_{ D_4}\big)=\big(88,47,2 \big)
 && {\rm rk}(\boldsymbol{{\mathcal W}}_{ D_4}\big)=137 \\ 
   \nonumber
    \rho^\bullet\big(\boldsymbol{{\mathcal W}}_{ D_5}\big)=& \, (125, 115,95,60,5) 
  &&     {\rm polrk}^\bullet\big(\boldsymbol{{\mathcal W}}_{ D_5}\big)=\big(235,149,8,\ldots \big)
 && 
\\
    \nonumber
      \rho^\bullet\big(\boldsymbol{{\mathcal W}}_{ D_6}\big)=&\,  (264, 249, 214, 144, 24) 
  &&     {\rm polrk}^\bullet\big(\boldsymbol{{\mathcal W}}_{ D_6}\big)=\big(504,344,20,\ldots  \big)\\
    \nonumber
      \rho^\bullet\big(\boldsymbol{{\mathcal W}}_{ D_7}\big)=&\,  (490, 469, 413, 287, 56)
  &&     {\rm polrk}^\bullet\big(\boldsymbol{{\mathcal W}}_{ D_7}\big)=\big(945, 671,40,\ldots  \big)\\
    \nonumber
         \rho^\bullet\big(\boldsymbol{{\mathcal W}}_{ D_8}\big)=&\,  (832, 804, 720, 510, 104)
  &&     {\rm polrk}^\bullet\big(\boldsymbol{{\mathcal W}}_{ D_8}\big)=\big( 1616,   1175 ,  70,\ldots \big)
\end{align}
(The ranks  $ {\rm polrk}^w(\boldsymbol{{\mathcal W}}_{ D_n}\big)$ for $w\geq 4$ and $ {\rm rk}(\boldsymbol{{\mathcal W}}_{ D_n}\big)$ are not given for $n$ from 5 to 8
 because we have been unable to compute them.   For $n=5$ already, the calculations were too huge to be carried out with the (yet powerful) machine we had at our disposal).
\mk 

Even if the table of ranks above is somehow incomplete, it is not too difficult to extrapolate from it what might be (most of) the ranks of $\boldsymbol{\mathcal X\hspace{-0.1cm}\mathcal W}_{\hspace{-0.05cm} D_n}$ for $n$ arbitrary. \sk 

Indeed, the values above satisfy the following relations:
\begin{align} \nonumber
\rho^1\big(
\boldsymbol{\mathcal X\hspace{-0.1cm}\mathcal W}_{\hspace{-0.05cm} D_n}\big)=&\, d_{D_n}^{\mathcal X}-n \\
\nonumber  
\rho^\sigma \big( \boldsymbol{\mathcal X\hspace{-0.1cm}\mathcal W}_{\hspace{-0.05cm} D_n} \big)=&\,
\rho^{\sigma-1} \big(\boldsymbol{\mathcal X\hspace{-0.1cm}\mathcal W}_{\hspace{-0.05cm} D_n}\big)-  { n-2+\sigma \choose \sigma}
\qquad \mbox{for } \, \sigma=2,3,4\, , \vspace{0.2cm}\\ 
\label{Ali:Ranks-XWDn}
{\rm polrk}^1\big(\boldsymbol{\mathcal X\hspace{-0.1cm}\mathcal W}_{\hspace{-0.05cm} D_n}\big)=&\,  2\,d_{D_n}^{\mathcal X}-\lvert (D_n)_{\geq -1}\lvert=2\,d_{D_n}^{\mathcal X}-n^2\\\nonumber
 {\rm polrk}^3\big( \boldsymbol{\mathcal X\hspace{-0.1cm}\mathcal W}_{\hspace{-0.05cm} D_n} \big)=&\, 2\, { n-1 \choose 3}\\
\mbox{and}
\quad {\rm polrk}^w\big(  \boldsymbol{\mathcal X\hspace{-0.1cm}\mathcal W}_{\hspace{-0.05cm} D_n}    \big)=&\, 0 
\quad \mbox{for any weight } \, w\geq 4\, .\nonumber
\end{align}
And we conjecture that this is the case in full generality: \mk 

{\bf Conjecture $\boldsymbol{{\mathcal X\hspace{-0.1cm} \mathcal W}_{D_n}}$.} 
{\it 
Let $n$ be an integer bigger than or equal to 4. 
\begin{enumerate}
\item 
 The ramification of $\boldsymbol{{\mathcal X\hspace{-0.1cm}\mathcal W}}_{\hspace{-0.05cm} D_n}$ is polylogarithmic and 
its  variety of common leaves  coincides with the corresponding 
 cluster arrangement:  one has 
 $
 \Sigma^c\Big(\boldsymbol{{\mathcal X\hspace{-0.1cm}\mathcal W}}_{\hspace{-0.05cm} D_n}\Big)=\boldsymbol{Arr}_{D_n}\, .
 $
\item The relations in \eqref{Ali:Ranks-XWDn} are satisfied;
\item 
  One has 
    $
  {\rm polrk}\big(\boldsymbol{\mathcal X\hspace{-0.1cm} \mathcal W}_{\hspace{-0.05cm}D_n}\big)={\rm rk}\big(\boldsymbol{\mathcal X\hspace{-0.1cm} \mathcal W}_{\hspace{-0.05cm}D_n}\big)< \rho\big(\boldsymbol{\mathcal X\hspace{-0.1cm} \mathcal W}_{\hspace{-0.05cm}D_n}\big)
  $  hence this web is not AMP but 
all its ARs  are polylogarithmic, of weight 1,2 or 3. 
\end{enumerate}}

Should it be satisfied, this conjecture does not tell us everything about the webs $\boldsymbol{{\mathcal X\hspace{-0.1cm} \mathcal W}_{D_n}}$'s, even 
for their very basic invariants. Indeed, the sequences corresponding to the fifth virtual ranks and the second polylogarithmic ranks of these webs start as follows
$$
 \rho^5\big(\boldsymbol{{\mathcal X}}_{ D_\bullet}\big)=  (1,5,24,56,\ldots) 
 \qquad \mbox{ and } \qquad 
 {\rm polrk}^2\big(\boldsymbol{{\mathcal X}}_{ D_\bullet}\big)=  (47, 149, 344, 671,\ldots)
$$
and both are quite mysterious. Even with the help of the OEIS, we haven't been able to guess a closed formula of its $n$-th term in function of $n$. 
\mk 

Note also that it would be a bit desapointing (at least for the author of theses lines) should these webs not be AMP.  Indeed, one would have expected a uniform behavior, similar to the corresponding webs in Dynkin type $A$, for cluster webs attached to simply laced Dynkin diagrams. \sk 

If  it turns out that this would have been too naive and that the $\boldsymbol{{\mathcal X}}_{ D_n}$'s are not AMP, these webs are nevertheless quite interesting since they carry trilogarithmic ARs.  It would be interesting to have a better understanding of such ARs for $n$ arbitrary.

\subsubsection{\bf Some remarks about   $\boldsymbol{{\mathcal X\hspace{-0.1cm} \mathcal W}_{\hspace{-0.05cm}D_4}}$ and $\boldsymbol{{\mathcal X\hspace{-0.1cm} \mathcal W}_{\hspace{-0.05cm}D_5}}$.}

The fact that cluster webs carry dilogarithmic ARs must be considered as  well established by now. What is more interesting and less well understood is the polylogarithmic ARs of higher weights that such webs could carry. \sk

The $\mathcal X$-cluster webs of type $D$ are interesting in this respect since, as suggests  the explicit polylogarithmics ranks of the cases $n=4,\ldots,7$ 
given at the begining of \S\ref{Par:Numerology-XWDn}, it could happen that 
each  $\boldsymbol{\mathcal X\hspace{-0.1cm} \mathcal W}_{\hspace{-0.05cm}D_n}$ carries many trilogarithmic abelian relations.

In the next two paragraphs, 
  we make a few preleminary remarks in this regard in case $D_4$ and $D_5$ respectively. 
We only touch upon these questions superficially. They definitely deserve to be studied more in depth.

\paragraph{\bf Some remarks about the trilogarithmic ARs of $\boldsymbol{{\mathcal X\hspace{-0.1cm} \mathcal W}_{\hspace{-0.05cm}D_4}}$.}
\label{Par:Remarks-trilog-AR-XWD4}

It would be interesting to know more about the two (linearly independant) trilogarithmic ARs of $\boldsymbol{\mathcal X\hspace{-0.1cm} \mathcal W}_{\hspace{-0.05cm}D_4}$.  
We introduce two trilogarithmic cluster symbols\footnote{We recall that  for $a,b,c\in \{0,1\}$, the notation  $abc$ stands for $\omega_a\otimes \omega_b\otimes \omega_c$ with $\omega_0=dz/z$ and $\omega_1=dz/(1+z)$.} by setting 
$$ 
\mathcal T_1=001-010 \qquad  \mbox{ and } \qquad 
\mathcal T_2=001-100\,.
$$

 By explicit computations, one can prove that 
\begin{enumerate}
\item  for $i=1,2$, there exists one AR of  $\boldsymbol{\mathcal X\hspace{-0.1cm} \mathcal W}_{\hspace{-0.05cm}D_4}$, 
 uniquely determined up to sign, whose non-trivial components all are trilogarithmic with symbol $\pm \mathcal T_i$;
 \item these two trilogarithmic ARs are irreducible but not complete. They share the same support which is  a certain 40-subweb of $\boldsymbol{\mathcal X\hspace{-0.1cm} \mathcal W}_{\hspace{-0.05cm}D_4}$, that we will denote by $\boldsymbol{\mathcal X\hspace{-0.1cm} \mathcal W}_{\hspace{-0.05cm}D_4}'$;
 \item  one has 
 $
\rho^\bullet\big(\boldsymbol{\mathcal X\hspace{-0.1cm} \mathcal W}_{\hspace{-0.05cm}D_4}'\big)=(36, 30, 20, 8, 1,0)$ and ${\rm polrk}^\bullet\big(\boldsymbol{\mathcal X\hspace{-0.1cm} \mathcal W}_{\hspace{-0.05cm}D_4}'\big)=(64,29,2)$  thus $
 \rho\big(\boldsymbol{\mathcal X\hspace{-0.1cm} \mathcal W}_{\hspace{-0.05cm}D_4}'\big)
 ={\rm polrk}\big(\boldsymbol{\mathcal X\hspace{-0.1cm} \mathcal W}_{\hspace{-0.05cm}D_4}'\big)=95$, hence $\boldsymbol{\mathcal X\hspace{-0.1cm} \mathcal W}_{\hspace{-0.05cm}D_4}'$ is AMP with only polylogarithmic ARs (of weight $\leq 3$).
\end{enumerate}

Actually, the web $\boldsymbol{\mathcal X\hspace{-0.1cm} \mathcal W}_{\hspace{-0.05cm}D_4}'$ as well as the fact that it carries a trilogarithmic identity, was known already: it is the web ${\boldsymbol{\mathcal W}}_{\mathcal G \mathcal G \mathcal S \mathcal V\mathcal V}$ considered in \cite{GSVV} and  discussed above in 
\S\ref{Par:GSVV}. \bk 

Viewed the properties it satisfies (carrying the two polylogarithmic ARs 
 of higher rank and being AMP), $\boldsymbol{\mathcal X\hspace{-0.1cm} \mathcal W}_{\hspace{-0.05cm}D_4}'$ appears as being more interesting than the complete $\mathcal X$-cluster web $\boldsymbol{\mathcal X\hspace{-0.1cm} \mathcal W}_{\hspace{-0.05cm}D_4}$ itself and for that reason, it is interesting to understand it better. \sk

It turns out that $\boldsymbol{\mathcal X\hspace{-0.1cm} \mathcal W}_{\hspace{-0.05cm}D_4}'$ can be constructed  from the complete $\boldsymbol{\mathcal X\hspace{-0.1cm} \mathcal W}_{\hspace{-0.05cm}D_4}$ independently of any consideration about abelian relations. Indeed, let $R$ be the cyclic (of order 3) symmetry of $D_4$ and denote by 
$\boldsymbol{\mathcal Y\hspace{-0.1cm} \mathcal W}_{\hspace{-0.05cm}D_4}^R$ the subweb of $\boldsymbol{\mathcal Y\hspace{-0.1cm} \mathcal W}_{\hspace{-0.05cm}D_4}$
whose first integrals are the $Y[\alpha]$'s for $\alpha\in (D_4)_{\geq -1}$ invariant by $R$.  
There are 4 such roots\footnote{These $R$-invariant roots  are the one with coordinates $(0,\pm 1,0,0)$, $(1,1,1,1)$ and $(1,2,1,1)$ relatively to the standard choice of the positive roots of the root system of type $D_4$.} hence $\boldsymbol{\mathcal Y\hspace{-0.1cm} \mathcal W}_{\hspace{-0.05cm}D_4}^R$ is a 4-web and it can be verified that 
\begin{equation}
\label{Eq:XWD4'}
\boldsymbol{\mathcal X\hspace{-0.1cm} \mathcal W}_{\hspace{-0.05cm}D_4}'=
\Big( \boldsymbol{\mathcal X\hspace{-0.1cm} \mathcal W}_{\hspace{-0.05cm}D_4} \setminus \boldsymbol{\mathcal Y\hspace{-0.1cm} \mathcal W}_{\hspace{-0.05cm}D_4}\Big)\sqcup  \boldsymbol{\mathcal Y\hspace{-0.1cm} \mathcal W}_{\hspace{-0.05cm}D_4}^R\, .
\end{equation}

However, this description of $
\boldsymbol{\mathcal X\hspace{-0.1cm} \mathcal W}_{\hspace{-0.05cm}D_4}'$ does not tell much and several natural questions do arise: 
\begin{questions}
\begin{enumerate}
\item Is there a way to generate $\boldsymbol{\mathcal X\hspace{-0.1cm} \mathcal W}_{\hspace{-0.05cm}D_4}'$ by means of certain mutations and/or choices of cluster variables?\footnote{Considering that $D_4$ is mutation equivalent to 
$A_2\boxtimes A_2$, what is discussed 
\S\ref{SS:PlabicWebsAA} further might be used to answer.} 
\item Is it possible to deduce from the description \eqref{Eq:XWD4'} of $\boldsymbol{\mathcal X\hspace{-0.1cm} \mathcal W}_{\hspace{-0.05cm}D_4}'$ that this web carries trilogarithmic ARs?
\item Does this generalizes to 
$\boldsymbol{\mathcal X\hspace{-0.1cm} \mathcal W}_{\hspace{-0.05cm}D_n}$
 for any higher $n$ ?\mk 
\end{enumerate}
\end{questions}


\paragraph{\bf About $\boldsymbol{{\mathcal X\hspace{-0.1cm} \mathcal W}_{\hspace{-0.05cm}D_5}}$.}
This web 
  is a  130-web in 5 variables,  
such that 
$$
    \rho^\bullet\big(\boldsymbol{{\mathcal X\hspace{-0.1cm} \mathcal W}}_{\hspace{-0.05cm}D_5}\big)= (125, 115,95,60,5)  \qquad  
\mbox{ and } 
 \qquad      {\rm polrk}^\bullet\big(\boldsymbol{{\mathcal X\hspace{-0.1cm} \mathcal W}}_{\hspace{-0.05cm}D_5}\big)=\big(235,149,8 \big)
$$
thus $\rho\big(\boldsymbol{{\mathcal X\hspace{-0.1cm} \mathcal W}}_{\hspace{-0.05cm}D_5}\big)=400$  and 
${\rm polrk}\big(\boldsymbol{{\mathcal X\hspace{-0.1cm} \mathcal W}}_{\hspace{-0.05cm}D_5}\big)=392$. We haven't been able to determine $
{\rm rk}\big(\boldsymbol{{\mathcal X\hspace{-0.1cm} \mathcal W}}_{\hspace{-0.05cm}D_5}\big)$ (the  computations were too heavy) so we do not know if this web is AMP, has only polylogarithmic ARs, or whether  something else happens. \mk

The space of trilogarithmic ARs of $\boldsymbol{{\mathcal X\hspace{-0.1cm} \mathcal W}}_{\hspace{-0.05cm}D_5}$ has dimension 8 and admits a basis $(\beta_i)_{i=1}^8$ which can be explicitly constructed. 
In particular, $\beta_1$ and $\beta_2$ can be chosen such that they form a basis of the trilog ARs of a 40-subweb of $\boldsymbol{{\mathcal X\hspace{-0.1cm} \mathcal W}}_{\hspace{-0.05cm}D_5}$, denoted by 
$\boldsymbol{{\mathcal X\hspace{-0.1cm} \mathcal W}}_{\hspace{-0.05cm}D_5}^{12}$ such that 
$$
\rho^\bullet\Big( \boldsymbol{{\mathcal X\hspace{-0.1cm} \mathcal W}}_{\hspace{-0.05cm}D_5}^{12}   \Big)=(36, 30, 20, 8, 1)\qquad 
\rho\Big(  \boldsymbol{{\mathcal X\hspace{-0.1cm} \mathcal W}}_{\hspace{-0.05cm}D_5}^{12}   \Big)=95 
\qquad 
{\rm polrk}^\bullet\Big(\boldsymbol{{\mathcal X\hspace{-0.1cm} \mathcal W}}_{\hspace{-0.05cm}D_5}^{12}\Big)=(64, 29, 2)\, .
$$
This shows that $\boldsymbol{{\mathcal X\hspace{-0.1cm} \mathcal W}}_{\hspace{-0.05cm}D_5}^{12}$ has AMP rank. Moreover, it can be verified that 
this web has intrinsic dimension equal to four hence 
it is likely equivalent to a pull-back of the  40-subweb $\boldsymbol{{\mathcal X\hspace{-0.1cm} \mathcal W}}_{\hspace{-0.05cm}D_4}'$ considered in the previous paragraph. \mk

More generally, it can be verified that the union of the supports of the height $\beta_i$'s is that of a strictly proper 100-subweb of 
$\boldsymbol{{\mathcal X\hspace{-0.1cm} \mathcal W}}_{\hspace{-0.05cm}D_5}$, 
denoted by $\boldsymbol{{\mathcal X\hspace{-0.1cm} \mathcal W}}_{\hspace{-0.05cm}D_5}'$, and such that 
\begin{align*}
\rho^\bullet\big(\boldsymbol{{\mathcal X\hspace{-0.1cm} \mathcal W}}_{\hspace{-0.05cm}D_5}'\big)=&\, (95, 85, 65, 35, 4)&& 
{\rm polrk}\big(  \boldsymbol{\mathcal X\hspace{-0.1cm} \mathcal W}_{\hspace{-0.05cm}D_5}'\big)= \big(175, 99, 8\big) \\ 
\rho\big(\boldsymbol{\mathcal X\hspace{-0.1cm} \mathcal W}_{\hspace{-0.05cm}D_5}'\big)= \, &284 && 
{\rm polrk}\big(  \boldsymbol{\mathcal X\hspace{-0.1cm} \mathcal W}_{\hspace{-0.05cm}D_5}'\big)=282 &&\mbox{ and }
{}^{} \hspace{0.7cm} {\rm rk}\big(  \boldsymbol{\mathcal X\hspace{-0.1cm} \mathcal W}_{\hspace{-0.05cm}D_5}'\big)=282\, . 
\end{align*}

Thus $\boldsymbol{\mathcal X\hspace{-0.1cm} \mathcal W}_{\hspace{-0.05cm}D_5}'$ is almost of maximal possible rank, but not quite. This contrasts with 
$\boldsymbol{\mathcal X\hspace{-0.1cm} \mathcal W}_{\hspace{-0.05cm}D_4}'$.

\begin{questions}
\begin{enumerate}
\item Is the space of trilogarithmic ARs of  $\boldsymbol{{\mathcal X\hspace{-0.1cm} \mathcal W}}_{\hspace{-0.05cm}D_5}$ spanned  
by ARs carried by subwebs isomorphic to $\boldsymbol{{\mathcal X\hspace{-0.1cm} \mathcal W}}_{\hspace{-0.05cm}D_4}'$?  More generally, are the 
general  trilogarithmic ARs of  $\boldsymbol{{\mathcal X\hspace{-0.1cm} \mathcal W}}_{\hspace{-0.05cm}D_5}$ accessible from the two trilogarithmic ARs of $\boldsymbol{{\mathcal X\hspace{-0.1cm} \mathcal W}}_{\hspace{-0.05cm}D_4}'$ of this type? 
\item  Is there a way to define/construct algebraically $\boldsymbol{{\mathcal X\hspace{-0.1cm} \mathcal W}}_{\hspace{-0.05cm}D_5}'$ as a subweb of $\boldsymbol{{\mathcal X\hspace{-0.1cm} \mathcal W}}_{\hspace{-0.05cm}D_5}$ without calling upon the notion of (polylogarithmic) abelian relation?
\end{enumerate}
\end{questions}
\begin{center}
$\star$
\end{center}

To finish our discussion about cluster webs in type $D$, let us say that they remain  still quite mysterious and that many elementary questions still arise about them. As an example, one can mention the secondary $\mathcal X$-cluster webs of this kind, namely the $\boldsymbol{{\mathcal U\hspace{-0.05cm}\mathcal X\hspace{-0.05cm} \mathcal W}}_{\hspace{-0.05cm}D_n}$'s. According to Theorem 
\ref{T:classical-cluster-webs}, $\boldsymbol{{\mathcal U\hspace{-0.05cm}\mathcal X\hspace{-0.05cm} \mathcal W}}_{\hspace{-0.05cm}D_4}$ is equivalent to Kummer's tetralogarithmic web hence carries tetralogarithmic ARs. It could be the case that this web is the first of a whole series of webs enjoying this interesting property. This question is wide open. Actually, nothing is known about $\boldsymbol{{\mathcal U\hspace{-0.05cm}\mathcal X\hspace{-0.05cm} \mathcal W}}_{\hspace{-0.05cm}D_n}$ for $n\geq 4$. For instance, having (even conjectural) statements about its numerology ({\it e.g.}\,closed formulas for its ranks) would already be  interesting.

 \subsection{\bf About the $\boldsymbol{\mathcal Y}$-cluster webs of  bi-Dynkin type.}
 \label{SS:YW-Delta-Delta'}

A central topic in this memoir is  constructing new AMP webs from cluster algebras. 
A relevant notion in this respect is that of cluster period, since to each period of a cluster algebra is associated a dilogarithmic AR of the associated web. 
And as the examples considered above show, numerous webs of this kind are AMP, with only polylogarithmic ARs, of weight 1 or 2. 
\mk 

In this subsection, we examine in more detail the $\mathcal Y$-cluster webs associated with pairs of Dynkin diagrams.  While it appears that not all of these webs are AMP, many of them are which leads us to formulate several conjectures about them.

 \subsubsection{\bf  Notations.}
We introduce/recall quickly some notations we will use in the lines below: 
\begin{itemize}
\item[$-$] $\Delta$ and $\Delta'$ are Dynkin diagrams, of  rank  $n$ and $n'$ respectively; 
\item[$-$]
To simplify the discussion, we will assume that $\Delta$ and $\Delta'$
are such that the type of the former is smaller than that of the latter diagram (with respect to the lexicographical order) and that $n\leq n'$ when both $\Delta$ and $\Delta'$ have the same type;
\item[$-$] $h$ and $h'$ stand respectively for the Coxeter numbers\footnote{The Coxeter numbers of Dynkin diagrams are given in Table \ref{Table:CoxeterNumbers} above.} of $\Delta$ and $\Delta'$;
\item[$-$] $\Delta \boxtimes \Delta'$ (resp.\,$\Delta \square \Delta'$) denotes the `triangle  product' (resp.\,the `square product') of the alternating Dynkin quivers $\vec{\Delta}$ and $\vec{\Delta}'$ (as defined in \cite[\S3.3]{KellerAnnals}, see also  \S\ref{Par:Interpretation-cluster-Y-systems-Dynkin-type} above);
\item[$-$]  We use the simpler notation $\mu_{\Delta,\Delta'}$ for   the period $i_{\Delta,\Delta'}^{half}$ defined in Remark 
\ref{Rem:Half-Periods};
\item[$-$] $\boldsymbol{{\mathcal Y\hspace{-0.1cm} \mathcal W}}_{\hspace{-0.05cm}\Delta,\Delta'}$ denotes the cluster web associated to the period $\mu_{\Delta,\Delta'}$;
\item[$-$] $\boldsymbol{{\mathcal Y\hspace{-0.1cm} \mathcal W}}_{\hspace{-0.05cm}\Delta,\Delta'}$ is   the cluster web associated to the period $i_{\Delta,\Delta'}^{half}$ defined in Remark 
\ref{Rem:Half-Periods};
\item[$-$] $d_{\Delta,\Delta'}$ stands for the degree of $\boldsymbol{{\mathcal Y\hspace{-0.1cm} \mathcal W}}_{\hspace{-0.05cm}\Delta,\Delta'}$;
\item[$-$] $\boldsymbol{F}_{\Delta,\Delta'}$ denotes the set of $F$-polynomials of 
$\boldsymbol{{\mathcal Y\hspace{-0.1cm} \mathcal W}}_{\hspace{-0.05cm}\Delta,\Delta'}$
({\it cf.}\, \eqref{Eq:F(W)} for a definition);
\item[$-$] $\boldsymbol{Arr}_{\Delta,\Delta'}$ is the arrangement of hypersurfaces 
in the initial torus
  cut out by 
   $\boldsymbol{F}_{\Delta,\Delta'}$.
\item[$-$] to save space, for $`{\rm rrr}'$ denoting $`\rho'$, $`{\rm polrk}'$ of $`{\rm rk}'$, we will sometimes denote them by ${\rm rrr}(\Delta,\Delta')$ (or ${\rm rrr}^\bullet(\Delta,\Delta')$) instead of 
${\rm rrr}\big(\boldsymbol{{\mathcal Y\hspace{-0.1cm} \mathcal W}}_{\hspace{-0.05cm}\Delta,\Delta'}\big)$.
\item[$-$] For $U,U'\in \{A,B,C,D,E,F,G\}$,  an object associated to a pair $(\Delta,\Delta')$ will be said of type $UU'$ if the Dynkin diagrams $\Delta$ and $\Delta'$ are of type $U$ and $U'$ respectively. 
\end{itemize}


We recall  that it has been conjectured ({\it cf.}\,Conjecture \ref{P:Y-web-type-Delta,Delta')}) that one has   $d_{\Delta,\Delta'}=nn'(h+h')/2$ for any pair $(\Delta,\Delta')$. This has been established   for any $\Delta'$   when $\Delta=A_1$ (see Corollary \ref{Coro:d-X-Delta--d-Y-Delta})   and in the case when both $\Delta$ and $\Delta'$ are of type $A$ (see Corollary \ref{Cor:Y-web-type-(Am,An)}).\mk 

Since the case when one of the two considered Dynkin diagrams is $A_1$ (or equivalently, has rank 1) has been extensively discussed in preceding sections, we assume from now on that both 
$\Delta$ and $\Delta'$ 
have rank bigger than 1:  one has 
$$
n={\rm rk}\big(\Delta\big)\geq 2 \quad \qquad \mbox{ and } \quad \qquad n'={\rm rk}\big(\Delta'\big)\geq 2\, .
$$

 \subsubsection{\bf  Some conjectural statements.}
 \label{SS:SomeConjecturalStatements}
Given the two following facts: 
\begin{enumerate}
\item  the  $\mathcal Y$-cluster webs 
 $\boldsymbol{{\mathcal Y\hspace{-0.1cm} \mathcal W}}_{\hspace{-0.03cm}A_n}$'s 
 all are AMP with only logarithmic and (one) dilogarithmic ARs (as established in \S\ref{SS:Y-cluster-web-An}); 
\item  the  $\mathcal Y$-cluster webs of bi-Dynkin type $AA$ (hence in particular the one in point 1.) can be described in a uniform manner  
 by means of projective geometry (more precisely, in terms of forgetful maps on 
 moduli spaces of projective configurations of points, {\it cf.}\,\S\ref{SS:Y-cluster-Web-type-AA}), 
 \end{enumerate}
it is natural to wonder whether  the $\boldsymbol{{\mathcal Y\hspace{-0.1cm} \mathcal W}}_{A_m,A_n}$'s are AMP or not, and by a straightforward generalization, to consider the same question for pairs with $\Delta,\Delta'$ of arbitrary types. \mk 

The explicit study on numerous cases that we have done (see below for some examples) leads us to the conclusion that the answer is different according to the bi-type $UU'$ of $(\Delta,\Delta')$, but  is independent of the ranks $n,n'$ for each fixed bi-type (possibly modulo a finite number of exceptions).  To provide a framework for the discussion to come, it is useful to consider the following series of statements which make sense mathematically for each Dynkin pair $(\Delta,\Delta')$ as above: \bk 

{\bf \underline{Statements }$\boldsymbol{{{\mathcal Y\hspace{-0.1cm} \mathcal W}}_{\Delta,\Delta'}}$.}   
{\it 
\begin{enumerate}
\item[{\bf 1.}]{\bf ${\big[}$Degree\big].} 
$\boldsymbol{\mathcal Y\hspace{-0.1cm} \mathcal W}_{\Delta,\Delta'}$ is a  
$d_{\Delta , \Delta'}$-web in $nn'$-variables with 
$$d_{\Delta ,\Delta'}=d_{\Delta ,\Delta'}^{\boldsymbol{\mathcal Y}}=\frac{1}{2}nn'\big(h+h'\big)\, .$$
\item[{\bf 2.}]{\bf ${\big[}${$\boldsymbol{F}$-polynomials}\big].} 
The set  
of  $F$-polynomials 
$\boldsymbol{F}_{\Delta,\Delta'}$ 
 has cardinality $d_{\Delta ,\Delta'}-nn'$. 
\sk 

Moreover, for any cluster variable $x$ element of  $\boldsymbol{\mathcal Y\hspace{-0.1cm} \mathcal W}_{\Delta,\Delta'}$, the set of irreducible factors of the non-monomial part of $1+x$ 
 is included in $\boldsymbol{F}_{\Delta,\Delta'} $.
\item[{\bf 3.}] {\bf ${\big[}$Ramification\big].}
 For any cluster first integral $x$   of $\boldsymbol{\mathcal Y\hspace{-0.1cm} \mathcal W}_{\Delta,\Delta'}$, one has
 $ \mathfrak B_x=\{0,-1,\infty\}$ and consequently, this web has polylogarithmic ramification. 
 Moreover, 
the  variety of common leaves $\Sigma^c\big(\boldsymbol{\mathcal Y\hspace{-0.1cm} \mathcal W}_{\Delta,\Delta'}\big)$ 
 coincides with the corresponding cluster arrangement:  one has 
 $$
\Sigma^c\Big(\boldsymbol{\mathcal Y\hspace{-0.1cm} \mathcal W}_{\Delta,\Delta'}\Big)=\boldsymbol{Arr}_{\Delta,\Delta'}\, .
 $$
\item[{\bf 4.}]  {\bf ${\big[}$Virtual rank(s)\big].} The virtual ranks  
of $\boldsymbol{\mathcal Y\hspace{-0.1cm} \mathcal W}_{\Delta,\Delta'}$ are  given by  
 $$
\rho^\bullet\Big(\boldsymbol{\mathcal Y\hspace{-0.1cm} \mathcal W}_{\Delta,\Delta'}  \Big)= \Big(\, d_{\Delta , \Delta'}-nn'\, , \,  nn'\, , \, 
1\, \Big)\qquad \mbox{ hence } \qquad 
\rho\Big(\boldsymbol{\mathcal Y\hspace{-0.1cm} \mathcal W}_{\Delta,\Delta'}  \Big)= d_{\Delta , \Delta'}+
1\, .  $$ 
 \item[{\bf 5.}] {\bf ${\big[}$Polylogarithmic rank(s)\big].}  Regarding the polylogarithmic ranks of 
 $\boldsymbol{\mathcal Y\hspace{-0.1cm} \mathcal W}_{\Delta,\Delta'}$, one has 
\begin{align*}
{\rm polrk}^\bullet\Big(\boldsymbol{\mathcal Y\hspace{-0.1cm} \mathcal W}_{\Delta,\Delta'}
\Big)= &\,  \Big( \, d_{\Delta , \Delta'}\, , \, 1 \, \Big)
 \qquad \mbox{ therefore } \qquad 
 {\rm polrk}\Big(\boldsymbol{\mathcal Y\hspace{-0.1cm} \mathcal W}_{\Delta,\Delta'}\Big)=d_{\Delta , \Delta'}+1\, . 
\end{align*}
\item[{\bf 6.}]  {\bf ${\big[}$Dilogarithmic AR{\big]}.} 
The space 
$
\boldsymbol{Dilo\hspace{-0.03cm}g\hspace{-0.03cm}A\hspace{-0.03cm}R}\big(\boldsymbol{\mathcal Y\hspace{-0.1cm} \mathcal W}_{\Delta,\Delta'} \big) 
$
  is 1-dimensional, spanned  by $ {\sf R}_{\Delta,\Delta'}$.  
\item[{\bf 7.}] {\bf ${\big[}$Rank and basis of ARs{\big]}.} 
  One has 
    $$
    {\rm polrk}\big(\boldsymbol{\mathcal Y\hspace{-0.1cm} \mathcal W}_{\Delta,\Delta'}\big)={\rm rk}\big(\boldsymbol{\mathcal Y\hspace{-0.1cm} \mathcal W}_{\Delta,\Delta'}\big)=\rho\big(\boldsymbol{\mathcal Y\hspace{-0.1cm} \mathcal W}_{\Delta,\Delta'}\big)
=d_{\Delta , \Delta'}+1
  $$
 hence 
$\boldsymbol{\mathcal Y\hspace{-0.1cm} \mathcal W}_{\Delta,\Delta'}$ is AMP 
 and 
 the following decomposition in direct sum holds true:
\begin{equation}
\label{Eq:AR-YWDeltaDelta'}
\boldsymbol{\mathcal A}\Big(\boldsymbol{\mathcal Y\hspace{-0.1cm} \mathcal W}_{\Delta,\Delta'}\Big) = 
\boldsymbol{Lo\hspace{-0.03cm}g\hspace{-0.03cm}A\hspace{-0.03cm}R}\Big(\boldsymbol{\mathcal Y\hspace{-0.1cm} \mathcal W}_{\Delta,\Delta'} \Big) 
\,  \oplus  \,
 \big\langle
\, {\sf R}_{\Delta,\Delta'}\, 
 \big\rangle\, . 
\end{equation}
\end{enumerate}}

 \subsubsection{\bf  When one of the Dynkin diagrams is exceptional.}
 \label{SS:OneDynkinDiagramsExceptional}
By considering several examples, we first are going to show that 
 the preceding statements are not all satisfied in full generality. It even seems to be the rule when 
one of the Dynkin diagrams is of exceptional
type : this is the case that we are going to discuss now. 
\sk 

If one of the Dynkin diagrams is of exceptional type, then necessarily  $\Delta'$ is of type $E,F$ or $G$ (since we assume $\Delta\leq \Delta'$). We have looked at the corresponding webs when $\Delta$ has rank 2. \sk 

Here are the results of some of our computations: 
\begin{itemize} 
\item 
The case when  $\Delta'=G_2$ and  $\Delta$ is $A_2$
 or $B_2$ have already be considered above, see  
\S\ref{Subpar:YWA2*G2} and 
\S\ref{Subpar:YWB2*G2} respectively. These webs are not AMP.  The same occurs when $\Delta=G_2$  as well since $\rho^\bullet(G_2,G_2)=(20,14,4)$ and 
${\rm polrk}^\bullet(G_2,G_2)=(24,1)$.
\item As  for $\boldsymbol{\mathcal Y\hspace{-0.04cm}\mathcal W}_{\hspace{-0.04cm}A_2,F_4}$, it is a 60-web in 8 variables such that  $\rho^\bullet= (52, 24, 1)$ and 
${\rm polrk}^\bullet=(60,1)$. 
\item For the case when $(\Delta,\Delta')=(A2,E_6)$:
$\boldsymbol{\mathcal Y\hspace{-0.04cm}\mathcal W}_{\hspace{-0.04cm}A_2,E_6}$
is a 90-web in 12 variables with 
$\rho^\bullet\big(A_2 ,  E_6\big)=\big( 78, 12,1\big)$ and 
${\rm polrk}^\bullet\big(A_2 ,  E_6\big)=(90,1)$. It follows that this web is AMP. 
\item  The web $\boldsymbol{\mathcal Y\hspace{-0.04cm}\mathcal W}_{\hspace{-0.04cm}A_2,E_7}$
is a 147-web in 14 variables such that  $\rho^\bullet(A_2,E_7)= (133, 42, 1)$. 
\end{itemize} 

Even if they are not many and incomplete (the computations to determine ${\rm polrk}^\bullet(A_2,E_6)$ or to investigate the case when $\Delta'=E_8$ were too heavy), the values above  strongly suggest that some of the 
{\bf {Statements }$\boldsymbol{{{\mathcal Y\hspace{-0.1cm} \mathcal W}}_{\Delta,\Delta'}}$} do not hold true when  $ \Delta'$ is exceptional and distinct from   $E_6$. In particular, these webs do not seem to be AMP. The case when $\Delta'=E_6$ is interesting since the seven {\bf {Statements }$\boldsymbol{{{\mathcal Y\hspace{-0.1cm} \mathcal W}}_{A_2,E_6}}$} are satisfied. It would be interesting to deal with the case $(A_3,E_6)$ to know if the same holds true. If yes, 
one might think that this is still the case for any 
web $\boldsymbol{{\mathcal Y\hspace{-0.1cm} \mathcal W}}_{A_n,E_6}$. \mk 

Anyway, since the case when one of the involved Dynkin diagrams is exceptional seems to be problematic in what concerns the validity of the preceding statements, we will leave this case aside from now on.




 \subsubsection{\bf  Both $\boldsymbol{\Delta}$ and $\boldsymbol{\Delta'}$ are of classical type.}
 \label{SS:BothDynkinDiagramsClassical}
We are considering here the case when $\boldsymbol{\Delta}$ and $\boldsymbol{\Delta'}$ are of type $A,D, C$ or $D$. \sk 

In  \S\ref{Subpar:YWB2*B2}, we considered the case of $\boldsymbol{\mathcal Y\hspace{-0.1cm} \mathcal W}_{B_2,B_2}$ and saw that if it is AMP, it admits two ARs not of polylogarithmic type and is such that $\rho^2\big( \boldsymbol{\mathcal Y\hspace{-0.1cm} \mathcal W}_{B_2,B_2}\big)=6>4=nn'$. The situation is even worse  with  $\boldsymbol{\mathcal Y\hspace{-0.1cm} \mathcal W}_{B_2,B_3}$ which is not AMP.\footnote{One has $\rho^\bullet(B_2,B_3)=(24, 9, 1)$, ${\rm polrk}^\bullet(B_2,B_3)=(30, 1)$ and 
${\rm rk}^\bullet(B_2,B_3)=32<34=\rho^\bullet(B_2,B_3)$.} These two examples might let us think that when both $\Delta$ and $\Delta'$ are multi-laced, then things are not as nice as they could be if the {\bf {Statements }$\boldsymbol{{{\mathcal Y\hspace{-0.1cm} \mathcal W}}_{\Delta,\Delta'}}$} were satisfied.\mk    

Letting this case aside too,  we then focus on the one when at most one of the two Dynkin diagrams is simply-laced.  More explicitly, we consider the following bi-Dynkin types: 
\begin{equation}
\label{Eq:AA->CD}
AA\,,  \quad  
AB\,,  \quad 
AC\,,  \quad 
AD\,, \quad 
BD \,, \quad 
CD 
\quad \mbox{ and } \quad 
DD\, .
\end{equation}

We have considered many cases of this kind, in particular the following:
$(A_n,A_{n'})$ for $n,n'\leq 4$,    $(A_n,B/C_{n'})$ with $n\leq 4$ and $n'\leq 3$, $(A_n,D_{n'})$ 
and $(B/C_n,D_{n'}) $
with $n\leq 3$ and $n'=4$ and $(n,n')=(2,5)$  in the  last two cases, and also $(D_4,D_4)$. In all these cases, we have verified that the seven {\bf {Statements }$\boldsymbol{{{\mathcal Y\hspace{-0.1cm} \mathcal W}}_{\Delta,\Delta'}}$} are satisfied. This leads us to state the 
\mk

{\bf Conjecture $\boldsymbol{{\mathcal Y\hspace{-0.1cm} \mathcal W}_{\hspace{-0.04cm} \Delta,\Delta'}}$.}   
{\it If both $\Delta$ and $\Delta'$ are of classical type with at least one of them 
simply-laced ({\it i.e.}\,we are in one of the  cases listed in \eqref{Eq:AA->CD}), then all the {\bf {`Statements }$\boldsymbol{{{\mathcal Y\hspace{-0.1cm} \mathcal W}}_{\hspace{-0.03cm}\Delta,\Delta'}}$'} are satisfied.}\mk

If true, this conjecture would tell us in particular that many of the series of $\mathcal Y$-cluster webs of bi-Dynkin types are AMP. \mk 

If we believe that the previous Conjecture indeed holds true, we must confess that we don't have any insight of any conceptual reason explaining why it is the case when the Dynkin diagrams are in the list  \eqref{Eq:AA->CD}, and why this (seemingly) does not happen for most of the other pairs. For instance, 
we have no idea  why  the  $
\boldsymbol{\mathcal Y\hspace{-0.1cm} \mathcal W}_{\hspace{-0.04cm}D_n}$'s (case $(A_1,D_n)$)  
might all be not AMP whereas  it may be the case for the $ \boldsymbol{\mathcal Y\hspace{-0.1cm} \mathcal W}_{\hspace{-0.04cm}A_2,D_n}$'s.\footnote{These two statements are conjectural in full generality when we write this (see {\bf Conjecture $\boldsymbol{{\mathcal Y\hspace{-0.1cm} \mathcal W}_{D_n}}$} in \S\ref{SSS:NumerologyClusterWebs-type-D} in case $(A_1,D_n)$) but we checked them for $n\leq 8$ in bi-type $(A_1,D_n)$ and for $n=5,6$ in type $(A_2,D_n)$.} This is unexpected and seems rather mysterious in our view.\mk 

We think that the cluster webs of bi-Dynkin types are very interesting webs and that there is still a lot to understand about them.   We discuss supplementary questions regarding them in the next subsection.

 \subsubsection{\bf  Complementary questions.}
 \label{SS:OneDynkinDiagramsExceptional}
We have seen that, in some cases, the 
{\bf {`Statements }$\boldsymbol{{{\mathcal Y\hspace{-0.1cm} \mathcal W}}_{\hspace{-0.03cm}\Delta,\Delta'}}$'} are not satisfied. But a closer look shows that this only concerns some of them, and not all.   
This observation led us to  state the following 
\begin{conjecture}For any pair of Dynkin diagrams $(\Delta,\Delta')$:
\begin{itemize}
\item  The three first  
{\bf {Statements }$\boldsymbol{{{\mathcal Y\hspace{-0.1cm} \mathcal W}}_{\hspace{-0.05cm}\Delta,\Delta'}}$}  
are satisifed.
\item  One has $\rho^1\big(\boldsymbol{\mathcal Y\hspace{-0.1cm} \mathcal W}_{\hspace{-0.05cm} \Delta,\Delta'})=d_{\Delta,\Delta'}-nn'$ and ${\rm polrk}^1\big(\boldsymbol{\mathcal Y\hspace{-0.1cm} \mathcal W}_{\hspace{-0.05cm} \Delta,\Delta' })=d_{\Delta,\Delta'}$. 
\end{itemize}
\end{conjecture}

The most important points in the previous conjecture are the fact that the ramification of the web is only polylogarithmic and that $\boldsymbol{F}_{\hspace{-0.04cm}\Delta,\Delta'}$ is of cardinality $d_{\Delta,\Delta'}$. For instance, the equality ${\rm polrk}^1\big(\boldsymbol{\mathcal Y\hspace{-0.1cm} \mathcal W}_{\hspace{-0.05cm} \Delta,\Delta' })=d_{\Delta,\Delta'}$ would follow quite easily from this.  \sk 

Since the degree $d_{\Delta,\Delta'}$ of the web  is conjecturally equal  to the length $nn'(h+h')/2$ of the period used to construct $\boldsymbol{\mathcal Y\hspace{-0.1cm} \mathcal W}_{\hspace{-0.05cm} \Delta,\Delta' }$, one can ask if the same statement holds true more generally

\begin{question}
Let $\boldsymbol{\mu}$ be
an irreducible cluster period of length $k$ such that the 
associated   $\mathcal P$-cluster web  $\boldsymbol{\mathcal W}_\mu= 
\boldsymbol{\mathcal W}\big(\, x_\mu(s)\, \lvert \, s=1,\ldots,k \, \big)
$ has degree $k$.\footnote{The assumption that the degree of $\boldsymbol{\mathcal W}_{\boldsymbol{\mu}} $  
coincides with the length of $\mu$ 
might be unnecessary, {\it cf.}\,Conjecture \ref{Conj:Diff-Indep-in-periods}.} 
 Has  $\boldsymbol{F}(\boldsymbol{\mathcal W}_{\boldsymbol{\mu}})$ cardinality $k$ and  does it contain the set of irreducible factors of $1+x$ for any cluster first integral $x$ of 
 $\boldsymbol{\mathcal W}_\mu$?
\end{question}

The first of the `{\bf {Statements }$\boldsymbol{{{\mathcal Y\hspace{-0.1cm} \mathcal W}}_{\hspace{-0.05cm}\Delta,\Delta'}}$}'
which is not satisfied in all the cases we have considered 
 is the formula $\rho^2(\Delta,\Delta)=nn'$ which does not always hold true, for  instance for pairs $(A_2,\Delta')$ with $\Delta'$ of  exceptional type.  
It would be interesting to have a general formula for the second virtual rank in any bi-type: 
\mk 

{\bf Problem.}
 {\it Find a closed formula for $\rho^2(\Delta,\Delta)$ in terms of the two involved Dynkin diagrams.}
\mk 

To finish our discussion of the $\mathcal Y$-webs of bi-Dynkin types, we list below some 
properties that a web $\boldsymbol{\mathcal Y\hspace{-0.1cm} \mathcal W}_{\Delta,\Delta'}$ may or may not satisfy,  which may be added to the seven {\bf \underline{Statements }$\boldsymbol{\mathcal Y\hspace{-0.1cm} \mathcal W}_{\Delta,\Delta'}$}. They are direct generalizations of some properties already considered before about the $\mathcal Y$-cluster webs $\boldsymbol{\mathcal Y\hspace{-0.1cm} \mathcal W}_{\Delta}$ associated to one Dynkin type (which can equally be seen as those of bi-Dynkin type $(A_1,\Delta)$).  We start by stating them in a row and comment on them right afterwards: 
\mk


\begin{enumerate}
\item[{\bf 8.}]  {\bf ${\big[}$Logarithmic ARs{\big]}.} {\it The logarithmic abelian relations $LogAR_{\boldsymbol{i}}(t)$ for $(\boldsymbol{i}, t)\in \mathfrak I^+_{
\hspace{-0.04cm}
\Delta ,  \Delta'  }$ form a basis of 
$\boldsymbol{Lo\hspace{-0.03cm}g\hspace{-0.03cm}A\hspace{-0.03cm}R}\big(\boldsymbol{\mathcal Y\hspace{-0.1cm} \mathcal W}_{\Delta,\Delta'} \big)$ 
 and consequently, one has ${\rm polrk}^1( \boldsymbol{\mathcal Y\hspace{-0.1cm} \mathcal W}_{\Delta,\Delta'})=
  d_{\Delta , \Delta'}$.} \sk
\item[{\bf 9.}]
{\bf ${\big[}$Monodromy of $ \boldsymbol{{\sf R}_{\Delta,\Delta'}}${\big]}.} 
{\it 
By monodromy,  the dilogarithmic AR $ {\sf R}_{\Delta,\Delta'}$ generates the whole space of polylogarithmic abelian relations which is 
$\boldsymbol{Lo\hspace{-0.03cm}g\hspace{-0.03cm}A\hspace{-0.03cm}R}\big(\boldsymbol{\mathcal Y\hspace{-0.1cm} \mathcal W}_{\Delta,\Delta'} \big)\oplus \big\langle
{\sf R}_{\Delta,\Delta'}
 \big\rangle $.}
\sk 
\item[{\bf 10.}]
{\bf ${\big[}$Differentiation of $ \boldsymbol{{\sf R}_{\Delta,\Delta'}}${\big]}.} 
{\it By differentiation,  $ {\sf R}_{\Delta,\Delta'}$ generates the whole space of polylogarithmic abelian relations of $\boldsymbol{\mathcal Y\hspace{-0.1cm} \mathcal W}_{\Delta,\Delta'}$.}
\sk 
\item[{\bf 11.}]
{\bf ${\big[}$Accessibility of $ \boldsymbol{{\sf R}_{\Delta,\Delta'}}${\big]}.} 
{\it The dilogarithmic identity $({\sf R}_{\Delta,\Delta'})$ is accessible from 
 $({\sf R}_{A_2})$. } 
\item[{\bf 12.}]
{\bf ${\big[}$Non linearizability{\big]}.}
{\it 
The web $\boldsymbol{\mathcal Y\hspace{-0.1cm} \mathcal W}_{\Delta,\Delta'}$ is not linearizable. }
\end{enumerate}

We do not have defined yet 
what the notations $LogAR_{\boldsymbol{i}}(t)$'s and 
$\mathfrak I^+_{\hspace{-0.04cm} \Delta ,  \Delta'  }$ of {\bf 8.} refer to, but these are just the straightforward generalizations  of those  
 of \S\ref{Par:LogAR-Y-method} to the bi-Dynkin cases.  Similarly,  
the ninth and tenth points just above are the direct generalization to pairs of Dynkin diagrams of the material discussed before in the case of a single diagram ({\it cf.}\,\S\ref{Par:Dilog-Generation-LogAR} regarding {\bf 10.} for instance).  
Notice that stating {\bf 9.}\,in a precise way would require beforehand to know what is the topology ( at least the fundamental group) of the complement of the set of common leaves 
of $\boldsymbol{\mathcal Y\hspace{-0.1cm} \mathcal W}_{\hspace{-0.04cm}\Delta,\Delta'}$ in the initial torus $\mathbf T^{nn'}$. According to the third of the {\bf {Statements }$\boldsymbol{\mathcal Y\hspace{-0.1cm} \mathcal W}_{\hspace{-0.04cm}\Delta,\Delta'}$}, it is conjectured that this complement  coincides with that of the cluster arrangement $\boldsymbol{Arr}_{\Delta,\Delta'}$. These arrangements of hypersurfaces appear as rather interesting (at least from this point of view) and would deserve to be studied more systematically.   \sk 

The accessibility of $ \big({{\sf R}_{\Delta,\Delta'}}\big)$ from the 5-terms relation is a natural questions stated at the end of the fifth section of \cite{KY}  (see also  \S\ref{Par:KimYamazaki} above). In this paper, the authors prove accessibility
 in finite type ({\it i.e.} when $\Delta=A_1$) and also in case $(A_3,A_3)$ which is treated as an example.\mk 
%

Finally, as for the  twelfth point above, it is a very natural question from the point of view of web-geometry to wonder about the linearizability of cluster webs 
$\boldsymbol{\mathcal Y\hspace{-0.1cm} \mathcal W}_{\hspace{-0.04cm}\Delta,\Delta'}$. 
We have proved that these webs are not linearizable when $\Delta=A_1$ (see Proposition \ref{Prop:YWD-non-linearizable}), it is very natural to conjecture that it is the case in full generality.   
Using the fact that a necessary condition for a web  to be linearizable is to 
be compatible with a projective connection ({\it cf.}\,\cite{PirioLin} for more details), we have verified this conjecture for many pairs $(\Delta,\Delta')$.


\newpage

%
%
%
%
%
%

\section{
Questions, problems and perspectives}
\label{S:Questions-Problems-Perspectives}
We believe that one of the main interests of the material and the results presented in the seven previous sections is the important number of questions and perspective they lead  to. \sk 

 Many questions and conjectures have already be stated before. We come back to some of them below.    But we also discuss  several new things. Some are precise questions rigorously stated, others just vague discussions about  general ideas which may (or may not) be fruitful. The unique general structure to be found in this section is that we have more or less tried to discuss different things according to their order of appearance in the preceding sections: first we discuss webs, then polylogarithms/iterated integrals and their functional equations, to eventually focus on questions related to cluster algebras.


\subsection{\bf In web geometry.}
\label{SS:InWebGeometry}

\subsubsection{\bf Characterisation of virtual ranks of AMP webs.}

A classical approach  to study the properties of a projective curve
$C\subset \mathbf P^n$  in algebraic geometry is to try to relate  its invariants to those of a general 
hyperplane section $\Gamma=C\cap H$.  In this regard, the property  for the curve to be ACM is particularly interesting since in this case  several of its numerical invariants  can  be recovered from those of $\Gamma$. This leads in particular to characterize the hyperplane section $\Gamma$ coming from an ACM curve $C$, for instance in terms of its Hilbert function $h_\Gamma$. 
The problem of characterizing and if possible of describing the Hilbert functions of general hyperplane sections of ACM curves is mentioned in  \cite[p.\,94-95]{HarrisEisenbud}  and has been touched in several papers. \sk

This problem about about algebraic curves has a straightforward generalization to webs: 
\mk 

{\bf Problem:}
{\it Characterize and describe the possible sequences of virtual ranks  
 of AMP webs.}\mk

Since AMP webs which are not algebraizable exist (the $\mathcal X$- or $\mathcal Y$-cluster webs in type $A_n$ for $n\geq 2$ for instance), this  problem is a strict generalization of the corresponding problem for projective curves and all the more so since it covers the case of non irreducible curves as well, a case left aside most of the time by algebraic geometers.\footnote{Algebraic geometers studying these questions often (if not always) assume that the general hyperplane section $\Gamma$ they are considering is formed by points satisfying what we call the `{\it Strong general position assumption}'.  As a rule, the corresponding condition is not satisfied by the cluster webs we are mainly dealing with in this paper.}
\sk

The problem above is easily answered for planar webs but we believe it is a very difficult one in full generality. Even the case of webs in three variables seems completely out of reach at the time of writing.  Even a partial answer would be interesting since it could give a way to discriminate non-linearizable AMP webs just by looking at their virtual ranks. The first questions in this direction which comes to our mind is the specific following one:

\begin{question} Given $n\geq 3$, does the sequence of virtual ranks $\rho^\bullet\big(\boldsymbol{\mathcal Y\hspace{-0.1cm} \mathcal W}_{\hspace{-0.04cm}A_n}\big)=\big( n(n+1)/2\, , \, n\, , \, 1\big)$ coincide with that of a projective curve
of degree $n(n+3)/2$ in  $ \mathbf P^n$?
\end{question}

We believe that the answer is negative for any $n\geq 3$. 
Knowing if it is indeed the  case when $n=3$ would  already be interesting.

\newpage



 \subsubsection{\bf Ranks of webs defined on configuration spaces}
 \label{SS:RanksWebsOnConfigurationSpaces}
 Given two positive integers $n$ and $m$ such that $m+n>2$, recall the web  in $mn$ variables
$\boldsymbol{\mathcal W}_{\hspace{-0.0cm}{\rm Conf}_{m+n+2}(\mathbf P^m)}$ introduced in \S\ref{Par:WebsOnConfugurationSpaces} which by definition is the 
 web on the moduli space of projective configurations ${\rm Conf}_{m+n+2}(\mathbf P^m)$
 defined by the  ${ m+n+2\choose n-1}{n+3 \choose 4}$ rational first integrals $\pi_{I,J}$  ({\it cf.}\,\eqref{Eq:pi-I-J}) for all pairs $(I,J)$ in $\mathfrak P(m,n)$. 
\mk 

For $m=1$, one recovers the web $\boldsymbol{\mathcal W}_{\hspace{-0.0cm}{\mathcal M}_{0,n+3}}$ studied in \cite{Pereira} where the author determines its virtual ranks, and shows that the total virtual rank coincides with the standard rank of this web which therefore is AMP. Considering this, a natural question is to compute the virtual ranks and the standard rank of $\boldsymbol{\mathcal W}_{\hspace{-0.0cm}{\rm Conf}_{m+n+2}(\mathbf P^m)}$ when $m\geq 2$. \sk 

More precisely, for $\sigma\geq 1$, let  $\ell^\sigma(m,n)$ stands for the dimension of the vector space spanned by the powers of order $\sigma$ of the differentials of the first integrals $\pi_{I,J}$'s at a generic point $\zeta\in \zeta\in {\rm Conf}_{m+n+2}(\mathbf P^m)$.  
The question which seems quite natural at this point is the following: 
\begin{question}
What is the 
value of  $\ell^\sigma(m,n)$ for $m,n,\sigma \in \mathbf N^*$ in general? Does it exist a closed formula for this quantity in terms of these three integers?
\end{question}

When $\sigma=1$, it is easy  to get that $\ell^1(m,n)=mn$ 
for every $m,n\geq 1$ but the computation of $\ell^\sigma(m,n)$ for $\sigma>1$ seems to be quite more difficult in general. This is a problem of 
  geometric and combinatorial  nature and we believe that being able to 
 give a formula for  $\ell^\sigma(m,n)$ would be meaningful regarding the (local or even global) geometry of the moduli space $ {\rm Conf}_{m+n+2}(\mathbf P^n)$ on which the web 
  $\boldsymbol{\mathcal W}_{\hspace{-0.0cm}{\rm Conf}_{m+n+2}(\mathbf P^m)}$ naturally lives.\mk 

A preliminary  exploration relying on brute-force computations leads us to conjecture that the following formulas hold true: 
\begin{enumerate}
\item[$\bullet$] Case $\sigma=2$ with $m,n\geq 1$ arbitrary:  
$$ \ell^2(m,n)=\frac{mn(mn+1)}{2}
 ; \hspace{10.5cm} {}^{} $$ 
\item[$\bullet$] Case $m=2$,  with $\sigma\geq 3$ and $n\geq 1$ arbitrary:   
$$ {}^{} 
\hspace{-0.3cm} 
\ell^3(2,n)=\frac{n(n+1)(7n+5)}{6}
\, , \hspace{0.2cm}
\ell^4(2,n)=\frac{n(n+1)(n+2) (11n+9)}{24}
\, ,  \hspace{0.2cm} \ell^\nu(2,n)=d_{2,n} \hspace{0.2cm}  \mbox{for } \nu\geq 5
$$ 
\item[$\bullet$] Case $m=3$, $\sigma=2,3$  and $m\geq 1$ arbitrary: 
\begin{equation}
\label{Eq:l3(m,n)}
 \ell^3(3,n)=\frac{n(23n^2+27n+10)}{6}\, . 
 \hspace{7.5cm} {}^{}
\end{equation}
\end{enumerate}

Although it would be interesting to know whether these formulas are indeed correct,  these are not satisfying in the form above since at the exception of the formula for $\ell^2(m,n)$, 
the symmetry between $m$ and $n$\footnote{Since both webs $\boldsymbol{\mathcal W}_{\hspace{-0.0cm}{\rm Conf}_{m+n+2}(\mathbf P^m)}$ and $\boldsymbol{\mathcal W}_{\hspace{-0.0cm}{\rm Conf}_{m+n+2}(\mathbf P^n)}$ are equivalent in a natural way (thanks to projective duality), one gets immediatly that 
$\ell^\sigma(m,n)=\ell^\sigma(n,m)$ for any $m,n,\sigma\geq 1$.}
is not apparent in any of them.   It would be more satisfying to have a disposal (even conjectural) formula for $\ell^\sigma(m,n)$ which be symmetric in $m$ and $n$.   \mk 

Even more interesting but less clear at this point, is the possibility to have geometric or/and combinatorial interpretations for the quantities $\ell^\sigma(m,n)$'s. For instance,  
formula \eqref{Eq:l3(m,n)} for $\ell^3(3,n)$ has been conjectured after  considering the first corresponding values for $n=1,2,3,4$, which are $10, 52, 149, 324$ respectively.  This string  of  four integers matches only one of the sequences of the OEIS, namely the sequence \href{https://oeis.org/A092966}{A092966}.  The $n$-th term of this sequence has a geometric interpretation (as the	number of interior balls in a truncated tetrahedral arrangement) as well as a combinatorial one (as the number 4-element subsets of $\{-n,\ldots ,0, \ldots, n\}$ having sum $n+1$).  It would be interesting to relate one or both of these interpretations for $A092966(n)$ to the geometric one corresponding to the very definition of $\ell^3(3,n)$. 
\begin{center}
$\star$
\end{center}

Regarding web geometry, the determinations of the polylogarithmic ranks (with ramification locus  $\{ 0, 1, \infty\}$) and of the actual (standard) ranks of the webs $\boldsymbol{\mathcal W}_{\hspace{-0.0cm}{\rm Conf}_{n+m+2}(\mathbf P^m)}$ is certainly as interesting as the one of their virtual ranks discussed above. \sk 

In \cite{Pereira}, Pereira has proved  that 
the  webs $\boldsymbol{\mathcal W}_{\hspace{-0.0cm}{\rm Conf}_{n+3}(\mathbf P^1)}$ (case $m=1$ with $n\geq 2$ arbitrary)  all are AMP, by computing explicitly their polylogarithmic ranks and total virtual ranks  and comparing theses quantities. This suggests immediately the following series of questions: 
\begin{questions}  Assume that both $m$ and $n$ are bigger than 1. 
\begin{enumerate}
\item  For $w\geq 3$ fixed, are there pairs $(m,n)$ such that $\boldsymbol{\mathcal W}_{\hspace{-0.0cm}{\rm Conf}_{m+n+2}(\mathbf P^m)}$ carries polylogarithmic ARs of weight $w$? If yes, determine the set of such pairs and the corresponding ARs. 
\item  Are there pairs $(m,n)$ such that $\boldsymbol{\mathcal W}_{\hspace{-0.0cm}{\rm Conf}_{m+n+2}(\mathbf P^m)}$ carries non polylogarithmic ARs? If yes, determine the set of such pairs and the ARs of this kind. 
\item Compute  the polylogarithmic ranks and the rank  of  
$\boldsymbol{\mathcal W}_{\hspace{-0.0cm}{\rm Conf}_{m+n+2}(\mathbf P^m)}$.
\item  Give a basis of the space of ARs of $\boldsymbol{\mathcal W}_{\hspace{-0.0cm}{\rm Conf}_{m+n+2}(\mathbf P^m)}$. 
\item  Determine the pairs $(m,n)$ such that  $\boldsymbol{\mathcal W}_{\hspace{-0.0cm}{\rm Conf}_{m+n+2}(\mathbf P^m)}$ is AMP.
\end{enumerate}
\end{questions}

In \S\ref{Par:W-Conf6(IP2)-is-AMP}, we have seen that if $\boldsymbol{\mathcal W}_{\hspace{-0.0cm}{\rm Conf}_{6}(\mathbf P^2)}$ is AMP  with only polylogarithmic ARs  (of weight 1 or 2), this is not the case for  
$\boldsymbol{\mathcal W}_{\hspace{-0.0cm}{\rm Conf}_{7}(\mathbf P^2)}$ which is even not AMP.  
This suggests that there is a priori no uniform answer to {\it 5.} and that 
even conjecturing what might be the set 
of  pairs $(m,n)$ such that  $\boldsymbol{\mathcal W}_{\hspace{-0.0cm}{\rm Conf}_{m+n+2}(\mathbf P^m)}$ is AMP  is non trivial.



 \subsubsection{\bf Webs on strata of degenerate configurations.}
The fact that  several interesting webs (such as Spence-Kummer's  9-web 
$\boldsymbol{\mathcal W}_{\hspace{0.0cm} \mathcal S\mathcal K}$
or 
Goncharov's 22-web $\boldsymbol{\mathcal W}_{\hspace{0.0cm} {\mathcal G}_{22}}$, which both are trilogarithmic) are obtained as restrictions of some {\it `configurational webs'} on some strata,  is a strong motivation to undertake  more systematically the study
 of the webs obtained in this way. \mk 
 
Via direct computations, we have already  looked at the case of seven points in the projective plane and have found several other interesting webs in addition to 
Spence-Kummer's 
$\boldsymbol{\mathcal W}_{\hspace{0.0cm} \mathcal S\mathcal K}$
 and Goncharov one $\boldsymbol{\mathcal W}_{\hspace{0.0cm} {\mathcal G}_{22}}$.  Below, we discuss the cases of  three examples associated to 
 the degenerate configurations of seven points in the projective plane pictured just below: \sk
 
\begin{figure}[h]
\begin{center}
\begin{tabular}{lcr}
\hspace{-2cm} 
\resizebox{1.1in}{0.9in}{
 \includegraphics{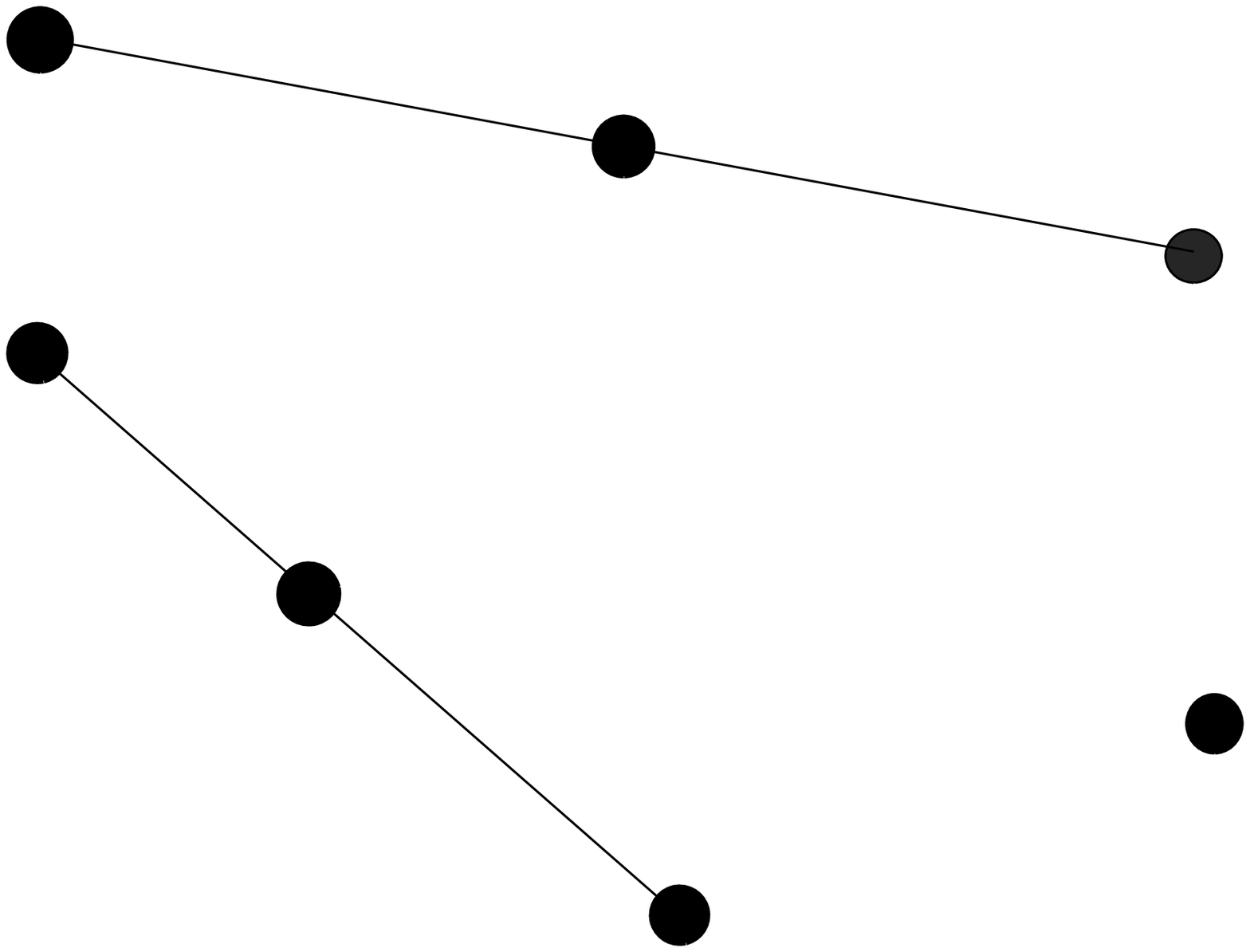}}
\hspace{1cm}  &   
 \resizebox{1.1in}{0.9in}{
 \includegraphics{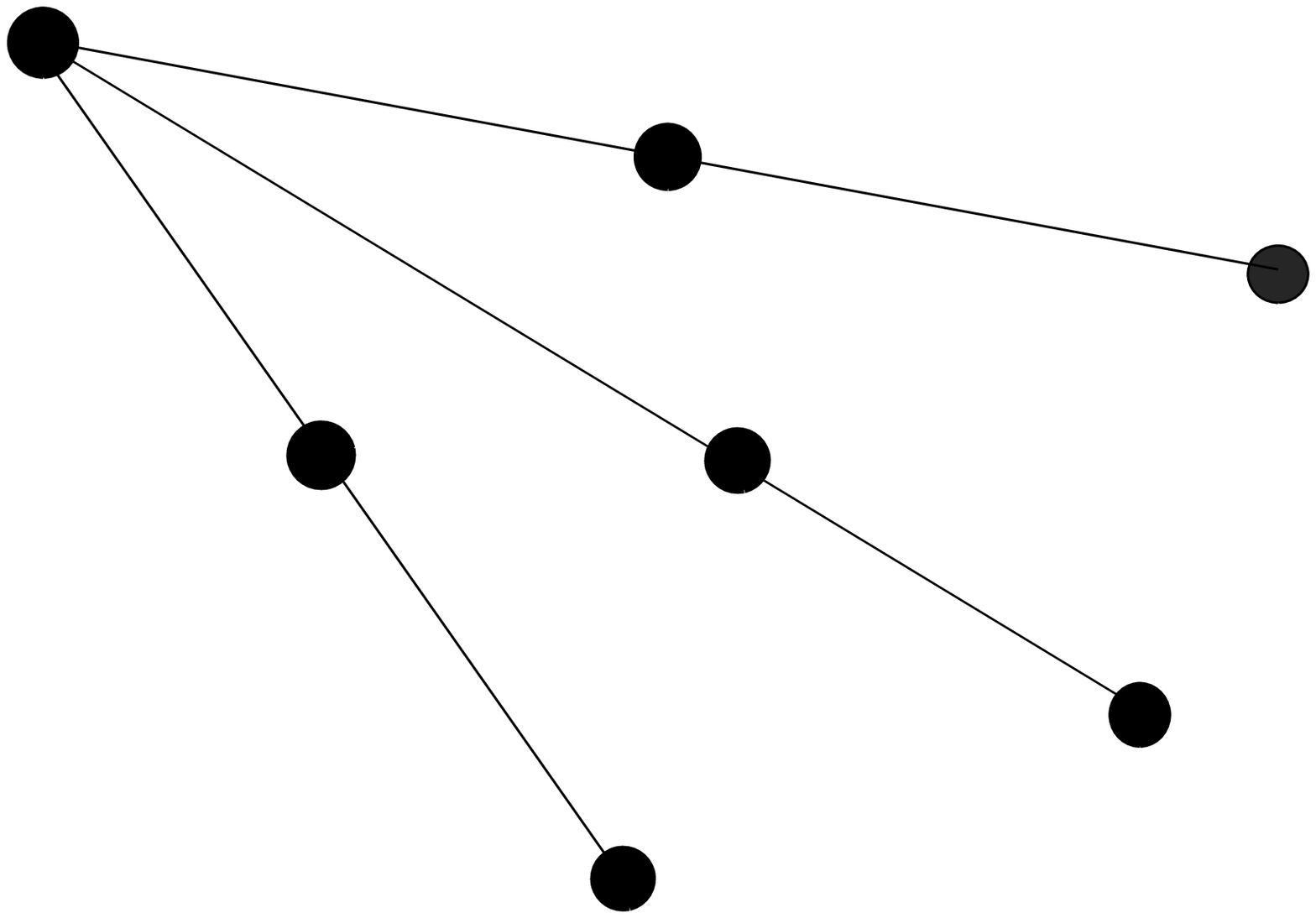}}
 & \hspace{1cm}    \resizebox{1.1in}{0.9in}{
 \includegraphics{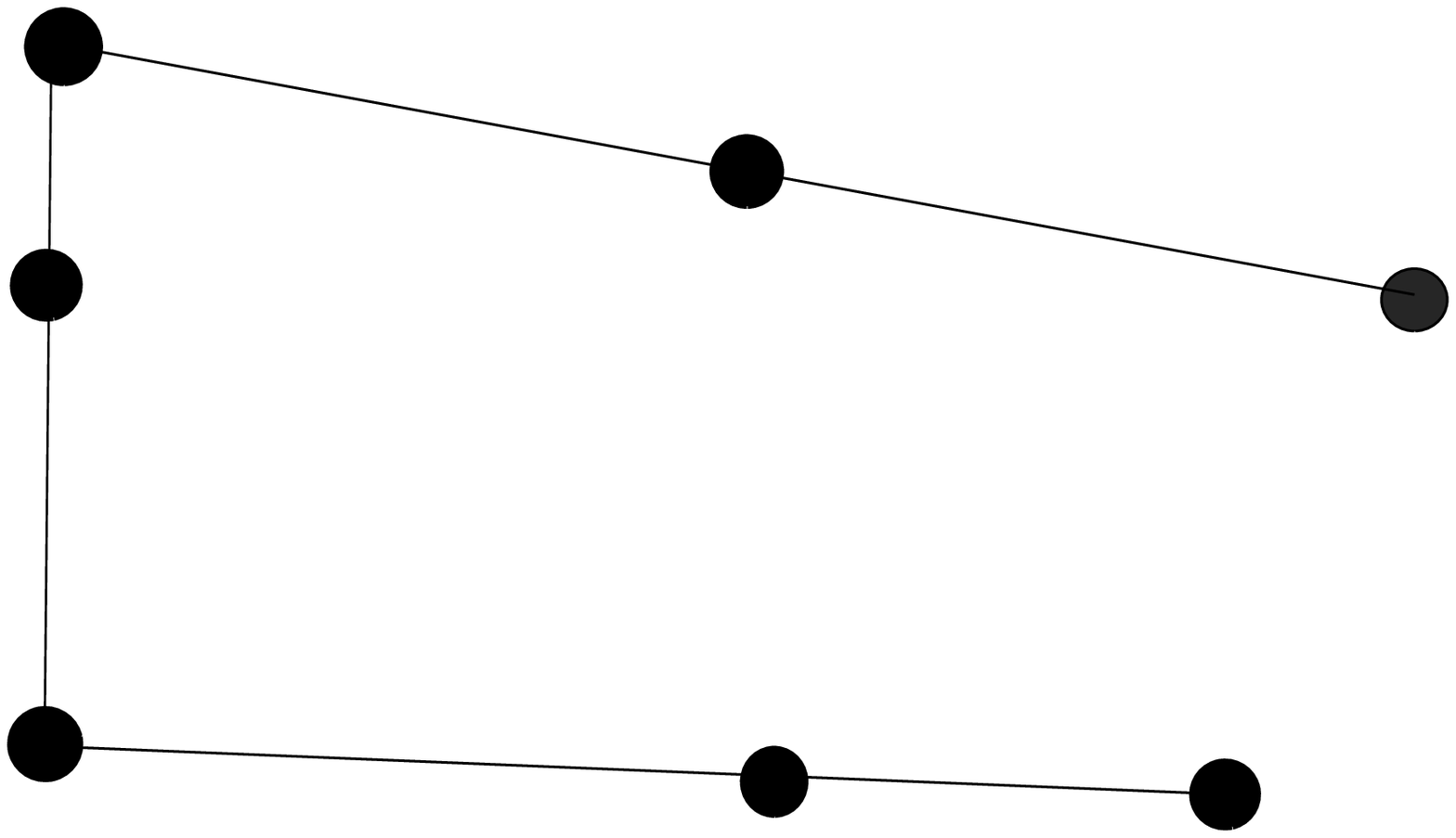}}
 \end{tabular}
\caption{three types of degenerate configurations of seven points on $\mathbf P^2$.}
\label{Fig:Strata}
\vspace{-0.5cm}
\end{center}
\end{figure}
 \begin{itemize}
 \item Let $\boldsymbol{S}_1$ be the 4-dimensional substratum of ${\rm Conf}_7(\mathbf P^2)$ formed by configurations of seven points similar to the one in the picture on the left of Figure \ref{Fig:Strata}.  
 Then  $\boldsymbol{\mathcal W}_{\boldsymbol{S}_1}=\boldsymbol{\mathcal W}_{{\rm Conf}_7(\mathbf P^2)} \lvert_{{\boldsymbol{S}_1}}$ is a 35-web in four variables, with $\rho^\bullet\big(\boldsymbol{\mathcal W}_{\boldsymbol{S}_1}\big)=(31,25,15)$ and ${\rm polrk}^\bullet\big(\boldsymbol{\mathcal W}_{\boldsymbol{S}_1}\big)=(56,16)$, and consequently, is AMP 
 with only polylogarithmic ARs (of weight at most 2). In this regard, it is similar to the cluster dilogarithmic webs associated to $Y$-systems.  
 \sk 
 \item Let $\boldsymbol{S}_2$ be the stratum of ${\rm Conf}_7(\mathbf P^2)$ formed by configurations of seven points similar to the one in the picture of the middle of  Figure \ref{Fig:Strata}. It has dimension 3 and one verifies that $\boldsymbol{\mathcal W}_{\boldsymbol{S}_2}=\boldsymbol{\mathcal W}_{{\rm Conf}_7(\mathbf P^2)} \lvert_{{\boldsymbol{S}_2}}$ is a 18-web in three variables, with $\rho^\bullet\big(\boldsymbol{\mathcal W}_{\boldsymbol{S}_2}\big)=(15,12,8,3)$  hence with total virtual rank equal to 38.  
Considering only  harmonic hyperlogarithmic ARs\footnote{That is, ARs whose components are iterated integrals with ramifications at the  four points
$0$, $-1$, $1$  and $\infty$ (these are four points on $\mathbf P^1$ in harmonic division, hence the name `{\it harmonic hyperlogarithms}').}  for this web, we have constructed 26 (resp.\,9) linearly independent such ARs of weight 1 (resp.\,of weight 2). Thus ${\rm rk}(\boldsymbol{\mathcal W}_{\boldsymbol{S}_2})\geq 35$. On the other hand, using the criterion of \S\ref{SubPar:CharacterizationWebsMaximalRank}, we have obtained that 
${\rm rk}(\boldsymbol{\mathcal W}_{\boldsymbol{S}_2})=
\rho(\boldsymbol{\mathcal W}_{\boldsymbol{S}_2})=
38$ hence this web is AMP. It would be interesting to know more about the 3 missing ARs which  may not be of hyperlogarithmic type (however, see our Conjecture \ref{Conj:Nature-Fi} above regarding them).
  \sk 
 \item Finally, let $\boldsymbol{S}_3$ be the stratum of dimension 3 
formed by configurations of the same type than the one in the picture on the right of Figure \ref{Fig:Strata}. 
 Then  $\boldsymbol{\mathcal W}_{\boldsymbol{S}_3}=\boldsymbol{\mathcal W}_{{\rm Conf}_7(\mathbf P^2)} \lvert_{{\boldsymbol{S}_3}}$ is a 15-web in 3  variables, with $\rho^\bullet=(12,9,5)$ and ${\rm polrk}^\bullet=(21,5)$.  Consequently, it is AMP 
 with only polylogarithmic ARs (of weight 2 at most).  Since it shares all these numerical properties with $\boldsymbol{\mathcal W}_{{\rm Conf}_6(\mathbf P^1)}
 =\boldsymbol{\mathcal W}_{\hspace{-0.1cm}{\mathcal M}_{0,6}}
 $, one can wonder whether  these two webs are equivalent or not, which is not clear a priori. 
\end{itemize}

\subsubsection{\bf Resonance webs associated to arrangements.}

One of the main features of the notions of virtual ranks and of being AMP is that 
they are flexible concepts which adapt to many different kinds of webs, but also robust, since they remain relevant in most situations.\mk 

In \cite{Pereira}, the author explains how to associate a web defined by rational maps to any arrangement by hyperplanes $\mathcal A\subset \mathbf P^n$ with complement $U_{\mathcal A}$. If    $\Sigma\subset H^1(M_{\mathcal{A}},\mathbf Z)$
is a subspace of dimension $k\geq 2$ which is totally isotropic for  the cup-product 
$H^1(U_{\mathcal{A}},\mathbf Z)\wedge H^1(U_{\mathcal{A}},\mathbf Z)\rightarrow H^2(U_{\mathcal{A}},\mathbf Z)$ and maximal for the inclusion among isotropic subspaces,  there exists a rational map $f_\Sigma: \mathbf P^n\dashrightarrow \mathbf P^1$ and a subset $\sigma_\Sigma\subset \mathbf P^1$ of cardinality $k+1$ such that
$f_\Sigma^{-1}(\sigma_\Sigma)\subset \mathcal A$ and 
 $\Sigma=(f_\Sigma)^*\big(H^1(\mathbf P^1\setminus \sigma_\Sigma, \mathbf Z )\big) $.  
Given $\Sigma$, the map $f_\Sigma$ is unique (up to post-composition by projective automorphisms). The number of such maximal isotropic subspaces being finite, one defines the {\bf resonance web $\boldsymbol{\mathcal W}(\mathcal A)$ associated to $\mathcal A$} as the web defined by these $f_\Sigma$'s: 
$$
\boldsymbol{\mathcal W}(\mathcal A)=
{\boldsymbol{\mathcal W}}\left(\,  f_\Sigma : \mathbf P^n\dashrightarrow \mathbf P^1 \hspace{0.1cm} 
 \Big\lvert 
 \hspace{0.00cm}  
 \begin{tabular}{l}
 {\it  $\Sigma\subset H^1(M_{\mathcal{A}},\mathbf Z)$ of dimension $\geq 2$,}
 \vspace{-0.07cm}
 \\
 {\it  totally $\wedge$-isotropic and maximal}
\end{tabular}
\right)\, .
$$

Some classical questions about arrangements concern the relations between the combinatorial properties of $\mathcal A$ and the algebraic and geometric properties of its complement $U_{\mathcal A}$. It is then natural to ask the following general questions when interested in web geometry: 
\begin{question} Let  $\mathcal A$ be an arrangement of hyperplanes in $\mathbf P^n$. 
\begin{enumerate}
\item  How much can be told about the web-theoretic invariants of 
$\boldsymbol{\mathcal W}(\mathcal A)$ (virtual, polylogarithmic, standard rank(s), (polylogarithmic) ARs, etc.) from the 
combinatorics of  $\mathcal A$ ?
\item Are there criterions on $\mathcal A$ ensuring that $\boldsymbol{\mathcal W}(\mathcal A)$ is AMP and/or with only polylogarithmic abelian relations?
\item Give examples of arrangements $\mathcal A$ such that  $\boldsymbol{\mathcal W}(\mathcal A)$ be AMP. 
\end{enumerate}
\end{question}

Of course, one does not expect nice or interesting answers to the previous question for any arrangement.  An exception may concern the polylogarithmic ranks of a resonance web $\boldsymbol{\mathcal W}(\mathcal A)$: we believe that it only depends on the combinatorics of $\mathcal A$.
\mk 


In \cite{Pereira},  Pereira first states several general results about the ARs of any resonance  web then he applies them to the web associated to 
an arrangement  ${\sf A}_n={\sf A}_{0,n+3}$  in $\mathbf P^n$ (see \S\ref{Par:BraidArr-2} above) related to the braid arrangement $B_{n+2}$ in $\mathbf P^{n+2}$ (for any $n\geq 2$).\footnote{See \cite[\S4]{Pereira} for some detaills about how   both arrangements ${\sf A}_n$ and $B_{n+2}$
are related.}
After identifying the web 
associated to ${\sf A}_n$ 
with 
$\boldsymbol{\mathcal W}_{\hspace{-0.05cm}\mathcal M_{0,n+3}}$, he computes its virtual ranks then proves that its rank is AMP ({\it cf.}\,Example \ref{Exm:Webs}.(4) above for precise statements).  In view to find other examples of webs enjoying this property, one can remark that the braid arrangements form just one class of an interesting class of arrangements, the Coxeter arrangements. \sk 

A \href{https://en.wikipedia.org/wiki/Coxeter_group}{Coxeter  group} $G$ is a finite subgroup of $GL_n(\mathbf R)$ generated by reflections and the associated {\bf Coxeter arrangement $\boldsymbol{\mathcal A(G)}$} is the complexification of the real projective arrangement in $\mathbf P(\mathbf R^n)$ whose affine cone is the union of 
 the reflecting hyperplanes in $\mathbf R^n$ of the reflections elements of $G$  ({\it cf.}\,\S6.2 in \cite{OrlikTerao}).\footnote{This terminology is not standard. In \cite{OrlikTerao}, 
 ${\mathcal A(G)}$ stands for the real arrangement in $\mathbf R^n$ formed by the reflecting hyperplanes of the reflections of $G$.} The (finite) Coxeter groups were classified by Coxeter in terms of Coxeter-Dynkin diagrams: there are three countable families of increasing rank corresponding to the Dynkin types $A,B$ and $D$,  
  the family of dihedral groups $I(p)$ with $p\geq 2$ (in dimension 2) and six exceptional examples, four of Lie type $E_6,E_7,E_8$, $F_4$, and two others, denoted by $H_3$ and $H_4$.

\begin{question} 
Let $G\subset GL_n(\mathbf R)$ be a finite Coxeter group. Assume that  $n\geq 2$. 
\begin{enumerate}
\item Describe the irreducible components $\Sigma$ of the first resonance variety $R^1(\mathcal A(G))$ and the corresponding rational maps $f_\Sigma$. In particular, give a formula for the degree of $\boldsymbol{\mathcal W}(\mathcal A(G))$. 
\item Does the set of common leaves of $\boldsymbol{\mathcal W}(\mathcal A(G))$ coincide with $\mathcal A(G)$? 
\item What is/are the   (standard, virtual, hyperlogarithmic) rank(s) of $\boldsymbol{\mathcal W}(\mathcal A(G))$?
\item  Give a basis of the space  of abelian relations $\mathcal A\big( \boldsymbol{\mathcal W}(\mathcal A(G))\big)$.  Is any AR of this web polylogarithmic?
\item Describe ({\it i.e.}\,monodromy, decomposition in irreducible factors, construction of these by extension, etc) the local system formed by the ARs of $\boldsymbol{\mathcal W}(\mathcal A(G))$ on the complement of its set of common leaves.  
\end{enumerate}
\end{question}





  All the Coxeter arrangements can be made explicit. For instance, here are defining equations in the cases $I(p)$ and Dynkin types $A$, $B$ and $D$ (see \cite{Coxeter}):
   \begin{align*}
 I(p)\, :&\qquad 0=  \prod_{0\leq i \leq p-1}\big(x-\omega^i y\big)=x^p-y^p \quad \mbox{\Big(with } \omega=\exp\big(2i\pi/p\big)\mbox{\Big)}\\
 A_n\, :&\qquad 0= \prod_{1\leq i <j\leq n}(x_i-x_j)\\
 B_n\, :&\qquad 0= \prod_{1\leq i <j\leq n}(x_i-x_j)(x_i+x_j) \prod_{1\leq k \leq n}x_k\\
 D_n\, :&\qquad 0=\prod_{1\leq i <j\leq n}(x_i-x_j)(x_i+x_j)\, .
 \end{align*}

 For the others types, explicit defining equations  are also available in the literature.\footnote{For instance, see  
 \cite[Table 3]{Iwasaki}  for the cases $H_3$, $H_4$, $F_4$ and 
 \cite{Metha} for the three exceptional types.} 
  There are effective methods to compute the resonance variety of any explicitly given arrangement, thus $R^1(\mathcal A(G))$ can be computed in all cases and consequently, the resonance web $\boldsymbol{\mathcal W}(\mathcal A(G))$ can be made explicit for any $G$.\footnote{Here we implicitly  assume that
 $G$ acts on $\mathbf R^n$ with $n\geq 3$, otherwise 
  the   web $\boldsymbol{\mathcal W}(\mathcal A(G))$ would not be defined. In particular,  there is no resonance web associated to the family of dihedral groups $I(p)$.}
Together with the algorithmic methods to determine basic invariants of webs considered in \S\ref{Para:DeterminingARs-and-Rank},  this opens a way to tackle the series of questions above by brute force computations. \sk 

In the recent paper \cite[\S6]{CohenSchenke}, the authors give descriptions of the resonance varieties of Coxeter arrangements for Dynkin types $B$ and $D$. It follows that 
the associated resonance web has degree $16 { n \choose 3} + 9 { n \choose 3}+{ n \choose 2}$  and $5{n \choose 3}+9 {n \choose 4}$  in types $B_n$ and $D_n$ respectively. It should not be very difficult to make these webs explicit and study them, at least  for $n$ small enough.  Comparing the degrees of these webs with those of the corresponding $\mathcal X$-cluster webs, it appears that 
the  Coxeter resonance webs and the $\mathcal X$-cluster webs do not coincide 
for Dynkin types $B$ and $D$, in contrast with what does occur in type $A$. 
\mk 

Given a Coxeter group $G$,  it acts as birational symmetries of $\boldsymbol{\mathcal W}(\mathcal A(G))$. Thus resonance webs associated to Coxeter arrangements admit many symmetries and one may think that this is the reason why they 
may carry many (poly- or hyperlogarithmic) ARs\footnote{This is what we believe. For now, this  has been proven to be true only in type $A$.}. But one should not think that the interest of the notion of `resonance web' introduced by Pereira only concerns very symmetrical arrangements. A striking example (with many symmetries though) is 
given by the `non-Fano arrangement' \cite[Example 5.4]{Pereira}: it is the arrangement ${\sf NF}$ of 7 lines in $\mathbf P^2$ pictured below:  
\begin{figure}[h]
\begin{center}
\resizebox{1.6in}{1.6in}{
 \includegraphics{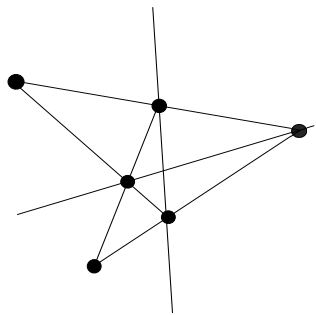}}
 \vspace{-0.05cm}
\caption{The non-Fano arrangement (see also Figure \ref{Fig:Spence-Kummer configuration}
above).} \label{Fig:Non-Fano}
\vspace{-0.6cm}
\end{center}
\end{figure}

(Observe that ${\sf NF}$ is precisely the set of common leaves of the web obtained by taking the restriction of  $\boldsymbol{\mathcal W}_{\hspace{-0.05cm} {\rm Conf}_7(\mathbf P^2)}$ to the stratum $\boldsymbol{S}({\hspace{-0.0cm}{\cal S}{\cal K}})$, see \S\ref{SubPar:GeomSK-Web}).\sk 

The resonance variety $R^1(${\sf NF}$)$ has been computed in \cite[Example 4.4]{CohenSuciu}: it has nine irreducible components, all of dimension 2. The pencils of curves corresponding to the associated rational fibrations  $ \mathbf P^2\setminus {\sf NF}\rightarrow \mathbf P^1\setminus \{Ê 0,1,\infty\}$ are easily described relatively to the six points represented by the black dots in Figure  \ref{Fig:Non-Fano}: there are six pencils of lines, formed by the line passing through one of the six points, and three pencils of conics, one for each 
subset of four points  in general position among the six. We deduce that the resonance web 
$\boldsymbol{\mathcal W}({\sf NF})$ coincides with Spence-Kummer web. 

\begin{question}
 In addition to Abel's 5-terms relations and Spence-Kummer identity, is there other classical polylogarithmic AFE whose associated web  is equivalent to the resonance 
 web of an arrangement? More specifically, what about the case of Kummer's identity $\boldsymbol{\mathcal K}_4$? Is the tetralogarithmic web $\boldsymbol{\mathcal W}_{\hspace{-0.05cm}\boldsymbol{\mathcal K}_4}$ equivalent to the resonance web of an arrangement of lines in $\mathbf P^2$?
\end{question}
\begin{center}
$\star$
\end{center}
To end this section, we mention that most of the material above generalises to complements of arrangements of hypersurfaces on a smooth (projective, or even just compact) manifold but 
the situation in general is   richer than in the case of arrangement by hyperplanes. 
For more details, we refer to  \cite{Dimca} and \cite{DPS} 
and the references therein. \sk

Say that  $X$ is a smooth projective variety  and let $A$ be the union of a finite number of irreducible hypersurfaces with complement $U=X\setminus A$. Then in addition to the first resonance variety $R^1(U)$ there is the (first) characteristic variety $V^1(U)$ which (as a set) is formed by the characters $\rho\in {\rm Hom}(\pi_1(U),\mathbf C^*)=H^1(U,\mathbf C^*)$ such that $H^1(U,\mathbf C_\rho)>0$, where $\mathbf C_\rho$ stands for the local system on $U$ with monodromy $\rho$.  One can thus attach two webs to $U$: its `resonance web 
$\boldsymbol{\mathcal W}_{R^1(U)}$' with one foliation defined by a holomorphic first integral for each irreducible component of $R^1(U)$ of dimension greater than 1 and its `{\bf characteristic web 
$\boldsymbol{\mathcal W}_{V^1(U)}$}', which is defined by holomorphic fibrations  as well. In general, the latter is a subweb of the former and both coincide under the condition that $U$ is 1-formal, a condition automatically satisfied when $U$ is the complement of an arrangement by hyperplanes, as proved by Brieskorn.\sk 

It would be interesting to know whether  
there exist  quasi-projective varieties $U$ as above such that 
one of the two webs $\boldsymbol{\mathcal W}_{R^1(U)}$ or $\boldsymbol{\mathcal W}_{V^1(U)}$ is AMP but is not  (equivalent to) a resonance web associated to a hyperplane arrangement.  We confess not having any  idea regarding this.
\mk 

Another interesting question, but to which specialists should be able to answer easily, is whether  classical algebraic webs associated to a projective curve are related to resonance or characteristic webs. More precisely, let $C\subset \mathbf P^n$ be a smooth (to simplify) projective curve  defining an algebraic web $\boldsymbol{\mathcal W}_C$ on $\check{\mathbf P}^n$.  The singular locus $\Sigma(\boldsymbol{\mathcal W}_C)$ of this web is a hypersurface containing the dual variety $C^\vee$ of $C$ (see \cite[\S2.2]{Nakai}). Is $\boldsymbol{\mathcal W}_C$ 
 related to a resonance or characteristic web constructed from the complement of $C^\vee$ or of $\Sigma(\boldsymbol{\mathcal W}_C)$ in $\check{\mathbf P}^n$? 




\subsection{\bf About polylogarithms and the functional equations they verify.}
\label{SS:AboutPolylog} 
Here we first discuss a question 
which came to our mind when elaborating the material of Section \S\ref{S:Polylogs-FEs}. 
Then the possibility that some famous polylogarithmic identities be of cluster type is discussed. Finally, we consider functional equations of multivariable polylogarithms and show that many aspects of the theory developed in the first two sections of this text 
can be generalized to webs of codimension higher than 1.


\subsubsection{\bf Iterated integrals, Schur functors and monodromy.}
\label{SS:II-Schur-Monodromy}
When looking for webs carrying many abelian relations (such as AMP webs), the ones with hyperlogarithmic  components are particularly interesting,  for two distinct reasons.
The first is that hyperlogarithms being multivalued, one can construct many other ARs of the same kind by analytic continuation along loops contained outside the singular set of the considered web. The second is that for any weight $w\geq 1$, the symmetric group $\mathfrak S_w$  acts naturally on the space of weight $w$  hyperlogarithmic 
ARs of a given web ({\it cf.}\,Lemma \ref{L:KWw-Sw-equivariant}) which can be used to get 
several new hyperlogarithmic ARs  from a given one (an illustration of this is given by  Theorem  \ref{T:FE-holomorphic}, especially its second part).
\sk 

The theory of representations of symmetric group actions on tensor products is by now a well established theory (notions of young symmetrizers, Spetch modules), as well as the study of the monodromy of single-variable hyperlogarithms. The problem of better understanding  the hyperlogarithmic ARs of webs leads to wonder how these two things are related.  Before considering this multivariable situation, a preliminary but already interesting case is the one-dimensional case about which we say a few words below. 
  \mk

Let 
$S=\{s_1,\ldots,s_n\}$ be a set of $n$ pairwise distinct points in $\mathbf C$ with $n\geq 2$  (the classical polylogarithmic case corresponding to $n=2$).  We set 
$S^\star=\{s_1,\ldots,s_n,\infty\}\subset \mathbf P^1$ and choose 
$\zeta\in \mathbf C\setminus S$, an arbitrary fixed base-point (what follows is independent of the choice of $\zeta$).  To the affine curve $\mathbf P^1\setminus S^\star$ is attached the vector space of $S$-logarithmic 1-forms 
$V= H^0\big(\mathbf P^1, \Omega_{\mathbf P^1}(\log S^\star)\big)$, 
 a basis of which is given by the set of  $d\log(z-s_i)=dz/(z-s_i)$'s for $i=1,\ldots,n$. 
 \sk 
 
 For any weight $w\geq 1$ and any partition $\lambda$ of $w$, let 
 $V_{\lambda}$ be the image of $V$ by the Schur functor $\boldsymbol{S}_{\hspace{-0.04cm}\lambda}$. \footnote{More concretely, for each partition $\lambda\vdash w$, $V_\lambda$
is the image in $V^{\otimes w}$ of  the \href{https://en.wikipedia.org/wiki/Young_symmetrizer}
{Young symmetrizer} $c_\lambda\in {\rm End}\big( 
V^{\otimes w} \big)$.} 
Since each $V_\lambda$  is an irreducible $\mathfrak S_w$-representation, 
we deduce that 
 \begin{equation}
 \label{Eq:Vw}
V^{\otimes w}=\bigoplus_{\lambda \vdash w}  V_{\lambda}
\end{equation}
 is the decomposition of the $w$-th tensor power of $V$ 
in irreducible $\mathfrak S_{ w}$-modules (it can be called the Schur or the Weyl decomposition of $V^{\otimes w}$ ({\it cf.}\,p.\,76 in \cite{FultonHarris}).  For any partition $\lambda$, we denote by 
$c_\lambda$ (resp.\,$\iota_\lambda$) the projection $V^{\otimes w}
\hspace{-0.05cm}
\twoheadrightarrow V_{\lambda}$  (resp.\,the injection $V_{\lambda}\hookrightarrow  V^{\otimes w}  $) induced by this decomposition.
\sk 

\vspace{-0.3cm}
We now use the terminology and the notations introduced in \S\ref{Par:II-AR} concerning iterated integrals. To a weight $w$ symbol  $\boldsymbol{\omega}
 \in V^{\otimes w}$, is associated the hyperlogarithm $L_{\boldsymbol{\omega},\zeta}={\rm II}_{S,\zeta}^w(\boldsymbol{\omega})\in \mathcal H_\zeta^w$, which is a well-defined germ of holomorphic function at $\zeta$. If $\mathcal C^\gamma$ stands for the analytic continuation along 
a loop $\gamma$ in  $\mathbf C\setminus S$, centered at $\zeta$, we have 
that $\mathcal C^\gamma\big(L_{\boldsymbol{\omega},\zeta}\big)
-L_{\boldsymbol{\omega},\zeta}
\in \oplus_{w'<w} \mathcal H^{w'}_\zeta$ (see \eqref{Eq:C-gamma}). 
Then considering the image of this by the symbol $\boldsymbol{\mathcal S}=({\rm II}_{S,\zeta}^w)^{-1}: \mathcal H^\bullet_\zeta\rightarrow V^{\otimes \bullet}$, one gets  a  $\mathbf C$-linear map 
$$
\rho^w(\gamma)=\rho_{S,\zeta}^w(\gamma) \, :\, V^{\otimes w} \longrightarrow \oplus_{\tilde w<w} V^{\otimes \tilde w}\, , \quad  \boldsymbol{\omega} \longmapsto \boldsymbol{\mathcal S}\circ (\mathcal C^\gamma- {\rm Id}\big)\circ {\rm II}_{S,\zeta}^w(\boldsymbol{\omega})\,. 
$$
Since $\rho(\gamma)$ only  depends on the homotopy class of $\gamma$ as a $\zeta$-based loop  in $\mathbf C\setminus S$, one obtains a map 
$$\rho^w=\rho^w_{S,\zeta} \, : \, \pi_1\big(\mathbf C\setminus S, \zeta\big)\longrightarrow 
{\rm End}_{\mathbf C}\Big(V^{\otimes w},\oplus_{\tilde w<w } V^{\otimes \tilde w}\Big)$$
which encodes the monodromy of  $S$-hyperlogarithms of weight $w$ 
 and has been described by several authors (for modern references, see 
 \cite[Thm.\,3]{HPV} in the polylogarithmic case ($n=2$) and the sixth section of \cite{Brown1} for the general case $n\geq 2$). 
\sk 

The question we are interested in is that of describing the decomposition of $\rho$ relatively to the Schur decompositions \eqref{Eq:Vw} of  $V^{\otimes k}$ at the source, 
and the decomposition of the many 
 blocks $V^{\otimes \tilde w}$ at the target. 
In order to be more precise,  we fix some loops $\gamma_i$ in  $\mathbf C\setminus S$, centered at $\zeta$, with index 1 with respect to $s_i$ and index 0 with respect to any $s_j$ with $j\neq i$, such that $\pi_1(\mathbf C\setminus S, \zeta)$ is the free group generated by the $\zeta$-based homotopy classes $[\gamma_1],\ldots,[\gamma_n]$.

\begin{question}  For any $i\in \{1,\ldots,n\}$, describe the decomposition of $\rho$ relatively to the Schur decompositions \eqref{Eq:Vw} of  $V^{\otimes w}$ at the source, and of the many blocks $V^{\otimes \tilde w}$ at the target. \sk 

 
More precisely,  for  $ \tilde w<w$ and any partitions $\lambda\vdash w$ and $\tilde \lambda \vdash \tilde w$, describe  
the $\mathbf C$-linear map 
  $$
 \rho^{w,\tilde w}_{\lambda, \tilde \lambda}(i)=c_{\tilde \lambda}\circ \rho^w(\gamma_i) \circ \iota_{\lambda}\,  : \,  V_\lambda\longrightarrow 
 V_{\tilde \lambda}\, .
 $$
\end{question}


As far we know, this question has not be considered before despite its natural and elementary character. 

 \subsubsection{\bf Cluster nature of classical polylogarithmic identities.}
 \label{SS:ClusterNature}
 The fact that there are links between the theory of  dilogarithmic identities  and that   of cluster algebras  is by now well-known. But the fact that some (classical or not) AFEs satisfied by polylogarithms of higher weight (3 and 4 for now) can also be realized using the cluster machinery may indicate that the `polylogarithmic cluster picture' may be bigger and more general than the dilogarithmic case alone.
 \sk 
 
 We have proved that Spence-Kummer trilogarithmic identity $(\mathcal S\mathcal K)$ and Kummer's tetralogarithmic one $({\mathcal K}_4)$ are of cluster type. A natural (but  clearly naive) consequence is to wonder about the next classical case: 
\begin{questions} 
Is Kummer's pentalogarithmic identity 
$({\mathcal K}_5)$  of cluster type? 
More precisely and in more invariant terms, is 
Kummer's pentalogarithmic 29-web  $ {\boldsymbol{\mathcal W}}_{\mathcal K_5}$  equivalent to a cluster web?  
\end{questions}

Of course the same question can be asked about any polylogarithmic AFE but it does not make sense in full generality  since one can obtain AFEs of great complexity and in particular not of cluster type from any given AFE.  Thus in order to consider cases for which an answer might be interesting, one has to focus on polylogarithmic AFEs which are irreducible, or more generally of special interest (being aware that the meaning of `special interest' is quite unprecise here). \sk 

Among the polylogarithmic identities which surely deserve to be considered, stands   trilogarithmic identity $(\boldsymbol{\mathcal G}_{22})$ which, according to Goncharov who discovered it, might play the same role for the trilogarithm as Abel's   
relation for the dilogarithm ({\it e.g.}\,see \cite[p.\,199]{Gangl0} or \cite[p.\,61]{Zagier}).

\begin{question} 
Is Goncharov's 22-web   $ {\boldsymbol{\mathcal W}}_{\boldsymbol{\mathcal G}_{22}}$   equivalent to a cluster web?  
\end{question}

This question is all the more interesting that Goncharov has given a geometric interpretation of the arguments appearing in the functional equation $(\boldsymbol{\mathcal G}_{22})$ by means of the 3-dimensional stratum  
$\boldsymbol{\Sigma}(\boldsymbol{{\mathcal G}}_{22})$
of  degenerate configurations of 7 points in $\mathbf P^2$ (see Figure \ref{Fig:CSKG}).
 The whole configuration space ${\rm Conf}_7(\mathbf P^2)$ naturally identifies itself (birationally) with  the cluster $\mathcal X$-manifold of type $A_2\boxtimes A_3$ (see 
 \S\ref{Parag:Conf-m+n+2(IPm)} and \cite{Volkov} or \cite{Weng} for more details). Since the latter quiver is mutation equivalent to $E_6$, it is of finite type hence  the set of  associated $\mathcal X$-cluster variables is finite.  Hence one  may wonder about  
 the possibility of describing Goncharov's web as a subweb of the restriction 
  of the whole cluster cluster web ${\boldsymbol{\mathcal X  \hspace{-0.04cm}\mathcal W}}_{\hspace{-0.03cm} A_2\boxtimes A_3}$ to Goncharov's stratum, which might be a way to answer the preceding question by the affirmative. \sk

However, some computations of ours lead us to the conclusion that 
$\boldsymbol{\Sigma}(\boldsymbol{\mathcal G}_{22})$ does not embed (even birationally) into the cluster $\mathcal X$-manifold of type $A_2\boxtimes A_3$. 
To explain this we use some material and notations from \cite{Weng} to which we refer for more details.  Given a certain plabic graph $P_0$, one can construct an explicit birational map 
$\eta^{P_0}: {\rm Conf}_7(\mathbf P^2)\longrightarrow \mathcal X_{A_2\boxtimes A_3}$. 
  Actually, the construction gives us the expression of this map 
  in the birational coordinates system $\boldsymbol{x}^{t_0}$  on the target given by a certain $\mathcal X$-seed  $S^{t_0}=(\boldsymbol{x}^{t_0}, Q^{t_0})$ for some vertex 
  $t_0$  of $ \mathbb T=\mathbb T^6$. We choose $S^{t_0}$  as initial  seed,  we set $\eta^{t_0}$ and for any $t\in \mathbb T$, we denote by 
  $\eta^{t} : {\rm Conf}_7(\mathbf P^2)\dashrightarrow \mathbf T^t$  the composition $x^t\circ \eta^{t_0}$ where $\mathbf T^t$ stands for the $\mathcal X$-cluster torus attached to the seed $S^{t}$.  Then  we have verified that: 
   \begin{itemize}
  \item[$-$]
    for most of the seeds $S^t$, the generic point of Goncharov's stratum is an indeterminacy point of $\eta^{t}$;  and

   \item[$-$] when it is not the case,  $\boldsymbol{\Sigma}(\boldsymbol{\mathcal G}_{22})$ is contracted by $\eta^t$ onto a subvariety  of $\mathbf T^t$ of dimension 2. 
  \end{itemize}
 
 These two possibilities prevent us from finding an avatar of  Goncharov's stratum in the cluster variety $\mathcal X_{A_2\boxtimes A_3}$ and suggest that, although it is naturally defined on the configuration space  ${\rm Conf}_7(\mathbf P^2)$ which is
 birational to $\mathcal X_{A_2\boxtimes A_3}$  in a nice way via $\eta^{t_0}$, Goncharov's web  $ {\boldsymbol{\mathcal W}}_{\boldsymbol{\mathcal G}_{22}}$ might not be of cluster type. 

\subsubsection{\bf Webs associated to multivariable polylogarithms.}
\label{SS:WebsMultivariablePolylogs}
 Here, we only give a few examples and sketch the main features of a possible theory of webs of codimension bigger than 1. Some of the notions and tools considered in this text for webs of codimension 1 can be generalized (rather straightforwardly to be honest)
 to these webs of high codimension. In particular, the idea that such a web has AMP  can make sense. We show that this is an interesting notion 
 by giving  some examples of such webs which are AMP, some of them being webs naturally associated to functional identities satisfied by polylogarithms of several variables.
 \mk

Classical polylogarithms  admit multivariable generalizations. A `{\it multiple polylogarithm}'  of  $m\geq 1$ variables   is defined as the sum  of a series of the type 
$$
{\bf L}{\rm i}_{\boldsymbol{n}}(z_1,\ldots,z_m)=\sum_{1\leq k_1<k_2<\cdots<k_m}\frac{z_1^{k_1}\cdots z_m^{k_m}}{k_1^{n_1}\cdots k_m^{n_m}}
$$
for some $m$-tuple $\boldsymbol{n}=(n_1,\ldots,n_m)\in (\mathbf N_{>0})^m$.   
Such a series  converges in the unit polydisk ({\it i.e.} when $\lvert z_i\lvert<1$  for all $i$) and defines a  multiple polylogarithm (of weigth $\lvert \boldsymbol{n}\lvert=n_1+\cdots+n_m$) which satisfies multivariable generalizations of the properties enjoyed by classical polylogarithms: ${\bf L}{\rm i}_{n_1,\ldots,n_m}$ admits an  iterated integral representation (see 
\cite[\S2]{GoncharovMultiple}), hence extends to $\mathbb C^m$ as a multivalued holomorphic functions with unipotent monodromy ({\it cf.}\,\cite{Zhao}).
What is interesting for us is that multiple polylogarithms  also satisfy some kind of `abelian functional equations'.  For instance, the weight 2 multiple logarithm 
$ {\bf L}{\rm i}_{1,1} $ satisfies the following identity:
$$ \big(L_{11}\big) \hspace{3cm} 
{\bf L}{\rm i}_{1,1} (x,y)=
{\bf L}{\rm i}_{2}\bigg( \frac{1-x}{1-1/y}\bigg)
-{\bf L}{\rm i}_{2}\bigg( \frac{1}{1-1/y}\bigg)-
{\bf L}{\rm i}_{2}( xy)\, .
\hspace{3cm} {}^{}
$$

The theory of functional equations satisfied by multiple polylogarithms is still largely unexplored but many (if not all) of these are of the following form  (see \cite{GanglAmplitude})
$$
\big( \bigstar\big)\hspace{4cm} 
 \sum_{i=1}^d c_i  \,{\bf L}{\rm i}_{\boldsymbol{n}_i}\big( u_i(x)\big)={P}(x) 
 \hspace{6cm}  {}^{}
 $$ 
where the $c_i$'s  are (rational) constants, the $\boldsymbol{n}_i$ are tuples 
 of multi-weight of the same total weight $\lvert\boldsymbol{n}_i\lvert=w\in \mathbf N_{>0}$,  the $u_i$'s rational  functions of $x\in \mathbf C^m$  (with values in $\mathbf  C^{m_{i}}$ for each $i$, where $m_i$ stands for the length of $\boldsymbol{n}_i$) and where $P$ is a rational expression of (possibly multiple) polylogarithms of weights $<w$ evaluated on rational functions as well.  To the  `AFE' $(\bigstar)$, one can associate the generalized `web'  $\boldsymbol{\mathcal W}_\bigstar$  of degree $d$ on $\mathbf C^r$,  possibly of mixed codimensions, defined by the rational first integrals $u_i: \mathbf C^m\dashrightarrow \mathbf C^{m_i}$ for $i=1,\ldots,d$: 
 $$\boldsymbol{\mathcal W}_\bigstar=\boldsymbol{\mathcal W}\big(u_1,\ldots,u_d\big)\, .$$
 
 Note that the term `web' is used in a very general and imprecise sense and that for this reason, one cannot expect to generalize for such a web the whole theory of webs presented for those of codimension 1 in the first section of this memoir. For instance, 
 if one can always define a finite invariant $\rho^\sigma(\boldsymbol{\mathcal W})$ for any $\sigma\geq 1$ to such a web $\boldsymbol{\mathcal W}$\footnote{The definition is a direct generalization of the one given for webs of codimension 1.}, it is quite possible that the total sum $\rho(\boldsymbol{\mathcal W})=\sum_{\sigma\geq 1 }\rho^\sigma(\boldsymbol{\mathcal W})$ be infinite (this already happens with $(L_{11})$).\sk

 However, when $\rho(\boldsymbol{\mathcal W})<+\infty$ holds true, the question of whether  or not $\boldsymbol{\mathcal W}$ is AMP makes sense and is relevant. An example of such a case is provided by a 5-web of codimension 2 in four variables which can be associated to a nice polylogarithmic functional equation  of weight 2 
 found in Rogers' paper \cite{Rogers}. Setting  
 $$
 \mathfrak R(x,a)={\mathcal R}(x) -
 \frac{1}{2}\log(a)\log(1-x)+\frac{1}{2}\log(x)\log(1-a)-
 {\mathcal R}(a)$$
 for $x,a\in ]0,1[$, where ${\mathcal R}(x)=\l {2} (x)+\frac{1}{2}\log(x)\log(1-x)$ stands for the original dilogarithm considered by Rogers (see \S\ref{Par:NotationDilogarithmicFunctions}), the following functional identity  $$ (\mathfrak R_2)\qquad \quad 
 \mathfrak R(x,a)+\mathfrak R(y,b)=\mathfrak R(xy,ab)+
 \mathfrak R\bigg(\, x\frac{1-y}{1-xy}\, , \, a\frac{1-b}{1-ab}\, \bigg)+
 \mathfrak R\bigg(\, y\frac{1-x}{1-xy}\, , \, b\frac{1-a}{1-ab}\, \bigg)
\, 
\hspace{2cm} {}^{}
$$ holds true
for all $x,y,a,b$ in some open domain of $\mathbf R^4$ (that needs not to be made explicit here).\footnote{When finishing the writing of this memoir, 
 has been released the preprint \cite{Monod} in which a functional identity in two variables very similar to 
 $ (\mathfrak R_2)$ appears, in relation with a particular case of a more general theory. It would be interesting to figure out whether  interesting general consequences in web geometry can be deduced from the results of \cite{Monod}.} 
 To this functional identity we associate the web   $\boldsymbol{\mathcal W}_{\mathfrak R_2}=\boldsymbol{\mathcal W}(u_1,\ldots,u_5)$ where the $u_i$'s stand for the following five $\mathbf C^2$-valued rational functions of $(x,y,a,b)\in \mathbf C^4$: 
  $$u_1=(x,a)\, , \hspace{0.16cm}
 u_2= (y,b)\, , \hspace{0.16cm}
 u_3=  (xy,ab)\, , \hspace{0.16cm}
   u_4=\bigg(\, x\frac{1-y}{1-xy}\, ,\,  a\frac{1-b}{1-ab}\, \bigg)\hspace{0.16cm} \mbox{ and }
  \hspace{0.16cm}
u_5=   \bigg(\, y\frac{1-x}{1-xy}\, ,\,  b\frac{1-a}{1-ab}\, \bigg)\, .
    $$
  
  The notion of abelian relation makes sense for this web of codimension 2: an AR for 
  $\boldsymbol{\mathcal W}_{\mathcal R_2}$ identifies itself with   a  non-trivial solution of the  following functional equation in four variables
$$
F_1(x,a)+F_2 (y,b)
+F_3(xy,ab)+F_4\bigg(\, x\frac{1-y}{1-xy}\, , \, a\frac{1-b}{1-ab}\, \bigg)+
F_5\bigg(\, y\frac{1-x}{1-xy}\, ,\,  b\frac{1-a}{1-ab}\, \bigg)\equiv 0
$$    
where the unknown $\big(F_i(\cdot,\cdot)\big)_{i=1}^5$ is a 5-tuple of (germs of) holomorphic functions in two variables.
  \sk

\vspace{-0.3cm}
One can easily construct several ARs for $\boldsymbol{\mathcal W}_{\mathcal R_2}$ by considering the two linear projections $\pi_{x,y}: (x,y,a,b)\mapsto (x,y)$ and  $\pi_{a,b}: (x,y,a,b)\mapsto (a,b)$.  Let $\boldsymbol{\mathcal B}_{x,y}=\boldsymbol{\mathcal W}\Big( x,y,xy,x\frac{1-y}{1-xy},y\frac{1-x}{1-xy}\Big)$ be Bol's 5-web in the variables $x,y$ and let us consider another copy $\boldsymbol{\mathcal B}_{a,b}$  of it,  but in the variables $a$ and $b$.
Then clearly $\boldsymbol{\mathcal W}_{\mathcal R_2}$ can be obtained from these two copies of Bol's web: in some sense (which can easily  be made precise), it is the intersection of the two pull-backs 
$\pi_{x,y}^*\big(\boldsymbol{\mathcal B}_{x,y}\big)$ and $\pi_{a,b}^*\big( \boldsymbol{\mathcal B}_{a,b} \big)$: 
$$
\boldsymbol{\mathcal W}_{\mathcal R_2} = 
\pi_{x,y}^*\big(\boldsymbol{\mathcal B}_{x,y}\big)\cap \pi_{a,b}^*\big( \boldsymbol{\mathcal B}_{a,b} \big)\, . 
$$
 
Clearly, by pull-back under $\pi_{x,y}$ and $\pi_{a,b}$, the spaces of classical ARs $\boldsymbol{\mathcal A}(\boldsymbol{\mathcal B}_{x,y}) $ and $\boldsymbol{\mathcal A}\big(\boldsymbol{\mathcal B}_{a,b}\big)$ respectively embed into $\boldsymbol{\mathcal A}\big(\boldsymbol{\mathcal W}_{\mathcal R_2}\big) $ and are obviously in direct sum in this space. This direct sum together with the AR associated to $(\mathfrak R_2)$ span a linear subspace of dimension 13 in 
$\mathcal A(\boldsymbol{\mathcal W}_{\mathfrak R_2}) $.\sk 

 On the other hand, by direct computations, one gets 
$$
\rho^\bullet\big(\boldsymbol{\mathcal W}_{\mathcal R_2}\big)=(6,5,2) \qquad \mbox{ hence }\qquad \rho\big(\boldsymbol{\mathcal W}_{\mathcal R_2}\big)=13\, , 
$$
from which it follows that the 2-codimensional web $\boldsymbol{\mathcal W}_{\mathcal R_2}$ is AMP  with polylogarithmic ARs.
\mk 

Actually, the previous construction just seems to be the second step of a whole series. Intersecting 3 copies of Bol's web in $\mathbf C^6$ with coordinates $x,y,a,b,u,v$, one obtains the following 3-codimensional 5-web 
\begin{align*}
3\boldsymbol{\mathcal B}
 = \boldsymbol{\mathcal W}\Bigg( 
 \, (x,a,u)\, , \, & 
  (y,b,v)\, , \, 
   (xy,ab,uv)\, , \, \\ 
   & 
   \bigg(x\frac{1-y}{1-xy}\, ,\,  a\frac{1-b}{1-ab}\,  , \,   u\frac{1-v}{1-uv}\bigg)\, , \, 
   \bigg(y\frac{1-x}{1-xy}\, , \, b\frac{1-a}{1-ab}, v\frac{1-u}{1-uv}\bigg)\, 
 \Bigg)
\end{align*}
which can be verified to be such that 
 $\rho^\bullet\big(3\boldsymbol{\mathcal B}\big)=(9,9,3)$.  Up to some pull-backs by obvious projections, the three 6-dimensional spaces  of ARs 
$\boldsymbol{\mathcal A}(\boldsymbol{\mathcal B}_{x,y}) $, $\boldsymbol{\mathcal A}(\boldsymbol{\mathcal B}_{a,b}) $ and $\boldsymbol{\mathcal A}(\boldsymbol{\mathcal B}_{u,v}) $ embed and are in direct sum in  
$\boldsymbol{\mathcal A}(3\boldsymbol{\mathcal B}) $.  On the other hand, 
$\boldsymbol{\mathcal B}^{\cap 3}$ carries one AR of type $(\mathcal R_2)$ for 
 each unordered pair of elements of $\{ (x,y), (a,b), (u,v)\}$ and the three corresponding  ARs are linearly independent modulo $\boldsymbol{\mathcal A}(\boldsymbol{\mathcal B}_{x,y}) \oplus \boldsymbol{\mathcal A}(\boldsymbol{\mathcal B}_{a,b}) \oplus \boldsymbol{\mathcal A}(\boldsymbol{\mathcal B}_{u,v}) $. This implies that 
  the space of polylogarithmic ARs of $3\boldsymbol{\mathcal B}$ has dimension $18+3=21$ at least. Since  this is equal to $  \rho\big(3\boldsymbol{\mathcal B} \big)=9+9+3$, one deduces that 
$3\mathcal B$ is AMP with only polylogarithmic ARs  (of weight 1 or 2).
\sk 

The construction above generalises in a straightforward way:  for any $k\geq 1$ one can define the web $k\boldsymbol{\mathcal B}$ which is a 5-web of codimension $k$ on $\mathbf C^{2k}$, which can be obtained by intersecting  $k$ distinct and independent  copies of Bol's web.  Then one can ask the 
\begin{question}
Is  $k\boldsymbol{\mathcal B}$  AMP with only 
 polylogarithmic ARs (of weight 1 or 2)  for any $k\geq 1$?
\end{question}
 
Since this is indeed the case for $k=1,2,3$, we believe that the answer is `yes' for any $k$.
\begin{center}
 $\star$
\end{center}

 The same constructions and questions seem to be relevant for other polylogarithmic webs as well, such as the Newman's dilogarithmic 6-web 
$\boldsymbol{\mathcal W}_{\hspace{-0.05cm}\mathcal N_6}$ of 
Spence-Kummer trilogarithmic web $\boldsymbol{\mathcal W}_{\hspace{-0.05cm}{\mathcal S}\mathcal{K}}$. 

For instance, intersecting two copies of the latter, each in two variables which are independant relatively to those involved in the other copy, one can construct the `double intersection'
$2\boldsymbol{\mathcal W}_{\hspace{-0.05cm}{\mathcal S}\mathcal{K}}$ which  is a 
 9-web of codimension 2 on $\mathbf C^4$. 
 One verifies that  $\rho^\bullet\big(2\boldsymbol{\mathcal W}_{\hspace{-0.05cm}{\mathcal S}\mathcal{K}}\big)=(14,17,16,10,6,4,2)$. So the total virtual rank of this web is finite. 
  By means of formal computations performed with the help of a computer algebra system, we have verified that $2\boldsymbol{\mathcal W}_{\hspace{-0.05cm}{\mathcal S}\mathcal{K}}$ is AMP to the 13th order.  Considering the approach described in \S\ref{SubPar:CharacterizationWebsMaximalRank}, this implies that $2\boldsymbol{\mathcal W}_{\hspace{-0.05cm}{\mathcal S}\mathcal{K}}$ is AMP.
%
%
%
%
 \begin{question}
  Give an explicit basis of  
 $\boldsymbol{\mathcal A} \big(2\boldsymbol{\mathcal W}_{\hspace{-0.05cm}{\mathcal S}\mathcal{K}}\big)${\rm \big)}. 
\end{question}
 
 %
 
What is potentially interesting is that this construction, 
which was suggested to us when  considering  the multivariable dilogarithmic identity 
$\big(\mathfrak R_2\big)$, seems to apply to non-polylogarithmic webs as well, such as the  web  $\boldsymbol{\mathcal W}_{\hspace{-0.05cm}xy}=\boldsymbol{\mathcal W}(x,y,x+y,x-y,xy)$. This 5-web is known to be exceptional despite the fact that it carries only rational abelian relations, see \cite{PirioChern}. Its `double' 
$$
2\boldsymbol{\mathcal W}_{\hspace{-0.05cm} xy}=
 \boldsymbol{\mathcal W}\Big( 
 \, (x,a)\, , \, 
  (y,b)\, , \, 
   (x+y,a+b)\, , \, 
   (x-y,a-b)\, , \, (xy,ab)\, 
 \Big)
$$ is such that 
$\rho^\bullet(2\boldsymbol{\mathcal W}_{\hspace{-0.05cm} xy})=(6,5,2)$ and therefore  
$\rho(2\boldsymbol{\mathcal W}_{\hspace{-0.05cm} xy})=13$. From the two copies $\boldsymbol{\mathcal W}_{\hspace{-0.05cm}xy}$ and $\boldsymbol{\mathcal W}_{\hspace{-0.05cm}ab}$ used to construct $
2\boldsymbol{\mathcal W}_{\hspace{-0.05cm} xy}$  one gets two subspaces of dimension 6 in direct sum in $\boldsymbol{\mathcal A}(2
\boldsymbol{\mathcal W}_{\hspace{-0.05cm}xy})$ which span a subspace of dimension 12. There  is another AR, namely the one corresponding to the following elementary polynomial identity 
 in the four variables $x,y,a$ and $b$:
$$
-2\,xa-2\,yb+(x+y)(a+b)+(x-y)(a-b)=0 \,. 
$$
It follows that,  as $\boldsymbol{\mathcal W}_{\hspace{-0.05cm}xy}$, its double $2
\boldsymbol{\mathcal W}_{\hspace{-0.05cm}xy}$ is AMP with only rational ARs. 
\begin{center}
 $\star$
\end{center}

The fact that the same quite simple construction applies to webs of different natures (polylogarithmic and rational) and gives interesting webs of codimension $>1$ suggests that it could be relevant for other type of webs as well.  Since our main interest lies in webs with many ARs, a natural class of webs to consider is that  of algebraic webs. 
 Let $\boldsymbol{\mathcal W}_{C}$ be an algebraic web associated to a projective curve 
$C\subset \mathbf P^n$. For any $k\geq 1$, the definition of $ k\boldsymbol{\mathcal W}_{C}$ is straighforward. 

 \begin{questions}
Let $k$ be bigger than 1. 
\begin{enumerate}
\item  Give a geometric description of $k\boldsymbol{\mathcal W}_{C}$  in terms of $C$. 
\item Describe the ARs of $\boldsymbol{\mathcal W}_{C}$.  
\item Can $k\boldsymbol{\mathcal W}_{C}$ be AMP? If yes, characterize the curves for which it is the case.
\end{enumerate}
\end{questions}

Answering the first question should be easy,  
but it could be much less so for the last two. \mk 

The general constructions and consideration as well as the examples considered above, show that it may be the case that an interesting and richer than expected  theory of webs of higher codimension does exist, provided that one relaxes the strong general position assumption which is classically assumed.  
We think that this research focus is promising and deserves further study.


\subsection{\bf About cluster variables and cluster webs}
\label{SS:AboutClusterVW}

Cluster webs are the central objects discussed in this text. However, despite its length, we are very far from having exhausted their study. Below, we list many questions left open which are interesting according to us.

 \subsubsection{\bf $\boldsymbol{\mathcal A}$-cluster webs and cluster webs with frozen variables.}
 Let $\Delta$ be a Dynkin diagram. We know that both cluster webs 
 $ {\boldsymbol{\mathcal A \hspace{-0.05cm}\mathcal W}}_{ \hspace{-0.05cm}\Delta}$ and 
  $ {\boldsymbol{\mathcal Y \hspace{-0.05cm}\mathcal W}}_{ \hspace{-0.05cm}\Delta}$
have the same degree (that is, are formed by the same number of foliations). But in all the cases 
 we have considered, we have noticed that these two cluster webs share  the same 
virtual ranks:   {\it i.e.}\,both sequences $\rho^\bullet\big( {\boldsymbol{\mathcal A \hspace{-0.05cm}\mathcal W}}_{ \hspace{-0.05cm}\Delta} \big) $ and 
$\rho^\bullet\big( {\boldsymbol{\mathcal Y \hspace{-0.05cm}\mathcal W}}_{ \hspace{-0.05cm}\Delta} \big) $ coincide. \sk

 This phenomenon seems also to hold true for webs associated to 
 pairs  of Dynkin diagrams. Given such a pair  $(\Delta,\Delta')$,  the associated $\mathcal T$-cluster web  (which is obviously the   $\mathcal A$-cluster avatar of
 the $\mathcal Y$-cluster web of bi-Dynkin type $(\Delta,\Delta')$)
 seems to have the same $\rho^\bullet$-sequence than 
 ${\boldsymbol{\mathcal Y \hspace{-0.05cm}\mathcal W}}_{ \hspace{-0.05cm} \Delta,\Delta'}$. 
  It would be interesting to know whether   this holds true in general and 
 to have an explanation of why if it is the case. \sk

On the other hand, again in all cases we considered, the  $\mathcal A$-cluster versions of the 
$\mathcal Y$-cluster webs of (bi-)Dynkin type have rather low logarithmic ranks and even trivial weight 2 hyperlogarithmic ranks.  As an example, we can consider the case of Dynkin type $A_3$. We have $\rho^\bullet\big( {\boldsymbol{\mathcal A \hspace{-0.05cm}\mathcal W}}_{ \hspace{-0.05cm}A_3} \big) =\rho^\bullet\big( {\boldsymbol{\mathcal Y \hspace{-0.05cm}\mathcal W}}_{ \hspace{-0.05cm}A_3} \big) =(6,3,1)$ and we know that ${\rm polrk}^\bullet\big( {\boldsymbol{\mathcal Y \hspace{-0.05cm}\mathcal W}}_{ \hspace{-0.05cm}A_3} \big) =(9,1)$. We have verified that  ${\boldsymbol{\mathcal A \hspace{-0.05cm}\mathcal W}}_{ \hspace{-0.05cm}A_3}$ has rank 6, with a subspace of logarithmic ARs of dimension 5 plus one extra AR, which is rational.  
Extrapolating to the case of an arbitrary Dynkin diagram, one can expect the  $\mathcal A$-cluster webs in finite type to have only rather elementary ARs (rational and logarithmic) and to be far from being 
of AMP rank. Despite this, it would be interesting to understand better the ARs of such webs since it would give us informations about  how the $\mathcal A$-cluster variables are related. 
\begin{center}
$\star$
\end{center}

The notion of frozen variable, which is important regarding the theory of cluster algebras, has not been consider in the whole text.  Allowing some variables to be frozen allows to get a more general class of cluster webs, and the fact that  some of them could be AMP and/or carry 
interesting (polylogarithmic?) ARs 
must not be excluded. \sk

As an example, let us consider the $\mathcal X$-cluster web associated to the 
  subgroup $N\subset {\rm SL}_4(\mathbf C)$  of unipotent upper triangular matrices.  In \cite[\S2.2]{GLS} the authors give an explicit seed for a 
structure of $\mathcal A$-cluster algebra on $\mathbf C[N]$. It is a cluster algebra of finite type $A_3$ in six variables, among which the last three are frozen. 
 By straightforward computations, we obtained  that there are 29 $\mathcal X$-variables associated to this cluster algebra, which give rise to a 
29-web in six variables, denoted here by ${\boldsymbol{\mathcal X \hspace{-0.05cm}\mathcal W}}_{\hspace{-0.04cm}N}$. 
It can be verified that: \sk 

${}^{}$ \hspace{0.3cm}$-$ the ramification of ${\boldsymbol{\mathcal X \hspace{-0.05cm}\mathcal W}}_{\hspace{-0.04cm}N}$ is $0,-1,\infty$;\vspace{0.06cm}\\ 
${}^{}$ \hspace{0.2cm} $-$  one has   
$\rho^\bullet({\boldsymbol{\mathcal X \hspace{-0.05cm}\mathcal W}}_{\hspace{-0.04cm}N})=(23,11,5)$, 
${\rm polrk}^\bullet({\boldsymbol{\mathcal X \hspace{-0.05cm}\mathcal W}}_{\hspace{-0.04cm}N})=(32,6)$
and  ${\rm rk}({\boldsymbol{\mathcal X \hspace{-0.05cm}\mathcal W}}_{\hspace{-0.04cm}N})=38$.\mk 

Hence one has  ${\rm rk}({\boldsymbol{\mathcal X \hspace{-0.05cm}\mathcal W}}_{\hspace{-0.04cm}N})= \rho({\boldsymbol{\mathcal X \hspace{-0.05cm}\mathcal W}}_{\hspace{-0.04cm}N})-1$ so this web can be said of `almost' AMP rank, which is a rather non-trivial property although not as pretty as being AMP. \sk 

This example naturally suggests the following 
\begin{questions}
 Let ${\boldsymbol{\mathcal X \hspace{-0.05cm}\mathcal W}}_{\hspace{-0.04cm}A}$ be the $\mathcal X$-web associated to a finite type  
 cluster algebra $A$ with frozen variables.  Then: \vspace{-0.2cm}
 \begin{enumerate}
 \item  is its ramification polylogarithmic (that is, included in $\{0, -1,\infty\}$)? 

 \item can  ${\boldsymbol{\mathcal X \hspace{-0.05cm}\mathcal W}}_{\hspace{-0.04cm}A}$ be of AMP rank?
 \item of which nature can  its abelian relations be?
 \end{enumerate}
\end{questions}

 \subsubsection{\bf Some questions about $\mathcal Y$-cluster webs and  webs associated to cluster periods.}

A substantial part of the present text is devoted to the study of cluster webs. We have obtained several results about them, but we are far  from having exhausted their study. Below we ask several questions about them. Some are rather precise, others much more vague and should better be described as indications to possibly relevant directions of research for future works.

  \paragraph{$\boldsymbol{\mathcal Y}$-cluster webs at a special point.}
Let $\Delta$ be a Dynkin diagram of rank $n\geq 2$  and Coxeter number $h$. Denote by $\boldsymbol{\mathcal W}$ the associated  $\mathcal X$- or $\mathcal Y$-cluster web. In this text, we gave formulas (which are conjectural 
in  many cases) for the virtual, polylogarithmic or standard rank(s)  for  the germ of $\boldsymbol{\mathcal W}$  at a generic point $\zeta$ of $(\mathbf C^*)^n$.  It is  natural to try  to be more precise by exhibiting such a generic point explicitely. By this, we mean an explicit point $\xi\in \mathbf C^n$ such that 
\begin{equation}
\label{Eq:208}
\scalebox{0.98}{
\begin{tabular}{l}
{\it $\bullet$ all the cluster variables of $\boldsymbol{\mathcal W}$  are defined at $\xi$ ({\it i.e.}\,$\xi$ is a regular point of $\boldsymbol{\mathcal W}$);}\vspace{0.15cm}\\
 {\it $\bullet$  if $x,\tilde x$ stand for two of them which are distinct, then $(dx\wedge d\tilde x)(\xi)\neq 0$  (in other terms:}\\ ${}^{}$ \hspace{0.15cm}{\it  the foliations of $\boldsymbol{\mathcal W}$ satisfy condition (wGP) at $\xi$);}\vspace{0.15cm} \\
 {\it $\bullet$ all the ranks (virtual, polylogarithmic, the standard) of $\boldsymbol{\mathcal W}$ at $\xi$ all coincide with}\\  ${}^{}$ \hspace{0.15cm}{\it  those at the generic point $\zeta$.}
\end{tabular}}
\end{equation}


An interesting candidate for such a base point $\xi$ is a special point considered in 
\cite{Mizuno}. Let $F_\Delta$ be the birational map defined in the last paragraph of \S\ref{SubPar:BipartiteBelt}:  it is of order $m_\Delta=2(h+2)/{\rm gcd}(h,2)$,   the iterate $(F_{\hspace{-0.05cm}\Delta})^{\circ( m_\Delta/2)}$ acts as a permutation of the coordinates and the $\mathcal Y$-cluster web of type $Y$ can be defined as the web whose first integrals are  the components of the iterates $(F_{\hspace{-0.05cm}\Delta})^{\circ s}$ for $s=0,\ldots, (m_\Delta/2)-1$.    \sk 

A natural choice for a  base point at which one could study $\boldsymbol{\mathcal Y\hspace{-0.05cm}
\mathcal W}_{\hspace{-0.05cm}\Delta}$ could be a fixed point for $F_\Delta$. 
Since the quiver $Q(\Delta,\ell)$ coincides with $\vec{\Delta}$ when $\Delta$ is simply laced and $\ell=2$, we deduce precisely from Proposition 3.6 on \cite{Mizuno} that there exists a unique point $\xi_\Delta\in (\mathbf R_{>0})^n$ such that $F_\Delta(\xi_\Delta)=\xi_\Delta$.

\begin{question}
\label{Question:812}
For $\Delta$ simply laced, does 
$\boldsymbol{\mathcal Y\hspace{-0.05cm}
\mathcal W}_{\hspace{-0.05cm}\Delta}$ satisfy
 the three conditions in \eqref{Eq:208} at $\xi_\Delta$? Is it possible to give closed (combinatorial) formulae for the differential  at this point of the cluster first integrals of this web?
\end{question}

 Actually, Mizuno proves that for any $\Delta$ and any level $\ell\geq 2$, the 
 finite order birational map  $F_{Q(\Delta,\ell)}$ associated to the quiver $Q(\Delta,\ell)$ admits a unique fixed point $\xi_{Q(\Delta,\ell)}$ with positive coordinates.  In  all the concrete $Y$-systems of bi-Dynkin type $\Delta\boxtimes \Delta'$ we have considered, we have verified that  the same holds true.  
 
\begin{questions}  Let 
$Q$ stand for a Dynkin diagram of the form $Q(\Delta,\ell)$ or $\Delta\boxtimes \Delta'$ (where $\Delta,\Delta'$ are Dynkin diagrams of rank $n$ and $n'$ respectively,  $\ell\geq 2$ being a level).
\begin{enumerate}
\item   For $Q=\Delta\boxtimes \Delta'$, is it true that  there exists a unique 
point $\xi_{Q}\in (\mathbf R_{>0})^{nn'}$ fixed by $F_{\hspace{-0.1cm} Q}$?
\item If yes: (a)  is it possible to give  a closed formula for $\xi_Q$?\sk \\
${}^{}$ \hspace{0.9cm} (b) what about the analogue for  
$\boldsymbol{\mathcal Y\hspace{-0.05cm}
\mathcal W}_{\hspace{-0.05cm}Q}$ ($= 
\boldsymbol{\mathcal Y\hspace{-0.05cm}
\mathcal W}_{\hspace{-0.05cm}\Delta\boxtimes \Delta'}$)
of   Question \ref{Question:812}?
\end{enumerate}
\end{questions}


 Some computations when $\Delta'=A_1$ and  with $\Delta$ of type $A$ or $D$ give explicit formulas for $n={\rm rk}(\Delta)$ small enough. In particular 
  if $\tilde \xi_{\Delta}$ stands for the fixed point $E(\xi_{\Delta})$ of the conjugation of $F_{\Delta}$ by the involution 
 $E: (u_i)_{i=1}^n\mapsto (u_i^{\epsilon(i)})_{i=1}^n$ where $\epsilon(i)=1$ (resp.\,-1) if $i$ is a sink (resp.\,a source) in $\vec{\Delta}$,  we have verified  for $n$ small enough ($n\leq 9$ say) that 
 \begin{itemize}
 \item  there exists   an  algebraic integer $\chi_n$ (which is quadratic for $n$ odd, cubic when  $n$ is even)  such that the $n$ coordinates of $\tilde \xi_{A_n}$  all belong to $\mathbf Q_{>0}+ \chi_n \mathbf Q_{>0}$;
 \item one has 
   $\tilde \xi_{D_n}\in (\mathbf N_{>0})^n$. 
    \end{itemize}
  It is these few facts that prompted us to ask  question {\it 2.(a)}  in the series of three above.

  \paragraph{About the dilogarithmic ARs of cluster webs.}
The dilogarithmic identities associated to cluster periods have not yet revealed all their secrets.

From \cite{KY}, we know that for any Dynkin diagram, the $\mathcal Y$-cluster dilogarithmic identity $({\sf R}_\Delta)$ is accessible from Abel's relation $({\sf R}_{A_2})$. 
To what extent does this generalize? More precisely:  
\begin{questions}  Let $\nu$ be a cluster period.
\begin{enumerate}
\item   Is the dilogarithmic identity $({\mathsf R}_\nu)$  accessible from $(\mathcal Ab)$? 
\item  In case $\nu$ is the $\mathcal Y$-period of bi-Dynkin type $(\Delta,\Delta')$ and when it can be done, 
find  an explicit writing of $({\mathsf R}_{\Delta,\Delta'})$ as a linear combination of Abel's relations,  if possible in an uniform manner with respect to the 
pair of Dynkin diagrams considered. 
\end{enumerate}
\end{questions}

In \cite{KY}, the authors conjecture that the answer to {\it 1.} when the cluster period is of bi-Dynkin type. As for {\it 2.}, an answer has been given only for the pairs $(A_n,A_1)$ with $n\geq 2$  ({\it cf.}\,the proof of Th\'eor\`eme 4.1 in \cite{Souderes}). 
\mk 

The fact that cluster periods give rise to dilogarithmic identities is now well-understood. 
On the other hand, nothing is really known with regard to a possible reciprocal:
\begin{question}
Let $\boldsymbol{\mathcal W}$ a cluster web carrying a non-trivial dilogarithmic identity. 
  Is it necessarily the one associated to a cluster period $\nu$ such that $\boldsymbol{\mathcal W}_\nu$ is a subweb of $\boldsymbol{\mathcal W}$? 
\end{question}

Another related but distinct question concerns the construction of the dilogarithmic identity $({\mathsf R}_\nu)$ when $\nu$ is  a given cluster period from the logarithmic ARs of $\boldsymbol{\mathcal W}_{\nu}$. 
In \S\ref{SS:cluster-webs-type-A2}, we have explained how to construct $({\sf R}_{A_2})$ from five logarithmic ARs of $\boldsymbol{\mathcal Y\hspace{-0.03cm} \mathcal W}_{\hspace{-0.04cm}A_2}$, considered both in functional and differential form.  One can wonder whether a similar approach is possible in full generality: 
\begin{questions} 
 Let $\nu$ be a cluster period.
\begin{enumerate}
\item Can a basis of $\boldsymbol{LogAR}(\boldsymbol{\mathcal W}_{\nu})$ be described in a uniform way (with respect to $\nu$)?
\item Can the computations in \eqref{Eq:RA2-from-LogARl} be generalized 
in order to get an effective construction of  $({\mathsf R}_\nu)$ from the (or some) logarithmic ARs of $\boldsymbol{\mathcal W}_{\nu}$?
\end{enumerate}
\end{questions}
An affirmative answer to the second question would give a constructive proof of 
Nakanishi's theorem (Theorem \ref{Thm:Nakanishis-identity} above)  using only elementary notions.

 \paragraph{\bf Arrangements, resonance and webs.}
\label{Parag:ArrangementsResonanceWebs}

In \S\ref{SubSec:ZigZag-Map}, we proved that the complement noted there $X_{A_n}$ of the braid arrangement in $\mathbf C^n$ is isomorphic to that of the cluster arrangement 
$\boldsymbol{Arr}_{A_n}$
of type $A_n$.  It follows that, despite this arrangement is cut out by polynomial equations of degree $>1$ (for most of them), its complement $U_{A_n}$ has the same topology as $X_{A_n}$, hence is well-understood (for instance, this gives for free that $U_{A_n}$ is 1-formal, which is not clear at first sight).  In particular, one can deduce from this an explicit description of the resonance variety $R^1(U_{A_n})$  then realize that the resonance web of 
$U_{A_n}$ coincides with the whole $\mathcal X$-cluster web 
$\boldsymbol{\mathcal X\hspace{-0.05cm}
\mathcal W}_{\hspace{-0.05cm}A_n}$. \mk 

The nice properties enjoyed by the complement of the cluster arrangement in type $A$ prompt to study more systematically the complement of the set of common leaves of cluster webs with AMP rank or carrying interesting ARs. \sk

As a first example, we consider $\boldsymbol{\mathcal U\hspace{-0.05cm}\mathcal W}_{\hspace{-0.05cm} A_3}$ which is equivalent  to Spence-Kummer web according to Theorem \ref{T:classical-cluster-webs}. Its set of common leaves $\Sigma^c(\boldsymbol{\mathcal U\hspace{-0.05cm}\mathcal W}_{\hspace{-0.05cm} A_3})$
in $\mathbf P^2$ is the union of seven hypersurfaces: the line at infinity and the closure of the affine hypersurfaces cut out by the six following `$F$-polynomials'
$$ 
u_1\, , \, u_2
\, , \, 1+u_1
\, , \,  1+u_2
\, , \,  1+u_1+u_2
\quad \mbox{ and }\quad   1+2u_1+u_2+u_1^2
\, . 
$$
 We set $U_{\boldsymbol{\mathcal U\hspace{-0.05cm}\mathcal W}_{\hspace{-0.05cm} A_3}}
 =\mathbf P^2_u \setminus \Sigma^c(\boldsymbol{\mathcal U\hspace{-0.05cm}\mathcal W}_{\hspace{-0.05cm} A_3})$. Since Spence-Kummer web is equivalent to the resonance web of the non-Fano arrangement ${\sf NF}$, the following question naturally arises: 
 \begin{question}
 Is the complement of cluster hypersurfaces $U_{\boldsymbol{\mathcal U\hspace{-0.05cm}\mathcal W}_{\hspace{-0.05cm} A_3}}$ linearizable ? More precisely, does it exist a birational map $\Psi: \mathbf P^2_u\dashrightarrow \mathbf P^2$  inducing an isomorphism from 
$U_{\boldsymbol{\mathcal U\hspace{-0.05cm}\mathcal W}_{\hspace{-0.05cm} A_3}}$ onto the complement of an arrangment of lines $\boldsymbol{A}_{\boldsymbol{\mathcal U\hspace{-0.05cm}\mathcal W}_{\hspace{-0.05cm} A_3}}\subset \mathbf P^2$?
\end{question}
Of course, we expect that the answer is affirmative and that 
one can take ${\sf NF}$ for 
 $\boldsymbol{A}_{\boldsymbol{\mathcal U\hspace{-0.05cm}\mathcal W}_{\hspace{-0.05cm} A_3}}$. \mk 
 
If the cluster complement $U_\Delta$ is linearizable when $\Delta$ is of type $A$, this is not clear for other Dynkin types and it would be interesting to understand better  these complements topologically.

 \begin{questions} Let $\Delta$ be a Dynkin diagram, not of type $A$. 
\begin{enumerate}
\item  Does the sets of common leaves  
of the web  $\boldsymbol{\mathcal Y\hspace{-0.05cm}\mathcal W}_{\hspace{-0.05cm} \Delta}$ coincide with the cluster arrangement $\boldsymbol{Arr}_\Delta$? 
Same question for $\boldsymbol{\mathcal X\hspace{-0.05cm}\mathcal W}_{\hspace{-0.05cm} \Delta}$.
\item What can be said about the combinatorics of $\boldsymbol{Arr}_\Delta$ and about the topology of the complement $U_\Delta$? 
\item In particular,  what are the irreducible components of the first resonance (resp.\,characteristic) variety $R^1(U_\Delta)$  (resp.\,$V^1(U_\Delta)$)? Is $U_\Delta$ 1-formal?  
\item Can the pair $(\boldsymbol{Arr}_\Delta,U_\Delta)$ be linearized? That is: 
does it exist $\Phi: \mathbf P^n\dashrightarrow \mathbf P^n$ 
inducing by restriction an isomorphism from $U_\Delta$ onto the complement of a 
 hyperplane arrangement?
\item How are related the cluster webs $\boldsymbol{\mathcal X\hspace{-0.05cm}\mathcal W}_{\hspace{-0.05cm} \Delta}$
and $\boldsymbol{\mathcal Y\hspace{-0.05cm}\mathcal W}_{\hspace{-0.05cm} \Delta}$ on the one hand, and the resonance and characteristic webs 
$\boldsymbol{\mathcal W}_{\hspace{-0.05cm} R^1(U_\Delta)}$ and 
$\boldsymbol{\mathcal W}_{\hspace{-0.05cm} V^1(U_\Delta)}$
on the other hand?
\end{enumerate}
\end{questions}

A cluster algebra of bi-Dynkin type $({\Delta,\Delta'})$ is not of finite type in general but 
even if it is  one can formulate similar questions as those just above for $\boldsymbol{\mathcal Y\hspace{-0.05cm}\mathcal W}_{\hspace{-0.05cm} \Delta,\Delta'}$. The interested readers, if any, will elaborate on this by themselves...

  \paragraph{About the characterization of $\boldsymbol{\mathcal Y}$-cluster webs.}
An interesting combinatorial invariant attached to a web ${\boldsymbol{\mathcal W}}$ 
is the {\bf rank function} $\boldsymbol{r_{\boldsymbol{\mathcal W}}}$ which is  the integer-valued function defined on the set of subwebs of ${\boldsymbol{\mathcal W}}$ by 
$r_{\boldsymbol{\mathcal W}}({\boldsymbol{\mathcal W}}')={\bf rk}({\boldsymbol{\mathcal W}}')$  for any  ${\boldsymbol{\mathcal W}}'\subset {\boldsymbol{\mathcal W}}$. \sk

According to Bol \cite{Bol} (see also \cite{Robert2} for a more recent document), Bol's web $\boldsymbol{\mathcal B}$ can be characterized (up to local equivalence) as the unique planar 5-web with at least nine of its 3-subwebs being hexagonal (which  is equivalent for them to have rank equal to 1). This implies in particular that this web is characterized by 
$r_{\boldsymbol{\mathcal B}}$. We recall that the $\mathcal Y$-cluster webs of type $A_n$ share for any $n\geq 2$ similar properties to those of  Bol's web (which corresponds to the case of $A_2$): AMP rank, all the ARs are logarithmic except the dilogarithmic one $({\sf R}_{A_n})$. For this reason, the question arises about the extension of the characterization of Bol's web by its rank function  to any $\mathcal Y$-cluster web, not necessarily of type $A$ but more generally,  of any given  Dynkin type $\Delta$: 
\begin{question}
Is ${\boldsymbol{\mathcal Y\hspace{-0.04cm}\mathcal W}}_{\hspace{-0.04cm}\Delta} $ characterized (up to local equivalence) by its rank function 
$r_{ {\boldsymbol{\mathcal Y\hspace{-0.04cm} \mathcal W}}_{\hspace{-0.04cm} \Delta}}$?
\end{question}

One expects the answer to be affirmative, at least for $\Delta$ of type $A$, 
but we do not have any serious argument on which to base this intuition. It should be possible to handle the $A_3$ case by direct calculations.

{{\paragraph{Webs associated to generalized cluster algebras.}
\label{Par:Generalized-cluster-Algebras}
Cluster algebras have been generalized by Chekov and Shapiro in \cite{ChekhovShapiro} where they introduce the notion of  `{\it generalized cluster algebra}' (GCA for short). 
The definition of such an algebra is similar to that of a cluster algebra in many ways: 
it is generated (as an algebra) by cluster variables, which come in 
`generalized clusters'. These are obtained from an initial `generalized seed' by means of `generalized mutations' which are generalizations of classical cluster mutations by means of formulas involving `exchange polynomials'. These polynomials  can be of arbitrary positive degree in the case of GCA's, contrarily to the case of classical  cluster algebras which corresponds to the one when all the exchange polynomials are of degree 1 (and actually coincide with $P=1+x$).\sk
%

{{What is interesting with regard to the geometry of webs is that Nakanishi has 
shown in {\rm \cite{Nakanishi2015}} (see also {\rm \cite{Nakanishi2018}}) that 
the notion of period generalizes in the realm of GCA and more interestingly for our purpose, that to each such `generalized period' is associated an AFE satisfied by what Nakanishi calls a `{generalized (Rogers) dilogarithm}', which can be seen as a weight 2 iterated  integral\footnote{Let $P$ be a non-constant polynomial whith no real root distinct from -1. The generalized dilogarithm associated to such a polynomial by Nakanishi is the function $L^P$ 
defined by  requiring that the identity $ L^P\left( {x^\delta}/{P(x)}\right)=\frac{1}{2}\int_{0}^x \left( 
{{\rm Log}\big(P(u)\big)}/{u}-
{{\rm Log}(u)}/{P(u)}
\right) du$ holds true for any $x>0$. When $P(x)=1+x$,  one has 
$L^{1+x}(y)={\sf R}(y/(1-y))={\mathcal R}(y)$ for every $y\in ]0,1[$, where 
${\sf R}$ and  $\mathcal R$  are respectively  the cluster and the Rogers dilogarithm  of \S\ref{Par:NotationDilogarithmicFunctions}.} (see \cite{Nakanishi2018}, in particular 
\S3 and Theorem 4.5 therein). If one is interested by constructing webs with AMP rank, it is natural to look at webs defined by the generalized cluster variables appearing as arguments of the generalized dilogarithms involved in Nakanishi's AFE associated to a generalized cluster period in a GCA. }}
\mk 

%
%

{{In order to deal with an explicit example, let us be a bit more precise and explicit about 
some basic notions of the theory of GCAs.  By definition, a generalized seed is a tuple  $\Sigma=(\boldsymbol{a},\boldsymbol{x}, B, \boldsymbol{\delta}, \boldsymbol{z})$ where: $(\boldsymbol{a},\boldsymbol{x}, B)$ is a seed in the classical sense; $\boldsymbol{\delta}=(\delta_1,\ldots,\delta_n)\in \mathbf N^n$ where $\delta_i$ specifies the degree of the $i$th exchange polynomial $P_i$ 
whose coefficients are encoded by $\boldsymbol{z}=(z_{i,s})$, which then is a  collection of scalars for $i=1,\ldots,n$ and $s=0,\ldots,\delta_i$, with $z_{i,0}=z_{i,\delta_i}=1$ for any $i=1,\ldots,n$, related to the $P_i$'s according to the relation 
$P_i(x)=\sum_{s=0}^{\delta_i} z_{i,s}x^s$ for each $i=1,\ldots,n$. }}
\bk 
\bk

 The associated generalized mutations are  birational transformations given by explicit formulas involving the degree $\delta_i$ and the coefficients of the $P_i$. For instance,  let $\Sigma'=(\boldsymbol{a}',\boldsymbol{x}', B', \boldsymbol{\delta}', \boldsymbol{z}')=\mu_k(\Sigma)$ be the mutation 
 in the $k$-th direction of a  generalized seed $\Sigma=(\boldsymbol{a},\boldsymbol{x}, B, \boldsymbol{\delta}, \boldsymbol{z})$.  Then the 
generalized  $\mathcal A$- and $\mathcal X$-cluster variables are given by the following formulas, where the $b_{kl}$'s stand for the coefficients of the exchange matrix $B$ and where 
we have set 
 $\hat a_k=\prod_{\ell=1}^n a_k^{b_{\ell k}}$:
 \begin{align*}
a_{j}'= & \, a_j  \hspace{-0.5cm} &\mbox{if }\, j\neq k   \qquad \quad    &\mbox{ and }\quad \qquad
a_k'=  \frac{1}{a_k}\left( \,   \prod_{\ell=1}^n  a_\ell^{[-b_{\ell k}]_+} \right)^{\delta_k}  
\hspace{-0.15cm}
P_k\big(\hat a_k\big) \hspace{0.9cm}  \mbox{if }\, j= k \sk \\
x_{j}'= & \, {x_j}^{-1}  \hspace{-0.5cm} &\mbox{ if }\, j= k   \qquad \quad    &\mbox{ and }\quad \qquad
x_j'=  {x_j}\left(   x_k^{[b_{kj}]_+} \right)^{\delta_k}  
P_k(x_k)^{-b_{ki}} 
\hspace{0.8cm}  \quad  \mbox{if }\, j\neq  k
\end{align*}
(these formulas should be compared to 
the classical mutation formulae  
\eqref{Eq:A-X-Mutation-formulae},  which  are the ones when all the  $P_k$'s coincide with the standard exchange polynomial $P(x)=1+x$).
\mk


We are going to consider the following example in rank 2, taken from 
\cite[2C\&3A]{Naka2015P} (see also \cite[\S6.2]{Nakanishi2018}),  which is the one when  the 
 initial exchange matrix is $B$ 
 and the corresponding exchange polynomials $P_1$ and $P_2$ are
  $$
 B =B_{A_2}=\Big[
\hspace{-0.2cm}
\scalebox{0.9}{
\begin{tabular}{c}
0 \hspace{0.05cm} $-1$\vspace{-0.1cm}\\
$1$  \hspace{0.16cm} 0
\end{tabular}} \hspace{-0.15cm} \Big]
\qquad \mbox{ and }\quad 
P_1(x)=1+s x+x^2,  \quad P_2(x)=1+x\, , $$ 
where  $s$ stands for a fixed parameter.  
One can verify that there is only a finite number of $\boldsymbol{\mathcal A}$- and $\boldsymbol{\mathcal X}$-cluster variables and that  they all    can  be computed explicitly (see Table 1 of  \cite{Naka2015P}).  For instance, up to inversion $f\leftrightarrow 1/f$, the set of these generalized $\boldsymbol{\mathcal A}$-cluster variables for this GCA is the following: 
$$
\Bigg\{\, {a_1}
\, ,\, {a_2}
\, ,\, {\frac {s{ a_2}+{{ a_2}}^{2}+1}{{a_1}}}\, ,\,{\frac {1+{ a_1}}{{ a_2}}}  \, ,\, {\frac {s{ a_2}+{{ a_2}}^{2}+{ a_1}+1}{{ a_1}\,{ a_2}}}\, ,\,
{\frac {s\, { a_1}{ a_2}+s\, { a_2}+{{ a_1}}^{2}+{{ a_2}}^{2}+2\,{ a_1}+1}{{ a_1}\,{{ a_2}}^{2}}}\, \Bigg\}\, 
$$
and there are similar, but a bit more involved,  formulas for the associated 
$\boldsymbol{\mathcal X}$-cluster variables. \sk 

The webs defined by the $\boldsymbol{\mathcal A}$-cluster variables or the $\boldsymbol{\mathcal X}$-cluster variables are equivalent\footnote{This is similar to what occurs in the classical case,  {\it cf.}\,Remark \ref{Rem:riri}. We believe that this is a general fact which holds true for all finite type  generalized cluster algebras of rank 2.} 
hence it is easier to deal with the former since they   have simpler expressions.  
By definition, the generalized $\boldsymbol{\mathcal A}$-cluster web associated to 
the chosen initial generalized seed, which we denote by $\Sigma(s)$, is the one admitting as first integrals these six generalized cluster variables
\begin{align*}
\boldsymbol{\mathcal A\hspace{-0.05cm}\mathcal W}_{\Sigma(s)}=
\boldsymbol{\mathcal W}
\Bigg(\, {a_1}\, ,\, {a_2}\, ,\, {\frac {s{ a_2}+{{ a_2}}^{2}+1}{{a_1}}}\, ,\,{\frac {1+{ a_1}}{{ a_2}}}  \, ,\, &{\frac {s{ a_2}+{{ a_2}}^{2}+{ a_1}+1}{{ a_1}\,{ a_2}}}\, ,\,
\\ & 
{\frac {\alpha\, { a_1}{ a_2}+s\, { a_2}+{{ a_1}}^{2}+{{ a_2}}^{2}+2\,{ a_1}+1}{{ a_1}\,{{ a_2}}^{2}}}\, \Bigg)\, .
\end{align*}

By direct computations, we have been able to establish the following facts: 
\begin{itemize}
\item  the web $\boldsymbol{\mathcal A\hspace{-0.05cm}\mathcal W}_{\Sigma(s)}$ 
has maximal rank 10   if and only if  $s\in \{ \pm 2\, , \, 0\}$. In this case, this web 
is equivalent to the $\boldsymbol{\mathcal X}$-cluster web of type $B_2$. 
\item the web $\boldsymbol{\mathcal A\hspace{-0.05cm}\mathcal W}_{\Sigma(s)}$ 
has  rank 9   if  $s\not \in \{ \pm 2\, , \, 0\}$, with polylogarithmic rank equal to $(6,2)$. 
\end{itemize}

Considering the invariants of the webs $\boldsymbol{\mathcal A\hspace{-0.05cm}\mathcal W}_{\Sigma(s)}$ and comparing them with those of the family of webs 
$\boldsymbol{\mathcal W}_{N_{6,\lambda}}$, $\lambda\in \mathbf C$  discussed in \S\ref{Subpar:N6-lambda} (see \eqref{Eq:Invariants-WN6}   and \eqref{Eq:Invariants-WN6-lambda} there), leads us to think  that these two families  may be isomorphic (we believe that this is rather easy to verify). If this were true,  it would be satisfying to have at disposal a general algebraic framework in which classical functional equations such as $(N_6,\lambda)$ of page \pageref{Page:N6lambda} can be understood. 
\begin{questions}
1. Are the two families of webs $\boldsymbol{\mathcal W}_{N_{6,\lambda}}$, $\lambda\in \mathbf C$ and $\boldsymbol{\mathcal A\hspace{-0.05cm}\mathcal W}_{\Sigma(s)}$, $s\in \mathbf C$ equivalent? If yes, find explicit equivalences between corresponding webs of each of these two families.
 \sk 
 
2. That the $A_2$-exchange matrix $B=B_{A_2}$ gives rise to a deformation of the cluster web  of type $B_2$ is a bit surprising. Is there a conceptual explanation of this?  
\sk

3.  Does this generalize to all classical rank 2 cluster webs $\boldsymbol{\mathcal X\hspace{-0.05cm} \mathcal W}_{\hspace{-0.05cm} \Delta}$ ? For $\Delta=A_2$ or $G_2$, does it exists an analytic deformation 
$\boldsymbol{\mathcal X\hspace{-0.05cm} \mathcal W}_{\hspace{-0.05cm} \Delta,\lambda}$ of $\boldsymbol{\mathcal X\hspace{-0.05cm} \mathcal W}_{\hspace{-0.05cm} \Delta,0}=\boldsymbol{\mathcal X\hspace{-0.05cm} \mathcal W}_{\hspace{-0.05cm} \Delta}$ such that all the polylogarithmic ARs of the latter web deform as well into hyperlogarithmic ARs for $\boldsymbol{\mathcal X\hspace{-0.05cm} \mathcal W}_{\hspace{-0.05cm} \Delta,\lambda}$? 
\end{questions}

A interesting candidate for possibly answering to the third question by the affirmative  in the case when $\Delta=G_2$ is given in  \cite[\S6.3]{Nakanishi2018}. 
The corresponding 2-dimensional family of generalized cluster 8-webs associated to the identities (6.15) in {\it loc.\,cit.}\,is quite interesting and deserves further study (in our opinion).  
\subsection{\bf Polylogarithmic AFEs of higher weight and cluster webs.}
\label{SS: Polylog-AFE-higher-weight-cluster-webs}

The fact that the theory of cluster algebras gives rise to many dilogarithmic identities (thanks to the work of many authors, culminating in Nakanishi's \cite{Nakanishi}, {\it cf.}\,also Theorem \ref{Thm:Nakanishis-identity} above) is by now common knowledge.  What is less well known, even by many specialists in these fields, is that one can get interesting AFEs satisfied by polylogarithms of weight higher than 2 from cluster algebras.  
For the moment, there are only a few examples of this (already discussed above in this text above,  see \S\ref{Par:GSVV} and also Theorem \ref{T:classical-cluster-webs}) 
and it is our opinion that this really deserves to be investigated further. \sk 

In this subsection, we discuss some points suggested by the few already known results about this as well as possible paths to construct new higher weight polylogarithmic identities from cluster algebras.
\bk

One of the main results of this text shows that many of the classical functional equations of polylogarithms (of low order) can be produced whithin the same framework, namely as cluster webs (for the dilogarithm) or as the restrictions of the latter to the secondary cluster manifolds of the corresponding cluster algebras.
 A natural question is then to ask if one can produce some functional equations for higher order polylogarithm using the same recipe.  
 This leads first to the question of understanding better the secondary cluster
manifold of the cluster algebras associated to quivers

 \subsubsection{\bf About the image of the cluster map $\boldsymbol{p: \mathcal A\rightarrow \mathcal X}$  in type $\boldsymbol{A}$.}

We think that our result that the secondary cluster web of type $A_3$ is equivalent to Spence-Kummer web hence carries trilogarithmic ARs is quite nice but we have to admit that we are not aware of any satisfying (conceptual) reason explaining it. We really believe that there is something to understand here. \sk 

In this perspective, a first step certainly consists in understanding better the secondary cluster variety $\boldsymbol{\mathcal U}_{A_3}$ which lies inside the cluster manifold  
$\boldsymbol{\mathcal X}_{A_3}$. Since the latter identifies naturally with $\mathcal M_{0,6}$ which has been intensively studied by algebraic geometers, one can expect to better understand $\boldsymbol{\mathcal U}_{A_3}$ by considering the corresponding surface in 
$\mathcal M_{0,6}$. 

\paragraph{} ${}^{}$\hspace{-0.6cm}
To explain what is to come, it is necessary to recall some (well-known) facts about the geometry of the moduli spaces $\mathcal M_{0,n+3}$'s. For any $n\geq 1$, this moduli space admits a nice (for instance smooth and proper) compactification $\overline{\mathcal M}_{0,n+3}$ constructed by Deligne, Knudsen and Mumford (hence called by their names) obtained by adding to $\mathcal M_{0,n+3}$ 
boundary strata  isomorphic to the product $\mathcal M_{0,n_1+3}\times \cdots\times 
\mathcal M_{0,n_s+3}
$  for some tuples $(n_i)_{i=1}^s\in \mathbf N^s$ such that $\sum_{i=1}^s n_i\leq n$. 
The study of the compactifications $\overline{\mathcal M}_{0,n+3}$ has been the source of a huge number of works. In particular,  there was a conjecture by Fulton about the cone ${\rm Eff}_k(\overline{\mathcal M}_{0,n+3})$ of effective algebraic cycles of codimension $k$ of $\overline{\mathcal M}_{0,n+3}$:  is it generated by boundary cycles?
This has been proved to be false for divisors ($k=1$) first in the case when $n=3$ in \cite{Vermeire} where the author proved that a certain  irreducible divisor in $\overline{\mathcal M}_{0,6}$
is effective but is not in the positive cone generated by boundary divisors. This divisor, denoted by $D_{KV}$ (where $KV$ is for `Keel-Vermeire') can be described geometrically as follows ({\it cf.}\,\cite{CT}): given a fixed Spence-Kummer configuration $C_{\hspace{-0.05cm}{\cal S}{\cal K}}\subset \mathbf P^2$ ({\it cf.} 
\S\ref{SubPar:GeomSK-Web}), $D_{KV}$ coincides with the image 
of the rational map $\mathbf P^2\dashrightarrow \overline{\mathcal M}_{0,6}$  associating to a generic point $p\in \mathbf P^2$ the image of 
$C_{\hspace{-0.05cm}{\cal S}{\cal K}}$ by the linear projection $\mathbf P^2\dashrightarrow \mathbf P^1$ from $p$.  Of course, by taking the image of $D_{KV}$ under the natural action of $\mathfrak S_6$ (by permuting the points), 
one gets other divisors that  will be called and denoted the same (a bit abusively).
\sk 

%
\begin{figure}[h]
\begin{center}
\resizebox{1.4in}{1.4in}{
 \includegraphics{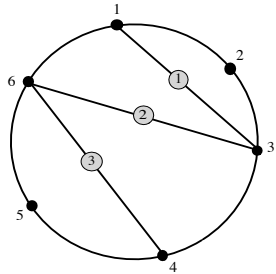}}
\vspace{-0.2cm}
\caption{The zig-zag triangulation $T_0$ of the hexagon.}
\label{Fig:T0-hexagon}
\vspace{-0.4cm}
\end{center}
\end{figure}

Let $x^{T_0}: \mathcal M_{0,6} \dashrightarrow \boldsymbol{\mathcal X}_{A_3}$
 be the birational identification associated to the zig-zag triangulation $T_0$ of the hexagon 
 ${\bf P}_6$ above  (see also   \S\ref{SSub:Type-An}).  One verifies that the strict transform of $\boldsymbol{\mathcal U}_{A_3}$ by $x^{T_0}$ is a surface hence gives rise to a divisor in 
 $\overline{\mathcal M}_{0,6}$ (after taking Zariski closure).

\begin{prop}
\label{Prop:UA3-KV}
The closure of  $(x^{T_0})^{-1}(\boldsymbol{\mathcal U}_{A_3})$  in $\overline{\mathcal M}_{0,6}$ is a Keel-Vermeire divisor.
\end{prop}
\begin{proof} We are going to compare the pull-backs of both divisors  in the statement under  the natural map $\varphi: \mathbf C^6\dashrightarrow  \mathcal M_{0,6},\, (z_i)_{i=1}^6\mapsto [z_1,\ldots,z_6]$. 
\mk

Let $X_1,X_2$ and $X_3$ be the pull-backs under $x^{T_0}\circ \varphi: \mathbf C^6\dashrightarrow \boldsymbol{\mathcal X}_{A_3}$ of the   $\mathcal X$-cluster variables
of the $A_3$-cluster associated  to the zig-zag triangulation $T_0$ of the hexagon (Fig.\,\ref{Fig:T0-hexagon}). One has 
$$ X_1=r(z_1,z_2,z_3,z_6)\, \qquad X_2=r(z_6,z_1,z_3,z_4)
\qquad \mbox{ and }\qquad 
X_3=r(z_6,z_3,z_4,z_5) \, $$
(where $r$ stands for Fock-Goncharov's cluster cross-ratio \eqref{Eq:cluster-r}).  On the other hand, the exchange matrix associated to $T_0$ is 
$$
B_{A_3}=\begin{bmatrix} 0 & 1 & 0\\
-1 & 0 & -1\\
0 & 1 & 0
\end{bmatrix}
$$
from which it comes immediately that 
the pull-back  of 
$ \boldsymbol{\mathcal U}_{A_3}$ under $x^{T_0}\circ \varphi $ is cut out by 
$ X_1=X_3$ or equivalently (but more explicitly) by 
 \begin{equation}
 \label{Eq:Keel-Vermeire}
 (12)(34)(65)+(16)(23)(45)=0
 \end{equation}
where we use the notation  $(ij)=z_i-z_j$ for any $i,j=1,\ldots,6$. \sk

From \cite[Remark\,9.2]{CT},
\eqref{Eq:Keel-Vermeire} can be seen as  an equation 
of the strict transform of $D_{KV}$ by $\varphi$ in $\mathbf C^6$ (modulo a permutation of the $z_i$'s). The proposition  follows. 
\end{proof}

One verifies easily that the cluster map $p: \boldsymbol{\mathcal A}_{A_n}\rightarrow \boldsymbol{\mathcal X}_{A_n}$ has rank $n$ for $n$ even  and rank $n-1$ when it is odd. Consequently, the secondary cluster variety $\boldsymbol{\mathcal U}_{A_n}$ is  a hypersurface in the latter case. Given the preceding proposition about the case when $n=3$, it is natural to try to understand better the corresponding irreducible divisor  of $\overline{\mathcal M}_{0,n+3}$.

\begin{questions} Assume that $n$ is odd.
What is the closure of $(x^{T_0})^{-1}(\boldsymbol{\mathcal U}_{A_n})$ 
in $\overline{\mathcal M}_{0,n+3}$?  Can this divisor be described geometrically? 
Is it one of the `hypertree divisors' defined in {\rm \cite{CT}}? If not, does it belong to the subcone of ${\rm Eff}_1(\overline{\mathcal M}_{0,n+3})$ spanned by boundary and hypertree divisors? 
\end{questions}

In addition of this algebraic geometry perspective,  the secondary cluster manifold  $\boldsymbol{\mathcal U}_{A_n}$ can also be considered from other points of view, such as Poisson geometry 
 as we explain quickly below.
\sk  

For $n$ odd, in the initial $\mathcal X$-cluster seed $\big( (u_i)_{i=1}^n, \vec{A}_n\big)$, the secondary manifold $\boldsymbol{\mathcal U}_{A_n}$ is cut out  by 
  \begin{equation}
 \label{Eq:Casimir}
 1=
 \frac{u_1}{u_3}\cdot \frac{u_5}{u_7}\cdots   
  \cdot \left(\frac{u_{n}}{u_{n-2}}\right)^{(-1)^{\lfloor n/2 \rfloor}}. 
 \end{equation}
 There is a more intrinsic `Poisson-theoretic' view on this equation: since the rank of $B_{\vec{A}_{n}}$ is $n-1$, 
 the cluster manifold $\boldsymbol{\mathcal X}_{A_{n}}$ is not symplectic but only Poisson. The field of Casimir rational functions is generated by one element $\mu$, which can be taken as the RHS of 
 \eqref{Eq:Casimir} (see the proof of Lemma 5.13 in \cite{GSV}).  One verifies that the 
pull-back  $(x^{T_0})^*\mu=\mu \circ x^{T_0}: \, 
 \mathcal M_{0,n+3}\rightarrow \mathbf C$ coincides up to sign with the classical multiratio denoted 
by  ${\sf M}_{n+3}$ in \cite{KSc} which has been considered by several authors, 
in relation with classical projective geometry or integrable discrete systems. 
 In particular, 
the relation ${\sf M}_6=-1$ (which corresponds to \eqref{Eq:Casimir} in the case when $n=3$) is related to Menelaus theorem and a geometric interpretation of ${\sf M}_8=1$ (that we believe to correspond to  \eqref{Eq:Casimir} for $n=4$) is given in \cite{KSc} to which we refer for details and references. 

\paragraph{} ${}^{}$\hspace{-0.6cm}
\label{Parag:Conf-m+n+2(IPm)}
It turns out that the questions above can be generalized to the bi-Dynkin type $(A_m,A_n)$.   The corresponding $\mathcal X$-cluster variety can be interpreted 
as  the  moduli space 
${\rm Conf}_{m+n+2}(\mathbf P^m)$ 
of projective classes of configurations of $m+n+2$ points  in $\mathbf P^m$. 
As explained in \cite{Weng} to which we refer for definitions and details, given any `minimal bipartite graph' $\Gamma$, 
one can construct a birational map $x^{\Gamma}: {\rm Conf}_{m+n+2}(\mathbf P^m)\dashrightarrow \boldsymbol{\mathcal X}_{A_m\boxtimes A_n}$ which corresponds to a $\mathcal X$-cluster associated to a quiver $Q_\Gamma$ of type $A_m\boxtimes A_n$ constructed from $\Gamma$ (when $m=1$, one has $A_n\simeq A_1\boxtimes A_n$ and  these maps specialize into the maps $x^T$ associated to triangulations of the $(n+3)$-gon described in \S\ref{SSub:Type-An}).  \sk

The map $p: \boldsymbol{\mathcal A}_{A_m\boxtimes A_n}\rightarrow \boldsymbol{\mathcal X}_{A_m\boxtimes A_n}$ is  dominant if and only if $\ell={\rm gcd}(m+1,n+1)=1$. Otherwise, 
 $\boldsymbol{\mathcal U}_{A_m\boxtimes A_n}$ is a subvariety of 
codimension $\ell-1$ in 
$\boldsymbol{\mathcal X}_{A_m\boxtimes A_n}$ which  it would be interesting to better understand. 
Let $x^{\Gamma_0}:  {\rm Conf}_{m+n+2}(\mathbf P^m)\dashrightarrow \boldsymbol{\mathcal X}_{A_m\boxtimes A_n}$ be the birational identification associated to the 
special minimal bipartite graph $\Gamma_0$ considered in \cite[\S2.2]{Weng} and set 
$ \boldsymbol{\mathcal U}^0_{m,n}=(x^{\Gamma_0})^{-1}(\boldsymbol{\mathcal U}_{A_m\boxtimes A_n})$. 

\begin{question} 
What is $ \boldsymbol{\mathcal U}^0_{m,n}$ as a subvariety of ${\rm Conf}_{m+n+2}(\mathbf P^m)$? Can it be defined geometrically (as for the case of $ \boldsymbol{\mathcal U}^0_{1,3}=(x^{T_0})^{-1}(\boldsymbol{\mathcal U}_{A_3})\subset {\rm Conf}_{6}(\mathbf P^1)=\mathcal M_{0,6}$ which can be constructed using Spence-Kummer configuration $C_{\hspace{-0.05cm}{\cal S}{\cal K}}$ as explained above)?
\end{question}

The geometry of configurations spaces ${\rm Conf}_{m+n+2}(\mathbf P^m)$ is a rich subject. In particular,  projective duality allows to consider the classical `Gale transform' $\mathcal G_{m,n}: {\rm Conf}_{m+n+2}(\mathbf P^m)\rightarrow  {\rm Conf}_{m+n+2}(\mathbf P^n)$ which has been the subject of many (classical or more recent) investigations (see \cite{EP}).

\begin{questions} 
1. How are $ \boldsymbol{\mathcal U}^0_{m,n}$ and $ \boldsymbol{\mathcal U}^0_{n,m}$ related with respect to the Gale transform?  
\sk 

${}^{}$\hspace{0.4cm} 2. When  $m=n$,  is $\boldsymbol{\mathcal U}^0_{m,m}$ formed of $\mathcal G_{m,m}$-fixed configurations ({\it cf.}\,{\rm \cite[Coro.\,8.4]{EP})}? \end{questions}

Some compactifications of the configuration space 
${\rm Conf}_{m+n+2}(\mathbf P^m)$ have been constructed as well, the most well-known being Kapranov's compactification in terms of Chow quotients (see \cite{KT}), denoted here by $\overline{\rm Conf}_{m+n+2}(\mathbf P^m)$.\footnote{Several other compactifications of ${\rm Conf}_{m+n+2}(\mathbf P^m)$ have been constructed recently, some of them only in special cases. The geometry of each and how these compactifications are related are interesting subjects  which have been little studied so far for $m\geq 2$. As interesting recent references on this, we mention \cite{ST} and \cite{GRou}.} When  $m$ is 1,  it coincides with Deligne-Knudsen-Mumford compactification $\overline{\mathcal M}_{0,n+3}$.  In view of Proposition \ref{Prop:UA3-KV}, a natural question  is the following: 


\begin{question} 
Does the closure of $ \boldsymbol{\mathcal U}^0_{m,n}$ in $\overline{\rm Conf}_{m+n+2}(\mathbf P^m)$ belong to 
the cone of effective cycles generated 
by the  boundary cycles\footnote{By definition, a boundary cycle is a cycle in $\overline{{\rm Conf}}_{m+n+2}(\mathbf P^m)$ obtained by intersecting some boundary divisors.} of the corresponding dimension?
\end{question}

Of course, we expect the answer to this question to be negative. 
More generally, one can wonder  to which extend the formalism of cluster algebras can give interesting informations about the algebraic geometry of classical moduli spaces such as the $\mathcal M_{0,n+3}$'s or the  more general configuration spaces ${\rm Conf}_{m+n+2}(\mathbf P^m)$'s. 

\paragraph{} ${}^{}$\hspace{-0.6cm}
\label{Par:UWXAA}
To finish this subsection, let us say a few words about the polylogarithmic ARs of the secondary cluster webs $\boldsymbol{\mathcal U\hspace{-0.04cm}\mathcal W}_{\hspace{-0.04cm} A_m,A_n}$.  
First, such a web is defined only when $A_m\boxtimes A_n$ is of finite type. If one requires that $m\leq n$, this does occur if and only if $(m,n)$ is of the form  $(1,n)$ with $n\geq 2$, or $(2,n)$ with $n=2,3,4$ (the square product $A_2\boxtimes A_n$ is mutation equivalent to $D_4$, $E_6$ and $E_8$ for $n=2,3,4$ respectively). Moreover, in order that $\boldsymbol{\mathcal U\hspace{-0.04cm}\mathcal W}_{\hspace{-0.04cm} A_m,A_n}$ differs from $\boldsymbol{\mathcal X\hspace{-0.04cm}\mathcal W}_{\hspace{-0.04cm} A_m,A_n}$, one has to assume in addition that $m+1$ and $n+1$ are not prime, which needs $(m,n)$ to be in the following list
\begin{equation*}
\label{Eq:listeAmAn}
(1,n)\quad \mbox{with } \, n\geq 2\, \mbox{ odd} \quad \qquad \mbox{ or }\qquad \quad  (m,n)=(2,2)\,. \end{equation*}

 For $(A_m,A_n)=(A_1,A_3)$, the corresponding  web carries trilogarithmic ARs since it is equivalent to Spence-Kummer web. Since $A_2\boxtimes A_2$ is mutation equivalent to $D_4$, one gets that $\boldsymbol{\mathcal U\hspace{-0.04cm}\mathcal W}_{\hspace{-0.04cm} A_2,A_2}$
 is isomorphic to Kummer's web $\boldsymbol{\mathcal W}_{\hspace{-0.04cm}\boldsymbol{\mathcal K}(4)}$ hence carries three tetralogarithmic ARs.  \sk 
 
 Considering these two examples, one might think that for any odd $n\geq 2$, 
 $\boldsymbol{\mathcal U\hspace{-0.04cm}\mathcal W}_{\hspace{-0.04cm} A_n}$ carries 
 polylogarithmic ARs of weight at least 3, which would be very interesting in order to produce (possibly new) polylogarithmic functional identities. However, this is too naive as the consideration of the case when $n=5$ shows. 
Indeed, contrarily to what happens for $A_1\boxtimes A_3$ or $A_2\boxtimes A_3$ both 
$ \boldsymbol{\mathcal X \hspace{-0.05cm} \mathcal W}_{\hspace{-0.04cm}A_5}$  
and 
$ \boldsymbol{\mathcal U \hspace{-0.05cm} \mathcal W}_{\hspace{-0.04cm}A_5}$ have the same number of foliations.\footnote{We have verified that this occurs for $A_7$ as well: both 
$ \boldsymbol{\mathcal X \hspace{-0.05cm} \mathcal W}_{\hspace{-0.04cm}A_7}$  
and 
$ \boldsymbol{\mathcal U \hspace{-0.05cm} \mathcal W}_{\hspace{-0.04cm}A_7}$ are formed by $210$ foliations.}
 Direct computations give us 
$$
\rho^\bullet\Big(\boldsymbol{\mathcal U \hspace{-0.05cm} \mathcal W}_{\hspace{-0.04cm}A_5}\Big)=\big( 
66, 60, 50, 35, 14\big)
\qquad 
{\rm polrk}^\bullet\Big(\boldsymbol{\mathcal U \hspace{-0.05cm} \mathcal W}_{\hspace{-0.04cm}A_5}\Big)=\big( 
121,55\big)
\qquad \mbox{and}\qquad 
{\rm rk}\Big(\boldsymbol{\mathcal U \hspace{-0.05cm} \mathcal W}_{\hspace{-0.04cm}A_5}\Big)\leq 211\, .
$$

In particular, one has  
$
{\rm polrk}\big(\boldsymbol{\mathcal U \hspace{-0.05cm} \mathcal W}_{\hspace{-0.04cm}A_5}\big)=176\leq 
{\rm rk}\big(\boldsymbol{\mathcal U \hspace{-0.05cm} \mathcal W}_{\hspace{-0.04cm}A_5}\big)\leq 211<225=
\rho\big(\boldsymbol{\mathcal U \hspace{-0.05cm} \mathcal W}_{\hspace{-0.04cm}A_5}\big)$
and ${\rm polrk}^3\Big(\boldsymbol{\mathcal U \hspace{-0.05cm} \mathcal W}_{\hspace{-0.04cm} A_5}\Big)=0$
hence 
$\boldsymbol{\mathcal U \hspace{-0.05cm} \mathcal W}_{\hspace{-0.04cm}  A_5}$ is not AMP and above all,  does not carry any polylogarithmic AR of weight higher than 2.

 \subsubsection{\bf About the notion of cluster subvariety.}

The most interesting cluster webs we have obtained so far are those obtained by restriction along the secondary cluster manifolds. In some cases,  this allows to get polylogarithmic ARs of weight $\geq 3$. The reason behind this phenomenon when it occurs is very unclear so far and is,  from our point of view at least, one of the most interesting things to be understood regarding the interplays between cluster algebras and polylogarithms. 

\paragraph{} ${}^{}$\hspace{-0.6cm}
Our experimentations with many concrete examples of webs have shown us that one can get new interesting webs by taking the restriction of a Dynkin-type cluster web 
$\boldsymbol{\mathcal X \hspace{-0.05cm} \mathcal W}_{\hspace{-0.04cm}\Delta}$ along a subvariety $\mathcal V$ of $
 \boldsymbol{\mathcal X}_{\hspace{-0.04cm}\Delta}$ which can be the secondary cluster manifold $\boldsymbol{\mathcal U}_{\hspace{-0.04cm}\Delta}$ or by identifying all the variables associated to the sources or/and to the sinks of the bipartite quiver $\vec{\Delta}$.\footnote{Note that in the case of $D_4$, these two examples actually coincide.} It would be interesting  to have at disposal a richer class of {\it `cluster subvarieties'} $\mathcal V$ along which taking the restriction of ($\mathcal X-$ or $\mathcal Y$-)cluster webs would give rise to new webs, in less variables, but possibly carrying new interesting ARS. 

More explicitly,  given a $\mathcal X$-initial seed $S^0=(\boldsymbol{x},Q)$ (of rank $n\geq 2$), we are looking for a way to associate to each vertex $t$ of $ \mathbb T^n$,  an irreducible affine subvariety $\mathcal V^t\subset \mathcal X\mathbf T^t$, in such a way that the associations $t\leadsto \mathcal V^t$  be  compatible with mutations, 
in the sense that the following condition would hold true:
 \begin{equation*}
{\rm (I)}
\hspace{1.8cm}
\begin{tabular}{l}
{\it for any edge  $ t\stackrel{k}{\line(1,0){15}} t'$ in $\mathbb T^n$, the corresponding mutation  $\mu_k$ }
\\ {\it     is defined at the generic point of $\mathcal V^t$ and $\mu_k(\mathcal V^t)=
\mathcal V^{t'}$.}
\end{tabular}
\hspace{4cm} {}^{}
\end{equation*}

As far as we know,  a  general theory of such cluster subvarieties has not be studied yet.  For the time being, we are aware of only two  specific classes of such varieties: 
\begin{itemize}
\item as explained by Fock and Goncharov in 
\cite{FockGoncharovENS} ({\it cf.}\,Corollary 2.8 and Lemma 2.10 therein), there exists an `exact sequence of cluster varieties' 
$ \boldsymbol{\mathcal A}\stackrel{p}{\longrightarrow} \boldsymbol{\mathcal X}\stackrel{\lambda}{\longrightarrow} H_{\mathcal X} \rightarrow 1$ where 
 $\lambda$ 
 stands for  a canonical fibration 
 of  the $\mathcal X$-cluster variety  over a certain torus $H_{\mathcal X}$.  In particular, the dimension of $H_{\mathcal X}$ is equal to the corank of  $p$ and  
 the secondary cluster variety  $\boldsymbol{\mathcal U}={\rm Im}(p)$ coincides with 
 the fiber of $\lambda$ over $1$.  In the case when $H_{\mathcal X}$ has positive dimension (which corresponds to the case when the rank of any exchange matrix 
of the cluster algebra  
 is not maximal), this torus parametrizes a family of cluster subarieties $\lambda^{-1}(t)\subset \boldsymbol{\mathcal X}$, $t\in H_{\mathcal X}$. These subarieties actually  are the symplectic leaves of $\boldsymbol{\mathcal X}$ for its canonical cluster Poisson structure and can be seen as symplectic deformations of $\boldsymbol{\mathcal U}={\lambda}^{-1}(1)$.
 \sk 
 
 For any Dynkin diagram $\Delta$ and for $\boldsymbol{Z}$ standing either for $\boldsymbol{\mathcal X}$ or for $\boldsymbol{\mathcal Y}$ here, one gets by restriction a family of webs $(\boldsymbol{Z \hspace{-0.05cm} \mathcal W}_{\hspace{-0.04cm}\Delta})\lvert_{\lambda=t}$ for $t\in H_{\mathcal X}$ which is a deformation of 
 $\boldsymbol{\mathcal U \hspace{-0.05cm}Z \hspace{-0.05cm} \mathcal W}_{\hspace{-0.04cm}\Delta}=\big(\boldsymbol{ Z \hspace{-0.05cm} \mathcal W}_{\hspace{-0.04cm}\Delta}\big)\lvert_{\lambda=1}$.  Since some webs  $\boldsymbol{\mathcal U \hspace{-0.05cm}Z \hspace{-0.05cm} \mathcal W}_{\hspace{-0.04cm}\Delta}$ carry interesting polylogarithmic ARs ({\it cf.}\,Theorem \ref{T:classical-cluster-webs}), it is interesting to look at their deformations $(\boldsymbol{Z \hspace{-0.05cm} \mathcal W}_{\hspace{-0.04cm}\Delta})\lvert_{\lambda=t}$: these webs could be related to interesting functional identities as well\footnote{By considering some specific cases, we already have verified that it is indeed the case.} hence deserve to be studied further.
\item 
Another notion which might be interesting for the purpose of constructing interesting cluster subvarieties is that of {`On-Shell subvarieties'} which originated in the theory of scattering amplitudes of high energy physics and  has been related to 
grassmannian varieties and the 
cluster algebras associated to them (see \cite{BFGW} for example).  It would be interesting to understand better these `On-Shell subvarieties' and to see if new interesting webs can be obtained from them. 
\end{itemize}

\paragraph{} ${}^{}$\hspace{-0.6cm}
Considering the examples of  \S\ref{Par:UWXAA}, it appears that the cases when  
the secondary cluster web carries new polylogarithmic ARs (of weight bigger than 2) seem to be in correspondance with those when some of the foliations of the initial cluster  web induce the same foliation after restriction to the secondary cluster manifold. 
These two conditions are not of the same nature (the first is web-theoretic, the second purely geometric), still it would be interesting to better understand how they are related. 
 We formalise this geometric condition in a more general context as follows: let $\mathcal V=\{ \, \mathcal V^t \, \}_{ t\in \mathbb T^n} $ be a cluster subvariety of  a $\mathcal X$-cluster variety $\boldsymbol{\mathcal X}$ as sketched in the previous paragraph
 and $\boldsymbol{\mathcal W}$ a $\mathcal X$-cluster web (for instance 
$\boldsymbol{\mathcal W}$ might be $\boldsymbol{\mathcal X \hspace{-0.05cm} \mathcal W}_{\hspace{-0.04cm}\Delta}$ or   $\boldsymbol{\mathcal Y \hspace{-0.05cm} \mathcal W}_{\hspace{-0.04cm}\Delta}$ when $\Delta$ and $\Delta'$ stand for Dynkin diagrams). 
One defines the restriction of $\boldsymbol{\mathcal W}$ along $\mathcal V$, denoted by 
$ \boldsymbol{\mathcal V \hspace{-0.05cm} \mathcal W}$,  as the web whose first integrals are the restriction on 
$\mathcal V$ of the cluster variables belonging  to $\boldsymbol{\mathcal W}$.   The geometric property we are  interested in is the following: 
\begin{equation*}
{\rm (II)}
\hspace{1.7cm}
\begin{tabular}{l}
{\it  one has ${\rm deg}\big( \boldsymbol{\mathcal V \hspace{-0.05cm} \mathcal W}\big)< {\rm deg}\big( \boldsymbol{\mathcal W} \big)
$, that is  there exist distinct cluster first }\\ {\it    integrals $x$ and $x'$ of $\boldsymbol{\mathcal W}$ such that $dx\wedge dx'=0$ in restriction along $\mathcal V$.}
\end{tabular}
\hspace{3cm}
\end{equation*}

\paragraph{} ${}^{}$\hspace{-0.6cm}
Now let  $\Delta$   be a Dynkin diagram of rank $n\geq 3$. 
The fact that  $\boldsymbol{\mathcal U}_{\hspace{-0.04cm}\Delta}$ is a proper subvariety of 
$\boldsymbol{\mathcal X}_{\hspace{-0.04cm}\Delta}$ is equivalent to the fact that the standard cluster Poisson structure of the latter cluster variety is degenerated.  Moreover in this case, 
the secondary cluster manifold $\boldsymbol{\mathcal U}_{\hspace{-0.04cm}\Delta}$ is a symplectif leaf.  We dont know if it is the case yet, but the fact that $\boldsymbol{\mathcal U}_{\hspace{-0.04cm}\Delta}$ is a symplectic manifold may play a role into the existence of supplementary polylogarithmic ARs carried by $\boldsymbol{\mathcal U \hspace{-0.05cm} \mathcal W}_{\hspace{-0.04cm}\Delta}$. For this reason, we consider the following property for a cluster subvariety $\mathcal V=\{ \, \mathcal V^t \, \}_{ t\in \mathbb T^n} $: 
\begin{equation*}
{\rm (III)}
\hspace{1.8cm}
\begin{tabular}{l}
{\it 
 $\mathcal V$   is a proper symplectic subvariety of $\boldsymbol{\mathcal X}_{\hspace{-0.04cm}\Delta}$.}
\end{tabular}
\hspace{6cm}
\end{equation*}

Let 
$\boldsymbol{\mathcal V}_{\hspace{-0.04cm}\Delta}$ be the class of cluster subvarieties
$\mathcal V=\{ \, \mathcal V^t \, \}_{ t\in \mathbb T^n} $ (as defined above)  satisfying conditions (I), (II) and (III).  Elements therein may give rise to interesting cluster webs.

\begin{questions}
1.  Is $\boldsymbol{\mathcal V}_{\hspace{-0.04cm}\Delta}$ empty? When it is not, what are (at least some of) its elements?\sk 

2.  For $\mathcal V$ in $ \boldsymbol{\mathcal V}_{\hspace{-0.04cm}\Delta}$, 
study the associated web $ \boldsymbol{\mathcal V \hspace{-0.05cm} \mathcal W}$ (ranks, ARs, etc). Are there new examples (in addition of those of Theorem \ref{T:classical-cluster-webs}) of such webs carrying polylogarithmic ARs of weight $\geq 3$? 
\end{questions}

To be honest, we have no idea whether  one may get examples of cluster webs carrying new  interesting polylogarithmic ARs that way. Even worse, we do not know in which way the conditions (I), (II) and (III) could be related to this when it happens. All of this seems very obscure for now, we hope that this will become conceptually clearer one day.

  \subsubsection{\bf Plabic webs of type $\boldsymbol{A_m\boxtimes A_{n}}$.} 
  \label{SS:PlabicWebsAA}
When considering secondary cluster webs, the case of bi-Dynkin type $(A_m,A_n)$ is interesting since  $\boldsymbol{\mathcal U}_{A_m\boxtimes A_n}$ is a proper (symplectic) cluster subvariety of  $\boldsymbol{\mathcal X}_{A_m\boxtimes A_n}$
for many pairs $(m,n)$ (namely, as soon as $m+1$ and $n+1$ are not relatively prime).  On the other hand,  agreeing that $m\leq n$, the cluster algebra with initial quiver $A_m\boxtimes A_n$ has finite type  only when $m=1$ or $m=2$ and $n=2,3, 4$. In our quest of new cluster webs which are interesting since they may carry polylogarithmic ARs, 
it appears natural to consider the restrictions to 
$\boldsymbol{\mathcal U}_{A_m\boxtimes A_n}$ of
webs defined by finite families $\Sigma_{m,n}$ of $\mathcal X$-cluster variables, even 
if the whole set of these variables  is infinite.  Then arises the problem of choosing/constructing 
 finite subsets  
$\Sigma_{m,n}$   of $\boldsymbol{\mathcal X var}\big( 
A_m\boxtimes A_n\big)$ in a relevant and/or natural way.  We describe such a construction below, essentially taken from  \cite[\S2.2]{PS}. 
\sk

A way to construct such a finite family of cluster variables for each type $A_m\times A_n$  relies on the description of some clusters  for the corresponding cluster algebra in terms of certain planar bipartite graphs embedded in a disk, the so-called {\it `plabic graphs'} introduced by  Postinov in his seminal unpublished text \cite{Postnikov} (graphs which are essentially equivalent to the `special minimal bipartite graph' used in \cite{Weng}).  The set 
of such graphs (of the corresponding bi-type, namely $(m+1,m+n+2)$ is finite and  given any two graphs in this family, one can be obtained from the other by a series of `plabic moves'. 
\sk

To each plabic graph $\mathcal P$ is associated a $\mathcal X$-cluster seed $(\boldsymbol{x}^{\mathcal P}, Q^{\mathcal P})$ (where $\boldsymbol{x}^{\mathcal P}$ stands for 
a $mn$-tuple of cluster variables) 
and  a plabic move giving a plabic graph from another is interpreted as 
 a mutation in terms of the corresponding clusters. 
Hence considering all the $\mathcal X$-cluster variables associated to the `plabic seeds' (which are those associated to plabic  graphs), one ends with a finite family of $\mathcal X$-cluster variable which are first integrals of the {\bf ($\boldsymbol{\mathcal X}$-)cluster  plabic web
of type} $\boldsymbol{A_m\boxtimes A_n}$, denoted by  $\boldsymbol{\mathcal X \hspace{-0.04cm}\mathcal W}_{\hspace{-0.05cm} A_m\boxtimes A_n}^{plabic}$. 
When $m=1$,  the plabic web of type  $ A_1\boxtimes A_n$ actually coincides with $\boldsymbol{\mathcal X \hspace{-0.04cm}\mathcal W}_{\hspace{-0.05cm} A_n}$ hence has already be considered before. But as soon $m,n\geq 2$, one gets new cluster webs despite the fact that 
the cluster algebra of type $A_m\boxtimes A_n$ is of infinite type, about which natural questions can be asked.

\begin{questions} 
\begin{enumerate}
\item  For $m,n$ both bigger than or equal to 2, study the plabic web 
$\boldsymbol{\mathcal X \hspace{-0.04cm}\mathcal W}_{\hspace{-0.05cm} A_m\boxtimes A_n}^{plabic}$ (degree, virtual and polylogarithmic ranks, ARs, etc). 

\item When ${\rm gcd}(m+1,n+1)\geq 2$, same questions for 
$\boldsymbol{\mathcal U \hspace{-0.04cm}\mathcal W}_{\hspace{-0.05cm} A_m\boxtimes A_n}^{plabic}$, the restriction of 
$\boldsymbol{\mathcal X \hspace{-0.04cm}\mathcal W}_{\hspace{-0.05cm} A_m\boxtimes A_n}^{plabic}$
along $\boldsymbol{\mathcal U}_{A_m\boxtimes A_n}$. Does any of these `secondary plabic cluster webs' carry new polylogarithmic identities of weight strictly bigger than 2? 
\end{enumerate}
\end{questions}


\subsection{\bf  Toric deformation of $\boldsymbol{\mathcal X}$-cluster varieties from the point of view of webs}

Clusters varieties are interesting objects from the  perspective of mirror symmetry, as Fock and Goncharov have conjectured that the $\mathcal A$- and the $\mathcal X$-cluster varieties 
(for a given cluster algebra) are mirror of each other.  Motivated by this conjecture\footnote{We refer to \cite{GHKK} for more background as well as results about this conjecture.}, several authors have considered the problem of constructing toric deformations of a given cluster variety. More precisely, given say a $\mathcal X$-cluster variety $\mathcal X$ with Fock-Goncharov's special completion $\widehat{\mathcal X}$ (see \cite{FockGoncharovXInfinity} and also \S\ref{SS:FockGoncharo-XInfinity}), 
the problem is about building  a family $\widehat{\mathscr X}\rightarrow \mathbf A^n$ over an affine space such that the fiber $\widehat{\mathscr X}_t$ over a generic $t\in  \mathbf A^n$ be isomorphic to $\widehat{\mathcal X}$, with the central (possibly singular) fiber $\widehat{\mathscr X}_0$ toric and satisfying also some other conditions.  \sk 

In case $\mathcal X=\mathcal X_\Delta$ is a cluster variety of finite type associated to a Dynkin diagram $\Delta$, each generic fiber $\widehat{\mathscr X}_t$ carries its own model of 
the cluster web $\boldsymbol{\mathcal X \hspace{-0.04cm}\mathcal W}_{\hspace{-0.05cm} \Delta}$, denoted by $\boldsymbol{\mathcal X \hspace{-0.04cm}\mathcal W}_{\hspace{-0.05cm} \Delta,t}$. We obtain a family of isomorphic webs parametrized by a Zariski-open subset of $\mathbf A^n$.
\begin{questions} 
\begin{enumerate}
\item   What about  the degeneration when $t$ goes at the origin of the $\boldsymbol{\mathcal X \hspace{-0.04cm}\mathcal W}_{\hspace{-0.05cm} \Delta,t}$'s: does this degeneration exist as a web? 
\item  In case this degeneration indeed exists, what are its properties as a web (ranks, ARs, etc)? 
What does this degeneration tell us about the original web $\boldsymbol{\mathcal X \hspace{-0.04cm}\mathcal W}_{\hspace{-0.05cm} \Delta}$?
\end{enumerate}
\end{questions}

As an explicit example that can be found in the literature, we mention the case when $\Delta=A_2$ considered in \cite{Boss}: in thise case, we have a deformation of $\mathcal X_{A_2}$ of dimension 2, and from \cite[Figure 1]{Boss}, we get that  in some affine coordinates  $(X_1,X_2)$,  one has 
$$
\boldsymbol{\mathcal X \hspace{-0.04cm}\mathcal W}_{\hspace{-0.05cm} \Delta,t}=
\boldsymbol{\mathcal W}\bigg( \,  X_1
\, , \,  X_2
\, , \, 
\frac{t_1+X_1}{X_2}
\, , \, 
\frac{1+t_2X_2}{X_1}
\, , \, 
\frac{t_1+X_1+t_1t_2X_2}{X_1X_2}
\, \bigg)
$$ 
for $t=(t_1,t_2)\in \mathbf A^2$ generic. From this explicit expression, we deduce that the degeneration at the origin $0\in \mathbf A^2$ exists and that it is the `logarithmic' 3-web 
$\boldsymbol{\mathcal X \hspace{-0.04cm}\mathcal W}_{\hspace{-0.05cm} A_2,0}=
\boldsymbol{\mathcal W}\Big( \,  X_1
\, , \,  X_2
\, , \, 
{X_1}/{X_2}
\, \Big)
$.  If this explicitly answers the first of the questions above, 
it is not clear how the quite elementary degeneration 
 $\boldsymbol{\mathcal X \hspace{-0.04cm}\mathcal W}_{\hspace{-0.05cm} A_2,0}$ could tell us anything about the original dilogarithmic web $\boldsymbol{\mathcal X \hspace{-0.04cm}\mathcal W}_{\hspace{-0.05cm} A_2}$.  \sk
 
 It would be interesting to work out other cases explicitly in order to get  other  examples of 
 degenerations $\boldsymbol{\mathcal X \hspace{-0.04cm}\mathcal W}_{\hspace{-0.05cm} \Delta,0}$ (under the assumption that these latter exist), hoping that these will be interesting.

 \subsection{\bf In connection with quantization and  categorification of cluster algebras} 
 \label{SS:In-connection-with-quantization}

There are interesting links between cluster algebras, quantification and category theory, which have given rise to a rich literature. For instance, quantization of cluster varieties is 
considered in \cite{FockGoncharovENS} and categorification of cluster algebra is one of the main tools used in \cite{KellerAnnals} to prove the periodicity property of $Y$-systems of bi-Dynkin type.\sk 

An interesting point is that many cluster dilogarithmic identities can be quantized as well and that in some cases, such quantized identities admit an interpretation within a certain categorification of the underlying cluster algebra.  For a geometer primarily interested on webs, this suggests many questions on possible generalizations of key notions of web geometry (the notion of web itself, but also those of abelian relation,  of rank, etc.) to other mathematical frameworks.  In return, considering some well-known results on the quantum/categorical side, one can wonder about possible extensions of classical web geometry to webs formed by infinitely many 
 foliations.  
\mk


%
%

Our main references for the material discussed in this subsection are
\cite{FockGoncharovENS}, 
\cite{KellerCTQDI} and 
\cite{KN}. More specific references are also given below.

\subsubsection{\bf The quantum pentagon identity.} 
 The most famous dilogarithmic identity is the `pentagon identity', namely Abel's five terms identity $({\sf R}_{A_2})$. 
 Quantum versions of it have been considered in many papers, from different  points of view. 
 Among many other references, we refer to  \cite{FaddeevKashaev}, \cite{Volkov2}, 
\cite{FockGoncharovENS}, \cite{KN} for the material discussed in this paragraph.
\sk 
 
There are (at least) two quantum versions of the dilogarithm and correspondingly two different quantum versions of 
 the  pentagon identity ({\it cf.}\,\cite[\S1.3]{FockGoncharovINV} or \cite{KN} and the references therein):
 \begin{enumerate}
 \item the first  {\it `quantum dilogarithm'} is a  formal series with coefficients in $\mathbf Q(q)$ where $q$ is a complex parameter known as the `quantum parameter'. There exist slightly different versions of it and correspondingly several essentially equivalent quantum pentagon identities. These  are 
multiplicative  identities between formal series in two non-commutative variables;
 \sk 
 \item there is also a {\it `non-compact quantum dilogarithm'} defined only for quantum parameters $q$ of modulus 1, by means of a planar integral along an infinite contour.  The corresponding quantum pentagon identity is equivalent to the fact that a certain operator acting on a Hilbert space is periodic, of period 5. 
\end{enumerate}
 
 The second version, of analytic and geometric nature,  is more in the spirit of classical quantum mechanics than the former which has a more formal and algebraic nature.  
How both  are related is not quite clear to us but  according  the abstract of \cite{Faddeev},  the first quantum pentagon identity can be derived from the non-compact one, which suggests that the latter may be  the most fundamental of both.
Despite this, we will leave aside all the non-compact side of the theory  to focus on the other one, which is more algebraic, easier to deal with and more easily connected with the theory of cluster algebras.
\mk 
 
 
As mentioned  above, there are several versions of the (compact) quantum dilogarithm and it is a bit tedious to specify how they are all connected. To simplify our exposition, we will allow ourselves to be not very rigorous at several places.  \sk

The most popular version of a (compact) quantum dilogarithm is the function 
$$\boldsymbol{\Psi}_q(x)
 = \sum_{n=0}^{+\infty}\frac{ q^{n^2} x^n } { \big(q-q^{-1}\big)\big(q^2-q^{-2}\big)\cdots \big(q^{n}-q^{-n}\big)}
 $$ 
 which, as a series of $x\in \mathbf C$,   converges for any $q$ such that $\lvert q\lvert <1$.  Its name is justified by its asymptotic behavior: assuming in addition that $\lvert x\lvert<1$, as 
 $q\rightarrow 1^-$  one has: 
\begin{equation}
\label{Eq=Psi-q-1}
\boldsymbol{\Psi}_q(x)\sim \exp\left( - \frac{ {\bf L}{\rm i}_2(-x)} { {\rm Log}(q^2) } 
\right)
\,.
\end{equation} 

Now let  $Y_1,Y_2$ be two  quasi-commuting variables 
  such that   $Y_1Y_2=q^2Y_2Y_1$. 
 Then setting 
 $$
 U_{1}=Y_1\, , \hspace{0.1cm} U_{2}=Y_2(1+qY_1)
 \, , \hspace{0.1cm} U_{3}=\big(1+qY_2+Y_1Y_2 \big)^{-1} Y_1
 \, , \hspace{0.1cm} U_{4}=q(1+qY_2)^{-1}Y_2Y_1
\hspace{0.2cm} \mbox{and}\hspace{0.2cm}
 U_{5}=Y_2\, ,  $$
  the following quantum multiplicative  identity  holds true: 
\begin{equation}
\label{Eq:Quantum-Pentagon}
\boldsymbol{\Psi}_q\big( U_{5} \big)^{-1}
\boldsymbol{\Psi}_q\big( U_{4}
\big)^{-1}
\boldsymbol{\Psi}_q\big( 
U_{3} \big)^{-1}
\boldsymbol{\Psi}_q\big( 
U_{2} \big)
\boldsymbol{\Psi}_q\big( 
U_{1} \big)=1\, . 
\end{equation}
 This relation is called the {\it `quantum pentagon identity'} (in `universal form' according to the terminology of \cite{KN}) which can be seen as a non-commutative deformation of the classical 5-terms dilogarithmic identity. Indeed, then 
 thanks to \eqref{Eq=Psi-q-1},  one obtains that up to some factor depending only on $q$,  the asymptotic development of the logarithm of
  the  LHS 
  of  \eqref{Eq:Quantum-Pentagon} is given formally by 
  \begin{equation}
\label{Eq:Dequantif-Quantum-Pentagon}
 {\bf L}{\rm i}_2\big(y_2\big)+
 {\bf L}{\rm i}_2\left( 
 \frac{y_1y_2}{1+y_2}
 \right)+
 {\bf L}{\rm i}_2
 \left( 
 \frac{y_1}{1+y_2+y_1y_2}
 \right)-
 {\bf L}{\rm i}_2\Big( 
 y_2\big(1+y_1\big)
 \Big)-
 {\bf L}{\rm i}_2(y_1)
\end{equation}
where $y_1,y_2$ are the classical limits of $Y_1$ and $Y_2$ hence are standard (commuting) 
 variables.  We arguments of the bilogarithm in the expression above are precisely the $\mathcal X$-cluster variables of type $A_2$ hence 
a  naive guess would be that 
\eqref{Eq:Dequantif-Quantum-Pentagon} is identically zero
since it is formally  obtained  from taking the asymptotic of the logarithm of 
the LHS of \eqref{Eq:Quantum-Pentagon}.  Actually, this is not the case 
 but (for some reasons investigated in \cite{KN}), 
 it becomes so if you replace the  classical bilogarithm  $ {\bf L}{\rm i}_2$  by the cluster dilogarithm ${\sf R}(x)=(1/2)\int_0^x \big( \log(1+u)/u-\log(u)/(1+u)\big) du$ since the 
 one has 
 $$
 ({\sf R}_{A_2})
 \hspace{2cm}
 {\sf R}\big(y_2\big)+
 {\sf R}\left( 
 \frac{y_1y_2}{1+y_2}
 \right)+
 {\sf R}
 \left( 
 \frac{y_1}{1+y_2+y_1y_2}
 \right)-
{\sf R}\Big( 
 y_2(1+y_1)
 \Big)-
{\sf R}\big(y_1\big)
\equiv 0\, .  
 \hspace{2cm} {}^{}
$$

 The formal derivation  sketched just above of the classical pentagon identity $({\sf R}_{A_2})$ from \eqref{Eq:Dequantif-Quantum-Pentagon} actually can be mathematically justified (but this is not trivial, {\it cf.}\,the appendix of \cite{KN})  which justifies the name `quantum pentagon equation' for \eqref{Eq:Dequantif-Quantum-Pentagon} which hence appears as a non-commutative deformation of the classical pentagon identity. 
 \mk 
 
The arguments $u_i=\lim_{q\rightarrow 1} U_i$ (for $i=1,\ldots,5$) of the cluster dilogarithm ${\sf R}$ in  $({\sf R}_{A_2})$ are the cluster first integrals of the 
 cluster web $ {\boldsymbol{\mathcal X\hspace{-0.04cm}\mathcal W}}_{\hspace{-0.04cm} A_2}$ on the initial cluster torus ${\mathbf T}_0={\rm Spec} \big( \mathbf C\big[ y_1^{\pm 1},
 y_2^{\pm 1}\big] \big) $.  These are the classical limits of the five arguments $U_i$ of $\boldsymbol{\Psi}_q$ in \eqref{Eq:Quantum-Pentagon}. They are elements of the non-commutative `initial cluster torus $\boldsymbol{q}{\mathbf T}_{0}$'\footnote{The `quantum torus $q{\mathbf T}_{0}$' does not exists {\it per se} as a geometric entity. What does exist however, is its putative  ring of regular functions 
 $\mathcal O(\boldsymbol{q}{\mathbf T}_{0})$
  which is the algebra over the field $\mathbf C(q)$, of non-commutative polynomials in $Y_1$, $Y_2$ and their inverses, subject to the relations $Y_1Y_2
 =q^2 Y_2Y_1$.} on which lives
$ \boldsymbol{q\mathcal X\hspace{-0.04cm}\mathcal W_{\hspace{-0.04cm} A_2}}$,  
    the {\bf `quantum $\boldsymbol{\mathcal X}$-cluster web of type $\boldsymbol{A_2}$'}, which is nothing but the 
 collection of the $U_i$'s appearing in \eqref{Eq:Quantum-Pentagon} (in the reverse order): 
$$
{\boldsymbol{q}\boldsymbol{\mathcal X\hspace{-0.04cm}\mathcal W}}_{\hspace{-0.04cm} A_2}=
{\boldsymbol{\mathcal W}}\bigg( \, 
Y_1\, , \, 
Y_2(1+qY_1)
\, , \, 
\big(1+qY_2+Y_1Y_2 \big)^{-1} Y_1
\, , \, 
q(1+qY_2)^{-1}Y_2Y_1
\, , \,  Y_2
\, \bigg)\, . 
$$
 
One would like to see this collection of five quantum cluster first integrals as a quantum deformation of its classical limit which is the usual cluster 5-web of type $A_2$: 
$$
{\boldsymbol{\mathcal X\hspace{-0.04cm}\mathcal W}}_{\hspace{-0.04cm} A_2}=\lim_{q\rightarrow 1^-} {\boldsymbol{q}\boldsymbol{\mathcal X\hspace{-0.04cm}\mathcal W}}_{\hspace{-0.05cm} A_2} . 
$$

At this point, an interesting although perhaps very naive attempt would be to investigate whether  some of the web-theoretic remarkable properties of ${\boldsymbol{\mathcal X\hspace{-0.04cm}\mathcal W}}_{\hspace{-0.04cm} A_2}$ can be quantized, in the sense that given a property $\mathcal P_1$ of this 5-web, there exists a property $P_q$ satisfied by ${\boldsymbol{q}\boldsymbol{\mathcal X\hspace{-0.04cm}\mathcal W}}_{\hspace{-0.04cm} A_2}$ say for $q$ such that $\lvert q\lvert <1$, in  such a way that $\mathcal P_1$ be the limit  of $P_q$ (in some sense to be made precise/rigorous). \mk 

Here is a non-exhaustive list of questions that may be asked about ${\boldsymbol{q}\boldsymbol{\mathcal X\hspace{-0.04cm}\mathcal W}}_{\hspace{-0.04cm} A_2}$ in this regards. 
\begin{questions}
\begin{enumerate}
\item  In addition to 
 $({\sf R}_{A_2})$, the 5-web ${\boldsymbol{\mathcal X\hspace{-0.04cm}\mathcal W}}_{\hspace{-0.04cm} A_2}$ carries 9 linearly independent logarithmic ARs (each of which can be chosen with only three terms). Does any one of these  admit 
a non-commutative deformation for the corresponding 3-subweb of ${\boldsymbol{q}\boldsymbol{\mathcal X\hspace{-0.04cm}\mathcal W}}_{\hspace{-0.04cm} A_2}$?\sk 

For a concrete example, let us take the logarithmic abelian relation   ${\rm Log}\big(  (1+y_2) /y_2\big)+
{\rm Log}\big(y_1y_2/(1+y_2)\big)-
{\rm Log}(y_1)=0$ of the  the 3-subweb $\boldsymbol{\mathcal W}(u_5,u_4,u_1)$. In multiplicative form, it can be written $\big(  (1+u_5) /u_5\big)\cdot u_4\cdot (u_1)^{-1}=1$, 
a rational identity which can be seen as  the classical limit (when $q\rightarrow 1^-$) of the quantum relation 
$$\Big(  U_5^{-1}(1+U_5) \Big)\cdot U_4\cdot (U_1)^{-1}=q$$
 which holds true in the corresponding algebra $\mathcal O(\boldsymbol{q}{\mathbf T}_{0})$ of non-commutative functions. We do believe that, similarly, any logarithmic AR of ${\boldsymbol{\mathcal X\hspace{-0.04cm}\mathcal W}}_{\hspace{-0.04cm} A_2}$ can be quantized in an elementary way so that the answer to the  question above is 'yes'.
\item  The importance of $({\sf R}_{A_2})$ among the ARs of ${\boldsymbol{\mathcal X\hspace{-0.04cm}\mathcal W}}_{\hspace{-0.04cm} A_2}$ is that all the other logarithmic ARs of this web can be derived from it (by action of the monodromy or by derivation, see \ref{Par:Dilog-Generation-LogAR}).  If 
the logarithmic ARs of  ${\boldsymbol{\mathcal X\hspace{-0.04cm}\mathcal W}}_{\hspace{-0.04cm} A_2}$ can all be quantized, can these quantizations be derived from 
\eqref{Eq:Quantum-Pentagon}, and if so, in which manner?
\item Can not only the ARs of ${\boldsymbol{\mathcal X\hspace{-0.04cm}\mathcal W}}_{\hspace{-0.04cm} A_2}$, but the `notion of AR' itself be quantized in a relevant way? 
The most naive attempt for a notion of `quantum abelian relation' 
would be to consider 5-tuples $(\boldsymbol{\Phi}_i)_{i=1}^5\in 
\mathbf C(q)[[z]]^5$ satisfying the following relation in $\mathcal O(\boldsymbol{q}{\mathbf T}_{0})$: 
\begin{equation}
\label{Eq:Quantum-AR-A2}
\boldsymbol{\Phi}_5\big( U_{5} \big)\, 
\boldsymbol{\Phi}_4\big( U_{4}
\big) \,
\boldsymbol{\Phi}_3\big( 
U_{3} \big) \,
\boldsymbol{\Phi}_2\big( 
U_{2} \big)\,
\boldsymbol{\Phi}_1\big( 
U_{1} \big)=1\, . 
\end{equation}
Can the space $\boldsymbol{\mathcal A}\big( {\boldsymbol{q\mathcal X\hspace{-0.04cm}\mathcal W}}_{\hspace{-0.04cm} A_2}\big)$ of such 5-tuples  be endowed with some (possibly non-commuta\- -tive) algebraic structure which would make it appear as a natural quantization of the space of abelian relations 
$\boldsymbol{\mathcal A}\big( {\boldsymbol{\mathcal X\hspace{-0.04cm}\mathcal W}}_{\hspace{-0.04cm} A_2}\big)$ 
 of ${\boldsymbol{\mathcal X\hspace{-0.04cm}\mathcal W}}_{\hspace{-0.04cm} A_2}$?
 \item   Since $\mathcal O(\boldsymbol{q}{\mathbf T}_{0})$ is not commutative, the order of the terms $\Phi_i(U_i)$'s in \eqref{Eq:Quantum-AR-A2}   is important. Taking another order gives another  notion of `quantum ARs' which may or may not be relevant. Such alternative notions of `quantum ARs' would deserve to  be investigated as well. For example: is there  $(\boldsymbol{\Phi}_i)_{i=1}^5$  as above satisfying 
 $\boldsymbol{\Phi}_5\big( U_{5} \big)\, 
\boldsymbol{\Phi}_4\big( U_{2}
\big) \,
\boldsymbol{\Phi}_3\big( 
U_{3} \big) \,
\boldsymbol{\Phi}_4\big( 
U_{2} \big)\,
\boldsymbol{\Phi}_1\big( 
U_{1} \big)=1$  from which $({\sf R}_{A_2})$ can be derived in the classical limit?
 \end{enumerate}
\end{questions}


   

\subsubsection{\bf Quantized cluster dilogarithmic identities and quantized cluster webs I.} 
 \label{Parag:QC-I}
An  interesting aspect of what has just been exposed above is that many things  actually generalize to any $\mathcal X$-cluster web associated to a cluster period.  This follows from results due to Reineke \cite{Reineke} for the $\mathcal Y$-periods associated to  (simply-laced) Dynkin diagrams and has been extended to (bi-simply-laced)  bi-Dynkin $\mathcal Y$-cluster periods by Keller in \cite{KellerCTQDI}.  
The case of a general cluster period discussed below is due to Kashaev and Nakanishi \cite{KN} (see also \cite{Nagao}).
\mk

We use the notation of \S\ref{Par:NakanishiDilogarithmicIdentity-Period} above. 
Let $\boldsymbol{\nu}=(\nu_1,\ldots,\nu_L)$ be a $\mathcal X$-cluster period of length $L$ of a cluster algebra with initial seed $(\boldsymbol{y}, B)$ with $\boldsymbol{y}=(y_i)_{i=1}^n$ and $B=(b_{ij})_{i,j=1}^n $ a skew-symmetric initial exchange matrix. Denote by $x_\ell=x_\ell(\boldsymbol{\nu})$ the corresponding cluster variables and  $\epsilon_\ell\in \{\pm 1\}$ the  associated tropical signs (with $\ell=1,\ldots,L$).  Then according to Nakanishi's theorem \cite{Nakanishi}, the following identity holds true identically:
\begin{equation*}
\Big(\boldsymbol{{\sf R}_{\boldsymbol{\nu}}}\Big)
\hspace{5cm}
\sum_{\ell=1}^L \epsilon_\ell \, {\sf R}\Big(  {x_\ell}^{\epsilon_\ell}
\Big) = 0 \, . 
\hspace{7cm}{}^{}
\end{equation*}

The initial $\mathcal X$-cluster torus $\mathbf T_0={\rm Spec}\big( \mathbf C[\boldsymbol{u}^{\pm 1}]\big)$ can be deformed into the quantum torus $\boldsymbol{q\mathbf T}_0$ associated to 
non-commutative variables $Y_i$'s which are quantum deformations  of the $y_i$'s satisfying the quasi-commuting relations 
$Y_iY_j=q^{2b_{ji}}Y_jY_i $.  The cluster period $\boldsymbol{\nu}$ and the corresponding cluster variables $x_\ell$ (for $\ell=1,\ldots,L$) can be quantized in a canonical way and give rise to a quantum cluster period with  associated $L$-tuple of quantum cluster variables $(X_\ell)_{\ell=1}^L$ (see \cite[Prop.\,3.4]{KN}) such that the following quantum dilogarithmic identity holds true (according to \cite[Cor.\,3.7]{KN}): 
\begin{equation*}
\Big(\boldsymbol{q{\sf R}_{\boldsymbol{\nu}}}\Big)
\hspace{3cm}
\boldsymbol{\Psi}_q\Big({X_L}^{\epsilon_L} \Big) ^{\epsilon_L} 
\cdots 
\boldsymbol{\Psi}_q\Big({X_{2}}^{\epsilon_2} \Big) ^{\epsilon_2}\,
\boldsymbol{\Psi}_q\Big({X_{1}}^{\epsilon_1} \Big) ^{\epsilon_1}
=1\, .
\hspace{5cm}{}^{}
\end{equation*}
Each classical cluster variable $x_\ell$ is obtained from $X_\ell$ by taking the limit $q\rightarrow 1$ and by using more sophisticated arguments ({\it cf.}\,Appendix A in \cite{KN}) one proves that 
$\big(\boldsymbol{{\sf R}_{\boldsymbol{\nu}}}\big)$ can be derived from  the quantum identity 
$\big(\boldsymbol{q{\sf R}_{\boldsymbol{\nu}}}\big)$
in the classical limit. \sk 

As in the $A_2$-case, discussed above, one defines the quantum cluster web associated to $\boldsymbol{\nu}$ as the set of $X_\ell$'s ordered according their 
 order of appearance in $\big(\boldsymbol{q{\sf R}_{\boldsymbol{\nu}}}\big)$: 
$$
{\boldsymbol{q\mathcal W}}_{\hspace{-0.04cm} \boldsymbol{\nu}}={\boldsymbol{\mathcal W}}
\Big( \,  X_L
\, , \,  X_{L-1}
\, , \, \ldots 
\, , \, 
\, , \,  X_2
\, , \,  X_1\, 
\Big)\, .
$$
(Rigorously,  the quantum cluster web ${\boldsymbol{q\mathcal W}}_{\hspace{-0.04cm} \boldsymbol{\nu}}$ is the ordered $L$-tuple $(X_L,\ldots,X_1)$. 
The `${\boldsymbol{\mathcal W}}$' appearing in the LHS is superfluous but is here to emphasize that our purpose is to deal with this ordered tuple of quantum cluster variables as if it was a usual web).  
\sk 

Clearly, one has $
{\boldsymbol{\mathcal W}}_{\hspace{-0.04cm} \boldsymbol{\nu}}=
{\boldsymbol{\mathcal W}}\big( \, x_{L-l}\, 
\lvert \, l=0,\dots,L-1
\big)=
\lim_{q\rightarrow 1^-} {\boldsymbol{q}\boldsymbol{\mathcal W}}_{\hspace{-0.04cm} \boldsymbol{\nu}} 
$ in a straightforward way hence one can see 
${\boldsymbol{q\mathcal W}}_{\hspace{-0.04cm} \boldsymbol{\nu}}$ as a quantum deformation of the cluster web ${\boldsymbol{\mathcal W}}_{\hspace{-0.04cm} \boldsymbol{\nu}}$ hence, 
similarly to what we did in the preceding paragraph in the $A_2$-case, several interesting 
web-theoretic questions can be asked about the former quantum web:

\begin{questions}
\label{Quest:qW-ARs}
\begin{enumerate}
\item  In addition to 
 $(\boldsymbol{{\sf R}_{\boldsymbol{\nu}}})$, can the other polylogarithmic ARs of 
 ${\boldsymbol{\mathcal W}}_{\hspace{-0.04cm} \boldsymbol{\nu}}$ be obtained  as limits from certain quantum identities which could be seen as quantum ARs for  ${\boldsymbol{q \mathcal W}}_{\hspace{-0.04cm} \boldsymbol{\nu}}$? And what about the other ARs of ${\boldsymbol{\mathcal W}}_{\hspace{-0.04cm} \boldsymbol{\nu}}$, not of polylogarithmic type? Do they admit natural or even canonical `quantum deformations'?
\item Does it exist a relevant notion of `quantum AR' for ${\boldsymbol{q\mathcal W}}_{\hspace{-0.04cm} \boldsymbol{\nu}}$? 
What about the most immediate one coming to mind, namely the one of  $L$-tuple $(\boldsymbol{\Phi}_\ell)_{\ell=1}^L\in \mathbf C(q)[[z]]^L$, such that 
$$\boldsymbol{\Phi}_L\Big({X_L} \Big)
\cdots 
\boldsymbol{\Phi}_1\Big({X_{1}} \Big) 
=1\, ?$$
\item
If the answer to the previous questions is affirmative, does the space $\boldsymbol{\mathcal A}\big( {\boldsymbol{q\mathcal W}}_{\hspace{-0.04cm} \boldsymbol{\nu}}\big)$ of quantum ARs carry a natural (possibly non-commutative) algebraic structure which would make it appear as a natural quantization of the space of abelian relations 
$\boldsymbol{\mathcal A}\big(
{\boldsymbol{\mathcal W}}_{\hspace{-0.04cm} \boldsymbol{\nu}}
\big)$ 
 of ${\boldsymbol{\mathcal W}}_{\hspace{-0.04cm} \boldsymbol{\nu}}$?
 \item   If 4. can be answered in a meaningful way, can the notion of `quantum rank' be defined? If so, what about the notion of being AMP? Does the expression `quantum cluster web with AMP rank' make sense?
 \end{enumerate}
\end{questions}

At this point, one can look for generalizations of the previous questions for 
cluster webs of 
other kinds.   
Below, we discuss  two interesting concrete examples: the first 
is the $\mathcal X$-cluster web ${\boldsymbol{\mathcal X\hspace{-0.04cm}\mathcal W}}_{\hspace{-0.04cm} B_2}$ equivalent to Newman's 6-web. Next 
we consider the case of the secondary cluster web ${\boldsymbol{\mathcal U\hspace{-0.04cm}\mathcal W}}_{\hspace{-0.04cm} A_3}$ which is a cluster model of Spence-Kummer 9-web.

\subsubsection{\bf Quantized cluster dilogarithmic identities and quantized cluster webs II.} 

%
%
%
%
%

%


The first examples 
beyond the skew-symmetric case
to be considered regarding the questions stated in the previous paragraph 
 certainly are those associated with cluster algebras of finite type $B_2$ or $G_2$. Both can be investigated quite  explicitly but to save space, we will only discuss the $B_2$ case. 

\paragraph{\bf About the quantum cluster web of type $\boldsymbol{B_2}$.} 
\label{SubPar:QuantumClusterWebB2}
From \cite[\S3.3]{FockGoncharovENS}, we know that the Dynkin cases $A_2$, $B_2$ and $G_2$ correspond to the initial exchange matrix  
$\Big[
\hspace{-0.2cm}
\scalebox{0.9}{
\begin{tabular}{c}
0 \hspace{0.05cm} $-1$\vspace{-0.1cm}\\
$c$  \hspace{0.16cm} 0
\end{tabular}} \hspace{-0.15cm} \Big]$
with $c=1,2$ and $c=3$ respectively and that the corresponding classical and quantum cluster mutations correspond respectively to the following classical and quantum recursion formulas (with $m\in \mathbf Z$):
$$
x_{m-1}x_{m+1}=
\begin{cases} \, 
(1+x_m)  \\
\, (1+x_m)^c 
\end{cases}
\quad \mbox{and} \qquad 
X_{m-1}X_{m+1}=
\begin{cases} \, 
\big(1+q^cX_m\big)  \\
\, \big(1+qX_m\big)\big(1+q^3 X_m\big)\cdots \big(1+q^{2c-1}X_m\big)
 \, , 
\end{cases}
$$
where the first line concerns the case when $m$ is even and the second the case when $m$ is odd.
\mk 

The case of $A_2$ ({\it i.e.}\,when $c=1$) has been considered above, here we deal with the one when $c=2$ (case $B_2$) which is even more interesting in what regards quantization (of the webs and of their ARs) since  if ${\boldsymbol{\mathcal X\hspace{-0.04cm}\mathcal W}}_{\hspace{-0.04cm} B_2}$ has maximal rank, not all its ARs are polylogarithmic.   
\sk 

From 
the proof of 
\cite[Lemma 3.9]{FockGoncharovENS}, we deduce the following expression for the quantum cluster web of type $B_2$: 
\begin{align}
\label{Eq:q-XW-B2}
{\boldsymbol{q\mathcal X\hspace{-0.04cm}\mathcal W}}_{\hspace{-0.04cm} B_2}=
&\, 
{\boldsymbol{\mathcal W}}\bigg( \, X_1
\, , \,  X_2^{-1}\big(1+q^{-1}X_1\big) \big(1+q^{-3}X_1\big) 
\, , \,
X_1^{-1}\big(1+q^2X_2 \big) 
\, , \, X_4
\, ,  \,  X_5
\, ,  \,
 X_2
\, \bigg)
\end{align}
with $X_4=X_2^{-1}X_1^{-2}\Big(
X_1+q(1+q^6X_2\big)
 \Big)
\Big(
X_1+q^3\big(1+q^2X_2\big)\Big)$ and a similar formula for $X_5$.\sk 

  Setting `$x_i=\lim_{q\rightarrow 1} X_i$' for $i=1,2$ for the two dequantized variables, one verifies 
   that 
$$
\lim_{q\rightarrow 1} 
{\boldsymbol{q\mathcal X\hspace{-0.04cm}\mathcal W}}_{\hspace{-0.04cm} B_2}=
{\boldsymbol{\mathcal W}}\Big( \, x_1
\, , \, \frac{(1+x_1)^2}{x_2}
\, , \,  \frac{1+x_2}{x_1}
\, , \,  \frac{(1+x_1+x_2)^2}{x_1^2x_2}
\, , \,   \frac{(1+x_1)^2+x_2}{x_1x_2}
\, , \, x_2
\, \Big)={\boldsymbol{\mathcal X\hspace{-0.04cm}\mathcal W}}_{\hspace{-0.04cm} B_2} \, .
$$
We recover precisely the 6-web studied in \S\ref{SS:cluster-webs:B2} where an explicit basis of its space of abelian relations is given.  \sk

Contrarily to the case of the dilogarithmic identity in the skew-symmetric case, there is no general result about quantum versions of dilog identities associated to cluster periods when the initial exchange matrix is just assumed to be skew-symmetrizable.  
 A conceptual proof that the quantum identities hold in the skew-symmetric case can be given in the realm of additive categorification of cluster algebras ({\it cf.}\,\cite{KellerCTQDI} or below for more details).  At the time writing, there is no general theory of categorification for skew-symmetrizable cluster algebras but one might expect that most of the properties obtained via categorical methods  in the skew-symmetric case actually also holds true for a skew-symmetrizable cluster algebra, under the assumption that it admits a suitable 
 categorification.  From \cite{Demonet}, it follows that the cluster algebras of finite type all admit a categorification hence one expects\footnote{The author of these lines, at least.} that the $\mathcal Y$-cluster identity $(\boldsymbol{\sf R}_\Delta)$ admits a quantum deformation $(\boldsymbol{q{\sf R}}_\Delta)$ for any Dynkin diagram $\Delta$. \sk

As far as we know, even if it is very likely the case, this has not been worked out 
in general yet  when $\Delta$ is multiply-laced (type $B,C,G_2$ and $F_4$). The single such cases considered in the literature we are aware of are the ones of rank 2 (namely $B_2$ and $G_2$): the corresponding quantum identities
$(\boldsymbol{q{\sf R}}_{B_2})$ and $(\boldsymbol{q{\sf R}}_{G_2})$
 indeed hold true, {\it cf.}\,\cite{Kuroki} and Proposition 6.1 of \cite{KY}.\footnote{Note however that no proof is given that 
 that these two quantum indentities hold true. 
According to the authors, these have been checked by explicit computations which are not given in \cite{KY}. We also note that the arguments of the quantum dilogarithm appearing in equation (13) of this paper (wich corresponds to the $B_2$ case)  do not formally correspond  to the quantum cluster variables $X_i$ we use in \eqref{Eq:q-XW-B2} to define 
${\boldsymbol{q\mathcal X\hspace{-0.04cm}\mathcal W}}_{\hspace{-0.04cm} B_2}$. 
However, we claim that there is a quantum deformation of 
$(\boldsymbol{{\sf R}}_{B_2})$, denoted by 
$(\boldsymbol{q{\sf R}}_{B_2})$, such that the arguments of the quantum dilogarithms appearing in it are precisely the elements of ${\boldsymbol{q\mathcal X\hspace{-0.04cm}\mathcal W}}_{\hspace{-0.04cm} B_2}$. Details are left to the reader.} 
\mk

Considering the above and \S\ref{SS:cluster-webs:B2}, some of the Questions \ref{Quest:qW-ARs} can be generalized to ${\boldsymbol{\mathcal X\hspace{-0.04cm}\mathcal W}}_{\hspace{-0.04cm} B_2}$ in a quite explicit form: 
\begin{questions}
\begin{enumerate}
\item[{\it 1.}]  {\it The web ${\boldsymbol{\mathcal X\hspace{-0.04cm}\mathcal W}}_{\hspace{-0.04cm} B_2}$ carries two linearly independant polylogarithmic ARs of weight 2: the standard dilogarithmic one $(\boldsymbol{{\sf R}_{B_2}})$, but also  the `symmetric' one 
$(\boldsymbol{{\sf S}_{B_2}})$. Since the former admits a quantum deformation (namely $\boldsymbol{q{\sf R}_{B_2}}$), what about the latter?  In the same vein, one can ask the same question for $(\boldsymbol{{\sf R}_{B_2}})-(\boldsymbol{{\sf S}_{B_2}})$ which is equivalent to Newman's identity $(\mathcal N_6)$ involving the sole bilogarithm ${\bf L}{\rm i}_2)$: it would be interesting to know whether  this classical identity admits a quantum deformation or not. }
\sk
\item[{\it 2.}]  {\it Does each of the non-polylogarithmic ARs of  ${\boldsymbol{\mathcal X\hspace{-0.04cm}\mathcal W}}_{\hspace{-0.04cm} B_2}$, namely the rational one  \eqref{E:EqXYZ=1} or the abelian relation $({\sf A}_2)$ (involving the function ${\sf A}(u)={\rm Arctan}(\sqrt{u})${\rm \big)}  admit a non-commutative deformation as a quantum AR for ${\boldsymbol{q\mathcal X\hspace{-0.04cm}\mathcal W}}_{\hspace{-0.04cm} B_2}$?}
\end{enumerate}
\end{questions}

It would be very nice and interesting if the questions above (or at least some of them) could be answered by the affirmative. We believe that this is the case.  However, we must confess that we have no conceptual argument to offer in support 
of this possibly too naive guess.
\sk

Finally, considering the explicit results in \S\ref{SS:cluster-webs:G2}, we mention that very similar questions to those above can be asked for the $G_2$ case as well.\sk

\paragraph{\bf A quantum deformation of Spence-Kummer web?} 
What  we want to discuss briefly here  can  be summarized by the following questions: 
\begin{equation}
\label{Question:kotok}
\begin{tabular}{l}
{\it Does the secondary cluster web 
${\boldsymbol{\mathcal U\hspace{-0.04cm}\mathcal W}}_{\hspace{-0.04cm} A_3}$ admit a quantum deformation and if so,}\\
{\it  what about its ARs and in particular its two linearly independent trilogarithmic}\\
{\it    ones? Do these identities admit quantum non-commutative deformations (possibly}\\
{\it  in multiplicative form)?}
\end{tabular}
\end{equation}
Affirmative answers would be very interesting. 
But there are several obstacles which can be identified  before answering the questions above. Here is a  list of the main ones: 
\begin{enumerate}
\item By definition, ${\boldsymbol{\mathcal U\hspace{-0.04cm}\mathcal W}}_{\hspace{-0.04cm} A_3}$ is  the web given by the restriction of  the $\mathcal X$-cluster web ${\boldsymbol{\mathcal X\hspace{-0.04cm}\mathcal W}}_{\hspace{-0.04cm} A_3}$ along the 
secondary cluster manifold $\mathcal U_{A_3}\subset \mathcal X_{A_3}$. In view of generalizing this to the quantum setting, the first step would be to define the `quantum deformation 
${\boldsymbol{q\mathcal X\hspace{-0.04cm}\mathcal W}}_{\hspace{-0.04cm} A_3}$'  of the full $\mathcal X$-cluster web of type $A_3$.  If the sets of classical or quantum $\mathcal X$-cluster variables are in bijection (a canonical bijection being induced by taking the classical limit $q\rightarrow 1$), there is a difficulty in order to define ${\boldsymbol{q\mathcal X\hspace{-0.04cm}\mathcal W}}_{\hspace{-0.04cm} A_3}$ since, as it follows from the considerations in the two previous paragraphs, the notion of `non-ordered web' does not make really sense in a non-commutative world. Having a relevant definition of ${\boldsymbol{q\mathcal X\hspace{-0.04cm}\mathcal W}}_{\hspace{-0.04cm} A_3}$ would require to have an ordering 
of the quantum cluster variables and since the classical version ${\boldsymbol{\mathcal X\hspace{-0.04cm}\mathcal W}}_{\hspace{-0.04cm} A_3}$ is not associated to a cluster period, we 
do not see from where such an ordering might come.
\item Another problem that arises, independently of the one that has just been evoked, is that 
of defining the (or a) `quantum version ${{q\mathcal U}}_{A_3}$' of the secondary cluster manifold $\mathcal U_{A_3}\subset \mathcal X_{A_3}$. Since ${\mathcal U}_{A_3}$ is the image of 
$p=p_{A_3}: \mathcal A_{A_3}\rightarrow \mathcal X_{A_3}$, a naive attempt would be to try to have a quantum version of this map to define ${{q\mathcal U}}_{A_3}$ similarly. This question has been investigated in full generality by several authors ({\it e.g.}\,see Appendix 2.7 of \cite{FockGoncharovINV} or \cite{GL}) from the work of which it follows that the notion of `cluster ensemble' $(\mathcal A,\mathcal X,p)$ can be generalized to the quantum setting if and only if when the initial exchange matrix has full rank, which is equivalent to $\mathcal U$ not being proper (hence coinciding with $\mathcal X$), which is precisely not the case under scrutiny.  The fact that  the quantum version of the $\mathcal A$-cluster variety $\mathcal A_{A_3}$ does not exist (and consequently, neither does a quantum version of $p_{A_3}$)  does not allow to treat the question of quantifying ${{\mathcal U}}_{A_3}$ via a similar approach to that of  the classical case.
\sk 

However,  the general theory tells us that ${{\mathcal U}}_{A_3}$ is a symplectic leaf of the cluster Poisson variety ${{\mathcal X}}_{A_3}$: the restriction along 
 ${{\mathcal U}}_{A_3}$
 of the cluster 
 Poisson bracket 
 $\{\cdot,\cdot\}$ of ${{\mathcal X}}_{A_3}$ (briefly discussed in \S\ref{Par:Cluster-varieties-cluster-ensembles}) is nondegenerate hence gives rise to a symplectic 2-form $\omega_{\mathcal U}$ on ${{\mathcal U}}_{A_3}$.  Hence the pair 
 $\big( {{\mathcal U}}_{A_3}, 
\omega_{\mathcal U} 
 \big)$ is a 
  holomorphic symplectic variety and as such, admits a canonical quantization.  Nevertheless making this quantization explicit and relating it to some `quantum cluster like variables'  does not seem obvious and would  probably  require  some extra work. 
\item Finally, the existence of a quantum deformation of  $({\cal S}{\cal K})$ ({\it cf.}\,\S\ref{SS:Spence-Kummer}) would in particular imply that a `quantum trilogarithm' does exist. This is still unknown at the time of writing although it is a very appealing question  that has certainly already caught the attention of specialists of the subject. 
\end{enumerate}

The existence of these three obstacles and the fact that they are 
of distinct natures might be the sign that the questions \eqref{Question:kotok} cannot be answered by the affirmative. It would be interesting to know more about this.


\subsubsection{\bf Categorification of cluster algebras and quantum dilogarithmic identities.} 
Our main reference for what is discussed in this paragraph is \cite{KellerCTQDI} to which we refer  for more background and details. In particular, below we will not work exactly with $\boldsymbol{\Psi}_q$ but with  $$\mathbb E=\boldsymbol{\Psi}_{q^{1/2}}\, ,$$ 
 where 
 $q^{1/2}$ stands 
  for an arbitrary (but fixed) square root  of $q$. \sk 
  
Let $\Delta$ stand for the bipartite quiver associated to a simply-laced  Dynkin diagram. 
Then building on some work by Kontsevich and Soibelman, 
Reineke gave in \cite{Reineke} a conceptual proof and a categorical interpretation of the quantum dilogarithmic identity, denoted by $(\mathbb E_\Delta)$,  associated to the cluster $Y$-period of type $\Delta$.  We review this very briefly below and state several questions all this suggests.
\sk 

 The main conceptual point here is the fact that cluster algebras can be categorified, at least in many cases, {\it e.g.}\,in the skew-symmetric case. If $Q$ stands for a quiver defining a 
(skew-symmetric) cluster algebra,  there is  a category $\mathscr C_Q$  within which most of the `cluster notions' of the cluster algebra associated to $Q$ can be interpreted conceptually.  
The construction of $\mathscr C_Q$, if now well-known by people familiar with this circle of ideas, is sophisticated as such hence will not be reviewed here (for more on this, see 
Keller's nice introductive papers  \cite{KellerCAQTC},\cite{KellerCat} or \cite{KellerCTQDI}). We will just say that in the simplest cases (such as $Q$ of Dynkin type $ADE$ with base field $k=\mathbf C$), it is constructed algebraically from the category 
${\sf mod }(kQ)$ of $k$-finite dimensional modules over the path algebra $kQ$.  
\sk

The cluster category satisfies several interesting properties (it is a triangulated 2-Calabi-Yau category, etc.) and in particular offers a setting in which the  quantum cluster identities can be understood conceptually.\footnote{See for instance the table p.\,124 of  \cite{KellerCAQTC}  which indicates some correspondances between the categorical and the classical settings of cluster algebras.}   A supplementary ingredient needed to achieve this is the notion of `stability condition' (also termed a `central charge' by physicits) which is a complex valued-map 
$Z: K_0\big({\sf mod }(kQ)\big)\rightarrow \mathbf C$ 
defined on the Grothendick group $K_0\big({\sf mod }(kQ)\big)$ satisfying several properties similar to those satisfied by the classical notion of stability of algebraic geometry. 
\sk 

We assume that $Q$ is of Dynkin type $ADE$ from now on. Denote this Dynkin type by $\Delta$, set $n$ as the rank of the latter and  let 
$\beta$ be the antisymmetric form defined by the initial exchange matrix associated to $Q$, namely  
$ \beta(e_i,e_j)= \# \{\,\mbox{arrows }  i\rightarrow j \,\mbox{ in } Q\,\}
-\# \{\,\mbox{arrows }  j\rightarrow i \,\mbox{ in } Q\,\}$ for any $i,j=1,\ldots,n$.  
If $(e_i)_{i=1}^n$ stands for the standard $\mathbf Z$-basis of $\mathbf Z^n$ ({\it i.e.}\, $e_i=(\delta_i^j)_{j=1}^n$ for any $i$) one sets $X^{e_i}=X_i$ for $i=1,\ldots,n$ 
and one defines 
 the {\bf  complete quantum affine space} $\widehat{\mathbb A}_Q$ as  the following non-commutative $\mathbf Q\big(q^{1/2}\big)$-algebra of formal series in the $X_i$'s:
 \begin{equation}
 \label{Eq:A-Q}
 \widehat{\mathbb A}_Q=\mathbf Q\big(q^{1/2}\big)
 \Big\langle  \Big\langle \, y^\alpha \,,\,  \alpha\in \mathbf N^n\, \Big\lvert \hspace{0.1cm}
 y^\alpha y^\beta=q^{\frac{1}{2}\beta(\alpha,\beta)}y^{\alpha+\beta}\,  \Big\rangle \Big\rangle \, . 
 \end{equation}

Then  assuming that $Z$ is a stability function of discrete type, one defines the {\bf refined DT-invariant}\footnote{Here `DT' stands for `Donaldson-Thomas', the two authors who first considered  this invariant (for the 
bounded derived category of coherent sheaves  on a CY 3-fold).} 
$$\boldsymbol{\mathbb E}_{Q,Z}\in \widehat{\mathbb A}_Q$$ as a certain ordered product
of the elements of $\widehat{\mathbb A}_Q$ 
of the form $\mathbb E\big(X^{\underline{\rm dim} M}\big)$ where $M$ ranges among the (finite) set of  $Z$-stable elements of ${\sf mod }(kQ)$ (with $\underline{\rm dim} M\in \mathbf N^n$ standing for the corresponding dimension vector for each such module $M$).  What justifies the use of the term `invariant' here is the following theorem due to Reineke \cite{Reineke} (see also Theorem 1.6 in \cite{KellerCTQDI}):
\begin{thm}
\label{Thm:Reineke}
The DT-invariant $\boldsymbol{\mathbb E}_{Q,Z}$ does not depend on the discrete stability condition $Z$. 
\end{thm}

Then using this theorem, it is possible to give a more conceptual proof of a quantum dilogarithmic identity equivalent to the quantum identity $(\boldsymbol{q{\sf R}}_{\Delta})$ 
of \S\ref{Parag:QC-I} associated to $\boldsymbol{i}_\Delta$, the cluster $Y$-period of Dynkin type $\Delta$ (which corresponds to the period \eqref{Eq:period-slice} specialized when $\Delta'=A_1$ and where $\Delta$ is the Dynkin type of $Q$). Indeed, one can prove (see after \cite[Corollary 1.7]{KellerCTQDI}) that there exist  two discrete stability conditions $Z_1$ and $Z_2$ such that 
\begin{itemize}
\item[$-$] \vspace{-0.15cm}the stable objects for $Z_1$, enumerated in some specific way, are precisely modules $S_{i_1},\ldots,S_{i_n}$ whose images in $\mathscr C_Q$ are the objects associated to the initial (quantum) cluster variables $X_{i_1},\ldots,X_{i_n}$ for a certain ordering $(i_1,\ldots,i_n)$ of $\{1,\ldots,n\}$;
\item[$-$] \vspace{-0.15cm} the $Z_2$-stable modules correspond to the indecomposable representations, 
with associated dimension vectors $\boldsymbol{\alpha}_1,\ldots,\boldsymbol{\alpha}_N \in \mathbf N^n$. 
\end{itemize}

From Reineke's theorem, it follows that both refined DT-invariants $\mathbb E_{Q,Z_1}$  and 
$\mathbb E_{Q,Z_2}$ coincide, which in explicit form can be written as the following 
quantum dilogarithmic identity 
$$
\big( \boldsymbol{\mathbb E_\Delta} \big)
\hspace{4cm}
\mathbb E(X_{i_1})\cdots \mathbb E(X_{i_n})=
\mathbb E\big(X^{\boldsymbol{\alpha}_1}\big)\cdots \mathbb E\big(X^{\boldsymbol{\alpha}_N}\big)\, ,
\hspace{6cm} {}^{}
 $$
an identity which is equivalent to $(\boldsymbol{q{\sf R}}_{\Delta})$ as mentioned above (see also \cite[\S3.5]{KN}).
 \sk
 
 The interpretation of $\big( \boldsymbol{\mathbb E_\Delta} \big)$ as an equality between two quantities attached to $\mathscr C_\Delta$ is very interesting and suggests whether  something similar happens for any other AR of the $\mathcal Y$-cluster web of type $\Delta$, 
 under the assumption that this AR admits a deformation quantization. 

Since there are many logarithmic ARs which can be constructed from $(\boldsymbol{{\sf R}}_{\Delta})$ (see \S\ref{Par:Dilog-Generation-LogAR}) and because we believe that any AR of this type admits a quantum deformation, it seems to us that the following question makes sense and is relevant:  
\begin{question} 
\label{Quest:cococo}
Assuming that the logarithmic ARs of ${\boldsymbol{\mathcal Y\hspace{-0.04cm}\mathcal W}}_{\hspace{-0.04cm} \Delta}$ admit deformation quantization, can the corresponding `quantum log-ARs' be understood within the categorification $\mathscr C_\Delta$ of the cluster algebra of Dynkin type $\Delta$?
\end{question}

But what would be the most interesting in this regard, would be to investigate similar questions for webs associated to multiply-laced Dynkin diagrams, since in many cases, these webs carry not only one, but two linearly independent polylogarithmic ARs of weight 2. To simplify the exposition, we only deal with the $B_2$ case below.\sk

If a quantum deformation of $(\boldsymbol{{\sf R}}_{B_2})$ is known to hold true (see \S\ref{SubPar:QuantumClusterWebB2} above), it is not 
known whether  this identity ({\it cf.}\,equation (14) in \cite{KY}) corresponds to the equality of two refined DT-invariants, for two distinct discrete stability functions. The main fact used in the simply-laced case  to argue  in this way and then prove the quantum identity $\big( \boldsymbol{\mathbb E_\Delta} \big)$ is the categorification $\mathscr C_\Delta$  of the underlying cluster algebra of type $\Delta$. Thanks to \cite{Demonet}, we know that the cluster algebra of type $B_2$ admits a categorification $\mathscr C_{B_2}$ as well although it is not of skew-symmetric (but only skew-symmetrizable)  type.  
It is then not unreasonable to 
expect that most of the things satisfied by categorifications $\mathscr C_\Delta$ with $\Delta$ simply laced also hold true, possibly under slight modifications, in the $B_2$ case as well. 
\sk


Considering the distinct kinds 
 of ARs of  
${\boldsymbol{\mathcal X\hspace{-0.04cm}\mathcal W}}_{\hspace{-0.04cm} B_2}$, 
one can ask the following questions which look even more interesting than the corresponding 
one stated just above  about the simply-laced case: 
\begin{questions} 
\label{Quest:cococo}
\begin{enumerate}
\item Assuming that a given AR of ${\boldsymbol{\mathcal Y\hspace{-0.04cm}\mathcal W}}_{\hspace{-0.04cm} B_2}$ admits a deformation quantization, can the corresponding `quantum AR' be understood within the categorification $\mathscr C_{B_2}$?
\item In particular, does the `symmetric dilogarithmic identity'  of ${\boldsymbol{\mathcal Y\hspace{-0.04cm}\mathcal W}}_{\hspace{-0.04cm} B_2}$
  (namely, the AR associated to identity $(\boldsymbol{{\sf S}}_{B_2})$ 
page \pageref{(S-B2)}) admit a quantum deformation and if so, can this quantum identity be interpreted in categorical terms within the associated cluster category? 
\end{enumerate}
\end{questions}

Until here in this subsection, we have discussed several questions about `quantum cluster webs' inspired from  classical notions of web geometry.  In the next paragraph, we will discuss briefly how, conversely, some results holding true in the non-commutative quantum world lead us  to 
speculate about a possible relevant generalization of the notion of webs, namely  webs made up of an infinite countable number of foliations and carrying non-trivial dilogarithmic ARs.


%
 
\subsubsection{\bf Quantum dilogarithmic identities with an infinite number of terms.}
 \label{Par:QdilogIdentInfiniteNumberOfTerms}
A very interesting upshot of \cite{Reineke} is that, at least in some situations,  Theorem \ref{Thm:Reineke} holds true for not necessarily discrete stability functions hence leads to identities with an infinite number of terms satisfied by the quantum dilogarithm $\mathbb E$.   The case of the general Kronecker quiver $K_m$ with $m\geq 1$ arbitrary is treated in 
  \cite[\S6.1]{Reineke} (see also \cite{Mozgovoy}). Another class of examples (namely, the $n$-cycle quivers) is considered in the preprint \cite{Allman}. 
 Here, to make the exposition as short as possible, we will only consider the case of the second Kronecker quiver $K_2$ as discussed in \cite[\S1.8]{KellerCTQDI}. 
 \mk

The second Kronecker quiver  and its associated exchange matrix are 
$$
\begin{tabular}{ccc}
$K_2$\,  : \, $\xymatrix{1 \ar@<0.5ex>[r] \ar@<-0.5ex>[r] & 2} $ &\quad   \qquad and \qquad  \quad  &  \qquad 
$B_{K_2}=\begin{bmatrix}
0 & 2\\
-2 & 0
\end{bmatrix}
$
\end{tabular}
$$
The associated cluster algebra is of infinite type, whose clusters are the pairs $(x_m,x_{m+1})$ for $m\in \mathbf Z$ where the $x_m$'s are constructed inductively from the initial cluster $(x_1,x_2)$ by means of the recurrence 
\begin{equation}
\label{Eq:classical K2}
x_{m-1}x_{m+1}=\big(1+x_m\big)^2 \quad \big(m\in \mathbf Z\big)\,. 
\end{equation}

Let $X^{\boldsymbol{\alpha}}$, for $\boldsymbol{\alpha}\in \mathbf Z^2$, be the elements  appearing in the definition  \eqref{Eq:A-Q} of  the associated quantum affine space 
$\widehat{\mathbb A}_{K_2}$.
In particular, one has $X^{(a,b)}= q^{ab} X_1^aX_2^b$ for any $(a,b)\in \mathbf Z^2$.
  For a first stability condition $Z_1$, the corresponding $\mathbb E_{K_2,Z_1}$ is equal to 
$\mathbb E(X_1)\mathbb E(X_2)$.  But for another non discrete stability condition $Z_2$, the corresponding DT-invariant is given by an infinite (countable) product of 
 terms. The equality between these two DT-invariants can be made entirely explicit: 
setting $\mathbb E(\boldsymbol{\alpha})=\mathbb E(X^{\boldsymbol{\alpha}})$ for any $\boldsymbol{\alpha}\in \mathbf Z^2$, the corresponding 
 quantum dilogarithmic identity is  
$$ \Big( 
\mathbb E_{K_2}
\Big) \hspace{0.5cm}
\mathbb E(X_1)\, \mathbb E(X_2)=\Big( \mathbb E(0,1)\,\mathbb E(1,2)\,\mathbb E(2,3)\cdots \Big)\,
\mathbb  E(1,1)^4\,\mathbb  E(2,2)^{-2}
\Big( \cdots \mathbb E(3,2)\,\mathbb E(2,1)\,\mathbb E(1,0\big) \Big)\, .
\hspace{5cm} {}^{}
$$
and one is interested by what might be a possible  classical counterpart to this quantum identity with countably many terms.  
\sk
 
In the case when $m=1$, the quiver $K_1$ is the one associated to the Dynkin diagram $A_2$ and the identity 
$(\mathbb E_{K_1})$ coincides with the first discovered form of a quantum dilogarithmic pentagon identity, namely $\mathbb E(X_1)\mathbb E(X_2)=\mathbb E(X_2)\mathbb E(q^{-1}X_1X_2)\mathbb E(X_1)$  which is equivalent to $(\boldsymbol{q{\sf R}}_{A_2})$. 
In the terminology of \cite{KN}, the former identity is said to be in `tropical form' whereas the second is in `universal form'.  If the classical pentagon identity can be recovered from 
$(\mathbb E_{K_1})$  ({\it cf.}\,\cite[\S2]{FaddeevKashaev}),  identities in universal form are much more suitable to get the  corresponding dilogarithmic identities in the classical limit. Indeed, for any cluster period $\boldsymbol{\nu}$ (of a skew-symmetric cluster algebra), the arguments of 
the cluster dilogarithm ${\sf R}$ in $(\boldsymbol{{\sf R}}_{\boldsymbol{\nu}})$ are the straightforward limits when $q^{1/2}\rightarrow 1$ of the quantum arguments of $
\mathbb E_q=\boldsymbol{\Psi}_{q^{1/2}}$ appearing in $(\mathbb E_{\boldsymbol{\nu}})=(\boldsymbol{q{\sf R}}_{\nu })$, 
 the quantum identity in universal form associated to $\boldsymbol{\nu}$. 
\sk 

In \cite[\S3.5]{KN}, the authors describe  Volkov's so-called `shuffle method' which allows to deduce the quantum identity in universal form  from the one in tropical form.  A priori, it only applies to  (quantum) dilogarithmic identities associated to (finite) cluster periods, but we believe that it could be adapted to handle the case of an identity with infinitely many terms such as $\big( 
\mathbb E_{K_2}
\big)$.  

The quantum versions  of the recurrence relations \eqref{Eq:classical K2}
are the following non-commutative ones 
\begin{equation}
\label{Eq:quantum K2}
X_{m-1}X_{m+1}=\big(1+q^{1/2}X_m\big)\big(1+q^{3/2}X_m\big) \quad \big(m\in \mathbf Z\big)\,. 
\end{equation}
and are precisely the $\mathcal X$-quantum variables associated to $K_2$ (up to inversion).  Indeed, from \eqref{Eq:quantum K2}, one deduces inductively a rational expression for $X_m$ in terms of $X_1$ and $X_2$ for all $m$, and setting 
$$
S_m^{q^{1/2}}=
\begin{cases}\, 
 \Big(\, \big(X_{m}^{-1}, X_{m+1}\big)\, , B_{K_2}\Big) \hspace{0.81cm}\mbox{for } m \mbox{ odd}; 
 \\
 \, 
\Big(\, \big(X_{m+1}, X_m^{-1}\big)\, , -B_{K_2}\Big) \hspace{0.6cm}\mbox{for } m \mbox{ even}\, , 
\end{cases}
$$
 it can be verified that 
$S_{m+1}^{q^{1/2}}$ is the quantum $\mathcal X$-seed obtained from $S_m^{q^{1/2}}$ by the quantum $\mathcal X$-mutation $\mu^q_{\varepsilon_m}$ of \cite[\S3]{FockGoncharovINV}, with $\varepsilon_m=1$ when $m$ is even and  $\varepsilon_m=2$ when it is  odd.  
\sk 

When $q^{1/2}$ goes to 1, the $X_m$'s converge to the quantities $x_m$ satisfying   \eqref{Eq:classical K2} which are the classical $\mathcal X$-cluster variables defined by the initial seed $S_0=\big( (x_1, x_2/(1+x_1)^2), -B_{K_2}\big)$. Hence $S_1=\mu_1(S_0)=\big( (x_1^{-1},x_2),B_{K_2}\big)$ and any $x_m$ is a positive Laurent polynomial in the variables $x_1,x_2$. Actually, closed formulas have been obtained for each $x_m$. In particular, for any $n\geq 0$, one has 
\begin{equation}
\label{Eq:xn+3-x-n}
x_{n+3}=\frac{\big(P_n\big)^2}{x_1^nx_2^{n-1}} \qquad 
\mbox{ and } \qquad 
x_{-n}=\frac{\big(P_{-n}\big)^2}{x_1^{n}x_2^{n+1}} \hspace{0.2cm} ,
\end{equation}
where $P_{\pm n}$ is an element of $\mathbf Z_{>0}[x_1,x_2]$ with constant term equal to 1 for any non negative $n$. \sk 

Instead of $(\mathbb E_{K_2})$,  we are going to work with the following identity 
$$ \Big( 
\mathbb E'_{K_2}
\Big) \hspace{0.4cm} 
 \cdots \mathbb E(3,2)\,\mathbb E(2,1)\,\mathbb E(1,0)
 \,
\mathbb  E(0,1)^{-1}\,\mathbb  E(1,0)^{-1}
 \mathbb E(0,1)\,\mathbb E(1,2)\,\mathbb E(2,3)\cdots  
=\mathbb E(X^{(2,2)})\mathbb E(X^{(1,1)})^2 
\hspace{5cm} {}^{}
$$
which is equivalent to $(\mathbb E_{K_2})$ (see equation (A.6) in \cite{GMN}). 
\sk 

From \eqref{Eq:xn+3-x-n}, it comes that relatively to the variables $x_1,x_2$ (1) the tropical signs of the $x_m$'s are all negative; the  $\boldsymbol{c}$-vector of $x_{n+3}$ (resp.\,of $x_{-n}$) is $(-n,-n-1)$ (resp.\,$(-n,-n+1)$) for any $n\geq 0$.\footnote{Note that if  $(x_1^{-1},x_2)$ is  a cluster (the one associated to the seed $S_1$), it is not the case for $(x_1,x_2)$. Hence the notions of tropical signs and of  $\boldsymbol{c}$-vector considered here for the $x_m$'s do not correspond {\it stricto sensu} to the standard notions, but are very slight modified versions of these.}.  Accordingly,  the LHS  of the quantum dilogarithmic identity just above  can be seen as the version 
written in tropical form of the expression in universal form obtained by taking the product of the  $\mathbb E(x_m)$'s for $m$ increasing from  $-\infty$ to $+\infty$. \sk 

On the other hand, since $X^{(2,2)}=X^{(1,1)}\cdot X^{(1,1)}$ (as it is easy to check), the RHS of 
$ \big( \mathbb E'_{K_2}\big)$ can be seen as a quantum dilogarithmic function of the single quantum expression $X^{(1,1)}$. 
Hence, under the assumption that Volkov's shuffle method applies to infinite product as well,  
 the version in universal form to   $\big( 
\mathbb E'_{K_2} \big)$ would be written as 
$$ \Big( 
\mathbb E_{K_2}^{univ}
\Big) \hspace{1.7cm}
\cdots 
\mathbb E(X_{-2})\, 
\mathbb E(X_{-1})\, \mathbb E(X_1)\, \mathbb E(X_2) \, \mathbb E(X_3)  \,  \mathbb E(X_4)   \, \mathbb E(X_5) \cdots  =\boldsymbol{\bf E}\Big(X^{(1,1)}\Big)\, .
\hspace{3cm}{}^{}
$$
for some   quantum dilogarithmic function $\boldsymbol{\bf E}(z)$ of 
$z\in \widehat{\mathbb A}_{K_2}$,  namely a finite  product 
of powers
  of terms of the form $\mathbb E(\varphi(z))$  for some rational functions  $\varphi\in \mathbf Q(q^{1/2})(z)$. 

\begin{question}
\label{Quest:kokoko}
\begin{enumerate}
\item Does a quantum identity of the form  $\big( 
\mathbb E_{K_2}^{univ}
\big) $ above hold true?  If yes, for which dilogarithmic right hand side  $\boldsymbol{\bf E}$  precisely?
\item  If the answer to 1.\,is affirmative, is 
the corresponding identity $\big( 
\mathbb E_{K_2}^{univ}
\big) $ equivalent (by means of formal and/or algebraic operations, possibly in infinite number) to $\big( 
\mathbb E_{K_2}
\big)$? 
\end{enumerate}
\end{question}

Since the arguments of $\mathbb E$ in the LHS of $\big( 
\mathbb E_{K_2}^{univ}
\big) $ all admit natural limits when $q^{1/2}\rightarrow 1$ contrarily to those in 
$\big( 
\mathbb E_{K_2}
\big) $ or $\big( 
\mathbb E_{K_2}'
\big) $, it is much  easier  to investigate the possible classical limit of the former identity than it is the case for the last two. As a preliminary step,  let us  discuss what could be the dequantization of the `{\it quantum web attached 
to $\big( 
\mathbb E_{K_2}^{univ}
\big) $}', denoted by 
$\boldsymbol{q\mathcal W}_{\hspace{-0.09cm}K_2}^{univ}$.  By definition, the latter is the infinite quantum web whose first integrals are the cluster variables  $X_m$ for $m\in \mathbf Z$ and $X^{(1,1)}$ (suitably ordered). 
Accordingly,  it is natural to require for its dequantization `$\boldsymbol{\mathcal W}_{\hspace{-0.09cm}K_2}^{univ}=\lim_{q^{1/2}\rightarrow 1} 
\boldsymbol{q\mathcal W}_{\hspace{-0.09cm}K_2}^{univ}$' 
 that it must have all the cluster variables 
$x_m=\lim_{q^{1/2}\rightarrow 1} X_m$ among its first integrals.  At the opposite, since $X^{(1,1)}$ is not a cluster variable, what might be its dequantization 
$x^{(1,1)}=\lim_{q^{1/2}\rightarrow 1} X^{(1,1)}$ is not straightforward.  \sk

Drawing inspiration from \cite[Appendix A]{GMN}, it is interesting to consider new (analytic) variables $a$ and $b$ such that $x_0=\cosh^2 ( b)/{\sinh^2 (a)}$ and 
$x_1=\cosh^2 ( a+b)/{\sinh^2 (a)}$. Hence for any $m\in \mathbf Z$, the quantity $x_m$ admits the following nice closed expression in terms of $a$ and $b$
$$x_m=\frac{\cosh^2 ( am+b)}{\sinh^2 (a)}$$
and considering the argumentation between (A.5) and (A.6) in  \cite[Appendix A]{GMN}, one deduces that in the coordinates $a,b$, one has $x^{(1,1)}=e^{-2a}$.   Hence  
\begin{equation}
\label{Eq:WK2-analytic-model}
\boldsymbol{\mathcal W}\bigg(\, 
\frac{\cosh^2 ( am+b)}{\sinh^2 (a)}
\, , \, e^{-2a}\hspace{0.2cm}  \Big\lvert \hspace{0.2cm} m\in \mathbf Z\, 
\bigg)\, 
\end{equation}
 is an analytic model of  $\boldsymbol{\mathcal W}_{\hspace{-0.09cm}K_2}^{univ}$.  Trying to express $a$ in terms of $x_1$ and $x_2$, one obtains that  $2a= {\rm Log}\, \xi$ where $\xi$ stands for a root of $z^2-\tau z+1$ where 
$$
\tau=\tau_{K_2}=\frac{1+2x_1+2x_2+x_1^2+x_2^2}{x_1x_2}\, .
$$  
Thus  $\xi=(\tau\pm \sqrt{\tau^2-4})/2$ from which it follows that, seen as a function of $x_1,x_2$,  $\tau$ can be taken as a first integral of the foliation corresponding to the one defined by $e^{-2a}$ in the coordinates $a,b$. 
This gives us that $\boldsymbol{\mathcal W}_{\hspace{-0.09cm}K_2}^{univ}$ is the following web with positive Laurent polynomials as first integrals: 
$$\boldsymbol{\mathcal W}_{\hspace{-0.09cm}K_2}^{univ}=
\boldsymbol{\mathcal W}\Big(
\hspace{0.15cm}
x_m 
\hspace{0.15cm} , \, 
\tau_{K_2} 
\hspace{0.15cm}
 \big\lvert \hspace{0.2cm} m\in \mathbf Z\, 
\Big)\, .
$$

As already mentioned above, for any cluster period $\boldsymbol{\nu}$, the classical dilogarithmic identity $(\boldsymbol{\sf R}_{\boldsymbol{\nu}})$ can be recovered from its quantum version 
$(\boldsymbol{q{\sf R}}_{\boldsymbol{\nu}})$ ({\it cf.}\,Appendix A of \cite{KN})\footnote{Actually,  \cite[Appendix A]{KN} is not about the identities 
$(\boldsymbol{\sf R}_{\boldsymbol{\nu}})$ and $(\boldsymbol{q{\sf R}}_{\boldsymbol{\nu}})$ we are dealing with here (which are said to be `{\it in universal form}' by Kashaev and Nakanishi),  but concerns other versions of them, namely dilogarithmic identities said in   `{\it in local form}'. We believe that the analysis in   
 \cite[Appendix A]{KN} can be adapted to the identities in universal form as well but a confirmation of that would be welcome.}.  Hence it is natural to wonder about the possible existence of a  semiclassical limit 
for 
$\big( 
\mathbb E_{K_2}^{univ}\big)$. 
when
$q^{1/2}\rightarrow 1$. 
\sk 

In case the answer to Question \ref{Quest:kokoko}.\,{1} is affirmative, 
 it is very natural to ask the following ones: 
\begin{question}
\begin{enumerate}
\item  Does it exist a functional identity  with an infinite number of dilogarithmic 
terms which can be obtained from $\big( 
\mathbb E_{K_2}
\big)$ in the semiclassical  limit? 
\item If yes, does it look like 
the identity which can be
derived   formally  in a straightforward way from  $\big( 
\mathbb E_{K_2}
\big)$, namely 
\begin{equation}
\label{Eq:dilog-infinite-terms}
 \sum_{m\in \mathbf Z} {\sf R}\big(x_m\big)={\bf R}(\, \tau_{K_2}\,)
\end{equation}
where ${\bf R}$ stands for a dilogarithmic function? If it is indeed the case, what is ${\bf R}$ precisely and does the LHS of \eqref{Eq:dilog-infinite-terms} converge as a function of $x_1,x_2$ on some domain of $\mathbf C^2$? 
\item  
The definition of $\boldsymbol{\mathcal W}_{\hspace{-0.09cm}K_2}^{univ}$ given above in terms of the 
  $\mathcal X$-cluster variables (of the cluster algebra determined by $K_2$) plus the extra first integral $\tau_{K_2}$ might be generalized to the case of any rank 2 cluster algebra  as follows: let $B$ be an initial exchange matrix whose antidiagonal coefficients are $b$ and $-c$ for any two given positive integers $b$ and $c$.  Denoting by $x^B_m$ the associated $\mathcal X$-cluster variables,\footnote{The $x_m^B$'s satisfy the following relations: $x^B_{m-1}x^B_{m+1}=(1+x_m^B)^b$ for $m$ odd and $x^B_{m-1}x^B_{m+1}=(1+x^B_m)^c$ if $m$ is even.} one can ask the following: are there coefficients  $n_m\in \mathbf Z$ and an at most countable set  $I$ indexing  positive Laurent polynomials $\tau_{b,c}^i \in \mathbf Z\big[x_1^{\pm 1}, x_2^{\pm 1}\big]$ and dilogarithmic functions ${\bf R}_{b,c}^i$ (for $i\in I$) such that the infinite dilogarithmic  relation below holds true?\footnote{Initially, this question was asked with the additional assumption that $I$ is a singleton (as conjectured in the case when $b=c=2$). It is after attending the seminar talk 
 \href{{https://cloud.math.univ-paris-diderot.fr/s/ZDWwLMfoNkqt8oY}}{\it `Wild quantum dilogarithm identities'}  by M.\,Reineke that we became aware  that when the exchange matrix $B$ is of wild type, we are more likely to expect that the infinite sum $  \sum_{m\in \mathbf Z} n_m\,{\sf R}\Big(x^B_m\Big)$ is equal to  a sum of countably many dilogarithmic terms of the form $ {\bf R}_{b,c}^i\big(\, \tau_{b,c}^i\,\big)$ as in \eqref{Eq:dilog-infinite-terms-b-c}.}
\begin{equation}
\label{Eq:dilog-infinite-terms-b-c}
 \sum_{m\in \mathbf Z} n_m\,{\sf R}\Big(x^B_m\Big)=
 \sum_{i \in I}
 {\bf R}_{b,c}^i\big(\, \tau_{b,c}^i\,\big)\, .
\end{equation}
 Moreover, when $B$ is skew-symmetric ({\it i.e.}\,when $b=c$ hence the initial quiver is the $b$-th Kronecker quiver $K_b$), one expects that
the coefficients $n_m$'s can all be taken equal to 1.
\end{enumerate}
\end{question}


Several preliminary computations lead us to think that the answers to the first two questions are affirmative.  In terms of the first integrals of the analytic model \eqref{Eq:WK2-analytic-model}
 of $\boldsymbol{\mathcal W}_{\hspace{-0.09cm}K_2}^{univ}$, the conjectural identity 
 \eqref{Eq:dilog-infinite-terms} 
 can be written 
\begin{equation}
\label{Eq:dilog-infinite-terms-II}
 \sum_{m\in \mathbf Z} {\sf } {\sf R}
 \left( \frac{\cosh^2 ( am+b)}{\sinh^2 (a)}
 \right)={\bf R}\left(
 e^{-2a}
 \right)
\end{equation}

Identities involving countably many dilogarithmic terms already appeared in the literature, in form very similar to the one just above, but involving the classical Rogers' dilogarithm $\mathcal R$ instead of ${\sf R}$: for instance see `Bridgeman orthospectrum identity' 
\eqref{Eq:Bridgeman-infinite-sum}
discussed in Remark \ref{Rem:Bridgeman}.  Another striking example is the infinite identity given in Theorem 2.1 of the recent preprint \cite{BridgemanArXiv} which looks like a one variable version of the conjectural identity \eqref{Eq:dilog-infinite-terms-II}. 
\begin{center}
$\star$
\end{center}

The material discussed above (some part of which is still conjectural) suggests a
possible 
 interesting generalization of many notions, concepts and results of classical web geometry to webs formed with an infinite (countable) number of foliations.  We brievly discuss this below. 
\mk

Let $I$ be an infinite countable set of indices. Assume that 
$\boldsymbol{u}=(u_i)_{i\in I}$ is a set of holomorphic functions at the origin of $\mathbf C^n$ (for some fixed $n\geq 2$), defining a {\it `web'} $\boldsymbol{\mathcal W}_{\hspace{-0.06cm}\boldsymbol{u}}={\boldsymbol{\mathcal W}}\big(\, u_i\, \lvert \, i\in I\,\big)$' (in the sense that $du_i\wedge du_j\neq 0$ for any $i,j\in I$ with $i\neq j$).  By definition, an AR for this web is 
a family $(F_i)_{i\in I}$ of holomorphic germs such that the infinite sum 
 $\sum_{i\in I} F_i(U_i)$ converges absolutely on a sufficiently small open neighbourhood of the origin and defines an element of $\mathcal O(\mathbf C^n,0)$ which actually is identically equal to a constant.  Working only with formal series, one can also consider the notion of formal abelian relation.  \sk 
 
 As in the  classical case, the space $\boldsymbol{\mathcal A}\big(\boldsymbol{\mathcal W}_{\hspace{-0.06cm}\boldsymbol{u}}\big)$ of ARs naturally carries a vector space structure.  
 As a subspace, it contains the space  of {\it `finite ARs'}, denoted by 
 $\boldsymbol{\mathcal A}^{f}\big(\boldsymbol{\mathcal W}_{\hspace{-0.06cm}\boldsymbol{u}}\big)$, which is the subpsace spanned by the family $\boldsymbol{\mathcal A}\big(\boldsymbol{\mathcal W}'\big)$'s for all finite subwebs $\boldsymbol{\mathcal W}'$ of $\boldsymbol{\mathcal W}_{\hspace{-0.06cm}\boldsymbol{u}}$. 
 This subspace (and consequently $\boldsymbol{\mathcal A}\big(\boldsymbol{\mathcal W}_{\hspace{-0.06cm}\boldsymbol{u}}\big)$ can have infinite dimension as easily shown by considering the case of $\boldsymbol{\mathcal W}_{\hspace{-0.09cm}K_2}^{univ}$: to each relation \eqref{Eq:classical K2} is associated the logarithmic abelian relation 
$${\rm Log}(x_{m-1})-
2\,{\rm Log}(1+x_m)+
{\rm Log}(x_{m+1})=0\, ,  $$
from which it follows that $\dim\big( \boldsymbol{\mathcal A}^{f}\big(
\boldsymbol{\mathcal W}_{\hspace{-0.09cm}K_2}^{univ} \big) \big)=\infty$. 
\sk

  An element of $\boldsymbol{\mathcal A}\big(\boldsymbol{\mathcal W}_{\hspace{-0.06cm}\boldsymbol{u}}\big) \setminus \boldsymbol{\mathcal A}^{f}\big(\boldsymbol{\mathcal W}_{\hspace{-0.06cm}\boldsymbol{u}}\big)$ is an {\it `infinite AR'}. 
 One defines the {\it `$\infty$-rank'}  $\rho^\infty\big(\boldsymbol{\mathcal W}_{\hspace{-0.06cm}\boldsymbol{u}}\big)$  
 of $\boldsymbol{\mathcal W}_{\hspace{-0.06cm}\boldsymbol{u}}$ as the (possibly infinite) dimension of the quotient vector space 
 $\boldsymbol{\mathcal A}^{\infty}\big(\boldsymbol{\mathcal W}_{\hspace{-0.06cm}\boldsymbol{u}}\big)=\boldsymbol{\mathcal A}\big(\boldsymbol{\mathcal W}_{\hspace{-0.06cm}\boldsymbol{u}}\big)/\boldsymbol{\mathcal A}^{f}\big(\boldsymbol{\mathcal W}_{\hspace{-0.06cm}\boldsymbol{u}}\big)$. 
Assuming that an identity \eqref{Eq:dilog-infinite-terms} indeed holds true, one expects that the corresponding dilogarithmic AR does not belong to 
$\boldsymbol{\mathcal A}^{f}\big(\boldsymbol{\mathcal W}_{\hspace{-0.09cm}K_2}^{univ} \big)$ which would implies that $\rho^\infty\big(
\boldsymbol{\mathcal W}_{\hspace{-0.09cm}K_2}^{univ} 
\big)$ is positive. 
 
 \begin{questions}
 \begin{enumerate}
 \item For classical ({\it i.e.}\,finite) webs, both notions of analytic and formal abelian relation actually coincide. Does the same hold true for infinite webs as well?
  \item  Is the $\infty$-rank of $\boldsymbol{\mathcal W}_{\hspace{-0.09cm}K_2}^{univ}$ finite? Can   a set of infinite ARs inducing a basis of  
   $\boldsymbol{\mathcal A}^{\infty}\big(\boldsymbol{\mathcal W}_{\hspace{-0.09cm}K_2}^{univ}\big)$ be explicited? 
   And if it is the case,  describe such a  set.
  \item  Same questions as those just above  for the web 
   defined by the first integrals $x_m^B$ ($m\in \mathbf Z$) and $\tau_{b,c}$ appearing in \eqref{Eq:dilog-infinite-terms-b-c}, under the assumption that such an identity indeed holds true. 
   \item  Does it exist  infinite webs carrying infinite polylogaritmic ARs of weight $w\geq 3$? 
  \end{enumerate}
 \end{questions}
 
 We are not aware of any trilogarithmic identity with infinitely many terms, hence we believe that the answer to the fourth question above might well be negative.  This constrasts with the case of infinite dilogarithmic identities, that we do know to exist thanks to the work of Bridgeman for instance.\sk 
 
  It is interesting to notice that infinite logarithmic identities also appear within the same context as `Bridgeman orthospectrum identity': among many examples, one can mention  `Basmajian identity' or the `generalized McShane-Mirzakhani identity', both for 
 finite area hyperbolic surfaces with geodesic boundary (see \cite{BridgemanTan}).  It would be interesting to consider such identities as ARs for infinite webs on certain moduli spaces of hyperbolic surfaces defined by first integrals corresponding to geometric quantities. \sk 
 
As seen above, when a finite web  carries a dilogarithmic AR, many other logarithmic ARs can be obtained from the latter (by monodoromy or derivation, see \S\ref{Par:Dilog-Generation-LogAR}). 
 One can wonder whether such a phenomenon also occurs for the infinite geometric webs just mentioned.  More explicitly, we wonder about the possibility to deduce in some way   Basmajian or McShane-Mirzakhani type logarithmic identities from dilogarithmic ones such as Bridgeman or Luo-Tan identities.


%

 \subsection{\bf Cluster algebras of finite type and moduli spaces of configurations.} 
 \label{SS:ClusterDelta-Configurations}

An interesting feature of finite type cluster algebras of type $A$ is that many things can be interpreted geometrically in terms of projective configurations of points on $\mathbf P^1$: 
for any $n\geq 1$,  the cluster variety $\boldsymbol{\mathcal X}_{A_n}$ identifies itself with $\mathcal M_{0,n+3}$, up to this identification, the cluster variables identify themselves with some cross-ratios of 4 points among the $n+3$, etc. (see \S\ref{SSub:Type-An}).  These identifications allow to use geometric arguments which proved to be crucial in several places in this text, in particular to establish Theorem \ref{T:YW-An} which is one of the most important results of this memoir.\sk 

In order to investigate further the  ($\mathcal X$- or  $\mathcal Y$-)cluster webs of Dynkin type distinct from $A$, it would be very valuable  to have similar geometric interpretations for the corresponding $\mathcal X$-cluster variety and cluster variables. 

 \begin{question}
 For any Dynkin diagram $\Delta$, is there a natural identification between the cluster variety 
 $\boldsymbol{\mathcal X}_{\Delta}$ and a moduli space of certain projective `configurations of type $\Delta$'? 
  \end{question}
As it is stated, this question is quite vague. 
A first step for answering it would be to make the term `natural' more precise and to define  the notion of  `configurations of type $\Delta$' rigorously.  \sk

In some recent preprints  \cite{AHHLT} and \cite{AHHL}, Arkani-Hamed and his coworkers have 
introduced some positive spaces which they call {\it `Cluster configuration spaces of finite type'}. 
If these spaces are certainly interesting (see the next subsection for instance), they do not provide the geometric interpretation we are looking for. Indeed, for each Dynkin diagram $\Delta$, the Cluster configuration space $\mathcal M_\Delta$ of \cite{AHHL} is defined as a {\it `smooth affine algebraic variety with a stratification in bijection with the faces of the Chapoton-Fomin-Zelevinsky generalized associahedron'} by means of explicit affine equations in coordinates.  
 But no interpretation of any kind of the points of $\mathcal M_\Delta$ in terms of geometric configurations appears in the paper hence we do not understand why the word `configuration' is used in the name given to these spaces. \mk 

Roughly, what we are looking for for the cluster variety $\boldsymbol{\mathcal X}_{\Delta}$
and the associated cluster variables for any $\Delta$  is a similar geometric/projective  description to the one in type $A_n$ in terms of $\mathcal M_{0,n+3}$ and of cross-ratios.  As far as we know, 
apart from this case, nothing has been published along this line  yet.  
On our part, we have rather precise ideas about how this question might be answered. This is the subject of an ongoing project \cite{Pirio-XDelta} that we hope to see completed one day.


 \subsection{\bf Cluster webs of higher codimension.} 
 \label{SS:ClusterWebs-higher-codim}
In type $A_n$, the geometric interpretations  
of the $\mathcal X$-cluster variables as cross-ratios on 
$\mathcal M_{0,n+3}$ tells interesting things about their anatomy.  In particular, they all are pull-backs of 
the cluster variables in type $A_{n-1}$ under one of the $n+3$ forgetful maps 
\begin{equation}
\label{Eq:varphi-i}
\varphi_i : \mathcal M_{0,n+3}
\rightarrow  \mathcal M_{0,n+2}
\qquad \mbox{ for }\,  i=1,\ldots,n+3\, 
\end{equation}
the $i$-th of which is given by forgetting the $i$-th point of the configuration.  This remark naturally leads to think that the $(n+3)$-web by curves  admitting the forgetful maps $f_i$ as first integrals is at least as important as the standard 1-codimensional cluster web  ${\boldsymbol{\mathcal X\hspace{-0.04cm}\mathcal W}}_{\hspace{-0.04cm}A_n} $. 
This curvilinear webs will be denoted
  $${\boldsymbol{\mathcal X\hspace{-0.04cm}\mathcal W}}_{\hspace{-0.04cm} A_n}^{(1)}
  = 
  {\boldsymbol{\mathcal W}}\Big(\varphi_1,\ldots,\varphi_n \Big)\, , 
  $$ the superscript with 1 between parentheses indicating that its leaves are of dimension 1.  Notice that when $n=2$, one gets nothing new since 
${\boldsymbol{\mathcal X\hspace{-0.04cm}\mathcal W}}_{\hspace{-0.04cm} A_2}^{(1)}$ coincides with ${\boldsymbol{\mathcal X\hspace{-0.04cm}\mathcal W}}_{\hspace{-0.04cm} A_2}$.

\subsubsection{}  
 Actually, the webs ${\boldsymbol{\mathcal X\hspace{-0.04cm}\mathcal W}}_{\hspace{-0.04cm} A_n}^{(1)} $ have already been considered by several authors and they appear as very similar to and as interesting as ${\boldsymbol{\mathcal X\hspace{-0.04cm}\mathcal W}}_{\hspace{-0.04cm} A_2}$: these webs have been considered first  in \cite{Burau} where Burau proves that there are projective models for which the leaves are lines (when $n$ is odd) or conics ($n$ even). These webs were studied from the point of view of their $(n-1)$-abelian relations and of their $(n-1)$-ranks more recently, by Damiano in \cite{Damiano}. He claimed  that he had proved that each web of this family is of maximal $(n-1)$-rank, with all its $(n-1)$-ARs of logarithmic type, except one which can be said of hyperlogarithmic nature and  is constructed by integrating the invariant volume form on the grassmannian $Gr_2^+(\mathbf R^{n+2})$ (of oriented $2$-planes in $\mathbf R^{n+2}$) along the fibers of the action of the Cartan subgroup of ${\rm GL}(\mathbf R^{n+2})$.  Some of the results of \cite{Damiano} are wrong and the proofs of some others are flawed but we believe that each web  ${\boldsymbol{\mathcal X\hspace{-0.04cm}\mathcal W}}_{\hspace{-0.04cm} A_n}^{(1)}$  indeed has  maximal $(n-1)$-rank.\footnote{We think that Damiano's proof of the maximality of the $(n-1)$-rank of ${\boldsymbol{\mathcal X\hspace{-0.04cm}\mathcal W}}_{\hspace{-0.04cm} A_n}^{(1)} $ is correct only  when $n$ is even. We know for sure that Damiano's proof is wrong for $n=3$ for instance, but we have a proof that the 2-rank of ${\boldsymbol{\mathcal X\hspace{-0.04cm}\mathcal W}}_{\hspace{-0.04cm} A_3}^{(1)}$ is maximal indeed. All this is the subject of a work in progress.}\mk 

The material evoked in the previous paragraph is quite interesting and if one believes (as the author of these lines) that for any Dynkin type $\Delta$, the  $\mathcal X$-cluster varieties/variables can be `geometrized' (see \S\ref{SS:ClusterDelta-Configurations}), one expects to be able to  factorize the cluster variables by maps 
$ \mathcal X_{\Delta}\rightarrow \mathcal X_{\Delta'}$
of geometric nature, for some Dynkin diagrams ${\Delta'}$ of rank  one less than that of $\Delta$. 
As explained above, we do not currently have the results that would allow us to approach this geometrically. However, it turns out that one can bypass this lack of geometric understanding of the cluster varieties in arbitrary type precisely relying on the cluster formalism. This is what we are going to explain now.  We will be rather allusive, claiming a lot without giving proofs, postponing a rigorous treatment to a possible future publication. 

%
%
%
%
%
%

 \subsubsection{}  
 Let $n$ be bigger than or equal to 2.  We set $\Delta=A_n$ in what follows, but being aware that everything below can be adapted to the case when $\Delta$ is any other Dynkin diagram of rank $n$. We set $I=\{1,\ldots,n\}$. We recall that $\vec{\Delta}$ stands for the bipartite quiver associated to $\Delta$ and we denote by $\epsilon: I\rightarrow \{\, \pm 1\,\}$ the map associating $-1$ (resp.\,1) to labels associated to sources (resp.\,sinks) in $\vec{\Delta}$.
 The initial seed we work with is $S_{\rm init}=(\boldsymbol{u},\vec{\Delta})$ with 
 $\boldsymbol{u}=(u_i^{\epsilon({i})})_{i\in I}$ (see Remark \ref{Rem:birat-E} for the reason of  that  choice of the initial cluster). 
 \sk

 The forgetful maps 
\eqref{Eq:varphi-i} 
 on 
$ \mathcal M_{0,n+3}$ correspond on the cluster side to rational dominant maps $\mathcal X_{A_n}\dashrightarrow 
\mathcal X_{A_{n-1}}$ or, up to the classical dual correspondance  spaces $\leftrightarrow$ algebras of functions, to cluster subalgebras of type $A_{n-1}$. 
A simple way to construct such a cluster subalgebra is as follows: 
we first consider a $\mathcal X$-seed $S^t=(\boldsymbol{x}^t, Q^t)$ of the initial $A_n$-cluster algebra with $\boldsymbol{x}^t=(x^t_i)_{i=1}^n$, and choose an element $k\in I=\{1,\ldots,n\}$.  We set $I_{\hat k}=I\setminus\{k\}$, $\pi_k: \boldsymbol{u}=(u_i)_{i\in I}\mapsto \boldsymbol{u}_{\hat k}=(u_i)_{i\in I_{\hat k}}$ and for any quiver $Q'$ whose vertices are labeled by elements of $I$ (such as $Q^t$), we denote by $Q'_{\hat k}$ the (possibly non connected) one obtained by removing from $Q'$ its $k$-th vertex and all the edges adjacent to it. 
Considering the variable $x_k^t$ and the $k$-th vertex of the quiver $Q^t$ as frozen, we perform all possible mutations from the seed $S^t$ under this restriction. 
We get a cluster pattern 
$$
\boldsymbol{\Sigma}_{\hat k}^t=\Big\{ 
\hspace{0.1cm}
S^{t'} =
\big(\boldsymbol{x}^{t'}, Q^{t'}\big) 
 \hspace{0.1cm}  \big\lvert  \hspace{0.1cm}
S^{t'}
=\mu_{i_1}\circ \ldots \circ \mu_{i_m}\big(S^t\big)
\, , \hspace{0.1cm}
 i_1,\ldots ,i_m\in I_{\hat k} 
 \hspace{0.1cm}
\Big\}
$$
 The images by the projection $\pi_k$ of all the clusters obtained from $S^t$ in this way form a $\mathcal X$-cluster pattern  
 which  coincides with the one associated with the initial seed  $\pi_k(S^t)= ( \boldsymbol{x}^t_{\hat k}, Q^t_{\hat k})$.   Given any  $\mathcal X$-variables  $x_i^{t'}$ of this cluster pattern with $i\neq k$, it does not depend on $x_k^t$ when seen as a rational function in the $x_i^t$'s hence 
 $$
{\boldsymbol{\mathcal X\hspace{-0.04cm}\mathcal W}}_{\hspace{-0.04cm} \hat k}^t= 
{\boldsymbol{\mathcal W}}\Big(
x_i^{t'}\hspace{0.1cm}\big\lvert \hspace{0.1cm}i\in I_{\hat k}\, , \hspace{0.1cm}
\big(\boldsymbol{x}^{t'}, Q^{t'}\big)  \in 
\boldsymbol{\Sigma}_{\hat k}^t\, 
\Big)  
$$
is a subweb of ${\boldsymbol{\mathcal X\hspace{-0.04cm}\mathcal W}}_{\hspace{-0.04cm} \Delta}$, 
 of intrinsic dimension $n-1$: the map $\pi_{\hspace{-0.04cm} \hat k}^t$ whose components are precisely the cluster variables $x_i^{t'}$ appearing in the definition of ${\boldsymbol{\mathcal X\hspace{-0.04cm}\mathcal W}}_{\hspace{-0.04cm} \hat k}^t$ has rank $n-1$ at the generic point hence $ \mathcal X^t_{\hat k}=\pi_{\hspace{-0.04cm} \hat k}^t(\mathcal X_\Delta)$ is an irreducible algebraic variety of dimension $n-1$. We get a rational map from $ \mathcal X_\Delta$ onto $ \mathcal X^t_{\hat k}$, again denoted by $\pi_{\hspace{-0.04cm} \hat k}^t$. These maps have to be seen as cluster analogues of the  maps  \eqref{Eq:varphi-i}. For instance, in type $A$, some of  these maps (the most interesting ones), are nothing else than the forgetful maps, but described in an alternative way.   Of course, 
  $\pi_{\hspace{-0.04cm} \hat k}^t: \mathcal X_\Delta\dashrightarrow \mathcal X^t_{\hat k}$ 
 admits an elementary model in the $x_i^t$'s: in these coordinates, it is given by the projection $\pi_k: (x_i^t)_{i\in I}\mapsto (x_i^t)_{i\in I_{\hat{k}}}$.\mk

Regarding the  subwebs 
${\boldsymbol{\mathcal X\hspace{-0.04cm}\mathcal W}}_{\hspace{-0.04cm} \hat k}^t$ and the associated maps $\pi_{\hspace{-0.04cm} \hat k}^t$, 
there are two natural questions: when are two  of them actually the same? And: what kind of cluster web  is it? This can be answered by noticing that for each initial data $(k,t)$ with $k\in I$ and $t$ in $\mathcal X\Gamma_\Delta$ (the $\mathcal X$-exchange graph in type $\Delta$), the cluster pattern $\boldsymbol{\Sigma}_{\hat k}^t$ meets the even bipartite belt  of the whole cluster algebra of type $\Delta$.
 Recall ({\it cf.}\,\S\ref{SubPar:BipartiteBelt}) that these $\mathcal X$-clusters are those 
obtained from the initial one by successive application of the compositions of mutations  $\mu_{\bullet\lvert\circ}=\mu_\bullet \circ \mu_\circ$ 
hence are of the form 
$$S_\ell=(\mu_{\bullet\lvert\circ})^\ell\big(S_{\rm init}\big)=\Big( \, \boldsymbol{Y}(\ell)
\, , \, \vec{\Delta} \, \Big)
\qquad \mbox{with} \quad 
\boldsymbol{Y}(\ell)=\big(Y_{\alpha(\ell,1)}, \ldots, Y_{\alpha(\ell,n)}\big)\, $$
where for any integer $\ell$ and any $i\in I$,  $Y_{\alpha(\ell,i)}$ stands for the $\mathcal Y$-cluster variables of Theorem 1.4 of \cite{FZ} associated to the root $\alpha(\ell,i)\in \Delta_{\geq -1}$. 
\sk

Assume that $h=h(\Delta)$ is even (when $h$ is odd, then $\Delta=A_{2m}$ for some $m\geq 1$ and this case, albeit slightly more subtle, can be treated similarly). Then $\mu_{\bullet\lvert\circ}$ has order $\kappa=(h+2)/2$ (modulo a permutation) hence any $\alpha\in \Delta_{\geq -1}$ can be written $\alpha=\alpha(\ell,i)$ for a unique pair $(\ell,i)$ with $\ell\in  \{0,\ldots,\kappa-1\}$ and $i\in  I$.  We denote by 
$ {\boldsymbol{\mathcal X\hspace{-0.04cm}\mathcal W}}_{\hspace{-0.04cm}\Delta,\alpha}$ 
the subweb obtained by the construction described above when starting from the cluster 
$S^\alpha=\big((Y_{\alpha(\ell,1)}, \ldots, Y_{\alpha(\ell,n)}), \vec{\Delta}\big)$ with the $i$-th coordinate (the one corresponding to $\alpha$) frozen. Then clearly a simple model for the rational associated map 
 in the initial coordinates $u_i$ is given by 
\begin{equation}
\label{Eq:pi-alpha}
\pi_\alpha=
\pi_{\alpha(\ell,i)}=
\Big( Y_{\alpha(\ell,1)},\ldots,  Y_{\alpha(\ell,i-1)}, Y_{\alpha(\ell,i+1)},\ldots, Y_{\alpha(\ell,n)} \Big)\, .\end{equation}

As for the $(n-1)$-dimensional target space of $\pi_\alpha$, one verifies easily that the following holds true:  if $i$ 
labels one of the (2 or 3) extremities of $\Delta$ then $\Delta_{\widehat \imath}$ 
is still a Dynkin diagram; otherwise, $\Delta_{\widehat \imath}$  admits two sub-Dynkin diagrams of $\Delta$ as connected components, denoted by $\Delta'_{\widehat \imath}$ and $\Delta''_{\widehat \imath}$. 
Accordingly,  it is natural to see 
$\pi_\alpha$ as a map  from $ \mathcal X_\Delta$ onto $\mathcal X_\alpha=\mathcal X_{\Delta_{\widehat \imath}}$
in the former case, and onto the product $\mathcal X_\alpha=\mathcal X_{\Delta'_{\widehat \imath}}\times 
\mathcal X_{\Delta''_{\widehat \imath}}$ in the latter case.\mk 

We thus have given a nice construction of certain subwebs of ${\boldsymbol{\mathcal X\hspace{-0.04cm}\mathcal W}}_{\hspace{-0.04cm} \Delta}$ which have intrinsic dimension $n-1$ and are maximal for the inclusion (as it can be verified). 
It is natural to ask if  this construction gives  all the subwebs satisfying the last two properties:
\begin{question} 
\label{Quest:tookoo}
Let $\boldsymbol{\mathcal W}$ be  a subweb  of ${\boldsymbol{\mathcal X\hspace{-0.04cm}\mathcal W}}_{\hspace{-0.04cm} \Delta}$ of intrinsic dimension $n-1$ which is also 
maximal for the inclusion.  Does it necessarily 
coincide with one of the  ${\boldsymbol{\mathcal X\hspace{-0.04cm}\mathcal W}}_{\hspace{-0.06cm} \Delta,\alpha}$'s? 
\end{question}
The answer is obvious (and affirmative) in rank 2. 
The same holds true for the two examples in rank 3 considered below (namely $A_3$ and $B_3$). 
We believe that the answer is affirmative in full generality.

\subsubsection{}  
\label{SS:FockGoncharo-XInfinity}
There is a nice geometric/topological way to understand the maps $F_\alpha$ all together. 
We now discuss this quickly without giving any justifications. For details/explanations/proofs, the reader can consult \cite{FockGoncharovXInfinity} as well as  \cite{AHHL}. 
\sk

The positive part $\mathcal X_\Delta(\mathbf R_{>0})$ can be identified with the interior 
of a convex polytope of dimension $n$,  the {\bf cluster associahedron   of type  $
\boldsymbol{\Delta}$}, denoted by $P_\Delta$. 
This polytope has been constructed combinatorically in \cite[\S3]{FZ}
(see also \cite{FR}), where it is proved that the facets (the faces  of codimension 1) of 
$P_\Delta$ are naturally in bijection with $\Delta_{\geq -1}$. 
For $\alpha\in \Delta_{\geq -1}$, we denote by  $P_\alpha$  the facet
of $P_\Delta$ 
 associated to it.  Then for each such root $\alpha$, the positive part $\mathcal X_\alpha(\mathbf R_{>0})$ 
of the image of 
\begin{equation}
\label{Eq:Pi-alpha}
\pi_\alpha:\mathcal X_\Delta\rightarrow \mathcal X_\alpha
\end{equation}
 (which, by the way,  is a positive map 
according to \eqref{Eq:pi-alpha}) identifies itself with the facet $P_\alpha$ of 
$\boldsymbol{\Delta}$.  Let $p_\alpha: P_\Delta\rightarrow P_\alpha$ be  (the restriction 
to $P_\Delta$ of the euclidean projection onto $P_\alpha$. Then the linear web  $\boldsymbol{\mathcal W}(\, p_\alpha\, \lvert \, \alpha\in \Delta_{\geq -1}\, \big)$ gives a combinatorial model on the $\Delta$-associahedron 
 of the web $
\boldsymbol{\mathcal W}(\, \pi_\alpha\, \lvert \, \alpha\in \Delta_{\geq -1}\, \big)$ 
 on  $\mathcal X_\Delta$ defined by all the $\pi_\alpha$'s.
 \mk

 An interesting question is whether  one can recover the poset of faces of the associahedron $P_\Delta$ from the web 
 ${\boldsymbol{\mathcal X\hspace{-0.04cm}\mathcal W}}_{\hspace{-0.04cm} \Delta}$.  
 Given a web ${\boldsymbol{\mathcal W}}$ of intrinsic dimension $n$, let ${C}^1({\boldsymbol{\mathcal W}})$ be the set of its subwebs with intrinsic dimension $n-1$ and which are maximal (with respect to inclusion), define ${C}^2({\boldsymbol{\mathcal W}})$ as the union of ${C}^1({\boldsymbol{\mathcal W}}')$'s for
 ${\boldsymbol{\mathcal W}}'$ in ${C}^1({\boldsymbol{\mathcal W}})$, etc.
  We define a notion of adjacency between two elements  ${\boldsymbol{\mathcal W}}^1,{\boldsymbol{\mathcal W}}^2 \in 
 {C}^{\geq 1}({\boldsymbol{\mathcal W}})=
 \cup_{c\geq 1}  {C}^{c}({\boldsymbol{\mathcal W}})$ by requiring inclusion 
 if the intrinsic dimensions of these two webs differ and otherwise  by the fact that 
 ${C}^{1}({\boldsymbol{\mathcal W}}^1)$ and ${C}^{1}({\boldsymbol{\mathcal W}}^2)$ have  (at least) an element in common.   Several interesting questions can be asked relatively to these notions: 

%
%
%
\begin{questions}
\begin{enumerate}
\item  
Can the poset of faces of $P_{\hspace{-0.05cm}\Delta}$ be reconstructed from  the set 
${C}^{\geq 1}\big({\boldsymbol{\mathcal X\hspace{-0.04cm}\mathcal W}}_{\hspace{-0.04cm} \Delta}\big)$ together with the notion of `adjacency' defined above?
\item
 In particular, can (and if yes, how) the $\mathcal X$-clusters  be obtained 
 from 
${\boldsymbol{\mathcal X\hspace{-0.04cm}\mathcal W}}_{\hspace{-0.04cm} \Delta}$ alone? 

\item More generally, what informations about $P_{\hspace{-0.05cm}\Delta}$ is it possible to reconstruct from  
${\boldsymbol{\mathcal X\hspace{-0.04cm}\mathcal W}}_{\hspace{-0.04cm} \Delta}$? 

\item   Does ${C}^{\geq 1}\big({\boldsymbol{\mathcal X\hspace{-0.04cm}\mathcal W}}_{\hspace{-0.04cm} \Delta}\big)$ together with the notion of `adjacency' characterize this web?

 More precisely, let ${\boldsymbol{\mathcal W}}$ be a $d^{\mathcal X}_{\Delta}$-web in $n$ variables such that  $({C}^{\geq 1}\big({\boldsymbol{\mathcal W}}\big), {adjacency}) $ is isomorphic 
to the corresponding combinatorial object for ${\boldsymbol{\mathcal X\hspace{-0.04cm}\mathcal W}}_{\hspace{-0.04cm}\Delta}$ (possibly up to relabelling the foliations). Is then 
${\boldsymbol{\mathcal W}}$ necessarily equivalent to ${\boldsymbol{\mathcal X\hspace{-0.04cm}\mathcal W}}_{\hspace{-0.04cm}\Delta}$?
\end{enumerate}
\end{questions}

Of course, ${C}^1({\boldsymbol{\mathcal W}})$ is trivial when $n=2$ hence these questions do not make sense in this case.  Note that answering to Question \ref{Quest:tookoo} is a first step for doing the same for the first of the four questions just above.  Since the answer to the former question is affirmative 
when $\Delta$ is  $A_3$ or $B_3$ (see below), we believe that the one to the latter is affirmative too for any Dynkin diagram $\Delta$.  



\subsubsection{}  
The reader will have noticed that we have not yet given any definition of a web of codimension 1  on $\mathcal X_\Delta$. A first definition could be to consider the one admitting   the face maps \eqref{Eq:Pi-alpha} for all $\alpha\in \Delta_{\geq -1}$ as first integrals, but this is actually not the  right one, as it will be clear from case $A_3$ that we consider below.\mk 

Let us then treat the case when  $\Delta=A_3$. Since $h(A_3)=4$ is even, one can construct the nine face maps $\pi_\alpha$ easily by following the recipe described above. 
The birational map corresponding to the composition of mutations $\mu_{\bullet\lvert \circ}$ is 

$$F_{A_3}\, :\, 
\big( u_1,u_2,u_3)\longmapsto 
\left( {\frac {1}{u_{{1}} ( 1+u_{{2}}) } } ,{\frac { \left( u_{{2}}u_{{3}}+u_{{3}}+1 \right)  \left( u_{{1}}u_{{2}}+u_{{1}}+1 \right) }{u_{{2}}}},{\frac {1}{u_{{3}} \left( 1+u_{{2}} \right)}} \right)
\, .
$$ 

We set $F_{[12]}=1+u_1+u_2$, $F_{[23]}=1+u_2+u_3$, $F_{[13]}=1+u_1+u_2+u_3+u_1u_3$  and 
$ \boldsymbol{Y}(0)= ( {1}/{u_{{1}}}\, ,\, u_{{2}}, \,  {1}/{u_{{3}}} )$. Then for 
$\boldsymbol{Y}(1)= F_{A_3}\big(\boldsymbol{Y}(0) \big)$ and 
$\boldsymbol{Y}(2)= F_{A_3}\big(\boldsymbol{Y}(1) \big)$, we have 
$$
 \boldsymbol{Y}(1)= 
 \left( 
 {\frac {u_{{1}}}{1+u_{{2}}}},{\frac { F_{[12]} F_{[23]} }{u_{{1}}u_{{2}}u_{{3}}}},{\frac {u_{{3}}}{1+u_{{2}}}}
  \right) 
\quad \mbox{ and } \quad 
  \boldsymbol{Y}(2)=
  \left( 
  {\frac {u_{{2}}u_{{3}}}{F_{[13]}}},{\frac { \left( u_{{3}}+1 \right)  \left( u_{{1}}+1 \right) }{u_{{2}}}},{\frac {u_{{1}}u_{{2}}}{F_{[13]}}}
  \right)\, .
$$

Consequently,  the three face maps of type $A_1\times A_1$ are 
$$ \pi_{-\alpha_2}=\big( \,  {u_{{1}}}\, ,\, {u_{{3}}}\, \big)\, , \quad 
\pi_{\alpha_1+\alpha_2+\alpha_3}=  \left( \, 
 {\frac {1+u_2}{u_1}} \, , \, {\frac {1+u_2}{u_3}}\, 
   \right) 
  \quad \mbox{ and }\quad 
  \pi_{\alpha_2}=
  \left( \, 
  \frac{F_{[13]}}{u_2u_3} \, ,\, \frac{  F_{[13]}}{u_2u_3}\, 
  \right)\, 
$$
and for the six ones of type $A_2$, one can take 
\begin{align*}
\pi_{-\alpha_1}=&\,  \big( u_2,u_3\big) 
& \pi_{\alpha_1}= \left( \,  {\frac { F_{[12]} F_{[23]} }{u_{{1}}u_{{2}}u_{{3}}}}\, ,\, {\frac {1+u_{{2}}}{u_{{3}}}} \, 
\right) && {}^{} \hspace{-0.1cm}
 \pi_{\alpha_1+\alpha_2}=  \left( \, 
  {\frac {F_{[13]}}{u_{{2}}u_{{3}}} }\, ,\, {\frac { \left( u_{{3}}+1 \right)  \left( u_{{1}}+1 \right) }{u_{{2}}}}
 \,  \right)
\\
 \pi_{-\alpha_3}= &\,   \big( u_1,u_2\big) 
 &
\pi_{\alpha_3}=  \left( \, 
 {\frac {1+u_{{2}}} {u_{{1}}}}\, ,\, {\frac { F_{[12]} F_{[23]} }{u_{{1}}u_{{2}}u_{{3}}}} \,  \right)  
 && {}^{} \hspace{0.1cm}  \pi_{\alpha_2+\alpha_3}=
  \left( \, {\frac { \left( u_{{3}}+1 \right)  \left( u_{{1}}+1 \right) }{u_{{2}}}}\, ,\, {\frac {F_{[13]}} {u_{{1}}u_{{2}}} }\, 
  \right)\, .
\end{align*}

We can then consider two curvilinear webs, the 9-web defined by considering all the face maps $\pi_\alpha$, or the 6-web by disregarding the face maps of non-irreducible type $A_1\times A_1$. 
Of course, the second option is the right one in order to recover the web noted by ${\boldsymbol{\mathcal X\hspace{-0.04cm}\mathcal W}}_{\hspace{-0.04cm} A_3}^{(1)}$  above. \mk 

But the fact that the 9-web 
$
{\boldsymbol{\mathcal W}}\big(\, \pi_\alpha \, \lvert \, \alpha\in (A_3)_{\geq -1}\, \big)$ is not the right one to be considered 
is really 
obvious from a web-theoretical point of view, if we notice  that 
it  does not satisfy the natural general position assumption usually required for curvilinear webs: if $\mathcal D_\alpha$ stands for the tangent distribution of dimension 1 defined by $d\pi_\alpha=0$ for any $\alpha$, there exist 3-tuples $(\alpha,\beta,\gamma)$ of pairwise distinct elements of $(A_3)_{\geq -1}$ such that $\mathcal D_\alpha,\mathcal D_\beta$ and $\mathcal D_\gamma$ are not in direct sum at the generic point of $\mathbf C^3$.   
\mk 

 As for webs of codimension 1, for any $\sigma\geq 0$, 
 one can define the $\sigma$-th virtual rank for  $(n-1)$-ARs of a curvilinear web 
 ${\boldsymbol{\mathcal W}}$
 in $n$ variables, denoted by $\rho_{(n-1)}^\sigma({\boldsymbol{\mathcal W}})$.\footnote{Stating a precise definition for $\rho_{(n-1)}^\sigma({\boldsymbol{\mathcal W}})$ is left as an exercise to the reader.} The fact that the foliations of $
{\boldsymbol{\mathcal W}}\big(\, \pi_\alpha \, \lvert \, \alpha\in (A_3)_{\geq -1}\, \big)$ are not in general position has the unfortunate consequence that 
$\rho_{(2)}^\sigma({\boldsymbol{\mathcal W}}\big(\, \pi_\alpha \, \lvert \, \alpha\in (A_3)_{\geq -1}\, \big)\big)=7$ for any $\sigma\geq 5$. Hence its total virtual 2-rank 
$$\rho_{(2)}\Big({\boldsymbol{\mathcal W}}\big(\, \pi_\alpha \, \lvert \, \alpha\in (A_3)_{\geq -1}\, \big)\Big)=\sum_{\sigma\geq 0}\rho_{(2)}^\sigma\Big({\boldsymbol{\mathcal W}}\big(\, \pi_\alpha \, \lvert \, \alpha\in (A_3)_{\geq -1}\, \big)\Big)$$ is infinite hence it makes no sense to ask whether  it is of AMP rank or not. \mk

The situation is quite the opposite for  the web 
$${\boldsymbol{\mathcal X\hspace{-0.04cm}\mathcal W}}_{\hspace{-0.04cm} A_3}^{(1)}=
{\boldsymbol{\mathcal W}}\Big(\, \pi_\alpha \, \lvert \, \alpha\in (A_3)_{\geq -1}\,\mbox{ of irreducible type}\, \Big)$$
which corresponds to ${\boldsymbol{\mathcal X\hspace{-0.04cm}\mathcal W}}_{\hspace{-0.04cm} \mathcal M_{0,6} }^{(1)}$ up to the identification $\mathcal X_{A_3}\simeq \mathcal M_{0,6}$. 
By direct computations, one verifies that $\rho_{(2)}^\bullet\big( \boldsymbol{\mathcal X\hspace{-0.05cm}\mathcal W}_{A_3}^{(1)}  \big)
=(3,4,3)$ hence 
$${\bf rk}_{(2)}\big( \boldsymbol{\mathcal X\hspace{-0.05cm}\mathcal W}_{A_3}^{(1)}
\big) \leq 
\rho_{(2)}\big( \boldsymbol{\mathcal X\hspace{-0.05cm}\mathcal W}_{A_3}^{(1)}  \big)=10\, , $$
 in accordance with Proposition 2.4 of \cite{Damiano}.   
It can be proved that  this bound actually is an equality hence 
 the web $\boldsymbol{\mathcal X\hspace{-0.05cm}\mathcal W}_{A_3}^{(1)} $ is AMP. \footnote{This would follow from some results of  \cite{Damiano} but some of the statements in it are wrong, in particular in the case of $\boldsymbol{\mathcal X\hspace{-0.05cm}\mathcal W}_{A_3}^{(1)}$. This web has maximal 2-rank but this has to be established following another approach than that of  \cite{Damiano}.}
\mk

The preceding example indicates that in order to get an interesting web from the face maps $\pi_\alpha$, one has to only consider those associated to the irreducible facets, that is the $P_\alpha$'s which are not a product of two associahedra of positive dimension. These are  the facets associated to the roots $\alpha(k,i)$ for some $k\in \mathbf Z$ where $i$ labels an extremal node of the Dynkin diagram.  Thus, for any Dynkin diagram, the curvilinear cluster web which it 
is  relevant to consider is the following one: 
$$
{\boldsymbol{\mathcal X\hspace{-0.04cm}\mathcal W}}_{\hspace{-0.06cm} \Delta}^{(1)}=
{\boldsymbol{\mathcal W}}\left(\, \pi_\alpha : \mathcal X_\Delta\rightarrow \mathcal X_{\Delta_\alpha} 
 \hspace{0.06cm} 
 \big\lvert 
 \hspace{0.05cm}  
 \begin{tabular}{l}
 {\it  $ \alpha=\alpha(k,i)$ for  $i$ labelling}
 \vspace{-0.07cm}
 \\
 {\it  an extremal node of $\Delta$}
\end{tabular}
\right)\, .
$$

\subsubsection{}  
We now consider another explicit example: $B_3$. Looking at the $B_3$-associahedron, which is called the 3-dimensional cyclohedron and is pictured in \cite[Fig.\,3]{CFZ}, it comes that 
$\boldsymbol{\mathcal X\hspace{-0.05cm}\mathcal W}_{\hspace{-0.03cm}B_3}^{(1)}$ is a 8-web, 
defined by four face maps of type $A_2$, and four others of type $B_2$. 
\sk

To simplify the notation, it is more convenient to label the face maps via the roots of the dual root system $C_3$: we set $p_{\alpha^\vee}=\pi_\alpha$ for any $\alpha\in \Delta_{\geq -1}$ with $\Delta=B_3$ here. Moreover,  the components of the $\pi_\alpha$  in \eqref{Eq:pi-alpha}
 can have quite complicate analytic expressions and for some almost-positive roots $\alpha^\vee\in \Delta^\vee$, it is more useful  to consider  instead a variant $\widetilde p_{\alpha^\vee}= (X[\alpha'],X[\alpha''])$ which defines the same foliation by rational curves but whose components $X[\alpha']$ and $X[\alpha'']$ are $\mathcal X$-variables which are formally simpler than the 
$\mathcal Y$-variables appearing in  \eqref{Eq:pi-alpha}.
Below, $\alpha_1,\alpha_2$ and $\alpha_3$ stand for the standard simple positive roots in $\Delta^\vee$ (= the root system of type $C_3$). \mk

The four 
$B_2$-face maps are the ones associated to the roots $\alpha_1$, $-\alpha_1$, $\alpha_1+\alpha_2$ and $\alpha_2+\alpha_3$ of $\Delta^\vee$. It can be verified that  the following cluster maps are first integrals for the corresponding foliations: 
\begin{align*}
p_{ -\alpha_1}= & \,  \left( \, \frac{1}{u_2}\, , \, \frac{1}{ u_3}  \, \right) \\
   \widetilde p_{ \alpha_1}= & \,  \left(   \, 
 \frac{1+u_1+u_2}{u_1u_2}
 \, ,\, {\frac { \left( 1+u_{{2}} \right) ^{2}}{u_{{3}}}}  \,  \right) \\
 \widetilde p_{\alpha_1+\alpha_2}= &\,\left(  \, 
  {\frac {(1+u_1)(1+u_3)}{u_{{2}}}} \, ,\,
   {\frac {1}{u_{{3}}}}  \,   \right)  \\
\mbox{and }\quad  \widetilde p_{\alpha_2+\alpha_3 }= & \,  \left(   \, 
\frac{\big(1+u_2+u_3\big)^2}{u_2^2u_3}
\, , \, 
\frac{\big(1+u_1+u_2+u_3+u_1u_3\big)^2}{u_2^2u_3}  \, 
\right) \, .
\end{align*}

As for the four 
$A_2$-face maps, these are the ones  associated to the roots $\alpha_3$, $-\alpha_3$, $2\alpha_2+\alpha_3$ and $\alpha_1+2\alpha_2+\alpha_3$ and the following cluster maps are first integrals for the corresponding foliations: 
\begin{align*}
p_{ -\alpha_3}= & \,  \left( \, \frac{1}{u_1}\, , \, \frac{1}{ u_2}\,  \right) \\
   \widetilde p_{\alpha_3 }= & \, 
   \left( \,  {\frac {1+u_{{2}}}{u_{{1}}}}\, ,\,  {\frac {
1+2\,u_{{2}}+u_{{3}}+
{u_{{2}}}^{2}}{u_2u_{{3}}}} \, \right)
 \\
 \widetilde p_{2\alpha_2+\alpha_3}= &\,
 \left(  \,  {\frac {1}{u_{{1}}}} \, ,\, 
  {\frac {(1+u_1)(1+u_3)}{u_{{2}}}} \, 
  \right)  \\  
\mbox{and }\quad  \widetilde p_{\alpha_1+2\alpha_2+\alpha_3 }= & \,  \left( 
\, \frac{1+u_1+u_2+u_3}{u_1u_2}
\, , \, 
\frac{1+u_1+2u_2+u_3+u_1u_2+u_1u_3+u_2^2}{u_2u_3}\, 
\right) \, .
\end{align*}




By direct computations, we have obtained that 
$$\rho_{(2)}^\bullet\Big( \boldsymbol{\mathcal X\hspace{-0.05cm}\mathcal W}_{\hspace{-0.05cm} B_3}^{(1)}  \Big)
=\big(5,8,9,8,5\big)
\qquad \mbox{ hence } \qquad 
\rho_{(2)}\Big( \boldsymbol{\mathcal X\hspace{-0.05cm}\mathcal W}_{\hspace{-0.05cm} B_3}^{(1)}  \Big)=35\, .
$$


Furthermore, by direct computations again, it can be verified that  the space of 
fifth-order jets of 2-ARs of 
 $\boldsymbol{\mathcal X\hspace{-0.05cm}\mathcal W}_{\hspace{-0.05cm} B_3}^{(1)}$ 
is indeed of dimension 35. But some obstructions appear  at the next two orders. 
Starting from order 8, the dimensions of the spaces of higher order jets of ARs  stabilize and all are equal to  21, which therefore is the 2-rank of  $\boldsymbol{\mathcal X\hspace{-0.05cm}\mathcal W}_{B_3}^{(1)}$.  Consequently,  this web is not 2-AMP contrarily to $\boldsymbol{\mathcal X\hspace{-0.05cm}\mathcal W}_{\hspace{-0.05cm} A_3}^{(1)}$ which is. \sk

The case $n=2$ is quite specific. First,  the three webs 
$\boldsymbol{\mathcal Y\hspace{-0.05cm}\mathcal W}_{\hspace{-0.05cm} B_2}$, 
$\boldsymbol{\mathcal X\hspace{-0.05cm}\mathcal W}_{\hspace{-0.05cm}B_2}$ and $\boldsymbol{\mathcal X\hspace{-0.05cm}\mathcal W}_{\hspace{-0.05cm}B_2}^{(1)}$ coincide hence the latter is AMP. But this has to be seen as an exceptional phenomenon resulting  from the isomorphism $B_2\simeq C_2$. 
Being both of type $B$ and of type $C$, 
 $\boldsymbol{\mathcal X\hspace{-0.05cm}\mathcal W}_{\hspace{-0.05cm}B_2}$ inherits from each type a  supplementary non-polylogarithmic AR (namely ${\sf A}_{B_2}$ for type $B_2$, ${\sf J}_{C_2}$ for type $C_2$) which both contribute to make its rank maximal.  For this reason, the $B_2$-case appears as exceptional and 
 and not considering it as the first member of the family of   $\boldsymbol{\mathcal X\hspace{-0.05cm}\mathcal W}_{\hspace{-0.05cm} B_n}^{(1)}$'s is more natural. \mk 

We assume $n\geq 3$. In \S\ref{SSS:ClusterWebs-type-B:Y-cluster web} and  \S\ref{SSS:ClusterWebs-type-B:X-cluster web},  we have verified for $n$ small and  conjectured for $n$ arbitray  that the standard ({\it i.e.}\,1-codimensional) cluster web of type $B_n$ is not AMP. Considering this, the fact that $\boldsymbol{\mathcal X\hspace{-0.05cm}\mathcal W}_{B_3}^{(1)}$ and possibly all the $\boldsymbol{\mathcal X\hspace{-0.05cm}\mathcal W}_{\hspace{-0.04cm}B_n}^{(1)}$'s
 for $n\geq 3$ are not AMP appears as predictable. However, even if not AMP, these curvilinear webs are interesting. And one can ask about them the same questions we asked about the webs 
 $\boldsymbol{\mathcal Y\hspace{-0.05cm}\mathcal W}_{\hspace{-0.04cm}B_n}$ and $\boldsymbol{\mathcal X\hspace{-0.05cm}\mathcal W}_{\hspace{-0.04cm} B_n}$: 

\begin{questions} 

1. Describe a basis of the spaces of $(n-1)$-abelian relations of $\boldsymbol{\mathcal X\hspace{-0.05cm}\mathcal W}_{\hspace{-0.04cm}B_n}^{(1)}$. \sk 

2.  In particular, does this web carry $(n-1)$-ARs 
which are  of polylogarithmic type or of another type, say of that 
of  the AR $({\sf A}_{B_n})$ of 
$\boldsymbol{\mathcal Y\hspace{-0.05cm}\mathcal W}_{\hspace{-0.04cm}B_n}$?\sk 
\end{questions}


\subsubsection{}  
We end our discussion about the webs $\boldsymbol{\mathcal X\hspace{-0.05cm}\mathcal W}_{\hspace{-0.05cm}\Delta}^{(1)}$ by a last question related  to a nice geometric result of Burau in the case when   $\Delta$ is of type $A$.\mk 

For any $\Delta$, since the components of the face maps $\pi_\alpha$ appearing as first integrals of 
$\boldsymbol{\mathcal X\hspace{-0.05cm}\mathcal W}_{\hspace{-0.05cm}\Delta}^{(1)}$ can be taken as elements of a same cluster, each foliation $\mathcal F_{\pi_\alpha}$ is birationally equivalent to the family of lines passing through  a point in $\mathbf P^n$.  Consequently, the leaves of 
$\boldsymbol{\mathcal X\hspace{-0.05cm}\mathcal W}_{\hspace{-0.05cm}\Delta}^{(1)}$ are rational curves.  
\sk


In \cite{Burau}, for any $n\geq 2$, Burau gave the construction of a projective variety in which the $(n+3)$-web $\boldsymbol{\mathcal X\hspace{-0.05cm}\mathcal W}_{\hspace{-0.05cm}{\mathcal M}_{0,n+3}}^{(1)}$ is realized by rational curves of  degree 1 or 2, depending on the parity of $n$:  
\begin{prop} We set $\delta_{A_n}=1$ when $n$ is odd, and $\delta_{A_n}=2$ when it is even.
\begin{enumerate}
\item  There exists a projective variety $V_{n}$ of dimension $n$ which 
\begin{itemize}
\item carries a $(n+3)$-web  by rational curves of degree $\delta_{A_n}$ which is isomorphic 
to $\boldsymbol{\mathcal X\hspace{-0.05cm}\mathcal W}_{\hspace{-0.05cm}{\mathcal M}_{0,n+3}}^{(1)}$;
\item is not covered by (rational) curves of degree less than $\delta_{A_n}$.\footnote{Of course,  this condition is empty when $n$ is odd since $\delta_{A_n}=1$ in this case, which is as minimal as it can be.}
\end{itemize}
\item  Consequently, $\boldsymbol{\mathcal X\hspace{-0.05cm}\mathcal W}_{\hspace{-0.05cm}A_{n}}^{(1)}$ can be linearized when $n$ is odd, and is equivalent to a web by conics (as leaves) when $n$ is even. 
\end{enumerate}
\end{prop}

Burau actually gave a construction of the variety $V_n$ by means of an explicit linear system: 
let $p_1,\ldots,p_{n+2}$ be $n+2$ points in general position in $\mathbf P^n$ and 
denote by 
$\mathcal L_{A_n}$ the linear system $\lvert (n+1)H-\sum_i (n-1)p_i\lvert$ when $n$ is even, and 
$\lvert ((n+1)/2)H-\sum_i ((n-1)/2)p_i\lvert$ when $n$ is odd (where $H$ stands for the classe of a hyperplane in ${\rm Pic}(\mathbf P^n)$).  Then the associated rational map $\varphi_{\mathcal L_{A_n}}$ is generically 1-1 onto its image 
$V_n=\varphi_{\mathcal L_{A_n}}(\mathbf P^n)$ which therefore has dimension $n$ in both cases.\footnote{In addition to Burau's paper \cite{Burau}, another classical reference 
where   the linear system $\mathcal L_{A_n}$ and the variety $V_n$ are discussed is  the treatise \cite{Room} by Room. For some recent references, see the papers \cite{Kumar1,Kumar2} by Kumar.}
\sk

Let $C$ be either a generic line passing through one of the $p_j$'s or a rational normal curve of degree $n$ passing through all of them. Then it can be verified that 
$\varphi_{\mathcal L_{A_n}}(C)$ has degree $\delta_{A_n}$. Since the $(n+3)$-web on $\mathbf P^n$  formed by the curves of this kind is a model of $\boldsymbol{\mathcal X\hspace{-0.05cm}\mathcal W}_{\hspace{-0.05cm}{\mathcal M}_{0,n+3}}^{(1)}$, it follows that the same holds for its 
 push-forward by $\varphi_{\mathcal L_{A_n}}$, which is a web by rational curves of degree $\delta_{A_n}$ on $V_n$. 
\sk

The cases when $n=2$ and $n=3$ in the previous description correspond to classical and  very well-known varieties and webs:  $V_2$ is the smooth quintic \href{https://en.wikipedia.org/wiki/Del_Pezzo_surface}{\it Del Pezzo surface} in $\mathbf P^5$ and the corresponding web on it is the one formed by the five pencils of conics included in it (it is a nice geometric model of Bol's web); and $V_3$ is the famous \href{https://en.wikipedia.org/wiki/Segre_cubic}{\it Segre cubic primal}\footnote{It is a cubic hypersurface in $\mathbf P^4$ first considered by \href{https://fr.wikipedia.org/wiki/Corrado_Segre}{C. Segre}, who proved that it can be characterized by the fact of having exactly 10 double points, which is the maximal possible number (see \cite{Dolgachev}).} in $\mathbf P^4$ 
 and the associated 6-web is the one formed by the lines contained in it.
\mk 

These results leads to wonder about the case of the curvilinear webs 
$\boldsymbol{\mathcal X\hspace{-0.05cm}\mathcal W}_{\hspace{-0.05cm}\Delta}^{(1)}$ for $\Delta$ arbitrary:
\begin{questions}
Let $\Delta$ be a Dynkin diagram of rank $n\geq 2$, not of type $A$. 
\begin{enumerate}
\item Is $\boldsymbol{\mathcal X\hspace{-0.05cm}\mathcal W}_{\hspace{-0.05cm}\Delta}^{(1)}$ linearizable? If yes, how?  If not, what is the smallest integer $\delta_\Delta \geq 2$ such that it can be realized as a web by rational curves of degree $\delta_\Delta$?
\item Is there a projective variety $V_\Delta$ of dimension $n$ carrying a web by rational curves of
 degree $\delta_\Delta$ equivalent to $\boldsymbol{\mathcal X\hspace{-0.05cm}\mathcal W}_{\hspace{-0.05cm}\Delta}^{(1)}$ but which is not covered by 
 rational curves of degree less than $\delta_\Delta$?
\end{enumerate}
\end{questions}

Answers to these questions are not known when $\Delta$ is not of type $A$, even in the simplest case when $\Delta$  has rank 2.  The case of $B_2$ should be the first  to be considered and already seems to be interesting.

\newpage

\bigskip\bigskip

{\small Luc Pirio\\
 LMV, UMR 8100,  CNRS -- Universit\'e Versailles--St-Quentin\\
Email: {\tt luc.pirio@uvsq.fr}
\end{document}